\documentclass[10pt,openany,oneside]{book}
\usepackage[textheight=538pt,
            textwidth=346pt,
            hmarginratio=1:1]{geometry}
\usepackage[utf8]{inputenc}
\usepackage{layout}
\usepackage{eucal}

\usepackage{etoolbox}

\usepackage{tocloft}

\usepackage{appendix}
\usepackage{trimspaces}
\usepackage{marginnote}

\usepackage{index}

\newindex{autind}{aidx}{aind}{Author index}
\newindex{notind}{nidx}{nind}{Notation index}
\newindex{default}{sidx}{sind}{Subject index}

\usepackage{multicol}



\usepackage[round]{natbib}
\usepackage{mathtools}
\usepackage{booktabs}
\usepackage{xspace}
\usepackage{makecell}
\usepackage{longtable}
\usepackage{relsize}
\usepackage[dvipsnames]{xcolor}
\usepackage{bigstrut}

\usepackage{tikz} %

\usepackage{subcaption} %
\usepackage{framed}

\usepackage{algorithm}
\usepackage[noend]{algpseudocode}

\algblockdefx{Block}{EndBlock}[1]{\textit{#1:}}{}
\algtext*{EndBlock}
\algrenewtext{Loop}[1]{#1}

\newcommand{\true}{\textbf{true}\xspace}
\newcommand{\false}{\textbf{false}\xspace}
\algrenewcommand\algorithmicif{\textrm{if}}
\algrenewcommand\algorithmicthen{\textrm{then}}
\algrenewcommand\algorithmicelse{\textrm{else}}
\algrenewcommand\Call[2]{\textup{\textsc{#1}(#2)}}
\algrenewcommand\algorithmiccomment[1]{\hfill\textit{(#1)}}

\makeatletter
\ifdefined\theHalgorithm
\else
\newcounter{Halgorithm}
\newcounter{HALG@line}
\fi
\renewcommand{\theHalgorithm}{\thechapter.\arabic{algorithm}}
\renewcommand{\theHALG@line}{\thechapter.\arabic{algorithm}.\arabic{ALG@line}}
\makeatother

\makeatletter
\begingroup
  \let\newcounter\@gobble
  \let\setcounter\@gobbletwo
  \globaldefs\@ne
  \let\c@loadepth\@ne
  \newlistof{algorithms}{loa}{\listalgorithmname}
\endgroup
\let\l@algorithm\l@algorithms
\makeatother

\def\listfigurealgoname{List of figures and algorithms}
\makeatletter
\def\listoffiguresalgos{%
   \chapter{\listfigurealgoname}
   \@mkboth{\MakeUppercase\listfigurealgoname}{\MakeUppercase\listfigurealgoname}
   \@starttoc{lof}
   \bigskip
   \@starttoc{loa}
}
\makeatother

\makeatletter
\patchcmd{\@chapter}{\addtocontents{lof}{\protect\addvspace{10\p@}}}{}{}{}
\makeatother

\usepackage{amsmath,amssymb,amsthm,bm,microtype,thmtools}
\usepackage{enumitem}

\cftpagenumbersoff{part}

\usepackage[pdfencoding=unicode,hypertexnames=true]{hyperref}
\usepackage{bookmark}

\usepackage{lineno}

\usepackage[nosort]{cleveref} %
\Crefname{equation}{Eq.}{Eqs.}
\Crefname{proposition}{Proposition}{Propositions}
\crefname{proposition}{proposition}{propositions}
\crefname{item}{part}{parts}
\Crefname{item}{Part}{Parts}
\Crefname{lemma}{Lemma}{Lemmas}
\crefname{lemma}{lemma}{lemmas}
\Crefname{figure}{Figure}{Figures}

\definecolor{cerulean}{rgb}{0.10, 0.58, 0.75}

\newcommand{\auxiliary}{magnitude\xspace}
\newcommand{\ndim}{$n$\nobreakdash-\hspace{0pt}dimensional\xspace}

\newcommand{\mathtogether}{%
\thickmuskip=5mu plus 5mu minus 1mu
\relpenalty=10000
\binoppenalty=10000
}
\newcommand{\squeezebin}{%
\medmuskip=1.5mu plus 2.5mu minus 1.5mu
}

\newcommand{\mathopfont}[1]{#1}%

\DeclareFontFamily{U}{stixex}{}
\DeclareFontShape{U}{stixex}{m}{n}{
    <-> stix-mathex
}{}
\DeclareSymbolFont{stixex}{U}{stixex}{m}{n}
\DeclareMathDelimiter{\lParen}{\mathopen}{stixex}{"DE}{stixex}{"02}
\DeclareMathDelimiter{\rParen}{\mathclose}{stixex}{"DF}{stixex}{"03}
\DeclareMathDelimiter{\lBrace}{\mathopen}{stixex}{"E8}{stixex}{"0E}
\DeclareMathDelimiter{\rBrace}{\mathclose}{stixex}{"E9}{stixex}{"0F}
\DeclareMathDelimiter{\lbrk}{\mathopen}{stixex}{"EE}{stixex}{"14}
\DeclareMathDelimiter{\rbrk}{\mathclose}{stixex}{"EF}{stixex}{"15}

\makeatletter
\newcommand{\pparens}[1]{
  \lParen\negKern#1\negKern\rParen%
  \@ifnextchar({}{}%
}
\newcommand{\bbraces}[1]{
  \lBrace\negKern#1\negKern\rBrace%
  \@ifnextchar({}{}%
}
\makeatother

\newcommand{\me}{e}

\newcommand{\pseudinvsym}{\dagger}

\newcommand{\zero}{\mathbf{0}}
\newcommand{\aaa}{\mathbf{a}}
\newcommand{\bb}{\mathbf{b}}

\newcommand{\tuu}{\tilde{\uu}}
\newcommand{\tu}{\tilde{u}}

\newcommand{\cc}{\mathbf{c}}
\newcommand{\dd}{\mathbf{d}}
\newcommand{\ee}{\mathbf{e}}
\newcommand{\xx}{\mathbf{x}}
\newcommand{\xing}{\hat{\xx}}

\newcommand{\ying}{\hat{y}}
\newcommand{\zing}{\hat{\zz}}
\newcommand{\uhat}{\hat{\uu}}
\newcommand{\vhat}{\hat{\vv}}
\newcommand{\ving}{\hat{v}}
\newcommand{\what}{\hat{\ww}}
\newcommand{\bhati}{\hat{b}}
\newcommand{\bhat}{{}\mkern-6mu\hat{\mkern6mu\mathbf{b}}}
\newcommand{\tbb}{{}\mkern-6mu\tilde{\mkern6mu\mathbf{b}}}
\newcommand{\chat}{\hat{\cc}}
\newcommand{\hatc}{\hat{c}}

\newcommand{\qhat}{{\hat{\bf q}}}
\newcommand{\xhat}{{\hat{\bf x}}}
\newcommand{\yhat}{{\hat{\bf y}}}
\newcommand{\zhat}{{\hat{\bf z}}}

\newcommand{\fstar}{{f^*}}
\newcommand{\fstarp}{(\fstar){\vphantom{\bigr)}}^{\mkern-1mu\prime}}
\newcommand{\fistar}{{f_i^*}}
\newcommand{\Fstar}{{F^*}}
\newcommand{\Gstar}{{G^*}}
\newcommand{\fpstar}{{f'^*}}
\newcommand{\gstar}{{g^*}}

\newcommand{\gpstar}{{g'^*}}
\newcommand{\hstar}{{h^*}}

\newcommand{\rstar}{{r^*}}
\newcommand{\sstar}{{s^*}}

\newcommand{\gtil}{{}\mkern2mu\tilde{\mkern-2mu g\mkern1mu}\mkern-1mu}%

\newcommand{\sbar}{\overline{\sS}}
\newcommand{\ubar}{\overline{\uu}}
\newcommand{\vbar}{\overline{\vv}}
\newcommand{\wbar}{\overline{\ww}}
\newcommand{\xbar}{\overline{\xx}}%
\newcommand{\bbar}{\overline{\bb}}
\newcommand{\cbar}{\overline{\cc}}
\newcommand{\dbar}{\overline{\dd}}
\newcommand{\ebar}{\overline{\ee}}

\newcommand{\barc}{\bar{c}}
\newcommand{\bare}{\bar{e}}

\newcommand{\barw}{\overlineKernIt{w}}
\newcommand{\barx}{\overlineKernIt{x}}
\newcommand{\bary}{\overlineKernIt{y}}
\newcommand{\barz}{\bar{z}}

\newcommand{\ybar}{\overline{\yy}}
\newcommand{\ybart}{\ybar_{\!t}}
\newcommand{\zbar}{\overline{\zz}}
\newcommand{\Cbar}{\overlineKernS{C}}

\newcommand{\Jbar}{\overlineKernf{J}}
\newcommand{\Kbar}{\overlineKernS{K}}
\newcommand{\Khat}{\hat{K}}
\newcommand{\Qbar}{\overlineKernIt{Q}}
\newcommand{\Lbar}{\overlineKernIt{L}}
\newcommand{\Mbar}{\overlineKernIt{M}}
\newcommand{\KLbar}{\overline{K\cap L}}
\newcommand{\Rbar}{\overlineKernR{R}}
\newcommand{\Rbari}{\overlineKernR{R}_i}

\newcommand{\AFstar}{(\A F)^{*}}
\newcommand{\AFdub}{(\A F)^{*\dualstar}}

\newcommand{\negKern}{\mkern-1.5mu}
\newcommand{\posKern}{\mkern1.5mu}
\newcommand{\posSkip}{\mskip1.5mu plus 2mu minus 1.5mu}
\newcommand{\overlineKernIt}[1]{{}\mkern2mu\overline{\mkern-2mu#1}}
\newcommand{\overlineKerng}[1]{{}\mkern2mu\overline{\mkern-2mu#1\mkern1mu}\mkern-1mu}
\newcommand{\overlineKerngspace}[1]{{}\mkern2mu\overline{\mkern-2mu#1\mkern1mu}}
\newcommand{\overlineKernp}[1]{{}\mkern1mu\overline{\mkern-1mu#1\mkern1mu}\mkern-1mu}
\newcommand{\overlineKernY}[1]{{}\overline{#1}}
\newcommand{\overlineKernR}[1]{{}\mkern3mu\overline{\mkern-3mu#1\mkern2mu}\mkern-2mu}
\newcommand{\overlineKernS}[1]{{}\mkern3mu\overline{\mkern-3mu#1\mkern1.5mu}\mkern-1.5mu}
\newcommand{\overlineKernO}[1]{{}\mkern3mu\overline{\mkern-3mu#1\mkern0.5mu}\mkern-0.5mu}
\newcommand{\overlineKernf}[1]{{}\mkern4.5mu\overline{\mkern-4.5mu#1\mkern1.5mu}\mkern-1.5mu}

\newcommand{\Abar}{\overlineKernf{A}}
\newcommand{\Abarinv}{\bigParens{\mkern-1.5mu\Abar\mkern2mu}^{\mkern-1.5mu-1}}%
\newcommand{\Alininv}{A^{-1}}
\newcommand{\Fbar}{\overlineKernR{F}}
\newcommand{\Sbar}{\overlineKernS{S}}

\newcommand{\Pbar}{\overlineKernR{P}}

\newcommand{\Ubar}{\overlineKernIt{U}}
\newcommand{\Nbar}{\overlineKernIt{N}}
\newcommand{\PNbar}{\overline{P\cap N}}

\newcommand{\Zbar}{\overlineKernS{Z}}
\newcommand{\Xbar}{\overlineKernIt{X}}
\newcommand{\Ybar}{\overlineKernY{Y}}

\newcommand{\Xbarp}{\overline{X'}}

\newcommand{\Xbarpi}{\overlineKernIt{X'_i}}

\newcommand{\sS}{{\bf s}}
\newcommand{\pp}{\mathbf{p}}
\newcommand{\qq}{{\bf q}}
\newcommand{\rr}{{\bf r}}
\newcommand{\uu}{\mathbf{u}}
\newcommand{\vv}{{\bf v}}
\newcommand{\ww}{{\bf w}}
\newcommand{\yy}{\mathbf{y}}

\newcommand{\zz}{{\bf z}}
\newcommand{\A}{\mathbf{A}}
\newcommand{\Ainv}{\A^{\negKern-1}}
\newcommand{\Ascalar}{A}
\newcommand{\transAscalar}{\trans{A}}
\newcommand{\B}{{\bf B}}

\newcommand{\PP}{\mathbf{P}}
\newcommand{\QQ}{\mathbf{Q}}
\newcommand{\QQy}{\QQ_2}
\newcommand{\QQz}{\QQ_1}
\newcommand{\RR}{\mathbf{R}}
\newcommand{\Rinv}{{\RR^{-1}}}
\newcommand{\Ss}{{\bf S}}
\newcommand{\VV}{\mathbf{V}}
\newcommand{\VVsub}[1]{\VV_{\negKern#1}}

\newcommand{\WW}{{\bf W}}
\newcommand{\ZZ}{{\bf Z}}
\newcommand{\Adag}{\A\negKern^\pseudinvsym}
\newcommand{\Bdag}{\B^\pseudinvsym}
\newcommand{\VVdag}{\VV^\pseudinvsym}
\newcommand{\VVtrans}{\trans{\VV}}
\newcommand{\VVdagtrans}{\trans{(\VV^\pseudinvsym)}}
\newcommand{\Bpseudoinv}{\B^\pseudinvsym}
\newcommand{\xpar}{\xx^L}%
\newcommand{\PPf}{\PP\mkern-1.5mu f}
\newcommand{\PPi}{\PPsub{i}}
\newcommand{\PPsub}[1]{\PP_{\mkern-2.5mu #1}}
\newcommand{\PPx}{{\PP\mkern-1.5mu}_1}%
\newcommand{\lmPPx}{P^{\vphantom{*}}_{\negKern 1}}

\newcommand{\lmPPy}{P^{\vphantom{*}}_{\negKern 2}}
\newcommand{\lmPPyinv}{P^{-1}_{\negKern 2}}
\newcommand{\PPy}{{\PP\mkern-1.5mu}_{2}}%
\newcommand{\PPxlow}{\PPx^{\vphantom{\top}}}
\newcommand{\PPylow}{\PPy^{\vphantom{\top}}}
\newcommand{\PPxE}{\PPx\mkern-0.5mu E}
\newcommand{\Qx}{\QQ_1}
\newcommand{\Qy}{\QQ_2}

\newcommand{\Pfext}{\overline{\PPf}}
\newcommand{\PfPext}{\overline{(\PPf)\PP}}

\newcommand{\zerov}[1]{\zero_{#1}}
\newcommand{\zerovec}{\zerov{0}}
\newcommand{\zeromat}[2]{\zero_{{#1}\times {#2}}}

\newcommand{\Iden}{\mathbf{I}}
\newcommand{\Idn}[1]{\Iden_{#1}}
\newcommand{\Idnn}{\Idn{n}}

\newcommand{\ellbar}{\overlineKernS{\ell}}

\newcommand{\irredIndexSet}{persistent index set\xspace}
\newcommand{\irredIndices}{persistent indices\xspace}
\newcommand{\irredIndicesFor}{persistent indices for\xspace}
\newcommand{\irredIndexSetFor}{persistent index set for\xspace}
\newcommand{\outA}{\mathsf{A}}
\newcommand{\outB}{\mathsf{B}}
\newcommand{\outC}{\mathsf{C}}
\newcommand{\outD}{\mathsf{D}}

\DeclareMathOperator{\irred}{\mathopfont{pers}}
\newcommand{\hardcore}[1]{\irred #1}
\newcommand{\indset}{I}
\newcommand{\uset}[1]{\uu[#1]}

\newcommand{\topo}{{\cal T}}

\newcommand{\expex}{\mathop{\overline{{\exp}\ministrut}}\nolimits}
\newcommand{\invex}{\mathop{\overline{\mathrm{\i nv}\ministrut}}\nolimits}
\newcommand{\logex}{\mathop{\overline{\ln}}\nolimits}
\DeclareMathOperator{\sqp}{\mathopfont{sqp}}

\newcommand{\Rtilde}{\tilde{\RR}}

\newcommand{\slinmapA}{A}%
\newcommand{\alinmapA}{\Abar}%
\newcommand{\slinmapAinv}{\Alininv}%
\newcommand{\alinmapAinv}{\Abarinv}%
\newcommand{\linmapalt}{F}

\newcommand{\Hcomp}{H^c}
\newcommand{\clHcomp}{\overline{H^c}}
\newcommand{\Ncomp}{{N^c}}
\newcommand{\Uv}{{U_{\vv}}}
\newcommand{\Uzero}{{U_{\zero}}}
\newcommand{\Unormv}{{U_{\vv/\norm{\vv}}}}

\newcommand{\US}{{U_{S}}}
\newcommand{\UV}{{U_{V}}}
\newcommand{\UVj}{{U_{V_j}}}

\newcommand{\rats}{{\mathbb{Q}}}
\newcommand{\nats}{{\mathbb{N}}}
\newcommand{\R}{{\mathbb{R}}}
\newcommand{\Rn}{\R^n}

\newcommand{\Rm}{\R^m}
\newcommand{\Rmn}{{\R^{m\times n}}}
\newcommand{\Rnm}{{\R^{n\times m}}}
\newcommand{\Rnk}{{\R^{n\times k}}}
\newcommand{\Rnp}{{\R^{n+1}}}
\newcommand{\Rnpnp}{{\R^{(n+1)\times (n+1)}}}

\newcommand{\Rk}{{\R^k}}

\newcommand{\Qn}{{\rats^n}}
\newcommand{\Rpos}{\mathbb{R}_{\geq 0}}

\newcommand{\Rstrictpos}{\mathbb{R}_{> 0}}
\newcommand{\Rneg}{\mathbb{R}_{\leq 0}}
\newcommand{\Rstrictneg}{\mathbb{R}_{< 0}}

\newcommand{\Rextpos}{{\Rpos\cup\{\oms\}}}
\newcommand{\Rextstrictpos}{{\Rstrictpos\cup\{\oms\}}}

\newcommand{\trans}[1]{#1^\top}
\newcommand{\transKern}[1]{#1^\top\!}
\newcommand{\transk}[1]{#1^\top\negKern}

\newcommand{\amatu}{\mathbf{A}}

\newcommand{\amatuinv}{\mathbf{B}}
\newcommand{\transamatuinv}{{\trans{\amatuinv}}}

\DeclareMathOperator{\resc}{\mathopfont{rec}}
\newcommand{\represc}[1]{{(\resc{#1})^{\triangle}}}
\newcommand{\rescperp}[1]{{(\resc{#1})^{\bot}}}

\newcommand{\rescbar}[1]{\overline{\resc{#1}}}%
\newcommand{\rescpol}[1]{{(\resc{#1})^{\circ}}}
\newcommand{\rescdubpol}[1]{{(\resc{#1})^{\circ\circ}}}

\DeclareMathOperator{\conssp}{\mathopfont{cons}}

\DeclareMathOperator{\ri}{\mathopfont{ri}}
\newcommand{\ric}[1]{{\ri(\resc{#1})}}

\DeclareMathOperator{\arescone}{\mathopfont{rec}}
\newcommand{\aresconeF}{{\arescone F}}
\newcommand{\aresconef}{{\arescone \fext}}
\newcommand{\aresconefp}{{\arescone \fpext}}
\newcommand{\aresconeg}{{\arescone \gext}}

\newcommand{\aresconegsub}[1]{{\arescone \gext_{#1}}}

\newcommand{\perpresf}{{(\aresconef)^{\bot}}}

\newcommand{\perpresg}{{(\aresconeg)^{\bot}}}
\newcommand{\perpresgsub}[1]{{(\aresconegsub{#1})^{\bot}}}
\newcommand{\perperpresf}{{(\aresconef)^{\bot\bot}}}

\newcommand{\fullshad}[1]{{#1^{\diamond}}}
\newcommand{\fullshadnoparen}[1]{{#1^{\diamond}}}

\newcommand{\fullshadfbar}{\overlineKernIt{f^{\diamond}}}

\DeclareMathOperator{\contset}{\mathopfont{cont}}
\DeclareMathOperator{\unimin}{\mathopfont{univ}}
\newcommand{\contsetf}{\contset{\fext}}

\DeclareMathOperator{\intr}{\mathopfont{int}}
\newcommand{\intdom}[1]{{\intr(\dom #1)}}
\newcommand{\ball}{B}%

\newcommand{\calB}{\mathcal{B}}

\newcommand{\calC}{{\cal C}}

\newcommand{\calS}{{\cal S}}

\newcommand{\polar}[1]{#1^{\circ}}
\newcommand{\dubpolar}[1]{{#1}^{\circ\circ}}

\newcommand{\apol}[1]{#1^{\bar{\circ}}}
\newcommand{\apolconedomfstar}{\apol{(\cone(\dom{\fstar}))}}
\newcommand{\apolconedomFstar}{\apol{(\cone(\dom{\Fstar}))}}
\newcommand{\apolconedomfpstar}{\apol{(\cone(\dom{\fpstar}))}}
\newcommand{\apolslopesf}{\apol{(\slopes{f})}}
\newcommand{\polapol}[1]{{#1}^{\circ\bar{\circ}}}
\newcommand{\polpol}[1]{{#1}^{\circ\circ}}
\newcommand{\apolpol}[1]{{#1}^{\bar{\circ}\circ}}

\newcommand{\Jpol}{\polar{J}}
\newcommand{\Kpol}{\polar{K}}
\newcommand{\Japol}{\apol{J}}
\newcommand{\Kapol}{\apol{K}}
\newcommand{\Kpolbar}{\clbar{\Kpol}}

\newcommand{\Jpolapol}{\polapol{J}}

\newcommand{\Kapolpol}{\apolpol{K}}
\newcommand{\Lpol}{\polar{L}}
\newcommand{\Lpolbar}{\clbar{\Lpol}}

\newcommand{\Khatpol}{\polar{\Khat}}
\newcommand{\Khatdubpol}{\dubpolar{\Khat}}

\newcommand{\aperp}[1]{#1^{\bar{\bot}}}
\newcommand{\perpaperp}[1]{#1^{\bot\bar{\bot}}}
\newcommand{\perperp}[1]{#1^{\bot\bot}}
\newcommand{\aperperp}[1]{#1^{\bar{\bot}\bot}}

\newcommand{\Lperpaperp}{\perpaperp{L}}
\newcommand{\Laperp}{\aperp{L}}
\newcommand{\Lapol}{\apol{L}}
\newcommand{\Uaperp}{\aperp{U}}
\newcommand{\Uperpaperp}{\perpaperp{U}}
\newcommand{\Saperp}{\aperp{S}}
\newcommand{\Uaperperp}{\aperperp{U}}

\newcommand{\ncone}[2]{N_{#2}({#1})}
\newcommand{\ancone}[2]{N\mkern-6.5mu\widetilde{\wtstrut}\mkern5mu_{#2}({#1})}

\DeclareMathOperator{\cl}{\mathopfont{cl}}
\DeclareMathOperator{\lsc}{\mathopfont{lsc}}

\newcommand{\clbar}{\overline}

\newcommand{\clmx}{\cl_{\mathrm{m}}}

\DeclareMathOperator{\cone}{\mathopfont{cone}}

\newcommand{\clcone}[1]{{\overline{\cone{#1}}}}

\newcommand{\repcl}[1]{{#1^{\triangle}}}

\newcommand{\rpair}[2]{\langle #1, #2\rangle}
\newcommand{\rpairf}[2]{#1, #2}
\newcommand{\mtuple}[1]{\langle #1\rangle}

\let\epi=\epiOp
\let\dom=\domOp
\let\relint=\relintOp

\DeclareMathOperator{\hypo}{\mathopfont{hypo}}
\DeclareMathOperator{\ray}{\mathopfont{ray}}

\DeclareMathOperator{\aray}{\mathopfont{ray}_{\star}}
\DeclareMathOperator{\acone}{\mathopfont{cone}_{\star}}
\DeclareMathOperator{\oconich}{\widetilde{\mathopfont{cone}}_{\star}}
\DeclareMathOperator{\acolspace}{\mathopfont{col}_{\star}}
\DeclareMathOperator{\aspan}{\mathopfont{span}_{\star}}
\DeclareMathOperator{\ahlineop}{\mathopfont{hline}_{\star}}

\DeclareMathOperator{\hlineop}{\mathopfont{hline}}
\newcommand{\hfline}[2]{\hlineop({#1},{#2})}
\newcommand{\ahfline}[2]{\ahlineop({#1},{#2})}

\newcommand{\epibar}[1]{\clmx(\epi #1)}
\newcommand{\epibarbar}[1]{\overline{\epi #1}}

\newcommand{\clepi}{\epibarbar}
\newcommand{\cldom}[1]{\overline{\dom #1}}
\newcommand{\cldomfext}{\cldom{\efspace}}
\newcommand{\clepifext}{\clepi{\efspace}}
\newcommand{\clepif}{\overline{\epi f}}

\DeclareMathOperator{\csmOp}{\mathopfont{csm}}
\newcommand{\cosmspn}{\csmOp \Rn}
\DeclareMathOperator{\rwdirOp}{\mathopfont{dir}}
\newcommand{\rwdir}{\rwdirOp}
\newcommand{\quof}{p}
\newcommand{\quofinv}{p^{-1}}
\newcommand{\xrw}{\hat{\xx}}
\newcommand{\ucball}{B}
\newcommand{\coshom}{\sigma}

\newcommand{\fext}{\ef}%
\newcommand{\fpext}{\efp}%
\newcommand{\fpextsub}[1]{\efpsub{#1}}%
\newcommand{\fextsub}[1]{\efsub{#1}}
\newcommand{\gext}{\eg}

\newcommand{\gpext}{\overline{g'}}
\newcommand{\hext}{\eh}%
\newcommand{\hpext}{\overline{h'}}
\newcommand{\dext}{\overline{d}}
\newcommand{\rext}{\overlineKernIt{r}}

\newcommand{\pext}{\ep}
\newcommand{\sext}{\overlineKerng{s}}
\newcommand{\lscfext}{\overline{\lsc f}}%
\newcommand{\lscfextf}{(\lscfext)}

\newcommand{\Af}{({\A f})}
\newcommand{\Afext}{(\overline{\A f})}
\newcommand{\Afextnp}{\overline{\A f}}
\newcommand{\fAext}{(\overline{\fA})}
\newcommand{\fAextnp}{\overline{\fA}}
\newcommand{\gAext}{(\overline{g \A})}

\newcommand{\ufext}{(\overline{\utransf})}
\newcommand{\ufextnp}{\overline{\utransf}}
\newcommand{\ufextstar}{\regParens{\overline{\utransf}}^{*}}
\newcommand{\ufstar}{(\utransf)^{*}}
\newcommand{\lamfext}{(\overline{\lambda f})}
\newcommand{\fplusgext}{(\overline{f+g})}

\newcommand{\Gofextnp}{\overline{G\circ f}}
\newcommand{\fgAext}{(\overline{f + g \A})}

\newcommand{\fAgstar}{\fstar[-\transA]+\gstar}

\newcommand{\vfext}{(\overline{v f})}
\newcommand{\sfprodext}[2]{(\overline{\sfprod{#1}{#2}})}

\newcommand{\fomin}{g}
\newcommand{\fominext}{\eg}
\newcommand{\normfcn}{\ell}
\newcommand{\normfcnext}{\ellbar}

\newcommand{\dualstar}{\bar{*}}
\newcommand{\phantomdualstar}{\vphantom{\dualstar}}

\newcommand{\fdubs}{f^{**}}
\newcommand{\fdub}{f^{*\dualstar}}
\newcommand{\fidub}{f_i^{*\dualstar}}

\newcommand{\gdub}{g^{*\dualstar}}
\newcommand{\hdub}{h^{*\dualstar}}

\newcommand{\psistar}{\psi^{*}}
\newcommand{\psistarb}{\psi^{\dualstar}}
\newcommand{\psineg}{\psi_{<0}}
\newcommand{\psipos}{\psi_{>0}}
\newcommand{\rhoneg}{\rho_{<0}}
\newcommand{\rhopos}{\rho_{>0}}
\newcommand{\rhostar}{\rho^{*}}
\newcommand{\rhostarext}{\overlineKernIt{\rho^*}}

\newcommand{\rhostarb}{\rho^{\dualstar}}
\newcommand{\Fdub}{F^{*\dualstar}}
\newcommand{\Gdub}{G^{*\dualstar}}
\newcommand{\affuv}{A_{\uu,v}}
\newcommand{\daffgen}[2]{\xi_{#1,#2}}
\newcommand{\daffsxy}{\daffgen{\xx}{y}}
\newcommand{\daffxy}{\daffgen{\xbar}{y}}
\newcommand{\fextstar}{{\efspace}^{\mkern0.5mu*}\negKern}
\newcommand{\gextstar}{{\egspace}^{\mkern0.5mu*}}
\newcommand{\fextdub}{{\efspace}^{\mkern0.5mu*\dualstar}}

\DeclareMathOperator{\fun}{\mathopfont{lenv}}%
\DeclareMathOperator{\funstar}{\mathopfont{lsup}}%
\newcommand{\settofcn}[1]{\fun #1}%
\newcommand{\settofcnE}{\settofcn{E}}
\newcommand{\settofcnEf}{(\settofcn{E})}

\newcommand{\settofcnEpf}{(\settofcn{E'})}
\newcommand{\settofcnEstar}{(\settofcnE)^*}%
\newcommand{\settofcnEdub}{(\settofcnE)^{*\dualstar}}%
\newcommand{\Esstar}{\funstar E}%
\newcommand{\Esstarf}{(\funstar E)}%
\newcommand{\Epsstar}{\funstar E'}%
\newcommand{\Epsstarf}{(\funstar E')}%
\newcommand{\Esdub}{(\Esstar)^{\dualstar}}%
\newcommand{\Esdubf}{(\Esstar)^{\dualstar}}%
\newcommand{\Epsdubf}{(\Epsstar)^{\dualstar}}%

\newcommand{\resfcn}[2]{{#1}|_{#2}}

\newcommand{\dlin}[1]{\xi_{#1}}
\newcommand{\dlinxbuz}{\dlin{\xbar,\uu_0,\beta}}
\newcommand{\dlinxbu}{\dlin{\xbar,\uu,\beta}}

\newcommand{\trivgen}[1]{\breve{#1}}
\newcommand{\ftriv}{\trivgen{f}}
\newcommand{\ftrivstar}{{\ftriv{}}^{\mkern1mu\smash{*}}\negKern\negKern}

\newcommand{\Fminussym}[1]{\raisebox{0.35pt}{$\scriptscriptstyle[$}%
                           #1%
                           \raisebox{0.35pt}{$\scriptscriptstyle]$}}
\newcommand{\Fminusgen}[2]{{#1}_{\Fminussym{#2}}}

\newcommand{\Fminusu}{\Fminusgen{F}{\uu}}
\newcommand{\Gminusu}{\Fminusgen{G}{\uu}}
\newcommand{\FminusU}{\Fminusgen{F}{u}}

\newcommand{\fminussym}[1]{\raisebox{0.35pt}{$\scriptscriptstyle[$}%
                           #1%
                           \raisebox{0.35pt}{$\scriptscriptstyle]$}}
\newcommand{\fminusgen}[2]{#1_{\fminussym{#2}}}

\newcommand{\aminusu}{\fminusgen{a}{\uu}}
\newcommand{\fminusu}{\fminusgen{f}{\uu}}
\newcommand{\fminusw}{\fminusgen{f}{\ww}}
\newcommand{\fminusa}{\fminusgen{f}{\aaa}}
\newcommand{\fminusU}{\fminusgen{f}{u}}
\newcommand{\fminusuext}{\efsub{\fminussym{\uu}}}
\newcommand{\fminuswext}{\efsub{\fminussym{\ww}}}
\newcommand{\fminusaext}{\efsub{\fminussym{\aaa}}}
\newcommand{\fextminusu}{\fminusgen{(\fext)}{\uu}}
\newcommand{\fextminusU}{\fminusgen{(\fext)}{u}}
\newcommand{\fminusUext}{\efsub{\fminussym{u}}}
\newcommand{\fminusutriv}{\fminusgen{\ftriv}{\uu}}
\newcommand{\ftrivminusu}{\fminusgen{(\ftriv)}{\uu}}
\newcommand{\lscftrivminusu}{\fminusgen{(\lsc\ftriv)}{\uu}}
\newcommand{\fminusustar}{\fminusu^{*}}

\newcommand{\fminusudub}{\fminusu^{*\dualstar}}

\newcommand{\fpminusminusu}{f'_{\fminussym{-\uu}}}
\newcommand{\fpminusminusuext}{\efpsub{\fminussym{-\uu}\kern-2.9ex}\kern2.9ex}
\newcommand{\gminusu}{\fminusgen{g}{\uu}}
\newcommand{\gminusU}{\fminusgen{g}{u}}
\newcommand{\gminusuext}{\fminusgen{\gext}{\uu}}
\newcommand{\gminusUext}{\fminusgen{\gext}{u}}
\newcommand{\hminusu}{\fminusgen{h}{\uu}}
\newcommand{\hminusw}{\fminusgen{h}{\ww}}
\newcommand{\hminusU}{\fminusgen{h}{u}}

\newcommand{\hminusuext}{\fminusgen{\hext}{\uu}}
\newcommand{\hminuswext}{\fminusgen{\hext}{\ww}}
\newcommand{\hminusUext}{\fminusgen{\hext}{u}}

\newcommand{\fminussub}[2]{f_{{#1},\fminussym{#2}}}
\newcommand{\fminussubext}[2]{\efsub{{#1},\fminussym{#2}}}
\newcommand{\fminususubd}[1]{\fminussub{#1}{\uu_{#1}}}
\newcommand{\fminususubdext}[1]{\fminussubext{#1}{\uu_{#1}}}
\newcommand{\fminususub}[1]{\fminussub{#1}{\uu}}
\newcommand{\fminususubext}[1]{\fminussubext{#1}{\uu}}
\newcommand{\fminusUsubext}[1]{\fminussubext{#1}{u}}

\newcommand{\fminusgenstar}[2]{{#1}^*_{\fminussym{#2}}}
\newcommand{\fminussubstar}[2]{f^*_{{#1},\fminussym{#2}}}
\newcommand{\fminusextsubstar}[2]{\efsubstar{#1,\fminussym{#2}}}%
\newcommand{\fiminusustar}{\fminussubstar{i}{\uu}}
\newcommand{\fiminusuextstar}{\fminusextsubstar{i}{\uu}}
\newcommand{\hminusustar}{\fminusgenstar{h}{\uu}}
\newcommand{\hminusuextstar}{\fminusgenstar{\hext}{\uu}}

\newcommand{\fVgen}[2]{{#1}^{\mkern1mu\vdash}_{#2}}

\newcommand{\fVgenstar}[2]{({#1}^{\mkern1mu\vdash}_{#2})^*}

\newcommand{\fV}{\fVgen{f}{\VV}}

\newcommand{\fVstar}{\fVgenstar{f}{\VV}}

\newcommand{\fv}{\fVgen{f}{\vv}}

\newcommand{\acolspaceVV}{\acolspace\!\VV}

\newcommand{\asubsym}{\partial}
\newcommand{\asubdifplain}[1]{\asubsym{#1}}
\newcommand{\asubdif}[2]{\asubsym{#1}({#2})}
\newcommand{\asubdiflagranglmext}[1]{\asubsym\lagranglmext{#1}}
\newcommand{\asubdiflagrangext}[2]{\asubsym\lagrangext{#1}{#2}}

\newcommand{\asubdiffext}[1]{\asubdif{\fext}{#1}}
\newcommand{\asubdifftriv}[1]{\asubdif{\ftriv}{#1}}

\newcommand{\asubdifAfext}[1]{\asubdif{\Afext}{#1}}
\newcommand{\asubdifgext}[1]{\asubdif{\gext}{#1}}
\newcommand{\asubdifhext}[1]{\asubdif{\hext}{#1}}
\newcommand{\asubdifrext}[1]{\asubdif{\rext}{#1}}
\newcommand{\asubdifsext}[1]{\asubdif{\sext}{#1}}
\newcommand{\asubdiffplusgext}[1]{\asubdif{\fplusgext}{#1}}
\newcommand{\asubdiffextsub}[2]{\asubdif{\fextsub{#1}}{#2}}

\newcommand{\asubdifsextsub}[2]{\asubdif{\sext_{#1}}{#2}}
\newcommand{\asubdiffminusuext}[1]{\asubdif{\fminusuext}{#1}}
\newcommand{\asubdiffminususubext}[2]{\asubdif{\fminususubext{#1}}{#2}}

\newcommand{\asubdifhminusuext}[1]{\asubdif{\hminusuext}{#1}}
\newcommand{\asubdifinddomfext}[1]{\asubdif{\inddomfext}{#1}}
\newcommand{\asubdifindaS}[1]{\asubdif{\indfa{S}}{#1}}
\newcommand{\asubdifindext}[2]{\asubdif{\indgenext{#1}}{#2}}
\newcommand{\asubdifindepifext}[1]{\asubdif{\indepifext}{#1}}
\newcommand{\asubdifindsext}[1]{\asubdif{\indsext}{#1}}

\newcommand{\asubdifgAext}[1]{\asubdif{\gAext}{#1}}
\newcommand{\asubdifF}[1]{\asubdif{F}{#1}}
\newcommand{\asubdifG}[1]{\asubdif{G}{#1}}

\newcommand{\asubdifpsistarb}[1]{\asubdif{\psistarb}{#1}}

\newcommand{\asubdifsfprodext}[3]{\asubdif{\sfprodext{#1}{#2}}{#3}}

\newcommand{\asubdifmuhext}[1]{\asubdif{\muhext}{#1}}

\newcommand{\asubdifNegogiext}[1]{\asubdif{\Negogiext}{#1}}

\newcommand{\asubdifindzohjext}[1]{\asubdif{\indzohjext}{#1}}
\newcommand{\asubdifindzext}[1]{\asubdif{\indzext}{#1}}
\newcommand{\asubdiflscF}[1]{\asubdif{(\lsc F)}{#1}}
\newcommand{\asubdiflscftriv}[1]{\asubdif{(\lsc\ftriv)}{#1}}

\newcommand{\asubdifpext}[1]{\asubdif{\pext}{#1}}
\newcommand{\asubdifvfext}[1]{\asubdif{\vfext}{#1}}

\newcommand{\wtstrut}{\vphantom{i}}
\newcommand{\adsubsym}{\partial\mkern-3.5mu\widetilde{\wtstrut}\mkern3.5mu}

\newcommand{\adsubdifplain}[1]{{\adsubsym{#1}}}
\newcommand{\adsubdif}[2]{{\adsubsym{#1}({#2})}}
\newcommand{\adsubdiffstar}[1]{\adsubdif{\fstar}{#1}}
\newcommand{\adsubdifgstar}[1]{\adsubdif{\gstar}{#1}}
\newcommand{\adsubdifFstar}[1]{\adsubdif{\Fstar}{#1}}

\newcommand{\adsubdifpsi}[1]{\adsubdif{\psi}{#1}}
\newcommand{\adsubdifpsidub}[1]{\adsubdif{\psidub}{#1}}
\newcommand{\adsubdiffminusustar}[1]{\adsubdif{\fminusustar}{#1}}
\newcommand{\adsubdifEpstar}[1]{\adsubdif{\Epsstarf}{#1}}
\newcommand{\psidubs}{\psi^{**}}
\newcommand{\psidub}{\psi^{\dualstar*}}

\newcommand{\rhodub}{\rho^{\dualstar*}}
\newcommand{\rhodubs}{\rho^{**}}
\newcommand{\psiuudub}{(\psi\uu)^{\dualstar*}}

\newcommand{\bpartial}{\asubsym}
\newcommand{\basubdifplain}[1]{{\bpartial{#1}}}
\newcommand{\basubdif}[2]{\bpartial{#1}({#2})}
\newcommand{\basubdiflagranglmext}[1]{\bpartial\lagranglmext{#1}}

\newcommand{\basubdiffext}[1]{\basubdif{\fext}{#1}}
\newcommand{\basubdiffAext}[1]{\basubdif{\fAext}{#1}}
\newcommand{\basubdifAfext}[1]{\basubdif{\Afext}{#1}}
\newcommand{\basubdifgext}[1]{\basubdif{\gext}{#1}}
\newcommand{\basubdifhext}[1]{\basubdif{\hext}{#1}}

\newcommand{\basubdiffextsub}[2]{\basubdif{\fextsub{#1}}{#2}}

\newcommand{\basubdifinddomfext}[1]{\basubdif{\inddomfext}{#1}}

\newcommand{\basubdifsfprodext}[3]{\basubdif{\sfprodext{#1}{#2}}{#3}}
\newcommand{\basubdiflamgiext}[1]{\basubdifsfprodext{\lambda_i\negKern}{\posKern g_i}{#1}}
\newcommand{\basubdiflamgext}[1]{\basubdifsfprodext{\lambda}{\posKern g}{#1}}
\newcommand{\basubdifmuhext}[1]{\basubdif{\muhext}{#1}}
\newcommand{\basubdifmuhjext}[1]{\basubdif{\muhjext}{#1}}

\newcommand{\basubdifnegext}[1]{\basubdif{\negfext}{#1}}
\newcommand{\basubdifNegogext}[1]{\basubdif{\Negogext}{#1}}
\newcommand{\basubdifNegogiext}[1]{\basubdif{\Negogiext}{#1}}

\newcommand{\basubdifindzohext}[1]{\basubdif{\indzohext}{#1}}
\newcommand{\basubdifindzohextTight}[1]{\basubdif{\indzohextTight}{#1}}
\newcommand{\basubdifindzohjext}[1]{\basubdif{\indzohjext}{#1}}
\newcommand{\basubdifindzext}[1]{\basubdif{\indzext}{#1}}
\newcommand{\basubdifpext}[1]{\basubdif{\pext}{#1}}

\newcommand{\gradf}{\nabla\negKern f}
\newcommand{\gradh}{\nabla h}

\newcommand{\dder}[3]{#1'(#2;#3)}
\newcommand{\dderf}[2]{\dder{f}{#1}{#2}}
\newcommand{\dderF}[2]{\dder{F}{#1}{#2}}
\newcommand{\dderfext}[2]{{\efspace}^{\smash{\prime}\vphantom{*}}\!(#1;#2)}
\newcommand{\dderh}[2]{\dder{h}{#1}{#2}}
\newcommand{\dderpsi}[2]{\dder{\psi}{#1}{#2}}
\newcommand{\rightlim}[2]{{#1}\rightarrow {#2}^{+}}

\newcommand{\shadexp}[1]{[{#1}]}
\newcommand{\genshad}[2]{{#1^{\shadexp{#2}}}}
\newcommand{\fshad}[1]{\genshad{f}{#1}}
\newcommand{\fpshad}[1]{\genshad{f'}{#1}}
\newcommand{\gshad}[1]{\genshad{g}{#1}}
\newcommand{\gishadvi}{\genshad{g_{i-1}}{\limray{\vv_i}}}

\newcommand{\fshadd}{\fshad{\ebar}}
\newcommand{\fshadv}{\fshad{\limray{\vv}}}

\newcommand{\fshadvb}{\fshad{\limray{\vbar}}}
\newcommand{\fpshadv}{\fpshad{\limray{\vv}}}

\newcommand{\fshadVo}{\fshad{\VV\omm}}
\newcommand{\fshadext}[1]{\overline{\fshad{#1}}}
\newcommand{\fshadextd}{\fshadext{\ebar}}
\newcommand{\fminusushadd}{\genshad{\fminusu}{\ebar}}

\newcommand{\finclset}{\mathcal{M}_{n+1}}
\newcommand{\homf}{\mu}
\newcommand{\homfinv}{\mu^{-1}}

\newcommand{\oms}{\omega}
\newcommand{\omsf}[1]{\oms {#1}}
\newcommand{\lmset}[1]{{\oms{#1}}}
\newcommand{\omm}{\boldsymbol{\omega}}%
\newcommand{\ommsub}[1]{\omm_{#1}}
\newcommand{\ommk}{\ommsub{k}}
\newcommand{\limray}[1]{\omsf{#1}}
\newcommand{\limrays}[1]{{[{#1}] \omm}}

\newcommand{\finv}{F^{-1}}

\newcommand{\uperp}{{\uu^\bot}}

\newcommand{\xperp}{{\xx^\bot}}

\newcommand{\xperpt}{{\xx_t^\bot}}

\newcommand{\ebarperp}{\ebar^{\bot}}

\newcommand{\xbarperp}{\xbar^{\bot}}

\newcommand{\xbarperpt}{\xbar_t^{\bot}}
\newcommand{\ybarperp}{\ybar^{\bot}}

\newcommand{\norm}[1]{\lVert#1\rVert}
\newcommand{\norms}[1]{\lVert#1\rVert_2}
\newcommand{\normp}[1]{\lVert#1\rVert_p}
\newcommand{\normq}[1]{\lVert#1\rVert_q}
\DeclareMathOperator{\colspace}{\mathopfont{col}}
\DeclareMathOperator{\spn}{\mathopfont{span}}
\newcommand{\spnfin}[1]{{\spn\{{#1}\}}}
\newcommand{\comColV}{(\colspace\VV)^\perp}

\DeclareMathOperator{\rspan}{\mathopfont{rspan}}
\newcommand{\rspanset}[1]{\rspan\set{#1}}
\newcommand{\rspanxbar}{\rspanset\xbar}
\newcommand{\rspanxbarperp}{\rspanset\xbarperp}

\DeclareMathOperator{\columns}{\mathopfont{columns}}

\newcommand{\Vxbar}{\VV_{\!\xbar}}
\newcommand{\qqxbar}{\qq_{\mkern1mu\xbar}}
\newcommand{\kxbar}{k_{\mkern1mu\xbar}}
\newcommand{\Rxbar}{R_{\mkern1mu\xbar}}
\newcommand{\Rxbari}{R_{\mkern1mu\xbar_i}}

\renewcommand{\eqref}[1]{Eq.~\textup{(\ref{#1})}}

\newcommand{\partsubref}[2]{(\ref{#1})(\ref{#2})}

\newcommand{\countsetgen}[1]{\set{#1}}
\newcommand{\countset}[1]{\countsetgen{#1_t}}

\newcommand{\refequiv}[3]{\ref{#1}(\ref{#2},\ref{#3})}
\newcommand{\Crefequiv}[3]{\Cref{#1}(\ref{#2},\ref{#3})}

\newcommand{\paren}[1]{\left({#1}\right)}
\newcommand{\brackets}[1]{\left[{#1}\right]}
\newcommand{\braces}[1]{\{#1\}}
\newcommand{\abs}[1]{\left\lvert#1\right\rvert}
\newcommand{\card}[1]{\left\lvert#1\right\rvert}

\newcommand{\regAbs}[1]{\lvert #1\rvert}

\newcommand{\extspac}[1]{\overline{\R^{#1}}}
\newcommand{\extspace}{\extspac{n}}
\newcommand{\extspacnp}{\extspac{n+1}}

\newcommand{\corez}[1]{\mathcal{E}_{#1}}
\newcommand{\corezn}{\corez{n}}
\newcommand{\coreznp}{\corez{n+1}}
\newcommand{\clcorezn}{\overlineKernIt{\mathcal{E}}_n}
\newcommand{\galax}[1]{{\cal G}_{{#1}}}
\newcommand{\galaxd}{\galax{\ebar}}
\newcommand{\galaxdp}{\galax{\ebar'}}
\newcommand{\galcl}[1]{\overlineKernIt{\cal G}_{{#1}}}
\newcommand{\galcld}{\galcl{\ebar}}
\newcommand{\galcldt}{\galcl{\ebar_t}}
\newcommand{\galcldp}{\galcl{\ebar'}}
\newcommand{\hogal}{\gamma}
\newcommand{\hogalinv}{\gamma'}%

\newcommand{\fcnspn}{{\cal F}}

\newcommand{\indfalet}{I}
\newcommand{\indflet}{\iota}

\newcommand{\indfa}[1]{\indfalet_{#1}}
\newcommand{\indfaAlt}[1]{\indfalet^{\vphantom{*}}_\subAlt{\mkern-1mu#1}}
\newcommand{\indfasing}[1]{\indfa{\mkern-0.5mu\set{#1}}}
\newcommand{\indfasingAlt}[1]{\indfaAlt{\set{#1}}}
\newcommand{\indfastar}[1]{\indfalet^*_{#1}}
\newcommand{\indfadub}[1]{\indfalet^{*\dualstar}_{#1}}

\newcommand{\indfdstar}[1]{\indflet^{\dualstar}_{#1}}
\newcommand{\indaS}{\indfa{S}}
\newcommand{\indaSstar}{\indfastar{S}}
\newcommand{\indaJstar}{\indfalet^*_{\negKern J}}
\newcommand{\indaJ}{\indfalet^{\vphantom{*}}_{\negKern J}}
\newcommand{\indaz}{\indfasing{0}}
\newcommand{\indazAlt}{\indfasingAlt{0}}

\newcommand{\indabetaAlt}{\indfasingAlt{\beta}}

\newcommand{\subAlt}[1]{{\smash[t]{#1}\vphantom{\le}}}
\newcommand{\lb}[2]{\seg(#1,#2)}
\newcommand{\indf}[1]{\indflet^{\vphantom{*}}_{\mkern-0.5mu#1}}
\newcommand{\indfAlt}[1]{\indflet^{\vphantom{*}}_\subAlt{\mkern-0.5mu#1}}
\newcommand{\indfstar}[1]{\indflet^{*}_{\mkern-0.5mu#1}}

\newcommand{\indfsing}[1]{\indfAlt{\set{#1}}}
\newcommand{\indfsingext}[1]{\indgenextAlt{\set{#1}}}
\newcommand{\indfsingstar}[1]{\indfstar{\set{#1}}}
\newcommand{\inds}{\indf{S}}
\newcommand{\indC}{\indf{C}}
\newcommand{\indK}{\indf{K}}
\newcommand{\indL}{\indf{L}}
\newcommand{\indKL}{\indf{K\cap L}}

\newcommand{\indz}{\indfsing{0}}

\newcommand{\indbeta}{\indfsing{\beta}}
\newcommand{\indJpol}{\indf{\negKern J^\circ}}
\newcommand{\indgenext}[1]{\ibar^{\vphantom{*}}_{\mkern-0.5mu#1}}
\newcommand{\indgenextAlt}[1]{\ibar^{\vphantom{*}}_\subAlt{\mkern-0.5mu#1}}

\newcommand{\indgenextstar}[1]{{\ibar}^{\mkern1mu*}_{#1}}
\newcommand{\indsext}{\indgenext{S}}
\newcommand{\indCext}{\indgenext{C}}
\newcommand{\indsextdub}{\ibar_S^{\mkern1mu*\dualstar}}

\newcommand{\indsextstar}{\indgenextstar{S}}
\newcommand{\indzext}{\indfsingext{0}}

\newcommand{\indbetaext}{\indfsingext{\beta}}
\newcommand{\indzoh}{\indz\negKern\circ h}
\newcommand{\indzohTight}{\indz\negKern\negKern\circ\negKern h}
\newcommand{\indzohj}{\indz\negKern\circ h_j}
\newcommand{\indzohext}{(\overline{\indzoh})}
\newcommand{\indzohextTight}{(\overline{\indzohTight})}
\newcommand{\indzohjext}{(\overline{\indzohj})}
\newcommand{\inddomf}{\indf{\dom{f}}}
\newcommand{\inddomfext}{\indgenext{\dom{f}}}

\newcommand{\indsdub}{\indflet^{*\dualstar}_S}
\newcommand{\indstar}[1]{{\indflet^*_{#1}}}
\newcommand{\indstars}{\indstar{S}}

\newcommand{\indzstar}{\indstar{0}}
\newcommand{\inddubstar}[1]{{\indflet^{**}_{#1}}}

\newcommand{\indepif}{\indf{\epi f}}
\newcommand{\indepifext}{\indgenext{\epi f}}
\newcommand{\indepifstar}{\indstar{\epi f}}
\newcommand{\inddomfstar}{\indstar{\dom f}}

\DeclareMathOperator{\conv}{\mathopfont{conv}}
\DeclareMathOperator{\seg}{\mathopfont{seg}}

\DeclareMathOperator{\ohull}{\widetilde{\mathopfont{conv}}}
\newcommand{\simplex}{\ohull}

\newcommand{\chs}[2]{{H}_{{#1},{#2}}}
\newcommand{\chsub}{\chs{\uu}{\beta}}
\newcommand{\chsua}{\chsub}
\newcommand{\chsuz}{\chs{\uu}{0}}

\newcommand{\lambar}{\hat{\lambda}}
\newcommand{\lamhat}{\hat{\lambda}}

\newcommand{\ahfsp}[3]{H_{{#1},{#2},{#3}}}

\DeclareMathOperator{\genhull}{\mathopfont{hull}}

\newcommand{\fmidpt}[3]{\mul{1-#1}#2\seqsum\mul{#1}#3}
\newcommand{\flammid}[2]{\fmidpt{\lambda}{#1}{#2}}
\newcommand{\fzrmid}[2]{\mul{1}#1\seqsum\mul{0}#2}

\newcommand{\mul}[2]{\pparens{#1}#2}
\newcommand{\almul}[1]{\mul{\alpha}{#1}}
\newcommand{\zmul}[1]{\mul{0}{#1}}

\newcommand{\lexless}{\mathrel{<_{L\negKern}}}

\newcommand{\sfprodskip}{\mskip0mu plus 0mu minus -2mu}
\newcommand{\sfprodsym}{\sfprodskip\mathbin{\mkern-1mu\varcirc\mkern-3mu}\sfprodskip}

\newcommand{\circkskip}{\mskip0mu plus 0mu minus -1mu}
\newcommand{\circk}{\circkskip\mathbin{\circ\mkern-2mu}\circkskip}

\newcommand{\sfprod}[2]{#1\sfprodsym#2}
\newcommand{\sfprodlamgi}{\sfprod{\lambda_i\negKern}{\posKern g_i}}

\newcommand{\sfprodstar}[2]{(\sfprod{#1}{#2})^*}

\newcommand{\rhoinv}{\rho^{-1}}
\newcommand{\vepsilon}{\boldsymbol{\epsilon}}

\newcommand{\lamvec}{\boldsymbol{\lambda}}
\newcommand{\muvec}{\boldsymbol{\mu}}
\newcommand{\lagrangfcn}{\ell}
\newcommand{\lagrangextsym}{\ellbar}
\newcommand{\lagrmid}{\mathbin{\mid}}
\newcommand{\lagrang}[2]{\lagrangfcn(#1 \lagrmid #2)}
\newcommand{\lagrangext}[2]{\lagrangextsym(#1 \lagrmid #2)}
\newcommand{\lagrangstar}[2]{\lagrangfcn^*(#1 \lagrmid #2)}
\newcommand{\lagranglm}[1]{\lagrang{#1}{\lamvec,\muvec}}
\newcommand{\lagranglmp}[1]{\lagrang{#1}{\lamvec',\muvec'}}
\newcommand{\lagranglmd}[1]{\lagrang{#1}{2\lamvec,2\muvec}}
\newcommand{\lagranglmz}[1]{\lagrang{#1}{\lamvec_0,\muvec_0}}
\newcommand{\lagranglmext}[1]{\lagrangext{#1}{\lamvec,\muvec}}
\newcommand{\lagrangflmext}{\lagranglmext{\posSkip\cdot}}
\newcommand{\lagranglmpext}[1]{\lagrangext{#1}{\lamvec',\muvec'}}
\newcommand{\lagranglmdext}[1]{\lagrangext{#1}{2\lamvec,2\muvec}}
\newcommand{\lagranglmzext}[1]{\lagrangext{#1}{\lamvec_0,\muvec_0}}
\newcommand{\lagranglmstar}[1]{\lagrangstar{#1}{\lamvec,\muvec}}

\newcommand{\lagrangflm}{\lagranglm{\posSkip\cdot}}

\newcommand{\lagpair}[2]{{\langle{#1}, {#2}\rangle}}
\newcommand{\lammupair}{\lagpair{\lamvec}{\muvec}}
\newcommand{\lammupairz}{\lagpair{\lamvec_0}{\muvec_0}}
\newcommand{\lammupairp}{\lagpair{\lamvec'}{\muvec'}}
\newcommand{\Lampairs}{\Gamma_{r,s}}
\newcommand{\sumitor}{\sum_{i=1}^r}
\newcommand{\sumjtos}{\sum_{j=1}^s}
\newcommand{\sumjtom}{\sum_{j=1}^m}
\newcommand{\negf}{\indf{\leq 0}}
\newcommand{\Negf}{\indfaAlt{\leq 0}}
\newcommand{\negfext}{\indgenext{\leq 0}}

\newcommand{\Negogiext}{(\overline{\Negf\circ g_i})}
\newcommand{\Negogext}{(\overline{\Negf\circ g})}

\newcommand{\lamgextgen}[1]{(\overline{\lambda_{#1} g_{#1}})}
\newcommand{\lamgiext}{\lamgextgen{i}}
\newcommand{\muhextgen}[1]{(\overline{\mu_{#1} h_{#1}})}
\newcommand{\muhext}{(\overline{\mu h})}
\newcommand{\muhjext}{\muhextgen{j}}

\newcommand{\Sperp}{{S^{\bot}}}
\newcommand{\Sperperp}{\perperp{S}}
\newcommand{\Sperpaperp}{\perpaperp{S}}
\newcommand{\Sbarperp}{{\Sbar}^{\posKern\bot}}
\newcommand{\Uperp}{{U^{\bot}}}
\newcommand{\Uperperp}{{U^{\bot\bot}}}
\newcommand{\Kperp}{{K^{\bot}}}
\newcommand{\Lperp}{{L^{\bot}}}
\newcommand{\Lperperp}{{L^{\bot\bot}}}
\newcommand{\Mperp}{{M^{\bot}}}

\DeclareMathOperator{\vertsl}{\mathopfont{vert}}
\DeclareMathOperator{\slopes}{\mathopfont{bar}}
\DeclareMathOperator{\barr}{\mathopfont{bar}}

\newcommand{\Rext}{\overline{\R}}

\newcommand{\ph}[1]{\phi_{#1}}
\newcommand{\phstar}[1]{{\phi^{\dualstar}_{#1}}}
\newcommand{\phstars}[1]{{\phi^{*}_{#1}}}
\newcommand{\phdub}[1]{{\phi^{\dualstar *}_{#1}}}
\newcommand{\phfcn}{{\varphi}}
\newcommand{\phfcninv}{\phfcn^{-1}}
\newcommand{\phx}{\ph{\xx}}
\newcommand{\phimg}{\Phi}%
\newcommand{\phimgcl}{\overline{\Phi}}%

\newcommand{\phimgA}{\Phi'}%

\newcommand{\calH}{{\cal H}}
\newcommand{\calQ}{{\cal Q}}

\newcommand{\Ep}{{E_{+}}}
\newcommand{\Em}{{E_{-}}}
\newcommand{\Hp}{H_{+}}
\newcommand{\Hm}{H_{-}}

\newcommand{\allseq}{{\cal S}}

\newcommand{\genexpfamtrans}[2]{{#1}_{#2}}

\newcommand{\distset}{I}
\newcommand{\distelt}{i}
\newcommand{\distalt}{j}
\newcommand{\distaltii}{k}%
\newcommand{\featmap}{\boldsymbol{\phi}}
\newcommand{\featmapsc}{\phi}
\newcommand{\featmapj}{\phi_j}
\newcommand{\featmapu}{\genexpfamtrans{\featmap}{\uu}}

\newcommand{\vparam}{\boldsymbol{\theta}}
\newcommand{\param}{\theta}
\newcommand{\barparam}{\overlineKernO{\param}}
\newcommand{\parambar}{\overlineKernO{\vparam}}
\newcommand{\qsub}[1]{p_{\mkern-1mu\vphantom{\barparam}#1}}
\newcommand{\qx}{\qsub{\vparam}}
\newcommand{\qsx}{\qsub{\param}}
\newcommand{\qbarx}{\qsub{\barparam}}
\newcommand{\qxbar}{\qsub{\parambar}}
\newcommand{\qxbarp}{\qsub{\parambar'\negKern}}
\newcommand{\qxt}{\qsub{\vparam_t}}
\newcommand{\qebar}{\qsub{\ebar}}

\newcommand{\medists}{Q}
\newcommand{\mldists}{P}

\newcommand{\sumex}{z}
\newcommand{\sumexu}{\genexpfamtrans{\sumex}{\uu}}

\newcommand{\sumexext}{\overlineKernIt{\sumex}}
\newcommand{\sumexextfi}{\genexpfamtrans{\sumexext}{\featmap(\distelt)}}
\newcommand{\sumexextu}{\genexpfamtrans{\sumexext}{\uu}}

\newcommand{\fullshadsumexu}{\fullshadnoparen{\sumexu}}

\newcommand{\lpart}{a}
\newcommand{\lpartu}{\lpart_{\uu}}
\newcommand{\lpartfi}{\lpart_{\featmap(\distelt)}}
\newcommand{\lpartext}{\overlineKernIt{\lpart}}
\newcommand{\lpartextu}{\lpartext_{\uu}}
\newcommand{\lpartextfi}{\lpartext_{\featmap(\distelt)}}
\newcommand{\lpartstar}{\lpart^*}
\newcommand{\lpartustar}{\lpartu^*}

\newcommand{\meanmap}{M}
\newcommand{\meanmapu}{\genexpfamtrans{\meanmap}{\uu}}

\newcommand{\entropy}{\mathrm{H}}

\newcommand{\convfeat}{\conv{\featmap(\distset)}}
\newcommand{\convfeatu}{\conv{\featmapu(\distset)}}

\newcommand{\subdiflpart}{\partial\lpart}

\newcommand{\subdiflpartstar}{\partial\lpartstar}
\newcommand{\asubdiflpart}{\asubdifplain{\lpartext}}
\newcommand{\asubdiflpartu}{\asubdifplain{\lpartextu}}
\newcommand{\basubdiflpart}{\basubdifplain{\lpartext}}

\newcommand{\adsubdiflpartstar}{\adsubdifplain{\lpartstar}}
\newcommand{\adsubdiflpartustar}{\adsubdifplain{\lpartustar}}

\newcommand{\Nj}[1]{m_{#1}}
\newcommand{\Nelt}{\Nj{\distelt}}
\newcommand{\numsamp}{m}
\newcommand{\loglik}{\ell}
\newcommand{\sfrac}[2]{\mbox{$\frac{#1}{#2}$}}

\newcommand{\phat}{\hat{p}}
\newcommand{\popdist}{\pi}

\newcommand{\regExp}[2]{\mathbf{E}_{#1}\regBracks{#2}}
\newcommand{\bigExp}[2]{\mathbf{E}_{#1}\negKern\bigBracks{#2}}
\DeclareMathOperator{\support}{\mathopfont{supp}}
\DeclareMathOperator{\affh}{\mathopfont{aff}}

\newcommand{\plusl}{\varplusl}
\newcommand{\plusd}{\varplusd}
\newcommand{\plusu}{\varplusu}

\newcommand{\seqeq}[2]{#1\sim #2}
\newcommand{\seqeqsymbol}{\sim}

\newcommand{\seqgt}{\gg}
\newcommand{\seqlt}{\ll}

\newcommand{\ci}[1]{i_{#1}}

\newcommand{\colset}{I}
\newcommand{\coldefset}{J}

\newcommand{\seqsum}{\mathbin{\varhash}}
\newcommand{\Seqsum}{\Sigma_{\seqsum}}

\newcommand{\glambdaletter}{\lambda}
\newcommand{\glambdait}{\glambdaletter_{it}}
\newcommand{\glambdapit}{\glambdaletter'_{it}}

\newcommand{\singerphi}{\varphi}

\newcommand{\scriptontop}[2]{\substack{#1 \\ #2}}

\newcommand{\InfseqLiminf}[4]{%
     \inf_{\substack{\vphantom{t\to\infty}#1 \textit{ in } #2:\\ #3}}
     \;
     \liminf_{\substack{t\to\infty\vphantom{#1 \textit{ in } #2:}\\ \vphantom{#3}}} {#4}}

\numberwithin{equation}{chapter}
\numberwithin{figure}{chapter}
\numberwithin{algorithm}{chapter}

\newtheoremstyle{mysubproof}
  {8pt plus 2pt minus 4pt} %
  {8pt plus 2pt minus 4pt} %
  {\normalfont}  %
  {0pt}       %
  {\itshape}  %
  {.}         %
  {.5em}      %
  {}          %

\newtheoremstyle{mysubclaim}
  {8pt plus 2pt minus 4pt} %
  {8pt plus 2pt minus 4pt} %
  {\normalfont}  %
  {0pt}       %
  {\itshape}  %
  {.}         %
  {.5em}      %
  {}          %

\theoremstyle{plain}
\newtheorem{theorem}{Theorem}[chapter]
\newtheorem{corollary}[theorem]{Corollary}
\newtheorem{proposition}[theorem]{Proposition}
\newtheorem{lemma}[theorem]{Lemma}

\theoremstyle{remark}
\newtheorem*{mainclaimp*}{Main claim}

\declaretheorem[title=Claim,within=theorem,style=mysubclaim]{claimpx}

\declaretheorem[title={\hskip0pt Proof},numbered=no,style=mysubproof,qed=$\lozenge$]{proofx}

\theoremstyle{definition}
\declaretheorem[title=Example,numberlike=theorem,style=definition,qed=$\blacksquare$]{example}
\newtheorem{definition}[theorem]{Definition}

\newcounter{tmpenumcounter}

\DeclareMathOperator*{\argmax*}{arg\,max}

\newcommand{\hpp}{\hat{\pp}}

\newcommand{\lsePP}{\pp}
\newcommand{\lseP}{p}
\newcommand{\lsePPh}{\hpp}
\newcommand{\lsePh}{\hat{p}}
\newcommand{\lsePPv}{\pp}
\newcommand{\lsePv}{p}

\newcommand{\probsimgen}[1]{\Delta_{#1}}
\newcommand{\probsim}{\probsimgen{n}}
\newcommand{\probsimm}{\probsimgen{m}}

\newcommand{\gamep}{p}
\newcommand{\ulam}{\uu_{\gamep}}
\newcommand{\ulamhalf}{\uu_{1/2}}

\newcommand{\eRn}{\overline{\R^n}}%
\newcommand{\eRm}{\overline{\R^m}}%
\newcommand{\eR}{\overline{\R}}
\newcommand{\eRf}[1]{\overline{\R^{#1}}}
\newcommand{\e}[1]{{}\mkern3mu\overline{\mkern-3mu#1}}
\newcommand{\ef}{\overlineKernf{f}}
\newcommand{\efp}{{}\mkern3mu\overline{\mkern-3mu f'\mkern1mu}\mkern-1mu}
\newcommand{\efpsub}[1]{{}\mkern3mu\overline{\mkern-3mu f'_{#1}\mkern1mu}\mkern-1mu}
\newcommand{\efspace}{{}\mkern4.5mu\overline{\mkern-4.5mu f\mkern1.5mu}}
\newcommand{\efsub}[1]{{}\mkern4.5mu\overline{\mkern-4.5mu f\mkern1.5mu}\mkern-1.5mu\mkern-2mu{}_{#1}}
\newcommand{\efsubstar}[1]{{}\mkern4.5mu\overline{\mkern-4.5mu f\mkern1.5mu}^{\mathrlap{*}}\mkern-1.5mu\mkern-3mu{\vphantom{f}}_{#1}}

\newcommand{\ep}{\overlineKernp{p}}
\newcommand{\eg}{\overlineKerng{g}}
\newcommand{\egspace}{\overlineKerngspace{g}}
\newcommand{\eh}{{}\mkern1mu\overline{\mkern-1mu h}}
\newcommand{\ibar}{\ibariialt}
\newcommand{\ministrut}{\vphantom{*}}
\newcommand{\ministrutalt}{\vphantom{\cdot}}

\newcommand{\ibariialt}{\mkern1mu\overline{\mkern-1mu\iota\ministrutalt\mkern1mu}\mkern-1mu}

\newcommand{\ex}{\overlineKernIt{x}}
\newcommand{\ey}{\e{y}}

\newcommand{\exx}{\overline{\xx}}
\newcommand{\eyy}{\overline{\yy}}

\newcommand{\veps}{\boldsymbol{\varepsilon}}

\newcommand{\unminus}{\mskip0mu plus 0mu minus -0.5mu}
\newcommand{\inprod}{\cdot}
\newcommand{\inprodk}{\unminus\cdot\unminus}

\newcommand{\wo}{\backslash}

\makeatletter
\newcommand{\varlozenge}{\mathpalette\v@rlozenge\relax}
\newcommand{\varcirc}{\mathpalette\v@rcirc\relax}
\newcommand{\varcircup}{\mathpalette\v@rcircup\relax}
\newcommand{\varsquare}{\mathpalette\v@rsquare\relax}
\newcommand\v@rlozenge[2]{%
  \setbox0=\hbox{$\m@th#1+$}%
  \setbox1=\hbox{$\m@th#1\vcenter{\hbox{\resizebox{!}{0.45\wd0}{$\m@th#1\blacklozenge\vphantom{+}$}}}$}%
  \mathbin{\ooalign{%
    \box1}}}
\newcommand\v@rcirc[2]{%
  \setbox0=\hbox{$\m@th#1+$}%
  \setbox1=\hbox{$\m@th#1\vcenter{\hbox{\resizebox{!}{0.5\wd0}{$\m@th#1\bullet\vphantom{+}$}}}$}%
  \mathbin{\ooalign{%
    \box1}}}
\newcommand\v@rcircup[2]{%
  \setbox0=\hbox{$\m@th#1+$}%
  \setbox1=\hbox{$\m@th#1\vcenter{\hbox{\resizebox{!}{0.35\wd0}{$\m@th#1\altcircup\vphantom{+}$}}}$}%
  \mathbin{\ooalign{%
    \box1}}}
\newcommand\v@rsquare[2]{%
  \setbox0=\hbox{$\m@th#1+$}%
  \setbox1=\hbox{$\m@th#1\vcenter{\hbox{\resizebox{!}{0.35\wd0}{$\m@th#1\blacksquare\vphantom{+}$}}}$}%
  \mathbin{\ooalign{%
    \box1}}}
\newcommand{\varplusl}{\mathpalette\v@rplusl\relax}
\newcommand{\varplusr}{\mathpalette\v@rplusr\relax}
\newcommand{\varpluslr}{\mathpalette\v@rpluslr\relax}
\newcommand{\varplusd}{\mathpalette\v@rplusd\relax}
\newcommand{\varplusu}{\mathpalette\v@rplusu\relax}
\newcommand{\varhash}{\mathpalette\v@rhash\relax}
\newlength{\v@rhashtot}
\newlength{\v@rhashraise}
\newcommand\v@rhashvar[2]{%
  \setbox0=\hbox{$\m@th#1=$}%
  \setbox1=\hbox{$\m@th#1\parallel$}%
  \setbox2=\hbox{$\m@th#1+$}%
  \setlength{\v@rhashtot}{\ht2}%
    \addtolength{\v@rhashtot}{\dp2}%
  \setbox3=\hbox{$\m@th#1\vbox{\hbox{\resizebox{\wd0}{\ht0}{$\m@th#1=$}}}$}%
  \setbox4=\hbox{$\m@th#1\vcenter{\hbox{\resizebox*{\wd1}{\v@rhashtot}{$\m@th#1\parallel$}}}$}%
  \mathbin{\ooalign{%
    \hidewidth\box4\hidewidth\cr
    \box3\cr}}}
\newcommand\v@rhashvarred[2]{%
  \setbox0=\hbox{$\m@th#1=$}%
  \setbox1=\hbox{$\m@th#1\parallel$}%
  \setbox2=\hbox{$\m@th#1+$}%
  \setlength{\v@rhashtot}{\ht2}%
    \addtolength{\v@rhashtot}{\dp2}%
  \setlength{\v@rhashraise}{0.1875\ht0}%
  \setbox3=\hbox{$\m@th#1\vbox{\hbox{\resizebox{\wd0}{0.75\ht0}{$\m@th#1=$}}}$}%
  \setbox4=\hbox{$\m@th#1\vcenter{\hbox{\resizebox*{0.75\wd1}{\v@rhashtot}{$\m@th#1\parallel$}}}$}%
  \mathbin{\ooalign{%
    \hidewidth\box4\hidewidth\cr
    \raise\v@rhashraise\box3\cr}}}
\newcommand\v@rhash[2]{%
  \setbox0=\hbox{$\m@th#1-$}%
  \setbox1=\hbox{$\m@th#1\mid$}%
  \setbox2=\hbox{$\m@th#1+$}%
  \setlength{\v@rhashtot}{\ht2}%
    \addtolength{\v@rhashtot}{\dp2}%
  \setbox3=\copy0
  \setbox4=\copy3
  \setbox5=\hbox{$\m@th#1\vcenter{\hbox{\resizebox*{\wd1}{\v@rhashtot}{$\m@th#1\mid$}}}$}%
  \setbox6=\copy5
  \mathbin{\ooalign{%
    \hidewidth\box5\kern.23\v@rhashtot\hidewidth\cr
    \hidewidth\kern.23\v@rhashtot\box6\hidewidth\cr
    \raise.115\v@rhashtot\box3\cr
    \lower.115\v@rhashtot\box4}}}
\newlength{\v@rplusraise}
\newcommand\v@rplusl[2]{%
  \setbox0=\hbox{$\m@th#1+$}%
  \setbox1=\hbox{$\m@th#1\vcenter{\hbox{\resizebox{!}{0.3\wd0}{$\m@th#1\bullet\vphantom{+}$}}}$}%
  \mathbin{\ooalign{%
    \box1\hidewidth\cr
    \box0\cr}}}
\newcommand\v@rplusr[2]{%
  \setbox0=\hbox{$\m@th#1+$}%
  \setbox1=\hbox{$\m@th#1\vcenter{\hbox{\resizebox{!}{0.3\wd0}{$\m@th#1\bullet\vphantom{+}$}}}$}%
  \mathbin{\ooalign{%
    \hidewidth\box1\cr
    \box0\cr}}}
\newcommand\v@rpluslr[2]{%
  \setbox0=\hbox{$\m@th#1+$}%
  \setbox1=\hbox{$\m@th#1\vcenter{\hbox{\resizebox{!}{0.3\wd0}{$\m@th#1\bullet\vphantom{+}$}}}$}%
  \setbox2=\copy1
  \mathbin{\ooalign{%
    \hidewidth\box1\cr
    \box2\hidewidth\cr
    \box0\cr}}}
\newcommand\v@rplusd[2]{%
  \setbox0=\hbox{$\m@th#1+$}%
  \setbox1=\hbox{\resizebox{!}{0.3\wd0}{$\m@th#1\bullet\vphantom{+}$}}%
  \setlength{\v@rplusraise}{-\dp0}%
    \addtolength{\v@rplusraise}{-0.25\ht1}%
  \mathbin{\ooalign{%
    \hidewidth\raise\v@rplusraise\box1\hidewidth\cr
    \box0\cr}}}
\newcommand\v@rplusu[2]{%
  \setbox0=\hbox{$\m@th#1+$}%
  \setbox1=\hbox{\resizebox{!}{0.3\wd0}{$\m@th#1\bullet\vphantom{+}$}}%
  \setlength{\v@rplusraise}{\ht0}%
    \addtolength{\v@rplusraise}{\dp1}%
    \addtolength{\v@rplusraise}{-0.75\ht1}%
  \mathbin{\ooalign{%
    \hidewidth\raise\v@rplusraise\box1\hidewidth\cr
    \box0\cr}}}
\makeatother

\newcommand{\bracks}[1]{[#1]}
\newcommand{\regBracks}[1]{[#1]}
\newcommand{\Bracks}[1]{\left[#1\right]}
\newcommand{\bigBracks}[1]{\bigl[#1\bigr]}
\newcommand{\BigBracks}[1]{\Bigl[#1\Bigr]}
\newcommand{\BigBraces}[1]{\Bigl\{#1\Bigr\}}
\newcommand{\biggBraces}[1]{\biggl\{#1\biggr\}}
\newcommand{\BiggBraces}[1]{\Biggl\{#1\Biggr\}}
\newcommand{\biggBracks}[1]{\biggl[#1\biggr]}
\newcommand{\BiggBracks}[1]{\Biggl[#1\Biggr]}
\newcommand{\bigBraces}[1]{\bigl\{#1\bigr\}}
\newcommand{\regBraces}[1]{\{#1\}}
\newcommand{\Braces}[1]{\left\{#1\right\}}
\newcommand{\parens}{\Parens}
\newcommand{\Parens}[1]{\left(#1\right)}
\newcommand{\regParens}[1]{(#1)}
\newcommand{\bigParens}[1]{\bigl(#1\bigr)}
\newcommand{\BigParens}[1]{\Bigl(#1\Bigr)}
\newcommand{\biggParens}[1]{\biggl(#1\biggr)}
\newcommand{\BiggParens}[1]{\Biggl(#1\Biggr)}
\newcommand{\set}[1]{\{#1\}}
\newcommand{\bigSet}[1]{\bigl\{#1\bigr\}}

\newcommand{\bigNorm}[1]{\bigl\lVert#1\bigr\rVert}
\newcommand{\BiggNorm}[1]{\Biggl\lVert#1\Biggr\rVert}

\newcommand{\seq}[1]{{(#1)}}

\newcommand{\bigAbs}[1]{\bigl\lvert#1\bigr\rvert}

\newcommand{\tupset}[2]{\langle#1\rangle_{#2}}
\newcommand{\tprojsym}{\pi}
\newcommand{\tproja}{\tprojsym_{\alpha}}
\newcommand{\tprojb}{\tprojsym_{\beta}}
\newcommand{\tprojinva}{\tprojsym^{-1}_{\alpha}}
\newcommand{\tprojz}{\tprojsym_{z}}
\newcommand{\tprojinvz}{\tprojsym^{-1}_{z}}

\newcommand{\RefsImplication}[2]{(\ref{#1})~$\Rightarrow$~(\ref{#2})}

\newcommand{\fA}{f\kernA}
\newcommand{\FA}{F\kernA}
\newcommand{\kernA}{\mkern-1mu\A}
\newcommand{\transA}{\trans{\A\mkern-3mu}}
\newcommand{\transAk}{\transA\mkern-1mu}
\newcommand{\ktransA}{\mkern-1mu\transA}

\newcommand{\utrans}{\uu^{\mkern-3mu\top}}
\newcommand{\utransA}{\trans{\uu}\mkern-3mu\A}
\newcommand{\xtransA}{\trans{\xx}\mkern-3mu\A}
\newcommand{\utransF}{\trans{\uu}\mkern-1mu F}
\newcommand{\utransf}{\trans{\uu}\mkern-1mu f}
\newcommand{\utransfext}{\trans{\uu}\mkern-1mu\ef}
\newcommand{\transB}{\trans{\B\mkern-1.5mu}}
\newcommand{\transBk}{\transB\mkern-1mu}

\newcommand{\pcite}[1]{p.~#1}
\newcommand{\ppcite}[1]{pp.~#1}

\newlength{\fencheldiffpairlen}
\newcommand{\fencheldiffpair}[2]{\makebox[\fencheldiffpairlen][l]{$#1$}\;\;and\;\;{$#2$}.}

\newlength{\tmpwidth}
\newcommand{\forcewidthof}[3][l]{%
  \settowidth{\tmpwidth}{#2}%
  \makebox[\tmpwidth][#1]{#3}%
}

\newcommand{\pmUni}{^^c2^^b1}%
\newcommand{\inftyUni}{^^e2^^88^^9e}%

\newlist{letter-compact}{enumerate}{1}
\setlist[letter-compact]{noitemsep, label={\upshape(\alph*)}, ref=\alph*}
\newlist{letter}{enumerate}{1}
\setlist[letter]{label={\upshape(\alph*)}, ref=\alph*}

\newlist{letter-compact-prime-simple}{enumerate}{1}
\setlist[letter-compact-prime-simple]{noitemsep, label={\upshape(\alph*$'$)}, ref=\alph*$'$}

\newlist{letter-compact-prime}{enumerate}{1}
\setlist[letter-compact-prime]{noitemsep, label={\upshape(\alph*)\hphantom{$'$}}, ref=\alph*}
\makeatletter
\newcommand{\itemprime}{%
\item[\upshape(\alph{letter-compact-primei}$'$)]\protected@edef\@currentlabel{\alph{letter-compact-primei}$'$}\ignorespaces%
}

\makeatletter
\newcommand{\itemx}[1]{%
\item[\upshape(#1)]\protected@edef\@currentlabel{#1}\ignorespaces%
}
\makeatother

\newcounter{tmplccounter}
\newcommand{\savelccounter}{
  \setcounter{tmplccounter}{\arabic{letter-compacti}}
}
\newcommand{\restorelccounter}{
  \setcounter{letter-compacti}{\thetmplccounter}
}

\newlist{roman-compact}{enumerate}{1}
\setlist[roman-compact]{noitemsep, label={\upshape(\roman*)}, ref=\roman*}

\newlist{item-compact}{itemize}{1}
\setlist[item-compact]{noitemsep, label=$\bullet$}

\newlist{proof-parts}{description}{2}
\setlist[proof-parts]{nosep, wide, font=\rmfamily\mdseries\itshape}
\makeatletter
\newcommand{\pfpart}[1]{\item[\trim@spaces{#1}]\hskip 0.1em minus 0.3em} %
\makeatother

\newcommand{\overviewsecheading}[1]{\paragraph{\scshape\mdseries\Cref{#1}:~~\nameref{#1}.}}

\newcommand{\figclifframp}[1]{%
  \begin{tikzpicture}
      \node at (0,0) {\includegraphics{\detokenize{#1}}};
      \node[xslant=-0.8, yscale=0.5, rotate=8, color=black] at (-1.8,1.85) {$+\infty$};
      \node[xslant=-0.35, yscale=0.5, rotate=11, color=black] at (-3.25,-1.2) {$-\infty$};
  \end{tikzpicture}%
}

\newcommand{\figentsubgrad}[1]{%
  \begin{tikzpicture}
      \node at (0,0) {\includegraphics{\detokenize{#1}}};
      \node[xslant=-0.8, yscale=0.5, rotate=13, color=black] at (-2.3,1.75) {$+\infty$};
      \node[xslant=-0.5, yscale=0.5, rotate=17, color=black] at (-4,-1) {$-\infty$};
  \end{tikzpicture}%
}

\newcommand{\mycaption}[2]{%
  \caption[#1]{\emph{#1.} #2}%
}

\newcommand{\nocomma}[1]{}

\newcommand{\indexintrotext}{\emph{References to figures and algorithms are marked, respectively, with \textup{f} and \textup{a}.}}

\newcommand{\indexa}{\indexgen{autind}}                %
\newcommand{\indexg}{\indexgen{default}}                %
\newcommand{\indexf}[1]{\indexgen{default}{#1|figidx}}  %
\newcommand{\indexalg}[1]{\indexgen{default}{#1|algidx}}  %
\newcommand{\indexm}[3]{%
  \indexgen{notind}{#1@\idxnotationentry{#2}#3}%
}

\newcommand{\indexgen}[1]{\index[#1]}

\newlength{\mathidxindent}
\newlength{\mathidxminspace}
\newlength{\tmpmathidxlen}

\newcommand{\idxnotationentry}[1]{%
  \settowidth{\tmpmathidxlen}{#1}%
  \addtolength{\tmpmathidxlen}{\mathidxminspace}%
  \ifdim \tmpmathidxlen<\mathidxindent%
      \setlength{\tmpmathidxlen}{\mathidxindent}%
  \fi%
  \makebox[\tmpmathidxlen][l]{#1}%
}

\newcommand{\idxroc}{\indexa{Rockafellar, R. T.}}
\newcommand{\idxwets}{\indexa{Wets, R. J-B}}
\newcommand{\idxroman}{\indexa{Roman, S.}}
\newcommand{\idxhj}{\indexa{Horn, R. A.}\indexa{Johnson, C. R.}}
\newcommand{\idxytt}{\indexa{Yanai, H.}\indexa{Takeuchi, K.}\indexa{Takane, Y.}}
\newcommand{\idxwainjord}{\indexa{Wainwright, M. J.}\indexa{Jordan, M. I.}}
\newcommand{\idxhitch}{\indexa{Aliprantis, C. D.}\indexa{Border, K. C.}}
\newcommand{\idxmunk}{\indexa{Munkres, J. R.}}
\newcommand{\idxmanetti}{\indexa{Manetti, M.}}
\newcommand{\idxborlew}{\indexa{Borwein, J. M.}\indexa{Lewis, A. S.}}
\newcommand{\idxhiriart}{\indexa{Hiriart-Urruty, J.-B.}\indexa{lemarechal, c.@Lemar\'echal, C.}}
\newcommand{\idxwaggoner}{\indexa{Waggoner, B.}}
\newcommand{\idxhandup}{\indexa{Hansen, G. L.}\indexa{Dupin, J.-C.}}
\newcommand{\idxsinger}{\indexa{Singer, I.}}
\newcommand{\idxschapire}{\indexa{Schapire, R. E.}}
\newcommand{\idxcollinsetal}{\indexa{Collins, M.}\idxschapire\indexa{Singer, Y.}}
\newcommand{\idxcoxetal}{\indexa{Cox, D. A.}\indexa{Little, J.}\indexa{oshea, d.@O'Shea, D.}}
\newcommand{\idxvandevel}{\indexa{Vel, M. L. J. van de}}
\newcommand{\idxmartinezsinger}{\idxsinger\indexa{martinez-legaz, j.-e.@Mart{\'i}nez-Legaz, J.-E.}}
\newcommand{\idxboydvand}{\indexa{Boyd, S.}\indexa{Vandenberghe, L.}}
\newcommand{\idxmatus}{\indexa{Telgarsky, M.}}
\newcommand{\idxneumann}{\indexa{von Neumann, J.}}
\newcommand{\idxziwei}{\indexa{Ji, Z.}}
\newcommand{\idxmiro}{\indexa{dudik, m.@Dud{\'i}k, M.}}
\newcommand{\idxfreund}{\indexa{Freund, Y.}}
\newcommand{\idxnesterov}{\indexa{Nesterov, Y.}}
\newcommand{\idxnemirovski}{\indexa{Nemirovski, A. S.}}
\newcommand{\idxtongbin}{\indexa{Yu, B.}\indexa{Zhang, T.}}

\begin{document}

\frontmatter
\pdfbookmark[1]{Astral Space: Convex Analysis at Infinity}{title}
\renewcommand*{\thefootnote}{\fnsymbol{footnote}}

\title{%
  \textbf{\Huge%
    Astral Space}
  \\[\bigskipamount]
  \textbf{\huge%
    Convex Analysis at Infinity}%
   \footnote{Final pre-publication draft.
             Book to be published by Princeton University Press in 2026.
   }
  \\[3\bigskipamount]
}

\author{%
  \Large%
  Miroslav Dud\'ik%
  \footnote{%
      Microsoft Research,
      New York City.
      \textit{\href{mailto:mdudik@microsoft.com}{mdudik@microsoft.com}},
      \textit{\href{mailto:schapire@microsoft.com}{schapire@microsoft.com}}.
  }%
  \and
  \Large%
  Robert E. Schapire%
  \footnotemark[2]%
  \and
  \Large%
  Matus Telgarsky%
  \footnote{%
      Courant Institute,
      New York University.
      \textit{\href{mailto:mjt10041@nyu.edu}{mjt10041@nyu.edu}}.
  }%
  \\[\bigskipamount]
}

\date{%
  \Large%
  December 2025
}
\maketitle

\clearpage
\thispagestyle{empty}

\vspace*{2in}
\begin{minipage}{4.05in}
\hfill\textbf{Abstract}\hfill
\medskip

Not all convex functions on $\R^n$ have finite minimizers;
some can only be minimized by a sequence as it heads to infinity.
In this work, we aim to develop a theory for understanding such
minimizers at infinity.
We study \emph{astral space}, a compact extension of $\R^n$
to which such points at infinity have been added.
Astral space is constructed to be as small as possible while still
ensuring that all linear functions can be continuously extended to the
new space.
Although astral space includes all of $\R^n$, it is not a vector space, nor even a metric space. However, it is sufficiently well-structured to allow useful and meaningful extensions of such concepts as convexity, conjugacy, and subdifferentials. We develop these concepts and analyze various properties of convex functions on astral space, including the detailed structure of their minimizers, exact characterizations of continuity, and convergence of descent algorithms.
\end{minipage}

\clearpage
\pdfbookmark[1]{\contentsname}{toc}
\tableofcontents

\listoffiguresalgos

\chapter{Preface}

The ideas for this book emerged over the course of our careers.
As machine-learning researchers, we have encountered optimization
problems that seemed very well-behaved and quite amenable to algorithmic approaches,
except for one problem: to fully solve them, one had to allow
the optimization parameter to go all the way ``to infinity,''
making proofs of convergence considerably more challenging.

For Rob and Matus,
the problem of optimizers at infinity came about in the study of AdaBoost and related algorithms. For Miro,
it came about in the study of maximum entropy, maximum likelihood, and
scoring rules in prediction markets. We all wrestled separately, and together,
as research collaborators, on analyzing such optimizers by means of sequences.
However, we all had a lurking suspicion that there should be a better way.
Perhaps, instead of reasoning about asymptotics of sequences, we could somehow add
new infinite points to
$n$-dimensional Euclidean space,
which could be treated
as mathematical objects in their own right and incorporated into the standard
tools of calculus and optimization.

Those were the seeds of what eventually grew into astral space. We started having
conversations and laying foundations for the basic theory of astral space in 2013, at a newly formed
Microsoft Research Lab in New York City, where Miro worked, Rob was on sabbatical,
and Matus was visiting.
Astral space was a curious new world that we explored
over several years, developing a theory that bit-by-bit came to
encompass extensions of many central topics of convex analysis, including 
convexity itself, conjugacy, separation theorems, and subgradients.

We originally planned to publish our theory in one or more research
articles, but the project grew to be too large,
so we thought it best to present a complete story as a book.
By the standards of our field,
this book is quite unusual since the ideas, notions, and results
here have never before been published, even partially. Yet, we believe that
the foundations belong together with their implications and
applications,
including those that first motivated us, which we return to toward the
end of the book.
We hope that you, dear reader, will likewise find value in this theory
for your own endeavors.

\paragraph{Acknowledgments.}
We are especially grateful to the anonymous reviewers of this book for
their time and invaluable feedback.
Thanks also to
Polina Baron,
Jonathan Eckstein,
Yoav Freund,
Raf Frongillo,
Surbhi Goel,
Zaid Harchaoui,
Terry Rockafellar,
and
Bo Waggoner
for their help, thoughts, comments, and suggestions.
Special thanks to Ziwei Ji for numerous ideas and discussions.

We are so very grateful to our editor, Hallie Stebbins,
for her patience, energy,
flexibility,
and for her vigorous and caring support of this project.
Thanks to everyone at Princeton University Press for
making this book a physical reality.

Miro and Rob thank their employer, Microsoft, for supporting
this work through many years of effort.
Matus was supported partially by NSF grant IIS-1750051,
and
part of this research was conducted while
visiting Microsoft Research, part while at
University of Illinois, Urbana-Champaign, Department of Computer
Science,
and now while at Courant Institute, New York University.

Finally, Miro thanks his husband, Jeremy, and children, Lili and Cate,
and Rob thanks his wife, Roberta, and children, Jeni and Zack,
for all their patience, love, and support.

\bigskip

\hfill Miroslav Dud\'ik, Robert Schapire, Matus Telgarsky

\hfill December 1, 2025

\mainmatter
\DeclareRobustCommand{\VAN}[3]{#2}

\chapter{Introduction}
\label{sec:intro}

Convex functions
play a fundamental role in
economics, operations
research, statistics, machine learning,
information theory, game theory,
and many other disciplines of engineering and science
\idxboydvand\citep{boyd_vandenberghe}.
For instance, many objective functions that are minimized in operations research and machine learning, such as those that come up in portfolio optimization and in fitting predictive models to data, are convex.
Convexity is an attractive structural property, because minimization of convex functions is algorithmically very well-behaved, with a vast array of efficient algorithms available
(%
\idxboydvand\citealp{boyd_vandenberghe};
\idxnemirovski\indexa{Yudin, D. B.}\citealp{nemirovsky_yudin};
\idxnesterov\citealp{nesterov-lectures-cvx-opt}%
).

\begin{figure}
  \centering
  \includegraphics{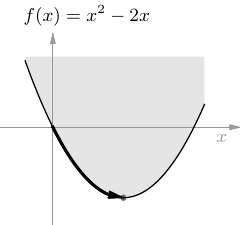}\hspace{0.15\textwidth}%
  \includegraphics{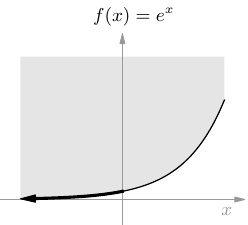}
  \mycaption{Minimization of convex functions in one dimension}{%
    \emph{Left:} An example where the approach of going down the slope (in the direction of negative gradient) converges to a minimizer.
    \emph{Right:} An example where the approach diverges to (negative) infinity.}
  \label{fig:1d-desc}%
\end{figure}

For instance,
for many convex functions, it is possible simply to ``go down the
slope'' of the function (that is, in the direction of the negative
gradient), which, for appropriately chosen step sizes,
is guaranteed to converge to a minimizer
(see \Cref{fig:1d-desc}, left, for an illustration in one dimension).

However, there are important examples of well-behaved functions (say, lower-bounded and differentiable on all of $\Rn$),
for which no algorithm can find a minimizer,
because no minimizer exists.
This can happen even in the simplest of cases,
like the exponential function $f(x)=e^x$ shown in \Cref{fig:1d-desc}~(right). In this case, going down the slope decreases function values,
but the algorithm diverges, with its iterates going to (negative) infinity. 
This function has
no finite minimizer,
and can only be minimized in a limit as
$x\to-\infty$.

Similarly, in multiple dimensions, there are examples of convex functions arising
in machine learning, statistics, and economics, which can only be
minimized by following a sequence to infinity. In these applications, we are often
interested in various properties of such minimizing sequences, but the analysis
of these sequences can be somewhat laborious and specialized to each case
(see, for instance,
\idxcollinsetal%
\citealp{collins_schapire_singer_adaboost_bregman};
\idxtongbin\citealp{zhang_yu_boosting};
and Section 3.6 of
\idxwainjord\citealp{WainwrightJo08}).

This book seeks to move past such specialized analyses.
Our overarching goal is to develop a more complete theory of convex
functions by extending $\Rn$ with suitable ``points at infinity.''
Ideally, with such an extension, there would be no need to analyze various types of divergent sequences; instead, we could work directly in the extended space, which would contain new ``infinite points''
to which these sequences converge.

In later chapters, we construct
a specific extension of $\Rn$ called \emph{astral space}, study its
properties, and develop an expansive theory of
convex functions on this space.
Before embarking on that journey, we aim in this introductory chapter
to bring out some key issues that arise in constructing such an
extension,
while also giving an overview of our approach
and introducing some of the themes (like differential theory and conjugacy) that will be examined in detail later in the book.

\section{Minimizers at infinity}
\label{sec:intro:sequences}

The task of extending $\Rn$ with points at infinity is deceptively simple in one
dimension, where we only need to add $+\infty$ and $-\infty$. The
resulting extension is called the
\indexg{extended reals}%
\indexm{r 600}{$\Rext$}{extended reals}%
\emph{extended real line} and denoted
$\eR=[-\infty,+\infty]$.
It is far from obvious how to generalize this concept to multiple dimensions. To develop some intuition about what is required of such an extension, we next look at several examples of convex functions and their minimizing sequences.

Before proceeding, some brief notational comments:
Vectors $\xx$ in $\Rn$ are usually written in bold, with components
$x_i$ given in italics so that
$\xx=\trans{[x_1,\dotsc,x_n]}$.
We use $t=1,2,\dotsc$ as the
index of sequences, unless noted otherwise,
and write $\seq{\xx_t}$ to mean
the sequence $\xx_1,\xx_2,\dotsc$. %
Limits and convergence are taken as $t\rightarrow+\infty$,
unless stated otherwise. For example, $\lim \xx_t$ means
$\lim_{t\to\infty}\xx_t$, and $\xx_t\rightarrow\xx$ means $\xx_t$ converges to~$\xx$ as $t\rightarrow+\infty$.
(More will be said about notation in
\Cref{sec:prelim-notation}.)

\begin{example}[Log-sum-exp]
\label{ex:log-sum-exp}
\indexg{minimizing sequences!examples|(}%
\indexg{Log-sum-exp|(}%
We first consider the log-sum-exp function, which comes up, for example, when fitting a multinomial model by maximum likelihood. Specifically,
let
\begin{equation}  \label{eqn:probsim-defn}
  \probsim = \BiggBraces{ \lsePPv\in [0,1]^n :\: \sum_{i=1}^n \lsePv_i = 1 }
\end{equation}
be the probability simplex in $\Rn$,
and
let $\lsePPh$ be any point in $\probsim$.
Consider minimization of the function $f:\Rn\rightarrow\R$ defined,
for $\xx\in\Rn$, by
\[
  f(\xx)=\ln\BiggParens{\sum_{i=1}^n e^{x_i}}-\lsePPh\inprod\xx.
\]
As is standard, $f$ can be rewritten as
\begin{equation}  \label{eqn:log-sum-exp-fcn}
  f(\xx)=-\sum_{i=1}^n \lsePh_i\ln\biggParens{\frac{e^{x_i}}{\sum_{j=1}^n e^{x_j}}}
        =-\sum_{i=1}^n \lsePh_i\ln \lseP_i(\xx),
\end{equation}
where $\lseP_i(\xx)=e^{x_i}/\sum_{j=1}^n e^{x_j}$,
so $\lsePP(\xx)=\trans{[\lseP_1(\xx),\ldots,\lseP_n(\xx)]}$
is a probability vector in $\probsim$. For this reformulation of $f$, it can be proved that
\begin{equation}  \label{eqn:log-sum-exp-inf}
  \inf_{\xx\in\Rn} f(\xx) = -\sum_{i=1}^n \lsePh_i\ln\lsePh_i,
\end{equation}
and that $\xx$ minimizes $f$ if and only if $\lsePP(\xx)=\lsePPh$
(see \Cref{pr:exp-log-props}\ref{pr:exp-log-props:a}).
So, if $\lsePh_i>0$ for all $i$, the minimum is attained by the vector $\xx$ with entries $x_i=\ln\lsePh_i$.

However, if $\lsePh_i=0$ for some $i$, then there is no $\xx\in\Rn$ that  attains the minimum.
In this case,
the infimum in \eqref{eqn:log-sum-exp-inf}
is instead reached by any sequence $\seq{\xx_t}$ for which
$\lsePP(\xx_t)\rightarrow\lsePPh$.
For example, if $n=3$ and
$\lsePPh=\trans{\bigBracks{0,\frac13,\frac23}}$, then the sequence
$\xx_t=\trans{\bigBracks{-t,\ln\frac13,\ln\frac23}}$
converges to the infimum as $t\to+\infty$.%
\indexg{minimizing sequences!examples|)}%
\indexg{Log-sum-exp|)}
\end{example}

This example suggests that, in extending $\Rn$,
it perhaps would suffice to allow
individual coordinates to take on values of $\pm\infty$. In other
words, we could consider a Cartesian product of the extended reals, $(\eR)^n$. The minimizer of \Cref{ex:log-sum-exp} would then be
$\trans{\Bracks{-\infty,\ln\frac13,\ln\frac23}}$. But the next example shows that this
is not enough because minimization may take us in a direction that is not aligned with coordinate axes.

\begin{figure}
  \centering
  \includegraphics{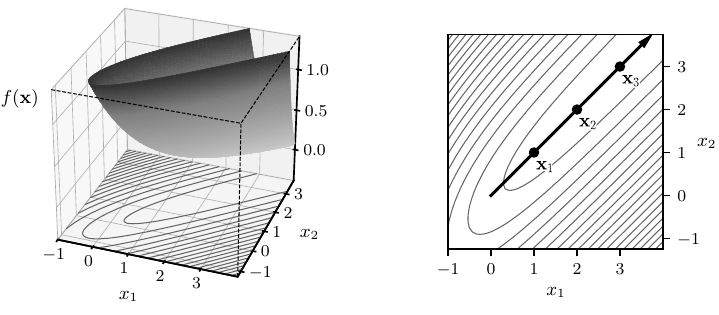}
  \mycaption{Diagonal valley}{%
  \indexf{Diagonal valley}%
    The function $f$ from \Cref{ex:diagonal-valley} (surface plot on the left, contour plot on the right).
    On the right, we also show the first few points of the sequence
    $\xx_t=\trans{[t,t]}$, which minimizes $f$ in the limit as it follows a halfline to infinity.}
  \label{fig:diag-valley}%
\end{figure}

\begin{example}[Diagonal valley]
\label{ex:diagonal-valley}
\indexg{Diagonal valley|(}%
\indexg{minimizing sequences!examples|(}%
Let $f:\R^2\to\R$
be defined, for $\xx\in\R^2$, as
\begin{equation}
\label{eqn:eg-diag-val}
  f(\xx)=f(x_1,x_2)=e^{-x_1}+(x_2-x_1)^2
\end{equation}
(see \Cref{fig:diag-valley}).
The infimum is obtained in the limit of any sequence $\seq{\xx_t}$ that satisfies $x_{t,1}\to+\infty$ while also $x_{t,2}-x_{t,1}\to 0$. One such sequence is $\xx_t=\trans{[t,t]}$, which goes to infinity in the direction of the vector $\vv=\trans{[1,1]}$
(see \Cref{fig:diag-valley}, right).
On this sequence, $f(\xx_t)=e^{-t}\rightarrow 0$.

If we were to just work in $(\Rext)^2$, we would find that $\xx_t\to\trans{[+\infty,+\infty]}$. However, to minimize $f$, the direction in which $\seq{\xx_t}$ goes to infinity is critical, and that direction is not represented by limit points in $(\Rext)^2$.
For instance,
the sequence $\xx'_t=\trans{[2t,t]}$, which goes to infinity in the
direction of $\trans{[2,1]}$,
also converges to $\trans{[+\infty,+\infty]}$ in $(\Rext)^2$, but it fails to minimize $f$ since $f(\xx'_t)\rightarrow+\infty$.%
\indexg{minimizing sequences!examples|)}%
\indexg{Diagonal valley|)}
\end{example}%

So perhaps we should extend $\R^n$ with ``limit points''
of sequences going to infinity in various directions in $\R^n$.
For instance, the sequences $\seq{\xx_t}$ and $\seq{\xx'_t}$, which
follow along halflines with different directions, would then also have
different limits.
Such limit points are considered
by
\idxroc\idxwets%
\citet[Section 3]{rock_wets}, who call them ``direction points''
and develop the topology and properties of the resulting ``cosmic
space.''
This concept was also studied by
\idxhandup\citet{hansen_dupin99}
under the name ``enlarged space.''
Another related concept is that of ``ideal points'' in the real projective plane, where a separate ideal point is introduced for each class of parallel lines
\idxcoxetal%
\citep[see][Chapter 8]{cox_little_oshea__ideals_varieties_algorithms}.

However, these abstractions are too coarse to characterize minimizing sequences in the diagonal valley example.
For example, in the cosmic space formalism, the sequence $\xx''_t=\trans{[t+1,t]}$ converges to the same limit point as $\xx_t=\trans{[t,t]}$, since they both go to infinity in the same direction, but the sequence $\seq{\xx''_t}$ fails to minimize $f$
since $f(\xx''_t)\rightarrow 1$.

Maybe we just need to consider an ``offset'' in addition to the direction of a sequence,
and consider limit points corresponding to following a halfline from a specified starting point in a specific direction.
All the minimizing sequences in the preceding examples can be written as $\xx_t=t\vv+\qq$ for a suitable choice of $\vv$ and $\qq$,
so each indeed proceeds along a halfline
with a starting point $\qq$ in the direction of $\vv$. 
It turns out that this is still not enough, as we show in the next example:

\begin{figure}
  \centering
    \includegraphics{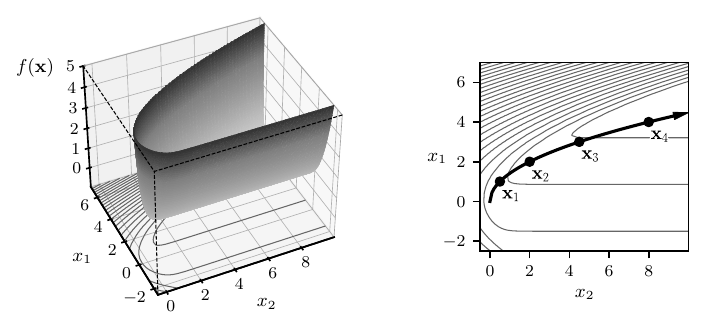}
  \mycaption{Two-speed exponential}{%
  \indexf{Two-speed exponential}%
    The function $f$ from \Cref{ex:two-speed-exp}.
    On the right, we show the first few points of the sequence
    $\smash{\xx_t=\trans{[t,t^2/2]}}$, which minimizes $f$ in the limit as it follows a parabolic curve to infinity.
  }
  \label{fig:twospeed-exp}
\end{figure}

\begin{example}[Two-speed exponential]
\label{ex:two-speed-exp}
\indexg{minimizing sequences!examples|(}%
\indexg{Two-speed exponential|(}%
Let $f:\R^2\to\R$
be defined, for $\xx\in\R^2$, as
\[
  f(\xx)=f(x_1,x_2)=e^{-x_1} + e^{-4x_2+x_1^2}
\]
(see \Cref{fig:twospeed-exp}).
The infimum is
obtained in the limit of any sequence $\seq{\xx_t}$ that satisfies $x_{t,1}\to+\infty$ while also $-4x_{t,2}+x_{t,1}^2\to-\infty$.
This means $f$ cannot be minimized along any halfline; rather, $x_{t,2}$ must go to $+\infty$ faster
than $x_{t,1}^2/4$.
One such sequence is $\smash{\xx_t=\trans{[t,t^2/2]}}$,
on which $f(\xx_t)\rightarrow 0$ (see \Cref{fig:twospeed-exp}, right).%
\indexg{minimizing sequences!examples|)}%
\indexg{Two-speed exponential|)}
\end{example}

The above examples show that the task of adding ``infinite points'' is subtle already in $\R^2$. So perhaps we should just stick to sequences to maintain the broadest flexibility.
The downside is that sequences seem to only be indirect proxies for the underlying ``infinite points,'' so working with sequences makes it harder to discover their structural properties.
Moreover, we hope that by adding the extra ``infinite points,''
the theory of convex functions will become more complete,
for example, by ensuring that all lower-bounded and continuous convex functions (like all those from our examples) have minimizers.

We aim
to extend $\R^n$ analogously to how the set of rational numbers is extended to obtain reals, or the set of reals is extended to obtain complex numbers. When moving from rational numbers to reals to complex numbers,
basic arithmetic operations naturally extend to the enlarged space,
and the enlarged space is more complete;
for example, bounded sets of reals have real-valued infima (which is not true for rational numbers), and all polynomials with complex coefficients have complex roots (whereas polynomials with real coefficients might have no real roots).

In the examples above,
we considered several ways of extending $\Rn$.
With each successive approach,
we saw that we could better capture the breadth of possible minimizers
of convex functions, but at the cost of the extended space becoming
larger and more intricate.
This tradeoff is inevitable, and so will require that we seek a middle
ground that balances the
extent to which we can characterize various minimizers at infinity
against our ability to generalize standard analytical tools to the
resulting space.

Specifically,
we seek an extension of $\R^n$ that will lend itself to natural extensions of the key concepts in convex analysis, like convex sets and functions, conjugacy, subdifferentials, and optimality conditions,
but which will also be more complete,
for example, when it comes to the existence of minimizers.
Our aim is to reveal the structure of how convex functions behave at
infinity.
We seek to build up the foundations of an expanded theory
of convex analysis that will, for instance, make
proving the convergence of optimization algorithms
rather routine, even if the function being optimized has no finite
minimizer.

\section{Extending first-order optimality conditions}
\label{sec:intro:first-order}

To get a taste of such an extended theory, we next examine how the first-order optimality
conditions, which characterize minimizers of convex functions, might generalize to include
points at infinity in just one dimension.

\begin{figure}
  \centering
  \includegraphics{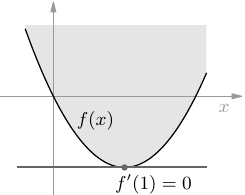}\hspace{0.15\textwidth}%
  \includegraphics{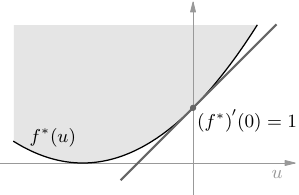}
  \mycaption{A quadratic function and its conjugate}{%
    The function $f(x)=x^2-2x$ and its conjugate $\fstar(u)=u^2/4+u+1$
    (calculated using Eq.~\ref{eq:intro-def-conj}).
    The function $f$ is minimized at $x=1$, with $f'(1)=0$.
    Dually, the derivative of $\fstar$ at $u=0$ is equal to~$1$. Thus,
    $1$ is the sole subgradient of $\fstar$
    at $0$, demonstrating that
    subgradients of $\fstar$ at the origin
    correspond to minimizers of $f$.}
  \label{fig:1d-quad}%
\end{figure}

\indexg{optimality conditions!standard|(}%
Let $f:\Rn\to\eR$ be convex.
One of the fundamental properties of convex functions
is that if
$\gradf(\xx)=\zero$ at a point $\xx\in\Rn$ where
$f$ is differentiable,
then $f$ must be (globally) minimized by $\xx$.
More generally, $\xx$ minimizes $f$
if and only if $\zero$ is a subgradient of $f$ at
$\xx$, meaning geometrically that there exists a
horizontal hyperplane that is tangent to $f$'s epigraph at $\xx$,
as illustrated in \Cref{fig:1d-quad}~(left), where $x=1$ is the minimizer
of the depicted function.
({Subgradients} are generalizations of gradients and will be
discussed further in
\Cref{ex:log-sum-exp:cont};
see also \Cref{sec:prelim:subgrads}.
And the {epigraph} of a function, defined formally in
Eq.~\ref{eqn:epi-def},
is the set of points above its graph,
for instance, the shaded regions in
\Cref{fig:1d-quad}.)

There is also a related ``dual'' characterization of the minimizers of
$f$.
Recall that every convex function $f$ has a
\indexg{conjugate (standard)}%
conjugate $f^*:\Rn\to\Rext$ defined as
\begin{equation}   \label{eq:intro-def-conj}
   f^*(\uu)=\sup_{\xx\in\Rn}\bigBracks{\xx\inprod\uu-f(\xx)},
\end{equation}
for $\uu\in\Rn$,
which is itself a convex function. The conjugate is closely related to
the problem of minimization of $f$:
If $\fstar(\zero)<+\infty$,
then $f$ is bounded below.
Further, under mild conditions on $f$,
it is known that
if the conjugate is differentiable at~$\zero$,
then $\gradf^*(\zero)$ must be the unique minimizer of $f$.
More generally, the subgradients of $f^*$ at $\zero$
are precisely the minimizers of $f$
(as follows from \Cref{pr:stan-subgrad-equiv-props}\ref{pr:stan-subgrad-equiv-props:a}\ref{pr:stan-subgrad-equiv-props:c}).
This is illustrated in \Cref{fig:1d-quad}, which shows a quadratic function $f$ and its conjugate $\fstar$. The function $f$
is minimized at $x=1$, with its derivative equal to~$0$; dually, the derivative of the conjugate at $u=0$ is equal to~$1$.
Since the subgradients of a function at any given point
are exactly the slopes of nonvertical tangents of the function's epigraph
at that point,
the analysis of minimizers of any convex function $f$ can be conducted by analyzing the tangents of the epigraph of $f^*$ at the origin.%
\indexg{optimality conditions!standard|)}%

Can this style of reasoning---a conjugate, subgradients at the origin,
tangents of the epigraph---be lifted to an extended
space that includes points at infinity?
And can we extend the space in such a way that the abundant structure
that these characterizations reflect is still retained?
The next example considers how all this might work out for the simple
case of the exponential function.

\begin{figure}
  \centering
  \includegraphics{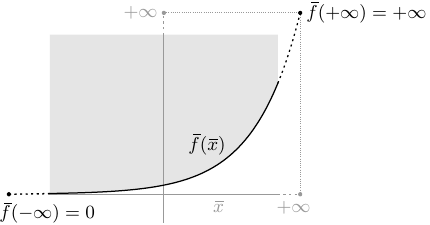}\hfill%
  \includegraphics{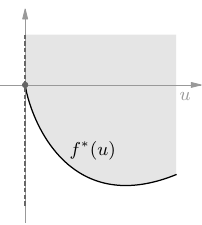}
  \mycaption{Extended exponential and its conjugate}{%
  \indexf{Extended exponential}%
  \indexf{exponential function!extension of}%
    Both the extension $\ef$ and
    the conjugate $\fstar$ are derived in \Cref{ex:exp:cont}.
    \emph{Left:}
    Although $f(x)=e^x$ has no minimizer, its extension $\ef$ is minimized at $\barx=-\infty$.
    \emph{Right:}
    The conjugate $\fstar$ is not differentiable at $0$,
    but its epigraph has a vertical tangent there, corresponding to the ``infinite
    rate of decrease'' and suggesting that it might be
    reasonable to define the derivative at $0$ to be $-\infty$, which would dually
    correspond to $-\infty$ being a minimizer of $\ef$.
  }%
  \label{fig:1d-exp}%
\end{figure}

\begin{example}[Extended exponential]   \label{ex:exp:cont}
\indexg{Extended exponential|(}%
\indexg{exponential function!extension of|(}%
Let $f(x)=e^x$ for $x\in\R$.
As already discussed,
this function has no finite minimizer.
Correspondingly, and
consistent with the characterizations just discussed, there
is no point $x\in\R$ where the function's
derivative is equal to $0$.
Further, the function's conjugate $\fstar$ can be computed to be, for
$u\in\R$,
\begin{equation}   \label{eq:ex:exp:cont:1}
 f^*(u)=
\begin{cases}
  u\ln u-u
    &\text{if $u\ge0$,}
\\
  +\infty
    &\text{if $u<0$,}
\end{cases}
\end{equation}
with the convention $0\ln 0=0$
(see \Cref{fig:1d-exp}, right). This function is
differentiable at all $u>0$, with derivative
$\fstarp(u)=\ln u$.
However, at $u=0$, even though the function is finite,
it is not differentiable since its
epigraph has no nonvertical tangents at that point.
This represents a certain incompleteness in the standard theory,
though consistent with $f$ having no finite minimizers.

Let us see how all this might change if we extend $f$ to $\Rext$.
A natural, continuous extension of $f$ to $\eR$ is the function
$\ef:\eR\to\eR$ defined as
\begin{align*}
  &\ef(x)=e^x\text{ when $x\in\R$},
\\
  &\textstyle
  \ef(-\infty)=\lim_{x\to-\infty}e^x=0,
\\
  &\textstyle
  \ef(+\infty)=\lim_{x\to+\infty}e^x=+\infty,
\end{align*}
(see \Cref{fig:1d-exp}, left).
Unlike $f$, which has no minimum,
the extension $\ef$
appears better-behaved
because it does attain its minimum (at $-\infty$).
Further, if we could somehow define the ``derivative''
of $\fext$ at $-\infty$,
we might reasonably expect that it should be equal
to $0$, since it would seem that an extended version of the horizontal
$x$-axis
is tangent to the epigraph of $\fext$ at $-\infty$
(as suggested by \Cref{fig:1d-exp}, left).

Next, let us attempt to generalize the notion of conjugacy to the extension $\ef$. 
Analogous to \eqref{eq:intro-def-conj},
this conjugate $\fextstar:\R\to\Rext$ could perhaps be defined,
for $u\in\R$, as
\[
  \fextstar(u)
  =\sup_{\ex}
    \bigBracks{\ex\cdot u - \ef(\ex)},
\]
where it is understood that the supremum here is only over those
$\ex\in\Rext$ for which the difference
$\ex\cdot u - \ef(\ex)$ is defined
(that is, not the difference of $+\infty$ with itself,
or $-\infty$ with itself).
This definition of $\smash{\fextstar}$
differs superficially from the one that we introduce
later in \Cref{sec:conjugacy-def}, but the two are nonetheless
equivalent.
It can be checked that
the function $\fextstar$ turns out to be identical to the
standard conjugate $f^*$ given in
\eqref{eq:ex:exp:cont:1}.
Therefore, as noted above, this function's epigraph has no nonvertical tangents at $0$.

Nonetheless, this function's epigraph does have a \emph{vertical}
tangent at $0$ (see \Cref{fig:1d-exp}, right). It would be natural to represent this vertical tangent
by positing the ``derivative'' at this point to be equal to $-\infty$,
corresponding to how $\fextstar$ is ``decreasing infinitely fast'' at $0$ (faster than any finite slope). 
If we could next apply a theorem saying that
the (suitably generalized) subgradients of $\fextstar$ at the origin are minimizers of $\ef$, we would
thereby obtain that $-\infty$ minimizes $\ef$, as is actually the case.%
\indexg{exponential function!extension of|)}%
\indexg{Extended exponential|)}
\end{example}

The above example suggests just one area, besides the existence of
minimizers, in which the theory of an extended space and extended
functions might give rise to a more complete convex analysis.
In this case,
we considered generalizing subgradients to take on
infinite values so that their extended form
might represent not only nonvertical, but also vertical tangents of
the epigraphs of convex functions.
If this can be done successfully,
we might expect all suitable convex functions to
have generalized subgradients everywhere that they are finite,
which is not the case in standard convex analysis.

In the example above,
we considered a function in one dimension and suggested
defining a generalized derivative that takes values in
the extended reals, $\Rext$,
and in particular, is equal to $\pm\infty$ at points
where the function is increasing or decreasing infinitely fast.
It is considerably less clear how gradients or subgradients of a multivariate
function might be similarly generalized since, at a given point,
the function may be increasing (or decreasing) infinitely fast in some
directions, but only increasing at a finite rate in other directions.
We will return to this question in Example~\ref{ex:log-sum-exp:cont}
where we will see how such generalized subgradients can themselves
be represented using points in an extension of $\Rn$ that includes
points at infinity.

To be clear,
our goal is not to merely extend subgradients to include vertical tangents;
rather,
we seek to develop an expansive theory
of convex analysis into which extended subgradients would fit like a
puzzle piece, together with
extended conjugates, optimality conditions, and an array of
other fundamental concepts.

\bigskip

In this book, we propose such a theory, extending $\R^n$ to an enlarged topological space called \emph{astral space}.
We study the new space's properties, and develop a theory of convex
functions on this space. Although astral space includes all of $\R^n$,
it is not a vector space, nor even a metric space. However, it is
sufficiently well-structured to allow useful and meaningful extensions
of such concepts as convexity, conjugacy, and subgradients. We develop these concepts and analyze various properties of convex functions on astral space, including the structure of their minimizers, characterization of continuity, and convergence of descent algorithms.

Although the conceptual underpinnings of astral space are simple, the
full formal development is somewhat involved since key topological properties need to be carefully established.
As a teaser of what the book is about, we next present a condensed development of astral space. This will allow us to revisit our earlier multidimensional examples
and situate them within the theory of functions on astral space.
We then finish this introductory chapter with a high-level
overview of the entire book.

\section{A quick introduction to astral space}
\label{sec:intro:astral}

In constructing astral space, our aim is to derive
a topological extension of $\R^n$ in which various ``points at infinity'' have been added, corresponding to limits of sequences such as those from
Examples~\ref{ex:log-sum-exp}, \ref{ex:diagonal-valley},
and~\ref{ex:two-speed-exp}.
\indexg{astral space!compactification@as compactification|(}%
In fact, we seek 
a \emph{compactification}%
\indexg{compactification}
of $\R^n$, in which every sequence has a convergent subsequence. There are many possible compactifications, which differ in how many new points they add and how finely they differentiate among different kinds of convergence to infinity.
We would like to add as few points as possible, but the fewer points we add, the fewer functions will be continuous in the new space, because there will be more sequences converging to each new point.%
\indexg{astral space!compactification@as compactification|)}%

\indexg{linear functions, continuous extension of|(}%
In this work,
we choose a tradeoff in which we add as few points as possible while still ensuring that all \emph{linear} functions, which are the bedrock of convex analysis, remain continuous.
As a consequence of this property, we will see that ``anything linear'' is likely to
behave nicely when extended to astral space,
including such notions as linear maps, hyperplanes, and halfspaces, as
well as convex sets, convex functions, conjugates, and subgradients, all
of which are
based, in their definitions, on linear functions.%
\indexg{linear functions, continuous extension of|)}

\paragraph{Convergence in all directions.}
To implement this idea, we focus on the limits of all
linear functions when evaluated on a given sequence.
Every such linear function takes the form $\xx\mapsto\xx\inprod\uu$
for some $\uu\in\R^n$, which,
when $\uu$ is a unit vector, can be viewed geometrically as
scalar projection in the direction of $\uu$.
\indexg{convergence in all directions|(}%
We define a notion of a well-behaved sequence
with respect to linear maps by saying that a sequence $\seq{\xx_t}$ in
$\R^n$ \emph{converges in all directions}
if, for all $\uu\in\Rn$,
its image $\seq{\xx_t\inprod\uu}$ under such a linear function
converges in $\eR$, meaning $\lim(\xx_t\inprod\uu)$ exists in
$\eR$ for all $\uu\in\Rn$.%
\indexg{convergence in all directions|)}

For example, every sequence that converges in $\R^n$ also converges in all directions, since $\xx_t\to\xx$ implies, for all $\uu\in\Rn$,
that $\xx_t\inprod\uu\to\xx\inprod\uu$.
There are additional sequences that converge in all directions like all those appearing in
Examples~\ref{ex:log-sum-exp}, \ref{ex:diagonal-valley}, and~\ref{ex:two-speed-exp}.\looseness=-1

\indexg{all-directions equivalence|(}%
If two sequences $\seq{\xx_t}$ and $\seq{\yy_t}$ both converge in all directions
and also ${\lim(\xx_t\inprod\uu)}={\lim(\yy_t\inprod\uu)}$ for all $\uu\in\R^n$, we say that
$\seq{\xx_t}$ and $\seq{\yy_t}$ are \emph{all-directions equivalent}.%
\indexg{all-directions equivalence|)}

\paragraph{Astral space.}
All-directions equivalence creates a partition into equivalence classes
of the set of sequences that converge in all directions.
As formally defined in
\Cref{sec:astral-space-intro},
to construct
\indexg{astral space!construction}%
\emph{astral space},
we associate an \emph{astral point} with each such equivalence class;
this point will then be exactly the common limit of every sequence in the associated equivalence class.
We write $\eRn$ for $n$-dimensional astral space, and use bar or overline notation to denote its
elements, such as $\exx$ or $\eyy$.\looseness=-1

\indexg{coupling function|(}%
This definition allows us to naturally extend the
inner product. Specifically, for all $\xbar\in\extspace$ and all
$\uu\in\R^n$, we define the \emph{coupling function}
\[
  \exx\inprod\uu=\lim(\xx_t\inprod\uu)\in\eR,
\]
where $\seq{\xx_t}$ is any sequence in the equivalence class associated with $\exx$. Note that the value of $\exx\inprod\uu$ does not depend on the choice of $\seq{\xx_t}$ because all the sequences
associated with $\exx$ are all-directions equivalent and therefore have identical limits $\lim(\xx_t\inprod\uu)$. In fact, the values of these limits
uniquely identify each astral point, so every astral point $\xbar$ is uniquely identified by the values of $\xbar\cdot\uu$ over all $\uu\in\R^n$.%
\indexg{coupling function|)}

The space $\R^n$ is naturally included in $\eRn$, coinciding with equivalence classes of sequences that converge to points in $\R^n$.
But there are additional elements, which are said to be
\emph{infinite}%
\indexg{astral points!infinite}%
\indexg{infinite points}
since they must satisfy
$\exx\inprod\uu\in\set{-\infty,+\infty}$ for at least one vector $\uu$.

\indexg{astrons|(}%
The simplest astral points (other than those in $\Rn$) are called \emph{astrons},
each obtained as the limit of a sequence of the form $\seq{t\vv}$ for some $\vv\in\R^n$.
The resulting astron is denoted $\limray{\vv}$.
By reasoning about the limit of the sequence
 $\seq{t\vv\inprod\uu}$, for $\uu\in\Rn$, it can be checked that
\[
  (\limray{\vv})\cdot\uu
       =
\begin{cases}
                 +\infty & \text{if $\vv\cdot\uu>0$,} \\
                 0       & \text{if $\vv\cdot\uu=0$,} \\
                 -\infty & \text{if $\vv\cdot\uu<0$.}
\end{cases}
\]
If $\vv\neq\zero$, then the associated astron $\limray{\vv}$ is infinite, and
is exactly the limit of a sequence that follows a halfline from the origin
in the direction of $\vv$.

In one-dimensional astral space, $\overline{\R^1}$, the only astrons
are $0$, and two others corresponding to $+\infty$ and $-\infty$.
In multiple dimensions, besides $\zero$,
there is a distinct astron $\omega\vv$ associated with every unit vector $\vv$.%
\indexg{astrons|)}

The astron construction can be generalized in a way that turns out to yield all additional astral points, including the limits in Examples~\ref{ex:log-sum-exp}, \ref{ex:diagonal-valley}, \ref{ex:two-speed-exp}.

\begin{example}[Polynomially graded sequence]
\label{ex:poly-speed-intro}
\indexg{Polynomially graded sequence|(}%
Let $\qq\in\Rn$, let
$\vv_1,\dotsc,\vv_k\in\Rn$, for some $k\geq 0$, and
consider the sequence
\begin{equation*}  %
  \xx_t = t^k \vv_1 + t^{k-1} \vv_2 + \dotsb + t \vv_k + \qq
    = \sum_{i=1}^k t^{k-i+1} \vv_i + \qq.
\end{equation*}
We can verify that this sequence converges in all directions, and therefore
corresponds to an astral point, by calculating the limit of $\seq{\xx_t\cdot\uu}$
for any vector $\uu\in\Rn$.

If $k=0$, then $\xx_t=\qq$ for all $t$, so trivially $\xx_t\cdot\uu\rightarrow\qq\cdot\uu$.
Otherwise, if $k\geq 1$, then
the evolution of the sequence $\seq{\xx_t}$ is dominated
by its overwhelmingly rapid growth in the direction of $\vv_1$.
As a result,
if $\vv_1\inprod\uu>0$ then $\xx_t\cdot\uu\rightarrow +\infty$,
and
if ${\vv_1\inprod\uu<0}$ then $\xx_t\cdot\uu\rightarrow -\infty$.
However, if $\vv_1\inprod\uu=0$,
then the term involving $\vv_1$ vanishes when considering
$\xx_t\cdot\uu$.
On the subspace orthogonal to $\vv_1$,
the direction of $\vv_2$ thus becomes dominant.
For such vectors $\uu$ that are perpendicular to $\vv_1$,
we find once again that $\xx_t\cdot\uu$ converges
to $+\infty$ or $-\infty$ if $\vv_2\inprod\uu>0$ or $\vv_2\inprod\uu<0$,
respectively.
This analysis can be continued, so that we next consider vectors
$\uu$ in the subspace orthogonal to both $\vv_1$ and $\vv_2$, where
$\vv_3$ is dominant.
And so on.
Eventually, for vectors $\uu$ that are orthogonal to all
$\vv_1,\dotsc,\vv_k$, we find that
$\xx_t\cdot\uu$ converges to the finite value $\qq\cdot\uu$.\looseness=-1

In summary, this argument shows that the sequence $\seq{\xx_t}$ converges in all directions, and its corresponding astral point $\exx$ is described by
\newcommand{\ttinprod}{\inprod}
\begin{align}
\label{eq:intro-astral-point}
\xbar\ttinprod\uu
&=
\begin{cases}
+\infty
&\text{if $\vv_1\ttinprod\uu=\dotsb=\vv_{i-1}\ttinprod\uu=0$ and $\vv_i\ttinprod\uu>0$ for some $i$,}
\\
-\infty
&\text{if $\vv_1\ttinprod\uu=\dotsb=\vv_{i-1}\ttinprod\uu=0$ and $\vv_i\ttinprod\uu<0$ for some $i$,}
\\
\qq\ttinprod\uu
&\text{if $\vv_1\ttinprod\uu=\dotsb=\vv_k\ttinprod\uu=0$.}%
\indexg{Polynomially graded sequence|)}
\end{cases}
\\[-6.5pt]
\tag*{\qedhere}
\end{align}%
\end{example}

\eqref{eq:intro-astral-point} implies that
if any of the vectors $\vv_j$ in \Cref{ex:poly-speed-intro} is a linear combination of the preceding vectors,
then it can be dropped without affecting the value of the limits
$\lim(\xx_t\inprod\uu)=\xbar\inprod\uu$, so that
the resulting sequence (after $\vv_j$ has been removed)
still converges to the same astral point $\xbar$.
In particular, 
this means,
for the purpose of constructing astral points as the
limits of sequences,
that it suffices
in
\Cref{ex:poly-speed-intro}
to consider only linearly
independent vectors $\vv_1,\dotsc,\vv_k$,
of which there can be at most $n$.

Although this example analyzes only one type of sequence,
we will later see
(by a different style of reasoning) that every sequence that converges in all directions must have a similar structure.
Namely, every
such sequence must have a most dominant direction $\vv_1$ in which it is tending to infinity most rapidly, followed by a second most dominant direction~$\vv_2$, and so on, with residual convergence to some finite point $\qq$.
As a result, as we prove
in \Cref{sec:astral:space:summary},
every astral point in $\extspace$
is characterized by \eqref{eq:intro-astral-point} for some choice of
$\vv_1,\dotsc,\vv_k,\qq\in\Rn$
(for some finite $k\geq 0$).

\paragraph{Astral topology.}
\indexg{astral topology!properties|(}%
In defining astral space, we also must specify a
{topology}, thus determining, for instance,
what it means for a sequence to converge in the space, or
for a function to be continuous on the space.
We specifically endow astral space with the
\emph{astral topology}, formally defined in
\Cref{subsec:astral-pts-as-fcns},
which turns out to be the ``coarsest'' (or ``smallest'') topology
under which all linear maps $\xx\mapsto\xx\inprod\uu$ defined over
$\R^n$ can be extended continuously to $\extspace$
(yielding the maps $\exx\mapsto\exx\inprod\uu$).
As proved in \Cref{thm:i:1} and \Cref{thm:first-count-and-conseq},
astral space has several important properties in this topology.

\indexg{compactness!astral space@of astral space|(}%
The first and perhaps most important property is
\emph{compactness},
which makes it often easier to work with astral space than with $\R^n$.
Because of compactness, every continuous function on astral space
attains its minimum.%
\indexg{compactness!astral space@of astral space|)}

\indexg{first-countability!astral space@of astral space|(}%
Another crucial property is called
\emph{first-countability}.
This property allows us to work mainly with sequences and their
limits, just as we would in $\Rn$, rather than the more involved
machinery afforded by general topology.
First-countability implies, for example, that a function $F$ on
$\extspace$ is continuous if and only if
$F(\xbar_t)\rightarrow F(\xbar)$ for every sequence
$\seq{\xbar_t}$ in $\extspace$ converging to some point
\indexg{first-countability!astral space@of astral space|)}%
$\xbar\in\extspace$.
Although all metric spaces (including $\Rn$) are also first-countable,
astral space for $n\ge 2$ is not metrizable, meaning
its topology cannot be obtained by any metric imposed on the space
(see \Cref{sec:not-second}).

\indexg{sequential compactness!astral space@of astral space|(}%
Since astral space is both compact and first-countable, it is also
\emph{sequentially compact},
meaning every sequence in $\eRn$ has a convergent subsequence,
greatly simplifying the analysis of optimization algorithms.
Many of our proofs depend critically on this property.%
\indexg{sequential compactness!astral space@of astral space|)}

Astral space, as already mentioned, is a compactification of $\Rn$.
In addition to being compact, this means also that $\Rn$ (in its
usual, familiar topology) is a topological subspace of astral space,
and that $\Rn$ is \emph{dense} in $\extspace$, implying
that every astral point is the limit of some sequence
of points in $\Rn$.

In one dimension, $\extspac{1}$ is homeomorphic with $\eR$
(that is, the same topologically).
In fact, we later define $\extspac{1}$ to be \emph{equal} to $\Rext$.

Convergence in astral space in the astral topology is closely related
to the notion of convergence in all directions, on which the space was
constructed.
In particular, a sequence $\seq{\xbar_t}$ in $\extspace$ converges to
a point $\xbar\in\extspace$ if and only if, for all $\uu\in\Rn$,
the sequence $\seq{\xbar_t\cdot\uu}$ converges in $\Rext$ to
$\xbar\cdot\uu$.
In this way, convergence
of the sequence $\seq{\xbar_t}$ in $n$-dimensional astral space
can be analyzed by considering convergence
of $\seq{\xbar_t\cdot\uu}$ in one-dimensional $\Rext$
(for all $\uu\in\Rn$).
Many proofs of convergence are based on this property.

Summarizing:

\begin{letter-compact}
\item \label{intro:prop:compact}
      $\extspace$ is compact.%
\item \label{intro:prop:first}
      $\extspace$ is first-countable.
\item \label{intro:prop:iad:1}
      $\xbar_t\rightarrow\xbar$ in $\extspace$ if and only if
      $\xbar_t\cdot\uu\rightarrow\xbar\cdot\uu$ in $\Rext$
  for all $\uu\in\Rn$.
\item \label{intro:prop:dense}
      $\Rn$ is a dense subspace of $\extspace$.
\item \label{intro:prop:R1}
      $\extspac{1}$ is homeomorphic with $\eR$.%
\indexg{astral topology!properties|)}
\end{letter-compact}

\paragraph{Representing astral points.}
Although astral space is a topological extension of the vector space $\R^n$,
it is not a vector space itself, and astral points cannot be
added.
The root problem is that the sum of $-\infty$ and $+\infty$ is undefined, making it impossible, for example, to establish the identity $(\xbar + \ybar)\inprod\uu=\xbar\inprod\uu+\ybar\inprod\uu$, which is meaningless when $-\infty$ and $+\infty$ appear on the right-hand side.

While standard addition does not generalize to astral space, a
noncommuta\-tive variant does. For $\barx,\bary\in\Rext$, we write this operation, called
\indexg{leftward addition!scalars@of scalars}%
\emph{leftward addition},
as
\indexm{x+y500}{$\barx\protect\plusl\bary$}{leftward sum (of scalars)}%
$\barx\plusl \bary$.                                  
It is the same as
ordinary addition except that,
when adding $-\infty$ and $+\infty$, the argument on the \emph{left} dominates.
Thus,
\begin{align}
  (+\infty) \plusl (-\infty) &= +\infty,
  \nonumber
  \\
  (-\infty) \plusl (+\infty) &= -\infty,
  \label{eqn:intro-left-sum-defn}
  \\
  \barx \plusl \bary &= \barx + \bary \;\;\mbox{in all other cases.}
  \nonumber
\end{align}
\indexg{leftward addition!astral points@of astral points|(}%
This operation can be extended from $\Rext$ to $\extspace$:
For $\xbar,\ybar\in\extspace$,
the leftward sum,
written
\indexm{x+y600}{$\xbar\protect\plusl\ybar$}{leftward sum (of astral points)}%
$\xbar\plusl\ybar$,
is defined to be that
unique point in $\extspace$
for which
\[ (\xbar\plusl\ybar)\cdot\uu = \xbar\cdot\uu \plusl \ybar\cdot\uu \]
for all $\uu\in\Rn$.
Such a point must always exist (see \Cref{pr:i:6}).
While leftward addition is not commutative, it is associative.
Scalar multiplication of a vector and standard matrix-vector multiplication also extend to astral space, and are distributive with leftward addition.%
\indexg{leftward addition!astral points@of astral points|)}

\indexg{representations of astral points!astrons@with astrons|(}%
All astral points can be decomposed as the leftward sum of astrons
and a finite part, corresponding exactly to the
progression of dominant directions that appeared in
\Cref{ex:poly-speed-intro}.
In that example,
\begin{equation}
\label{eq:intro-sum-astrons}
     \xbar = \limray{\vv_1}\plusl\dotsb\plusl\limray{\vv_k}\plusl\qq.
\end{equation}
This can be verified by comparing the definitions and observations above with \eqref{eq:intro-astral-point}.
Such a representation of an astral point is not unique.
\indexg{canonical representation|(}%
Nevertheless, every astral point does have a
unique \emph{canonical representation}
in which the vectors $\vv_1,\dotsc,\vv_k$ are orthonormal,
and $\qq$ is orthogonal to all of them
(see \Cref{pr:uniq-canon-rep}).%
\indexg{representations of astral points!astrons@with astrons|)}%
\indexg{canonical representation|)}

\indexg{rank, astral|(}%
Every astral point
has an intrinsic \emph{astral rank},
which is the
smallest number of astrons needed to represent it,
and which is also equal to the number
of astrons appearing in its canonical representation.
In particular, this means that every astral point in $\extspace$ has
astral rank at most $n$.
Points in $\Rn$ have astral rank~$0$.
Astral points of the form $\omega\vv\plusl\qq$ with $\vv\ne\zero$ have astral rank $1$;
these are the points obtained as limits of sequences that go to infinity along a halfline
with starting point $\qq$ 
in the direction of $\vv$, such as $\xx_t= t\vv+\qq$.%
\indexg{rank, astral|)}

\section{Minimization revisited}
\label{sec:intro-min-revisit}

In each of the examples from \Cref{sec:intro:sequences},
the function being considered was
minimized only via a sequence of points going to infinity.
The preceding development allows us now
to take the limit of each of those sequences in astral space.
The next step is to extend the functions themselves to astral space so
that they can be evaluated at such an infinite limit point.
\indexg{lower semicontinuous extension|(}%
To do so, we focus especially on a natural extension of a function $f:\R^n\to\eR$ to astral space called the \emph{lower semicontinuous extension}
 (or simply,
\emph{extension}) of $f$ to $\eRn$,
which we introduce in \Cref{sec:lsc:ext}.
This function, written $\fext:\extspace\rightarrow\Rext$, is defined at any $\xbar\in\eRn$ as
\[
\fext(\xbar) = \InfseqLiminf{\seq{\xx_t}}{\Rn}{\xx_t\rightarrow \xbar}
                            {f(\xx_t)},
\]
where the infimum is taken over all sequences $\seq{\xx_t}$ in $\Rn$ that converge to $\xbar$.
In words, $\ef(\xbar)$ is the infimum across all possible limit values achievable by any sequence $\bigParens{f(\xx_t)}$ for which $\xx_t\to\xbar$.
For example, this definition yields the extension of the exponential
function given in Example~\ref{ex:exp:cont}.

In general,
the extension $\ef$ is a lower semicontinuous function on a compact
space, and therefore always has a minimizer.
Moreover, by $\fext$'s construction,
an astral point~$\xbar$ minimizes the extension $\fext$ if and only if
there exists a sequence in $\Rn$ converging to~$\xbar$ that minimizes
the original function $f$.
Additionally,
if $\ef$ is not just lower semicontinuous, but actually continuous at
$\xbar$,
and if $\xbar$ minimizes $\fext$, then
\emph{every} sequence in~$\Rn$ converging to $\xbar$ minimizes $f$.
This makes $\fext$'s continuity properties algorithmically appealing; understanding when and where the extension $\fext$ is continuous is therefore
one important focus of this book.%
\indexg{lower semicontinuous extension|)}

\bigskip

Having extended both $\Rn$ and functions on $\Rn$,
we are now ready to return to the
minimization problems from \Cref{sec:intro:sequences}
as an
illustration of the notions discussed above, as well as
a few of the general results appearing later in the book
(though somewhat specialized for simpler presentation).

In what follows, $\ee_1,\dotsc,\ee_n$ denote the standard basis
vectors in $\Rn$
(so $\ee_i$ is all $0$'s except for a $1$ in the $i$-th component).

\begin{example}[Diagonal valley, continued]
\label{ex:diagonal-valley:cont}
\indexg{minimizing sequences!examples|(}%
\indexg{Diagonal valley|(}%
We derive the extension $\fext$ of the diagonal valley function
from \Cref{ex:diagonal-valley}, which turns out to be continuous everywhere. First, we rewrite the function using inner products, since these can be continuously extended to astral space:
\[
  f(\xx)=e^{-x_1}+(x_2-x_1)^2=e^{\xx\inprod(-\ee_1)}+\bigBracks{\xx\inprod(\ee_2-\ee_1)}^2.
\]
Now to obtain the continuous extension, we can just rely on the
continuity of the coupling function to obtain
\[
  \ef(\exx)=e^{\exx\inprod(-\ee_1)}+\bigBracks{\exx\inprod(\ee_2-\ee_1)}^2.
\]
Here, implicitly, the functions $e^x$ and $x^2$ have been extended in the natural way to~$\Rext$ according to their limits as $x\rightarrow\pm\infty$.
Note that both terms of the summation are nonnegative (though possibly $+\infty$), so the sum is always defined.

The minimizing sequence $\xx_t=\trans{[t,t]}$ from \Cref{ex:diagonal-valley} follows a halfline,
and converges in astral space to the point $\xbar=\limray{\vv}$ where $\vv=\trans{[1,1]}$.
This point is an astron and its astral rank is one.
That $\fext$ is both continuous everwhere and has a rank-one minimizer is not a coincidence:
In Sections~\ref{subsec:rank-one-minimizers}
and~\ref{subsec:cond-for-cont},
we prove that if the extension~$\fext$ is continuous everywhere, then it must have a minimizer of astral rank at most one.%
\indexg{Diagonal valley|)}
\looseness=-1
\end{example}

\begin{example}[Two-speed exponential, continued]
\label{ex:two-speed-exp:cont}
\indexg{Two-speed exponential|(}%
Recall the two-speed exponential function from
\Cref{ex:two-speed-exp}:
\[
  f(\xx)=e^{-x_1} + e^{-4x_2+x_1^2}.
\]
Unlike the previous example, this function's extension $\fext$ is not continuous everywhere.
Earlier, we argued that the sequence $\xx_t=\trans{[t,t^2/2]}$ minimizes $f$, satisfying $f(\xx_t)\to 0$.
The sequence $(\xx_t)$ converges to the astral point $\exx=\omega\ee_2\plusl\omega\ee_1$, with astral rank~2.
On the other hand, the sequence $\xx'_t=\trans{[t,t^2/5]}$
also converges to $\omega\ee_2\plusl\omega\ee_1$, but $f(\xx'_t)\to +\infty$.
This shows that $\fext$ is not continuous at $\xbar$, and means, more specifically, that the extension $\ef$ satisfies $\ef(\xbar)=0$, but not all sequences converging to $\xbar$ minimize $f$.\looseness=-1

It turns out that $\xbar$ is the only minimizer of $\ef$, so $\ef$ is also an example of a function that does not have a rank-one minimizer,
and so cannot be minimized by following a halfline to infinity.
As discussed in the previous example, the fact that $\fext$ has no rank-one minimizer implies generally that it cannot be continuous everywhere.%
\indexg{Two-speed exponential|)}
\end{example}

\begin{example}[Log-sum-exp, continued]
\label{ex:log-sum-exp:cont}
\indexg{Log-sum-exp|(}%
The log-sum-exp function from
Example~\ref{ex:log-sum-exp}
has a continuous extension similar to the diagonal valley function, but a direct construction is slightly more involved.  Instead, we take a different route that showcases an exact dual characterization of continuity in terms of a geometric property related to $f$'s conjugate, $\fstar$.
In \Cref{subsec:cond-for-cont}, we will see
that when $f$ is finite everywhere, its extension $\ef$ is continuous everywhere if and only if
the effective domain of $\fstar$
has a conic hull that is \emph{polyhedral}, meaning that it is equal to the intersection of finitely many halfspaces in $\Rn$.
(The \emph{effective domain} of a function that takes values in $\Rext$
is the set of points where the function is not $+\infty$.
The \emph{conic hull} of a set $S\subseteq\Rn$ is obtained by taking
all of the nonnegative combinations of points in $S$.)

For the log-sum-exp function, we have, for $\uu\in\Rn$,
\[
  \fstar(\uu)=\begin{cases}
  \sum_{i=1}^n (u_i+\lsePh_i)\ln(u_i+\lsePh_i)
  &\text{if $\uu+\lsePPh\in\probsim$,}
\\
  +\infty
  &\text{otherwise,}
  \end{cases}
\]
where, as before, $\probsim$ is the probability simplex in $\Rn$.
Hence, the effective domain of $\fstar$ is $\probsim-\lsePPh$, which, being a translation of a simplex, is polyhedral. Its conic hull is therefore also polyhedral,
so $\ef$ is continuous everywhere
(by \Cref{thm:cont-conds-finiteev}\ref{thm:cont-conds-finiteev:e}\ref{thm:cont-conds-finiteev:a}).

Earlier, we considered the particular case that $n=3$ and
$\lsePPh=\trans{[0,\frac13,\frac23]}$.
We saw, in this case, that $f$ is minimized by the sequence
$\xx_t=t\vv+\qq$, where
$\vv=\trans{[-1,0,0]}$ and
$\qq=\trans{[0,\ln\frac13,\ln\frac23]}$.
In astral terms, the sequence $\seq{\xx_t}$ converges to
$\xbar=\limray{\vv}\plusl\qq$, which does indeed minimize $\fext$
(as does every sequence converging to $\xbar$,
since $\fext$ is continuous everywhere).%
\indexg{minimizing sequences!examples|)}%

These facts can further be related to differential properties of $f^*$
(which, as in 
\Cref{ex:exp:cont}, is the same as $\fext$'s conjugate,
as will be defined in \Cref{sec:conjugacy-def}).
As discussed in \Cref{sec:intro:first-order},
every subgradient of $\fstar$ at $\zero$
is a minimizer of $f$.
By its standard definition, $\xx\in\Rn$ is a
\emph{subgradient} of $\fstar$ at $\uu\in\Rn$ if
\begin{equation}  \label{eqn:intro-stand-subgrad}
  \fstar(\uu') \geq \fstar(\uu) + \xx\cdot (\uu'-\uu)
\end{equation}
for all $\uu'\in\Rn$, so that the epigraph of $\fstar$ is supported at $\uu$ by the
affine function (in~$\uu'$) on the right-hand side of the inequality.
In this case,
as in \Cref{ex:exp:cont},
$\fstar$ has \emph{no} subgradients at $\zero$, even though
$\zero$ is in the effective domain of $\fstar$, corresponding to $f$
having no finite minimizers.

This incompleteness goes away when working in astral space.
As hinted in the discussion of Example~\ref{ex:exp:cont},
we can define generalized subgradients which are themselves astral points
and which so include ``infinite subgradients'' corresponding to vertical
tangents to the function's epigraph,
in addition to all standard (finite) subgradients.
Although our formulation in \Cref{sec:astral-dual-subgrad}
is more general, when $\fstar(\uu)\in\R$,
we can define these simply by replacing
the inner product in \eqref{eqn:intro-stand-subgrad} with the
astral coupling function;
thus, an astral point
$\xbar\in\extspace$ is an \emph{astral dual subgradient} of $\fstar$
at $\uu$ if
\begin{equation}   \label{eq:ex:log-sum-exp:cont:1}
  \fstar(\uu') \geq \fstar(\uu) + \xbar\cdot (\uu'-\uu)
\end{equation}
for all $\uu'\in\Rn$.
Whereas, as we just noted, it is possible for a convex function to have no
standard subgradients at a particular point, in astral space, \emph{every}
point must have an astral dual subgradient (\Cref{thm:adsubdiff-nonempty}),
even points not in the function's effective domain (using
a generalized version of the definition in Eq.~\ref{eq:ex:log-sum-exp:cont:1}).
Moreover, the astral dual subgradients of a convex function at a
particular point fully determine the function's instantaneous rate of
change in every direction, even if infinite in some directions and
finite in others (\Cref{thm:adsub-gives-dderiv}).

In this particular case,
from the preceding development, it can be checked that
$\xbar=\limray{\vv}\plusl\qq$, as defined above, does indeed satisfy the
condition in \eqref{eq:ex:log-sum-exp:cont:1}
at $\uu=\zero$, and so is an astral dual subgradient.
By general results given in
\Cref{chp:char-astral-subgrads},
this implies
that $\zero$ is
a so-called
astral \emph{primal} subgradient of $\fext$ at $\xbar$
(by \Cref{cor:strict-adif-fext-inverses}),
which in turn implies that $\xbar$ minimizes $\fext$
(by \Cref{pr:asub-zero-is-min}).

Thus, whereas the conjugate of the original function $f$ has no
subgradient at $\zero$, corresponding to $f$ having no finite minimizer,
in astral space, the conjugate does have an astral dual subgradient at
$\zero$ which, correspondingly, is also a minimizer of
the function's extension $\fext$.%
\indexg{Log-sum-exp|)}%
\end{example}

These examples illustrate some of the key questions and topics that we study in
astral space, including continuity, conjugacy, convexity, the structure of
minimizers, differential theory,
and optimality.

\section{Overview of the book}
\label{sec:intro:overview}

We end this chapter with a brief summary of the different parts of the book.
At the highest level, this book defines astral space, studies its
properties, then develops the theory of convex sets and
convex functions on astral space.
We follow the conceptual framework developed for convex analysis on $\R^n$ by \idxroc\citet{ROC},
which in turn grew out of
\indexa{Fenchel, W.}%
Fenchel's lecture notes~\citeyearpar{fenchel1953convex}.

\overviewsecheading{part:prelim}
We begin with preliminaries and a review of relevant topics in
linear algebra, topology, and convex analysis.
This review is meant to make the book self-contained, and to provide
a convenient reference, not as a substitute for a beginning
introduction to these areas.

\overviewsecheading{part:astral-space}
In the first main part of the book,
we formally construct astral space along the lines
given in \Cref{sec:intro:astral}.
We then define its topology, and prove its main properties, including
compactness and first-countability.

Next, we extend linear maps to astral space and study in detail how
astral points can be represented using astrons and leftward addition,
or even more conveniently using matrices.
We directly link these representations to the sequences converging to
a particular astral point.

We then explore how astral points can be decomposed using their
representations in ways that are not only useful, but also revealing
of the structure of astral space itself.
We end this part with a comparison to the cosmic space of
\citet[Section 3]{rock_wets}.%
\idxroc\idxwets

\overviewsecheading{part:extending-functions}
In the next part, we explore how functions $f$ on $\Rn$ can be
extended to astral space.
We define the lower semicontinuous extension $\fext$, briefly introduced in
\Cref{sec:intro-min-revisit}, and study its properties.

We then define astral conjugates.
We especially focus on the astral \emph{biconjugate} of a function
$f$ on $\Rn$, which
can be seen as an alternative way of extending $f$ to
astral space.
We give precise conditions for when this biconjugate is the same as
the lower semicontinuous extension.

In the last chapter of this part, we derive rules for computing lower
semicontinuous extensions (for instance, for the sum of two functions,
or for the composition of a function with an affine map).

\overviewsecheading{part:convex-sets}
The next part of the book explores convexity in astral space.
We begin by defining astral convex sets, along with such related
notions as astral halfspaces and astral convex hull.
We then look at how astral convex sets can be constructed or operated
upon, for instance, by applying an affine map.
We introduce here some particularly useful operations on sets, such as
\emph{sequential sum}, which is the astral analogue of taking the sum of
two sets in $\Rn$.

Astral cones are another major topic of this part, which we here
define and study in detail.
We especially look at useful operations which relate standard
cones in $\Rn$ to astral cones,
such as various types of closure and polarity.
This study of astral cones turns out to be particularly relevant when
later characterizing continuity of functions.
We also define and explore astral linear subspaces,
which are a special case of astral cones.

With these notions in hand, we are able to prove separation theorems
for astral space, for instance,
showing that any two disjoint, closed, convex sets in
astral space can be strongly separated by an astral hyperplane.
As in standard convex analysis, such theorems are fundamental, with
many consequences.

\overviewsecheading{part:astral-functions}
We next explore properties of functions defined on astral space.
To begin, we define and study astral convex functions in general
(including, but not limited to, the extensions of convex functions on
$\Rn$).
We give precise characterizations for when an astral function is
convex and when it is equal to its own biconjugate.
We also study operations on astral functions that preserve their
convexity.

Next, we explore in detail the structure of
minimizers of astral convex functions, especially the extensions of
convex functions on $\Rn$.
We relate the ``infinite part'' of all such minimizers to a set called
the \emph{astral recession cone},
and their ``finite part'' to a particularly well-behaved convex function on $\Rn$
called the \emph{universal reduction}.
We give characterizations and a kind of procedure for enumerating all
minimizers of a given function.

We then explore the structure of minimizers in specific cases,
focusing especially on the astral rank of minimizers, which turns out
to be related to continuity properties of the function.

In the final chapter of this part, we precisely characterize the set of
points where the extension of a convex function on $\Rn$ is
continuous, and also
characterize when that extension is continuous
everywhere.

\overviewsecheading{part:differential-theory}
Next, we study differential theory extended to
astral space.
We define two forms of astral subgradient:
the first, primal form assigns (finite) subgradients to points at infinity;
the second, dual form assigns possibly infinite subgradients to finite points
(for instance, to points with vertical tangents).
We give several characterizations for working with astral
subgradients,
including one that relates primal and dual subgradients via conjugacy.
Among other properties,
we show that every convex function has astral dual subgradients at
every point, which is not the case for standard subgradients.

Finally,
we derive calculus rules for astral subgradients, for instance, for the
sum of two functions.

\overviewsecheading{part:optimality}
In the last part, we explore how astral space is relevant for solving a convex
optimization problem or for characterizing its solutions.
We take a close look at how Fenchel duality and 
important optimality conditions, such as the KKT
conditions, extend to astral space, making for a more complete
characterization of the solutions of a convex program, including when
those solutions might be at infinity.

We also study iterative methods for minimizing a convex function which
might have no finite minimizer, and apply these to a range of commonly
encountered methods and problems.
In the last chapter, we explore how the many ideas and concepts
developed throughout the book can be applied to the
study of exponential-family distributions, which are of fundamental
interest in statistics and machine learning.

\part{Preliminaries and Background}
\label{part:prelim}

\chapter{General preliminaries and review of linear algebra}
\chaptermark{General preliminaries}
\label{sec:more-prelim}

This chapter provides some notational conventions that we
adopt, and some general background, including a review of a few topics
from linear algebra.

\section{Notational conventions}
\label{sec:prelim-notation}

We write $\R$%
\indexm{r 100}{$\R$}{reals}
for the set of reals,
$\rats$%
\indexm{q}{$\rats$}{rationals}
for the set of rational numbers,
and $\nats$%
\indexm{n}{$\nats$}{natural numbers}%
\indexg{natural numbers}
for the set of (strictly) positive integers.
We let $\eR$ denote the set of extended reals:%
\indexg{extended reals}%
\indexm{r 600}{$\Rext$}{extended reals}
\[\Rext=[-\infty,+\infty]=\R\cup\set{-\infty,+\infty}.\]
Also,
$\Rpos=[0,+\infty)$ is the set of all nonnegative reals,
and
$\Rstrictpos=(0,+\infty)$ is the set of all strictly positive reals;
likewise, $\Rneg=(-\infty,0]$ and $\Rstrictneg=(-\infty,0)$.

In general, we mostly adhere to the following conventions:
Scalars are denoted like~$x$, in italics.
Vectors in $\Rn$ are denoted like $\xx$, in bold.
Points in astral space are denoted like $\xbar$, in bold
with a bar.
Matrices in $\R^{m\times n}$ are denoted like $\A$, in bold
uppercase, with the transpose written as $\trans{\A\!}$%
\indexm{A 300}{$\transA$}{transpose}.
Scalars in $\Rext$ are sometimes written like $\barx$, with a bar.
Greek letters, like $\alpha$ and $\lambda$, are also used for scalars.

Vectors $\xx\in\Rn$ are usually understood to have components
$x_i$, and are taken to have a column shape, so
$\xx=\trans{[x_1,\dotsc,x_n]}$.
\indexg{matrices!notation|(}%
A matrix in $\R^{n\times k}$ with columns $\vv_
1,\dotsc,\vv_k\in\R^n$ is written
\indexm{[...]}{$[\ldots]$}{matrix with given columns}%
$[\vv_1,\dotsc,\vv_k]$.
We also use this notation to piece together a larger matrix from other
matrices and vectors in the natural way.
For instance, if $\VV\in\R^{n\times k}$, $\ww\in\Rn$, and
$\VV'\in\R^{n\times k'}$, then
$[\VV,\ww,\VV']$ is a matrix in $\R^{n\times (k+1+k')}$
whose first $k$ columns are a copy of $\VV$, whose next column is
$\ww$, and whose last $k'$ columns are a copy of $\VV'$.
For a matrix $\VV\in\Rnk$, we write
\indexm{columns}{$\columns{\A}$}{columns of matrix}%
$\columns{\VV}$
for the set of
columns of $\VV$; thus, if $\VV=[\vv_1,\ldots,\vv_k]$ then
$\columns{\VV}=\{\vv_1,\ldots,\vv_k\}$.%
\indexg{matrices!notation|)}

We write
\indexm{e i}{$\ee_i$}{standard basis vector}%
$\ee_1,\dotsc,\ee_n\in\Rn$ for the
\indexg{standard basis vectors}%
\emph{standard basis vectors} (where the dimension $n$ is provided or
implied by the context).
Thus, $\ee_i$ has $i$-th component equal to $1$, and all other
components equal to $0$.
The $n\times n$ identity matrix is written
\indexm{i n, i}{$\Idnn,\,\Iden$}{identity matrix}%
$\Idnn$, or without the
subscript when clear from context.
We write $\zerov{n}$%
\indexm{0 n, 0}{$\zerov{n},\,\zero$}{zero vector}
for the all-zeros vector in $\Rn$, dropping
the subscript when clear from context, and also
write $\zeromat{m}{n}$%
\indexm{0 mxn}{$\zeromat{m}{n}$}{zero matrix}
for the $m\times n$ matrix with entries
that are all equal to zero.

The inner product%
\indexg{inner product}
of vectors $\xx,\yy\in\Rn$ is defined as
\indexm{x s100 y}{$\xx\cdot\yy$}{inner product}%
$\xx\inprod\yy=\sum_{i=1}^n x_i y_i$;
the corresponding
\indexg{norm (Euclidean)}%
\indexg{Euclidean norm}%
(Euclidean) norm
is $\norm{\xx}=\sqrt{\xx\inprod\xx}$,%
\indexm{x 100}{$\norm{\xx}$}{Euclidean norm}
also called the length of $\xx$.
For $\xx,\yy\in\Rn$, the \emph{Cauchy-Schwarz inequality}%
\indexg{Cauchy-Schwarz inequality}
states that $|\xx\cdot\yy|\leq \norm{\xx}\cdot\norm{\yy}$,
and the \emph{triangle inequality}%
\indexg{triangle inequality}
states that
$\norm{\xx+\yy}\leq\norm{\xx}+\norm{\yy}$.
The
\indexg{ball (open Euclidean)}%
\emph{(Euclidean) open ball}
with center $\xx\in\Rn$ and radius $\epsilon>0$ is defined as
\begin{equation}
\label{eqn:open-ball-defn}
\indexm{ball}{$\ball(\xx,\epsilon)$}{Euclidean open ball}
 \ball(\xx,\epsilon) = \Braces{ \zz\in\Rn :\: \norm{\zz-\xx} < \epsilon }.
\end{equation}

\indexg{sequences|(}%
Throughout this book, the variable $t$%
\indexm{t}{$t$}{usual sequence index}
has special meaning as
the usual index of all sequences.
Thus, we write \emph{$\seq{x_t}$ in $X$} for the sequence
$x_1,x_2,\ldots$ with all elements in $X$.
Limits and convergence are taken as $t\to +\infty$,
unless stated otherwise. For example,
$\lim f(x_t)$ means $\lim_{t\to\infty} f(x_t)$.
We say a property of a sequence holds for all $t$ to mean that it
holds for all $t\in \nats$.
We say that it holds for all sufficiently large $t$ to mean that it
holds for all but finitely many values of $t$ in $\nats$.

It is sometimes convenient to allow some elements of a sequence to be
undefined.
For instance, the sequence elements may themselves be quotients, some
of which have denominators equal to $0$.
In such cases, provided only finitely many of the elements are
undefined, we define the sequence's limit and other convergence
properties to be the same as those of the sequence obtained by deleting all of
the (finitely many) undefined elements.%
\indexg{sequences|)}

For a function $f:X\to Y$, the
\indexg{image (of a set)}%
\emph{image of a set $A\subseteq X$ under $f$}
is the set
\[ f(A)=\set{f(x):\:{x\in A}},\]
while the \emph{inverse image of a set $B\subseteq Y$ under $f$}%
\indexg{inverse image (of a set)}
is the set
\[f^{-1}(B)=\set{x\in X:\:f(x)\in B}.\]

Let $f:X\rightarrow Y$.
The function $f$ is \emph{injective}%
\indexg{injective (function)}
(or \emph{one-to-one})%
\indexg{one-to-one (function)}
if for all $x,x'\in X$, if $x\neq x'$ then $f(x)\neq f(x')$.
The function is \emph{surjective}%
\indexg{surjective (function)}
(or is said to map \emph{onto} $Y$)%
\indexg{onto (map)}
if $f(X)=Y$.
The function is \emph{bijective}, or is a \emph{bijection},%
\indexg{bijection}
if it is both injective and surjective.
Every bijective function has an inverse $f^{-1}:Y\rightarrow X$,
meaning
$f^{-1}(f(x))=x$ for all $x\in X$
and
$f(f^{-1}(y))=y$ for all $y\in Y$.
If $f$ is bijective and $B\subseteq Y$, then the inverse image
$f^{-1}(B)$ is the same as the image of $B$ under the inverse map $f^{-1}$, so our two
uses of $f^{-1}$ are consistent with each other.

A family of sets $\{S_1,S_2,\ldots\}$ indexed by $t\in\nats$
is written more succinctly as
\indexm{S t}{$\countset{S}$}{countable family of sets}%
$\countset{S}$.
The \emph{set difference}%
\indexg{set difference}
of sets $X$ and $Y$, written
$X\setminus Y$, is the set
\[
  X\setminus Y
  =
  \Braces{x \in X :\: x\not\in Y}.
\]
Let $X,Y\subseteq\Rn$ and $\Lambda\subseteq\R$.
\indexg{set sum (standard)}%
We define
\[
  X+Y
  =
  \Braces{\xx + \yy :\: \xx \in X,\, \yy\in Y},
  \quad
  \text{and}
  \quad
  \Lambda X
  =
  \Braces{\lambda \xx :\: \lambda\in\Lambda,\, \xx \in X}.
\]
For $\xx,\yy\in\Rn$ and $\lambda\in\R$, we also define
$\xx+Y=\{\xx\}+Y$,
$X+\yy=X+\{\yy\}$,
$\lambda X=\{\lambda\}X$,
$-X=(-1)X$,
and
$X-Y = X+(-Y)$.
For a matrix $\A\in\Rmn$, we likewise define
\[
  \A X = \Braces{ \A \xx :\: \xx\in X}.
\]

\indexg{tuples|(}%
We write tuples as
\indexm{<...>}{$\mtuple{\ldots}$}{tuple with given entries}%
$\mtuple{x_1,\ldots,x_m}$.
The set of all such tuples with
${x_i\in X_i}$, for $i=1,\ldots,m$,
is the set $X_1\times \dotsb\times X_m$.
In the special case that
${\xx_i\in\R^{k_i}}$ for $i=1,\ldots,m$,
$\mtuple{\xx_1,\ldots,\xx_m}$ is a point in
$\R^{k_1}\times\dotsb\times\R^{k_m}$,
which we also interpret
as the equivalent (column) vector in $\R^{k_1+\dotsb+k_m}$
whose first $k_1$ components are given by $\xx_1$,
whose next $k_2$ components by $\xx_2$, and so on.
Thus, this notation can be used to concatenate real vectors of
possibly varying dimension in the natural way.

In particular, (ordered) pairs are written $\rpair{x}{y}$
(which should not be confused with similar notation commonly used
by other authors for inner product).
If, as is often the case, $\xx\in\Rn$ and $y\in\R$, then
$\rpair{\xx}{y}$ can be viewed either as a point in $\Rn\times\R$ or
as a vector in $\R^{n+1}$
whose first $n$ components are given by $\xx$ and last
component by $y$.

As is customary,
for functions $f$ defined on tuples, 
we usually write
$f(x_1,\ldots,x_m)$ as shorthand for $f(\mtuple{x_1,\ldots,x_m})$.
Specifically,
functions $f$ that are defined on $\R^{n+1}$ can be viewed
equivalently as functions on $\Rn\times\R$, so
in such cases, 
for $\xx\in\Rn$ and $y\in\R$,
we write simply $f(\xx,y)$ rather than
$f(\rpair{\xx}{y})$.%
\indexg{tuples|)}

The \emph{epigraph}%
\indexg{epigraph!defined}
 of a function
$f:X\rightarrow \Rext$,
denoted $\epi{f}$,
is the set of pairs $\rpair{x}{y}$, with $x\in X$ and $y\in\R$,
for which $f(x)\leq y$:
\begin{equation}
\label{eqn:epi-def}
\indexm{epi f}{$\epi{f}$}{epigraph}%
  \epi{f} = \regBraces{ \rpair{x}{y} \in X\times\R:\: y \geq f(x) }.
\end{equation}
Likewise, the \emph{hypograph}%
\indexg{hypograph}
 of $f$,
denoted $\hypo{f}$,
is the set 
\begin{equation}
\label{eqn:hypo-def}
\indexm{hypo f}{$\hypo{f}$}{hypograph}%
  \hypo{f} = \regBraces{ \rpair{x}{y} \in X\times\R:\: y \leq f(x) }.
\end{equation}
\indexg{domain, effective|(}%
The \emph{effective domain}
 of $f$
(or simply its \emph{domain}, for short)
is the set of points where $f(x)<+\infty$:
\[
\indexm{domf}{$\dom{f}$}{effective domain}%
  \dom{f} = \regBraces{ x \in X:\: f(x) < +\infty }.%
\indexg{domain, effective|)}%
\]
We write $\inf f$ for the function's infimum,
$\inf_{x\in X} f(x)$, and similarly define
$\sup f$, as well as $\min f$ and $\max f$,
when these are attained.

For any $\alpha\in\Rext$, we write
$f\equiv \alpha$
to mean $f(x)=\alpha$ for all $x\in X$.
Likewise, $f>\alpha$ means
$f(x)>\alpha$ for all $x\in X$, with
$f\geq \alpha$, $f<\alpha$, $f\leq \alpha$ defined similarly.
For any
$f:X\rightarrow \Rext$
and
$g:X\rightarrow \Rext$,
we write $f=g$ to mean $f(x)=g(x)$ for all $x\in X$, and similarly
define $f<g$, $f\leq g$, etc.
We say that $f$ \emph{majorizes}%
\indexg{majorization}
$g$ if $f\geq g$.

We often extend scalar operations to functions in a straightforward
pointwise fashion.
For instance, if
$f:X\rightarrow \Rext$
and
$g:X\rightarrow \Rext$,
then $f+g$ denotes the function $(f+g):X\rightarrow\Rext$ defined by
$(f+g)(x)=f(x)+g(x)$
for $x\in X$ (assuming this sum is defined).
Similarly, $\max\{f,g\}$ is the pointwise maximum of $f$ and $g$,
and if $f_i:X\rightarrow\Rext$ for $i\in\indset$ (where $\indset$ is
any index set), then
$\sup_{i\in\indset} f_i$ is the pointwise supremum of the functions
$f_i$.

\section{Working with \texorpdfstring{$\pm\infty$}{\pmUni\inftyUni}}
\label{sec:prelim-work-with-infty}

\indexg{extended reals!arithmetic|(}%
The sum of $-\infty$ and $+\infty$ is undefined, but other sums and
products involving $\pm\infty$ are defined as usual
(see, for example, \idxroc\citealp[Section~4]{ROC}):
\begin{gather*}
\begin{aligned}
&
 \alpha+(+\infty)=(+\infty)+\alpha=+\infty
&&\text{if $\alpha\in(-\infty,+\infty]$,}
\\
&
 \alpha+(-\infty)=(-\infty)+\alpha=-\infty
&&\text{if $\alpha\in[-\infty,+\infty)$,}
\end{aligned}
\\[\medskipamount]
\begin{aligned}
&
 \alpha\cdot(+\infty)=(+\infty)\cdot\alpha=(-\alpha)\cdot(-\infty)=(-\infty)\cdot(-\alpha)=+\infty
&&\text{if $\alpha\in(0,+\infty]$,}
\\
&
 \alpha\cdot(-\infty)=(-\infty)\cdot\alpha=(-\alpha)\cdot(+\infty)=(+\infty)\cdot(-\alpha)=-\infty
&&\text{if $\alpha\in(0,+\infty]$,}
\end{aligned}
\\[\medskipamount]
 0\cdot(+\infty)=(+\infty)\cdot0=0\cdot(-\infty)=(-\infty)\cdot0=0.
\end{gather*}
Note importantly that $0\cdot(\pm\infty)$ is defined to be $0$.%
\indexg{extended reals!arithmetic|)}

We also define the symbol $\oms$%
\indexm{omega}{$\oms$}{infinity}
to be equal to $+\infty$.
We introduce this synonymous notation for readability and to emphasize
$\oms$'s primary role as a scalar multiplier.
Evidently, for $\barx\in\Rext$,
\[
  \omsf{\barx} =
\begin{cases}
                            +\infty & \text{if $\barx>0$,}\\
                            0       & \text{if $\barx=0$,}\\
                            -\infty & \text{if $\barx<0$.}
\end{cases}
\]
Both $\oms$ and its negation, $-\oms=-\infty$, are of course included
in $\Rext$.
Later, we will extend this notation along the lines
discussed in \Cref{sec:intro:astral}.

\indexg{summability|(}%
We say that extended reals $\alpha,\beta\in\eR$ are \emph{summable}
if their sum is defined, that is,
if it is not the case that one of them is $+\infty$ and the other is $-\infty$.
More generally, extended reals $\alpha_1,\dotsc,\alpha_m\in\eR$ are
{summable} if their sum $\alpha_1+\dotsb+\alpha_m$ is defined,
that is, if $-\infty$ and $+\infty$ are not both included in
$\{\alpha_1,\dotsc,\alpha_m\}$.
Also,
letting $f_i:X\rightarrow\Rext$ for $i=1,\ldots,m$,
we say that $f_1,\ldots,f_m$ are summable
if $f_1(x),\ldots,f_m(x)$ are
summable for all $x\in X$
(so that the function $f_1+\dotsb+f_m$ is defined).%
\indexg{summability|)}

Summable scalars have the following monotonicity property:

\begin{proposition}
\label{pr:summable:inc}
Let $\barx,\bary,\barx',\bary'\in\eR$ be such that $\barx\le\barx'$ and
$\bary\le\bary'$.
If $\barx$ and $\bary$ are summable and $-\infty<\barx+\bary$, then $\barx'$ and $\bary'$
are also summable and $\barx+\bary\le\barx'+\bary'$.
\end{proposition}

Since the ordinary sum of $-\infty$ and $+\infty$ is undefined,
we will make frequent use of modified versions of addition
in which that sum is defined.
The modifications that we make are meant only to resolve
the ``tie'' that occurs when combining $-\infty$ and
$+\infty$;
these operations are the same as ordinary
addition in all other situations.
We define three different such operations,
corresponding to the different ways of resolving such ties.

\indexg{leftward addition!scalars@of scalars|(}%
The first and most important of these,
called
\emph{leftward addition},
 was previously introduced in
\eqref{eqn:intro-left-sum-defn}.
For this operation,
when adding $-\infty$ and $+\infty$, the leftward sum is defined to be
equal to the argument on the left.
Thus,
the leftward sum of extended reals $\barx,\bary\in\Rext$, written
$\barx\plusl \bary$,%
\indexm{x+y500}{$\barx\protect\plusl\bary$}{leftward sum (of scalars)}
is defined as
\begin{equation}
\label{eqn:left-sum-alt-defn}
 \barx \plusl \bary =
\begin{cases}
                            \barx & \text{if $\barx\in\{-\infty,+\infty\}$,}\\
                            \barx + \bary  & \text{otherwise.}
\end{cases}
\end{equation}%
In the special case when the arguments are summable (and in particular, if either is in~$\R$),
this operation is commutative and equivalent to ordinary addition. In the general
case, this operation is not commutative, but it is always associative and distributive.
These and other basic properties are summarized in the next
proposition, whose proof is straightforward.

\begin{proposition}  \label{pr:i:5}
  For all $\barx,\bary,\barz,\barx',\bary'\in\Rext$,
  the following hold:
  \begin{letter-compact}
  \item  \label{pr:i:5a}
    $(\barx\plusl \bary)\plusl \barz=\barx\plusl (\bary\plusl \barz)$.
  \item  \label{pr:i:5b}
    $\lambda(\barx\plusl \bary)=\lambda\barx\plusl \lambda\bary$,
    for $\lambda\in\R$.
  \item  \label{pr:i:5bb}
    $\alpha\barx\plusl \beta\barx=\alpha\barx+\beta\barx
    =(\alpha+\beta)\barx$, for $\alpha,\beta\in\Rext$
    with $\alpha\beta\geq 0$.
  \item  \label{pr:i:5c}
    If $\barx$ and $\bary$ are summable then
    $\barx\plusl\bary=\barx+\bary=\bary+\barx=\bary\plusl\barx$.
  \item  \label{pr:i:5ord}
    If $\barx\leq \barx'$ and $\bary\leq \bary'$ then $\barx\plusl \bary \leq \barx'\plusl \bary'$.
  \item  \label{pr:i:5lim}
    $\barx\plusl\bary_t\to\barx\plusl\bary$, for any sequence $\seq{\bary_t}$ in $\eR$
     with $\bary_t\to\bary$.
  \end{letter-compact}
\end{proposition}

Note that part~(\ref{pr:i:5b}) of \Cref{pr:i:5} holds for
$\lambda\in\R$, but does not
hold in general if $\lambda=\oms$.
For instance, if $x=1$ and $y=-2$ then
$\oms(x\plusl y)=-\infty$
but $\oms x \plusl \oms y = +\infty$.

Also,
part~(\ref{pr:i:5lim}) does not hold in general if the arguments are
reversed.
That is, $\bary_t\plusl\barx$ need not converge to
$\bary\plusl\barx$ (where $\bary_t\rightarrow\bary$).
For instance, suppose
$\barx=-\infty$, $\bary=+\infty$, and $\bary_t=t$ for all $t$.
Then $\bary_t\rightarrow\bary$, $\bary_t\plusl\barx=-\infty$
for all $t$, but $\bary\plusl\barx=+\infty$.%
\indexg{leftward addition!scalars@of scalars|)}

\indexg{downward addition|(}%
\indexg{lower addition|(}%
Next, we define the \emph{downward sum}
(or \emph{lower sum})
of $\ex, \bary\in\Rext$, denoted
$\ex \plusd \bary$,%
\indexm{x+y100}{$\barx\protect\plusd\bary$}{downward sum}
to be equal to $-\infty$ if
\emph{either} $\ex$ or $\bary$ is $-\infty$
(and is like ordinary addition otherwise).
Thus,
\begin{equation}  \label{eq:down-add-def}
 \ex \plusd \bary =
 \begin{cases}
   -\infty
   & \text{if $\ex=-\infty$ or $\bary=-\infty$,}
 \\
   \ex + \bary
   & \text{otherwise.}
 \end{cases}
\end{equation}
The next proposition summarizes some properties of this operation:

\begin{proposition}  \label{pr:plusd-props}
  For all $\barx,\bary,\barz,\barx',\bary'\in\Rext$,
  the following hold:
  \begin{letter-compact}
  \item  \label{pr:plusd-props:a}
    $\barx\plusd \bary = \bary\plusd \barx$.
  \item  \label{pr:plusd-props:b}
    $(\barx\plusd \bary) \plusd \barz = \barx\plusd (\bary \plusd \barz)$.
  \item  \label{pr:plusd-props:lambda}
    $\lambda(\barx\plusd \bary) = \lambda\barx \plusd \lambda\bary$, for $\lambda\in\Rpos$.
  \item  \label{pr:plusd-props:c}
    If $\barx$ and $\bary$ are summable,
    then
    $\barx\plusd \bary = \barx+\bary$.
  \item  \label{pr:plusd-props:d-gen}
    For $S\subseteq\Rext$,
    $
    \sup \braces{\barx \plusd \barw:\:\barw\in S}
    =
    \barx\plusd \sup \braces{\barw:\:\barw\in S}
    $.
  \item  \label{pr:plusd-props:d}
    $\sup \braces{\barx - w:\:w\in\R,\,w\geq \bary} = -\bary \plusd \barx$.
  \item  \label{pr:plusd-props:e}
    $\barx \geq \bary \plusd \barz$
    if and only if
    $-\bary \geq -\barx \plusd \barz$.
  \item  \label{pr:plusd-props:f}
    If $\barx\leq\barx'$ and $\bary\leq\bary'$
    then $\barx\plusd\bary\leq\barx'\plusd\bary'$.
  \end{letter-compact}
\end{proposition}

\begin{proof}
Part~(\ref{pr:plusd-props:d-gen}) can be checked by separately
considering the following cases:
$\sup S=-\infty$
(which includes the case $S=\emptyset$);
$\sup S>-\infty$ and $\barx\in\{-\infty,+\infty\}$;
$\sup S>-\infty$ and $\barx\in\R$.
Part~(\ref{pr:plusd-props:d}) is a special case of
part~(\ref{pr:plusd-props:d-gen}) with $S=-\{w\in\R :\: w\geq\bary\}$.

The other parts can be checked in a similar fashion.%
\indexg{lower addition|)}%
\indexg{downward addition|)}%
\end{proof}

\indexg{upward addition|(}%
Finally, we define the
\emph{upward sum}
(or \emph{upper sum})
of $\barx,\bary\in\Rext$, denoted
\indexm{x+y200}{$\barx\protect\plusu\bary$}{upward sum}%
$\barx\plusu\bary$,
to be
\begin{equation}  \label{eq:up-add-def}
 \barx \plusu \bary =
 \begin{cases}
   +\infty
   & \text{if $\barx=+\infty$ or $\bary=+\infty$,}
 \\
   \barx + \bary
   & \text{otherwise.}
 \end{cases}
\end{equation}
Analogous properties to those given in
\Cref{pr:plusd-props}
can be proved for upward sum.
Here are some facts relating upward and downward addition:

\begin{proposition}  \label{pr:plusd-plusu-props}
  For all $\barx,\bary,\barz\in\Rext$,
  the following hold:
  \begin{letter-compact}
  \item  \label{pr:plusd-plusu-props:a}
    $-(\barx\plusd \bary) = (-\barx)\plusu(-\bary)$.
  \item  \label{pr:plusd-plusu-props:b}
    $\barx \geq \bary \plusd \barz$
    if and only if
    $\barx \plusu (-\bary) \geq \barz$.
  \end{letter-compact}
\end{proposition}

\begin{proof}
For part~(\ref{pr:plusd-plusu-props:b}), we have
\[
  \barx \geq \bary \plusd \barz
  \;\Leftrightarrow\;
  -\barz \geq \bary \plusd (-\barx)
  \;\Leftrightarrow\;
  \barz \leq (-\bary)\plusu \barx,
\]
where the first equivalence is by
\Cref{pr:plusd-props}(\ref{pr:plusd-props:e}),
and the second by
part~(\ref{pr:plusd-plusu-props:a}).%
\indexg{upward addition|)}
\end{proof}

\indexg{extended reals!convergence in|(}%
A sequence $\seq{\alpha_t}$ in $\eR$ converges to $+\infty$ if for every $M\in\R$
the sequence eventually stays in $(M,+\infty]$. It converges to $-\infty$
if for every $M\in\R$ the sequence eventually stays in $[-\infty,M)$. Finally, it converges to $x\in\R$ if for every $\epsilon>0$ the sequence eventually stays in
\indexg{extended reals!convergence in|)}%
$(x-\epsilon,x+\epsilon)$.
Limits in $\eR$ satisfy the
following:

\begin{proposition}
\label{prop:lim:eR}
\indexg{limit inferior and superior|(}%
Let $\seq{\alpha_t}$ and $\seq{\beta_t}$ be sequences in $\eR$. Then:
\begin{letter-compact}
\item \label{i:liminf:eR:sum}
      $\liminf(\alpha_t\plusd\beta_t)\ge(\liminf\alpha_t)\plusd(\liminf\beta_t)$.
\item \label{i:liminf:eR:mul}
      $\liminf(\lambda\alpha_t)=\lambda(\liminf\alpha_t)$ for all $\lambda\in\Rpos$.
\item \label{i:liminf:eR:min}
      $\liminf\bigParens{\min\set{\alpha_t,\beta_t}}
       =
       \min\bigBraces{\liminf\alpha_t,\,\liminf\beta_t}$.
\medskip
\item \label{i:limsup:eR:sum}
      $\limsup(\alpha_t\plusu\beta_t)\le(\limsup\alpha_t)\plusu(\limsup\beta_t)$.
\item \label{i:limsup:eR:mul}
      $\limsup(\lambda\alpha_t)=\lambda(\limsup\alpha_t)$ for all $\lambda\in\Rpos$.
\item \label{i:limsup:eR:max}
      $\limsup\bigParens{\max\set{\alpha_t,\beta_t}}
       =
       \max\bigBraces{\limsup\alpha_t,\,\limsup\beta_t}$.%
\indexg{limit inferior and superior|)}%
\end{letter-compact}
\indexg{extended reals!continuity of operations|(}%
Furthermore, if $\alpha_t\to\alpha$ and $\beta_t\to\beta$
where $\alpha,\beta\in\eR$, then:
\begin{letter-compact}[resume]
\item \label{i:lim:eR:sum}
  If $\alpha$ and $\beta$ are summable,
  then $\alpha_t$ and $\beta_t$ are summable for all sufficiently
  large $t$, and $\alpha_t+\beta_t\to\alpha+\beta$.

\item \label{i:lim:eR:genmul}
  If $\alpha\beta\neq 0$ or if $\alpha$ and $\beta$ are both finite,
  then $\alpha_t \beta_t \rightarrow \alpha \beta$.

\item \label{i:lim:eR:mul}
      $\lambda\alpha_t\to\lambda\alpha$ for all $\lambda\in\R$.
\end{letter-compact}
\end{proposition}

\begin{proposition} \label{prop:limsup:eR:conv-comb}
  Let $\gamma\in\Rext$,
  and, for $i=1,\dotsc,m$, let $\seq{\gamma_{it}}$ be a sequence
  in~$\Rext$ 
  and $\seq{\lambda_{it}}$ be a sequence in $\Rpos$.
  Assume
   $\gamma_{it}\to\gamma$ for all $i$, and that
   $\sum_{i=1}^m\lambda_{it}\rightarrow 1$.
   Then $\sum_{i=1}^m \lambda_{it}\gamma_{it}\rightarrow\gamma$.%
\indexg{extended reals!continuity of operations|)}%
\end{proposition}

\section{Linear algebra}
\label{sec:prelim:lin-alg}

\indexg{linear maps (standard)|(}%
Every matrix $\A\in\Rmn$ is associated with a linear map
$A:\Rn\rightarrow\Rm$ defined by $A(\xx)=\A\xx$ for $\xx\in\Rn$.%
\indexg{linear maps (standard)|)}%

A nonempty subset $L\subseteq\Rn$ is a \emph{linear subspace}%
\indexg{linear subspaces (standard)}
(also called simply a \emph{subspace})
if it is closed
under vector addition and multiplication by any scalar, that is,
if $\xx+\yy\in L$ and $\lambda\xx\in L$ for all $\xx,\yy\in L$ and
$\lambda\in\R$.
The intersection of an arbitrary collection of linear subspaces is
also a linear subspace.
The \emph{span}%
\indexg{span (standard)}
of a set $S\subseteq\Rn$, denoted
\indexm{span s300}{$\spn{S}$}{span (standard)}%
$\spn{S}$,
is the smallest linear
subspace that includes $S$, that is, the intersection of all linear
subspaces that include $S$.

A point $\zz\in\Rn$ is a \emph{linear combination}%
\indexg{linear combinations}
of points
$\xx_1,\ldots,\xx_k\in\Rn$, if
$\zz=\lambda_1\xx_1+\dotsb+\lambda_k\xx_k$
for some $\lambda_1,\ldots,\lambda_k\in\R$.

\begin{proposition}  \label{pr:span-is-lin-comb}
  Let $S\subseteq\Rn$.
  Then $\spn{S}$ consists of all linear combinations of zero or more
  points in $S$.
\end{proposition}

A set of vectors $\xx_1,\ldots,\xx_k\in\Rn$ is
\emph{linearly dependent}%
\indexg{linear dependence}
if there exist
$\lambda_1,\ldots,\lambda_k\in\R$, not all $0$, such that
$\lambda_1\xx_1+\dotsb+\lambda_k\xx_k=\zero$,
or equivalently, if one of the vectors is a linear combination of the
others.
If $\xx_1,\ldots,\xx_k$ are not linearly dependent, then they are
\indexg{linear independence}%
\emph{linearly independent}.
A \emph{basis} for a linear subspace $L\subseteq\Rn$ is a
linearly independent subset $B\subseteq\Rn$ that spans $L$
(meaning $L=\spn{B}$).
Every basis for $L$ is of the same cardinality, called the
\indexg{linear subspaces (standard)!dimension of}%
\indexg{dimension!linear subspace@of linear subspace}%
\emph{dimension} of $L$, and denoted
\indexm{dim l}{$\dim{L}$}{dimension (standard)}%
$\dim{L}$.
The standard basis vectors, $\ee_1,\ldots,\ee_n$,
form the
\indexg{standard basis}%
\emph{standard basis}
for $\Rn$.

\indexg{orthogonality|(}%
We say that vectors $\xx,\yy\in\Rn$ are \emph{orthogonal} and write $\xx\perp\yy$%
\indexm{x s200 y}{$\xx\perp\yy$}{orthogonal to}
 if $\xx\inprod\yy=0$.
We say that $\xx\in\Rn$ is orthogonal to a matrix
$\A\in\R^{n\times k}$
and write $\xx\perp\A$ if $\xx$ is orthogonal to all columns of
$\A$, that is, if $\xtransA=\zero$.
Still more generally, two matrices $\A\in\R^{n\times k}$
and $\A'\in\R^{n\times k'}$ are said to be orthogonal,
written $\A\perp\A'$,
if every column of $\A$ is orthogonal to every column of $\A'$.

A set $B\subseteq\Rn$ is orthogonal if every vector in $B$ is
orthogonal to every other vector in
\indexg{orthogonality|)}%
$B$.
The set is \emph{orthonormal}%
\indexg{orthonormality}
if it is orthogonal and if also all vectors
in the set have unit length.

The \emph{column space}%
\indexg{column space!standard}
of a matrix $\A\in\R^{n\times k}$ is the span of its columns,
denoted
\indexm{col}{$\colspace\A$}{column space (standard)}%
$\colspace\A$.
\indexg{rank (of matrix)|(}%
The \emph{rank}
of $\A$
is the dimension of its column space.

\begin{proposition}  \label{pr:row-rank-is-col-rank}
  The rank of matrix $\A\in\R^{n\times k}$ 
  is equal to the rank of its transpose,~$\transA$.
\end{proposition}

\begin{proof}
  See \citet[Theorem~1.16]{roman-lin-alg}\idxroman.
\end{proof}

\indexg{full column or row rank|(}%
We say $\A$ has
\emph{full column rank} if its columns are linearly
independent,
and that it has
\emph{full row rank} if its transpose, $\transA$,
has full column
\indexg{rank (of matrix)|)}%
\indexg{full column or row rank|)}%
rank.
We say $\A$ is \emph{column-orthogonal}%
\indexg{column-orthogonal (matrix)}
if its columns are orthonormal, that is, if the columns all have unit length and are orthogonal to one
another.

If $L_1$ and $L_2$ are linear subspaces of $\Rn$, then
$L_1+L_2$ and $L_1\cap L_2$ are as well.

\begin{proposition}  \label{pr:subspace-intersect-dim}
  Let $L_1$ and $L_2$ be linear subspaces of $\Rn$.
  Then
  \[
     \dim(L_1+L_2)
     +
     \dim(L_1\cap L_2)
     =
     \dim{L_1}
     +
     \dim{L_2}.
  \]
\end{proposition}

\begin{proof}
  See \citet[Theorem~1.14]{roman-lin-alg}\idxroman.
\end{proof}

A square matrix $\RR\in\R^{n\times n}$ with entries $r_{ij}$, where $i$ indexes
rows and $j$ columns,
is \emph{upper triangular}%
\indexg{matrices!upper triangular}%
\indexg{upper triangular matrices}
 if $r_{ij} = 0$
whenever $j<i$.
\indexg{positive upper triangular matrices|(}%
We say $\RR$ is \emph{positive upper triangular}
if it is upper
triangular and if all of its diagonal entries are strictly positive.

\begin{proposition}
\label{prop:pos-upper}
~
  \begin{letter-compact}
  \item  \label{prop:pos-upper:prod}
    The product of two positive upper triangular matrices is also
    positive upper triangular.
  \item  \label{prop:pos-upper:inv}
    Every positive upper triangular matrix $\RR$ is invertible,
    and its inverse $\Rinv$ is also positive upper triangular.
  \end{letter-compact}
\end{proposition}

\begin{proof}
  ~

\begin{proof-parts}
\pfpart{Part~(\ref{prop:pos-upper:prod}):}
This is straightforward to check.
(See also \citealp[Section~0.9.3]{horn_johnson_2nd}\idxhj).

\pfpart{Part~(\ref{prop:pos-upper:inv}):}
Let $\RR\in\R^{n\times n}$ be positive upper triangular
with entries $r_{ij}$
and columns $\rr_1,\ldots,\rr_n$.
By induction on $j=1,\ldots,n$, we claim that
each standard basis vector
$\ee_j$ can be written as
\begin{equation}   \label{eq:prop:pos-upper:1}
  \ee_j=\sum_{i=1}^j  s_{ij} \rr_i
\end{equation}
for some $s_{1j},\ldots,s_{jj}\in\R$ with
$s_{jj}>0$.
This is because
$\rr_j=\sum_{i=1}^j r_{ij} \ee_i$,
so
\[
  \ee_j
   =
   \frac{1}{r_{jj}}
             \brackets{ \rr_j - \sum_{i=1}^{j-1} r_{ij}\ee_i }
   =
   \frac{1}{r_{jj}}
             \brackets{ \rr_j - \sum_{i=1}^{j-1} r_{ij} \sum_{k=1}^i  s_{ki} \rr_k },
\]
where the second equality is by inductive hypothesis.
Thus, $\ee_j$ can be written as a linear combination of $\rr_1,\ldots,\rr_j$,
and therefore as in \eqref{eq:prop:pos-upper:1};
in particular, $s_{jj}=1/r_{jj}>0$.
Setting $s_{ij}=0$ for $i>j$ and letting $\Ss$
denote the $n\times n$ matrix with entries $s_{ij}$,
it follows that $\Ss$ is positive upper
triangular, and that $\RR\Ss=\Iden$, that is, $\Ss=\Rinv$.
\qedhere
\end{proof-parts}
\end{proof}

\begin{proposition}[QR factorization]
\label{prop:QR}
\indexg{QR factorization|(}%
Let $\A\in \R^{m\times n}$ have full column rank.
Then there exists a column-orthogonal matrix $\QQ\in\R^{m\times n}$
and a positive upper triangular matrix
$\RR\in\R^{n\times n}$ such that $\A=\QQ\RR$.
Furthermore, the factors $\QQ$ and $\RR$ are uniquely determined.
\end{proposition}
\begin{proof}
See Theorem~2.1.14(a,b,e) of
\indexg{QR factorization|)}%
\indexg{positive upper triangular matrices|)}%
\citet{horn_johnson_2nd}\idxhj.
\end{proof}

\section{Orthogonal complement, projection matrices, pseudoinverse}
\label{sec:prelim:orth-proj-pseud}

\indexg{orthogonal complement (primal)!standard sets@for standard sets|(}%
For a set $S\subseteq\Rn$, we define the
\emph{orthogonal complement} of $S$, denoted $\Sperp$,
to be the set of vectors
orthogonal to all points in $S$, that is,
\begin{equation}  \label{eq:std-ortho-comp-defn}
\indexm{s 500}{$\Sperp$}{(primal) orthogonal complement}%
  \Sperp
  =
  \Braces{ \xx\in\Rn :\: \xx\cdot\uu=0 \text{ for all }\uu\in S }.
\end{equation}
Note that this definition is sometimes only applied when $S$ is a linear
subspace, but we allow $S$ to be an arbitrary subset.
We write $\Sperperp$ for $(\Sperp)^\bot$.

\begin{proposition}  \label{pr:std-perp-props}
  Let $S,U\subseteq\Rn$.
  \begin{letter-compact}
  \item  \label{pr:std-perp-props:a}
    $\Sperp$ is a linear subspace.
  \item  \label{pr:std-perp-props:b}
    If $S\subseteq U$ then $\Uperp\subseteq\Sperp$.
  \item  \label{pr:std-perp-props:c}
    $\Sperperp = \spn{S}$.
    In particular, if $S$ is a linear subspace then $\Sperperp=S$.
  \item  \label{pr:std-perp-props:d}
    $(\spn{S})^{\bot} = \Sperp$.
  \end{letter-compact}
\end{proposition}

If $L\subseteq\Rn$ is a linear subspace of dimension $k$,
then its orthogonal complement $L^\perp$ is a linear subspace of
dimension
\indexg{orthogonal complement (primal)!standard sets@for standard sets|)}%
$n-k$.
\indexg{orthogonal projection|(}%
\indexg{projection (orthogonal)|(}%
Furthermore, every point can be uniquely decomposed as the sum of some
point in $L$ and some point in $\Lperp$:

\begin{proposition}   \label{pr:lin-decomp}
  Let $L\subseteq\Rn$ be a linear subspace, and let $\xx\in\Rn$.
  Then there exist unique vectors $\xpar\in L$ and $\xperp\in \Lperp$ such
  that $\xx=\xpar+\xperp$.
\end{proposition}

The vector $\xpar$ appearing in this proposition is called the
\emph{orthogonal projection of $\xx$ onto $L$}.
The mapping $\xx\mapsto\xpar$, called
\emph{orthogonal projection onto $L$},
is a linear map described by a unique
matrix $\PP\in\R^{n\times n}$ called the
\indexg{matrices!projection (orthogonal)|(}%
\indexg{projection matrices!orthogonal|(}%
\indexg{orthogonal projection matrices|(}%
\emph{orthogonal projection matrix onto $L$};
that is, $\xpar=\PP\xx$ for all $\xx\in\Rn$.
(When clear from context, we often omit ``orthogonal'' when referring
to each of these.)
Specifically, we have:

\begin{proposition}  \label{pr:basis-to-proj-mat}
  Let $L\subseteq\Rn$ be a linear subspace, 
  let $\set{\vv_1,\dotsc,\vv_k}$ be an orthonormal basis for $L$,
  and let $\VV=[\vv_1,\dotsc,\vv_k]$.
  Let
  $\PP=\VV\trans{\VV}$,
  and let $\PP'=\Iden-\PP$ (where $\Iden$ is the $n\times n$ identity
  matrix).
  Then $\PP$ is the orthogonal projection matrix onto $L$,
  and $\PP'$ is the orthogonal projection matrix onto $\Lperp$.
\end{proposition}

\begin{proof}
That $\PP$ is the orthogonal projection matrix onto $L$
follows as a special case of
\citet[Theorem~2.6]{yanai-takeuchi-takane}\idxytt.
Combined with \Cref{pr:lin-decomp}, this then also implies
that $\PP'$ is the orthogonal projection matrix onto $\Lperp$.
\end{proof}

As in the \namecref{pr:basis-to-proj-mat},
projection onto $L^\perp$, the orthogonal complement of linear
subspace $L$,
is called \emph{projection orthogonal to~$L$}.
In particular,
when $L=\spn\set{\vv}$ for some $\vv\in\Rn$,
we call this \emph{projection orthogonal to $\vv$}.

Here are properties of orthogonal projection matrices:

\begin{proposition}  \label{pr:proj-mat-props}
  Let $\PP\in\R^{n\times n}$ be the orthogonal projection matrix onto
  some linear subspace $L\subseteq\Rn$.
  Let $\xx\in\Rn$.
  Then:
  \begin{letter-compact}
  \item     \label{pr:proj-mat-props:a}
    $\PP$ is symmetric.
  \item     \label{pr:proj-mat-props:b}
    $\PP^2=\PP$.
  \item     \label{pr:proj-mat-props:c}
    $\PP\xx\in L$.
  \item     \label{pr:proj-mat-props:d}
    If $\xx\in L$ then $\PP\xx=\xx$.
  \item     \label{pr:proj-mat-props:e}
    If $\xx\in\Lperp$ then $\PP\xx=\zero$.
  \end{letter-compact}
\end{proposition}

\begin{proof}
For
  parts~(\ref{pr:proj-mat-props:a})
  and~(\ref{pr:proj-mat-props:b}),
see \citet[Theorem~2.5]{yanai-takeuchi-takane}\idxytt.
For
  parts~(\ref{pr:proj-mat-props:d})
  and~(\ref{pr:proj-mat-props:e}),
see \citet[Theorem~2.2]{yanai-takeuchi-takane}\idxytt.
Part~(\ref{pr:proj-mat-props:c})
follows immediately from how projection is
\indexg{orthogonal projection|)}%
\indexg{projection (orthogonal)|)}%
\indexg{matrices!projection (orthogonal)|)}%
\indexg{orthogonal projection matrices|)}%
\indexg{projection matrices!orthogonal|)}%
defined.
\end{proof}

From \Cref{pr:lin-decomp}, it follows that,
for any matrix $\VV\in\Rnk$,
every vector in~$\Rn$ can be decomposed
as a linear combination of the columns of $\VV$
plus a vector that is orthogonal to all of those columns:

\begin{proposition}  \label{pr:lin-decomp-rel-vecs}
  Let $\VV\in\Rnk$,
  and
  let $\xx\in\Rn$.
  Then
    there exists $\bb\in\Rk$ and unique $\qq\in\Rn$ such that
    $\xx=\VV\bb + \qq$ and $\qq\perp\VV$.
    Specifically, $\qq=\PP\xx$ where $\PP$ is the orthogonal
    projection matrix onto $(\colspace \VV)^{\perp}$.
  Furthermore,
    if $\VV$ has full column rank,
    then the vector $\bb$ is also uniquely determined.
  In particular,
    if $\VV$ is column-orthogonal, then
    $\bb=\trans{\VV}\xx$.
\end{proposition}

\begin{proposition}  \label{pr:pseudoinv-dfn}
\indexg{pseudoinverse|(}%
\indexg{Moore-Penrose generalized inverse|(}%
For every matrix $\A\in\R^{m\times n}$, there exists
a unique matrix $\Adag\in\R^{n\times m}$%
\indexm{A 700}{$\Adag$}{pseudoinverse}
satisfying the following conditions:
\begin{letter-compact}
\item \label{pr:pseudoinv-dfn:1}
  $\A\Adag\A=\A$.
\item \label{pr:pseudoinv-dfn:2}
  $\Adag\A\Adag=\Adag$.
\item \label{pr:pseudoinv-dfn:3}
  $\trans{(\A\Adag)}=\A\Adag$.
\item \label{pr:pseudoinv-dfn:4}
  $\trans{(\Adag\A)}=\Adag\A$.
\end{letter-compact}
\end{proposition}

\begin{proof}
  See \citet[Theorem~1]{penrose-gen-inv}\indexa{Penrose, R.}.
\end{proof}

The matrix $\Adag\in\R^{n\times m}$
given in \Cref{pr:pseudoinv-dfn} is
called $\A\negKern$'s \emph{pseudoinverse}
(or \emph{Moore-Penrose generalized inverse}).

\begin{proposition}   \label{pr:pseudoinv-props}
Let $\A\in\R^{m\times n}$.
Then:
\begin{letter-compact}
\item   \label{pr:pseudoinv-props:dagdag}
  $(\Adag)^\pseudinvsym=\A$
  and
  $\trans{(\Adag)}=(\transA)^\pseudinvsym$.
\item   \label{pr:pseudoinv-props:b}
 $\A\Adag$ is the orthogonal projection matrix onto $\colspace\A$.
\item   \label{pr:pseudoinv-props:a}
 $\colspace\regParens{\trans{(\Adag)}}=\colspace\A$.
\item   \label{pr:pseudoinv-props:c}
 If $\A$ has full column rank, then $\Adag\A=\Idn{n}$.
\item   \label{pr:pseudoinv-props:d}
 If $\A$ is column-orthogonal, then $\Adag=\transA$.
\end{letter-compact}
\end{proposition}

\begin{proof}
  ~

\begin{proof-parts}
\pfpart{Part~(\ref{pr:pseudoinv-props:dagdag}):}
See \citet[Lemma~1.1 and~1.2]{penrose-gen-inv}\indexa{Penrose, R.}.

\pfpart{Part~(\ref{pr:pseudoinv-props:b}):}
See \citet[Eq.~3.88]{yanai-takeuchi-takane}\idxytt.

\pfpart{Part~(\ref{pr:pseudoinv-props:a}):}
Let $\B=\trans{(\Adag)}$;
we aim to show $\colspace\A=\colspace\B$.
Suppose $\xx\in\colspace{\A}$.
Then $\xx=\A\zz$ for some $\zz\in\Rn$,
so
\[
  \xx
  =
  \A\zz
  =
  \A\Adag\A\zz
  =
  \trans{(\A\Adag)} \A\zz
  =
  \B \bigParens{\transAk\A\zz},
\]
where the second and third equalities are
by \Cref{pr:pseudoinv-dfn}(\ref{pr:pseudoinv-dfn:1},\ref{pr:pseudoinv-dfn:3}).
Thus, $\xx\in\colspace\B$,
so
$\colspace\A\subseteq\colspace\B$.
Applied with
$\A$ replaced by $\B$, this further implies
$\colspace\B\subseteq\colspace(\trans{(\Bdag)})=\colspace\A$,
noting for the equality that
$\A=\trans{(\Bdag)}$ by
part~(\ref{pr:pseudoinv-props:dagdag}).

\pfpart{Part~(\ref{pr:pseudoinv-props:c}):}
Suppose $\A$ has full column rank.
Then
\[
  \dim(\colspace\Adag)
  =
  \dim\bigParens{\colspace\regParens{\trans{(\Adag)}}}
  =
  \dim(\colspace\A)
  =
  n,
\]
where the first equality is by
\Cref{pr:row-rank-is-col-rank},
the second is by part~(\ref{pr:pseudoinv-props:a}),
and the third is because the $n$ columns of 
$\A$ are linearly independent.
Thus, $\colspace\Adag$ is a subset of $\Rn$ of dimension $n$, and so
must be equal to $\Rn$.

Let $\xx\in\Rn=\colspace\Adag$.
Then $\xx=\Adag\zz$ for some $\zz\in\Rm$, and
\[
  \Adag\A\xx
  =
  \Adag\A\Adag\zz
  =
  \Adag\zz
  =
  \xx,
\]
where the second equality is by
property~(\ref{pr:pseudoinv-dfn:2})
of \Cref{pr:pseudoinv-dfn}.
Since this holds for all $\xx\in\Rn$, it follows that
$\Adag\A=\Idn{n}$.

\pfpart{Part~(\ref{pr:pseudoinv-props:d}):}
If $\A$ is column-orthogonal,
then $\transAk\A=\Idn{n}$, so $\transA$, if substituted for $\Adag$,
satisfies the four properties of 
\Cref{pr:pseudoinv-dfn}, implying $\Adag=\transA$
(since $\Adag$ is unique).%
\indexg{Moore-Penrose generalized inverse|)}%
\indexg{pseudoinverse|)}
\qedhere
\end{proof-parts}
\end{proof}

\begin{proposition}  \label{pr:lim-vec-proj}
  Let $\vv_1,\ldots,\vv_k\in\Rn$,
  let $\PP$ be the orthogonal projection matrix onto
  $\spnfin{\vv_1,\ldots,\vv_k}$,
  let $\seq{\xx_t}$ be a sequence in $\Rn$ and let $\xx\in\Rn$.
  Suppose $\xx_t\cdot\vv_i\rightarrow\xx\cdot\vv_i$ for
  $i=1,\ldots,k$.
  Then $\PP\xx_t\rightarrow\PP\xx$.
  In particular, if $\vv_1,\ldots,\vv_k$ span $\Rn$, then
  $\xx_t\rightarrow\xx$.
\end{proposition}

\begin{proof}
Let $\VV=[\vv_1,\ldots,\vv_k]$, with pseudoinverse $\VVdag$.
The proposition's assumption means that
$\VVtrans \xx_t \rightarrow \VVtrans \xx$,
implying
that
\[
  \PP \xx_t
  = 
  \VVdagtrans \VVtrans \xx_t
  \rightarrow
  \VVdagtrans \VVtrans \xx
  =
  \PP \xx.
\]
The convergence is by continuity of matrix multiplication.
The equalities are because
$\PP = \VV \VVdag$
by \Cref{pr:pseudoinv-props}(\ref{pr:pseudoinv-props:b}),
implying
$\PP = \trans{(\VV \VVdag)}$
since $\PP$ is symmetric.
\end{proof}

\section{Zero vectors, zero matrices, zero-dimensional Euclidean space}
\label{sec:prelim-zero-dim-space}

We will occasionally find use for vectors and matrices with zero rows
or columns.
We here briefly outline the algebra associated with these.

As noted already, $\zerov{n}$
denotes the all-zeros vector in $\Rn$,
and $\zeromat{m}{n}$ denotes the $m\times n$ all-zeros matrix,
corresponding to the linear map mapping all of $\R^n$ to $\zerov{m}$.

\indexg{Euclidean space!zero-dimensional|(}%
We write $\R^0$ for zero-dimensional Euclidean space, consisting of a single
point, the origin, denoted $\zerov{0}$ (or simply $\zero$ when clear
from context). As the zero of a vector space,
this point satisfies the standard identities $\zerovec+\zerovec=\zerovec$ and $\lambda\zerovec=\zerovec$
for all $\lambda\in\R$. The inner product is also defined as
$\zerovec\cdot\zerovec=0$.
The vector $\zerovec$ is the unique zero-dimensional tuple (or an ``empty'' tuple)
with entries (vacuously) from $\R$;
this tuple corresponds to the unique map from the empty set
to $\R$. 

\indexg{Euclidean space!zero-dimensional|)}%
Further,
for any subset $S\subseteq\R$, we
have $S^0=\{\zerov{0}\}$
(where $S^0$ is the zero-fold Cartesian product of $S$ with itself)
since $\zerovec$ is also the unique empty tuple
with entries (vacuously) from $S$.
Thus, $S^0=\R^0=\{\zerov{0}\}$.

\indexg{matrices!zero rows or columns@with zero rows or columns|(}%
Interpreting matrices in $\Rmn$ as linear maps from $\R^n$ to $\R^m$, for any $m,n\ge 0$,
we find that $\R^{m\times 0}$ contains only a single matrix, denoted
$\zeromat{m}{0}$, corresponding to the only linear map from $\R^0$ to
$\R^m$,
which maps $\zerov{0}$ to $\zerov{m}$.
Similarly, $\R^{0\times n}$ contains only one matrix, denoted $\zeromat{0}{n}$,
corresponding to the linear map from $\R^n$ to $\R^0$
mapping all of $\R^n$ to $\zerov{0}$.
The matrix $\zeromat{0}{0}$ is both a zero matrix as well as an
identity matrix, because it is the identity map on $\R^0$;
thus, $\Idn{0}=\zeromat{0}{0}$.
Vacuously, $\zeromat{0}{0}$ is positive upper triangular,
and $\zeromat{m}{0}$ is column-orthogonal for $m\geq 0$.

Since $\R^{m\times n}$, for any $m,n\ge 0$, is itself a vector space
with $\zeromat{m}{n}$ as its zero, we obtain standard identities
$\zeromat{m}{n}+\zeromat{m}{n}=\zeromat{m}{n}$ and $\lambda\zeromat{m}{n}=\zeromat{m}{n}$ for all $\lambda\in\R$.
Interpreting matrix product as composition of linear maps, we further
obtain the following identities for $m,n,k\geq 0$,
and $\A\in \R^{m\times k}$ and $\B\in\R^{k\times n}$:
\[
  \A \zeromat{k}{n} = \zeromat{m}{k} \B = \zeromat{m}{n}.
\]

Finally, just as vectors in $\Rn$ are interpreted as $n\times 1$
matrices, we identify $\zerovec$ with $\zeromat{0}{1}$ so that, for example, we can write
$\zerov{n}=\zeromat{n}{1}=\zeromat{n}{0}\zeromat{0}{1}=\zeromat{n}{0}\zerov{0}$.%
\indexg{matrices!zero rows or columns@with zero rows or columns|)}

\chapter{Review of topology}
\label{sec:prelim:topology}

Convex analysis in $\Rn$ uses the standard Euclidean norm to define topological
concepts such as convergence, closure, and continuity.
We will see that astral space is not a normed space, and
is not even metrizable, so instead of building the topology of astral
space from a metric,
we will need to use the framework of general topology,
which we briefly review in this \namecref{sec:prelim:topology}.
Depending on background, this \namecref{sec:prelim:topology} can be skimmed or skipped,
and revisited only as needed.
A more complete introduction to topology can be found, for instance,
in
\idxhitch\citet{hitchhiker_guide_analysis}
or
\idxmunk\citet{munkres}.

\section{Topological space, open sets, base, subbase}

\begin{definition}
\label{def:open}
\indexg{topology!defined|(}%
A \emph{topology} $\topo$ on a set $X$ is a collection of subsets of $X$, called
\emph{open sets},%
\indexg{open sets}
that satisfy the following conditions:
\begin{letter-compact}
\item\label{i:open:a} $\emptyset\in\topo$ and $X\in\topo$.
\item\label{i:open:b} $\topo$ is closed under finite intersections, that is,
$U\cap V\in\topo$ for all $U,V\in\topo$.
\item\label{i:open:c} $\topo$ is closed under arbitrary unions; that is,
  if $U_\alpha\in\topo$ for all $\alpha\in\indset$
  (where $\indset$ is any index set), then
  $\bigcup_{\alpha\in\indset} U_\alpha\in\topo$.
\end{letter-compact}
The set~$X$ with topology $\topo$ is called a \emph{topological space}.%
\indexg{topological spaces}%
\indexg{topology!defined|)}
\end{definition}

\begin{example}[Topology on $\R$]
\indexg{topology!reals@on reals|(}%
The standard topology on $\R$ consists of all sets $U$ such that for every $x\in U$ there exists $\epsilon>0$ such that $(x-\epsilon,x+\epsilon)\subseteq U$.%
\indexg{topology!reals@on reals|)}%
\end{example}

\begin{example}[Euclidean topology on $\Rn$]
\label{ex:topo-on-rn}
\indexg{Euclidean topology|(}%
\indexg{topology!Euclidean|(}%
\indexg{Euclidean space!topology|(}%
The Euclidean topology on $\Rn$ consists
of all sets $U$ such that for every $\xx\in U$, there exists $\epsilon>0$ such that
$\ball(\xx,\epsilon)\subseteq U$.
Except where explicitly stated otherwise, we will always assume this
topology on $\Rn$.%
\indexg{Euclidean space!topology|)}%
\indexg{topology!Euclidean|)}%
\indexg{Euclidean topology|)}
\end{example}

\indexg{topology!generated by base or subbase|(}%
Instead of specifying the topology $\topo$ directly, it is usually more convenient
to specify a suitable subfamily of $\topo$ and use it to generate $\topo$ by means of unions and intersections.
As such,
a subfamily $\calB$ of $\topo$ is called a \emph{base}%
\indexg{base}
for $\topo$ if every
element of $\topo$ can be written as a union of elements of $\calB$.
A subfamily $\calS$ of $\topo$ is called
a \emph{subbase}%
\indexg{subbase}
for~$\topo$ if
the collection of all finite intersections of elements of $\calS$ forms a base for $\topo$.
(Some authors, including
\citealp{munkres}, instead use the terms \emph{basis}
and \emph{subbasis}.)
Note that
any collection $\calS$ of subsets of $X$
whose union is all of $X$
is a subbase for a topology on $X$ obtained by taking all unions of all finite
intersections of elements of $\calS$. In this case, we say that
$\calS$ \emph{generates}
the topology $\topo$.

\begin{proposition}   \label{pr:base-equiv-topo}
  Let $X$ be a topological space with topology $\topo$,
  let $\calB$ be a base for~$\topo$,
  and let $U\subseteq X$.
  Then $U$ is open in $\topo$ if and only if
  for all $x\in U$, there exists a set $B\in\calB$
  such that $x\in B\subseteq U$.
\end{proposition}

\begin{proof}
  See \idxmunk\citet[Lemma~13.1]{munkres}.
\end{proof}

\begin{example}[Base and subbase for topology on $\R$]
\indexg{topology!reals@on reals|(}%
A base for the standard topology on $\R$ consists of all open intervals
$(a,b)$, where $a,b\in\R$. A subbase consists of sets of the form
$(-\infty,b)$ and $(b,+\infty)$ for $b\in\R$.%
\indexg{topology!reals@on reals|)}
\end{example}

\begin{example}[Base and subbase for topology on $\eR$]
\label{ex:topo-rext}
\indexg{topology!extended reals@on extended reals|(}%
On $\eR$, we will always assume the topology with a
base consisting of all intervals
$(a,b)$, $[-\infty,b)$, and $(b,+\infty]$,
for $a,b\in\R$. This topology is generated by a subbase consisting of sets $[-\infty,b)$ and $(b,+\infty]$
for all $b\in\R$.%
\indexg{topology!extended reals@on extended reals|)}
\end{example}

\begin{example}[Base for the Euclidean topology on $\Rn$]
\indexg{topology!Euclidean|(}%
\indexg{Euclidean space!topology|(}%
\indexg{Euclidean topology|(}%
A base for the Euclidean topology on $\Rn$ consists of open balls
$\ball(\xx,\epsilon)$
for all $\xx\in\Rn$ and $\epsilon>0$.%
\indexg{topology!Euclidean|)}%
\indexg{Euclidean space!topology|)}%
\indexg{Euclidean topology|)}%
\indexg{topology!generated by base or subbase|)}
\end{example}%

For a metric $d:X\times X\rightarrow\R$,
the \emph{metric topology}%
\indexg{topology!metric}\indexg{metric topology}
on $X$ induced by $d$
is that topology whose base consists of all balls
$\{ y\in X: d(x,y)<\epsilon \}$,
for $x\in X$ and $\epsilon>0$.
Thus, Euclidean topology on $\Rn$ is the same as the metric topology
on $\Rn$ with metric given by the Euclidean distance,
$d(\xx,\yy)=\norm{\xx-\yy}$.

A topological space $X$ is
\emph{metrizable}%
\indexg{topology!metrizable}\indexg{metrizability}
if there exists a metric that induces its topology.

\section{Closed sets, neighborhoods, dense sets, separation properties}
\label{sec:prelim:topo:closed-sets}

\indexg{closed sets|(}%
Complements of open sets are called \emph{closed sets}. Their basic properties
are symmetric to the properties of open sets from \Cref{def:open}:
\begin{letter-compact}
\item $\emptyset$ and $X$ are closed in any topology on $X$.
\item Finite unions of closed sets are closed.
\item Arbitrary intersections of closed sets are closed.%
\indexg{closed sets|)}
\end{letter-compact}%

The \emph{closure}%
\indexg{closure (topological)}
of a set $A\subseteq X$, denoted $\Abar$,%
\indexm{s 100}{$\Sbar$}{closure (general topological)}
 is the intersection of all closed sets containing~$A$.
The
\indexg{interior}%
\emph{interior of $A$},
denoted $\intr A$,%
\indexm{int}{$\intr S$}{interior}
is the union of all
open sets contained in $A$. A~\emph{neighborhood}%
\indexg{neighborhood}
of a point $x\in X$ is an open set that contains $x$.
This definition of neighborhood follows
\idxmunk\citet{munkres} and others, but
note that some authors,
including \idxhitch\citet{hitchhiker_guide_analysis},
instead define a neighborhood of $x$ as any
set that contains an open set containing $x$, that is, any set that is a superset of our notion of a neighborhood.
They would refer to our neighborhoods as open neighborhoods.

\begin{proposition}   \label{pr:closure:intersect}
  Let $X$ be a topological space, let $A\subseteq X$, and let $x\in X$.
  \begin{letter-compact}
  \item   \label{pr:closure:intersect:s1}
    $A=\Abar$ if and only if $A$ is closed.
  \item   \label{pr:closure:intersect:s2}
    For every closed set $C$ in $X$,
    if $A\subseteq C$ then $\Abar\subseteq C$.
  \item   \label{pr:closure:intersect:a}
    $x\in\Abar$ if and only if $A\cap U\neq\emptyset$ for every
    neighborhood $U$ of $x$.
  \item   \label{pr:closure:intersect:b}
    $x\in\intr{A}$ if and only if $U\subseteq A$ for some neighborhood
    $U$ of $x$.
  \item   \label{pr:closure:intersect:comp}
    $\intr{A}=X\setminus(\clbar{X\setminus A})$.
  \end{letter-compact}
\end{proposition}

\begin{proof}
  ~

\begin{proof-parts}
\pfpart{Parts~(\ref{pr:closure:intersect:s1})
  and~(\ref{pr:closure:intersect:s2}):}
These are straightforward from the definition of closure.

\pfpart{Part~(\ref{pr:closure:intersect:a}):}
See \idxmunk\citet[Theorem 17.5]{munkres}.

\pfpart{Part~(\ref{pr:closure:intersect:b}):}
Let $x\in X$.
Suppose $x\in U$ where $U$ is open and $U\subseteq A$.
Then $U$ is also in $\intr A$, by its definition, so $x\in\intr A$.
Conversely, if $x\in\intr A$ then
$x$ is in an open set that is included in $A$, namely,
$\intr A$.

\pfpart{Part~(\ref{pr:closure:intersect:comp}):}
See \idxhitch\citet[Lemma~2.4]{hitchhiker_guide_analysis}.
\qedhere
\end{proof-parts}
\end{proof}

A subset $Z$ of a topological space $X$ is said to be \emph{dense}%
\indexg{dense (set)}
in $X$ if $\Zbar=X$.
By \Cref{pr:closure:intersect}(\ref{pr:closure:intersect:a}),
this means that $Z$ is dense in $X$ if and only if
every nonempty open set in $X$ has a nonempty intersection with $Z$.

Two nonempty subsets $A$ and $B$ of a topological space are said to be \emph{separated by open sets}%
\indexg{separation by open sets}
if they are included in disjoint open sets,
that is, if there exist disjoint open sets $U$ and $V$
such that $A\subseteq U$ and $B\subseteq V$.
A topological space is called:
\begin{item-compact}
\item \emph{Hausdorff}\indexg{Hausdorff (space)}
    if any two distinct points are separated by open sets;
\item \emph{regular}\indexg{regular (space)}
    if any nonempty closed set and any point disjoint from it are separated by open sets;
\item \emph{normal}\indexg{normal (space)}
   if any two disjoint nonempty closed sets are separated by open sets.
\end{item-compact}

\begin{proposition}
\label{prop:hausdorff:normal}
In a Hausdorff space, all singleton sets are closed, so every normal Hausdorff space is regular.
\end{proposition}

\begin{proof}
That singletons in a Hausdorff space are closed follows from
\idxmunk\citet[Theorem~17.8]{munkres}.
That a normal Hausdorff space is regular then follows from the definitions
above.
\end{proof}

\section{Subspace topology}
\label{sec:subspace}

\indexg{subspace (topological)|(}%
Let $X$ be a space with topology $\topo$, and let
$Y\subseteq X$.
Then
the collection $\topo_Y=\set{U\cap Y:\:U\in\topo}$ defines a topology on $Y$,
called the
\indexg{topology!subspace}%
\emph{subspace topology}.
 The set $Y$ with topology $\topo_Y$
is called a \emph{(topological) subspace} of $X$.
Furthermore, for any base $\calB$
for $\topo$, the collection $\set{U\cap Y:\:U\in\calB}$ is a base for $\topo_Y$.
Likewise, for any subbase $\calS$ for $\topo$, the collection $\set{V\cap Y:\:V\in\calS}$ is a subbase for $\topo_Y$.

\begin{example}[Topology on \mbox{$[0,1]$}]
\label{ex:topo-zero-one}
The topology on $[0,1]$, as a subspace of $\R$, has as subbase
all sets $[0,b)$ and $(b,1]$ for $b\in [0,1]$.
\end{example}

Note that, when working with subspaces, it can be important to specify
which topology we are referencing since, for instance, a set might be
closed in subspace $Y$ (meaning in the subspace topology $\topo_Y$),
but not closed in $X$ (meaning in the original topology $\topo$).

\begin{proposition}
\label{prop:subspace}
Let $X$ be a topological space, $Y$ its subspace, and $Z\subseteq Y$.
\begin{letter-compact}
\item
\label{i:subspace:closed}
  $Z$ is closed in $Y$ if and only if there exists a set $C$ closed in $X$ such that $Z=C\cap Y$.
\item
\label{i:subspace:closure}
  The closure of $Z$ in $Y$ equals $\Zbar\cap Y$, where $\Zbar$ is the closure of $Z$ in $X$.
\item
\label{i:subspace:Hausdorff}
  If $X$ is Hausdorff, then so is $Y$.
\item
\label{i:subspace:dense}
  If $Z$ is dense in $Y$ and $Y$ is dense in $X$, then $Z$ is dense in $X$.
\item
\label{i:subspace:subspace}
  The topology on $Z$ as a subspace of $Y$ coincides with the topology on $Z$ as a subspace of $X$.
\end{letter-compact}
\end{proposition}
\begin{proof}
Parts~(\ref{i:subspace:closed}, \ref{i:subspace:closure}, \ref{i:subspace:Hausdorff}) follow by \idxmunk\citet[Theorems 17.2, 17.4, 17.11]{munkres}.
To prove part~(\ref{i:subspace:dense}), note that for
every nonempty open set $U\subseteq X$, the intersection $U\cap Y$ is
nonempty and open in $Y$, and therefore, the intersection
$(U\cap Y)\cap Z=U\cap Z$ is
nonempty, meaning that $Z$ is dense in $X$.
Finally, part~(\ref{i:subspace:subspace}) follows from the definition
of subspace topology.%
\indexg{subspace (topological)|)}
\end{proof}

\section{Continuity, homeomorphism, compactness}

\indexg{continuity|(}%
A function $f:X\to Y$ between topological spaces is said to be
\emph{continuous}
 if $f^{-1}(V)$ is open in $X$ for every $V$ that is open in~$Y$. 
We say that $f$ is \emph{continuous at a point $x\in X$} if for every
neighborhood $V$ of $f(x)$, there exists a neighborhood $U$ of $x$
such that $f(U)\subseteq V$.

\begin{proposition}
\label{prop:cont}
For a function $f:X\to Y$ between topological spaces the following are equivalent:
\begin{letter-compact}
\item  \label{prop:cont:a}
   $f$ is continuous.
\item  \label{prop:cont:b}
   $f$ is continuous at every point in $X$.
\item  \label{prop:cont:inv:closed}
   $f^{-1}(C)$ is closed in $X$ for every $C$ that is a closed subset of $Y$.
\item  \label{prop:cont:sub}
   $f^{-1}(V)$ is open in $X$ for every element $V$ of some subbase for the topology on $Y$.
\item  \label{prop:cont:c}
   $f(\Abar)\subseteq\overlineKernIt{f(A)}$ for every $A\subseteq X$.
\end{letter-compact}
\end{proposition}

\begin{proof}
  See \idxhitch\citet[Theorem 2.27]{hitchhiker_guide_analysis},
  or Theorem~18.1 and the preceding discussion of
  \idxmunk\citet{munkres}.
\end{proof}

\begin{proposition}  \label{pr:sup-f-Abar}
  Let $X$ be a topological space,
  let $A\subseteq X$,
  and let $f:X\rightarrow\Rext$ be continuous.
  Then $\sup f(A)=\sup f(\Abar)$;
  that is,
  $
     \sup_{x\in A} f(x)
     =
     \sup_{x\in \Abar} f(x)
  $.
\end{proposition}

\begin{proof}
First, $\sup f(A)\leq\sup f(\Abar)$ since $A\subseteq\Abar$.
For the reverse inequality, let $\alpha=\sup f(A)$,
and let $U=[-\infty,\alpha]$.
Then
$f(\Abar)\subseteq\overlineKernIt{f(A)}\subseteq U$,
where the first inclusion is by
\Cref{prop:cont}(\ref{prop:cont:a},\ref{prop:cont:c})
since $f$ is continuous, and the second is
because $U$, which is closed, includes $f(A)$, and therefore its
closure as well.
Thus, $\sup f(\Abar)\leq \alpha$, completing the proof.%
\indexg{continuity|)}
\end{proof}

\indexg{homeomorphism|(}%
Let $X$ and $Y$ be topological spaces.
A function $f:X\rightarrow Y$ is a \emph{homeomorphism} if
$f$ is a bijection, and if both $f$ and $f^{-1}$ are continuous.
When such a function exists, we say that
$X$ and $Y$ are \emph{homeomorphic}.
From the topological perspective, spaces $X$ and $Y$ are identical;
only the names of the points have been changed. Any topological property of $X$, meaning a property derived from the topology alone, is also a property of $Y$ (and vice versa).

\begin{example}[Homeomorphism of $\Rext$ with \mbox{$[0,1]$}]
\label{ex:homeo-rext-finint}
We assume topologies as in
Examples~\ref{ex:topo-rext}
and~\ref{ex:topo-zero-one}.
Let $f:\Rext \rightarrow [0,1]$ be defined,
for $\barx\in\Rext$, by
\[
   f(\barx)
   =
   \begin{cases}
     0               & \text{if $\barx=-\infty$}, \\
     e^{\barx}/(1+e^{\barx}) & \text{if $\barx\in\R$}, \\
     1               & \text{if $\barx=+\infty$.}
   \end{cases}
\]
This function is bijective and continuous with a continuous inverse.
Therefore, it is a homeomorphism, so $\Rext$ and $[0,1]$ are
homeomorphic.

This means, for instance, that $\Rext$ is metrizable,
since $[0,1]$ is metrizable (using ordinary distance
between points).%
\indexg{homeomorphism|)}
\end{example}

Let $X$ be a topological space.
An \emph{open cover}%
\indexg{open cover}
 of a set $K\subseteq X$
is any family of open sets whose union includes $K$, that is,
any family $\set{U_\alpha}_{\alpha\in\indset}$ of open sets such that
$K\subseteq\bigcup_{\alpha\in\indset} U_\alpha$.
\indexg{compactness|(}%
We say that a set $K$ is \emph{compact}
if every open cover of $K$ has a
finite open subcover, that is, if
for every open cover
$\set{U_\alpha}_{\alpha\in\indset}$ of $K$,
there exist $\alpha_1,\dotsc,\alpha_m\in\indset$ such that
$K\subseteq U_{\alpha_1}\cup\dotsb\cup U_{\alpha_m}$.
If $X$ is itself compact, then we say that $X$ is
a \emph{compact space}.

A set $S\subseteq\Rn$ is \emph{bounded} if
$S\subseteq\ball(\zero,R)$ for some $R\in\R$.

\begin{proposition}  \label{pr:compact-in-rn}
A set in $\Rn$ is compact
if and only if it is closed and bounded.
\end{proposition}

\begin{proof}
  See \idxmunk\citet[Theorem~27.3]{munkres}.
\end{proof}

In particular, this proposition implies that $\Rn$ itself is not compact
(for $n\geq 1$).
Although $\R$ is not compact, its extension $\Rext$ is compact:

\begin{example}[Compactness of $\Rext$]
\label{ex:rext-compact}
\indexg{compactness!extended reals@of extended reals|(}%
From Example~\ref{ex:homeo-rext-finint},
$\Rext$ is homeomorphic with $[0,1]$, which is compact, being closed
and bounded.
Therefore, $\Rext$ is also compact.%
\indexg{compactness!extended reals@of extended reals|)}%
\end{example}

Here are some implications of compactness:

\begin{proposition}~
\label{prop:compact}
\begin{letter-compact}
\item
\label{prop:compact:cont-compact}
The image of a compact set under a continuous function is compact.
\item
\label{prop:compact:closed-subset}
Every closed subset of a compact set is compact.
\item
\label{prop:compact:closed}
Every compact subset of a Hausdorff space is closed.
\item
\label{prop:compact:subset-of-Hausdorff}
Every compact Hausdorff space is normal, and therefore regular.
\end{letter-compact}
\end{proposition}
\begin{proof}
See
\idxmunk\citet[Theorems 26.5, 26.2, 26.3, 32.3]{munkres}. Regularity
in part (\ref{prop:compact:subset-of-Hausdorff}) follows by \Cref{prop:hausdorff:normal}.
\end{proof}

\begin{proposition}   %
\label{pr:cont-compact-attains-max}
\indexg{compactness!extremum attained}%
Let $f:X\to\eR$ be continuous and $X$ compact. Then $f$ attains both its minimum and maximum on $X$.
\end{proposition}

\begin{proof}
  This is a special case of \idxmunk\citet[Theorem~27.4]{munkres}.
\end{proof}

\begin{proposition}
\label{pr:cont-from-compact}
Let $f:X\to Y$ be continuous, $X$ compact, $Y$ Hausdorff,
and let $S\subseteq X$.
Then:
\begin{letter-compact}
  \item   \label{pr:cont-from-compact:a}
    If $S$ is closed in $X$, then $f(S)$ is closed in $Y$.
  \item   \label{pr:cont-from-compact:b}
    $f(\Sbar)=\overline{f(S)}$.
\end{letter-compact}
\end{proposition}

\begin{proof}
  ~
\begin{proof-parts}
\pfpart{Part~(\ref{pr:cont-from-compact:a}):}
Suppose $S$ is closed in $X$,
implying it is also compact,
being a closed subset of a compact space
(\Cref{prop:compact}\ref{prop:compact:closed-subset}).
Since $f$ is continuous, $f(S)$ is compact as well
(\Cref{prop:compact}\ref{prop:compact:cont-compact}),
and therefore closed
(\Cref{prop:compact}\ref{prop:compact:closed}).

\pfpart{Part~(\ref{pr:cont-from-compact:b}):}
Since $f$ is continuous,
$f(\Sbar)\subseteq\clbar{f(S)}$ by
\Cref{prop:cont}(\ref{prop:cont:a},\ref{prop:cont:c}).
Also, $f(S)\subseteq f(\Sbar)$ (since $S\subseteq\Sbar$),
implying $\clbar{f(S)}\subseteq f(\Sbar)$
since $f(\Sbar)$ is closed by part~(\ref{pr:cont-from-compact:a}).
\qedhere
\end{proof-parts}
\end{proof}

For topological spaces $X$ and $Y$, we say that $Y$ is a
\emph{compactification}%
\indexg{compactification}
of $X$ if $Y$ is compact,
$X$ is a subspace of $Y$,
and $X$ is dense in $Y$.
(In defining compactification, others, including \idxmunk\citealp{munkres},
also require that $Y$ be Hausdorff.)%
\indexg{compactness|)}%

\section{Sequences, first-countability, second-countability}
\label{sec:prelim:countability}

\indexg{sequence convergence|(}%
Let $X$ be a topological space.
We say that
a sequence $\seq{x_t}$ in $X$ \emph{converges}
to a point $x$ in $X$
if for every neighborhood $U$ of $x$,
$x_t\in U$ for all $t$ sufficiently large
(that is, there exists an index $t_0$ such that
$x_t\in U$ for $t\ge t_0$).
When this occurs,
we write $x_t\to x$ and refer to $x$ as a \emph{limit} of $\seq{x_t}$.
We say that a sequence
$\seq{x_t}$ \emph{converges in $X$} if it has a limit $x\in X$. A
sequence in a Hausdorff space can have at most one limit. In that
case, we write $\lim x_t=x$ when $x_t\to x$.%
\indexg{sequence convergence|)}

A \emph{subsequence}%
\indexg{subsequences}
 of a sequence $\seq{x_t}$ is any sequence
$\seq{x_{s(t)}}$ where $s:\nats\rightarrow\nats$, and
$s(1)<s(2)<\dotsb$.
We say that $x$ is a \emph{subsequential limit}%
\indexg{subsequential limit} 
of a sequence
if $x$ is the limit of one of its subsequences.

\indexg{countability|(}%
A set $S$ is \emph{countably infinite} if there exists a bijection
between $S$ and $\nats$.
It is \emph{countable} if it is either finite or countably infinite.

\begin{proposition}  \label{pr:count-equiv}
  Let $S$ be a nonempty set.
  Then the following are equivalent:
  \begin{letter-compact}
  \item    \label{pr:count-equiv:a}
    $S$ is countable.
  \item    \label{pr:count-equiv:b}
    There exists a surjective function from $\nats$ to $S$.
  \item    \label{pr:count-equiv:c}
    There exists an injective function from $S$ to $\nats$.
  \end{letter-compact}
\end{proposition}

\begin{proof}
  See \idxmunk\citet[Theorem~7.1]{munkres}.
\end{proof}

\begin{proposition}  \label{pr:uncount-interval}
Let $a,b\in\R$ with $a<b$.
  Then the intervals $[a,b]$, $[a,b)$, $(a,b]$, and $(a,b)$ are all
  uncountable.
\end{proposition}

\begin{proof}
See \idxhitch\citet[Section 1.9]{hitchhiker_guide_analysis}.%
\indexg{countability|)}
\end{proof}

\indexg{countable neighborhood base|(}%
\indexg{first-countability|(}%
A \emph{neighborhood base}%
\indexg{neighborhood base}
at a point $x$ is a collection $\calB_x$ of
neighborhoods of $x$ such that each neighborhood of $x$ contains at
least one of the elements of $\calB_x$. A topological space $X$ is said to
be
\emph{first-countable}%
\indexg{first-countability!defined}
if there exists a countable neighborhood base at every point $x\in X$.

\begin{example}[First-countability of $\Rn$]
\label{ex:rn-first-count}
\indexg{first-countability!Euclidean space@of Euclidean space|(}%
\indexg{Euclidean space!first-countability of|(}%
Let $\xx\in\Rn$, and for each $t$, let
$B_t=\ball(\xx,1/t)$.
Then $\countset{B}$ is a countable neighborhood base for $\xx$.
This is because if $U$ is a neighborhood of $\xx$, then there exist
$\yy\in\Rn$ and $\epsilon>0$ such that
$\xx\in \ball(\yy,\epsilon)\subseteq U$.
Letting $t$ be so large that
$\norm{\xx-\yy}+1/t<\epsilon$,
it then follows that
$B_t\subseteq  \ball(\yy,\epsilon)\subseteq U$.
Thus, in the Euclidean topology, $\Rn$ is first-countable.%
\indexg{first-countability!Euclidean space@of Euclidean space|)}%
\indexg{Euclidean space!first-countability of|)}%
\end{example}

\indexg{countable neighborhood base!nested|(}%
A countable neighborhood base $\countset{B}$ at a point $x$
is \emph{nested}
if $B_1\supseteq B_2\supseteq\dotsb$.
For instance, the countable neighborhood base in
Example~\ref{ex:rn-first-count} is nested.
Any countable neighborhood base $\countset{B}$ can be turned into
a nested one, $\countset{B'}$, by setting $B'_t=B_1\cap\dotsb\cap B_t$
for all $t$.%
\indexg{countable neighborhood base!nested|)}

If $\countset{B}$ is a nested countable neighborhood base at $x$,
then every sequence $\seq{x_t}$ with each point $x_t$ in $B_t$
must converge to $x$:

\begin{proposition}
\label{prop:nested:limit}
Let $x$ be a point in a topological space $X$. Let $\countset{B}$ be a nested countable neighborhood base at $x$, and let $x_t\in B_t$
for all $t$. Then $x_t\to x$.
\end{proposition}
\begin{proof}
Let $U$ be a neighborhood of $x$. Then there exists $t_0$ such that $B_{t_0}\subseteq U$,
implying, since the base sets are nested,
that $x_t\in B_t\subseteq B_{t_0}\subseteq U$
for $t\ge t_0$. Thus, $x_t\to x$.%
\indexg{countable neighborhood base|)}
\end{proof}

First-countability allows us to work with sequences as we would in
$\Rn$, for instance, for characterizing the closure of a set or the
continuity of a function:

\begin{proposition}
\label{prop:first:properties}
Let $X$ and $Y$ be topological spaces, and let $x\in X$.
\begin{letter-compact}
\item
\label{prop:first:closure}
Let $A\subseteq X$.
If there exists a sequence $\seq{x_t}$ in $A$ with $x_t\rightarrow x$
then $x\in\Abar$.
If $X$ is first-countable, then the converse holds as well.
\item
\label{prop:first:cont}
Let $f:X\to Y$.
If $f$ is continuous at $x$ then
$f(x_t)\rightarrow f(x)$
for every sequence $\seq{x_t}$ in $X$ that converges to $x$.
If $X$ is first-countable, then the converse holds as well.
\end{letter-compact}
\end{proposition}%

\begin{proof}
  ~

\begin{proof-parts}
\pfpart{Part~(\ref{prop:first:closure}):}
See \citet[Theorem~30.1a]{munkres}.

\pfpart{Part~(\ref{prop:first:cont}):}
Assume first that $f$ is continuous at $x$,
and let $\seq{x_t}$ be any sequence in $X$ that converges to $x$.
Then for any neighborhood $V$ of $f(x)$, there exists a neighborhood $U$ of $x$ such that $f(U)\subseteq V$. Since $x_t\to x$, there exists $t_0$ such that $x_t\in U$ for all $t\ge t_0$, and so also $f(x_t)\in V$ for all $t\ge t_0$. Thus, $f(x_t)\to f(x)$.

For the converse, assume $X$ is first-countable, and
that $x_t\to x$ implies $f(x_t)\to f(x)$.
Let $\countset{B}$ be a nested countable neighborhood base at $x$, which exists by first-countability. For contradiction,
suppose $f$ is not continuous at $x$,
and hence that there exists a neighborhood $V$ of $f(x)$ such that $f^{-1}(V)$ does not contain any
neighborhood of $x$, and so in particular does not contain any $B_t$.
Then for all $t$, there exists $x_t\in B_t\setminus f^{-1}(V)$.
Thus, $f(x_t)\not\in V$ for all $t$, but also $x_t\to x$, by
\Cref{prop:nested:limit},
implying $f(x_t)\to f(x)$, by assumption.
The set $V^c=Y\setminus V$ is closed, and $f(x_t)\in V^c$ for all $t$, so by
part~(\ref{prop:first:closure}), we also have that $f(x)\in V^c$.
This contradicts
that $V$ is a neighborhood of $f(x)$.
\qedhere
\end{proof-parts}
\end{proof}

\indexg{countable neighborhood base|(}%
The next proposition proves a kind of converse
to \Cref{prop:nested:limit}.
In particular, it shows that if, in some topological space $X$,
a family of neighborhoods
$\countset{B}$ of a point~$x$ has the property given in
\Cref{prop:nested:limit}
(that is, every sequence $\seq{x_t}$ with $x_t\in B_t$ converges to $x$),
then those neighborhoods must constitute a neighborhood base at $x$.
The proposition shows more generally that this holds even if
we only assume the property from
\Cref{prop:nested:limit}
to be true for sequences in a dense subset of the space $X$
(rather than all of $X$).
Note that the proposition does not assume that $X$ is first-countable.\looseness=-1

\begin{proposition}
\label{thm:first:fromseq}
  Let $X$ be a regular topological space with a dense subset $Z$.
  Let $x\in X$, and let $\countset{B}$ be a countable
  family of neighborhoods of $x$.
  Suppose that for every sequence $\seq{z_t}$ with each $z_t\in Z\cap B_t$,
  we have $z_t\rightarrow x$.
  Then the family $\countset{B}$ is a countable neighborhood base at $x$.
\end{proposition}

\begin{proof}
Contrary to the claim, suppose
$\countset{B}$ is not a countable neighborhood base at~$x$,
meaning there exists a
neighborhood $N$ of $x$ that does not contain any $B_t$.
Then for each $t$, there is a point
$x_t\in B_t\cap \Ncomp$
where $\Ncomp=X \setminus N$ is the complement of $N$.

Since $\Ncomp$ is closed (being the complement of an open set) and
$x\not\in\Ncomp$, the fact that the space $X$ is regular implies that there must
exist disjoint open sets $U$ and $V$ such that $\Ncomp\subseteq U$ and
$x\in V$.

For each $t$, $x_t\in B_t\cap\Ncomp \subseteq B_t\cap U$,
meaning $B_t\cap U$ is a neighborhood of~$x_t$.
Therefore, since $Z$ is dense in $X$,
there exists a point
$z_t\in B_t\cap U\cap Z$.

By assumption, the resulting sequence $\seq{z_t}$ converges to $x$.
Since $V$ is a neighborhood of $x$, this means that all but
finitely many of the points $z_t$ are in $V$.
But this is a contradiction since every $z_t\in U$, and $U$ and $V$
are disjoint.%
\indexg{countable neighborhood base|)}
\end{proof}

\indexg{sequential compactness|(}%
A topological space $X$ is \emph{sequentially compact}
if every sequence in $X$ has a subsequence converging to an element of $X$.%
\indexg{sequential compactness|)}%

\begin{proposition}
\label{prop:first:subsets}
~
\begin{letter-compact}
\item
\label{i:first:subspace}
  Every subspace of a first-countable space is first-countable.
\item
\label{i:first:compact}
  Every first-countable, compact space is sequentially compact.
\end{letter-compact}
\end{proposition}

\begin{proof}
~

\begin{proof-parts}
\pfpart{Part~(\ref{i:first:subspace}):}
See \idxmunk\citet[Theorem~30.2]{munkres}.

\pfpart{Part~(\ref{i:first:compact}):}
See \idxmanetti\citet[Lemma~6.21]{manetti-topology}.%
\indexg{first-countability|)}
\qedhere
\end{proof-parts}
\end{proof}

\indexg{second-countability|(}%
A topological space is said to be \emph{second-countable} if it has a
countable base.

\begin{proposition}
\label{prop:sep:metrizable}
\indexg{metrizability!second-countability and}%
Let $X$ be a topological space.
If $X$ is metrizable and if $X$ includes a countable dense subset,
then $X$ is second-countable.
\end{proposition}

\begin{proof}
  See \idxhitch\citet[Theorem~3.40]{hitchhiker_guide_analysis},
  noting that if $X$ is metrizable, then it is Hausdorff.%
\indexg{second-countability|)}
\end{proof}

\begin{example}[Countable dense subset and countable base for $\Rn$]
\label{ex:2nd-count-rn}
\indexg{second-countability!Euclidean space@of Euclidean space|(}%
\indexg{Euclidean space!second-countability of|(}%
The set $\rats$ is countable, as is its $n$-fold Cartesian product
$\rats^n$
\idxmunk\citep[Example~7.3 and Theorem~7.6]{munkres}.
Further, $\rats^n$ is dense in $\Rn$ (as can be seen, for instance, by
constructing, for each $\xx\in\Rn$, a sequence in $\rats^n$ that
converges to $\xx$).
Since $\Rn$ is a metric space, \Cref{prop:sep:metrizable} then implies
that it is second-countable, and so has a countable base.
Explicitly, such a base is given by all sets
$\ball(\yy,1/k)$ for all $\yy\in\rats^n$ and all $k\in\nats$.%
\indexg{second-countability!Euclidean space@of Euclidean space|)}%
\indexg{Euclidean space!second-countability of|)}%
\end{example}

\section{Product topology, Tychonoff's theorem, weak topology}
\label{sec:prod-top}

The \emph{Cartesian product}%
\indexg{Cartesian product}
of an indexed family of sets $\{X_\alpha\}_{\alpha\in\indset}$,
where $\indset$ is some index set,
is the set
\begin{equation}   \label{eqn:cart-prod-notation}
  \prod_{\alpha\in\indset} X_\alpha
  =
  \bigBraces{\tupset{x_\alpha}{\alpha\in I} :\:
    x_\alpha \in X_\alpha \mbox{ for all } \alpha\in I
  }.
\end{equation}
\indexg{tuples|(}%
Here,
\indexm{x alpha alpha i}{$\tupset{x_\alpha}{\alpha\in\indset}$}{tuple (over general index set)}%
$x=\tupset{x_\alpha}{\alpha\in\indset}$
is a \emph{tuple},
which, as it appears in this equation, is formally the function
$x:\indset\rightarrow \bigcup_{\alpha\in\indset} X_\alpha$
with $x(\alpha)=x_\alpha$ for
\indexg{tuples|)}%
$\alpha\in\indset$.
The product
$\prod_{\alpha\in\indset} X_\alpha$ thus consists of all such tuples
with $x_\alpha\in X_\alpha$ for $\alpha\in\indset$.

Since tuples are functions, we sometimes write the components of a
tuple $x$ as $x(\alpha)$ rather than $x_\alpha$, depending on context.
When $X_\alpha=X$ for all $\alpha\in\indset$,
the Cartesian product in
\eqref{eqn:cart-prod-notation}
is written $X^\indset$, and every tuple in this set is simply
a function mapping $\indset$ into $X$.
Thus, in general, $X^\indset$ denotes the set of all functions from
$\indset$ to $X$.

When $\indset=\{1,2\}$, the Cartesian product in
\eqref{eqn:cart-prod-notation}
is simply
$X_1\times X_2$, with elements given by all pairs
$\rpair{x_1}{x_2}$ with $x_1\in X_1$ and $x_2\in X_2$.

\indexg{topology!product|(}%
\indexg{product topology|(}%
Suppose each $X_\alpha$ is a topological space with topology
$\topo_\alpha$, for all $\alpha\in\indset$.
Let $X=\prod_{\alpha\in\indset} X_\alpha$.
For each $\beta\in\indset$, let
$\tprojb:X\rightarrow X_\beta$ denote the
\indexg{projection map}%
\emph{projection map}
defined by $\tprojb(x)=\tprojb(\tupset{x_\alpha}{\alpha\in\indset})=x_\beta$
for all tuples $x\in X$.
The \emph{product topology} on $X$
is then the topology generated by a subbase consisting of all sets
\begin{equation}   \label{eqn:prod-topo-subbase}
  \tprojinva(U_\alpha)
  =
  \Braces{ x \in X
    :\:
    x_\alpha \in U_\alpha
  },
\end{equation}
for all $\alpha\in\indset$ and all $U_\alpha\in\topo_\alpha$, that is,
all sets $U_\alpha$ that are open in $X_\alpha$.
The resulting topological space is called a
\indexg{product space}%
\emph{product space}.
Unless explicitly stated otherwise, we always assume this product
topology when working with Cartesian products.

If each topology $\topo_\alpha$ is generated by a subbase
$\calS_\alpha$, then
the same product topology is generated if only subbase elements are
used, that is, if a subbase for the product space is constructed
using only sets
$\tprojinva(U_\alpha)$ for $\alpha\in\indset$ and
$U_\alpha\in\calS_\alpha$.

Equivalently, the product topology is generated by a base consisting
of all sets $\prod_{\alpha\in\indset} U_\alpha$ where each $U_\alpha$
is open in $X_\alpha$ for all $\alpha\in\indset$, and where
$U_\alpha=X_\alpha$ for all but finitely many values of $\alpha$
\idxmunk\citep[Theorem~19.1]{munkres}.
In particular, $X\times Y$, where $X$ and $Y$ are topological spaces, has as
base all sets $U\times V$ where $U$ and $V$ are open in $X$ and $Y$,
respectively.

Here are some properties of product spaces in the product topology.
Note especially that the product of compact spaces is compact.

\begin{proposition}  \label{pr:prod-top-props}
  Let $\indset$ be an index set, and for each $\alpha\in\indset$,
  let $X_\alpha$ be a topological space.
  Let $X=\prod_{\alpha\in\indset} X_\alpha$ be endowed with the product
  topology.
  Then:
  \begin{letter-compact}
  \item  \label{pr:prod-top-props:a}
    If each $X_\alpha$ is Hausdorff, then $X$ is also Hausdorff.
  \item  \label{pr:prod-top-props:b}
    \indexg{Tychonoff's theorem}%
    \indexg{compactness!product space@of product space}%
    \indexg{product space!compactness of}%
    \textup{(Tychonoff's theorem)~}
    If each $X_\alpha$ is compact, then $X$ is also compact.
  \item  \label{pr:prod-top-props:c}
    If each $X_\alpha$ is first-countable, and if $\indset$ is
    countable,
    then $X$ is also first-countable.
  \item  \label{pr:prod-top-props:d}
    Let $A_\alpha\subseteq X_\alpha$, for each $\alpha\in\indset$.
    Then
    \[
       \clbar{\prod_{\alpha\in\indset} A_\alpha}
       =
       \prod_{\alpha\in\indset} \Abar_\alpha,
    \]
    where the closure on the left is in $X$,
    and
    where $\Abar_\alpha$ denotes closure of $A_\alpha$ in $X_\alpha$.
  \item  \label{pr:prod-top-props:e}
    Let $\seq{x_t}$ be a sequence in $X$, and let $x\in X$.
    Then $x_t\rightarrow x$ in $X$
    if and only if $x_t(\alpha)\rightarrow x(\alpha)$ in $X_\alpha$
    for all $\alpha\in\indset$.
  \end{letter-compact}
\end{proposition}

\begin{proof}
For the first four parts,
see, respectively, Theorems~19.4,~37.3,~30.2,~19.5
of \idxmunk\citet{munkres}.
Part~(\ref{pr:prod-top-props:e}) follows by
Lemma 2.52 of \idxhitch\citet{hitchhiker_guide_analysis}.
\end{proof}

In particular, if $X$ and $Y$ are topological spaces, then
\Cref{pr:prod-top-props}(\ref{pr:prod-top-props:a},\ref{pr:prod-top-props:b},\ref{pr:prod-top-props:c})
shows that if both $X$ and $Y$ are, respectively,
Hausdorff or compact or first-countable, then so is $X\times Y$.
Further,
\Cref{pr:prod-top-props}(\ref{pr:prod-top-props:e})
shows that for all sequences $\seq{x_t}$ in $X$
and $\seq{y_t}$ in $Y$, and for all $x\in X$ and $y\in Y$,
the sequence $\seq{\rpair{x_t}{y_t}}$ converges to
$\rpair{x}{y}$ (in $X\times Y$) if and only if
$x_t\rightarrow x$ and $y_t\rightarrow y$.%
\indexg{product topology|)}%
\indexg{topology!product|)}

If $\topo$ and $\topo'$ are topologies on some space $X$ with
$\topo\subseteq\topo'$, then
we say that $\topo$ is
\emph{coarser} (or \emph{smaller} or \emph{weaker})
than $\topo'$, and likewise, that $\topo'$ is
\emph{finer} (or \emph{larger} or \emph{stronger})
than $\topo$.
The intersection $\bigcap_{\alpha\in\indset}\topo_\alpha$ of a family of topologies
$\set{\topo_\alpha}_{\alpha\in\indset}$ on $X$
is also a topology on~$X$.
Consequently,
for any family $\calS$ of sets whose union is $X$, there exists a unique coarsest
topology that contains $\calS$, namely,
the intersection of all topologies
that contain $\calS$. This is precisely the topology generated by
$\calS$ as a subbase.

\indexg{topology!weak|(}%
\indexg{weak topology|(}%
Let $\calH$ be a family of functions from a set $X$ to a topological
space $Y$.
The \emph{weak topology}
with respect to $\calH$ is the coarsest topology on $X$
under which all the functions in $\calH$ are continuous.
By \Cref{prop:cont}(\ref{prop:cont:sub}), the weak topology with respect to $\calH$ has a subbase
consisting of all sets $h^{-1}(U)$
for all $h\in\calH$ and all $U$ open in $Y$
(or, alternatively, all elements of a subbase for $Y$).
In particular,
the product topology on $X^\indset$ is the same as
the weak topology on this set with respect to
the set of all projection maps $\tproja$ for $\alpha\in\indset$.%
\indexg{topology!weak|)}%
\indexg{weak topology|)}

Let $Z$ be a set, and let $Y$ be a topological space.
As noted earlier, $Y^Z$ is the set of all functions mapping $Z$ into
$Y$, that is, all functions $f:Z\rightarrow Y$.
This set then is an important special case of a Cartesian product.
In this case,
the projection maps $\tprojz$ simply become function evaluation maps
so that $\tprojz(f)=f(z)$,
for $z\in Z$ and $f:Z\rightarrow Y$.
As such, a subbase for the product topology on $Y^Z$ consists of all
sets
\begin{equation}   \label{eq:fcn-cl-gen-subbase}
   \tprojinvz(U)
   =
   \Braces{ f\in Y^Z :\: f(z)\in U },
\end{equation}
for all $z\in Z$ and $U$ open in $Y$, or alternatively, all $U$ that
are elements of a subbase for~$Y$.\looseness=-1

\indexg{topology!pointwise convergence@of pointwise convergence|(}%
The product topology on the function space $Y^Z$ is also called
the \emph{topology of pointwise convergence}
due to the following property:

\begin{proposition}  \label{pr:prod-top-ptwise-conv}
  Let $Z$ be a set, let $Y$ be a topological space,
  let $f:Z\rightarrow Y$, and
  let $\seq{f_t}$ be a sequence of functions in $Y^Z$.
  Then $f_t\rightarrow f$ in the product topology on $Y^Z$
  if and only if
  $f_t(z)\rightarrow f(z)$ for all $z\in Z$.
\end{proposition}

\begin{proof}
  This
  follows by
  \idxhitch\citet[Lemma~2.52]{hitchhiker_guide_analysis}.%
\indexg{topology!pointwise convergence@of pointwise convergence|)}
\end{proof}

\section{Lower semicontinuity on first-countable spaces}
\label{sec:prelim:lower-semicont}

\indexg{lower semicontinuity|(}%
A function $f:X\to\eR$ on a first-countable space $X$ is said to be
\indexg{lower semicontinuity!defined|(}%
\emph{lower semicontinuous at $x\in X$}
if $f(x)\le\liminf f(x_t)$
for every sequence $\seq{x_t}$ in~$X$ with $x_t\to x$.
Likewise, $f$ is
\indexg{upper semicontinuity}%
\emph{upper semicontinuous at $x\in X$}
if $f(x)\ge\limsup f(x_t)$
for every sequence $\seq{x_t}$ in $X$ with $x_t\to x$.
We only define semicontinuity on
first-countable spaces, but it is also possible to give definitions for general
topological spaces (see \idxhitch\citealp[Section~2.10]{hitchhiker_guide_analysis}).
By \Cref{prop:first:properties}(\ref{prop:first:cont}), $f$ is
continuous at~$x$ if and only if it is both lower semicontinuous and
upper semicontinuous at $x$.
Below, we focus on lower semicontinuity, but symmetric results hold
for upper semicontinuous functions as well since
$f$ is upper semicontinuous at $x$ if and only if $-f$ is lower
semicontinuous at~$x$.\looseness=-1

A function $f:X\rightarrow\Rext$ is \emph{lower semicontinuous} if it
is lower semicontinuous at every point in
\indexg{lower semicontinuity!defined|)}%
$X$.
\indexg{sublevel sets!lower semicontinuity and|(}%
These are precisely the functions whose epigraphs are closed in $X\times\R$, or equivalently, the functions whose sublevel sets are closed:

\begin{proposition}
\label{prop:lsc}
Let $f:X\to\eR$ where $X$ is a first-countable space.
Then the following are equivalent:
\begin{letter-compact}
\item $f$ is lower semicontinuous.
  \label{prop:lsc:a}
\item The epigraph of $f$ is closed in $X\times\R$.%
\indexg{epigraph!lower semicontinuity and}
  \label{prop:lsc:b}
\item The sublevel set $\set{x\in X:\:f(x)\le\alpha}$ is closed for every $\alpha\in\eR$.
  \label{prop:lsc:c}
\end{letter-compact}
\end{proposition}%

\begin{proof}
~
\begin{proof-parts}
\pfpart{%
  (\ref{prop:lsc:a})
  $\Rightarrow$
  (\ref{prop:lsc:b}):}
Note first that $X\times\R$ is first-countable by
\Cref{pr:prod-top-props}(\ref{pr:prod-top-props:c})
since both $X$ and $\R$ are first-countable.

Assume $f$ is lower semicontinuous, and let
$\rpair{x}{y}\in X\times\R$ be in the closure of $\epi f$.
We aim to show that $\rpair{x}{y}\in\epi f$, proving that
$\epi f$ is closed.
By
\Cref{prop:first:properties}(\ref{prop:first:closure}),
there exists a sequence $\seq{\rpair{x_t}{y_t}}$ in $\epi f$ converging to $\rpair{x}{y}$,
implying $x_t\rightarrow x$ and $y_t\rightarrow y$
by \Cref{pr:prod-top-props}(\ref{pr:prod-top-props:e}).
By lower semicontinuity, and the fact that $f(x_t)\le y_t$, we have
\[
  f(x)\le\liminf f(x_t)\le\liminf y_t=y.
\]
Thus, $\rpair{x}{y}\in\epi f$, so $\epi f$ is closed.

\pfpart{%
  (\ref{prop:lsc:b})
  $\Rightarrow$
  (\ref{prop:lsc:c}):}
Assume $\epi{f}$ is closed, and
for $\alpha\in\eR$, let
\[ S_{\alpha}=\set{x\in X:\:f(x)\le\alpha}. \]
We need to show that $S_\alpha$ is closed for every
$\alpha\in\eR$. This is trivially true for $\alpha=+\infty$, since
$S_{+\infty}=X$.

For $\alpha\in\R$, let $x\in X$ be in the closure of $S_\alpha$,
implying, since $X$ is first-countable, that there exists a sequence
$\seq{x_t}$ in $S_\alpha$ with $x_t\rightarrow x$.
Then $\rpair{x_t}{\alpha}\in\epi f$ for all $t$,
and furthermore, the sequence $\seq{\rpair{x_t}{\alpha}}$ in $\epi f$
converges to $\rpair{x}{\alpha}$
(by \Cref{pr:prod-top-props}\ref{pr:prod-top-props:e}).
Since $\epi f$ is closed, this implies that $\rpair{x}{\alpha}$ is in $\epi f$. Thus, $x\in S_\alpha$, and therefore $S_\alpha$ is closed.

The remaining case, $\alpha=-\infty$, follows because $S_{-\infty}=\bigcap_{\alpha\in\R} S_\alpha$, which must be closed as an intersection of closed sets.\looseness=-1

\pfpart{%
  (\ref{prop:lsc:c})
  $\Rightarrow$
  (\ref{prop:lsc:a}):}
Assume all sublevel sets are closed.
Let $\seq{x_t}$ be a sequence in $X$ converging to some $x\in X$. We will
argue that $f(x)\le\liminf f(x_t)$. If $\liminf f(x_t)=+\infty$, then
this trivially holds. Otherwise, let $\alpha\in\R$ be such that
$\alpha>\liminf f(x_t)$. Then the sequence $\seq{x_t}$ has infinitely
many elements $x_t$ with $f(x_t)\le \alpha$. Let $\seq{x'_t}$ be the
subsequence consisting of those elements. The set
$S=\set{x\in X:\: f(x)\le \alpha}$ is closed, by assumption,
and includes every $x'_t$. Since $x'_t\to x$, we must have $x\in S$ by
\Cref{prop:first:properties}(\ref{prop:first:closure}), so $f(x)\le \alpha$. Since this is true for all $\alpha>\liminf f(x_t)$, we obtain $f(x)\le\liminf f(x_t)$.%
\indexg{sublevel sets!lower semicontinuity and|)}%
\qedhere
\end{proof-parts}
\end{proof}

\indexg{lower semicontinuity!minimum attained|(}%
\indexg{compactness!extremum attained|(}%
Lower semicontinuous functions on compact sets are particularly well-behaved because they always attain a minimum:

\begin{proposition}   %
\label{thm:weierstrass}
Let $f:X\to\eR$ be a lower semicontinuous function on a first-countable compact space $X$. Then $f$ attains a minimum on $X$.
\end{proposition}
\begin{proof}
Let $\alpha=\inf f$. If $\alpha=+\infty$, then every $x\in X$ is a minimizer,
so we assume henceforth that $\alpha<+\infty$. We construct a sequence $\seq{x_t}$ such that $f(x_t)\to\alpha$:
If $\alpha=-\infty$, then let each $x_t$ be such that $f(x_t)<-t$.
Otherwise,
if $\alpha\in\R$, then let $x_t$ be such that $f(x_t)<\alpha+1/t$. In
either case, such points $x_t$ must exist since $\alpha=\inf f$.

By compactness (and first-countability),
the sequence $\seq{x_t}$ has a convergent subsequence, which we denote $\seq{x'_t}$. Let $x$ be a limit of $\seq{x'_t}$. Then by lower semicontinuity we have
\[
  f(x)\le\liminf f(x'_t)=\lim f(x'_t)=\alpha.
\]
Since $\alpha=\inf f$, the point $x$ must be a minimizer of~$f$.%
\indexg{lower semicontinuity!minimum attained|)}%
\indexg{compactness!extremum attained|)}
\end{proof}

Here are some operations that preserve lower semicontinuity:

\begin{proposition}   \label{pr:lsc-res-domain}
  Let $f:X\to\eR$ be a lower semicontinuous function on a
  first-countable space $X$.
  \begin{letter-compact}
  \item   \label{pr:lsc-res-domain:a}
    Let $Y$ be first-countable, and suppose $g:Y\rightarrow X$ is
    continuous.
    Then $f\circ g$ is lower semicontinuous.
  \item   \label{pr:lsc-res-domain:b}
    Let $Z\subseteq X$.
    Then $\resfcn{f}{Z}$ is lower semicontinuous.
  \item   \label{pr:lsc-res-domain:c}
    Let $g:X\to\eR$ be lower semicontinuous,
    and assume $f$ and $g$ are summable.
    Then $f+g$ is lower semicontinuous.
  \end{letter-compact}
\end{proposition}

\begin{proof}
  ~

\begin{proof-parts}
\pfpart{Part~(\ref{pr:lsc-res-domain:a}):}
Let $h=f\circ g$, and 
let $\seq{y_t}$ be any sequence in $Y$ converging to some point
$y\in Y$.
Since $g$ is continuous, $g(y_t)\rightarrow g(y)$.
Therefore, since $f$ is lower semicontinuous,
\[
  \liminf h(y_t)
  =
  \liminf f\bigParens{g(y_t)}
  \geq
  f\bigParens{g(y)}
  =
  h(y),
\]
proving that $h$ is lower semicontinuous.

\pfpart{Part~(\ref{pr:lsc-res-domain:b}):}
Let $i:Z\rightarrow X$ be the inclusion map $i(z)=z$,
which is continuous.
Then by part~(\ref{pr:lsc-res-domain:a}),
$\resfcn{f}{Z}=f\circ i$ is lower semicontinuous.

\pfpart{Part~(\ref{pr:lsc-res-domain:c}):}
Let $h=f+g$, and let $\seq{x_t}$ be any sequence in $X$ converging
to some point $x\in X$. Then
\begin{align*}
  \liminf h(x_t)
  &=\liminf [f(x_t)+g(x_t)]
\\
  &
  \ge [\liminf f(x_t)] \plusd [\liminf g(x_t)]
  \ge f(x) + g(x)
  = h(x).
\end{align*}
The first inequality is by \Cref{prop:lim:eR}(\ref{i:liminf:eR:sum});
the second inequality is by lower semicontinuity and summability of $f$ and $g$.
This proves that $h$ is lower semicontinuous.
\qedhere
\end{proof-parts}
\end{proof}

\begin{proposition}
\label{pr:lsc-sup}
Let $f_i:X\to\eR$ for all $i\in I$, where $I$ is any index set,
and $X$ is a first-countable space. Let $h=\sup_{i\in I}f_i$
be their pointwise supremum. Let $x\in X$ be a point where
each $f_i$ is lower semicontinuous. Then $h$ is also lower
semicontinuous at $x$. Consequently, if each $f_i$ is
lower semicontinuous then so is $h$.
\end{proposition}
\begin{proof}
  Let $\seq{x_t}$ be a sequence in $X$ such that $x_t\to x$.
  Then for each $i\in I$,
  \[
    \liminf h(x_t)\ge\liminf f_i(x_t)\ge f_i(x).
  \]
  Taking the supremum over $i$ yields
  $\liminf h(x_t)\ge\sup_{i\in I} f_i(x)=h(x)$. Thus, $h$
  is lower semicontinuous at $x$.%
\indexg{lower semicontinuity|)}%
\end{proof}

\indexg{lower semicontinuous hull|(}%
The \emph{lower semicontinuous hull}%
\indexg{lower semicontinuous hull!defined}
of $f:X\to\eR$, where $X$ is
first-countable, is the function $(\lsc f):X\to\eR$
defined, for $x\in X$, as
\begin{equation}
\label{eq:lsc:liminf:X:prelims}
\indexm{lsc}{$\lsc f$}{lower semicontinuous hull}%
 (\lsc f)(x)
 = \InfseqLiminf{\seq{x_t}}{X}{x_t\rightarrow x}
                {f(x_t)},
\end{equation}
where the notation means that the infimum is over all sequences $\seq{x_t}$ in $X$ converging to~$x$.

\begin{proposition}  \label{pr:lsc-is-idempotent}
  Let $f:X\to\eR$ where $X$ is a first-countable space,
  and let $x\in X$.
  Then
  $f(x) = (\lsc f)(x)$ if and only if $f$ is lower semicontinuous
  at $x$.
\end{proposition}

\begin{proof}
From definitions, 
$f$ is lower semicontinuous at $x$
if and only if
$f(x) \leq (\lsc f)(x)$.
Since $(\lsc f)(x)\leq f(x)$ always
(as witnessed by the trivial sequence $x_t=x$),
this proves the claim.
\end{proof}

The function $\lsc f$ can be characterized as the greatest lower semicontinuous function that is majorized by $f$, and its epigraph as the closure of $\epi f$ in $X\times\R$:

\begin{proposition}   \label{prop:lsc:characterize}
  Let $f:X\to\eR$ where $X$ is a first-countable space.
  Then:
  \begin{letter-compact}
  \item   \label{prop:lsc:characterize:a}
    $\lsc f$ is lower semicontinuous and majorized by $f$.
  \item   \label{prop:lsc:characterize:b}
    For every lower semicontinuous function $h:X\to\eR$, if
    $h\leq f$ then $h\leq\lsc f$.
  \item   \label{prop:lsc:characterize:c}
    $\epi(\lsc f)$ is the closure of $\epi f$ in $X\times\R$.
  \end{letter-compact}
\end{proposition}

\begin{proof}
Let $E=\clbar{\epi f}$, the closure of $\epi{f}$ in $X\times\R$.
For $x\in X$, let
\[ Y_x=\set{y\in\R :\: \rpair{x}{y}\in E}, \]
and let $g(x)=\inf Y_x$.
We prove several claims regarding $g$, and then finally show that it
is equal to $\lsc{f}$, which will prove the proposition.

\begin{claimpx}   \label{cl:prop:lsc:characterize:1}
  $E=\epi g$.
\end{claimpx}

\begin{proofx}
Let $x\in X$.
We prove the claim by showing that
\begin{equation}   \label{eq:prop:lsc:characterize:1}
  Y_x=\set{y\in\R:\: y\ge g(x)}.
\end{equation}

If $g(x)=+\infty$ then $Y_x=\emptyset$, vacuously proving
\eqref{eq:prop:lsc:characterize:1}
in this case.

Assume therefore that $g(x)<+\infty$.
If $y\in Y_x$, then $y\geq \inf Y_x = g(x)$;
thus, $Y_x\subseteq [g(x),+\infty)$.
For the reverse inclusion, let $y\in\R$ be such that $y>g(x)$. Since $g(x)=\inf Y_x$, there exists $y'\in Y_x$, such that $g(x)\le y'< y$. This means that $\rpair{x}{y'}\in E$, so $\rpair{x_t}{y'_t}\to\rpair{x}{y'}$ for some sequence $\seq{\rpair{x_t}{y'_t}}$ in $\epi f$.
For any $s\in\Rstrictpos$, we then also have $\rpair{x}{y'+s}\in E$ since the points $\rpair{x_t}{y'_t+s}$ are also in $\epi f$ and converge to $\rpair{x}{y'+s}$. In particular, setting $s=y-y'$, we obtain that $\rpair{x}{y}\in E$, and hence also $y\in Y_x$. Thus, we have shown that $(g(x),+\infty)\subseteq Y_x$. 
If $g(x)=-\infty$, this completes the proof of
\eqref{eq:prop:lsc:characterize:1}.
It therefore remains only to argue that $g(x)\in Y_x$ when $g(x)\in\R$.

Suppose then that $g(x)\in\R$, and let $y_t=g(x)+1/t$ for all $t$.
Then $y_t\in Y_x$, so $\rpair{x}{y_t}\in E$. And since $E$ is closed, the point $\rpair{x}{g(x)}$,
which is a limit of $\seq{\rpair{x}{y_t}}$, is also in $E$. Thus, $g(x)\in Y_x$.
This completes the proof of
\eqref{eq:prop:lsc:characterize:1}
in all cases, proving the claim.
\end{proofx}

\begin{claimpx}   \label{cl:prop:lsc:characterize:2}
  $g$ is lower semicontinuous and $g\leq f$.
\end{claimpx}
  
\begin{proofx}
Since $\epi g$ is closed, $g$ is lower semicontinuous by
\Cref{prop:lsc}(\ref{prop:lsc:b},\ref{prop:lsc:a}).
Also, $\epi f\subseteq E = \epi g$, so $g$ is majorized by $f$. 
\end{proofx}

\begin{claimpx}   \label{cl:prop:lsc:characterize:3}
Suppose $h:X\rightarrow\Rext$
is lower semicontinuous and that $h\leq f$.
Then $h\leq g$.
\end{claimpx}

\begin{proofx}
Since $h$ is lower semicontinuous, $\epi h$ is closed
(by \Cref{prop:lsc}\ref{prop:lsc:a}\ref{prop:lsc:b}).
And since $h\leq f$,
$\epi f\subseteq\epi h$.
Therefore,
$\epi g=\clbar{\epi f}\subseteq\epi h$, so $h\leq g$.
\end{proofx}

\begin{claimpx}   \label{cl:prop:lsc:characterize:4}
  $g=\lsc f$.
\end{claimpx}

\begin{proofx}
Let $x\in X$. Then
\[
  g(x)
  \le \InfseqLiminf{\seq{x_t}}{X}{x_t\rightarrow x}
                   {g(x_t)}
  \le \InfseqLiminf{\seq{x_t}}{X}{x_t\rightarrow x}
                   {f(x_t)}
  =
  (\lsc f)(x)
\]
with both inequalities following from
\Cref{cl:prop:lsc:characterize:2}.
It remains to prove $g(x)\ge(\lsc f)(x)$.
This is trivial if $g(x)=+\infty$, so we assume
$g(x)<+\infty$.
Let $y\in\R$ be such that $y>g(x)$.
Then $\rpair{x}{y}\in E$ (by \Cref{cl:prop:lsc:characterize:1}),
meaning that
$\rpair{x_t}{y_t}\to\rpair{x}{y}$ for some sequence
$\seq{\rpair{x_t}{y_t}}$ in $\epi f$. Therefore,
\[
  y
  =
  \liminf y_t
  \ge
  \liminf f(x_t)
  \ge
  (\lsc f)(x),
\]
with the last inequality following from $\lsc f$'s definition
(Eq.~\ref{eq:lsc:liminf:X:prelims}).
Since this holds for all $y>g(x)$, we obtain $g(x)\ge(\lsc f)(x)$, finishing the proof.
\end{proofx}

Combining \Cref{cl:prop:lsc:characterize:4}
with
Claims~\ref{cl:prop:lsc:characterize:2},~\ref{cl:prop:lsc:characterize:3},
and~\ref{cl:prop:lsc:characterize:1}
now yields,
respectively,
parts~(\ref{prop:lsc:characterize:a}),~(\ref{prop:lsc:characterize:b}),
and~(\ref{prop:lsc:characterize:c})
of the proposition.
\end{proof}

\begin{proposition}  \label{pr:lsc-pt-in-neigh}
  Let $f:X\to\eR$ where $X$ is a first-countable space,
  and let $x\in X$.
  Let $U$ be a neighborhood of $x$, and suppose
  $(\lsc f)(x)<\beta$ for some $\beta\in\R$.
  Then there exists a point $x'\in U$ with $f(x')<\beta$.
\end{proposition}

\begin{proof}
Since $(\lsc f)(x)<\beta$,
by $\lsc f$'s definition (Eq.~\ref{eq:lsc:liminf:X:prelims}),
there exists a sequence $\seq{x_t}$ in $X$
converging to $x$ such that
$f(x_t) < \beta$  for infinitely many values of $t$.
Since $U$ is a neighborhood of $x$,
we have
$x_t\in U$ for all $t$ sufficiently large.
So, among the sequence elements with $f(x_t) < \beta$,
it suffices to pick $x_t$ with a sufficiently large index~$t$.
\end{proof}

\begin{proposition}   \label{pr:lsc-seq-lim-exists}
  Let $f:X\to\eR$ where $X$ is a first-countable space,
  and let $x\in X$.
  Then there exists a sequence $\seq{x_t}$ in $X$
  such that $x_t\rightarrow x$ and $f(x_t)\rightarrow(\lsc f)(x)$.
\end{proposition}

\begin{proof}
If $(\lsc f)(x)=+\infty$ then also $f(x)=+\infty$
(by \Cref{prop:lsc:characterize}\ref{prop:lsc:characterize:a}),
so the trivial
sequence $x_t=x$ satisfies the claim.
In the remainder of the proof, we consider the case $(\lsc f)(x)<+\infty$.

Let $\countset{B}$ be a nested countable neighborhood base
for $x$.
For each $t$, let
\[
  \beta_t =
  \begin{cases}
    -t & \text{if $(\lsc f)(x) = -\infty$,}\\
    (\lsc f)(x) + 1/t
       & \text{if $(\lsc f)(x) \in \R$.}\\
  \end{cases}
\]
Then by \Cref{pr:lsc-pt-in-neigh},
there must exist a point $x_t\in B_t$ for which
$f(x_t) < \beta_t$.
Thus,
\[
  (\lsc f)(x)
  \leq
  \liminf f(x_t)
  \leq
  \limsup f(x_t)
  \leq
  \limsup \beta_t
  = (\lsc f)(x),
\]
where the first inequality is by $\lsc f$'s definition
(Eq.~\ref{eq:lsc:liminf:X:prelims}) since
the sequence $\seq{x_t}$ converges to $x$
by \Cref{prop:nested:limit}.
Hence, $f(x_t) \rightarrow (\lsc f)(x)$.
\end{proof}

\begin{proposition}   \label{pr:lsc-min-exists-equals-sup}
  Let $f:X\to\eR$ where $X$ is a compact and first-countable space.
  Then $\lsc f$ attains its minimum over $X$.
  Moreover, $\min(\lsc f) = \inf f$.
\end{proposition}

\begin{proof}
The function $\lsc f$ is lower semicontinuous
(by \Cref{prop:lsc:characterize}\ref{prop:lsc:characterize:a}),
so it attains its minimum by \Cref{thm:weierstrass}.

\Cref{prop:lsc:characterize}(\ref{prop:lsc:characterize:a})
also implies that $\lsc f \leq f$, so
$\min(\lsc f)\leq \inf f$.
On the other hand, $\inf f\leq f(x)$ for all $x\in X$, meaning
the constant function $x\mapsto\inf f$
is lower semicontinuous and majorized by $f$.
Therefore, this same function is majorized by $\lsc f$ by
\Cref{prop:lsc:characterize}(\ref{prop:lsc:characterize:b}).
Thus, $\inf f\leq (\lsc f)(x)$ for all $x\in X$,
so $\inf f\leq \min(\lsc f)$.%
\indexg{lower semicontinuous hull|)}
\end{proof}

\chapter{Review of convex analysis}
\label{sec:rev-cvx-analysis}

We briefly review some fundamentals of convex analysis that
we use throughout this work.
For a more complete introduction,
the reader is referred to a standard text, such
as \idxroc\citet{ROC}, on which this review is mainly based.
As with the previous chapter,
depending on background, readers may prefer to skip or skim
this chapter and only refer back to it as needed.

\section{Convex sets}  \label{sec:prelim:convex-sets}

The
\indexg{line segments}%
\emph{line segment joining points $\xx$ and $\yy$} in $\Rn$ is the
set $\braces{(1-\lambda)\xx+\lambda\yy :\: \lambda\in [0,1]}$.
For any $\xx\in\Rn$ and $\vv\in\Rn\wo\set{\zero}$,
the
\indexg{halflines (standard)}%
\emph{halfline with endpoint $\xx$ in the direction of $\vv$} is the set
$\set{\xx+\lambda\vv :\: \lambda\in \Rpos}$. Halflines with the origin as endpoint
are called
\indexg{rays (standard)}\emph{rays}.\looseness=-1

\indexg{convex sets (standard)|(}%
A set $S\subseteq\Rn$ is \emph{convex}
if the line segment joining any two points in $S$ is also
entirely included in $S$.
If $\lambda_1,\ldots,\lambda_m$ are nonnegative and sum to $1$,
then the point
$\sum_{i=1}^m \lambda_i \xx_i$
is said to be a
\indexg{convex combinations|(}%
\emph{convex combination}
of $\xx_1,\ldots,\xx_m\in\Rn$.

The arbitrary intersection of a family of convex sets is also
\indexg{convex sets (standard)|)}%
convex.
\indexg{convex hull (standard)|(}%
As such, for any set $S\subseteq\Rn$, there exists a smallest convex
set that includes $S$, namely, the intersection of all convex sets
that include $S$, which is called the \emph{convex hull} of $S$,
and denoted
\indexm{conv}{$\conv S$}{convex hull}%
$\conv{S}$.
Here are some facts about convex hulls:

\begin{proposition}  \label{roc:thm2.3}
Let $S\subseteq\Rn$.
Then $\conv{S}$ consists of all convex combinations of points in $S$.
\end{proposition}

\begin{proof}
  See \idxroc\citet[Theorem~2.3]{ROC}.
\end{proof}

\indexg{caratheodorys theorem@Carath\'{e}odory's theorem|(}%
Carath\'{e}odory's theorem is more specific:

\begin{proposition}  \label{roc:thm17.1}
Let $S\subseteq\Rn$.
If $\xx\in\conv{S}$ then $\xx$ is a convex combination of at most
$n+1$ points in $S$.
\end{proposition}

\begin{proof}
  See \idxroc\citet[Theorem~17.1]{ROC}.%
\indexg{convex hull (standard)|)}%
\indexg{convex combinations|)}%
\indexg{caratheodorys theorem@Carath\'{e}odory's theorem|)}
\end{proof}

\begin{proposition}  \label{roc:thm3.3}
Let $C_1,\ldots,C_m\subseteq\Rn$ be convex and nonempty, and let
\[
   C
   =
   \conv\BiggParens{\bigcup_{i=1}^m C_i}.
\]
Then $C$ consists exactly of all points of the form
$\sum_{i=1}^m \lambda_i \xx_i$
for some $\lambda_1,\ldots,\lambda_m\in [0,1]$ which sum
to $1$, and
some $\xx_i\in C_i$ for $i=1,\ldots,m$.
\end{proposition}

\begin{proof}
This follows as a special case of \idxroc\citet[Theorem~3.3]{ROC}
(with $I=\{1,\ldots,m\}$).
\end{proof}

\begin{proposition}  \label{roc:thm3.1}
\indexg{convex sets (standard)!operations on|(}%
  Let $C,D\subseteq\Rn$ be convex.
  Then $C+D$ is also convex.
\end{proposition}

\begin{proof}
  See \idxroc\citet[Theorem~3.1]{ROC}.
\end{proof}

\begin{proposition}  \label{pr:aff-preserves-cvx}
  Let $\A\in\Rmn$ and $\bb\in\Rm$,
  and let $F:\Rn\rightarrow\Rm$
  be defined by
  $F(\xx)=\A\xx+\bb$
  for $\xx\in\Rn$.
  Let $C\subseteq\Rn$ and $D\subseteq\Rm$ both be convex.
  Then $F(C)$ and $\finv(D)$ are also both convex.
\end{proposition}

\begin{proof}
  Let $A:\Rn\rightarrow\Rm$ be the linear map associated with $\A$ so
  that $A(\xx)=\A\xx$ for $\xx\in\Rn$.
  Then $F(C)=A(C)+\bb$ and $\finv(D)=\Alininv(D-\bb)$.
  From the definition of convexity, it is straightforward to show that
  the translation of a convex set is also convex.
  Combined with \idxroc\citet[Theorem~3.4]{ROC}, the claim follows.%
\indexg{convex sets (standard)!operations on|)}
\end{proof}

\section{Hull operators}
\label{sec:prelim:hull-ops}

\indexg{hull operator|(}%
The convex hull operation is an example of a more generic hull operation,
defined as follows: Let $X$ be any set and let $\calC$ be a collection of subsets of
$X$ that includes $X$ itself and that is closed under
arbitrary intersection
(meaning that if $\calS\subseteq\calC$ then
$\bigcap_{S\in\calS} S$ is also in $\calC$).
For a set $S\subseteq X$,
we then define $\genhull S$ to be the smallest set in $\calC$ that
includes $S$, that is, the intersection of all sets in $\calC$ that
include $S$:
\[
  \genhull S
  =
  \bigcap\Braces{C\in\calC:\: S\subseteq C}.
\]
The mapping $S\mapsto \genhull S$ is called
the \emph{hull operator for $\calC$}. A hull operator for any collection of subsets of a set $X$ (that satisfies the conditions above) is referred to as a \emph{hull operator on $X$}.

For example, the convex hull operation $S\mapsto\conv{S}$ for
$S\subseteq\Rn$ is the hull operator for the set of all convex sets in
$\Rn$.
It is a hull operator on $\Rn$.
As another example, for any topological space $X$, the closure operation
$S\mapsto \Sbar$ is the hull operator for the set of all closed sets,
as we saw in \Cref{sec:prelim:topo:closed-sets}.
It is a hull operator on $X$.
We will shortly see several other examples.

Here are some general properties of hull operators:

\begin{proposition}  \label{pr:gen-hull-ops}
  Let $\calC$ be a collection of subsets of $X$ that includes
  $X$ itself and that is closed under arbitrary intersection.
  Let $\genhull$ be the hull operator for $\calC$.
  Let $S,U\subseteq X$.
  Then the following hold:
  \begin{letter-compact}
  \item  \label{pr:gen-hull-ops:a}
    $S\subseteq \genhull S$ and $\genhull S\in\calC$.
    Also,
    $S= \genhull S$ if and only if $S\in\calC$.
  \item  \label{pr:gen-hull-ops:b}
    If $S\subseteq U$ and $U\in\calC$,
    then $\genhull S \subseteq U$.
  \item  \label{pr:gen-hull-ops:c}
    If $S\subseteq U$, then $\genhull S\subseteq \genhull U$.
  \item  \label{pr:gen-hull-ops:d}
    If $S\subseteq U\subseteq \genhull S$,
    then $\genhull U = \genhull S$.
  \end{letter-compact}
\end{proposition}

\begin{proof}
  ~

\begin{proof-parts}
\pfpart{Parts~(\ref{pr:gen-hull-ops:a})
  and~(\ref{pr:gen-hull-ops:b}):}
These are immediate from definitions.

\pfpart{Part~(\ref{pr:gen-hull-ops:c}):}
Since $S\subseteq U\subseteq\genhull{U}$, and since
$\genhull{U}\in\calC$,
$\genhull{S}\subseteq\genhull{U}$
by part~(\ref{pr:gen-hull-ops:b}).

\pfpart{Part~(\ref{pr:gen-hull-ops:d}):}
$\genhull{S}\subseteq\genhull U$ by part~(\ref{pr:gen-hull-ops:c}).
And since $U\subseteq\genhull{S}$ and $\genhull{S}\in\calC$,
$\genhull U\subseteq\genhull{S}$
by part~(\ref{pr:gen-hull-ops:b}).
\qedhere
\end{proof-parts}
\end{proof}

\begin{proposition}
\label{pr:hull:union}
Let $\genhull$ be a hull operator on $X$ and let $\set{S_i}_{i\in I}$ be an indexed family of sets $S_i\subseteq X$ for some arbitrary index set $I$. Then
\[
  \genhull\BiggParens{\bigcup_{i\in I} \genhull S_i}
  =
  \genhull\BiggParens{\bigcup_{i\in I} S_i}.
\]
\end{proposition}
\begin{proof}
Let $S=\bigcup_{i\in I} S_i$ and $U=\bigcup_{i\in I}\genhull S_i$. Since $S_i\subseteq\genhull S_i$ (\Cref{pr:gen-hull-ops}\ref{pr:gen-hull-ops:a}), we have $S\subseteq U$. Moreover, $\genhull S_i\subseteq\genhull S$ (\Cref{pr:gen-hull-ops}\ref{pr:gen-hull-ops:b}) and so also $U=\bigcup_{i\in I}\genhull S_i\subseteq\genhull S$. Thus, $S\subseteq U\subseteq\genhull S$,
and so by \Cref{pr:gen-hull-ops}(\ref{pr:gen-hull-ops:d}), $\genhull U=\genhull S$.%
\indexg{hull operator|)}
\end{proof}

\section{Affine sets}   \label{sec:prelim:affine-sets}

\indexg{affine sets|(}%
A set $A\subseteq\Rn$ is \emph{affine} if it includes
$(1-\lambda)\xx+\lambda\yy$
for all $\xx,\yy\in A$ and for all $\lambda\in\R$.
The linear subspaces of $\Rn$ are
precisely the affine sets that include the origin \citep[Theorem 1.1]{ROC}.
Affine sets $A$ and $A'$ in $\Rn$ are said to be \emph{parallel} if
$A'=A+\uu$ for some $\uu\in\Rn$.

\begin{proposition}  \label{roc:thm1.2}
  Let $A\subseteq\Rn$ be a nonempty affine set.
  Then there exists a unique linear subspace parallel to $A$.
\end{proposition}

\begin{proof}
  See \idxroc\citet[Theorem~1.2]{ROC}.
\end{proof}

The
\indexg{dimension!affine set@of affine set}%
\emph{dimension} of a nonempty affine set $A$ is defined to be the
dimension of the unique linear subspace that is parallel to it.%
\indexg{affine sets|)}

\indexg{affine hull|(}%
The \emph{affine hull} of a set
$S\subseteq\Rn$, denoted $\affh{S}$,%
\indexm{aff}{$\affh{S}$}{affine hull}
is the smallest affine set that includes~$S$;
that is, $\affh{S}$ is the intersection of all affine sets that
include $S$.
Equivalently,
$\affh{S}$ consists of all \emph{affine combinations}%
\indexg{affine combinations}
of finitely many points in $S$, that is,
all combinations
$\sum_{i=1}^m \lambda_i \xx_i$ where
$\xx_1,\dotsc,\xx_m\in S$,
$\lambda_1,\dotsc,\lambda_m\in\R$,
and $\sum_{i=1}^m \lambda_i = 1$.

An expression for
the linear subspace parallel to $\affh{S}$ is given in the next
proposition:

\begin{proposition}
\label{pr:lin-aff-par}
  Let $S\subseteq\Rn$,
  and let $\uu\in\affh S$.
  Then $\spn(S-\uu)$ is the linear subspace parallel to
  $\affh S$.
\end{proposition}

\begin{proof}
  See \idxroman\citet[Theorem~16.4]{roman-lin-alg}.%
\indexg{affine hull|)}
\end{proof}

\section{Closure, interior, relative interior}

When working in $\Rn$, it is understood that we are always using the
standard Euclidean topology, as defined in
Example~\ref{ex:topo-on-rn},
unless explicitly stated otherwise.

\indexg{closure (of set in rn)@closure (of set in $\Rn$)|(}%
Let $S\subseteq\Rn$.
The closure of $S$ in $\Rn$
(that is, the smallest closed set in $\Rn$ that includes $S$)
is denoted
\indexg{closure (of set in rn)@closure (of set in $\Rn$)|)}%
\indexm{cls300}{$\cl S$}{closure (of set) in $\Rn$}%
$\cl S$.
As in \Cref{sec:prelim:topo:closed-sets},
the interior of $S$
(that is, the largest open set in $\Rn$ included in $S$)
is denoted $\intr S$.
The \emph{boundary}%
\indexg{boundary}
of $S$ is the set difference
$(\cl{S})\setminus(\intr S)$.

\indexg{relative interior|(}%
\indexg{interior!relative|(}%
Often, it is more useful to consider $S$ in the subspace topology
induced by its affine hull, $\affh{S}$.
The interior of $S$ in $\affh S$ (when viewed as a topological
subspace of $\Rn$) is called its
\emph{relative interior}, and is denoted
\indexm{ri s}{$\ri S$}{relative interior}%
$\ri S$.
Thus, a point $\xx\in\Rn$ is in $\ri S$ if and only if
$\ball(\xx,\epsilon)\cap(\affh{S})\subseteq S$ for some $\epsilon\in\Rstrictpos$
(where $\ball(\xx,\epsilon)$ is as in Eq.~\ref{eqn:open-ball-defn}).

From definitions, it follows that
\begin{equation}  \label{eq:ri-in-cl}
   \intr S \subseteq \ri S \subseteq S \subseteq \cl S
\end{equation}
for all $S\subseteq\Rn$.
The set difference $(\cl S)\setminus(\ri S)$ is called the
\emph{relative boundary}%
\indexg{boundary!relative}%
\indexg{relative boundary}
of $S$.
If $S=\ri S$, then $S$ is said to be
\indexg{relatively open (set)}%
\emph{relatively open}.
Every singleton set is relatively open and
its relative boundary is empty.
Every affine set is both closed and relatively open.

Here are some facts about the relative interior of a convex set:

\begin{proposition}   \label{pr:ri-props}
\indexg{closure (of set in rn)@closure (of set in $\Rn$)|(}%
  Let $C,D\subseteq\Rn$ be convex.
  \begin{letter-compact}
  \item   \label{pr:ri-props:roc-thm6.2a}
    $\ri C$ and $\cl C$ are convex.
  \item   \label{pr:ri-props:roc-thm6.2b}
    If $C\neq\emptyset$ then $\ri{C} \neq\emptyset$.
  \item   \label{pr:ri-props:roc-thm6.3}
    $\ri(\cl C)=\ri C$ and $\cl(\ri C) = \cl C$.
  \item   \label{pr:ri-props:roc-cor6.6.2}
    $\ri(C+D)=(\ri C) + (\ri D)$.
  \item   \label{pr:ri-props:intC-nonemp-implies-eq-riC}
    If $\intr C\neq\emptyset$
    then $\intr C = \ri C$.
  \item   \label{pr:ri-props:intC-D-implies-riC-riD}
    If $(\intr C)\cap(\cl D)\neq\emptyset$
    then $(\ri C)\cap(\ri D)\neq\emptyset$.
  \item   \label{pr:ri-props:roc-cor6.3.1}
    The following are equivalent:
    \begin{roman-compact}
    \item
      $\ri C = \ri D$.
    \item
      $\cl C = \cl D$.
    \item
      $\ri C \subseteq D \subseteq \cl C$.
    \end{roman-compact}
  \end{letter-compact}
\end{proposition}

\begin{proof}
~

\begin{proof-parts}
\pfpart{Parts~(\ref{pr:ri-props:roc-thm6.2a})
          and~(\ref{pr:ri-props:roc-thm6.2b}):}
See \idxroc\citet[Theorem~6.2]{ROC}.

\pfpart{Part~(\ref{pr:ri-props:roc-thm6.3}):}
See \idxroc\citet[Theorem~6.3]{ROC}.

\pfpart{Part~(\ref{pr:ri-props:roc-cor6.6.2}):}
See \idxroc\citet[Corollary~6.6.2]{ROC}.

\pfpart{Part~(\ref{pr:ri-props:intC-nonemp-implies-eq-riC}):}
Suppose $\xx\in\intr C$.
Then, for some $\epsilon\in\Rstrictpos$, the ball $\ball(\xx,\epsilon)$ is
included in $C$, implying that $\affh{C}=\Rn$, and so that
$\ri C = \intr C$ by definition of relative interior.\looseness=-1

\pfpart{Part~(\ref{pr:ri-props:intC-D-implies-riC-riD}):}
Suppose
$(\intr C)\cap(\cl D)\neq\emptyset$.
Since $\intr C$ is open, it follows that
$(\intr C)\cap(\ri D)\neq\emptyset$, by
\idxroc\citet[Corollary 6.3.2]{ROC}.
Therefore,
${(\ri C)\cap(\ri D)\neq\emptyset}$, by
\eqref{eq:ri-in-cl}.

\pfpart{Part~(\ref{pr:ri-props:roc-cor6.3.1}):}
See \idxroc\citet[Corollary~6.3.1]{ROC}.%
\indexg{closure (of set in rn)@closure (of set in $\Rn$)|)}%
\qedhere
\end{proof-parts}
\end{proof}

\begin{proposition}   \label{roc:thm6.1}
  Let $C\subseteq\Rn$ be convex,
  let $\xx\in\ri C$,
  let $\yy\in\cl C$,
  and let $\lambda\in (0,1]$.
  Then $\lambda\xx + (1-\lambda)\yy \in \ri C$.
\end{proposition}

\begin{proof}
See \idxroc\citet[Theorem~6.1]{ROC}.
\end{proof}

\begin{proposition}   \label{roc:thm6.4}
  Let $C\subseteq\Rn$ be convex and nonempty,
  and let $\xx\in\Rn$.
  Then $\xx\in\ri C$ if and only if
  for all $\yy\in C$ there exists $\delta\in\Rstrictpos$ such that
  $\xx+\delta (\xx-\yy) \in C$.
\end{proposition}

\begin{proof}
See \idxroc\citet[Theorem~6.4]{ROC}.
\end{proof}

\indexg{relative boundary|(}%
\indexg{boundary!relative|(}%
Relative boundary points of a convex set $C$ are arbitrarily
close to points that are not in $C$'s closure:

\begin{proposition}   \label{pr:bnd-near-cl-comp}
  Let $C\subseteq\Rn$ be convex,
  let $\xx\in(\cl C)\setminus(\ri C)$,
  and let $\epsilon\in\Rstrictpos$.
  Then there exists a point $\zz\in\Rn\setminus(\cl C)$
  such that $\norm{\xx-\zz}<\epsilon$.
\end{proposition}

\begin{proof}
Let $C'=\cl C$.
By \Cref{pr:ri-props}(\ref{pr:ri-props:roc-thm6.3}),
$\ri C' = \ri C$, so
$\xx\not\in\ri C'$.
By \Cref{roc:thm6.4} (applied to $C'$),
it follows that there must exist a point
$\yy\in C'$ such that
for all $\delta\in\Rstrictpos$,
$\zz_{\delta}\not\in C'$,
where
$\zz_{\delta}=\xx+\delta(\xx-\yy)$.
Letting $\zz=\zz_{\delta}$ for sufficiently small $\delta\in\Rstrictpos$
(so that $\norm{\xx-\zz}=\delta \norm{\xx-\yy} < \epsilon$),
the claim now follows.%
\indexg{boundary!relative|)}%
\indexg{relative boundary|)}
\end{proof}

\begin{proposition}   \label{roc:thm6.5}
\indexg{closure (of set in rn)@closure (of set in $\Rn$)|(}%
  Let $C_1,\ldots,C_m\subseteq\Rn$ be convex, and assume
  $\bigcap_{i=1}^m (\ri C_i)\neq\emptyset$.
  Then
  \[
     \cl\Parens{\bigcap_{i=1}^m C_i}
     =
     \bigcap_{i=1}^m (\cl C_i),
  \]
  and
  \[
     \ri\Parens{\bigcap_{i=1}^m C_i}
     =
     \bigcap_{i=1}^m (\ri C_i).
  \]
\end{proposition}

\begin{proof}
See \idxroc\citet[Theorem~6.5]{ROC}.%
\indexg{closure (of set in rn)@closure (of set in $\Rn$)|)}%
\end{proof}

\begin{proposition}   \label{roc:thm6.7}
  Let $C\subseteq\Rm$ be convex,
  let $\A\in\Rmn$, and let
  $A:\Rn\rightarrow\Rm$ be defined by
  $A(\xx)=\A\xx$ for $\xx\in\Rn$.
  Assume there exists $\xing\in\Rn$ for which $\A\xing\in\ri C$.
  Then
  $\Alininv(\ri C) = \ri(\Alininv(C))$.
\end{proposition}

\begin{proof}
See \idxroc\citet[Theorem~6.7]{ROC}.
\end{proof}

\begin{proposition}  \label{pr:ri-conv-finite}
  Let
  $\xx_1,\ldots,\xx_m\in\Rn$.
  Then
  \[
     \ri\bigParens{\conv \{\xx_1,\ldots,\xx_m\}}
     =
     \Braces{
       \sum_{i=1}^m \lambda_i \xx_i
       :\:
       \lambda_1,\ldots,\lambda_m \in\Rstrictpos,\,
       \sum_{i=1}^m \lambda_i = 1
     }.
  \]
\end{proposition}

\begin{proof}
This follows from
\idxroc\citet[Theorem~6.9]{ROC}
(with each set $C_i$, in his notation, set to the singleton
$\{\xx_i\}$).%
\indexg{interior!relative|)}%
\indexg{relative interior|)}
\end{proof}

\section{Cones}
\label{sec:prelim:cones}

\indexg{cones (standard)|(}%
A set $K\subseteq\Rn$ is a \emph{cone} if $\zero\in K$
and if $K$ is closed under
multiplication by positive scalars, that is, if $\lambda \xx\in K$ for
all $\xx\in K$ and $\lambda\in\Rstrictpos$.
Equivalently, $K$ is a cone if it is nonempty and closed under
multiplication by nonnegative scalars.
Note importantly that some authors, such as
\idxroc\idxwets\citet{rock_wets} and
\idxborlew\citet{borwein_lewis_06},
require that a cone include the origin, as we have done here,
but others, including \idxroc\citet[\pcite{13}]{ROC}, do not.

Every linear subspace is a closed convex cone.
\indexg{conic hull (standard)|(}%
The arbitrary intersection of cones is itself a cone, and therefore
the same holds for convex cones.
Consequently,
for any set $S\subseteq\Rn$, there exists a smallest convex cone that
includes $S$, which is the intersection of all convex cones
containing $S$.
The resulting convex cone, denoted $\cone{S}$,%
\indexm{cone s300}{$\cone{S}$}{conic hull (standard)}
is called the
\emph{conic hull of $S$} or the
\emph{convex cone generated by $S$}.

\indexg{conic combinations|(}%
For $m\geq 0$ and $\lambda_1,\ldots,\lambda_m\in\Rpos$,
the point $\sum_{i=1}^m \lambda_i \xx_i$ is said to be
a \emph{conic combination} of the points $\xx_1,\ldots,\xx_m\in\Rn$.

\begin{proposition}   \label{pr:scc-cone-elts}
  Let $S\subseteq\Rn$.
  \begin{letter-compact}
  \item   \label{pr:scc-cone-elts:b}
    The conic hull
    $\cone{S}$ consists of all conic combinations of
    points in $S$.%
\indexg{conic combinations|)}
  \item   \label{pr:scc-cone-elts:c:new}
    Suppose $S$ is convex.
    Then
    \[
       \cone{S} = \{\zero\} \cup
             \Braces{ \lambda\xx :\: \xx\in S,\,\lambda\in\Rstrictpos }.
    \]
  \item   \label{pr:scc-cone-elts:d}
    Suppose $S$ is a cone.
    Then $S$ is convex if and only if it is closed under vector
    addition
    (so that $\xx+\yy\in S$ for all $\xx,\yy\in S$).
  \item   \label{pr:scc-cone-elts:span}
    $\spn{S}=\cone(S \cup -S)$.
  \end{letter-compact}
\end{proposition}

\begin{proof}
  ~
  
\begin{proof-parts}
\pfpart{%
  Parts~(\ref{pr:scc-cone-elts:b})
  and~(\ref{pr:scc-cone-elts:c:new}):
}
See Corollaries~2.6.2 and~2.6.3 of \idxroc\citet{ROC}
as well as the immediately following discussion
(keeping in mind, as mentioned above, that
Rockafellar's definition of a cone differs slightly from ours).

\pfpart{Part~(\ref{pr:scc-cone-elts:d}):}
See \idxroc\citet[Theorem~2.6]{ROC}.

\pfpart{Part~(\ref{pr:scc-cone-elts:span}):}
By \Cref{pr:span-is-lin-comb},
the set $\spn{S}$ consists of all linear combinations
$\sum_{i=1}^m \lambda_i \xx_i$
where each $\lambda_i\in\R$ and
each $\xx_i\in S$.
From part~(\ref{pr:scc-cone-elts:b}), $\cone{(S\cup -S)}$ consists of
all such linear combinations with each $\lambda_i\in\Rpos$ and each
$\xx_i\in S\cup -S$.
Using the simple fact that $\lambda_i\xx_i=(-\lambda_i)(-\xx_i)$,
it follows that these two sets are the same.%
\indexg{conic hull (standard)|)}
\qedhere
\end{proof-parts}
\end{proof}

\begin{proposition}
\label{prop:cone-linear}
\indexg{cones (standard)!under linear map|(}%
Let $K\subseteq\Rn$ be a convex cone and let $\A\in\Rmn$. Then $\A K$ is also a convex cone.
\end{proposition}

\begin{proof}
Let $\yy\in \A K$ and $\lambda\in\Rstrictpos$. Then we
must have $\yy=\A\xx$ for some $\xx\in K$, so also $\lambda\xx\in K$,
and thus $\lambda\yy=\A(\lambda\xx)\in\A K$.
Also, $\A K$ includes the origin since $\A\zerov{n}=\zerov{m}$.
Therefore, $\A K$ is a cone, and is convex by
\Cref{pr:aff-preserves-cvx}.%
\indexg{cones (standard)!under linear map|)}%
\end{proof}

\indexg{polar (primal)!defined for standard cones|(}%
Let $K\subseteq\Rn$ be a convex cone.
Its \emph{polar},
denoted $\Kpol$,
is the set of points whose inner product with every
point in $K$ is nonpositive, that is,
\begin{equation}  \label{eqn:polar-def}
\indexm{k 300}{$\Kpol$}{(primal) polar}%
  \Kpol = \Braces{\uu\in\Rn :\:
                   \xx\cdot\uu\leq 0 \mbox{ for all } \xx\in K}.
\end{equation}
We write $\dubpolar{K}$ for the polar of $\Kpol$; that is,
$\dubpolar{K}=\polar{(\Kpol)}$.

\begin{proposition}  \label{pr:polar-props}
  Let $K,J\subseteq\Rn$ be convex cones.
  \begin{letter-compact}
  \item  \label{pr:polar-props:clK-cvx-cone}
    $\cl K$ is also a convex cone.
  \item  \label{pr:polar-props:a}
    $\Kpol=\polar{(\cl K)}$.
  \item  \label{pr:polar-props:b}
    $\Kpol$ is a closed (in $\Rn$), convex cone.
  \item  \label{pr:polar-props:c}
    $\dubpolar{K}=\cl K$.
  \item  \label{pr:polar-props:d}
    If $J\subseteq K$ then $\Kpol\subseteq\Jpol$.
  \item  \label{pr:polar-props:e}
    $J+K$ is also a convex cone, and
    $\polar{(J+K)} = \Jpol \cap \Kpol$.
  \item  \label{pr:polar-props:f}
    If $K$ is a linear subspace then
    $\Kpol = \Kperp$.
  \item  \label{pr:polar-props:coneSpol}
    Let $S\subseteq\Rn$.
    Then
    $
      \polar{(\cone{S})}
      =
      \Braces{\uu\in\Rn :\:
        \xx\cdot\uu\leq 0 \textup{ for all } \xx\in S
      }
    $.
  \end{letter-compact}
\end{proposition}

\begin{proof}
  ~

\begin{proof-parts}
\pfpart{%
  Parts~(\ref{pr:polar-props:clK-cvx-cone})
  and~(\ref{pr:polar-props:d}):
}
These are routine to argue from definitions.

\pfpart{%
  Parts~(\ref{pr:polar-props:a})
  and~(\ref{pr:polar-props:f}):
}
See the discussion preceding and following Theorem~14.1 of \idxroc\citet{ROC}.

\pfpart{%
  Parts~(\ref{pr:polar-props:b})
  and~(\ref{pr:polar-props:c}):
}
These follow from \idxroc\citet[Theorem~14.1]{ROC}
applied to $\cl K$
(and using part~(\ref{pr:polar-props:a})).

\pfpart{Part~(\ref{pr:polar-props:e}):}
From definitions, it is straightforward to show that $J+K$ is a cone,
and therefore a convex cone by \Cref{roc:thm3.1}.
For the expression for $\polar{(J+K)}$,
see \idxroc\citet[Corollary~16.4.2]{ROC}.

\pfpart{Part~(\ref{pr:polar-props:coneSpol}):}
Let
$U=\Braces{\uu\in\Rn :\: \xx\cdot\uu\leq 0 \text{ for all } \xx\in S}$.
We aim to show that a point $\uu\in\Rn$ is in
$\polar{(\cone{S})}$ if and only if it is in $U$.
If $\uu\in\polar{(\cone{S})}$, then $\xx\cdot\uu\leq 0$ for all
$\xx$ in $\cone{S}$, and so also for all $\xx$ in $S$;
hence, $\uu\in U$.
For the converse, suppose $\uu\in U$.
Then $\xx\cdot\uu\leq 0$ for all $\xx\in S$,
implying that the same holds for every conic combination $\xx$ of
points in $S$, and so for all $\xx\in\cone{S}$ by
\Cref{pr:scc-cone-elts}(\ref{pr:scc-cone-elts:b}).
Thus, $\uu\in\polar{(\cone{S})}$.%
\indexg{cones (standard)|)}%
\indexg{polar (primal)!defined for standard cones|)}
\qedhere
\end{proof-parts}
\end{proof}

\section{Separation theorems}
\label{sec:prelim-sep-thms}

\indexg{hyperplanes (standard)|(}%
A \emph{hyperplane}
$J$ in $\Rn$ is defined by a vector
$\vv\in\Rn\wo\{\zero\}$ and
scalar $\beta\in\R$, and consists of the set of points
\begin{equation}  \label{eqn:std-hyp-plane}
 J=\braces{\xx\in\Rn :\: \xx\cdot\vv = \beta}.
\end{equation}
\indexg{halfspaces (standard)|(}%
The hyperplane is associated with two \emph{closed halfspaces}
consisting of
those points $\xx\in\Rn$ for which $\xx\cdot\vv\leq\beta$, and those
for which $\xx\cdot\vv\geq\beta$.
It is similarly associated with two \emph{open halfspaces} defined by
corresponding strict inequalities.
If the hyperplane or the boundary of the halfspace includes the
origin, we say that it is \emph{homogeneous}%
\indexg{homogeneous (halfspace or hyperplane)}
(that is, if $\beta=0$).%
\indexg{hyperplanes (standard)|)}%
\indexg{halfspaces (standard)|)}

\indexg{separation by hyperplane (standard)|(}%
Let $C$ and $D$ be nonempty subsets of $\Rn$.
We say that a hyperplane $J$ \emph{separates} $C$ and $D$ if $C$ is
included in one of the closed halfspaces associated with $J$, and $D$
is included in the other; thus, for $J$ as defined above,
$\xx\cdot\vv\leq\beta$ for all $\xx\in C$ and
$\xx\cdot\vv\geq\beta$ for all $\xx\in D$.
(This is without loss of generality since if the reverse inequalities
hold, we can simply negate $\vv$ and $\beta$.)

\indexg{proper separation|(}%
We say that $J$ \emph{properly}
 separates $C$ and $D$ if, in addition,
$C$ and $D$ are not both included in $J$
(that is, $C\cup D \not\subseteq J$).%
\indexg{proper separation|)}%

\indexg{strong separation (standard)|(}%
We say that $J$ \emph{strongly} separates $C$ and $D$ if
for some $\epsilon>0$,
$C+\ball(\zero,\epsilon)$ is included in one of the open halfspaces
associated with $J$, and
$D+\ball(\zero,\epsilon)$ is included in the other
(where $\ball(\zero,\epsilon)$ is as defined in
Eq.~\ref{eqn:open-ball-defn}).
Equivalently, with $J$ as defined in \eqref{eqn:std-hyp-plane},
this means that
$\xx\cdot\vv<\beta-\epsilon$ for all $\xx\in C$,
and
$\xx\cdot\vv>\beta+\epsilon$ for all $\xx\in D$.

Here are several facts about separating convex sets:

\begin{proposition}   \label{roc:thm11.1}
  Let $C,D\subseteq\Rn$ be nonempty.
  Then there exists a hyperplane that strongly separates $C$ and $D$
  if and only if
  there exists $\vv\in\Rn$ such that
  \[
     \sup_{\xx\in C} \xx\cdot \vv
     <
     \inf_{\xx\in D} \xx\cdot \vv.
  \]
\end{proposition}

\begin{proof}
See \idxroc\citet[Theorem~11.1]{ROC}.
\end{proof}

\begin{proposition}   \label{roc:cor11.4.2}
  Let $C,D\subseteq\Rn$ be convex and nonempty.
  Assume $(\cl{C})\cap(\cl{D})=\emptyset$ and that
  $C$ is bounded.
  Then there exists a hyperplane that strongly separates $C$ and $D$.
\end{proposition}

\begin{proof}
See \idxroc\citet[Corollary~11.4.2]{ROC}.%
\indexg{strong separation (standard)|)}
\end{proof}

\begin{proposition}   \label{roc:thm11.2}
  Let $C\subseteq\Rn$ be nonempty, convex and relatively open,
  and let $A\subseteq\Rn$ be nonempty and affine.
  Assume $C\cap A=\emptyset$.
  Then there exists a hyperplane $J$ that includes $A$
  and such that $C$ is included in one of the open halfspaces
  associated with $J$.
  That is, there exists $\vv\in\Rn$ and $\beta\in\R$ such that
  $\xx\cdot\vv=\beta$ for all $\xx\in A$,
  and
  $\xx\cdot\vv<\beta$ for all $\xx\in C$.
\end{proposition}

\begin{proof}
See \idxroc\citet[Theorem~11.2]{ROC}.
\end{proof}

\begin{proposition}   \label{roc:thm11.3} 
\indexg{proper separation|(}%
 Let $C,D\subseteq\Rn$ be nonempty and convex.
  Then there exists a hyperplane separating $C$ and $D$ properly
  if and only if
  $(\ri C)\cap(\ri D)=\emptyset$.
\end{proposition}

\begin{proof}
See \idxroc\citet[Theorem~11.3]{ROC}.%
\indexg{separation by hyperplane (standard)|)}%
\indexg{proper separation|)}
\end{proof}

\begin{proposition}   \label{pr:con-int-halfspaces}
\indexg{convex hull (standard)!as intersection of standard halfspaces|(}%
\indexg{conic hull (standard)!as intersection of standard halfspaces|(}%
  Let $S\subseteq\Rn$.
  \begin{letter-compact}
  \item   \label{roc:cor11.5.1}
    The closure of the convex hull of $S$,
    $\cl(\conv{S})$, is equal to the intersection of all closed
    halfspaces in $\Rn$ that include~$S$.
  \item   \label{roc:cor11.7.2}
    The closure of the conic hull of $S$,
    $\cl(\cone{S})$, is equal to the intersection of all
    homogeneous closed
    halfspaces in $\Rn$ that include $S$.
  \end{letter-compact}
\end{proposition}

\begin{proof}
  ~

\begin{proof-parts}
\pfpart{Part~(\ref{roc:cor11.5.1}):}
See \idxroc\citet[Corollary~11.5.1]{ROC}.

\pfpart{Part~(\ref{roc:cor11.7.2}):}
See \idxroc\citet[Corollary~11.7.2]{ROC}.%
\indexg{convex hull (standard)!as intersection of standard halfspaces|)}%
\indexg{conic hull (standard)!as intersection of standard halfspaces|)}%
\qedhere
\end{proof-parts}
\end{proof}

\indexg{sharp discrimination|(}%
\indexg{discrimination, sharp|(}%
Let $K_1$ and $K_2$ be convex cones in $\Rn$,
and let $L_i=K_i\cap -K_i$, for $i\in\{1,2\}$.
Then $K_1$ and $K_2$ are said to be
\emph{sharply discriminated} by a vector $\uu\in\Rn$ if:
\[
   \xx\cdot\uu
   \begin{cases}
     < 0 & \text{for all $\xx\in K_1\setminus L_1$,} \\
     = 0 & \text{for all $\xx\in L_1 \cup L_2$,} \\
     > 0 & \text{for all $\xx\in K_2\setminus L_2$.}
   \end{cases}
\]
\indexg{separation by hyperplane (standard)|(}%
Except when $\uu=\zero$, this condition means that
$K_1$ and $K_2$ are separated by the hyperplane
$J=\{\xx\in\Rn :\: \xx\cdot\uu=0\}$
(since $\xx\cdot\uu\leq 0$ for $\xx\in K_1$, and
$\xx\cdot\uu\geq 0$ for $\xx\in K_2$),
and moreover that
$K_1\setminus L_1$ is included in one of the open halfspaces
associated with $J$, while
$K_2\setminus L_2$ is included in the other.
(Nonetheless,
we use the phrase ``discriminated by a vector'' rather than
``separated by a hyperplane''
since in fact we do allow $\uu=\zero$, in which case the set $J$
above is not technically a hyperplane.
If $K_1\neq L_1$ or $K_2\neq L_2$, then this possibility is ruled out.)%
\indexg{separation by hyperplane (standard)|)}

The next proposition and the main ideas of its proof are due to
\indexa{Klee, V. L.}%
\citet[Theorem~2.7]{Klee55}.

\begin{proposition}  \label{pr:cones-sharp-sep}
  Let $K_1,K_2\subseteq\Rn$ be closed convex cones, and assume
  $K_1\cap K_2=\{\zero\}$.
  Then there exists a vector $\uu\in\Rn$ that sharply discriminates
  $K_1$ and $K_2$.
\end{proposition}

\begin{proof}
Let $K=K_1-K_2$, and let $L=K\cap-K$.
Then $K$ is a closed convex cone.
To see this, note first that
$-K_2$ is a convex cone (by \Cref{prop:cone-linear}),
and is closed (by \Cref{prop:cont}\ref{prop:cont:inv:closed},
since the map $\xx\mapsto-\xx$ is continuous and is its
own inverse).
Therefore, $K=K_1+(-K_2)$ is a convex cone by
\Cref{pr:polar-props}(\ref{pr:polar-props:e}),
and is closed by
\idxroc\citet[Corollary~9.1.3]{ROC} since $K_1$ and $K_2$ are closed convex
cones whose intersection is $\{\zero\}$.

Let $\uu$ be any point in $\ri \Kpol$
(which must exist by
Propositions~\ref{pr:polar-props}\ref{pr:polar-props:b}
and~\ref{pr:ri-props}\ref{pr:ri-props:roc-thm6.2b}).
We first show that $\xx\cdot\uu\leq 0$ for $\xx\in K$
and that $\xx\cdot\uu<0$ for $\xx\in K\setminus L$
(which is equivalent to $\uu$ sharply discriminating
$K$ and $\{\zero\}$).
Then below,
we show this implies that $\uu$ sharply discriminates $K_1$ and
$K_2$.

First, $\uu\in\Kpol$, meaning $\xx\cdot\uu\leq 0$ for all $\xx\in K$.
Let $\xx\in K\setminus L$; we aim to show $\xx\cdot\uu<0$.
Then $\xx\not\in -K$, so $-\xx\not\in K = \dubpolar{K}$,
where the equality is by
\Cref{pr:polar-props}(\ref{pr:polar-props:c}).
Therefore, there exists $\ww\in\Kpol$ such that
$-\xx\cdot\ww>0$, that is, $\xx\cdot\ww<0$.
Since $\uu\in\ri{\Kpol}$,
by \Cref{roc:thm6.4},
there exists $\delta\in\Rstrictpos$ such that
the point $\vv=\uu+\delta(\uu-\ww)$ is in $\Kpol$.
We then have
\[
  \xx\cdot\uu
  =
  \frac{\xx\cdot\vv + \delta \xx\cdot\ww}
       {1+\delta}
  <
  0.
\]
The equality is by rearranging $\vv$'s definition.
The inequality is because
$\xx\cdot\ww<0$
and
$\xx\cdot\vv\leq 0$
(since $\xx\in K$ and $\vv\in \Kpol$).

Next,
we show $\uu$ sharply discriminates $K_1$ and $K_2$.
Observe first that
$K_1=K_1-\{\zero\}\subseteq K_1-K_2=K$,
since $\zero\in K_2$.

We claim $K_1\cap L\subseteq L_1$.
Let $\xx\in K_1\cap L$, which we aim to show is in $L_1$.
Then $\xx\in -K$ so $\xx=\yy_2-\yy_1$ for some $\yy_1\in K_1$ and
$\yy_2\in K_2$.
This implies $\yy_2=\xx+\yy_1$.
Since $\xx$ and $\yy_1$ are both in the convex cone $K_1$,
their sum $\xx+\yy_1$ is as well.
Therefore, since $K_1\cap K_2=\{\zero\}$, we must have
$\yy_2=\zero=\xx+\yy_1$.
Thus, $\xx=-\yy_1\in -K_1$, so $\xx\in L_1$, as claimed.

Consequently,
\begin{equation}  \label{eq:pr:cones-sharp-sep:1}
   K_1\setminus L_1
   \subseteq
   K_1\setminus (K_1\cap L)
   =
   K_1 \setminus L
   \subseteq
   K\setminus L,
\end{equation}
with the first inclusion from the preceding argument, and the second from
$K_1\subseteq K$.

Thus, if $\xx\in K_1$ then $\xx\in K$ so $\xx\cdot\uu\leq 0$,
from the argument above.
Hence,
if $\xx\in L_1$, then $\xx\in K_1$ and $-\xx\in K_1$, together
implying that $\xx\cdot\uu=0$.
And if $\xx\in K_1\setminus L_1$ then
$\xx\in K\setminus L$
(from Eq.~\ref{eq:pr:cones-sharp-sep:1}),
implying $\xx\cdot\uu<0$, as argued earlier.

By similar arguments, $K_2\subseteq -K$ and
$K_2\setminus L_2 \subseteq (-K)\setminus L$.
So if $\xx\in K_2$ then $-\xx\in K$ implying $\xx\cdot\uu\geq 0$.
Thus, if $\xx\in L_2$ then $\xx\cdot\uu=0$.
And if $\xx\in K_2\setminus L_2$ then
$-\xx\in K\setminus L$ implying $\xx\cdot\uu>0$.%
\indexg{sharp discrimination|)}%
\indexg{discrimination, sharp|)}
\end{proof}

\section{Faces and polyhedral sets}
\label{sec:prelim:faces}

\indexg{faces|(}%
Let $C\subseteq\Rn$ be convex.
A convex subset $F\subseteq C$ is said to be \emph{a face of $C$} if
for all points $\xx,\yy\in C$, and for all $\lambda\in (0,1)$,
if $(1-\lambda)\xx+\lambda\yy$ is in $F$, then $\xx$ and $\yy$ are
also in~$F$.
For instance, the faces of the cube $[0,1]^3$ in $\R^3$ are the cube's
eight vertices, twelve edges, six square faces, the entire cube, and
the empty set.

\begin{proposition}  \label{pr:face-props}
  Let $C\subseteq\Rn$ be convex, and let $F$ be a face of $C$.
  \begin{letter-compact}
  \item  \label{pr:face-props:cor18.1.1}
    If $C$ is closed (in $\Rn$), then so is $F$.
  \item  \label{pr:face-props:thm18.1}
    Let $D\subseteq C$ be convex with $(\ri D)\cap F\neq\emptyset$.
    Then $D\subseteq F$.
  \item  \label{pr:face-props:cor18.1.2}
    Let $E$ be a face of $C$ with $(\ri E)\cap (\ri F)\neq\emptyset$.
    Then $E=F$.
  \item  \label{pr:face-props:intersect}
    An arbitrary intersection of faces of $C$ is also a face of $C$.
  \end{letter-compact}
\end{proposition}

\begin{proof}
~

\begin{proof-parts}
\pfpart{Part~(\ref{pr:face-props:cor18.1.1}):}
See \idxroc\citet[Corollary~18.1.1]{ROC}.

\pfpart{Part~(\ref{pr:face-props:thm18.1}):}
See \idxroc\citet[Theorem~18.1]{ROC}.

\pfpart{Part~(\ref{pr:face-props:cor18.1.2}):}
See \idxroc\citet[Corollary~18.1.2]{ROC}.

\pfpart{Part~(\ref{pr:face-props:intersect}):}
This is straightforward from the definition.
\qedhere
\end{proof-parts}
\end{proof}

\begin{proposition}   \label{roc:thm18.2}
  Let $C\subseteq\Rn$ be convex.
  Then the relative interiors of the nonempty faces of $C$ form a
  partition of $C$.
  That is, for all $\xx\in C$, there exists a unique face $F$ of~$C$
  whose relative interior, $\ri F$, includes $\xx$.
\end{proposition}

\begin{proof}
See \idxroc\citet[Theorem~18.2]{ROC}.
\end{proof}

\begin{proposition}   \label{roc:thm18.3}
  Let $S\subseteq\Rn$,
  let $C=\conv{S}$,
  and
  let $F$ be a face of $C$.
  Then $F=\conv(S\cap F)$.
\end{proposition}

\begin{proof}
See \idxroc\citet[Theorem~18.3]{ROC}.%
\indexg{faces|)}
\end{proof}

\indexg{polyhedral sets (standard)|(}%
A convex set $C$ is \emph{polyhedral} if it is equal to the
intersection of finitely many halfspaces.
\indexg{finitely generated (set)|(}%
The set is \emph{finitely generated} if there exist vectors
$\vv_1,\ldots,\vv_m\in\Rn$ and $k\in\{1,\ldots,m\}$
such that
\begin{equation}   \label{eq:fin-gen-form}
   C
   =
   \Braces{ \sum_{i=1}^m \lambda_i \vv_i :\:
            \lambda_1,\ldots,\lambda_m\in\Rpos,\,
            \sum_{i=1}^k \lambda_i = 1 }.
\end{equation}
(Note in this expression that all $m$ coefficients $\lambda_i$ must be
nonnegative, but only the first $k$ need add up to $1$.)
A convex set $C$ is both bounded and finitely generated
if and only if it is a \emph{polytope}, that is, the convex hull of
finitely many points.

\begin{proposition}   \label{roc:thm19.1}
  Let $C\subseteq\Rn$ be convex.
  Then the following are equivalent:
  \begin{letter-compact}
  \item   \label{roc:thm19.1:a}
    $C$ is polyhedral.
  \item   \label{roc:thm19.1:b}
    $C$ is finitely generated.
  \item   \label{roc:thm19.1:c}
    $C$ is closed (in $\Rn$) and has finitely many faces.
  \end{letter-compact}
\end{proposition}

\begin{proof}
See \idxroc\citet[Theorem~19.1]{ROC}.
\end{proof}

\begin{proposition}  \label{pr:fin-gen-cvx-cone}
  Let $K\subseteq\Rn$.
  Then the following are equivalent:
  \begin{letter-compact}
  \item  \label{pr:fin-gen-cvx-cone:a}
    $K$ is a finitely generated
    (or equivalently, polyhedral)
    convex cone.
  \item  \label{pr:fin-gen-cvx-cone:b}
    $K=\cone{V}$ for some finite $V\subseteq\Rn$.
  \item  \label{pr:fin-gen-cvx-cone:c}
    $K$ is the intersection of finitely many homogeneous closed
    halfspaces; that is,
    \[
      K
      =
      \Braces{
        \xx\in\Rn
        :\:
        \xx\cdot\uu_i\leq 0
        \textup{ for } i=1,\ldots,m
      }
    \]
    for some $\uu_1,\ldots,\uu_m\in\Rn\wo\{\zero\}$.
  \end{letter-compact}
\end{proposition}

\begin{proof}
  See discussion at beginning of
  \idxroc\citet[Section~19]{ROC}.%
\indexg{finitely generated (set)|)}
\end{proof}

\begin{proposition}   \label{roc:cor19.2.2}
  Let $K\subseteq\Rn$ be a polyhedral convex cone.
  Then $\Kpol$ is also polyhedral.
\end{proposition}

\begin{proof}
See \idxroc\citet[Corollary~19.2.2]{ROC}.
\end{proof}

\begin{proposition}   \label{roc:thm19.6}
  Let $C_1,\ldots,C_m\subseteq\Rn$ be convex and polyhedral.
  Then
  \[  \cl\BiggParens{\conv \bigcup_{i=1}^m C_i}  \]
  is also polyhedral.
\end{proposition}

\begin{proof}
See \idxroc\citet[Theorem~19.6]{ROC}.%
\indexg{polyhedral sets (standard)|)}%
\end{proof}

\section{Convex functions}
\label{sec:prelim:cvx-fcns}

Let $f:\Rn\rightarrow\Rext$.
The epigraph of $f$,
as defined in \eqref{eqn:epi-def},
is a subset of $\Rn\times\R$, which
can thus be viewed as a subset of $\R^{n+1}$.
We say that $f$ is \emph{convex}%
\indexg{convex functions (standard)!defined}
if its epigraph, $\epi{f}$, is
a convex set in $\R^{n+1}$.
The function is \emph{concave}%
\indexg{concave functions}
if $-f$ is convex,
or equivalently, if $f$'s hypograph
(defined in Eq.~\ref{eqn:hypo-def})
is convex as a subset of $\R^{n+1}$.
The function is \emph{finite everywhere}%
\indexg{finite everywhere (function)}
if $f>-\infty$ and $f<+\infty$.
The function
is \emph{proper}%
\indexg{proper functions}
if $f > -\infty$ and $f\not\equiv+\infty$.
The function is \emph{positively homogeneous}%
\indexg{positively homogeneous functions}
if $f(\lambda\xx)=\lambda f(\xx)$ for all $\xx\in\Rn$ and all
$\lambda\in\Rstrictpos$.

\indexg{affine functions (standard)|(}%
Every affine
function (of the form $\xx\mapsto\xx\cdot\uu+\beta$ for
$\xx\in\Rn$,
where $\uu\in\Rn$ and $\beta\in\R$) is both convex and concave.%
\indexg{affine functions (standard)|)}%

\indexg{convex functions (standard)!sublevel sets of|(}%
\indexg{convex functions (standard)!domain|(}%
If $f$ is convex, then its domain, $\dom{f}$, is also convex, as are
all its sublevel sets:

\begin{proposition}   \label{roc:thm4.6}
  Let $f:\Rn\rightarrow\Rext$ be convex, and let $\alpha\in\Rext$.
  Then the sets
  \[
     \{\xx\in\Rn :\: f(\xx)\leq \alpha\}
     \;\;\mbox{ and }\;\;
     \{\xx\in\Rn :\: f(\xx) < \alpha\}
  \]
  are both convex.
  In particular, $\dom{f}$ is convex.
\end{proposition}

\begin{proof}
See \idxroc\citet[Theorem~4.6]{ROC}.%
\indexg{convex functions (standard)!domain|)}%
\indexg{convex functions (standard)!sublevel sets of|)}
\end{proof}

\indexg{convex functions (standard)!characterized|(}%
Here is a useful characterization for when a function is convex.
It is similar to other standard characterizations,
such as \idxroc\citet[Theorems~4.1 and~4.2]{ROC}.
The one given here is applicable even when the function is improper.

\begin{proposition}   \label{pr:stand-cvx-fcn-char}
  Let $f:\Rn\rightarrow\Rext$.
  Then the following are equivalent:
  \begin{letter-compact}
    \item
    \label{pr:stand-cvx-fcn-char:a}
      $f$ is convex.
    \item
    \label{pr:stand-cvx-fcn-char:b}
      For all $\xx_0,\xx_1\in\dom f$
      and all $\lambda\in [0,1]$,
      \begin{equation}
      \label{eq:stand-cvx-fcn-char:b}
        f\bigParens{(1-\lambda)\xx_0 + \lambda \xx_1}
        \leq
        (1-\lambda) f(\xx_0) + \lambda f(\xx_1).
      \end{equation}
  \end{letter-compact}
\end{proposition}

\begin{proof}
~
\begin{proof-parts}
\pfpart{%
  (\ref{pr:stand-cvx-fcn-char:a})
  $\Rightarrow$
  (\ref{pr:stand-cvx-fcn-char:b}):
}
Suppose $f$ is convex, and let $\xx_0,\xx_1\in\dom f$
and $\lambda\in [0,1]$.
We aim to show that
\eqref{eq:stand-cvx-fcn-char:b}
holds.
For $i\in\{0,1\}$,
let $y_i\in\R$ with $f(\xx_i)\leq y_i$
so that $\rpair{\xx_i}{y_i}\in\epi{f}$.
Let $\xx=(1-\lambda)\xx_0+\lambda\xx_1$
and $y=(1-\lambda)y_0+\lambda y_1$.
Since $f$ is convex, its epigraph is convex, so the point
\begin{equation}
\label{eq:pr:stand-cvx-fcn-char:2}
  \rpair{\xx}{y}
  =
  (1-\lambda)\rpair{\xx_0}{y_0}+\lambda\rpair{\xx_1}{y_1}
\end{equation}
is in $\epi{f}$ as well.
Hence, $f(\xx)\leq y$, that is,
$f((1-\lambda)\xx_0 + \lambda \xx_1)\leq (1-\lambda) y_0 + \lambda y_1$.
Since this holds for all $y_0\geq f(\xx_0)$ and $y_1\geq f(\xx_1)$,
this implies
\eqref{eq:stand-cvx-fcn-char:b}.

\pfpart{%
  (\ref{pr:stand-cvx-fcn-char:b})
  $\Rightarrow$
  (\ref{pr:stand-cvx-fcn-char:a}):
}
Suppose statement~(\ref{pr:stand-cvx-fcn-char:b}) holds.
Let $\rpair{\xx_i}{y_i}\in\epi{f}$, for $i\in\{0,1\}$,
implying $\xx_i\in\dom{f}$ and $f(\xx_i)\leq y_i$.
Let
$\lambda\in [0,1]$,
and let $\xx=(1-\lambda)\xx_0+\lambda\xx_1$
and $y=(1-\lambda)y_0+\lambda y_1$,
implying
\eqref{eq:pr:stand-cvx-fcn-char:2}.
Then, by our assumption,
\[
  f(\xx)
  \leq
  (1-\lambda)f(\xx_0) + \lambda f(\xx_1)
  \leq
  (1-\lambda) y_0 + \lambda y_1
  =
  y.
\]
Therefore, $\rpair{\xx}{y}\in\epi{f}$, so
$\epi{f}$ is convex, and so is $f$.%
\indexg{convex functions (standard)!characterized|)}
\qedhere
\end{proof-parts}
\end{proof}

\indexg{convex functions (standard)!strictly|(}%
\indexg{strictly convex functions|(}%
A function $f:\Rn\rightarrow\Rext$ is
\emph{strictly convex on $C$}, where $C$ is a convex subset of $\dom f$, if
\eqref{eq:stand-cvx-fcn-char:b} holds with strict inequality
for all $\lambda\in (0,1)$ and all $\xx_0,\xx_1\in C$ with
$\xx_0\neq\xx_1$.
The function is \emph{strictly convex}
if it is strictly convex on
\indexg{convex functions (standard)!strictly|)}%
\indexg{strictly convex functions|)}%
$\dom f$.

We next review some natural ways of constructing functions that are
convex.

\indexg{convex functions (standard)!constructing|(}%
\indexg{indicator functions (standard)|(}%
For a set $S\subseteq\Rn$,
the \emph{indicator function} $\inds:\Rn\rightarrow\Rext$
is defined, for $\xx\in\Rn$, by
\begin{equation} \label{eq:indf-defn}
\indexm{i s 200}{$\inds$}{indicator function (standard)}
  \inds(\xx) =
  \begin{cases}
    0 & \text{if $\xx\in S$,} \\
    +\infty
      & \text{otherwise.}
  \end{cases}
\end{equation}
If $S$ is a convex set, then $\inds$ is a convex function.
Also, when $S$ is a singleton $\{\zz\}$, we sometimes write
$\indf{\zz}$ as shorthand for $\indfsing{\zz}$,
and when working in $\R$, we sometimes write
$\indf{\leq\beta}$ for $\indf{(-\infty,\beta]}$, where $\beta\in\R$
(and similarly for other inequality relations).%
\indexg{indicator functions (standard)|)}

For a function $f:X\rightarrow\Rext$ and $\lambda\in\R$,
we define the function $\lambda f:X\rightarrow\Rext$ by
$(\lambda f)(x)=\lambda [f(x)]$ for $x\in X$.
If $f:\Rn\rightarrow\Rext$ is convex and $\lambda\in\Rpos$ then
$\lambda f$ is also convex.

It is sometimes useful to consider a slight variant in which
only function values for points in the effective domain are
scaled.
Thus, for a function $f:\Rn\rightarrow\Rext$ and
a nonnegative scalar $\lambda\in\Rpos$,
we define the function
$(\sfprod{\lambda}{f}):\Rn\rightarrow\Rext$,
for $\xx\in\Rn$, by
\begin{equation}
\label{eq:sfprod-defn}
\indexm{lambda f}{$\sfprod{\lambda}{f}$}{restricted scalar multiple}%
  (\sfprod{\lambda}{f})(\xx)
  =
  \begin{cases}
    \lambda f(\xx)
      & \text{if $f(\xx) < +\infty$,}
    \\
    +\infty
      & \text{if $f(\xx) = +\infty$.}
  \end{cases}
\end{equation}
If $\lambda>0$ then $\sfprod{\lambda}{f}$ is the same as the standard scalar multiple
$\lambda f$.
However, if $\lambda = 0$ then $\sfprod{0}{f}$ zeros out only the set
$\dom{f}$, resulting in an indicator function on that set,
whereas $0 f\equiv 0$.
Thus,
\begin{equation}
\label{eq:sfprod-identity}
    \sfprod{\lambda}{f}
    =
    \begin{cases}
      \lambda f
        & \text{if $\lambda>0$,} \\
      \inddomf
        & \text{if $\lambda=0$.}
    \end{cases}
\end{equation}
If $f$ is convex, then $\sfprod{\lambda}{f}$ is as well, for $\lambda\in\Rpos$.

\indexg{sum of functions (standard)!convexity of|(}%
The sum of summable convex functions is also convex
(slightly generalizing \idxroc\citealp[Theorem~5.2]{ROC}):

\begin{proposition}  \label{pr:std-sum-fcns-cvx}
  Let $f:\Rn\rightarrow\Rext$ and $g:\Rn\rightarrow\Rext$ be convex
  and summable.
  Then $f+g$ is convex.
\end{proposition}

\begin{proof}
Let $h=f+g$,
let $\xx_0,\xx_1\in\dom{h}$, and let $\lambda\in[0,1]$.
Then $\xx_0$ and $\xx_1$ are also in $\dom{f}$ and $\dom{g}$, so,
by \Cref{pr:stand-cvx-fcn-char},
\eqref{eq:stand-cvx-fcn-char:b} holds for both $f$ and $g$.
Adding these inequalities and again applying that proposition yields
the claim.%
\indexg{sum of functions (standard)!convexity of|)}
\end{proof}

\indexg{linear map in composition with standard function|(}%
For a function $f:\Rm\rightarrow\Rext$ and matrix $\A\in\Rmn$,
the function $\fA:\Rn\rightarrow\Rext$ is defined, for $\xx\in\Rn$, by
\begin{equation}  \label{eq:fA-defn}
\indexm{f a 100}{$\fA$}{composition with linear map (standard)}%
  (\fA)(\xx)=f(\A\xx).
\end{equation}
Thus, $\fA$ is the composition of $f$ with the linear map associated with
$\A$.
\indexg{affine map in composition with standard function|(}%
If $f$ is convex, then so is $\fA$, as is $f$'s composition with any
affine map:

\begin{proposition}   \label{roc:thm5.7:fA}
  Let $f:\Rm\rightarrow\Rext$ be convex,
  let $\A\in\Rmn$, and let $\bb\in\Rm$.
  Let $h:\Rn\rightarrow\Rext$ be defined by
  $h(\xx)=f(\A\xx+\bb)$ for $\xx\in\Rn$.
  Then $h$ is convex.
\end{proposition}

\begin{proof}
This follows straightforwardly using
\Cref{pr:stand-cvx-fcn-char}.%
\indexg{linear map in composition with standard function|)}%
\indexg{affine map in composition with standard function|)}
\end{proof}

\indexg{linear image of standard function|(}%
For a function $f:\Rn\rightarrow\Rext$ and matrix $\A\in\Rmn$,
the function $\A f:\Rm\rightarrow\Rext$ is defined,
for $\xx\in\Rm$, by
\begin{equation}  \label{eq:lin-image-fcn-defn}
\indexm{a f 100}{$\A f$}{linear image of function (standard)}%
  (\A f)(\xx) =  \inf\,\bigBraces{f(\zz):\:\zz\in\Rn,\,\A\zz=\xx}.
\end{equation}
The function $\A f$ is called
the \emph{image of $f$ under $\A$}.

\begin{proposition}   \label{roc:thm5.7:Af}
  Let $f:\Rn\rightarrow\Rext$ be convex
  and let $\A\in\Rmn$.
  Then $\A f$ is convex.
\end{proposition}

\begin{proof}
See \idxroc\citet[Theorem~5.7]{ROC}.
\end{proof}

\begin{proposition}  \label{pr:std-lin-img-props}
  Let $f:\Rn\rightarrow\Rext$ and let $\A\in\Rmn$.
  \begin{letter-compact}
  \item   \label{pr:std-lin-img-props:inf-fA}
    $f\ktransA\geq\inf f$.
  \item   \label{pr:std-lin-img-props:a}
    $(\A f) \A \leq f$.
  \item   \label{pr:std-lin-img-props:inf-Af}
    $\inf (\A f) = \inf f$.
  \item   \label{pr:std-lin-img-props:b}
    Let $g:\Rn\rightarrow\Rext$ and suppose $f\leq g$.
    Then
    $\A f \leq \A g$
    and
    $f \ktransA \leq g \transA$.
  \item   \label{pr:std-lin-img-props:c}
    Let $\B\in\R^{k\times m}$.
    Then $\B(\A f) = (\B\A) f$.
  \end{letter-compact}
\end{proposition}

\begin{proof}
  ~
  
\begin{proof-parts}
\pfpart{%
  Parts~(\ref{pr:std-lin-img-props:inf-fA}),
  (\ref{pr:std-lin-img-props:a}),
  and (\ref{pr:std-lin-img-props:b}):
}
These follow directly
from definitions
(Eqs.~\ref{eq:fA-defn} and~\ref{eq:lin-image-fcn-defn}).

\pfpart{Part~(\ref{pr:std-lin-img-props:inf-Af}):}
From $\A f$'s definition
(Eq.~\ref{eq:lin-image-fcn-defn}), we have
$\A f\geq\inf f$, so $\inf(\A f)\geq\inf f$.
For the reverse inequality, we have
$f\geq (\A f) \A \geq \inf (\A f)$,
with the first inequality from 
part~(\ref{pr:std-lin-img-props:a}),
and the second from
part~(\ref{pr:std-lin-img-props:inf-fA})
(applied with $f$ and $\A$ set to
$\A f$ and $\transA$).
Thus, $\inf f\geq\inf(\A f)$.

\pfpart{Part~(\ref{pr:std-lin-img-props:c}):}
For $\xx\in\R^{k}$,
\begin{align*}
  \bigParens{\B(\A f)}(\xx)
  &=
  \inf\,\bigBraces{ (\A f)(\yy) :\: \yy\in\Rm,\, \B\yy=\xx }
  \\
  &=
  \inf\,\bigBraces{
    f(\zz)
    :\:
    \zz\in\Rn,\, \yy\in\Rm,\, \A\zz=\yy,\, \B\yy=\xx
  }
  \\
  &=
  \inf\,\bigBraces{ f(\zz) :\: \zz\in\Rn,\, \B\A\zz=\xx }
  =
  \bigParens{(\B\A) f}(\xx),
\end{align*}
where the first, second and last equalities are by definition of
image of a function
(Eq.~\ref{eq:lin-image-fcn-defn}).%
\indexg{linear image of standard function|)}
\qedhere
\end{proof-parts}
\end{proof}

\begin{proposition}   \label{roc:thm5.5}
\indexg{suprema and maxima, pointwise!standard convexity of|(}%
  Let $f_{i}:\Rn\rightarrow\Rext$ be convex
  for all $i\in\indset$, where
  $\indset$ is any index set.
  Let $h=\sup_{i\in I} f_i$ be their pointwise supremum.
  Then $h$ is convex.
\end{proposition}

\begin{proof}
See \idxroc\citet[Theorem~5.5]{ROC}.%
\indexg{suprema and maxima, pointwise!standard convexity of|)}
\end{proof}

\indexg{nondecreasing function, composition with!standard convexity of|(}%
Let $g:\R\to\eR$ be nondecreasing. We say that $G:\eR\to\eR$ is a
\emph{monotone extension}%
\indexg{monotone extension}
of~$g$ if $G$ is nondecreasing and agrees with $g$ on $\R$,
that is, if $G(-\infty)\le\inf g$,
$G(x)=g(x)$ for $x\in\R$, and $G(+\infty)\ge\sup g$.

The next proposition shows that
a function obtained by composing a convex function
$f:\Rn\to\eR$ with a nondecreasing convex function
$g:\R\to\Rext$,
appropriately extended to $\eR$, is also convex.
This slightly generalizes \idxroc\citet[Theorem~5.1]{ROC}
by now allowing both $f$ and $g$
to include $-\infty$ in their ranges.

\begin{proposition}
\label{prop:nondec:convex}
  Let $f:\Rn\rightarrow\Rext$ be convex, let
  $g:\R\to\eR$ be convex and nondecreasing,
  and let $G:\eR\to\eR$ be a monotone extension of $g$ such that $G(+\infty)=+\infty$.
  Then the function $h=G\circ f$ is convex.
\end{proposition}

\begin{proof}
Let $\xx_0,\xx_1\in\dom{h}$ and let $\lambda\in[0,1]$.
We will show
\begin{equation}   \label{eq:prop:nondec:convex:1}
  h\bigParens{(1-\lambda)\xx_0 + \lambda \xx_1}
  \leq
  (1-\lambda) h(\xx_0) + \lambda h(\xx_1),
\end{equation}
and so that $h$ satisfies condition (\ref{pr:stand-cvx-fcn-char:b}) of
\Cref{pr:stand-cvx-fcn-char}.
\eqref{eq:prop:nondec:convex:1} holds trivially if $\lambda\in\set{0,1}$,
so we assume henceforth that $\lambda\in(0,1)$.
Since $G(+\infty)=+\infty$
and
$\xx_0,\xx_1\in\dom f$, we have
\begin{align}
  h\bigParens{(1-\lambda)\xx_0+\lambda\xx_1}
&=
  G\bigParens{f\bigParens{(1-\lambda)\xx_0+\lambda\xx_1}}
\nonumber
\\
&\le
  G\bigParens{(1-\lambda)f(\xx_0)+\lambda f(\xx_1)},
\label{eq:prop:nondec:convex:2}
\end{align}
where the inequality follows by monotonicity of $G$ and convexity of $f$ (using \Cref{pr:stand-cvx-fcn-char}\ref{pr:stand-cvx-fcn-char:b}). Now, if either $f(\xx_0)=-\infty$
or $f(\xx_1)=-\infty$, then
\begin{align*}
  G\bigParens{(1-\lambda)f(\xx_0)+\lambda f(\xx_1)}
=
  G(-\infty)
&=
  (1-\lambda)G(-\infty)+\lambda G(-\infty)
\\
&\le
  (1-\lambda)G\bigParens{f(\xx_0)}+\lambda G\bigParens{f(\xx_1)},
\end{align*}
where the inequality follows by monotonicity of $G$.
Combined with \eqref{eq:prop:nondec:convex:2} and since
$h=G\circ f$, this proves
\eqref{eq:prop:nondec:convex:1} in this case.

In the remaining case, $f(\xx_0)$ and $f(\xx_1)$ are in $\R$,
implying $g(f(\xx_i))=G(f(\xx_i))=h(\xx_i)<+\infty$ for $i\in\{0,1\}$,
and so that $f(\xx_i)\in\dom{g}$.
Thus,
\begin{align*}
  G\bigParens{(1-\lambda)f(\xx_0)+\lambda f(\xx_1)}
&=
  g\bigParens{(1-\lambda)f(\xx_0)+\lambda f(\xx_1)}
\\
&\le
  (1-\lambda)g\bigParens{f(\xx_0)}+\lambda g\bigParens{f(\xx_1)},
\end{align*}
where the inequality follows by convexity of $g$ (using \Cref{pr:stand-cvx-fcn-char}\ref{pr:stand-cvx-fcn-char:b}).
Again combining with \eqref{eq:prop:nondec:convex:2}, this proves
\eqref{eq:prop:nondec:convex:1},
completing the proof.%
\indexg{nondecreasing function, composition with!standard convexity of|)}%
\indexg{convex functions (standard)!constructing|)}
\end{proof}

Here are some facts about relative interiors as they relate to a
function's domain and epigraph:

\begin{proposition}  \label{roc:lem7.3}
\indexg{epigraph!relative interior of|(}%
  Let $f:\Rn\rightarrow\Rext$ be convex.
  Then $\ri(\epi{f})$ consists of all pairs
  $\rpair{\xx}{y}\in\Rn\times\R$ with
  $\xx\in\ri(\dom{f})$ and $y>f(\xx)$.
\end{proposition}

\begin{proof}
See \idxroc\citet[Lemma~7.3]{ROC}.%
\indexg{epigraph!relative interior of|)}
\end{proof}

\begin{proposition}   \label{pr:Ax-ridomf-ridomfA}
\indexg{affine map in composition with standard function!relative interior of domain|(}%
  Let $f:\Rm\rightarrow\Rext$ be convex,
  let $\A\in\Rmn$, and let $\bb\in\Rm$.
  Let $h:\Rn\rightarrow\Rext$ be defined by
  $h(\xx)=f(\A\xx+\bb)$ for $\xx\in\Rn$.
  Suppose there exists $\xing\in\Rn$ for which
  $\A\xing+\bb\in\ri(\dom f)$.
  Then for all $\xx\in\Rn$,
  $\A\xx+\bb\in\ri(\dom f)$
  if and only if
  $\xx\in\ri(\dom h)$.
  In particular, $\xing\in\ri(\dom h)$.
\end{proposition}

\begin{proof}
Let $g:\Rm\rightarrow\Rext$ be defined by $g(\xx)=f(\xx+\bb)$ for
$\xx\in\Rn$.
Then $\dom g = (\dom f) - \bb$ so
\begin{equation}   \label{eq:pr:Ax-ridomf-ridomfA:1}
  \ri(\dom g)=\ri(\dom f) + \ri\{-\bb\}=\ri(\dom f) - \bb,
\end{equation}
with the first equality from
\Cref{pr:ri-props}(\ref{pr:ri-props:roc-cor6.6.2})
(with $C=\dom f$ and $D=\{-\bb\}$)
and the second because
$\ri\{-\bb\}=\{-\bb\}$.

Let $\slinmapA:\Rn\rightarrow\Rm$ be the linear
map associated with $\A$, meaning
$\slinmapA(\xx)=\A\xx$ for $\xx\in\Rn$.
We then have
\begin{equation}   \label{eq:pr:Ax-ridomf-ridomfA:2}
  \ri\regParens{\dom h}
  =
  \ri\bigParens{\slinmapAinv(\dom{g})}
  =
  \slinmapAinv\bigParens{\ri(\dom{g})}
  =
  \slinmapAinv\bigParens{\ri(\dom{f}) - \bb}.
\end{equation}
The first equality is because
$\dom{h}=\{\xx\in\Rn :\: \A\xx\in\dom{g}\} = \slinmapAinv(\dom{g})$.
The second equality is by
\Cref{roc:thm6.7} (with $C=\dom{g}$),
noting that
$\A\xing\in\ri(\dom g)$ by 
\eqref{eq:pr:Ax-ridomf-ridomfA:1}
since $\A\xing+\bb\in\ri(\dom f)$.
The last equality is also by
\eqref{eq:pr:Ax-ridomf-ridomfA:1}.
Since, for all $\xx\in\Rn$,
$\A\xx+\bb\in\ri(\dom f)$ if and only if $\xx$ is included in the set
on the right-hand side of
\eqref{eq:pr:Ax-ridomf-ridomfA:2},
this proves the claim.%
\indexg{affine map in composition with standard function!relative interior of domain|)}%
\end{proof}

\begin{proposition}  \label{pr:Ax-in-ri-dom-Af}
\indexg{linear image of standard function!relative interior of domain|(}%
  Let $f:\Rn\rightarrow\Rext$ be convex with $f\not\equiv+\infty$,
  and let $\A\in\Rmn$.
  Then there exists $\xhat\in\Rn$ such that $\A\xhat\in\ri(\dom \A f)$.
\end{proposition}

\begin{proof}
By assumption, there exists some point $\xx\in\dom{f}$,
implying $(\A f)(\A\xx)\leq f(\xx)<+\infty$
by \Cref{pr:std-lin-img-props}(\ref{pr:std-lin-img-props:a}).
Thus, $\dom{\A f}$ is nonempty, and
is also convex, being the domain of a convex function
(Propositions~\ref{roc:thm4.6} and~\ref{roc:thm5.7:Af}).
Therefore, there exists a point $\zz$ in its relative interior,
$\ri(\dom{\A f})$
(\Cref{pr:ri-props}\ref{pr:ri-props:roc-thm6.2b}).
Since $(\A f)(\zz)<+\infty$, this further implies
$\zz=\A\xhat$ for some $\xhat\in\Rn$
(from $\A f$'s definition in Eq.~\ref{eq:lin-image-fcn-defn}).
This proves the claim.%
\indexg{linear image of standard function!relative interior of domain|)}%
\end{proof}

\section{Lower semicontinuity and continuity}
\label{sec:prelim:lsc}

As was defined more generally in \Cref{sec:prelim:lower-semicont},
a function $f:\Rn\rightarrow\Rext$ is
{lower semicontinuous at a point $\xx\in\Rn$} if
\[
  \liminf f(\xx_t)\geq f(\xx)
\]
for every sequence $\seq{\xx_t}$ in $\Rn$ that converges to $\xx$.
The function is
{lower semicontinuous}
if it is lower semicontinuous at every point $\xx\in\Rn$.
Lower semicontinuity of $f$ is equivalent to the epigraph of $f$ being
closed as a subset of $\R^{n+1}$
(see \Cref{prop:lsc}).

\indexg{lower semicontinuity!convex functions@of convex functions|(}%
If $f$ is lower semicontinuous, convex and improper, then it is
infinite at every point:

\begin{proposition}   \label{pr:improper-vals}
  Let $f:\Rn\rightarrow\Rext$ be convex and improper.
  \begin{letter-compact}
  \item   \label{pr:improper-vals:thm7.2}
    For all $\xx\in\ri(\dom{f})$,
    $f(\xx)=-\infty$.
  \item   \label{pr:improper-vals:cor7.2.1}
    Suppose, in addition, that $f$ is lower semicontinuous.
    Then $f(\xx)\in\{-\infty,+\infty\}$ for all $\xx\in\Rn$.
  \end{letter-compact}
\end{proposition}

\begin{proof}
~

\begin{proof-parts}
\pfpart{Part~(\ref{pr:improper-vals:thm7.2}):}
See \idxroc\citet[Theorem~7.2]{ROC}.

\pfpart{Part~(\ref{pr:improper-vals:cor7.2.1}):}
See \idxroc\citet[Corollary~7.2.1]{ROC}.%
\indexg{lower semicontinuity!convex functions@of convex functions|)}
\qedhere
\end{proof-parts}
\end{proof}

\indexg{lower semicontinuous hull!convex function@of convex function|(}%
\indexg{convex functions (standard)!lower semicontinuous hull of|(}%
As a special case of the more general definition
from \Cref{sec:prelim:lower-semicont},
the {lower semicontinuous hull} of a function
$f:\Rn\rightarrow\Rext$ is the function
$(\lsc f):\Rn\rightarrow\Rext$
defined, for $\xx\in\Rn$, by
\begin{equation*}  %
   (\lsc f)(\xx)
    = \InfseqLiminf{\seq{\xx_t}}{\Rn}{\xx_t\rightarrow \xx}{f(\xx_t)},
\end{equation*}
where the infimum is over all sequences $\seq{\xx_t}$ in $\Rn$
converging to $\xx$.
As we saw in \Cref{prop:lsc:characterize},
$\lsc f$ is the greatest lower semicontinuous function
that is majorized by $f$, and its epigraph is the closure of
the epigraph of $f$ in $\R^{n+1}$.

\indexg{closure (of convex function)|(}%
\indexg{convex functions (standard)!closure of|(}%
For a convex function $f:\Rn\rightarrow\Rext$, its \emph{closure}
$(\cl f):\Rn\rightarrow\Rext$%
\indexm{cl f}{$\cl f$}{closure (of function)}
 is defined to be the same as its lower
semicontinuous hull
if $f>-\infty$,
and is identically $-\infty$ otherwise;
that is, $\cl f=\lsc f$ if $f>-\infty$, and
$\cl f\equiv-\infty$ if $f(\xx)=-\infty$ at any point $\xx\in\Rn$.
A convex function $f$ is \emph{closed}%
\indexg{closed (convex function)}%
\indexg{convex functions (standard)!closed}
 if $f=\cl f$, that is,
if it is lower
semicontinuous and either $f > -\infty$ or $f\equiv-\infty$.
Thus, if $f$ is proper, then it is closed if and only if it is
lower semicontinuous;
if it is improper, then it is closed if and only if either
$f\equiv+\infty$ or $f\equiv-\infty$.

This definition of a closed convex function follows
\idxroc\citet[Section~7]{ROC}, but note that there is not full
agreement in the literature on this terminology
(see, for instance,
\indexa{Bertsekas, D. P.}\citealp[Section~1.1.2]{Bertsekas2009}
or
\idxborlew\citealp[Section~4.2]{borwein_lewis_06}).

\begin{proposition}  \label{pr:lsc-props}
  Let $f:\Rn\rightarrow\Rext$ be convex.
  \begin{letter-compact}
  \item  \label{pr:lsc-props:a}
    $\lsc f$ is convex and lower semicontinuous.
  \item  \label{pr:lsc-props:b}
    If $\xx\in\Rn$ is not a relative boundary point of $\dom{f}$,
    then $(\lsc f)(\xx)=f(\xx)$.
  \item  \label{pr:lsc-props:c}
    $\dom(\lsc f)$ and $\dom{f}$ have the same closure and relative
    interior.
  \item  \label{pr:lsc-props:d}
    If, in addition, $f$ is proper, then
    $\cl f (=\lsc f)$ is closed and proper.
  \item \label{pr:lsc-props:e}
    If $f$ is lower semicontinuous, but not closed, then it
    takes on only the values $+\infty$ and $-\infty$,
    but is not identically equal to either.
  \end{letter-compact}
\end{proposition}

\begin{proof}
~

\begin{proof-parts}
\pfpart{Part~(\ref{pr:lsc-props:a}):}
That $\lsc f$ is lower semicontinuous follows from
\Cref{prop:lsc:characterize}(\ref{prop:lsc:characterize:a}).

By definition, a function is convex if and only if its epigraph is
convex, so $\epi{f}$ is convex, and so also
is its closure
(by \Cref{pr:ri-props}\ref{pr:ri-props:roc-thm6.2a}),
which is the same as $\epi(\lsc f)$
by
\Cref{prop:lsc:characterize}(\ref{prop:lsc:characterize:c}).
Therefore, $\lsc f$ is convex.

\pfpart{Part~(\ref{pr:lsc-props:b}):}
If $f$ is proper, then this follows from
\idxroc\citet[Theorem~7.4]{ROC}.

Suppose instead that $f$ is not proper and consider any $\xx\in\Rn$
that is not a relative boundary point of $\dom f$.
If $\xx\in\ri(\dom f)$, then
$f(\xx)=-\infty$
by \Cref{pr:improper-vals}(\ref{pr:improper-vals:thm7.2}),
implying $(\lsc f)(\xx)=-\infty$ as well,
by
\Cref{prop:lsc:characterize}(\ref{prop:lsc:characterize:a}).

Otherwise, we must have $\xx\not\in\cl(\dom f)$,
implying $f(\xx)=+\infty$ and further that there
must exist a neighborhood $U\subseteq\Rn$ of $\xx$ that is disjoint
from $\dom f$.
Let $\seq{\xx_t}$ be any sequence in $\Rn$ converging to $\xx$.
Then all but finitely many elements of the sequence must be in $U$,
being a neighborhood of $\xx$, and none of these are in $\dom{f}$.
Thus, $f(\xx_t)\rightarrow+\infty$.
Since this holds for all such sequences, this shows that
$(\lsc f)(\xx)=+\infty=f(\xx)$.

\pfpart{Part~(\ref{pr:lsc-props:c}):}
Let $\xx\in\Rn$.
If $\xx\in\ri(\dom f)$ then by
part~(\ref{pr:lsc-props:b}),
$(\lsc f)(\xx)=f(\xx)<+\infty$.
Similarly, if $\xx\not\in\cl(\dom f)$ then
$(\lsc f)(\xx)=f(\xx)=+\infty$.
Thus,
\[
   \ri(\dom f)
   \subseteq
   \dom(\lsc f)
   \subseteq
   \cl(\dom f).
\]
The claim then follows from
\Cref{pr:ri-props}(\ref{pr:ri-props:roc-cor6.3.1}).

\pfpart{Part~(\ref{pr:lsc-props:d}):}
See \idxroc\citet[Theorem~7.4]{ROC}.
\pfpart{Part~(\ref{pr:lsc-props:e}):}
By definition of closedness, since $f$ is lower semicontinuous,
it must be improper. Therefore, it only takes on the values of $+\infty$ and
$-\infty$
(by \Cref{pr:improper-vals}\ref{pr:improper-vals:cor7.2.1}), but is not identically equal to either
(which would mean that it is closed).%
\indexg{lower semicontinuous hull!convex function@of convex function|)}%
\indexg{convex functions (standard)!lower semicontinuous hull of|)}%
\indexg{closure (of convex function)|)}%
\indexg{convex functions (standard)!closure of|)}
\qedhere
\end{proof-parts}
\end{proof}

\begin{proposition}  \label{roc:thm7.6-mod}
\indexg{convex functions (standard)!sublevel sets of|(}%
\indexg{sublevel sets!standard convex function@of standard convex function|(}%
  Let $f:\Rn\rightarrow\Rext$ be convex,
  and let $\beta\in\R$ with $\inf f < \beta$.
  Let
  \begin{align*}
    L &= \Braces{\xx\in\Rn :\: f(\xx) \leq \beta}, \\
    M &= \Braces{\xx\in\Rn :\: f(\xx) < \beta}.
  \end{align*}
  Then
  \begin{align*}
    \cl L &= \cl M
       = \Braces{\xx\in\Rn :\: (\lsc f)(\xx) \leq \beta},
    \\
    \ri L &= \ri M
       = \Braces{\xx\in \ri(\dom f) :\: f(\xx) < \beta}.
  \end{align*}
\end{proposition}

\begin{proof}
This follows by the same proof as in
\idxroc\citet[Theorem~7.6]{ROC} with $\cl f$ replaced by
$\lsc f$.%
\indexg{convex functions (standard)!sublevel sets of|)}%
\indexg{sublevel sets!standard convex function@of standard convex function|)}
\end{proof}

\indexg{convex functions (standard)!continuity of|(}%
\indexg{continuity!standard convex functions@of standard convex functions|(}%
By \Cref{prop:first:properties}(\ref{prop:first:cont})
(and since $\Rn$ is first-countable),
a function $f:\Rn\rightarrow\Rext$ is continuous at a point
$\xx\in\Rn$ if and only if $f(\xx_t)\rightarrow f(\xx)$ for every sequence
$\seq{\xx_t}$ in $\Rn$ converging to $\xx$.
The function is continuous (or continuous everywhere) if it is
continuous at every point $\xx\in\Rn$.
Convex functions are continuous at all points in $\Rn$ except possibly
the boundary points of $\dom{f}$:

\begin{proposition}   \label{pr:stand-cvx-cont}
  Let $f:\Rn\rightarrow\Rext$ be convex,
  and let $\xx\in\Rn$.
  If $\xx$ is not a boundary point of $\dom{f}$, then
  $f$ is continuous at $\xx$.
  Consequently, if $f$ is finite everywhere, then $f$ is continuous
  everywhere.
\end{proposition}

\begin{proof}
Suppose $\xx$ is not a boundary point of $\dom{f}$.
Let $\seq{\xx_t}$ be any sequence in $\Rn$ converging to $\xx$.
We aim to show $f(\xx_t)\rightarrow f(\xx)$.

As a first case, suppose $\xx$ is in $\intdom{f}$, which so must be equal to
$\ri(\dom{f})$ by
\Cref{pr:ri-props}(\ref{pr:ri-props:intC-nonemp-implies-eq-riC}).
Since $\intdom{f}=\ri(\dom{f})$ is a neighborhood of $\xx$, it must
include all but finitely many sequence elements $\xx_t$.
Therefore, $f(\xx_t)\rightarrow f(\xx)$
by \idxroc\citet[Theorem~10.1]{ROC}.

In the alternative case, $\xx$ is in $\Rn\setminus\cl(\dom{f})$,
which thus is a neighborhood of $\xx$, implying that it includes all
but finitely many elements $\xx_t$, and
so that $f(\xx_t)=+\infty$ for all $t$ sufficiently large.
Thus, $f(\xx_t)\rightarrow +\infty=f(\xx)$.%
\indexg{convex functions (standard)!continuity of|)}%
\indexg{continuity!standard convex functions@of standard convex functions|)}
\end{proof}

\section{Conjugacy}
\label{sec:prelim:conjugate}

\indexg{conjugate (standard)|(}%
The \emph{conjugate} of a function $f:\Rn\rightarrow\Rext$ is the
function
$\fstar:\Rn\rightarrow\Rext$ defined, for $\uu\in\Rn$, by
\begin{equation}  \label{eq:fstar-def:intro}
\indexm{f 400}{$\fstar$}{conjugate (standard)}%
\indexg{conjugate (standard)!defined}%
  \fstar(\uu) = \sup_{\xx\in\Rn} \bigBracks{\xx\cdot\uu - f(\xx)}.
\end{equation}
We can interpret $\fstar$ as encoding, via its epigraph,
exactly those affine functions $\xx\mapsto \xx\cdot\uu - v$ that are
majorized by $f$:

\begin{proposition}  \label{pr:epi-fstar-maj-affn}
  Let $f:\Rn\rightarrow\Rext$ and let $\rpair{\uu}{v}\in\Rn\times\R$.
  Then $\rpair{\uu}{v}\in\epi{\fstar}$ if and only if
  $f(\xx)\geq \xx\cdot\uu - v$ for all $\xx\in\Rn$.
\end{proposition}

\begin{proof}
From \eqref{eq:fstar-def:intro},
the pair $\rpair{\uu}{v}$ is in
$\epi{\fstar}$, meaning $\fstar(\uu)\leq v$,
if and only if $\xx\cdot\uu - f(\xx)\leq v$ for all $\xx\in\Rn$,
or equivalently,
if and only if $f(\xx)\geq \xx\cdot\uu - v$
for all $\xx\in\Rn$.
\end{proof}

We write $\fdubs$ for the \emph{biconjugate}%
\indexg{biconjugate (standard)}
 of $f$,
the conjugate of its conjugate; that is,
$\fdubs=(\fstar)^*$.
From the preceding proposition,
it can be argued that the biconjugate of $f$ is equal to the pointwise
supremum of all affine functions that are majorized by $f$.

Here are more facts about conjugates:

\begin{proposition}  \label{pr:conj-props}
  Let $f:\Rn\rightarrow\Rext$.
  \begin{letter-compact}
  \item  \label{pr:conj-props:a}
    $\inf f = -\fstar(\zero)$.
  \item  \label{pr:conj-props:b}
    Let $g:\Rn\rightarrow\Rext$.
    If $f\geq g$ then $\fstar \leq \gstar$.
  \item  \label{pr:conj-props:c1}
    If $f\equiv+\infty$ then $\fstar\equiv-\infty$;
    otherwise,
    if $f\not\equiv+\infty$ then $\fstar>-\infty$.
  \item  \label{pr:conj-props:c2}
    If $f(\xx)=-\infty$ for some $\xx\in\Rn$ then $\fstar\equiv+\infty$.
  \item  \label{pr:conj-props:d}
    $\fstar$ is convex and closed.
  \item  \label{pr:conj-props:e}
    $\fstar = (\lsc{f})^* = (\cl{f})^*$.
  \end{letter-compact}
\end{proposition}

\begin{proof}
~

\begin{proof-parts}
\pfpart{Parts~(\ref{pr:conj-props:a}), (\ref{pr:conj-props:b}),
  (\ref{pr:conj-props:c1}), (\ref{pr:conj-props:c2}):}
These are straightforward from definitions.

\pfpart{Part~(\ref{pr:conj-props:d}):}
If $\fstar(\uu)=-\infty$ for some $\uu\in\Rn$, then
$\xx\cdot\uu-f(\xx) = -\infty$ for all $\xx\in\Rn$, implying
$f\equiv+\infty$ and so that $\fstar\equiv-\infty$
(by part~(\ref{pr:conj-props:c1})).
Therefore, $\fstar$ is convex and closed in this case.

Assume then that $\fstar>-\infty$.
For each $\xx\in\Rn$, let
$\psi_{\xx}:\Rn\rightarrow\Rext$ be defined by
$\psi_{\xx}(\uu)=\xx\cdot\uu-f(\xx)$ for $\uu\in\Rn$.
Each of these functions is convex and continuous,
including when $f(\xx)\in\{-\infty,+\infty\}$.
Therefore, $\fstar$, which is the pointwise supremum of all
$\psi_\xx$, is convex
(by \Cref{roc:thm5.5})
and lower semicontinuous
(by \Cref{pr:lsc-sup}).
Since $\fstar>-\infty$, this also shows that $\fstar$ is closed.

\pfpart{Part~(\ref{pr:conj-props:e}):}
We first show $\epi\fstar=\epi(\lsc f)^*$.
Let $\rpair{\uu}{v}\in\Rn\times\R$, and let
$h:\Rn\rightarrow\R$ be defined by
$h(\xx)=\xx\cdot\uu-v$
for $\xx\in\Rn$.
Then
\[
  \rpair{\uu}{v} \in \epi\fstar
  \;\;\Leftrightarrow\;\;
  h \leq f
  \;\;\Leftrightarrow\;\;
  h \leq \lsc f
  \;\;\Leftrightarrow\;\;
  \rpair{\uu}{v} \in \epi(\lsc f)^*.
\]  
The first and third equivalences are by
\Cref{pr:epi-fstar-maj-affn}
(applied, respectively, to $f$ and $\lsc f$).
The second equivalence is by
\Cref{prop:lsc:characterize}(\ref{prop:lsc:characterize:a},\ref{prop:lsc:characterize:b}),
since $h$ is continuous.
Thus, $\fstar$ and $(\lsc{f})^*$ have the same epigraphs, and
therefore are equal.

If $f>-\infty$, then $\cl{f}=\lsc{f}$, implying their conjugates are
also equal.
Otherwise, if $f\not>-\infty$, then $\cl{f}\equiv-\infty$ implying
both $\fstar$ and $(\cl f)^*$ are identically $+\infty$
(by part~(\ref{pr:conj-props:c2})),
and so are equal.
\qedhere
\end{proof-parts}
\end{proof}

\begin{proposition}  \label{pr:conj-props-cvx}
  Let $f:\Rn\rightarrow\Rext$ be convex.
  \begin{letter-compact}
  \item  \label{pr:conj-props-cvx:a}
    $f$ is proper if and only if  $\fstar$ is proper.
  \item  \label{pr:conj-props-cvx:b}
    $\fdubs=\cl{f}$.
  \item  \label{pr:conj-props-cvx:c}
    For $\uu\in\Rn$,
    \[
       \fstar(\uu)
       =
       \sup_{\xx\in\ri(\dom f)} [\xx\cdot\uu - f(\xx)].
    \]
  \end{letter-compact}
\end{proposition}

\begin{proof}
~

\begin{proof-parts}
\pfpart{%
  Parts~(\ref{pr:conj-props-cvx:a})
  and~(\ref{pr:conj-props-cvx:b}):
}
See \idxroc\citet[Theorem~12.2]{ROC}.

\pfpart{Part~(\ref{pr:conj-props-cvx:c}):}
See \idxroc\citet[Corollary~12.2.2]{ROC}.%
\indexg{conjugate (standard)|)}
\qedhere
\end{proof-parts}
\end{proof}

We next look at the conjugates of specific types of functions and functions
obtained by various functional operations.

As defined in \eqref{eq:indf-defn},
let $\inds$ be the indicator function for
some set $S\subseteq\Rn$.
\indexg{indicator functions (standard)!conjugate of|(}%
\indexg{support functions (standard)|(}%
\indexg{conjugate (standard)!indicator function@of indicator function|(}%
Its conjugate, $\indstars$, is called
the \emph{support function} for $S$, and is equal,
for $\uu\in\Rn$, to
\begin{equation} \label{eq:e:2}
  \indstars(\uu) = \sup_{\xx\in S} \xx\cdot\uu.
\end{equation}

\begin{proposition}  \label{roc:thm14.1-conj}
  Let $K\subseteq\Rn$ be a closed convex cone.
  Then $\indstar{K}=\indf{\Kpol}$
  and ${\indstar{\Kpol}=\indf{K}}$.
\end{proposition}

\begin{proof}
See \idxroc\citet[Theorem~14.1]{ROC}.%
\indexg{conjugate (standard)!indicator function@of indicator function|)}%
\indexg{indicator functions (standard)!conjugate of|)}%
\indexg{support functions (standard)|)}
\end{proof}

\begin{proposition}
\label{pr:conj-shift}
  Let $f:\Rn\to\eR$ and let $\bb\in\Rn$.
\begin{letter-compact}
  \item \label{pr:conj-shift:a}
    Let
    $g(\xx)=f(\xx+\bb)$ for $\xx\in\Rn$.
    Then $\gstar(\uu)=\fstar(\uu)-\bb\inprod\uu$ for $\uu\in\Rn$.
  \item \label{pr:conj-shift:b}
    Let
    $h(\xx)=f(\xx)-\xx\inprod\bb$ for $\xx\in\Rn$.
    Then $\hstar(\uu)=\fstar(\uu+\bb)$ for $\uu\in\Rn$.
\end{letter-compact}
\end{proposition}

\begin{proof}
Let $\uu\in\Rn$.

\begin{proof-parts}
\pfpart{Part~(\ref{pr:conj-shift:a}):}
  From the definition of conjugate (Eq.~\ref{eq:fstar-def:intro}),
  \begin{align*}
  \gstar(\uu)
  &= \sup_{\xx\in\Rn} \bigBracks{\xx\cdot\uu - f(\xx+\bb)}
  \\
  &= \sup_{\xx\in\Rn} \bigBracks{(\xx+\bb)\cdot\uu - f(\xx+\bb)} - \bb\cdot\uu
   = \fstar(\uu) - \bb\cdot\uu.
\end{align*}

\pfpart{Part~(\ref{pr:conj-shift:b}):}
  Similarly,
\begin{align*}
  \hstar(\uu)
  &=
  \sup_{\xx\in\Rn} \bigBracks{\xx\cdot\uu - \bigParens{f(\xx)-\xx\cdot\bb}}
\\
  &=
  \sup_{\xx\in\Rn} \bigBracks{\xx\cdot(\uu+\bb) - f(\xx)}
  =
  \fstar(\uu+\bb).
\qedhere
\end{align*}
\end{proof-parts}
\end{proof}

\begin{proposition}  \label{roc:thm16.4}
\indexg{conjugate (standard)!sum of functions@of sum of functions|(}%
\indexg{sum of functions (standard)!conjugate of|(}%
  Let $f_i:\Rn\rightarrow\Rext$ be convex and proper,
  for $i=1,\ldots,m$.
  Assume $\bigcap_{i=1}^m \ri(\dom{f})\neq\emptyset$.
  Let $\uu\in\Rn$.
  Then
  \[
     (f_1+\dotsb+f_m)^*(\uu)
     =
     \inf\Braces{\sum_{i=1}^m f^*_i(\uu_i):\:
                 \uu_1,\ldots,\uu_m\in\Rn,\,
                 \sum_{i=1}^m \uu_i = \uu}.
  \]
  Furthermore, the infimum is attained.
\end{proposition}

\begin{proof}
See \idxroc\citet[Theorem~16.4]{ROC}.%
\indexg{sum of functions (standard)!conjugate of|)}%
\indexg{conjugate (standard)!sum of functions@of sum of functions|)}%
\end{proof}

\begin{proposition}   \label{roc:thm16.3:fA}
\indexg{conjugate (standard)!composition with linear map@of composition with linear map|(}%
\indexg{linear map in composition with standard function!standard conjugate of|(}%
  Let $f:\Rm\rightarrow\Rext$ be convex,
  and let $\A\in\Rmn$.
  Assume there exists $\xing\in\Rn$ such that
  $\A\xing\in\ri(\dom f)$.
  Let $\uu\in\Rn$.
  Then
  \[
     (\fA)^*(\uu)
     =
     \inf\Braces{ \fstar(\ww):\: \ww\in\Rm,\,\transA\ww = \uu }.
  \]
  Furthermore, if $\uu\in\colspace{\transA}$, then
  the infimum is attained.
  (Otherwise, if $\uu\not\in\colspace{\transA}$, then the infimum
  is vacuously equal to $+\infty$.)
  Thus, $(\fA)^*=\transAk\fstar$.
\end{proposition}

\begin{proof}
See \idxroc\citet[Theorem~16.3]{ROC}.%
\indexg{conjugate (standard)!composition with linear map@of composition with linear map|)}%
\indexg{linear map in composition with standard function!standard conjugate of|)}
\end{proof}

\begin{proposition}   \label{roc:thm16.3:Af}
\indexg{conjugate (standard)!linear image@of linear image|(}%
\indexg{linear image of standard function!standard conjugate of|(}%
  Let $f:\Rn\rightarrow\Rext$ be convex
  and let $\A\in\Rmn$.
  Then
  \[ (\A f)^* = \fstar \ktransA. \]
\end{proposition}

\begin{proof}
See \idxroc\citet[Theorem~16.3]{ROC}.%
\indexg{conjugate (standard)!linear image@of linear image|)}%
\indexg{linear image of standard function!standard conjugate of|)}
\end{proof}

\begin{proposition}   \label{pr:prelim-conj-max}
\indexg{conjugate (standard)!pointwise maximum@of pointwise maximum|(}%
\indexg{suprema and maxima, pointwise!conjugate of|(}%
  Let $f_i:\Rn\rightarrow\R$ be convex, for $i=1,\dotsc,m$,
  and let
  $  h=\max\Braces{f_1,\dotsc,f_m} $.
  Let $\uu\in\Rn$.
  Then
  \[
     \hstar(\uu)
     =
     \inf\Braces{ \sum_{i=1}^m \lambda_i \fistar(\uu_i):\:
       \sum_{i=1}^m \lambda_i \uu_i = \uu
     },
  \]
  where the infimum is over all $\uu_1,\ldots,\uu_m\in\Rn$ and all
  nonnegative $\lambda_1,\ldots,\lambda_m$ which sum to~$1$,
  subject to $\sum_{i=1}^m \lambda_i \uu_i = \uu$.
  Furthermore, the infimum is attained.
\end{proposition}

\begin{proof}
See \idxhiriart\citet[Theorem~X.2.4.7]{HULL-big-v2}.%
\indexg{conjugate (standard)!pointwise maximum@of pointwise maximum|)}%
\indexg{suprema and maxima, pointwise!conjugate of|)}%
\end{proof}

\indexg{conjugate (standard)!composition with nondecreasing function@of composition with nondecreasing function|(}%
\indexg{nondecreasing function, composition with!standard conjugate of|(}%
We next give an expression for the conjugate of
the composition of a nondecreasing convex function
on $\R$, suitably extended to $\Rext$, with a convex function
$f:\Rn\rightarrow\Rext$.
This expression is in terms of $\indepifstar$,
the support function for the epigraph of $f$, for which
an explicit expression in terms of $\fstar$ is given below in
\Cref{pr:support-epi-f-conjugate}.
This material, including the proof, is based closely on
\idxhiriart\citet[Theorem~X.2.5.1]{HULL-big-v2}.

\begin{proposition}   \label{thm:conj-compose-our-version}
  Let $f:\Rn\rightarrow\Rext$ be convex, and
  let $g:\R\rightarrow\Rext$ be convex, proper and nondecreasing.
  Assume there exists $\xing\in\Rn$ such that
  $f(\xing)<\sup(\dom g)$.
  Let $G:\eR\to\eR$ be the monotone extension of $g$ with $G(-\infty)=\inf g$
  and $G(+\infty)=+\infty$.
  Let $\uu\in\Rn$.
  Then
  \[
     (G\circ f)^*(\uu)
     =
     \min_{v\in\R} \bigBracks{g^*(v)+\indepifstar(\rpairf{\uu}{-v})},
  \]
  with the minimum attained.
\end{proposition}

\begin{proof}
In the proof, we consider points $\zz\in\R^{n+1}$ of the form $\rpair{\xx}{y}$ where $\xx\in\Rn$ and $y\in\R$. To extract $\xx$ and $y$ from $\zz$, we introduce two matrices:
The first matrix, in $\R^{n\times(n+1)}$, is
$\PPx=[\Idnn,\zerov{n}]$,
where $\Idnn$ is the $n\times n$ identity matrix and $\zerov{n}$ is
the all-zeros vector in $\Rn$.
The other matrix, in $\R^{1\times(n+1)}$, is
$\PPy=[\trans{\zerov{n}},1]=\trans{\ee_{n+1}}$,
where $\ee_{n+1}=\rpair{\zerov{n}}{1}$ is the $(n+1)$-st standard basis vector in
$\R^{n+1}$.
Thus, for $\zz=\rpair{\xx}{y}$, we obtain $\PPx\zz=\xx$ and $\PPy\zz=y$.

Let $h=G\circ f$.
We aim to derive an expression for $\hstar(\uu)$.

Because $g$ is nondecreasing, we can rewrite $h$, for $\xx\in\Rn$, as
\begin{align}
  h(\xx)
  \nonumber
  &=
  \inf\Braces{ g(y):\:y\in\R,\, y\geq f(\xx) }
  \\
  \notag
  &=
  \inf\Braces{ g(y) + \indepif(\rpairf{\xx}{y}):\:y\in\R }
  \\
  \label{eq:thm:conj-compose-our-version:1}
  &=
  \inf\Braces{ g(\PPy\zz) + \indepif(\zz):\:\zz\in\R^{n+1},\, \PPx\zz=\xx }.
\end{align}
Note that the definition of $G$ at $\pm\infty$ implies
that
the first equality
is valid even if $f(\xx)$ is $\pm\infty$.

Let $r=g\PPy$; in
other words, for
$\xx\in\Rn$ and $y\in\R$,
$r(\rpairf{\xx}{y})=g(y)$.
Furthermore, let $s=r+\indepif=g\PPy+\indepif$.
Then \eqref{eq:thm:conj-compose-our-version:1} shows that
$h=\PPx s$.
Consequently,
by \Cref{roc:thm16.3:Af},
$\hstar=\sstar \trans{\PPx}$, that is,
\begin{equation}  \label{eq:thm:conj-compose-our-version:4}
  \hstar(\uu) = \sstar(\trans{\PPx} \uu) = \sstar(\rpairf{\uu}{0}).
\end{equation}
It only remains then to compute $\sstar(\rpairf{\uu}{0})$.
For this, we prove the following facts which are needed for
applying some of the preceding conjugacy tools.
We state these as a lemma for later reference.

\begin{lemma}   \label{lem:thm:conj-compose-our-version:1}
  Let $f,g,\xing$ be as stated in
  \Cref{thm:conj-compose-our-version},
  and let $r=g\PPy$ where $\PPy=[\trans{\zero_n},1]$.
  Then:
  \begin{letter-compact}
  \item   \label{lem:thm:conj-compose-our-version:1:a}
    There exists $\zing\in\R^{n+1}$ such that
    $\PPy\zing\in \ri(\dom{g})$.    
  \item   \label{lem:thm:conj-compose-our-version:1:b}
    $\ri(\dom{r})\cap\ri(\epi f)\neq\emptyset$.
  \end{letter-compact}
\end{lemma}

\begin{proofx}
Since $f(\xing)<\sup(\dom g)$, there exists
$z\in\dom g$ with $z>f(\xing)$.
Let $U=(-\infty,z)$, which is open (in $\R$) and which is included in
$\dom g$ since $g$ is nondecreasing and $g(z)<+\infty$.
Also, let $\ying\in\R$ be such that $f(\xing)\leq\ying<z$,
implying $\ying\in U\subseteq\dom g$,
and thus that $\ying\in\intdom{g}$,
and so also that $\ying\in\ri(\dom g)$
(\Cref{pr:ri-props}\ref{pr:ri-props:intC-nonemp-implies-eq-riC}).

\begin{proof-parts}
\pfpart{Part~(\ref{lem:thm:conj-compose-our-version:1:a}):}
By the foregoing,
$\PPy\rpair{\zero}{\ying}=\ying\in\ri(\dom{g})$,
proving the claim.

\pfpart{Part~(\ref{lem:thm:conj-compose-our-version:1:b}):}
The set $\Rn\times U$ is open in $\Rn\times\R=\R^{n+1}$ and includes the
point $\rpair{\xing}{\ying}$.
Furthermore,
$\Rn\times U\subseteq\dom{r}$
since
$r(\rpairf{\xx}{y})=g(y)\leq g(z)<+\infty$
for $\xx\in\Rn$ and $y\in U$,
since $g$ is nondecreasing.
Therefore, $\rpair{\xing}{\ying}\in\intdom{r}$.
Also, $\rpair{\xing}{\ying}\in\epi{f}$.
Thus, $\intdom{r}\cap\epi{f}\neq\emptyset$.
The claim now follows
by \Cref{pr:ri-props}(\ref{pr:ri-props:intC-D-implies-riC-riD}).
\qedhere
\end{proof-parts}
\end{proofx}

Since $s=r+\indepif$ and
in light of
\Cref{lem:thm:conj-compose-our-version:1}(\ref{lem:thm:conj-compose-our-version:1:b}),
we can now compute $\sstar(\rpairf{\uu}{0})$ using
\Cref{roc:thm16.4},
yielding
(with Eq.~\ref{eq:thm:conj-compose-our-version:4}) that
\begin{equation}  \label{eq:thm:conj-compose-our-version:3}
  \hstar(\uu)
  =
  \sstar(\rpairf{\uu}{0})
  =
  \inf\,\BigBraces{\rstar(\ww) + \indepifstar\bigParens{\rpair{\uu}{0} - \ww}
                 :\:\ww\in\R^{n+1} },
\end{equation}
and furthermore that the infimum is attained.

Next,
since $r = g\PPy$ and
in light of
\Cref{lem:thm:conj-compose-our-version:1}(\ref{lem:thm:conj-compose-our-version:1:a}),
we can apply
\Cref{roc:thm16.3:fA}
to compute $\rstar$, yielding, for $\ww\in\R^{n+1}$, that
\[
  \rstar(\ww)
  =
  \inf\,\BigBraces{ \gstar(v) :\: v\in\R,\,\trans{\PPy} v  = \ww }.
\]
This means that if $\ww$ has the form
$\ww= \trans{\PPy} v=\rpair{\zero}{v}$, for some $v\in\R$,
then this $v$ must be unique and
we must have
$\rstar(\ww)=\gstar(v)$.
Otherwise, if $\ww$ does not have this form, then vacuously
$\rstar(\ww)=+\infty$.

This last fact implies that, in calculating the infimum in
\eqref{eq:thm:conj-compose-our-version:3},
we need only consider $\ww$ of the form
$\ww=\rpair{\zero}{v}$, for $v\in\R$.
Thus, that equation becomes
\begin{align}
\notag
  \hstar(\uu)
  &=
  \inf\,\BigBraces{\rstar(\rpairf{\zero}{v}) + \indepifstar\bigParens{\rpair{\uu}{0} - \rpair{\zero}{v}}
                 :\: v\in\R }
\\
\label{eq:thm:conj-compose-our-version:6}
  &=
  \inf\,\BigBraces{\gstar(v) + \indepifstar(\rpairf{\uu}{-v})
                 :\: v\in\R }.
\end{align}

It remains to show that the infimum is always attained.
If $\hstar(\uu)<+\infty$, then, as previously noted, the infimum in
\eqref{eq:thm:conj-compose-our-version:3}
is attained by some $\ww\in\R^{n+1}$, which must be of the form
$\ww=\rpair{\zero}{v}$, for some $v\in\R$
(since otherwise $\rstar(\ww)=+\infty$),
implying that the infimum in
\eqref{eq:thm:conj-compose-our-version:6}
is attained as well.
Otherwise, if $\hstar(\uu)=+\infty$, then
\eqref{eq:thm:conj-compose-our-version:6}
implies that the infimum is attained at every point $v\in\R$.
\end{proof}

\indexg{epigraph!support function of|(}%
\indexg{support functions (standard)!epigraph@of epigraph|(}%
The support function $\indepifstar$ that appears in
\Cref{thm:conj-compose-our-version}
can be expressed as follows (where $\sfprod{v}{f}$ is as defined in
Eq.~\ref{eq:sfprod-defn}):

\begin{proposition}
\label{pr:support-epi-f-conjugate}
  Let $f:\Rn\rightarrow\Rext$,
  let $\uu\in\Rn$, and let $v\in\R$.
  If $v\geq 0$ then
  \[
    \indepifstar(\rpairf{\uu}{-v})
    =
    \sfprodstar{v}{f}(\uu)
    =
    \begin{cases}
      v \fstar(\uu/v)
        & \text{if $v>0$,} \\
      \inddomfstar(\uu)
        & \text{if $v=0$.}
    \end{cases}
  \]
  Otherwise, if $v<0$ then
  \[
    \indepifstar(\rpairf{\uu}{-v})
    =
    \begin{cases}
      -\infty
        & \text{if $f\equiv+\infty$,} \\
      +\infty
        & \text{otherwise.}
    \end{cases}
  \]
\end{proposition}

\begin{proof}
Suppose first that $v\geq 0$.
Then
\begin{align*}
  \indepifstar(\rpairf{\uu}{-v})
  &=
  \sup_{\rpair{\xx}{y}\in\epi f}
      \bigBracks{ \xx\cdot\uu - yv }
  \\
  &=
  \adjustlimits\sup_{\xx\in\dom f} \sup_{\,\,y\in\R:\:y\geq f(\xx)\,\,}
      \!\bigBracks{ \xx\cdot\uu - yv }
  \\
  &=
  \sup_{\xx\in\dom f} \bigBracks{ \xx\cdot\uu - v f(\xx) }
  \\
  &=
  \sup_{\xx\in\Rn} \bigBracks{ \xx\cdot\uu - (\sfprod{v}{f})(\xx) }
  \\
  &=
  \sfprodstar{v}{f}(\uu).
\end{align*}
The first equality is by definition of support function
(Eq.~\ref{eq:e:2}).
The fourth equality is because $(\sfprod{v}{f})(\xx)$ equals
$v f(\xx)$ if $\xx\in\dom{f}$, and is $+\infty$ otherwise.
And the last equality is by definition of conjugate
(Eq.~\ref{eq:fstar-def:intro}).

If $v=0$ then $\sfprodstar{0}{f}=\inddomfstar$
by \eqref{eq:sfprod-identity}.
And if $v>0$ then
\[
  \sfprodstar{v}{f}(\uu)
  =
  (v f)^*(\uu)
  =
  \sup_{\xx\in\Rn} \bigBracks{ \xx\cdot\uu - v f(\xx) }
  =
  v\!\sup_{\xx\in\Rn} \BigBracks{ \xx\cdot\frac{\uu}{v} -  f(\xx) }
  =
  v \fstar\Parens{\frac{\uu}{v}},
\]
using \eqref{eq:fstar-def:intro}.

Finally, suppose $v<0$.
If $f\equiv+\infty$ then $\indepif\equiv+\infty$
so $\indepifstar\equiv-\infty$.
Otherwise, if $f\not\equiv+\infty$,
let $\xing$ be any point with $f(\xing)<+\infty$.
Then
\[
  \indepifstar(\rpairf{\uu}{-v})
  =
  \sup_{\rpair{\xx}{y}\in\epi f} \bigBracks{ \xx\cdot\uu - yv }
  \geq
  \sup_{y\in\R:\:y\geq f(\xing)}
  \bigBracks{ \xing\cdot\uu - yv }
  =
  +\infty
\]
since $\rpair{\xing}{y}\in\epi f$ for all $y\geq f(\xing)$.%
\indexg{conjugate (standard)!composition with nondecreasing function@of composition with nondecreasing function|)}%
\indexg{nondecreasing function, composition with!standard conjugate of|)}%
\indexg{epigraph!support function of|)}%
\indexg{support functions (standard)!epigraph@of epigraph|)}%
\end{proof}

\section{Subgradients, subdifferentials, and directional derivatives}
\label{sec:prelim:subgrads}

\indexg{subdifferentials, standard|(}%
We say that $\uu\in\Rn$ is a \emph{subgradient} of a function
$f:\Rn\rightarrow\Rext$ at a point $\xx\in\Rn$ if
\begin{equation}  \label{eqn:prelim-standard-subgrad-ineq}
\indexg{subdifferentials, standard!defined}%
   f(\xx') \geq f(\xx) + (\xx'-\xx)\cdot\uu
\end{equation}
for all $\xx'\in\Rn$.
The \emph{subdifferential} of $f$ at $\xx$, denoted
\indexm{d fx200}{$\partial f(\xx)$}{subdifferential (standard)}%
$\partial f(\xx)$, is the set of all subgradients of $f$ at $\xx$.
\indexg{optimality conditions!standard|(}%
It is immediate from this definition that
$\zero\in\partial f(\xx)$ if and only if $\xx$ minimizes
\indexg{optimality conditions!standard|)}%
$f$.
Subdifferentials and their extension to astral space will be explored
in considerable detail beginning in
\Cref{sec:gradients}.

Here are some properties and characterizations of subgradients:

\begin{proposition}  \label{roc:thm23.4}
  Let $f:\Rn\rightarrow\Rext$ be convex and proper,
  and let $\xx\in\Rn$.
  \begin{letter-compact}
  \item  \label{roc:thm23.4:a}
    If $\xx\in\ri(\dom{f})$ then $\partial f(\xx)\neq\emptyset$.
  \item  \label{roc:thm23.4:b}
    If $\partial f(\xx)\neq\emptyset$ then $\xx\in\dom{f}$.
  \end{letter-compact}
\end{proposition}

\begin{proof}
See \idxroc\citet[Theorem~23.4]{ROC}.
\end{proof}

\begin{proposition}  \label{pr:stan-subgrad-equiv-props}
  Let $f:\Rn\rightarrow\Rext$,
  and let $\xx,\uu\in\Rn$.
  Then the following are equivalent:
  \begin{letter-compact}
  \item  \label{pr:stan-subgrad-equiv-props:a}
    $\uu\in\partial f(\xx)$.
  \item  \label{pr:stan-subgrad-equiv-props:b}
    $\fstar(\uu) = \xx\cdot\uu - f(\xx)$.
  \savelccounter
  \end{letter-compact}
  In addition,
  if $f$ is convex and proper, then the statements above and the
  statements below are also all equivalent to one another:
  \begin{letter-compact}
  \restorelccounter
  \item  \label{pr:stan-subgrad-equiv-props:c}
    $(\cl{f})(\xx) = f(\xx)$ and $\xx\in\partial\fstar(\uu)$.
  \item  \label{pr:stan-subgrad-equiv-props:d}
    $(\cl{f})(\xx) = f(\xx)$ and $\uu\in\partial(\cl{f})(\xx)$.
  \end{letter-compact}
\end{proposition}

\begin{proof}
By simple rearranging,
the definition of standard subgradient given in
\eqref{eqn:prelim-standard-subgrad-ineq}
holds if and only if
$\xx\cdot\uu-f(\xx) \geq \xx'\cdot\uu-f(\xx')$
for all $\xx'\in\Rn$, that is, if and only if
\[
  \xx\cdot\uu-f(\xx)
  =
  \sup_{\xx'\in\Rn} [\xx'\cdot\uu-f(\xx')].
\]
Since the term on the right is exactly $\fstar(\uu)$, this
proves the equivalence of
(\ref{pr:stan-subgrad-equiv-props:a})
and
(\ref{pr:stan-subgrad-equiv-props:b}).

When $f$ is convex and proper, the remaining equivalences
follow directly from
Theorem~23.5 and Corollary~23.5.2 of
\idxroc\citet{ROC}.
\end{proof}

\begin{proposition}  \label{roc:thm24.4}
\indexg{subdifferentials, standard!continuity|(}%
\indexg{continuity!standard subgradients@of standard subgradients|(}%
  Let $f:\Rn\rightarrow\Rext$ be closed, proper and convex.
  Let $\seq{\xx_t}$ and $\seq{\uu_t}$ be sequences in $\Rn$,
  and let $\xx,\uu\in\Rn$.
  Suppose $\xx_t\rightarrow\xx$, $\uu_t\rightarrow\uu$,
  and that $\uu_t\in\partial f(\xx_t)$ for all $t$.
  Then $\uu\in\partial f(\xx)$.
\end{proposition}

\begin{proof}
See \idxroc\citet[Theorem~24.4]{ROC}.%
\indexg{subdifferentials, standard!continuity|)}%
\indexg{continuity!standard subgradients@of standard subgradients|)}%
\end{proof}

\begin{proposition}  \label{roc:thm25.1}
\indexg{differentiability|(}%
\indexg{subdifferentials, standard!gradient and|(}%
\indexg{gradients (standard)|(}%
  Let $f:\Rn\rightarrow\Rext$ be convex,
  and let $\xx\in\Rn$ with $f(\xx)\in\R$.
  \begin{letter-compact}
  \item  \label{roc:thm25.1:a}
    If $f$ is differentiable at $\xx$, then $\gradf(\xx)$ is $f$'s
    only subgradient at $\xx$; that is,
    $\partial f(\xx) = \{ \gradf(\xx) \}$.
  \item  \label{roc:thm25.1:b}
    Conversely, if $\partial f(\xx)$ is a singleton, then
    $f$ is differentiable at $\xx$.
  \end{letter-compact}
\end{proposition}

\begin{proof}
See \idxroc\citet[Theorem~25.1]{ROC}.%
\indexg{gradients (standard)|)}%
\indexg{differentiability|)}%
\indexg{subdifferentials, standard!gradient and|)}%
\indexg{subdifferentials, standard|)}
\end{proof}

Here are some rules for computing subdifferentials, which extend
standard calculus rules for finding ordinary gradients:

If $f:\Rn\rightarrow\Rext$ and $\lambda\in\Rstrictpos$,
then for $\xx\in\Rn$,
$\partial (\lambda f)(\xx) = \lambda[\partial f(\xx)]$
(as can be shown straightforwardly from definitions).

As is often the case, the next rule, for sums of functions,
depends partially on a topological condition:

\begin{proposition}   \label{roc:thm23.8}
\indexg{sum of functions (standard)!subgradients of|(}%
\indexg{subdifferentials, standard!sum of functions@of sum of functions|(}%
  Let $f_i:\Rn\rightarrow\Rext$ be convex and proper, for
  $i=1,\ldots,m$.
  Let $h=f_1+\dotsb+f_m$, and let $\xx\in\Rn$.
  \begin{letter-compact}
  \item   \label{roc:thm23.8:a}
    $\partial f_1(\xx)+\dotsb+\partial f_m(\xx) \subseteq \partial h(\xx)$.
  \item   \label{roc:thm23.8:b}
    If, in addition,
    $\bigcap_{i=1}^m \ri(\dom{f_i})\neq\emptyset$,
    then
    $\partial f_1(\xx)+\dotsb+\partial f_m(\xx) = \partial h(\xx)$.
  \end{letter-compact}
\end{proposition}

\begin{proof}
See \idxroc\citet[Theorem~23.8]{ROC}.%
\indexg{subdifferentials, standard!sum of functions@of sum of functions|)}%
\indexg{sum of functions (standard)!subgradients of|)}%
\end{proof}

\begin{proposition}   \label{roc:thm23.9}
\indexg{subdifferentials, standard!composition with linear map@of composition with linear map|(}%
\indexg{linear map in composition with standard function!standard subgradients of|(}%
  Let $f:\Rn\rightarrow\Rext$ be convex and proper,
  and let $\A\in\Rmn$.
  Let $\xx\in\Rn$.
  \begin{letter-compact}
  \item   \label{roc:thm23.9:a}
    $\transA\partial f(\A \xx) \subseteq \partial (\fA)(\xx)$.
  \item   \label{roc:thm23.9:b}
    If, in addition,
    there exists $\xing\in\Rn$ such that $\A \xing\in\ri(\dom{f})$,
    then
    $\transA\partial f(\A \xx) = \partial (\fA)(\xx)$.
  \end{letter-compact}
\end{proposition}

\begin{proof}
See \idxroc\citet[Theorem~23.9]{ROC}.%
\indexg{linear map in composition with standard function!standard subgradients of|)}%
\indexg{subdifferentials, standard!composition with linear map@of composition with linear map|)}%
\end{proof}

\indexg{subdifferentials, standard!linear image of function@of linear image of function|(}%
\indexg{linear image of standard function!standard subgradients of|(}%
The following is based on
\idxhiriart\citet[Theorem~VI.4.5.1]{HULL-big-v1}:

\begin{proposition}   \label{pr:stan-subgrad-lin-img}
  Let $f:\Rn\rightarrow\Rext$
  and $\A\in\Rmn$.
  Let $\xx,\uu\in\Rm$.
  Then the following are equivalent:
  \begin{letter-compact}
  \item   \label{pr:stan-subgrad-lin-img:a}
    $\uu\in\partial (\A f)(\xx)$,
    and also
    the infimum defining $(\A f)(\xx)$ in
    \eqref{eq:lin-image-fcn-defn} is attained by some $\zz\in\Rn$
    with $\A\zz=\xx$.
  \item   \label{pr:stan-subgrad-lin-img:b}
    There exists $\zz\in\Rn$ such that $\A \zz = \xx$
    and $\transA \uu \in \partial f(\zz)$.
  \end{letter-compact}
\end{proposition}

\begin{proof}
Let $h=\A f$.

\begin{proof-parts}
\pfpart{%
  (\ref{pr:stan-subgrad-lin-img:a})
  $\Rightarrow$
  (\ref{pr:stan-subgrad-lin-img:b}):
}
Suppose $\uu\in\partial h(\xx)$ and that
there exists $\zz\in\Rn$
attaining the infimum in \eqref{eq:lin-image-fcn-defn},
implying $\A\zz=\xx$ and $h(\xx)=f(\zz)$.
Then for all $\zz'\in\Rm$,
\[
  f(\zz')
  \geq
  h(\A \zz')
  \geq
  h(\xx) + \uu\cdot(\A\zz' - \xx)
  =
  f(\zz) + (\transA\uu)\cdot(\zz'-\zz).
\]
The first inequality is by $\A f$'s definition.
The second is because $\uu\in\partial h(\xx)$.
The equality
follows because
$h(\xx)=f(\zz)$ and $\xx=\A\zz$.
Thus, $\transA\uu\in\partial f(\zz)$, as claimed.

\pfpart{%
  (\ref{pr:stan-subgrad-lin-img:b})
  $\Rightarrow$
  (\ref{pr:stan-subgrad-lin-img:a}):
}
Suppose there exists $\zz\in\Rn$ with $\A \zz = \xx$
and $\transA \uu \in \partial f(\zz)$.
Let $\zz'\in\Rn$, and let $\xx'=\A\zz'$.
Then
\[
  f(\zz')
  \geq
  f(\zz)
  +
  (\transA \uu)\cdot (\zz' - \zz)
  =
  f(\zz)
  +
  \uu\cdot (\xx' - \xx).
\]
The inequality is because
$\transA \uu \in \partial f(\zz)$.
The equality is
because $\xx'=\A\zz'$ and $\xx=\A\zz$.
Since
this
holds for all $\zz'$ with $\A\zz'=\xx'$, it follows that
\[
  h(\xx')
  \geq
  f(\zz)
  +
  \uu\cdot (\xx' - \xx)
  \geq
  h(\xx)
  +
  \uu\cdot (\xx' - \xx),
\]
where the second inequality is because $f(\zz)\geq h(\xx)$ by $\A f$'s
definition.
This shows that
$\uu\in\partial h(\xx)$.
Furthermore, applied with $\xx'=\xx$, it shows that
$f(\zz)=h(\xx)$ and thus that the infimum
in \eqref{eq:lin-image-fcn-defn} is attained (by $\zz$).%
\indexg{subdifferentials, standard!linear image of function@of linear image of function|)}%
\indexg{linear image of standard function!standard subgradients of|)}%
\qedhere
\end{proof-parts}
\end{proof}

\indexg{suprema and maxima, pointwise!standard subgradients of|(}%
\indexg{subdifferentials, standard!pointwise supremum@of pointwise supremum|(}%
For subgradients of a pointwise supremum, we have the
following general inclusion, followed by a special case in which the
inclusion holds with equality.
A more general result is given by
\idxhiriart\citet[Theorem~VI.4.4.2]{HULL-big-v1}.

\begin{proposition}   \label{pr:stnd-sup-subgrad}
  Let $f_i:\Rn\rightarrow\Rext$ for all $i\in\indset$,
  where $\indset$ is an index set, and let
  $ h = \sup_{i\in\indset} f_i $.
  Let $\xx\in\Rn$, and let
  $J = \Braces{i \in \indset :\: h(\xx)=f_i(\xx)}$.
  Then
  \[
         \cl\conv\Parens{ \bigcup_{i\in J} \partial f_i(\xx) }
           \subseteq
           \partial h(\xx).
  \]
\end{proposition}

\begin{proof}
This can be proved as in Lemma~VI.4.4.1 of
\idxhiriart\citet{HULL-big-v1} (even though that lemma itself is stated for a
more limited case).%
\end{proof}

\begin{proposition}   \label{pr:std-max-fin-subgrad}
  Let $f_i:\Rn\rightarrow\R$ be convex, for $i=1,\dotsc,m$,
  and let
  $ h=\max\regBraces{f_1,\dotsc,f_m} $.
  Let $\xx\in\Rn$, and let
  $ J = \bigBraces{i \in \{1,\ldots,m\} :\: h(\xx)=f_i(\xx)}$.
  Then
  \[
     \partial h(\xx)
     =
     \conv\Parens{ \bigcup_{i\in J} \partial f_i(\xx) }.
  \]
\end{proposition}

\begin{proof}
See \idxhiriart\citet[Corollary~VI.4.3.2]{HULL-big-v1}.%
\indexg{subdifferentials, standard!pointwise supremum@of pointwise supremum|)}%
\indexg{suprema and maxima, pointwise!standard subgradients of|)}%
\end{proof}

\begin{proposition}  \label{pr:std-subgrad-comp-inc}
\indexg{subdifferentials, standard!composition with nondecreasing function@of composition with nondecreasing function|(}%
\indexg{nondecreasing function, composition with!standard subgradients of|(}%
  Let $f:\Rn\rightarrow\R$ be convex,
  and let $g:\R\rightarrow\R$ be convex and nondecreasing.
  Let $h=g\circ f$ and let $\xx\in\Rn$.
  Then
  \[
     \partial h(\xx)
     =
     \Braces{ v \uu :
       \uu\in\partial f(\xx),
       v \in \partial g(f(\xx))
     }.
  \]
\end{proposition}

\begin{proof}
See
\idxhiriart\citet[Theorem~VI.4.3.1]{HULL-big-v1}.%
\indexg{nondecreasing function, composition with!standard subgradients of|)}%
\indexg{subdifferentials, standard!composition with nondecreasing function@of composition with nondecreasing function|)}%
\end{proof}

\indexg{one-sided directional derivatives|(}%
Let $f:\Rn\to\eR$.
The \emph{one-sided directional derivative} of $f$
at a point $\xx\in\Rn$ where $f(\xx)\in\R$,
with respect to a vector $\yy\in\Rn$,
is
defined to be
\begin{equation}   \label{eq:direc-deriv-dfn}
\indexm{f x y}{$\dderf{\xx}{\yy}$}{one-sided directional derivative}
  \dderf{\xx}{\yy}
  =
  \lim_{\rightlim{\lambda}{0}} \frac{f(\xx+\lambda\yy) - f(\xx)}
                                      {\lambda},
\end{equation}
where the limit is taken from the right,
that is, as $\lambda>0$ approaches~$0$.
Thus, $\dderf{\xx}{\yy}$ captures the instantaneous rate of change of
$f$ at $\xx$ along the vector $\yy$.

For convex functions, the limit in the definition in
\eqref{eq:direc-deriv-dfn} can be replaced by an infimum
(and so always exists).

\begin{proposition}   \label{roc:thm23.1}
  Let $f:\Rn\rightarrow\Rext$ be convex, and let $\xx\in\Rn$.
  Assume $f(\xx)\in\R$.
  Then for all $\yy\in\Rn$,
  \[
    \dderf{\xx}{\yy}
    =
    \inf_{\lambda\in\Rstrictpos} \frac{f(\xx+\lambda\yy) - f(\xx)}
                                      {\lambda}.
  \]
  Furthermore, $\dderf{\xx}{\yy}$, as a function of $\yy\in\Rn$,
  is a convex, positively homogeneous function
  with $\dderf{\xx}{\zero}=0$.
\end{proposition}

\begin{proof}
  See \idxroc\citet[Theorem~23.1]{ROC}.
\end{proof}

Directional derivatives are closely related to subgradients.
In particular, a function $f$'s subdifferential at every point $\xx$
where $f(\xx)\in\R$ is fully determined by $f$'s directional
derivatives at $\xx$, according to the following relationship:

\begin{proposition}   \label{roc:thm23.2}
  Let $f:\Rn\rightarrow\Rext$ be convex, and let $\xx,\uu\in\Rn$.
  Assume $f(\xx)\in\R$.
  Then $\uu\in\partial f(\xx)$ if and only if
  $\yy\cdot\uu\leq \dderf{\xx}{\yy}$ for all $\yy\in\Rn$.
\end{proposition}

\begin{proof}
  See \idxroc\citet[Theorem~23.2]{ROC}.
\end{proof}

For proper convex functions $f$, and at points $\xx$ in the relative
interior of $\dom f$, the directional derivatives are likewise
determined by the subdifferential at $\xx$:

\begin{proposition}   \label{roc:thm23.4-dd}
  Let $f:\Rn\rightarrow\Rext$ be convex and proper, let
  $\xx\in\ri(\dom f)$, and let $\yy\in\Rn$.
  Then
  \[
     \dderf{\xx}{\yy}     
     =
     \indfstar{\partial f(\xx)}(\yy)
     =
     \sup\bigBraces{
       \yy\cdot\uu
       :\:
       \uu\in\partial f(\xx)
     }.
  \]
\end{proposition}

\begin{proof}
  See \idxroc\citet[Theorem~23.4]{ROC}.%
\indexg{one-sided directional derivatives|)}
\end{proof}

\section{Recession cone and constancy space}
\label{sec:prelim:rec-cone}

\indexg{recession cone (standard)|(}%
The \emph{recession cone} of a function $f:\Rn\rightarrow\Rext$,
denoted $\resc{f}$,
is the set of directions in which the function never increases:
\begin{equation}  \label{eqn:resc-cone-def}
\indexg{recession cone (standard)!defined}%
\indexm{rec f200}{$\resc f$}{recession cone (standard)}%
  \resc{f}
  =
  \Braces{\vv\in\Rn :\:
    \forall \xx\in\Rn,\,
    \forall \lambda\in\Rpos,\,
    f(\xx+\lambda\vv)\leq f(\xx)
  }.
\end{equation}
Note that \idxroc\citet[\pcite{69}]{ROC} only defines recession cones for convex functions $f\not\equiv+\infty$, whereas we
define them for general functions, including $f\equiv+\infty$ (in which case $\resc{f}=\Rn$).

\begin{proposition}   \label{pr:resc-cone-basic-props}
  Let $f:\Rn\rightarrow\Rext$.
  Then $f$'s recession cone, $\resc{f}$, is a convex cone.
  In addition,
  if $f$ is lower semicontinuous, then $\resc{f}$ is
  closed in~$\Rn$.
\end{proposition}

\begin{proof}
  ~

\begin{proof-parts}

\pfpart{Cone:}
This is
immediate from definitions.

\pfpart{Convex:}
Let $\vv$ and $\ww$ be in $\resc{f}$, and let $\lambda\in\Rpos$.
Then for all $\xx\in\Rn$,
\[
   f\bigParens{\xx+\lambda(\vv+\ww)}
   =
   f(\xx+\lambda\vv+\lambda\ww)
   \leq
   f(\xx+\lambda\vv)
   \leq
   f(\xx).
\]
The first inequality is because $\ww\in\resc{f}$, and the second is
because $\vv\in\resc{f}$.
Thus, $\vv+\ww\in\resc{f}$.
Therefore, $\resc{f}$ is convex by
\Cref{pr:scc-cone-elts}(\ref{pr:scc-cone-elts:d}).

\pfpart{Closed:}
Assume $f$ is lower semicontinuous.
Let $\seq{\vv_t}$ be any convergent sequence in $\resc{f}$, and
suppose its limit is $\vv$.
Then for any $\xx\in\Rn$ and $\lambda\in\Rpos$,
\[
  f(\xx)
  \geq
  \liminf  f(\xx+\lambda\vv_t)
  \geq
  f(\xx+\lambda\vv)
\]
since $\xx+\lambda\vv_t\rightarrow\xx+\lambda\vv$
and $f$ is lower semicontinuous.
Thus, $\vv\in\resc{f}$, so $\resc{f}$ is closed.
\qedhere
\end{proof-parts}
\end{proof}

\begin{proposition}  \label{pr:stan-rec-equiv}
  Let $f:\Rn\rightarrow\Rext$ be convex and lower semicontinuous,
  and let $\vv\in\Rn$.
  Then the following are equivalent:
  \begin{letter-compact}
  \item  \label{pr:stan-rec-equiv:a}
    $\vv\in\resc{f}$.
  \item  \label{pr:stan-rec-equiv:b}
    For all $\xx\in\Rn$, $f(\xx+\vv)\leq f(\xx)$.
  \item  \label{pr:stan-rec-equiv:c}
    Either $f\equiv+\infty$ or
    \begin{equation}  \label{eq:pr:stan-rec-equiv:1}
      \liminf_{\lambda\rightarrow+\infty}  f(\xx+\lambda\vv)
      <
      +\infty
    \end{equation}
    for some $\xx\in\Rn$.
  \end{letter-compact}
\end{proposition}

\begin{proof}
  ~

\begin{proof-parts}
\pfpart{%
  (\ref{pr:stan-rec-equiv:a})
  $\Rightarrow$
  (\ref{pr:stan-rec-equiv:b}):
}
This is immediate from the definition of $\resc{f}$
(Eq.~\ref{eqn:resc-cone-def}).

\pfpart{%
  (\ref{pr:stan-rec-equiv:b})
  $\Rightarrow$
  (\ref{pr:stan-rec-equiv:c}):
}
Suppose statement~(\ref{pr:stan-rec-equiv:b}) holds
and that $f\not\equiv+\infty$.
Let $\xx$ be any point in $\dom{f}$.
Then for all $t$,
$f(\xx+t\vv)\leq f(\xx+(t-1)\vv)$,
so by an inductive argument,
$f(\xx+t\vv)\leq f(\xx)$ for all $t$.
Therefore,
$\liminf_{\lambda\rightarrow+\infty} f(\xx+\lambda\vv) \leq f(\xx) < +\infty$,
proving
\eqref{eq:pr:stan-rec-equiv:1}.

\pfpart{%
  (\ref{pr:stan-rec-equiv:c})
  $\Rightarrow$
  (\ref{pr:stan-rec-equiv:a}):
}
If $f$ is proper, then this follows directly from
\idxroc\citet[Theorem~8.6]{ROC}.
Also, if $f\equiv+\infty$, then $\resc{f}=\Rn$, again proving the claim.
So assume henceforth that $f$ is improper and that
\eqref{eq:pr:stan-rec-equiv:1}
holds for some $\xx\in\Rn$.
Then $f(\zz)\in\{-\infty,+\infty\}$ for all $\zz\in\Rn$, by
\Cref{pr:improper-vals}(\ref{pr:improper-vals:cor7.2.1}).
Consequently, there exists a sequence $\seq{\lambda_t}$ in $\R$ with
$\lambda_t\rightarrow+\infty$ and $f(\xx+\lambda_t\vv)=-\infty$
for all $t$.

Let $S=\dom{f}$, which is convex
(\Cref{roc:thm4.6})
and closed
(by \Cref{prop:lsc}\ref{prop:lsc:a}\ref{prop:lsc:c},
applied with $\alpha=-\infty$).
Let $f'=\inds$, the indicator function for $S$
(Eq.~\ref{eq:indf-defn}), which is convex, closed and proper.
Then $f'(\xx+\lambda_t\vv)=0$ for all $t$, so
\eqref{eq:pr:stan-rec-equiv:1}
holds for $f'$, implying $\vv\in\resc{f'}$ by the argument above.
Since $f(\yy)\leq f(\zz)$ if and only if $f'(\yy)\leq f'(\zz)$ for all
$\yy,\zz\in\Rn$,
it follows from the definition of recession cone that
$\resc{f'}=\resc{f}$.
Therefore, $\vv\in\resc{f}$.
\qedhere
\end{proof-parts}
\end{proof}

The recession cone of a closed, proper convex function
$f:\Rn\rightarrow\Rext$
can be expressed as a polar of
$\cone(\dom{\fstar})$,
the cone generated by the effective domain of its conjugate:

\begin{proposition}  \label{pr:rescpol-is-con-dom-fstar}
  Let $f:\Rn\rightarrow\Rext$ be closed and convex,
  with $f\not\equiv+\infty$.
  Then
  \begin{equation}
    \label{eq:pr:rescpol-is-con-dom-fstar:1}
    \rescpol{f}
    =
    \cl\bigParens{\cone(\dom{\fstar})}.
  \end{equation}
  Consequently,
  \begin{equation}
    \label{eq:pr:rescpol-is-con-dom-fstar:2}
    \resc{f}
    =
    \polar{\bigParens{\cone(\dom{\fstar})}}
    =
    \Braces{\vv\in\Rn :\:  \uu\cdot\vv\leq 0 \textup{ for all } \uu\in\dom{\fstar}}.
    \end{equation}
\end{proposition}

\begin{proof}
For
\eqref{eq:pr:rescpol-is-con-dom-fstar:1},
see
\idxroc\citet[Theorem~14.2]{ROC}
for when $f$ is proper.
Otherwise, if $f$ is improper, then we must have $f\equiv-\infty$
since $f$ is closed and $f\not\equiv+\infty$.
In this case, $\fstar\equiv+\infty$ and $\resc{f}=\Rn$, implying
$\rescpol{f}=\set{\zero}=\cl\bigParens{\cone(\dom{\fstar})}$.

For the first equality of
\eqref{eq:pr:rescpol-is-con-dom-fstar:2}, we
then have
\[
  \resc{f}
  =
  \rescdubpol{f}
  =
  \polar{\bigBracks{\cl\bigParens{\cone(\dom{\fstar})}}}
  =
  \polar{\bigParens{\cone(\dom{\fstar})}}.
\]
The first equality is from
\Cref{pr:polar-props}(\ref{pr:polar-props:c})
and since $\resc{f}$ is a closed (in $\Rn$) convex cone
(\Cref{pr:resc-cone-basic-props}).
The second equality is from
\eqref{eq:pr:rescpol-is-con-dom-fstar:1}.
The last equality is by
\Cref{pr:polar-props}(\ref{pr:polar-props:a}).

Finally, the second equality of
\eqref{eq:pr:rescpol-is-con-dom-fstar:2}
follows from
\Cref{pr:polar-props}(\ref{pr:polar-props:coneSpol}).%
\indexg{recession cone (standard)|)}
\end{proof}

\indexg{constancy space|(}%
The \emph{constancy space} of a function $f:\Rn\rightarrow\Rext$,
denoted $\conssp{f}$,
consists of those directions in which the value of $f$ remains constant:
\begin{equation}  \label{eq:conssp-defn}
\indexm{cons}{$\conssp{f}$}{constancy space}%
  \conssp{f} =
    \Braces{\vv\in\Rn :\: \forall \xx\in\Rn,\,
                        \forall \lambda\in\R,
                        f(\xx+\lambda\vv) = f(\xx)
    }.
\end{equation}
As was the case with recession cones,
\idxroc\citet[\pcite{69}]{ROC} only defines the constancy space for convex functions $f\not\equiv+\infty$, whereas we
define it for general functions, including $f\equiv+\infty$ (in which case $\conssp{f}=\Rn$).

\begin{proposition}  \label{pr:prelim:const-props}
  Let $f:\Rn\rightarrow\Rext$.
  \begin{letter-compact}
  \item  \label{pr:prelim:const-props:a}
    $\conssp{f}=(\resc{f}) \cap (-\resc{f})$.
    Consequently,
    $\zero\in\conssp{f}\subseteq\resc{f}$.
  \item  \label{pr:prelim:const-props:b}
    $\conssp{f}$ is a linear subspace.
  \end{letter-compact}
\end{proposition}

\begin{proof}
  ~

\begin{proof-parts}

\pfpart{Part~(\ref{pr:prelim:const-props:a}):}
If $\vv\in\conssp{f}$ then from definitions, both $\vv$ and $-\vv$ are
in $\resc{f}$.
Thus,
$\conssp{f}\subseteq(\resc{f}) \cap (-\resc{f})$.

For the reverse inclusion, suppose
$\vv\in(\resc{f}) \cap (-\resc{f})$.
Then for all $\xx\in\Rn$ and all $\lambda\in\Rpos$,
$f(\xx+\lambda\vv)\leq f(\xx)$
and
$f(\xx)=f((\xx+\lambda\vv)-\lambda\vv)\leq f(\xx+\lambda\vv)$
since $-\vv\in\resc{f}$.
Thus,
$f(\xx+\lambda\vv) = f(\xx)$.
Therefore, $\vv\in\conssp{f}$.

\pfpart{Part~(\ref{pr:prelim:const-props:b}):}
Since $\resc{f}$ is a convex cone
(\Cref{pr:resc-cone-basic-props}),
this follows from
part~(\ref{pr:prelim:const-props:a})
and
\idxroc\citet[Theorem~2.7]{ROC}.
\qedhere
\end{proof-parts}
\end{proof}

\begin{proposition}
\label{pr:cons-PP}
  Let $f:\Rn\to\eR$. Let $S\subseteq\Rn$, and let $\PP$ be the orthogonal projection matrix onto $\Sperp$. Then the following are equivalent:
  \begin{letter-compact}
  \item \label{pr:cons-PP:a}
     $S\subseteq\conssp{f}$.
  \item \label{pr:cons-PP:b}
     $f(\xx)=f(\PP\xx)$ for all $\xx\in\Rn$.
  \end{letter-compact}
\end{proposition}

\begin{proof}
Let $L=\Sperp$, so that $\PP$ is the projection matrix onto $L$
and implying $\Lperp=\spn{S}$ by
\Cref{pr:std-perp-props}(\ref{pr:std-perp-props:c}).

\begin{proof-parts}
\pfpart{%
  (\ref{pr:cons-PP:a})
  $\Rightarrow$
  (\ref{pr:cons-PP:b}):
}
Suppose $S\subseteq\conssp{f}$.
Then also $\Lperp=\spn{S}\subseteq\conssp{f}$ (by \Cref{pr:prelim:const-props}\ref{pr:prelim:const-props:b}). 
Let $\xx\in\Rn$.
By \Cref{pr:lin-decomp}, there exists $\zz\in \Lperp$ such that
$\xx=\PP\xx+\zz$.
Then
$f(\xx)=f(\PP\xx+\zz)=f(\PP\xx)$ since $\zz\in\Lperp\subseteq\conssp{f}$.

\pfpart{%
  (\ref{pr:cons-PP:b})
  $\Rightarrow$
  (\ref{pr:cons-PP:a}):
}
Suppose $f=f\PP$, and let $\vv\in S$.
Then for all $\xx\in\Rn$ and $\lambda\in\R$,
\[
  f(\xx+\lambda\vv)=f\bigParens{\PP(\xx+\lambda\vv)}
  =f(\PP\xx)=f(\xx),
\]
where the first and last equalities follow by assumption,
and the second is because
$\PP\vv=\zero$ by \Cref{pr:proj-mat-props}(\ref{pr:proj-mat-props:e})
since $\vv\in S\subseteq \Lperp$.
Thus, $\vv\in\conssp{f}$.
\qedhere
\end{proof-parts}
\end{proof}

\begin{proposition}  \label{pr:cons-equiv}
  Let $f:\Rn\rightarrow\Rext$ be convex and lower semicontinuous,
  and let $\vv\in\Rn$.
  Then the following are equivalent:
  \begin{letter-compact}
  \item  \label{pr:cons-equiv:a}
    $\vv\in\conssp{f}$.
  \item  \label{pr:cons-equiv:b}
    For all $\xx\in\Rn$, $f(\xx+\vv) = f(\xx)$.
  \item  \label{pr:cons-equiv:c}
    Either $f\equiv+\infty$ or
    \begin{equation}  \label{eq:pr:cons-equiv:1}
      \sup_{\lambda\in\R}  f(\xx+\lambda\vv)
      <
      +\infty
    \end{equation}
    for some $\xx\in\Rn$.
  \end{letter-compact}
\end{proposition}

\begin{proof}
  ~

\begin{proof-parts}
\pfpart{%
  (\ref{pr:cons-equiv:a})
  $\Rightarrow$
  (\ref{pr:cons-equiv:b}):
}
This is immediate from the definition of $\conssp{f}$
(Eq.~\ref{eq:conssp-defn}).

\pfpart{%
  (\ref{pr:cons-equiv:b})
  $\Rightarrow$
  (\ref{pr:cons-equiv:a}):
}
Suppose~(\ref{pr:cons-equiv:b}) holds, implying
$\vv\in\resc{f}$ by
\Cref{pr:stan-rec-equiv}(\ref{pr:stan-rec-equiv:a},\ref{pr:stan-rec-equiv:b}).
Then also $f(\xx-\vv)=f(\xx)$ for all $\xx\in\Rn$,
similarly implying $-\vv\in\resc{f}$.
Thus, $\vv\in(\resc{f})\cap(-\resc{f})=\conssp{f}$
by
\Cref{pr:prelim:const-props}(\ref{pr:prelim:const-props:a}).

\pfpart{%
  (\ref{pr:cons-equiv:a})
  $\Rightarrow$
  (\ref{pr:cons-equiv:c}):
}
Suppose $\vv\in\conssp{f}$ and $f\not\equiv+\infty$.
Let $\xx$ be any point in $\dom{f}$.
Then
$f(\xx+\lambda\vv)=f(\xx)<+\infty$ for all $\lambda\in\R$, implying
\eqref{eq:pr:cons-equiv:1}.

\pfpart{%
  (\ref{pr:cons-equiv:c})
  $\Rightarrow$
  (\ref{pr:cons-equiv:a}):
}
If $f\equiv+\infty$ then $\conssp{f}=\Rn$, implying the claim.
Otherwise, suppose there exist $\xx\in\Rn$ and $\beta\in\R$ such that
$f(\xx+\lambda \vv) \leq \beta$ for all $\lambda\in\R$.
Then
\eqref{eq:pr:stan-rec-equiv:1} holds for both $\vv$ and $-\vv$,
so both points are in $\resc{f}$, by
\Cref{pr:stan-rec-equiv}(\ref{pr:stan-rec-equiv:a},\ref{pr:stan-rec-equiv:c}).
Thus, as above,
$\vv\in(\resc{f})\cap(-\resc{f})=\conssp{f}$.%
\indexg{constancy space|)}
\qedhere
\end{proof-parts}
\end{proof}

\part{Astral Space and Astral Points}
\label{part:astral-space}

\chapter{Constructing astral space}  \label{sec:astral-space-intro}

\indexg{astral space!construction|(}%
We begin our development with the formal construction
of {astral space}.
Earlier, in \Cref{sec:intro:astral}, we gave a quick
introduction to astral space, including the main ideas of its
construction.
Nevertheless, in this chapter, we start again fresh so as to bring out more
fully the intuition and rationale for our construction, and to fill in
its details more formally and precisely while also
establishing some first basic properties of astral points.

At a high level,
our aim is to construct astral space as
a compact extension of Euclidean space
in which various points at infinity have been added.
We require compactness because it is such a powerful property, and
because the lack of compactness is at the heart of what makes $\Rn$
difficult to work with when sequences to infinity are involved.
There are numerous ways of compactifying $\Rn$.
Our purpose here is to construct a compactification that will be
favorably compatible with convex analysis, which we also hope to
extend to the enlarged space.
For that reason, we specifically aim for a compactification
that allows the continuous extension of all linear functions, since linear
functions are so fundamental to convex analysis.
We will see through the course of our development how this property
helps to ensure that notions and properties that are
built on linearity tend to extend reasonably to astral space.

\section{The construction}
\label{sec:astral:construction}
\label{subsec:astral-intro-motiv}

Unbounded sequences in $\Rn$ cannot converge because there is
nothing ``out there'' for them to converge to.
The main idea in constructing astral space is to add to $\Rn$ various
points at infinity so that such sequences can have limits.
Astral points thus represent
limits of sequences in $\Rn$, including those that are unbounded.
To construct the space then, we need to answer two key questions:
\begin{item-compact}
\item[(1)]
  Which sequences in $\Rn$ should have limits in the new space?
\item[(2)]
  When should two such sequences have the same limit in the new space?
\end{item-compact}
With answers to these questions, constructing the space will be fairly
straightforward:
We will just need to add one ``new'' point for every group of
sequences that are all meant to be convergent according to the first
question, and that are all meant to have the same limit according to
the second question.

In this work, as a fundamental, guiding principle, we base our answers
to these questions on the limits of all linear functions
when evaluated on a given sequence.
In other words, for a sequence $\seq{\xx_t}$ in $\Rn$, we focus,
for all vectors $\uu\in\Rn$, on
the limit of the associated sequence $\seq{\xx_t\cdot\uu}$,
which is the image of $\seq{\xx_t}$ under the linear function
$\xx\mapsto\xx\cdot\uu$.

\begin{figure}
  \centering
  \includegraphics{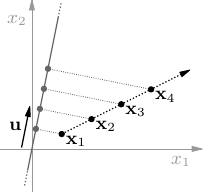}\hfill%
  \includegraphics{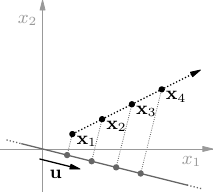}\hfill%
  \includegraphics{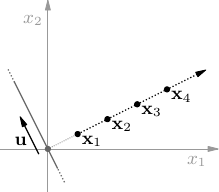}
\bigskip

  \includegraphics{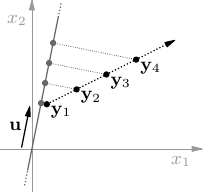}\hfill%
  \includegraphics{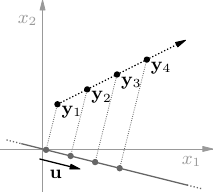}\hfill%
  \includegraphics{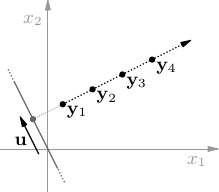}
\bigskip

  \includegraphics{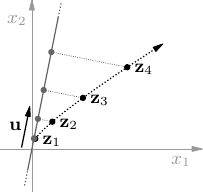}\hfill%
  \includegraphics{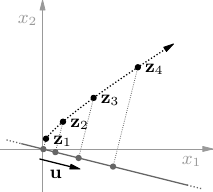}\hfill%
  \includegraphics{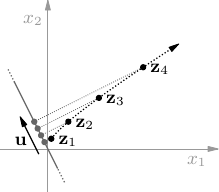}
  \caption[All-directions convergence]{%
\indexf{all-directions equivalence}%
\indexf{convergence in all directions}%
    \emph{All-directions convergence of sequences $\seq{\xx_t}$, $\seq{\yy_t}$, and~$\seq{\zz_t}$
    from Examples~\ref{ex:basic-lim-ray}, \ref{ex:basic-with-offset}, and~\ref{ex:lim-quadratic},
    respectively.}
    \emph{First row, sequence $\seq{\xx_t}$:}
    Convergence of scalar projections
    $\xx_t\inprod\uu$ for three choices of vector $\uu$. For the first and second choice of $\uu$,
    the sequence $\xx_t\inprod\uu$ converges to $+\infty$. For the last choice of $\uu$, the sequence converges
    to $0$.
    \emph{Second and third row, sequences $\seq{\yy_t}$ and $\seq{\zz_t}$:}
    Convergence of $\yy_t\inprod\uu$ and $\zz_t\inprod\uu$ for the same choices of~$\uu$ as in the first row.
    For the first and second choice of~$\uu$, sequences $\yy_t\inprod\uu$ and $\zz_t\inprod\uu$ converge to $+\infty$ just like $\xx_t\inprod\uu$, but for the last choice,
    $\yy_t\inprod\uu$ converges to a finite value different from $0$ and $\zz_t\inprod\uu$ converges to $+\infty$.
    Thus, the three sequences are not all-directions equivalent, although each of them converges in all directions
    (as we argue in the respective examples).
  }
  \label{fig:seq-ex5.234}%
\end{figure}

Geometrically, if $\uu$ is a unit vector in $\Rn$, then this is
equivalent to considering the sequence's projection onto 
the line passing through the origin in the direction
of $\uu$.
This is because
$\xx_t\cdot\uu$ is the signed distance from the origin to $\xx_t\negKern$'s
projection along this line, which thus tracks the progress of
these projected points along the line.
If $\uu$ is not a unit vector (and is not $\zero$), then
the distance is simply
scaled by the constant $\norm{\uu}$. For an illustration, see
\Cref{fig:seq-ex5.234}, where we show projections of three different
sequences (one in each row)
onto three lines corresponding to three different vectors $\uu$.

\indexg{convergence in all directions|(}%
\indexg{convergence in all directions!defined|(}%
In accord with the guiding principle above, we say that a sequence $\seq{\xx_t}$ in
$\Rn$ \emph{converges in all directions} if
for all $\uu\in\Rn$, the sequence $\seq{\xx_t\cdot\uu}$ converges to a
limit in $\Rext$.
In our construction, we will then answer the first question above by
defining limits in the new space exactly for those sequences that
converge in all directions.%
\indexg{convergence in all directions!defined|)}

\indexg{all-directions equivalence|(}%
\indexg{all-directions equivalence!defined|(}%
Similarly, we say that a pair of sequences
$\seq{\xx_t}$ and $\seq{\yy_t}$ in
$\Rn$ 
are \emph{all-directions equivalent}
(or just \emph{equivalent}, when the context is clear)
if both sequences converge in all directions,
and
if for all $\uu\in\Rn$,
$\lim (\xx_t\cdot\uu) = \lim (\yy_t\cdot\uu)$,
in other words, if they have the same limits when projected in every
direction.
For the second question above, we will require in our construction
that two convergent sequences have the same limit in the new space if and only if
they are all-directions equivalent.%
\indexg{all-directions equivalence!defined|)}

As a first observation, we note that sequences that are convergent in
$\Rn$ also converge in all directions:

\begin{proposition}  \label{pr:conv-rn-all-dir-conv}
  Let $\seq{\xx_t}$ be a sequence in $\Rn$ and let $\xx\in\Rn$.
  Then $\xx_t\rightarrow\xx$
  if and only if
  for all $\uu\in\Rn$, $\xx_t\cdot\uu\rightarrow\xx\cdot\uu$.
\end{proposition}

\begin{proof}
If $\xx_t\rightarrow\xx$ then $\xx_t\cdot\uu\rightarrow\xx\cdot\uu$,
for all $\uu\in\Rn$, since every linear function is continuous.
Conversely, if $\xx_t\cdot\uu\rightarrow\xx\cdot\uu$
for all $\uu\in\Rn$,
then in particular,
$\xx_t\cdot\ee_i\rightarrow\xx\cdot\ee_i$, for
$i=1,\ldots,n$, implying
$\xx_t\rightarrow\xx$
by \Cref{pr:lim-vec-proj}.
\end{proof}

Here are some examples involving unbounded sequences that tend to infinity:

\begin{example}  \label{ex:basic-lim-ray}
In $\R^2$, let $\xx_t=\trans{[2t,\,t]}=t\vv$, where
$\vv=\trans{[2,1]}$.
The resulting sequence $\seq{\xx_t}$
follows a ray from the origin in the direction of~$\vv$.
Even though this sequence does not converge in $\R^2$,
it does converge in all directions since, for all $\uu\in\R^2$,
\begin{equation}  \label{eq:ex:basic-lim-ray:1}
   \xx_t\cdot\uu
   =
   (t \vv)\cdot\uu
   =
   t (\vv\cdot\uu)
   \rightarrow
   \begin{cases}
     +\infty & \text{if $\vv\cdot\uu>0$,} \\
     -\infty & \text{if $\vv\cdot\uu<0$,} \\
     0       & \text{if $\vv\cdot\uu=0$.}
   \end{cases}
\end{equation}
  The projections of $\seq{\xx_t}$ in three sample
  directions are shown in \Cref{fig:seq-ex5.234} (first row).
  Because this sequence converges in all directions, it will have a
  limit in the new space,
  which will be denoted $\limray{\vv}$ (we define this expression formally in \Cref{sec:astral:space:summary}).
  
Intuitively, ``similar'' sequences that ``eventually follow''
this same ray should have the same limit,
such as the sequences of points
$\xx'_t=\sqrt{t}\vv$ or $\xx''_t=\trans{[2t,\,t+1/t]}$.
Indeed, the images of these sequences under every linear
function have the same limits as those given in
\eqref{eq:ex:basic-lim-ray:1} for $\seq{\xx_t}$.
Therefore, sequences $\seq{\xx'_t}$ and $\seq{\xx''_t}$ are all-directions equivalent to
$\seq{\xx_t}$, and so will have the same limit $\limray{\vv}$.
\end{example}

\begin{example}  \label{ex:basic-with-offset}
Continuing the last example,
let us next consider the sequence $\seq{\yy_t}$ where
$\yy_t=\trans{[2t-1,\,t+2]}=t\vv+\ww$,
where $\ww=\trans{[-1,2]}$.
This sequence is very similar to the sequence $\seq{\xx_t}$.
Both follow along parallel halflines in the direction of $\vv$;
however, these halflines differ in their starting points,
with the $\xx_t\negKern$'s halfline beginning at $\zero$, and the
$\yy_t\negKern$'s at $\ww$.
Does this small difference matter as $t$ gets large?

Note first that $\seq{\yy_t}$ converges in all directions, and so
will have a limit.
Specifically, for $\uu\in\R^2$, we can compute that
\begin{equation}  \label{eq:ex:basic-with-offset:1}
   \yy_t\cdot\uu
   =
   (t \vv + \ww)\cdot\uu
   =
   t (\vv\cdot\uu) + \ww\cdot\uu
   \rightarrow
   \begin{cases}
     +\infty     & \text{if $\vv\cdot\uu>0$,} \\
     -\infty     & \text{if $\vv\cdot\uu<0$,} \\
     \ww\cdot\uu & \text{if $\vv\cdot\uu=0$.}
   \end{cases}
\end{equation}
Sample projections of the sequence $\seq{\yy_t}$
are shown in \Cref{fig:seq-ex5.234} (second row).
For both sequences $\seq{\xx_t}$ and $\seq{\yy_t}$, the growth in the
direction of $\vv$ is so overpowering that the respective limits given
in Eqs.~(\ref{eq:ex:basic-lim-ray:1})
and~(\ref{eq:ex:basic-with-offset:1})
will be the same and infinite for all vectors $\uu\in\R^2$ for which
$\vv\cdot\uu\neq 0$.
However, when $\uu$ is orthogonal to $\vv$, the projections of the
two sequences have different limiting behaviors.
For instance, if
$\uu=\ww=\trans{[-1,2]}$,
then
$\xx_t\cdot\uu\rightarrow 0$
but
$\yy_t\cdot\uu\rightarrow 5$.
In particular, this shows that the sequences are not all-directions
equivalent and so will have different limits in our construction.
\end{example}

\indexg{linear functions, continuous extension of!all-directions equivalence and|(}%
Our requirement that only sequences which are all-directions
equivalent should have the same limit is directly related to the
continuity of linear functions when extended to the new space.
For instance, in the last example, if we were to construct an
extended space in which the sequences $\seq{\xx_t}$ and
$\seq{\yy_t}$ had the same limit, then it would
not be possible to extend the simple function
$f(x_1,x_2)= -x_1 + 2 x_2$
(that is, $f(\xx)=\xx\cdot\uu$, with $\uu=\ww=\trans{[-1,2]}$
as in the example)
continuously to this space.
This is because
if we were to have $\lim \xx_t = \lim \yy_t$
and defined the value of an extended version of $f$ to be equal to any $\alpha\in\eR$
at the common limit of $\seq{\xx_t}$ and $\seq{\yy_t}$,
then at least one of the sequences
$\smash{\seq{f(\xx_t)}}$ or $\smash{\seq{f(\yy_t)}}$ would fail to converge
to $\alpha$, regardless of our choice of~$\alpha$,
since, as we showed in \Cref{ex:basic-with-offset},
$\lim f(\xx_t)=0$, but $\lim f(\yy_t)=5$. Thus, the extended function
cannot be continuous.
Therefore, and more generally, if we want linear functions to have
continuous extensions in the extended space we are constructing, then we need to
treat sequences as having distinct limits if they differ in their
limit in any direction, that is, whenever they are not all-directions
equivalent.%
\indexg{linear functions, continuous extension of!all-directions equivalence and|)}

The next example shows that
astral space will include not only the limits
of sequences along rays or halflines, but also sequences which are
growing in multiple directions at varying rates.

\begin{example}
\label{ex:lim-quadratic}
Continuing the preceding examples,
consider finally the sequence $\seq{\zz_t}$
with elements
$\zz_t= t^2 \vv + t \ww = \trans{\regBracks{2 t^2 - t,\,t^2 + 2t}}$.
Like the sequences $\seq{\xx_t}$ and $\seq{\yy_t}$,
this one is moving to infinity most rapidly in the direction of~$\vv$.
But this sequence is also growing to infinity, though at a slower
rate, in the direction of~$\ww$. (See \Cref{fig:seq-ex5.234}, third row,
where $\uu$ in the last column points in the direction of $\ww$.)
As before, we can compute, for $\uu\in\R^2$, that
\begin{equation*}
   \zz_t\cdot\uu
   =
   (t^2 \vv + t \ww)\cdot\uu
   =
   t^2 (\vv\cdot\uu) + t (\ww\cdot\uu)
   \rightarrow
   \begin{cases}
     +\infty & \text{if $\vv\cdot\uu>0$,} \\
     -\infty & \text{if $\vv\cdot\uu<0$,} \\
     +\infty & \text{if $\vv\cdot\uu=0$ and $\ww\cdot\uu>0$,} \\
     -\infty & \text{if $\vv\cdot\uu=0$ and $\ww\cdot\uu<0$,} \\
     0       & \text{if $\vv\cdot\uu=0$ and $\ww\cdot\uu=0$.}
   \end{cases}
\end{equation*}
Thus, the sequence converges in all directions and so will have a
limit in the new space.
Furthermore, this sequence is not all-directions equivalent to either
$\seq{\xx_t}$ or $\seq{\yy_t}$ since, for instance, if
$\uu=\ww$, then
$\zz_t\cdot\uu\rightarrow+\infty$ which differs from either of the
limits for the other sequences given in
Example~\ref{ex:basic-with-offset}.
Therefore, all three sequences will have distinct limits.%
\indexg{all-directions equivalence|)}%
\indexg{convergence in all directions|)}
\end{example}

\indexg{astral space!defined|(}%
With these notions, we are ready to formally construct astral space.
The space will consist of $\Rn$ together with newly constructed points
at infinity, providing limits for all sequences which converge in all
directions,
while ensuring that two convergent sequences will have the same
limit if and only if they are all-directions equivalent.

Let $\allseq$ denote the set of all sequences in $\Rn$ that converge
in all directions.
Then all-directions equivalence defines an equivalence relation on
$\allseq$.
We therefore can use all-directions equivalence to partition $\allseq$ into equivalence classes, and use those classes to define astral
space.

Let $\Pi$ be the resulting collection of equivalence classes;
that is, the sets in $\Pi$ are all nonempty, their union is all of
$\allseq$, and two sequences are in the same set in $\Pi$ if and only if
they are all-directions equivalent
(implying that all sets in $\Pi$ are disjoint from one another).
Each equivalence class will correspond to a common limit of the
sequences included in that class.
As such, we define astral space so that every point is
effectively identified with one of the equivalence classes of $\Pi$.

For every point $\xx\in\Rn$, note that there must exist a unique
equivalence class in $\Pi$ consisting exactly of all sequences in
$\Rn$ that converge to $\xx$.
This is because
the trivial sequence whose every element is equal to
$\xx$ converges in all directions and so must be included in
some equivalence class.
Moreover, by \Cref{pr:conv-rn-all-dir-conv},
this sequence is all-directions equivalent to any sequence
$\seq{\xx_t}$ in $\Rn$ if and only if $\xx_t\rightarrow\xx$.
Naturally then, we will want to identify $\xx$ with this
class so that $\Rn$ is included in astral space.

Combining these ideas yields the following formal definition:

\begin{definition}  \label{def:astral-space}
\emph{Astral space} is a set denoted $\extspace$%
\indexm{r n 700}{$\extspace$}{astral space}
such that the following hold:
\begin{item-compact}
\item
  There exists a bijection
  \indexm{pi x}{$\pi(\xbar)$}{all-directions equivalence class}%
  $\pi:\extspace\rightarrow\Pi$ identifying
  each element of $\extspace$ with an equivalence class in $\Pi$
  (where $\Pi$ is defined above).
\item
  $\Rn\subseteq\extspace$.
\item
  For all $\xx\in\Rn$, $\pi(\xx)$ is the equivalence class consisting
  of all sequences that converge to $\xx$, establishing the natural
  correspondence discussed above.%
\indexg{astral space!defined|)}
\end{item-compact}
\end{definition}

\indexg{astral space!0- and 1-dimensional@$0$- and $1$-dimensional|(}%
In the special case that $n=1$, we choose $\extspac{1}=\Rext$.
This is possible because, for every $\barx\in\Rext$ (including
$\pm\infty$), there is one equivalence class consisting of all sequences
$\seq{x_t}$ in $\R$ that converge to $\barx$;
naturally, we define $\pi(\barx)$ to be equal to this class.
Furthermore, these are the only equivalence classes in $\Pi$.

When $n=0$, it follows from definitions that
$\extspac{0}=\R^0=\{\zerovec\}$ since the only possible sequence
has every element equal to $\R^0$'s only point, $\zerovec$.%
\indexg{astral space!0- and 1-dimensional@$0$- and $1$-dimensional|)}%
\indexg{astral space!construction|)}

Later, in \Cref{sec:astral-as-fcns},
we will define a natural topology for $\extspace$.
For every point $\xbar\in\extspace$, we will see that, in this
topology, every sequence $\seq{\xx_t}$ in the associated equivalence
class $\pi(\xbar)$ converges to
$\xbar$, so that the astral point $\xbar$ truly can be understood as
the limit of sequences, as previously discussed.
We will also see that $\extspace$ is
indeed a
compactification of~$\Rn$.

\section{The coupling function}

\indexg{coupling function|(}%
Let $\xbar\in\extspace$, and let $\uu\in\Rn$.
By construction,
all of the sequences $\seq{\xx_t}$ in $\pi(\xbar)$
have the same limit $\lim (\xx_t\cdot\uu)$ (and in particular, this limit exists).
\indexg{coupling function!defined|(}%
We use the notation $\xbar\cdot\uu$ to denote this common limit.
That is, we define
\begin{equation}   \label{eqn:coupling-defn}
\indexm{x s100 u}{$\xbar\cdot\uu$}{coupling function}%
  \xbar\cdot\uu
  = \lim (\xx_t\cdot\uu),
\end{equation}
where $\seq{\xx_t}$ is any sequence in $\pi(\xbar)$, noting that the
same value will result regardless of which one is selected.
The operation $\xbar\cdot\uu$,
when viewed as a function mapping $\eRn\times\Rn$ to $\eR$,
is called the \emph{(astral) coupling function}, or simply the \emph{coupling}.%
\indexg{coupling function!defined|)}

Note that $\xbar\cdot\zero=0$ for all $\xbar\in\extspace$.
Furthermore,
if $\xbar=\xx$ for some $\xx\in\Rn$,
then $\xbar\cdot\uu=\xx\cdot\uu$
since
$\xx_t\rightarrow\xx$ for every $\seq{\xx_t}$ in $\pi(\xx)$.
In other words, the coupling notation is compatible with the usual definition
of $\xx\cdot\uu$ as the inner product between $\xx$ and $\uu$.
For similar reasons, if $n=1$ so that $\extspac{1}=\Rext$,
then $\barx\cdot u$, the coupling of $\barx\in\Rext$ and $u\in\R$,
is the same as $\barx u$, their ordinary product, even if
$\barx\in\{-\infty,+\infty\}$.

The coupling function
is critically central to our development;
indeed, we will see that all properties of astral points
$\xbar\in\extspace$ can be expressed in terms of the values of
$\xbar\cdot\uu$ over all $\uu\in\Rn$.
As a start,
these values uniquely determine $\xbar$'s identity:

\begin{proposition}  \label{pr:i:4}
  Let $\xbar$ and $\xbar'$ be in $\extspace$.
  Then $\xbar=\xbar'$ if and only if
  $\xbar\cdot\uu=\xbar'\cdot\uu$
  for all $\uu\in\Rn$.
\end{proposition}

\begin{proof}
Suppose $\xbar\cdot\uu=\xbar'\cdot\uu$ for all $\uu\in\Rn$.
If $\seq{\xx_t}$ is a sequence in $\pi(\xbar)$, then
$\xx_t\cdot\uu\rightarrow \xbar\cdot\uu=\xbar'\cdot\uu$,
for all $\uu\in\Rn$, implying $\seq{\xx_t}$ is also in $\pi(\xbar')$.
Thus, $\pi(\xbar)$ and $\pi(\xbar')$ are equivalence classes with a
nonempty intersection, which implies that they must actually be equal.
Therefore, $\xbar=\xbar'$ since $\pi$ is a bijection.

The reverse implication is immediate.
\end{proof}

Despite the suggestiveness of the notation,
the coupling function
is not actually an inner product.
Nonetheless, it does have some similar properties, as
we show in the next two propositions.
\indexg{coupling function!distributivity|(}%
In the first, we show that it is partially distributive, except when
adding $-\infty$ with $+\infty$ would be involved.

\begin{proposition}  \label{pr:i:1}
  Let $\xbar\in\extspace$ and $\uu,\vv\in\Rn$.
  Suppose
  $\xbar\cdot\uu$ and $\xbar\cdot\vv$
  are summable.
  Then
  \[ \xbar\cdot(\uu+\vv) = \xbar\cdot\uu + \xbar\cdot\vv. \]
\end{proposition}

\begin{proof}
Let $\seq{\xx_t}$ in $\Rn$ be any sequence in $\pi(\xbar)$.
Then
\begin{align*}
  \xbar\cdot\uu + \xbar\cdot\vv
  =
  \lim (\xx_t\cdot\uu) + \lim (\xx_t\cdot\vv)
  &=
  \lim \regBracks{\xx_t\cdot\uu+\xx_t\cdot\vv}
  \\
  &=
  \lim \regBracks{\xx_t\cdot(\uu+\vv)}
  =
  \xbar\cdot(\uu+\vv).
\end{align*}
The first and last equalities are because $\seq{\xx_t}$ is in
$\pi(\xbar)$, and
the second equality is by continuity of addition in $\eR$
(\Cref{prop:lim:eR}\ref{i:lim:eR:sum}).%
\indexg{coupling function!distributivity|)}%
\end{proof}

\indexg{scalar multiples (astral)|(}%
For any point $\xbar\in\extspace$ and scalar $\lambda\in\R$, we define the
\emph{scalar product} $\lambda\xbar$ to be the unique point in $\extspace$ for which
$(\lambda\xbar)\cdot\uu=\lambda(\xbar\cdot\uu)$ for all $\uu\in\Rn$.
The next proposition proves that such a point exists.
Note that when $\xbar=\xx$ is in $\Rn$,
$\lambda\xbar$ is necessarily
equal to the usual product $\lambda\xx$ of scalar $\lambda$ with vector $\xx$.
For the case $\lambda=0$, this proposition (combined with
\Cref{pr:i:4})
 implies $0\xbar=\zero$ for all $\xbar\in\extspace$
(keeping in mind that $0\cdot(\pm\infty)=0$).
By definition, we let $-\xbar=(-1)\xbar$.

\begin{proposition}  \label{pr:i:2}
  Let $\xbar\in\extspace$ and let $\lambda\in\R$.
  Then there exists a unique point in $\extspace$, henceforth denoted
  \indexm{alpha xbar300}{$\alpha\xbar$}{scalar multiple}%
  $\lambda\xbar$, for which
  \[ (\lambda\xbar)\cdot\uu = \lambda (\xbar\cdot\uu) = \xbar\cdot (\lambda\uu) \]
  for all $\uu\in\Rn$.
\end{proposition}

\begin{proof}
Let $\seq{\xx_t}$ in $\Rn$ be any sequence in $\pi(\xbar)$.
Then for all $\uu\in\Rn$,
\begin{align}
  \xbar\cdot (\lambda\uu)
  =
  \lim \bigBracks{\xx_t\cdot (\lambda\uu)}
  &=
  \lim \bigBracks{(\lambda\xx_t) \cdot \uu}
  \nonumber \\
  &=
  \lim \bigBracks{\lambda (\xx_t\cdot\uu)}
  =
  \lambda \lim (\xx_t\cdot\uu)
  =
  \lambda (\xbar\cdot\uu),
  \label{eq:pr:i:2:b}
\end{align}
where the first and last equalities are because $\seq{\xx_t}$ is in 
$\pi(\xbar)$, and the fourth equality
is by continuity of scalar multiplication
(\Cref{prop:lim:eR}\ref{i:lim:eR:mul}).
Thus, 
$\seq{(\lambda \xx_t) \cdot \uu}$ has a limit
for all $\uu\in\Rn$; therefore,
the sequence $\seq{\lambda \xx_t}$ is in $\allseq$.
Since $\Pi$ is a partition of $\allseq$,
this means this sequence is in the equivalence class
$\pi(\ybar)$ for some unique $\ybar\in\extspace$,
implying $\ybar\cdot\uu=\lim[(\lambda\xx_t)\cdot\uu]$
for $\uu\in\Rn$.
Combined with \eqref{eq:pr:i:2:b}
and
defining $\lambda\xbar$ to be $\ybar$, this proves the result.%
\indexg{scalar multiples (astral)|)}
\end{proof}

\indexg{astral points!finite|(}%
Points in $\Rn$ are characterized by being finite in every direction,
or equivalently, by being finite along every coordinate axis.

\begin{proposition}  \label{pr:i:3}
  Let $\xbar\in\extspace$. Then the following are equivalent:
\begin{letter-compact}
\item\label{i:3a}
  $\xbar\in\Rn$.
\item\label{i:3b}
  $\xbar\cdot\uu\in\R$ for all $\uu\in\Rn$.
\item\label{i:3c}
  $\xbar\cdot\ee_i\in\R$ for $i=1,\dotsc,n$.
\end{letter-compact}
\end{proposition}
\begin{proof}
Implications (\ref{i:3a})$\,\Rightarrow\,$(\ref{i:3b}) and (\ref{i:3b})$\,\Rightarrow\,$(\ref{i:3c}) are
immediate. It remains to prove (\ref{i:3c})$\,\Rightarrow\,$(\ref{i:3a}).
Suppose that $\xbar\cdot\ee_i\in\R$ for $i=1,\dotsc,n$. Let $x_i=\xbar\cdot\ee_i$, and
$\xx=\trans{[x_1,\dotsc,x_n]}$.
Then for $\uu\in\Rn$, 
\[
  \xbar\cdot\uu
  =
  \xbar\cdot\Parens{
    \sum_{i=1}^n u_i \ee_i
  }
  =
  \sum_{i=1}^n u_i (\xbar\cdot\ee_i)
  =
  \sum_{i=1}^n u_i x_i = \xx\cdot\uu,
\]
where the second equality is by repeated application of
Propositions~\ref{pr:i:1} and~\ref{pr:i:2}.
Since this holds for all $\uu\in\Rn$, $\xbar$ must be equal to $\xx$,
by \Cref{pr:i:4},
and therefore is in~$\Rn$.%
\indexg{astral points!finite|)}%
\indexg{coupling function|)}%
\end{proof}

\indexg{extended reals!Cartesian product of|(}%
\indexg{astral space!r n@versus $(\Rext)^n$|(}%
Although projections along coordinate axes fully characterize
points in $\Rn$, they are not sufficient to characterize
astral points outside $\Rn$.
In other words, astral space $\extspace$ is distinct from
$(\Rext)^n$, the $n$-fold Cartesian product of $\Rext$ with itself.
Points in either space can be regarded as the limits of possibly
unbounded sequences in $\Rn$.
But astral points embody far more information.
To see this, suppose some sequence $\seq{\xx_t}$ in $\Rn$
converges to a point $\xhat$ in $(\Rext)^n$.
Then $\xhat$, by its form, encodes the limit of
$(\xx_t\cdot\ee_i)$, for each standard basis vector $\ee_i$,
since this limit is exactly the $i$-th component of $\xhat$.
In comparison, if instead the sequence $\seq{\xx_t}$ is in the equivalence
class of an
astral point $\xbar$ in $\extspace$, then $\xbar$ encodes the limit
of $(\xx_t\cdot\uu)$ for \emph{all} vectors $\uu\in\Rn$, not just
the standard basis vectors.
For instance, if $n=2$ and $\xhat=\trans{[+\infty,+\infty]}$
in $(\Rext)^2$,
then $\xx_t\cdot\ee_i\rightarrow+\infty$, for $i=1,2$.
From this information,
if $\uu=\ee_1-\ee_2$,
for example,
there is no way to deduce the limit of
$\seq{\xx_t\cdot\uu}$, or even if the limit exists.
On the other hand, this limit, or the limit for any other vector $\uu\in\R^2$,
would be readily available from $\xbar\in\extspac{2}$.

Thus, as limits of sequences,
astral points $\xbar\in\extspace$
retain all of the information embodied
by points $\xhat$ in $(\Rext)^n$, and
usually far more.
This results in astral space having a remarkably rich, powerful structure,
as will be developed throughout this work.%
\indexg{extended reals!Cartesian product of|)}%
\indexg{astral space!r n@versus $(\Rext)^n$|)}

\section{Astral points}
\label{subsec:astral-pt-form}

What kinds of points comprise astral space?
The space includes all of $\Rn$, of course.
All of the other ``new'' points correspond to equivalence classes
of unbounded sequences
$\seq{\xx_t}$ in $\Rn$
(that is, for which $\norm{\xx_t}\rightarrow+\infty$).
This follows from
\Cref{pr:i:3}(\ref{i:3a},\ref{i:3b})
which shows that if
$\xbar\in\extspace\setminus\Rn$ then $\xbar\cdot\uu$ must be infinite
for some $\uu\in\Rn$, implying $\xx_t\cdot\uu$ converges to
$\pm\infty$ for any sequence $\seq{\xx_t}$ in $\pi(\xbar)$.
For this reason, we say astral
points in $\Rn$ are
\indexg{astral points!finite}%
\emph{finite}, and that all others
are \emph{infinite}.%
\indexg{astral points!infinite}%
\indexg{infinite points}

In later chapters,
we will see that all points in $\extspace$ have a
specific structure, as suggested by the examples in
\Cref{subsec:astral-intro-motiv}.
Every astral point $\xbar$, outside those in $\Rn$, corresponds
to the equivalence class of
sequences which have a particular dominant direction $\vv$ in which they
grow to infinity most rapidly.
In addition, these sequences may be growing to infinity in other
directions that are secondary, tertiary,~etc.
These sequences may also have a finite part in the sense of converging
to a finite value in some directions.
Importantly, the details of this structure are entirely determined by
the point $\xbar$ itself so that every sequence that converges to
$\xbar$ will have the same dominant direction, same finite part,~etc.

\indexg{Polynomially graded sequence|(}%
To demonstrate this structure,
we introduced, in Example~\ref{ex:poly-speed-intro}, the polynomially graded sequence
\begin{equation}
\label{eq:h:6}
  \xx_t = t^k \vv_1 + t^{k-1} \vv_2 + \dotsb + t \vv_k + \qq,
\end{equation}
where
$\qq,\vv_1,\dotsc,\vv_k\in\Rn$. We argued that $\seq{\xx_t}$ converges in all
directions and is in the equivalence class of the unique astral point $\exx$ that satisfies, for $\uu\in\Rn$,
\begin{equation}
\label{eqn:pf:2}
\newcommand{\ttinprod}{\inprod}
\xbar\ttinprod\uu
=
\begin{cases}
+\infty
&\text{if $\vv_1\ttinprod\uu=\dotsb=\vv_{i-1}\ttinprod\uu=0$ and $\vv_i\ttinprod\uu>0$ for some $i$,}
\\
-\infty
&\text{if $\vv_1\ttinprod\uu=\dotsb=\vv_{i-1}\ttinprod\uu=0$ and $\vv_i\ttinprod\uu<0$ for some $i$,}
\\
\qq\ttinprod\uu
&\text{if $\vv_1\ttinprod\uu=\dotsb=\vv_k\ttinprod\uu=0$.}
\end{cases}
\end{equation}
The same astral point $\xbar$ would be obtained if polynomial coefficients in
\eqref{eq:h:6} were replaced by any coefficient sequences going to infinity at decreasing rates and $\qq$ were replaced
by any sequence of $\qq_t\negKern$'s converging to $\qq$ (see \Cref{thm:i:seq-rep}). Under natural conditions, this turns out to characterize all the sequences in the equivalence class of $\xbar$ (see \Cref{thm:seq-rep}).

Moreover, as we show in \Cref{sec:char:func},
every point in astral space must have exactly the form of \eqref{eqn:pf:2}, for some
$\qq,\vv_1,\dotsc,\vv_k\in\Rn$.
This specific structure is fundamental to our understanding of $\extspace$,
and will be central to the foundation on which all of the later results
are based.%
\indexg{Polynomially graded sequence|)}%

\chapter{Astral topology}
\label{sec:astral-as-fcns}

We next study a topology for astral space.
We will see that astral space, in this topology, is compact,
as was the goal of our construction,
and we will also see that it is first-countable, which
is another useful property, ensuring
that
we can reason about convergence and continuity in this space using sequences, similar to $\Rn$.
Nonetheless, we will see that astral space is not
second-countable and not metrizable (for $n\geq 2$).

In \Cref{sec:astral:construction},
astral space was constructed to consist of points representing
equivalence classes of sequences.
We begin this chapter by presenting
an alternative view of the space in which
astral points are in one-to-one correspondence with certain functions.
This perspective will allow us
to derive a topology for the space, and to prove properties like
compactness.
It will also allow us to develop some fundamental operations on astral
space, such as leftward addition.

\section{Astral points as functions}
\label{subsec:astral-pts-as-fcns}

\indexg{functional representation of astral points|(}%
Every point $\xbar\in\extspace$
defines
a function $\ph{\xbar}:\Rn\rightarrow\Rext$
corresponding to the evaluation of the coupling function with the first argument
fixed to $\xbar$, namely,
\begin{equation}  \label{eq:ph-xbar-defn}
\indexm{phi x}{$\ph{\xbar}$}{functional representation of a point}%
\indexg{functional representation of astral points!defined}%
  \ph{\xbar}(\uu) = \xbar\cdot\uu
\end{equation}
for $\uu\in\Rn$.
Every such function is included in
the set of all
functions mapping $\Rn$ to $\Rext$,
\[
\indexm{f 100}{$\fcnspn$}{all functions on $\Rn$ to $\Rext$}
  \fcnspn=\Rext^{\Rn},
\]
which we endow with the product topology,
that is, the topology of pointwise convergence,
as in \Cref{sec:prod-top}.
Since $\Rext$ is compact
(Example~\ref{ex:rext-compact}),
$\fcnspn$ is also compact
by Tychonoff's theorem
(\Cref{pr:prod-top-props}\ref{pr:prod-top-props:b}).

We consider two topological subspaces of $\fcnspn$. The first consists of the functions $\ph{\xx}$ defined by $\xx\in\Rn$, while the second, larger one, consists of the functions $\ph{\xbar}$ defined by $\xbar\in\eRn$:
\begin{align*}
  \phimg&=
  \set{\ph{\xx}:\: \xx\in\Rn},
\\
  \phimgA&=
  \set{\ph{\xbar}:\: \xbar\in\extspace}.
\end{align*}

\begin{figure}
  \centering
  \includegraphics{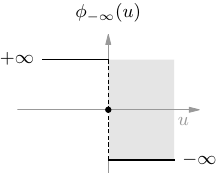}\hfill%
  \includegraphics{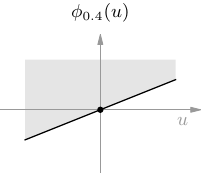}\hfill%
  \includegraphics{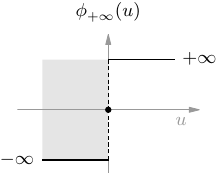}
  \mycaption{%
    Functions $\ph{\barx}$ in one dimension}{%
    \indexf{functional representation of astral points}%
    For $\barx\in\eR$, we define $\ph{\barx}(u)=\barx u$ for all $u\in\R$.
    If $\barx=x\in\R$, then $\ph{\barx}$ is a standard linear function
    (center),
    but if $\barx\in\set{-\infty,+\infty}$, then $\ph{\barx}$ is improper
    (left and right).
  }%
  \label{fig:1d-linear}%
\end{figure}

\begin{example}  \label{ex:fcn-classes-n=1}
When $n=1$, $\fcnspn$ consists of all functions
$\psi:\R\rightarrow\Rext$.
The set $\phimg$ includes all functions $\ph{x}$, for $x\in\R$,
where $\ph{x}(u)=xu$ for $u\in\R$, in other words, all lines that pass
through the origin with arbitrary (finite) slope $x$.
The set $\phimgA$ includes all of these functions as well as
$\ph{\barx}$ for $\barx\in\{-\infty,+\infty\}$.
When $\barx=+\infty$, this function is given, for $u\in\R$, by
\begin{equation}  \label{eq:ex:fcn-classes-n=1:1}
  \ph{+\infty}(u)
  =
  (+\infty)\cdot u
  =
  \begin{cases}
               +\infty & \text{if $u>0$,}\\
               0       & \text{if $u=0$,}\\
               -\infty & \text{if $u<0$.}
  \end{cases}
\end{equation}
Symmetrically, $\ph{-\infty}=(-\infty)\cdot u$.
These two functions can be viewed as having an infinite
(positive or negative) slope at~$0$. They are depicted in \Cref{fig:1d-linear}.%
\indexg{functional representation of astral points|)}
\end{example}

\indexg{astral topology!construction|(}%
Astral space, $\extspace$, is in natural correspondence with
$\phimgA$, as given by the map
$\phfcn:\extspace\to\phimgA$
where
\[ \phfcn(\xbar)=\ph{\xbar} \]
for $\xbar\in\extspace$.
This map is a bijection by \Cref{pr:i:4}.
Furthermore,
we will soon define a topology on $\extspace$ under which $\phfcn$ is
also a homeomorphism, so topological properties proved for $\phimgA$
will apply to $\extspace$ as well.
In particular, we will prove below that $\phimgA$ is closed, and
therefore compact (since $\fcnspn$ is compact), which will show that
$\extspace$ is compact as well.

Along the way to deriving topological properties of~$\phimgA$,
we will also prove that functions in $\phimgA$ have specific
convexity properties and a particular functional form.
When translated back to $\extspace$, this will characterize the structural
form of all points in $\extspace$ as discussed
in \Cref{subsec:astral-pt-form}.

As already mentioned, we assume a product
topology on $\fcnspn$.
Following the definition of product topology from \Cref{sec:prod-top}
(specifically, Eq.~\ref{eq:fcn-cl-gen-subbase}),
and because the topology on $\Rext$ is generated by
a subbase consisting of sets of the form
$[-\infty,b)$ and $(b,+\infty]$ for $b\in\R$
(\Cref{ex:topo-rext}),
we obtain that a subbase for the product
topology on $\fcnspn$ consists of all sets of the form
\[
   \set{\psi\in\fcnspn:\: \psi(\uu) < b}
   \mbox{~~or~~}
   \set{\psi\in\fcnspn:\: \psi(\uu) > b},
\]
where $b\in\R$, $\uu\in\Rn$.
Or, more succinctly, this subbase consists of all sets
\begin{equation}  \label{eq:fcnspn-subbase}
   \set{\psi\in\fcnspn:\: s\psi(\uu) < b},
\end{equation}
where $s\in\set{-1,+1}$, $b\in\R$, $\uu\in\Rn$.
In this topology, as is generally true for the product topology,
if $\seq{\psi_t}$ is any sequence in $\fcnspn$,
and $\psi\in\fcnspn$, then
$\psi_t\rightarrow\psi$ if and only if
$\psi_t(\uu)\rightarrow\psi(\uu)$ for all $\uu\in\Rn$
(\Cref{pr:prod-top-ptwise-conv}).

From \eqref{eq:fcnspn-subbase},
the topology on $\phimgA$ as a subspace of $\fcnspn$ is
generated by a subbase consisting of all sets of the form
\begin{align}
\notag
  \regBraces{\psi\in\phimgA:\: s \psi(\uu) < b}
  &=
  \regBraces{\ph{\xbar}:\: \xbar\in\extspace,\, s\ph{\xbar}(\uu) < b}
\\
\notag
  &=
  \regBraces{\ph{\xbar}:\: \xbar\in\extspace,\, s(\xbar\cdot\uu) < b}
\\
\label{eq:h:3:sub:phimgA}
  &=
  \phfcn\bigParens{
    \regBraces{\xbar\in\extspace :\: s(\xbar\cdot\uu) < b}
  },
\end{align}
where $s\in\set{-1,+1}$, $b\in\R$, $\uu\in\Rn$.

We noted that $\phfcn$ defines a bijection between
$\extspace$ and $\phimgA$.
To ensure that $\phfcn$ defines a homeomorphism as well,
we define our topology for astral space
so that a set
$U\subseteq\extspace$ is open if and only if its image $\phfcn(U)$ is
open in~$\phimgA$. This can be ensured by defining a subbase
for $\extspace$ with those elements $S\subseteq\extspace$
whose images $\phfcn(S)$ form a subbase for~$\phimgA$.
From \eqref{eq:h:3:sub:phimgA}, such a subbase consists of elements of
the form given in the next definition.
For the entirety of this work, we will always assume this
topology for astral space, unless explicitly stated otherwise.

\begin{definition}  \label{def:astral-topology}
\indexg{astral topology!defined|(}%
The \emph{astral topology} on $\extspace$ is the topology generated
by a subbase consisting of all sets
\begin{equation}
\label{eq:h:3a:sub}
  \set{\xbar\in\extspace:\: s (\xbar\cdot\uu) < b},
\end{equation}
where $s\in\{-1,+1\}$, $\uu\in\Rn$, and $b\in\R$.%
\indexg{astral topology!defined|)}
\end{definition}

By \Cref{pr:i:2},
we have $s(\xbar\cdot\uu)=\xbar\cdot(s\uu)$, so, when convenient,
we can take $s=+1$ when working with sets as in
\eqref{eq:h:3a:sub}.
Furthermore, when $\uu=\zero$, the set in
\eqref{eq:h:3a:sub}
is either the empty set or the
entire space $\extspace$.
Both of these can be discarded from the
subbase without changing the generated topology
since $\emptyset$ and $\extspace$ will be
included anyway in the topology as, respectively,
the empty union and empty intersection of subbase elements.
Thus, simplifying \Cref{def:astral-topology},
we can say that the astral topology is equivalently generated by a
subbase consisting of all sets
\begin{equation}
\label{eq:h:3a:sub-alt}
  \set{\xbar\in\extspace:\: \xbar\cdot\uu < b},
\end{equation}
where $\uu\in\Rn\wo\{\zero\}$ and $b\in\R$.

Moreover, by taking finite intersections,
this subbase generates a base consisting of all sets
\begin{equation}
\label{eq:h:3a}
  \set{\xbar\in\extspace:\: \xbar\cdot\uu_i < b_i \text{ for all }
                             i=1,\dotsc,k},
\end{equation}
for some finite $k\geq 0$,
and some $\uu_i\in\Rn\wo\{\zero\}$ and $b_i\in\R$,
for $i=1,\dotsc,k$.

When $n=1$,
the astral topology for $\extspac{1}$ is the same as the
usual topology for $\Rext$ since, in this case,
the subbase elements given in
\eqref{eq:h:3a:sub} coincide with those
in Example~\ref{ex:topo-rext}
for $\Rext$ (along with $\emptyset$ and $\Rext$).

\indexg{linear functions, continuous extension of!astral space@in astral space|(}%
\indexg{weak topology!astral space and|(}%
\indexg{astral topology!weak topology@as weak topology|(}%
\Cref{def:astral-topology} says that astral topology is precisely the weak topology
(as defined in \Cref{sec:prod-top})
on $\eRn$ with respect to the family of all
maps $\xbar\mapsto\xbar\inprod\uu$, where $\uu\in\Rn$.
In other words, astral topology is the coarsest topology on
$\extspace$ under which all such maps are continuous.
Since these maps are the natural astral extensions of standard
linear functions on $\Rn$,
this choice of topology
is in line with our goal
of constructing astral space
in a way that allows all linear functions to be extended continuously.%
\indexg{linear functions, continuous extension of!astral space@in astral space|)}%
\indexg{astral topology!weak topology@as weak topology|)}%
\indexg{weak topology!astral space and|)}%
\indexg{astral topology!construction|)}%

\begin{proposition}  \label{pr:astral-topo-homeo}
  The map $\phfcn:\extspace\to\phimgA$ is a homeomorphism.
\end{proposition}

\begin{proof}
\Cref{pr:i:4} directly implies
that $\phfcn$ is a bijection.

To see that $\phfcn$ is continuous, let $V$ be a
subbase element of $\phimgA$, as in \eqref{eq:h:3:sub:phimgA},
so that $V=\phfcn(S)$ for some subbase element $S$ of $\extspace$
as in \eqref{eq:h:3a:sub}.
Then $\phfcninv(V)=S$, which is open.
Therefore, $\phfcn$ is continuous by
\Cref{prop:cont}(\ref{prop:cont:sub},\ref{prop:cont:a}).

Similarly, for every subbase element $S$ of $\extspace$, as in
\eqref{eq:h:3a:sub}, the set $\phfcn(S)$ is a
subbase element of $\phimgA$, according to \eqref{eq:h:3:sub:phimgA},
and is therefore open.
Thus, $\phfcninv$ is continuous as well
(again by
\Cref{prop:cont}\ref{prop:cont:sub}\ref{prop:cont:a}).
\end{proof}

\indexg{astral topology!Euclidean space as subspace of|(}%
\indexg{Euclidean space!subspace of astral space@as subspace of astral space|(}%
As we show next,
the topology on $\Rn$ as a subspace of $\extspace$
(in the astral topology)
is the same as the standard Euclidean topology on $\Rn$.
Furthermore, any set $U\subseteq\Rn$ is open in the Euclidean topology
on~$\Rn$ if and only if it is open in the astral topology on~$\eRn$.
For the purposes of this proposition, we say that a set $U$ is open in
$\extspace$ if it is open in the astral topology, and open in
$\Rn$ if it is open in the Euclidean topology.

\begin{proposition}   \label{pr:open-sets-equiv-alt}
  Let $U\subseteq\Rn$.
  Then the following are equivalent:
  \begin{letter-compact}
  \item     \label{pr:open-sets-equiv-alt:a}
    $U$ is open in $\Rn$.
  \item     \label{pr:open-sets-equiv-alt:b}
    $U$ is open in $\extspace$.
  \item     \label{pr:open-sets-equiv-alt:c}
    $U=V\cap\Rn$ for some set $V\subseteq\extspace$ that is open in
    $\extspace$.
    (That is, $U$ is open in the topology on $\Rn$ as a subspace of
    $\extspace$.)
  \end{letter-compact}
  Consequently,
  the topology on $\Rn$ as a subspace of $\extspace$ is the same as
  the Euclidean topology on $\Rn$.
\end{proposition}

\begin{proof}
We assume $U\neq\emptyset$ since otherwise
the claimed equivalence is trivial.

\begin{proof-parts}
\pfpart{%
  (\ref{pr:open-sets-equiv-alt:a})
  $\Rightarrow$
  (\ref{pr:open-sets-equiv-alt:b}):
}
Suppose $U$ is open in $\Rn$.
Let $\yy\in U$.
Then there exists $\epsilon>0$ such that
$\ball(\yy,\epsilon) \subseteq U$.
Let $C$ be the following base element (in $\extspace$):
\[
   C = \BigBraces{
     \xbar\in\extspace :\:
     |\xbar\cdot\ee_i - \yy\cdot\ee_i| < \delta
     \text{ for } i=1,\dotsc,n
   }
\]
where
$\delta = \epsilon/\sqrt{n}$.
Then $C\subseteq\Rn$ since
if $\xbar\in C$, then we must have $\xbar\cdot\ee_i\in\R$, for
$i=1,\dotsc,n$, implying $\xbar\in\Rn$ by
\Crefequiv{pr:i:3}{i:3c}{i:3a}.
Further, $C\subseteq\ball(\yy,\epsilon)$ since if
$\xx\in C\subseteq\Rn$ then
$|(\xx-\yy)\cdot\ee_i|<\delta$, for $i=1,\dotsc,n$, implying
$\norm{\xx-\yy}<\delta\sqrt{n} = \epsilon$.
Thus, $\yy\in C\subseteq\ball(\yy,\epsilon)\subseteq U$,
proving $U$ is open in $\extspace$
(by \Cref{pr:base-equiv-topo}).

\pfpart{%
  (\ref{pr:open-sets-equiv-alt:b})
  $\Rightarrow$
  (\ref{pr:open-sets-equiv-alt:c}):
}
Since $U=U\cap\Rn$, if $U$ satisfies
(\ref{pr:open-sets-equiv-alt:b}),
then it also satisfies
(\ref{pr:open-sets-equiv-alt:c}).

\pfpart{%
  (\ref{pr:open-sets-equiv-alt:c})
  $\Rightarrow$
  (\ref{pr:open-sets-equiv-alt:a}):
}
Suppose $U=V\cap\Rn$ for some $V\subseteq\extspace$ that is open in
$\extspace$.
Let $\yy\in U$.
Then $\yy$ is also in $V$, which is open in $\extspace$, so
by \Cref{pr:base-equiv-topo},
there exists a base element $C$ as in \eqref{eq:h:3a} such that
$\yy\in C \subseteq V$.
Since each $\uu_i$ in that equation is not $\zero$, we can assume
further without loss of generality that $\norm{\uu_i}=1$
(since otherwise we can divide both sides of the inequality by
$\norm{\uu_i}$).
Let $\epsilon\in\Rstrictpos$ be such that
\[   \epsilon < \min \bigBraces{b_i - \yy\cdot\uu_i :\: i=1,\dotsc,k},  \]
which must exist since $\yy\in C$.

Suppose $\xx\in \ball(\yy,\epsilon)$.
Then
for $i=1,\dotsc,k$,
\[
   \xx\cdot\uu_i
   =
   \yy\cdot\uu_i + (\xx-\yy)\cdot\uu_i
   <
   \yy\cdot\uu_i + \epsilon
   <
   \yy\cdot\uu_i + (b_i - \yy\cdot\uu_i)
   =
   b_i.
\]
The first inequality is by the Cauchy-Schwarz inequality
since $\xx\in \ball(\yy,\epsilon)$,
and the second is by our choice of $\epsilon$.
Thus, $\xx\in C$.
Therefore,
$\ball(\yy,\epsilon)\subseteq C\cap\Rn \subseteq V\cap\Rn = U$,
proving $U$ is open in $\Rn$
(again by \Cref{pr:base-equiv-topo}).

\pfpart{Equivalent topologies:}
The equivalence of~(\ref{pr:open-sets-equiv-alt:a})
and~(\ref{pr:open-sets-equiv-alt:c})
shows that
these two topologies are identical.%
\indexg{astral topology!Euclidean space as subspace of|)}%
\indexg{Euclidean space!subspace of astral space@as subspace of astral space|)}
\qedhere
\end{proof-parts}
\end{proof}

\Cref{pr:open-sets-equiv-alt}(\ref{pr:open-sets-equiv-alt:a},\ref{pr:open-sets-equiv-alt:b})
shows that a set in $\Rn$ is
open in the astral topology on $\eRn$
if and only if it
is open in the standard Euclidean topology on $\Rn$.
For this reason, we can simply refer to a set in $\Rn$ as being open
without specifying which of these two topologies we mean.

The same is \emph{not} true for closed sets:
A set in $\Rn$ can be closed with respect to
the Euclidean topology, but not closed in
the astral topology.
For instance, $\Rn$ itself is closed in the former,
but not in the latter.
When necessary, to make clear which topology is being
referred to, 
we therefore use the phrase \emph{closed in $\Rn$} to refer to Euclidean
topology,
or \emph{closed in $\extspace$} for astral topology.
Often, however, the topology will be clear from context, especially
since this ambiguity only arises when considering sets that are
necessarily in $\Rn$.

\indexg{closure, astral|(}%
\indexg{closure (topological)!astral versus Euclidean|(}%
In the same way, the closure of a set in $\Rn$ can in general be
different if the closure is with respect to the Euclidean topology on
$\Rn$ or the astral topology on $\extspace$.
Note importantly that we therefore adopt the notational convention of
writing $\Sbar$%
\indexm{s 300}{$\Sbar$}{closure in $\eRn$}
for the closure in $\eRn$ of any set $S\subseteq\eRn$
(including subsets of $\Rn$),
while continuing to write
$\cl S$ for the closure in $\Rn$ of a set $S\subseteq\Rn$.

Here are some useful facts about closures and closed sets, all
standard results specialized to the current setting:

\begin{proposition}  \label{pr:closed-set-facts}
  \mbox{}
  \begin{letter-compact}
  \item  \label{pr:closed-set-facts:a}
    Let $S \subseteq \Rn$.
    Then $\cl S = \Sbar \cap \Rn$.
  \item  \label{pr:closed-set-facts:aa}
    Let $S \subseteq \Rn$.
    Then $\clbar{\cl S} = \Sbar$.
  \item  \label{pr:closed-set-facts:b}
    Let $U \subseteq \extspace$ be open.
    Then $\clbar{(U \cap \Rn)} = \Ubar$.
  \end{letter-compact}
\end{proposition}

\begin{proof}~
\begin{proof-parts}
\pfpart{Part~(\ref{pr:closed-set-facts:a}):}
Since $\Rn$ in the Euclidean topology is a subspace of $\extspace$
(by
\Cref{pr:open-sets-equiv-alt}\ref{pr:open-sets-equiv-alt:a}\ref{pr:open-sets-equiv-alt:c}),
this is a special case of
\Cref{prop:subspace}(\ref{i:subspace:closure}).

\pfpart{Part~(\ref{pr:closed-set-facts:aa}):}
On the one hand, $\Sbar\subseteq \clbar{\cl S}$
since $S\subseteq \cl S$.
On the other hand,
$\cl S = \Sbar\cap\Rn\subseteq \Sbar$
by part~(\ref{pr:closed-set-facts:a}).
Since $\Sbar$ is closed in $\extspace$, this implies
$\clbar{\cl S} \subseteq \Sbar$.

\pfpart{Part~(\ref{pr:closed-set-facts:b}):}
First, $\clbar{(U \cap \Rn)} \subseteq \Ubar$
since $U \cap \Rn \subseteq U$.
To prove the reverse inclusion, suppose $\xbar\in\Ubar$.
Let $V$ be any neighborhood of $\xbar$.
Then by
\Cref{pr:closure:intersect}(\ref{pr:closure:intersect:a}),
there exists a point $\ybar$ in the open set $V\cap U$,
which in turn implies, since
$\Rn$ is dense in $\extspace$, that
$V\cap U\cap \Rn\neq\emptyset$.
Thus, $U\cap\Rn$ intersects every neighborhood $V$ of $\xbar$.
Therefore,
$\xbar\in \clbar{(U \cap \Rn)}$,
again by \Cref{pr:closure:intersect}(\ref{pr:closure:intersect:a}).%
\indexg{closure, astral|)}%
\indexg{closure (topological)!astral versus Euclidean|)}
\qedhere
\end{proof-parts}
\end{proof}

\section{Characterizing functional form}
\label{sec:char:func}

Let $\phimgcl$ denote the closure of $\phimg$ in $\fcnspn$.
We will study the properties of $\phimgcl$ and of the
functions comprising it.
Eventually, we will prove that $\phimgcl=\phimgA$, which will
imply that $\phimgA$ is closed and
therefore compact (so that $\extspace$ is as well).

\indexg{functional representation of astral points!characterizations|(}%
Additionally, we will show that $\phimgA$
consists exactly of those functions
$\psi:\Rn\rightarrow\Rext$
that are convex, concave,
and that vanish at the origin (meaning $\psi(\zero)=0$).
Indeed, for every $\xx\in\Rn$,
the linear function $\ph{\xx}(\uu)=\xx\cdot\uu$
is convex, concave, and vanishes at the origin,
which shows that every function in
$\phimg$
has these
properties.
Nonetheless, these are not the only ones;
there are other functions, all improper,
that have these properties as well.
For instance,
the functions $\phi_{+\infty}$ and $\phi_{-\infty}$ from
\Cref{ex:fcn-classes-n=1} (shown in \Cref{fig:1d-linear})
are examples of such functions.
Here is another:

\begin{figure}
  \centering
  \figclifframp{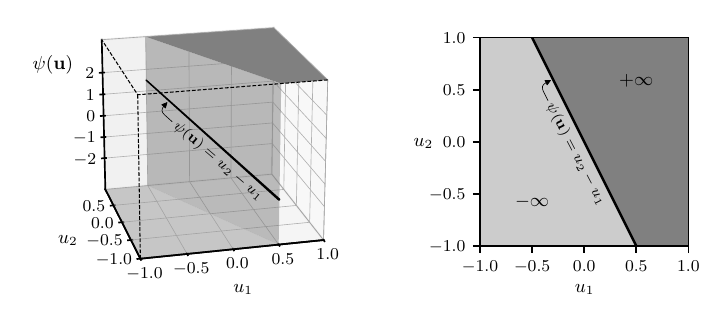}
  \mycaption{Cliff ramp}{%
    \indexf{Cliff ramp}%
    The function $\psi$ from \Cref{ex:cliff-ramp}.
    This function coincides with
    the linear function $u_2-u_1$
    along the line $2u_1+u_2=0$,
    but equals
    $-\infty$ on one side of that line, and $+\infty$ on the other side.
  }%
  \label{fig:cliff-ramp}%
\end{figure}

\begin{example}[Cliff ramp]
\label{ex:cliff-ramp}
\indexg{Cliff ramp|(}%
In $\R^2$, let
\begin{align}
\nonumber
  \psi(\uu)
  =
  \psi(u_1,u_2)
  &=
  \begin{cases}
             -\infty & \text{if $2u_1+u_2<0$},
  \\
             u_2-u_1 & \text{if $2u_1+u_2=0$},
  \\
             +\infty & \text{if $2u_1+u_2>0$}.
  \end{cases}
\end{align}
As shown in \Cref{fig:cliff-ramp},
this function coincides with
the linear function $u_2-u_1$
along the line $2u_1+u_2=0$
(the ``ramp'' portion of the function),
but equals
$-\infty$ on one side of that line, and $+\infty$ on the other side
(the ``cliff'').
It can be checked that the epigraphs of both $\psi$ and $-\psi$ are
convex sets.
Therefore, the function
is convex, concave, and satisfies $\psi(\zero)=0$.%
\indexg{Cliff ramp|)}
\end{example}

We begin our development with the following inclusion:

\begin{proposition}  \label{pr:h:7}
  $\phimgA\subseteq\phimgcl$.
\end{proposition}

\begin{proof}
Let $\xbar\in\extspace$.
We aim to show that $\ph{\xbar}$, which is an arbitrary element of
$\phimgA$, is included in $\phimgcl$.
Let $\seq{\xx_t}$ be any sequence in $\pi(\xbar)$. Then,
from the definition of coupling function (Eq.~\ref{eqn:coupling-defn}),
we have
$\xx_t\cdot\uu\rightarrow\xbar\cdot\uu$ for all $\uu\in\Rn$,
or equivalently,
$\ph{\xx_t}(\uu)\rightarrow\ph{\xbar}(\uu)$ for all $\uu\in\Rn$.
By \Cref{pr:prod-top-ptwise-conv},
this implies that
$\ph{\xx_t}\rightarrow\ph{\xbar}$.
Since $\ph{\xx_t}\in\phimg$ for all $t$, this shows that
$\ph{\xbar}\in\phimgcl$
(by \Cref{prop:first:properties}\ref{prop:first:closure}).
\end{proof}

The next theorem gives some properties of
functions in $\phimgcl$, in particular, showing that they all vanish
at the origin, and that the set is closed under negation:

\begin{theorem}  \label{pr:h:2}
  Let $\psi\in\phimgcl$.
  Then:
  \begin{letter-compact}
  \item  \label{pr:h:2:a}
    $\psi(\zero)=0$.
  \item  \label{pr:h:2:b}
    $-\psi\in\phimgcl$.
  \end{letter-compact}
\end{theorem}

\begin{proof}
  ~

\begin{proof-parts}
\pfpart{Part~(\ref{pr:h:2:a}):}
Let $C=\{\xi\in\fcnspn :\: \xi(\zero)=0\}$.
This set is closed,
since it is the complement of the set
\[
  \{\xi\in\fcnspn:\:\xi(\zero) < 0\}
   \,\cup\,
  \{\xi\in\fcnspn:\:\xi(\zero) > 0\},
\]
which is a union of two subbase elements and therefore open.
Since $\phx\in C$ for all $\xx\in\Rn$, we have
$\phimg\subseteq C$, implying
$\phimgcl\subseteq C$ since $C$ is closed.
Therefore,
$\psi(\zero)=0$
since $\psi\in\phimgcl$.

\pfpart{Part~(\ref{pr:h:2:b}):}
Let $\eta:\fcnspn\rightarrow\fcnspn$ be defined by
$\eta(\xi)=-\xi$ for $\xi\in\fcnspn$.
Then $\eta$ is continuous since, for any subbase element
$S = \set{\xi\in\fcnspn:\: s\xi(\uu) < b}$
as in \eqref{eq:fcnspn-subbase},
where
$s\in\set{-1,+1}$, $b\in\R$, $\uu\in\Rn$,
we have $\eta^{-1}(S)={\set{\xi\in\fcnspn:\: -s\xi(\uu) < b}}$,
which is open, being also a subbase element.
Therefore, $\eta$ is continuous by
\Cref{prop:cont}(\ref{prop:cont:sub},\ref{prop:cont:a}).

We then have
\[
  -\psi
  =
  \eta(\psi)
  \in
  \eta(\phimgcl)
  \subseteq
  \clbar{\eta(\phimg)}
  \subseteq
  \phimgcl.
\]
The first inclusion is because $\psi\in\phimgcl$.
The second is by 
\Cref{prop:cont}(\ref{prop:cont:a},\ref{prop:cont:c})
since $\eta$ is continuous.
The last is because $\eta(\phimg)\subseteq\phimg$
since for all $\xx\in\Rn$, $\eta(\phx)=-\phx=\ph{-\xx}\in\phimg$.
\qedhere
\end{proof-parts}
\end{proof}

Next, we show that all functions in $\phimgcl$
are both convex and concave.

\begin{example}[Cliff ramp, continued]
\label{ex:cliff-ramp-cont}
\indexg{Cliff ramp|(}%
The cliff ramp function $\psi$ from
Example~\ref{ex:cliff-ramp}
can be written as
\[
  \psi(\uu)
  =\begin{cases}
              -\infty     & \text{if $\vv\cdot\uu<0$,} \\
              \qq\cdot\uu & \text{if $\vv\cdot\uu=0$,} \\
              +\infty     & \text{if $\vv\cdot\uu>0$,}
   \end{cases}
\]
where $\vv=\trans{[2,1]}$ and $\qq=\trans{[-1,1]}$. This
is a special case of the right-hand side of \eqref{eqn:pf:2},
which resulted from the analysis
of the polynomially graded
sequence in \eqref{eq:h:6}. Following that example,
we obtain that for the sequence $\xx_t=t\vv+\qq$, we must have $\xx_t\cdot\uu\rightarrow\psi(\uu)$
for all $\uu\in\R^2$.
Thus, $\ph{\xx_t}\to\psi$, and therefore, $\psi\in\phimgcl$.
The next theorem shows this implies that $\psi$
is convex and concave, as was previously noted in Example~\ref{ex:cliff-ramp}.%
\indexg{Cliff ramp|)}
\end{example}

\begin{theorem}  \label{thm:h:3}
  Let $\psi\in\phimgcl$.
  Then $\psi$ is both convex and concave.
\end{theorem}

\begin{proof}
To show that $\psi$ is convex, we prove that it satisfies the condition in
\Cref{pr:stand-cvx-fcn-char}(\ref{pr:stand-cvx-fcn-char:b}).
Let $\uu,\vv\in\dom\psi$ and let $\lambda\in [0,1]$.
We aim to show that
\begin{equation}  \label{eq:thm:h:3:1}
  \psi(\ww)\leq (1-\lambda)\psi(\uu) + \lambda\psi(\vv)
\end{equation}
where
$\ww = (1-\lambda)\uu + \lambda\vv$.

Let $\alpha,\beta\in\R$ be such that $\psi(\uu)<\alpha$ and
$\psi(\vv)<\beta$ (which exist since $\uu,\vv\in\dom\psi$).
We claim that
$\psi(\ww)\leq \gamma$ where
$\gamma = (1-\lambda)\alpha + \lambda\beta$.

Suppose to the contrary that $\psi(\ww)>\gamma$.
Let
\[
  U
  =
  \Braces{
    \xi\in\fcnspn :\: \xi(\ww)>\gamma,\, \xi(\uu)<\alpha,\, \xi(\vv)<\beta
  }.
\]
Then $U$ is open (being a finite intersection of subbase elements)
and $\psi\in U$, so $U$ is a neighborhood of $\psi$.
Therefore, since $\psi\in\phimgcl$,
there exists a function $\phx\in U\cap\phimg$,
for some $\xx\in\Rn$
(by \Cref{pr:closure:intersect}\ref{pr:closure:intersect:a}),
implying $\phx(\ww)>\gamma$, $\phx(\uu)<\alpha$ and $\phx(\vv)<\beta$.
By definition of $\phx$, it follows that
\[
  \gamma
  <
  \xx\cdot\ww
  =
  (1-\lambda) \xx\cdot\uu + \lambda \xx\cdot\vv
  \leq (1-\lambda) \alpha + \lambda \beta = \gamma,
\]
with the first equality following from $\ww$'s definition.
This is clearly a contradiction.

Thus,
\[ \psi(\ww)\leq (1-\lambda) \alpha + \lambda \beta. \]
Since this holds for all $\alpha>\psi(\uu)$ and $\beta>\psi(\vv)$,
this proves \eqref{eq:thm:h:3:1}, and
therefore that $\psi$ is convex
by
\Cref{pr:stand-cvx-fcn-char}.

Since this holds for all of $\phimgcl$,
this also means that $-\psi$ is convex
by \Cref{pr:h:2}(\ref{pr:h:2:b}).
Therefore, $\psi$ is concave as well.
\end{proof}

We come next to a central theorem showing that any function
$\psi:\Rn\rightarrow\Rext$ that is convex, concave
and that vanishes at the origin must have a particular structural form,
which we express using the $\oms$ notation and leftward sum operation
that were introduced in \Cref{sec:prelim-work-with-infty}.
By \Cref{thm:h:3}, this will apply to every function in
$\phimgcl$, and so also every function in $\phimgA$ (by \Cref{pr:h:7}).
Since $\phimgA$ consists of all functions $\ph{\xbar}$, we obtain that
every function $\ph{\xbar}(\uu)=\xbar\inprod\uu$ must have the form given
in \Cref{thm:h:4}, which is exactly the form we
derived in Example~\ref{ex:poly-speed-intro} for a specific sequence,
though expressed more succinctly in \Cref{thm:h:4}.

For instance,
the cliff ramp function $\psi$ from Examples~\ref{ex:cliff-ramp} and~\ref{ex:cliff-ramp-cont}, which we have already seen is convex, concave and vanishes at the origin,
can be written as
\[
  \psi(\uu)
  =
  \omsf{(\vv\inprod\uu)} \plusl \qq\inprod\uu
  =
  \omsf{(2u_1+u_2)} \plusl (u_2 - u_1).
\]

\begin{theorem}  \label{thm:h:4}
  Let $\psi:\Rn\rightarrow\Rext$ be convex and
  concave with $\psi(\zero)=0$.
  Then, for some $k\ge 0$, there exist orthonormal vectors $\vv_1\dotsc,\vv_k\in\Rn$
  and $\qq\in\Rn$ that is orthogonal to them
  such that for all $\uu\in\Rn$,
  \begin{equation}  \label{thm:h:4:eq1}
  \psi(\uu) = \omsf{(\vv_1\cdot\uu)} \plusl \dotsb \plusl
                          \omsf{(\vv_k\cdot\uu)} \plusl \qq\cdot\uu.
  \end{equation}
\end{theorem}

\begin{proof}
We prove the following \emph{main claim} by induction on $d=0,\dotsc,n$:

\begin{mainclaimp*}
For every linear subspace $L$ of $\Rn$ of dimension $d$,
there exist orthonormal vectors $\vv_1,\dotsc,\vv_k\in L$ and
$\qq\in L$ that is orthogonal to them such that
\eqref{thm:h:4:eq1} holds for all $\uu\in L$.
\end{mainclaimp*}

Once proved,
the theorem will follow by letting $L=\Rn$.

In the base case that $d=0$, we must have $L=\{\zero\}$.
Since $\psi(\zero)=0$, 
we can simply let $k=0$ and $\qq=\zero$ in this case.

For the inductive step, let $d>0$, and assume the main claim holds for
all subspaces of dimension $d-1$.
Let $L\subseteq\Rn$ be any linear subspace of dimension $d$.

Since $\psi$ is both convex and concave, it must have both a convex
epigraph and a convex {hypograph}.
Let us define the sets corresponding to $\psi$'s
epigraph and hypograph on the set $L$:
\begin{align*}
  \Ep &= \{\rpair{\uu}{z} \in L\times\R : z\geq \psi(\uu)\}
  = (L\times\R) \cap (\epi\psi),
\\
  \Em &= \{\rpair{\uu}{z} \in L\times\R : z\leq \psi(\uu)\}
  = (L\times\R) \cap (\hypo\psi).
\end{align*}
Both $\Ep$ and $\Em$ are convex, being the respective intersections of
$\psi$'s epigraph and hypograph with the convex set
$L\times\R$.
They are also
nonempty since they both contain
$\rpair{\zero}{0}$.
The main idea of the proof is to separate these sets with a
hyperplane, which, together with our inductive hypothesis, will allow
us to derive a representation of $\psi$.

To start, we claim that
$\Ep$ and $\Em$ have disjoint relative interiors, that is,
that $(\relint{\Ep})\cap(\relint{\Em})=\emptyset$.
Note that $\Ep$ is itself an epigraph of a convex function, equal to $\psi$ on $L$ and
$+\infty$ outside~$L$, and so, by
\Cref{roc:lem7.3},
any point $\rpair{\uu}{z}\in\relint{\Ep}$
must satisfy $z>\psi(\uu)$, and so cannot be in
$\Em$, nor its relative interior.

Since $\Ep$ and $\Em$ are nonempty convex subsets of
$\Rnp$ with disjoint relative interiors,
by \Cref{roc:thm11.3},
there exists a hyperplane that properly separates them,
meaning that $\Ep$ is included in the closed halfspace $\Hp$ on one side of
the hyperplane, $\Em$ is included in the opposite closed halfspace $\Hm$, and
the sets $\Ep$ and $\Em$ are not both entirely included in the separating
hyperplane itself.
Thus, there exist
$\rpair{\vv}{b}\in\R^{n+1}$ and $c\in\R$,
with $\rpair{\vv}{b}\ne\rpair{\zero}{0}$,
such that
\begin{gather}
\label{eq:sep:1}
  \Ep
  \subseteq\Hp = \set{ \rpair{\uu}{z}\in \Rnp:\:  \vv\cdot\uu + bz\le c },
\\
\label{eq:sep:2}
  \Em
  \subseteq\Hm = \set{ \rpair{\uu}{z}\in \Rnp:\: \vv\cdot\uu + bz\ge c },
\\
\label{eq:sep:3}
\Ep\cup\Em\not\subseteq\Hp\cap\Hm.\vphantom{\Rnp}
\end{gather}

We assume, without loss of generality, that $\vv\in L$. Otherwise, we can replace $\vv$ by $\vv^\prime=\PP\vv$ where $\PP$ is the projection matrix onto the subspace $L$. To see this, let $\Hp'$ and $\Hm'$ be defined like $\Hp$ and $\Hm$, but with $\vv$ replaced by $\vv'$. For all $\uu\in L$, we have
\[
  \vv'\!\inprod\uu=\trans{\vv}\PP\uu=\trans{\vv}\uu=\vv\inprod\uu
\]
(using \Cref{pr:proj-mat-props}\ref{pr:proj-mat-props:a}\ref{pr:proj-mat-props:d}).
Therefore, if $\rpair{\uu}{z}\in\Ep$, then $\vv'\!\inprod\uu+bz=\vv\inprod\uu+bz\le c$,
so $\Ep\subseteq\Hp'$. Symmetrically, $\Em\subseteq\Hm'$. Also, since $\Ep$ and $\Em$ are separated properly, there must exist $\rpair{\uu}{z}\in\Ep\cup\Em$ such that $\rpair{\uu}{z}\not\in\Hp\cap\Hm$, that is, $\vv\inprod\uu+bz\ne c$. For this same
$\rpair{\uu}{z}$ we also have $\vv'\!\inprod\uu+bz=\vv\inprod\uu+bz\ne c$ and thus
$\rpair{\uu}{z}\not\in\Hp'\cap\Hm'$.
Hence, Eqs.~(\ref{eq:sep:1}--\ref{eq:sep:3}) continue to hold when $\vv$ is replaced by $\vv'$.

Also, we must have $c=0$ in Eqs.~(\ref{eq:sep:1}) and~(\ref{eq:sep:2})
since
$\rpair{\zero}{0}\in\Ep\cap\Em\subseteq \Hp\cap\Hm$.
Furthermore,
$b\le c=0$ since
$\rpair{\zero}{1}\in\Ep\subseteq\Hp$.

The remainder of the proof considers separately the cases when $b<0$, corresponding to a nonvertical separating hyperplane, and $b=0$, corresponding to a vertical
separating hyperplane.

\begin{proof-parts}
\pfpart{Case $b<0$:} In this case, the inequalities defining the halfspaces $\Hp$ and $\Hm$ can be rearranged as
\begin{align*}
  \Hp &= \set{ \rpair{\uu}{z}\in \Rnp:\: z \ge -(\vv/b)\cdot\uu}
  = \epi\xi,
\\
  \Hm &= \set{ \rpair{\uu}{z}\in \Rnp:\: z \le -(\vv/b)\cdot\uu}
  = \hypo\xi,
\end{align*}
where $\xi:\Rn\rightarrow\R$
is the linear function $\xi(\uu)=-(\vv/b)\cdot\uu$
for $\uu\in\Rn$.
We next claim that $\psi(\uu)=\xi(\uu)$ for all $\uu\in L$.

Let $\uu\in L$.
We first show $\psi(\uu)\geq\xi(\uu)$.
This is immediate if $\psi(\uu)=+\infty$.
Otherwise, for all $z\in\R$, if $z\geq \psi(\uu)$ then
$\rpair{\uu}{z}\in\Ep\subseteq\Hp$, so
$z\geq \xi(\uu)$.
Since $z\geq\xi(\uu)$ for all $z\geq \psi(\uu)$, it follows that
$\psi(\uu)\geq \xi(\uu)$.

By a symmetric argument, using $\Em\subseteq\Hm$,
we obtain that $\psi(\uu)\leq \xi(\uu)$.

Thus, $\psi(\uu)=\xi(\uu)$ for all $\uu\in L$. Since $\xi(\uu)=\qq\inprod\uu$
for $\uu\in\Rn$,
where $\qq=-(\vv/b)\in L$,
this proves the main claim with $k=0$.

\pfpart{Case $b=0$:}
In this case, we have $\vv\ne\zero$
(since $\rpair{\vv}{b}\ne\rpair{\zero}{0}$),
and
the halfspaces $\Hp$ and $\Hm$ can be written as
\begin{align*}
  \Hp&= \set{ \rpair{\uu}{z}\in \Rnp:\: \vv\cdot\uu\le 0 },
\\
  \Hm&= \set{ \rpair{\uu}{z}\in \Rnp:\: \vv\cdot\uu\ge 0 }.
\end{align*}
 We assume, without loss of generality, that $\norm{\vv}=1$;
otherwise we can
replace $\vv$ by $\vv/\norm{\vv}$.
We examine the function $\psi(\uu)$ for $\uu\in L$, separately
for when $\vv\inprod\uu>0$, $\vv\inprod\uu<0$, and
$\vv\inprod\uu=0$.

First, suppose $\uu\in L$ with $\vv\inprod\uu>0$. 
Then we must have $\psi(\uu)=+\infty$, because if $\psi(\uu)<+\infty$,
then there exists $z\in\R$ with $z\geq\psi(\uu)$, implying
$\rpair{\uu}{z}\in\Ep\subseteq\Hp$, thus
contradicting $\vv\inprod\uu>0$.

By a symmetric argument,
if $\uu\in L$ with $\vv\inprod\uu<0$ then $\psi(\uu)=-\infty$.

It remains to consider the case $\uu\in L$ with $\vv\inprod\uu=0$, that is, $\uu\in L\cap M$, where $M=\set{\uu\in \Rn :\: \vv\cdot\uu=0}$. The separating hyperplane can be expressed in terms of~$M$ as
\[
  \Hp\cap\Hm=\{\rpair{\uu}{z}\in \Rnp :\: \vv\cdot\uu=0\} = M\times\R.
\]
Because $\Ep$ and $\Em$ are \emph{properly} separated, they are not
both entirely contained in this hyperplane;
that is,
\[
  L\times\R = \Ep\cup\Em \not\subseteq M\times\R,
\]
so
$L\not\subseteq M$.

We next apply our inductive hypothesis to the linear subspace
$L' = L\cap M $.
Since $L\not\subseteq M$,
this space has dimension $d-1$.
Thus, there exist orthonormal vectors $\vv_1\dotsc,\vv_k\in L'$ and
$\qq\in L'$ that is orthogonal to them 
such that $\psi(\uu)=\xi'(\uu)$ for all $\uu\in L'$, where
$\xi':\Rn\to\eR$ is the function
\[ \xi'(\uu) =
             \omsf{(\vv_1\cdot\uu)} \plusl \dotsb \plusl
                          \omsf{(\vv_k\cdot\uu)} \plusl \qq\cdot\uu
\]
for $\uu\in\Rn$.

Let us define $\xi:\Rn\rightarrow\Rext$ by
\begin{align*}
  \xi(\uu)
  &= \omsf{(\vv\cdot\uu)} \plusl \xi'(\uu)
\\
  &= \omsf{(\vv\cdot\uu)} \plusl
               \omsf{(\vv_1\cdot\uu)} \plusl \dotsb \plusl
                          \omsf{(\vv_k\cdot\uu)} \plusl \qq\cdot\uu
\end{align*}
for $\uu\in\Rn$.
Then for $\uu\in L$,
\[
  \xi(\uu) =
\begin{cases}
                    +\infty & \text{if $\vv\cdot\uu>0$,}\\
                    \xi'(\uu) & \text{if $\vv\cdot\uu=0$,}\\
                    -\infty & \text{if $\vv\cdot\uu<0$.}
\end{cases}
\]
We claim $\psi(\uu)=\xi(\uu)$ for all $\uu\in L$.
As we argued above, if $\vv\cdot\uu>0$, then $\psi(\uu)=+\infty=\xi(\uu)$,
and if $\vv\cdot\uu<0$, then $\psi(\uu)=-\infty=\xi(\uu)$. Finally,
if $\vv\cdot\uu=0$, then $\uu\in L'$ so
$\psi(\uu)=\xi'(\uu)=\xi(\uu)$ by
inductive hypothesis.
  Also by inductive hypothesis, vectors $\vv_1,\dotsc,\vv_k$ are
  orthonormal and in $L'\subseteq L$, and $\qq$ is orthogonal to them and in $L'\subseteq L$. The vector $\vv$ is of unit length, in~$L$, and orthogonal to $L'$, because $L'\subseteq M$.
Thus,
$\vv,\vv_1,\ldots,\vv_k$ and $\qq$ satisfy all the requirements of the
main claim, completing the induction and the proof.
\qedhere
\end{proof-parts}
\end{proof}

The form of the function $\psi$ in \Cref{thm:h:4} is the same as the one
that arose in the analysis of the polynomially graded sequence
in Example~\ref{ex:poly-speed-intro}:

\begin{proposition}
\label{lemma:chain:psi}
Let $\vv_1,\dotsc,\vv_k,\qq\in\Rn$, and for $\uu\in\Rn$, let
\begin{align*}
  \psi(\uu) &= \omsf{(\vv_1\cdot\uu)} \plusl \dotsb \plusl
                          \omsf{(\vv_k\cdot\uu)}\plusl \qq\cdot\uu,
\\
  \xi(\uu) &=
\newcommand{\ttinprod}{\inprod}
\begin{cases}
+\infty
&\text{if $\vv_1\ttinprod\uu=\dotsb=\vv_{i-1}\ttinprod\uu=0$ and $\vv_i\ttinprod\uu>0$ for some $i$,}
\\
-\infty
&\text{if $\vv_1\ttinprod\uu=\dotsb=\vv_{i-1}\ttinprod\uu=0$ and $\vv_i\ttinprod\uu<0$ for some $i$,}
\\
\qq\ttinprod\uu
&\text{if $\vv_1\ttinprod\uu=\dotsb=\vv_k\ttinprod\uu=0$.}
\end{cases}
\end{align*}
Then $\psi(\uu)=\xi(\uu)$ for all $\uu\in\Rn$.
\end{proposition}

\begin{proof}
The proof follows by verifying that $\xi(\uu)=\psi(\uu)$ in
each of the cases in the definition of $\xi(\uu)$.
\end{proof}

Pulling these results together, we can now conclude that
the sets $\phimgA$ and $\phimgcl$ are actually identical, each
consisting exactly of all of the functions on $\Rn$ that are convex,
concave and that vanish at the origin, and furthermore, that these are
exactly the functions of the form given in \Cref{thm:h:4}.

\begin{theorem}  \label{thm:h:5}
  Let $\psi\in\fcnspn$, that is, $\psi:\Rn\rightarrow\Rext$.
  Then the following are equivalent:
  \begin{letter-compact}
  \item  \label{thm:h:5a0}
    $\psi\in\phimgA$.
    That is, for some $\xbar\in\extspace$,
    $\psi(\uu)=\ph{\xbar}(\uu)=\xbar\cdot\uu$ for all $\uu\in\Rn$.
  \item  \label{thm:h:5a}
    $\psi\in\phimgcl$.
  \item  \label{thm:h:5b}
    $\psi$ is convex, concave and $\psi(\zero)=0$.
  \item  \label{thm:h:5c}
    There exist $\qq,\vv_1\dotsc,\vv_k\in\Rn$, for some $k\geq 0$,
    such that for all $\uu\in\Rn$,
    \[
    \psi(\uu) = \omsf{(\vv_1\cdot\uu)} \plusl \dotsb \plusl
                          \omsf{(\vv_k\cdot\uu)}\plusl \qq\cdot\uu.
    \]
\end{letter-compact}
Furthermore, the same equivalence holds if in
part~(\ref{thm:h:5c}) we additionally require that the vectors
$\vv_1,\dotsc,\vv_k$ are orthonormal and that
$\qq$ is orthogonal to them.
\end{theorem}

\begin{proof}
  ~

\begin{proof-parts}
\pfpart{%
  (\ref{thm:h:5a0})
  $\Rightarrow$
  (\ref{thm:h:5a}):
}
This is immediate from \Cref{pr:h:7}.

\pfpart{%
  (\ref{thm:h:5a})
  $\Rightarrow$
  (\ref{thm:h:5b}):
}
This follows
from
Theorems~\ref{pr:h:2}(\ref{pr:h:2:a}) and~\ref{thm:h:3}.

\pfpart{%
  (\ref{thm:h:5b})
  $\Rightarrow$
  (\ref{thm:h:5c}):
}
This was shown in \Cref{thm:h:4}, including with the more stringent
requirement on the form of $\qq$ and the $\vv_i$'s.

\pfpart{%
  (\ref{thm:h:5c})
  $\Rightarrow$
  (\ref{thm:h:5a0}):
}
Suppose $\psi$ has the form given in part~(\ref{thm:h:5c}).
As in Example~\ref{ex:poly-speed-intro}, let
$\seq{\xx_t}$ be the polynomially graded sequence
given in \eqref{eq:h:6}.
Then for $\uu\in\Rn$,
\[
  \xx_t\cdot\uu
  =
  t^k (\vv_1\cdot\uu) + t^{k-1} (\vv_2\cdot\uu)
  + \dotsb + t (\vv_k\cdot\uu) + (\qq\cdot\uu).
\]
This expression has a limit, namely, the
expression for $\xi(\uu)$ in \Cref{lemma:chain:psi},
which, according to that proposition, is equal to $\psi(\uu)$.
Therefore, the sequence $\seq{\xx_t}$ converges in all directions and
so is in $\pi(\xbar)$ for some $\xbar\in\extspace$
(where $\pi$ is as in \Cref{def:astral-space}).
For this $\xbar$, we have $\xbar\cdot\uu=\lim(\xx_t\cdot\uu)=\xi(\uu)=\psi(\uu)$
for all $\uu\in\Rn$.
That is, $\ph{\xbar}=\psi$, proving $\psi\in\phimgA$.%
\indexg{functional representation of astral points!characterizations|)}
\qedhere
\end{proof-parts}
\end{proof}

\idxwaggoner\citet{waggoner25}
defines
\indexg{linear extended functions}%
\emph{linear extended functions}
as functions $\psi:\Rn\rightarrow\Rext$
that satisfy $\psi(\lambda\uu)=\lambda \psi(\uu)$ for all $\uu\in\Rn$, $\lambda\in\R$,
and also
$\psi(\uu+\vv)=\psi(\uu)+\psi(\vv)$
for all $\uu,\vv\in\Rn$ for which
$\psi(\uu)$ and $\psi(\vv)$ are summable.
\idxwaggoner\citet[Proposition~2.7]{waggoner25}
shows essentially
that a function $\psi:\Rn\rightarrow\Rext$ is linear extended if
and only if it has the form given in part~(\ref{thm:h:5c}) of
\Cref{thm:h:5} (though stated somewhat differently).
As a result,
the set $\phimgA$ of all functions $\ph{\xbar}$, for $\xbar\in\extspace$,
comprises exactly all linear extended functions.
Indeed, Propositions~\ref{pr:i:1} and~\ref{pr:i:2}
show that every function $\ph{\xbar}$ has the two properties defining
such functions.

\section{Astral space, astrons, leftward addition}
\label{sec:astral:space:summary}

\indexg{astral topology!properties|(}%
As already discussed, we can translate topological facts about
$\phimg$ and $\phimgA$ back to $\Rn$ and $\extspace$ via the
homeomorphism $\phfcn$.
We summarize these implications in the following theorem.
Most importantly, we have shown that $\extspace$ is a compactification
of $\Rn$:

\begin{theorem}  \label{thm:i:1}
~
  \begin{letter-compact}
  \item  \label{thm:i:1a}
\indexg{compactness!astral space@of astral space|(}%
\indexg{astral space!compactness of|(}%
    $\extspace$ is compact
    and Hausdorff.
  \item  \label{thm:i:1aa}
    $\extspace$ is normal and regular.
  \item  \label{thm:i:1b}
    $\Rn$ is dense in $\extspace$.
  \item  \label{thm:i:1-compactification}
\indexg{astral space!compactification@as compactification|(}%
\indexg{compactification!astral space as|(}%
    $\extspace$ is a compactification of $\Rn$.
  \end{letter-compact}
Furthermore, let $\seq{\xbar_t}$ be a sequence in $\extspace$ and let $\xbar\in\extspace$. Then:
  \begin{letter-compact}[resume]
  \item  \label{thm:i:1c}
\indexg{coupling function!continuity of|(}%
\indexg{linear functions, continuous extension of!astral space@in astral space|(}%
\indexg{sequence convergence!astral points@of astral points|(}%
    $\xbar_t\rightarrow\xbar$
    if and only if
    for all $\uu\in\Rn$,
    $\xbar_t\cdot\uu\rightarrow\xbar\cdot\uu$.
  \item  \label{thm:i:1e}
    $\seq{\xbar_t}$ converges in $\extspace$
    if and only if
    for all $\uu\in\Rn$,
    $\seq{\xbar_t\cdot\uu}$ converges in $\eR$.
  \item  \label{thm:i:1d}
    There exists a sequence $\seq{\xx_t}$ in $\Rn$ with
    $\xx_t\rightarrow\xbar$.
  \end{letter-compact}
\end{theorem}

\begin{proof}~
\begin{proof-parts}

\pfpart{Part~(\ref{thm:i:1a}):}
As earlier noted, the product space $\fcnspn$ is compact by
Tychonoff's theorem
(\Cref{pr:prod-top-props}\ref{pr:prod-top-props:b}), since
$\Rext$ is compact
(Example~\ref{ex:rext-compact}).
By \Crefequiv{thm:h:5}{thm:h:5a0}{thm:h:5a},
$\phimgA$ equals $\phimgcl$ and is therefore compact,
being a closed subset of the compact space $\fcnspn$
(\Cref{prop:compact}\ref{prop:compact:closed-subset}).
Since $\Rext$ is Hausdorff, $\fcnspn$ and its
subspace $\phimgA$ are also Hausdorff by
Propositions~\ref{pr:prod-top-props}(\ref{pr:prod-top-props:a})
and~\ref{prop:subspace}(\ref{i:subspace:Hausdorff}).
Therefore, $\extspace$, which is homeomorphic with $\phimgA$, is
also compact and Hausdorff.%
\indexg{compactness!astral space@of astral space|)}%
\indexg{astral space!compactness of|)}%

\pfpart{Part~(\ref{thm:i:1aa}):}
Since $\extspace$ is compact and Hausdorff,
by \Cref{prop:compact}(\ref{prop:compact:subset-of-Hausdorff})
it is also normal and regular.

\pfpart{Part~(\ref{thm:i:1b}):}
That $\phimgA=\phimgcl$ means exactly that $\phimg$ is dense in
$\phimgA$.
So, by our homeomorphism, $\Rn$ is dense in $\extspace$.

\pfpart{Part~(\ref{thm:i:1-compactification}):}
This follows from
parts~(\ref{thm:i:1a}) and~(\ref{thm:i:1b}),
and from
\Crefequiv{pr:open-sets-equiv-alt}{pr:open-sets-equiv-alt:a}{pr:open-sets-equiv-alt:c}
which shows that $\Rn$, in the Euclidean topology, is a subspace of
$\extspace$.%
\indexg{astral space!compactification@as compactification|)}%
\indexg{compactification!astral space as|)}%

\pfpart{Part~(\ref{thm:i:1c}):}
By our homeomorphism,
${\xbar_t}\rightarrow{\xbar}$ if and only if
$\ph{\xbar_t}\rightarrow\ph{\xbar}$,
which,
because we are using product topology in $\fcnspn$,
in turn holds if and only if
$\xbar_t\cdot\uu=\ph{\xbar_t}(\uu)\rightarrow\ph{\xbar}(\uu)=\xbar\cdot\uu$
for all $\uu\in\Rn$
(by \Cref{pr:prod-top-ptwise-conv}).%
\indexg{sequence convergence!astral points@of astral points|)}%
\indexg{coupling function!continuity of|)}%
\indexg{linear functions, continuous extension of!astral space@in astral space|)}%

\pfpart{Part~(\ref{thm:i:1e}):}
If ${\xbar_t}\rightarrow{\xbar}$ then part~(\ref{thm:i:1c}) implies
that $\seq{\xbar_t\cdot\uu}$ converges in $\eR$. For the converse,
let $\psi:\Rn\rightarrow\Rext$ be defined by
$\psi(\uu)=\lim (\xbar_t\cdot\uu)$
for $\uu\in\Rn$.
Then $\ph{\xbar_t}(\uu)\rightarrow\psi(\uu)$ for all $\uu\in\Rn$,
so $\ph{\xbar_t}\rightarrow\psi$.
Since $\ph{\xbar_t}\in\phimgA$, for all $t$, and since $\phimgA$ is closed,
$\psi$ must also be in $\phimgA$
(\Cref{prop:first:properties}\ref{prop:first:closure}).
Thus, $\psi=\ph{\xbar}$ for some $\xbar\in\extspace$, so
$\ph{\xbar_t}\rightarrow\ph{\xbar}$, implying,
by our homeomorphism, that
${\xbar_t}\rightarrow{\xbar}$.

\pfpart{Part~(\ref{thm:i:1d}):}
Let $\seq{\xx_t}$ in $\Rn$ be any sequence in $\pi(\xbar)$
(where $\pi$ is as given in \Cref{def:astral-space}).
Then for all $\uu\in\Rn$, $\xx_t\cdot\uu\rightarrow\xbar\cdot\uu$
by definition of the coupling function (Eq.~\ref{eqn:coupling-defn}),
so $\xx_t\rightarrow\xbar$ by part~(\ref{thm:i:1c}).%
\indexg{astral topology!properties|)}
\qedhere
\end{proof-parts}
\end{proof}

Returning to the construction of astral space given in
\Cref{def:astral-space}, for every $\xbar\in\extspace$,
\Cref{thm:i:1}(\ref{thm:i:1c}) shows
that a sequence $\seq{\xx_t}$ in $\Rn$
converges to $\xbar$ if and only if
$\xx_t\cdot\uu\rightarrow\xbar\cdot\uu$ for all $\uu\in\Rn$,
and so if and only if
$\seq{\xx_t}$ is in $\pi(\xbar)$.
Thus, $\pi(\xbar)$ consists of exactly those sequences in $\Rn$ that
converge to $\xbar$.

\indexg{astrons|(}%
In \Cref{sec:intro:astral}
and Example~\ref{ex:basic-lim-ray},
we encountered
an important example of an astral point not in $\Rn$:
For $\vv\in\Rn$, the point $\limray{\vv}$, which is the product of the
scalar $\oms=+\infty$ with the vector $\vv$, is
called an \emph{astron}, and is defined to be that point in
$\extspace$ that is the limit of the sequence $\seq{t\vv}$.
\indexg{scalar multiples (astral)|(}%
\indexg{omega-multiples@$\oms$-multiples|(}%
More generally, for $\xbar\in\extspace$, we define the product
$\limray{\xbar}$ analogously as the limit of the sequence
$\seq{t\xbar}$.
(Note that a point $\limray{\xbar}$ is \emph{not} called an astron
unless it is equal to $\limray{\vv}$ for some $\vv\in\Rn$.)
The next simple proposition proves the existence and general form of
such points, including all astrons.

As we will see shortly, astrons act as fundamental building blocks
which, together with real vectors, generate all of astral space using
a generalized form of leftward addition.

\begin{proposition}  \label{pr:astrons-exist}
  Let $\xbar\in\extspace$.
  Then the sequence $\seq{t\xbar}$ has a limit in $\extspace$,
  henceforth denoted
\indexm{alpha xbar300}{$\alpha\xbar$}{scalar multiple}%
$\limray{\xbar}$.
  Furthermore, for all $\uu\in\Rn$,
  \[
     (\limray{\xbar})\cdot\uu
     = \omsf{(\xbar\cdot\uu)}
           =
     \begin{cases}
                 +\infty & \text{if $\xbar\cdot\uu>0$,} \\
                 0       & \text{if $\xbar\cdot\uu=0$,} \\
                 -\infty & \text{if $\xbar\cdot\uu<0$.}
     \end{cases}
  \]
\end{proposition}

\begin{proof}
For $\uu\in\Rn$,
\[
  (t\xbar)\cdot\uu
  =
  t(\xbar\cdot\uu)
  \rightarrow
  \omsf{(\xbar\cdot\uu)}
\]
(with the equality following from \Cref{pr:i:2}).
Since this limit exists for each $\uu\in\Rn$,
by
\Cref{thm:i:1}(\ref{thm:i:1e}),
the sequence
$\seq{t\xbar}$ must have a limit in $\extspace$,
namely, $\limray{\xbar}$.
This also proves $(\limray{\xbar})\cdot\uu$ has the form given
in the proposition.
\end{proof}

When multiplying a scalar $\lambda\in\R$ by a real vector
$\xx=\trans{[x_1,\ldots,x_n]}$, the product $\lambda\xx$ can of course
be computed simply by multiplying each component $x_i$ by $\lambda$
to obtain $\trans{[\lambda x_1,\ldots,\lambda x_n]}$.
The same is \emph{not} true when multiplying by $\oms$.
Indeed, multiplying each component $x_i$ by $\oms$ would yield
a vector in $(\Rext)^n$, not an astral point in $\extspace$.
For instance, if we multiply each component of
$\xx=\trans{[2,0,-1]}$ by $\oms$, we obtain the point
$\trans{[+\infty,0,-\infty]}$ in $(\Rext)^3$, which is very different
from the astron $\limray{\xx}$ in $\extspac{3}$
(see the discussion of the distinction between $(\eR)^n$ and $\eRn$ after \Cref{pr:i:3}).%
\indexg{omega-multiples@$\oms$-multiples|)}

In general, in $\extspace$,
each of the unit vectors $\vv\in\Rn$, with $\norm{\vv}=1$,
gives rise to a distinct astron $\limray{\vv}$.
Together with the origin,
$\limray{\zero}=\zero$, these comprise all of the
astrons in $\extspace$.
For example, $\extspac{1}=\Rext$ has three astrons:
$\limray{(1)}=\oms=+\infty$,
$\limray{(-1)}=-\oms=-\infty$,
and
$\limray{(0)}=0$.%
\indexg{astrons|)}

Although used less often, for completeness, we also define scalar
multiplication by $-\oms$ as
$(-\oms)\xbar=\limray{(-\xbar)}$ for $\xbar\in\extspace$.

The next proposition summarizes properties of scalar multiplication:

\begin{proposition}  \label{pr:scalar-prod-props}
  Let $\xbar\in\extspace$, and let $\alpha,\beta\in\Rext$.
  Then the following hold:
  \begin{letter-compact}
  \item  \label{pr:scalar-prod-props:a}
    $(\alpha\xbar)\cdot\uu = \alpha(\xbar\cdot\uu)$
    for all $\uu\in\Rn$.
  \item  \label{pr:scalar-prod-props:b}
    $\alpha(\beta\xbar)=(\alpha \beta) \xbar$.
  \item  \label{pr:scalar-prod-props:c}
    $0 \xbar = \zero$.
  \item  \label{pr:scalar-prod-props:d}
    $1 \xbar = \xbar$.
  \item  \label{pr:scalar-prod-props:e}
    $\lambda_t \xbar_t \rightarrow \lambda \xbar$ for any sequence
    $\seq{\xbar_t}$ in $\extspace$ and any sequence
    $\seq{\lambda_t}$ in $\R$ such that
    $\xbar_t\to\xbar$
    and
    $\lambda_t\to\lambda$ for some $\lambda\in\R\wo\set{0}$.
  \end{letter-compact}
\end{proposition}

\begin{proof}
~

\begin{proof-parts}
\pfpart{Part~(\ref{pr:scalar-prod-props:a}):}
If $\alpha\in\R$, then the claim follows from
\Cref{pr:i:2}.
If $\alpha=\oms$, then it follows from
\Cref{pr:astrons-exist}.
If $\alpha=-\oms$, then for $\uu\in\Rn$,
\[
  [(-\oms)\xbar]\cdot\uu
  =
  [\limray{(-\xbar)}]\cdot\uu
  =
  \limray{(-\xbar\cdot\uu)}
  =
  (-\oms)(\xbar\cdot\uu).
\]
The first equality is by definition of multiplication by $-\oms$.
The second is by
\Cref{pr:astrons-exist}.
And the last is by commutativity of multiplication over $\Rext$.

\pfpart{Part~(\ref{pr:scalar-prod-props:b}):}
For all $\uu\in\Rn$,
\[
  [\alpha(\beta\xbar)]\cdot\uu
  =
  \alpha [(\beta\xbar)\cdot\uu]
  =
  \alpha [\beta(\xbar\cdot\uu)]
  =
  (\alpha\beta)(\xbar\cdot\uu)
  =
  [(\alpha\beta)\xbar]\cdot\uu,
\]
where the first, second and last equalities are all 
by part~(\ref{pr:scalar-prod-props:a}).
The claim then follows by \Cref{pr:i:4}.

The proofs of the remaining parts are similar.%
\indexg{scalar multiples (astral)|)}
\qedhere
\end{proof-parts}
\end{proof}

\indexg{leftward addition!astral points@of astral points|(}%
Next, we extend
leftward addition, previously defined
in
Eqs.~(\ref{eqn:intro-left-sum-defn})
and~(\ref{eqn:left-sum-alt-defn})
for scalars in
$\Rext$, to all points in $\extspace$.
For $\xbar,\ybar\in\extspace$, we define $\xbar\plusl\ybar$ to be that
unique point in $\extspace$ for which
\[ (\xbar\plusl\ybar)\cdot\uu = \xbar\cdot\uu \plusl \ybar\cdot\uu \]
for all $\uu\in\Rn$.
The next proposition shows that such a point must exist.
This operation will turn out to be useful for describing and
working with points in $\extspace$.
It is much like taking the vector sum of two points in $\Rn$, except
that for those vectors $\uu$ for which this would lead to adding
$+\infty$ and $-\infty$, the left summand dominates.

\begin{proposition}  \label{pr:i:6}
  Let $\xbar$ and $\ybar$ be in $\extspace$.
  Then there exists a unique point in $\extspace$, henceforth denoted
\indexm{x+y600}{$\xbar\protect\plusl\ybar$}{leftward sum (of astral points)}%
  $\xbar\plusl\ybar$,
  for which
  \[ (\xbar\plusl\ybar)\cdot\uu = \xbar\cdot\uu \plusl \ybar\cdot\uu \]
  for all $\uu\in\Rn$.
\end{proposition}

\begin{proof}
By \Crefequiv{thm:h:5}{thm:h:5a0}{thm:h:5c},
since $\ph{\xbar}$ and $\ph{\ybar}$ are both
in $\phimgA$,
there exist $\qq,\vv_1,\dotsc,\vv_k\in\Rn$, for some $k\geq 0$,
and $\qq',\vv'_1,\dotsc,\vv'_{k'}\in\Rn$, for some $k'\geq 0$,
such that for all $\uu\in\Rn$,
\begin{align}
\SwapAboveDisplaySkip
  \ph{\xbar}(\uu) &=
  \omsf{(\vv_1\cdot\uu)} \plusl \dotsb \plusl
                          \omsf{(\vv_k\cdot\uu)} \plusl \qq\cdot\uu,
\notag
  \\
  \ph{\ybar}(\uu) &=
  \omsf{(\vv'_1\cdot\uu)} \plusl \dotsb \plusl
                          \omsf{(\vv'_{k'}\cdot\uu)} \plusl \qq'\cdot\uu.
  \label{eq:i:3}
\end{align}
Let $\psi:\Rn\rightarrow\Rext$ be defined by
\begin{align*}
  \psi(\uu)
  &=
  \xbar\cdot\uu \plusl \ybar\cdot\uu
  \\
  &=
  \ph{\xbar}(\uu)\plusl \ph{\ybar}(\uu)
  \\
  &=
     \omsf{(\vv_1\cdot\uu)} \plusl \dotsb \plusl \omsf{(\vv_k\cdot\uu)}
     \plusl
     \omsf{(\vv'_1\cdot\uu)} \plusl \dotsb \plusl \omsf{(\vv'_{k'}\cdot\uu)}
     \plusl (\qq+\qq')\cdot\uu,
\end{align*}
for $\uu\in\Rn$.
where the last equality uses \eqref{eq:i:3} and
\Cref{pr:i:5}(\ref{pr:i:5c}).
Because $\psi$ has this form, \Crefequiv{thm:h:5}{thm:h:5c}{thm:h:5a0}
implies
it is in $\phimgA$ and therefore
equal to $\ph{\zbar}$ for some $\zbar\in\extspace$.
Setting $\xbar\plusl\ybar$ equal to $\zbar$ proves the proposition,
with uniqueness following from \Cref{pr:i:4}.
\end{proof}

Just as was the case for scalars, leftward addition of astral
points is associative, distributive and partially commutative
(specifically, if either of the summands is in $\Rn$).
It is also the same as vector addition when both summands are in
$\Rn$.
We summarize these and other properties in the next proposition.

\begin{proposition}  \label{pr:i:7}
  Let $\xbar,\ybar,\zbar\in\extspace$,
  and $\xx,\yy\in\Rn$.
  Then:
  \begin{letter-compact}
  \item  \label{pr:i:7a}
    $(\xbar\plusl \ybar)\plusl \zbar=\xbar\plusl (\ybar\plusl \zbar)$.
  \item  \label{pr:i:7b}
    $\lambda(\xbar\plusl \ybar)=\lambda\xbar\plusl \lambda\ybar$,
    for $\lambda\in\R$.
  \item  \label{pr:i:7c}
    For all $\alpha,\beta\in\Rext$, if $\alpha\beta\geq 0$
    then
    $\alpha\xbar\plusl \beta\xbar=(\alpha+\beta)\xbar$.
  \item  \label{pr:i:7d}
    $\xbar\plusl \yy = \yy\plusl \xbar$.
    In particular,
    $\xbar\plusl \zero = \zero\plusl \xbar = \xbar$.
  \item  \label{pr:i:7e}
    $\xx\plusl \yy = \xx + \yy$.
  \item  \label{pr:i:7f}
    $\xbar\plusl \ybart\to\xbar\plusl\ybar$, for any sequence
    $\seq{\ybart}$ in $\eRn$ with $\ybart\to\ybar$.
  \item  \label{pr:i:7g}
    $\xbar_t\plusl \yy_t\to\xbar\plusl\yy$, for any sequences
    $\seq{\xbar_t}$ in $\eRn$ and $\seq{\yy_t}$ in $\Rn$
    with $\xbar_t\to\xbar$ and $\yy_t\rightarrow\yy$.
  \end{letter-compact}
\end{proposition}

\begin{proof}
~
\begin{proof-parts}
\pfpart{Part~(\ref{pr:i:7a}):}
For $\uu\in\Rn$,
\begin{align*}
  \bigParens{(\xbar\plusl \ybar)\plusl \zbar}\cdot\uu
  &=
  (\xbar\plusl \ybar)\cdot\uu \plusl \zbar\cdot\uu
  \\
  &=
  (\xbar\cdot\uu\plusl \ybar\cdot\uu) \plusl \zbar\cdot\uu
  \\
  &=
  \xbar\cdot\uu\plusl (\ybar\cdot\uu \plusl \zbar\cdot\uu)
  \\
  &=
  \xbar\cdot\uu\plusl (\ybar \plusl \zbar) \cdot\uu
  =
  \bigParens{\xbar \plusl (\ybar \plusl \zbar)} \cdot\uu,
\end{align*}
where the third equality is by
\Cref{pr:i:5}(\ref{pr:i:5a}),
and each of the others 
by \Cref{pr:i:6}.
The claim now follows by \Cref{pr:i:4}.
\end{proof-parts}

The proofs of the other parts are similar.%
\indexg{leftward addition!astral points@of astral points|)}
\end{proof}

\indexg{leftward addition!sets@of sets|(}%
The leftward sum of two sets $X$ and $Y$ in
$\extspace$ is defined similarly as for ordinary addition:
\[
  X\plusl Y = \{ \xbar\plusl\ybar:\: \xbar\in X,\, \ybar\in Y \}.
\]
We also define
$X\plusl\ybar = X\plusl \{\ybar\}$
and
$\ybar\plusl X =  \{\ybar\} \plusl X$
for $X\subseteq\extspace$ and
\indexg{leftward addition!sets@of sets|)}%
$\ybar\in\extspace$.
\indexg{scalar multiples (astral)!sets@of sets|(}%
Analogously, if $A\subseteq\Rext$ and $X\subseteq\extspace$,
then we define the product
\[
    A X = \{ \alpha \xbar:\: \alpha\in A,\, \xbar\in X \},
\]
with $\alpha X = \{\alpha\} X$ and
$A \xbar = A \{\xbar\}$
for $\alpha\in\Rext$ and $\xbar\in\extspace$.%
\indexg{scalar multiples (astral)!sets@of sets|)}%

\indexg{decompositions, astral|(}%
\indexg{astral points!decomposition of|(}%
\indexg{representations of astral points!astrons@with astrons|(}%
Using astrons and leftward addition,
and invoking \Cref{thm:h:5},
we obtain
the following succinct representation for all points in $\extspace$:

\begin{corollary}
\label{cor:h:1}
  Astral space $\extspace$ consists of all points $\xbar$
  of the form
  \begin{equation}   \label{eq:cor:h:1:0}
    \xbar=\limray{\vv_1} \plusl \dotsb \plusl \limray{\vv_k} \plusl \qq
  \end{equation}
  for $\qq,\vv_1,\dotsc,\vv_k\in\Rn$, $k\geq 0$.
  Furthermore, the same statement holds if we require
  the vectors $\vv_1,\dotsc,\vv_k$ to be
  orthonormal and $\qq$ to be orthogonal to them.
\end{corollary}

\begin{proof}
Every expression of the form given in the corollary is in $\extspace$
since astral space includes all astrons and all of $\Rn$, and
since the space is closed under leftward addition
(by Propositions~\ref{pr:astrons-exist} and~\ref{pr:i:6}).

To show every point has this form,
suppose $\xbar\in\extspace$.
Then $\ph{\xbar}\in\phimgA$, so
by \Crefequiv{thm:h:5}{thm:h:5a0}{thm:h:5c},
there exist orthonormal $\vv_1,\dotsc,\vv_k\in\Rn$ and $\qq\in\Rn$ that is orthogonal to them,
such that, for $\uu\in\Rn$,
\begin{align*}
  \xbar\cdot\uu
  =
  \ph{\xbar}(\uu)
  &=
  \omsf{(\vv_1\cdot\uu)} \plusl \dotsb \plusl
                          \omsf{(\vv_k\cdot\uu)} \plusl \qq\cdot\uu
  \\
  &=
  (\limray{\vv_1})\cdot\uu \plusl \dotsb \plusl
                          (\limray{\vv_k})\cdot\uu \plusl \qq\cdot\uu
  \\
  &=
  \bigParens{\limray{\vv_1} \plusl \dotsb \plusl \limray{\vv_k} \plusl\qq}\cdot\uu,
\end{align*}
where
the third and fourth equalities are by
Propositions~\ref{pr:astrons-exist}
and~\ref{pr:i:6}, respectively.
By \Cref{pr:i:4}, this proves
\eqref{eq:cor:h:1:0}.
\end{proof}

We refer to the expression $\xbar=\limray{\vv_1} \plusl \dotsb \plusl \limray{\vv_k} \plusl \qq$
as an \emph{(astral) decomposition} of $\xbar$; it is not necessarily unique.
When vectors $\vv_1,\dotsc,\vv_k$ are orthonormal
and $\qq$ is orthogonal to them, then the astral decomposition is said to be
\indexg{decompositions, astral!orthogonal}%
\emph{orthogonal}.
In \Cref{sec:canonical:rep}, we show that each astral point has a unique orthogonal decomposition.

Using \Cref{lemma:chain:psi} we can also characterize the value of $\xbar\inprod\uu$
based on the astral decomposition of $\xbar$:

\begin{lemma}[Case Decomposition Lemma]
\label{lemma:case}
\indexg{Case Decomposition Lemma|(}%
Let $\xbar=\limray{\vv_1} \plusl \dotsb \plusl \limray{\vv_k} \plusl \qq$
\ for some $k\ge 0$ and $\vv_1,\dotsc,\vv_k,\qq\in\Rn$,
and let $\uu\in\Rn$.
Then
  \[
  \xbar\inprod\uu =
  \newcommand{\ttinprod}{\inprod}
  \begin{cases}
  +\infty
  &\text{if $\vv_1\ttinprod\uu=\dotsb=\vv_{i-1}\ttinprod\uu=0$ and $\vv_i\ttinprod\uu>0$ for some $i$,}
  \\
  -\infty
  &\text{if $\vv_1\ttinprod\uu=\dotsb=\vv_{i-1}\ttinprod\uu=0$ and $\vv_i\ttinprod\uu<0$ for some $i$,}
  \\
  \qq\ttinprod\uu
  &\text{if $\vv_1\ttinprod\uu=\dotsb=\vv_k\ttinprod\uu=0$.}
  \end{cases}
  \]
\end{lemma}
\begin{proof}
Let $\psi(\uu)$ and $\xi(\uu)$ be as given in
\Cref{lemma:chain:psi}.
Then from Propositions~\ref{pr:astrons-exist}
and~\ref{pr:i:6}, it follows that $\xbar\cdot\uu=\psi(\uu)$.
Since $\psi(\uu)=\xi(\uu)$
by \Cref{lemma:chain:psi}, this proves the claim.%
\indexg{Case Decomposition Lemma|)}%
\indexg{decompositions, astral|)}%
\indexg{astral points!decomposition of|)}%
\indexg{representations of astral points!astrons@with astrons|)}
\end{proof}

\Cref{thm:h:5} and \Cref{cor:h:1}
tell us that
elements of astral space can be viewed
from two main perspectives.
First, as discussed in \Cref{sec:astral-space-intro},
we can think of astral points $\xbar\in\extspace$
as limits of sequences in~$\Rn$.
Such sequences can have a finite limit $\qq$ in~$\Rn$, exactly when they converge
to $\qq$ in the Euclidean topology.
But more interestingly, \Cref{cor:h:1} expresses precisely the form of
every astral point, and so also
all of the ways in which convergent sequences in astral space can
tend to infinity:
Every such sequence has a primary dominant direction~$\vv_1$, and may also have a secondary direction~$\vv_2$, a tertiary
direction~$\vv_3$, etc.,
and a remaining finite part
$\qq\in\Rn$.
This is exactly the form discussed in
\Cref{subsec:astral-pt-form};
the polynomially graded sequence constructed in Example~\ref{ex:poly-speed-intro} is an example of a
sequence with just these properties.

Alternatively, as we saw in
Sections~\ref{subsec:astral-pts-as-fcns}
and~\ref{sec:char:func},
astral points can be represented
as functions mapping $\Rn$ to
$\Rext$, and indeed their topology and other properties are largely
derived from this correspondence.
Astral space, we proved, is a homeomorphic copy of the function space
$\phimgA=\phimgcl$, a space that
consists
of those
functions on $\Rn$ that are convex, concave and that vanish at the
origin.
We also proved that these functions have a very particular functional
form.

In the foregoing,
we followed an approach in which astral space was constructed
based on the first perspective.
But it would have been possible to instead construct this same space
based on the second perspective.
Such an approach would begin
with $\phimg$, which is exactly the space of all
linear functions, and which is
equivalent topologically and as a vector space to $\Rn$.
The next step would be to form the closure of this space,~$\phimgcl$, which
is compact as a closed subset of the compact space~$\fcnspn$.
The actual astral space~$\extspace$ could then be constructed as a
homeomorphic copy of $\phimgcl$.
This is very much like the approach used, for instance,
in constructing the
\indexg{Stone-Cech compactification@Stone-\v{C}ech compactification}%
Stone-\v{C}ech compactification
\idxmunk\citep[Section~38]{munkres}.

\indexg{astral space!compactification@as compactification|(}%
More significantly, this view
shows that astral space is a special case of the
\indexg{Q-compactification@$\calQ$-compactification}%
$\calQ$-compacti\-fications studied by
\indexa{Loeb, P. A.}%
\citet{q_compactification}, in which we set
$\calQ$ to be the set of all linear functions on $\Rn$.
\indexg{linear functions, continuous extension of!astral space@in astral space|(}%
In particular, this connection implies that astral space is the
smallest compactification of $\Rn$ under which all linear functions can be
extended continuously, which was a primary aim for its construction.%
\indexg{astral space!compactification@as compactification|)}%
\indexg{linear functions, continuous extension of!astral space@in astral space|)}%

While our approach to extending $\Rn$ is motivated by compactness and continuity,
\indexa{Hamel, A. H.}\indexa{Schrage, C.}%
\citet{hamel-schrage}
use order theory to
develop a different approach,
motivated by the goal of resolving ambiguity when adding $+\infty$ and
$-\infty$. Their approach also yields extensions of various notions from convex analysis,
including subgradients, conjugacy, and duality.

\section{A sufficient condition for convergence of sequences}
\label{sec:suf-cond-seq-conv}

\indexg{representations of astral points!sequences converging to|(}%
\indexg{sequence convergence!decomposition@to decomposition|(}%
The next theorem proves
a very general sufficient condition for when a sequence in $\Rn$
converges to a particular point $\xbar\in\extspace$.
This generalizes the convergence of a polynomially graded
sequence explored in
Example~\ref{ex:poly-speed-intro},
here replacing polynomial coefficients with arbitrary
sequences going to infinity at successively decreasing rates.
Later, in
\Cref{sec:seqs-to-matrix-rep},
we will see the sense in which
these same conditions are both necessary and sufficient for such
convergence.

In the theorem,
as usual, $b_{t,i}$ denotes the $i$-th component of vector $\bb_t$.

\begin{theorem}
\label{thm:i:seq-rep}
Let
$\xbar=\limray{\vv_1} \plusl \dotsb \plusl \limray{\vv_k} \plusl \qq$
\ for some $k\ge 0$ and $\vv_1,\dotsc,\vv_k,\qq\in\Rn$.
Let $\seq{\bb_t}$ be a sequence in $\Rk$ and $\seq{\qq_t}$ a
sequence in $\Rn$.
For each $t$, let
\begin{equation*}  %
  \xx_t
  =
  b_{t,1} \vv_1
  + \dotsb + b_{t,k}\vv_k + \qq_t
  =
  \sum_{i=1}^k b_{t,i} \vv_i + \qq_t.
\end{equation*}
Assume all of the following hold:
\begin{letter-compact}
  \item  \label{thm:i:seq-rep:a}
    $b_{t,i}\rightarrow+\infty$, for $i=1,\dotsc,k$.
  \item  \label{thm:i:seq-rep:b}
    $b_{t,i+1}/b_{t,i}\rightarrow 0$, for $i=1,\dotsc,k-1$.
  \item  \label{thm:i:seq-rep:c}
    $\qq_t\rightarrow\qq$.
\end{letter-compact}
Then $\xx_t\to\xbar$.%
\indexg{representations of astral points!sequences converging to|)}%
\indexg{sequence convergence!decomposition@to decomposition|)}
\end{theorem}

Condition~(\ref{thm:i:seq-rep:b}) of the theorem means that the
sequence $\seq{b_{t,i}}$ is growing to $+\infty$ much faster than
$\seq{b_{t,i+1}}$.
\indexg{asymptotic dominance|(}%
Before proving the theorem, we focus for a moment on the relation
between sequences that this is an example of, and which we now define
more generally:

\begin{definition}
Let $\seq{x_t}$ and $\seq{y_t}$ be sequences in $\R$,
both converging to $+\infty$.
We say that $\seq{x_t}$
\emph{asymptotically dominates}
$\seq{y_t}$,
written $x_t \seqgt y_t$%
\indexm{xt>yt}{$x_t \seqgt y_t$}{asymptotically dominates}
(or $y_t \seqlt x_t$),
if $y_t/x_t \rightarrow 0$.
\end{definition}

Thus, in the theorem,
condition~(\ref{thm:i:seq-rep:b})
means that
\begin{equation}  \label{eq:b-seq-asymp-dom}
   b_{t,1}
   \seqgt 
   b_{t,2}
   \seqgt 
   \dotsb
   \seqgt
   b_{t,k}.
\end{equation}
When this condition is satisfied by a sequence $\seq{\bb_t}$ in $\Rk$,
in other words, when both conditions~(\ref{thm:i:seq-rep:a})
and~(\ref{thm:i:seq-rep:b}) of the theorem hold,
we say that
\indexg{convergence of entries at decreasing rates}%
\emph{the entries of $\bb_t$ converge to $+\infty$ at decreasing rates}.
Note that for $k=0$, there exists only one sequence in $\R^0$,
namely, the constant sequence
with every element equal to $\zerovec$.
Vacuously, this sequence satisfies both
conditions~(\ref{thm:i:seq-rep:a}) and~(\ref{thm:i:seq-rep:b}) of
\Cref{thm:i:seq-rep}, so its entries converge to $+\infty$ at
decreasing rates, though in a rather degenerate sense.

Here are some order properties of asymptotic dominance:

\begin{proposition}   \label{pr:asymp-dom-props}
  Let $\seq{x_t}$, $\seq{y_t}$, and $\seq{z_t}$ be 
  sequences in $\R$, all converging to $+\infty$.
  \begin{letter-compact}
  \item   \label{pr:asymp-dom-props:a}
    If $x_t\seqgt y_t$ and $y_t\seqgt z_t$, then
    $x_t\seqgt z_t$.
  \item   \label{pr:asymp-dom-props:b}
    If $x_t\seqgt y_t$, then $x_t/y_t\rightarrow+\infty$
    (and therefore $y_t\not\seqgt x_t$).
  \item   \label{pr:asymp-dom-props:c}
    There exist sequences $\seq{x'_t}$ and $\seq{z'_t}$
    in $\R$, both converging to $+\infty$, such that
    $x'_t\seqgt y_t\seqgt z'_t$.
  \item   \label{pr:asymp-dom-props:d}
    If $x_t\seqgt z_t$, then there exists
    a sequence $\seq{y'_t}$ in $\R$ converging to $+\infty$
    such that
    $x_t\seqgt y'_t\seqgt z_t$.
  \item   \label{pr:asymp-dom-props:e}
    If $x_t\seqgt y_t \seqgt z_t$,
    then
    $x_t/z_t \seqgt y_t/z_t$.
  \end{letter-compact}
\end{proposition}

\begin{proof}
  ~
  
\begin{proof-parts}
\pfpart{Part~(\ref{pr:asymp-dom-props:a}):}
If $x_t\seqgt y_t$ and $y_t\seqgt z_t$,
then
\[
  \frac{z_t}{x_t}
  =
  \frac{z_t}{y_t}\cdot\frac{y_t}{x_t}\rightarrow 0,
\]
where the equality holds for all $t$ sufficiently large
(so that $x_t$ and $y_t$ are positive).
Thus, $x_t\seqgt z_t$.

\pfpart{Part~(\ref{pr:asymp-dom-props:b}):}
This follows from definitions and since 
$x_t$ and $y_t$ are positive for all $t$ sufficiently large.

\pfpart{Part~(\ref{pr:asymp-dom-props:c}):}
For all $t$,
let $x'_t=ty_t$ and $z'_t=\sqrt{|y_t|}$.
Then $\seq{x'_t}$ and $\seq{z'_t}$ both converge to $+\infty$
(since $\seq{y_t}$ does).
Furthermore,
$y_t / x'_t = 1/t\rightarrow 0$,
and
$z'_t / y_t = 1/\sqrt{y_t} \rightarrow 0$
(with the equalities holding for $t$ so large that $y_t>0$).
Thus, $x'_t\seqgt y_t\seqgt z'_t$.

\pfpart{Part~(\ref{pr:asymp-dom-props:d}):}
For all $t$,
let $y'_t=\sqrt{|x_t z_t|}$.
Then $y'_t\rightarrow+\infty$.
Also,
\[
   \frac{y'_t}{x_t}
   =
   \frac{z_t}{y'_t}
   =
   \sqrt{\frac{z_t}{x_t}}
   \rightarrow 0,
\]
with the equalities holding for all $t$ sufficiently large
(so that $x_t$ and $z_t$ are positive),
and with convergence from $x_t\seqgt z_t$.
Thus, $x_t\seqgt y'_t\seqgt z_t$.

\pfpart{Part~(\ref{pr:asymp-dom-props:e}):}
Suppose $x_t\seqgt y_t \seqgt z_t$.
Then $x_t/z_t\rightarrow+\infty$
and $y_t/z_t\rightarrow+\infty$
by parts~(\ref{pr:asymp-dom-props:a})
and~(\ref{pr:asymp-dom-props:b}),
and
$(y_t/z_t)/(x_t/z_t)=y_t/x_t\rightarrow+\infty$
(with equality holding for $t$ sufficiently large, so that $x_t,z_t>0$).
Thus, $x_t/z_t \seqgt y_t/z_t$.%
\indexg{asymptotic dominance|)}
\qedhere
\end{proof-parts}
\end{proof}

The next simple observation follows immediately:

\begin{proposition}
\label{prop:dec:trans}
\indexg{convergence of entries at decreasing rates|(}%
Let $\seq{\bb_t}$ be a sequence in $\Rk$ whose entries converge to $+\infty$ at decreasing rates.
Suppose $i,j\in\{1,\ldots,k\}$, with $i<j$.
Then $b_{t,j}/b_{t,i}\to 0$ and $b_{t,i}/b_{t,j}\to+\infty$.
\end{proposition}

\begin{proof}
By \eqref{eq:b-seq-asymp-dom} and
since asymptotic dominance is transitive
(\Cref{pr:asymp-dom-props}\ref{pr:asymp-dom-props:a}),
$b_{t,i}\seqgt b_{t,j}$, implying also that
$b_{t,i}/b_{t,j}\to+\infty$
(by \Cref{pr:asymp-dom-props}\ref{pr:asymp-dom-props:b}).%
\indexg{convergence of entries at decreasing rates|)}
\end{proof}

\indexg{sequence convergence!decomposition@to decomposition|(}%
We now return to the proof of \Cref{thm:i:seq-rep}.

\begin{proof}[Proof of \Cref{thm:i:seq-rep}:]
Condition~(\ref{thm:i:seq-rep:a})
implies, for $i=1,\dotsc,k$, that
$b_{t,i}>0$ for all but finitely many values of $t$.
By discarding all other sequence elements, we therefore can assume
that $b_{t,i}>0$ for all $t$ and all $i\in\{1,\dotsc,k\}$.

We need to show $\xx_t\cdot\uu\rightarrow\xbar\cdot\uu$ for all
$\uu\in\Rn$, which, by
\Cref{thm:i:1}(\ref{thm:i:1c}), will imply $\xx_t\rightarrow\xbar$.
Let $\uu\in\Rn$.
Suppose first that $\vv_i\inprod\uu=0$ for $i=1,\dotsc,k$.
Then
\[
   \xx_t\cdot\uu
   =
   \qq_t\cdot\uu
   \rightarrow
   \qq\cdot\uu
   =
   \xbar\cdot\uu,
\]
where the convergence in the middle follows from
condition~(\ref{thm:i:seq-rep:c}), and
the last equality from the case decomposition of $\xbar\inprod\uu$ in
\Cref{lemma:case}.

Otherwise, let $j\in\set{1,\dotsc,k}$ be the first index for which
$\vv_j\inprod\uu\neq 0$, so that
$\vv_i\inprod\uu=0$ for $i=1,\dotsc,j-1$.
Assume further that $\vv_j\inprod\uu>0$; the case
$\vv_j\inprod\uu<0$ can be proved symmetrically.
Then $\xbar\inprod\uu=+\infty$ by
\Cref{lemma:case}.

Also,
\begin{align}
  \xx_t\cdot\uu
  &=
  \sum_{i=1}^k b_{t,i} \vv_i \cdot \uu   + \qq_t\cdot\uu
\label{eqn:thm:i:seq-rep:1}
  =
  b_{t,j} \biggBracks{ \vv_j\cdot\uu
                +
                \sum_{i=j+1}^k \frac{b_{t,i}}{b_{t,j}} \vv_i\cdot\uu
                  }
         + {\qq_t\cdot\uu}.
\end{align}
By \Cref{prop:dec:trans},
$b_{t,i}/b_{t,j}\rightarrow 0$
for $i=j+1,\dotsc,k$.
Therefore, the bracketed expression in
\eqref{eqn:thm:i:seq-rep:1} is converging to $\vv_j\cdot\uu>0$.
Since $b_{t,j}\rightarrow+\infty$, by condition~(\ref{thm:i:seq-rep:a}),
it follows by continuity of multiplication
(\Cref{prop:lim:eR}\ref{i:lim:eR:genmul})
that the first term on the right-hand side of
\eqref{eqn:thm:i:seq-rep:1} is converging to $+\infty$, while
$\qq_t\cdot\uu$ is converging to $\qq\cdot\uu\in\R$
from condition~(\ref{thm:i:seq-rep:c}).
Therefore,
$\xx_t\cdot\uu\rightarrow+\infty=\xbar\cdot\uu$.

Thus, $\xx_t\cdot\uu\rightarrow\xbar\cdot\uu$ for all
$\uu\in\Rn$, so $\xx_t\rightarrow\xbar$.%
\indexg{sequence convergence!decomposition@to decomposition|)}
\end{proof}

\section{First-countability}
\label{sec:first:count}

\indexg{first-countability!astral space@of astral space|(}%
\indexg{astral space!first-countability of|(}%
\indexg{countable neighborhood base!astral point@for astral point|(}%
We next prove that astral space is \emph{first-countable}, meaning
that every point $\xbar\in\extspace$ has a {countable neighborhood base} (see
\Cref{sec:prelim:countability}).
This fundamental topological property allows us to work with
sequential characterizations of closure and continuity (see \Cref{prop:first:properties}),
and also implies that astral space is {sequentially compact}
(\Cref{prop:first:subsets}\ref{i:first:compact}).

For the proof, we use the vectors from $\xbar$'s orthogonal
decomposition (which must exist by \Cref{cor:h:1})
to construct a family of neighborhoods $\countset{B}$.
We then prove that any sequence with elements
$\xx_t\in B_t\cap\Rn$ will satisfy the
conditions of \Cref{thm:i:seq-rep}, and therefore must converge to
$\xbar$.
By \Cref{thm:first:fromseq}, this will prove that $\countset{B}$ is
a countable neighborhood base at $\xbar$.

\begin{figure}
  \centering
  \includegraphics{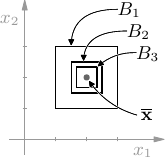}\hfill%
  \includegraphics{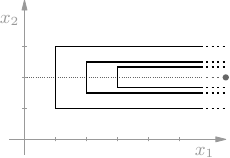}\hfill%
  \includegraphics{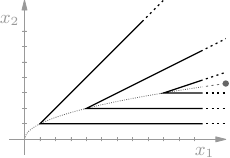}
  \mycaption{Countable neighborhood base}{%
    \indexf{countable neighborhood base!astral point@for astral point}%
      Each panel depicts the boundaries of the first three elements of
      the countable neighborhood base from \Cref{thm:first:local} for
      a specific choice of $\xbar$ in $\extspac{2}$.
      \emph{Left:}
      $\xbar=\trans{[2,2]}$.
      \emph{Center:}
      $\xbar=\limray{\ee_1}\protect\plusl 2\ee_2$.
      \emph{Right:}
      $\xbar=\limray{\ee_1}\protect\plusl\limray{\ee_2}$.}
  \label{fig:base}%
\end{figure}

\Cref{thm:first:local} explicitly lists the sets composing the
countable base for $\xbar$.
For illustration,
the first three elements of that neighborhood base are plotted
in \Cref{fig:base}
for three different points in $\extspac{2}$.

\begin{theorem}
\label{thm:first:local}
Let
$\xbar\in\extspace$ with an orthogonal decomposition
$\xbar=\limray{\vv_1} \plusl \dotsb \plusl \limray{\vv_k} \plusl \qq$,
and let $\vv_{k+1},\dotsc,\vv_{n}$ be orthonormal vectors that are orthogonal to $\vv_1,\dotsc,\vv_k$, so that $\vv_1,\dotsc,\vv_n$ is an orthonormal basis for $\Rn$.

For each $t$, let
$B_t$ be the set of all points $\zbar\in\extspace$ satisfying
all of the following conditions:
\begin{letter-compact}
\item
\label{it:count-base-cond:a}
  $\zbar\cdot\vv_i > t \;$ for $i=1,\dotsc,k$.
\item
\label{it:count-base-cond:b}
  $\zbar\cdot(\vv_i - t\vv_{i+1}) > 0 \;$ for $i=1,\dotsc,k-1$.
\item
\label{it:count-base-cond:c}
  $|\zbar\cdot\vv_j - \qq\cdot\vv_j| < 1/t \;$ for $j=k+1,\dotsc,n$.
\end{letter-compact}
Then the collection $\countset{B}$ is a nested countable neighborhood base for $\xbar$.
\end{theorem}

\begin{proof}
By \Cref{thm:i:1}(\ref{thm:i:1aa},\ref{thm:i:1b}), $\eRn$ is a regular topological space
and $\Rn$ is dense in~$\eRn$. Therefore, in order to use
\Cref{thm:first:fromseq} to show that $\countset{B}$ is a countable neighborhood base, it suffices to check that the sets
$B_t$ are neighborhoods of $\xbar$ and then verify that any sequence with elements $\xx_t\in B_t\cap\Rn$ converges to $\xbar$. We will then separately argue that the collection $\countset{B}$ is in fact nested.

To start, note that each $B_t$ is open since it is a finite intersection of subbase elements.
Also, using the orthogonality of
$\vv_1,\dotsc,\vv_k,\qq$, the case decomposition of $\xbar\inprod\uu$ in
\Cref{lemma:case} implies that:
\begin{itemize}[noitemsep]
\item
  $\xbar\cdot\vv_i=+\infty$ for $i=1,\dotsc,k$.
\item
  $\xbar\cdot(\vv_i - t\vv_{i+1}) = +\infty$ for $i=1,\dotsc,k-1$.
\item
  $\xbar\cdot\vv_j = \qq\cdot\vv_j$ for $j=k+1,\dotsc,n$.
\end{itemize}
Thus, each $B_t$ is a neighborhood of $\xbar$.

Now let $\xx_t\in B_t\cap\Rn$ for all $t$. We
prove that the sequence $\seq{\xx_t}$ satisfies the
conditions of \Cref{thm:i:seq-rep} and hence
converges to $\xbar$.
Since $\vv_1,\dotsc,\vv_k$ are orthonormal,
we can write each $\xx_t$ as
\[
  \xx_t = \sum_{i=1}^k b_{t,i} \vv_i + \qq_t,
\]
where $b_{t,i}=\xx_t\cdot\vv_i$ and $\qq_t$ is orthogonal to
$\vv_1,\dotsc,\vv_k$
(\Cref{pr:lin-decomp-rel-vecs}).

The fact that $\xx_t\in B_t$ implies that
$b_{t,i}=\xx_t\cdot\vv_i > t$,
for $i=1,\dotsc,k$,
so $b_{t,i}\to+\infty$,
and hence
condition~(\ref{thm:i:seq-rep:a}) of
\Cref{thm:i:seq-rep} is satisfied.

We also have
$\xx_t\cdot(\vv_i - t\vv_{i+1}) > 0$, for $i=1,\dotsc,k-1$,
so
\[
   b_{t,i}=\xx_t\cdot\vv_i > t(\xx_t\cdot\vv_{i+1}) = t b_{t,i+1}.
\]
As a result,
$0<b_{t,i+1}/b_{t,i} < 1/t$, and thus,
$b_{t,i+1}/b_{t,i}\rightarrow 0$.
Therefore,
condition~(\ref{thm:i:seq-rep:b}) of
\Cref{thm:i:seq-rep} is satisfied as well.

Finally, for $j=k+1,\dotsc,n$, we have
$\xx_t\cdot\vv_j=\qq_t\cdot\vv_j$.
Since $\xx_t\in B_t$, this implies
\[
   \bigAbs{\qq_t\cdot\vv_j - \qq\cdot\vv_j}
   =
   \bigAbs{\xx_t\cdot\vv_j - \qq\cdot\vv_j}
   <
   \frac{1}{t},
\]
so $\qq_t\cdot\vv_j \to \qq\cdot\vv_j$.
Since $\qq_t$ and $\qq$ are orthogonal to $\vv_1,\dotsc,\vv_k$,
we also have $\qq_t\cdot\vv_i=\qq\cdot\vv_i=0$ for $i=1,\dotsc,k$.
Thus, $\qq_t\cdot\vv_i\to \qq\cdot\vv_i$ for all basis vectors $\vv_i$, $i=1,\dotsc,n$, so $\qq_t\to\qq$
(by \Cref{pr:lim-vec-proj}),
thereby satisfying condition~(\ref{thm:i:seq-rep:c}) of
\Cref{thm:i:seq-rep}.

Having satisfied \Cref{thm:i:seq-rep}'s
three conditions,
we conclude that $\xx_t\to\xbar$. Thus, by \Cref{thm:first:fromseq}, the collection $\countset{B}$ is a countable neighborhood base at $\xbar$.

It remains to show that the collection $\countset{B}$ is nested. Let
$\zbar\in B_{t+1}$ for any $t\in\nats$.
We will show that $\zbar\in B_t$. First, note that conditions (\ref{it:count-base-cond:a}) and (\ref{it:count-base-cond:c}) for membership in~$B_t$ immediately follow from the corresponding conditions for $B_{t+1}$. 
Furthermore,
for $i=1,\ldots,k-1$,
by conditions (\ref{it:count-base-cond:b}) and (\ref{it:count-base-cond:a}) for $B_{t+1}$, the expressions   $\zbar\cdot\regParens{\vv_i - (t+1)\vv_{i+1}}$ and $\zbar\cdot\vv_{i+1}$ are
both positive and therefore summable, so by \Cref{pr:i:1},
\[
\zbar\cdot(\vv_i - t\vv_{i+1})
   =
   \zbar\cdot\BigParens{\vv_i - (t+1)\vv_{i+1}}
   +
   \zbar\cdot\vv_{i+1}
   >
   t+1
   >
   0.
\]
Hence, condition (\ref{it:count-base-cond:b}) for
$B_t$ also holds.
Thus, $\zbar\in B_t$, finishing the proof.%
\indexg{countable neighborhood base!astral point@for astral point|)}
\end{proof}

As an immediate consequence, astral space is first-countable, and thanks
to its compactness, it is also sequentially compact:

\begin{corollary}
\label{thm:first-count-and-conseq}
\indexg{sequential compactness!astral space@of astral space|(}%
\indexg{astral space!sequential compactness of|(}%
$\eRn$ is first-countable and sequentially compact.
\end{corollary}
\begin{proof}
By \Cref{thm:first:local},
every point in $\extspace$
has a countable neighborhood base, so
$\eRn$ is first-countable. Since $\eRn$ is compact and first-countable,
it is also sequentially compact, by \Cref{prop:first:subsets}(\ref{i:first:compact}).%
\indexg{first-countability!astral space@of astral space|)}%
\indexg{astral space!first-countability of|)}%
\indexg{sequential compactness!astral space@of astral space|)}%
\indexg{astral space!sequential compactness of|)}
\end{proof}

\section{Not second-countable and not metrizable}
\label{sec:not-second}

Although $\extspace$ is first-countable, we show next that it is not
second-countable (for $n\geq 2$),
meaning its topology does not have a countable base.
By \Cref{prop:sep:metrizable}, this will further imply
that it also is not metrizable.

\indexg{astrons!topological isolation of|(}%
To prove this,
for each $\vv\in\Rn$, we will define in the next theorem
an open set $\Uv$ with
the property that $\Uv$ includes the astron $\limray{\vv}$, but does
not include any other astron.
More specifically, we show that all of the points in $\Uv$ must either
be in $\Rn$ or have the form $\limray{\vv}\plusl\qq$ for some
$\qq\in\Rn$, thereby excluding all other nonzero astrons.
These properties will imply that a countable base for $\extspace$
must include a distinct open set for every astron, which is
impossible since the set of all astrons is uncountable when $n\geq 2$. 

\begin{figure}
  \centering
  \includegraphics{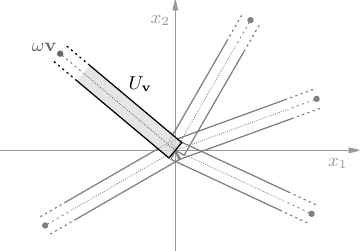}
  \mycaption{%
    Construction from \Cref{thm:formerly-lem:h:1:new} in two dimensions}{%
    \indexf{astrons!topological isolation of}%
    Open sets $\Uv$ for several choices of $\vv\in\R^2$, $\norm{\vv}=1$.
    Each $\Uv$ contains the corresponding astron~$\limray{\vv}$, but no other astron.
  }
  \label{fig:astrons}%
\end{figure}

\begin{theorem}  \label{thm:formerly-lem:h:1:new}
  For every $\vv\in\Rn$,
  there exists an open set $\Uv\subseteq\extspace$
  that includes
  $\limray{\vv}$ and for which the following hold:
  \begin{letter-compact}
  \item  \label{thm:formerly-lem:h:1:b:new}
    $\Uv\subseteq \Rn \cup [\limray{\vv}\plusl\Rn]$.
  \item  \label{thm:formerly-lem:h:1:c:new}
     If $\ww\in\Rn$ and $\norm{\vv}=\norm{\ww}=1$, then
     $\limray{\ww}\in\Uv$ if and only if $\ww=\vv$.
  \end{letter-compact}
\end{theorem}

\begin{proof}
Let $\vv\in\Rn$.
If $\vv=\zero$, we can choose $\Uzero=\Rn$, which
satisfies
part~(\ref{thm:formerly-lem:h:1:b:new}) trivially
and
part~(\ref{thm:formerly-lem:h:1:c:new}) vacuously.
When $\vv\neq\zero$,
it suffices to consider only the case that $\vv$ is a unit
vector since
if $\norm{\vv}\neq 1$, then
we can choose $\Uv$ to be the same as the corresponding set
for a normalized version of $\vv$; that is, we can choose
$\Uv=\Unormv$
(noting that $\limray{(\lambda\vv)}=(\oms\lambda)\vv=\limray{\vv}$
for $\lambda\in\Rstrictpos$).
Therefore, we assume henceforth that
$\norm{\vv}=1$.

Let $\uu_1,\dotsc,\uu_{n-1}$ be
any orthonormal basis for the linear space orthogonal to $\vv$;
thus, $\vv,\uu_1,\dotsc,\uu_{n-1}$ form an orthonormal basis for all
of $\Rn$.
We define $\Uv$ to be the set of all $\xbar\in\extspace$ that satisfy both of the following conditions:
\begin{item-compact}
\item
  $\xbar\cdot\vv>0$.
\item
  $|\xbar\cdot\uu_j| < 1$ for $j=1,\dotsc,n-1$.
\end{item-compact}
Then $\Uv$ is open, since it is a finite intersection of subbase elements.
Furthermore,
$\limray{\vv}\cdot\vv=+\infty$ and
$\limray{\vv}\cdot\uu_j=0$ for $j=1,\dotsc,n-1$;
therefore, $\limray{\vv}\in\Uv$.
See \Cref{fig:astrons} for an illustration of this contruction in $\extspac{2}$.

Both parts of the theorem will follow from the next claim:

\begin{claimpx}  \label{cl:thm:formerly-lem:h:1:new:1}
Let $\xbar\in\Uv\setminus\Rn$, and suppose that
$\xbar=\limray{\vv_1}\plusl\dotsc\plusl\limray{\vv_k}\plusl\qq$
is an orthogonal decomposition of $\xbar$.
Then $k=1$ and $\vv_1=\vv$ (so $\xbar=\limray{\vv}\plusl\qq$).
\end{claimpx}

\begin{proofx}
Since $\xbar\not\in\Rn$, $k>0$.

For $j\in\{1,\dotsc,n-1\}$, if
$\vv_i\cdot\uu_j\neq 0$ for any $i\in\{1,\dotsc,k\}$,
then the case decomposition of $\xbar\inprod\uu_j$
from \Cref{lemma:case}
yields
$\xbar\cdot\uu_j\in\{-\infty,+\infty\}$,
contradicting that $|\xbar\cdot\uu_j|<1$
since $\xbar\in\Uv$.

Therefore,
$\vv_i\cdot\uu_j=0$ for all $i$, $j$.
Since $\vv_i$ is a unit vector, and
since $\vv,\uu_1,\dotsc,\uu_{n-1}$ form an orthonormal basis for
$\Rn$, the only possibility then
is that each $\vv_i$ is either equal to $\vv$ or $-\vv$.
Furthermore,
since the $\vv_i$'s are orthogonal to one another, this further
implies that $k=1$, with $\vv_1$ either $\vv$ or $-\vv$.
But if $\vv_1=-\vv$, then
$\xbar\cdot\vv=\limray{\vv_1}\cdot\vv\plusl\qq\cdot\vv=-\infty$,
contradicting that
$\xbar\cdot\vv>0$ since $\xbar\in\Uv$.
Therefore, $\vv_1=\vv$ and
$\xbar=\limray{\vv}\plusl\qq$.
\end{proofx}

We now prove the two parts of the theorem:

\begin{proof-parts}
\pfpart{Part~(\ref{thm:formerly-lem:h:1:b:new}):}
Let $\xbar$ in $\Uv$.
If $\xbar\not\in\Rn$, then we can apply
\Cref{cl:thm:formerly-lem:h:1:new:1} to an orthogonal
decomposition of $\xbar$
(which must exist by \Cref{cor:h:1}), implying
$\xbar\in\limray{\vv}\plusl\Rn$, proving the claim.

\pfpart{Part~(\ref{thm:formerly-lem:h:1:c:new}):}
Let $\ww\in\Rn$ be such that $\norm{\ww}=1$.
If $\ww=\vv$ then $\limray{\ww}\in\Uv$, as observed above.

For the converse,
let $\xbar=\limray{\ww}$ and assume that $\xbar\in\Uv$.
Then $\xbar\not\in\Rn$ since
$\xbar\cdot\ww=\limray{\ww}\cdot\ww=+\infty$.
Therefore, applying
\Cref{cl:thm:formerly-lem:h:1:new:1}
to $\limray{\ww}$, which is an orthogonal decomposition of $\xbar$,
yields that $\ww=\vv$.%
\indexg{astrons!topological isolation of|)}
\qedhere
\end{proof-parts}
\end{proof}

\begin{theorem} \label{thm:h:8:new}
\indexg{second-countability!astral space not second-countable|(}%
\indexg{astral space!not second-countable|(}%
  For $n\geq 2$, astral space $\extspace$ is not second-countable.
\end{theorem}

\begin{proof}
Let $n\geq 2$, and define the set
\[
  C
  =
  \bigBraces{
    \transk{[\,\cos\theta,\, \sin\theta,\, \underbrace{0, \dotsc, 0}_{n-2}\,]}
    :\:
    \theta\in[0,2\pi)
    },
\]  
which is a circle of radius~1 in $\Rn$.
Then $C$ is uncountable, being in bijection with the interval
$[0,2\pi)$
(which is itself uncountable by \Cref{pr:uncount-interval}).

For each $\vv\in C$,
let $\Uv$ be as in \Cref{thm:formerly-lem:h:1:new}.
Suppose, contrary to the theorem's claim,
that there exists a countable base
$\countset{B}$
for $\extspace$.
For each $\vv\in C$, $\Uv$ is a neighborhood of $\limray{\vv}$,
which implies,
by \Cref{pr:base-equiv-topo},
that there exists
an index $i(\vv)\in\nats$ such that
$\limray{\vv}\in B_{i(\vv)} \subseteq \Uv$.
The resulting function $i: C\rightarrow\nats$ is then injective, because if
$i(\vv)=i(\ww)$ for some $\vv,\ww\in C$, then
\[ \limray{\ww}\in B_{i(\ww)} = B_{i(\vv)} \subseteq \Uv, \]
which implies $\vv=\ww$ by
\Cref{thm:formerly-lem:h:1:new}(\ref{thm:formerly-lem:h:1:c:new}).
Therefore,
by \Cref{pr:count-equiv}(\ref{pr:count-equiv:c},\ref{pr:count-equiv:a}),
$C$ must be countable,
which is a contradiction
with $C$ being uncountable as established above,
thus completing the
\indexg{second-countability!astral space not second-countable|)}%
\indexg{astral space!not second-countable|)}%
proof.
\end{proof}

\indexg{metrizability!astral space not metrizable|(}%
\indexg{astral space!not metrizable|(}%
As a corollary,
$\extspace$ also cannot be metrizable when $n\ge 2$:

\begin{corollary}
For $n\geq 2$, astral space, $\extspace$, is not metrizable.
\end{corollary}

\begin{proof}
Assume $n\geq 2$.
We previously saw that
$\Qn$ is a countable dense subset of~$\Rn$
(Example~\ref{ex:2nd-count-rn}),
and that $\Rn$ is dense in $\extspace$
(\Cref{thm:i:1}\ref{thm:i:1b}).
Therefore,
$\Qn$ is a countable dense subset of $\eRn$
by \Cref{prop:subspace}(\ref{i:subspace:dense})
(and since $\Rn$ is a subspace of $\extspace$ by
\Cref{pr:open-sets-equiv-alt}).
Every metrizable space with a countable dense subset is second-countable
(\Cref{prop:sep:metrizable}),
but $\extspace$ is not second-countable
(by \Cref{thm:h:8:new}).
Therefore, $\extspace$ is not metrizable.%
\indexg{metrizability!astral space not metrizable|)}%
\indexg{astral space!not metrizable|)}
\end{proof}

\indexg{astrons!topological isolation of|(}%
\Cref{thm:formerly-lem:h:1:new}
highlights an important
topological property regarding astrons.
While the set of infinite astrons is in bijection with the unit sphere,
it is topologically quite different. On the unit sphere, any unit vector
$\vv$ is arbitrarily close to other unit vectors, meaning that any of its neighborhoods,
no matter how small, will include other unit vectors.
\Cref{thm:formerly-lem:h:1:new}
shows that this is not true for infinite astrons.
Rather, every astron $\limray{\vv}$ (with $\norm{\vv}=1$) has a neighborhood $\Uv$ that excludes
\emph{all} other astrons $\limray{\ww}$ (with $\norm{\ww}=1$, $\ww\ne\vv$), no matter how tiny the
distance between $\vv$ and $\ww$. Thus, unlike points on the unit sphere,
infinite astrons are topologically isolated from one another.
We will consider related issues again in \Cref{sec:conv-in-dir}.%
\indexg{astrons!topological isolation of|)}

\chapter{Linear maps and the matrix representation of astral points}
\chaptermark{Linear maps and matrix representation}
\label{sec:representing}

Astral space was specifically constructed to allow continuous extensions
of linear functions mapping $\Rn$ to $\R$, each corresponding 
to the operation of taking an inner product
$\xx\mapsto\xx\inprod\uu$ for some $\uu\in\Rn$. In this chapter,
we will see more broadly that all linear maps that map $\Rn$ to $\Rm$ for any $m\ge 0$
can be extended continuously to astral space. Such linear maps correspond to
matrix-vector multiplication $\xx\mapsto\A\xx$ for some $\A\in\Rnm$.
With the use of extended linear maps, we can adapt various notions
and operations from linear algebra to
all of astral space, including to the points at infinity.
Indeed, we will see in this chapter that every astral point can be represented using
standard matrices, with properties of those matrices related directly
to properties of the represented astral point.
Although such matrix representation of a point is not unique,
it is always possible to pick a matrix that satisfies additional orthonormality conditions
and obtain a canonical representation, which is unique.

\section{Linear and affine maps}
\label{sec:linear-maps}

{\mathtogether%
\indexg{linear maps, astral|(}%
\indexg{matrix product, astral|(}%
The
coupling function, $\xbar\inprodk\uu$, viewed as a function
of~$\xbar$, was constructed to be a continuous extension of the
standard inner product,
that is, of the linear function $\xx\mapsto\xx\inprodk\uu$
(see \Cref{thm:i:1}\ref{thm:i:1c}).
In this section, we will see that it is similarly possible to extend
all linear maps from $\Rn$ to $\Rm$ to maps from $\eRn$ to $\eRm$.
To begin, for any matrix $\A\in\Rmn$ and point $\xbar\in\extspace$,
we consider in the next \namecref{thm:mat-mult-def:alt}
sequences of the form $\seq{\A\xx_t}$, where $\seq{\xx_t}$ is a
sequence in $\Rn$ converging to $\xbar$, and show that
all such sequences must have a common, uniquely determined limit
in $\extspac{m}$.}

\begin{theorem}  \label{thm:mat-mult-def:alt}
  Let $\A\in\Rmn$, and let $\xbar\in\extspace$.
  Then there exists a unique point $\zbar\in\extspac{m}$
  such that
  for every sequence $\seq{\xx_t}$ in $\Rn$,
  if $\xx_t\rightarrow\xbar$, then $\A\xx_t\rightarrow \zbar$.
  Furthermore, if $\xbar=\xx\in\Rn$ then $\zbar=\A\xx$.
\end{theorem}
\begin{proof}
Let $\seq{\xx_t}$ be any sequence in $\Rn$ that converges to $\xbar$,
which must exist by \Cref{thm:i:1}(\ref{thm:i:1d}).
Then for all $\uu\in\Rm$,
\begin{equation} \label{eq:i:5:alt}
  (\A\xx_t)\cdot\uu
  = \trans{(\A\xx_t)}\uu
  = \trans{\xx}_t \transAk\uu
  = \xx_t\cdot(\transAk\uu)
  \rightarrow \xbar\cdot (\transAk\uu),
\end{equation}
where the last step follows by continuity of the
coupling function in its first argument
(\Cref{thm:i:1}\ref{thm:i:1c}).
Thus, for all $\uu\in\Rm$, the sequence $(\A\xx_t)\cdot\uu$ converges in $\eR$, and therefore,
by \Cref{thm:i:1}(\ref{thm:i:1e}), the sequence
$\seq{\A\xx_t}$ has a limit $\zbar\in\eRm$.

Since $\A\xx_t\to\zbar$, we must have
$(\A\xx_t)\cdot\uu\rightarrow\zbar\cdot\uu$ for all $\uu\in\Rm$
(again by \Cref{thm:i:1}\ref{thm:i:1c}).
As we saw in \eqref{eq:i:5:alt},
this same sequence, $(\A\xx_t)\cdot\uu$, also converges
to $\xbar\cdot (\transAk\uu)$, so
we have $\zbar\cdot\uu=\xbar\cdot (\transAk\uu)$ for all $\uu\in\Rm$.

If $\seq{\xx'_t}$ is any other sequence in $\Rn$ converging to $\xbar$, then
this same reasoning shows that $\seq{\A\xx'_t}$ must have some limit
$\zbar'\in\eRm$ satisfying
$\zbar'\cdot\uu=\xbar\cdot (\transAk\uu)$
for all $\uu\in\Rm$. Thus, $\zbar\cdot\uu=\zbar'\cdot\uu$ for all $\uu\in\Rm$,
so $\zbar=\zbar'$ (by \Cref{pr:i:4}),
proving that $\zbar$ is unique.

Moreover, if $\xbar=\xx\in\Rn$, then $\A\xx_t\to\A\xx$ by continuity of standard 
linear maps, so $\zbar=\A\xx$ as claimed.
\end{proof}

In light of this theorem, we can define the action of the matrix $\A$ on any point $\xbar\in\extspace$ as follows:

\begin{definition}
\label{def:astral-linear-map}
\indexg{linear maps, astral!defined|(}%
\indexg{matrix product, astral!defined|(}%
  Let $\A\in\Rmn$. For any $\xbar\in\extspace$,
  we define
\indexm{A x}{$\A\xbar$}{multiplication by matrix}%
  $\A\xbar$ to be the unique point in $\extspac{m}$
  such that for every sequence $\seq{\xx_t}$ in $\Rn$,
  if $\xx_t\rightarrow\xbar$, then $\A\xx_t\rightarrow \A\xbar$.
  The mapping $\xbar\mapsto\A\xbar$ from $\eRn$ to $\eRm$
  is called the \emph{astral linear map} associated with $\A$.\looseness=-1
\end{definition}

\indexg{matrix product, astral!defined|)}%
\indexg{linear maps, astral!defined|)}%
Succinctly,
\[
   \A\xbar = \lim \A\xx_t
\]
for any (and every) sequence $\seq{\xx_t}$ in $\Rn$
converging to $\xbar$.

By \Cref{thm:mat-mult-def:alt}, astral linear maps are well-defined and agree
with standard linear maps on $\Rn$, so if $\xbar=\xx$ for some $\xx\in\Rn$ then
$\A\xbar$ is equal to the standard matrix-vector product $\A\xx$.
The following property of astral linear maps was established in the proof of \Cref{thm:mat-mult-def:alt}:

\begin{theorem}
\label{thm:Ax-dot-u}
\indexg{matrix product, astral!coupling function of|(}%
  Let $\A\in\Rmn$ and let $\xbar\in\extspace$.
  Then for all $\uu\in\Rm$,
  \[
    (\A\xbar)\cdot\uu = \xbar\cdot(\transAk\uu).
  \]
\end{theorem}

\begin{proof}
Let $\seq{\xx_t}$ be any sequence in $\Rn$ that converges to $\xbar$
(which exists by \Cref{thm:i:1}\ref{thm:i:1d}). Following similar reasoning
as in the proof of \Cref{thm:mat-mult-def:alt}, we obtain that for all $\uu\in\Rm$,
\[
  (\A\xx_t)\inprod\uu
  =\trans\xx_t\transAk\uu
  =\xx_t\inprod(\transAk\uu)
  \to
  \xbar\inprod(\transAk\uu),
\]
where the last step follows by \Cref{thm:i:1}(\ref{thm:i:1c}).
Since $\A\xx_t\to\A\xbar$
(by \Cref{def:astral-linear-map}), we must also have $(\A\xx_t)\cdot\uu\to(\A\xbar)\cdot\uu$
(by \Cref{thm:i:1}\ref{thm:i:1c}). Thus, $(\A\xbar)\cdot\uu=\xbar\cdot(\transAk\uu)$,
proving the \namecref{thm:Ax-dot-u}.%
\indexg{matrix product, astral!coupling function of|)}%
\end{proof}

\indexg{linear maps, astral!continuity of|(}%
\indexg{continuity!astral linear maps@of astral linear maps|(}%
The next theorem shows that astral linear maps are
unique continuous extensions of the corresponding standard linear maps to astral space:

\begin{theorem}   \label{thm:linear:cont}
  Let $\A\in\Rmn$, and let
  $\slinmapA:\Rn\rightarrow\Rm$ and $\alinmapA:\extspace\rightarrow\extspac{m}$
  be the associated standard and
  astral linear maps
  defined by $\slinmapA(\xx)=\A\xx$ for $\xx\in\Rn$
  and $\alinmapA(\xbar)=\A\xbar$ for $\xbar\in\extspace$.
  Then:
  \begin{letter-compact}
  \item   \label{thm:linear:cont:a}
    $\alinmapA(\xx)=\slinmapA(\xx)$ for $\xx\in\Rn$.
  \item   \label{thm:linear:cont:b}
    $\alinmapA$ is continuous.
  \end{letter-compact}
  Moreover, $\alinmapA$ is the only function mapping $\eRn$ to $\eRm$
  with both of these properties.
\end{theorem}

\begin{proof}
  ~

\begin{proof-parts}  
\pfpart{Part~(\ref{thm:linear:cont:a}):}
Immediate by \Cref{thm:mat-mult-def:alt}.

\pfpart{Part~(\ref{thm:linear:cont:b}):}
Let $\seq{\xbar_t}$ be any
sequence in $\extspace$ converging to some point $\xbar\in\extspace$.
Then for all $\uu\in\Rm$,
\[
   \alinmapA(\xbar_t)\cdot\uu
   =
   (\A\xbar_t)\cdot\uu
   =
   \xbar_t\cdot(\transAk\uu)
   \rightarrow
   \xbar\cdot(\transAk\uu)
   =
   (\A\xbar)\cdot\uu
   =
   \alinmapA(\xbar)\cdot\uu,
\]
where
the second and third equalities are by
\Cref{thm:Ax-dot-u}, and
the convergence is by \Cref{thm:i:1}(\ref{thm:i:1c}).
It follows that
$\alinmapA(\xbar_t)\rightarrow \alinmapA(\xbar)$
(again by \Cref{thm:i:1}\ref{thm:i:1c}).
Thus, $\alinmapA$ is continuous
(by \Cref{prop:first:properties}\ref{prop:first:cont}
and first-countability).

\pfpart{Uniqueness:}
Suppose $\linmapalt:\eRn\rightarrow\eRm$ is another function with these same
properties, and let $\xbar\in\extspace$.
Let $\seq{\xx_t}$ be any sequence in $\Rn$ converging to $\xbar$
(which exists by \Cref{thm:i:1}\ref{thm:i:1d}).
Then
\[
   \linmapalt(\xbar)
   =
   \lim \linmapalt(\xx_t)
   =
   \lim \slinmapA(\xx_t)
   =
   \lim \alinmapA(\xx_t)
   =
   \alinmapA(\xbar),
\]
where the first and last equalities are by
property~(\ref{thm:linear:cont:b}),
and the second and third are by
property~(\ref{thm:linear:cont:a}).
Thus, $\linmapalt=\alinmapA$.%
\indexg{linear maps, astral!continuity of|)}%
\indexg{continuity!astral linear maps@of astral linear maps|)}
\qedhere
\end{proof-parts}
\end{proof}

\indexg{column space!astral|(}%
The (standard) column space of $\A$ is the span of its columns, or
equivalently, the image of $\Rn$ under the linear map $\xx\mapsto\A\xx$:
\[
  \colspace\A
  =
  \set{\A\xx :\: \xx\in\Rn}.
\]
Analogously, we define the \emph{astral column space} of $\A$,
denoted $\acolspace\A$,%
\indexm{colstar}{$\acolspace\A$}{astral column space}
to be the image of $\eRn$ under the astral linear map $\xbar\mapsto\A\xbar$:
\[
  \acolspace\A
  =
  \set{\A\xbar:\:\xbar\in\eRn}.%
\indexg{column space!astral|)}
\]

\indexg{coupling function!matrix product@as matrix product|(}%
\indexg{linear maps, astral!coupling function as|(}%
If $\uu\in\Rn$ and $\xbar\in\extspace$,
then \Cref{thm:linear:cont}(\ref{thm:linear:cont:b}) implies that
$\trans{\uu}\xbar$, which is a point in $\extspac{1}=\Rext$, is the
same as $\xbar\cdot\uu$:

\begin{proposition}   \label{pr:trans-uu-xbar}
  Let $\xbar\in\extspace$ and $\uu\in\Rn$.
  Then $\trans{\uu}\xbar=\xbar\cdot\uu$.
\end{proposition}

\begin{proof}
Let $\seq{\xx_t}$ be any sequence in $\Rn$ converging to $\xbar$.
Then
\[
   \trans{\uu}\xbar
   = \lim (\trans{\uu} \xx_t)
   = \lim (\xx_t\cdot\uu)
   = \xbar\cdot\uu,
\]
with the first equality from
\Cref{thm:linear:cont}(\ref{thm:linear:cont:b})
(with $\A=\trans{\uu}$), and the third from
\Cref{thm:i:1}(\ref{thm:i:1c}).%
\indexg{coupling function!matrix product@as matrix product|)}%
\indexg{linear maps, astral!coupling function as|)}
\end{proof}

Here are some basic properties of astral linear maps. Note that the first three generalize analogous properties of the coupling function
(Propositions~\ref{pr:i:6},~\ref{pr:i:2} and~\ref{pr:astrons-exist}).

\begin{proposition} \label{pr:h:4}
  Let $\A\in\Rmn$, and let $\xbar\in\extspace$.
  Then:
  \begin{letter-compact}
  \item  \label{pr:h:4c}
    $\A(\xbar\plusl\ybar) = \A\xbar \plusl \A\ybar$, for $\ybar\in\extspace$.
  \item  \label{pr:h:4lam-mat}
    $\A(\lambda\xbar) = \lambda(\A\xbar) = (\lambda \A) \xbar$,
    for $\lambda\in\R$.
  \item  \label{pr:h:4g}  \label{pr:h:4e}
    $\A(\alpha\xbar) = \alpha(\A\xbar)$,
    for $\alpha\in\Rext$.
  \item  \label{pr:h:4d}
    $\B(\A\xbar) = (\B\A)\xbar$, for $\B\in\R^{\ell\times m}$.
  \item  \label{pr:h:4:iden}
    $\Iden\,\xbar = \xbar$, where $\Iden$ is the $n\times n$ identity matrix.
  \end{letter-compact}
\end{proposition}

\begin{proof}
~
\begin{proof-parts}
\pfpart{Part~(\ref{pr:h:4c}):}
By \Cref{thm:Ax-dot-u} and \Cref{pr:i:6},
for all $\uu\in\Rn$,
\begin{align*}
  \smash{\bigParens{\A(\xbar\plusl\ybar)}\cdot\uu}
  &= (\xbar\plusl\ybar)\cdot(\transAk\uu)
\\
  &= \xbar\cdot(\transAk\uu)
    \plusl
    \ybar\cdot(\transAk\uu)
\\
  &= (\A\xbar)\cdot\uu
    \plusl
    (\A\ybar)\cdot\uu
\\
  &= (\A\xbar \plusl \A\ybar)\cdot\uu.
\end{align*}
The claim now follows by \Cref{pr:i:4}.
\end{proof-parts}

The proofs of the other parts are similar.%
\indexg{linear maps, astral|)}
\end{proof}

\indexg{matrix product, astral!row-by-row calculation of|(}%
When multiplying a vector $\xx\in\Rn$ by a matrix $\A\in\Rmn$,
the result $\A\xx$ is a vector in $\Rm$ whose $j$-th component is
computed by multiplying $\xx$ by the $j$-th row of~$\A$.
In the same way, for an astral point $\xbar\in\extspace$,
we might expect that $\A\xbar$ in $\extspac{m}$
can be viewed as the result of multiplying $\xbar$ by the
individual rows of~$\A$.
However, this is not the case.
Matrix multiplication of astral points is far more holistic,
retaining much more information about the astral points $\xbar$
being multiplied, as shown in the next example:

\begin{example}   \label{ex:mat-prod-not-just-row-prods}
With $m=n=2$, let
$\zbar=\A\xbar$ where
\[
  \A =
  \begin{bmatrix*}[r]
                  1 & 0 \\
                 -1 & 1 \\
  \end{bmatrix*}
\]
and $\xbar=\limray{\ee_1}\plusl \beta \ee_2$,
for some $\beta\in\R$.
Let $\trans{\aaa}_1=[1, 0]$ and $\trans{\aaa}_2=[-1,1]$
denote the rows of $\A$.
Multiplying these rows separately by $\xbar$,
we obtain by the Case Decomposition Lemma~\ref{lemma:case}
that $\trans{\aaa}_1\xbar=\xbar\cdot\aaa_1=+\infty$
and
$\trans{\aaa}_2\xbar=\xbar\cdot\aaa_2=-\infty$.
From this,
it might appear that information about $\beta$ has been
entirely erased by the process of multiplying by $\A$.
But this is incorrect.
Indeed, using \Cref{pr:h:4}, we can compute
\[
  \zbar=\A\xbar=\limray{(\A\ee_1)}\plusl\beta(\A\ee_2)
      =\limray{\vv}\plusl\beta\ee_2
\]
where $\vv=\trans{[1, -1]}$.
As a result,
the value of $\beta$ can be readily extracted from $\zbar$
by coupling it with $\uu=\trans{[1,1]}$ since
\[
  \zbar\cdot\uu
  =
  \limray{(\vv\cdot\uu)}\plusl\beta\ee_2\cdot\uu
  =
  \beta.
\]
Alternatively, since this matrix $\A$ is invertible,
$\xbar=\Ainv \zbar$, by
\Cref{pr:h:4}(\ref{pr:h:4d},\ref{pr:h:4:iden}),
which means $\xbar$ (including the value of $\beta$)
can be entirely recovered from $\zbar$, though clearly not from
the row products,
$\trans{\aaa}_1\xbar$
and
$\trans{\aaa}_2\xbar$.
\end{example}

Thus, in general, $\A\xbar$ need not be fully determined by the
corresponding row products, as described above.
Nevertheless, as shown next,
$\A\xbar$ \emph{is} fully determined by those row products if,
unlike in Example~\ref{ex:mat-prod-not-just-row-prods},
they are all finite, or if $\A\xbar$ is finite.

\begin{proposition}   \label{pr:mat-prod-is-row-prods-if-finite}
  Let $\A\in\Rmn$, let $\xbar\in\extspace$, and let $\bb\in\Rm$.
  For $j=1,\dotsc,m$, let $\trans{\aaa}_j$ be the $j$-th row of
  $\A$
  (so that $\transA=[\aaa_1,\dotsc,\aaa_m]$ and each $\aaa_j\in\Rn$).
  Then $\A\xbar=\bb$ if and only if $\xbar\cdot\aaa_j=b_j$ for
  $j=1,\dotsc,m$.
\end{proposition}

\begin{proof}
Suppose first that $\A\xbar=\bb$.
Then for $j=1,\dotsc,m$,
\[
  b_j
  =
  \bb\cdot\ee_j
  =
  (\A\xbar)\cdot\ee_j
  =
  \xbar\cdot(\transAk\ee_j)
  =
  \xbar\cdot\aaa_j
\]
where the third equality is by
\Cref{thm:Ax-dot-u}
(and $\ee_1,\dotsc,\ee_m$ are the standard basis vectors in $\Rm$).

For the converse, suppose now that $\xbar\cdot\aaa_j=b_j$ for
$j=1,\dotsc,m$.
Let $\uu=\trans{[u_1,\dotsc,u_m]}\in\Rm$.
Then for all $j$,
\begin{equation}   \label{eq:pr:mat-prod-is-row-prods-if-finite:1}
  \xbar\cdot(u_j \aaa_j) = u_j (\xbar\cdot\aaa_j) = u_j b_j.
\end{equation}
Since these values are all in $\R$,
$\xbar\cdot(u_1 \aaa_1),\dotsc,\xbar\cdot(u_m \aaa_m)$
are summable.
Consequently,
\begin{align}
\notag  
  (\A\xbar)\cdot\uu
  =
  \xbar\cdot(\transAk\uu)
  &=
  \xbar\cdot\BiggParens{\,\sum_{j=1}^m u_j \aaa_j }
\\
  &=
  \sum_{j=1}^m \xbar\cdot(u_j \aaa_j)
  =
  \sum_{j=1}^m u_j b_j
  =
  \bb\cdot\uu.
  \label{eq:pr:mat-prod-is-row-prods-if-finite:2}
\end{align}
The first equality is by
\Cref{thm:Ax-dot-u},
the third by \Cref{pr:i:1}, and the fourth by
\eqref{eq:pr:mat-prod-is-row-prods-if-finite:1}.
Since
\eqref{eq:pr:mat-prod-is-row-prods-if-finite:2}
holds for all $\uu\in\Rm$, it follows that $\A\xbar=\bb$
(by \Cref{pr:i:4}).%
\indexg{matrix product, astral!row-by-row calculation of|)}%
\indexg{matrix product, astral|)}%
\end{proof}

\indexg{linear maps, astral!characterization|(}%
\Cref{pr:h:4}
lists some basic properties of the astral linear map
$\xbar\mapsto \A\xbar$.
As we show next, if a function $F:\Rn\rightarrow\Rm$ satisfies
a subset of these properties, then actually $F$ must be an astral linear map,
that is, we must have $F(\xbar)=\A\xbar$ for some matrix $\A\in\Rmn$.
Thus, the properties listed in the next theorem are both necessary
and sufficient for a map to be astral linear.

\begin{theorem}  \label{thm:linear-trans-is-matrix-map}
  Let $F:\eRn\to\eRm$ satisfy the following:
  \begin{letter-compact}
  \item  \label{thm:linear-trans-is-matrix-map:a}
    $F(\xbar\plusl\ybar)=F(\xbar)\plusl F(\ybar)$
    for all $\xbar,\ybar\in\extspace$.
  \item  \label{thm:linear-trans-is-matrix-map:b}
    $F(\alpha \xx) = \alpha F(\xx)$ for all
    $\alpha\in\Rext$ and $\xx\in\Rn$.
  \end{letter-compact}
  Then there exists a unique matrix $\A\in\Rmn$ for which
  $F(\xbar)=\A\xbar$
  for all $\xbar\in\extspace$.
\end{theorem}

\begin{proof}
We prove first that $F(\Rn)\subseteq \Rm$.
Suppose to the contrary that $\ybar=F(\xx)$ for some $\xx\in\Rn$,
and some $\ybar\in\extspac{m}\setminus\Rm$.
Then
\begin{equation}
\label{eqn:thm:linear-trans-is-matrix-map:2}
  \zerov{m}
  =
  F(\zerov{n})
  =
  F\bigParens{\xx\plusl (-\xx)}
  =
  F(\xx) \plusl F(-\xx)
  =
  \ybar \plusl F(-\xx).
\end{equation}
The first equality is because, by
condition~(\ref{thm:linear-trans-is-matrix-map:b}),
${F(\zerov{n})=F(0\cdot\zerov{n})=0\cdot F(\zerov{n}) = \zerov{m}}$.
The
second equality follows because leftward addition extends standard addition (\Cref{pr:i:7}\ref{pr:i:7e}), and
the third equality follows from
condition~(\ref{thm:linear-trans-is-matrix-map:a}).

Since $\ybar\not\in\Rm$, we must have
$\ybar\cdot\uu\not\in\R$
for some $\uu\in\Rm$
(by \Cref{pr:i:3}).
Without loss of generality, we assume $\ybar\cdot\uu=+\infty$
(since if $\ybar\cdot\uu=-\infty$, we can replace $\uu$ with $-\uu$).
Thus, by
\eqref{eqn:thm:linear-trans-is-matrix-map:2},
\[
   0
   =
   \zerov{m} \cdot\uu
   =
   \ybar\cdot\uu \plusl F(-\xx)\cdot\uu
   =
   +\infty,
\]
a clear contradiction.

Having established  $F(\Rn)\subseteq \Rm$,
let $\A\in\Rmn$ be the matrix
\[ \A=[F(\ee_1),\dotsc,F(\ee_n)], \]
that is, the matrix whose $i$-th column is $F(\ee_i)$,
where $\ee_1,\dotsc,\ee_n$ are the standard basis vectors in $\Rn$.
Then for any vector
$\xx=\trans{[x_1,\dotsc,x_n]}$
in $\Rn$, we have
\begin{equation}
\label{eqn:thm:linear-trans-is-matrix-map:1}
  F(\xx)
  =
  F\paren{\sum_{i=1}^n x_i \ee_i}
  =
  \sum_{i=1}^n F(x_i \ee_i)
  =
  \sum_{i=1}^n x_i F(\ee_i)
  =
  \A\xx
\end{equation}
with the second equality following from repeated application
of condition~(\ref{thm:linear-trans-is-matrix-map:a}),
the third from
condition~(\ref{thm:linear-trans-is-matrix-map:b}),
and the last equality from $\A$'s definition.

Next, let $\xbar\in\extspace$.
By \Cref{cor:h:1}, we can write
$\xbar= \limray{\vv_1} \plusl \dotsb \plusl \limray{\vv_k} \plusl \qq$
for some $\qq,\vv_1,\dotsc,\vv_k\in\Rn$.
Then
\begin{align*}
  F(\xbar)
  &=
  F(\limray{\vv_1} \plusl \dotsb \plusl \limray{\vv_k} \plusl \qq)
  \\
  &=
  F(\limray{\vv_1}) \plusl \dotsb \plusl F(\limray{\vv_k}) \plusl F(\qq)
  \\
  &=
  \limray{\bigParens{F(\vv_1)}} \plusl \dotsb \plusl \limray{\bigParens{F(\vv_k)}} \plusl F(\qq)
  \\
  &=
  \limray{(\A\vv_1)} \plusl \dotsb \plusl \limray{(\A\vv_k)} \plusl \A\qq
  \\
  &=
  \A(\limray{\vv_1}) \plusl \dotsb \plusl \A(\limray{\vv_k}) \plusl \A\qq
  \\
  &=
  \A\regParens{
    \limray{\vv_1} \plusl \dotsb \plusl \limray{\vv_k} \plusl \qq
  }
  \\
  &=
  \A\xbar.
\end{align*}
The second and third equalities are by
conditions~(\ref{thm:linear-trans-is-matrix-map:a})
and~(\ref{thm:linear-trans-is-matrix-map:b}), respectively.
The fourth equality is by
\eqref{eqn:thm:linear-trans-is-matrix-map:1}.
The fifth and sixth equalities are by
\Cref{pr:h:4}(\ref{pr:h:4e},\ref{pr:h:4c}).

Thus, $F(\xbar)=\A\xbar$ for all $\xbar\in\extspace$.
If, for some matrix $\B\in\Rmn$, we also have
$F(\xbar)=\B\xbar$ for all $\xbar\in\extspace$,
then $\A\ee_i=F(\ee_i)=\B\ee_i$ for $i=1,\dotsc,n$, meaning
the $i$-th columns of $\A$ and $\B$ are identical,
and therefore $\A=\B$.
Hence, $\A$ is also unique.
\end{proof}

The two conditions of
\Cref{thm:linear-trans-is-matrix-map} are roughly analogous to
the standard definition of a linear transformation between two vector
spaces, with noncommutative leftward addition replacing vector
addition in
condition~(\ref{thm:linear-trans-is-matrix-map:a}),
and the scalar $\alpha$ in condition~(\ref{thm:linear-trans-is-matrix-map:b})
ranging over all of $\Rext$,
rather than only $\R$ (albeit the condition is only required to hold for
points $\xx$ in $\Rn$ rather than all of $\extspace$).

If we instead only require that
condition~(\ref{thm:linear-trans-is-matrix-map:b})
hold for finite $\alpha$ (in $\R$, not $\Rext$),
even while considering all of $\xbar\in\eRn$,
then
\Cref{thm:linear-trans-is-matrix-map} is no longer true, in
general, as we show in the next example:

\begin{example}
With $m=n=1$, suppose
$F:\Rext\rightarrow\Rext$ is defined, for $\barx\in\Rext$, by
\[
   F(\barx) =
   \begin{cases}
                    +\infty & \text{if $\barx=-\infty$,}\\
                    \barx   & \text{if $\barx\in\R$,}\\
                    -\infty & \text{if $\barx=+\infty$.}
   \end{cases}
\]
Then it can be checked that $F$ satisfies
condition~(\ref{thm:linear-trans-is-matrix-map:a})
of \Cref{thm:linear-trans-is-matrix-map},
and also condition~(\ref{thm:linear-trans-is-matrix-map:b}), but only
for $\alpha\in\R$, though for all $\barx\in\eR$.
That is, $F(\lambda\barx) = \lambda F(\barx)$ for all
$\lambda\in\R$ and $\barx\in\eR$, but, for instance,
$F(\limray{(1)}) \neq \limray{(F(1))}$.
Furthermore, we cannot write $F$ as a linear map
$\barx\mapsto A\barx$ for any matrix (scalar, actually)
$A\in\R^{1\times 1}$.%
\indexg{linear maps, astral!characterization|)}%
\end{example}

\indexg{affine maps, astral|(}%
Using astral linear maps as building blocks, we define \emph{astral affine maps} as functions that map each
point $\xbar$ in $\extspace$ to $\bbar\plusl \A\xbar$ in
$\extspac{m}$, for some matrix $\A\in\R^{m\times n}$ and point
$\bbar\in\extspac{m}$.
\indexg{continuity!astral affine maps@of astral affine maps|(}%
As a corollary of \Cref{thm:linear:cont}(\ref{thm:linear:cont:b}), all such maps are continuous:

\begin{corollary} \label{cor:aff-cont}
  Let $\A\in\Rmn$ and $\bbar\in\extspac{m}$,
  and let $F:\extspace\rightarrow\extspac{m}$
  be defined by
  $F(\xbar)=\bbar\plusl\A\xbar$
  for $\xbar\in\extspace$.
  Then
  $F$ is continuous.
\end{corollary}

\begin{proof}
Let $\seq{\xbar_t}$ be a sequence in $\extspace$ that converges to
some point $\xbar\in\extspace$.
Then $\A\xbar_t\to\A\xbar$ by \Cref{thm:linear:cont}(\ref{thm:linear:cont:b}),
and thus $\bbar\plusl\A\xbar_t\to\bbar\plusl\A\xbar$ by
\Cref{pr:i:7}(\ref{pr:i:7f}). The continuity of $F$ now follows from
$\extspace$ being first-countable
(and \Cref{prop:first:properties}\ref{prop:first:cont}).
\end{proof}

Note that the map formed by adding a point $\bbar$ on the \emph{right},
that is, a map of the form $\xbar\mapsto\A\xbar\plusl\bbar$, need not be continuous.
For example, with $n=m=1$, the function
$F(\barx)=\barx\plusl(-\infty)$, where $\barx\in\Rext$, is not continuous:
On the sequence $x_t=t$, whose limit is $+\infty$, we see that
$F(x_t)=-\infty$ for all $t$, but $F(+\infty)=+\infty$.%
\indexg{continuity!astral affine maps@of astral affine maps|)}

Astral affine maps are precisely the maps $F:\eRn\to\extspac{m}$
that satisfy $F(\xbar\plusl\ybar)=F(\xbar)\plusl L(\ybar)$ for some astral linear map $L$.
(Astral affine maps satisfy this condition immediately from their definition, and any map that satisfies the condition can be seen to be astral affine by setting $\bbar=F(\zero)$ and noting that $L(\xbar)=\A\xbar$ for a suitable $\A$.) This generalizes an analogous characterization of affine maps from $\Rn$ to $\Rm$. Astral affine maps include continuous extensions of the standard affine maps, corresponding to $\bbar\in\R^{m}$, but also additional maps, corresponding to infinite $\bbar$.

As an immediate corollary of continuity of astral affine maps,
the affine image of any closed set is also closed:

\begin{corollary}  \label{cor:aff-img-closed-is-closed}
  Let $F:\eRn\to\eRm$ be an affine map as defined in \Cref{cor:aff-cont},
  and let $S\subseteq\extspace$.
  Then:
  \begin{letter-compact}
  \item   \label{cor:aff-img-closed-is-closed:a}
    If $S$ is closed in $\eRn$, then $F(S)$ is closed in $\eRm$.
  \item   \label{cor:aff-img-closed-is-closed:b}
    $\clbar{F(S)}=F(\Sbar)$.
  \end{letter-compact}
\end{corollary}

\begin{proof}
Since $F$ is continuous (\Cref{cor:aff-cont})
and astral space is compact and Hausdorff (\Cref{thm:i:1}\ref{thm:i:1a}), both
claims of the
\namecref{cor:aff-img-closed-is-closed}
follow
by \Cref{pr:cont-from-compact}(\ref{pr:cont-from-compact:a},\ref{pr:cont-from-compact:b}).%
\indexg{affine maps, astral|)}
\end{proof}

\section{Astral points in matrix form}
\label{sec:matrix:form}

\indexg{representations of astral points!matrix form|(}%
We have seen already, for instance in \Cref{cor:h:1}, that
leftward sums of the form
\begin{equation}
\label{eqn:pf:4}
 \limray{\vv_1} \plusl \dotsb \plusl \limray{\vv_k}
\end{equation}
emerge naturally in studying astral space.
We will next see how such leftward sums can be expressed using matrices, resulting in more compact expressions that can be analyzed and manipulated with tools from linear algebra.
\indexg{icons!omega@$\omm$|(}%
Specifically, we will
rewrite \eqref{eqn:pf:4} in a compact form as $\VV\omm$ where
$\VV=[\vv_1,\dotsc,\vv_k]$ denotes the matrix in $\R^{n\times k}$
whose $i$-th column is $\vv_i$, and where $\omm$ is a suitable
astral point in $\extspac{k}$.

To motivate our choice of $\omm$,
note first that,
by \Cref{thm:i:seq-rep},
the astral point appearing in \eqref{eqn:pf:4} is the limit of any
sequence whose elements have the form
\[
  \VV\bb_t = b_{t,1} \vv_1 + \dotsb + b_{t,k} \vv_k,
\]
where
$\seq{\bb_t}$ is a sequence in $\Rk$,
and, for $i=1,\dotsc,k$, the sequence $\seq{b_{t,i}}$ grows
to $+\infty$, with each successive sequence
growing slower than the preceding one.
Thus, to express \eqref{eqn:pf:4} as $\VV\omm$, the limit of a
sequence $\seq{\VV\bb_t}$ as above,
we intuitively want to think of $\omm$ as being like a column vector,
all of whose elements are infinite,
with the first component being the ``most infinite'' and each
following component being ``less infinite'' than the previous one.

Formally, we define $\ommk$ to be the point in $\extspac{k}$
defined by
\begin{equation}
\label{eq:i:7}
\indexm{omega k, omega}{$\ommk,\,\omm$}{icon used in matrix representation}%
\ommk = \limray{\ee_1} \plusl \dotsb \plusl \limray{\ee_k},
\end{equation}
where $\ee_i$ is the $i$-th standard basis vector
in $\R^k$.
Then by \Cref{pr:h:4}(\ref{pr:h:4c},\ref{pr:h:4g}),
\begin{align}
\notag
\VV \ommk
    =
    \VV \bigParens{ \limray{\ee_1} \plusl \dotsb \plusl \limray{\ee_k} }
    &=
    \VV(\limray{\ee_1}) \plusl \dotsb \plusl \VV(\limray{\ee_k})
\\
\notag
    &=
    \limray{(\VV\ee_1)} \plusl \dotsb \plusl \limray{(\VV\ee_k)}
\\
\label{eqn:mat-omm-defn}
    &=
    \limray{\vv_1} \plusl \dotsb \plusl \limray{\vv_k}.
\end{align}
In other words, $\VV \ommk$, with $\ommk$ defined as in
\eqref{eq:i:7}, is
equal to \eqref{eqn:pf:4}, now stated in
matrix form.
When clear from context,
we
omit $\ommk\negKern$'s subscript and write simply $\omm$.

When $k=0$, we define
$\ommsub{0}=\zerovec$, the only element of
$\extspac{0}=\R^0$.
In this case, $\VV=\zeromat{n}{0}$, the only matrix
in $\R^{n\times 0}$.
Thus, $\VV\ommsub{0}=\zeromat{n}{0} \zerovec = \zerov{n}$,
or more simply,
$\VV\omm=\zero$.
(See
\Cref{sec:prelim-zero-dim-space}.)%
\indexg{icons!omega@$\omm$|)}%

In matrix notation, \Cref{cor:h:1} states that astral space
consists exactly of all points of the form
$\VV\omm\plusl \qq$ for some matrix $\VV$ and vector $\qq$.
That is,
\begin{equation}   \label{eq:i:6}
 \extspace = \{ \VV\omm\plusl\qq :\: \VV\in\R^{n\times k},\,
                                \qq\in\Rn,\, k\geq 0 \}.
\end{equation}
Thus, every point
$\xbar\in\extspace$ can be written as
$\xbar=\VV\omm\plusl\qq$.
We refer to $\VV\omm\plusl\qq$, or more formally, the pair
$(\VV,\qq)$, as a
\emph{matrix representation}, or
\emph{matrix form}, of the point $\xbar$.%
\indexg{representations of astral points!matrix form|)}

\indexg{representations of astral points!manipulating|(}%
\indexg{representations of astral points!matrix product of|(}%
In matrix form, matrix-vector product is simple to compute:
If $\A\in\Rmn$ and $\xbar\in\extspace$, then we can write $\xbar$ in
the form
$\xbar=\VV\omm\plusl\qq$ as above,
and can then immediately compute the corresponding form of
$\A\xbar$ using
\Cref{pr:h:4}(\ref{pr:h:4c},\ref{pr:h:4d})
as
\[
  \A\xbar = (\A\VV)\omm \plusl (\A\qq).%
\indexg{representations of astral points!matrix product of|)}%
\]

We next explore in detail various properties of matrix representations
and rules for manipulating them.
We will see how these can be used, for instance, to determine if two
representations are for the same point, or to manipulate a given
representation into a more favorable form.

\indexg{representations of astral points!scalar multiple of|(}%
We begin with how such representations are affected when multiplied by
a positive scalar:

\begin{proposition}   \label{pr:mult-rep-by-scalar}
  Let
  $\alpha\in\Rextstrictpos$,
  and let
  $\xbar=\VV\omm\plusl \qq$
  where $\VV\in\R^{n\times k}$, for some $k\geq 0$, and $\qq\in\Rn$.
  Then $\alpha\xbar = \VV\omm\plusl \alpha\qq$.
\end{proposition}

\begin{proof}
Suppose first that $\alpha\in\Rstrictpos$.
Then
\begin{align}
  \alpha\ommk
  &=
  \alpha (\limray{\ee_1}\plusl\cdots\plusl\limray{\ee_k})
  \nonumber
  \\
  &=
  (\alpha \oms){\ee_1}\plusl\cdots\plusl(\alpha\oms){\ee_k}
  \nonumber
  \\
  &=
  \limray{\ee_1}\plusl\cdots\plusl\limray{\ee_k}
  =
  \ommk,
  \label{eq:pr:mult-rep-by-scalar:1}
\end{align}
where
the first and last equalities are by \eqref{eq:i:7},
the second is by
Propositions~\ref{pr:i:7}(\ref{pr:i:7b})
and
\ref{pr:scalar-prod-props}(\ref{pr:scalar-prod-props:b}),
and the third is simply because
$\alpha\oms=+\infty=\oms$.
Thus,
\[
  \alpha \xbar
  =
  \alpha(\VV\omm\plusl \qq)
  =
  \alpha(\VV \omm) \plusl \alpha\qq
  =
  \VV (\alpha\omm) \plusl \alpha\qq
  =
  \VV\omm \plusl \alpha\qq,
\]
where the second, third and fourth equalities follow respectively from
\Cref{pr:i:7}(\ref{pr:i:7b}),
\Cref{pr:h:4}(\ref{pr:h:4g}),
and
\eqref{eq:pr:mult-rep-by-scalar:1}.

For the remaining case that $\alpha=\oms$, we have
\[
  \limray{\xbar}
  =
  \lim (t \xbar)
  =
  \lim (\VV\omm \plusl t\qq)
  =
  \VV\omm\plusl \limray{\qq}
\]
with the second equality by the foregoing, and the last by
\Cref{pr:i:7}(\ref{pr:i:7f}).%
\indexg{representations of astral points!scalar multiple of|)}%
\end{proof}

Next, we give a simple generalization of \eqref{eqn:mat-omm-defn} which shows how to decompose $\VV\omm$ when $\VV$ is written as a concatenation of matrices rather than just columns:

\begin{proposition}
\label{prop:split}
Consider matrices $\VVsub{i}\in\R^{n\times k_i}$ for $i=1,\dots,\ell$, with $k_i\ge 0$,
and let $\VV=[\VVsub{1},\dotsc,\VVsub{\ell}]\in\R^{n\times k}$
be their concatenation, where $k=k_1+\dotsb+k_\ell$.
Then
\[
  \VV\omm=\VVsub{1}\omm\plusl\dotsb\plusl\VVsub{\ell}\omm.
\]
\end{proposition}

\begin{proof}
This follows by expanding both sides of the equality as
(possibly empty) leftward sums
of astrons using \eqref{eqn:mat-omm-defn}.
\end{proof}

\indexg{representations of astral points!coupling function and|(}%
\indexg{coupling function!matrix representation and|(}%
For the next few results, it will be useful to explicitly characterize $\omm\inprod\uu$:
\begin{proposition}
\label{prop:omm:inprod}
Let $\uu\in\R^k$. Then
\[
\omm\inprod\uu=
\begin{cases}
  0&\text{if $\uu=\zero$,}
\\
  +\infty&\text{if the first nonzero coordinate of $\uu$ is positive,}
\\
  -\infty&\text{if the first nonzero coordinate of $\uu$ is negative.}
\end{cases}
\]
\end{proposition}
\begin{proof}
Immediate from the definition of $\omm$ in \eqref{eq:i:7} and the Case Decomposition Lemma~\ref{lemma:case}.
\end{proof}

A key consequence of \Cref{prop:omm:inprod} is that $\omm\inprod\uu\in\set{\pm\infty}$ whenever $\uu\ne\zero$. 
We use this to show next how to express the value of $\xbar\cdot\uu$ for $\xbar=\VV\omm\plusl\qq$, depending on whether $\uu$ is orthogonal to the columns of $\VV$.
In particular, we show that $\xbar\cdot\uu\in\R$ if and only if $\uu\perp\VV$.
\begin{proposition}
\label{pr:vtransu-zero}
  Let $\xbar=\VV\omm\plusl \qq$
  where $\VV\in\R^{n\times k}$, for some $k\geq 0$, and $\qq\in\Rn$.
  For all $\uu\in\Rn$,
  if $\uu\perp\VV$ then
  $\xbar\cdot\uu=\qq\cdot\uu\in\R$;
  otherwise,
  $\xbar\cdot\uu=(\VV\omm)\cdot\uu \in \{-\infty,+\infty\}$.
\end{proposition}

\begin{proof}
If $\uu\perp\VV$, then $\trans{\VV}\uu=\zero$, and hence
\begin{align*}
  \xbar\cdot\uu
  &
  =(\VV\omm)\cdot\uu\plusl\qq\cdot\uu
  =\omm\cdot(\trans{\VV}\uu) \plusl \qq\cdot\uu
  =\omm\cdot\zero \plusl \qq\cdot\uu
  =\qq\cdot\uu,
\intertext{%
where the second equality is
by \Cref{thm:Ax-dot-u}.
Otherwise, if $\uu\not\perp\VV$, then
$\trans{\VV}\uu\ne\zero$, implying, by
\Cref{prop:omm:inprod}, that
$(\VV\omm)\cdot\uu=\omm\cdot(\trans{\VV}\uu)\in\set{-\infty,+\infty}$,
and so that
}
  \xbar\cdot\uu
  &
  =(\VV\omm)\cdot\uu\plusl\qq\cdot\uu
  =(\VV\omm)\cdot\uu.
\qedhere
\end{align*}
\end{proof}

In words, \Cref{pr:vtransu-zero} states that if $\uu$
has a nonzero inner product with any column of $\VV$,
then the value of $(\VV\omm\plusl\qq)\inprod\uu$ is already determined by $\VV\omm$.
Only when $\uu$ is orthogonal to all of $\VV$'s columns
does the evaluation proceed to the next term, which is $\qq$.%
\indexg{representations of astral points!coupling function and|)}%
\indexg{coupling function!matrix representation and|)}%

\indexg{Projection Lemma|(}%
A similar case analysis underlies the following important property of leftward sums:
\begin{lemma}[Projection Lemma]
\label{lemma:proj}
Let $\VV\in\R^{n\times k}$, $\zbar\in\eRn$, and let $\PP$ be the projection matrix onto
$(\colspace\VV)^\perp$. Then $\VV\omm\plusl\zbar=\VV\omm\plusl\PP\zbar$.
\end{lemma}

\begin{proof}
Let $\uu\in\Rn$. If $\uu\perp\VV$, then
$\uu\in(\colspace\VV)^\perp$ so $\PP\uu=\uu$
(by \Cref{pr:proj-mat-props}\ref{pr:proj-mat-props:d}).
Hence,
\begin{align*}
  \bigParens{\VV\omm\plusl\zbar}\cdot\uu
  =
  (\VV\omm)\cdot\uu\plusl\zbar\cdot\uu
  &=
  (\VV\omm)\cdot\uu\plusl\zbar\cdot(\PP\uu)
  \\
  &=
  (\VV\omm)\cdot\uu\plusl(\PP\zbar)\cdot\uu
  =
  \bigParens{\VV\omm\plusl\PP\zbar}\cdot\uu,
\intertext{%
where the third equality is by
\Cref{thm:Ax-dot-u}
(and \Cref{pr:proj-mat-props}\ref{pr:proj-mat-props:a}).
Otherwise,
if $\uu\not\perp\VV$, then
 $\trans{\VV}\uu\ne\zero$ so that, by \Cref{prop:omm:inprod},
$(\VV\omm)\cdot\uu=\omm\cdot(\trans{\VV}\uu)\in\set{-\infty,+\infty}$, implying
}
  \bigParens{\VV\omm\plusl\zbar}\cdot\uu
  =
  (\VV\omm)\cdot\uu\plusl\zbar\cdot\uu
  &=
  (\VV\omm)\cdot\uu
  \\
  &=
  (\VV\omm)\cdot\uu\plusl(\PP\zbar)\cdot\uu
  =
  \bigParens{\VV\omm\plusl\PP\zbar}\cdot\uu.
\end{align*}
Thus, for all $\uu\in\Rn$,
$(\VV\omm\plusl\zbar)\inprod\uu=(\VV\omm\plusl\PP\zbar)\inprod\uu$, so
$\VV\omm\plusl\zbar=\VV\omm\plusl\PP\zbar$
(by \Cref{pr:i:4}).
\end{proof}

The lemma says that leftward sums of the form $\VV\omm\plusl\zbar$ are
unchanged if we replace $\zbar$ by its projection onto $\comColV$,
or, for that matter, by any other astral point whose projection is the same as that of $\zbar$.
For example, if $\xbar=[\vv_1,\dotsc,\vv_k]\omm\plusl\qq$, then the
lemma implies that we can add any multiple of $\vv_1$ to the vectors $\vv_2,\dotsc,\vv_k$ or
to $\qq$ without affecting the value of $\xbar$. This follows by the Projection Lemma~\ref{lemma:proj}
with $\VV=[\vv_1]$ and $\zbar=[\vv_2,\dotsc,\vv_k]\omm\plusl\qq$. More generally, we can add any linear combination of $\vv_1,\dotsc,\vv_{i-1}$ 
to $\vv_i,\dotsc,\vv_k$ or $\qq$ without affecting the value of $\xbar$.%
\indexg{Projection Lemma|)}

\indexg{positive upper triangular matrices!matrix representation and|(}%
To study such equivalent transformations more generally, it will be
convenient to represent them using {positive upper triangular}
matrices, as defined in
\Cref{sec:prelim:lin-alg}.
The
next proposition shows that, when applied to a matrix
$\VV=[\vv_1,\dotsc,\vv_k]$,
positive upper triangular matrices
represent the operation of adding a linear combination of
$\vv_1,\dotsc,\vv_{i-1}$ to the vectors $\vv_i,\dotsc,\vv_k$,
while possibly also multiplying each $\vv_i$ by a positive scalar.
We will see soon that such operations have no effect on $\VV\omm$.

\begin{proposition}
\label{prop:pos-upper:descr}
Let $k\ge 0$ and let $\VV=[\vv_1,\dotsc,\vv_k]$ and $\VV'=[\vv'_1,\dotsc,\vv'_k]$, where $\vv_i,\vv'_i\in\R^n$ for $i=1,\dotsc,k$.
Then the following are equivalent:
\begin{letter-compact}
\item
\label{prop:pos-upper:descr:a}
  There exists a positive upper triangular matrix $\RR\in\R^{k\times k}$ such that $\VV'=\VV\RR$.
\item
\label{prop:pos-upper:descr:b}
  For $j=1,\dotsc,k$,
  each column $\vv'_{\negKern j}$ is the sum of a positive finite multiple of $\vv_j$ and a linear combination
  of the vectors $\vv_1,\dotsc,\vv_{j-1}$.
\end{letter-compact}
\end{proposition}

\begin{proof}
  ~
  
\begin{proof-parts}
\pfpart{%
  (\ref{prop:pos-upper:descr:a})
  $\Rightarrow$
  (\ref{prop:pos-upper:descr:b}):
}
Suppose $\VV'=\VV\RR$ where
$\RR\in\R^{k\times k}$ is a positive upper triangular matrix with entries $r_{ij}$
so that $r_{jj}>0$ for all $j$,
and $r_{ij}=0$ for $i>j$.
Then from the definition of matrix multiplication,
\begin{equation}   \label{eq:prop:pos-upper:descr:1}
  \vv'_{\negKern j}
  =
  \sum_{i=1}^k r_{ij}\vv_i
  =
  r_{jj} \vv_j + \sum_{i=1}^{j-1} r_{ij}\vv_i
\end{equation}
for $j=1,\dotsc,k$, implying
statement~(\ref{prop:pos-upper:descr:b}).

\pfpart{%
  (\ref{prop:pos-upper:descr:b})
  $\Rightarrow$
  (\ref{prop:pos-upper:descr:a}):
}
If statement~(\ref{prop:pos-upper:descr:b}) holds,
then there exist $r_{ij}\in\R$ for which
\eqref{eq:prop:pos-upper:descr:1}
holds with $r_{jj}>0$ for all $j$,
and $r_{ij}=0$ for $i>j$.
Setting $\RR$ to be the corresponding
positive upper triangular
matrix with entries $r_{ij}$,
this means $\VV'=\VV\RR$.
\qedhere
\end{proof-parts}
\end{proof}

\indexg{representations of astral points!positive upper triangular matrices and|(}%
We next show that the astral point $\omm$
is unchanged when multiplied by a positive upper
triangular matrix $\RR$, and when any finite vector $\bb$ is
leftwardly added to it.
\indexg{Push Lemma|(}%
This will immediately imply that an astral point with representation
$\VV\omm\plusl\qq$
can be expressed alternatively with $\VV$ replaced by $\VV\RR$ and $\qq$ replaced by $\qq+\VV\bb$,
corresponding to the operation of adding linear combinations of vectors appearing earlier in the matrix representation to those appearing later in the matrix representation.
The fact that this operation has no effect on the value of an astral point is stated formally below as the Push Lemma.

\begin{lemma}
\label{lemma:RR_omm}
  Let $k\ge 0$, let $\RR\in\R^{k\times k}$ be a positive upper
  triangular matrix, and let $\bb\in\R^k$.
  Then $\RR\omm\plusl\bb=\omm$.
\end{lemma}

\begin{proof}
First note that $\omm\plusl\bb=\Iden\omm\plusl\bb$, so by the Projection Lemma~\ref{lemma:proj}, $\Iden\omm\plusl\bb=\Iden\omm\plusl\zero=\omm$. Thus, $\omm\plusl\bb=\omm$. It remains to show that $\RR\omm=\omm$.

Denote the entries of $\RR$ as $r_{ij}$ and its columns as $\rr_j$.
To prove that $\RR\omm=\omm$, we need to show that $(\RR\omm)\cdot\uu = \omm\cdot\uu$ for all $\uu\in\R^{k}$.
As such,
let $\uu\in\R^k$, with entries~$u_j$, and let $\uu'=\trans{\RR}\uu$,
with entries $u'_j$.
Then by matrix algebra, and since $\RR$ is upper triangular,
\begin{equation}
\label{eq:RR_omm:2}
  u'_j=\sum_{i=1}^j  u_i r_{ij}.
\end{equation}

We claim that $\omm\cdot\uu' = \omm\cdot\uu$.
If $\uu=\zero$ then also $\uu'=\zero$, implying the claim.
Otherwise, let $\ell$ be the smallest index for which $u_{\ell}\neq 0$,
so that $u_j=0$ for $j<\ell$. By \eqref{eq:RR_omm:2},
$u'_j=0$ for $j<\ell$ and $u'_{\ell}=u_{\ell}r_{\ell\ell}$.
Thus, $\ell$ is also the smallest index for which $u'_{\ell}\neq 0$, and $u'_{\ell}$ has the same sign as $u_{\ell}$.
Therefore, by \Cref{prop:omm:inprod},
$\omm\cdot\uu = \omm\cdot\uu'$.

Thus, in all cases,
\[
  \omm\cdot\uu
  =
  \omm\cdot\uu'
  =
  \omm\cdot(\trans{\RR}\uu)
  =
  (\RR\omm)\cdot\uu,
\]
with the last equality from \Cref{thm:Ax-dot-u}.
Since this holds for all $\uu\in\R^k$,
we conclude that $\omm=\RR\omm$
(by \Cref{pr:i:4}).
\end{proof}

\begin{lemma}[Push Lemma]
\label{pr:g:1}
  Let $\VV\in\R^{n\times k}$ and $\qq\in\Rn$. Further, let
  $\RR\in\R^{k\times k}$ be a positive
  upper triangular matrix and $\bb\in\R^k$. Then
  $\VV\omm\plusl\qq=(\VV\RR)\omm\plusl(\qq+\VV\bb)$.
\end{lemma}

\begin{proof}
We have
\[
  (\VV\RR)\omm\plusl(\qq+\VV\bb)
=
  \VV(\RR\omm)\plusl\VV\bb\plusl\qq
=
  \VV(\RR\omm\plusl\bb)\plusl\qq
=
  \VV\omm\plusl\qq,
\]
with the last equality from \Cref{lemma:RR_omm}.%
\indexg{Push Lemma|)}%
\indexg{representations of astral points!positive upper triangular matrices and|)}%
\indexg{positive upper triangular matrices!matrix representation and|)}
\end{proof}

\indexg{representations of astral points!full-rank matrix@with full-rank matrix|(}%
The next proposition, following directly from
the Projection Lemma~\ref{lemma:proj}, shows that any point $\VV\omm$ can be
rewritten using a linearly independent subset of the columns of $\VV$.
(An empty set of columns is vacuously linearly independent.)

\begin{proposition}
\label{prop:V:indep}
Let
$\xbar= [\vv_1,\dotsc,\vv_k]\omm$
with $k\ge 0$. If $\vv_i\in\spn\set{\vv_1,\dotsc,\vv_{i-1}}$,
for some $i\ge 1$,
then $\xbar=[\vv_1,\dotsc,\vv_{i-1},\vv_{i+1},\dotsc,\vv_k]\omm$. Therefore,
there exists $\ell\ge 0$ and a subset of indices
$ 1\leq i_1 < i_2 < \dotsb < i_{\ell} \leq k $
such that
$\xbar=[\vv_{i_1},\dotsc,\vv_{i_{\ell}}]\omm$
and
$\vv_{i_1},\dotsc,\vv_{i_{\ell}}$
are linearly independent.
\end{proposition}

\begin{proof}
Suppose $\vv_i\in\spnfin{\vv_1,\dotsc,\vv_{i-1}}$, and let
$\VV=[\vv_1,\dotsc,\vv_{i-1}]$ and $\VV'=[\vv_{i+1},\dotsc,\vv_k]$.
Let $\PP$ be the projection matrix onto
$(\colspace\VV)^\perp$, the orthogonal complement of the
column space of $\VV$.
We then have
\begin{align*}
  \xbar
  =[\VV,\vv_i,\VV']\omm
&
  =\VV\omm\plusl\vv_i\omm\plusl\VV'\omm
\\
&
  =\VV\omm\plusl(\PP\vv_i)\omm\plusl\VV'\omm
  \\
  &
  =\VV\omm\plusl\zero\plusl\VV'\omm
  =[\VV,\VV']\omm
  =
  [\vv_1,\dotsc,\vv_{i-1},\vv_{i+1},\dotsc,\vv_k]\omm.
\end{align*}
The second and fifth equalities are by \Cref{prop:split}.
The third is by the Projection Lemma~\ref{lemma:proj}.
The fourth is because $\PP\vv_i=\zero$
(by \Cref{pr:proj-mat-props}\ref{pr:proj-mat-props:e})
since $\vv_i$ is in $\colspace{\VV}$.\looseness=-1

By repeatedly removing linearly dependent vectors in this way, we eventually must end up
with a linearly independent subset, as claimed in the proposition.%
\indexg{representations of astral points!full-rank matrix@with full-rank matrix|)}%
\indexg{representations of astral points!manipulating|)}%
\end{proof}

\begin{definition}
\indexg{rank, astral|(}%
The \emph{astral rank} of a point $\xbar\in\extspace$
is
the smallest number of columns $k$ for which there exist a matrix
$\VV\in\R^{n\times k}$
and $\qq\in\Rn$ with $\xbar=\VV\omm\plusl\qq$.
\end{definition}

Astral rank captures the
dimensionality of the set of directions in which a sequence with limit
$\xbar$ is going to infinity.
The set of points with astral rank $0$ is exactly~$\Rn$.
Points with astral rank $1$ take the form $\limray{\vv}\plusl\qq$ with $\vv\ne\zero$
and arise as limits of sequences that go to infinity along a
halfline.
Points with higher astral ranks arise as limits of sequences that go to infinity in one dominant direction, but also in some secondary direction, and possibly a
tertiary direction, and so on, as discussed in
\Cref{subsec:astral-pt-form}.

\begin{theorem}  \label{thm:ast-rank-is-mat-rank} %
\indexg{rank (of matrix)!astral rank and|(}%
  Let $\xbar=\VV\omm\plusl\qq$ for some
  $\VV\in\R^{n\times k}$ and $\qq\in\Rn$.
  Then the astral rank of $\xbar$ is equal to the
  matrix rank of $\VV$.
\end{theorem}

\begin{proof}
Let $r$ be the astral rank of $\xbar$, which means that there exist a matrix $\VV'\in\R^{n\times r}$ and a vector $\qq'\in\R^n$ such that $\xbar=\VV'\omm\plusl\qq'$. The columns of $\VV'$ must be linearly independent,
since otherwise we could remove some of them (by \Cref{prop:V:indep}),
which would contradict the fact that $r$ is the astral rank of
$\xbar$. Since $\VV'$ consists of $r$ linearly independent columns,
its matrix rank is $r$. To prove the proposition, it suffices to show that
the matrix rank of $\VV$ is equal to the matrix rank of~$\VV'$.

Let $L$ be the set of vectors $\uu\in\Rn$ such that $\xbar\inprod\uu\in\R$. Since $\xbar=\VV\omm\plusl\qq=\VV'\omm\plusl\qq'$, we obtain by \Cref{pr:vtransu-zero}
that $L=(\colspace\VV)^\perp=(\colspace\VV')^\perp$.
Thus, $L^\perp=\colspace\VV=\colspace\VV'$
(by \Cref{pr:std-perp-props}\ref{pr:std-perp-props:c}).
Therefore, the matrix rank of $\VV$ is equal to the matrix rank of $\VV'$.%
\indexg{rank (of matrix)!astral rank and|)}
\end{proof}

As a corollary, we can show that astral rank is subadditive in the
following sense:

\begin{corollary}
\label{cor:rank:subadd}
Let $\xbar,\xbar_1,\xbar_2\in\eRn$ be such that $\xbar=\xbar_1\plusl\xbar_2$, and let
$r,r_1,r_2$ be their respective astral ranks.
Then
\[
   \max\{r_1, r_2\} \leq r\le r_1+r_2.
\]
\end{corollary}

\begin{proof}
Let
$\xbar_1=\VVsub{1}\omm\plusl\qq_1$ and
$\xbar_2=\VVsub{2}\omm\plusl\qq_2$. By Propositions~\ref{pr:i:7}(\ref{pr:i:7d}) and~\ref{prop:split},
\[
  \xbar
  =
  \xbar_1\plusl\xbar_2
  =
  (\VVsub{1}\omm\plusl\qq_1)
  \plusl
  (\VVsub{2}\omm\plusl\qq_2)
  =
  [\VVsub{1},\VVsub{2}]\omm\plusl(\qq_1+\qq_2).
\]
The column space of $[\VVsub{1},\VVsub{2}]$ is equal to the set sum
$(\colspace{\VVsub{1}})+(\colspace{\VVsub{2}})$.
Therefore,
\[
  r
  =
  \dim(\colspace[\VVsub{1},\VVsub{2}])
  \leq
  \dim(\colspace{\VVsub{1}})
  +
  \dim(\colspace{\VVsub{2}})
  =
  r_1 + r_2,
\]
where
the inequality is by \Cref{pr:subspace-intersect-dim},
and the two equalities are both by
\Cref{thm:ast-rank-is-mat-rank}
(and definition of matrix rank as the dimension of its column
space).

Similarly, the linear subspace $\colspace[\VVsub{1},\VVsub{2}]$ includes
$\colspace{\VVsub{i}}$, for $i=1,2$, so
$r_i=\dim(\colspace{\VVsub{i}})\leq \dim(\colspace[\VVsub{1},\VVsub{2}])=r$.%
\indexg{rank, astral|)}
\end{proof}

\section{Canonical representation}
\label{sec:canonical:rep}

\indexg{canonical representation|(}%
Although the same astral point may have multiple representations, it
is possible to obtain a unique representation of a particularly
natural form:

\begin{definition}
  For a matrix $\VV\in\Rnk$ and $\qq\in\Rn$, we say that
  $\VV\omm\plusl\qq$ is a \emph{canonical representation}
  if $\VV$ is column-orthogonal (so that the columns of $\VV$ are
  orthonormal), and if $\qq\perp\VV$.
\end{definition}

\Cref{cor:h:1} showed that every astral point can be represented in
this canonical form.
We show next that the canonical representation of every point is in
fact unique.

\begin{theorem}  \label{pr:uniq-canon-rep}
  Every point in $\extspace$ has a unique canonical representation.
\end{theorem}

Before proving the theorem, we give a lemma that will be one of the
main steps in proving uniqueness.
We state the lemma in a form that is a bit more general than needed
at this point since it will get used again shortly. To state the lemma,
we need some additional terminology. We say that a matrix $\VV'\in\R^{n\times k'}$
is a \emph{prefix} of a matrix $\VV\in\R^{n\times k}$ if $k'\le k$ and
if the first $k'$ columns of $\VV$ are identical to $\VV'$, that is,
if $\VV=[\VV',\VV'']$ for some matrix $\VV''$.

\begin{lemma}  \label{lem:uniq-ortho-prefix}
  Let $\VV\in\R^{n\times k}$ and $\VV'\in\R^{n\times k'}$ be column-orthogonal matrices,
  and let $\zbar,\,\zbar'\in\extspace$.
  If $\VV\omm\plusl\zbar=\VV'\omm\plusl\zbar'$
  then one of the matrices must be a prefix of the other.
  Consequently:
  \begin{letter-compact}
  \item  \label{lem:uniq-ortho-prefix:a}
    If either $k\geq k'$ or $\zbar\in\Rn$, then
    $\VV'$ is a prefix of $\VV$.
  \item  \label{lem:uniq-ortho-prefix:b}
    If either $k=k'$ or $\zbar,\zbar'\in\Rn$, then
    $\VV=\VV'$.
  \end{letter-compact}
\end{lemma}

\begin{proof}
Let
$\xbar = \VV\omm\plusl\zbar$
and
$\xbar' = \VV'\omm\plusl\zbar'$.
To prove the lemma, we assume that neither of the matrices $\VV$ and $\VV'$ is a prefix of the other, and show this implies $\xbar\neq\xbar'$.

Since neither is a prefix of the other, we can write $\VV=[\VVsub{0},\vv,\VVsub{1}]$ and
$\VV'=[\VVsub{0},\vv',\VVsub{1}']$ where $\VVsub{0}$ is the initial (possibly empty) set of columns where $\VV$ and $\VV'$ agree, and $\vv,\vv'\in\Rn$
is the first column where $\VV$ and $\VV'$ disagree
(so $\vv\ne\vv'$).

Let $\uu=\vv-\vv'$. Note that $\vv,\vv'\perp\VVsub{0}$ since the matrices $\VV$ and
$\VV'$ are column-orthogonal, and therefore also $\uu\perp\VVsub{0}$.
We will show that $\xbar\inprod\uu\ne\xbar'\inprod\uu$.
We have
\begin{align*}
  \xbar\inprod\uu
  &=
  (\VV\omm\plusl\zbar)\inprod\uu
  =
  ([\VVsub{0},\vv,\VVsub{1}]\omm\plusl\zbar)\inprod\uu
\\
  &=
  (\VVsub{0}\omm)\inprod\uu
  \plusl
  \limray{\vv}\inprod\uu
  \plusl
  (\VVsub{1}\omm)\inprod\uu
  \plusl
  \zbar\inprod\uu
\\
  &=
  0\plusl(+\infty)
  \plusl
  (\VVsub{1}\omm)\inprod\uu
  \plusl
  \zbar\inprod\uu
  =
  +\infty.
\end{align*}
The third equality is by \Cref{prop:split}
(and since $\vv\omm=\limray{\vv}$).
In the fourth equality, we first used $\uu\perp\VVsub{0}$,
implying $(\VVsub{0}\omm)\inprod\uu=0$
(by \Cref{pr:vtransu-zero}).
We then
used that $\vv\inprod\uu=\vv\inprod(\vv-\vv')=1-\vv\inprod\vv'>0$,
since $\norm{\vv}=\norm{\vv'}=1$ and $\vv\ne\vv'$,
thereby implying that $\limray{\vv}\inprod\uu=+\infty$.

Similarly,
\begin{align*}
  \xbar'\inprod\uu
  &=
  (\VVsub{0}\omm)\inprod\uu
  \plusl
  \limray{\vv'}\inprod\uu
  \plusl
  (\VV'_1\omm)\inprod\uu
  \plusl
  \zbar'\inprod\uu
\\
  &=
  0\plusl(-\infty)
  \plusl
  (\VV'_1\omm)\inprod\uu
  \plusl
  \zbar'\inprod\uu
  =
  -\infty.
\end{align*}
This time, in the second equality,
we used that $\vv'\inprod\uu=\vv'\inprod(\vv-\vv')=\vv'\inprod\vv-1<0$,
implying
$\limray{\vv'}\inprod\uu=-\infty$.
Thus, $\xbar\inprod\uu\ne\xbar'\inprod\uu$;
hence, $\xbar\neq\xbar'$.

Therefore, one of the matrices $\VV$ and $\VV'$ must be a prefix of
the other.

\begin{proof-parts}
\pfpart{Consequence~(\ref{lem:uniq-ortho-prefix:a}):}
If $k\geq k'$ then the only possibility is that $\VV'$ is a prefix of
$\VV$.
If $\zbar=\zz\in\Rn$, then, by
\Cref{thm:ast-rank-is-mat-rank},
the astral rank of $\VV\omm\plusl\zz$ is $k$, while
the astral rank of $\VV'\omm\plusl\zbar'$ is at least $k'$
by \Cref{cor:rank:subadd}.
Being equal, these two points must have the same astral rank, so
$k\geq k'$, implying, as just noted, that $\VV'$ is a prefix of $\VV$.

\pfpart{Consequence~(\ref{lem:uniq-ortho-prefix:b}):}
If $k=k'$ or $\zbar,\zbar'\in\Rn$, then applying
(\ref{lem:uniq-ortho-prefix:a}) in both directions yields
that $\VV$ and $\VV'$ are prefixes of each other, and so are equal.
\qedhere
\end{proof-parts}
\end{proof}

\begin{proof}[Proof of \Cref{pr:uniq-canon-rep}]
Existence was proved in \Cref{cor:h:1}. To prove uniqueness,
consider canonical representations of two points,
$\xbar = \VV\omm\plusl\qq$
and
$\xbar' = \VV'\omm\plusl\qq'$,
and assume
that
$\xbar=\xbar'$. We show that this implies
that the two representations are identical.

To start,
\Cref{lem:uniq-ortho-prefix}(\ref{lem:uniq-ortho-prefix:b})
implies that $\VV=\VV'$.
Next, note that $\qq$ and $\qq'$ are both orthogonal to all the columns
of $\VV=\VV'$, so \Cref{pr:vtransu-zero} yields
\begin{align*}
  \xbar\cdot(\qq-\qq')  &= \qq \cdot(\qq - \qq') \\
  \xbar'\cdot(\qq-\qq')  &= \qq' \cdot(\qq - \qq').
\end{align*}
Since $\xbar=\xbar'$, these quantities must be equal.
Taking their difference, we obtain
$\norm{\qq-\qq'}^2=0$, and therefore $\qq=\qq'$.
Thus, the representations are identical.%
\indexg{canonical representation|)}
\end{proof}

\indexg{representations of astral points!equivalence of|(}%
Building on the uniqueness of the canonical representation, we next
characterize fully when two representations
$\VV\omm\plusl \qq$ and $\VV'\omm\plusl \qq'$ represent the same
astral point, assuming that $\VV$ and $\VV'$ both have full column
rank.
This characterization further implies that any such representation
for a given point can be obtained
from any other such representation for that same point
by means of a push operation
(including from its canonical representation).

\begin{theorem}  \label{thm:i:2}
  Let
  $\xbar=\VV\omm\plusl \qq$
  and
  $\xbar'=\VV'\omm\plusl \qq'$
  be two points in $\extspace$, for some $\VV,\VV'\in\R^{n\times k}$,
  $k\geq 0$, and $\qq,\qq'\in\Rn$.
  Assume $\VV$ and $\VV'$ both have full column rank.
  Then $\xbar=\xbar'$ if and only if
  $\VV'=\VV\Rtilde$ and $\qq'=\qq+\VV\tbb$
  for some positive upper triangular matrix
  $\Rtilde\in\R^{k\times k}$ and some $\tbb\in\R^{k}$.
\end{theorem}

\begin{proof}
That the given conditions imply $\xbar=\xbar'$ is proved by
the Push Lemma~\ref{pr:g:1}.

To prove the converse, assume $\xbar=\xbar'$. By \Cref{prop:QR}, we can write $\VV=\WW\RR$
for some column-orthogonal $\WW\in\Rnk$ and positive upper triangular
$\RR\in\R^{k\times k}$,
and similarly decompose $\VV'$ as $\VV'=\WW'\RR'$. Thus, 
\begin{align}
  \xbar&=\VV\omm\plusl\qq=\WW(\RR\omm)\plusl\qq=\WW\omm\plusl\qq,
  \nonumber
  \\
  \xbar'&=\VV'\omm\plusl\qq'=\WW'(\RR'\omm)\plusl\qq'=\WW'\omm\plusl\qq',
  \label{eq:thm:i:2:e1}
\end{align}
where the last equality on each line is
by \Cref{lemma:RR_omm}.
Since $\xbar=\xbar'$,
\Cref{lem:uniq-ortho-prefix}(\ref{lem:uniq-ortho-prefix:b})
implies that $\WW=\WW'$. Thus,
\[
  \VV'=\WW\RR'=(\VV\Rinv)\RR'=\VV(\Rinv\RR')
  =\VV\Rtilde,
\]
where we let $\Rtilde=\Rinv\RR'$, which
is positive upper triangular, as follows, along with $\RR$'s
invertibility,
from \Cref{prop:pos-upper}.

Next, let $\PP$ be the projection matrix onto
$(\colspace\WW)^\perp$,
the orthogonal complement of the column space of $\WW$.
By \eqref{eq:thm:i:2:e1}
and the Projection Lemma~\ref{lemma:proj},
\begin{align*}
  \xbar&=\WW\omm\plusl\qq=\WW\omm\plusl\PP\qq,
\\
  \xbar'&=\WW\omm\plusl\qq'=\WW\omm\plusl\PP\qq'.
\end{align*}
Note that both of the representations on the right are canonical since
$\WW$ is column-orthogonal, and since $\PP\qq$ and $\PP\qq'$ are both orthogonal
to the columns of $\WW$.
Since the canonical representation of $\xbar=\xbar'$ 
is unique
(\Cref{pr:uniq-canon-rep}),
this implies $\PP\qq=\PP\qq'$.

By \Cref{pr:lin-decomp-rel-vecs}, we can write
$\qq=\PP\qq+\WW\cc$ and $\qq'=\PP\qq'+\WW\cc'$
for suitable $\cc,\cc'\in\R^k$.
Then
\[
  \qq
  =
  \PP\qq + \WW\cc
  =
  (\PP\qq' + \WW\cc') + \WW(\cc-\cc')
  =
  \qq' + (\VV\Rinv)(\cc-\cc')
  =
  \qq' + \VV\tbb,
\]
where the second equality is because $\PP\qq=\PP\qq'$,
and where we define $\tbb=\Rinv(\cc-\cc')$.
This completes the proof.%
\indexg{representations of astral points!equivalence of|)}
\end{proof}

\chapter{Representation and sequences}
\label{sec:rep-seq}

Astral space was constructed from sequences, with each astral point
the common limit of a set of sequences in $\Rn$.
From a different perspective,
we saw in the last chapter how all astral points
can be represented using matrices.
In this chapter, we bring these two views of astral space closer by
characterizing
the sequences in $\Rn$ that converge to a
specific astral point in terms of that point's matrix representation.
Earlier, in \Cref{thm:i:seq-rep},
which generalized
the polynomially graded sequences from
Example~\ref{ex:poly-speed-intro},
we gave sufficient conditions for a sequence to converge to a point
with a particular representation.
Here, we prove the sense in which those conditions are necessary as
well, thereby providing a complete characterization of sequences
converging to a point.

We also introduce some useful notions for describing certain
properties of sequences.
We then apply these to further characterize the sequences
converging to points
that can be represented
by specific kinds of expressions,
such as $\A\xbar$.
The tools developed here will allow us to easily move back and forth
between astral points and the sequences that converge to them.

\section{Characterizing all sequences converging to an astral point}
\label{sec:seqs-to-matrix-rep}

\indexg{representations of astral points!sequences converging to|(}%
\indexg{sequence convergence!matrix representation@to matrix representation|(}%
We first study sequences converging to a point
$\xbar\in\eRn$ with representation $\xbar=\VV\omm\plusl\qq$
where $\VV\in\R^{n\times k}$ has linearly independent columns and
$\qq\in\Rn$ is orthogonal to all the columns of $\VV$.
Let $\seq{\xx_t}$ be any sequence in $\Rn$.
Then each point $\xx_t$ can be uniquely
represented in the form $\xx_t=\VV \bb_t + \qq_t$ for some $\bb_t\in\Rk$ and
some $\qq_t\in\Rn$ with $\qq_t\perp\VV$
(\Cref{pr:lin-decomp-rel-vecs}).
The next theorem provides necessary and sufficient conditions for
when the sequence $\seq{\xx_t}$ will converge to $\xbar$,
in terms of this representation.
As usual, $b_{t,i}$ denotes the $i$-th component of $\bb_t$.

\begin{theorem}  \label{thm:seq-rep}
  Let $\VV\in\R^{n\times k}$ be a matrix with $k\geq 0$ columns that are
  linearly independent.
  Let $\xbar=\VV\omm\plusl\qq$ where $\qq\in\Rn$ and
  $\qq\perp \VV$.
  Let $\seq{\xx_t}$ be a sequence in~$\Rn$ with
  $\xx_t=\VV \bb_t + \qq_t$ for some $\bb_t\in\Rk$ and $\qq_t\in\Rn$
  with $\qq_t \perp \VV$, for all $t$.
  Then $\xx_t\rightarrow\xbar$ if and only if all of the following
  hold:
  \begin{letter-compact}
  \item  \label{thm:seq-rep:a}
    $b_{t,i}\rightarrow+\infty$, for $i=1,\dotsc,k$.
  \item  \label{thm:seq-rep:b}
    $b_{t,i+1}/b_{t,i}\rightarrow 0$, for $i=1,\dotsc,k-1$.
  \item  \label{thm:seq-rep:c}
    $\qq_t\rightarrow\qq$.
  \end{letter-compact}
\end{theorem}

\begin{proof}
\Cref{thm:i:seq-rep}, applied to $\qq$ and the columns of $\VV$,
shows that if the sequence
$\seq{\xx_t}$ satisfies conditions
(\ref{thm:seq-rep:a}),
(\ref{thm:seq-rep:b}), and
(\ref{thm:seq-rep:c}),
then $\xx_t\rightarrow\xbar$. It remains to prove the converse.

Assume that $\xx_t\rightarrow\xbar$.
Let $\VVdag$ be $\VV$'s pseudoinverse
(as discussed in \Cref{sec:prelim:orth-proj-pseud}),
and let $\ZZ=\trans{(\VVdag)}$.
Since $\VV$ has full column rank,
$\trans{\ZZ}\VV=\VVdag \VV = \Iden$
(where $\Iden$ is the $k\times k$ identity matrix)
by \Cref{pr:pseudoinv-props}(\ref{pr:pseudoinv-props:c}).

Also, the column space of $\ZZ$ is the same as that
of $\VV$
by \Cref{pr:pseudoinv-props}(\ref{pr:pseudoinv-props:a}).
In particular, this means that every column of $\ZZ$ is a linear
combination of the columns of~$\VV$.
Since $\qq$ is orthogonal to all columns of $\VV$, it follows that
$\qq$ is also orthogonal to all columns of $\ZZ$, that is,
$\qq\perp\ZZ$.
Likewise, $\qq_t\perp\ZZ$ for all $t$.

Thus, we have
\[
  \trans{\ZZ} \xx_t
  =
  \trans{\ZZ} (\VV \bb_t + \qq_t)
  =
  (\trans{\ZZ}\VV) \bb_t + \trans{\ZZ}\qq_t
  =\bb_t
\]
where the last equality follows because $\trans{\ZZ}\VV=\Iden$
and $\qq_t\perp\ZZ$.
Similarly,
\[
  \trans{\ZZ} \xbar
  =
  \trans{\ZZ} (\VV \omm \plusl \qq)
  =
  (\trans{\ZZ}\VV) \omm \plusl \trans{\ZZ}\qq
  =\omm,
\]
with the second equality by
\Cref{pr:h:4}(\ref{pr:h:4c},\ref{pr:h:4d}).
Thus,
\begin{equation}
\label{eqn:thm:seq-rep:1}
  \bb_t
  =
  \trans{\ZZ}\xx_t
  \rightarrow
  \trans{\ZZ}\xbar
  =
  \omm
\end{equation}
by \Cref{thm:linear:cont}(\ref{thm:linear:cont:b}), since $\xx_t\rightarrow\xbar$.

For $i=1,\dotsc,k$,
since $\bb_t\to\omm$,
and with $\ee_1,\dotsc,\ee_k$ the standard basis vectors in~$\R^k$,
we must also have
\[b_{t,i} = \bb_t \cdot \ee_i \rightarrow \omm \cdot \ee_i = +\infty,\]
with convergence from \Cref{thm:i:1}(\ref{thm:i:1c}) and
the last equality from \Cref{prop:omm:inprod}
(or \Cref{lemma:case}).
This proves part~(\ref{thm:seq-rep:a}).

Next, let $\epsilon\in\Rstrictpos$.
Then, by similar reasoning, for $i=1,\dotsc,k-1$,
\[
  \epsilon b_{t,i} - b_{t,i+1}
  =
  \bb_t \cdot (\epsilon \ee_i - \ee_{i+1})
  \rightarrow
  \omm \cdot (\epsilon \ee_i - \ee_{i+1})
  =
  +\infty,
\]
where the last step again follows
from \Cref{prop:omm:inprod}.
Thus, for all $t$ sufficiently large,
$\epsilon b_{t,i} - b_{t,i+1} > 0$,
and also $b_{t,i}>0$ and $b_{t,i+1}>0$ by
part~(\ref{thm:seq-rep:a}).
That is,
$0 < b_{t,i+1}/b_{t,i} < \epsilon$.
Since this holds for all $\epsilon\in\Rstrictpos$,
this proves
part~(\ref{thm:seq-rep:b}).

Finally,
let $\PP\in\R^{n\times n}$ be the projection matrix onto
$(\colspace\VV)^\perp$.
Then
\[
  \PP\xx_t
  =
  \PP (\VV \bb_t + \qq_t)
  =
  (\PP\VV) \bb_t + \PP\qq_t
  =\qq_t,
\]
where the last equality follows because,
by
\Cref{pr:proj-mat-props}(\ref{pr:proj-mat-props:d},\ref{pr:proj-mat-props:e}),
$\PP\VV=\zeromat{n}{k}$
and $\PP\qq_t=\qq_t$ since $\qq_t\perp\VV$. Similarly,
$\PP\xbar=\qq$, and therefore,
$\qq_t=\PP\xx_t\rightarrow\PP\xbar=\qq$, proving
part~(\ref{thm:seq-rep:c}).
\end{proof}

Recall from \Cref{sec:suf-cond-seq-conv}
that a sequence
$\seq{\bb_t}$ in $\Rk$ that satisfies conditions~(\ref{thm:seq-rep:a})
and~(\ref{thm:seq-rep:b}) of \Cref{thm:seq-rep} is said to have
{entries converging to $+\infty$ at decreasing rates}.
The next corollary shows
that this property characterizes all sequences converging to $\omm$:

\begin{corollary}
\label{cor:omm:seq}
Let $\seq{\bb_t}$ be a sequence in $\R^k$. Then $\bb_t\to\omm$
if and only if the entries of $\bb_t$ converge to $+\infty$ at decreasing rates.
\end{corollary}

\begin{proof}
This follows directly from \Cref{thm:seq-rep} with $n=k$,
$\VV=\Iden$, and $\qq=\zero$.
\end{proof}

\Cref{thm:seq-rep} characterizes convergence of a sequence to a
point $\xbar=\VV\omm\plusl\qq$ when the columns of $\VV$ are linearly
independent and when $\qq\perp\VV$.
If those columns are \emph{not} linearly independent, then a sequence element
$\xx_t\in\Rn$ can still be expressed as
$\xx_t=\VV \bb_t + \qq_t$ for some $\bb_t\in\Rk$ and $\qq_t\in\Rn$,
but this way of writing $\xx_t$ is no longer {unique}, even with
the additional requirement that
$\qq_t\perp\VV$.
Nevertheless, as the next theorem shows, there always exists a
choice of $\bb_t$ and $\qq_t$ so that the conditions of
\Cref{thm:seq-rep} hold.
Combined with
\Cref{thm:i:seq-rep},
this then provides both a necessary and sufficient condition for a
sequence to converge to a point with a representation as above, even
when $\VV$'s columns are not linearly independent.
Note that we also no longer require $\qq\perp \VV$, nor that
$\qq_t\perp \VV$.

\begin{theorem}  \label{thm:seq-rep-not-lin-ind}
  Let $\VV\in\R^{n\times k}$ with $k\geq 0$, let $\qq\in\Rn$,
  and
  let $\xbar=\VV\omm\plusl\qq$.
  Let $\seq{\xx_t}$ be a sequence in $\Rn$.
  Then $\xx_t\rightarrow\xbar$ if and only if
  there exist sequences $\seq{\bb_t}$ in $\Rk$ and
  $\seq{\qq_t}$ in $\Rn$ such that
  $\xx_t=\VV\bb_t + \qq_t$ for all $t$,
  and such that conditions
  (\ref{thm:seq-rep:a}),
  (\ref{thm:seq-rep:b}), and
  (\ref{thm:seq-rep:c})
  of \Cref{thm:seq-rep} are satisfied.

  Furthermore, this equivalence still holds if we additionally
  require that ${\bb_t\in\Rstrictpos^k}$ for all $t$.
\end{theorem}

To begin, in the following lemma,
we prove one of the theorem's implications
for the special case that the nonzero columns of
$\VV$ are linearly independent and that $\qq\perp\VV$.
We then show how the general case can be reduced to this special case.

\begin{lemma}   \label{lem:seq-rep-nz-lin-ind}
  Let $\VV$, $k$, $\qq$, $\xbar$, and $\seq{\xx_t}$
  be as in
  \Cref{thm:seq-rep-not-lin-ind}.
  Assume the nonzero columns of $\VV$ are linearly independent
  and that $\qq\perp\VV$,
  and suppose further that $\xx_t\rightarrow\xbar$.
  Then there exist sequences $\seq{\bb_t}$ in $\Rk$ and
  $\seq{\qq_t}$ in $\Rn$
  such that $\xx_t=\VV\bb_t + \qq_t$ for all $t$,
  and such that conditions
  (\ref{thm:seq-rep:a}),
  (\ref{thm:seq-rep:b}), and
  (\ref{thm:seq-rep:c})
  of \Cref{thm:seq-rep}
  are satisfied.
\end{lemma}

\begin{proof}
Let $\vv_1,\dotsc,\vv_k\in\Rn$ be the columns of $\VV$ so that
$\VV=[\vv_1,\dotsc,\vv_k]$.
Let $\ci{1},\dotsc,\ci{r}$ be the indices of all of $\VV$'s
nonzero columns, with
$1\leq\ci{1}<\cdots<\ci{r}\leq k$.

Let $\VV'=[\vv_{\ci{1}},\dotsc,\vv_{\ci{r}}]$.
Then $\qq\perp\VV'$ (since $\qq\perp\VV$)
and $\xbar=\VV'\omm\plusl\qq$ since $\VV\omm$ is unchanged when
omitting columns that are $\zero$.
Further, by assumption, the columns of~$\VV'$ are linearly
independent.
Thus, by \Cref{thm:seq-rep}, there exist sequences
$\seq{\bb'_t}$ in~$\R^r$ and $\seq{\qq_t}$ in $\R^n$
such that: $\xx_t=\VV'\bb'_t+\qq_t$ for all $t$;
$b'_{t,j}\rightarrow+\infty$ for $j=1,\dotsc,r$;
$b'_{t,j} \seqgt b'_{t,j+1}$
for $j=1,\dotsc,r-1$;
and
$\qq_t\rightarrow\qq$
(satisfying condition~(\ref{thm:seq-rep:c})
of \Cref{thm:seq-rep}).

To complete the proof, we need to construct a sequence $\seq{\bb_t}$
in $\Rk$ whose entries converge to $+\infty$ at
decreasing rates.
For entries $i$ corresponding to nonzero columns of $\VV$, so that
$i=\ci{j}$ for some $j$,
we can simply set $b_{t,i}=b'_{t,j}$.
For all other entries, we can then ``pad'' $\bb_t$ with any
sequences satisfying the required asymptotic-dominance relationships.

In more detail, we construct the sequences $\seq{b_{t,i}}$ in steps.
We let $\coldefset$ denote the set of indices $i$ for which
$\seq{b_{t,i}}$ has already been defined.
Naturally, this set will 
evolve through the course of the
construction; our aim is for it eventually to include all
of $\colset$, where $\colset=\{1,\dotsc,k\}$ is the set of all column
indices.

Throughout the construction, we maintain the following invariants:
\begin{item-compact}
  \item
    For all $i\in\coldefset$, $b_{t,i}\rightarrow+\infty$.
  \item
    For all $i,j\in\coldefset$,
    if $i<j$ then $b_{t,i}\seqgt b_{t,j}$.
\end{item-compact}
Once $\colset\subseteq\coldefset$, these will imply that $\seq{\bb_t}$
satisfies conditions~(\ref{thm:seq-rep:a})
and~(\ref{thm:seq-rep:b})
of \Cref{thm:seq-rep}.

To begin the construction, for $j=1,\dotsc,r$, we define
$b_{t,\ci{j}}=b'_{t,j}$ for all $t$.
Correspondingly, we let $\coldefset=\{\ci{1},\dotsc,\ci{r}\}$.
Since $\seq{\bb'_t}$ has entries converging to $+\infty$ at decreasing
rates, the invariants given above are satisfied.

Next, for convenience, we define two ``extra'' sequences
with indices $0$ and $k+1$:
If $r>0$, we define $\seq{b_{t,0}}$ and $\seq{b_{t,k+1}}$ 
to be any sequences converging to $+\infty$ for which
$b_{t,0}\seqgt b'_{t,1}$
and
$b'_{t,r}\seqgt b_{t,k+1}$;
such sequences exist by 
\Cref{pr:asymp-dom-props}(\ref{pr:asymp-dom-props:c}).
Otherwise, if $r=0$, we simply let $b_{t,0}=t^2$ and
$b_{t,k+1}=t$
for all $t$, implying $b_{t,0}\seqgt b_{t,k+1}$.
In either case, we then add $0$ and $k+1$ to $\coldefset$, noting that
the two invariants are still maintained.

Finally, we fill in all the rest of the sequences as follows:
If $\colset\subseteq\coldefset$, then we are finished.
Otherwise, let $i$ be any element of
$\colset\setminus\coldefset$, that is, any index whose sequence
$\seq{b_{t,i}}$ has not yet been defined.
Let
\[
  j_0
  =
  \max\Braces{ j\in\coldefset :\: j<i },
  \qquad
  j_1
  =
  \min\Braces{ j\in\coldefset :\: j>i },
\]
which are both attained since $0$ and $k+1$ are included in
$\coldefset$.
This means: $j_0,j_1\in\coldefset$;
$j_0<i<j_1$; and
$\coldefset$ does not include any index strictly between $j_0$ and
$j_1$.
We next define $\seq{b_{t,i}}$ to be any sequence converging to
$+\infty$ and for which
$b_{t,j_0}\seqgt b_{t,i}\seqgt b_{t,j_1}$,
which must exist by
\Cref{pr:asymp-dom-props}(\ref{pr:asymp-dom-props:d})
(since, by our invariants,
$b_{t,j_0}\seqgt b_{t,j_1}$).
We then add $i$ to $\coldefset$, noting that the invariants above are
still maintained.

This process continues until $\colset\subseteq\coldefset$
(which must happen eventually since
an element of the finite set $\colset$ is added to
$\coldefset$ on each iteration).
At that point, $\seq{\bb_t}$ and $\seq{\qq_t}$ satisfy
conditions
  (\ref{thm:seq-rep:a}),
  (\ref{thm:seq-rep:b}), and
  (\ref{thm:seq-rep:c})
of \Cref{thm:seq-rep}
(where, to be clear,
$\bb_t=\trans{[b_{t,1},\dotsc,b_{t,k}]}$
does not include the ``extra'' components
$b_{t,0}$ and $b_{t,k+1}$ defined above).
Since,
$b_{t,\ci{j}}=b'_{t,j}$ for $j=1,\dotsc,r$,
and since all other columns of $\VV$ are $\zero$,
we also have that
$\xx_t=\VV'\bb'_t+\qq_t=\VV\bb_t+\qq_t$,
completing the proof.
\end{proof}

\begin{proof}[Proof of \Cref{thm:seq-rep-not-lin-ind}:]
Let $\vv_1,\dotsc,\vv_k\in\Rn$ be the columns of $\VV$, so
$\VV=[\vv_1,\dotsc,\vv_k]$.

If there exist sequences $\seq{\bb_t}$ and
$\seq{\qq_t}$ satisfying the conditions of the theorem
then
\Cref{thm:i:seq-rep}
(applied to $\qq$ and $\vv_1,\dotsc,\vv_k$)
shows that $\xx_t\rightarrow\xbar$.

To prove the converse, we first use $\VV$ to construct a matrix $\VV'$
to which \Cref{lem:seq-rep-nz-lin-ind} can be applied.
For $i=1,\dotsc,k$, let
\[
  \vv'_i
  =
  \begin{cases}
    \zero & \text{if $\vv_i \in \spnfin{\vv_1,\dotsc,\vv_{i-1}}$,} \\
    \vv_i & \text{otherwise,}
  \end{cases}
\]
and let $\VV'=[\vv'_1,\dotsc,\vv'_k]$.
Then in either of the cases above,
$\vv'_i$ is equal to $\vv_i$ minus a
(possibly empty) linear combination of
$\vv_1,\dotsc,\vv_{i-1}$.
Therefore, by
\Cref{prop:pos-upper:descr},
$\VV'=\VV\RR$ for some positive upper triangular matrix
$\RR\in\R^{k\times k}$.
Also, we can write $\qq=\qq'+\VV\cc$ where $\qq'$ is the projection of
$\qq$ onto $(\colspace{V})^{\perp}$, and for some $\cc\in\Rk$
(\Cref{pr:lin-decomp-rel-vecs}).

Thus,
\[
  \xbar
  =
  \VV\omm\plusl\qq
  =
  \VV\RR\omm\plusl (\qq-\VV\cc)
  =
  \VV'\omm\plusl \qq',
\]
where the second equality is by the Push Lemma~\ref{pr:g:1}.
Further, $\qq'$ is orthogonal to $\vv_1,\dotsc,\vv_k$,
implying $\qq'\perp\VV'$.
Also, the nonzero columns of $\VV'$ are linearly independent since, by
construction, each nonzero column $\vv'_i$ is outside the span of 
$\vv'_1,\dotsc,\vv'_{i-1}$.

Therefore, by \Cref{lem:seq-rep-nz-lin-ind},
there exist sequences $\seq{\bb'_t}$ and $\seq{\qq'_t}$ such that:
$\xx_t=\VV'\bb'_t+\qq'_t$ for all $t$;
$\qq'_t\rightarrow\qq'$;
and
$\bb'_t\rightarrow\omm$ (using \Cref{cor:omm:seq}).
For each $t$, let $\bhat_t=\RR\bb'_t - \cc$
and $\qhat_t=\qq'_t+\VV\cc$.
Then for all $t$,
\begin{equation}  \label{eq:thm:seq-rep-not-lin-ind:1}
  \xx_t
  =
  \VV'\bb'_t+\qq'_t
  =
  \VV(\RR\bb'_t - \cc) + (\qq'_t + \VV\cc)
  =
  \VV \bhat_t + \qhat_t.
\end{equation}
Further,
$\qhat_t=\qq'_t+\VV\cc\rightarrow\qq'+\VV\cc=\qq$,
and
\[
  \bhat_t
  =
  \RR\bb'_t - \cc
  \rightarrow
  \RR\omm \plusl (-\cc)
  =
  \omm,
\]
where convergence is by
\Cref{thm:linear:cont}(\ref{thm:linear:cont:b})
and
\Cref{pr:i:7}(\ref{pr:i:7g}),
and the last equality is by
\Cref{lemma:RR_omm}.
Combined with \Cref{cor:omm:seq},
this shows that
$\seq{\bhat_t}$ and $\seq{\qhat_t}$ satisfy conditions
(\ref{thm:seq-rep:a}),
(\ref{thm:seq-rep:b}), and
(\ref{thm:seq-rep:c})
of \Cref{thm:seq-rep}.

Finally, for each $t$, if $\bhati_{t,i}>0$ for all $i\in\{1,\dotsc,k\}$,
then we let $\bb_t=\bhat_t$ and $\qq_t=\qhat_t$.
Otherwise, if this is not the case, then
we choose $\bb_t$ to be an arbitrary vector in $\Rstrictpos^k$,
and we let $\qq_t=\qhat_t+\VV(\bhat_t-\bb_t)$.
Thus, in either case,
\eqref{eq:thm:seq-rep-not-lin-ind:1} implies that
$\xx_t=\VV\bhat_t+\qhat_t=\VV\bb_t+\qq_t$ for all $t$.
Moreover, since all entries of $\bhat_t$ converge to $+\infty$,
this modification can only affect
finitely many sequence elements; therefore, the modified sequences
$\seq{\bb_t}$ and $\seq{\qq_t}$ have the same convergence properties
as $\seq{\bhat_t}$ and $\seq{\qhat_t}$,
implying in particular that $\seq{\bb_t}$ and $\seq{\qq_t}$
satisfy conditions
(\ref{thm:seq-rep:a}),
(\ref{thm:seq-rep:b}), and
(\ref{thm:seq-rep:c})
of \Cref{thm:seq-rep}
with each $\bb_t\in\Rstrictpos^k$.
This completes the proof.%
\indexg{representations of astral points!sequences converging to|)}%
\indexg{sequence convergence!matrix representation@to matrix representation|)}
\end{proof}

\section{Strong equivalence and span-bound sequences}
\label{sec:strong-equiv-span-bound}

We next define and study two properties for describing sequences.
Both of these will appear in various characterizations in \Cref{sec:conv-pts-part-forms}.
\indexg{strong equivalence|(}%
We begin with a notion of asymptotic equivalence between sequences
called strong equivalence which, as we discuss below, is stronger than the all-directions equivalence introduced in \Cref{sec:astral:construction}.

\begin{definition}
  Two sequences $\seq{\xx_t}$ and $\seq{\yy_t}$ in $\Rn$
  are said to be \emph{strongly equivalent},
  written $\seqeq{\xx_t}{\yy_t}$,%
  \indexm{xt~yt}{$\seqeq{\xx_t}{\yy_t}$}{strongly equivalent}
  if $\xx_t - \yy_t \rightarrow\zero$.
\end{definition}

It can be checked that
strong equivalence is an equivalence relation, and so induces a
(rather fine) partition of all sequences in $\Rn$ into equivalence classes.
As shown next,
a strongly equivalent pair of sequences must have the same
convergence properties in astral space; that is, they
either both converge or both diverge, and if convergent, must
also have the same astral limit.

\begin{proposition}  \label{pr:eq-in-lim-same-lim}
  Let $\seq{\xx_t}$ and $\seq{\yy_t}$ be strongly equivalent sequences in $\Rn$.
  Let $\xbar\in\extspace$.
  Then $\xx_t\rightarrow\xbar$
  if and only if $\yy_t\rightarrow\xbar$.
\end{proposition}

\begin{proof}
It suffices to prove one implication since the other then follows by symmetry.
Suppose $\xx_t\rightarrow\xbar$.
Then
\[
  \yy_t
  =
  \xx_t + (\yy_t - \xx_t)
  \rightarrow
  \xbar \plusl \zero
  =
  \xbar,
\]
with the convergence following from
\Cref{pr:i:7}(\ref{pr:i:7g}).
\end{proof}

\Cref{pr:eq-in-lim-same-lim} shows that
if two convergent sequences are strongly equivalent, then they must
have the same limit.
The converse is false, in general.
In other words, two sequences with the same astral limit need not be
strongly equivalent.
For instance, if $x_t=t+1$ and $y_t=t$, then both the sequences
$\seq{x_t}$ and $\seq{y_t}$ converge to $+\infty$, but these sequences
are not strongly equivalent since $x_t-y_t=1\not\rightarrow 0$.
Thus, strong equivalence induces a finer partition of sequences than equivalence in all directions and is therefore a stronger notion of equivalence.

Strong equivalence is preserved under various operations that are continuous on $\Rn$:

\begin{proposition}
\label{prop:strong:eq:propties}
Let $\seq{\xx_t}, \seq{\xx'_t}, \seq{\yy_t}, \seq{\yy'_t}$ be sequences in $\Rn$ such that
$\seqeq{\xx_t}{\xx'_t}$ and $\seqeq{\yy_t}{\yy'_t}$. Furthermore, let $\A\in\Rmn$ and
let $\seq{\lambda_t}$ be a bounded sequence in $\R$. Then:
\begin{letter-compact}
\item \label{prop:strong:eq:propties:a}
  $\seqeq{\xx_t+\yy_t}{\xx'_t+\yy'_t}$.
\item \label{prop:strong:eq:propties:b}
  $\seqeq{\A\xx_t}{\A\xx'_t}$.
\item \label{prop:strong:eq:propties:c}
  $\seqeq{\lambda_t\xx_t}{\lambda_t\xx'_t}$.
\end{letter-compact}
\end{proposition}
\begin{proof}
The statements follow from the definition of strong equivalence and continuity properties of standard vector addition, linear maps, and scalar multiplication.%
\indexg{strong equivalence|)}
\end{proof}

\indexg{representational span!singleton@of singleton|(}%
Let $\xbar=\VV\omm\plusl\qq$ where $\VV=[\vv_1,\dotsc,\vv_k]$ for some
$\vv_1,\dotsc,\vv_k,\qq\in\Rn$.
In considering sequences $\seq{\xx_t}$ in $\Rn$ that converge to
$\xbar$, it is sometimes convenient to require that
each element $\xx_t$ is constructed as a linear combination of only
$\vv_1,\dotsc,\vv_k,\qq$,
that is, so that each $\xx_t$ is in the column space of
$[\VV,\qq]$.
This would seem natural since, in all other directions, the sequence
is converging to zero.
This was the case, for instance, for the polynomially graded sequence of
Example~\ref{ex:poly-speed-intro}.

To study such sequences, we give a name to
the column space just mentioned:
\begin{definition}  \label{def:rep-span-snglton}
Let $\xbar=\VV\omm\plusl\qq$ where $\VV\in\Rnk$, $k\geq 0$, and $\qq\in\Rn$.
Then the \emph{representational span} of the singleton $\{\xbar\}$,
denoted $\rspanxbar$,%
\indexm{rspan200}{$\rspanxbar$}{representational span (of singleton)}
is equal to
$\colspace[\VV,\qq]$, the column space of the matrix $[\VV,\qq]$,
or equivalently, the
span of $(\columns{\VV})\cup\{\qq\}$.
\end{definition}
For now, we only define representational span for singletons;
later, in \Cref{sec:astral-span}, we extend this definition for
arbitrary astral sets.

Superficially, this definition would appear to depend on a
particular representation of~$\xbar$, which might seem problematic
since such representations are not unique.
However, the next \namecref{pr:rspan-sing-equiv-dual}
implies that this set is
determined by $\xbar$ alone, without any dependence on the specific
representation chosen for it.

\begin{proposition}  \label{pr:rspan-sing-equiv-dual}
  Let $\xbar=\VV\omm\plusl\qq$ where $\VV\in\Rnk$, $k\geq 0$,
  and $\qq\in\Rn$.
  Let $\zz\in\Rn$.
  Then the following are equivalent:
  \begin{letter-compact}
  \item  \label{pr:rspan-sing-equiv-dual:a}
    $\zz\in\colspace[\VV,\qq]$ (that is, $\zz\in\rspanxbar$).
  \item  \label{pr:rspan-sing-equiv-dual:b}
    For all $\uu\in\Rn$, if $\xbar\cdot\uu=0$ then $\zz\cdot\uu=0$.
  \end{letter-compact}
\end{proposition}

\begin{proof}
  ~

\begin{proof-parts}
\pfpart{%
  (\ref{pr:rspan-sing-equiv-dual:a})
  $\Rightarrow$
  (\ref{pr:rspan-sing-equiv-dual:b}):
}
Suppose $\zz\in\colspace[\VV,\qq]$.
Let $\uu\in\Rn$ be such that $\xbar\cdot\uu=0$.
Then, by \Cref{pr:vtransu-zero}, $\uu\perp\VV$
and $\qq\cdot\uu=0$.
Consequently, $\zz\cdot\uu=0$ since
$\zz$ is a linear combination of $\qq$ and the columns of $\VV$.

\pfpart{%
  (\ref{pr:rspan-sing-equiv-dual:b})
  $\Rightarrow$
  (\ref{pr:rspan-sing-equiv-dual:a}):
}
Suppose~(\ref{pr:rspan-sing-equiv-dual:b}) holds.
Let $L=\colspace[\VV,\qq]$, which is a linear subspace.
By linear algebra, we can write $\zz=\yy+\uu$ where
$\yy\in L$ and $\uu$ is in $\Lperp$, the set of vectors orthogonal to every point in $L$.
In particular, $\uu\perp\VV$ and $\uu\perp\qq$ so
$\xbar\cdot\uu=0$
by \Cref{pr:vtransu-zero}.
Therefore,
\[
  0
  =
  \zz\cdot\uu
  =
  \yy\cdot\uu
  +
  \uu\cdot\uu
  =
  \norm{\uu}^2,
\]
where the first equality is by assumption (since $\xbar\cdot\uu=0$),
and the third is because $\yy\perp\uu$.
Thus, $\uu=\zero$ so $\zz=\yy\in L$.%
\indexg{representational span!singleton@of singleton|)}
\qedhere
\end{proof-parts}
\end{proof}

\begin{definition}
\indexg{span-boundness|(}%
Let $\seq{\xx_t}$ be a sequence in $\Rn$ that converges in astral
space, and let $\xbar\in\extspace$ denote its limit.
We say that $\seq{\xx_t}$ is \emph{span-bound} if 
$\xx_t\in\rspanxbar$ for all $t$.
\end{definition}

Not every convergent sequence is span-bound.
For instance, the sequence in $\R$ with elements $x_t=1/t$ converges
to $0$, but it is not span-bound since $\rspanset{0}=\{0\}$.
Nevertheless, as we show next,
every convergent sequence is strongly equivalent to some
span-bound sequence:

\begin{theorem}   \label{thm:spam-limit-seqs-exist}
  Let $\seq{\xx_t}$ be a sequence in $\Rn$ converging to some point
  $\xbar\in\extspace$.
  Then there exists a span-bound sequence $\seq{\xx'_t}$ in $\Rn$
  that is strongly equivalent to $\seq{\xx_t}$
  (implying that also $\xx'_t\rightarrow\xbar$).
\end{theorem}

\begin{proof}
We can write $\xbar=\VV\omm\plusl\qq$ for some $\VV\in\Rnk$, $k\geq 0$,
and $\qq\in\Rn$.
Let $L=\rspanxbar=\colspace[\VV,\qq]$.
Let $\PP\in\R^{n\times n}$ be the projection matrix onto $L$,
and let $\Iden$ be the $n\times n$ identity matrix.
Then for all $\ww\in L$,
$\PP\ww=\ww$
(\Cref{pr:proj-mat-props}\ref{pr:proj-mat-props:d}),
so $(\Iden-\PP)\ww=\zero$,
implying
\begin{equation}    \label{eq:thm:spam-limit-seqs-exist:1}
  (\Iden-\PP)\xbar
  =
  (\Iden-\PP)\VV\omm\plusl(\Iden-\PP)\qq
  =
  \zero.
\end{equation}
For each $t$, let $\xx'_t=\PP \xx_t$.
Then
$\xx_t-\xx'_t=(\Iden-\PP)\xx_t\rightarrow(\Iden-\PP)\xbar=\zero$,
by \Cref{thm:linear:cont}(\ref{thm:linear:cont:b}) and \eqref{eq:thm:spam-limit-seqs-exist:1},
so $\seqeq{\xx_t}{\xx'_t}$,
implying $\xx'_t\rightarrow\xbar$ by
\Cref{pr:eq-in-lim-same-lim}.
Furthermore, $\seq{\xx'_t}$ is span-bound
since each $\xx'_t\in L$.%
\indexg{span-boundness|)}
\end{proof}

\section{Convergence to points represented by specific expressions}
\label{sec:conv-pts-part-forms}

We next apply the tools and notions developed so far,
especially \Cref{thm:seq-rep-not-lin-ind},
to characterize sequences that converge
to limits
that can be represented by specific kinds of expressions.
\indexg{sequence convergence!matrix product@to matrix product|(}%
\indexg{matrix product, astral!sequences converging to|(}%
The first such result is a characterization of sequences
$\seq{\yy_t}$ in $\Rn$
that converge to limits of the form $\A\xbar$ for some $\xbar\in\eRn$
and some matrix $\A\in\Rmn$. By continuity of linear maps (\Cref{thm:linear:cont}\ref{thm:linear:cont:b}), we
know that if $\seq{\xx_t}$ is a
sequence in $\Rn$ that converges to $\xbar\in\extspace$, then
the sequence
$\yy_t=\A\xx_t$ must converge to $\A\xbar$.
Using \Cref{thm:seq-rep-not-lin-ind}, we prove a type of converse:
We show that any
sequence $\seq{\yy_t}$ in $\Rm$ that converges to $\A\xbar$ must be strongly equivalent
to a sequence of the form $\seq{\A\xx_t}$ such that $\xx_t\to\xbar$.

\begin{theorem}   \label{thm:inv-lin-seq}
  Let $\A\in\Rmn$, let $\seq{\yy_t}$ be a sequence in $\Rm$,
  and let $\xbar\in\extspace$.
  Then the following are equivalent:
  \begin{letter-compact}
  \item     \label{thm:inv-lin-seq:a}
    $\yy_t\rightarrow\A\xbar$.
  \item     \label{thm:inv-lin-seq:b}
    There exists a sequence $\seq{\xx_t}$ in $\Rn$ such that
    $\xx_t\rightarrow\xbar$ and
    $\seqeq{\yy_t}{\A\xx_t}$.
  \end{letter-compact}
\end{theorem}

\begin{proof}
  ~

\begin{proof-parts}
\pfpart{%
  (\ref{thm:inv-lin-seq:b})
  $\Rightarrow$
  (\ref{thm:inv-lin-seq:a}):
}
If there exists a sequence $\seq{\xx_t}$ as in
part~(\ref{thm:inv-lin-seq:b}),
then $\A\xx_t\rightarrow\A\xbar$
by \Cref{thm:linear:cont}(\ref{thm:linear:cont:b}),
so $\yy_t\rightarrow\A\xbar$ by
\Cref{pr:eq-in-lim-same-lim}.

\pfpart{%
  (\ref{thm:inv-lin-seq:a})
  $\Rightarrow$
  (\ref{thm:inv-lin-seq:b}):
}
Suppose $\yy_t\rightarrow\A\xbar$.
We can write $\xbar=\VV\omm\plusl\qq$ for some $\VV\in\Rnk$ and
$\qq\in\Rn$, implying
$\A\xbar=\A\VV\omm\plusl\A\qq$.
Since $\yy_t\rightarrow\A\xbar$, by
\Cref{thm:seq-rep-not-lin-ind},
there exist sequences $\seq{\bb_t}$ in $\Rk$
and $\seq{\rr_t}$ in $\Rm$ such that:
$\yy_t=\A\VV\bb_t + \rr_t$;
conditions
(\ref{thm:seq-rep:a}) and (\ref{thm:seq-rep:b})
of \Cref{thm:seq-rep}
hold; and
$\rr_t\rightarrow\A\qq$.

Let $\xx_t=\VV\bb_t + \qq$.
Then $\xx_t\rightarrow\xbar$ by
\Cref{thm:seq-rep-not-lin-ind}
(or \Cref{thm:i:seq-rep}).
Furthermore,
\[
  \A\xx_t - \yy_t
  =
  (\A\VV\bb_t + \A\qq) - (\A\VV\bb_t + \rr_t)
  =
  \A\qq - \rr_t
  \rightarrow
  \zero
\]
since $\rr_t\rightarrow\A\qq$,
so $\seqeq{\yy_t}{\A\xx_t}$.
\qedhere
\end{proof-parts}
\end{proof}

\Cref{thm:inv-lin-seq} is especially useful when working with
projection matrices, as shown in the
next \namecref{cor:inv-proj-seq}.
For instance, if $\PP$ is a projection matrix
onto some subspace $L$, then the \namecref{cor:inv-proj-seq}
shows that any sequence in $\Rn$ that converges to
the projection $\PP\xbar$ of an astral point $\xbar$
can be turned into a sequence that converges to~$\xbar$ by adding
to it a sequence of vectors that are each orthogonal to $L$.

\begin{corollary}  \label{cor:inv-proj-seq}
  Let $\WW\in\Rnk$, and let $\PP\in\R^{n\times n}$ be the projection
  matrix
  onto $(\colspace{\WW})^\perp$.
  Let
  $\seq{\yy_t}$ be a sequence in $\Rn$ converging to $\PP\xbar$
  for some $\xbar\in\extspace$.
  Then there exists a sequence $\seq{\cc_t}$ in $\Rk$ such that
  $\yy_t + \WW\cc_t \rightarrow \xbar$.
\end{corollary}

\begin{proof}
By \Cref{thm:inv-lin-seq} (applied with $\A=\PP$),
there exists a sequence $\seq{\xx_t}$ in $\Rn$ with
$\xx_t\rightarrow\xbar$ and $\seqeq{\yy_t}{\PP\xx_t}$.
Decomposing $\xx_t$ as in
\Cref{pr:lin-decomp-rel-vecs},
for each $t$,
we can write
$\xx_t=\WW\cc_t+\PP\xx_t$ for some sequence $\seq{\cc_t}$ in $\Rk$.
Thus,
by \Cref{prop:strong:eq:propties}(\ref{prop:strong:eq:propties:a}),
$\yy_t+\WW\cc_t\seqeqsymbol\PP\xx_t+\WW\cc_t=\xx_t$,
implying
$\yy_t+\WW\cc_t\to\xbar$
by \Cref{pr:eq-in-lim-same-lim}.%
\indexg{matrix product, astral!sequences converging to|)}%
\indexg{sequence convergence!matrix product@to matrix product|)}
\end{proof}

\indexg{sequence convergence!omega x+y@to $\limray{\xbar}\plusl\ybar$|(}%
Next, we characterize
sequences that converge to limits of the form
$\limray{\xbar}\plusl\ybar$,
for some points $\xbar,\ybar\in\extspace$.
Points of this form will appear later in various contexts.
For instance, by \Cref{cor:h:1}, any infinite astral point $\xbar$
can be written
as $\xbar=\limray{\vv}\plusl\ybar$ for
some $\vv\in\Rn\wo\set{\zero}$ and some $\ybar\in\extspace$.
This form of $\xbar$ is a special case of the general form we now consider.

Intuitively, if $\seq{\xx_t}$ is a sequence in $\Rn$ converging to
$\xbar$
and $\seq{\lambda_t}$ is a sequence in $\R$ converging to
$\oms=+\infty$,
then we might expect that $\lambda_t \xx_t\rightarrow\limray{\xbar}$.
Indeed, this intuition is correct, but only if we further
stipulate that $\seq{\xx_t}$ is span-bound.
(That this is necessary is shown momentarily in
Example~\ref{ex:no-span-bnd-no-conv}.)
If, in addition, the sequence $\seq{\yy_t}$ in $\Rn$ converges to
$\ybar$, then we might further expect that
$\lambda_t\xx_t+\yy_t\rightarrow\limray{\xbar}\plusl\ybar$.
This is also correct,
provided that $\lambda_t$ is growing fast enough to
ensure the dominance of $\lambda_t\xx_t$ over $\yy_t$.
The next theorem makes all this precise:

\begin{theorem}   \label{thm:lim-plusl}
  Let $\seq{\xx_t}$ be a span-bound sequence in $\Rn$ converging to
  $\xbar\in\eRn$,
  let $\seq{\yy_t}$ be a sequence in $\Rn$ converging to
  $\ybar\in\eRn$,
  and let $\seq{\lambda_t}$ be a sequence in $\R$ converging to $+\infty$ such that
  $\norm{\yy_t}/\lambda_t\to 0$. Then
  $\lambda_t\xx_t+\yy_t\to\limray{\xbar}\plusl\ybar$.
\end{theorem}

\begin{proof}
Let $\zbar = \limray{\xbar}\plusl\ybar$,
and for each $t$,
let $\zz_t=\lambda_t \xx_t + \yy_t$.
Let $\uu\in\Rn$.
We aim to show $\zz_t\cdot\uu\rightarrow\zbar\cdot\uu$, which will
prove the claim by
\Cref{thm:i:1}(\ref{thm:i:1c}).

Suppose first that $\xbar\cdot\uu>0$,
implying $\zbar\cdot\uu=\limray{\xbar}\cdot\uu\plusl\ybar\cdot\uu=+\infty$.
Also, for all $t$ sufficiently large (so that $\lambda_t>0$),
\begin{equation}    \label{eq:thm:lim-plusl:1}
  \zz_t\cdot\uu
  =
  \lambda_t \biggParens{ \xx_t\cdot \uu + \frac{\yy_t\cdot\uu}{\lambda_t} }
  \geq
  \lambda_t \biggParens{ \xx_t\cdot \uu - \frac{\norm{\yy_t} \norm{\uu}}
                                        {\lambda_t}
                  }.
\end{equation}
The equality is from the definition of $\zz_t$.
The inequality is by the Cauchy-Schwarz inequality.
Since $\norm{\yy_t}/\lambda_t\rightarrow 0$,
and since $\xx_t\cdot\uu\rightarrow\xbar\cdot\uu$
(by \Cref{thm:i:1}\ref{thm:i:1c}),
the parenthesized expression on the right-hand side of
\eqref{eq:thm:lim-plusl:1}
is converging to $\xbar\cdot\uu>0$.
Since $\lambda_t\rightarrow+\infty$, this shows
(by \Cref{prop:lim:eR}\ref{i:lim:eR:genmul})
that $\zz_t\cdot\uu\rightarrow+\infty=\zbar\cdot\uu$.

The case $\xbar\cdot\uu<0$ can be handled symmetrically, or by applying
the preceding argument to $-\uu$.

Finally, suppose $\xbar\cdot\uu=0$,
implying $\xx_t\cdot\uu=0$ for all $t$
by \Cref{pr:rspan-sing-equiv-dual} since $\seq{\xx_t}$ is span-bound.
Then
\[
   \zz_t\cdot\uu
   =
   \lambda_t \xx_t\cdot\uu + \yy_t\cdot\uu
   =
   \yy_t\cdot\uu
   \rightarrow
   \ybar\cdot\uu
   =
   \limray{\xbar}\cdot\uu \plusl \ybar\cdot\uu
   =
   \zbar\cdot\uu
\]
(with convergence by \Cref{thm:i:1}\ref{thm:i:1c}),
completing the proof.
\end{proof}

For example, given $\vv\in\Rn$ and a sequence $\seq{\yy_t}$ in $\Rn$
converging to some point ${\ybar\in\extspace}$, we can
construct a sequence $\seq{\zz_t}$ in $\Rn$ converging to
$\limray{\vv}\plusl\ybar$ using
\Cref{thm:lim-plusl} 
(applied with $\xbar=\xx_t=\vv$)
by letting
$\zz_t=\lambda_t \vv + \yy_t$ with $\lambda_t$ chosen so that
$\lambda_t\rightarrow+\infty$ and
$\norm{\yy_t}/\lambda_t\rightarrow 0$, for instance, by
letting $\lambda_t=t (1+\norm{\yy_t})$ for all $t$.

Without the span-boundness assumption on $\seq{\xx_t}$, \Cref{thm:lim-plusl}
does not hold in general:

\begin{example}  \label{ex:no-span-bnd-no-conv}
In one dimension,  let
$x=y=0$, and for all $t$, let
$x_t=1/t$, $y_t=0$, and $\lambda_t=t$.
Then
$x_t\to x$, $y_t\to y$, $\lambda_t\to+\infty$ and $\abs{y_t}/\lambda_t\to 0$,
but $\lambda_t x_t+y_t\to 1\neq 0 = \limray{x}\plusl y$.
\end{example}

We next prove a kind of converse to \Cref{thm:lim-plusl}
showing that a sequence $\seq{\zz_t}$ in $\Rn$
converges to a point of the form $\limray{\xbar}\plusl\ybar$
if and only if it
can be written
as $\zz_t=\lambda_t\xx_t+\yy_t$ for some sequences
$\seq{\xx_t}$, $\seq{\yy_t}$ and $\seq{\lambda_t}$ satisfying the
conditions of \Cref{thm:lim-plusl}.

\begin{theorem}   \label{thm:lim-plusl-inv}
  Let $\xbar,\ybar\in\extspace$ and
  let $\seq{\zz_t}$ be a sequence in $\Rn$.
  Then the following are equivalent:
  \begin{letter-compact}
  \item   \label{thm:lim-plusl-inv:a}
    $\zz_t\to\limray{\xbar}\plusl\ybar$.
  \item   \label{thm:lim-plusl-inv:b}
    There exist
    a span-bound sequence $\seq{\xx_t}$ in $\Rn$,
    a sequence $\seq{\yy_t}$ in $\Rn$,
    and a sequence $\seq{\lambda_t}$ in $\R$
    such that
    $\zz_t=\lambda_t\xx_t+\yy_t$ for all $t$,
    $\xx_t\rightarrow\xbar$,
    $\yy_t\rightarrow\ybar$,
    $\lambda_t\rightarrow+\infty$,
    and
    $\norm{\yy_t}/\lambda_t\to 0$.
  \end{letter-compact}
  This equivalence still holds if,
  in statement~(\ref{thm:lim-plusl-inv:b}),
  we also require $\lambda_t>0$ for all $t$.
\end{theorem}

\begin{proof}
That
(\ref{thm:lim-plusl-inv:b})
$\Rightarrow$
(\ref{thm:lim-plusl-inv:a})
was proved in \Cref{thm:lim-plusl}.

For the converse,
suppose $\zz_t\to\limray{\xbar}\plusl\ybar$.
Let $\xbar=\VV\omm\plusl\qq$ and $\ybar=\VV'\omm\plusl\qq'$ be
representations of $\xbar$ and~$\ybar$, respectively, and let $\VV''=[\VV,\qq,\VV']$.
Then
\[
  \limray{\xbar}\plusl\ybar=(\VV\omm\plusl\limray{\qq})\plusl(\VV'\omm\plusl\qq')
  =[\VV,\qq,\VV']\omm\plusl\qq'
  =\VV''\omm\plusl\qq',
\]
with the first equality by
\Cref{pr:mult-rep-by-scalar},
and the second by
\Cref{prop:split}.
Let $k$, $k'$, and $k''=k+1+k'$ denote, respectively, the number of columns of $\VV$, $\VV'$, and $\VV''$.

Since $\zz_t\to\VV''\omm\plusl\qq'$,
\Cref{thm:seq-rep-not-lin-ind} implies that $\zz_t=\VV''\bb''_t+\qq'_t$
for some sequences $\seq{\bb''_t}$ in $\Rstrictpos^{k''}$ and $\seq{\qq'_t}$ in $\R^n$ such that
the entries of $\bb''_t$ converge to $+\infty$ at decreasing rates and $\qq'_t\to\qq'$.
For all $t$, write $\bb''_t=\mtuple{\bb_t,\lambda_t,\bb'_t}$ where $\bb_t\in\Rstrictpos^{k}$,
$\lambda_t\in\Rstrictpos$, and $\bb'_t\in\Rstrictpos^{k'}$,
implying, in particular, that $\lambda_t\rightarrow+\infty$.

Let
\[
  \xx_t=\frac{\VV\bb_t}{\lambda_t}+\qq,
  \qquad
  \yy_t=\VV'\bb'_t+\qq'_t,
\]
which we now argue converge to $\xbar$ and $\ybar$, respectively.
Since the entries of $\bb''_t$
converge to~$+\infty$ at decreasing rates, the same is also true for entries of
$\bb'_t$.
Further, for $i=1,\dotsc,k$,
$b_{t,i}/\lambda_t\rightarrow+\infty$
(by \Cref{prop:dec:trans}, applied to $\bb''_t$),
and if $i<k$ then
${b_{t,i}\seqgt b_{t,i+1}\seqgt \lambda_t}$,
implying
$b_{t,i}/\lambda_t \seqgt b_{t,i+1}/\lambda_t$
by \Cref{pr:asymp-dom-props}(\ref{pr:asymp-dom-props:e}).
Thus, the entries of~$\bb_t/\lambda_t$ also converge to $+\infty$ at
decreasing rates.
Therefore, by \Cref{thm:seq-rep-not-lin-ind},
\begin{equation*}
\xx_t\to\VV\omm\plusl\qq=\xbar,
\qquad
\yy_t\to\VV'\omm\plusl\qq'=\ybar.
\end{equation*}
Moreover, 
the sequence $\seq{\xx_t}$ is span-bound
since $\xx_t\in\colspace[\VV,\qq]$ for all $t$.

Next, letting $\vv'_1,\dotsc,\vv'_{k'}\in\Rn$ denote the columns of $\VV'$,
we have for all $t$,
\begin{equation}   \label{eq:thm:lim-plusl-inv:1}
  0
  \leq
  \frac{\norm{\yy_t}}{\lambda_t}
  =
  \frac{1}{\lambda_t}
  \cdot
  \BiggNorm{\sum_{i=1}^{k'} b'_{t,i} \vv'_i + \qq'_t}
  \leq
  \sum_{i=1}^{k'} \biggParens{\frac{b'_{t,i}}{\lambda_t}} \norm{\vv'_i}
  +
  \frac{\norm{\qq'_t}}{\lambda_t},
\end{equation}
where the second inequality is by the triangle inequality.
Applying \Cref{prop:dec:trans} to $\bb''_t$
implies that $b'_{t,i}/\lambda_t\to 0$ for $i=1,\dotsc,k'$.
Also,
$\norm{\qq'_t}/\lambda_t\to 0$
since $\qq'_t\to\qq'$ and $\lambda_t\rightarrow+\infty$.
Therefore, the right-hand side of \eqref{eq:thm:lim-plusl-inv:1}
converges to $0$, so $\norm{\yy_t}/\lambda_t$ does as well.
  
Thus, for all $t$,
\begin{align*}
  \zz_t
  =
  \VV''\bb''_t+\qq'_t
  &=
  \Bracks{\VV,\qq,\VV'}\mtuple{\bb_t,\lambda_t,\bb'_t}+\qq'_t
  \\
  &=
  \lambda_t\biggParens{\frac{\VV\bb_t}{\lambda_t}+\qq}
  + (\VV'\bb'_t+\qq'_t)
  =
  \lambda_t\xx_t+\yy_t,
\end{align*}
with the sequences $\seq{\xx_t}$, $\seq{\yy_t}$, and $\seq{\lambda_t}$ satisfying all the conditions
of the \namecref{thm:lim-plusl-inv}
(including that $\lambda_t>0$ for all $t$).
\end{proof}

\indexg{omega-multiples@$\oms$-multiples!sequences converging to|(}%
As a special case of \Cref{thm:lim-plusl,thm:lim-plusl-inv},
we obtain the following characterization of sequences in $\Rn$ that converge to
points of the form $\limray{\xbar}$. Note that even in this special case,
the span-boundness requirement cannot be dropped
(as follows again from Example~\ref{ex:no-span-bnd-no-conv}).

\begin{corollary}   \label{thm:seq-to-limray:0}
  Let $\seq{\xx_t}$ be a span-bound sequence in $\Rn$ converging to $\xbar\in\eRn$, and let $\seq{\lambda_t}$ be a sequence in $\R$ converging to $+\infty$. Then $\lambda_t\xx_t\to\limray{\xbar}$.
\end{corollary}

\begin{proof}
This is immediate from \Cref{thm:lim-plusl}
(with $\ybar=\zero$ and $\yy_t=\zero$ for all~$t$).
\end{proof}

\begin{corollary}   \label{thm:seq-to-limray}
  Let $\xbar\in\extspace$ and
  let $\seq{\zz_t}$ be a sequence in $\Rn$.
  Then the following are equivalent:
  \begin{letter-compact}
  \item     \label{thm:seq-to-limray:a}
    $\zz_t\rightarrow\limray{\xbar}$.
  \item     \label{thm:seq-to-limray:b}
    There exist a sequence $\seq{\lambda_t}$ in $\R$
    and a span-bound sequence $\seq{\xx_t}$ in $\Rn$
    such that $\lambda_t\rightarrow+\infty$,
    $\xx_t\rightarrow\xbar$,
    and
    $\seqeq{\zz_t}{\lambda_t\xx_t}$.
  \end{letter-compact}
  This equivalence still holds if,
  in statement~(\ref{thm:seq-to-limray:b}),
  we also require $\lambda_t>0$ for all $t$.
\end{corollary}

\begin{proof}
  ~

\begin{proof-parts}
\pfpart{%
  (\ref{thm:seq-to-limray:b})
  $\Rightarrow$
  (\ref{thm:seq-to-limray:a}):
}
Suppose $\seq{\xx_t}$ and $\seq{\lambda_t}$
are as in
statement~(\ref{thm:seq-to-limray:b}).
Then
by \Cref{thm:seq-to-limray:0}, $\lambda_t\xx_t\to\limray{\xbar}$, and since
$\seqeq{\zz_t}{\lambda_t\xx_t}$, we obtain $\zz_t\to\limray{\xbar}$
(by \Cref{pr:eq-in-lim-same-lim}).

\pfpart{%
  (\ref{thm:seq-to-limray:a})
  $\Rightarrow$
  (\ref{thm:seq-to-limray:b}):
}
Suppose $\zz_t\rightarrow\limray{\xbar}$.
Then by \Cref{thm:lim-plusl-inv},
$\zz_t=\lambda_t\xx_t+\qq_t$ for some sequence $\seq{\qq_t}$ in $\Rn$ converging to $\zero$,
some span-bound sequence $\seq{\xx_t}$ in $\Rn$ converging to $\xbar$, and some sequence
$\seq{\lambda_t}$ in $\Rstrictpos$ converging to $+\infty$.
Thus, $\zz_t-\lambda_t\xx_t=\qq_t\to\zero$, meaning that $\seqeq{\zz_t}{\lambda_t\xx_t}$, and
that the sequences $\seq{\xx_t}$ and $\seq{\lambda_t}$ satisfy the conditions of statement~(\ref{thm:seq-to-limray:b}).%
\indexg{omega-multiples@$\oms$-multiples!sequences converging to|)}%
\indexg{sequence convergence!omega x+y@to $\limray{\xbar}\plusl\ybar$|)}%
\qedhere
\end{proof-parts}
\end{proof}

\chapter{Astral-point decompositions and the structure of astral space}
\chaptermark{Astral-point decompositions}
\label{sec:ast:decomp}

As we saw in
\Cref{cor:h:1},
every astral point in $\extspace$ can be fully decomposed into a
leftward sum of astrons, together with a finite vector in $\Rn$.
Based on that decomposition,
this chapter explores how astral points can be split apart in useful
ways.
Specifically, we explore how every astral point can be split
into a finite part and another component, called an icon,
which, if the point is infinite, can be viewed as its
``purely infinite'' part.
We will see moreover how this decomposition reveals the structure of
astral space as a whole, yielding a natural
partitioning of the space into subregions called galaxies,
each a topological copy of Euclidean space,
though of varying dimensions.

Next, for any infinite astral point, we focus on the first,
or most dominant, astron in its decomposition (where, without loss of
generality, we assume all astrons in that decomposition are nonzero).
We will see how that first astron captures the strongest direction in which
every sequence converging to the point is going to infinity.
This will lay the groundwork for a technique in which an astral point
can be pulled apart and analyzed, one astron at a time, which will be
used frequently going forward, particularly in analyzing astral
functions.

\section{Icons}

\indexg{icons|(}%
We begin by introducing icons and how they can be used to decompose
any astral point.

\indexg{astral space!algebraic structure|(}%
In algebraic terms, we have seen that astral space
is
closed under leftward addition (\Cref{pr:i:6}), and that
this operation is associative
(\Cref{pr:i:7}\ref{pr:i:7a}).
This shows that astral space is a semigroup under this operation, and
furthermore is a monoid since $\zero$ is an identity
\indexg{astral space!algebraic structure|)}%
element.
\indexg{idempotent|(}%
We will especially be interested in idempotent elements of this
semigroup, that is, those points that are unchanged when leftwardly
added to themselves.
Such points are called icons:

\begin{definition}
\indexg{icons!defined|(}%
A point $\ebar\in\extspace$ is an \emph{icon},
or is said to be \emph{iconic}, if
$\ebar\plusl\ebar=\ebar$.
The set of all icons is denoted
\[
\indexg{icons!defined|)}%
\indexm{en}{$\corezn$}{set of all icons}%
  \corezn = \regBraces{ \ebar\in\extspace:\:
                       \ebar\plusl\ebar = \ebar }.
\]
\end{definition}

For example, in $\Rext=\extspac{1}$, the only icons are
$\corez{1}=\{-\infty,0,+\infty\}$.

The term ``icon'' is
derived as a contraction of ``idempotent'' and ``cone,'' with the
latter referring to a cone-like property of such points,
specifically, that $\ebar=\lambda\ebar$ for any $\lambda\in\Rstrictpos$
(see \Cref{pr:i:8}\ref{pr:i:8d}
\indexg{idempotent|)}%
below).
The next proposition gives three equivalent characterizations
of icons:

\begin{proposition}  \label{pr:icon-equiv}
  Let $\ebar\in\extspace$.
  Then the following are equivalent:
  \begin{letter-compact}
  \item  \label{pr:icon-equiv:a}
    $\ebar$ is an icon; that is, $\ebar\plusl\ebar=\ebar$.
  \item  \label{pr:icon-equiv:b}
\indexg{coupling function!icons and|(}%
    For all $\uu\in\Rn$, $\ebar\cdot\uu\in\{-\infty,0,+\infty\}$.
  \item  \label{pr:icon-equiv:c}
\indexg{representations of astral points!icons@of icons|(}%
    $\ebar=\VV\omm$ for some
    $\VV\in\R^{n\times k}$, $k\geq 0$.
  \end{letter-compact}
Furthermore, the same equivalence holds if the matrix in
part~(\ref{pr:icon-equiv:c}) is required to be column-orthogonal.
\end{proposition}

\begin{proof}
~
\begin{proof-parts}

\pfpart{\RefsImplication{pr:icon-equiv:a}{pr:icon-equiv:b}:}
Suppose $\ebar\plusl\ebar=\ebar$.
Then for all $\uu\in\Rn$,
by \Cref{pr:i:6},
$\ebar\cdot\uu\plusl\ebar\cdot\uu=\ebar\cdot\uu$, which is impossible
if $\ebar\cdot\uu$ is a positive or negative real number, that is,
unless $\ebar\cdot\uu$ is in $\{-\infty,0,+\infty\}$.

\pfpart{\RefsImplication{pr:icon-equiv:b}{pr:icon-equiv:c}:}
Suppose $\ebar$ satisfies (\ref{pr:icon-equiv:b}).
Let
$\ebar=\VV\omm\plusl\qq$
be its canonical representation
(which exists by \Cref{pr:uniq-canon-rep}).
Then because $\qq$ is orthogonal to all the columns of $\VV$,
by \Cref{pr:vtransu-zero},
$\ebar\cdot\qq=\qq\cdot\qq=\norm{\qq}^2$.
Since $\ebar\cdot\qq\in\{-\infty,0,+\infty\}$, this quantity must be
$0$, implying $\qq=\zero$.

Note that the matrix $\VV$ is column-orthogonal, thereby
satisfying the additional
requirement stated at the end of the proposition.

\pfpart{\RefsImplication{pr:icon-equiv:c}{pr:icon-equiv:a}:}
Suppose $\ebar=\VV\omm$ for some $\VV\in\R^{n\times k}$, $k\geq 0$,
and let $\PP\in\R^{n\times n}$ be the projection matrix onto $(\colspace\VV)^\perp$.
Then
\[
   \VV\omm \plusl \VV\omm
   =
   \VV\omm \plusl \PP\VV\omm
   =
   \VV\omm\plusl\zero  %
   =
   \VV\omm,
\]
where the first equality is by the Projection Lemma
(\Cref{lemma:proj}),
and the second because $\PP\VV=\zeromat{n}{k}$
(by
\Cref{pr:proj-mat-props}\ref{pr:proj-mat-props:e}).
Therefore, $\ebar$ is
\indexg{coupling function!icons and|)}%
\indexg{representations of astral points!icons@of icons|)}%
an icon.
\qedhere
\end{proof-parts}
\end{proof}

\indexg{astral points!iconic decomposition of|(}%
\indexg{representations of astral points!icon and finite part@as icon and finite part|(}%
\indexg{decompositions, astral!icon and finite part@as icon and finite part|(}%
\indexg{icons!decomposition using|(}%
We saw in \Cref{cor:h:1}
that every point $\xbar\in\extspace$ can be represented as
$\VV\omm\plusl\qq$.
As a result of
\Crefequiv{pr:icon-equiv}{pr:icon-equiv:c}{pr:icon-equiv:a},
$\VV\omm$ is an icon,
which means that $\xbar$ can be written $\xbar=\ebar\plusl\qq$ for
some icon $\ebar\in\corezn$ and some finite vector $\qq\in\Rn$.
\indexg{iconic part|(}%
\indexg{finite part|(}%
In other words, $\xbar$ can be decomposed into an
\emph{iconic part}, $\ebar$, and a
\emph{finite part}, $\qq$.
Furthermore, as shown next,
the iconic part is uniquely determined by $\xbar$.
The finite part, on the other hand, is not in general uniquely
determined.
For example, in $\Rext$,
$+\infty=+\infty\plusl q$ for all $q\in\R$,
so $+\infty$'s decomposition is not unique.%
\indexg{finite part|)}

\begin{theorem}  \label{thm:icon-fin-decomp}
  Let $\xbar\in\extspace$.
  Then
  $\xbar=\ebar\plusl\qq$ for some $\ebar\in\corezn$ and $\qq\in\Rn$.
  Thus, $\extspace = \corezn\plusl\Rn$.
  Furthermore, $\ebar$ is uniquely determined by $\xbar$; that is,
  if it also holds that $\xbar=\ebar'\plusl\qq'$ for some
  $\ebar'\in\corezn$ and $\qq'\in\Rn$, then $\ebar=\ebar'$.
\end{theorem}

\begin{proof}
By \Cref{cor:h:1}, $\xbar$ can be written as
$\xbar=\VV\omm\plusl\qq$ for some $\qq\in\Rn$ and some $\VV\in\Rnk$,
$k\geq 0$.
By \Crefequiv{pr:icon-equiv}{pr:icon-equiv:c}{pr:icon-equiv:a},
$\ebar=\VV\omm$ is an icon, so indeed $\xbar=\ebar\plusl\qq$ for $\ebar\in\corezn$ and $\qq\in\Rn$.

To show uniqueness of $\xbar$'s iconic part, suppose
$\xbar=\ebar\plusl\qq=\ebar'\plusl\qq'$ for some
$\ebar,\ebar'\in\corezn$ and $\qq,\qq'\in\Rn$.
Then
by
\Crefequiv{pr:icon-equiv}{pr:icon-equiv:a}{pr:icon-equiv:c},
there are column-orthogonal
matrices $\VV\in\R^{n\times k}$ and $\VV'\in\R^{n\times k'}$,
for some $k,k'\geq 0$,
such that $\ebar=\VV\omm$ and $\ebar'=\VV'\omm$.
Further, by \Cref{thm:ast-rank-is-mat-rank}, $\VV$ and $\VV'$ must both have
matrix rank equal to $\xbar$'s astral rank.
Since both matrices are column-orthogonal, they must both have full
column rank, implying that $\xbar$'s astral rank must be equal to both
$k$ and $k'$, so $k=k'$.
Since $\xbar=\VV\omm\plusl\qq=\VV'\omm\plusl\qq'$, it then follows by
\Cref{lem:uniq-ortho-prefix}(\ref{lem:uniq-ortho-prefix:b})
that $\VV=\VV'$.
Hence,
$\ebar=\VV\omm=\VV'\omm=\ebar'$.%
\indexg{iconic part|)}%
\indexg{astral points!iconic decomposition of|)}%
\indexg{representations of astral points!icon and finite part@as icon and finite part|)}%
\indexg{decompositions, astral!icon and finite part@as icon and finite part|)}%
\indexg{icons!decomposition using|)}
\end{proof}

Here are additional properties of icons.
Note specifically that all astrons are icons.
Note also that \Cref{pr:mult-rep-by-scalar} is a special case of
part~(\ref{pr:i:8d}) below (when combined with
\Cref{pr:icon-equiv}\ref{pr:icon-equiv:a}\ref{pr:icon-equiv:c}).

\begin{proposition}  \label{pr:i:8}
  Let $\ebar\in\corezn$ be an icon.

  \begin{letter-compact}
  \item  \label{pr:i:8b}
    The only icon in $\Rn$ is $\zero$; that is, $\corezn\cap\Rn=\{\zero\}$.
  \item  \label{pr:i:8d}
    Let $\alpha\in\Rextstrictpos$ and $\zbar\in\extspace$.
    Then
    $\alpha (\ebar \plusl \zbar) = \ebar \plusl \alpha \zbar$.
    In particular,
    $\alpha\ebar=\ebar$.
  \item  \label{pr:i:8-infprod}
    $\limray{\xbar}$ is an icon for all $\xbar\in\extspace$.
    In particular, all astrons are icons.
  \item  \label{pr:i:8-matprod}
    $\A \ebar$ is an icon for all $\A\in\Rmn$,
    and $\alpha \ebar$ is an icon for all $\alpha\in\Rext$.
  \item  \label{pr:i:8-leftsum}
    The set of icons is closed under leftward addition; that is,
    if $\dbar$ is some icon in $\corezn$, then so is
    $\ebar\plusl\dbar$.
  \item  \label{pr:i:8e}
    The set of all icons, $\corezn$, is closed in $\extspace$.
  \end{letter-compact}
\end{proposition}

\begin{proof}
~
\begin{proof-parts}
\pfpart{Part~(\ref{pr:i:8b}):}
Let $\ee\in\Rn$.
Then $\ee$ is an icon if and only if
$\ee=\ee\plusl\ee=\ee+\ee$,
which, by subtracting $\ee$, holds if and only if $\ee=\zero$.

\pfpart{Part~(\ref{pr:i:8d}):}
We first argue that $\alpha\ebar=\ebar$.
For all $\uu\in\Rn$,
$(\alpha\ebar)\cdot\uu = \alpha(\ebar\cdot\uu) = \ebar\cdot\uu$
by \Cref{pr:scalar-prod-props}(\ref{pr:scalar-prod-props:a}),
and since
$\alpha>0$ and
$\ebar\cdot\uu\in\{-\infty,0,+\infty\}$ 
(by
\Cref{pr:icon-equiv}\ref{pr:icon-equiv:a}\ref{pr:icon-equiv:b}).
Therefore, $\alpha\ebar=\ebar$
(by \Cref{pr:i:4}).

Next, suppose $\alpha\in\Rstrictpos$.
Then
$\alpha(\ebar\plusl\zbar)=\alpha\ebar\plusl\alpha\zbar=\ebar\plusl\alpha\zbar$
by \Cref{pr:i:7}(\ref{pr:i:7b})
and the preceding argument.

Finally, for the case $\alpha=\oms$, we now have
\[
  \limray{(\ebar\plusl\zbar)}
  =
  \lim [t(\ebar\plusl\zbar)]
  =
  \lim [\ebar\plusl t\zbar]
  =
  \ebar\plusl\limray{\zbar}
\]
where the second equality applies the argument just given,
and the third is by
\Cref{pr:i:7}(\ref{pr:i:7f}).

\pfpart{Part~(\ref{pr:i:8-infprod}):}
Let $\xbar\in\extspace$ and let $\uu\in\Rn$.
Then $(\limray{\xbar})\cdot\uu\in\{-\infty,0,+\infty\}$
by \Cref{pr:astrons-exist}.
Therefore, $\limray{\xbar}$ is an icon by
\Crefequiv{pr:icon-equiv}{pr:icon-equiv:b}{pr:icon-equiv:a}.

\pfpart{Part~(\ref{pr:i:8-matprod}):}
Let $\A\in\Rmn$.
Since $\ebar$ is an icon,
by
\Crefequiv{pr:icon-equiv}{pr:icon-equiv:a}{pr:icon-equiv:c},
$\ebar=\VV\omm$ for some $\VV\in\Rnk$, $k\geq 0$.
Then $\A\ebar = \A(\VV\omm) = (\A\VV)\omm$.
Thus, $\A\ebar$ is also an icon, again by
\Crefequiv{pr:icon-equiv}{pr:icon-equiv:c}{pr:icon-equiv:a}.

If $\alpha\in\R$, then
$\alpha\ebar=\alpha(\Iden\ebar)=(\alpha\Iden)\ebar$ is an icon, as
just argued (where $\Iden$ is the $n\times n$ identity matrix).

Finally, $\limray{\ebar}$ and $(-\oms)\ebar=\limray{(-\ebar)}$ are
icons by part~(\ref{pr:i:8-infprod}).

\pfpart{Part~(\ref{pr:i:8-leftsum}):}
Let $\dbar\in\corezn$,
and let $\uu\in\Rn$.
Then $\ebar\cdot\uu$ and $\dbar\cdot\uu$ are both
in $\{-\infty,0,+\infty\}$, by
\Crefequiv{pr:icon-equiv}{pr:icon-equiv:a}{pr:icon-equiv:b}.
Combining with
\Cref{pr:i:6},
this implies that
$(\ebar\plusl\dbar)\cdot\uu=\ebar\cdot\uu\plusl\dbar\cdot\uu$
is in $\{-\infty,0,+\infty\}$ as well.
Therefore, $\ebar\plusl\dbar\in\corezn$, again by
\Crefequiv{pr:icon-equiv}{pr:icon-equiv:b}{pr:icon-equiv:a}.

\pfpart{Part~(\ref{pr:i:8e}):}
Let $\seq{\ebar_t}$ be a sequence in $\corezn$ with some limit $\xbar\in\eRn$. Let $\uu\in\Rn$.
Then by
\Crefequiv{pr:icon-equiv}{pr:icon-equiv:a}{pr:icon-equiv:b},
$\ebar_t\inprod\uu\in\set{-\infty,0,+\infty}$ for all $t$, and therefore also
$\xbar\inprod\uu=\lim(\ebar_t\inprod\uu)\in\set{-\infty,0,+\infty}$
(with convergence by \Cref{thm:i:1}\ref{thm:i:1c}).
Since this holds for all $\uu\in\Rn$,
we obtain $\xbar\in\corezn$, again by
\Crefequiv{pr:icon-equiv}{pr:icon-equiv:b}{pr:icon-equiv:a}.
Thus, $\corezn=\clcorezn$, so $\corezn$ is closed.
\qedhere
\end{proof-parts}
\end{proof}

\Cref{pr:i:8}(\ref{pr:i:8d}) shows in particular that for
any icon $\ebar\in\corezn$, the singleton set~$\{\ebar\}$ acts
somewhat like a cone
in the sense of being closed under multiplication by a
positive scalar.
Furthermore,
icons are the only points with this property
(since if $\xbar\in\extspace$ is not an icon, then
by definition and \Cref{pr:i:7}\ref{pr:i:7c},
$2\xbar=\xbar\plusl\xbar\neq\xbar$).
We will have much more to say about cones in astral space in
\Cref{sec:cones}.%
\indexg{icons|)}

\section{Galaxies and prefix structure}
\label{sec:galaxies}

\indexg{galaxies|(}%
\Cref{thm:icon-fin-decomp}
shows that
every astral point has a uniquely determined iconic part.
\indexg{galaxies!defined|(}%
As a result, astral space itself can be naturally partitioned into
disjoint sets called \emph{galaxies}
consisting of all points
with an identical iconic part,
that is, into sets
\[
\indexm{G e}{$\galax{\ebar}$}{galaxy}%
   \galax{\ebar}
      = \regBraces{ \ebar \plusl \qq :\: \qq\in\Rn }
      = \ebar \plusl \Rn,
\]
defined for each
\indexg{galaxies!defined|)}%
$\ebar\in\corezn$.
The set of finite vectors, $\Rn$, is the galaxy $\galax{\zero}$ with icon $\zero$.
\indexg{astral space!algebraic structure|(}%
In algebraic terms, each galaxy $\galax{\ebar}$ is a
commutative subgroup of $\eRn$ under leftward addition,
which acts like
vector addition on elements of the same galaxy; this is because
if $\ebar\in\corezn$ and $\qq,\qq'\in\Rn$ then
$(\ebar\plusl\qq)\plusl(\ebar\plusl\qq') = \ebar\plusl(\qq+\qq')$
using
\indexg{astral space!algebraic structure|)}%
\Cref{pr:i:7}.
We next explore
properties of galaxies, their
closures and relationships to one another, and how these relate to the
overall topological structure of astral space.

We first show that
the closure of a galaxy $\galaxd$ includes exactly those
astral points that can be written in the form $\ebar\plusl\zbar$, for
some $\zbar\in\extspace$.
For readability, we write the closure of galaxy
$\galaxd$ as $\galcld$, rather than
$\clbar{\galaxd}$.

\begin{proposition}  \label{pr:galaxy-closure}
  Let $\ebar\in\corezn$ be an icon.
  Then the closure of galaxy $\galaxd$ is
  \[
      \galcld
      = \regBraces{ \ebar \plusl \zbar :\: \zbar\in\extspace }
      = \ebar \plusl \extspace.
  \]
\end{proposition}

\begin{proof}
Let $\xbar\in\extspace$.
We prove $\xbar\in\galcld$
if and only if $\xbar\in\ebar\plusl\extspace$.

Suppose $\xbar\in\ebar\plusl\extspace$, meaning
$\xbar=\ebar\plusl\zbar$ for some $\zbar\in\extspace$.
Then by
\Cref{thm:i:1}(\ref{thm:i:1d}),
there exists a sequence $\seq{\zz_t}$ in $\Rn$ that converges to
$\zbar$.
By \Cref{pr:i:7}(\ref{pr:i:7f}),
$\ebar\plusl\zz_t\rightarrow\ebar\plusl\zbar=\xbar$.
Since each $\ebar\plusl\zz_t$ is in $\galaxd$, we must have
$\xbar\in\galcld$.

Conversely, suppose $\xbar\in\galcld$. Then there exists a sequence
$\seq{\xbar_t}$ of elements of $\galaxd$ converging to $\xbar$. 
For each $t$,
we can
write $\xbar_t=\ebar\plusl\zz_t$ for some $\zz_t\in\Rn$. By sequential
compactness
of $\eRn$, there exists a subsequence of $(\zz_t)$ that converges to
some element $\zbar\in\eRn$. Discarding all other sequence elements,
and the corresponding elements in the sequence $\seq{\xbar_t}$,
we then obtain $\xbar_t=\ebar\plusl\zz_t\to\ebar\plusl\zbar$.
Thus, $\xbar=\ebar\plusl\zbar\in\ebar\plusl\extspace$.
\end{proof}

A point $\xbar\in\extspace$ is in the closure $\galcld$ if and only
if, for all $\uu\in\Rn$, $\xbar\cdot\uu$ and $\ebar\cdot\uu$ are equal
whenever $\ebar\cdot\uu\neq 0$
(or equivalently, whenever $\ebar\cdot\uu$ is infinite):

\begin{proposition}   \label{pr:cl-gal-equiv}
  Let $\ebar\in\corezn$ be an icon, and let $\xbar\in\extspace$.
  Then the following are equivalent:
  \begin{letter-compact}
  \item   \label{pr:cl-gal-equiv:a}
    $\xbar\in\ebar\plusl\extspace$
    (that is, $\xbar\in\galcld$).
  \item   \label{pr:cl-gal-equiv:c}
    $\xbar=\ebar\plusl\xbar$.
  \item   \label{pr:cl-gal-equiv:b}
    For all $\uu\in\Rn$, if $\ebar\cdot\uu\neq 0$ then
    $\xbar\cdot\uu=\ebar\cdot\uu$.
  \end{letter-compact}
\end{proposition}

\begin{proof}
  ~

\begin{proof-parts}
\pfpart{%
  (\ref{pr:cl-gal-equiv:a})
  $\Rightarrow$
  (\ref{pr:cl-gal-equiv:b}):
}
Suppose $\xbar=\ebar\plusl\zbar$ for some $\zbar\in\extspace$.
Let $\uu\in\Rn$, and suppose $\ebar\cdot\uu\neq 0$.
Since $\ebar$ is an icon, we then must have
$\ebar\cdot\uu\in\{-\infty,+\infty\}$ by
\Crefequiv{pr:icon-equiv}{pr:icon-equiv:a}{pr:icon-equiv:b}.
Therefore, $\xbar\cdot\uu=\ebar\cdot\uu\plusl\zbar\cdot\uu=\ebar\cdot\uu$.

\pfpart{%
  (\ref{pr:cl-gal-equiv:b})
  $\Rightarrow$
  (\ref{pr:cl-gal-equiv:c}):
}
Suppose statement~(\ref{pr:cl-gal-equiv:b}) holds.
For all $\uu\in\Rn$,
we claim
\begin{equation}    \label{eq:pr:cl-gal-equiv:1}
  \xbar\cdot\uu=\ebar\cdot\uu \plusl \xbar\cdot\uu.
\end{equation}
This is immediate if $\ebar\cdot\uu= 0$.
Otherwise, if $\ebar\cdot\uu\neq 0$, then
$\ebar\cdot\uu\in\{-\infty,+\infty\}$ since $\ebar$ is an icon,
and also
$\xbar\cdot\uu=\ebar\cdot\uu$ by assumption.
Together, these imply \eqref{eq:pr:cl-gal-equiv:1}.
Therefore,
$\xbar=\ebar\plusl\xbar$
by \Cref{pr:i:4}.

\pfpart{%
  (\ref{pr:cl-gal-equiv:c})
  $\Rightarrow$
  (\ref{pr:cl-gal-equiv:a}):
}
This is immediate.
\qedhere
\end{proof-parts}
\end{proof}

\indexg{astral space!galactic structure|(}%
For any icons $\ebar,\,\ebar'\in\corezn$,
the galaxies
$\galaxd$ and $\galaxdp$ are disjoint, unless $\ebar=\ebar'$,
as noted already.
Nevertheless, the \emph{closures} of these galaxies might not be disjoint.
Indeed,
the next \namecref{pr:galaxy-closure-inclusion} shows that
if $\galcld$ and $\galcldp$ intersect at even a single point, then one
must fully contain the other.
\indexg{icons!prefix of|(}%
\indexg{prefix (of icon)|(}%
Furthermore, $\galcld$ is entirely included in $\galcldp$
if and only if $\ebar'$ is a \emph{prefix} of $\ebar$, meaning
$\ebar=\ebar'\plusl\dbar$ for some icon $\dbar\in\corezn$.
This can be restated equivalently
in terms of canonical representations
since
if $\ebar=\VV\omm$ and $\ebar'=\VV'\omm$
are canonical representations,
then
\Cref{lem:uniq-ortho-prefix}(\ref{lem:uniq-ortho-prefix:a})
implies
that $\ebar'$ is a prefix of $\ebar$ if and only if $\VV'$ is a prefix of $\VV$.%
\indexg{prefix (of icon)|)}%
\indexg{icons!prefix of|)}

\begin{theorem}  \label{pr:galaxy-closure-inclusion}
  Let $\ebar,\,\ebar'\in\corezn$ be icons.
  Then:
  \begin{letter}
  \item
  \label{pr:galaxy-closure-inclusion:a}
    $\galcld\subseteq\galcldp$
    if and only if
    $\ebar=\ebar'\plusl\dbar$ for some $\dbar\in\corezn$.
  \item
  \label{pr:galaxy-closure-inclusion:b}
    If $\galcld\cap\galcldp\neq\emptyset$
    then
    either
    $\galcld\subseteq\galcldp$
    or
    $\galcldp\subseteq\galcld$.
  \end{letter}
\end{theorem}

\begin{proof}
~
\begin{proof-parts}
\pfpart{Part~(\ref{pr:galaxy-closure-inclusion:a}):}
Suppose first that
$\ebar=\ebar'\plusl\dbar$ for some $\dbar\in\corezn$.
Then, by \Cref{pr:galaxy-closure},
\[
  \galcld
  =\ebar\plusl\eRn
  =\ebar'\plusl\dbar\plusl\eRn
  \subseteq\ebar'\plusl\eRn
  =\galcldp.
\]

Conversely, suppose now that
$\galcld\subseteq\galcldp$.
Then, in particular, $\ebar\in\galcldp$, so
$\ebar=\ebar'\plusl\zbar$ for some $\zbar\in\extspace$,
by \Cref{pr:galaxy-closure}.
By \Cref{thm:icon-fin-decomp},
$\zbar=\dbar\plusl\qq$ for some $\dbar\in\corezn$ and $\qq\in\Rn$,
so
$\ebar=(\ebar'\plusl\dbar)\plusl\qq$.
Since the iconic part of this point is unique, it follows that
$\ebar=\ebar'\plusl\dbar$
(by \Cref{thm:icon-fin-decomp}
and \Cref{pr:i:8}\ref{pr:i:8-leftsum}).

\pfpart{Part~(\ref{pr:galaxy-closure-inclusion:b}):}
Suppose there exists a point $\xbar$ in $\galcld\cap\galcldp$.
Then by \Cref{pr:galaxy-closure},
$\xbar=\ebar\plusl\zbar=\ebar'\plusl\zbar'$,
for some $\zbar,\,\zbar'\in\extspace$.

By
\Crefequiv{pr:icon-equiv}{pr:icon-equiv:a}{pr:icon-equiv:c},
$\ebar=\VV\omm$ and $\ebar'=\VV'\omm$
for some column-orthogonal matrices
$\VV\in\R^{n\times k}$ and
$\VV'\in\R^{n\times k'}$.
By \Cref{lem:uniq-ortho-prefix}, one of the matrices must be a prefix of the other.
Without loss of generality, assume that
$\VV'$ is a prefix of $\VV$, so $\VV=[\VV',\VV'']$ for some matrix $\VV''$, and thus by \Cref{prop:split},
\[
\ebar=[\VV',\VV'']\omm=\VV'\omm\plusl\VV''\omm=\ebar'\plusl\VV''\omm,
\]
so
$\galcld\subseteq\galcldp$
by part~(\ref{pr:galaxy-closure-inclusion:a}).
\qedhere
\end{proof-parts}
\end{proof}

\indexg{galaxies!tree structure of|(}%
As a result, we can arrange the galaxies in a directed rooted tree
capturing these inclusion relationships.
The vertices of the tree
correspond to icons in $\corezn$
(representing galaxies),
with $\zero$ as the root.
For all icons $\ebar\in\corezn$ and for all $\vv\in\Rn$,
an edge is directed from $\ebar$ to
$\ebar\plusl\limray{\vv}$ (which is also an icon,
by \Cref{pr:i:8}\ref{pr:i:8-leftsum}), unless this would
result in a self-loop, that is, unless
$\ebar=\ebar\plusl\limray{\vv}$.
That the resulting structure is a tree rather than a directed acyclic
graph can be proved using \Cref{lem:uniq-ortho-prefix}.
Equivalently, the tree can be formulated in terms of canonical
representations, so that, using
\Crefequiv{pr:icon-equiv}{pr:icon-equiv:a}{pr:icon-equiv:c},
for column-orthogonal matrices $\VV\in\R^{n\times k}$
and $\VV'\in\R^{n\times (k+1)}$,
an edge is directed from $\VV\omm$ to $\VV'\omm$ if and only if
$\VV$ is a prefix of $\VV'$.

Thus, the children of $\zero$, the root, are exactly the
nonzero astrons.
Their children are the icons of astral rank~2, and so on.
In general,
the depth (distance from the root) of every icon $\ebar$ is exactly
equal to its astral rank, so
the height of the tree is $n$.
Further, the tree captures all galactic inclusions since
a path exists from $\ebar$ to $\ebar'$ if and only if $\ebar$
is a prefix of $\ebar'$, that is,
if and only if
$\galcld\supseteq\galcldp$
(by
\Cref{pr:galaxy-closure-inclusion}\ref{pr:galaxy-closure-inclusion:a}).%
\indexg{galaxies!tree structure of|)}%
\indexg{astral space!galactic structure|)}%
\indexg{galaxies|)}

\indexg{commutativity (of astral points)|(}%
\indexg{astral points!commutativity of|(}%
\indexg{leftward addition!commutativity of|(}%
The prefix structure also plays a role in characterizing when two elements $\xbar,\ybar\in\eRn$ \emph{commute} with respect to leftward addition, that is, when they satisfy $\xbar\plusl\ybar=\ybar\plusl\xbar$:

\begin{theorem}
\label{prop:commute}
Let $\xbar,\ybar\in\eRn$. Let $\xbar=\ebar\plusl\qq$ and $\ybar=\ebar'\plusl\qq'$ for some $\ebar,\ebar'\in\corezn$ and $\qq,\qq'\in\Rn$. Then the following are equivalent:
\begin{letter-compact}
\item  \label{prop:commute:a}
  $\xbar\plusl\ybar=\ybar\plusl\xbar$.
\item  \label{prop:commute:b}
  Either $\ebar$ is a prefix of $\ebar'$, or $\ebar'$ is a prefix of $\ebar$.
\item  \label{prop:commute:c}
  For all $\uu\in\Rn$, $\xbar\inprod\uu$ and $\ybar\inprod\uu$ are summable.
\end{letter-compact}
\end{theorem}
\begin{proof}
~
\begin{proof-parts}

\pfpart{\RefsImplication{prop:commute:a}{prop:commute:b}:}
  Let $\zbar=\xbar\plusl\ybar=\ebar\plusl(\qq\plusl\ybar)$, so $\zbar\in\galcl{\ebar}$ by \Cref{pr:galaxy-closure}. Since also $\zbar=\ybar\plusl\xbar$, we similarly obtain $\zbar\in\galcl{\ebar'}$. Thus, $\galcl{\ebar}\cap\galcl{\ebar'}\ne\emptyset$, so by \Cref{pr:galaxy-closure-inclusion},
  either $\ebar$ is a prefix of $\ebar'$, or $\ebar'$ is a prefix of $\ebar$.

\pfpart{\RefsImplication{prop:commute:b}{prop:commute:c}:} Without loss of generality,
  assume that $\ebar$ is a prefix of $\ebar'$, so $\ebar'=\ebar\plusl\dbar$ for some $\dbar\in\corezn$. Thus, $\ybar=\ebar\plusl\zbar'$ with $\zbar'=\dbar\plusl\qq'$. Let $\uu\in\Rn$. We will show that $\xbar\inprod\uu$ and $\ybar\inprod\uu$ are summable.

  If $\ebar\inprod\uu=0$ then $\xbar\inprod\uu=\qq\inprod\uu\in\R$, so $\xbar\inprod\uu$ and $\ybar\inprod\uu$ are summable regardless of the value of $\ybar\inprod\uu$.
Otherwise, if $\ebar\inprod\uu\ne 0$ then
$\xbar\cdot\uu=\ebar\cdot\uu=\ybar\cdot\uu$ by
\Cref{pr:cl-gal-equiv}(\ref{pr:cl-gal-equiv:a},\ref{pr:cl-gal-equiv:b}).
Since $\xbar\inprod\uu$ and $\ybar\inprod\uu$ are equal, they are also summable, finishing the proof.

\pfpart{\RefsImplication{prop:commute:c}{prop:commute:a}:}
  Let $\uu\in\Rn$. Then
\[
  (\xbar\plusl\ybar)\inprod\uu=
  (\xbar\inprod\uu)\plusl(\ybar\inprod\uu)
  =
  (\ybar\inprod\uu)\plusl(\xbar\inprod\uu)
  =
  (\ybar\plusl\xbar)\inprod\uu,
\]
  where the middle equality follows by \Cref{pr:i:5}(\ref{pr:i:5c}). Since this holds
  for all $\uu\in\Rn$, we obtain $\xbar\plusl\ybar=\ybar\plusl\xbar$
  (by \Cref{pr:i:4}).
\qedhere
\end{proof-parts}
\end{proof}

Here are some consequences of commutativity.
The first of these generalizes
\Cref{pr:i:7}(\ref{pr:i:7g}):

\begin{proposition}   \label{pr:sum-seq-commuting-limits}
  Let $\xbar,\ybar\in\extspace$, and let
  $\seq{\xbar_t},\seq{\ybar_t}$ be sequences in $\extspace$
  with $\xbar_t\rightarrow\xbar$ and
  $\ybar_t\rightarrow\ybar$.
  Assume
  $\xbar\plusl\ybar=\ybar\plusl\xbar$.
  Then
  $\xbar_t\plusl\ybar_t \rightarrow\xbar\plusl\ybar$.
\end{proposition}

\begin{proof}
Let $\uu\in\Rn$.
Then by \Cref{thm:i:1}(\ref{thm:i:1c}),
$\xbar_t\cdot\uu\rightarrow\xbar\cdot\uu$
and
$\ybar_t\cdot\uu\rightarrow\ybar\cdot\uu$.
Also, by
\Cref{prop:commute}(\ref{prop:commute:a},\ref{prop:commute:c}),
$\xbar\cdot\uu$ and $\ybar\cdot\uu$ are summable.
By \Cref{prop:lim:eR}(\ref{i:lim:eR:sum}), this implies,
for all sufficiently large $t$, that
$\xbar_t\cdot\uu$ and $\ybar_t\cdot\uu$ are summable,
so that
$(\xbar_t\plusl\ybar_t)\cdot\uu = \xbar_t\cdot\uu + \ybar_t\cdot\uu$.
This latter expression converges to
$\xbar\cdot\uu + \ybar\cdot\uu = (\xbar\plusl\ybar)\cdot\uu$
(again by \Cref{prop:lim:eR}\ref{i:lim:eR:sum}).
Therefore,
$(\xbar_t\plusl\ybar_t)\cdot\uu \rightarrow (\xbar\plusl\ybar)\cdot\uu$.
Since this holds for all $\uu\in\Rn$,
the claim follows
by \Cref{thm:i:1}(\ref{thm:i:1c}).
\end{proof}

\Cref{pr:sum-seq-commuting-limits} need not hold if $\xbar$
and $\ybar$ do not commute.
For instance, in $\Rext$, if $x_t=t$ and $y_t=-t$ then
$x_t\rightarrow+\infty$,
$y_t\rightarrow-\infty$,
but
$x_t+y_t\rightarrow 0\neq (+\infty)\plusl(-\infty)$.

\indexg{matrix product, astral!commutativity and|(}%
We next use commutativity to
derive a sufficient condition for when $\A\xbar\plusl\B\xbar=(\A+\B)\xbar$.

\begin{proposition}
\label{prop:commute:AB}
Let $\xbar\in\eRn$, let $\A,\B\in\Rmn$, and assume that
$\A\xbar\plusl\B\xbar = \B\xbar\plusl\A\xbar$.
Then $\A\xbar\plusl\B\xbar=(\A+\B)\xbar$.
\end{proposition}

\begin{proof}
Let $\uu\in\Rm$.
Then by
\Cref{prop:commute}(\ref{prop:commute:a},\ref{prop:commute:c}),
$(\A\xbar)\cdot\uu$ and $(\B\xbar)\cdot\uu$ are summable. Therefore,
\begin{align*}
  (\A\xbar\plusl\B\xbar)\cdot\uu
  &=
  (\A\xbar)\cdot\uu + (\B\xbar)\cdot\uu
  \\
  &=
  \xbar\cdot\bigParens{\transAk\uu}
  +
  \xbar\cdot\bigParens{\transBk\uu}
  \\
  &=
  \xbar\cdot\bigParens{\transAk\uu + \transBk\uu}
  =
  \xbar\cdot\BigParens{\bigParens{\transA + \transB}\uu}
  =
  \BigParens{(\A+\B)\xbar}\cdot\uu.
\end{align*}
The first equality is by
Propositions~\ref{pr:i:6}
and~\ref{pr:i:5}(\ref{pr:i:5c}).
The second and fifth are by \Cref{thm:Ax-dot-u}.
The third is by
\Cref{pr:i:1}, since the terms are summable.

Since this holds for all $\uu\in\Rn$,
the claim follows
by \Cref{pr:i:4}.%
\indexg{matrix product, astral!commutativity and|)}%
\indexg{leftward addition!commutativity of|)}%
\indexg{commutativity (of astral points)|)}%
\indexg{astral points!commutativity of|)}
\end{proof}

\begin{figure}[t]
  \centering
  \includegraphics{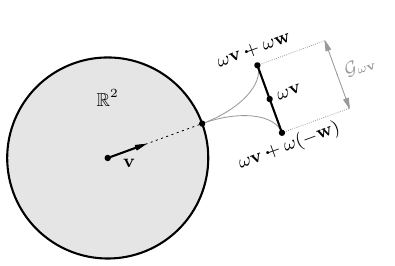}
  \mycaption{Galaxies associated with astrons in two dimensions}{%
    \indexf{galaxies!homeomorphic to Euclidean space}%
    \indexf{astrons!topological isolation of}%
    In two dimensions, any nonzero astron $\limray{\vv}$ is reached
    by following a ray in the direction of $\vv$ to infinity. This astron plays the role of
    an ``origin'' of the galaxy $\galax{\limray{\vv}}$
    (depicted as the segment on the right of the figure),
    which is homeomorphic to~$\R$,
    and of its closure $\galcl{\limray{\vv}}$, which is homeomorphic to $\eR$.
    If $\ww$
    is a vector orthogonal to $\vv$, then $\galax{\limray{\vv}}$ consists of all points of the form
    $\limray{\vv}\protect\plusl\lambda\ww$, for $\lambda\in\R$, and its closure additionally includes the
    points $\limray{\vv}\protect\plusl\limray{\ww}$ and $\limray{\vv}\protect\plusl\limray{(-\ww)}$.
    This figure also illustrates how each astron is topologically isolated from all other astrons
    (see \Cref{sec:not-second}).
    The zero astron~$\zero$ is surrounded by its galaxy~$\R^2$, depicted as an open disc.
    The boundary of~$\R^2$ is depicted as a circle, but it is not 
    a homeomorphic copy of a circle. Instead,
    if we ``zoom in,''  %
    each of the points on the circle
    actually contains a homeomorphic copy of $\eR$, corresponding to a closed galaxy~$\galcl{\limray{\vv}}$, with the astron
    $\limray{\vv}$ at its ``center.'' Thus, each of the infinite astrons
    $\limray{\vv}$
    is separated from the ``nearby'' infinite astrons by its own galaxy.
  }
  \label{fig:galaxy}%
\end{figure}

\begin{figure}[t]
  \centering
  \includegraphics{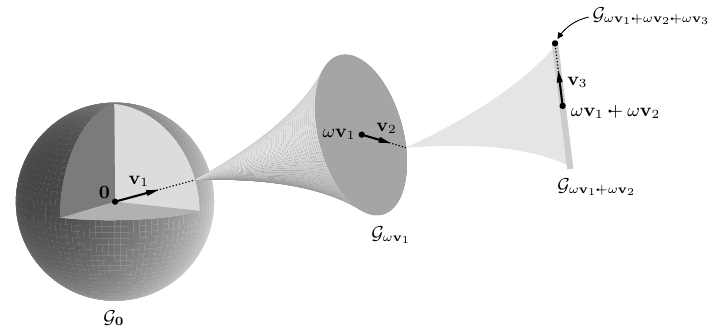}
  \mycaption{Galaxies in three dimensions}{%
    \indexf{galaxies!homeomorphic to Euclidean space}%
    There are four kinds of galaxies in three dimensions, corresponding to icons
    of rank $k=0,1,2,3$. The icon $\zero$ is the only icon of rank $k=0$. It is at the center of the galaxy $\galax{\zero}$,
    which coincides with $\R^3$. Picking any $\vv_1\ne\zero$ and following the ray in the direction
    of $\vv_1$ to infinity, we reach the icon $\limray{\vv_1}$ of rank $k=1$; the associated galaxy
    $\galax{\limray{\vv_1}}$ is homeomorphic to $\R^2$.
    Next, for any $\vv_2\ne\zero$, $\vv_2\perp\vv_1$, we can follow the sequence of points $\limray{\vv_1}\plusl t\vv_2$ to
    reach the icon $\limray{\vv_1}\plusl\limray{\vv_2}$ of rank $k=2$; the associated galaxy
    is homeomorphic to $\R^1$.
    Finally, for any $\vv_3\ne\zero$, $\vv_3\perp\vv_1$, $\vv_3\perp\vv_2$, we can follow the sequence of points
    $\limray{\vv_1}\plusl\limray{\vv_2}\plusl t\vv_3$ to reach the icon $\limray{\vv_1}\plusl\limray{\vv_2}\plusl\limray{\vv_3}$
    of rank $k=3$; the associated galaxy is a singleton set.}
  \label{fig:galaxy3d}%
\end{figure}

\indexg{galaxies!homeomorphic to Euclidean space|(}%
Returning to galaxies, we next consider their topology
and the topology of their closures.
The next theorem shows that the galaxy $\galaxd$ of an
icon $\ebar$ of astral rank $k$
is homeomorphic to $(n-k)$-dimensional Euclidean space,
$\R^{n-k}$, and
furthermore, the closure of that galaxy is
homeomorphic to
$(n-k)$-dimensional astral space,
$\extspac{n-k}$.
In other words, all galaxies and their closures
are topological copies of
lower-dimensional Euclidean spaces and astral spaces
(respectively). See Figures~\ref{fig:galaxy} and~\ref{fig:galaxy3d} for illustrations in two and three dimensions.

\begin{theorem}  \label{thm:galaxies-homeo}
  Let $\ebar\in\corezn$ be an icon of astral rank $k$.
  Then:
  \begin{letter-compact}
  \item  \label{thm:galaxies-homeo:a}
    $\galcld$ is homeomorphic to $\extspac{n-k}$.
  \item  \label{thm:galaxies-homeo:b}
    $\galaxd$ is homeomorphic to $\R^{n-k}$.
  \end{letter-compact}
\end{theorem}

\begin{proof}
~
\begin{proof-parts}
\pfpart{Part~(\ref{thm:galaxies-homeo:a}):}
By \Cref{thm:ast-rank-is-mat-rank} and
\Crefequiv{pr:icon-equiv}{pr:icon-equiv:a}{pr:icon-equiv:c},
$\ebar=\VV\omm$, for
some column-orthogonal matrix $\VV\in\R^{n\times k}$.
Let $\ww_1,\dotsc,\ww_{n-k}\in\Rn$ be an orthonormal basis for $(\colspace\VV)^\perp$,
and let $\WW=[\ww_1,\dotsc,\ww_{n-k}]$. We construct continuous
maps $\hogal: \galcld\rightarrow\extspac{n-k}$ and
$\hogalinv: \extspac{n-k}\rightarrow\galcld$ and show they are
inverses of each other, thus establishing that they are in fact homeomorphisms.

Specifically, we let
\begin{equation}
\label{eq:galaxies-homeo:0}
\hogal(\xbar) = \trans{\WW} \xbar,
\qquad
\hogalinv(\ybar)=\ebar\plusl \WW\ybar
\end{equation}
for $\xbar\in\galcld$ and $\ybar\in\extspac{n-k}$.
Both maps are affine and therefore continuous by \Cref{cor:aff-cont}.

Let $\xbar\in\galcld$, so $\xbar=\ebar\plusl\zbar$ for some $\zbar\in\extspace$
by \Cref{pr:galaxy-closure}. We will show that $\hogalinv(\hogal(\xbar))=\xbar$
and thus establish that $\hogalinv$ is a left inverse of $\hogal$.
Let $\PP=\WW\trans{\WW}$,
which is the projection matrix onto
$(\colspace\VV)^\perp$
(by \Cref{pr:basis-to-proj-mat}).
We calculate:
\begin{align*}
\notag
  \hogalinv(\hogal(\xbar))
=
  \ebar\plusl\WW\trans{\WW}(\ebar\plusl\zbar)
\notag
&=
  \ebar\plusl\WW(\trans{\WW}\ebar)\plusl\PP\zbar
\\
&=
  \ebar\plusl\PP\zbar
\notag
=
  \VV\omm\plusl\PP\zbar
=
  \VV\omm\plusl\zbar = \ebar\plusl\zbar = \xbar.
\end{align*}
The fifth equality is by the Projection Lemma~\ref{lemma:proj},
and the third is because
\begin{equation}
\label{eq:galaxies-homeo:3}
  \trans{\WW}\ebar=(\trans{\WW}\VV)\omm=\zeromat{(n-k)}{k}\omm=\zero
\end{equation}
since $\WW\perp\VV$.
Thus, $\hogalinv(\hogal(\xbar))=\xbar$,
proving that $\hogalinv$ is a left inverse of $\hogal$.

Next, let $\ybar\in\extspac{n-k}$.
Then
\[
  \hogal(\hogalinv(\ybar))
  =
  \trans{\WW}(\ebar\plusl\WW\ybar)
  =
  \trans{\WW}\ebar\plusl \trans{\WW}\WW\ybar
  =
  \ybar
\]
with the last equality following from
\eqref{eq:galaxies-homeo:3} and that $\trans{\WW}\WW=\Idn{n-k}$ since $\WW$ is column-orthogonal.
Thus, $\hogalinv$ is also a right inverse of $\hogal$, completing the proof.

\pfpart{Part~(\ref{thm:galaxies-homeo:b}):}
We redefine
$\hogal$ and $\hogalinv$ according to the same rules given in
\eqref{eq:galaxies-homeo:0}, but now with restricted
domain and range so that
$\hogal: \galaxd\rightarrow\R^{n-k}$
and
$\hogalinv: \R^{n-k}\rightarrow\galaxd$.
If $\yy\in\R^{n-k}$, then
$\hogalinv(\yy)=\ebar\plusl \WW\yy$, which is in $\galaxd$.
And
if $\xbar\in\galaxd$ then $\xbar=\ebar\plusl\zz$ for some
$\zz\in\Rn$, so
$\hogal(\xbar)=\trans{\WW}\ebar\plusl\trans{\WW}\zz=\trans{\WW}\zz$,
by
\eqref{eq:galaxies-homeo:3};
thus,
$\hogal(\xbar)$ is
in $\R^{n-k}$.\looseness=-1

Even with restricted domain and range,
the arguments used in
part~(\ref{thm:galaxies-homeo:a}) can be applied, proving that
$\hogal$ and $\hogalinv$ are continuous maps and
inverses of each other and therefore homeomorphisms.%
\indexg{galaxies!homeomorphic to Euclidean space|)}
\qedhere
\end{proof-parts}
\end{proof}

\section{Dominant directions}
\label{sec:dom:dir}

\indexg{dominant direction|(}%
When an infinite astral point
is decomposed into astrons and a finite part, as in
\Cref{cor:h:1},
the astrons are ordered by dominance, as previously discussed.
In this section, we focus on the first, most dominant,
astron appearing in such a representation,
as well as the natural decomposition of
any infinite astral point into that first astron and ``everything else.''
(Here, without loss of generality, we are only considering
representations in which all astrons are nonzero.)
This dominant astron has an important interpretation in terms of sequences:
We will see that for every sequence $\seq{\xx_t}$ in
$\Rn$ converging to an infinite point
$\xbar\in\extspace\setminus\Rn$, the directions of the sequence
elements, $\xx_t/\norm{\xx_t}$, must have a limit $\vv\in\Rn$
whose associated astron, $\limray{\vv}$, is always the first in any
representation of $\xbar$.

\indexg{direction|(}%
In more detail,
a \emph{direction} in $\Rn$ is represented by a unit vector, that is,
a vector $\vv\in\Rn$ with $\norm{\vv}=1$.
The direction of a nonzero vector $\xx\in\Rn$ is
\indexg{direction|)}%
$\xx/\norm{\xx}$.
An infinite astral point's dominant direction is the direction
associated with its first astron:

\begin{definition}
\indexg{dominant direction!defined|(}%
A vector $\vv\in\Rn$
is a \emph{dominant direction} of an infinite point
$\xbar\in\extspace\setminus\Rn$ if
$\norm{\vv}=1$ and if
$\xbar=\limray{\vv}\plusl\zbar$
for some $\zbar\in\extspace$.%
\indexg{dominant direction!defined|)}
\end{definition}

The next theorem proves the fact just mentioned, that
if $\seq{\xx_t}$ is any sequence in $\Rn$ converging to $\xbar$, then
the directions of the vectors $\xx_t$ must converge to
$\xbar$'s dominant direction.
Moreover, every infinite astral point has a unique dominant
direction.

As a preliminary step, we give the following proposition, to be
used in the proof, regarding the limit of $\norm{\xx_t}$ for an
astrally convergent sequence.

\begin{proposition}  \label{pr:seq-to-inf-has-inf-len}
  Let $\seq{\xx_t}$ be a sequence in $\Rn$ that converges to some
  point $\xbar\in\extspace$.
  Then
  \[
    \norm{\xx_t}
    \rightarrow
    \begin{cases}
                    \norm{\xbar} & \text{if $\xbar\in\Rn$,}\\
                    +\infty      & \text{otherwise.}
    \end{cases}
  \]
\end{proposition}

\begin{proof}
If $\xbar=\xx\in\Rn$, then $\xx_t\rightarrow\xx$, implying
$\norm{\xx_t}\rightarrow\norm{\xx}<+\infty$
by continuity.

Otherwise, if $\xbar\not\in\Rn$, then
$\xbar\cdot\uu\not\in\R$ for some $\uu\in\Rn$
(by \Cref{pr:i:3}\ref{i:3a}\ref{i:3b}),
implying
$|\xx_t\cdot\uu|\rightarrow|\xbar\cdot\uu|=+\infty$
(by \Cref{thm:i:1}\ref{thm:i:1c}).
Since $|\xx_t\cdot\uu|\leq\norm{\xx_t}\norm{\uu}$
by the Cauchy-Schwarz inequality,
this proves the claim in this case.
\end{proof}

\begin{theorem}  \label{thm:dom-dir}
  Let $\xbar\in\extspace\setminus\Rn$, and let $\vv\in\Rn$
  with $\norm{\vv}=1$.
  Also,
  let $\seq{\xx_t}$ and $\seq{\dd_t}$ be
  sequences in $\Rn$ such that $\xx_t\rightarrow\xbar$,
  and $\dd_t=\xx_t / \norm{\xx_t}$
  whenever $\xx_t\neq \zero$
  (or equivalently,
  $\xx_t =\norm{\xx_t}\dd_t$ for all $t$).
  Then the following are equivalent:
  \begin{letter-compact}
  \item  \label{thm:dom-dir:a}
    $\xbar=\limray{\vv}\plusl\zbar$ for some $\zbar\in\extspace$.
    That is, $\vv$ is a dominant direction of $\xbar$.
  \item  \label{thm:dom-dir:b}
    For all $\uu\in\Rn$,
    if $\vv\cdot\uu>0$ then $\xbar\cdot\uu=+\infty$,
    and
    if $\vv\cdot\uu<0$ then $\xbar\cdot\uu=-\infty$.
  \item  \label{thm:dom-dir:c}
    $\dd_t \rightarrow \vv$.
  \end{letter-compact}
  Furthermore, every point $\xbar\in\extspace\setminus\Rn$ has a
  unique dominant direction.
\end{theorem}

\begin{proof}
We begin by establishing the existence and uniqueness of the dominant direction for any
point $\xbar\in\extspace\setminus\Rn$, and then show the equivalence
of conditions (\ref{thm:dom-dir:a}), (\ref{thm:dom-dir:b}), (\ref{thm:dom-dir:c}).

\begin{proof-parts}
\pfpart{Existence:}
Being in $\extspace$, $\xbar$ must have the form given in
\Cref{cor:h:1}.
In the notation of that corollary, $k\geq 1$ (since
$\xbar\not\in\Rn$),
and
$\vv_1$ must be a dominant direction (assuming, without loss of
generality, that $\norm{\vv_1}=1$).

\pfpart{Uniqueness:}
Suppose both $\vv$ and $\vv'$ are dominant directions of
$\xbar$.
Then $\norm{\vv}=\norm{\vv'}=1$, and
$\xbar=\limray{\vv}\plusl\zbar=\limray{\vv'}\plusl\zbar'$,
for some $\zbar,\,\zbar'\in\extspace$.
That is,
$\VV\omm\plusl\zbar=\VV'\omm\plusl\zbar'$,
where $\VV=[\vv]$ and $\VV'=[\vv']$.
Since these matrices are column-orthogonal,
we can apply
\Cref{lem:uniq-ortho-prefix}(\ref{lem:uniq-ortho-prefix:b}),
yielding $\VV=\VV'$.
Thus, $\vv=\vv'$.

\pfpart{(\ref{thm:dom-dir:b}) $\Rightarrow$ (\ref{thm:dom-dir:a}):}
Suppose (\ref{thm:dom-dir:b}) holds.
Let $\uu\in\Rn$.
If $\limray{\vv}\cdot\uu=+\infty$ then $\vv\cdot\uu>0$ so
$\xbar\cdot\uu=+\infty$ by assumption.
Likewise,
if $\limray{\vv}\cdot\uu=-\infty$ then
$\xbar\cdot\uu=-\infty$.
Thus, $\xbar\cdot\uu=\limray{\vv}\cdot\uu$ whenever
$\limray{\vv}\cdot\uu\neq 0$.
Therefore, we can apply
\Crefequiv{pr:cl-gal-equiv}{pr:cl-gal-equiv:a}{pr:cl-gal-equiv:b}
with
$\ebar=\limray{\vv}$, implying (\ref{thm:dom-dir:a}).

\pfpart{(\ref{thm:dom-dir:c}) $\Rightarrow$ (\ref{thm:dom-dir:b}):}
Suppose (\ref{thm:dom-dir:c}) holds.
Let $\uu\in\Rn$, and suppose $\vv\cdot\uu>0$.
Then
$\xx_t\cdot\uu = \norm{\xx_t}(\dd_t\cdot\uu) \rightarrow +\infty$
since $\dd_t\cdot\uu\rightarrow \vv\cdot\uu >0$
and since $\norm{\xx_t}\rightarrow+\infty$
by \Cref{pr:seq-to-inf-has-inf-len}.
Therefore, $\xbar\cdot\uu=+\infty$,
by \Cref{thm:i:1}(\ref{thm:i:1c}),
since $\xx_t\rightarrow\xbar$.
By a symmetric argument, if $\vv\cdot\uu<0$ then
$\xbar\cdot\uu=-\infty$.

\pfpart{(\ref{thm:dom-dir:a}) $\Rightarrow$ (\ref{thm:dom-dir:c}):}
Suppose $\xbar=\limray{\vv}\plusl\zbar$ for some $\zbar\in\extspace$.
Then
by \Cref{thm:lim-plusl-inv}, $\xx_t=\lambda_t\vv_t+\zz_t$ for some span-bound sequence $\seq{\vv_t}$ in $\Rn$ and sequences $\seq{\zz_t}$ in $\Rn$ and $\seq{\lambda_t}$ in $\Rstrictpos$
with $\vv_t\to\vv$, $\zz_t\to\zbar$, $\lambda_t\to+\infty$, and
$\norm{\zz_t}/\lambda_t\rightarrow 0$.
Thus, $\zz_t/\lambda_t\to\zero$, so
\begin{equation}
\label{eq:dom-dir:1}
  \vv_t + \zz_t/\lambda_t \to \vv.
\end{equation}
Therefore, for all $t$ sufficiently large (so that $\xx_t\neq\zero$),
\[
  \dd_t
  =
  \frac{\xx_t}{\norm{\xx_t}}
  =
  \frac{\lambda_t(\vv_t+\zz_t/\lambda_t)}{\lambda_t\norm{\vv_t+\zz_t/\lambda_t}}
  =
  \frac{\vv_t+\zz_t/\lambda_t}{\norm{\vv_t+\zz_t/\lambda_t}}
  \to
  \vv,
\]
where the last step follows from \eqref{eq:dom-dir:1} and because
$\norm{\vv_t + \zz_t/\lambda_t} \to \norm{\vv}=1$ by continuity of the norm.
\qedhere
\end{proof-parts}
\end{proof}

Thus, every infinite astral point $\xbar\in\extspace\setminus\Rn$ has
a unique dominant direction associated with the first
(nonzero)
astron in any representation for
\indexg{dominant direction|)}%
$\xbar$.
We next turn our attention to the rest of $\xbar$'s representation,
the part following that first astron.

\indexg{dominant direction!decomposition using|(}%
\indexg{decompositions, astral!dominant direction@with dominant direction|(}%
\indexg{projection orthogonal to vector|(}%
Let $\vv$ be $\xbar$'s dominant direction.
By definition, this means that $\xbar=\limray{\vv}\plusl\zbar$ for
some $\zbar\in\extspace$, implying, by
\Crefequiv{pr:cl-gal-equiv}{pr:cl-gal-equiv:a}{pr:cl-gal-equiv:c},
that actually $\xbar=\limray{\vv}\plusl\xbar$.
Applying the Projection Lemma~\ref{lemma:proj},
it then follows that $\xbar=\limray{\vv}\plusl\PP\xbar$
where $\PP$ is the projection matrix onto the linear subspace
orthogonal to $\vv$.
Thus, in general, $\xbar$ can be decomposed as
\begin{equation}   \label{eqn:domdir-xperp-decomp}
   \xbar
   =
   \limray{\vv}\plusl\xbarperp,
\end{equation}
where $\xbarperp$,%
\indexm{x 700}{$\xbarperp$}{projection orthogonal to vector}
both here and in what follows, is shorthand
for $\PP\xbar$
whenever the vector $\vv$ (and hence also the matrix $\PP$) is clear
from context
(for any point $\xbar$, including points in~$\Rn$).\looseness=-1

This decomposition
will be used in many of our proofs, which are often by induction
on astral rank, since
$\xbarperp$ has lower astral rank than $\xbar$.
The next proposition provides
the basis for such proofs:

\begin{proposition} \label{pr:h:6}
  Let $\xbar\in\extspace\setminus\Rn$.
  Let $\vv$ be its dominant direction, and
  let $\xbarperp$ be $\xbar$'s projection
  orthogonal to $\vv$.
  Then $\xbar=\limray{\vv}\plusl\xbarperp$.
  Further, suppose $\xbar$'s canonical representation is
  $\xbar=[\vv_1,\dotsc,\vv_k]\omm\plusl\qq$,
  so that its astral rank is $k\ge 1$.
  Then $\xbarperp\negKern$'s canonical representation is
  $\xbarperp=[\vv_2,\dotsc,\vv_k]\omm\plusl\qq$,
  implying its astral rank is $k-1$.
\end{proposition}

\begin{proof}
That $\xbar=\limray{\vv}\plusl\xbarperp$ was argued above.

Let $\VV'=[\vv_2,\dotsc,\vv_k]$ so that
$\xbar=[\vv_1,\VV']\omm\plusl\qq$, and let $\PP$ be the projection
matrix onto the space orthogonal to $\vv_1$. Since $\vv$ is the unique
dominant direction of $\xbar$, we must have $\vv=\vv_1$. Thus,
\[
  \xbarperp
  =\PP\bigParens{[\vv_1,\VV']\omm\plusl\qq}
  =[\PP\vv_1,\PP\VV']\omm\plusl\PP\qq
  =[\zero_n,\VV']\omm\plusl\qq
  =\VV'\omm\plusl\qq,
\]
where the second equality follows from \Cref{pr:h:4}(\ref{pr:h:4c},\ref{pr:h:4d}),
the third from
\Cref{pr:proj-mat-props}(\ref{pr:proj-mat-props:d},\ref{pr:proj-mat-props:e})
since $\VV'\perp\vv_1$ and $\qq\perp\vv_1$, and the last
from \Cref{prop:split}.%
\indexg{dominant direction!decomposition using|)}%
\indexg{decompositions, astral!dominant direction@with dominant direction|)}
\end{proof}

\indexg{canonical representation!derived from sequence|(}%
Let $\seq{\xx_t}$ in $\Rn$ be any sequence converging to
$\xbar\in\extspace$, and let $\xbar$'s canonical representation be as
given in \Cref{pr:h:6}.
As we know from
\Crefequiv{thm:dom-dir}{thm:dom-dir:a}{thm:dom-dir:c},
if $\xbar\not\in\Rn$,
we can determine its dominant direction, $\vv_1$, from the
sequence as the limit of $\xx_t/\norm{\xx_t}$.
In fact, from any such sequence,
using the decomposition given in
\eqref{eqn:domdir-xperp-decomp}, we can determine the rest of
$\xbar$'s canonical representation as well.
To do so, we next project the sequence using $\PP$,
forming the sequence with elements
$\xperpt=\PP\xx_t$, which, by \Cref{thm:linear:cont}(\ref{thm:linear:cont:b}),
converges to $\xbarperp=\PP\xbar$.
From this sequence, we can determine $\xbarperp$'s dominant direction
(as the limit of the directions of the projected points $\xperpt$),
which is $\vv_2$ by \Cref{pr:h:6}, thus
yielding the second astron in $\xbar$'s canonical representation.
This process can be repeated to determine one astron after another,
until these repeated projections finally
yield a sequence converging to a finite point $\qq\in\Rn$.

Thus, any astral point $\xbar$'s entire canonical representation is
determined by this process for any sequence converging to $\xbar$.%
\indexg{canonical representation!derived from sequence|)}

As already noted,
the projections $\xbar\mapsto\xbarperp$ are special
cases of astral linear maps.
As such,
the next proposition summarizes some properties of those projections.

\begin{proposition} \label{pr:h:5}
  Let $\vv\in\Rn$,
  let $\xbar\in\extspace$, and let $\xbarperp$ be its projection
  orthogonal to~$\vv$.
  Then:
  \begin{letter-compact}
  \item  \label{pr:h:5a}
    $\xbarperp\cdot\uu=\xbar\cdot\uperp$ for all $\uu\in\Rn$.
  \item  \label{pr:h:5b}
    For any sequence $\seq{\xbar_t}$ in $\extspace$,
    if $\xbar_t\rightarrow\xbar$ then
    $\xbarperpt\rightarrow \xbarperp$.
  \item  \label{pr:h:5c}
    $(\xbar\plusl\ybar)^\bot = \xbarperp \plusl \ybarperp$ for $\ybar\in\extspace$.
  \item  \label{pr:h:5e}
    $(\alpha\ww)^\bot = \alpha(\ww^\bot)$ for $\ww\in\Rn$ and $\alpha\in\Rext$.
  \end{letter-compact}
\end{proposition}

\begin{proof}
Let $\PP$ be the projection matrix onto the linear space orthogonal to $\vv$.
Then for part~(\ref{pr:h:5a}), we have
$\xbarperp\cdot\uu=(\PP\xbar)\cdot\uu=\xbar\cdot(\PP\uu)=\xbar\cdot\uperp$
by \Cref{thm:Ax-dot-u} and since $\trans{\PP}=\PP$
(\Cref{pr:proj-mat-props}\ref{pr:proj-mat-props:a}).
The other parts follow respectively from
\Cref{thm:linear:cont}(\ref{thm:linear:cont:b})
and  \Cref{pr:h:4}(\ref{pr:h:4c},\ref{pr:h:4e}).%
\indexg{projection orthogonal to vector|)}
\end{proof}

\section{Convergence in direction}
\label{sec:conv-in-dir}

\indexg{astrons!sequence of|(}%
\indexg{sequence convergence!astrons@of astrons|(}%
Suppose $\seq{\vv_t}$ is a sequence in $\Rn$ converging to some point
$\vv\in\Rn$.
What can we say about the limit of the associated sequence of astrons
$\limray{\vv_t}$
(assuming, for this discussion, that the limit exists)?
We might be tempted to guess that
$\limray{\vv_t}\to\limray{\vv}$ since $\vv_t\to\vv$.
But this holds only in very degenerate cases.
This is because, as we saw in \Cref{sec:not-second},
astrons are isolated from one another,
so for any astron $\limray{\vv}$,
there exists a neighborhood of $\limray{\vv}$ that excludes
all other astrons (see \Cref{thm:formerly-lem:h:1:new}).
As a result, a sequence of astrons can converge to an astron $\limray{\vv}$
only if all but finitely many of the sequence elements are
acually equal to $\limray{\vv}$.

However, there are many examples of convergent squences of astrons that
do not follow this pattern. Since astrons are icons and the set
of icons is closed (\Cref{pr:i:8}\ref{pr:i:8e}),
such sequences must always converge to an icon.
Here is an example:

\begin{example}[Limits of astrons]
\label{ex:lim-astrons}
\indexg{Limits of astrons|(}%
In $\R^2$, let
$\vv_t = \ee_1 + (1/t) \ee_2$.
Then $\vv_t\rightarrow\ee_1$, but
$\limray{\vv_t}\not\rightarrow\limray{\ee_1}$.
Rather,
$\limray{\vv_t}\rightarrow\zbar$
where $\zbar=\limray{\ee_1}\plusl\limray{\ee_2}$.
To see this, let $\uu\in\Rn$.
If $\ee_1\cdot\uu>0$, then
$\vv_t\cdot\uu\rightarrow\ee_1\cdot\uu$, so for all
sufficiently large $t$, $\vv_t\cdot\uu>0$,
implying $\limray{\vv_t}\cdot\uu=+\infty=\zbar\cdot\uu$.
The case $\ee_1\cdot\uu<0$ is symmetric.
And if $\ee_1\cdot\uu=0$, then
$\vv_t\cdot\uu=(1/t)\ee_2\cdot\uu$,
so in this case,
$\limray{\vv_t}\cdot\uu=\limray{\ee_2}\cdot\uu=\zbar\cdot\uu$ for all $t$.
In every case, $\limray{\vv_t}\cdot\uu\rightarrow\zbar\cdot\uu$,
so $\limray{\vv_t}\rightarrow\zbar$
(by \Cref{thm:i:1}\ref{thm:i:1c}).%
\indexg{Limits of astrons|)}%
\end{example}

We can try to give some (very) informal intuition for what is
happening in this example.
A naive guess might have been that
$\seq{\limray{\vv_t}}$ would converge to $\limray{\ee_1}$,
since $\vv_t\rightarrow\ee_1$, but that
would be impossible, as already discussed.
We might then expect the sequence to get ``close'' to
$\limray{\ee_1}$, for instance, to converge to a point in its galaxy,
$\galax{\limray{\ee_1}}$, but this also does not happen
(and could not happen since, as already mentioned, the limit must be
an icon).
The galaxy
$\galax{\limray{\ee_1}}$
is topologically a line, consisting of all points
$\limray{\ee_1}\plusl \lambda\ee_2$, for $\lambda\in\R$.
The sequence includes points
$\limray{\vv_t}=\limray{(\ee_1+(1/t)\ee_2)}$ that
get ever ``closer'' to the galaxy, approaching it from
``above'' (that is, in the direction of $\ee_2$).
Eventually, the sequence converges to
$\limray{\ee_1}\plusl\limray{\ee_2}$, a point in the galaxy's closure
at its ``upper'' edge.
So the sequence approaches but never enters the galaxy, rather
converging to a point in its closure.

In this example, although the sequence does not converge to
$\limray{\ee_1}$, as might have been naively expected, we did see that the
sequence converges to a point $\zbar$ whose dominant direction is $\ee_1$,
that is, whose first astron is $\limray{\ee_1}$.
This kind of convergence holds always, and in much greater generality,
as shown in the next theorem.
For instance, if $\vv_t\rightarrow\vv$, as in the discussion above,
this theorem implies that $\limray{\vv_t}$ must converge to a point
of the form $\limray{\vv}\plusl\ybar$ for some $\ybar\in\extspace$.

\begin{theorem}   \label{thm:gen-dom-dir-converg}
\indexg{omega-multiples@$\oms$-multiples!sequence of|(}%
\indexg{sequence convergence!omega-multiples@of $\oms$-multiples|(}%
  Let $\seq{\vbar_t}$ be a sequence in $\extspace$ converging to some
  point $\vbar\in\extspace$.
  Let $\seq{\xbar_t}$ be a sequence with each element $\xbar_t$ in
  $\limray{\vbar_t}\plusl\extspace$ and which converges to some point
  $\xbar\in\extspace$.
  Then $\xbar\in\limray{\vbar}\plusl\extspace$.
\end{theorem}

\begin{proof}
We prove the theorem using
\Cref{pr:cl-gal-equiv}.
As such, let $\uu\in\Rn$.
If $\limray{\vbar}\cdot\uu=+\infty$ then
$\vbar\cdot\uu>0$.
Since $\vbar_t\cdot\uu\rightarrow\vbar\cdot\uu$
(by \Cref{thm:i:1}\ref{thm:i:1c}),
this implies that for all $t$ sufficiently large,
$\vbar_t\cdot\uu>0$, and so that
$\limray{\vbar_t}\cdot\uu=+\infty$,
implying $\xbar_t\cdot\uu=+\infty$
(by \Cref{pr:cl-gal-equiv}\ref{pr:cl-gal-equiv:a}\ref{pr:cl-gal-equiv:b}).
Since $\xbar_t\cdot\uu\rightarrow\xbar\cdot\uu$, it follows that
$\xbar\cdot\uu=+\infty$.
Likewise, if $\limray{\vbar}\cdot\uu=-\infty$ then
$\xbar\cdot\uu=-\infty$.
Since $\limray{\vbar}$ is an icon
(by \Cref{pr:i:8}\ref{pr:i:8-infprod}),
we conclude that $\xbar\cdot\uu=\limray{\vbar}\cdot\uu$ whenever
$\limray{\vbar}\cdot\uu\neq 0$.
By
\Crefequiv{pr:cl-gal-equiv}{pr:cl-gal-equiv:b}{pr:cl-gal-equiv:a},
this proves the theorem.%
\indexg{omega-multiples@$\oms$-multiples!sequence of|)}%
\indexg{sequence convergence!astrons@of astrons|)}%
\indexg{sequence convergence!omega-multiples@of $\oms$-multiples|)}%
\indexg{astrons!sequence of|)}
\end{proof}

In the special case that each $\vbar_t$ is equal to an icon
$\ebar_t\in\corezn$, \Cref{thm:gen-dom-dir-converg} implies
that if $\ebar_t\rightarrow\ebar$ for some $\ebar\in\corezn$
and each $\xbar_t\in\galcldt$ then their limit $\xbar$ is in
$\galcld$.

\indexg{sequence convergence!dominant directions@of dominant directions|(}%
\indexg{dominant direction!sequence of|(}%
Using \Cref{thm:gen-dom-dir-converg}, we can now prove
that the dominant direction of the limit of a sequence is identical to the limit of the dominant directions of sequence elements:

\begin{theorem}  \label{thm:dom-dirs-continuous}
  Let $\seq{\xbar_t}$ be a sequence in $\extspace\setminus\Rn$ that
  converges to some point $\xbar\in\extspace$
  (which cannot be in $\Rn$).
  For each $t$,
  let $\vv_t$ be the dominant direction of $\xbar_t$,
  and let $\vv\in\Rn$.
  Then
    $\vv$ is the dominant direction of $\xbar$
  if and only if
    $\vv_t\rightarrow\vv$.
\end{theorem}

\begin{proof}
First, since $\Rn$ is open in $\eRn$,
the set $\eRn\setminus\Rn$ is closed in $\eRn$, so $\xbar\in\eRn\setminus\Rn$.

\begin{proof-parts}
\pfpart{``If'' ($\Leftarrow$):}
Suppose $\vv_t\rightarrow\vv$.
Then
\Cref{thm:gen-dom-dir-converg}, applied with $\vbar_t=\vv_t$
and $\vbar=\vv$, implies that
$\xbar\in\limray{\vv}\plusl\extspace$, that is, that
$\vv$ is the dominant direction of $\xbar$.

\pfpart{``Only if'' ($\Rightarrow$):}
Suppose $\vv$ is the dominant direction of $\xbar$, and that,
contrary to the theorem's claim,
$\vv_t\not\rightarrow\vv$.
Then there exists a neighborhood $U\subseteq\Rn$ of $\vv$ that
excludes infinitely many $\vv_t$.
By discarding all other sequence elements, we can assume
$\vv_t\not\in U$
for all $t$.
Furthermore, because each $\vv_t$ is on the unit sphere in $\Rn$,
which is compact, the sequence must have a subsequence that converges
to some unit vector $\ww\in\Rn$.
By again discarding all other sequence elements, we can assume the
entire sequence converges so that $\vv_t\rightarrow\ww$, and
still $\xbar_t\rightarrow\xbar$.
Since $U$'s complement, $\Rn\setminus U$, is closed (in $\Rn$)
and includes each
$\vv_t$, it must also include their limit $\ww$, but not $\vv$;
thus, $\ww\neq \vv$.

Since $\vv_t\rightarrow\ww$, as argued above, $\ww$ must be the
dominant direction of $\xbar$.
However, by assumption, $\vv$ is also $\xbar$'s dominant direction,
a contradiction since $\xbar$'s dominant direction is unique
(by \Cref{thm:dom-dir}), but $\vv\neq\ww$.

Having reached a contradiction, we conclude that $\vv_t\rightarrow\vv$.%
\indexg{sequence convergence!dominant directions@of dominant directions|)}%
\indexg{dominant direction!sequence of|)}
\qedhere
\end{proof-parts}
\end{proof}

\section{Comparison to cosmic space}
\label{sec:cosmic}

\indexg{cosmic space|(}%
As mentioned earlier,
\idxroc\idxwets\citet[Section~3A]{rock_wets}
study a different compactification of
$\Rn$ called \emph{cosmic space},
also studied by \idxhandup\citet{hansen_dupin99}
under the name
\indexg{enlarged space}\emph{enlarged space}.
Here, we explore how cosmic space and astral space
are related.
As we will see, the previously introduced concepts of galaxies and
dominant directions play a key role in understanding this
relationship.

Cosmic space consists of $\Rn$ together with
\emph{direction points}, one for every ray from the origin.
Rockafellar\idxroc\ and Wets\idxwets\
denote such points as
$\rwdir{\vv}$,%
\indexm{dir v}{$\rwdir{\vv}$}{direction point}
for $\vv\in\Rn\wo\{\zero\}$.
Thus, \ndim cosmic space,
which they write as $\cosmspn$,%
\indexm{csmrn}{$\cosmspn$}{cosmic space}
is $\Rn$ together with all direction points:
\[
  \cosmspn
  =
  \Rn \cup \bigBraces{\rwdir{\vv} :\: \vv\in\Rn\wo\{\zero\}}.
\]
Note that for vectors $\vv,\vv'\in\Rn\wo\{\zero\}$,
$\rwdir{\vv}=\rwdir{\vv'}$
if and only if $\vv$ and $\vv'$ have the same direction
(so that $\vv'/\norm{\vv'}=\vv/\norm{\vv}$).

As discussed by
Rockafellar\idxroc\ and Wets,\idxwets\
the topology on cosmic space is homeomorphic to the closed unit
ball in $\Rn$,
with $\Rn$ itself mapped to the interior of the
ball and the direction points mapped to its surface.
Thus, we can picture cosmic space being formed by shrinking $\Rn$ down
to the open unit ball in $\Rn$, and then taking its closure so that
the points on the surface of the ball correspond exactly to the
``new'' direction points.\looseness=-1

In more detail, let
$\ucball=\set{\xx\in\Rn :\: \norm{\xx}\leq 1}$
be the closed unit ball in $\Rn$,
and let $\coshom:\cosmspn\rightarrow\ucball$ be defined,
for $\xrw\in\cosmspn$, by
\begin{equation*}
  \coshom(\xrw)
  =
  \begin{cases}
    \xx / (1+\norm{\xx}) & \text{if $\xrw=\xx\in\Rn$,}\\
    \vv / \norm{\vv}
      & \text{if $\xrw=\rwdir{\vv}$ for some $\vv\in\Rn\wo\{\zero\}$.}
  \end{cases}
\end{equation*}
This function is a bijection between $\cosmspn$ and $\ucball$,
and it also maps (bijectively) $\Rn$
onto $\ucball$'s interior, and the set of direction
points in $\cosmspn$ onto $\ucball$'s surface.
The topology on $\cosmspn$ is
defined to ensure
that $\coshom$ is a homeomorphism. Thus, a set
$U\subseteq\cosmspn$ is defined to be open if and only if its image, $\coshom(U)$,
is open in $\ucball$ (as a subspace of $\Rn$).

Cosmic space is compact (because its homeomophic image, the closed unit ball,
is compact) and is the closure of $\Rn$ (because the closed unit ball is the closure
of its interior).
Furthermore, $\Rn$, in its usual Euclidean topology,
is a topological subspace of cosmic space.
(This is because $\coshom$, if restricted to $\Rn$, can be shown to
define a
homeomorphism between $\Rn$ in the Euclidean topology and
$\intr\ucball$.)
Thus, cosmic space is a compactification of $\Rn$.%
\indexg{cosmic space|)}%

\indexg{cosmic space!compared to astral space|(}%
\indexg{astral space!cosmic space comparison|(}%
However, while both cosmic space and astral space are compactifications of $\Rn$,
their topological properties are quite different.
In particular, any direction point $\rwdir{\vv}$ is arbitrarily
close to other direction points in cosmic space in the sense that every neighborhood
of $\rwdir{\vv}$ includes other direction points
(since the set of direction points is homeomorphic to the unit sphere in~$\Rn$).
This is not the case for astrons, which are
the natural counterparts of direction points in astral space:
Just like direction points in cosmic space,
astrons arise as limits of rays in astral space,
but unlike direction points,
astrons are topologically isolated from one another,
with each possessing a neighborhood that excludes all other astrons
(\Cref{thm:formerly-lem:h:1:new}).

\indexg{linear functions, continuous extension of!cosmic space@in cosmic space|(}%
As another example,
we know from \Cref{thm:i:1}(\ref{thm:i:1c}) that
for every $\uu\in\Rn$,
the linear function $f(\xx)=\xx\cdot\uu$, for $\xx\in\Rn$,
can be extended continuously to astral space
to obtain the function
$\xbar\mapsto\xbar\cdot\uu$, for $\xbar\in\extspace$.
The same is \emph{not} true for cosmic space.
To see this, suppose $n\geq 2$ and $\uu\neq\zero$.
Let $\ww\in\Rn$ be any nonzero vector with
$\ww\perp\uu$,
and let
$\alpha\in\R$.
Consider the sequence of points
$\xx_t=\alpha \uu + t \ww$.
Then regardless of $\alpha$, this sequence converges to $\rwdir{\ww}$ in
$\cosmspn$
(since
$\coshom(\xx_t)=\xx_t/(1+\norm{\xx_t})\rightarrow\ww/\norm{\ww}=\coshom(\rwdir{\ww})$,
implying $\xx_t\rightarrow\rwdir{\ww}$).
On the other hand, $f(\xx_t)=\alpha \norm{\uu}^2$ for all $t$,
implying that also $\lim f(\xx_t)=\alpha\norm{\uu}^2$,
which is not independent of $\alpha$ (since $\uu\ne\zero$).
Thus, for different choices of $\alpha$, the sequence
$\seq{\xx_t}$ always has
the same limit $\rwdir{\ww}$ in cosmic space, but the function values $f(\xx_t)$
have different limits.
As a result, no extension of $f$ to cosmic space can
be continuous at $\rwdir{\ww}$.
Indeed,
this argument shows that \emph{no} linear function on $\Rn$
can be extended
continuously to cosmic space if $n\ge 2$, except for the identically zero
function.%
\indexg{linear functions, continuous extension of!cosmic space@in cosmic space|)}%

We next look at the particular relationship between the topologies on
cosmic space and astral space.
To explain this connection,
let us first define the map
$\quof:\extspace\rightarrow\cosmspn$ as follows:
For $\xx\in\Rn$, the map is simply the identity, so
$\quof(\xx)=\xx$.
For all other points $\xbar\in\extspace\setminus\Rn$,
we let $\quof(\xbar)=\rwdir{\vv}$
where $\vv\in\Rn$ is $\xbar$'s dominant direction
(which exists and is unique by \Cref{thm:dom-dir}).
Thus, $\quof$ maps all infinite points $\xbar\in\extspace$
with the same dominant direction $\vv$
to the same direction point $\rwdir{\vv}\in\cosmspn$.
In other words, $\quofinv(\rwdir{\vv})$ consists exactly of those
astral points of the form $\limray{\vv}\plusl\zbar$, for some
$\zbar\in\extspace$, which means
$\quofinv(\rwdir{\vv})$ is
$\galcl{\limray{\vv}}$, the closure of $\limray{\vv}$'s galaxy
(\Cref{pr:galaxy-closure}).
In this sense, applying the map $\quof$ to astral space causes every such set
$\galcl{\limray{\vv}}$ to ``collapse'' down to a single point, namely,
$\rwdir{\vv}$.

We claim that the topology on $\cosmspn$ inherited from $\extspace$ as a result of
this collapsing
operation is
precisely the cosmic topology defined earlier.
\indexg{quotient topology|(}%
\indexg{topology!quotient|(}%
Formally, as shown in the next theorem,
we are claiming
that $\quof$ is a \emph{quotient map}, meaning
that it is surjective, and that,
for all subsets $U\subseteq\cosmspn$,
$\quofinv(U)$ is open in $\extspace$ if and only if
$U$ is open in $\cosmspn$.
As a consequence of this property, the topology on $\cosmspn$ is said to be the
\indexg{quotient topology|)}%
\indexg{topology!quotient|)}%
\emph{quotient topology} induced by $\quof$, and
$\cosmspn$ is an instance of a \emph{quotient space} of $\extspace$.

\begin{theorem}
\label{thm:quotient-map}
  Let $\quof:\extspace\rightarrow\cosmspn$ be the map
  defined above.
  Then $\quof$ is a quotient map.
  Therefore, $\cosmspn$ is a quotient space of $\extspace$.
\end{theorem}

\begin{proof}
We first show that $\quof$ is continuous.
Let $\seq{\xbar_t}$ be a sequence in $\extspace$ that converges to
some $\xbar\in\extspace$.
We aim to show
$\quof(\xbar_t)\rightarrow\quof(\xbar)$ in $\cosmspn$
(which, by
\Cref{prop:first:properties}\ref{prop:first:cont},
is sufficient
for proving continuity since $\extspace$ is
first-countable).

If $\xbar=\xx\in\Rn$, then, because $\Rn$ is a neighborhood of $\xx$,
all but finitely many of the elements $\xbar_t$ must also be in $\Rn$.
Since $\quof$ is the identity function on $\Rn$, and since the
topologies on $\cosmspn$ and $\extspace$ are the same when restricted
to $\Rn$, the claim follows directly in this case.

Suppose then that $\xbar\not\in\Rn$, and therefore
has some dominant direction $\vv\in\Rn$  %
so that $\quof(\xbar)=\rwdir{\vv}$.
We consider cases based on the elements of the sequence
$\seq{\xbar_t}$.

First, suppose the sequence includes at most finitely many elements
not in $\Rn$.
By discarding %
such elements, which does not
affect the sequence's convergence properties, we can assume
that the entire sequence is in $\Rn$
so that
$\xbar_t=\xx_t\in\Rn$ for all $t$,
and $\xx_t\rightarrow\xbar$ in $\extspace$.
For each $t$, let $\dd_t=\xx_t/\norm{\xx_t}$ if $\xx_t\neq\zero$, and
otherwise let $\dd_t=\zero$.
Then
\[
   \coshom\bigParens{\quof(\xbar_t)}
   =
   \coshom(\xx_t)
   =
   \frac{\xx_t}{1+\norm{\xx_t}}
   =
   \dd_t \cdot \frac{\norm{\xx_t}}{1+\norm{\xx_t}}
   \rightarrow
   \vv
   =
   \coshom(\rwdir{\vv})
   =
   \coshom\bigParens{\quof(\xbar)},
\]
where convergence is because $\dd_t\rightarrow\vv$ by
\Crefequiv{thm:dom-dir}{thm:dom-dir:a}{thm:dom-dir:c},
and
$\norm{\xx_t}/(1+\norm{\xx_t})\rightarrow 1$
since $\norm{\xx_t}\rightarrow+\infty$
(by \Cref{pr:seq-to-inf-has-inf-len}).
Since $\coshom$ is a homeomorphism, it follows that
$\quof(\xbar_t)\rightarrow\quof(\xbar)$.

Next, suppose the sequence includes at most finitely many elements
that are in $\Rn$.
Then, as before, we can discard these so
that none of the sequence elements $\xbar_t$ are in~$\Rn$.
For each $t$,
let $\vv_t$ be $\xbar_t$'s dominant direction
so that $\quof(\xbar_t)=\rwdir{\vv_t}$.
We then have that
\[
   \coshom\bigParens{\quof(\xbar_t)}
   =
   \coshom(\rwdir{\vv_t})
   =
   \vv_t
   \rightarrow
   \vv
   =
   \coshom(\rwdir{\vv})
   =
   \coshom\bigParens{\quof(\xbar)},
\]
with convergence following from \Cref{thm:dom-dirs-continuous}.
Therefore, $\quof(\xbar_t)\rightarrow\quof(\xbar)$ in this case as
well.

If the sequence $\seq{\xbar_t}$ is a mix of infinitely many elements
in $\Rn$ and infinitely many elements not
in $\Rn$, then we can treat the two subsequences of elements
in or not in $\Rn$ separately.
The arguments above show that the images of each of these subsequences
under $\quof$ converge to $\quof(\xbar)$.
Therefore, the image of the entire sequence converges to
$\quof(\xbar)$.
(This is because, for any neighborhood $U$ of $\quof(\xbar)$ in
$\cosmspn$, all elements of each subsequence must eventually be in
$U$; therefore, all elements of the entire sequence must eventually be
in $U$.)

Thus, in all cases, $\quof(\xbar_t)\rightarrow\quof(\xbar)$.
Therefore, $\quof$ is continuous.

We can now show that $\quof$ is a quotient map:
First, $\quof$ is surjective since
$\quof(\xx)=\xx$ for all $\xx\in\Rn$
and
$\quof(\limray{\vv})=\rwdir{\vv}$ for all
$\vv\in\Rn\wo\{\zero\}$.

Let $U\subseteq\cosmspn$.
It remains to show that 
$U$ is open in $\cosmspn$
if and only if
$\quofinv(U)$ is open in $\extspace$.
If $U$ is open, then $\quofinv(U)$ is open, since $\quof$ is
continuous.
Conversely, if $\quofinv(U)$ is open, then
$\extspace\setminus\quofinv(U)$ is closed.
Its image is
$\quof(\extspace\setminus\quofinv(U))=(\cosmspn)\setminus U$,
since $\quof$ is surjective.
This image must be closed by
\Cref{pr:cont-from-compact}(\ref{pr:cont-from-compact:a})
since $p$ is continuous, $\extspace$ is compact, and
since $\cosmspn$ is Hausdorff (being homeomorphic to a subspace of
Euclidean space).
Therefore, $U$ is open, completing the proof.%
\indexg{cosmic space!compared to astral space|)}%
\indexg{astral space!cosmic space comparison|)}%
\end{proof}

\part{Extending Functions to Astral Space}
\label{part:extending-functions}

\chapter{Lower semicontinuous extension}
\label{sec:functions}

We are now ready to begin the study of functions
that have been extended to astral space.
We are especially motivated by the fundamental problem of minimizing a
convex function $f$ on $\Rn$.
In general, such a function might not be minimized at any
finite point in its domain and its minimizers might only be attained
at infinity, by following a sequence.
\indexg{lower semicontinuous extension|(}%
To study this situation within our framework,
we focus particularly on an extension $\fext$ of $f$ to $\extspace$,
which is constructed in such a way that $\ef$ is lower semicontinuous,
and so that $f$'s minimum over sequences in
$\Rn$ coincides with $\fext$'s minimum
over points in $\extspace$.
Further, $\fext$ always attains its minimum thanks to its lower
semicontinuity and the compactness of $\eRn$.
Much of the rest of this book studies $\fext$'s properties,
for example, its continuity and the structure of its minimizers.

\section{Definition and basic properties}
\label{sec:lsc:ext}

To begin,
we define the extension $\fext$, prove its lower semicontinuity, and derive
several related properties.

As discussed in \Cref{sec:prelim:lower-semicont},
a function $f:X\to\eR$ on a first-countable space~$X$ is
{lower semicontinuous} at $x\in X$ if
$f(x)\le\liminf f(x_t)$ whenever $x_t\to x$.
The function is lower semicontinuous if it is lower
semicontinuous at every point in~$X$.
Lower semicontinuous functions are precisely those whose epigraphs are closed in $X\times\R$
(\Cref{prop:lsc}\ref{prop:lsc:a}\ref{prop:lsc:b}).
On compact sets, such functions always attain their
minimum~(\Cref{thm:weierstrass}).

We also defined the lower semicontinuous hull of a function $f:X\to\eR$
(Eq.~\ref{eq:lsc:liminf:X:prelims})
and showed that its epigraph is the closure of $\epi f$ in $X\times\R$.
Specializing this definition to the case when $X=\Rn$,
we noted in \Cref{sec:prelim:lsc} that
the lower semicontinuous hull of a function $f:\Rn\rightarrow\Rext$
is given by
\begin{equation}
\label{eq:lsc:liminf}
 (\lsc f)(\xx) = \InfseqLiminf{\seq{\xx_t}}{\Rn}{\xx_t\rightarrow \xx}
                              {f(\xx_t)}
\end{equation}
for $\xx\in\Rn$.
This operation is particularly natural
for convex functions $f$ since they are already lower semicontinuous on the relative interior of $\dom{f}$, so only the function values on the relative boundary of $\dom{f}$ may need to be adjusted;
further,
the resulting function, $\lsc f$, remains convex
(\Cref{pr:lsc-props}\ref{pr:lsc-props:a}\ref{pr:lsc-props:b}).

To extend $f$ to astral space, we use the same idea but now
considering all sequences converging to points in $\extspace$, not
just $\Rn$:

\begin{definition}  \label{def:lsc-ext}
\indexg{lower semicontinuous extension!defined|(}%
The \emph{lower semicontinuous extension} of a function
$f:\Rn\rightarrow\Rext$
(or simply the \emph{extension} of $f$)
is the function $\fext:\eRn\to\eR$%
\indexm{f 200}{$\fext$}{lower semicontinuous extension}
 defined by
\begin{equation}
\label{eq:e:7}
\fext(\xbar) = \InfseqLiminf{\seq{\xx_t}}{\Rn}{\xx_t\rightarrow \xbar}
                            {f(\xx_t)},
\end{equation}
for all $\xbar\in\extspace$,
where (as usual for such notation)
the infimum is over all sequences $\seq{\xx_t}$ in $\Rn$
converging to $\xbar$.
(The set of such sequences is nonempty by
\indexg{lower semicontinuous extension!defined|)}%
\Cref{thm:i:1}\ref{thm:i:1d}.)
\end{definition}

Clearly,
$\ef(\xx)=(\lsc f)(\xx)$ for all $\xx\in\Rn$, but $\ef$ is also defined at points $\xbar\in\eRn\setminus\Rn$.
Here are some examples:

\begin{example}[Extension of an affine function]
\label{ex:ext-affine}
\mathtogether%
\indexg{lower semicontinuous extension!affine function@of affine function|(}%
\indexg{affine functions (standard)!extension of|(}%
Let $\uu\in\Rn$, $b\in\R$, and let $f:\Rn\rightarrow\R$ be defined by
$f(\xx)=\xx\cdot\uu+b$ for $\xx\in\Rn$.
In this simple case, \eqref{eq:e:7} yields
$\fext(\xbar)=\xbar\cdot\uu+b$ since if
$\seq{\xx_t}$ is any sequence in $\Rn$ that converges to $\xbar\in\extspace$,
then
$f(\xx_t)=\xx_t\cdot\uu+b\rightarrow\xbar\cdot\uu+b$ by continuity (\Cref{thm:i:1}\ref{thm:i:1c}).
In fact, $\ef$ is an instance of an astral affine map from \Cref{sec:linear-maps}.%
\indexg{lower semicontinuous extension!affine function@of affine function|)}%
\indexg{affine functions (standard)!extension of|)}
\end{example}

\begin{figure}
  \centering
  \includegraphics{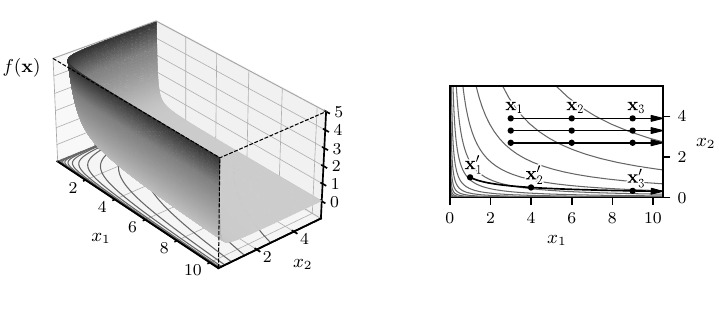}
  \mycaption{Product of hyperbolas}{%
    \indexf{Product of hyperbolas}%
    The function $f$ from \Cref{ex:recip-fcn-eg}. On the right,
    we show three example sequences of the form
    $\xx_t=\trans{[3t, \beta]}$ with $\beta\in\Rstrictpos$,
    converging
    to points of the form $\xbar=\limray{\ee_1}\plusl \beta\ee_2$ and satisfying $\lim f(\xx_t)=0$.
    We also show the sequence
    $\xx'_t=\trans{[t^2,1/t]}$, which converges 
    to $\xbar'=\limray{\ee_1}$ with $\lim f(\xx'_t)=0$.
    Since $f\ge 0$, these sequences demonstrate that $\ef(\xbar)=0$
    for all $\xbar=\limray{\ee_1}\plusl\beta\ee_2$ with $\beta\in\Rpos$.}
  \label{fig:prod-hyp}%
\end{figure}

\begin{example}[Product of hyperbolas]   \label{ex:recip-fcn-eg}
\indexg{Product of hyperbolas|(}%
Let $f:\R^2\rightarrow\Rext$ be defined, for $\xx\in\R^2$, as
\begin{equation}  \label{eqn:recip-fcn-eg}
  f(\xx) = f(x_1, x_2)
  =
\begin{cases}
    \dfrac{1}{x_1 x_2}
    & \text{if $x_1>0$ and $x_2>0$,}
\\[2ex]
    +\infty
    & \text{otherwise.}
\end{cases}
\end{equation}
(see \Cref{fig:prod-hyp}).
This function is convex, closed, proper and continuous everywhere.
Suppose $\beta\in\R$ and $\xbar=\limray{\ee_1}\plusl\beta\ee_2$
(where $\ee_1$ and $\ee_2$ are standard basis
\indexg{Product of hyperbolas|)}%
vectors).
\indexg{Product of hyperbolas!extension of|(}%
If $\beta>0$,
then $\fext(\xbar)=0$ since on any sequence $\seq{\xx_t}$ converging
to $\xbar$, the first component $\xx_t\cdot\ee_1=x_{t1}$ converges
to $\xbar\cdot\ee_1=+\infty$,
while the second component $\xx_t\cdot\ee_2=x_{t2}$
converges to $\xbar\cdot\ee_2=\beta>0$, implying that
$f(\xx_t)\rightarrow 0$ (see \Cref{fig:prod-hyp}, right).
If $\beta<0$, then a similar argument shows that
$\fext(\xbar)=+\infty$.
And if $\beta=0$, so that
$\xbar=\limray{\ee_1}$, then $\fext(\xbar)$ is again equal to $0$,
although more care is now needed in finding a sequence that shows this. One example
is the sequence $\xx_t=t^2\ee_1 + (1/t)\ee_2$, for which
$f(\xx_t)=1/t\rightarrow 0$ (see \Cref{fig:prod-hyp}, right).
This implies $\fext(\xbar)\leq 0$,
and since $f$ is nonnegative everywhere,
$\fext$ is as well, so $\fext(\xbar)=0$.%
\indexg{Product of hyperbolas!extension of|)}%
\end{example}

\indexg{trivial extension|(}%
Another way of extending a function $f$ on $\Rn$ to astral space,
which will turn out to be useful in proofs,
is to simply define the extended function to be
$+\infty$ at all infinite points:\looseness=-1

\begin{definition}  \label{def:triv-ext}
The \emph{trivial extension} of a function
$f:\Rn\rightarrow\Rext$
is the function $\ftriv:\eRn\to\eR$ defined, for $\xbar\in\eRn$, as
\begin{equation}   \label{eq:def:triv-ext:1}
\indexm{f 300}{$\ftriv$}{trivial extension}%
  \ftriv(\xbar)
  =
  \begin{cases}
    f(\xbar)    &  \text{if $\xbar\in\Rn$,}  \\
    +\infty     &  \text{otherwise.}
  \end{cases}
\end{equation}
\end{definition}
Except for being defined over all of $\extspace$,
the trivial extension $\ftriv$ is very much the same as the original
function $f$;
in particular, they both have the same epigraphs.
\indexg{lower semicontinuous extension!trivial extension and|(}%
The usefulness of a trivial extension comes from the following
proposition which shows that its lower semicontinuous hull,
$\lsc\ftriv$, is the same as $f$'s lower semicontinuous
extension,~$\fext$.
This will allow us to apply general facts about lower semicontinuous
hulls, as in \Cref{sec:prelim:lower-semicont},
to prove properties of $\fext$.

\begin{proposition}  \label{pr:lsc-ftriv-is-fext}
  Let $f:\Rn\rightarrow\Rext$ with trivial extension $\ftriv$ and
  lower semicontinuous extension $\fext$.
  Then $\fext=\lsc\ftriv$.
\end{proposition}

\begin{proof}
Let $\xbar\in\extspace$.
Suppose $\seq{\xx_t}$ is a sequence in $\Rn$ converging to $\xbar$.
Then
\[
  \liminf f(\xx_t)
  =
  \liminf \ftriv(\xx_t)
  \geq
  (\lsc\ftriv)(\xbar),
\]
with the inequality by definition of lower semicontinuous hull
(Eq.~\ref{eq:lsc:liminf:X:prelims}).
Since this holds for all such sequences, it follows that
$\fext(\xbar)\geq(\lsc\ftriv)(\xbar)$ by $\fext$'s definition
(Eq.~\ref{eq:e:7}).

For the reverse inequality,
$\fext(\xbar)\leq(\lsc\ftriv)(\xbar)$,
we assume $(\lsc\ftriv)(\xbar)<+\infty$ since the inequality is
immediate otherwise.
By \Cref{pr:lsc-seq-lim-exists},
there exists a sequence $\seq{\xbar_t}$ in $\extspace$ converging to
$\xbar$ with $\ftriv(\xbar_t)\rightarrow(\lsc\ftriv)(\xbar)$.
Since $(\lsc\ftriv)(\xbar)<+\infty$, this implies that
$\ftriv(\xbar_t)$ can be $+\infty$ for only finitely many sequence elements
$\xbar_t$; discarding these, we can assume henceforth that
$\ftriv(\xbar_t)<+\infty$ for all $t$.
By $\ftriv$'s definition, this implies, for all $t$, that
$\xbar_t\in\Rn$ and so that $\ftriv(\xbar_t)=f(\xbar_t)$.
Consequently,
\[
  (\lsc\ftriv)(\xbar)
  =
  \lim \ftriv(\xbar_t)
  =
  \lim f(\xbar_t)
  \geq
  \fext(\xbar),
\]  
with the inequality from $\fext$'s definition.%
\indexg{lower semicontinuous extension!trivial extension and|)}%
\indexg{trivial extension|)}
\end{proof}

The next theorem characterizes $\ef$ as the greatest lower
semicontinuous function that, when restricted to $\Rn$, is majorized by
$f$.
Moreover, the epigraph of $\ef$ is the closure of the epigraph of
$f$ in $\eRn\times\R$.

\begin{theorem}   \label{prop:ext:F}
  Let $f:\Rn\to\eR$.
  Then:
  \begin{letter-compact}
  \item   \label{prop:ext:F:a}
    $\fext$ is lower semicontinuous, and its restriction to $\Rn$ is
    majorized by $f$ (that is, $\resfcn{\fext}{\Rn}\leq f$).
  \item   \label{prop:ext:F:b}
    For every lower semicontinuous function $G:\extspace\to\eR$, if
    $\resfcn{G}{\Rn}\leq f$ then $G\leq\fext$.

  \item   \label{prop:ext:F:c}
\indexg{epigraph!extension@of extension|(}%
\indexg{lower semicontinuous extension!epigraph of|(}%
    $\epi\ef$ is the closure of $\epi f$ in $\eRn\times\R$.
  \end{letter-compact}
\end{theorem}

\begin{proof}
  ~

\begin{proof-parts}
\pfpart{Part~(\ref{prop:ext:F:a}):}
That $\fext$ is lower semicontinuous is immediate from
Propositions~\ref{prop:lsc:characterize}(\ref{prop:lsc:characterize:a})
and~\ref{pr:lsc-ftriv-is-fext}.
Further,
\[
  \resfcn{\fext}{\Rn}
  =
  \resfcn{(\lsc\ftriv)}{\Rn}
  \leq
  \resfcn{\ftriv}{\Rn}
  =
  f,
\]  
with the first equality from
\Cref{pr:lsc-ftriv-is-fext},
the inequality from
\Cref{prop:lsc:characterize}(\ref{prop:lsc:characterize:a}),
and the second equality from $\ftriv$'s definition.

\pfpart{Part~(\ref{prop:ext:F:b}):}
Let $G:\extspace\to\eR$ be lower semicontinuous with
$\resfcn{G}{\Rn}\leq f$.
From $\ftriv$'s definition (Eq.~\ref{eq:def:triv-ext:1}), this implies
that $G\leq \ftriv$, and so that
$G\leq\lsc\ftriv=\fext$
by
Propositions~\ref{prop:lsc:characterize}(\ref{prop:lsc:characterize:b})
and~\ref{pr:lsc-ftriv-is-fext}.

\pfpart{Part~(\ref{prop:ext:F:c}):}
By \Cref{pr:lsc-ftriv-is-fext},
$\epi\fext=\epi(\lsc\ftriv)$, which is equal to
the closure of $\epi\ftriv$ in $\extspace\times\R$ by
\Cref{prop:lsc:characterize}(\ref{prop:lsc:characterize:c}).
Since $\epi\ftriv=\epi f$ by $\ftriv$'s definition, this proves the
claim.%
\indexg{epigraph!extension@of extension|)}%
\indexg{lower semicontinuous extension!epigraph of|)}
\qedhere
\end{proof-parts}
\end{proof}

The next proposition shows that the infimum appearing in \eqref{eq:e:7}
must be realized by some sequence for which
$f(\xx_t)$ converges to $\fext(\xbar)$:

\begin{proposition}  \label{pr:d1}
  Let $f:\Rn\rightarrow\Rext$, and let $\xbar\in\extspace$.
  Then there exists a sequence $\seq{\xx_t}$ in $\Rn$ such that
  $\xx_t\rightarrow\xbar$ and $f(\xx_t)\rightarrow \fext(\xbar)$.
\end{proposition}

\begin{proof}
If $\fext(\xbar)=+\infty$ then let $\seq{\xx_t}$ be any sequence in
$\Rn$ converging to $\xbar$
(which exists by \Cref{thm:i:1}\ref{thm:i:1d}).
Then, by $\fext$'s definition,
$\liminf f(\xx_t)\geq\fext(\xbar)=+\infty$,
implying $f(\xx_t)\rightarrow+\infty$ and proving the claim in this
case.

Otherwise, $\fext(\xbar)<+\infty$.
Then by \Cref{pr:lsc-seq-lim-exists}, there exists a sequence
$\seq{\xbar_t}$ in $\extspace$ such that $\xbar_t\rightarrow\xbar$ and
$\ftriv(\xbar_t)\rightarrow (\lsc\ftriv)(\xbar) = \fext(\xbar)$,
with equality from \Cref{pr:lsc-ftriv-is-fext}.
Since $\fext(\xbar)<+\infty$, $\ftriv(\xbar_t)$ can be $+\infty$ for
only finitely many elements~$\xbar_t$; discarding these, we can assume
henceforth that $\ftriv(\xbar_t)<+\infty$ for all $t$.
By $\ftriv$'s definition,
this implies that each $\xbar_t\in\Rn$ and that
$\ftriv(\xbar_t)=f(\xbar_t)$.
This then proves the claim.
\end{proof}

Here are some other basic properties of $\fext$:

\begin{proposition}  \label{pr:h:1}
  Let $f:\Rn\rightarrow\Rext$.
  Then the following hold for $f$'s extension, $\fext$:
  \begin{letter-compact}
  \item  \label{pr:h:1a}
    For all $\xx\in\Rn$,
    $\fext(\xx)=(\lsc f)(\xx)\leq f(\xx)$.
    Thus, $\fext(\xx)=f(\xx)$ if $f$ is already lower
    semicontinuous at $\xx$.
  \item  \label{pr:h:1:geq}
    Let $g:\Rn\rightarrow\Rext$.
    If $f\geq g$ then $\fext\geq\gext$.
  \item  \label{pr:h:1aa}
    The extension of $f$ is the same as that of its lower
    semicontinuous hull.
    That is, $\fext=\lscfext$.
  \item  \label{pr:h:1b-neigh}
    Let $U\subseteq\extspace$ be a neighborhood of some
    $\xbar\in\extspace$, and
    let $V\subseteq\Rext$ be a neighborhood (in $\Rext$) of
    $\fext(\xbar)$.
    Then there exists $\xx\in U\cap \Rn$ with $f(\xx)\in V$.
  \item  \label{pr:h:1c}
    In $\extspace$, the closures of the effective domains of $f$ and
    $\fext$ are identical.
    That is,
    $\cldom{f}=\cldomfext$.
  \end{letter-compact}
\end{proposition}

\begin{proof}
~
\begin{proof-parts}
\pfpart{Part~(\ref{pr:h:1a}):}
The equality follows from the definitions of $\lsc f$ and $\ef$ (Eqs.~\ref{eq:lsc:liminf} and~\ref{eq:e:7}). 
The inequality was proved in \Cref{prop:ext:F}(\ref{prop:ext:F:a}).

\pfpart{Part~(\ref{pr:h:1:geq}):}
Let $\xbar\in\extspace$, and suppose $f\geq g$.
Let $\seq{\xx_t}$ be any sequence converging to $\xbar$.
Then
\[
  \liminf f(\xx_t)
  \geq
  \liminf g(\xx_t)
  \geq
  \gext(\xbar),
\]
where the last inequality is by definition of $\gext$
(\Cref{def:lsc-ext}).
Since this holds for all such sequences, the claim follows
(by definition of $\fext$).

\pfpart{Part~(\ref{pr:h:1aa}):}
We need to show that $\ef=\eg$ where $g=\lsc f$. Since $g\le f$, we
have $\eg\le\ef$ (by parts~(\ref{pr:h:1a}) and~(\ref{pr:h:1:geq})).
To show that $\ef\le\eg$, note that $\ef$ is lower semicontinuous (by \Cref{prop:ext:F}\ref{prop:ext:F:a}) and majorized by $g$ on $\Rn$ by part~(\ref{pr:h:1a}).
Therefore, it must be majorized by $\eg$ (by \Cref{prop:ext:F}\ref{prop:ext:F:b} applied to $g$).

\pfpart{Part~(\ref{pr:h:1b-neigh}):}
By \Cref{pr:d1}, there exists a sequence $\seq{\xx_t}$ in $\Rn$
converging to $\xbar$ with $f(\xx_t)\rightarrow\fext(\xbar)$.
Hence, for all $t$ sufficiently large,
$\xx_t\in U$, since $U$ is a neighborhood of $\xbar$,
and $f(\xx_t)\in V$,
since $V$ is a neighborhood of $\fext(\xbar)$.
Therefore, the claim holds for any $\xx_t$ with a sufficiently large
index $t$.

\pfpart{Part~(\ref{pr:h:1c}):}
By part~(\ref{pr:h:1a}),
$\fext(\xx)\leq f(\xx)$ for all $\xx\in\Rn$.
Therefore,
$\dom f \subseteq \dom \fext$,
implying
$\cldom{f} \subseteq \cldomfext$.

For the reverse inclusion, suppose $\xbar\in\dom{\fext}$.
By \Cref{pr:d1}, there exists a sequence $\seq{\xx_t}$ in $\Rn$
converging to $\xbar$ with $f(\xx_t)\rightarrow\fext(\xbar)$.
Since $\fext(\xbar)<+\infty$, this implies that $f(\xx_t)$ can be
$+\infty$ for finitely many sequence elements $\xx_t$;
discarding these, we can assume $\xx_t\in\dom{f}$ for all $t$.
Thus, $\xbar\in\cldom{f}$.
Since $\dom{\fext}\subseteq\cldom{f}$, it follows that
$\cldomfext\subseteq\cldom{f}$
since $\cldom{f}$ is closed (in $\extspace$).
\qedhere
\end{proof-parts}
\end{proof}

\indexg{lower semicontinuous extension!minimum attained|(}%
Because astral space is compact and $\fext$ is lower semicontinuous,
the minimum of $\fext$ is always realized at some point
$\xbar\in\extspace$. The next proposition shows
that this minimum is equal to the infimum of $f$.
Thus, minimizing a function $f$ on $\Rn$ is equivalent to finding a
minimizer of its extension $\fext$:

\begin{proposition}  \label{pr:fext-min-exists}
  Let $f:\Rn\rightarrow\Rext$.
  Then $\fext$ attains its minimum at some point in $\extspace$.
  Moreover, $\min\fext=\inf f$.
\end{proposition}

\begin{proof}
By \Cref{pr:lsc-ftriv-is-fext}, $\fext=\lsc\ftriv$,
implying, by \Cref{pr:lsc-min-exists-equals-sup},
that $\fext$ attains its minimum.
These propositions further show that
$\min\fext=\min(\lsc\ftriv)=\inf\ftriv=\inf f$,
with the last equality from $\ftriv$'s definition.%
\indexg{lower semicontinuous extension!minimum attained|)}%
\indexg{lower semicontinuous extension|)}
\end{proof}

\indexg{indicator functions (astral)|(}%
Analogous to the standard indicator function $\inds$ defined in
\eqref{eq:indf-defn},
for an astral set $S\subseteq\eRn$, we define
the \emph{astral indicator function}
$\indfa{S}:\extspace\rightarrow\Rext$, for $\xbar\in\eRn$, as
\begin{equation} \label{eq:indfa-defn}
\indexm{i s 600}{$\indfa{S}$}{astral indicator function}
  \indfa{S}(\xbar) =
  \begin{cases}
    0 & \text{if $\xbar\in S$,} \\
    +\infty
      & \text{otherwise.}
  \end{cases}
\end{equation}
Also, as with standard indicator functions, for a single point
$\zbar\in\extspace$,
we sometimes write $\indfa{\zbar}$ as shorthand for
$\indfa{\{\zbar\}}$,
and when working in $\Rext$, we sometimes write
$\indfa{\leq\beta}$ for $\indfa{[-\infty,\beta]}$, where $\beta\in\Rext$
(and likewise for other inequality relations).

\indexg{lower semicontinuous extension!indicator function@of indicator function|(}%
\indexg{indicator functions (standard)!extension of|(}%
As shown in the next proposition,
for $S\subseteq\Rn$,
the extension of a standard indicator function $\inds$ is the astral
indicator function for $\Sbar$, the closure of $S$ in~$\extspace$:
\begin{proposition}  \label{pr:inds-ext}
  Let $S\subseteq\Rn$.
  Then $\indsext=\indfa{\Sbar}$.
\end{proposition}
\begin{proof}
By \Cref{prop:ext:F}(\ref{prop:ext:F:c}), the epigraph of $\indsext$ is the closure of
$\epi\inds$ in $\eRn\times\R$. Since $\epi\inds=S\times\Rpos$, and the
closure of a product is the product of the closures
(\Cref{pr:prod-top-props}\ref{pr:prod-top-props:d}),
we obtain $\epi\indsext=\Sbar\times\Rpos=\epi{\indfa{\Sbar}}$.%
\indexg{indicator functions (standard)!extension of|)}%
\indexg{lower semicontinuous extension!indicator function@of indicator function|)}%
\indexg{indicator functions (astral)|)}
\end{proof}

\section{Continuity}

\indexg{continuity of extensions|(}%
We next look more closely at what continuity means for an extension
$\fext$.
Although $\fext$ must always be lower semicontinuous
(\Cref{prop:ext:F}\ref{prop:ext:F:a}), it need not be continuous everywhere.
Later, we will exactly characterize at what points $\fext$ is
continuous and when it is continuous everywhere
(see especially \Cref{sec:continuity}).

Because $\extspace$ is first-countable, an astral function
$F:\extspace\rightarrow\Rext$, including an extension $\fext$,
is {continuous} at some point $\xbar$ if and only if, for every sequence of
points $\seq{\xbar_t}$ in $\extspace$, if $\xbar_t\rightarrow\xbar$ then
$F(\xbar_t)\rightarrow F(\xbar)$
\indexg{continuity of extensions|)}%
(\Cref{prop:first:properties}\ref{prop:first:cont}).
For extensions, there is also another natural notion of
continuity that we now define:

\begin{definition}  \label{dfn:extens-cont}
\indexg{extensible continuity|(}%
\indexg{extensible continuity!defined|(}%
  Let $f:\Rn\rightarrow\Rext$, and let $\xbar\in\extspace$.
  We say that $f$ is \emph{extensibly continuous} at $\xbar$ if
  for every sequence $\seq{\xx_t}$ in $\Rn$,
  if $\xx_t\rightarrow\xbar$
  then $f(\xx_t)\rightarrow\fext(\xbar)$.%
\indexg{extensible continuity!defined|)}
\end{definition}
Thus, extensible continuity only involves sequences in $\Rn$
and the values on this sequence of the original function $f$
(rather than sequences in $\extspace$ and the values of $\fext$,
as is the case for the more general notion of continuity applied to
$\fext$).

\indexg{continuity of extensions!extensible continuity and|(}%
The next \namecref{thm:ext-cont-f}
shows in part that extensible continuity, which is often easier to
work with, implies ordinary continuity of $\fext$ at $\xbar$.
For instance, we earlier argued that the function $f$
in \Cref{ex:recip-fcn-eg}
is extensibly continuous at
$\xbar=\limray{\ee_1}\plusl\beta\ee_2$ if $\beta\neq 0$,
implying that $\fext$ is continuous at this same point.

In general, continuity of $\fext$ at a point~$\xbar$ need not imply
extensible continuity of $f$ at~$\xbar$.
Here is an example:

\begin{example}
Let $f:\R\to\R$ be the nonconvex function defined, for $x\in\R$, by
\[
  f(x)=
  \begin{cases}
    1
      & \text{if $x=0$,}
    \\
    0
      & \text{otherwise.}
  \end{cases}
\]
Then $\fext\equiv 0$, which is continuous everywhere, but
$f$ is not extensibly continuous at $0$ since, for instance,
the sequence with elements $x_t=0$ converges trivially to $0$,
but $\lim f(x_t) = 1 \neq \fext(0)$.
\end{example}

Nevertheless, in general, as shown further in the next
\namecref{thm:ext-cont-f},
for a function $f$ that is either convex or lower semicontinuous,
extensible continuity of $f$ at $\xbar$ is equivalent to
ordinary continuity of $\fext$ at $\xbar$:

\begin{theorem}  \label{thm:ext-cont-f}
  Let $f:\Rn\rightarrow\Rext$, and let $\xbar\in\extspace$.
  \begin{letter-compact}
  \item  \label{thm:ext-cont-f:a}
    If $f$ is extensibly continuous at $\xbar$ then
    $\fext$ is continuous at $\xbar$.
  \item  \label{thm:ext-cont-f:b}
    Suppose $f$ is either convex or lower semicontinuous (or both).
    Then
    $f$ is extensibly continuous at $\xbar$ if and only if
    $\fext$ is continuous at $\xbar$.
  \end{letter-compact}
\end{theorem}

\begin{proof}
~
\begin{proof-parts}

\pfpart{Part~(\ref{thm:ext-cont-f:a}):}
Suppose $f$ is extensibly continuous at $\xbar$.
Let $\seq{\xbar_t}$ be any sequence in $\extspace$ converging to
$\xbar$.
Since $\ef$ is lower semicontinuous (\Cref{prop:ext:F}\ref{prop:ext:F:a}),
$\liminf \fext(\xbar_t) \geq \fext(\xbar)$.
So, to prove continuity of $\fext$ at $\xbar$, we only need
to show
\begin{equation}  \label{eqn:thm:ext-cont-f:1}
  \limsup \fext(\xbar_t) \leq \fext(\xbar) .
\end{equation}
Assume that $\limsup\fext(\xbar_t)>-\infty$ (since otherwise
Eq.~\ref{eqn:thm:ext-cont-f:1}
is immediate), and let ${\beta\in\R}$ be such that
$\limsup \fext(\xbar_t) > \beta$.
Let $\countset{B}$ be a nested countable neighborhood base
for $\xbar$.
For each $t$,
there must exist some
$s$ with $\xbar_{s}\in B_t$ and $\fext(\xbar_{s})>\beta$
(since
all but finitely many of the sequence elements $\xbar_s$
must be included in $B_t$,
and since
$\fext(\xbar_s)>\beta$ for infinitely many values of
$s$).
Therefore, by
\Cref{pr:h:1}(\ref{pr:h:1b-neigh})
(applied with $U=B_t$ and $V=(\beta,+\infty]$),
there exists
$\xx_t\in B_t\cap\Rn$ with $f(\xx_t)>\beta$.

By \Cref{prop:nested:limit}, the resulting sequence
$\seq{\xx_t}$ converges to $\xbar$.
Therefore, by extensible continuity,
$\fext(\xbar)=\lim f(\xx_t) \geq \beta$,
since $f(\xx_t)>\beta$ for all $t$.
Since this holds for all $\beta<\limsup \fext(\xbar_t)$,
this proves \eqref{eqn:thm:ext-cont-f:1}, completing the proof.

\pfpart{Part~(\ref{thm:ext-cont-f:b}):}
Suppose that $f$ is either convex or lower semicontinuous.
In light of part~(\ref{thm:ext-cont-f:a}), it suffices to prove,
under these assumptions, that
if $\fext$ is continuous at $\xbar$ then $f$ is extensibly continuous
at $\xbar$.
Therefore, we assume henceforth that $\fext$ is continuous at~$\xbar$.\looseness=-1

First, suppose $f$ is lower semicontinuous.
Let $\seq{\xx_t}$ be any sequence in $\Rn$ converging to $\xbar$.
Then $f(\xx_t)=\fext(\xx_t)$, for all $t$, by
\Cref{pr:h:1}(\ref{pr:h:1a}), since $f$ is lower
semicontinuous.
Also, $\fext(\xx_t)\rightarrow\fext(\xbar)$ since $\fext$ is
continuous at $\xbar$.
Therefore, $f(\xx_t)\rightarrow\fext(\xbar)$, so
$f$ is extensibly continuous at $\xbar$.

In the remaining case, we assume
$f$ is convex (and not necessarily lower
semicontinuous).
Suppose,
by way of contradiction, that contrary to the claim, there
exists a sequence $\seq{\xx_t}$ in $\Rn$ that converges to $\xbar$,
but for which $f(\xx_t)\not\rightarrow\fext(\xbar)$.
This implies $\fext(\xbar)<+\infty$ since otherwise we would have
$\liminf f(\xx_t)\geq\fext(\xbar)=+\infty$, implying
$f(\xx_t)\rightarrow+\infty=\fext(\xbar)$.

Also, since $\fext(\xx_t)\rightarrow\fext(\xbar)$ (by continuity of $\fext$
at $\xbar$), we must have $f(\xx_t)\neq\fext(\xx_t)$ for infinitely
many values of $t$.
By discarding all other elements of the sequence, we assume
henceforth that $f(\xx_t)\neq\fext(\xx_t)$ for all $t$, while
retaining the property that $\xx_t\rightarrow\xbar$.
This implies further that, for all $t$,
$f(\xx_t)\neq(\lsc f)(\xx_t)$
(by \Cref{pr:h:1}\ref{pr:h:1a}),
and so, by \Cref{pr:lsc-props}(\ref{pr:lsc-props:b}),
that
$\xx_t\in\cl(\dom{f})\setminus\ri(\dom{f})$,
the relative boundary of $\dom{f}$.

Consequently,
by \Cref{pr:bnd-near-cl-comp} (applied to $\dom{f}$),
for each $t$, there
exists a point $\xx'_t\in\Rn\setminus\regParens{\cl(\dom{f})}$ with
$\norm{\xx'_t-\xx_t} < 1/t$.
Thus,
$\xx'_t=\xx_t+\vepsilon_t$ for some $\vepsilon_t\in\Rn$
with $\norm{\vepsilon_t} < 1/t$.
Since $\vepsilon_t\rightarrow\zero$ and $\xx_t\rightarrow\xbar$,
we also have $\xx'_t\rightarrow\xbar$ (by \Cref{pr:i:7}\ref{pr:i:7g}).

We have
\[
   \fext(\xx'_t) = (\lsc f)(\xx'_t) = f(\xx'_t) = +\infty,
\]
where the first equality is from
\Cref{pr:h:1}(\ref{pr:h:1a}),
the second is from
\Cref{pr:lsc-props}(\ref{pr:lsc-props:b}) since
$\xx'_t\not\in\cl(\dom{f})$,
and the last follows also from this latter fact.
Thus, by continuity of $\ef$ at $\xbar$,
$\fext(\xbar)=\lim \fext(\xx'_t)=+\infty$,
which contradicts $\fext(\xbar)<+\infty$.%
\indexg{continuity of extensions!extensible continuity and|)}%
\indexg{extensible continuity|)}%
\qedhere
\end{proof-parts}
\end{proof}

\indexg{nondecreasing functions, extension of|(}%
\indexg{lower semicontinuous extension!nondecreasing functions@of nondecreasing functions|(}%
\indexg{continuity of extensions!nondecreasing functions@of nondecreasing functions|(}%
As a simple illustration,
we next study properties of an extension of a convex nondecreasing function on $\R$.
(In the proposition, for $x\in\R$,
we use $\sup_{y<x} g(y)$ as shorthand for
$\sup\Braces{g(y) :\: y\in (-\infty,x)}$.)

\begin{proposition}
\label{pr:conv-inc:prop}
Let $g:\R\to\eR$ be convex and nondecreasing. Then:
\begin{letter-compact}
\item
  \label{pr:conv-inc:infsup}
  $\eg$ is continuous at $\pm\infty$ with $\eg(-\infty)=\inf g$ and $\eg(+\infty)=\sup g$. Moreover,
  if $g$ is continuous at a point $x\in\R$ then $\eg(x)=g(x)$ and $\eg$ is continuous at $x$.
\item
  \label{pr:conv-inc:nondec}
  For $x\in\R$, $\eg(x)=\sup_{y<x} g(y)$, so $\eg$ is nondecreasing.
\item
  \label{pr:conv-inc:nonconst}
  If $g$ is not constant then neither is $\eg$, and $\eg(+\infty)=+\infty$.
\item
  \label{pr:conv-inc:strictly}
  If $g$ is strictly increasing then so is $\eg$, and $\eg(+\infty)=+\infty$.
\item
  \label{pr:conv-inc:discont}
  Suppose $g$ is not continuous.
  Then there must exist a unique $z\in\R$ where
  $g$ is not continuous,
  namely, $z=\sup(\dom g)$.
  Further, $\eg(x)=g(x)$ for $x\in\R\wo\set{z}$, $\eg(z)=\sup_{x<z} g(x)$, and
  $\dom \eg=[-\infty,z]$.
  Thus, $\eg$ is continuous on $\eR\wo\set{z}$, but not at $z$.
\end{letter-compact}
\end{proposition}
\begin{proof}
~
\begin{proof-parts}
\pfpart{Part~(\ref{pr:conv-inc:infsup}):}
Let $\seq{x_t}$ be a sequence in $\R$ converging to $+\infty$. Then for any $y\in\R$, we have $x_t>y$ for all $t$ sufficiently large, so $\liminf g(x_t)\ge g(y)$. Therefore,
\[
  \sup_{y\in\R} g(y)\le\liminf g(x_t)\le\limsup g(x_t)\le\sup_{y\in\R} g(y),
\]
showing that $g(x_t)\to\sup g$. Thus, $\eg(+\infty)=\sup g$ and $g$ is extensively continuous at $+\infty$, so $\eg$ is continuous at $+\infty$ (by \Cref{thm:ext-cont-f}\ref{thm:ext-cont-f:a}). By a symmetric argument, $\eg(-\infty)=\inf g$ and $\eg$ is continuous at $-\infty$.

For the second part of the claim, let $x\in\R$ be a point where $g$ is continuous. Then $g$ is both lower semicontinuous and extensibly continuous at $x$, so $\eg(x)=g(x)$ (by \Cref{pr:h:1}\ref{pr:h:1a}) and $\eg$ is continuous at~$x$ (by \Cref{thm:ext-cont-f}\ref{thm:ext-cont-f:a}).

\pfpart{Part~(\ref{pr:conv-inc:nondec}):}
Let $x\in\R$ and let $\seq{x_t}$
be a sequence in $\R$ such that $x_t\to x$ and $g(x_t)\to\eg(x)$
(which exists by \Cref{pr:d1}). 
Then for all $y\in (-\infty,x)$, we have $x_t>y$ for all $t$ sufficiently large, so $\lim g(x_t)\ge g(y)$. Taking supremum
then yields
\[
  \sup_{y<x} g(y)\le\lim g(x_t)=\eg(x)
  \le\liminf g(x-1/t)
  \le\sup_{y<x} g(y),
\]
where the second inequality follows from $\gext$'s definition
(Eq.~\ref{eq:e:7}) since $x-1/t\to x$. Thus,
$\eg(x)=\sup_{y<x} g(y)$ as claimed.
Since, by part~(\ref{pr:conv-inc:infsup}),
$\eg(-\infty)=\inf g$ and $\eg(+\infty)=\sup g$,
we obtain that $\eg$ is nondecreasing.

\pfpart{Part~(\ref{pr:conv-inc:nonconst}):}
Suppose $g$ is not constant, implying,
since it is also nondecreasing, that there exist
$x,x'\in\R$ with $x<x'$ and $g(x)<g(x')$.
Let $y\in(x',+\infty)$.
Then $\gext(x)\leq g(x)<g(x')\leq \gext(y)$
by part~(\ref{pr:conv-inc:nondec}), so $\gext$ is not constant.

Let $g'=\lsc g$, which is convex and lower semicontinuous
(\Cref{pr:lsc-props}\ref{pr:lsc-props:a})
and equal to $\resfcn{\gext}{\R}$
by \Cref{pr:h:1}(\ref{pr:h:1a}), implying $g'(x)<g'(y)$.
Therefore, $1$ is not in
$\resc{g'}$
(as defined in Eq.~\ref{eqn:resc-cone-def}),
implying
\[
  \gext(+\infty)
  =
  \lim_{z\rightarrow+\infty} \gext(z)
  =
  \lim_{z\rightarrow+\infty} g'(z)
  =
  +\infty
\]
where the first equality is by part~(\ref{pr:conv-inc:infsup})
and the third by
\Cref{pr:stan-rec-equiv}(\ref{pr:stan-rec-equiv:a},\ref{pr:stan-rec-equiv:c})
(applied to $g'$ with $\vv=1$).

\pfpart{Part~(\ref{pr:conv-inc:strictly}):}
By
parts~(\ref{pr:conv-inc:infsup}) and~(\ref{pr:conv-inc:nondec}),
$\eg(-\infty)=\inf g$, $\eg(x)=\sup_{y<x} g(y)$ for $x\in\R$, and $\eg(+\infty)=\sup g$.
Suppose $\barx,\barx'\in\Rext$ with $\barx<\barx'$.
Then there exist $z,z'\in\R$ such that $\barx<z<z'<\barx'$,
implying $\eg(\barx)\leq g(z)<g(z')\leq\eg(\barx')$ since $g$ is strictly
increasing.
Thus, $\eg$ is strictly increasing as well.

Since $g$ is strictly increasing on $\R$, it is not constant, so $\eg(+\infty)=+\infty$ by part~(\ref{pr:conv-inc:nonconst}).

\pfpart{Part~(\ref{pr:conv-inc:discont}):}
Assume that $g$ is not continuous and let $z=\sup(\dom g)$. We claim that $z\in\R$.
If $z\in\{-\infty,+\infty\}$ then $\dom g=\R$ or $\dom g=\emptyset$, so in both cases the boundary of $\dom g$ is empty and $g$ is continuous everywhere (by \Cref{pr:stand-cvx-cont}). Thus, we must have $z\in\R$. By monotonicity of $g$, $\dom g$ is either $(-\infty,z)$ or $(-\infty,z]$.
Its boundary consists of the single point $z$, so $z$ must be the sole point where $g$ is discontinuous (by \Cref{pr:stand-cvx-cont}), and $g$ is continuous everywhere else. Therefore, by part~(\ref{pr:conv-inc:infsup}), $\eg(x)=g(x)$ for $x\in\R\wo\set{z}$ and $\eg$ is continuous on $\eR\wo\set{z}$.

By part~(\ref{pr:conv-inc:nondec}), $\eg(z)=\sup_{x<z} g(x)$. Moreover, we must have $\eg(z)<+\infty$. Otherwise, $g(z)\ge\eg(z)=+\infty$ (by \Cref{pr:h:1}\ref{pr:h:1a}),
which would imply continuity of $g$ at~$z$, because any sequence $z_t\to z$ would satisfy $\liminf g(z_t)\ge\eg(z)=+\infty$, and hence $g(z_t)\to+\infty=g(z)$. Thus, $\eg(z)<+\infty$, but $\lim\eg(z+1/t)=+\infty$, so $\eg$ is not continuous at $z$.%
\indexg{continuity of extensions!nondecreasing functions@of nondecreasing functions|)}%
\indexg{lower semicontinuous extension!nondecreasing functions@of nondecreasing functions|)}%
\indexg{nondecreasing functions, extension of|)}%
\qedhere
\end{proof-parts}
\end{proof}

\section{Working with epigraphs}
\label{sec:work-with-epis}

When working with a function $f:\Rn\rightarrow\Rext$, we have seen
that a point in the function's epigraph can be regarded either as a
pair in $\Rn\times\R$ or equivalently as a vector in $\R^{n+1}$.
In this section, we derive a
similar
equivalence
for epigraphs of astral functions.

\indexg{xy pairs@$\rpair{\xbar}{y}$ pairs!embedded in rnplusone@embedded in $\extspac{n+1}$|(}%
\indexg{embedding of rnbartimesr in rnplusone@embedding of $\extspace\times\R$ in $\extspac{n+1}$|(}%
\indexg{epigraph!astral function@of astral function|(}%
Let $F:\extspace\rightarrow\Rext$.
Then $F$'s epigraph
is a subset of $\extspace\times\R$, namely, the set of pairs
$\rpair{\xbar}{y}\in\extspace\times\R$ for which $F(\xbar)\leq y$.
It will often be beneficial
to regard $\extspace\times\R$ as a subset of
$\extspac{n+1}$, so that $\epi{F}$ also becomes a subset of
$\extspac{n+1}$, an astral space with numerous favorable properties,
such as compactness.
Although $\extspace\times\R$ is formally not a subset of
$\extspac{n+1}$,
we show next that it can be embedded in $\extspac{n+1}$,
that is, shown to be homeomorphic to a subset of the larger space.
Going forward, this will allow us to treat $\extspace\times\R$,
and so also the epigraph of any astral function,
as a subset of $\extspac{n+1}$.%
\indexg{epigraph!astral function@of astral function|)}

\indexg{component projection matrices!pairs in rntimesr@for pairs in $\Rn\times\R$|(}%
To develop this embedding, first consider
a point $\zz=\rpair{\xx}{y}$ in $\R^{n+1}$, where $\xx\in\Rn$ and $y\in\R$.
Then we can extract $\xx$ and $y$ from $\zz=\rpair{\xx}{y}$, by
applying the linear maps
associated with the matrices
\begin{equation}
\label{eq:PPx:PPy}
\indexm{p 1, p 2}{$\PPx$, $\PPy$}{component projection matrices (for pairs in $\Rn\times\R$)}%
   \PPx=[\Idnn,\zero_n]
   \qquad
   \text{and}
   \qquad
   \PPy=[\trans{\zero_n},1].
\end{equation}
Here, as usual,
$\Idnn$ denotes the $n\times n$ identity matrix and $\zerov{n}$ the
$n\times 1$ all-zeros (column) vector.
Thus, $\PPx$ is an $n\times(n+1)$ matrix whose first $n$ columns are
the $n\times n$ identity matrix and whose last column is all zeros,
while $\PPy$ is a $1\times (n+1)$ row vector
that is all zeros except
the last entry which is equal to $1$.
When applied to $\zz$,
we obtain
$\PPx\zz=\xx$ and $\PPy\zz=y$;
in this way, $\PPx$ extracts the first $n$ coordinates of $\zz$,
while $\PPy$ extracts $\zz$'s last coordinate.

Note that $\PPy=\trans{\ee_{n+1}}$, where $\ee_{n+1}$ is the last
standard basis vector in $\R^{n+1}$.
Consequently, $\PPy\zz=\trans{\ee_{n+1}}\zz=\zz\cdot\ee_{n+1}$,
yielding, as just noted, $\zz$'s last coordinate.
Likewise, when $\PPy$ is applied to an astral point
$\zbar\in\extspacnp$,
we obtain
$\PPy\zbar=\trans{\ee_{n+1}}\zbar=\zbar\cdot\ee_{n+1}$
(by \Cref{pr:trans-uu-xbar}), whose value in $\Rext$
can also informally be regarded as the ``last coordinate'' of $\zbar$.%
\indexg{component projection matrices!pairs in rntimesr@for pairs in $\Rn\times\R$|)}%

When constructing the embedding of $\eRn\times\R$ in $\extspac{n+1}$,
we posit that  $\extspace\times\R$ should correspond to those points $\zbar$ in $\extspac{n+1}$
whose ``last coordinate''
(obtained by applying $\PPy$) is in $\R$. We denote the set of all such
points
\begin{equation}  \label{eq:finclset-defn}
\indexm{m n+1}{$\finclset$}{all points in $\extspacnp$ with finite last coordinate}%
   \finclset = \{ \zbar\in \extspac{n+1} :\:
                     \PPy\zbar\in \R \}.
\end{equation}
The next theorem shows that $\extspace\times\R$ is homeomorphic to $\finclset$. In its proof,
we use the following orthogonality identities, which follow from
$\PPx$ and $\PPy\negKern$'s definitions:
\begin{equation}  \label{eq:PPx:PPy:orth}
   \PPxlow\trans{\PPx} = \Idnn,\quad
   \PPylow\trans{\PPy} = 1,\quad
   \PPxlow\trans{\PPy} = \zero_n,\quad
   \PPylow\trans{\PPx} = \trans{\zero_n}.
\end{equation}

\begin{theorem}  \label{thm:homf}
  Define
  $\homf:\extspace\times\R\rightarrow\finclset$ to be the function
  \begin{equation}   \label{eq:thm:homf:1}
    \indexm{mu}{$\homf$}{embedding of $\extspace\times\R$}%
    \homf(\xbar,y) = \trans{\PPx} \xbar \plusl \trans{\PPy}y
  \end{equation}
  for $\xbar\in\extspace$, $y\in\R$,
  where 
  $\PPx=[\Idnn,\zero_n]$ and $\PPy=[\trans{\zero_n},1]$.
  Then:
  \begin{letter-compact}
  \item  \label{thm:homf:aa}
    $\homf$ is bijective with inverse
    \begin{equation}  \label{eqn:homfinv-def}
      \homfinv(\zbar) = \rpair{\PPx\zbar}{\,\PPy\zbar}
    \end{equation}
    for $\zbar\in\finclset$.
  \item  \label{thm:homf:b}
    $\homf$ is a homeomorphism (that is, both
    $\homf$ and its inverse are continuous).
  \item  \label{thm:homf:a}
    For all $\xbar\in\extspace$, $y\in\R$,
    \[
      \homf(\xbar,y) \cdot \rpair{\uu}{v} =  \xbar\cdot\uu + y v
    \]
    for all $\rpair{\uu}{v}\in\Rn\times\R=\Rnp$.
  \item  \label{thm:homf:c}
    $\homf(\xx,y)=\rpair{\xx}{y}$ for all $\rpair{\xx}{y}\in\Rn\times\R=\Rnp$.
  \end{letter-compact}
\end{theorem}

\begin{proof}
~

\begin{proof-parts}
\pfpart{Part~(\ref{thm:homf:aa}):}
Let $\gamma:\finclset\rightarrow\extspace\times\R$ be the function
given in \eqref{eqn:homfinv-def}, which we aim to show is the
functional inverse of $\homf$.
Let $\xbar\in\extspace$ and $y\in\R$.
We first show
$ \gamma(\homf(\xbar,y)) = \rpair{\xbar}{y} $.
Let $\zbar=\homf(\xbar,y)$.
Then by $\homf$'s definition, properties of astral linear maps
(\Cref{pr:h:4}\ref{pr:h:4c}), and orthogonality identities (Eq.~\ref{eq:PPx:PPy:orth}),
\begin{align*}
    \PPx\zbar
    =\PPxlow\trans{\PPx}\xbar\plusl\PPxlow\trans{\PPy}y
    =\xbar,
\\
    \PPy\zbar
    =\PPylow\trans{\PPx}\xbar\plusl\PPylow\trans{\PPy}y
    =y.
\end{align*}
This shows both that $\zbar\in\finclset$ and that $\gamma$ is a left inverse of $\homf$.

Next, let $\zbar\in\finclset$, so $\PPy\zbar\in\R$. We aim to show that
$\homf(\gamma(\zbar))=\zbar$. We calculate:
\begin{align*}
    \homf\bigParens{\gamma(\zbar)}
    = \homf\bigParens{\rpair{\PPx\zbar}{\,\PPy\zbar}}
    &= \trans{\PPx}\PPxlow\zbar
       \plusl
       \trans{\PPy}\PPylow\zbar
\\
    &= \bigParens{\trans{\PPx}\PPxlow+\trans{\PPy}\PPylow}\zbar
     = \Idn{n+1}\zbar
     = \zbar.
\end{align*}
The third equality is by \Cref{prop:commute:AB}
(noting that
$\PPy\zbar\in\R$ and so that $\trans{\PPy}\PPy\zbar\in\R^{n+1}$).
The fourth
equality is from the definition of $\PPx$ and $\PPy$.
Thus, $\gamma$ is also a right inverse of $\homf$.
Therefore, $\homf$ is bijective with inverse
$\homfinv=\gamma$, as claimed.

\pfpart{Part~(\ref{thm:homf:b}):}
The function $\homf$ is continuous since if
$\seq{\rpair{\xbar_t}{y_t}}$ is any sequence in $\extspace\times\R$
that converges to $\rpair{\xbar}{y}\in\extspace\times\R$,
then
\[
  \homf(\xbar_t,y_t)
  =\trans{\PPx}\xbar_t\plusl\trans{\PPy}y_t
  \to
  \trans{\PPx}\xbar\plusl\trans{\PPy}y
  =\homf(\xbar,y)
\]
by continuity of astral linear maps
(\Cref{thm:linear:cont}\ref{thm:linear:cont:b})
and \Cref{pr:i:7}(\ref{pr:i:7g}).
The function $\homfinv$ is continuous by
continuity of astral linear maps
(\Cref{thm:linear:cont}\ref{thm:linear:cont:b}).
Thus, $\homf$ is a homeomorphism.

\pfpart{Part~(\ref{thm:homf:a}):}
Let $\xbar\in\extspace$, $y\in\R$, and
$\rpair{\uu}{v}\in\Rn\times\R$.
Then
\begin{align*}
  \homf(\xbar,y) \cdot \rpair{\uu}{v}
  &=
  (\trans{\PPx} \xbar) \cdot \rpair{\uu}{v}
  \plusl
  (\trans{\PPy} y) \cdot \rpair{\uu}{v}
  \\
  &=
  \xbar \cdot (\PPx \rpair{\uu}{v})
  \plusl
  y \cdot (\PPy \rpair{\uu}{v})
  \\
  &=
  \xbar \cdot \uu
  +
  y v.
\end{align*}
The first equality is by $\homf$'s definition and
\Cref{pr:i:6};
the second by \Cref{thm:Ax-dot-u};
the third by definition of $\PPx$ and $\PPy$ (and since $y v\in\R$).

\pfpart{Part~(\ref{thm:homf:c}):}
Let $\xx\in\Rn$, $y\in\R$, and $\zbar=\homf(\xx,y)$.
Then
$\zbar\cdot\rpair{\uu}{v}=\rpair{\xx}{y}\cdot\rpair{\uu}{v}$
for all $\rpair{\uu}{v}\in\Rn\times\R$, by
part~(\ref{thm:homf:a}).
Combined with \Cref{pr:i:4},
this implies $\zbar=\rpair{\xx}{y}$.
(Alternatively, this could be proved directly from the definition of
$\homf$.)
\qedhere
\end{proof-parts}
\end{proof}

The function $\homf$ from \Cref{thm:homf}
is called the
\indexg{natural embedding of $\eRn\times\R$}%
\emph{natural embedding of $\eRn\times\R$ in $\extspacnp$}.
The
\namecref{thm:homf}
shows that
$\homf$ is a homeomorphism between the set of pairs
$\extspace\times\R$ and the set
$\finclset\subseteq\extspacnp$ given in \eqref{eq:finclset-defn},
with each point $\rpair{\xbar}{y}\in\eRn\times\R$ mapped to
$\homf(\xbar,y)\in\finclset$.
Moreover, by \Cref{thm:homf}(\ref{thm:homf:c}), $\mu$ is
an extension of the natural isomorhism between $\Rn\times\R$ and $\Rnp$ (as topological vector
spaces) that we have been tacitly using throughout this book.
The existence of $\mu$ allows us to treat
$\extspace\times\R$ as a subset of $\extspacnp$, and to continue to identify $\Rn\times\R$
with $\Rnp$.%
\indexg{embedding of rnbartimesr in rnplusone@embedding of $\extspace\times\R$ in $\extspac{n+1}$|)}

To simplify notation, for the rest of this book, we identify each pair
$\rpair{\xbar}{y}$ with its homeomorphic image $\homf(\xbar,y)$ so that,
when clear from context,
$\rpair{\xbar}{y}$ may denote either the given pair in
$\extspace\times\R$ or the point
$\homf(\xbar,y)\in\finclset\subseteq\extspacnp$.
Note importantly that this convention only applies when $y$ is finite
(in $\R$, not $\pm\infty$).
In some cases, it may be that the interpretation cannot be determined
from context; usually, this means that either interpretation can be
used.
Nonetheless, in the vast majority of cases, $\rpair{\xbar}{y}$ should
be regarded as the point $\homf(\xbar,y)$ in $\extspacnp$.

We also apply this simplification to subsets
$S\subseteq\extspace\times\R$,
such as the epigraph $\epi F$ of an astral function
$F:\extspace\rightarrow\Rext$,
writing simply $S$, when clear from context, to denote
its homeomorphic image, $\homf(S)\subseteq\finclset$.
For instance, we can now write
$\extspace\times\R\subseteq\extspacnp$, which formally means that
$\homf(\extspace\times\R)\subseteq\extspacnp$.%
\indexg{xy pairs@$\rpair{\xbar}{y}$ pairs!embedded in rnplusone@embedded in $\extspac{n+1}$|)}%

\indexg{xy pairs@$\rpair{\xbar}{y}$ pairs!operations on|(}%
The next proposition summarizes properties of such pairs
$\rpair{\xbar}{y}$, and is largely a mere restatement of
\Cref{thm:homf}.
In this proposition, the notation $\rpair{\xbar}{y}$ always refers to
a point in $\extspacnp$.

\begin{proposition}   \label{pr:xy-pairs-props}
  Let $\xbar,\xbar'\in\extspace$, let $y,y'\in\R$,
  let $\zbar\in\extspacnp$,
  let $\uu\in\Rn$, and let $v\in\R$.
  Also, let $\PPx=[\Idnn,\zero_n]$ and $\PPy=[\trans{\zero_n},1]$.
  Then:
  \begin{letter-compact}
  \item    \label{pr:xy-pairs-props:a}
    $\rpair{\xbar}{y} = \trans{\PPx} \xbar \plusl \trans{\PPy}y$.
  \item    \label{pr:xy-pairs-props:b}
    $\rpair{\xbar}{y}\cdot\rpair{\uu}{v} = \xbar\cdot\uu + yv$.
  \item    \label{pr:xy-pairs-props:b-new}
    $\zbar\cdot\rpair{\uu}{0}=(\PPx\zbar)\cdot\uu$
    and
    $\zbar\cdot\rpair{\zero}{v}=(\PPy\zbar)v$.
    (In particular,
    $\zbar\cdot\rpair{\zero}{1}=\PPy\zbar$.)
  \item    \label{pr:xy-pairs-props:c}
    $\PPx\rpair{\xbar}{y} = \xbar$ and $\PPy\rpair{\xbar}{y} = y$.
  \item    \label{pr:xy-pairs-props:d}
    If $\PPy\zbar\in\R$ then
    $\zbar = \rpair{\PPx\zbar}{\,\PPy\zbar}$.
  \item    \label{pr:xy-pairs-props:h}
    $\rpair{\xbar}{y}=\rpair{\xbar'}{y'}$
    if and only if
    $\xbar=\xbar'$ and $y=y'$.
  \item    \label{pr:xy-pairs-props:e}
    $\rpair{\xbar}{y}\plusl\rpair{\xbar'}{y'} = \rpair{\xbar\plusl\xbar'}{y+y'}$.
  \item    \label{pr:xy-pairs-props:f}
    $\lambda \rpair{\xbar}{y} = \rpair{\lambda\xbar}{\lambda y}$
    for all $\lambda\in\R$.
  \item    \label{pr:xy-pairs-props:g}
    Let $\seq{\xbar_t}$ be a sequence in $\extspace$,
    and
    let $\seq{y_t}$ be a sequence in $\R$.
    Then $\rpair{\xbar_t}{y_t}\rightarrow\rpair{\xbar}{y}$
    if and only if
    $\xbar_t\rightarrow\xbar$
    and
    $y_t\rightarrow y$.
  \end{letter-compact}
\end{proposition}

\begin{proof}
Part~(\ref{pr:xy-pairs-props:a})
follows from \eqref{eq:thm:homf:1} of
\Cref{thm:homf}.
Part~(\ref{pr:xy-pairs-props:b}) is from
\Cref{thm:homf}(\ref{thm:homf:a}).

For part~(\ref{pr:xy-pairs-props:b-new}), we have
\[
  \zbar\cdot\rpair{\uu}{0}
  =
  \zbar\cdot(\trans{\PPx}\uu)
  =
  (\PPx\zbar)\cdot\uu,
\]
where the first equality is from the definition of $\PPx$,
and the second by
\Cref{thm:Ax-dot-u}.
Similarly,
$
  \zbar\cdot\rpair{\zero}{v}
  =
  \zbar\cdot(\trans{\PPy}v)
  =
  (\PPy\zbar)v
$.

Parts~(\ref{pr:xy-pairs-props:c})
and~(\ref{pr:xy-pairs-props:d}) follow from
\Cref{thm:homf}(\ref{thm:homf:aa}).
Part~(\ref{pr:xy-pairs-props:h}) is because $\homf$ is a bijection
(\Cref{thm:homf}\ref{thm:homf:aa}).

For part~(\ref{pr:xy-pairs-props:e}), we have
\begin{align*}
  \rpair{\xbar}{y}\plusl\rpair{\xbar'}{y'}
  &=
  \bigParens{\trans{\PPx} \xbar \plusl \trans{\PPy}y}
  \plusl
  \bigParens{\trans{\PPx} \xbar' \plusl \trans{\PPy}y'}
  \\
  &=
  \trans{\PPx} (\xbar\plusl\xbar') \plusl \trans{\PPy}(y+y')
  \\
  &=
  \rpair{\xbar\plusl\xbar'}{y+y'}.
\end{align*}
The first and third equalities are by
part~(\ref{pr:xy-pairs-props:a}).
The second equality uses
\Cref{pr:h:4}(\ref{pr:h:4c}).
The proof of
part~(\ref{pr:xy-pairs-props:f})
is similar.

Part~(\ref{pr:xy-pairs-props:g}) is because $\homf$ is a homeomorphism
(\Cref{thm:homf}\ref{thm:homf:b}).%
\indexg{xy pairs@$\rpair{\xbar}{y}$ pairs!operations on|)}
\end{proof}

In the next two propositions,
we write $\clmx{S}$%
\indexm{cls700}{$\clmx{S}$}{closure in $\extspace\times\R$} 
for the closure in $\extspace\times\R$
of any set $S$ in that space.
(For a set
$S\subseteq\extspacnp$, we continue to write $\Sbar$ for the
closure of $S$ in $\extspacnp$.)
\indexg{lower semicontinuous hull!epigraph of|(}%
\indexg{epigraph!astral function@of astral function|(}%
The first of these propositions relates the epigraph
of a general astral function to that of its lower semicontinuous hull.
The second applies these to relate the epigraph of a
function $f:\Rn\rightarrow\Rext$ to that of its extension $\fext$.

\begin{proposition}  \label{pr:epiF-epi-lscF}
  Let $F:\extspace\rightarrow\Rext$,
  and let $G=\lsc{F}$.
  Then:
  \begin{letter}
  \item  \label{pr:epiF-epi-lscF:a}
    $\epi G = \epibar F$.
  \item  \label{pr:epiF-epi-lscF:b}
    $\epi G = \clepi F \cap \finclset$.
  \item  \label{pr:epiF-epi-lscF:c}
    $\clepi G = \clepi F$.
  \end{letter}
\end{proposition}

\begin{proof}
~

\begin{proof-parts}
\pfpart{Part~(\ref{pr:epiF-epi-lscF:a}):}
This is immediate from
\Cref{prop:lsc:characterize}(\ref{prop:lsc:characterize:c}).

\pfpart{Part~(\ref{pr:epiF-epi-lscF:b}):}
By \Cref{thm:homf}, $\eRn\times\R$ is homeomorphic with
$\finclset$, which is a topological subspace of $\extspacnp$, so
by part (\ref{pr:epiF-epi-lscF:a}) and \Cref{prop:subspace}(\ref{i:subspace:closure}),
$\epi G = \epibar{f} = \clepi{F} \cap \finclset$.

\pfpart{Part~(\ref{pr:epiF-epi-lscF:c}):}
By parts~(\ref{pr:epiF-epi-lscF:a})
and~(\ref{pr:epiF-epi-lscF:b}),
$ \epi{F}\subseteq\epi{G}\subseteq\clepi{F} $,
so by \Cref{pr:gen-hull-ops}(\ref{pr:gen-hull-ops:d}),
$\clepi{G}=\epibarbar{F}$ (since topological closure is a hull operator).%
\indexg{lower semicontinuous hull!epigraph of|)}%
\indexg{epigraph!astral function@of astral function|)}
\qedhere
\end{proof-parts}
\end{proof}

\begin{proposition}  \label{pr:wasthm:e:3}
\indexg{lower semicontinuous extension!epigraph of|(}%
\indexg{epigraph!extension@of extension|(}%
  Let $f:\Rn\rightarrow\Rext$.
  Then:
  \begin{letter}
  \item  \label{pr:wasthm:e:3a}
    $\epi{\fext}=\epibar{f}$.
  \item  \label{pr:wasthm:e:3b}
    $\epi{\fext} = \epibarbar{f} \cap \finclset$.
  \item  \label{pr:wasthm:e:3c}
    $\clepifext = \epibarbar{f}$.
  \end{letter}
\end{proposition}

\begin{proof}
Let $F=\ftriv$ and let $G=\lsc{F}$.
Then $\epi{F}=\epi f$ by $\ftriv$'s definition,
and $\epi{G}=\epi{\fext}$ since $G=\fext$
by \Cref{pr:lsc-ftriv-is-fext}.
The claims therefore all follow from
\Cref{pr:epiF-epi-lscF}.%
\indexg{lower semicontinuous extension!epigraph of|)}%
\indexg{epigraph!extension@of extension|)}
\end{proof}

\section{Astronic reductions}
\label{sec:shadow}

We next begin the study of reductions, a core technique for
analyzing astral functions.

In \Cref{pr:h:6}, we saw how every infinite astral point $\xbar$ can
be decomposed into the astron associated with its dominant direction
$\vv$ and the
projection $\exx^\perp$ orthogonal to~$\vv$, whose astral rank is
lower than that of $\exx$. This decomposition forms the basis of proofs by induction on astral rank.
For the purpose of applying this technique when analyzing an extension $\ef$,
we next introduce a kind of projection operation, which effectively reduces the dimensionality of the domain of $\ef$ while preserving its key properties.

\begin{definition}   \label{def:astron-reduction}
\indexg{reductions, astronic|(}%
\indexg{reductions, astronic!defined|(}%
Let $f:\Rn\rightarrow\Rext$, and let $\vv\in\Rn$.
The \emph{reduction of $f$ at astron $\limray{\vv}$}
is the function
$\fshadv:\Rn\rightarrow\Rext$ defined, for $\xx\in\Rn$, by
\[
\indexm{f e 300}{$\fshadv$}{reduction (astronic)}%
  \fshadv(\xx) = \fext(\limray{\vv} \plusl \xx).
\]
\end{definition}
We refer to this type of reduction as
\indexg{reductions, astronic!defined|)}%
\emph{astronic};
more general reductions will be
introduced later in \Cref{sec:ent-closed-fcn}.

Let $g=\fshadv$ be such a reduction.
This function is a constant in the
direction of
$\vv$ (that is,
$g(\xx)=g(\xx+\lambda\vv)$ for all $\xx\in\Rn$ and $\lambda\in\R$), which
means that the reduction~$g$ can be regarded informally as a function
only over the space orthogonal to $\vv$, even if it is formally
defined over all of $\Rn$.
In this sense, $f$ has been ``reduced'' in forming $g$.

\begin{example}[Reduction of the product of hyperbolas]
\label{ex:recip-fcn-eg-cont}
\indexg{Product of hyperbolas!reduction of|(}%
Suppose that $f$ is the product of hyperbolas
from \Cref{ex:recip-fcn-eg}, and let $\vv=\ee_1$.
Then, for $\xx\in\R^2$,
\begin{equation} \label{eqn:recip-fcn-eg:g}
  g(\xx)
  =
  g(x_1,x_2)
  =
  \fext(\limray{\ee_1}\plusl\xx)
  =
  \begin{cases}
    0
      & \text{if $x_2\geq 0$,}
    \\
    +\infty
      & \text{otherwise,}
  \end{cases}
\end{equation}
as can be seen
by plugging in the values $\fext(\limray{\ee_1}\plusl\xx)$ derived in \Cref{ex:recip-fcn-eg}.%
\indexg{Product of hyperbolas!reduction of|)}
\end{example}

When forming a reduction $\fshadv$,
the vector $\vv$ is usually (but not always) assumed to be
in the {recession cone} of $f$, $\resc{f}$, as defined in
\Cref{sec:prelim:rec-cone}.
For instance, for the product of hyperbolas
(\Cref{ex:recip-fcn-eg,ex:recip-fcn-eg-cont}),
the recession cone consists of all vectors $\vv=\trans{[v_1,v_2]}$
with $v_1,v_2\geq 0$, that is, vectors in $\Rpos^2$.

If $\vv\in\resc{f}$, then, as shown next,
the minimum of the reduction $g$
is the same as the minimum of~$f$.
This
suggests that $f$ can be minimized by first minimizing $g$ and
then adjusting the resulting solution appropriately.
Later, in \Cref{sec:minimizers},
we will develop ideas along these lines
which constructively characterize the minimizers of
$f$ by defining and recursively minimizing an astronic reduction.

Here are some simple properties of reductions:

\begin{proposition}  \label{pr:d2}
Let $f:\Rn\rightarrow\Rext$,
let $\vv\in\resc{f}$,
and let
$g=\fshadv$
be the reduction of $f$ at $\limray{\vv}$.
Let $\xx\in\Rn$,
and let $\xperp$ denote its projection orthogonal to~$\vv$.
Then:
\begin{letter-compact}
\item  \label{pr:d2:a}
  $g(\xperp) = g(\xx)$.
\item  \label{pr:d2:b}
  $g(\xx) \leq f(\xx)$.
\item  \label{pr:d2:c}
  $\inf g = \inf f$.
  Consequently, $g\equiv+\infty$ if and only if $f\equiv+\infty$.
\end{letter-compact}
\end{proposition}

\begin{proof}
~
\begin{proof-parts}
\pfpart{Part~(\ref{pr:d2:a}):}
The statement follows because
$\limray{\vv}\plusl\xx=\limray{\vv}\plusl\xperp$
by the Projection Lemma (\Cref{lemma:proj}).

\pfpart{Part~(\ref{pr:d2:b}):}
Since $\vv\in\resc{f}$, for all $t$,
$ f(\xx) \geq f(\xx + t \vv) $.
Thus,
\[ f(\xx) \geq \liminf f(\xx+t\vv)
          \geq \fext(\limray{\vv}\plusl\xx) = g(\xx),
\]
with the second inequality by $\fext$'s definition
(Eq.~\ref{eq:e:7}) and since
$\xx+t\vv\rightarrow\limray{\vv}\plusl\xx$
by \Cref{pr:i:7}(\ref{pr:i:7f}).

\pfpart{Part~(\ref{pr:d2:c}):}
Since $g(\xx)=\ef(\limray{\vv}\plusl\xx)$, we have $\inf
g\ge\inf\ef=\inf f$, where the equality follows by \Cref{pr:fext-min-exists}.
On the other hand, by part~(\ref{pr:d2:b}),
$\inf g \leq \inf f$.
\qedhere
\end{proof-parts}
\end{proof}

\indexg{reductions, astronic!extension of|(}%
The reduction $g=\fshadv$ has its own extension $\gext$.
In a moment,
in \Cref{thm:d4},
we will show that the basic properties for $g$'s behavior
given in
\Cref{pr:d2} and
\Cref{def:astron-reduction}
carry over in a natural way to its extension.

Before proving the theorem, we give the following more general
proposition which will be used in the proof, specifically in showing
that the kind of property
given in \Cref{pr:d2}(\ref{pr:d2:a}) generally carries over
to extensions.

\begin{proposition}
\label{pr:PP:ef}
Let $f:\Rn\to\eR$ and $\PP\in\R^{n\times n}$ be an orthogonal projection matrix such that
$f(\xx)=f(\PP\xx)$ for all $\xx\in\Rn$. Then $\ef(\xbar)=\ef(\PP\xbar)$ for all $\xbar\in\eRn$.
\end{proposition}
\begin{proof}
Let $\xbar\in\eRn$. By \Cref{pr:d1}, there exists a sequence $\seq{\xx_t}$ in $\Rn$ with $\xx_t\to\xbar$
and $f(\xx_t)\to\ef(\xbar)$. Then
\[
  \ef(\xbar)=\lim f(\xx_t)=\lim f(\PP\xx_t)\ge\ef(\PP\xbar)
\]
where the inequality follows from the definition of $\ef$ and the fact
that $\PP\xx_t\to\PP\xbar$ by continuity of linear maps
(\Cref{thm:linear:cont}\ref{thm:linear:cont:b}).

It remains to prove the reverse inequality. By \Cref{pr:d1}, there exists a sequence $\seq{\yy_t}$ in $\Rn$ with $\yy_t\to\PP\xbar$ and $f(\yy_t)\to\ef(\PP\xbar)$. Let $\ww_1,\dotsc,\ww_k$, for some $k\ge 0$, be an orthonormal basis for $(\colspace\PP)^\perp$ and let $\WW=[\ww_1,\dotsc,\ww_k]$. By
\Cref{cor:inv-proj-seq}, there exists a sequence $\seq{\cc_t}$ in $\R^k$ such that $\yy_t+\WW\cc_t\to\xbar$. Therefore,
\begin{align}
\notag
  \ef(\PP\xbar)=\lim f(\yy_t)
&               =\lim f(\PP\yy_t)
\\
\notag
&
               =\lim f\bigParens{\PP(\yy_t+\WW\cc_t)}
               =\lim f(\yy_t+\WW\cc_t)
                             \ge\ef(\xbar)
\end{align}
where the third equality is because
$\PP\WW=\zero_{n\times k}$
(by \Cref{pr:proj-mat-props}\ref{pr:proj-mat-props:e}).
Thus,
$\ef(\xbar)=\ef(\PP\xbar)$ for all $\xbar\in\eRn$.%
\end{proof}

\begin{theorem}  \label{thm:d4}
Let $f:\Rn\rightarrow\Rext$,
let $\vv\in\resc{f}$,
and let
$g=\fshadv$
be the reduction of $f$ at $\limray{\vv}$.
Let $\xbar\in\extspace$,
and let $\xbarperp$ denote its projection orthogonal to~$\vv$.
Then:
\begin{letter-compact}
\item  \label{thm:d4:a}
  $\gext(\xbarperp) = \gext(\xbar)$.
\item  \label{thm:d4:b}
  $\gext(\xbar)\leq \fext(\xbar)$.
\item  \label{thm:d4:c}
  $\gext(\xbar)=\fext(\limray{\vv}\plusl\xbar)$.
\end{letter-compact}
\end{theorem}

\begin{proof}
~
\begin{proof-parts}
\pfpart{Part~(\ref{thm:d4:a}):}
This follows directly from \Cref{pr:d2}(\ref{pr:d2:a}) and
\Cref{pr:PP:ef}
(since $\xbarperp=\PP\xbar$ where $\PP$ is the projection
matrix orthogonal to $\vv$).

\pfpart{Part~(\ref{thm:d4:b}):}
By \Cref{pr:d2}(\ref{pr:d2:b}), we have $g\le f$, so also $\eg\le\ef$
(\Cref{pr:h:1}\ref{pr:h:1:geq}).

\pfpart{Part~(\ref{thm:d4:c}):}
First,
\begin{equation}   \label{eq:thm:d4:1}
  \gext(\xbar)
  =
  \gext(\xbarperp)
  =
  \gext\paren{(\limray{\vv}\plusl\xbar)^{\bot}}
  =
  \gext(\limray{\vv}\plusl\xbar)
  \leq
  \fext(\limray{\vv}\plusl\xbar).
\end{equation}
The first and third equalities are by
part~(\ref{thm:d4:a}), and
the second is because
$(\limray{\vv}\plusl\xbar)^{\bot} = \xbarperp$
by \Cref{pr:h:5}(\ref{pr:h:5c},\ref{pr:h:5e}).
The inequality is by part~(\ref{thm:d4:b}).

To show the reverse inequality, by \Cref{pr:d1},
there exists a sequence $\seq{\xx_t}$
in~$\Rn$ with $\xx_t\rightarrow\xbar$ and
$g(\xx_t)\rightarrow\gext(\xbar)$.
For each $t$, let $\ybar_t=\limray{\vv}\plusl\xx_t$.
Then $\ybar_t\rightarrow\limray{\vv}\plusl\xbar$
(by \Cref{pr:i:7}\ref{pr:i:7f}).
Thus,
\[
  \gext(\xbar)
  =
  \lim g(\xx_t)
  =
  \lim \fext(\limray{\vv}\plusl \xx_t)
  =
  \lim \fext(\ybar_t)
  \geq
  \fext(\limray{\vv}\plusl \xbar)
\]
where the last inequality follows by lower semicontinuity of $\ef$
(\Cref{prop:ext:F}\ref{prop:ext:F:a}).%
\indexg{reductions, astronic!extension of|)}
\qedhere
\end{proof-parts}
\end{proof}

In a moment, we will see in
\Cref{thm:a10-nunu}
that, under the assumption that $f$ is convex
and $\vv$ is in its recession cone,
the reduction $g=\fshadv$ is convex and lower semicontinuous.
\indexg{shadow (of function)|(}%
This theorem further relates $g$ to another function
$\gtil:\Rn\rightarrow\Rext$ given, for $\xx\in\Rn$, by
\begin{equation}  \label{eqn:gtil-defn}
  \gtil(\xx)=\inf_{\lambda\in\R}  f(\xx+\lambda \vv).
\end{equation}
This function can be viewed as a kind of ``shadow'' of $f$ in the
direction of~$\vv$.
Specifically, the theorem shows that $g$ is the
{lower semicontinuous hull} of $\gtil$.

\begin{example}[Shadow of the product of hyperbolas]
\indexg{Product of hyperbolas!shadow of|(}%
Consider again the product of hyperbolas
from \Cref{ex:recip-fcn-eg} and let $\vv=\ee_1$.
Then, for $\xx\in\R^2$,
\[
  \gtil(\xx)
  =
  \gtil(x_1,x_2)
  =
  \begin{cases}
    0
      & \text{if $x_2 > 0$,}
    \\
    +\infty
      & \text{otherwise,}
  \end{cases}
\]
which differs from the reduction $g$ derived in \Cref{ex:recip-fcn-eg-cont} only
when $x_2=0$.
Evidently, $g=\lsc\gtil$, as is true in general.%
\indexg{Product of hyperbolas!shadow of|)}%
\end{example}

Before proving \Cref{thm:a10-nunu}, we first consider functions of a
form that generalizes that given in \eqref{eqn:gtil-defn},
noting some useful properties of such functions:
\begin{definition}  \label{dfn:shadow-fcns}
  Let $f:\Rn\rightarrow\Rext$, and let $\VV\in\Rnk$.
  The \emph{shadow of $f$ along $\VV$} is the function
  $\fV:\Rn\rightarrow\Rext$ defined, for $\xx\in\Rn$, by
  \begin{equation}   \label{eqn:fV-init-defn}
\indexm{f v}{$\fV$}{shadow (of function)}%
    \fV(\xx)
    =
    \inf_{\yy\in\colspace\VV} f(\xx+\yy)
    =
    \inf_{\bb\in\Rk} f(\xx+\VV\bb).
  \end{equation}
\end{definition}
Thus, $\gtil$, as given in
\eqref{eqn:gtil-defn},
is equal to $\fv$, the shadow of $f$ along $\VV=[\vv]=\vv$,
the matrix with the single column $\vv$.
The next proposition
provides a
succinct and convenient representation for all such functions.
The proposition uses the notation for the composition of a function with a linear map
and the image of a function under a linear map,
which were introduced in
\Cref{sec:prelim:cvx-fcns}
(Eqs.~\ref{eq:fA-defn} and~\ref{eq:lin-image-fcn-defn}).

\begin{proposition}     \label{pr:gtil-is-PPfP-gen}
  Let $f:\Rn\rightarrow\Rext$, and let $\VV\in\Rnk$.
  Let $\PP$ be the projection matrix onto
  $(\colspace{\VV})^\perp$.
  Then $\fV=(\PPf)\PP$.
  Consequently, if $f$ is convex then so is $\fV$.
\end{proposition}

\begin{proof}
Let $\xx\in\Rn$.
Then
\begin{align*}
  [(\PPf)\PP](\xx)
  =
  (\PPf)(\PP\xx)
  &=
  \inf\Braces{ f(\zz) :\: \zz\in\Rn,\, \PP\zz = \PP\xx }
  \\
  &=
  \inf\Braces{ f(\zz) :\: \zz\in\xx + \colspace{\VV} }
  =
  \fV(\xx).
\end{align*}  
The first two equalities are by
Eqs.~(\ref{eq:fA-defn}) and~(\ref{eq:lin-image-fcn-defn}),
respectively.
The fourth equality is by $\fV$'s definition.
The third equality is because, for $\zz\in\Rn$,
$\PP\zz=\PP\xx$ if and only if
$\zz\in\xx+\colspace{\VV}$.
This is because if $\zz=\xx+\VV\bb$
for some $\bb\in\Rk$ then $\PP\zz=\PP\xx+\PP\VV\bb=\PP\xx$
(since $\PP\VV=\zeromat{n}{k}$
by \Cref{pr:proj-mat-props}\ref{pr:proj-mat-props:e}).
Conversely, by \Cref{pr:lin-decomp-rel-vecs},
$\xx=\PP\xx+\VV\cc$ and $\zz=\PP\zz+\VV\dd$ for some $\cc,\dd\in\Rk$,
so if $\PP\xx=\PP\zz$ then $\zz=\xx+\VV(\dd-\cc)$.
Thus, $\fV=(\PPf)\PP$ as claimed.

If $f$ is convex, then so is $\PPf$ by
\Cref{roc:thm5.7:Af}, and so also is $\fV$ by
\Cref{roc:thm5.7:fA}.%
\indexg{shadow (of function)|)}
\end{proof}

\begin{theorem}  \label{thm:a10-nunu}
\indexg{reductions, astronic!lower semicontinuous hull@as lower semicontinuous hull|(}%
\indexg{shadow (of function)!astronic reductions and|(}%
Let $f:\Rn\rightarrow\Rext$,
let $\vv\in\resc{f}$,
let $g=\fshadv$ be the reduction of $f$ at~$\limray{\vv}$,
and let $\gtil=\fv$ be the shadow of $f$ along $\vv$.
Then $g=\lsc\gtil$.
Therefore, $g$ is lower semicontinuous,
and furthermore, if $f$ is convex then so is $g$.
\end{theorem}

\begin{proof}
For $\xx\in\Rn$, we have
\[
  g(\xx)
  =
  \fext(\limray{\vv}\plusl\xx)
  =
  \gext(\xx)
  =
  (\lsc g)(\xx),
\]
where the first equality is by $g$'s definition,
the second is by \Cref{thm:d4}(\ref{thm:d4:c}),
and the third is by \Cref{pr:h:1}(\ref{pr:h:1a}).
Thus, $g=\lsc g$, so $g$ is lower semicontinuous.

Since $\vv\in\resc{f}$,
$f(\xx+\lambda\vv)$ is nonincreasing as a function of $\lambda\in\R$,
for any fixed $\xx\in\Rn$.
Therefore, the expression for $\gtil$ from \eqref{eqn:gtil-defn} can be strengthened to
\begin{equation*}   %
  \gtil(\xx)=
  \lim_{\lambda\rightarrow+\infty} f(\xx+\lambda\vv),
\end{equation*}
implying
\[
   \gtil(\xx)=\lim f(\xx+t\vv)\ge\ef(\limray{\vv}\plusl\xx)=g(\xx),
\]
where the inequality is by definition of $\ef$
(Eq.~\ref{eq:e:7}) and since $\xx+t\vv\rightarrow\limray{\vv}\plusl\xx$
(by \Cref{pr:i:7}\ref{pr:i:7f}).
Thus, $g\leq \gtil$.
Since $g$ is lower semicontinuous, this implies
$g\leq \lsc \gtil$ by
\Cref{prop:lsc:characterize}(\ref{prop:lsc:characterize:b}).

For the reverse inequality, let $\xx\in\Rn$.
By \Cref{pr:d1},
there exists a sequence $\seq{\yy_t}$ in $\Rn$ such that
$\yy_t\rightarrow\limray{\vv}\plusl\xx$
and
$f(\yy_t)\rightarrow\fext(\limray{\vv}\plusl\xx)=g(\xx)$.
Then by \Cref{thm:lim-plusl-inv} (applied with $\xbar=\vv$ and
$\ybar=\xx$),
there exist a span-bound sequence $\seq{\vv_t}$ in $\Rn$,
a sequence $\seq{\xx_t}$ in $\Rn$, and a sequence $\seq{\lambda_t}$ in
$\R$ such that
$\yy_t=\lambda_t\vv_t+\xx_t$ for all $t$,
$\vv_t\rightarrow\vv$, and $\xx_t\rightarrow\xx$.
Further, since $\seq{\vv_t}$ is span-bound with limit $\vv$,
for each $t$, $\vv_t=\gamma_t \vv$ for some $\gamma_t\in\R$.
Thus,
\begin{align*}
  g(\xx)
  =
  \lim f(\yy_t)
  \geq
  \liminf \gtil(\yy_t)
  &=
  \liminf \gtil(\lambda_t \gamma_t \vv + \xx_t)
  \\
  &=
  \liminf \gtil(\xx_t)
  \geq
  (\lsc \gtil)(\xx).
\end{align*}
The first inequality is because $\gtil\leq f$, by $\gtil$'s definition.
The third equality is because $\gtil(\zz+\lambda\vv)=\gtil(\zz)$
for all $\zz\in\Rn$ and $\lambda\in\R$,
again by $\gtil$'s definition.
The last inequality is by definition of lower semicontinuous hull
(Eq.~\ref{eq:lsc:liminf}),
and since $\xx_t\rightarrow\xx$.

Thus, $g=\lsc\gtil$.
Therefore, if $f$ is convex, then $\gtil$ is as well
by \Cref{pr:gtil-is-PPfP-gen},
and so also is $g$ by \Cref{pr:lsc-props}(\ref{pr:lsc-props:a}).%
\indexg{reductions, astronic|)}%
\indexg{reductions, astronic!lower semicontinuous hull@as lower semicontinuous hull|)}%
\indexg{shadow (of function)!astronic reductions and|)}%
\end{proof}

\chapter{Astral conjugates and biconjugates}
\label{sec:conjugacy}

As discussed briefly in \Cref{sec:prelim:conjugate},
the {conjugate} of a function
$f:\Rn\rightarrow\Rext$ is the function
$\fstar:\Rn\rightarrow\Rext$ given, for $\uu\in\Rn$, by
\begin{equation}  \label{eq:fstar-def}
  \fstar(\uu) = \sup_{\xx\in\Rn} \bracks{\xx\cdot\uu - f(\xx)}.
\end{equation}
This is a centrally important notion in standard convex
analysis.
As we saw in \Cref{pr:epi-fstar-maj-affn},
the conjugate $\fstar$ encodes exactly those affine functions that are
majorized by~$f$.
Moreover, by
\Cref{pr:conj-props-cvx}(\ref{pr:conj-props-cvx:b}),
functions $f$ that are closed and convex are equal to their own biconjugate, that is, the conjugate of $\fstar$,
so that
$f=(\fstar)^* = \fdubs$.
This implies, in such cases,
that $f$ can be represented as the supremum over all affine
functions that it majorizes, and furthermore,
that the original function $f$ can be fully reconstructed from the dual
representation afforded by~$\fstar$.

In this chapter,
we extend these fundamental concepts to
astral space.

\section{Definitions and basic properties}
\label{sec:conjugacy-def}

\indexg{conjugate, primal (astral)|(}%
Let $F:\extspace\rightarrow\Rext$.
For now, we allow $F$ to be any function defined on
astral space, although later we will often focus on when
$F$ is the extension $\fext$ of some convex
function $f:\Rn\rightarrow\Rext$.

How can the definition of the conjugate $\fstar$ given in
\eqref{eq:fstar-def} be extended to a function $F$?
A natural idea is to simply replace $f$ by $F$ and $\xx$ by
$\xbar$ so that the expression inside the supremum
(now over $\xbar\in\extspace$) becomes
$\xbar\cdot\uu - F(\xbar)$.
The problem, of course, is that
$\xbar\cdot\uu$ and $F(\xbar)$ might both be $+\infty$ (or both
$-\infty$)
so the question arises how to resolve the ``ties'' between
infinite values in those cases.

To address this, we can re-express the conjugate in a way
that generalizes more easily.
\indexg{conjugate (standard)!re-expressed|(}%
In particular, since the epigraph of $f$ consists of all pairs
$\rpair{\xx}{y}$ in $\Rn\times\R$ with $y\geq f(\xx)$,
we can rewrite \eqref{eq:fstar-def} as
\begin{equation}  \label{eq:fstar-mod-def}
  \fstar(\uu) = \sup_{\rpair{\xx}{y}\in\epi f} \Bracks{\xx\cdot\uu - y}.%
\indexg{conjugate (standard)!re-expressed|)}
\end{equation}
This expression generalizes directly to
astral space by making the simple substitutions suggested earlier, and
specifically replacing $\xx\cdot\uu$ with its astral
analogue $\xbar\cdot\uu$.
Note that
although $\xbar\cdot\uu$ might be infinite, $y$ is always in $\R$, so
the earlier issue is sidestepped.

\begin{definition}  \label{def:ast-conjugate}
\indexg{conjugate, primal (astral)!defined|(}%
  The \emph{conjugate} of a function $F:\extspace\rightarrow\Rext$ is
  the function $\Fstar:\Rn\rightarrow\Rext$ defined,
  for $\uu\in\Rn$, as
\begin{equation}  \label{eq:Fstar-def}
\indexm{f 500}{$\Fstar$}{primal (astral) conjugate}%
  \Fstar(\uu)
  =
  \sup_{\rpair{\xbar}{y}\in\epi F} \Bracks{\xbar\cdot\uu - y}.
\end{equation}
\end{definition}
We sometimes also refer to $\Fstar$ as a \emph{primal} conjugate, to
emphasize its distinction from dual conjugates, which will be
introduced shortly.

For any $F:\extspace\rightarrow\Rext$,
the resulting conjugate function $\Fstar$ is always convex:

\begin{proposition}  \label{pr:conj-is-convex}
  Let $F:\extspace\rightarrow\Rext$.
  Then its conjugate $\Fstar$ is convex.
\end{proposition}

\begin{proof}
For any fixed $\xbar\in\extspace$ and $y\in\R$,
$\xbar\cdot\uu$, viewed as a function of $\uu\in\Rn$,
is convex by
\Crefequiv{thm:h:5}{thm:h:5a0}{thm:h:5b},
so $\xbar\cdot\uu-y$ is convex as well.
Therefore, $\Fstar$ is convex since it is the pointwise supremum of
convex functions
(\Cref{roc:thm5.5}).
\end{proof}

The definition for $\Fstar$ can be
rewritten using the downward sum operation
(as defined in Eq.~\ref{eq:down-add-def}).
In particular,
using \Cref{pr:plusd-props}(\ref{pr:plusd-props:d}),
we can rewrite \eqref{eq:Fstar-def} as
\begin{align}
  \Fstar(\uu)
  &=
  \sup_{\xbar\in\extspace}
       \BigBracks{
         \sup \bigBraces{\xbar\cdot\uu - y:\:y\in\R,\,y\geq F(\xbar)}
       }
  \notag
  \\
  &=
  \sup_{\xbar\in\extspace}
       \BigBracks{ - F(\xbar) \plusd \xbar\cdot\uu }.%
\indexg{conjugate, primal (astral)!defined|)}
  \label{eq:Fstar-down-def}
\end{align}
This answers the earlier question of how to resolve the ``ties'' between
infinite values in the standard definition.
Note also that $-F(\xbar)\plusd\xbar\cdot\uu$ is equal to $-\infty$ if
$-F(\xbar)$ and $\xbar\cdot\uu$ are not summable
(since then one of the terms must be $-\infty$),
and otherwise is equal to the ordinary
difference $\xbar\cdot\uu-F(\xbar)$.
Thus, \eqref{eq:Fstar-down-def} shows that $\Fstar(\uu)$ could
equivalently be expressed as the supremum of
$\xbar\cdot\uu-F(\xbar)$ over all
$\xbar\in\extspace$ for which this difference is defined.

In this form, the conjugate $\Fstar$
(as well as the dual conjugate to be defined shortly)
is the same as that defined by
\idxsinger\citet[Definition~8.2]{singer_book}
in his abstract treatment of convex analysis,
which can be instantiated to the astral setting
by letting his variables
$X$ and $W$ be equal to $\extspace$ and $\Rn$, respectively,
and defining his coupling function~$\singerphi$ to coincide
with ours, $\singerphi(\xbar,\uu)=\xbar\cdot\uu$ for $\xbar\in\extspace$ and
$\uu\in\Rn$.
The equivalence of Eqs.~(\ref{eq:Fstar-down-def}) and~(\ref{eq:Fstar-def}) is mentioned by
\idxsinger\citet[Eq.~8.34]{singer_book}.
This general form of conjugate was originally proposed by
\indexa{Moreau, J. J.}%
\citet[Eq.~14.7]{moreau__convexity}.

\begin{example}[Conjugate of an astral affine function]
\indexg{conjugate, primal (astral)!affine function@of affine function|(}%
\indexg{affine functions (astral)!primal conjugate of|(}%
Consider the affine function
$F(\xbar)=\xbar\cdot\ww+b$, for
$\xbar\in\extspace$, where $\ww\in\Rn$ and $b\in\R$.
Then
\[
  \Fstar(\uu)
  =
  \begin{cases}
    -b
    & \text{if $\uu=\ww$,}
  \\
    +\infty
    & \text{otherwise.}
  \end{cases}
\]
This can be shown directly using
\eqref{eq:Fstar-down-def}.
The function $\Fstar$ coincides with the conjugate of the
standard affine function
$f(\xx)=\xx\cdot\ww+b$, whose extension
is $\ef=F$.
Thus, $\fextstar=\fstar$.
We will show that this holds more generally in
\Cref{pr:fextstar-is-fstar} below.%
\indexg{conjugate, primal (astral)!affine function@of affine function|)}%
\indexg{affine functions (astral)!primal conjugate of|)}%
\end{example}

\begin{example}[Conjugate of an astral point indicator]
\label{ex:ind-astral-point}
\indexg{indicator functions (astral)!conjugate for singleton|(}%
Let $\zbar\in\extspace$.
As mentioned earlier, $\indfa{\zbar}$ is shorthand for
$\indfa{\set{\zbar}}$, the
astral indicator of the singleton $\set{\zbar}$, as defined in
\eqref{eq:indfa-defn}.
This function's conjugate is
$\indfastar{\zbar}(\uu)=\zbar\cdot\uu$, for $\uu\in\Rn$;
that is, $\indfastar{\zbar}=\ph{\zbar}$, where $\ph{\zbar}$ is the
function corresponding to evaluation of the coupling function in the second argument, introduced earlier in the construction of astral space (see Eq.~\ref{eq:ph-xbar-defn}).
Thus,
whereas the standard conjugate of any function
on $\Rn$
is always closed and convex,
this example shows that the astral conjugate
of a function
defined on $\extspace$ might not be closed or even lower
semicontinuous (although it must be convex by
\Cref{pr:conj-is-convex}).
For example, if $n=1$ and $\barz=+\infty$, then
$\indfastar{\barz}(u)=\limray{u}$, for $u\in\R$,
which
is not proper or closed or lower semicontinuous.%
\indexg{indicator functions (astral)!conjugate for singleton|)}%
\end{example}

The operation of taking the conjugate,
referred to as \emph{conjugation},
maps a function $F$ defined on $\extspace$ to the function
$\Fstar$ defined on
\indexg{conjugate, primal (astral)|)}%
$\Rn$.
\indexg{conjugate, dual (astral)|(}%
We next define
a dual operation
that maps in the
reverse direction, from functions on $\Rn$ to functions on
$\extspace$.
In standard convex analysis, both a function $f$ and its conjugate
$\fstar$ are defined on $\Rn$ so the same
conjugation operation
can be used in either direction.
But in the astral setting, as is the case more generally in the
abstract setting of
\idxsinger\citet{singer_book}, a different dual
conjugation
is required.
This asymmetry reflects the fact that the coupling function
$\xbar\cdot\uu$ is defined with its two arguments
$\xbar\in\extspace$ and $\uu\in\Rn$
belonging to different spaces.

Let $\psi: \Rn\rightarrow\Rext$ be any function on $\Rn$.
We use a Greek letter to emphasize
that we think of $\psi$ as operating on the dual variable $\uu\in\Rn$;
later, we will often take $\psi$ itself to
be the conjugate of some other function.
By direct analogy with the preceding definition of $\Fstar$,
we have:

\begin{definition}
\indexg{conjugate, dual (astral)!defined|(}%
The \emph{dual conjugate} of a function $\psi:\Rn\rightarrow\Rext$ is the function
$\psistarb:\extspace\rightarrow\Rext$
defined,
for $\xbar\in\extspace$,
as
\begin{equation}  \label{eq:psistar-def}
  \psistarb(\xbar)
  =
  \sup_{\rpair{\uu}{v}\in\epi \psi} \Bracks{\xbar\cdot\uu - v}.
\end{equation}
\end{definition}

By \Cref{pr:plusd-props}(\ref{pr:plusd-props:d}), this
definition can be stated equivalently as
\begin{equation}  \label{eq:psistar-def:2}
\indexm{psistar}{$\psistarb$}{dual (astral) conjugate}%
  \psistarb(\xbar)
  =
  \sup_{\uu\in\Rn} \bigBracks{ - \psi(\uu) \plusd \xbar\cdot\uu }.%
\indexg{conjugate, dual (astral)!defined|)}%
\end{equation}
We use the notation $\psistarb$ rather than $\psistar$ because
the latter denotes the standard conjugate of $\psi$,
so $\psistar$ is a function defined on $\Rn$ while $\psistarb$ is
defined on $\extspace$. These two kinds of conjugates
($\psistarb$ and $\psistar$) agree on $\Rn$:

\begin{proposition}
\label{pr:psistarb:psistar}
  Let $\psi:\Rn\rightarrow\Rext$.
  Then
  $\psistarb(\xx)=\psi^*(\xx)$ for all $\xx\in\Rn$.
\end{proposition}
\begin{proof}
The statement follows immediately from the definitions
of standard conjugate (Eq.~\ref{eq:fstar-def})
and dual conjugate (Eq.~\ref{eq:psistar-def:2}).
\end{proof}

As seen in Example~\ref{ex:ind-astral-point}, the astral (primal)
conjugate of an astral function need not be lower semicontinuous.
Nonetheless, the astral dual conjugate of any function is always lower
semicontinuous:

\begin{proposition}  \label{pr:dual-conj-lsc}
  Let $\psi:\Rn\rightarrow\Rext$.
  Then $\psistarb$ is lower semicontinuous.
\end{proposition}

\begin{proof}
By \Cref{eq:psistar-def}, $\psistarb$ is the pointwise supremum
of a family of functions of the form $\xbar\mapsto\xbar\inprod\uu-v$,
defined over $\xbar\in\eRn$,
with $\uu\in\Rn$ and $v\in\R$. Since these functions are lower semicontinuous
(by \Cref{ex:ext-affine}
combined with \Cref{prop:ext:F}\ref{prop:ext:F:a}),
so is $\psistarb$ (by \Cref{pr:lsc-sup}).%
\indexg{conjugate, dual (astral)|)}
\end{proof}

\indexg{biconjugate (astral)|(}%
We will be especially interested in the \emph{biconjugate}
of a function $F:\extspace\rightarrow\Rext$,
that is, $\Fdub=(\Fstar)^{\dualstar}$,
the dual conjugate of the conjugate of $F$.
In standard convex analysis, the biconjugate $\fdubs$ for a
function $f:\Rn\rightarrow\Rext$ is equal to its closure, $\cl f$, if
$f$ is convex
(\Cref{pr:conj-props-cvx}\ref{pr:conj-props-cvx:b}).
Thus, as already discussed, if $f$ is closed and convex, then
$f=\fdubs$.
Furthermore, $\fdubs$ is in general
equal to the pointwise supremum over all
affine functions that are majorized by~$f$.
An analogous result holds in the astral setting, as we show now.

For $\uu\in\Rn$ and $v\in\R$, let
$\affuv(\xbar)=\xbar\cdot\uu-v$, for $\xbar\in\extspace$.
As we saw in \Cref{ex:ext-affine}, this is the extension of the standard affine
function $\xx\mapsto\xx\cdot\uu-v$,
and is also
an instance of an astral affine map defined in
\Cref{sec:linear-maps}.
Then $\Fdub$ is exactly the pointwise supremum over all such functions
that are majorized by $F$, as we show in the next
theorem, which is a special case of
\idxsinger\citet[Theorem~8.5]{singer_book}.

\begin{theorem}  \label{thm:fdub-sup-afffcns}
  Let $F:\extspace\rightarrow\Rext$, and
  for $\uu\in\Rn$ and $v\in\R$, let $\affuv:\extspace\rightarrow\Rext$
  be defined by
  $\affuv(\xbar)=\xbar\cdot\uu-v$ for $\xbar\in\extspace$.
  \begin{letter-compact}
  \item    \label{thm:fdub-sup-afffcns:a}
    For all $\uu\in\Rn$ and $v\in\R$,
    we have
    $\rpair{\uu}{v}\in\epi\Fstar$
    if and only if
    $F\geq\affuv$.
  \item    \label{thm:fdub-sup-afffcns:b}
    For all $\xbar\in\extspace$,
    \begin{equation*} %
      \Fdub(\xbar) = \sup\bigBraces{\affuv(\xbar) :\: \uu\in\Rn,\,v\in\R,\,\affuv\leq F}.
    \end{equation*}
  \item    \label{thm:fdub-sup-afffcns:c}
    $F\geq \Fdub$.
  \item    \label{thm:fdub-sup-afffcns:d}
    $F=\Fdub$ if and only if $F$ is the pointwise supremum
    of some collection of affine functions $\affuv$.
  \end{letter-compact}
\end{theorem}

\begin{proof}
  ~
  
\begin{proof-parts}
\pfpart{Part~(\ref{thm:fdub-sup-afffcns:a}):}
Let $\uu\in\Rn$ and $v\in\R$.
Then
\begin{align*}
  \rpair{\uu}{v}\in\epi \Fstar
  \;\Leftrightarrow\;
  \Fstar(\uu)\leq v
  &\;\Leftrightarrow\;
  \xbar\cdot\uu - y \leq v
  \text{ for all } \rpair{\xbar}{y}\in \epi F
  \\
  &\;\Leftrightarrow\;
  \xbar\cdot\uu - v \leq y
  \text{ for all } \rpair{\xbar}{y}\in \epi F
  \\
  &\;\Leftrightarrow\;
  \xbar\cdot\uu - v \leq F(\xbar)
  \text{ for all } \xbar\in\extspace
  \\
  &\;\Leftrightarrow\;
  \affuv\leq F,
\end{align*}
where the second equivalence is  
by definition of conjugate (Eq.~\ref{eq:Fstar-def}).

\pfpart{Part~(\ref{thm:fdub-sup-afffcns:b}):}
This follows directly from
part~(\ref{thm:fdub-sup-afffcns:a})
combined with the definition of dual conjugate (Eq.~\ref{eq:psistar-def}) with
$\psi=\Fstar$.

\pfpart{Part~(\ref{thm:fdub-sup-afffcns:c}):}
This is immediate from
part~(\ref{thm:fdub-sup-afffcns:b}).

\pfpart{Part~(\ref{thm:fdub-sup-afffcns:d}):}
By part~(\ref{thm:fdub-sup-afffcns:b}), if
$F=\Fdub$ then $F$ is the pointwise supremum over a collection of
affine functions.
For the converse, suppose, for some set
$E\subseteq \Rn\times\R$, that
$F(\xbar) = \sup_{\rpair{\uu}{v}\in E} \affuv(\xbar)$
for $\xbar\in\extspace$.
If $\rpair{\uu}{v}\in E$, then this implies $\affuv\leq F$,
so that $\affuv\leq\Fdub$ by
part~(\ref{thm:fdub-sup-afffcns:b}).
Since this is true for every such pair in $E$, it follows that
$F\leq\Fdub$.
Combined with 
part~(\ref{thm:fdub-sup-afffcns:c}), this implies $F=\Fdub$.%
\indexg{biconjugate (astral)|)}%
\qedhere
\end{proof-parts}
\end{proof}

\indexg{biconjugate (astral)|(}%
Although we largely focus on the biconjugates $\Fdub$
(and later, $\fdub$, where $f$ is defined over~$\Rn$),
it is also possible to form a dual form of biconjugate
$\psidub=(\psi^{\bar{*}})^*$ from a function
$\psi:\Rn\rightarrow\Rext$ by applying
the conjugations
in the opposite order.
Analogous properties to those shown in
\Cref{thm:fdub-sup-afffcns} apply
to $\psidub$, as stated in the next \namecref{thm:psi-geq-psidub}.
The proof is exactly analogous to that of
\Cref{thm:fdub-sup-afffcns}, and is therefore omitted.

\begin{theorem}  \label{thm:psi-geq-psidub}
  Let $\psi:\Rn\rightarrow\Rext$, and
  for $\xbar\in\extspace$ and $y\in\R$,
  let $\daffxy:\Rn\rightarrow\Rext$
  be defined by
  $\daffxy(\uu)=\xbar\cdot\uu-y$ for $\uu\in\Rn$.
  \begin{letter-compact}
  \item    \label{thm:psi-geq-psidub:a}
    For all $\xbar\in\extspace$ and $y\in\R$,
    we have
    $\rpair{\xbar}{y}\in\epi\psistarb$
    if and only if
    $\psi\geq\daffxy$.

  \item    \label{thm:psi-geq-psidub:b}
    For all $\uu\in\Rn$,
    \begin{equation*}
      \psidub(\uu)
      =
      \sup\bigBraces{
        \daffxy(\uu)
        :\:
        \xbar\in\extspace,\,y\in\R,\,\daffxy\leq \psi
      }.
    \end{equation*}
  \item    \label{thm:psi-geq-psidub:c}
    $\psi\geq \psidub$.
  \item    \label{thm:psi-geq-psidub:d}
    $\psi=\psidub$ if and only if $\psi$
    is the pointwise supremum
    of some collection of functions $\daffxy$.
  \end{letter-compact}
\end{theorem}

\begin{example}[Dual conjugate of $\ph{\zbar}$.]
\label{ex:dual:ph:zbar}
\indexg{functional representation of astral points!dual conjugate of|(}%
Let $\zbar\in\extspace$, and let $\ph{\zbar}:\Rn\to\eR$
be defined by
$\ph{\zbar}(\uu)=\zbar\cdot\uu$
for $\uu\in\Rn$, as
introduced earlier in the construction of astral space (see Eq.~\ref{eq:ph-xbar-defn}). Its dual conjugate is
\[
  \phstar{\zbar}(\xbar)
  =
  \sup_{\uu\in\Rn} \Bracks{ - \zbar\cdot\uu \plusd \xbar\cdot\uu }
  =
  \indfa{\zbar}^{\phantomdualstar}(\xbar),
\]
with $\indfa{\zbar}$ as in \Cref{ex:ind-astral-point}.
To see this, note that $-\zbar\cdot\uu\plusd\zbar\cdot\uu$ is equal
to~$0$ if $\zbar\cdot\uu\in\R$ (including if $\uu=\zero$), and otherwise
is equal to $-\infty$;
thus, $\phstar{\zbar}(\zbar)=0$.
And if $\xbar\neq\zbar$ then, for some $\uu\in\Rn$,
$\xbar\cdot\uu\neq\zbar\cdot\uu$
(by \Cref{pr:i:4}), which means that $\xbar\inprod\uu$
and $-\zbar\inprod\uu$ must be summable,
and $-\zbar\inprod\uu+\xbar\inprod\uu\ne 0$.
Therefore, the expression
$- \zbar\cdot(\lambda\uu) \plusd \xbar\cdot(\lambda\uu)
= \lambda(- \zbar\cdot\uu +\xbar\cdot\uu )$
can be made arbitrarily large for an appropriate choice of $\lambda\in\R$,
implying $\phstar{\zbar}(\xbar)=+\infty$.
Thus, $\phstar{\zbar}=\indfa{\zbar}^{\phantomdualstar}$. Also, as we saw in \Cref{ex:ind-astral-point}, $\indfastar{\zbar}=\ph{\zbar}^{\phantomdualstar}$, and hence
$\indfadub{\zbar}=\indfa{\zbar}^{\phantomdualstar}$
and $\phdub{\zbar}=\ph{\zbar}^{\phantomdualstar}$.

If $\zbar$ is not in $\Rn$,
then $\ph{\zbar}$'s standard conjugate is
$\phstars{\zbar}\equiv+\infty$ (since it can be shown using
\Cref{pr:i:3}
that $\ph{\zbar}(\uu)=-\infty$ for some $\uu\in\Rn$).
Thus, standard
conjugation
entirely erases the identity of $\zbar$.
On the other hand, from the astral dual conjugate
$\phstar{\zbar}$, we have just seen that
it is possible to reconstruct
$\ph{\zbar}^{\phantomdualstar}=\phdub{\zbar}$, and thus recover $\zbar$.
This shows that the astral dual conjugate $\psistarb$ can
retain more information about the original function
$\psi:\Rn\rightarrow\Rext$
than the standard conjugate $\psistar$.%
\indexg{functional representation of astral points!dual conjugate of|)}%
\indexg{biconjugate (astral)|)}%
\end{example}

\indexg{closedness, astral (primal)|(}%
When a function $F:\extspace\rightarrow\Rext$ is equal to its own
biconjugate, so that $F=\Fdub$, we say that $F$ is
\emph{astral closed} (or sometimes just closed).
We also say that $F$ is astral closed at a point $\xbar\in\extspace$
if $F(\xbar)=\Fdub(\xbar)$.
\indexg{closedness, astral dual|(}%
Likewise, for a function $\psi:\Rn\rightarrow\Rext$, we say that
$\psi$ is \emph{astral dual closed}
(or just dual closed) if $\psi=\psidub$, and that $\psi$ is astral
dual closed at a point $\uu\in\Rn$ if $\psi(\uu)=\psidub(\uu)$.%
\indexg{closedness, astral dual|)}

Here are some straightforward properties of biconjugates, including
simple characterizations of when a function is (dual) closed.
We state these
in separate propositions for the primal and dual forms.
Nevertheless, they should be regarded as
if proved in parallel, beginning with a proof of
part~(\ref{pr:primal-biconj-props:a})
for both propositions,
followed by a proof of
part~(\ref{pr:primal-biconj-props:b}) for both,
and so on.
In this way, the seemingly forward reference in the proof of
\Cref{pr:primal-biconj-props}(\ref{pr:primal-biconj-props:c})
to
\Cref{pr:dual-biconj-props}(\ref{pr:dual-biconj-props:b})
can be seen to be instead to a part that was already proved.

\begin{proposition}  \label{pr:primal-biconj-props}
\indexg{conjugate, primal (astral)|(}%
  Let $F:\extspace\rightarrow\Rext$.
  \begin{letter-compact}
  \item   \label{pr:primal-biconj-props:a}
    Let $G:\extspace\rightarrow\Rext$ and suppose $F\leq G$.
    Then $\Fstar\geq\Gstar$.
  \item   \label{pr:primal-biconj-props:b}
    $F^{*\dualstar *} = \Fstar$.
  \item   \label{pr:primal-biconj-props:c}
    $F=\Fdub$ if and only if $F=\psistarb$ for some
    $\psi:\Rn\rightarrow\Rext$.
  \end{letter-compact}
\end{proposition}

\begin{proof}
~
\begin{proof-parts}
\pfpart{Part~(\ref{pr:primal-biconj-props:a}):}
This follows directly from the definition of astral conjugate given
in
\eqref{eq:Fstar-down-def}
(or from Eq.~\ref{eq:Fstar-def} since $\epi{G}\subseteq\epi{F}$).

\pfpart{Part~(\ref{pr:primal-biconj-props:b}):}
By
\Cref{thm:fdub-sup-afffcns}(\ref{thm:fdub-sup-afffcns:c}),
$\Fdub\leq F$, implying that
$F^{*\dualstar *} = (\Fdub)^* \geq \Fstar$
by part~(\ref{pr:primal-biconj-props:a}).
On the other hand, by
\Cref{thm:psi-geq-psidub}(\ref{thm:psi-geq-psidub:c})
(with $\psi=\Fstar$),
$F^{*\dualstar *} = (\Fstar)^{\dualstar *} \leq \Fstar$.

\pfpart{Part~(\ref{pr:primal-biconj-props:c}):}
If $F=\Fdub$ then $F=\psistarb$ with $\psi=\Fstar$.
Conversely, if $F=\psistarb$ for some $\psi:\Rn\rightarrow\Rext$, then
${\Fdub=(\psistarb)^{*\dualstar}=\psi^{\dualstar * \dualstar}=\psistarb=F}$
with the third equality from
\Cref{pr:dual-biconj-props}(\ref{pr:dual-biconj-props:b}).%
\indexg{closedness, astral (primal)|)}%
\indexg{conjugate, primal (astral)|)}
\qedhere
\end{proof-parts}
\end{proof}

\begin{proposition}  \label{pr:dual-biconj-props}
\indexg{closedness, astral dual|(}%
  Let $\psi:\Rn\rightarrow\Rext$.
  \begin{letter-compact}
  \item   \label{pr:dual-biconj-props:a}
    Let $\rho:\Rn\rightarrow\Rext$.
    If $\psi\leq \rho$ then $\psistarb\geq\rhostarb$.
  \item   \label{pr:dual-biconj-props:b}
    $\psi^{\dualstar * \dualstar} = \psistarb$.
  \item   \label{pr:dual-biconj-props:c}
    $\psi=\psidub$ if and only if $\psi=\Fstar$ for some
    $F:\extspace\rightarrow\Rext$.
  \end{letter-compact}
\end{proposition}

\begin{proof}
  Analogous to \Cref{pr:primal-biconj-props}.%
\indexg{closedness, astral dual|)}
\end{proof}

\indexg{conjugate, primal (astral)|(}%
\indexg{conjugate, dual (astral)|(}%
Propositions~\ref{pr:primal-biconj-props}(\ref{pr:primal-biconj-props:c})
and~\ref{pr:dual-biconj-props}(\ref{pr:dual-biconj-props:c})
together imply that the conjugation operation $F\mapsto\Fstar$
maps astral closed functions $F:\extspace\rightarrow\Rext$
(for which $F=\Fdub$) to
astral dual closed functions $\psi:\Rn\rightarrow\Rext$
(for which $\psi=\psidub$).
Moreover, this map is a bijection, with inverse given by the dual
conjugation operation, $\psi\mapsto\psistarb$.
Thus, there is a one-to-one correspondence, defined by the conjugation
operations, between astral closed and astral dual closed functions.%
\indexg{conjugate, primal (astral)|)}%
\indexg{conjugate, dual (astral)|)}%

\indexg{closedness, astral dual!standard closedness and|(}%
We also note that astral dual closedness is implied by standard
closedness:

\begin{proposition}  \label{pr:std-clos-implies-ast-dual-clos}
  Let $\psi:\Rn\rightarrow\Rext$, and suppose
  $\psi=\psidubs$.
  Then
  $\psi=\psidub$.
\end{proposition}

\begin{proof}
By assumption, and
from \eqref{eq:fstar-mod-def}, we have, for $\uu\in\Rn$, that
\[
  \psi(\uu)
  =
  \sup_{\rpair{\xx}{y}\in\epi\psistar} \Bracks{\xx\cdot\uu - y}.
\]
Therefore, $\psi$ is the pointwise supremum of a collection of
functions $\daffsxy$ as in \Cref{thm:psi-geq-psidub};
hence, $\psi=\psidub$ by part~(\ref{thm:psi-geq-psidub:d})
of that \namecref{thm:psi-geq-psidub}.%
\indexg{closedness, astral dual!standard closedness and|)}%
\end{proof}

\indexg{conjugate, primal (astral)|(}%
The next proposition relates a function's lower
semicontinuous hull and its conjugate to the function's conjugate and
biconjugate:

\begin{proposition}  \label{pr:conj-lsc-props}
  Let $F:\extspace\rightarrow\Rext$.
  Then
  $F\geq (\lsc F) \geq \Fdub$
  and
  $\Fstar = (\lsc F)^*$.
\end{proposition}    

\begin{proof}
\Cref{thm:fdub-sup-afffcns}(\ref{thm:fdub-sup-afffcns:c})
shows that $F\geq\Fdub$,
and \Cref{pr:dual-conj-lsc} implies that $\Fdub$ is lower
semicontinuous.
Therefore, $F\geq (\lsc F) \geq \Fdub$
by
\Cref{prop:lsc:characterize}(\ref{prop:lsc:characterize:a},\ref{prop:lsc:characterize:b}).

Consequently,
$\Fstar\leq (\lsc F)^* \leq F^{*\dualstar *} = \Fstar$
by \Cref{pr:primal-biconj-props}(\ref{pr:primal-biconj-props:a},\ref{pr:primal-biconj-props:b}).%
\indexg{conjugate, primal (astral)|)}%
\end{proof}

\section{Conjugates and biconjugates of extensions}
\label{sec:conj-exts}

\indexg{conjugate, primal (astral)!extensions@of extensions|(}%
\indexg{lower semicontinuous extension!conjugate of|(}%
We next briefly consider conjugates and biconjugates of extensions of a
function $f:\Rn\rightarrow\Rext$.
These will be treated in much greater depth in
\Cref{sec:ast-close-extensions}.
For now, we establish some basic facts.
First, we show that the astral conjugate of $\fext$, denoted
$\fextstar=\bigParens{\fext\,}^{\!*}$,
is equal to the (standard)
conjugate of $f$,
and the same is true for the conjugate of $f$'s trivial extension.

\begin{proposition}  \label{pr:fextstar-is-fstar}
\indexg{trivial extension!conjugate of|(}%
  Let $f:\Rn\rightarrow\Rext$.
  Then
  $\fextstar=\ftrivstar=\fstar$.
\end{proposition}

\begin{proof}
First, $\fext=\lsc\ftriv$
(by \Cref{pr:lsc-ftriv-is-fext}),
implying
$\fextstar = (\lsc\ftriv)^* = \ftrivstar$,
with the second equality from
\Cref{pr:conj-lsc-props} (applied to $F=\ftriv$).

That $\ftrivstar(\uu)=\fstar(\uu)$, for $\uu\in\Rn$, then
follows by comparison
of conjugate definitions, Eqs.~(\ref{eq:fstar-mod-def})
and~(\ref{eq:Fstar-def}), and since $\epi\ftriv=\epi f$.%
\indexg{conjugate, primal (astral)!extensions@of extensions|)}%
\indexg{lower semicontinuous extension!conjugate of|)}%
\indexg{trivial extension!conjugate of|)}%
\end{proof}

\indexg{biconjugate (astral)!extension@of extension|(}%
\indexg{lower semicontinuous extension!biconjugate of|(}%
Applied to the extension $\fext$ of a function
$f:\Rn\rightarrow\Rext$,
this shows
that $\fext$'s
biconjugate is
\begin{equation}  \label{eqn:fdub-form-lower-plus}
  \fextdub(\xbar)
  =
  \fdub(\xbar)
  =
  \sup_{\uu\in\Rn} \bigBracks{-\fstar(\uu) \plusd \xbar\cdot\uu}
\end{equation}
for $\xbar\in\extspace$, where the second equality is from \eqref{eq:psistar-def:2}.
This expression is very close in form to the standard biconjugate
$\fdubs$, and shows that $\fdubs(\xx)=\fdub(\xx)$ for all $\xx\in\Rn$.

In several later proofs, we decompose $\xbar\in\extspace$ in
\Cref{eq:psistar-def:2,eqn:fdub-form-lower-plus} as
$\xbar={\limray{\vv}\plusl\zbar}$ for some $\vv\in\Rn$ and
$\zbar\in\extspace$.
In such cases,
it will be convenient to use variants of those equations with
leftward addition instead of downward addition.
\indexg{conjugate, dual (astral)!expressed with leftward addition|(}%
The next proposition provides sufficient conditions for rewriting \eqref{eq:psistar-def:2} using leftward addition,
which will be satisfied, for instance, when $\psi$ is convex and closed.

\begin{proposition}  \label{pr:psi-with-plusl}
  Let $\psi:\Rn\rightarrow\Rext$, and assume that either
  $\psi>-\infty$ or $\psi(\zero)=-\infty$.
  Then for $\xbar\in\extspace$,
  \begin{equation}    \label{eq:pr:psi-with-plusl:1}
    \psistarb(\xbar)
    =
    \sup_{\uu\in\Rn} \bigBracks{ - \psi(\uu) \plusl \xbar\cdot\uu }.
  \end{equation}
\end{proposition}

\begin{proof}
Note that $\barx\plusd \bary = \barx\plusl \bary$ for all $\barx,\bary\in\Rext$, except if
$\barx=+\infty$ and $\bary=-\infty$.
Therefore, if $\psi>-\infty$ then
$ - \psi(\uu) \plusd \xbar\cdot\uu = - \psi(\uu) \plusl \xbar\cdot\uu $
for all $\uu\in\Rn$ and $\xbar\in\extspace$, proving the claim in this
case.

Otherwise, if $\psi(\zero)=-\infty$ then for all $\xbar\in\extspace$,
\[
  -\psi(\zero)\plusd \xbar\cdot\zero
  = +\infty
  = -\psi(\zero)\plusl \xbar\cdot\zero,
\]
implying that both $\psistarb(\xbar)$ and
the right-hand side of \eqref{eq:pr:psi-with-plusl:1} are equal to
$+\infty$, proving the claim in this case as well.%
\indexg{conjugate, dual (astral)!expressed with leftward addition|)}%
\end{proof}

The next theorem summarizes results for $\fext$ and its biconjugate:

\begin{theorem}  \label{thm:fext-dub-sum}
  Let $f:\Rn\rightarrow\Rext$.
  Then
  \begin{letter-compact}
    \item  \label{thm:fext-dub-sum:a}
      $\fext\geq\fextdub = \fdub$.
    \item  \label{thm:fext-dub-sum:b}
      For all $\xbar\in\extspace$,
      \[
       \fdub(\xbar)
       =
       \sup_{\uu\in\Rn} \bigBracks{-\fstar(\uu) \plusd \xbar\cdot\uu}
       =
       \sup_{\uu\in\Rn} \bigBracks{-\fstar(\uu) \plusl \xbar\cdot\uu}.
      \]
    \item  \label{thm:fext-dub-sum:fdubs}
      For all $\xx\in\Rn$, $\fdub(\xx)=\fdubs(\xx)$.
  \end{letter-compact}
\end{theorem}

\begin{proof}
~
\begin{proof-parts}
\pfpart{Part~(\ref{thm:fext-dub-sum:a}):}
This is direct from
\Cref{thm:fdub-sup-afffcns}(\ref{thm:fdub-sup-afffcns:c})
(applied to $\fext$) and
\Cref{pr:fextstar-is-fstar}.

\pfpart{Part~(\ref{thm:fext-dub-sum:b}):}
The first equality is from
\eqref{eq:psistar-def:2}.
The second equality holds by \Cref{pr:psi-with-plusl}
since $f^*$ is closed (by \Cref{pr:conj-props}\ref{pr:conj-props:d}).

\pfpart{Part~(\ref{thm:fext-dub-sum:fdubs}):}
This is immediate from 
\Cref{pr:psistarb:psistar}
(with $\psi=\fstar$).%
\indexg{biconjugate (astral)!extension@of extension|)}%
\indexg{lower semicontinuous extension!biconjugate of|)}%
\qedhere
\end{proof-parts}
\end{proof}

\section{Functional operations involving linear maps}
\label{sec:conjs-lin-fctl-ops}

We next explore properties of conjugates of functions constructed in
particular ways using linear maps.

Recall from \eqref{eq:fA-defn} that
for a function $f:\Rm\rightarrow\Rext$ and matrix $\A\in\Rmn$,
we write $\fA$ for the composition of $f$ with the linear map
associated with $\A$, that is, for the function given by
$(\fA)(\xx)=f(\A\xx)$ for $\xx\in\Rn$.
We extend this notation straightforwardly to astral space:
\indexg{linear map in composition with astral function|(}%
For a function
$F:\extspac{m}\rightarrow\Rext$, we define the function
$\FA:\extspace\rightarrow\Rext$ by
\begin{equation}  \label{eq:FA-dfn}
\indexm{f a 300}{$\FA$}{composition with linear map (astral)}%
  (\FA)(\xbar)=F(\A\xbar)
\end{equation}
for $\xbar\in\extspace$.
(For clarity, we sometimes also write $\fA$ or $\FA$
as $f[\A]$ or $F[\A]$.)

Likewise,
recall from \eqref{eq:lin-image-fcn-defn} that
for a matrix $\A\in\Rmn$ and a function $f:\Rn\rightarrow\Rext$,
the image of $f$ under $\A$, denoted $\A f$, is defined, for any $\xx\in\Rm$, as
\begin{equation*}  %
  (\A f)(\xx) =  \inf\,\bigBraces{f(\zz):\:\zz\in\Rn,\,\A\zz=\xx}.
\end{equation*}
This operation also generalizes straightforwardly
to astral functions:
\indexg{linear image of astral function|(}%
Thus,
for $F:\extspace\rightarrow\Rext$ and $\A\in\Rmn$,
we define the
\emph{image of $F$ under $\A$}
to be the function $\A F:\extspac{m}\rightarrow\Rext$ defined,
for $\xbar\in\extspac{m}$, as
\begin{equation}  \label{eqn:image-F-dfn}
\indexm{a f 300}{$\A F$}{linear image of function (astral)}%
  (\A F)(\xbar) = \inf\,\bigBraces{F(\zbar) :\: \zbar\in\extspace,\, \A\zbar=\xbar}.
\end{equation}

In this section, we consider properties of functions constructed in
this way, focusing especially on connections to conjugates.

First,
here are some straightforward facts, analogous to
\Cref{pr:std-lin-img-props}:

\begin{proposition}  \label{pr:ast-lin-img-props}
  Let $F:\extspace\rightarrow\Rext$ and let $\A\in\Rmn$.
  \begin{letter-compact}
  \item   \label{pr:ast-lin-img-props:inf-FA}
    $F\ktransA\geq\inf F$.
  \item   \label{pr:ast-lin-img-props:a}
    $(\A F) \A \leq F$.
  \item   \label{pr:ast-lin-img-props:inf-AF}
    $\inf (\A F) = \inf F$.
  \item   \label{pr:ast-lin-img-props:b}
    Let $G:\extspace\rightarrow\Rext$ and suppose $F\leq G$.
    Then
    $\A F \leq \A G$
    and
    $F \ktransA \leq G \transA$.
  \item   \label{pr:ast-lin-img-props:c}
    Let $\B\in\R^{k\times m}$.
    Then $\B(\A F) = (\B\A) F$.
  \end{letter-compact}
\end{proposition}

\begin{proof}
Analogous to the proof of \Cref{pr:std-lin-img-props},
but now using the definitions in
Eqs.~(\ref{eq:FA-dfn}) and~(\ref{eqn:image-F-dfn}).%
\indexg{linear image of astral function|)}%
\indexg{linear map in composition with astral function|)}%
\end{proof}

\indexg{linear map in composition with astral function!conjugate identities and inequalities|(}%
\indexg{linear map in composition with astral function!astral closedness preserved|(}%
\indexg{linear image of astral function!conjugate identities|(}%
\indexg{closedness, astral (primal)!preserved by composition with linear map|(}%
\indexg{conjugate, primal (astral)!image under linear map@of image under linear map|(}%
The next proposition gives three conjugate identities.
Part~(\ref{pr:prml-lin-op-conj-idens:a})
shows that the formula for the (standard) conjugate of the linear
image of a standard function given in \Cref{roc:thm16.3:Af} holds also
for astral conjugates.
As a consequence,
in part~(\ref{pr:prml-lin-op-conj-idens:b}),
we obtain formulas for the conjugate of a function
(or really the biconjugate of the function)
composed with a linear map.
Finally,
in part~(\ref{pr:prml-lin-op-conj-idens:c}),
we see that if a function is astral closed, then so is its composition
with a linear map.

\Cref{pr:dual-lin-op-conj-idens} gives 
a dual analogue of
\Cref{pr:prml-lin-op-conj-idens}.
These two propositions should 
be regarded as if proved in parallel,
similar to \Cref{pr:primal-biconj-props,pr:dual-biconj-props}.
The same also applies to
\Cref{pr:conj-AF-props,pr:dual-conj-Apsi-props},
which follow afterwards.

\begin{proposition}   \label{pr:prml-lin-op-conj-idens}
  Let $F:\extspac{m}\rightarrow\Rext$
  and let $\A\in\Rmn$.
  Then:
  \begin{letter-compact}
  \item    \label{pr:prml-lin-op-conj-idens:a}
    $(\transAk F)^* = \Fstar\kernA$.
  \item    \label{pr:prml-lin-op-conj-idens:b}
    $(\Fdub \kernA)^* = (\transAk \Fstar)^{\dualstar *}$.
  \item    \label{pr:prml-lin-op-conj-idens:c}
    If $F=\Fdub$ then $\FA = (\FA)^{* \dualstar}$.
  \end{letter-compact}
\end{proposition}

\begin{proof}
  ~

\begin{proof-parts}
\pfpart{Part~(\ref{pr:prml-lin-op-conj-idens:a}):}
Let $\uu\in\Rn$.
Then
\begin{align*}
  \Fstar(\A\uu)
  &=
  \sup_{\zbar\in\extspac{m}}
      \BigBracks{-F(\zbar) \plusd \zbar\cdot(\A\uu)}
  \\  
  &=
  \sup_{\zbar\in\extspac{m}}
      \BigBracks{-F(\zbar) \plusd (\transA\zbar)\cdot\uu}
  \\  
  &=
  \sup_{\xbar\in\extspace} \BigBracks{
    \sup \bigBraces{-F(\zbar) \plusd \xbar\cdot\uu :\:
                   \zbar\in\extspac{m},\, \transA\zbar=\xbar}
  }
  \\  
  &=
  \sup_{\xbar\in\extspace} \BigBracks{
    \sup \bigBraces{-F(\zbar) :\: \zbar\in\extspac{m},\, \transA\zbar=\xbar}
    \plusd \xbar\cdot\uu
    }
  \\  
  &=
  \sup_{\xbar\in\extspace} \BigBracks{
    -\inf\,\bigBraces{F(\zbar) :\: \zbar\in\extspac{m},\, \transA\zbar=\xbar}
    \plusd \xbar\cdot\uu
    }
  \\  
  &=
  \sup_{\xbar\in\extspace}
    \BigBracks{ - (\transAk F)(\xbar) \plusd \xbar\cdot\uu}
  =
  (\transAk F)^*(\uu).
\end{align*}  
The first and last equalities are by definition of conjugate
(Eq.\ref{eq:Fstar-down-def}).
The second equality is by
\Cref{thm:Ax-dot-u}.
The fourth is by
\Cref{pr:plusd-props}(\ref{pr:plusd-props:d-gen}).
The sixth is by $\transAk F$'s definition
(Eq.~\ref{eqn:image-F-dfn}).

\pfpart{Part~(\ref{pr:prml-lin-op-conj-idens:b}):}
From
\Cref{pr:dual-lin-op-conj-idens}(\ref{pr:dual-lin-op-conj-idens:a}),
applied with $\psi=\Fstar$, we have
\begin{equation}   \label{eq:pr:prml-lin-op-conj-idens:1}
  (\transAk\Fstar)^{\dualstar} = \Fdub\kernA.
\end{equation}
Taking (primal) conjugate of both sides then yields the claim.

\pfpart{Part~(\ref{pr:prml-lin-op-conj-idens:c}):}
Suppose $F=\Fdub$.
By \eqref{eq:pr:prml-lin-op-conj-idens:1},
$\FA = \Fdub\kernA = (\transAk \Fstar)^{\dualstar}$,
which shows that $\FA$ is the astral dual conjugate of some function
(namely, $\transAk \Fstar$).
Therefore,
$\FA = (\FA)^{* \dualstar}$
by
\Cref{pr:primal-biconj-props}(\ref{pr:primal-biconj-props:c}).%
\indexg{linear map in composition with astral function!conjugate identities and inequalities|)}%
\indexg{linear map in composition with astral function!astral closedness preserved|)}%
\indexg{linear image of astral function!conjugate identities|)}%
\indexg{closedness, astral (primal)!preserved by composition with linear map|)}%
\indexg{conjugate, primal (astral)!image under linear map@of image under linear map|)}%
\qedhere
\end{proof-parts}
\end{proof}

\begin{proposition}   \label{pr:dual-lin-op-conj-idens}
\indexg{linear map in composition with standard function!astral conjugate identities and inequalities|(}%
\indexg{linear map in composition with standard function!astral dual closedness preserved|(}%
\indexg{linear image of standard function!astral conjugate identities|(}%
\indexg{closedness, astral dual!preserved by composition with linear map|(}%
\indexg{conjugate, dual (astral)!image under linear map@of image under linear map|(}%
  Let $\psi:\Rm\rightarrow\Rext$
  and let $\A\in\Rmn$.
  Then:
  \begin{letter-compact}
  \item    \label{pr:dual-lin-op-conj-idens:a}
    $(\transA \psi)^{\dualstar} = \psistarb\kernA$.
  \item    \label{pr:dual-lin-op-conj-idens:b}
    $(\psidub \kernA)^{\dualstar} = (\transA \psistarb)^{* \dualstar}$.
  \item    \label{pr:dual-lin-op-conj-idens:c}
    If $\psi=\psidub$ then $\psi\A = (\psi\A)^{\dualstar *}$.
  \end{letter-compact}
\end{proposition}

\begin{proof}
  Analogous to \Cref{pr:prml-lin-op-conj-idens}.%
\indexg{linear map in composition with standard function!astral conjugate identities and inequalities|)}%
\indexg{linear map in composition with standard function!astral dual closedness preserved|)}%
\indexg{linear image of standard function!astral conjugate identities|)}%
\indexg{closedness, astral dual!preserved by composition with linear map|)}%
\indexg{conjugate, dual (astral)!image under linear map@of image under linear map|)}%
\end{proof}

The analogues of
Propositions~\ref{pr:prml-lin-op-conj-idens}(\ref{pr:prml-lin-op-conj-idens:c})
and~\ref{pr:dual-lin-op-conj-idens}(\ref{pr:dual-lin-op-conj-idens:c}) do not hold
for $\transAk F$ and $\transA\psi$ without further assumptions, as we
investigate in more detail in
\Cref{sec:ast-close-primal-biconj}
(in \Cref{cor:AF-ast-closed}
and in \Cref{ex:AF-not-ast-closed,ex:Af:not-closed}).

\indexg{linear map in composition with astral function!conjugate identities and inequalities|(}%
The next two propositions give useful inequalities relating various
conjugates and biconjugates.
The first proposition is for primal conjugates, the other
for dual conjugates.

\begin{proposition}   \label{pr:conj-AF-props}
  Let $F:\extspac{m}\rightarrow\Rext$ and let $\A\in\Rmn$.
  Let $G:\extspace\rightarrow\Rext$ and assume that $G\geq \FA$.
  Then:
  \begin{letter-compact}
  \item    \label{pr:conj-AF-props:b}
    $\Gstar\leq\transAk \Fstar$.
  \item    \label{pr:conj-AF-props:c}
    $\Fstar\geq\Gstar\ktransA$.
  \item    \label{pr:conj-AF-props:d}
    $\Gdub\geq\Fdub\kernA$.
  \end{letter-compact}
\end{proposition}

\begin{proof}
  ~

\begin{proof-parts}
\pfpart{Part~(\ref{pr:conj-AF-props:b}):}
First,
\begin{equation}   \label{eq:pr:conj-AF-props:1}
  G\geq \FA \geq \Fdub\kernA,
\end{equation}
where the second inequality is by
\Cref{pr:ast-lin-img-props}(\ref{pr:ast-lin-img-props:b})
since $F\geq\Fdub$
(by \Cref{thm:fdub-sup-afffcns}\ref{thm:fdub-sup-afffcns:c}).
We then have
\[
   \Gstar
   \leq
   (\Fdub \kernA)^*
   =
   (\transAk \Fstar)^{\dualstar *}
   \leq
   \transAk \Fstar.
\]
The first inequality is by
\eqref{eq:pr:conj-AF-props:1} combined with
\Cref{pr:primal-biconj-props}(\ref{pr:primal-biconj-props:a}).
The equality is by
\Cref{pr:prml-lin-op-conj-idens}(\ref{pr:prml-lin-op-conj-idens:b}).
The last inequality is by
\Cref{thm:psi-geq-psidub}(\ref{thm:psi-geq-psidub:c}).

\pfpart{Part~(\ref{pr:conj-AF-props:c}):}
We have
\[
  \Gstar \ktransA
  \leq
  (\transAk \Fstar) \transA
  \leq
  \Fstar.
\]
The first inequality is by
part~(\ref{pr:conj-AF-props:b})
combined with
\Cref{pr:std-lin-img-props}(\ref{pr:std-lin-img-props:b}).
The second is by
\Cref{pr:std-lin-img-props}(\ref{pr:std-lin-img-props:a}).

\pfpart{Part~(\ref{pr:conj-AF-props:d}):}
By part~(\ref{pr:conj-AF-props:c}),
$\Fstar\geq\Gstar\ktransA$, implying
$\Gdub\geq\Fdub\kernA$ by 
\Cref{pr:dual-conj-Apsi-props}(\ref{pr:dual-conj-Apsi-props:c})
(with $\rho$, $\psi$ and $\A$, as they appear in that proposition,
set to $\Fstar$, $\Gstar$ and $\transA$,
\indexg{linear map in composition with astral function!conjugate identities and inequalities|)}%
respectively).
\qedhere
\end{proof-parts}
\end{proof}

\begin{proposition}   \label{pr:dual-conj-Apsi-props}
\indexg{linear map in composition with standard function!astral conjugate identities and inequalities|(}%
  Let $\psi:\Rm\rightarrow\Rext$ and let $\A\in\Rmn$.
  Let $\rho:\Rn\rightarrow\Rext$ and assume that $\rho\geq \psi\A$.
  Then:
  \begin{letter-compact}
  \item    \label{pr:dual-conj-Apsi-props:b}
    $\rhostarb\leq\transA \psistarb$.
  \item    \label{pr:dual-conj-Apsi-props:c}
    $\psistarb\geq\rhostarb\ktransA$.
  \item    \label{pr:dual-conj-Apsi-props:d}
    $\rhodub\geq\psidub\kernA$.
  \end{letter-compact}
\end{proposition}

\begin{proof}
  Analogous to \Cref{pr:conj-AF-props}.%
\indexg{linear map in composition with standard function!astral conjugate identities and inequalities|)}%
\end{proof}

\chapter{Astral closedness of extensions}
\label{sec:ast-close-extensions}

Beginning with a function $f:\Rn\rightarrow\Rext$,
we saw in \Cref{sec:lsc:ext} how $f$'s lower semicontinuous
extension $\fext$ provides one natural way of extending $f$ to all of
astral space.
The biconjugate $\fdub$, considered in
\Cref{sec:conj-exts}, can be viewed as providing another way of
extending $f$ to astral space.
In this chapter, we derive necessary and sufficient conditions
characterizing when these two ways of extending $f$ to astral space
yield the same function, so that $\fext=\fdub$.
Since $\smash{\fextdub}=\fdub$
(by \Cref{pr:fextstar-is-fstar}),
this is the same as characterizing when
$\fext$ is astral closed.
Along the way to proving these characterizations,
we also develop a number of techniques that will be
more broadly useful in studying astral functions, especially
extensions.

The methods in this chapter build particularly on properties of
reductions.
As we saw in
\Cref{sec:shadow},
the form and properties of a reduction $\fshadv$
can depend on whether $\vv\in\resc f$. For that and other reasons that will
become apparent later,
the analysis of $\resc f$ is at the core of our development.

However, standard characterizations of $\resc f$ focus on closed, proper and convex functions.
One typical result in this area is \Cref{pr:rescpol-is-con-dom-fstar}, which states that
the recession cone of a closed, proper convex function
$f:\Rn\rightarrow\Rext$
can be expressed as
\begin{equation}
\label{eq:rescpol-is-con-dom-fstar}
    \resc{f}
    = \polar{\bigParens{\cone(\dom{\fstar})}}.
\end{equation}
Results of this type are
insufficient for our needs because functions that arise in the analysis of astral space, such as astronic reductions, can be improper
even when they are derived from functions on $\Rn$ that are closed, convex and proper.
Therefore, before we can characterize when $\fext=\fdub$, we need to generalize results like the one in \eqref{eq:rescpol-is-con-dom-fstar} to improper convex functions.

There are two components in our generalization approach. First, we replace $\cone(\dom{\fstar})$
with a slightly larger set, called the {barrier cone of~$f$}, which we introduce and analyze in \Cref{sec:slopes}.
Second, we give a technique for
replacing any improper function $f$ by a proper function $f'$
using a monotone transformation (namely, composition with the
exponential function) which preserves convexity and many other of $f$'s
properties.
We can then apply existing results to the modified function $f'$ to obtain results for the original function $f$. This technique is developed in \Cref{sec:exp-comp}.

\section{Barrier cone of a function}
\label{sec:slopes}

\indexg{barrier cone (of set)|(}%
The \emph{barrier cone} of a set $S\subseteq\Rn$,
denoted $\barr S$,
is the set
\begin{equation}    \label{eqn:barrier-cone-defn}
\indexm{bars}{$\barr S$}{barrier cone of set}%
   \barr S
   =
   \Braces{ \uu\in\Rn :\:
           \sup_{\xx\in S} \xx\cdot\uu < +\infty
          },
\end{equation}
which can be viewed as the set of vectors in whose direction
the set $S$ is eventually bounded.
The barrier cone can also be expressed in terms of
the support function for $S$
(as defined in
Eq.~\ref{eq:e:2}),
specifically, as
the effective domain of $\indstars$:
\begin{equation}  \label{eqn:bar-cone-defn}
   \barr S
   =
   \bigBraces{ \uu\in\Rn :\:
           \indstar{S}(\uu) < +\infty
         }
   =
   \dom{\indstar{S}}.
\end{equation}
Note that the barrier cone of the empty set is $\Rn$ since then the
supremum appearing in
\eqref{eqn:barrier-cone-defn} is over an empty set and so equals
$-\infty$.
(This also follows from
\eqref{eqn:bar-cone-defn} since
$\indstar{\emptyset}\equiv-\infty$.)

Here are some properties of barrier cones:

\begin{proposition}  \label{pr:bar-cone-props}
  Let $S\subseteq\Rn$.
  Then:
  \begin{letter-compact}
  \item      \label{pr:bar-cone-props:a}
    $\barr S$ is a convex cone.
  \item      \label{pr:bar-cone-props:b}
    $\barr{S}=\barr(\cl S)$.
  \item      \label{pr:bar-cone-props:c}
    $\barr{S}=\Rn$ if and only if $S$ is bounded.
  \end{letter-compact}
\end{proposition}

\begin{proof}
Each of the claims holds when $S=\emptyset$ since
$\barr \emptyset=\Rn$, which is a convex cone,
since $\cl\emptyset=\emptyset$, and since
$\emptyset$ is bounded.
We therefore assume henceforth that $S$ is nonempty.

\begin{proof-parts}
\pfpart{Part~(\ref{pr:bar-cone-props:a}):}
From \eqref{eqn:barrier-cone-defn}, it is clear that $\barr{S}$
includes $\zero$ and is closed under multiplication by a positive
scalar.
Thus, $\barr{S}$ is a cone.
Moreover, by \eqref{eqn:bar-cone-defn}, $\barr{S}$ is the effective domain of the convex function $\indstar{S}$, so it must be convex
(by Propositions~\ref{roc:thm4.6}
and~\ref{pr:conj-props}\ref{pr:conj-props:d}).

\pfpart{Part~(\ref{pr:bar-cone-props:b}):}
Let $\uu\in\Rn$.
Then by \Cref{pr:sup-f-Abar},
$
  \sup_{\xx\in S} \xx\cdot\uu
  =
  \sup_{\xx\in \cl{S}} \xx\cdot\uu
$,
since $\xx\mapsto\xx\cdot\uu$ is continuous.
Therefore, $\uu\in\barr{S}$ if and only if
$\uu\in\barr{(\cl{S})}$, proving the claim.

\pfpart{Part~(\ref{pr:bar-cone-props:c}):}
If $S$ is bounded, then there exists $M\in\Rpos$ such that
$\norm{\xx}\le M$ for all $\xx\in S$. Thus,
by the Cauchy-Schwarz inequality,
for any $\uu\in\Rn$,
$\sup_{\xx\in S} \xx\cdot\uu \le M\norm{\uu}<+\infty$.
Hence, $\uu\in\barr{S}$.

If $S$ is not bounded, then there exists a sequence $\seq{\xx_t}$ in $S$ such that $\norm{\xx_t}\to+\infty$.
By sequential compactness,
this sequence must have a subsequence that converges in $\eRn$;
by discarding all other elements, we can
assume that $\seq{\xx_t}$ has a limit, say $\xbar\in\extspace$
(which cannot be finite since $\seq{\xx_t}$ is unbounded).
Let $\vv\in\Rn$ be the dominant direction of $\xbar$. Then $\xx_t\inprod\vv\to\xbar\inprod\vv=+\infty$
(by Theorems~\ref{thm:i:1}\ref{thm:i:1c} and~\ref{thm:dom-dir}\ref{thm:dom-dir:a}\ref{thm:dom-dir:b}),
so $\vv\not\in\barr{S}$.
Hence, $\barr{S}\ne\Rn$.%
\indexg{barrier cone (of set)|)}%
\qedhere
\end{proof-parts}
\end{proof}

\indexg{barrier cone (of set)!epigraph@of epigraph|(}%
\indexg{epigraph!barrier cone of|(}%
We will especially be interested in the
barrier cone of the epigraph of a function $f:\Rn\rightarrow\Rext$,
which is the set
\[
   \barr(\epi f)
   =
   \biggBraces{ \rpair{\uu}{v}\in\Rn\times\R :
           \sup_{\rpair{\xx}{y}\in\epi f} (\xx\cdot\uu + yv) < +\infty
         }
   =
   \dom{\indstar{\epi f}}.
\]
The form of $\indstar{\epi f}$
was derived in
\Cref{pr:support-epi-f-conjugate}.
Using that proposition, we can
relate the set $\cone(\dom{\fstar})$, appearing in the standard characterization of $\resc{f}$ in~\eqref{eq:rescpol-is-con-dom-fstar}, to the barrier cone $\barr(\epi f)=\dom{\indepifstar}$:
\begin{proposition}  \label{pr:conedomfstar-and-barrepif}
  Let $f:\Rn\rightarrow\Rext$.
  Let $\uu\in\Rn\wo\{\zero\}$.
  Then $\uu\in\cone(\dom{\fstar})$ if and only if
  $\rpair{\uu}{-v}\in\barr(\epi f)$ for some $v\in\Rstrictpos$.
\end{proposition}

\begin{proof}
By \Cref{pr:support-epi-f-conjugate},
for any $v\in\Rstrictpos$,
\[
    \indepifstar(\rpairf{\uu}{-v})
    =
    v\fstar(\uu/v).
\]
Therefore, $\rpair{\uu}{-v}\in\barr(\epi f)=\dom{\indepifstar}$ if and only if $(\uu/v)\in\dom\fstar$.

This proves that if $\rpair{\uu}{-v}\in\barr(\epi f)$ for some
$v\in\Rstrictpos$
then $\uu\in\cone(\dom\fstar)$.
For the reverse implication, suppose $\uu\in\cone(\dom{\fstar})$. Since $\uu\ne\zero$ and $\dom\fstar$ is convex,
there exists
$v\in\Rstrictpos$
such that $(\uu/v)\in\dom\fstar$ (by \Cref{pr:scc-cone-elts}\ref{pr:scc-cone-elts:c:new}); hence, $\rpair{\uu}{-v}\in\barr(\epi f)$.
\end{proof}

Thus, except possibly for the origin,
$\cone(\dom{\fstar})$ consists of
vectors $\uu$ that correspond to the elements of $\barr(\epi f)$
taking the
form $\rpair{\uu}{v}$ with $v<0$.
Nevertheless, to analyze improper functions,
we need to consider the \emph{entire} set
$\barr(\epi f)$, including
the elements $\rpair{\uu}{v}$ with $v=0$.

\indexg{barrier cone (of function)|(}%
Accordingly, we define the barrier cone of a function:

\begin{definition}
\indexg{barrier cone (of function)!defined|(}%
Let $f:\Rn\rightarrow\Rext$.
Then the \emph{barrier cone of $f$},
denoted $\slopes{f}$, is the set
\begin{equation}  \label{eq:bar-f-defn}
\indexm{barf}{$\slopes{f}$}{barrier cone of function}
  \slopes{f}
  =
  \BigBraces{ \uu\in\Rn :\:
    \exists v\in\R,\,
    \rpair{\uu}{v}\in\barr(\epi f)
  }.%
\indexg{barrier cone (of function)!defined|)}
\end{equation}
\end{definition}

In other words, $\slopes{f}$
is the projection of $\barr(\epi f)$
onto $\Rn$;
that is,
\begin{equation}  \label{eq:bar-f-defn-alt}
  \slopes{f}
  =
  \lmPPx\bigParens{\barr(\epi f)},
\end{equation}
where $\lmPPx:\R^{n+1}\rightarrow\Rn$ is the projection map
$\lmPPx(\uu,v)=\uu$ for
$\rpair{\uu}{v}\in\Rn\times\R$.
Note in particular that if $f\equiv+\infty$ then
$\barr{f}=\Rn$.

\indexg{vertical barrier cone|(}%
The pairs $\rpair{\uu}{v}$ in
$\barr(\epi f)$ with $v=0$ are of particular interest, as they
do not generally correspond to vectors in $\cone(\dom\fstar)$.
We call the set of all such vectors $\uu\in\Rn$
with $\rpair{\uu}{0}\in\barr(\epi f)$
the \emph{vertical barrier cone of $f$}, denoted
$\vertsl{f}$:
\begin{equation}  \label{eq:vert-bar-cone-def}
\indexm{vertf}{$\vertsl{f}$}{vertical barrier cone}
  \vertsl{f}
  =
  \bigBraces{ \uu\in\Rn :\: \rpair{\uu}{0}\in\barr(\epi f) }.%
\indexg{barrier cone (of set)!epigraph@of epigraph|)}%
\indexg{epigraph!barrier cone of|)}%
\end{equation}
This set is so named because
each vector $\uu$ in the set defines a vertical halfspace
in $\R^{n+1}$ that includes all of $\epi f$.
This set can also be expressed as the barrier cone of $\dom f$:

\begin{proposition}  \label{pr:vert-bar-is-bar-dom}
  Let $f:\Rn\rightarrow\Rext$.
  Then $\vertsl{f}$ is a convex cone,
  and
  \[  \vertsl{f}=\dom{\indstar{\dom{f}}}=\barr(\dom{f}). \]
\end{proposition}

\begin{proof}
We have
\[
  \vertsl{f}
  =
  \bigBraces{ \uu\in\Rn :\: \indstar{\epi f}(\uu,0) < +\infty }
  =
  \dom{\indstar{\dom{f}}}
  =
  \barr(\dom{f}).
\]
The first equality is by
\eqref{eqn:bar-cone-defn} combined with
\eqref{eq:vert-bar-cone-def}.
The second is because
$\indepifstar(\rpairf{\uu}{0})=\inddomfstar(\uu)$
by \Cref{pr:support-epi-f-conjugate}.
The third is also by
\eqref{eqn:bar-cone-defn}.
In particular, this shows that $\vertsl{f}$ is a convex cone,
by \Cref{pr:bar-cone-props}(\ref{pr:bar-cone-props:a}).
\end{proof}

When $\vertsl{f}$ is combined with
$\cone(\dom{\fstar})$, we obtain exactly the barrier cone of $f$,
as shown in the next proposition.
Thus, $\cone(\dom{\fstar})$ is always included in the barrier cone of
$f$, but this latter set might also include additional elements from
the vertical barrier cone.

\begin{theorem}  \label{thm:slopes-equiv}  %
  Let $f:\Rn\rightarrow\Rext$.
  Then $\slopes{f}$ is a convex cone, and
  \[
  \slopes{f} =
  \cone(\dom{\fstar}) \cup (\vertsl{f}).
  \]
\end{theorem}

\begin{proof}
If $f\equiv+\infty$, then $\fstar\equiv-\infty$ so
$\slopes{f}=\cone(\dom{\fstar})=\Rn$, implying the claim.
We therefore assume henceforth that $f\not\equiv+\infty$.

As a projection of the convex cone $\barr(\epi f)$, the set $\slopes{f}$ is also a convex cone (by \Cref{prop:cone-linear}).
By
\Cref{pr:support-epi-f-conjugate},
all the points $\rpair{\uu}{v}\in\barr(\epi f)=\dom{\indstar{\epi f}}$
satisfy $v\le 0$. Let $B_0$ denote the set of points $\rpair{\uu}{v}\in\barr(\epi f)$
with $v=0$ and $B_{-}$ the set of points with $v<0$; thus, $\barr(\epi f)=B_0\cup B_{-}$.
We then have
(with $\lmPPx$ as used in Eq.~\ref{eq:bar-f-defn-alt}):
\begin{align*}
  \slopes{f}
  =
  \lmPPx\bigParens{\barr(\epi f)}
  &=
  \lmPPx(B_0) \cup \lmPPx(B_{-})
  \\
  &=
  (\vertsl{f}) \cup \lmPPx(B_{-})
\\
  &= (\vertsl{f}) \cup \bigBracks{\lmPPx(B_{-})\setminus\set{\zero}}
\\
  &= (\vertsl{f}) \cup \bigBracks{\cone(\dom{\fstar})\setminus\set{\zero}}
\\
  &= (\vertsl{f}) \cup\cone(\dom{\fstar}).
\end{align*}
The first and third equalities are by the definitions of
$\slopes{f}$ and $\vertsl{f}$
(Eqs.~\ref{eq:bar-f-defn-alt} and~\ref{eq:vert-bar-cone-def}).
The fourth and last equalities are because $\vertsl{f}$
includes $\zero$, being a cone
(\Cref{pr:vert-bar-is-bar-dom}).
And the fifth equality is by
\Cref{pr:conedomfstar-and-barrepif}.%
\indexg{vertical barrier cone|)}%
\end{proof}

The barrier cone of a function is the same as that of its
lower semicontinuous hull:

\begin{proposition}  \label{pr:slopes-same-lsc}
  Let $f:\Rn\rightarrow\Rext$.
  Then $\slopes(\lsc f)=\slopes{f}$.
\end{proposition}

\begin{proof}
We have
$\barr(\epi f)=\barr\regParens{\cl(\epi f)}=\barr\regParens{\epi(\lsc f)}$,
where the first equality is by
\Cref{pr:bar-cone-props}(\ref{pr:bar-cone-props:b}),
and the second is from
\Cref{prop:lsc:characterize}(\ref{prop:lsc:characterize:c}).
From the definition of the barrier cone of a function
(Eq.~\ref{eq:bar-f-defn}),
this proves the claim.
\end{proof}

In general, $\slopes{f}$ need not be the same as
$\cone(\dom{\fstar})$, even when $f$ is convex, closed, and proper,
as shown by the next example:

\begin{example}[Restricted linear function]
\label{ex:negx1-else-inf}
\indexg{Restricted linear function|(}%
\indexg{Restricted linear function!barrier cone of|(}%
Let $f:\R^2\rightarrow\Rext$ be defined, for $\xx\in\R^2$, by
\begin{equation}  \label{eqn:ex:1}
  f(\xx) = 
   f(x_1,x_2) =
   \begin{cases}
     -x_1
     & \text{if $x_2\geq 0$,}
   \\
     +\infty
     & \text{otherwise.}
   \end{cases}
\end{equation}
This function is convex, closed, and proper.
\indexg{Restricted linear function|)}%
Its conjugate, for $\uu\in\R^2$, can be computed to be
\begin{equation}  \label{eqn:ex:1:conj}
  \fstar(\uu) =
   \fstar(u_1,u_2) =
   \begin{cases}
             0   & \mbox{if $u_1=-1$ and $u_2\leq 0$,} \\
          +\infty & \mbox{otherwise.}
   \end{cases}
\end{equation}
Thus, using \Cref{pr:scc-cone-elts}(\ref{pr:scc-cone-elts:c:new}),
\begin{equation}  \label{eqn:ex:1:cone-dom-fstar}
  \cone(\dom{\fstar})=
    \{\zero\} \,\cup\,
    \bigBraces{\uu\in\R^2:\: u_1<0 \text{ and } u_2\leq 0}.
\end{equation}
From \Cref{pr:vert-bar-is-bar-dom}, it can further be calculated that
\[
  \vertsl{f} = \{\uu\in\R^2:\: u_1=0 \text{ and } u_2\leq 0\},
\]
so by \Cref{thm:slopes-equiv},
\[\slopes{f} = \{\uu\in\R^2:\: u_1\leq 0 \text{ and } u_2\le 0\}.\]
Hence, in this case $\slopes{f}\neq\cone(\dom\fstar)$.%
\indexg{barrier cone (of function)|)}%
\indexg{Restricted linear function!barrier cone of|)}%
\end{example}

Later, in \Cref{sec:dual-char-ent-clos},
we will study when
$\slopes{f}=\cone(\dom{\fstar})$, a property that we will see
characterizes
when a function's biconjugate $\fdub$ is the same as its
extension~$\fext$.

\section{Exponential composition}
\label{sec:exp-comp}

We next introduce a simple trick that will allow us to translate many results for closed, proper convex functions to analogous results for improper lower semicontinuous functions. Specifically, we transform a possibly improper convex function $f\not\equiv+\infty$ into a proper function~$f'$ via a composition with the exponential function. We show that $f'$ retains many of $f$'s properties, so we can apply existing results to the modified function $f'$ to obtain results for the original function $f$.

\indexg{strictly increasing function, composition with|(}%
\indexg{exponential function!composition with|(}%
Let $g:\R\to\R$ be any strictly increasing convex function that is bounded below, and let $\eg:\eR\to\eR$
be its extension. Then $\eg$ must be continuous and strictly increasing over its entire domain (by \Cref{pr:conv-inc:prop}\ref{pr:conv-inc:infsup}\ref{pr:conv-inc:strictly}).
\indexg{exponential function!extension of|(}%
For example, when $g(x)=e^x=\exp(x)$, 
the extension $\eg=\expex$ is given, for $\ex\in\Rext$,
by
\begin{equation}  \label{eq:expex-defn}
\indexm{exp}{$\expex$}{exponential's extension}
   \expex(\ex) =
   \begin{cases}
     0
       & \text{if $\ex=-\infty$,} \\
     \exp(\ex)
       & \text{if $\ex\in\R$,} \\
     +\infty
       & \text{if $\ex=+\infty$.}%
\indexg{exponential function!extension of|)}%
   \end{cases}
\end{equation}
Any convex function $f:\Rn\rightarrow\Rext$ can be composed with $\eg$ to derive
a new convex function $f'=\eg\circ f$, which is bounded below.
This and other properties of $f'$ are given in
the next proposition. Canonically, we will invoke this proposition with $\eg=\expex$.
\begin{proposition}  \label{pr:j:2}
  Let $f:\Rn\rightarrow\Rext$ be convex, let
  $g:\R\to\R$ be a strictly increasing convex function
  with $\inf g>-\infty$,
  and let $f'=\eg\circ f$, that is, $f'(\xx)=\eg(f(\xx))$ for $\xx\in\Rn$.
  Then the following hold:
  \begin{letter-compact}
  \item  \label{pr:j:2a}
    $f'$ is convex and bounded below, with $\inf f'\ge \inf g$.
  \item  \label{pr:j:2c}
    $\fpext=\overline{\eg\circ f}=\eg\circ\ef$.
  \item  \label{pr:j:2lsc}
    $\lsc{f'}=\lsc(\eg\circ f)=\eg\circ(\lsc f)$.
    Therefore,
    if $f$ is lower semicontinuous, then $f'$ is as well.
  \item  \label{pr:j:2b}
    For all $\xx,\yy\in\Rn$,
    $f'(\xx)\leq f'(\yy)$ if and only if $f(\xx)\leq f(\yy)$.
    So, $\resc{f'}=\resc{f}$.
  \item  \label{pr:j:2d}
    For all $\xbar,\ybar\in\extspace$,
    $\fpext(\xbar)\leq \fpext(\ybar)$ if and only if
    $\fext(\xbar)\leq \fext(\ybar)$.
  \item \label{pr:j:2lim}
    Let $\xbar\in\eRn$ and let $\seq{\xx_t}$ be a sequence in $\Rn$ such that $\xx_t\to\xbar$. Then $f(\xx_t)\to\ef(\xbar)$ if and only if $f'(\xx_t)\to\fpext(\xbar)$.
  \end{letter-compact}
  In particular, the above statements hold when $g=\exp$ and $\eg=\expex$.
\end{proposition}

\begin{proof}
As a real-valued convex function, $g$ must be continuous (\Cref{pr:stand-cvx-cont}), so
\Cref{pr:conv-inc:prop}(\ref{pr:conv-inc:infsup}) implies that $\eg$ is continuous on $\eR$ with the
values $\eg(-\infty)=\inf g$, $\eg(x)=g(x)$ for $x\in\R$,
and $\eg(+\infty)=\sup g$. Moreover, by \Cref{pr:conv-inc:prop}(\ref{pr:conv-inc:strictly}),
$\eg$ is strictly increasing and $\eg(+\infty)=+\infty$.

\begin{proof-parts}
\pfpart{Part~(\ref{pr:j:2a}):}
From the definition of $f'$, we have
$\inf f'\geq\min\eg=\inf g>-\infty$
(using \Cref{pr:fext-min-exists}).
Convexity of $f'$ follows by \Cref{prop:nondec:convex}.

\pfpart{Part~(\ref{pr:j:2c}):}
Let $\xbar\in\eR$; we aim to show
$\eg\regParens{\fext(\xbar)}=\fpext(\xbar)$.
By \Cref{pr:d1},
there exists a sequence $\seq{\xx_t}$ in $\Rn$ converging to $\xbar$
with $f(\xx_t)\rightarrow \fext(\xbar)$.
This implies
\[
  \eg\bigParens{\fext(\xbar)}
  =
  \eg\bigParens{\lim f(\xx_t)}
  =
  \lim \eg\bigParens{f(\xx_t)}
  =
  \lim f'(\xx_t)
  \geq
  \fpext(\xbar),
\]
where the second equality is by continuity of $\eg$,
and the inequality is by definition of extensions
(Eq.~\ref{eq:e:7}).

By a similar argument,
there exists a sequence $\seq{\xx'_t}$ in $\Rn$ converging to $\xbar$
with $f'(\xx'_t)\rightarrow \fpext(\xbar)$.
Further, since $\Rext$ is sequentially compact, there exists a
subsequence of the sequence $\seq{f(\xx'_t)}$ that converges; by
discarding all other elements, we can assume that the entire sequence
converges in $\Rext$.
Then
\[
  \fpext(\xbar)
  =
  \lim f'(\xx'_t)
  =
  \lim \eg\bigParens{f(\xx'_t)}
  =
  \eg\bigParens{\lim f(\xx'_t)}
  \geq
  \eg\bigParens{\fext(\xbar)}.
\]
The third equality is because $\eg$ is continuous.
The inequality is because $\lim f(\xx'_t)\geq \fext(\xbar)$ by
definition of $\fext$, and because $\eg$ is strictly increasing.

\pfpart{Part~(\ref{pr:j:2lsc}):}
For all $\xx\in\Rn$,
\[
  (\lsc f')(\xx)
  =
  \fpext(\xx)
  =
  \eg\bigParens{\fext(\xx)}
  =
  \eg\bigParens{(\lsc f)(\xx)},
\]
where the first and third equalities are by
\Cref{pr:h:1}(\ref{pr:h:1a}), and the second is
by part~(\ref{pr:j:2c}).
If $f$ is lower semicontinuous then
the right-hand side is equal to $\eg(f(\xx))=f'(\xx)$,
proving $f'$ is lower semicontinuous as well.

\pfpart{Part~(\ref{pr:j:2b}):}
The first claim is immediate from the fact that $\eg$ is strictly increasing.
The second claim follows from
the first claim and
the definition of the recession cone
(Eq.~\ref{eqn:resc-cone-def}).

\pfpart{Part~(\ref{pr:j:2d}):} This follows from
part~(\ref{pr:j:2c}) since $\eg$ is strictly increasing.

\pfpart{Part~(\ref{pr:j:2lim}):} If $f(\xx_t)\to\ef(\xbar)$,
then
\[
  \fpext(\xbar)
  =
  \eg\bigParens{\ef(\xbar)}
  =
  \eg\bigParens{\lim f(\xx_t)}
  =
  \lim \eg\bigParens{f(\xx_t)}
  =
  \lim f'(\xx_t).
\]
The first equality is by part~(\ref{pr:j:2c}).
The third is by continuity of $\eg$.
Thus, $f'(\xx_t)\rightarrow  \fpext(\xbar)$.

For the converse, suppose $f'(\xx_t)\to\fpext(\xbar)$.
From the definition of $\ef$
(Eq.~\ref{eq:e:7}), $\liminf f(\xx_t)\ge\ef(\xbar)$.
For a contradiction, suppose $\limsup f(\xx_t)>\ef(\xbar)$.
Then there exists $\beta\in\R$ such that $\beta>\ef(\xbar)$ and
$f(\xx_t)\geq\beta$ for infinitely many values of $t$.
For all such $t$, by the strict monotonicity of $\eg$,
we have
$ f'(\xx_t) = \eg\regParens{f(\xx_t)} \geq \eg(\beta)$,
and furthermore that
$ \eg(\beta) > \eg\regParens{\ef(\xbar)} = \fpext(\xbar)$,
with the last equality from part~(\ref{pr:j:2c}).
This contradicts that
\indexg{strictly increasing function, composition with|)}%
$f'(\xx_t)\to\fpext(\xbar)$.
\qedhere
\end{proof-parts}
\end{proof}

\indexg{barrier cone (of function)!exponential composition and|(}%
Now focusing on the special case when $g=\exp$, we analyze the conjugate
of $f'=\expex\circ f$ and show
that the effective domain of $\fpstar$ turns out to be exactly equal to $f$'s barrier cone, $\slopes{f}$:

\begin{theorem}
\label{thm:dom-fpstar}
  Let $f:\Rn\rightarrow\Rext$ be convex,
  and let $f'=\expex\circ f$.
  Then
  \begin{equation}  \label{eq:thm:dom-fpstar:1}
     \slopes{f}
     =
     \dom{\fpstar}
     =
     \cone(\dom{\fpstar}).
  \end{equation}
\end{theorem}

\begin{proof}
Since $\slopes{f}$ is a convex cone (by \Cref{thm:slopes-equiv}), it
suffices to prove $\slopes{f}=\dom{\fpstar}$,
since this will then imply the second equality of
\eqref{eq:thm:dom-fpstar:1}.
If $f\equiv+\infty$, then $f'\equiv+\infty$ and $\fpstar\equiv-\infty$,
so $\barr{f}=\Rn=\dom{\fpstar}$, proving the claim in this case.
We therefore assume henceforth that $f\not\equiv+\infty$.

Let $\uu\in\Rn$.
Setting $g=\exp$ and $G=\expex$,
\Cref{thm:conj-compose-our-version} then yields
\begin{equation}
\label{eq:from:conj-compose}
     \fpstar(\uu)=(G\circ f)^*(\uu)
     =
     \min_{v\in\R} \bigBracks{g^*(v)+\indepifstar(\rpairf{\uu}{-v})},
\end{equation}
noting that,
since $f\not\equiv+\infty$, there exists $\xing\in\Rn$ such
that
$f(\xing)<+\infty=\sup(\dom g)$.
By a standard calculation, for $v\in\R$,
\begin{equation}
\label{eq:expstar}
  g^*(v)=
  \begin{cases}
     v\ln v - v &\text{if $v\ge 0$,}
   \\
     +\infty &\text{if $v<0$,}
   \end{cases}
\end{equation}
with the convention $0\ln 0=0$.
We then have
\begin{align*}
 \uu\in\dom{\fpstar}
 &\;\Leftrightarrow\;
 \exists v \in \R
   \text{ such that }
   v\in\dom{\gstar}
   \text{ and }
   \rpair{\uu}{-v}\in\dom{\indfstar{\epi f}}
\\
&\;\Leftrightarrow\;
   \exists v \in \R
   \text{ such that }
   \rpair{\uu}{-v}\in\dom{\indfstar{\epi f}}
\\
 &\;\Leftrightarrow\;
   \exists v\in\R
   \text{ such that }
   \rpair{\uu}{-v}\in\barr(\epi f)
\\
 &\;\Leftrightarrow\;
   \uu\in\barr f.
\end{align*}
The first equivalence is by \eqref{eq:from:conj-compose}.
The second equivalence is
because, for $v\in\R$, if 
$\indfstar{\epi f}(\uu,-v)<+\infty$, then
$v\geq 0$ by
\Cref{pr:support-epi-f-conjugate},
implying
$\gstar(v)<+\infty$ by \eqref{eq:expstar}.
The third equivalence is because
$\dom{\indfstar{\epi f}}=\barr(\epi f)$ by
\eqref{eqn:bar-cone-defn}.
The final equivalence is from the definition of $\barr f$
(Eq.~\ref{eq:bar-f-defn}).
Hence, $\dom{\fpstar}=\barr f$.%
\indexg{exponential function!composition with|)}%
\indexg{barrier cone (of function)!exponential composition and|)}%
\end{proof}

\section{Reductions and conjugacy}
\label{sec:conjugacy:reductions}

\indexg{reductions, astronic|(}%
In \Cref{sec:shadow} we studied properties of astronic reductions $g=\fshadv$
under the assumption that $\vv$ is in $f$'s recession cone, $\resc{f}$.
We next use exponential composition
to show that if $\vv$ is \emph{not}
in $\resc{f}$, then $g$ is identically equal to $+\infty$:

\begin{theorem}  \label{thm:i:4}
Let $f:\Rn\rightarrow\Rext$ be convex and lower semicontinuous.
Let $\vv\in\Rn$
and let
$g=\fshadv$
be the reduction of $f$ at $\limray{\vv}$.
Assume $\vv\not\in\resc{f}$.
Then
\[
  \gext(\xbar)
  =\fext(\limray{\vv}\plusl\xbar)
  =+\infty
\]
for all $\xbar\in\extspace$.
If, in addition, $f$ is closed, then
$\fdub(\limray{\vv}\plusl\xbar) =+\infty$
for all $\xbar\in\extspace$.
\end{theorem}

\begin{proof}
The cases $f\equiv +\infty$ or $f\equiv -\infty$
are both impossible since then we would have
$\resc{f}=\Rn$, but $\vv\not\in\resc{f}$.
Consequently, being lower semicontinuous,
$f$ is closed if and only if it is proper.

Let us assume for the moment that $f$ is closed, and so also proper.
Then
by \Cref{pr:rescpol-is-con-dom-fstar},
\[
  \resc f
  =
    \bigSet{\vv\in\Rn:\: \forall \uu\in\dom{\fstar},\, \uu\cdot\vv\leq 0}.
\]
Since $\vv\not\in\resc{f}$,
there must exist a point
$\uu\in\dom{\fstar}$ with $\vv\cdot\uu>0$.
For all $\xbar\in\extspace$, we then have
\begin{align*}
  \fext(\limray{\vv}\plusl\xbar)
  &\geq
  \fdub(\limray{\vv}\plusl\xbar)
  \\
  &\geq
  -\fstar(\uu)
  \plusl (\limray{\vv}\plusl \xbar)\cdot\uu
  \\
  &=
  -\fstar(\uu)
  \plusl \limray{\vv}\cdot\uu
  \plusl \xbar\cdot\uu
  = +\infty.
\end{align*}
The inequalities are by
\Cref{thm:fext-dub-sum}(\ref{thm:fext-dub-sum:a},\ref{thm:fext-dub-sum:b}).
The last equality is because $-\fstar(\uu)>-\infty$
(since $\uu\in\dom{\fstar}$),
and $\limray{\vv}\cdot\uu=+\infty$
(since $\vv\cdot\uu>0$).
Thus,
$\fext(\limray{\vv}\plusl\xbar)=\fdub(\limray{\vv}\plusl\xbar)=+\infty$
for all $\xbar\in\extspace$.

Returning to the general case,
suppose $f\not\equiv +\infty$ but $f$ is not
necessarily closed and proper.
Let $f'=\expex\circ f$.
Then $f'$ is convex, lower-bounded, and lower semicontinuous
by \Cref{pr:j:2}(\ref{pr:j:2a},\ref{pr:j:2lsc}),
so $f'$ is also proper and closed.
Also, $\resc{f'}=\resc{f}$
(by \Cref{pr:j:2}\ref{pr:j:2b}),
so $\vv\not\in\resc{f'}$.
Therefore,
for all $\xx\in\Rn$,
\[
  \expex\bigParens{g(\xx)}
  =
  \expex\bigParens{\fext(\limray{\vv}\plusl\xx)}
  =
  \fpext(\limray{\vv}\plusl\xx)=+\infty,
\]
with the second equality from
\Cref{pr:j:2}(\ref{pr:j:2c}),
and the third from the argument above, applied to $f'$.
This implies
$g(\xx)=+\infty$.
Thus, $g\equiv+\infty$ so
$\gext\equiv+\infty$.
\end{proof}

As a summary, the next corollary combines \Cref{thm:i:4} with
Theorems~\ref{thm:d4} and~\ref{thm:a10-nunu}, making no assumption
about $\vv$:

\begin{corollary}  \label{cor:i:1}
  Let $f:\Rn\rightarrow\Rext$ be convex, let $\vv\in\Rn$,
  and let
  $g=\fshadv$
  be the reduction of $f$ at $\limray{\vv}$.
  Then:
  \begin{letter-compact}
  \item \label{cor:i:1a}
    $g$ is convex and lower semicontinuous.
  \item \label{cor:i:1b}
    $\gext(\xbarperp)=\gext(\xbar)=\fext(\limray{\vv}\plusl\xbar)$
    for $\xbar\in\extspace$,
    with $\xbarperp$ denoting its projection orthogonal to~$\vv$.
  \item \label{cor:i:1c}
    If $\vv\in\resc{f}$, then $\gext\leq\fext$.
  \end{letter-compact}
\end{corollary}

\begin{proof}
Part~(\ref{cor:i:1c}) is exactly what is stated in
\Cref{thm:d4}(\ref{thm:d4:b}).

Assume for the remainder of the proof that $f$ is lower
semicontinuous.
This is without loss of generality since if it is not, we can replace
$f$ with $\lsc f$, which does not affect either $\fext$ or $g$
(by \Cref{pr:h:1}\ref{pr:h:1aa}).

If $\vv\in\resc{f}$, then parts~(\ref{cor:i:1a}) and~(\ref{cor:i:1b})
follow from \Cref{thm:d4}(\ref{thm:d4:a},\ref{thm:d4:c})
and \Cref{thm:a10-nunu}.

If $\vv\not\in\resc{f}$, then
$\gext(\xbar)=\fext(\limray{\vv}\plusl\xbar)=+\infty$
for all $\xbar\in\extspace$,
by \Cref{thm:i:4}.
This implies part~(\ref{cor:i:1b}).
This also implies that $g$ is the constant function $+\infty$, which is
convex and lower semicontinuous, proving part~(\ref{cor:i:1a}).%
\indexg{reductions, astronic|)}%
\end{proof}

We next explore how conjugacy and astronic reductions relate to one
another.
As we showed in \Cref{thm:a10-nunu}, for a function $f:\Rn\rightarrow\Rext$
and for $\vv\in\resc{f}$, the reduction $\fshadv$ is the lower
semicontinuous hull of the shadow $\fv$.
\indexg{shadow (of function)!conjugate of|(}%
\indexg{conjugate (standard)!shadow@of shadow|(}%
As a first step toward computing the conjugate of reductions,
we therefore next calculate, for any matrix $\VV\in\Rnk$,
the conjugate of the shadow $\fV$ of $f$ along $\VV$:

\begin{theorem}  \label{thm:fV-shad-conj}
  Let $f:\Rn\rightarrow\Rext$ be convex,
  $f\not\equiv+\infty$,
  and let $\VV\in\Rnk$.
  Let $L=\{\uu\in\Rn :\: \uu\perp\VV\}$.
  Then $\fVstar = \fstar \plusu \indf{L}$.
  That is, for $\uu\in\Rn$,
  \begin{equation*}  %
    \fVstar(\uu)
    =
    \begin{cases}
      \fstar(\uu)
      & \text{if $\uu\perp\VV$,}
      \\
      +\infty
      & \text{otherwise.}
    \end{cases}
  \end{equation*}
\end{theorem}

\begin{proof}
Let $\PP$ be the projection matrix onto $L$.
Then by \Cref{pr:gtil-is-PPfP-gen},
$\fV=(\PPf)\PP$
(noting that $L=(\colspace{\VV})^\perp$ as follows from
\Cref{pr:std-perp-props}\ref{pr:std-perp-props:d},
applied with $S=\columns{\VV}$).
Thus,
\begin{equation}  \label{eq:thm:fV-shad-conj:1}
  \fVstar=[(\PPf)\PP]^*=\PP[(\PPf)^*]=\PP(f^*\PP).
\end{equation}
The second equality is from
\Cref{roc:thm16.3:fA}, noting that,
by \Cref{pr:Ax-in-ri-dom-Af}, there exists a point
$\xhat\in\Rn$ for which
$\PP\xhat\in\ri(\dom{\PPf})$.
The last equality is from
\Cref{roc:thm16.3:Af}.
(In both steps, we also used that $\PP$ is symmetric by 
\Cref{pr:proj-mat-props}\ref{pr:proj-mat-props:a}.)

Hence, for $\uu\in\Rn$,
\begin{align*}
  \fVstar(\uu)
  =
  [\PP(f^*\PP)](\uu)
  &=
  \inf\,\bigBraces{ f^*(\PP\ww) :\: \ww\in\Rn,\, \PP\ww=\uu }
  \nonumber
  \\
  &=
  \begin{cases}
    f^*(\uu) &\text{if $\uu\in\colspace\PP$,}
    \\
    +\infty  &\text{otherwise,}
  \end{cases}
\end{align*}
where the first equality is by \eqref{eq:thm:fV-shad-conj:1},
and
in the second, we expanded definitions
(Eqs.~\ref{eq:fA-defn} and~\ref{eq:lin-image-fcn-defn}).
This
implies further that
$\fVstar=\fstar\plusu\indf{L}$
since $\colspace\PP=L$, as follows from
\Cref{pr:proj-mat-props}(\ref{pr:proj-mat-props:c},\ref{pr:proj-mat-props:d}),
completing the proof.%
\indexg{shadow (of function)!conjugate of|)}%
\indexg{conjugate (standard)!shadow@of shadow|)}%
\end{proof}

\indexg{conjugate (standard)!astronic reduction@of astronic reduction|(}%
\indexg{reductions, astronic!conjugate of|(}%
Next, we apply \Cref{thm:fV-shad-conj} to
relate the conjugate of a convex function $f$ to that of its
reduction $g=\fshadv$, provided $\vv\in\resc{f}$
(and $f\not\equiv+\infty$).

\begin{theorem}  \label{thm:e1}
Let $f:\Rn\rightarrow\Rext$ be convex, $f\not\equiv+\infty$, and
let $\vv\in\resc{f}$.
Let
$g=\fshadv$
be the reduction of $f$ at $\limray{\vv}$,
and let
$L=\{\vv\}^{\perp}$
be the linear subspace orthogonal to~$\vv$. Then:
\begin{letter-compact}
\item  \label{thm:e1:a}
  $g^*=f^* \plusu\indf{L}$; that is, for $\uu\in\Rn$,
  \[
    \gstar(\uu) =
    \begin{cases}
      \fstar(\uu)
      & \text{if $\uu\perp\vv$,}
      \\
      +\infty
      & \text{otherwise.}
    \end{cases}
  \]
\item  \label{thm:e1:b}
  $\dom{\gstar} = (\dom{\fstar}) \cap L$.
\item  \label{thm:e1:c}
  $\cone(\dom{\gstar}) = \cone(\dom{\fstar}) \cap L$.
\item  \label{thm:e1:d}
  $\slopes{g} = (\slopes{f}) \cap L$.
\end{letter-compact}
\end{theorem}

\begin{proof}
  ~

\begin{proof-parts}
\pfpart{Part~(\ref{thm:e1:a}):}
Let $\gtil=\fv$.
Then $g=\lsc\gtil$ by \Cref{thm:a10-nunu}.
Thus,
\[
  \gstar
  =
  (\lsc\gtil)^{*}
  =
  \gtil^*
  =
  \fstar\plusu\indf{L}.
\]
The second equality is by
\Cref{pr:conj-props}(\ref{pr:conj-props:e}).
The third is by 
\Cref{thm:fV-shad-conj}
(with $\VV=\vv$).

\pfpart{Part~(\ref{thm:e1:b}):}
This follows directly from part~(\ref{thm:e1:a}).

\pfpart{Part~(\ref{thm:e1:c}):}
Since $\dom{\gstar}\subseteq \dom{\fstar}$, we have
$\cone(\dom{\gstar})\subseteq\cone(\dom{\fstar})$, and since $L$ is a convex cone
and $\dom{\gstar}\subseteq L$, we have $\cone(\dom{\gstar})\subseteq L$. Thus,
$
  \cone(\dom{\gstar})\subseteq\cone(\dom{\fstar})\cap L
$.

For the reverse inclusion, let
$\ww\in\cone(\dom{\fstar})\cap L$.
Clearly, $\zero$ is in $\cone(\dom{\gstar})$.
Otherwise, if $\ww\neq\zero$, then since
$\ww\in \cone(\dom{\fstar})$,
by \Cref{pr:scc-cone-elts}(\ref{pr:scc-cone-elts:c:new}),
we can write $\ww=\lambda\uu$ for some $\lambda\in\Rstrictpos$ and $\uu\in\dom{\fstar}$.
Moreover, since $\ww\in L$ and $L$ is a linear subspace, we have $\uu=\ww/\lambda\in L$, implying
$\uu\in (\dom{\fstar})\cap L = \dom{\gstar}$.
Thus,
$\ww=\lambda\uu\in\cone(\dom{\gstar})$,
and so
$\cone(\dom{\fstar})\cap L\subseteq \cone(\dom{\gstar})$.

\pfpart{Part~(\ref{thm:e1:d}):}
Let $f'=\expex\circ f$, which is convex by \Cref{pr:j:2}(\ref{pr:j:2a}), and $f'\not\equiv+\infty$.
In addition,
let $g'=\fpshadv$ be the reduction of $f'$ at
$\limray{\vv}$,
noting that $\vv\in\resc{f'}$ by
\Cref{pr:j:2}(\ref{pr:j:2b}).
Then for all $\xx\in\Rn$,
\[
  g'(\xx)
  =
  \fpext(\limray{\vv}\plusl\xx)
  =
  \expex\bigParens{\fext(\limray{\vv}\plusl\xx)}
  =
  \expex\bigParens{g(\xx)},
\]
with the second equality following from
\Cref{pr:j:2}(\ref{pr:j:2c}).
The reductions $g$ and $g'$ are convex by \Cref{thm:a10-nunu},
and not identically $+\infty$ (by \Cref{pr:d2}\ref{pr:d2:c}).
Therefore,
\[
  \slopes{g}
  =
  \dom{\gpstar}
  =
  (\dom{\fpstar}) \cap L
  =
  (\slopes{f}) \cap L,
\]
where the second equality is
by part~(\ref{thm:e1:b}) applied to $f'$,
and the first and third equalities are
by \Cref{thm:dom-fpstar}.%
\indexg{conjugate (standard)!astronic reduction@of astronic reduction|)}%
\indexg{reductions, astronic!conjugate of|)}%
\qedhere
\end{proof-parts}
\end{proof}

\indexg{biconjugate (astral)!astronic reduction@of astronic reduction|(}%
\indexg{reductions, astronic!astral biconjugate of|(}%
We can also relate the biconjugate $\gdub$ of a
reduction $g=\fshadv$ to that of the function $f$ from which it was
derived.
This fact will be used soon in characterizing when
$\fext=\fdub$.

\begin{theorem}  \label{thm:i:5}
  Let $f:\Rn\rightarrow\Rext$ be convex and closed, let $\vv\in\Rn$,
  and let
  $g=\fshadv$
  be the reduction of $f$ at $\limray{\vv}$.
  Then $\gdub(\xbar)=\fdub(\limray{\vv}\plusl\xbar)$
  for all $\xbar\in\extspace$.
\end{theorem}

\begin{proof}
We proceed in cases.

Suppose first that $f$ is not proper.
In this case, since $f$ is convex, closed and improper, it must be either
identically $+\infty$ or
identically $-\infty$.
If $f\equiv +\infty$, then it can be checked that
$g\equiv +\infty$, $\fstar=\gstar\equiv -\infty$
and $\fdub=\gdub\equiv +\infty$.
Similarly, if $f\equiv -\infty$ then
$\fdub=\gdub\equiv -\infty$.
Either way, the theorem's claim holds.

For the next case, suppose $\vv\not\in\resc{f}$.
Then
\Cref{thm:i:4} implies that
${\fdub(\limray{\vv}\plusl\xbar)} =+\infty$
for all $\xbar\in\extspace$,
and also that $g\equiv +\infty$, so $\gdub\equiv +\infty$.
Thus, the claim holds in this case as well.

We are left only with the case that $f$ is closed, convex and proper, and that
$\vv\in\resc{f}$, which we assume for the remainder of the proof.

Let $\uu\in\Rn$.
We argue next that
\begin{equation} \label{eq:e:9}
  {-\gstar(\uu)}
  =
  -\fstar(\uu) \plusl \limray{\vv}\cdot\uu.
\end{equation}
If $\uu\cdot\vv=0$ then by \Cref{thm:e1}(\ref{thm:e1:a}),
$\gstar(\uu)=\fstar(\uu)$ implying
\eqref{eq:e:9} in this case.
If $\uu\cdot\vv<0$, then $\gstar(\uu)=+\infty$ by
\Cref{thm:e1}(\ref{thm:e1:a}),
$\limray{\vv}\cdot\uu=-\infty$,
and $\fstar(\uu)>-\infty$ since $f$ is proper (so
that $\fstar$ is as well, by
\Cref{pr:conj-props-cvx}\ref{pr:conj-props-cvx:a}).
These imply that both sides of \eqref{eq:e:9} are equal to $-\infty$
in this case.
And if $\uu\cdot\vv>0$, then $\uu\not\in\dom{\fstar}$
by \Cref{pr:rescpol-is-con-dom-fstar},
since $\vv\in\resc{f}$,
so $\uu\not\in\dom{\gstar}$ by
\Cref{thm:e1}(\ref{thm:e1:b}).
Therefore, in this case as well,
both sides of \eqref{eq:e:9} are equal to $-\infty$.

Thus, for all $\xbar\in\extspace$,
\begin{align*}
  \gdub(\xbar)
  &=
  \sup_{\uu\in\Rn} \bigBracks{- \gstar(\uu) \plusl \xbar\cdot\uu }
  \\
  &=
  \sup_{\uu\in\Rn} \bigBracks{- \fstar(\uu) \plusl \limray{\vv}\cdot\uu \plusl \xbar\cdot\uu }
  \\
  &=
  \sup_{\uu\in\Rn} \bigBracks{- \fstar(\uu) \plusl (\limray{\vv} \plusl \xbar)\cdot\uu }
  \\
  &=
  \fdub(\limray{\vv}\plusl \xbar),
\end{align*}
where
the first and fourth equalities are
by \Cref{thm:fext-dub-sum}(\ref{thm:fext-dub-sum:b}),
and the second equality is by
\eqref{eq:e:9}.
This completes the proof.%
\indexg{biconjugate (astral)!astronic reduction@of astronic reduction|)}%
\indexg{reductions, astronic!astral biconjugate of|)}%
\end{proof}

\section{Reduction-closedness}
\label{sec:ent-closed-fcn}

We next develop a general and precise characterization
for when $\fext=\fdub$.

\indexg{reductions, iconic|(}%
As a step in that direction,
we first generalize the reductions introduced in
\Cref{sec:shadow} in a
straightforward way, allowing reductions now to be at any icon, not
just astrons.

\begin{definition}   \label{def:iconic-reduction}
\indexg{reductions, iconic!defined|(}%
Let $f:\Rn\rightarrow\Rext$, and let
$\ebar\in\corezn$ be an icon.
The
\emph{reduction of $f$ at icon $\ebar$}
is the function $\fshadd:\Rn\rightarrow\Rext$
defined, for $\xx\in\Rn$, by
\begin{equation}   \label{eqn:def:iconic-reduction}
\indexm{f e 700}{$\fshadd$}{reduction (iconic)}
  \fshadd(\xx) = \fext(\ebar\plusl \xx).
\end{equation}
Such a reduction is said to be
\indexg{reductions, iconic!defined|)}%
\emph{iconic}.
\end{definition}
Clearly, a reduction $\fshadv$ at astron $\limray{\vv}$
is a special case of such an iconic reduction
in which $\ebar=\limray{\vv}$.
When $\ebar=\zero$, the resulting reduction at $\zero$ is
$\fshad{\zero}=\lsc f$,
the lower semicontinuous hull of $f$,
by
\Cref{pr:h:1}(\ref{pr:h:1a}).

The function $\fshadd$
captures the behavior of $\fext$ on the single galaxy
$\galax{\ebar}=\ebar\plusl\Rn$,
but it can also be viewed as a kind of composition of multiple
astronic reductions, as shown in the next proposition.
It is the closedness of all reductions at
all icons in $\corezn$
that will characterize when $\fext=\fdub$.

\begin{proposition}  \label{pr:icon-red-decomp-astron-red}
  Let $f:\Rn\rightarrow\Rext$ be convex, and let
  $\vv_1,\dotsc,\vv_k\in\Rn$, for some $k\ge 0$.
  Let $g_0 = \lsc f$, and for $i=1,\dotsc,k$,
  let $g_i = \gishadvi$.
  Then for $i=0,\dotsc,k$, each function $g_i$ is convex and lower
  semicontinuous, and furthermore,
  \begin{align}
    g_i(\xx)
    &=
    \fext\bigParens{[\vv_1,\dotsc,\vv_i]\omm\plusl\xx},
    \label{eq:pr:icon-red-decomp-astron-red:1}
    \\
    \intertext{for all $\xx\in\Rn$, and}
    \gext_i(\xbar)
    &=
    \fext\bigParens{[\vv_1,\dotsc,\vv_i]\omm\plusl\xbar}.
    \label{eq:pr:icon-red-decomp-astron-red:2}
  \end{align}
  for all $\xbar\in\extspace$.
  In particular, $g_k = \fshadd$
  where $\ebar=\limrays{\vv_1,\dotsc,\vv_k}$.
\end{proposition}

\begin{proof}
Proof is by induction on $i=0,\dotsc,k$.
In the base case that $i=0$,
we have,
by \Cref{pr:h:1}(\ref{pr:h:1a},\ref{pr:h:1aa}),
that
$g_0(\xx)=(\lsc{f})(\xx)=\fext(\xx)$
for $\xx\in\Rn$,
and that
$\gext_0=\lscfext=\fext$.
Further, $g_0$ is convex and lower semicontinuous by
\Cref{pr:lsc-props}(\ref{pr:lsc-props:a}).

For the inductive step when $i>0$, suppose the claim holds for $i-1$.
Then $g_i$ is convex and lower semicontinuous by
\Cref{cor:i:1}(\ref{cor:i:1a}).
Further, for $\xbar\in\extspace$,
\begin{align*}
  \gext_i(\xbar)
  =
  \gext_{i-1}(\limray{\vv_i}\plusl\xbar)
  &=
  \fext\bigParens{[\vv_1,\dotsc,\vv_{i-1}]\omm
        \plusl \limray{\vv_i}\plusl\xbar}
  \\
  &=
  \fext\bigParens{[\vv_1,\dotsc,\vv_i]\omm\plusl\xbar}.
\end{align*}
The first equality is by \Cref{cor:i:1}(\ref{cor:i:1b}).
The second is by inductive hypothesis.
This proves \eqref{eq:pr:icon-red-decomp-astron-red:2}.
Since $g_i$ is lower semicontinuous, it follows then
by \Cref{pr:h:1}(\ref{pr:h:1a}),
that
$g_i(\xx)=(\lsc{g_i})(\xx)=\gext_i(\xx)$
for $\xx\in\Rn$,
proving \eqref{eq:pr:icon-red-decomp-astron-red:1}.
\end{proof}

\Cref{pr:icon-red-decomp-astron-red} then implies the following simple
properties of iconic reductions:

\begin{proposition}  \label{pr:i:9}
  Let $f:\Rn\rightarrow\Rext$ be convex, and
  let $\ebar\in\corezn$.
  Then:
  \begin{letter}
  \item  \label{pr:i:9a}
    $\fshadd$ is convex and lower semicontinuous.
  \item  \label{pr:i:9b}
    $\fshadextd(\xbar) = \fext(\ebar\plusl\xbar)$
    for all $\xbar\in\extspace$.
  \item  \label{pr:i:9c}
    Either $\fshadextd\leq \fext$ or $\fshadd\equiv +\infty$.
  \item  \label{pr:i:9cons}
    Let $\PP$ be the projection matrix onto
    $(\rspanset{\ebar})^\perp$.
    Then
    $\fshadd(\xx)=\fshadd(\PP\xx)$
    for all $\xx\in\Rn$.
    Hence,
    $\rspanset{\ebar}\subseteq\conssp{\fshadd}$.
  \end{letter}
\end{proposition}

\begin{proof}
By
\Crefequiv{pr:icon-equiv}{pr:icon-equiv:a}{pr:icon-equiv:c},
$\ebar=\VV\omm$
where $\VV=[\vv_1,\dotsc,\vv_k]$ for some
$\vv_1,\dotsc,\vv_k\in\Rn$.
Let $g_0 = \lsc f$ and let
$g_i = \gishadvi$
for $i=1,\dotsc,k$.
By
\Cref{pr:icon-red-decomp-astron-red},
$\fshadd=g_k$ and each $g_i$ is convex and lower semicontinuous,
for $i=0,\dotsc,k$.

\begin{proof-parts}
\pfpart{Part~(\ref{pr:i:9a}):}
This is immediate from the preceding remarks.

\pfpart{Part~(\ref{pr:i:9b}):}
By \Cref{pr:icon-red-decomp-astron-red},
each $\gext_i$ is as given in
\eqref{eq:pr:icon-red-decomp-astron-red:2}.
Thus,
for $\xbar\in\extspace$,
$\fshadextd(\xbar)=\gext_k(\xbar)=\fext(\ebar\plusl\xbar)$.

\pfpart{Part~(\ref{pr:i:9c}):}
Suppose first that
$\vv_i\not\in\resc{g_{i-1}}$ for some $i\in\{1,\dotsc,k\}$.
Then $g_i\equiv+\infty$, by \Cref{thm:i:4}, so $g_j\equiv+\infty$
for $j\geq i$,
and in particular $\fshadd=g_k\equiv+\infty$.

Otherwise,
$\vv_i\in\resc{g_{i-1}}$ for all $i\in\{1,\dotsc,k\}$, so
\[
  \fshadextd
  =
  \gext_k
  \leq
  \gext_{k-1}
  \leq \dotsb \leq
  \gext_0 = \fext,
\]
with each inequality following from
\Cref{cor:i:1}(\ref{cor:i:1c})
(and the last equality from
\Cref{pr:h:1}\ref{pr:h:1aa}).

\pfpart{Part~(\ref{pr:i:9cons}):}
Note that $\rspanset{\ebar}=\colspace{\VV}$ so
$\PP$ is the projection matrix onto $(\colspace{\VV})^\perp$.
Thus,
for $\xx\in\Rn$,
$
  \fshadd(\xx)=\ef(\ebar\plusl\xx)=\ef(\ebar\plusl\PP\xx)=\fshadd(\PP\xx)
$,
with the second equality following from
the Projection Lemma~\ref{lemma:proj}.
That $\rspanset{\ebar}\subseteq\conssp{\fshadd}$
then follows from \Cref{pr:cons-PP}.
\qedhere
\end{proof-parts}
\end{proof}

So all of $f$'s iconic reductions are lower semicontinuous.
However, they are not necessarily closed.
\indexg{reduction-closedness|(}%
The property of \emph{all} of the iconic reductions being closed,
called reduction-closedness, turns
out to exactly characterize when $\fext=\fdub$, as we will see shortly.

\begin{definition}
\indexg{reduction-closedness!defined|(}%
A convex function $f:\Rn\rightarrow\Rext$
is \emph{reduction-closed} if
all of its reductions are closed at all icons, that is,
if $\fshadd$ is closed for every icon
$\ebar\in\corezn$.%
\indexg{reduction-closedness!defined|)}%
\end{definition}

Here are some useful facts about this property:

\begin{proposition} \label{pr:j:1}
  Let $f:\Rn\rightarrow\Rext$ be convex.
  Then the following hold:
  \begin{letter-compact}
  \item \label{pr:j:1a}
    Let $\ebar\in\corezn$, and
    let $g=\fshadd$
    be the reduction of $f$ at $\ebar$.
    Then $\gshad{\ebar'}=\fshad{\ebar\plusl\ebar'}$
    for all $\ebar'\in\corezn$.
    Consequently,
    if $f$ is reduction-closed, then so is $g$.

  \item \label{pr:j:1b}
    $f$ is \emph{not} reduction-closed if and only if
    there exists
    an icon
    $\ebar\in\corezn$ and $\qq,\qq'\in\Rn$ such that
    $\fshadd(\qq)=\fext(\ebar\plusl\qq)=+\infty$
    and
    $\fshadd(\qq')=\fext(\ebar\plusl\qq')=-\infty$.
  \item \label{pr:j:1c}
    If $\inf f > -\infty$ then $f$ is reduction-closed.
  \item \label{pr:j:1d}
    If $f < +\infty$ then $f$ is reduction-closed.
  \end{letter-compact}
\end{proposition}

\begin{proof}
~
\begin{proof-parts}
\pfpart{Part~(\ref{pr:j:1a}):}
For all icons $\ebar'\in\corezn$ and for all $\xx\in\Rn$,
\[
  \gshad{\ebar'}(\xx)
  =
  \gext(\ebar'\plusl\xx)
  =
  \fext(\ebar\plusl\ebar'\plusl\xx)
  =
  \fshad{\ebar\plusl\ebar'}(\xx),
\]
with the second equality
by \Cref{pr:i:9}(\ref{pr:i:9b}).
Thus, $\gshad{\ebar'}=\fshad{\ebar\plusl\ebar'}$.
Therefore, if $f$ is reduction-closed, then for all $\ebar'\in\corezn$,
$\gshad{\ebar'}=\fshad{\ebar\plusl\ebar'}$ is closed,
implying that $g$ is also reduction-closed.

\pfpart{Part~(\ref{pr:j:1b}):}
By definition, $f$ is not reduction-closed if and only if
$\fshadd$ is not closed for some $\ebar\in\corezn$.
Since $\fshadd$ is convex lower semicontinuous (by
\Cref{pr:i:9}\ref{pr:i:9a}),
it is not closed if and only if it is equal to $-\infty$ at some
point, and equal to $+\infty$ at some other point
(by
\Cref{pr:lsc-props}\ref{pr:lsc-props:e}).

\pfpart{Part~(\ref{pr:j:1c}):}
If $\inf f>-\infty$
then, by \Cref{pr:fext-min-exists},
\[ \fshadd(\xx)=\fext(\ebar\plusl\xx)\geq\inf f>-\infty \]
for all
$\ebar\in\corezn$ and $\xx\in\Rn$.
Therefore, by part~(\ref{pr:j:1b}), $f$ must be reduction-closed.

\pfpart{Part~(\ref{pr:j:1d}):}
If $f<+\infty$, then for all $\ebar\in\corezn$,
by \Cref{pr:i:9}(\ref{pr:i:9c}),
either $\fshadd\equiv+\infty$ or
$\fshadd\leq\fext$, implying,
$\fshadd(\xx)\leq\fext(\xx)\leq f(\xx)<+\infty$
for all $\xx\in\Rn$
(with the second inequality from \Cref{pr:h:1}\ref{pr:h:1a}).
In either case, the condition of part~(\ref{pr:j:1b}) is ruled out;
therefore, $f$ is reduction-closed.%
\indexg{reductions, iconic|)}%
\qedhere
\end{proof-parts}
\end{proof}

\indexg{closedness, astral (primal)!extension@of extension|(}%
\indexg{lower semicontinuous extension!astral closedness of|(}%
\indexg{closedness, astral (primal)!reduction-closedness and|(}%
\indexg{reduction-closedness!astral closedness and|(}%
We are now ready to show that reduction-closedness characterizes when $\fext=\fdub$.

\begin{theorem}  \label{thm:dub-conj-new}
  Let $f:\Rn\rightarrow\Rext$ be convex.
  Then $f$ is reduction-closed if and only if
  $\fext=\fdub$.
\end{theorem}

\begin{proof}
Since $\fext=\lscfext$
(\Cref{pr:h:1}\ref{pr:h:1aa}),
since reduction-closedness is actually a property of $\fext$,
and since
$\fstar=(\lsc f)^*$
(\Cref{pr:conj-props}\ref{pr:conj-props:e}),
we can assume without loss of generality that $f$
is lower semicontinuous, replacing it with $\lsc f$ if it is
not.\looseness=-1

We first prove that if $f$ is reduction-closed then
$\fext(\xbar)=\fdub(\xbar)$, by induction on the astral rank of
$\xbar$.
More precisely, we show by induction on $k=0,\dotsc,n$ that for every
lower semicontinuous, reduction-closed convex function
$f:\Rn\rightarrow\Rext$,
and for all $\xbar\in\extspace$, if $\xbar$ has astral rank $k$ then
$\fext(\xbar)=\fdub(\xbar)$.

So suppose that $f$ is convex, lower semicontinuous, and
reduction-closed.
In particular, this implies
that $\fshad{\zero}=\lsc f = f$ is closed.
Also, let $\xbar\in\extspace$ have astral rank $k$.\looseness=-1

In the base case that $k=0$, $\xbar$ must be some point
$\xx\in\Rn$.
Since $f$ is closed, $\fdubs=f$
(\Cref{pr:conj-props-cvx}\ref{pr:conj-props-cvx:b}),
so
\begin{equation}  \label{eq:thm:dub-conj-new:1}
  \fdub(\xx)=\fdubs(\xx)=f(\xx)=\fext(\xx),
\end{equation}
where the first equality follows from
\Cref{thm:fext-dub-sum}(\ref{thm:fext-dub-sum:fdubs}) and the last
from lower semicontinuity of $f$
(by \Cref{pr:h:1}\ref{pr:h:1a}).

For the inductive step when $k>0$, we can write
$\xbar=\limray{\vv}\plusl\xbarperp$ where $\vv$ is $\xbar$'s dominant
direction, and $\xbarperp$, its projection orthogonal to $\vv$,
has astral rank $k-1$
(\Cref{pr:h:6}).
Let $g=\fshadv$ be the reduction of $f$ at $\limray{\vv}$.
Then $g$ is convex and lower semicontinuous
by
\Cref{cor:i:1}(\ref{cor:i:1a}),
and is reduction-closed by \Cref{pr:j:1}(\ref{pr:j:1a}).
Thus, we can apply our inductive hypothesis to $g$
at $\xbarperp$, yielding
\[
  \fext(\xbar)
  =
  \gext(\xbarperp)
  =
  \gdub(\xbarperp)
  =
  \fdub(\xbar),
\]
where the %
equalities follow, respectively, from
\Cref{cor:i:1}(\ref{cor:i:1b}),
our inductive hypothesis,
and
\Cref{thm:i:5}.
This completes the induction.

Conversely, suppose now that $\fext=\fdub$.
Further, suppose by way of contradiction that $f$ is not
reduction-closed.
Then by \Cref{pr:j:1}(\ref{pr:j:1b}), there exists an
icon $\ebar\in\extspace$ and $\qq,\qq'\in\Rn$ such that
$\fext(\ebar\plusl\qq)=+\infty$
and
$\fext(\ebar\plusl\qq')=-\infty$.
Then for all $\xx\in\Rn$,
\begin{equation}  \label{eq:i:8}
  \fext(\ebar\plusl\xx)
  =
  \fdub(\ebar\plusl\xx)
  =
  \sup_{\uu\in\Rn} \bigBracks{-\fstar(\uu) \plusl \ebar\cdot\uu \plusl \xx\cdot\uu},
\end{equation}
where the first equality is by assumption, and the second
by \Cref{thm:fext-dub-sum}(\ref{thm:fext-dub-sum:b}).
In particular, when $\xx=\qq'$, the left-hand side is equal to
$-\infty$, so the expression inside the supremum is
equal to $-\infty$ for all $\uu\in\Rn$.
Since $\qq'\cdot\uu\in\R$, this means that
$-\fstar(\uu) \plusl \ebar\cdot\uu = -\infty$
for all $\uu\in\Rn$.
In turn, this implies, by \eqref{eq:i:8}, that actually
$\fext(\ebar\plusl\xx)=-\infty$ for all $\xx\in\Rn$.
But this contradicts that
$\fext(\ebar\plusl\qq)=+\infty$.

Thus, $f$ is reduction-closed.%
\indexg{reduction-closedness|)}%
\indexg{closedness, astral (primal)!reduction-closedness and|)}%
\indexg{reduction-closedness!astral closedness and|)}%
\end{proof}

Combined with \Cref{pr:j:1}(\ref{pr:j:1c},\ref{pr:j:1d}),
we immediately obtain the following corollary:

\begin{corollary}  \label{cor:all-red-closed-sp-cases}
  Let $f:\Rn\rightarrow\Rext$ be convex.
  If either $\inf f > -\infty$ or $f<+\infty$
  then $\fext=\fdub$.
\end{corollary}

The next example shows that it is indeed possible that
$\fext\neq\fdub$,
even if $f$ is convex, closed, and proper.

\begin{example}[Astral biconjugate of the restricted linear function]
\label{ex:biconj:notext}
\indexg{Restricted linear function!astral biconjugate of|(}%
Consider the restricted linear function from
\Cref{ex:negx1-else-inf},
\[
   f(x_1,x_2) =
   \begin{cases}
     -x_1
     & \text{if $x_2\geq 0$,}
   \\
     +\infty
     & \text{otherwise,}
   \end{cases}
\]
and its conjugate $\fstar$, which is the indicator of the set $\set{\uu\in\R^2:\:u_1=-1,u_2\le 0}$.
Let $\ebar=\limray{\ee_1}$.
By \Cref{thm:fext-dub-sum}(\ref{thm:fext-dub-sum:b}),
for all $\xx\in\R^2$,
\[
   \fdub(\ebar\plusl\xx)
   =
   \sup_{\uu\in\R^2:\:u_1=-1,u_2\le 0} (\ebar\plusl\xx)\cdot\uu
   =-\infty.
\]
On the other hand,
\begin{align}
\label{eq:j:2}
   \fext(\ebar\plusl\xx)
   &=
   \begin{cases}
     -\infty & \text{if $x_2\geq 0$,}
     \\
     +\infty & \text{otherwise.}
   \end{cases}
\end{align}
Thus, $+\infty=\fext(\ebar\plusl\xx)\neq\fdub(\ebar\plusl\xx)=-\infty$
if $x_2<0$.
\eqref{eq:j:2} also shows that the reduction
$\fshadd$ is not closed;
thus, consistent with \Cref{thm:dub-conj-new}, $f$ is not
reduction-closed.%
\indexg{Restricted linear function!astral biconjugate of|)}%
\end{example}

In general, as we show next,
even if $f$ is not reduction-closed, $\fext$
and $\fdub$ must agree at every
point in the closure of the effective domain of
$\fext$ (which is the same as the closure of $\dom f$,
by \Cref{pr:h:1}\ref{pr:h:1c}) and
at every point $\xbar\in\extspace$ with $\fdub(\xbar)>-\infty$.
This means that
at points $\xbar$ where $\fext$ and $\fdub$ differ, we can
say exactly what values each function will take, namely,
$\fext(\xbar)=+\infty$ and $\fdub(\xbar)=-\infty$,
as was the case in \Cref{ex:biconj:notext}.

\begin{theorem}  \label{thm:fext-neq-fdub}
  Let $f:\Rn\rightarrow\Rext$ be convex, and let $\xbar\in\extspace$.
  If $\xbar\in\cldom{f}$ or $\fdub(\xbar)>-\infty$
  then $\fext(\xbar)=\fdub(\xbar)$.

  Consequently,
  if $\fext(\xbar)\neq\fdub(\xbar)$ then
  $\fext(\xbar)=+\infty$ and $\fdub(\xbar)=-\infty$.
\end{theorem}

\begin{proof}
The proof is similar to the first part of the proof of
\Cref{thm:dub-conj-new}.

By Propositions~\ref{pr:h:1}(\ref{pr:h:1c})
and~\ref{pr:fextstar-is-fstar},
$\cldom{f}=\cldomfext$ and
$\fstar=\fextstar$,
which means the theorem's claim can be stated entirely in terms of
$\fext$.
Thus,
since $\fext=\lscfext$
(by \Cref{pr:h:1}\ref{pr:h:1aa}),
it suffices to prove the result for $\lsc f$.
Therefore, without loss of generality, we assume henceforth that $f$
is lower semicontinuous.

First, suppose that $f$ is not closed,
and that either
$\xbar\in\cldom{f}$ or $\fdub(\xbar)>-\infty$.
Under these assumptions, we show that $\fext(\xbar)=\fdub(\xbar)$.
Since $f$ is lower semicontinuous but not closed,
$f$ must be infinite everywhere so that $f(\xx)\in\{-\infty,+\infty\}$
for all $\xx\in\Rn$, and furthermore, $f\not\equiv+\infty$
(by \Cref{pr:lsc-props}\ref{pr:lsc-props:e}).
These facts imply $\fstar\equiv+\infty$, so $\fdub\equiv-\infty$.
From our assumptions, this further implies that $\xbar\in\cldom{f}$,
which means there exists a sequence $\seq{\xx_t}$ in $\dom{f}$ with
$\xx_t\rightarrow\xbar$.
For all $t$,
since $f(\xx_t)<+\infty$, the foregoing implies we must actually have
$f(\xx_t)=-\infty$, so
\[
   -\infty
   =
   \liminf f(\xx_t)
   \geq
   \fext(\xbar)
\]
by definition of $\fext$.
Therefore, as claimed,
$\fext(\xbar)=-\infty=\fdub(\xbar)$ in this case.

As in the proof of
\Cref{thm:dub-conj-new}, we prove the theorem for the case when $f$ is closed by induction
on the astral rank of $\xbar$.
Specifically,
we prove by induction on $k=0,\dotsc,n$ that for every
lower semicontinuous, convex function $f$,
and for all $\xbar\in\extspace$,
if $\xbar$ has astral rank $k$,
and if either
$\xbar\in\cldom{f}$ or $\fdub(\xbar)>-\infty$,
then
$\fext(\xbar)=\fdub(\xbar)$.

Let $f$ be such a function and let $\xbar$ have astral rank $k$
with
$\xbar\in\cldom{f}$ or $\fdub(\xbar)>-\infty$.
We further assume that $f$ is closed, since the case that $f$ is not
closed was handled above.

In the base case that $k=0$, $\xbar$ is a point $\xx\in\Rn$.
Since $f$ is closed, $\ef(\xx)=f(\xx)=\fdubs(\xx)=\fdub(\xx)$,
where the first equality follows from lower semicontinuity of $f$
(by \Cref{pr:h:1}\ref{pr:h:1a}) and the last from
\Cref{thm:fext-dub-sum}(\ref{thm:fext-dub-sum:fdubs}).

For the inductive step when $k>0$,
by \Cref{pr:h:6}, we can write
$\xbar=\limray{\vv}\plusl\xbarperp$ where $\vv\in\Rn$ is $\xbar$'s
dominant direction, and $\xbarperp$, the projection of $\xbar$
orthogonal to~$\vv$, has astral rank $k-1$.
If $\vv\not\in\resc{f}$, then \Cref{thm:i:4}
immediately implies that
$\fext(\xbar)=+\infty=\fdub(\xbar)$ in this case.
Therefore, we assume henceforth that $\vv\in\resc{f}$.

Let
$g=\fshadv$ be the reduction of $f$ at $\limray{\vv}$,
which is convex and lower semicontinuous
(\Cref{cor:i:1}\ref{cor:i:1a}).
Then
$ \fext(\xbar) = \gext(\xbarperp)$
by \Cref{cor:i:1}(\ref{cor:i:1b}),
and
$ \gdub(\xbarperp) = \fdub(\xbar) $
by \Cref{thm:i:5}.
Therefore, to prove
$\fext(\xbar)=\fdub(\xbar)$, it suffices to show
$\gext(\xbarperp)=\gdub(\xbarperp)$.

If $\fdub(\xbar)>-\infty$ then
$ \gdub(\xbarperp) > -\infty $,
so that, by inductive hypothesis,
$\gext(\xbarperp)=\gdub(\xbarperp)$.

Otherwise, we must have $\xbar\in\cldom{f}$,
meaning there exists a sequence $\seq{\xx_t}$ in $\dom f$ with
$\xx_t\rightarrow\xbar$.
Thus, for each $t$,
by \Cref{pr:d2}(\ref{pr:d2:a},\ref{pr:d2:b}),
$g(\xperpt)\leq f(\xx_t)<+\infty$,
so $\xperpt\in\dom g$.
Since $\xperpt\rightarrow\xbarperp$
(by \Cref{pr:h:5}\ref{pr:h:5b}),
this means $\xbarperp\in\cldom{g}$.
Therefore, by inductive hypothesis,
$\gext(\xbarperp)=\gdub(\xbarperp)$
in this case as well.

This completes the induction and the proof.%
\indexg{closedness, astral (primal)!extension@of extension|)}%
\indexg{lower semicontinuous extension!astral closedness of|)}%
\end{proof}

\section{A dual characterization of reduction-closedness}
\label{sec:dual-char-ent-clos}

\indexg{reduction-closedness!barrier cone and|(}%
\indexg{barrier cone (of function)!reduction-closedness and|(}%
We next give a dual characterization of when a function $f$ is reduction-closed,
and thus when $\fext=\fdub$.
In \Cref{sec:slopes}, we defined $\slopes{f}$, the
barrier cone of $f$, which
includes all of $\cone(\dom{\fstar})$
as well as $f$'s vertical barrier cone, $\vertsl{f}$.
We show now that $f$ is reduction-closed if and only if
$\slopes{f}=\cone(\dom{\fstar})$, that is,
if and only if $f$'s vertical barrier cone is already entirely included
in $\cone(\dom{\fstar})$.
Whereas reduction-closedness would appear to be an astral property
involving the behavior of $\fext$ on every galaxy, this
shows that actually the property can be precisely characterized just in
terms of $f$ and its conjugate.%
\indexg{reduction-closedness!barrier cone and|)}%
\indexg{barrier cone (of function)!reduction-closedness and|)}

\indexg{barrier cone (of function)!recession cone and|(}%
\indexg{recession cone (standard)!barrier cone and|(}%
We begin by generalizing the standard result relating $\resc f$ and $\cone(\dom\fstar)$
when $f$ is closed, proper, and convex
(\Cref{pr:rescpol-is-con-dom-fstar}) to also cover improper lower
semicontinuous functions that are not identically $+\infty$:

\begin{theorem}  \label{thm:rescpol-is-slopes}
  Let $f:\Rn\rightarrow\Rext$ be convex and lower semicontinuous
  with $f\not\equiv+\infty$.
  Then
  \[ \resc{f} = \polar{(\slopes{f})} \]
  and
  \[ \rescpol{f} = \cl(\slopes{f}). \]
\end{theorem}

\begin{proof}
Let $f'=\expex\circ f$.
Then
\[
   \resc{f}
   =
   \resc{f'}
   =
   \polar{\bigParens{\cone(\dom{\fpstar})}}
   =
   \polar{(\slopes{f})}.
\]
The first equality is by
\Cref{pr:j:2}(\ref{pr:j:2b}).
The second is by \Cref{pr:rescpol-is-con-dom-fstar}
applied to $f'$
(which is convex, lower semicontinuous and lower-bounded by
\Cref{pr:j:2}\ref{pr:j:2a}\ref{pr:j:2lsc}, and so also closed and
proper).
The third is
by \Cref{thm:dom-fpstar}.

This proves the first claim of the
\namecref{thm:rescpol-is-slopes}.
Taking the polar of both sides then yields the second claim
by \Cref{pr:polar-props}(\ref{pr:polar-props:c}),
noting that $\slopes{f}$ is a convex cone
(by \Cref{thm:slopes-equiv}).%
\indexg{barrier cone (of function)!recession cone and|)}%
\indexg{recession cone (standard)!barrier cone and|)}%
\end{proof}

In particular, combined with
\Cref{pr:rescpol-is-con-dom-fstar},
when $f$ is closed, proper and convex,
this shows that
$\slopes{f}$ and $\cone(\dom{\fstar})$ have the same closures (in
$\Rn$) and therefore can differ only in their relative boundaries.
\indexg{reduction-closedness!barrier cone and|(}%
\indexg{barrier cone (of function)!reduction-closedness and|(}%
Our dual characterization of reduction-closedness states that $f$ is reduction-closed
precisely when these two sets coincide.

We prove each direction of the characterization separately.

\begin{theorem}  \label{thm:slopescone-implies-entclosed}
  Let $f:\Rn\rightarrow\Rext$ be convex and lower semicontinuous.
  If $\slopes{f}=\cone(\dom{\fstar})$
  then $f$ is reduction-closed.
\end{theorem}

\begin{proof}
We need to show that
$\fshadd$ is closed for all $\ebar\in\corezn$.
The proof is by induction on the astral rank of $\ebar$ (which, like that of all points in $\extspace$, cannot exceed $n$).
More precisely, we prove the following by induction on $k=0,\dotsc,n$:
For every convex, lower semicontinuous function
$f:\Rn\rightarrow\Rext$, and
for every $\ebar\in\corezn$ of astral rank at most $k$,
if $\slopes{f}=\cone(\dom{\fstar})$ then $\fshadd$ is closed.

In the base case that $k=0$,
the only icon of astral rank $0$ is $\ebar=\zero$
(by \Cref{pr:i:8}\ref{pr:i:8b}). Consider the corresponding
reduction $\fshad{\zero}=\lsc f=f$. If $f>-\infty$ then $f$ is closed,
because it is lower semicontinuous. Otherwise,
$f(\qq)=-\infty$ for some $\qq\in\Rn$.
In that case, $\fstar\equiv+\infty$
(by \Cref{pr:conj-props}\ref{pr:conj-props:c2}),
so $\dom{\fstar}=\emptyset$ and
$\cone(\dom{\fstar})=\set{\zero}$. By the theorem's assumption,
$\slopes{f}=\cone(\dom{\fstar})=\set{\zero}$, implying that
$\resc{f}=\polar{(\slopes{f})}=\Rn$,
by \Cref{thm:rescpol-is-slopes}.
By definition of recession cone (Eq.~\ref{eqn:resc-cone-def}),
this in turn implies that $f(\xx)\leq f(\qq)=-\infty$ for all
$\xx\in\Rn$, and so that
$f\equiv-\infty$, which is closed, finishing the proof of the base case.

For the inductive step, we assume $k>0$ and that the inductive
hypothesis holds for $k-1$.
Consider $\ebar$ of rank $k$, and write $\ebar=\limray{\vv}\plusl\ebarperp$, where $\vv$ is
the dominant direction of $\ebar$, and $\ebarperp$ is
the projection of $\ebar$ orthogonal to $\vv$.
Let
$g=\fshadv$
be the reduction of $f$ at $\limray{\vv}$,
which is convex and lower semicontinuous
(\Cref{thm:a10-nunu}).
Note that $\gshad{\ebarperp}=\fshad{\ebar}$ by
\Cref{pr:j:1}(\ref{pr:j:1a}).
If $\vv\not\in\resc{f}$ then
$g\equiv+\infty$ (by \Cref{thm:i:4}),
so $\gshad{\ebarperp}\equiv+\infty$;
therefore, $\fshad{\ebar}=\gshad{\ebarperp}\equiv+\infty$ is closed.

In the remaining case, $\vv\in\resc{f}$.
If $f\equiv+\infty$ then $\fshad{\ebar}\equiv+\infty$, which is closed.
Otherwise, $f\not\equiv+\infty$ so $g\not\equiv+\infty$ by
\Cref{pr:d2}(\ref{pr:d2:c}).
Moreover,
\[
  \slopes{g}=(\slopes{f})\cap L=\bigParens{\cone(\dom{\fstar})}\cap L=\cone(\dom{\gstar}),
\]
where $L=\{\vv\}^\perp$, and
where the first and last equalities are by
\Cref{thm:e1}(\ref{thm:e1:c},\ref{thm:e1:d}).
Since $\ebarperp$ has astral rank $k-1$
(by~\Cref{pr:h:6}), we can apply our
inductive hypothesis to $g$ and $\ebarperp$, and obtain that
$\fshad{\ebar}=\gshad{\ebarperp}$ is closed,
completing the induction and the proof.
\end{proof}

Next, we prove the converse:

\begin{theorem}  \label{thm:entclosed-implies-slopescone}
  Let $f:\Rn\rightarrow\Rext$ be convex and lower semicontinuous.
  If $f$ is reduction-closed
  then $\slopes{f}=\cone(\dom{\fstar})$.
\end{theorem}

\begin{proof}
If $f\equiv+\infty$ then $\fstar\equiv-\infty$ so
$\slopes{f}=\Rn=\cone(\dom{\fstar})$.
We therefore assume henceforth that $f\not\equiv+\infty$.

We proceed by induction on the dimension of $\cone(\dom{\fstar})$.
That is, we prove by induction on $k=0,\dotsc,n$ that
for every convex, lower semicontinuous function $f$
with $f\not\equiv+\infty$,
if $f$ is reduction-closed and
$\cone(\dom{\fstar})\subseteq L$ for some linear
subspace $L\subseteq\Rn$ of dimension at most $k$, then
$\slopes{f}=\cone(\dom{\fstar})$.

We begin with the base case that $k=0$.
Then
$\cone(\dom{\fstar})\subseteq L =\{\zero\}$,
which means $\dom{\fstar}$ is either $\emptyset$ or $\set{\zero}$.
We need to show that in both cases, $\barr f=\set{\zero}$.
If $\dom{\fstar}=\{\zero\}$ then
$\inf f = -\fstar(\zero) > -\infty$, implying that $f$ is closed and
proper.
Therefore,
\[
  \{\zero\}
  =
  \cl\bigParens{\cone(\dom{\fstar})}
  =
  \rescpol{f}
  =
  \cl(\slopes{f}),
\]
where the second and third equalities are
by \Cref{pr:rescpol-is-con-dom-fstar}
and \Cref{thm:rescpol-is-slopes}, respectively.
This implies that $\slopes{f}=\set{\zero}$.

In the remaining case, $\dom{\fstar}=\emptyset$, so $\fstar\equiv+\infty$
and $\fdubs\equiv-\infty$.
Since $f$ is reduction-closed, $f=\lsc f=\fshad{\zero}$ is closed, and
therefore $f=\fdubs\equiv-\infty$. Hence, $\resc{f}=\Rn$, and $\barr f=\set{\zero}$
by \Cref{thm:rescpol-is-slopes}, completing the proof of the base case.

For the inductive step, let $f$ and $L$ satisfy the conditions of
the induction with $k>0$.

If $\cone(\dom{\fstar})$ is contained in a linear space of dimension $k-1$,
then $\barr f=\cone(\dom{\fstar})$ by inductive hypothesis, so we assume henceforth that this
is not the case. In particular, this means that $\cone(\dom{\fstar})\ne\emptyset$ and thus
$\dom{\fstar}\ne\emptyset$, which implies that $f>-\infty$
(by \Cref{pr:conj-props}\ref{pr:conj-props:c2}),
and therefore $f$ is closed and proper.

For the sake of contradiction, suppose that $f$ is reduction-closed,
but that $\barr f\ne\cone(\dom{\fstar})$. Since $\cone(\dom{\fstar})\subseteq\barr f$
(by \Cref{thm:slopes-equiv}), this means that there exists some point $\uu$
in $(\slopes{f})\setminus\cone(\dom{\fstar})$.

Since $f$ is closed and proper,
\[
  \uu
  \in
  \slopes{f}
  \subseteq
  \cl(\slopes{f})
  =
  \rescpol{f}
  =
  \cl\bigParens{\cone(\dom{\fstar})},
\]
where the two equalities are respectively
from
\Cref{thm:rescpol-is-slopes} and
\Cref{pr:rescpol-is-con-dom-fstar}.
On the other hand, $\uu\not\in\cone(\dom{\fstar})$.
Therefore, $\uu$ is a relative boundary point of
$\cone(\dom{\fstar})$.

Before continuing, we pause to prove the following
lemma regarding convex cones and their relative boundary points:

\begin{lemma}  \label{lem:rel-bnd-convcone}
\indexg{separation by hyperplane (standard)|(}%
  Let $K\subseteq\Rn$ be a convex cone, and
  let $\uu\in(\cl K) \setminus(\ri K)$ be a relative boundary point of
  $K$.
  Then there exists a vector $\vv\in\Rn$ such that all of the
  following hold:
  \begin{letter-compact}
  \item
    $\uu\cdot\vv = 0$.
  \item
    $\xx\cdot\vv < 0$ for all $\xx\in\ri K$
    (implying that there exists $\yy\in K$ with $\yy\cdot\vv<0$).
  \item
    $\xx\cdot\vv \leq 0$ for all $\xx\in K$.
    That is, $\vv\in \Kpol$.
  \end{letter-compact}
\end{lemma}

\begin{proofx}
The set $\ri K$ is nonempty
(by \Cref{pr:ri-props}\ref{pr:ri-props:roc-thm6.2b})
and relatively open, and the singleton
$\set{\uu}$ is an affine set disjoint from $\ri K$.
Therefore,
by \Cref{roc:thm11.2},
there exists a hyperplane
that contains $\set{\uu}$, and such that $\ri K$ is contained in one of the
open halfspaces associated with the hyperplane.
In other words, there exist $\vv\in\Rn$ and $\beta\in\R$ such that
$\uu\cdot\vv=\beta$ and $\xx\cdot\vv<\beta$ for all $\xx\in\ri K$.

We claim that $\xx\cdot\vv\leq\beta$ for all $\xx\in K$.
This is because if $\xx\in K$ then it is also in
$\cl{K} = \cl(\ri{K})$
(\Cref{pr:ri-props}\ref{pr:ri-props:roc-thm6.3}),
implying there exists a sequence $\seq{\xx_t}$ in $\ri{K}$ with
$\xx_t\rightarrow\xx$.
Since $\xx_t\cdot\vv<\beta$ for all $t$,
it follows that
$\xx\cdot\vv=\lim(\xx_t\cdot\vv)\leq\beta$.

Next, we claim that $\beta=0$.
To see this, let $\epsilon>0$.
Then the open set
$\set{\xx\in\Rn :\: \xx\cdot\vv > \beta - \epsilon}$
is a neighborhood of $\uu$.
Since $\uu\in\cl{K}$, this neighborhood must intersect~$K$ at some
point $\zz$
(\Cref{pr:closure:intersect}\ref{pr:closure:intersect:a});
that is, $\zz\in K$ and $\zz\cdot\vv>\beta-\epsilon$.
Since $K$ is a cone, $2\zz$ is also in $K$, implying
$2(\beta-\epsilon) < (2\zz)\cdot\vv \leq \beta$, and therefore that
$\beta < 2\epsilon$. Similarly, $\zz/2$ is in~$K$, implying that
$(\beta-\epsilon)/2 < (\zz/2)\cdot\vv \leq \beta$, and therefore that
$-\epsilon < \beta$. Thus, $-\epsilon<\beta<2\epsilon$ for all $\epsilon>0$, which implies that $\beta=0$.

We have shown that $\uu\cdot\vv=0$;
$\xx\cdot\vv < 0$ for all $\xx\in\ri K$;
and
$\xx\cdot\vv \leq 0$ for all $\xx\in K$.
By definition, this last fact means that $\vv\in\Kpol$.
And since $\ri K$ is not empty,
there exists $\yy\in \ri K \subseteq K$ with $\yy\cdot\vv<0$.%
\indexg{separation by hyperplane (standard)|)}%
\end{proofx}

Since $\uu$ is a relative boundary point of
$\cone(\dom{\fstar})$,
\Cref{lem:rel-bnd-convcone} implies that
there exists $\vv\in\Rn$ with $\uu\cdot\vv=0$ and $\yy\cdot\vv<0$ for
some $\yy\in\cone(\dom{\fstar})$.
Furthermore,
$\vv\in \polar{\regParens{\cone(\dom{\fstar})}} = \resc{f}$
(by \Cref{pr:rescpol-is-con-dom-fstar}).

Let
$g=\fshadv$
be the reduction of $f$ at $\limray{\vv}$.
Also, let
$M=\{\vv\}^\perp$
be the linear subspace orthogonal to $\vv$,
and let $L'=L\cap M$.
Note that $\yy\in\cone(\dom{\fstar})\subseteq L$ but
$\yy\not\in M$, because $\yy\inprod\vv<0$, and so $\yy\not\in L'$.
Thus, $L'\subseteq L$ but $L'\neq L$ so $L'$ has dimension strictly
less than $L$.

The function $g$ is convex and lower semicontinuous (by
\Cref{thm:a10-nunu})
with $g\not\equiv+\infty$ (by \Cref{pr:d2}\ref{pr:d2:c}),
and is also reduction-closed (by \Cref{pr:j:1}\ref{pr:j:1a}).
Moreover, 
\[
   \cone(\dom{\gstar})
   =
   \cone(\dom{\fstar}) \cap M
   \subseteq
   (L \cap M) = L',
\]
where the first equality is by \Cref{thm:e1}(\ref{thm:e1:c}).
Therefore, 
\[
   \cone(\dom{\gstar})
   =
   \barr g
   =
   (\barr f) \cap M,
\]
where the first equality is
by inductive hypothesis applied to $g$ and $L'$, and the second by
\Cref{thm:e1}(\ref{thm:e1:d}).
Since
$\uu\in\barr f$
and $\uu\in M$, we obtain
\[
  \uu\in(\barr f) \cap M=\cone(\dom{\gstar})\subseteq\cone(\dom{\fstar}),
\]
which contradicts the assumption that $\uu\not\in\cone(\dom{\fstar})$. Thus, we must have $\barr f=\cone(\dom{\fstar})$.
This completes the induction and the proof.
\end{proof}

Combined with \Cref{thm:dub-conj-new}, we thus have proved:

\begin{corollary}  \label{cor:ent-clos-is-slopes-cone}
\indexg{closedness, astral (primal)!extension@of extension|(}%
\indexg{lower semicontinuous extension!astral closedness of|(}%
\indexg{closedness, astral (primal)!barrier cone and|(}%
\indexg{barrier cone (of function)!astral closedness and|(}%
  Let $f:\Rn\rightarrow\Rext$ be convex.
  Then the following are equivalent:
  \begin{letter-compact}
  \item  \label{cor:ent-clos-is-slopes-cone:b}
    $\fext=\fdub$.
  \item  \label{cor:ent-clos-is-slopes-cone:a}
    $f$ is reduction-closed.
  \item  \label{cor:ent-clos-is-slopes-cone:c}
    $\slopes{f}=\cone(\dom{\fstar})$.
  \end{letter-compact}
\end{corollary}

\begin{proof}
By Propositions~\ref{pr:h:1}(\ref{pr:h:1aa})
and~\ref{pr:slopes-same-lsc},
and since
$\fstar=(\lsc f)^*$
(by \Cref{pr:conj-props}\ref{pr:conj-props:e}),
we can assume without loss of generality that $f$
is lower semicontinuous, replacing it with $\lsc f$ if it is not.

In this case,
the equivalence of
(\ref{cor:ent-clos-is-slopes-cone:b})
and
(\ref{cor:ent-clos-is-slopes-cone:a})
was proved in \Cref{thm:dub-conj-new}.

The equivalence of
(\ref{cor:ent-clos-is-slopes-cone:a})
and
(\ref{cor:ent-clos-is-slopes-cone:c})
was proved by
\Cref{thm:slopescone-implies-entclosed,thm:entclosed-implies-slopescone}.%
\indexg{reduction-closedness!barrier cone and|)}%
\indexg{barrier cone (of function)!reduction-closedness and|)}%
\indexg{closedness, astral (primal)!extension@of extension|)}%
\indexg{lower semicontinuous extension!astral closedness of|)}%
\indexg{closedness, astral (primal)!barrier cone and|)}%
\indexg{barrier cone (of function)!astral closedness and|)}%
\end{proof}

\chapter{Calculus rules for extensions}
\label{sec:calculus:extensions}

This chapter develops some simple tools for computing
extensions of functions constructed using standard operations like
addition of two functions or the composition of a function with a
linear map.

\section{Scalar multiple}

\indexg{scalar multiples (astral)!extension of|(}%
\indexg{lower semicontinuous extension!scalar multiple@of scalar multiple|(}%
The extension of a nonnegative scalar multiple $\lambda f$ of a
function $f:\Rn\rightarrow\Rext$ is straightforward to compute:

\begin{proposition}  \label{pr:scal-mult-ext}
  Let $f:\Rn\rightarrow\Rext$, let $\lambda\in\Rpos$, and let
  $\xbar\in\extspace$.
  Then $\lamfext(\xbar)=\lambda\fext(\xbar)$.
  That is, $\lamfext=\lambda \fext$.
\end{proposition}

\begin{proof}
If $\lambda=0$ then $\lambda f\equiv 0$ so
$\lamfext(\xbar)=0=\lambda \fext(\xbar)$.

Otherwise, $\lambda>0$.
From the definition of extension (Eq.~\ref{eq:e:7}), and
since $\ex\mapsto \lambda\ex$ is strictly increasing over
$\ex\in\Rext$,
\begin{align*}
  \lamfext(\xbar)
  &=
  \InfseqLiminf{\seq{\xx_t}}{\Rn}{\xx_t\rightarrow \xbar}
               {\lambda f(\xx_t)}
  \\
  &=
  \lambda
  \InfseqLiminf{\seq{\xx_t}}{\Rn}{\xx_t\rightarrow \xbar}
               {f(\xx_t)}
  =
  \lambda\fext(\xbar).%
\indexg{scalar multiples (astral)!extension of|)}%
\indexg{lower semicontinuous extension!scalar multiple@of scalar multiple|)}%
\qedhere
\end{align*}
\end{proof}

\section{Sums of functions}
\label{sec:calc-ext-sum-fcns}

\indexg{sum of functions (standard)!extension of|(}%
\indexg{lower semicontinuous extension!sum of functions@of sum of functions|(}%
Suppose $h=f+g$, where
$f:\Rn\rightarrow\Rext$ and $g:\Rn\rightarrow\Rext$.
We might expect $h$'s extension to be the sum of
the extensions of $f$ and $g$, so that $\hext=\fext+\gext$
(provided, of course, that the sums involved are defined).
We will see that
this is indeed the case under various conditions.
First, however, we give two examples showing that,
in general,
$\hext(\xbar)$ need not be equal to $\fext(\xbar)+\gext(\xbar)$,
even when $\fext(\xbar)$ and $\gext(\xbar)$ are summable:

\begin{example}[Positives and negatives]
\label{ex:pos-neg}
\indexg{Positives and negatives|(}%
Let $P=\R_{>0}$ and $N=\R_{<0}$ denote, respectively, the sets of positive and negative reals, and let $f=\indf{P}$ and $g=\indf{N}$ be the associated indicator functions (see Eq.~\ref{eq:indf-defn}). Then $P\cap N=\emptyset$, so $h=f+g\equiv+\infty$, and therefore also $\eh\equiv+\infty$. On the other hand, $\Pbar=[0,+\infty]$ and $\Nbar=[-\infty,0]$, so by \Cref{pr:inds-ext},
\[
  \ef+\eg=\indfa{\Pbar}+\indfa{\Nbar}=\indfa{\Pbar\cap\Nbar}=\indfa{\set{0}},
\]
so $\ef(0)+\eg(0)=0$ but $\eh(0)=+\infty$.%
\indexg{Positives and negatives|)}%
\end{example}

In the previous example, the functions $f$ and $g$ are not lower semicontinuous at $0$, which is also the point where $\ef+\eg$ disagrees with $\eh$.
The next example shows that even when $f$ and $g$ are lower semicontinuous, we do not necessarily have $\eh=\ef+\eg$.

\begin{figure}
  \centering
  \includegraphics{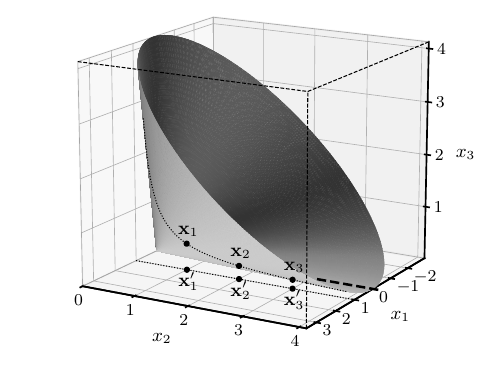}
  \mycaption{Sideways cone}{%
    \indexf{Sideways cone}%
    The cone $K$ from \Cref{ex:KLsets:extsum-not-sum-exts}.
    The cone is tangent to the horizontal plane $L$ corresponding to the ``bottom wall''
    of the figure. The intersection of $K$ and $L$ is
    the ray $K\cap L =\set{\lambda\ee_2:\:\lambda\in\Rpos}$ depicted as the dashed line.
    The closure of the ray is the set $\KLbar=\set{\alpha\ee_2:\:\alpha\in\Rextpos}$.
    The sequence of points $\xx_t=\trans{[1,\,t,\,1/(2t)]}$ is entirely in $K$;
    the sequence of points $\xx'_t=\trans{[1,\,t,\,0]}$ is entirely in $L$.
    Both sequences converge to
    $\xbar=\limray{\ee_2}\plusl\ee_1$, so $\xbar\in\Kbar\cap\Lbar$, but
    $\xbar\not\in\KLbar$.
  }%
  \label{fig:sideways-cone}%
\end{figure}

\begin{example}[Sideways cone]
\label{ex:KLsets:extsum-not-sum-exts}
\indexg{Sideways cone|(}%
The standard
\indexg{second-order cone}%
second-order cone in $\R^3$ is the set
\begin{equation}  \label{eq:stan-2nd-ord-cone}
  \biggBraces{ \zz\in\R^3 :\: \sqrt{z_1^2 + z_2^2} \leq z_3 },
\end{equation}
which is a classic, upward-oriented ``ice-cream cone''
in which every horizontal slice, with $z_3$ held to any nonnegative
constant, is a disc of radius $z_3$.

Let
$K$ be this same cone rotated so that its axis is instead pointing in the direction of the vector $\uu=\trans{[0,1,1]}$ (see \Cref{fig:sideways-cone}).
To derive an analytic description of $K$, note that it consists of all
points $\xx\in\R^3$ whose angle $\theta$ with $\uu$ is at
most $45^\circ$ so that
\[
  \xx\inprod\uu
  = \norm{\xx}\norm{\uu} \cos \theta
  \ge \norm{\xx}\norm{\uu} \cos\BigParens{\frac{\pi}{4}}
  =   \norm{\xx}.
\]
Squaring, rearranging, and restricting to $\xx$
with $\xx\inprod\uu\ge 0$ then yields
\begin{equation}  \label{eqn:bad-set-eg-K}
  K
  =\set{\xx\in\R^3:\: x_1^2\le 2 x_2 x_3\text{ and }x_2,x_3\ge 0}.
\end{equation}
(In the same way,
Eq.\ref{eq:stan-2nd-ord-cone} can be derived by instead setting
$\uu=\ee_3$.)%
\indexg{Sideways cone|)}%

\indexg{Sideways cone!extension of sum and|(}%
Next, let $L\subseteq\R^3$ be the plane
\begin{equation}
\label{eqn:bad-set-eg-L}
  L
  =
  \Braces{
    \xx\in\R^3 :\: x_3 = 0
  },
\end{equation}
and let $f=\indK$ and $g=\indL$ be the corresponding indicator
functions (as defined in Eq.~\ref{eq:indf-defn}).
Then $h=f+g$ is equal to $\indf{K\cap L}$,
the indicator function of the ray
\begin{equation}  \label{eqn:bad-set-eg-KL}
  K\cap L
  =
  \Braces{
    \xx\in\R^3 :\: x_1 = x_3 = 0,\, x_2 \geq 0
  },
\end{equation}
depicted as a dashed line in \Cref{fig:sideways-cone}.

Let $\xbar=\limray{\ee_2}\plusl\ee_1$.
Then $\xbar\in\Kbar$ since,
as shown in the figure,
for each $t$, the point
$\xx_t=\transKern{[1,\,t,\,1/(2t)]}\in K$,
so
$\xbar=\lim \xx_t$ is in~$\Kbar$.
Likewise, $\xbar\in\Lbar$ since,
for each $t$, the point
$\xx'_t=\transKern{[1,\,t,\,0]}\in L$,
so
$\xbar=\lim \xx'_t$ is in~$\Lbar$.
On the other hand,
$\xbar\not\in \KLbar$ since if $\seq{\zz_t}$ is any sequence
in $K\cap L$ with a limit $\zbar\in\extspac{3}$, then $\zz_t\cdot\ee_1=0$ for
all $t$, so $\zbar\cdot\ee_1=\lim(\zz_t\cdot\ee_1)=0$
(by \Cref{thm:i:1}\ref{thm:i:1c}),
whereas $\xbar\cdot\ee_1=1$.
Therefore, by
\Cref{pr:inds-ext},
$\fext(\xbar)=\gext(\xbar)=0$,
but $\hext(\xbar)=+\infty$.%
\indexg{Sideways cone!extension of sum and|)}%
\end{example}

With these counterexamples in mind, we proceed to give sufficient
conditions for when in fact
$\hext(\xbar)=\fext(\xbar)+\gext(\xbar)$.
The next proposition provides such conditions based on
extensible continuity (see \Cref{dfn:extens-cont})
and lower semicontinuity.

\begin{proposition}  \label{pr:ext-sum-fcns}
  Let $f:\Rn\rightarrow\Rext$ and $g:\Rn\rightarrow\Rext$,
  and assume $f$ and $g$ are summable.
  Let $h=f+g$,
  let $\xbar\in\extspace$, and assume
  $\fext(\xbar)$ and $\gext(\xbar)$ are summable.
  Then the following hold:
  \begin{letter-compact}
  \item    \label{pr:ext-sum-fcns:a}
    $\hext(\xbar)\geq\fext(\xbar)+\gext(\xbar)$.
  \item    \label{pr:ext-sum-fcns:b}
    If either $f$ or $g$ is extensibly continuous at $\xbar$,
    then $\hext(\xbar)=\fext(\xbar)+\gext(\xbar)$.
  \item    \label{pr:ext-sum-fcns:c}
    If both $f$ and $g$ are extensibly continuous at $\xbar$,
    then $h$ is as well.
  \item    \label{pr:ext-sum-fcns:d}
    If $\xbar=\xx\in\Rn$ and both $f$ and $g$ are lower semicontinuous at $\xx$, then
    $h$ is as well, and
    $\hext(\xx)=\fext(\xx)+\gext(\xx)$.
  \end{letter-compact}
\end{proposition}

\begin{proof}
~
\begin{proof-parts}

\pfpart{Part~(\ref{pr:ext-sum-fcns:a}):}
By \Cref{pr:d1}, there exists a sequence $\seq{\xx_t}$ in $\Rn$ converging to $\xbar$ such that
$h(\xx_t)\to\eh(\xbar)$.
Then
\begin{align*}
  \fext(\xbar) + \gext(\xbar)
  &\le
  \liminf f(\xx_t) \plusd \liminf g(\xx_t)
\\
  &\le
  \liminf \bigBracks{f(\xx_t) + g(\xx_t)}
\\
  &=
  \liminf h(\xx_t)
  =
  \eh(\xbar).
\end{align*}
The first inequality follows by monotonicity of
downward addition (\Cref{pr:plusd-props}\ref{pr:plusd-props:f}),
since $\fext(\xbar)\le\liminf f(\xx_t)$
and similarly for~$\eg(\xbar)$
(by Eq.~\ref{eq:e:7}).
The second inequality follows by superadditivity of $\liminf$ (\Cref{prop:lim:eR}\ref{i:liminf:eR:sum}), and the final equality is by our choice of $\seq{\xx_t}$.

\pfpart{Part~(\ref{pr:ext-sum-fcns:b}):}
Without loss of generality, assume $f$ is extensibly continuous at
$\xbar$ (swapping $f$ and $g$ otherwise).
By \Cref{pr:d1}, there exists a
sequence $\seq{\xx_t}$ in $\Rn$ that converges to $\xbar$ and for
which $g(\xx_t)\rightarrow\gext(\xbar)$.
Also, since $f$ is extensibly continuous at~$\xbar$,
$f(\xx_t)\rightarrow\fext(\xbar)$.
Thus,
$h(\xx_t) = f(\xx_t) + g(\xx_t)\rightarrow \fext(\xbar) + \gext(\xbar)$
(\Cref{prop:lim:eR}\ref{i:lim:eR:sum}),
so, by \eqref{eq:e:7},
$\hext(\xbar)\leq \lim h(\xx_t) = \fext(\xbar) + \gext(\xbar)$.
Combined with
part~(\ref{pr:ext-sum-fcns:a}), this proves the claim.

\pfpart{Part~(\ref{pr:ext-sum-fcns:c}):}
Suppose $f$ and $g$ are extensibly continuous at $\xbar$,
and let $\seq{\xx_t}$ be any sequence in $\Rn$ converging to $\xbar$.
Then $f(\xx_t)\rightarrow\fext(\xbar)$ and
$g(\xx_t)\rightarrow\gext(\xbar)$, so
\[
  h(\xx_t)
  =
  f(\xx_t) + g(\xx_t)
  \rightarrow
  \fext(\xbar) + \gext(\xbar)
  =
  \hext(\xbar)
\]
with the last equality following from
part~(\ref{pr:ext-sum-fcns:b}).

\pfpart{Part~(\ref{pr:ext-sum-fcns:d}):}
Suppose $\xbar=\xx\in\Rn$ and $f$ and $g$ are lower semicontinuous at $\xx$.
Let $\seq{\xx_t}$ be any sequence in $\Rn$ converging to $\xx$. Then
\begin{align*}
  h(\xx)=
  f(\xx)+g(\xx)
  &\le
  \liminf f(\xx_t)
  \plusd
  \liminf g(\xx_t)
\\
  &\le
  \liminf \bigBracks{f(\xx_t) + g(\xx_t)}
  =
  \liminf h(\xx_t),
\end{align*}
where the first inequality is by lower semicontinuity of $f$ and $g$
at $\xx$ and monotonicity of downward addition
(\Cref{pr:plusd-props}\ref{pr:plusd-props:f}),
and the second inequality is by \Cref{prop:lim:eR}(\ref{i:liminf:eR:sum}). Thus, $h$ is lower semicontinuous at $\xx$ as claimed. Since $\ef$, $\eg$, $\eh$ are lower semicontinuous at $\xx$, we have, by \Cref{pr:h:1}(\ref{pr:h:1a}),
\[
 \eh(\xx)=h(\xx)=f(\xx)+g(\xx)=\ef(\xx)+\eg(\xx).
\qedhere
\]
\end{proof-parts}
\end{proof}

In particular, \Cref{pr:ext-sum-fcns}(\ref{pr:ext-sum-fcns:d}) shows that if $f$ and $g$ are lower semicontinuous
then $\eh$ can differ from $\ef+\eg$ only at infinite points, as was the case in \Cref{ex:KLsets:extsum-not-sum-exts}.

In \Cref{pr:ext-sum-fcns}, we made continuity and lower semicontinuity assumptions, but did
not assume convexity.
We will next give different sufficient conditions for when the extension of a sum is the sum of the extensions,
based on convexity and the assumption that the relative interiors of the
domains of the functions being added have at least one point in
common.
We first prove a lemma for when all of the functions are nonnegative,
and later, in \Cref{thm:ext-sum-fcns-w-duality},
give a more general result.
The proof is based on conjugacy.

\begin{lemma}  \label{lem:ext-sum-nonneg-fcns-w-duality}
  For $i=1,\dotsc,m$,
  let $f_i:\Rn\rightarrow\Rext$ be convex
  with $f_i\geq 0$.
  Assume
  $ \bigcap_{i=1}^m \ri(\dom{f_i}) \neq \emptyset $.
  Let $h=f_1+\dotsb+f_m$, and let $\xbar\in\extspace$.
  Then
  \[ \hext(\xbar)=\efsub{1}(\xbar)+\dotsb+\efsub{m}(\xbar). \]
\end{lemma}

\begin{proof}
It suffices to prove
\begin{equation}  \label{eqn:thm:ext-sum-fcns-w-duality:1}
  \hext(\xbar)
  \leq
  \sum_{i=1}^m \efsub{i}(\xbar)
\end{equation}
since the reverse inequality follows from
\Cref{pr:ext-sum-fcns}(\ref{pr:ext-sum-fcns:a}),
applied inductively.
For $i=1,\dotsc,m$,
we assume henceforth that
$\efsub{i}(\xbar)<+\infty$, since otherwise
\eqref{eqn:thm:ext-sum-fcns-w-duality:1} is immediate.
This implies $f_i\not\equiv+\infty$ and so that
$f_i$ is proper.

Let $\uu\in\Rn$.
We will prove
\begin{equation}  \label{eqn:thm:ext-sum-fcns-w-duality:2}
  {-\hstar(\uu)}\plusd \xbar\cdot\uu
  \leq
  \sum_{i=1}^m \efsub{i}(\xbar).
\end{equation}
Once proved, this will imply
\eqref{eqn:thm:ext-sum-fcns-w-duality:1} since then
\[
  \hext(\xbar)
  =
  \hdub(\xbar)
  =
  \sup_{\uu\in\Rn} \bigBracks{ -\hstar(\uu)\plusd \xbar\cdot\uu }
  \leq
  \sum_{i=1}^m \efsub{i}(\xbar),
\]
with the equalities following respectively from
\Cref{cor:all-red-closed-sp-cases}
(since $h\geq 0$)
and \Cref{thm:fext-dub-sum}(\ref{thm:fext-dub-sum:b}).

We assume $\hstar(\uu)<+\infty$ since otherwise
\eqref{eqn:thm:ext-sum-fcns-w-duality:2}
is immediate.
By a standard characterization of the conjugate of a sum of proper
convex functions (\Cref{roc:thm16.4}) and since
$ \bigcap_{i=1}^m \ri(\dom{f_i}) \neq \emptyset $,
there exist $\uu_1,\dotsc,\uu_m\in\Rn$ with
\begin{equation}  \label{eq:lem:ext-sum-nonneg-fcns-w-duality:L1}
  \sum_{i=1}^m \uu_i = \uu,
  \quad\text{ and }\quad
  \sum_{i=1}^m \fistar(\uu_i) = \hstar(\uu).
\end{equation}
Since $\hstar(\uu)<+\infty$, it follows that
each $\fistar(\uu_i)<+\infty$, and so actually
$\fistar(\uu_i)\in\R$,
since each $f_i$ is proper, so each $\fistar$ is as well
(by
\Cref{pr:conj-props-cvx}\ref{pr:conj-props-cvx:a}).

By \Cref{thm:fext-dub-sum}(\ref{thm:fext-dub-sum:a},\ref{thm:fext-dub-sum:b}),
\begin{equation}  \label{eqn:thm:ext-sum-fcns-w-duality:3}
  \efsub{i}(\xbar)
  \geq
  \fidub(\xbar)
  \geq
  -\fistar(\uu_i)\plusd \xbar\cdot\uu_i,
\end{equation}
for $i=1,\dotsc,m$.
Since $\fistar(\uu_i)\in\R$ and since we assumed
$\efsub{i}(\xbar) < +\infty$,
it follows that $\xbar\cdot\uu_i<+\infty$.

Thus, $\xbar\cdot\uu_1,\dotsc,\xbar\cdot\uu_m$ are summable.
By \Cref{pr:i:1}, applied repeatedly, this implies
\begin{equation}  \label{eq:lem:ext-sum-nonneg-fcns-w-duality:L2}
  \xbar\cdot\uu
  =
  \xbar\cdot\BiggParens{\sum_{i=1}^m \uu_i}
  =
  \sum_{i=1}^m \xbar\cdot\uu_i.
\end{equation}

Combining now yields
\begin{align*}
  -\hstar(\uu)\plusd\xbar\cdot\uu
  &=
  -\BiggParens{\sum_{i=1}^m \fistar(\uu_i)}
  +
  \sum_{i=1}^m \xbar\cdot\uu_i
  \\
  &=
  \sum_{i=1}^m \bigBracks{ -\fistar(\uu_i) + \xbar\cdot\uu_i }
  \leq
  \sum_{i=1}^m \efsub{i}(\xbar).
\end{align*}
The first equality is by
Eqs.~(\ref{eq:lem:ext-sum-nonneg-fcns-w-duality:L1})
and~(\ref{eq:lem:ext-sum-nonneg-fcns-w-duality:L2}),
with summability ensured since $\fistar(\uu_i)\in\R$
and $\xbar\cdot\uu_i<+\infty$ for $i=1,\dotsc,m$.
The inequality follows from
\eqref{eqn:thm:ext-sum-fcns-w-duality:3}.
This proves
\eqref{eqn:thm:ext-sum-fcns-w-duality:2},
completing the lemma.
\end{proof}

\indexg{simultaneous convergence to multiple extensions|(}%
\indexg{lower semicontinuous extension!simultaneous convergence to multiple|(}%
In \Cref{pr:d1}, we saw that for any function
$f:\Rn\rightarrow\Rext$ and point $\xbar\in\extspace$,
there must exist a sequence $\seq{\xx_t}$ in $\Rn$ that converges to
$\xbar$ and for which $f(\xx_t)\rightarrow\fext(\xbar)$.
Using
\Cref{lem:ext-sum-nonneg-fcns-w-duality},
we show next that,
for any finite collection of convex functions
$f_1,\dotsc,f_m$ with effective domains whose relative interiors
overlap,
there must exist a \emph{single} sequence $\seq{\xx_t}$ converging to
$\xbar$ for which this property holds simultaneously
for each of the functions, so that
$f_i(\xx_t)\rightarrow\efsub{i}(\xbar)$ for $i=1,\dotsc,m$.
This will also help in proving the more general form of
\Cref{lem:ext-sum-nonneg-fcns-w-duality}
that we present in
\Cref{thm:ext-sum-fcns-w-duality}.

\begin{theorem}  \label{thm:seq-for-multi-ext}
  Let $f_i:\Rn\rightarrow\Rext$ be convex, for $i=1,\dotsc,m$,
  and let $\xbar\in\extspace$.
  Assume
  $\bigcap_{i=1}^m \ri(\dom{f_i}) \neq \emptyset$.
  Then there exists a sequence $\seq{\xx_t}$ in $\Rn$ such that
  $\xx_t\rightarrow\xbar$
  and
  $f_i(\xx_t)\rightarrow\efsub{i}(\xbar)$
  for $i=1,\dotsc,m$.
\end{theorem}

\begin{proof}
First, by possibly rearranging the indices of the functions, we can
assume without loss of generality that
$\efsub{i}(\xbar)<+\infty$ for $i=1,\dotsc,\ell$
and
$\efsub{i}(\xbar)=+\infty$ for $i=\ell+1,\dotsc,m$,
for some $\ell\in\{0,\dotsc,m\}$.
Let $f'_i=\expex\circ f_i$, for $i=1,\dotsc,m$
(where $\expex$ is as in Eq.~\ref{eq:expex-defn}),
and let $h=f'_1+\dotsb + f'_\ell$.

Then for $i=1,\dotsc,\ell$,
each $f'_i$ is
convex and nonnegative
(by \Cref{pr:j:2}\ref{pr:j:2a}).
Also, $\dom{f'_i}=\dom{f_i}$, so
$\bigcap_{i=1}^\ell \ri(\dom{f'_i}) \neq \emptyset$.
Therefore, by
\Cref{lem:ext-sum-nonneg-fcns-w-duality},
\begin{equation}  \label{eqn:thm:seq-for-multi-ext:1}
  \hext(\xbar)
  =
  \sum_{i=1}^\ell \fpextsub{i}(\xbar).
\end{equation}

By \Cref{pr:d1}, there exists a sequence $\seq{\xx_t}$ in
$\Rn$ converging to $\xbar$ with $h(\xx_t)\rightarrow\hext(\xbar)$.
For $i=1,\dotsc,\ell$, we consider each sequence $\seq{f'_i(\xx_t)}$,
one by one.
By sequential compactness of $\eR$, this sequence must have a
convergent subsequence; discarding all other elements
then yields one on which the entire sequence
converges.
Repeating for each $i$,
we obtain a sequence $\seq{\xx_t}$ that converges to
$\xbar$ and for which $h(\xx_t)\rightarrow\hext(\xbar)$ and
$\seq{f'_i(\xx_t)}$ converges in $\eR$ for all $i=1,\dotsc,\ell$.

We then have
\begin{equation}  \label{eq:seq-multi-ext}
  \sum_{i=1}^\ell
  \fpextsub{i}(\xbar)
  =
  \hext(\xbar)
  =
  \lim h(\xx_t)
  =\lim\Bracks{\sum_{i=1}^\ell f'_i(\xx_t)}
  =\sum_{i=1}^\ell\BigBracks{\lim f'_i(\xx_t)},
\end{equation}
where the last equality is
by continuity (\Cref{prop:lim:eR}\ref{i:lim:eR:sum}), noting that
$f'_i(\xx_t)\ge 0$, so also
$\lim f'_i(\xx_t)\ge 0$, implying that the limits are summable.

Furthermore, for $i=1,\dotsc,\ell$,
$\lim f'_i(\xx_t)\ge\fpextsub{i}(\xbar)$
(by definition of extensions, Eq.~\ref{eq:e:7}) and also
$\fpextsub{i}(\xbar)\in\R$.
Combined with
\eqref{eq:seq-multi-ext},
these facts imply that actually
$\lim f'_i(\xx_t)=\fpextsub{i}(\xbar)$
for all $i=1,\dotsc,\ell$,
since otherwise the rightmost expression in
\eqref{eq:seq-multi-ext} would be strictly greater than
the leftmost expression in that equality.

By \Cref{pr:j:2}(\ref{pr:j:2lim}), this implies that
$
  \lim f_i(\xx_t)
  =\efsub{i}(\xbar)
$
for $i=1,\dotsc,\ell$.
Also, for $i=\ell+1,\dotsc,m$,
we have $\liminf f_i(\xx_t)\geq\efsub{i}(\xbar)=+\infty$
(by Eq.~\ref{eq:e:7}), implying
$f_i(\xx_t)\rightarrow+\infty$.
Thus,
$f_i(\xx_t) \rightarrow  \efsub{i}(\xbar) $
for $i=1,\dotsc,m$.%
\indexg{simultaneous convergence to multiple extensions|)}%
\indexg{lower semicontinuous extension!simultaneous convergence to multiple|)}%
\end{proof}

We can now prove the more general form of
\Cref{lem:ext-sum-nonneg-fcns-w-duality}:

\begin{theorem}  \label{thm:ext-sum-fcns-w-duality}
  Let $f_i:\Rn\rightarrow\Rext$ be convex,
  for $i=1,\dotsc,m$.
  Assume $f_1,\dotsc,f_m$ are summable,
  and that
  $\bigcap_{i=1}^m \ri(\dom{f_i}) \neq \emptyset$.
  Let $h=f_1+\dotsb+f_m$, and let $\xbar\in\extspace$.
  Assume also that
  $\efsub{1}(\xbar),\dotsc,\efsub{m}(\xbar)$ are summable.
  Then
  \[ \hext(\xbar)=\efsub{1}(\xbar)+\dotsb+\efsub{m}(\xbar). \]
\end{theorem}

\begin{proof}
By \Cref{thm:seq-for-multi-ext},
there exists a sequence $\seq{\xx_t}$ in $\Rn$
with $\xx_t\rightarrow\xbar$ and
$f_i(\xx_t)\rightarrow\efsub{i}(\xbar)$
for $i=1,\dotsc,m$.
For this sequence,
\[
  \hext(\xbar)
  \ge
  \efsub{1}(\xbar)+\dotsb+\efsub{m}(\xbar)
  =
  \lim\bigBracks{
    f_1(\xx_t)+\dotsb+f_m(\xx_t)
  }
  =
  \lim h(\xx_t)
  \ge
  \hext(\xbar).
\]
The first inequality is
by repeated application of
\Cref{pr:ext-sum-fcns}(\ref{pr:ext-sum-fcns:a}),
and the first equality is by \Cref{prop:lim:eR}(\ref{i:lim:eR:sum}).
The last inequality is by definition of extensions
(Eq.~\ref{eq:e:7}).
This proves the theorem.%
\indexg{sum of functions (standard)!extension of|)}%
\indexg{lower semicontinuous extension!sum of functions@of sum of functions|)}%
\end{proof}

\indexg{closure, astral!intersection of sets@of intersection of sets|(}%
When applied to indicator functions,
\Cref{thm:ext-sum-fcns-w-duality}
(or \Cref{lem:ext-sum-nonneg-fcns-w-duality}) implies that for any
finite collection of convex sets in $\Rn$ with overlapping relative
interiors,
the intersection of their (astral) closures is equal to the closure
of their intersection:

\begin{corollary}  \label{cor:clos-int-sets}
  Let $S_1,\dotsc,S_m\subseteq\Rn$ be convex,
  and let $U=\bigcap_{i=1}^m S_i$.
  Assume
  $\bigcap_{i=1}^m \ri{S_i}\neq\emptyset$.
  Then
  \[
     \Ubar
     =
     \bigcap_{i=1}^m \Sbar_i.
  \]
\end{corollary}

\begin{proof}
Let $f_i=\indf{S_i}$
be the indicator function of the set $S_i$, for $i=1,\dotsc,m$.
These functions are convex and nonnegative.
Let $h= f_1 + \dotsb + f_m$, implying
$h=\indf{U}$.

Since
$\bigcap_{i=1}^m \ri(\dom{f_i})= \bigcap_{i=1}^m (\ri{S_i})\neq\emptyset$,
we can apply
\Cref{lem:ext-sum-nonneg-fcns-w-duality}
yielding
$\hext= \efsub{1} + \dotsb + \efsub{m}$.
We then have
\[
  \indfa{\Ubar}=\hext = \efsub{1} + \dotsb + \efsub{m}
           =\indfa{\Sbar_1}+\dotsb+\indfa{\Sbar_m}
           =\indfa{\cap_{i=1}^m\Sbar_i},
\]
where the first and third equalities are by
\Cref{pr:inds-ext}.
This proves the claim.
\end{proof}

Without assuming overlapping relative interiors,
\Cref{cor:clos-int-sets} need not hold in general.
For instance, the sets $P=\R_{>0}$ and $N=\R_{<0}$ from \Cref{ex:pos-neg} do not intersect, but their astral closures do, so $\PNbar\ne\Pbar\cap\Nbar$.
Similarly, in \Cref{ex:KLsets:extsum-not-sum-exts}, we constructed sets $K,L\subseteq\Rn$ and
a point $\xbar\in\eRn$, such that $\xbar\in\Kbar\cap\Lbar$, but $\xbar\not\in\KLbar$, so $\Kbar\cap\Lbar\ne\KLbar$.%
\indexg{closure, astral!intersection of sets@of intersection of sets|)}

\section{Composition with an affine map}

\indexg{linear map in composition with standard function!extension of|(}%
\indexg{affine map in composition with standard function!extension of|(}%
\indexg{lower semicontinuous extension!composition with affine map@of composition with affine map|(}%
We next consider functions $h:\Rn\rightarrow\Rext$ of the form
$h(\xx)=f(\bb+\A\xx)$, where $f:\Rm\rightarrow\Rext$,
$\A\in\Rmn$ and $\bb\in\Rm$, so that $h$ is the composition of $f$
with an affine map.
In this case, we might expect $h$'s extension to be
$\hext(\xbar)=\fext(\bb\plusl\A\xbar)$.
As before, this will be the case under some fairly mild conditions.
First, however, we show that this need not hold in
general:

\begin{example}[Positives and negatives, continued]
\indexg{Positives and negatives|(}%
Let $f=i_P$ where $P=\R_{>0}$ as in \Cref{ex:pos-neg},
let $A=[0]$ (the $1\times 1$ zero matrix),
and let $h(x)=f(A x)$ for $x\in\R$. Thus, $h(x)=f(0)=+\infty$ for all $x\in\R$, so also $\eh(\ex)=+\infty$ for all $\ex\in\eR$. However, as we saw in \Cref{ex:pos-neg}, $\ef=\indfa{\Pbar}$ where $\Pbar=[0,+\infty]$, so $\ef(A\ex)=\ef(0)=0$ for all~$\ex\in\Rext$.%
\indexg{Positives and negatives|)}%
\end{example}

The problem here is that $f$ is not lower semicontinuous at $0$.
The next example shows that lower semicontinuity of $f$ is not sufficient to guarantee that $\eh(\xbar)=\ef(\bb\plusl\A\xbar)$ for all $\xbar$.

\begin{example}[Sideways cone, continued]
\label{ex:fA-ext-countereg}
\indexg{Sideways cone!extension of composition with linear map and|(}%
Let $f=\indK$ where $K\subseteq\R^3$ is the sideways cone from \Cref{ex:KLsets:extsum-not-sum-exts},
let
\[
   \A
   =
   \begin{bmatrix}
        1 & 0 & 0 \\
        0 & 1 & 0 \\
        0 & 0 & 0
   \end{bmatrix},
\]
and let $h(\xx)=f(\A\xx)$ for $\xx\in\R^3$.
Let $\xbar=\limray{\ee_2}\plusl\ee_1$,
where
$\ee_1,\ee_2,\ee_3$ are the standard basis vectors in $\R^3$.
Then $\A\xbar=\limray{\ee_2}\plusl\ee_1$
(using \Cref{pr:h:4}).
In Example~\ref{ex:KLsets:extsum-not-sum-exts},
we argued that $\fext(\A\xbar)=0$.

For $\xx\in\R^3$,
$h(\xx) = h(x_1,x_2,x_3)  =  f(x_1,x_2,0)$,
which is an indicator function~$\indf{M}$ for the set
$M=\{\xx\in\R^3 :\: x_1=0,\,x_2\geq 0\}$.
Therefore,
$\hext=\indfa{\Mbar}$
by \Cref{pr:inds-ext},
so
if $\hext(\xbar)=0$ then there exists a sequence
$\seq{\xx_t}$ in $M$ that converges to $\xbar$,
implying
$0=\xx_t\cdot\ee_1\rightarrow\xbar\cdot\ee_1$,
a contradiction since $\xbar\cdot\ee_1=1$.
Thus, $\hext(\xbar)=+\infty$ but
$\fext(\A\xbar)=0$.%
\indexg{Sideways cone!extension of composition with linear map and|)}%
\end{example}

Here are some conditions
for when such extensions have the expected form:

\begin{proposition}  \label{pr:ext-affine-comp}
  Let $f:\Rm\rightarrow\Rext$,
  $\A\in\Rmn$, and $\bb\in\Rm$.
  Let $h:\Rn\rightarrow\Rext$ be defined by
  $h(\xx) = f(\bb + \A\xx)$
  for $\xx\in\Rn$.
  Let $\xbar\in\extspace$.
  Then the following hold:
  \begin{letter-compact}
  \item  \label{pr:ext-affine-comp:a}
    $\hext(\xbar)\geq\fext(\bb\plusl\A\xbar)$.
  \item  \label{pr:ext-affine-comp:b}
    If $m=n$ and $\A$ is invertible then
    $\hext(\xbar)=\fext(\bb\plusl\A\xbar)$.
  \item  \label{pr:ext-affine-comp:c}
    If $f$ is extensibly continuous at $\bb\plusl\A\xbar$
    then
    $\hext(\xbar)=\fext(\bb\plusl\A\xbar)$
    and
    $h$ is extensibly continuous at $\xbar$.
  \item    \label{pr:ext-affine-comp:d}
    If $\xbar=\xx\in\Rn$ and $f$ is lower semicontinuous at $\bb+\A\xx$
    then
    $\hext(\xx)=\fext(\bb+\A\xx)$
    and
    $h$ is lower semicontinuous at $\xx$.
  \end{letter-compact}
\end{proposition}

\begin{proof}
~
\begin{proof-parts}

\pfpart{Part~(\ref{pr:ext-affine-comp:a}):}
Suppose $\seq{\xx_t}$ is a sequence in $\Rn$ that converges to
$\xbar$.
Then $\bb+\A\xx_t=\bb\plusl\A\xx_t\rightarrow \bb\plusl\A\xbar$
by continuity (\Cref{cor:aff-cont}).
Therefore, 
\[
  \liminf h(\xx_t)
  =
  \liminf f(\bb+\A\xx_t)
  \geq
  \fext(\bb\plusl\A\xbar),
\]
with the inequality by definition of extension
(Eq.~\ref{eq:e:7}).
Since this holds for all such sequences, the claim follows
(again by Eq.~\ref{eq:e:7}).

\pfpart{Part~(\ref{pr:ext-affine-comp:b}):}
Suppose $m=n$ and $\A$ is invertible.
Then by part~(\ref{pr:ext-affine-comp:a}),
$\hext(\xbar)\geq\fext(\bb\plusl\A\xbar)$.

For the reverse inequality,
note that, for $\xx,\yy\in\Rn$,
$\yy=\bb+\A\xx$ if and only if
$\xx=\bb'+\Ainv\yy$,
where $\bb' = -\Ainv\bb$.
Therefore,
$f(\yy)=h(\bb'+\Ainv\yy)$
for all $\yy\in\Rn$.
Applying part~(\ref{pr:ext-affine-comp:a}) with $f$ and $h$ swapped,
it follows that
$\fext(\ybar)\geq\hext(\bb'\plusl\Ainv\ybar)$
for $\ybar\in\extspace$.
In particular, setting $\ybar=\bb\plusl\A\xbar$, this yields
$\fext(\bb\plusl\A\xbar)\geq\hext(\xbar)$.

\pfpart{Part~(\ref{pr:ext-affine-comp:c}):}
Suppose $f$ is extensibly continuous at $\bb\plusl\A\xbar$.
Let $\seq{\xx_t}$ be any sequence in~$\Rn$ converging to
$\xbar$, implying
$\bb+\A\xx_t\rightarrow \bb\plusl\A\xbar$
by continuity (\Cref{cor:aff-cont}).
Since $f$ is extensibly continuous at $\bb\plusl\A\xbar$,
it follows that
$h(\xx_t)=f(\bb+\A\xx_t)\rightarrow\fext(\bb\plusl\A\xbar)$.

Thus, $\lim h(\xx_t)$ exists and has the same value
$\fext(\bb\plusl\A\xbar)$ for every such sequence.
Furthermore, this limit must be equal to $\hext(\xbar)$ for at least
one such sequence, by \Cref{pr:d1}.
This proves both parts of the claim.

\pfpart{Part~(\ref{pr:ext-affine-comp:d}):}
Suppose $\xbar=\xx\in\Rn$ and $f$ is lower semicontinuous at $\bb+\A\xx$. Let $\seq{\xx_t}$ be any sequence in $\Rn$ converging to $\xx$. Then $\bb+\A\xx_t\to\bb+\A\xx$ by continuity, so
\[
  h(\xx)
  =
  f(\bb+\A\xx)
  \le
  \liminf f(\bb+\A\xx_t)
  =
  \liminf h(\xx_t)
\]
by lower semicontinuity of $f$. Thus, $h$ is lower semicontinuous at $\xx$. Therefore,
by \Cref{pr:h:1}(\ref{pr:h:1a}),
\[
 \eh(\xx)
 =
 h(\xx)
 =
 f(\bb+\A\xx)
 =
 \ef(\bb+\A\xx).%
\indexg{affine map in composition with standard function!extension of|)}%
\qedhere
\]
\end{proof-parts}
\end{proof}

Henceforth, we consider only composition with a linear map, that is,
with $\bb=\zero$.
For a matrix $\A\in\Rmn$ and a function $f:\Rm\rightarrow\Rext$
or $F:\extspace\rightarrow\Rext$,
we use the notation $\fA$ or $\FA$ for such a composition,
as was defined in Eqs.~(\ref{eq:fA-defn})
and~(\ref{eq:FA-dfn}).

We show next that if $f$ is convex and if the column space of the
matrix $\A$ intersects the relative interior of $\dom f$, then
$\fAext(\xbar)=\fext(\A\xbar)$; that is, $\fAextnp=\fext \A$.

\begin{theorem}  \label{thm:ext-linear-comp}
  Let $f:\Rm\rightarrow\Rext$ be convex
  and
  let $\A\in\Rmn$.
  Assume there exists $\xing\in\Rn$ such that
  $\A\xing\in\ri(\dom f)$.
  Then  $\fAextnp = \fext \A$; that is,
  for $\xbar\in\extspace$,
  \[ \fAext(\xbar)=\fext(\A\xbar). \]
\end{theorem}

\begin{proof}
Let $\xbar\in\eRn$, and let $h=\fA$.
It suffices to prove $\hext(\xbar)\leq\fext(\A\xbar)$
since the reverse inequality was proved in
\Cref{pr:ext-affine-comp}(\ref{pr:ext-affine-comp:a}).

Let us assume for now that $f\geq 0$,
later returning to the general case.
Then $h\geq 0$, so both $f$ and $h$ are reduction-closed
(\Cref{pr:j:1}\ref{pr:j:1c}). 
Furthermore,
both $f$ and $h$ are proper
since $h(\xing)=f(\A\xing)<+\infty$.

Let $\uu\in\Rm$.
We will show
$-\hstar(\uu)\plusd\xbar\cdot\uu \leq \fext(\A\xbar)$,
which we will see is sufficient to prove the main claim.
This is immediate if $\hstar(\uu)=+\infty$, so we assume henceforth
that $\hstar(\uu)<+\infty$. Since $h$ is proper, so is $\hstar$
(by
\Cref{pr:conj-props-cvx}\ref{pr:conj-props-cvx:a});
hence, $\hstar(\uu)\in\R$.

Using the characterization of the conjugate $(\fA)^*$
(\Cref{roc:thm16.3:fA})
and the fact that $\hstar(\uu)\in\R$, we obtain that $\hstar(\uu)=(\fA)^*(\uu)=\fstar(\ww)$ for some $\ww\in\Rm$ such that $\transA\ww=\uu$.
Thus,
\begin{align*}
  -\hstar(\uu)\plusd\xbar\cdot\uu
  &=
  -\fstar(\ww)\plusd\xbar\cdot(\transA\ww)
\\
  &=
  -\fstar(\ww)\plusd(\A\xbar)\cdot\ww
  \leq
  \fdub(\A\xbar)
  \leq
  \fext(\A\xbar).
\end{align*}
The second equality is by
\Cref{thm:Ax-dot-u},
and the inequalities are by
\Cref{thm:fext-dub-sum}(\ref{thm:fext-dub-sum:a}\ref{thm:fext-dub-sum:b}).

Taking supremum over $\uu\in\Rm$ then yields
\begin{equation}
\label{eq:ext-linear-comp:final}
  \hext(\xbar)
  =
  \hdub(\xbar)
  =
  \sup_{\uu\in\Rm}\bigBracks{-\hstar(\uu)\plusd\xbar\cdot\uu}
  \leq
  \fext(\A\xbar),
\end{equation}
with the equalities following from
\Cref{cor:all-red-closed-sp-cases}
and
\Cref{thm:fext-dub-sum}(\ref{thm:fext-dub-sum:b}).

This proves the claim when $f\geq 0$.
For the general case, we apply the exponential composition trick
(\Cref{sec:exp-comp}).
Let $f$ be as stated in the theorem (not necessarily nonnegative).
Let $f'=\expex\circ f$ and $h'=\expex\circ h$.
Then for $\xx\in\Rn$,
$h'(\xx)=\expex(h(\xx))=\expex(f(\A\xx))=f'(\A\xx)$.
Also, $\dom f'=\dom f$ so if $\A\xing\in\ri(\dom f)$ then also
$\A\xing\in\ri(\dom f')$.

Thus, for $\xbar\in\extspac{n}$,
\[
  \expex\bigParens{\hext(\xbar)}
  =
  \hpext(\xbar)
  \leq
  \fpext(\A\xbar)
  =
  \expex\bigParens{\fext(\A\xbar)}.
\]
The inequality follows from the argument above (namely, Eq.~\ref{eq:ext-linear-comp:final}) applied to $f'$
(which is nonnegative).
The two equalities are
from \Cref{pr:j:2}(\ref{pr:j:2c}).
Since $\expex$ is strictly increasing, it follows that
$\hext(\xbar)\leq\fext(\A\xbar)$ as claimed.%
\indexg{linear map in composition with standard function!extension of|)}%
\indexg{lower semicontinuous extension!composition with affine map@of composition with affine map|)}%
\end{proof}

\indexg{closure, astral!inverse image under linear map@of inverse image under linear map|(}%
\indexg{linear maps (standard)!astral closure of inverse image under|(}%
Similar to \Cref{cor:clos-int-sets},
we can apply this theorem to an indicator function, yielding
the next corollary which shows that, under suitable conditions,
the astral closure of the inverse
image of a convex set $S\subseteq\Rm$ under a linear map
is the same as the inverse image of the
corresponding astral closure $\Sbar$ under the corresponding
astral linear map.

\begin{corollary}  \label{cor:lin-map-inv}
  Let $\A\in\Rmn$ and let
  $\slinmapA:\Rn\rightarrow\Rm$ and $\alinmapA:\extspace\rightarrow\extspac{m}$
  be the associated standard and
  astral linear maps
  (so that $\slinmapA(\xx)=\A\xx$ for $\xx\in\Rn$
  and $\alinmapA(\xbar)=\A\xbar$ for $\xbar\in\extspace$).
  Let $S\subseteq\Rm$ be convex, and assume there exists some point
  $\xing\in\Rn$ with $\A\xing\in\ri{S}$.
  Then
  \[ \clbar{\slinmapAinv(S)}=\alinmapAinv(\Sbar). \]
\end{corollary}

\begin{proof}
Let $f=\indf{S}$ be $S$'s indicator function, which is convex,
and let
$h=\fA$.
Then $h(\xx)=0$ if $\A\xx\in S$ (and otherwise is
$+\infty$);
that is, $h=\indf{U}$ where
\[
  U
  =
  \{ \xx\in\Rn :\: \A\xx\in S \}
  =
  \slinmapAinv(S).
\]

Since $f$ and $h$ are indicators of convex sets $S$ and $U$, their extensions are astral indicators of the closures of $S$ and $U$ (by \Cref{pr:inds-ext}),
that is,
$\fext=\indset_{\Sbar}$ and
$\hext=\indset_{\Ubar}$.

By assumption, there exists a point $\xing$ with
$\A\xing\in\ri{S}=\ri(\dom f)$, so we can apply
\Cref{thm:ext-linear-comp}
yielding $\hext(\xbar)=\fext(\A\xbar)$ for $\xbar\in\extspace$.
Since $\fext$ is the indicator of $\Sbar$ and $\hext$ is the indicator of $\Ubar$, the set $\Ubar$ satisfies
\[
  \Ubar
  =
  \{ \xbar\in\extspace :\: \A\xbar\in\Sbar \}
  =
  \alinmapAinv(\Sbar),
\]
completing the proof.%
\indexg{linear maps (standard)!astral closure of inverse image under|)}%
\indexg{closure, astral!inverse image under linear map@of inverse image under linear map|)}%
\end{proof}

\indexg{closure, astral!affine set@of affine set|(}%
\indexg{affine sets!astral closure of|(}%
As an additional corollary, we can compute the astral closure of any
affine set:

\begin{corollary}  \label{cor:closure-affine-set}
  Let $\A\in\Rmn$ and
  $\bb\in\Rm$,
  and let
  $M = \{ \xx\in\Rn :\: \A\xx = \bb \}$.
  Then
  \begin{equation}  \label{eq:cor:closure-affine-set:1}
    \Mbar
    =
    \Braces{ \xbar\in\extspace :\: \A\xbar = \bb }.
  \end{equation}
\end{corollary}

\begin{proof}
Let $U$ denote the set on the right-hand side of
\eqref{eq:cor:closure-affine-set:1}.

Suppose first that $M$ is empty, so that $\Mbar$ is empty as well.
To prove the corollary in this case,
we need to show that $U$ is empty.
Suppose to the contrary that
$\A\xbar=\bb$ for some $\xbar\in\extspace$.
We can write $\xbar=\ebar\plusl\qq$ for some $\ebar\in\corezn$ and $\qq\in\Rn$ (\Cref{thm:icon-fin-decomp}). Then
\begin{equation}  \label{eq:cor:closure-affine-set:2}
  \bb = \A(\ebar\plusl\qq)
      = (\A\ebar)\plusl(\A\qq),
\end{equation}
where the second equality is by \Cref{pr:h:4}(\ref{pr:h:4c}).
Since $\A\ebar$ is an icon (by \Cref{pr:i:8}\ref{pr:i:8-matprod}),
and since the icons for the expressions on the left- and right-hand
sides of
\eqref{eq:cor:closure-affine-set:2}
must be the same
(by \Cref{thm:icon-fin-decomp}),
it follows that $\A\ebar=\zero$.
Thus, $\A\qq=\bb$,
contradicting that $M=\emptyset$.

In the alternative case, $M$ includes some point $\xing\in\Rn$.
We apply
\Cref{cor:lin-map-inv}
with $S=\{\bb\}$, noting that
$\A\xing=\bb\in\{\bb\}=\ri S$
and that
$\Sbar=\{\bb\}$.
Thus,
in the notation of that corollary,
$M = \slinmapAinv(S)$
and
$U = \smash{\alinmapAinv(\Sbar)}$.
Therefore, $\Mbar=U$.%
\indexg{closure, astral!affine set@of affine set|)}%
\indexg{affine sets!astral closure of|)}%
\end{proof}

\section{Linear image of a function}
\label{sec:ext-lin-img-fcn}

The image of a function $f:\Rn\rightarrow\Rext$ under
matrix $\A\in\Rmn$, denoted $\A f$, was defined in
\eqref{eq:lin-image-fcn-defn}.
\indexg{linear image of astral function|(}%
Generalizing to astral space, we also defined in
\eqref{eqn:image-F-dfn} the image of an astral function
$F:\extspace\rightarrow\Rext$ under $\A$, denoted $\A F$, to be the
function given by
\begin{equation}  \label{eqn:image-F-dfn2}
  (\A F)(\xbar) = \inf\,\bigBraces{F(\zbar) :\: \zbar\in\extspace,\, \A\zbar=\xbar}
\end{equation}
for $\xbar\in\extspac{m}$.
In a moment, we will see that $\Afextnp=\A \fext$.
To prove this, we first give a general result showing that if $F$ is
lower semicontinuous then so is $\A F$;
furthermore, the infimum appearing in \eqref{eqn:image-F-dfn2}
is always attained (unless it is vacuous).

\begin{theorem}   \label{thm:AF-lsc}
  Let $\A\in\Rmn$, and let $F:\extspace\rightarrow\Rext$ be lower
  semicontinuous.
  Then:
  \begin{letter-compact}
  \item   \label{thm:AF-lsc:a}
    $\A F$ is lower semicontinuous.
  \item   \label{thm:AF-lsc:b}
    If $\xbar$ is in the astral column space of $\A$,
    then the infimum appearing in
    \eqref{eqn:image-F-dfn2}
    is attained by some $\zbar\in\extspace$ with $\A\zbar=\xbar$
    (and otherwise is vacuous).
  \end{letter-compact}
\end{theorem}

\begin{proof}
We prove the parts of the theorem in the opposite
order of how they were stated.

\begin{proof-parts}
\pfpart{Part~(\ref{thm:AF-lsc:b}):}
Suppose $\xbar\in\acolspace \A$.
Let $\gamma=(\A F)(\xbar)$.
If $\gamma=+\infty$, then $F(\zbar)=+\infty$ for all
$\zbar\in\extspace$ with $\A\zbar=\xbar$; therefore, the infimum
in \eqref{eqn:image-F-dfn2} is attained by any of these.

Otherwise, for all $t$, let
\[
  \beta_t
  =
  \begin{cases}
    \gamma + 1/t
    & \text{if $\gamma\in\R$,}
  \\
    -t
    & \text{if $\gamma=-\infty$.}
  \end{cases}
\]
For each $t$,
by $\A F$'s definition and since $\gamma<+\infty$,
there exists $\zbar_t\in\extspace$ with
$\A\zbar_t=\xbar$ and $F(\zbar_t)<\beta_t$.
By sequential compactness,
the resulting sequence $\seq{\zbar_t}$ must have a convergent subsequence.
By discarding all other elements, we can assume the entire sequence
converges to some point $\zbar\in\extspace$.
Note that $\xbar=\A\zbar$ since
$\xbar=\A\zbar_t\rightarrow\A\zbar$
(by \Cref{thm:linear:cont}\ref{thm:linear:cont:b}).

Thus,
\[
  \gamma
  \leq
  F(\zbar)
  \leq
  \liminf F(\zbar_t)
  \leq
  \liminf \beta_t
  =
  \gamma.
\]
The first inequality is by $\gamma$'s definition, since
$\A\zbar=\xbar$.
The second is because $F$ is lower semicontinuous
and $\zbar_t\rightarrow\zbar$.
The third is because $F(\zbar_t)<\beta_t$.
Thus, $(\A F)(\xbar)=\gamma=F(\zbar)$, so $\zbar$ attains
the infimum in \eqref{eqn:image-F-dfn2}.

\pfpart{Part~(\ref{thm:AF-lsc:a}):}
Let $\xbar\in\extspac{m}$, and let $\seq{\xbar_t}$ be any sequence in
$\extspac{m}$ converging to $\xbar$.
We aim to show that
\begin{equation}   \label{eq:thm:AF-lsc:1}
  (\A F)(\xbar)
  \leq
  \liminf\,(\A F)(\xbar_t).
\end{equation}
Let $\alpha=  \liminf\,(\A F)(\xbar_t)$.
If $\alpha=+\infty$ then
\eqref{eq:thm:AF-lsc:1} holds trivially.
Otherwise, let $\beta\in\R$ be such that $\beta>\alpha$.
Then there exists a sequence of indices $\seq{s(t)}$ in $\nats$
such that
$s(1)<s(2)<\dotsb$
and for which $(\A F)(\xbar_{s(t)})<\beta$ for all $t$.
For each $t$, by part~(\ref{thm:AF-lsc:b}), there exists
$\zbar_t\in\extspace$ such that $\A\zbar_t=\xbar_{s(t)}$
and $(\A F)(\xbar_{s(t)})=F(\zbar_t)$.
By sequential compactness, the sequence $\seq{\zbar_t}$ has a
convergent subsequence; by discarding all other elements, as well as
the corresponding elements of the sequence $\seq{s(t)}$, we can
assume that $\seq{\zbar_t}$ converges to some point
$\zbar\in\extspace$.
Since $\xbar_{s(t)}\rightarrow\xbar$ and
$\xbar_{s(t)}=\A\zbar_t\rightarrow\A\zbar$
(by \Cref{thm:linear:cont}\ref{thm:linear:cont:b}), it follows that $\A\zbar=\xbar$.
Thus,
\[
  (\A F)(\xbar)
  \leq
  F(\zbar)
  \leq
  \liminf F(\zbar_t)
  =
  \liminf\,(\A F)(\xbar_{s(t)})
  \leq
  \beta.
\]
The first inequality is by $\A F$'s definition, and since
$\A\zbar=\xbar$.
The second is because $F$ is lower semicontinuous and
$\zbar_t\rightarrow\zbar$.
The third is because $(\A F)(\xbar_{s(t)})<\beta$ for all $t$.
Since this holds for all $\beta>\alpha$, it follows that
$(\A F)(\xbar)\leq \alpha$, proving
\eqref{eq:thm:AF-lsc:1}.%
\indexg{linear image of astral function|)}%
\qedhere
\end{proof-parts}
\end{proof}

\indexg{linear image of standard function!extension of|(}%
\indexg{lower semicontinuous extension!linear image of function@of linear image of function|(}%
Using this theorem, we can now show that the extension of $\A f$
always has the expected form:

\begin{theorem}   \label{thm:inf-lin-ext}
  Let $f:\Rn\rightarrow\Rext$ and $\A\in\Rmn$,
  and let $\xbar\in\extspac{m}$.
  Then
  \begin{equation}   \label{eqn:thm:inf-lin-ext:1}
     \Afext(\xbar)
     =
     \inf\,\bigBraces{\fext(\zbar):\:\zbar\in\extspace,\,\A\zbar=\xbar}.
  \end{equation}
  That is, $\Afextnp = \A \fext$.

  Furthermore, if $\xbar$ is in the astral column space of $\A$,
  then the infimum appearing in \eqref{eqn:thm:inf-lin-ext:1}
  is attained by some $\zbar\in\extspace$ with $\A\zbar=\xbar$
  (and otherwise is vacuous).
\end{theorem}

\begin{proof}
We first show that $(\A\fext)(\xbar)\geq \Afext(\xbar)$.
This is trivial if $(\A\fext)(\xbar)=+\infty$.
Otherwise, $(\A\fext)(\xbar)<+\infty$, implying
the set on the right-hand side of \eqref{eqn:thm:inf-lin-ext:1} is not
empty.
Consider any $\zbar\in\extspace$ such that $\A\zbar=\xbar$, and let $\seq{\zz_t}$ be a sequence in $\Rn$ converging to $\zbar$
with $f(\zz_t)\rightarrow\fext(\zbar)$
(which exists by \Cref{pr:d1}).
Then
\[
  \fext(\zbar)
  =
  \lim f(\zz_t)
  \geq
  \liminf\,(\A f)(\A\zz_t)
  \geq
  \Afext(\xbar),
\]
where the first inequality is because $f(\zz_t)\geq (\A f)(
\A\zz_t)$
by \Cref{pr:std-lin-img-props}(\ref{pr:std-lin-img-props:a}),
and the second inequality is because $\A\zz_t\rightarrow\A\zbar=\xbar$
by continuity (\Cref{thm:linear:cont}\ref{thm:linear:cont:b}). Taking infimum over all
$\zbar\in\extspace$ such that $\A\zbar=\xbar$ then yields
$(\A \fext)(\xbar)\geq \Afext(\xbar)$, as claimed.

It remains to show $\A\fext\leq \Afextnp$.
Let $\xx\in\Rm$.
Then
\begin{align*}
  (\A\fext)(\xx)
  &=
  \inf\,\bigBraces{\fext(\zbar) :\: \zbar\in\extspace,\, \A\zbar=\xx}
  \\
  &\leq
  \inf\,\bigBraces{\fext(\zz) :\: \zz\in\Rn,\, \A\zz=\xx}
  \\
  &\leq
  \inf\,\bigBraces{f(\zz) :\: \zz\in\Rn,\, \A\zz=\xx}
  =
  (\A f)(\xx).
\end{align*}
The first inequality is because $\Rn\subseteq\extspace$,
and the second is by \Cref{pr:h:1}(\ref{pr:h:1a}).

Thus, $(\A\fext)(\xx)\leq(\A f)(\xx)$ for all $\xx\in\Rm$;
that is, $\A\fext$ is majorized by $\A f$ when restricted to $\Rm$.
Furthermore, $\A\fext$ is lower semicontinuous by
\Cref{thm:AF-lsc}(\ref{thm:AF-lsc:a})
since $\fext$ is lower semicontinuous
(by \Cref{prop:ext:F}\ref{prop:ext:F:a}).
Therefore, $\A\fext\leq\Afextnp$
by \Cref{prop:ext:F}(\ref{prop:ext:F:b}),
proving $\A\fext=\Afextnp$.
That the infimum in \eqref{eqn:thm:inf-lin-ext:1} is attained
(unless vacuous) now follows directly from
\Cref{thm:AF-lsc}(\ref{thm:AF-lsc:b}).
\end{proof}

The next \namecref{thm:inf-lin-ext-seq}
provides a different expression for the extension
of $\A f$, along the lines of the original definition of $\fext$ given
in \eqref{eq:e:7}:

\begin{theorem}   \label{thm:inf-lin-ext-seq}
  Let $f:\Rn\rightarrow\Rext$, and
  $\A\in\Rmn$.
  Let $\xbar\in\extspac{m}$.
  Then
  \begin{equation}   \label{eqn:thm:inf-lin-ext-seq:1}
     \Afext(\xbar)
     =
     \InfseqLiminf{\seq{\zz_t}}{\Rn}{\A\zz_t\rightarrow \xbar}{f(\zz_t)}.
  \end{equation}
  Furthermore, if $\xbar$ is in the astral column space of $\A$,
  then there exists a sequence $\seq{\zz_t}$ in $\Rn$
  with $\A\zz_t\rightarrow\xbar$ and $f(\zz_t)\rightarrow\Afext(\xbar)$.
\end{theorem}

\begin{proof}
Let $R$ denote the value of the expression on the right-hand side of
\eqref{eqn:thm:inf-lin-ext-seq:1}, which we aim to show is equal to
$\Afext(\xbar)$.

We first show $R\geq \Afext(\xbar)$. If the infimum in \eqref{eqn:thm:inf-lin-ext-seq:1} is over an empty set then $R=+\infty$ and the claim holds. Otherwise, consider any sequence $\seq{\zz_t}$ in $\Rn$
such that $\A\zz_t\to\xbar$. Then
\[
  \liminf f(\zz_t)
  \ge \liminf\,(\A f)(\A\zz_t)
  \ge \Afext(\xbar),
\]
where the first inequality is by
\Cref{pr:std-lin-img-props}(\ref{pr:std-lin-img-props:a}),
and the second by definition of extensions
(Eq.~\ref{eq:e:7})
since $\A\zz_t\to\xbar$.
Taking infimum over all sequences $\seq{\zz_t}$ in $\Rn$ such that $\A\zz_t\to\xbar$ then yields $R\ge\Afext(\xbar)$.

For the reverse inequality, note first that
if $\xbar$ is not in $\A\negKern$'s astral column space, then
$\Afext(\xbar)=+\infty$ by
\Cref{thm:inf-lin-ext},
implying $R=+\infty$ as well, so $\Afext(\xbar)=R$.
We therefore assume henceforth that $\xbar$ is in $\A$'s astral column
space.

By \Cref{thm:inf-lin-ext}, there exists
$\zbar\in\extspace$ with $\A\zbar=\xbar$ that attains the infimum in
\eqref{eqn:thm:inf-lin-ext:1}
so that $\Afext(\xbar)=\fext(\zbar)$.
By \Cref{pr:d1},
there exists a sequence $\seq{\zz_t}$ in~$\Rn$
with $\zz_t\rightarrow\zbar$ and
$f(\zz_t)\rightarrow\fext(\zbar)$.
Also, $\A\zz_t\rightarrow\A\zbar=\xbar$
(by \Cref{thm:linear:cont}\ref{thm:linear:cont:b}), so
\[
  \Afext(\xbar)
  =
  \fext(\zbar)
  =
  \lim f(\zz_t)
  \geq
  R
  \geq
  \Afext(\xbar).
\]
The first inequality is by $R$'s definition, and the second was shown
above.
Thus, $\Afext(\xbar)=R$ and $\seq{\zz_t}$ has the stated properties.%
\indexg{linear image of standard function!extension of|)}%
\indexg{lower semicontinuous extension!linear image of function@of linear image of function|)}%
\end{proof}

\indexg{infimal convolution|(}%
The \emph{infimal convolution} of two proper functions
$f:\Rn\rightarrow\Rext$
and
$g:\Rn\rightarrow\Rext$
is the function
$h:\Rn\rightarrow\Rext$ defined, for $\zz\in\Rn$, by
\[
  h(\zz)
  =
  \inf\,\BigBraces{f(\xx) + g(\zz-\xx)
         :\:
        \xx\in\Rn
      }.
\]
As an application of the foregoing results, we can
characterize the extension of such a function.
To compute $\hext$, define the function $s:\R^{2n}\to\eR$ by
$s(\rpairf{\xx}{\yy})=f(\xx)+g(\yy)$
for $\xx,\yy\in\Rn$,
and let $\A=[\Idnn,\Idnn]$ be the $n\times 2n$ matrix
such that
$\A\rpair{\xx}{\yy}=\xx+\yy$ for $\xx,\yy\in\Rn$. We can then rewrite $h$ as
\begin{align*}
  h(\zz)
  &=
  \inf\,\BigBraces{s(\rpairf{\xx}{\yy})
         :\:
        \xx,\yy\in\Rn,\,\A\rpair{\xx}{\yy}=\zz
      }
  \\
  &=
  \inf\,\BigBraces{s(\ww)
         :\:
        \ww\in\R^{2n},\,\A\ww=\zz
      }.
\intertext{%
By \Cref{thm:inf-lin-ext,thm:inf-lin-ext-seq},
it then follows that,
for $\zbar\in\extspace$,}
  \hext(\zbar)
  &=
  \inf\,\BigBraces{\sext(\wbar)
         :\:
        \wbar\in\extspac{2n},\,\A\wbar=\zbar
      }
  \\[4pt]
  &=
  \InfseqLiminf{\seq{\xx_t},\seq{\yy_t}}{\Rn}{\xx_t+\yy_t\rightarrow \zbar}
               {\bigBracks{f(\xx_t) + g(\yy_t)}}.%
\indexg{infimal convolution|)}%
\end{align*}

\section{Pointwise supremum of a collection of functions}

\indexg{suprema and maxima, pointwise!extension of|(}%
\indexg{lower semicontinuous extension!pointwise supremum@of pointwise supremum|(}%
We next consider the extension of
the pointwise maximum or supremum of a collection of
functions:
\[
  h(\xx) = \sup_{i\in\indset} f_i(\xx),
\]
defined for $\xx\in\Rn$,
where $f_i:\Rn\rightarrow\Rext$ for all $i$ in some
index set $\indset$.
If each $f_i$ is convex, then $h$ is as well
(\Cref{roc:thm5.5}).
In this case, we expect $h$'s extension to be
\begin{equation}  \label{eq:hext-sup-fexti}
  \hext(\xbar) = \sup_{i\in\indset} \efsub{i}(\xbar)
\end{equation}
for $\xbar\in\extspace$.
Below, we establish sufficient conditions for when this holds.

First note that this is not always the case.
For instance, for $f$ and $g$ as defined in
\Cref{ex:pos-neg} (or as in \Cref{ex:KLsets:extsum-not-sum-exts}),
we can set $h=\max\set{f,g}=f+g$, and obtain $\hext\ne\ef+\eg=\max\set{\ef,\eg}$.

Nevertheless, we can show that $\hext(\xbar) \geq \sup_{i\in\indset} \efsub{i}(\xbar)$ for all $\xbar\in\eRn$, and that the equality holds, under lower semicontinuity assumptions, for $\xbar=\xx\in\Rn$.

\begin{proposition}  \label{pr:ext-sup-of-fcns-bnd}
  Let $f_i:\Rn\rightarrow\Rext$ for $i\in\indset$, and let
  $h = \sup_{i\in\indset} f_i$.
  Then:
  \begin{letter-compact}
  \item \label{pr:ext-sup:a}
    For $\xbar\in\extspace$,
    $\hext(\xbar) \geq \sup_{i\in\indset} \efsub{i}(\xbar)$.
  \item \label{pr:ext-sup:b}
    For $\xx\in\Rn$,
    if $f_i$ is lower semicontinuous at $\xx$ for all $i\in\indset$,
    then $h$ is as well,
    and $\hext(\xx)=\sup_{i\in\indset}\efsub{i}(\xx)$.
  \end{letter-compact}
\end{proposition}

\begin{proof}
~
\begin{proof-parts}
\pfpart{Part~(\ref{pr:ext-sup:a}):}
For all $i\in\indset$,
$h\geq f_i$, implying $\hext\geq\efsub{i}$
by \Cref{pr:h:1}(\ref{pr:h:1:geq}), proving the claim.

\pfpart{Part~(\ref{pr:ext-sup:b}):}
Let $\xx\in\Rn$ and suppose each $f_i$ is lower semicontinuous at
$\xx$.
Then
\[
  h(\xx)
  \geq
  (\lsc h)(\xx)
  =
  \hext(\xx)
  \ge
  \sup_{i\in\indset} \efsub{i}(\xx)
  =
  \sup_{i\in\indset} f_i(\xx)
  =
  h(\xx).
\]  
The first inequality and first equality are
by \Cref{pr:h:1}(\ref{pr:h:1a}).
The second inequality is by part~(\ref{pr:ext-sup:a}).
The second equality is also
by \Cref{pr:h:1}(\ref{pr:h:1a})
since each $f_i$ is lower semicontinuous at $\xx$.
This proves both parts of the claim.
\qedhere
\end{proof-parts}
\end{proof}

When $\indset$ is finite, we can prove sufficient conditions for when
the equality holds for general $\xbar\in\eRn$.
These will all follow from the next lemma:

\begin{lemma}  \label{lem:ext-sup-fcns-seq}
  Let $f_i:\Rn\rightarrow\Rext$ for $i=1,\dotsc,m$, and let
  \[ h=\max\regBraces{f_1,\dotsc,f_m}. \]
  Let $\xbar\in\extspace$, and
  let $\seq{\xx_t}$ be a sequence in $\Rn$ that converges to $\xbar$.
  Suppose $f_i(\xx_t)\rightarrow\efsub{i}(\xbar)$ for all $i=1,\dotsc,m$.
  Then
  $h(\xx_t)\rightarrow\hext(\xbar)$
  and
  \begin{equation*}  %
     \hext(\xbar)
     =
     \max\regBraces{\efsub{1}(\xbar),\dotsc,\efsub{m}(\xbar)}.
  \end{equation*}
\end{lemma}

\begin{proof}
Let $\indset=\set{1,\dots,m}$ be the index set for $i$. Then
\begin{align}
\notag
  \limsup h(\xx_t)
  &=
  \limsup \BigBracks{\max_{i\in\indset} f_i(\xx_t)}
\\
\label{eq:ext-sup:lemma}
  &=
  \max_{i\in\indset} \BigBracks{\limsup f_i(\xx_t)}
  =
  \max_{i\in\indset}\efsub{i}(\xbar),
\end{align}
where the second equality is by repeated application of \Cref{prop:lim:eR}(\ref{i:limsup:eR:max}),
and the last equality is by the assumption in the lemma. Therefore,
\[
  \hext(\xbar)
  \leq
  \liminf h(\xx_t)
  \leq
  \limsup h(\xx_t)
  =
  \max_{i\in\indset}\efsub{i}(\xbar)
  \leq
  \hext(\xbar).
\]
The first inequality is by definition of extensions
(Eq.~\ref{eq:e:7})
since $\xx_t\rightarrow\xbar$.
The equality is by \eqref{eq:ext-sup:lemma}.
The last inequality follows from \Cref{pr:ext-sup-of-fcns-bnd}(\ref{pr:ext-sup:a}).
This proves both of the lemma's claims.
\end{proof}

Focusing now on the case when $m=2$, so $h=\max\{f,g\}$, we use \Cref{lem:ext-sup-fcns-seq}
to show how continuity properties of $f$ and $g$ at $\xbar$ affect the form of $\hext(\xbar)$.

\begin{proposition}  \label{pr:ext-sup-fcns-cont}
  Let $f:\Rn\rightarrow\Rext$ and $g:\Rn\rightarrow\Rext$,
  and let
  $h=\max\{f,g\}$.
  Let $\xbar\in\eRn$. Then
  \begin{letter-compact}
  \item  \label{pr:ext-sup-fcns-cont:a}
    If either $f$ or $g$ is extensibly continuous at $\xbar$,
    then $\hext(\xbar)=\max\{\fext(\xbar),\gext(\xbar)\}$.
  \item  \label{pr:ext-sup-fcns-cont:b}
    If both $f$ and $g$ are extensibly continuous at $\xbar$,
    then $h$ is as well.
  \end{letter-compact}
\end{proposition}

\begin{proof}
~
\begin{proof-parts}
\pfpart{Part~(\ref{pr:ext-sup-fcns-cont:a}):}
Suppose $g$ is extensibly continuous at $\xbar$.
By \Cref{pr:d1},
there exists a sequence $\seq{\xx_t}$ in $\Rn$ converging to
$\xbar$ and for which
$f(\xx_t)\rightarrow\fext(\xbar)$.
Since $g$ is extensibly continuous, it also follows that
$g(\xx_t)\rightarrow\gext(\xbar)$.
Therefore, the claim follows by
\Cref{lem:ext-sup-fcns-seq}.

\pfpart{Part~(\ref{pr:ext-sup-fcns-cont:b}):}
Suppose $f$ and $g$ are extensibly continuous at $\xbar$.
Let $\seq{\xx_t}$ be any sequence in $\Rn$ converging to
$\xbar$.
Then
$f(\xx_t)\rightarrow\fext(\xbar)$
and
$g(\xx_t)\rightarrow\gext(\xbar)$,
since the two functions are extensibly continuous.
Therefore, $h(\xx_t)\rightarrow\hext(\xbar)$
by
\Cref{lem:ext-sup-fcns-seq}.
Since this holds for every such sequence, the claim follows.
\qedhere
\end{proof-parts}
\end{proof}

We next give a condition based on convexity that holds when the relative interiors of the functions' effective domains have a point in common:

\begin{theorem}   \label{thm:ext-finite-max-convex}
  Let $f_i:\Rn\rightarrow\Rext$ be convex, for $i=1,\dotsc,m$,
  and let
  \[  h=\max\Braces{f_1,\dotsc,f_m}. \]
  Assume $\bigcap_{i=1}^m \ri(\dom{f_i}) \neq \emptyset$.
  Then
  \[
    \hext
    =
    \max\Braces{\efsub{1},\dotsc,\efsub{m}}.
  \]
\end{theorem}

\begin{proof}
Let $\xbar\in\extspace$.
By \Cref{thm:seq-for-multi-ext},
there exists a sequence $\seq{\xx_t}$ in $\Rn$ converging to $\xbar$
and with
$f_i(\xx_t)\rightarrow\efsub{i}(\xbar)$
for $i=1,\dotsc,m$.
Therefore, by \Cref{lem:ext-sup-fcns-seq},
$\hext(\xbar) = \max\regBraces{\efsub{1}(\xbar),\dotsc,\efsub{m}(\xbar)}$.
\end{proof}

\Cref{pr:ext-sup-fcns-cont} and \Cref{thm:ext-finite-max-convex} cannot be generalized to work with
infinite collections of functions without making further assumptions. In fact, any
example of a closed proper convex function $h:\Rn\to\eRn$ such that $\eh\ne\hdub$
can be turned into a counterexample, because $h(\xx)=\sup_{\uu\in\dom h^*}[-h^*(\uu)+\xx\inprod\uu]$ and
$\hdub(\xbar)=\sup_{\uu\in\dom h^*}[-h^*(\uu)+\xbar\inprod\uu]$
(by \Cref{thm:fext-dub-sum}\ref{thm:fext-dub-sum:b}). Here, $h$ is expressed as a supremum over affine
functions $f_\uu(\xx)=-h^*(\uu)+\xx\inprod\uu$ for $\xx\in\Rn$, with $\uu\in\dom h^*$.

\indexg{Restricted linear function!extension of supremum and|(}%
We demonstrate this on an example based on the restricted linear function from
\Cref{ex:negx1-else-inf}, whose extension does not coincide with its astral biconjugate,
as we saw in \Cref{ex:biconj:notext}.

\begin{example}
Let $U=\set{\uu\in\R^2:\:u_1=-1,\,u_2\le 0}$ and consider the collection of linear functions
$f_\uu(\xx)=\xx\inprod\uu$ for $\xx\in\R^2$, with $\uu\in U$. Define $h:\R^2\to\eR$, for $\xx\in\Rn$, as
\[
  h(\xx)=\sup_{\uu\in U} f_\uu(\xx)
        =\sup_{u_2\in\R_{\le 0}} [-x_1 + u_2 x_2]
        =
   \begin{cases}
     -x_1
     & \text{if $x_2\geq 0$,}
   \\
     +\infty
     & \text{otherwise,}
   \end{cases}
\]
so $h$ is the restricted linear function from
\Cref{ex:negx1-else-inf,ex:biconj:notext}.
Let $\ebar=\limray{\ee_1}$.
As we saw in \Cref{ex:biconj:notext}, for all $\xx\in\R^2$,
\[
   \sup_{\uu\in U} \efsub{\uu}(\ebar\plusl\xx)
   =
   \sup_{\uu\in\R^2:\:u_1=-1,\,u_2\le 0} (\ebar\plusl\xx)\cdot\uu
   =-\infty,
\]
but
\[
   \eh(\ebar\plusl\xx)
   =
   \begin{cases}
     -\infty & \text{if $x_2\geq 0$,}
     \\
     +\infty & \text{otherwise.}
   \end{cases}
\]
Thus, $\eh(\ebar\plusl\xx)\ne\sup_{\uu\in U} \fextsub{\uu}(\ebar\plusl\xx)$ if $x_2<0$.%
\indexg{suprema and maxima, pointwise!extension of|)}%
\indexg{lower semicontinuous extension!pointwise supremum@of pointwise supremum|)}%
\indexg{Restricted linear function!extension of supremum and|)}%
\end{example}

\indexg{constraint region, extension when restricted to|(}%
\indexg{lower semicontinuous extension!restricted function@of restricted function|(}%
As a simple corollary of \Cref{thm:ext-finite-max-convex}, we can compute the extension of a function
obtained by restricting a function to a convex subregion:

\begin{corollary}   \label{cor:ext-restricted-fcn}
  Let $f:\Rn\rightarrow\Rext$ be convex, and
  let $S\subseteq\Rn$ be convex.
  Assume
  $\ri(\dom{f})\cap (\ri S)\neq\emptyset$.
  Let $h:\Rn\rightarrow\Rext$ be defined, for $\xx\in\Rn$, by
  \begin{align*}
     h(\xx)
     &=
     f(\xx)\plusu\indf{S}(\xx)
     =
        \begin{cases}
              f(\xx)   & \text{if $\xx\in S$,} \\
              +\infty  & \text{otherwise.}
        \end{cases}
  \end{align*}
  Let $\xbar\in\extspace$.
  Then
  \begin{align*}
     \hext(\xbar)
     &=
     \fext(\xbar)\plusu\indfa{\Sbar}(\xbar)
     =
        \begin{cases}
              \fext(\xbar)  & \text{if $\xbar\in \Sbar$,} \\
              +\infty       & \text{otherwise.}
        \end{cases}
  \end{align*}
\end{corollary}

\begin{proof}
Let $g:\Rn\rightarrow\Rext$ be defined,
for $\xx\in\Rn$, by
\[
   g(\xx)
     =
        \begin{cases}
              -\infty   & \text{if $\xx\in S$,} \\
              +\infty   & \text{otherwise,}
        \end{cases}
\]
which is convex.
To compute $g$'s extension, let $g'=\expex\circ\posKern g$.
Then $g'=\inds$ so $\gpext=\indfa{\Sbar}$ by
\Cref{pr:inds-ext}.
Further, by
\Cref{pr:j:2}(\ref{pr:j:2c}),
$\gpext=\expex\circ\posKern\gext$.
Hence, 
\begin{equation}   \label{eq:cor:ext-restricted-fcn:1}
   \gext(\xbar)
     =
        \begin{cases}
              -\infty   & \text{if $\xbar\in\Sbar$,} \\
              +\infty   & \text{otherwise.}
        \end{cases}
\end{equation}

Note that $h(\xx)=\max\{f(\xx),g(\xx)\}$ for $\xx\in\Rn$.
Also,
$\ri(\dom{f})\cap\ri(\dom{g})\neq\emptyset$, by assumption,
since $\dom{g}=S$.
Therefore, we can apply
\Cref{thm:ext-finite-max-convex}, yielding
$\hext(\xbar)=\max\{\fext(\xbar),\gext(\xbar)\}$.
Combined with \eqref{eq:cor:ext-restricted-fcn:1},
this proves the claim.%
\indexg{constraint region, extension when restricted to|)}%
\indexg{lower semicontinuous extension!restricted function@of restricted function|)}%
\end{proof}

\section{Composition with a nondecreasing function}
\label{sec:calc-ext-comp-inc-fcn}

\indexg{lower semicontinuous extension!composition with nondecreasing function@of composition with nondecreasing function|(}%
\indexg{nondecreasing function, composition with!extension of|(}%
Next, we consider the composition
$G\circ f$
of
a function $f:\Rn\rightarrow\Rext$ with a nondecreasing function $G:\eR\rightarrow\eR$.
In \Cref{sec:exp-comp}, we saw that
when $G$ is the extension of a suitable convex function (say $G=\expex$)
then $\Gofextnp=G\circ\ef$. Here we generalize that result to a wider range of functions $G$.

First, we give two examples showing that it is not enough to assume that $G$ is lower semicontinuous:

\begin{example}
\label{ex:Gof:pos-neg}
Let $G=\indfa{\le 0}$ be the indicator of the set $[-\infty,0]$, and
define $f:\R\rightarrow\Rext$, for $x\in\R$, as
\[
  f(x)=
  \begin{cases}
    x &\text{if $x>0$,}
  \\
    +\infty &\text{otherwise.}
  \end{cases}
\]
Let $h=G\circ f$. Then $h\equiv+\infty$, so also
$\eh\equiv+\infty$. However, for $\ex\in\eR$,
\[
  \ef(\ex)=
  \begin{cases}
    \ex &\text{if $\ex\ge0$,}
  \\
    +\infty &\text{otherwise,}
  \end{cases}
\]
so $G\circ\ef=\indfa{\set{0}}$. Thus, $\eh\ne G\circ\ef$.
\end{example}

The issue above is that $f$ is not lower semicontinuous at $0$.
The next example shows that even if $f$ is closed,
convex, and proper, and even if $G$ is an extension of a closed,
proper convex function,
we might still have $\Gofextnp\ne G\circ\ef$.

\begin{example}
\label{ex:Gof:sideways}
For $\xx\in\R^3$, let
\[
   f(\xx) = f(x_1,x_2,x_3)=
   \begin{cases}
       x_3     & \text{if $\xx\in K$,}
   \\
       +\infty & \text{otherwise,}
   \end{cases}
\]
where $K$ is the sideways cone from
\Cref{ex:KLsets:extsum-not-sum-exts}
(Eq.~\ref{eqn:bad-set-eg-K}),
and let $G=\indfa{\le 0}$,
which is the extension of the
standard indicator function $\indf{\leq 0}$
(the indicator for the set $\Rneg$).
The composition $h=G\circ f$ is then the indicator
$\indf{K\cap L}$, where
$L$ is the plane defined in
\Cref{ex:KLsets:extsum-not-sum-exts}
(Eq.\ref{eqn:bad-set-eg-L}).

As in \Cref{ex:KLsets:extsum-not-sum-exts},
let $\xbar=\limray{\ee_2}\plusl\ee_1$, and
let $\xx_t=\trans{[1,\,t,\,1/(2t)]}$, which converges to $\xbar$.
Then $f(\xx_t)=1/(2t)\rightarrow 0$.
Since $\inf f\geq 0$, it follows that $\fext(\xbar)=0$.
Thus, $G\regParens{\fext(\xbar)}=G(0)=0$.
On the other hand, as argued
in \Cref{ex:KLsets:extsum-not-sum-exts},
$\hext(\xbar)=+\infty$.
\end{example}

Although lower semicontinuity of $G$ does not suffice to
ensure that $\Gofextnp=G\circ\ef$, it does suffice to ensure an inequality, and
the continuity of $G$ suffices to ensure equality:

\begin{proposition}
\label{pr:Gf-cont}
  Let $f:\Rn\rightarrow\eR$,
  let $G:\eR\rightarrow\eR$ be nondecreasing,
  and let $h=G\circ f$.
  Let $\xbar\in\extspace$.
  Then:
  \begin{letter-compact}
  \item
  \label{pr:Gf-cont:a}
    If $G$ is lower semicontinuous at $\ef(\xbar)$ then
    $\eh(\xbar)\geq G\regParens{\ef(\xbar)}$.
  \item
  \label{pr:Gf-cont:b}
    If $G$ is continuous at $\ef(\xbar)$ then
    $\eh(\xbar)= G\regParens{\ef(\xbar)}$. If additionally $f$ is extensibly continuous at $\xbar$
    then $h$ is extensibly continuous at $\xbar$ as well.
  \end{letter-compact}
\end{proposition}

\begin{proof}
~
\begin{proof-parts}
\pfpart{Part (\ref{pr:Gf-cont:a}):}
Suppose $G$ is lower semicontinuous at $\ef(\xbar)$.
Let $\seq{\xx_t}$ be a sequence in $\Rn$ such that $\xx_t\to\xbar$ and $h(\xx_t)\to\eh(\xbar)$
(which exists by \Cref{pr:d1}).
Further, by sequential compactness of $\eR$, there exists a subsequence of
$\seq{f(\xx_t)}$ that converges in $\eR$; by discarding all other
sequence elements, we can assume the entire sequence converges to
some limit $\alpha\in\Rext$.
Note that $\alpha=\lim f(\xx_t)\geq\fext(\xbar)$
by $\fext$'s definition (Eq.~\ref{eq:e:7}).

We have
$\eh(\xbar)=\lim h(\xx_t) = \lim G\regParens{f(\xx_t)}$.
Therefore, to complete the proof, it suffices to show that
\begin{equation}   \label{eq:pr:Gf-cont:r1}
  \lim G\bigParens{f(\xx_t)}\geq G\bigParens{\fext(\xbar)}.
\end{equation}
If $\alpha=\fext(\xbar)$ then this follows from our assumption that
$G$ is lower semicontinuous at $\fext(\xbar)$.
Otherwise, if $\alpha>\fext(\xbar)$, then
for all $t$ sufficiently large,
$f(\xx_t)>\fext(\xbar)$, implying,
since $G$ is nondecreasing, that
$G\regParens{f(\xx_t)}\geq G\regParens{\fext(\xbar)}$.
Therefore, \eqref{eq:pr:Gf-cont:r1} holds in this case as well,
completing the proof.

\pfpart{Part (\ref{pr:Gf-cont:b}):}
Suppose $G$ is continuous (so also lower semicontinuous) at $\ef(\xbar)$,
and let $\seq{\xx_t}$ be a sequence in $\Rn$ such that $\xx_t\to\xbar$ and $f(\xx_t)\to\ef(\xbar)$.
Then by sequential compactness of $\Rext$, the sequence
$\seq{h(\xx_t)}$ has a convergent subsequence; by discarding all other
sequence elements, we can assume the entire sequence converges in $\Rext$.
Then
\[
 \eh(\xbar)
 \le
 \lim
 h(\xx_t)
 =
 \lim
 G\bigParens{f(\xx_t)}
 =
 G\bigParens{\fext(\xbar)},
\]
where the inequality is by $\eh$'s definition,
and the last equality follows by continuity of $G$
at $\fext(\xbar)$ since
$f(\xx_t)\to\ef(\xbar)$.
Combining with part~(\ref{pr:Gf-cont:a})
then yields $\eh(\xbar)= G\regParens{\ef(\xbar)}$.

For the second part of the claim, suppose also that $f$ is extensibly
continuous at $\xbar$.
Let $\seq{\xx_t}$ be any sequence in $\Rn$ such that $\xx_t\to\xbar$.
Then
\[
  \eh(\xbar)
  =
  G\bigParens{\ef(\xbar)}
  =
  G\bigParens{\lim f(\xx_t)}
  =
  \lim G\bigParens{f(\xx_t)}
  =
  \lim h(\xx_t).
\]
The first equality was just proved.
The second is by extensible continuity of $f$ at $\xbar$.
The third is by continuity of $G$ at $\fext(\xbar)$.
Thus, $h$ is extensibly continuous at $\xbar$.
\qedhere
\end{proof-parts}
\end{proof}

Next we focus on compositions of the form $h=\eg\circ\negKern f$ where $f:\Rn\rightarrow\Rext$ is convex and
$\eg$ is the extension of a nondecreasing convex function
$g:\R\rightarrow\Rext$.
We will show that $\eh=\eg\circ\negKern\ef$ without assuming that $\eg$
is continuous, as long as $f$ attains a value in the interior of $\dom g$
(which was not the case in \Cref{ex:Gof:pos-neg,ex:Gof:sideways}).

\begin{theorem}
\label{thm:Gf-conv}
  Let $f:\Rn\rightarrow\Rext$ be convex,
  let $g:\R\rightarrow\Rext$ be convex and nondecreasing,
  and assume there exists $\xing\in\Rn$ such that
  $f(\xing)<\sup(\dom g)$.
  Let $h=\eg\circ\negKern f$. Then
  $\eh=\eg\circ\negKern\ef$.
\end{theorem}

\begin{proof}
We note first that $\gext$ is nondecreasing by
\Cref{pr:conv-inc:prop}(\ref{pr:conv-inc:nondec}).

Let $\xbar\in\eRn$. We will show that $\eh(\xbar)=\eg(\ef(\xbar))$. This holds if $\eg$
is continuous at~$\ef(\xbar)$ (by \Cref{pr:Gf-cont}\ref{pr:Gf-cont:b}), so we
assume henceforth that $\eg$ is \emph{not}
continuous at $\ef(\xbar)$.
By \Cref{pr:conv-inc:prop}(\ref{pr:conv-inc:infsup},\ref{pr:conv-inc:discont}), this can only happen when $g$ and  $\eg$ are both discontinuous at a single point $z\in\R$,
namely, $z=\sup(\dom g)$.
Since $\eg$ is not continuous at $\ef(\xbar)$,
this implies $z=\ef(\xbar)$.

We show next that $f>-\infty$:
By \Cref{pr:h:1}(\ref{pr:h:1aa}),
$\lscfextf(\xbar)=\fext(\xbar)=z\in\R$,
implying,
by \Cref{pr:h:1}(\ref{pr:h:1b-neigh})
(applied to $\lsc{f}$ with $U=\extspace$ and $V=\R$),
that there exists $\xx'\in\Rn$ with
$(\lsc f)(\xx')\in\R$.
Therefore, $\lsc{f}$, which is convex and lower semicontinuous, is
also proper by
\Cref{pr:improper-vals}(\ref{pr:improper-vals:cor7.2.1}).
Hence, $f\geq\lsc{f}>-\infty$.

In particular, this shows that
$f(\xing)\in(-\infty,z)$.

\begin{claimpx}
  There exists a sequence $\seq{\xx_t}$ in $\Rn$ that
  converges to $\xbar$ and with
  $f(\xx_t)\in (-\infty,z)$ for all $t$.
\end{claimpx}

\begin{proofx}
By \Cref{pr:d1}, there exists a sequence $\seq{\xx'_t}$ in
$\Rn$ converging to $\xbar$ and with
$f(\xx'_t)\rightarrow\fext(\xbar)=z$.
Since $z\in\R$, there can be at most finitely many values of $t$ for
which $f(\xx'_t)$ is infinite; discarding these, we can assume
$f(\xx'_t)\in\R$ for all $t$.

If $f(\xx'_t)<z$ for infinitely many values of $t$, then discarding all
other elements results in a sequence with the properties stated
in the claim.
Therefore, we focus henceforth on the alternative case in which
$f(\xx'_t)<z$ for only finitely many sequence elements.
Further,
by discarding these,
we can assume henceforth that
$f(\xx'_t) \geq z$ for all $t$.

For each $t$, let
\[
  \lambda_t
  =
  \frac{\bigParens{1-\frac{1}{t}} \bigParens{z - f(\xing)}}
       {f(\xx'_t) - f(\xing)}.
\]
Then $0\leq \lambda_t < 1$ since $f(\xing)<z\leq f(\xx'_t)$.
Also, $\lambda_t\rightarrow 1$ since $f(\xx'_t)\rightarrow z$.
Let
$ \xx_t = \lambda_t \xx'_t + (1-\lambda_t) \xing $.
Then
\[
  f(\xx_t)
  \leq
  \lambda_t f(\xx'_t) + (1-\lambda_t) f(\xing)
  =
  \bigParens{1-\sfrac{1}{t}} z + \sfrac{1}{t} f(\xing)
  <
  z.
\]
The first inequality is because $f$ is convex
(\Cref{pr:stand-cvx-fcn-char}).
The equality is by our choice of $\lambda_t$.
The last inequality is because $f(\xing)<z$.

Further, since $\lambda_t\xx'_t\to\xbar$ (by \Cref{pr:scalar-prod-props}\ref{pr:scalar-prod-props:e}) and $(1-\lambda_t) \xing\to\zero$, we must have $\xx_t\rightarrow\xbar$ (by \Cref{pr:i:7}\ref{pr:i:7g}).
Since $f>-\infty$, as shown above,
the resulting sequence $\seq{\xx_t}$ satisfies all
the claimed properties.
\end{proofx}

Let $\seq{\xx_t}$ be a sequence with the properties stated in the
preceding claim.
Then
\begin{align*}
  \hext(\xbar)
  \leq
  \liminf h(\xx_t)
  &=
  \liminf \gext(f(\xx_t))
  \\
  &=
  \liminf g(f(\xx_t))
  \leq
  \sup_{y<z} g(y)
  =
  \gext(z)
  =
  \gext(\fext(\xbar))
  \leq
  \hext(\xbar).
\end{align*}
The first inequality is by definition of extensions
(Eq.~\ref{eq:e:7}).
The second equality is because $\eg$ agrees with $g$ on $(-\infty,z)$
(\Cref{pr:conv-inc:prop}\ref{pr:conv-inc:discont}),
and since $f(\xx_t)\in(-\infty,z)$ for all $t$.
The second inequality is also because of this latter fact.
The third equality is by
\Cref{pr:conv-inc:prop}(\ref{pr:conv-inc:nondec}).
And the last inequality is by
\Cref{pr:Gf-cont}(\ref{pr:Gf-cont:a})
since $\gext$ is nondecreasing and lower semicontinuous
(\Cref{prop:ext:F}\ref{prop:ext:F:a}).
This completes the proof.%
\indexg{lower semicontinuous extension!composition with nondecreasing function@of composition with nondecreasing function|)}%
\indexg{nondecreasing function, composition with!extension of|)}%
\end{proof}

\indexg{sublevel sets!astral closure of|(}%
\indexg{closure, astral!sublevel sets@of sublevel sets|(}%
As an application, we can use
\Cref{thm:Gf-conv}
to relate the sublevel sets of a function and its extension:

\begin{theorem}
\label{thm:closure-of-sublev-sets}
  Let $f:\Rn\rightarrow\Rext$ be convex,
  let $\beta\in\R$ be such that
  $\beta>\inf f$, and let
  \begin{align*}
    L &= \set{\xx\in\Rn :\: f(\xx) \leq \beta},
    \\
    M &= \set{\xx\in\Rn :\: f(\xx) < \beta}.
  \end{align*}
  Then
  \begin{equation}
  \label{eq:closure-of-sublev-sets}
    \Lbar = \Mbar = \set{\xbar\in\extspace :\: \fext(\xbar) \leq \beta}.
  \end{equation}
\end{theorem}

\begin{proof}
Let $g=\indf{\le\beta}$ be the indicator of the set $(-\infty,\beta]$, so
$\eg=\indfa{\le\beta}$ is the indicator of the set $[-\infty,\beta]$
(by \Cref{pr:inds-ext}).
Let $R$ denote the set
on the right-hand side of \eqref{eq:closure-of-sublev-sets}
and let $h=\eg\circ\negKern f$.
Then $h$ is exactly $\indf{L}$, the indicator for~$L$, and
$\eh=\indfa{\Lbar}$ (by \Cref{pr:inds-ext}).

Since $\beta>\inf f$, there exists a point $\xing$ with
$f(\xing)<\beta=\sup(\dom g)$.
Therefore,
\[
  \indfa{\Lbar}
  =
  \eh=\eg\circ\negKern\ef=\indfa{\le\beta}\circ\ef=\indfa{R},
\]
where the second equality is by \Cref{thm:Gf-conv}.
Thus, $\Lbar=R$, proving the claim
for $\Lbar$.

This also proves the claim for $\Mbar$ since
\[
  \Mbar
  =
  \clbar{\cl M}
  =
  \clbar{\cl L}
  =
  \Lbar,
\]
where the first and third equalities are by
\Cref{pr:closed-set-facts}(\ref{pr:closed-set-facts:aa}),
and the second equality is because
$\cl M = \cl L$
(\Cref{roc:thm7.6-mod}).%
\indexg{sublevel sets!astral closure of|)}%
\indexg{closure, astral!sublevel sets@of sublevel sets|)}%
\end{proof}

\indexg{closure, astral!halfspaces@of halfspaces|(}%
\indexg{halfspaces (standard)!astral closure of|(}%
As a corollary, we can compute the closures of closed and
open halfspaces:

\begin{corollary}
\label{cor:halfspace-closure}
  Let $\uu\in\Rn\wo\{\zero\}$, let $\beta\in\R$,
  and let
  \begin{align*}
    L &= \set{\xx\in\Rn :\: \xx\cdot\uu \leq \beta},
    \\
    M &= \set{\xx\in\Rn :\: \xx\cdot\uu < \beta}.
  \end{align*}
  Then
  \[
  \Lbar = \Mbar = \set{\xbar\in\extspace :\: \xbar\cdot\uu \leq \beta}.
  \]
\end{corollary}

\begin{proof}
Define $f:\Rn\rightarrow\Rext$ to be
$f(\xx)=\xx\cdot\uu$ for $\xx\in\Rn$,
which is convex, and whose extension is
$\fext(\xbar)=\xbar\cdot\uu$
for $\xbar\in\extspace$
(Example~\ref{ex:ext-affine}).
In addition,
\[
  \inf f
  \leq
  \lim_{\lambda\rightarrow-\infty} f(\lambda\uu)
  =
  -\infty
  <
  \beta.
\]
Therefore, we can apply
\Cref{thm:closure-of-sublev-sets}
to $f$,
yielding the claim.%
\indexg{closure, astral!halfspaces@of halfspaces|)}%
\indexg{halfspaces (standard)!astral closure of|)}%
\end{proof}

\part{Convex Sets and Cones}
\label{part:convex-sets}

\chapter{Convex sets}
\label{sec:convexity}

We next study how the notion of convexity can be extended to
astral space.

In $\Rn$, a set is convex if for every pair of points $\xx,\yy$ in the
set, their convex combination
$(1-\lambda)\xx+\lambda\yy$ is also in the set, for all
$\lambda\in [0,1]$. Said differently, a set is convex
if it includes the line segment connecting every pair of points in the set.
As a natural first attempt at extending this notion to astral space,
we might try to define line segments and
convexity in astral space using some kind
of convex combination
$\zbar=(1-\lambda)\xbar+\lambda\ybar$
of two points $\xbar,\ybar\in\extspace$.
We have not actually defined the ordinary sum of two astral points,
as in this expression, nor will we define it.
But if we could,
then for any $\uu\in\Rn$, we would naturally want it
to be the case that
\[
 \zbar\cdot\uu = (1-\lambda) \xbar\cdot\uu +
                 \lambda \ybar\cdot\uu.
\]
Such an expression, however, would be problematic
since $\xbar\inprod\uu$ and $\ybar\inprod\uu$ are not
necessarily summable across all $\uu$ (see \Cref{prop:commute}), so
that the right-hand side might be undefined.

To avoid this difficulty, we take a different approach below.
Rather than generalizing the notion of a convex combination, we will
instead
directly generalize the notion of a line segment and define convexity using
such generalized line segments.
This kind of convexity is called \emph{interval convexity}%
\indexg{interval convexity}
by
\idxvandevel\citet[Chapter 4]{vandeVel}
and an \emph{inner approach}%
\indexg{inner approach}
by
\idxsinger\citet[Section 0.1a]{singer_book}.
In \Cref{sec:sequential}, and in particular in
\Cref{sec:seq-cvx-conic-combs},
we return to the question of how to
define astral analogues of convex combinations.

\section{Halfspaces, segments, outer hull, convexity}
\label{sec:def-convexity}

In this section,
we develop several fundamental notions in a progression that eventually culminates
with a definition of convexity for astral sets.
\indexg{halfspaces, astral|(}%
\indexg{hyperplanes, astral|(}%
We begin by extending the standard definitions of hyperplanes and
halfspaces to astral space.

A (standard) closed halfspace is the set of points $\xx\in\Rn$ for which
$\xx\cdot\uu\leq \beta$,
for some $\uu\in\Rn\wo\{\zero\}$ and $\beta\in\R$.
We can immediately extend this to astral space, defining a
\emph{closed astral halfspace}, denoted $\chsua$, to be the set
\begin{equation}
\label{eq:chsua-defn}
\indexm{h u beta}{$\chsua$}{closed astral halfspace}%
  \chsua = \Braces{\xbar\in\extspace :\: \xbar\cdot\uu \leq \beta}.
\end{equation}
In the same way, \emph{open astral halfspaces} and
\emph{astral hyperplanes} are sets of the form
\begin{equation}   \label{eqn:open-hfspace-defn}
 \Braces{\xbar\in\extspace :\: \xbar\cdot\uu<\beta}
\end{equation}
and
\begin{equation}   \label{eqn:hyperplane-defn}
 \Braces{\xbar\in\extspace :\: \xbar\cdot\uu=\beta},
\end{equation}
respectively.
(When clear from context, we sometimes drop the astral modifier
when referring to such sets.)
As usual, these forms accommodate halfspaces defined by reverse inequalities by substituting
$-\uu$ for $\uu$ and $-\beta$ for $\beta$.%
\indexg{hyperplanes, astral|)}

In fact, the subbase for the astral topology, given in
\eqref{eq:h:3a:sub-alt},
consists exactly of all open astral halfspaces.
Likewise, the base elements for the astral topology
(given in Eq.~\ref{eq:h:3a})
consist of all finite intersections of open astral
halfspaces.
In particular, this means that all
open astral halfspaces are indeed open,
and that all closed astral halfspaces are closed,
being complements of open halfspaces.

The intersection of a closed or open astral halfspace with $\Rn$ is
just the corresponding standard closed or open halfspace in $\Rn$;
for instance,
\begin{equation}  \label{eqn:chsua-cap-rn}
 \chsua\cap\Rn = \regBraces{\xx\in\Rn :\: \xx\cdot\uu \leq \beta}.
\end{equation}
The closure (in $\extspace$) of a closed or open halfspace like this one in
$\Rn$ is exactly the closed astral halfspace $\chsua$
(\Cref{cor:halfspace-closure}).
Similar facts hold for hyperplanes.

As expected, the interior of a closed astral halfspace is the
corresponding open astral halfspace:

\begin{proposition}
\label{prop:intr:H}
Let $\uu\in\Rn\wo\set{\zero}$, $\beta\in\R$, and $H = \braces{\xbar\in\extspace :\: \xbar\cdot\uu \leq \beta}$.
Then
\[
  \intr H=\braces{\xbar\in\extspace :\: \xbar\cdot\uu < \beta}.
\]
\end{proposition}
\begin{proof}
Let $\Hcomp=\eRn\setminus H$. Since $H$ is closed, $\Hcomp$ is open. Let
\[
 M=\Hcomp\cap\Rn=\braces{\xx\in\Rn :\: \xx\cdot\uu > \beta}.
\]
By \Cref{pr:closed-set-facts}(\ref{pr:closed-set-facts:b}) and \Cref{cor:halfspace-closure}, $\clHcomp=\Mbar=\braces{\xbar\in\extspace :\: \xbar\cdot\uu\ge\beta}$.
Hence, by \Cref{pr:closure:intersect}(\ref{pr:closure:intersect:comp}),
$\intr H=\extspace\setminus\clHcomp=\braces{\xbar\in\extspace :\: \xbar\cdot\uu < \beta}$.%
\indexg{halfspaces, astral|)}%
\end{proof}

The starting point for standard convex analysis is the line segment
joining two points $\xx$ and $\yy$.
As discussed above,
in $\Rn$, this is the set of all convex combinations of the two
points, a perspective that does not immediately
generalize to astral space.
However, there is another way of thinking about the line segment
joining $\xx$ and $\yy$, namely, as the intersection of all
halfspaces that include both of the endpoints.
This interpretation generalizes directly to astral space, leading, in
a moment, to a definition of the segment joining two astral points
as the intersection of all astral halfspaces that include both of the points.

\indexg{outer convex hull|(}%
To state this more formally, we first give a more general definition
for the intersection of all closed astral halfspaces that include an arbitrary set
$S\subseteq\extspace$, called the outer convex hull.
A segment will then be a special case in which $S$ has two elements:

\begin{definition}  \label{def:outer-cvx-hull}
\indexg{outer convex hull!defined|(}%
Let $S\subseteq\extspace$.
The \emph{outer convex hull} of $S$
(or \emph{outer hull}, for short),
denoted $\ohull S$, is the
intersection of all closed astral halfspaces that include $S$; that
is,
\begin{equation}    \label{eqn:ohull-defn}
\indexm{conv s}{$\ohull S$}{outer convex hull}%
  \ohull S
  =
  \bigcap{\BigBraces{ \chsua:\:
      \uu\in\Rn\wo\{\zero\},\, \beta\in\R,\, S\subseteq \chsua}}.
\end{equation}
\end{definition}
Note that if $S$ is not included in any halfspace $\chsua$, then the
intersection on the right-hand side of
\eqref{eqn:ohull-defn}
is vacuous and therefore equal to $\extspace$.
Thus, in this case, $\ohull S = \extspace$.
(Alternatively, we could allow $\uu=\zero$ in
Eq.~\ref{eqn:ohull-defn} so that this intersection is never empty since
$\chs{\zero}{0}=\extspace$ will always be included.)%
\indexg{outer convex hull!defined|)}

\indexg{segments, astral!defined|(}%
\indexg{segments, astral|(}%
In these terms, we can now define astral segments:

\begin{definition}
Let $\xbar,\ybar\in\extspace$.
The \emph{segment joining} $\xbar$ and $\ybar$,
denoted $\lb{\xbar}{\ybar}$, is the intersection of all closed astral
halfspaces that include both $\xbar$ and $\ybar$;
thus,
\[
\indexm{seg x y}{$\lb{\xbar}{\ybar}$}{segment joining points}%
\lb{\xbar}{\ybar} = \ohull\{\xbar,\ybar\}.%
\indexg{segments, astral|)}%
\indexg{segments, astral!defined|)}%
\]
\end{definition}

Here are some properties of outer convex hull:

\begin{proposition}
\label{pr:ohull:hull}
Outer convex hull is a hull operator. Moreover, for any $S\subseteq\eRn$, the
set $\ohull S$ is closed (in $\eRn$).
Consequently, $\ohull\Sbar=\ohull S$.
\end{proposition}

\begin{proof}
Let $\calC$ be the collection consisting of all possible intersections
of closed astral halfspaces
in $\eRn$
(including the empty intersection, which is all of $\eRn$).
Then $\ohull$ is the hull operator for $\calC$
(as defined in \Cref{sec:prelim:hull-ops}).

Let $S\subseteq\extspace$.
Then $\ohull S$ is an intersection of closed astral halfspaces and is
therefore closed.
As such, $S\subseteq\Sbar\subseteq\ohull S$.
Therefore,
since $\ohull$ is a hull operator, 
$\ohull\Sbar=\ohull S$ by
\Cref{pr:gen-hull-ops}(\ref{pr:gen-hull-ops:d}).
\end{proof}

The outer hull of any set $S\subseteq\extspace$ can be
characterized in a way that is often more useful, as shown next:

\begin{proposition}  \label{pr:ohull-simplify}
  Let $S\subseteq\extspace$, and let $\zbar\in\extspace$.
  Then the following are equivalent:
  \begin{letter-compact}
  \item  \label{pr:ohull-simplify:a}
    $\zbar\in\ohull{S}$.
  \item  \label{pr:ohull-simplify:b}
    For all $\uu\in\Rn$,
    \begin{equation}   \label{eq:pr:ohull-simplify:1}
       \zbar\cdot\uu
       \leq
       \sup_{\xbar\in S} \xbar\cdot\uu.
    \end{equation}
  \item  \label{pr:ohull-simplify:c}
    For all $\uu\in\Rn$,
    \begin{equation}   \label{eq:pr:ohull-simplify:2}
       \inf_{\xbar\in S} \xbar\cdot\uu
       \leq
       \zbar\cdot\uu
       \leq
       \sup_{\xbar\in S} \xbar\cdot\uu.
    \end{equation}
  \end{letter-compact}
\end{proposition}

\begin{proof}
  ~

\begin{proof-parts}
\pfpart{%
  (\ref{pr:ohull-simplify:a})
  $\Rightarrow$
  (\ref{pr:ohull-simplify:b}):
}
Suppose $\zbar\in\ohull{S}$,
which implies $S$ is not empty
(since $\ohull\emptyset=\emptyset$).
Let $\uu\in\Rn$, and
let $\gamma$ be equal to the right-hand side of
\eqref{eq:pr:ohull-simplify:1};
we aim to show that $\zbar\cdot\uu\leq\gamma$.
We assume $\uu\neq\zero$ and that $\gamma<+\infty$
since otherwise the claim is immediate.
Let $\beta\in\R$ be such that $\beta\geq\gamma$.
Then $\sup_{\xbar\in S} \xbar\cdot\uu = \gamma\leq\beta$,
implying that $S\subseteq\chsua$, and so that
$\zbar\in\chsua$ by definition of outer hull
(Eq.~\ref{eqn:ohull-defn}).
Therefore, $\zbar\cdot\uu\leq\beta$.
Since this holds for all $\beta\geq\gamma$, it follows that
$\zbar\cdot\uu\leq\gamma$.

\pfpart{%
  (\ref{pr:ohull-simplify:b})
  $\Rightarrow$
  (\ref{pr:ohull-simplify:a}):
}
Assume \eqref{eq:pr:ohull-simplify:1} holds for all $\uu\in\Rn$.
Suppose $S\subseteq\chsua$ for some $\uu\in\Rn\wo\{\zero\}$ and
$\beta\in\R$.
This means $\xbar\cdot\uu\leq\beta$ for all $\xbar\in S$,
so, combined with
\eqref{eq:pr:ohull-simplify:1},
$\zbar\cdot\uu\leq\sup_{\xbar\in S}\xbar\cdot\uu\leq\beta$.
Therefore, $\zbar\in\chsua$.
Since this holds for all closed astral halfspaces that include $S$,
$\zbar$ must be in $\ohull{S}$, by its definition.

\pfpart{%
  (\ref{pr:ohull-simplify:c})
  $\Rightarrow$
  (\ref{pr:ohull-simplify:b}):
}
This is immediate.

\pfpart{%
  (\ref{pr:ohull-simplify:b})
  $\Rightarrow$
  (\ref{pr:ohull-simplify:c}):
}
\eqref{eq:pr:ohull-simplify:1}
immediately implies the second inequality of
\eqref{eq:pr:ohull-simplify:2},
and also implies the first inequality
by substituting $\uu$ with $-\uu$.
\qedhere
\end{proof-parts}
\end{proof}

\indexg{segments, astral|(}%
As an immediate consequence, we obtain the following characterization
of segments, which, by definition, are outer hulls of pairs of points:

\begin{proposition}  \label{pr:seg-simplify}
  Let $\xbar,\ybar,\zbar\in\extspace$.
  Then the following are equivalent:
  \begin{letter-compact}
  \item  \label{pr:seg-simplify:a}
    $\zbar\in\lb{\xbar}{\ybar}$.
  \item  \label{pr:seg-simplify:b}
    For all $\uu\in\Rn$,
    $\zbar\cdot\uu \leq \max\{\xbar\cdot\uu,\, \ybar\cdot\uu\}$.
  \item  \label{pr:seg-simplify:c}
    For all $\uu\in\Rn$,
    $
      \min\{\xbar\cdot\uu,\, \ybar\cdot\uu\}
      \leq
      \zbar\cdot\uu
      \leq
      \max\{\xbar\cdot\uu,\, \ybar\cdot\uu\}
    $.
  \end{letter-compact}
\end{proposition}

For example,
segments in $\Rext$ are closed intervals, as expected:

\begin{example}[Segments in $\Rext$]  \label{ex:segs-in-Rext}
Let $\barx,\bary\in\Rext$ with $\barx\leq\bary$.
Then $\lb{\barx}{\bary}=[\barx,\bary]$.
This is because, by
\Crefequiv{pr:seg-simplify}{pr:seg-simplify:a}{pr:seg-simplify:c},
a point $\barz\in\Rext$ is in $\lb{\barx}{\bary}$ if and only if
$\min\{\barx u,\, \bary u\} \leq \barz u \leq \max\{\barx u,\, \bary u\}$
for all $u\in\R$.
By considering separately the cases that $u$ is positive, negative or
zero, it follows that this will hold if and only if
$\barx\leq\barz\leq\bary$.%
\indexg{segments, astral|)}%
\end{example}

\indexg{support functions, astral|(}%
\indexg{indicator functions (astral)!conjugate and biconjugate of|(}%
Analogous to the standard support function
$\indstars$ given in \eqref{eq:e:2},
the \emph{astral support function} for a set $S\subseteq\extspace$ is
the conjugate $\indaSstar$ of the astral indicator function
$\indaS$ (defined in Eq.~\ref{eq:indfa-defn}).
From \eqref{eq:Fstar-down-def}, this is the function
\begin{equation}  \label{eqn:astral-support-fcn-def}
  \indaSstar(\uu)
  =
  \sup_{\xbar\in S} \xbar\cdot\uu,
\end{equation}
for
\indexg{support functions, astral|)}%
$\uu\in\Rn$.
\indexg{support functions, astral!dual conjugate of|(}%
\indexg{outer convex hull!biconjugate of indicator@as biconjugate of indicator|(}%
Note that this function is the same as
the expression that appears on the right-hand side
of \eqref{eq:pr:ohull-simplify:1}.
As a result,
the outer convex hull of $S$ is directly linked to
the biconjugate $\indfadub{S}$:

\begin{theorem}  \label{thm:ohull-biconj}
  Let $S\subseteq\extspace$.
  Then
  $
    \indfadub{S}
    =
    \indfa{\ohull S}
  $.
\end{theorem}

\begin{proof}
Let $\zbar\in\eRn$.
Then we have
\begin{align}
  \zbar\in\ohull S
  &\;\Leftrightarrow\;
  \forall\uu\in\Rn,\,\zbar\inprod\uu\le\indfastar{S}(\uu)
  \notag
\\
  &\;\Leftrightarrow\;
  \forall\uu\in\Rn,\,-\indfastar{S}(\uu)\plusd\zbar\inprod\uu\le 0
  \notag
\\
  &\;\Leftrightarrow\;
  \indfadub{S}(\zbar)\le 0.
  \label{eq:thm:ohull-biconj:1}
\end{align}
The first equivalence is by
\Cref{pr:ohull-simplify}(\ref{pr:ohull-simplify:a},\ref{pr:ohull-simplify:b})
and
\eqref{eqn:astral-support-fcn-def}.
The second equivalence is by \Cref{pr:plusd-props}(\ref{pr:plusd-props:e}).
The last equivalence is because
\begin{equation}  \label{eq:thm:ohull-biconj:2}
  \indfadub{S}(\zbar)
  =
  \sup_{\uu\in\Rn}\bigBracks{-\indfastar{S}(\uu)\plusd\zbar\inprod\uu}
\end{equation}
by definition of the dual conjugate
(Eq.~\ref{eq:psistar-def:2}).

We claim that $\indfadub{S}(\zbar)\in\set{0,+\infty}$, which,
in light of \eqref{eq:thm:ohull-biconj:1},
will suffice to prove the proposition.

If $S=\emptyset$ then
$\indaS\equiv+\infty$ implying $\indaSstar\equiv-\infty$ and
$\indfadub{S}\equiv+\infty$, so the claim holds.
Otherwise, $\indfastar{S}(\zero)=0$, so
$\indfadub{S}(\zbar)\ge-\indfastar{S}(\zero)\plusd\zbar\inprod\zero=0$. If
$\indfadub{S}(\zbar)=0$ then the claim holds, so in the rest of the
proof we consider the case $\indfadub{S}(\zbar)>0$.

Since $\indfadub{S}(\zbar)>0$,
\eqref{eq:thm:ohull-biconj:2}
implies there exists $\uhat\in\Rn$ such that
$-\indfastar{S}(\uhat)\plusd\zbar\inprod\uhat>0$.
For any $\lambda\in\Rstrictpos$,
\[
  \indfastar{S}(\lambda\uhat)=
  \sup_{\xbar\in S}\bigBracks{\xbar\cdot(\lambda\uhat)}=
  \lambda\sup_{\xbar\in S}\bigBracks{\xbar\cdot\uhat}=
  \lambda\indfastar{S}(\uhat),
\]
so
\begin{align*}
  \indfadub{S}(\zbar)
  &\ge
  \bigBracks{-\indfastar{S}(\lambda\uhat)\plusd\zbar\inprod(\lambda\uhat)}
  =
  \lambda\bigBracks{-\indfastar{S}(\uhat)\plusd\zbar\inprod\uhat}.
\end{align*}
Taking $\lambda\to+\infty$, we obtain $\indfadub{S}(\zbar)=+\infty$,
completing the proof.%
\indexg{outer convex hull|)}%
\indexg{indicator functions (astral)!conjugate and biconjugate of|)}%
\indexg{support functions, astral!dual conjugate of|)}%
\indexg{outer convex hull!biconjugate of indicator@as biconjugate of indicator|)}%
\end{proof}

\indexg{segments, astral!examples|(}%
To get a sense of what astral segments are like,
we give a few examples in $\eRf{2}$, thereby also demonstrating how
\Cref{pr:seg-simplify} can be used for this purpose.
Later, we will develop general tools
(such as
Theorems~\ref{thm:lb-with-zero} and~\ref{thm:conv-lmset-char})
which can be applied more
directly for these same examples,
greatly simplifying calculations.

\begin{example}
\label{ex:seg-zero-e1}
First, consider $S=\lb{\zero}{\ee_1}$. We instantiate
\Cref{pr:seg-simplify}(\ref{pr:seg-simplify:c})
with values $\uu$ of the form $\transk{[1,\alpha]}$ and $\transk{[0,1]}$, for $\alpha\in\R$.
Then any $\zbar\in S$ must satisfy
\begin{gather*}
  0
  =
  \min\set{0,\,\ee_1\inprod\trans{[1,\alpha]}}
  \,\le\,
  \zbar\inprod\trans{[1,\alpha]}\!
  \le\,\max\set{0,\,\ee_1\inprod\trans{[1,\alpha]}}=1,
\\
  0
  =
  \min\set{0,\,\ee_1\inprod\trans{[0,1]}}
  \,\le\,
  \zbar\inprod\trans{[0,1]}\!
  \le\,\max\set{0,\,\ee_1\inprod\trans{[0,1]}}=0.
\end{gather*}
Every vector $\uu\in\R^2$ is a scalar multiple of one of the ones used
above.
Consequently, every constraint as given in
\Cref{pr:seg-simplify}(\ref{pr:seg-simplify:c})
must be equivalent to one of the constraints above.
Therefore, $\zbar\in S$ if and only if it satisfies these two
constraints. We summarize these as
\[
\textup{(C1)}\quad
  \zbar\inprod\trans{[1,\alpha]}\!\!\in[0,1] \text{ for all }\alpha\in\R,
\qquad
\textup{(C2)}\quad
  \zbar\inprod\trans{[0,1]}\!\!=0.
\]
By (C2), $\zbar$ must be orthogonal to $\ee_2$ (that is, if $\zbar=\VV\omm\plusl\qq$, then $\VV\perp\ee_2$ and $\qq\perp\ee_2$),
and so must take the form $\zbar=\beta\ee_1$ for some $\beta\in\eR$.
By (C1), we in fact must have $\beta\in[0,1]$. Thus,
the set of points $\zbar$ that satisfy (C1) and (C2) is the standard line segment connecting $\zero$ and $\ee_1$:
\[
  \lb{\zero}{\ee_1} = \set{\beta\ee_1:\:\beta\in[0,1]}.
\qedhere
\]
\end{example}

\begin{example}
\label{ex:seg-zero-oe2}
Now consider $S=\lb{\zero}{\limray{\ee_2}}$, and instantiate
\Cref{pr:seg-simplify}(\ref{pr:seg-simplify:c})
with $\uu$ of the form $\transk{[\alpha,1]}$ and $\transk{[1,0]}$, for $\alpha\in\R$. Scalar multiples of these vectors include all of~$\R^2$,
so, by similar reasoning as in the last example,
$\zbar\in S$ if and only if it satisfies
\[
\textup{(C1)}\quad
  \zbar\inprod\trans{[\alpha,1]}\!\!\in[0,+\infty] \text{ for all }\alpha\in\R,
\qquad
\textup{(C2)}\quad
  \zbar\inprod\trans{[1,0]}\!\!=0.
\]
Similar to the last example, (C2) implies that $\zbar=\beta\ee_2$ for some $\beta\in\eR$, and
(C1) implies that $\beta\in[0,+\infty]$. Thus, the set of points that satisfy
(C1) and (C2)
is the astral closure of the ray in the direction of $\ee_2$:
\[
  \lb{\zero}{\limray{\ee_2}} = \set{\beta\ee_2:\:\beta\in[0,+\infty]}.
\qedhere
\]
\end{example}

\begin{example}
\label{ex:seg-zero-oe2-plus-e1}
Next, consider $S=\lb{\zero}{\limray{\ee_2}\plusl\ee_1}$. As in the previous example, instantiate
\Cref{pr:seg-simplify}(\ref{pr:seg-simplify:c})
with $\uu$ of the form $\transk{[\alpha,1]}$ and $\transk{[1,0]}$, for $\alpha\in\R$. Scalar multiples of these vectors include all of~$\R^2$,
so $\zbar\in S$ if and only if it satisfies
\[
\textup{(C1)}\quad
  \zbar\inprod\trans{[\alpha,1]}\!\!\in[0,+\infty] \text{ for all }\alpha\in\R,
\qquad
\textup{(C2)}\quad
  \zbar\inprod\trans{[1,0]}\!\!\in[0,1].
\]
This means that if $\zbar$ is infinite, its dominant direction must be $\ee_2$,
and its projection orthogonal to $\ee_2$ must be of the form $\beta\ee_1$ with $\beta\in[0,1]$. On the other hand,
if $\zbar$ is finite, say $\zbar=\trans{[z_1,z_2]}\negKern$, then (C1) implies
that $\alpha z_1+z_2\ge 0$ for all $\alpha\in\R$, which is only possible if $z_1=0$ and $z_2\ge 0$. Thus, the set of points that satisfy (C1) and (C2) is
\[
  \lb{\zero}{\limray{\ee_2}\plusl\ee_1} = \set{\beta\ee_2:\:\beta\in[0,+\infty)}
  \,\cup\,
  \set{\omega\ee_2\plusl\beta\ee_1:\:\beta\in[0,1]}.
\qedhere
\]
\end{example}

\begin{example}
\label{ex:seg-oe1-oe2}
The final example in $\eRf{2}$ is
$S=\lb{\limray{\ee_1}}{\limray{\ee_2}}$. We instantiate
\Cref{pr:seg-simplify}(\ref{pr:seg-simplify:c})
with $\uu$ of the form
  $\transk{[1,\alpha]}$,
  $\transk{[1,-\alpha]}$,
  $\transk{[1,0]}$,
  and
  $\transk{[0,1]}$,
where $\alpha\in\R_{>0}$. Scalar multiples of these vectors include all of~$\R^2$,
so $\zbar\in S$ if and only if it satisfies
\begin{gather*}
\begin{aligned}
&
\textup{(C1)}\quad
  \zbar\inprod\trans{[1,\alpha]}\!\!=+\infty
  \text{ for all }\alpha\in\R_{>0},
\\&
\textup{(C2)}\quad
  \zbar\inprod\trans{[1,-\alpha]}\!\!\in[-\infty,+\infty]
  \text{ for all }\alpha\in\R_{>0},
\end{aligned}
\\
\textup{(C3)}\quad
  \zbar\inprod\trans{[1,0]}\!\!\in[0,+\infty],
\qquad
\textup{(C4)}\quad
  \zbar\inprod\trans{[0,1]}\!\!\in[0,+\infty].
\end{gather*}
By (C1), $\zbar$ must be infinite.
Let $\vv=\trans{[v_1,v_2]}$ be its dominant direction
(implying $\vv\neq\zero$).
By~(C1), $v_1+\alpha v_2\ge 0$ for all $\alpha\in\R_{>0}$, so $v_1,v_2\ge 0$.

If $v_1,v_2>0$, then $\zbar=\limray{\vv}\plusl\zbar'$ satisfies (C1--C4) for any $\zbar'\in\eRf{2}$.
By the Projection Lemma~\ref{lemma:proj},
all such points $\zbar$ can be expressed as $\limray{\vv}\plusl\beta\ww$,
where $\beta\in\eR$ and $\ww=\transk{[v_2,-v_1]}$ is the vector that spans the linear space orthogonal to $\vv$.

If $v_1=0$ then $\zbar$ can be written as $\limray{\ee_2}\plusl\beta\ee_1$ with $\beta\in\eR$. By (C3), we must actually have $\beta\in[0,+\infty]$. Symmetrically, if $v_2=0$ then $\zbar=\limray{\ee_1}\plusl\beta\ee_2$ where $\beta\in[0,+\infty]$.

Altogether, we thus obtain
\begin{align}
  \lb{\limray{\ee_1}}{\limray{\ee_2}}
  &=
  \bigSet{\limray{\vv}\plusl\beta\ww:\:
          \vv=\trans{[v_1,v_2]}\!\!\in\R^2_{>0},\,
          \ww=\trans{[v_2,-v_1]}\negKern,\,
          \beta\in\eR}
\nonumber
\\
  &\qquad{}
   \cup\,
   \bigSet{\limray{\ee_1}\plusl\beta\ee_2:\:\beta\in[0,+\infty]}
\nonumber
\\
  &\qquad{}
   \cup\,
   \bigSet{\limray{\ee_2}\plusl\beta\ee_1:\:\beta\in[0,+\infty]}.
\label{eq:ex:seg-oe1-oe2:1}
\end{align}
This set can be visualized as the concatenation of closed galaxies associated with astrons in the directions $\vv=\trans{[\cos\varphi,\sin\varphi]}\negKern$ with the polar angle $\varphi\in(0,\pi/2)$. The two ``boundary astrons'' $\limray{\ee_1}$ and $\limray{\ee_2}$ ($\varphi=0$ and $\varphi=\pi/2$) only contribute a closed ``half-galaxy'' that attaches to the other included galaxies.
\end{example}

\begin{example}
\label{ex:seg-negiden-iden}
We next consider a less intuitive example in $\eRn$.
Let $\Iden$ be the $n\times n$ identity matrix.
Then the segment joining the points
$-\Iden \omm$ and $\Iden \omm$
turns out to be all of~$\extspace$;
that is,
$  \lb{-\Iden \omm}{\Iden \omm}=\extspace  $.
To see this, let $\zbar\in\extspace$ and let $\uu\in\Rn$.
If $\uu\neq\zero$, then
$\Iden\omm\cdot\uu\in\{-\infty,+\infty\}$
by \Cref{pr:vtransu-zero},
so
\[
  \zbar\cdot\uu\leq +\infty=\max\{-\Iden\omm\cdot\uu,\,\Iden\omm\cdot\uu\}.
\]
Otherwise,
if $\uu=\zero$, then
$\zbar\cdot\zero=0=\max\{-\Iden\omm\cdot\zero,\,\Iden\omm\cdot\zero\}$.
Thus, in all cases, the inequality appearing in
\Cref{pr:seg-simplify}(\ref{pr:seg-simplify:b})
is satisfied, so
$\zbar\in\lb{-\Iden \omm}{\Iden \omm}$.
\end{example}

This last example shows, in an extreme way,
that the ``segment'' joining two infinite astral
points can be very different from the standard, one-dimensional
segment joining points in $\Rn$.%
\indexg{segments, astral!examples|)}

\indexg{convex sets, astral|(}%
Continuing our general development,
we are now ready to define convexity for subsets of astral space.
Having defined segments, this definition is entirely analogous to the
standard one for subsets of $\Rn$:
\begin{definition}
\indexg{convex sets, astral!defined|(}%
  A set $S\subseteq\extspace$ is
  \emph{astrally convex} (or simply \emph{convex})
  if the segment joining every pair of points in $S$ is
  also included in $S$, that is, if
  \[
    \lb{\xbar}{\ybar} \subseteq S
    \quad
    \text{for all $\xbar,\ybar\in S$.}%
\indexg{convex sets, astral!defined|)}%
  \]
\end{definition}

For example, as expected,
the convex sets in $\Rext$ are all intervals, whether open, closed, or
half-open:

\begin{example}[Convex sets in $\Rext$]  \label{ex:cvx-sets-in-Rext}
Let $S\subseteq\Rext$ be nonempty.
Then $S$ is convex if and only if $S$ has the form
$[\alpha,\beta]$,
$(\alpha,\beta)$,
$[\alpha,\beta)$,
or
$(\alpha,\beta]$,
for some $\alpha,\beta\in\Rext$.
That sets of these forms are convex follows straightforwardly from
Example~\ref{ex:segs-in-Rext} and definition of convexity.
For the converse, suppose $S$ is convex, and let
$\alpha=\inf{S}$ and $\beta=\sup{S}$.
Suppose $\barz\in(\alpha,\beta)$.
Then there exists $\barx,\bary\in S$ such that
$\alpha\leq\barx<\barz<\bary\leq\beta$, implying
$\barz\in[\barx,\bary]=\lb{\barx}{\bary}\subseteq S$,
again by Example~\ref{ex:segs-in-Rext} and definition of convexity.
Thus, $(\alpha,\beta)\subseteq S\subseteq[\alpha,\beta]$,
proving the claim.

This also shows that the closure of any convex set $S$ in $\Rext$
is a closed interval, namely,
$\Sbar=[\inf{S},\sup{S}]$.
\end{example}

We next show that several important examples of sets are astrally convex, and derive some fundamental
properties of astral convex sets.

One of these concerns {updirected} families of sets, where an indexed
collection of sets $\set{S_i:\:i\in\indset}$ is said to be
\indexg{updirected (family of sets)|(}%
\emph{updirected} if for any $i,j\in\indset$ there exists
$k\in\indset$ such that $S_i\cup S_j\subseteq S_k$.
The union of an updirected family of convex sets in $\Rn$ is convex.
We show the same is true for an updirected family of astrally convex
sets.
We also show that intersections of arbitrary collections of astrally
convex sets are astrally convex, and that the empty set and the entire space $\eRn$ are astrally convex. Altogether, these properties, stated as \Cref{pr:e1}(\ref{pr:e1:univ},\ref{pr:e1:b},\ref{pr:e1:union}) below, imply that astral convexity is an instance of abstract convexity as defined by
\idxvandevel\citet{vandeVel}.

\begin{proposition}  \label{pr:e1}
~
  \begin{letter-compact}
  \item  \label{pr:e1:a}
\indexg{convex sets (standard)!astral convexity and|(}%
    If $\xx,\yy\in\Rn$, then $\lb{\xx}{\yy}$ is the standard line segment
    joining $\xx$ and $\yy$ in $\Rn$:
    \[
       \lb{\xx}{\yy} = \bigBraces{(1-\lambda)\xx + \lambda\yy :\:
                                   \lambda\in [0,1]}.
    \]
    Therefore, a subset of $\Rn$ is astrally convex if
    and only if it is convex in the standard sense.\looseness=-1
  \item  \label{pr:e1:univ}
    The empty set $\emptyset$ and the full astral space $\eRn$ are astrally convex.
  \item  \label{pr:e1:b}
    The intersection of an arbitrary collection of astrally convex
    sets is astrally convex.
  \item  \label{pr:e1:int-rn}
    If $S\subseteq\extspace$ is astrally convex, then so is
    $S\cap\Rn$.
  \item  \label{pr:e1:union}
    Let $\set{S_i:\:i\in\indset}$ be an {updirected} family of
    astrally convex sets $S_i\subseteq\eRn$
    (where $\indset$ is any index set).
    Then their union $\bigcup_{i\in\indset} S_i$ is also astrally convex.
  \item  \label{pr:e1:c}
\indexg{halfspaces, astral!convexity of|(}%
\indexg{hyperplanes, astral!convexity of|(}%
    Every set of the form given in
    Eqs.~(\ref{eq:chsua-defn}),~(\ref{eqn:open-hfspace-defn})
    or~(\ref{eqn:hyperplane-defn}), with $\uu\in\Rn$ and
    $\beta\in\Rext$, is astrally convex.
    Therefore,
    every astral hyperplane and every closed or open astral halfspace is
    astrally convex.
  \item  \label{pr:e1:ohull}
\indexg{outer convex hull!convexity of|(}%
\indexg{segments, astral!convexity of|(}%
    For any $S\subseteq\extspace$,
    the outer hull, $\ohull S$, is astrally convex.
    Also, for all $\xbar\in\extspace$,
    $\ohull \{\xbar\} = \lb{\xbar}{\xbar} = \{\xbar\}$.
    Therefore, all segments and singletons are astrally convex.

  \item  \label{pr:e1:d}
    Every base element of the astral topology
    given in \eqref{eq:h:3a}
    is astrally convex.

  \item  \label{pr:e1:base}
\indexg{countable neighborhood base!convexity of|(}%
    Every point $\xbar\in\extspace$ has a nested countable neighborhood base consisting
    of astrally convex sets.%
\indexg{countable neighborhood base!convexity of|)}
  \end{letter-compact}
\end{proposition}

\begin{proof}
~
\begin{proof-parts}
\pfpart{Part~(\ref{pr:e1:a}):}
Suppose $\xx,\yy\in\Rn$.
If $\zz\in\lb{\xx}{\yy}$,
then for all $\uu\in\Rn$, by
\Cref{pr:seg-simplify}(\ref{pr:seg-simplify:a},\ref{pr:seg-simplify:c}),
\[
  -\infty
  <
  \min\{\xx\cdot\uu,\,\yy\cdot\uu\}
  \leq
  \zz\cdot\uu
  \leq
  \max\{\xx\cdot\uu,\,\yy\cdot\uu\}
  <
  +\infty,
\]
implying that $\zz\in\Rn$
(by \Cref{pr:i:3}\ref{i:3b}\ref{i:3a}).
Thus, $\lb{\xx}{\yy}\subseteq\Rn$, so
$\lb{\xx}{\yy}$ is included in some closed astral halfspace $\chsua$ in $\extspace$
if and only if it is included in the corresponding
(standard) closed halfspace in $\Rn$
given in \eqref{eqn:chsua-cap-rn}.
It follows that $\lb{\xx}{\yy}$ is the intersection of all closed
halfspaces in $\Rn$ containing both $\xx$ and $\yy$.
Therefore,
\[
   \lb{\xx}{\yy}
   =
   \cl\bigParens{\conv\{\xx,\yy\}}
   =
   \conv\{\xx,\yy\}
   =
   \bigBraces{
     (1-\lambda)\xx + \lambda\yy :\:
     \lambda\in [0,1]
   },
\]
where the first equality is
by \Cref{pr:con-int-halfspaces}(\ref{roc:cor11.5.1}),
the second because $\conv\{\xx,\yy\}$ is closed
by
\Cref{roc:thm19.1}(\ref{roc:thm19.1:b},\ref{roc:thm19.1:c}),
and the third by \Cref{roc:thm2.3}.

This implies that, when restricted to subsets of $\Rn$, astral
convexity coincides exactly with the
standard definition of convexity in $\Rn$.%
\indexg{convex sets (standard)!astral convexity and|)}

\pfpart{Part~(\ref{pr:e1:univ}):}
Immediate from the definition of astral convexity.

\pfpart{Part~(\ref{pr:e1:b}):}
Let
\[
  M
  = \bigcap_{i\in\indset} S_i
\]
where each $S_i\subseteq\extspace$ is convex, and $\indset$ is an
arbitrary index set.
Let $\xbar$, $\ybar$ be in $M$.
Then for all $i\in\indset$, $\xbar,\ybar\in S_i$, so
$\lb{\xbar}{\ybar}\subseteq S_i$.
Since this holds for all $i\in\indset$,
$\lb{\xbar}{\ybar}\subseteq M$, and $M$ is convex.

\pfpart{Part~(\ref{pr:e1:int-rn}):}
Suppose $S\subseteq\extspace$ is astrally convex.
By part~(\ref{pr:e1:a}), $\Rn$ is astrally convex, being convex in the
standard sense.
Therefore, $S\cap\Rn$ is as well by
part~(\ref{pr:e1:b}).

\pfpart{Part~(\ref{pr:e1:union}):}
Let $\set{S_i:\:i\in\indset}$ be an updirected collection of astrally convex sets in $\eRn$, and let
$M=\bigcup_{i\in\indset} S_i$.
We will show that $M$ is convex.
Let $\xbar$, $\ybar\in M$. Then there exist $i,j\in\indset$ such that
$\xbar\in S_i$ and $\ybar\in S_j$, and so also $k\in\indset$ such that
$\xbar,\ybar\in S_k$. By astral convexity of $S_k$,
$\lb{\xbar}{\ybar}\subseteq S_k\subseteq M$, so $M$ is convex.%
\indexg{updirected (family of sets)|)}

\pfpart{Part~(\ref{pr:e1:c}):}
Let $\uu\in\Rn$ and $\beta\in\Rext$, and let
$H=\set{\xbar\in\extspace :\: \xbar\cdot\uu \leq \beta}$,
which we aim to show is convex.
Suppose $\xbar,\ybar\in H$, and
that $\zbar\in\lb{\xbar}{\ybar}$.
Then
\[
  \zbar\cdot\uu \leq \max\{\xbar\cdot\uu,\,\ybar\cdot\uu\} \leq \beta,
\]
where
the first inequality is from
\Cref{pr:seg-simplify}(\ref{pr:seg-simplify:a},\ref{pr:seg-simplify:b}),
and the
second is because $\xbar,\ybar\in H$.
Thus, $\zbar\in H$, so $H$ is convex.
Therefore, sets of the form given in \eqref{eq:chsua-defn}, including
closed astral halfspaces, are convex.

The proof is the same for sets as given in
\eqref{eqn:open-hfspace-defn},
but using strict inequalities; these
include open astral halfspaces.
The set of points for which
$\xbar\cdot\uu=\beta$, as in \eqref{eqn:hyperplane-defn},
is convex by the above and
part~(\ref{pr:e1:b}), since this set is the
intersection of the two convex sets
defined by $\xbar\cdot\uu\le\beta$ and $\xbar\cdot\uu\ge\beta$,
respectively.
Thus, astral hyperplanes are also convex.%
\indexg{halfspaces, astral!convexity of|)}%
\indexg{hyperplanes, astral!convexity of|)}

\pfpart{Part~(\ref{pr:e1:ohull}):}
An outer hull is an intersection of closed
astral halfspaces and therefore convex by parts~(\ref{pr:e1:b}) and~(\ref{pr:e1:c}).
Segments, as outer hulls, are thus convex. It remains to argue that $\ohull \{\xbar\} = \{\xbar\}$
for all $\xbar\in\eRn$, which will also establish that singletons are convex.

Let $\xbar,\zbar\in\extspace$.
Then by
\Cref{pr:ohull-simplify}(\ref{pr:ohull-simplify:a},\ref{pr:ohull-simplify:c}),
$\zbar\in\ohull\{\xbar\}$ if and only if, for all $\uu\in\Rn$,
$\xbar\cdot\uu\leq\zbar\cdot\uu\leq\xbar\cdot\uu$,
which,
by \Cref{pr:i:4}, holds if and only if $\zbar=\xbar$.
Therefore, $\ohull\{\xbar\}=\{\xbar\}$.%
\indexg{outer convex hull!convexity of|)}%
\indexg{segments, astral!convexity of|)}

\pfpart{Part~(\ref{pr:e1:d}):}
Since base elements, as in \eqref{eq:h:3a},
are intersections of open astral halfspaces,
they are convex by parts~(\ref{pr:e1:b}) and~(\ref{pr:e1:c}).

\pfpart{Part~(\ref{pr:e1:base}):}
\indexg{countable neighborhood base!convexity of|(}%
Elements of the nested countable neighborhood base from
\Cref{thm:first:local} are also base elements in the astral topology
(of the form given in Eq.~\ref{eq:h:3a}),
and are therefore convex by part~(\ref{pr:e1:d}).%
\indexg{countable neighborhood base!convexity of|)}
\qedhere
\end{proof-parts}
\end{proof}

Thanks to \Cref{pr:e1}(\ref{pr:e1:a}), there is no ambiguity
in referring to an astrally convex set $S\subseteq\eRn$ simply
as convex; we use this terminology from now on.

\indexg{closure, astral!convexity of|(}%
\indexg{outer convex hull!sets in $\Rn$@for sets in $\Rn$|(}%
We next show that if $S$ is any convex subset of $\Rn$, then its
astral closure $\Sbar$ in $\extspace$ is also convex, and more specifically, is
exactly equal to the outer hull of $S$, the intersection of all
closed astral halfspaces that contain $S$.
The key step of the proof is in applying the results of
\Cref{sec:ent-closed-fcn} to $\inds$, the indicator function
for $S$, to show that $\indsext = \indsdub$,
thereby connecting $\Sbar$ and $\ohull S$
via \Cref{pr:inds-ext} and \Cref{thm:ohull-biconj}.

\begin{theorem}  \label{thm:e:6}
  Let $S\subseteq\Rn$ be convex.
  Then $\Sbar$, its closure in $\extspace$,
  is exactly equal to its outer hull;
  that is, $\Sbar=\ohull S$.
  Consequently, $\Sbar$ is convex.
\end{theorem}

\begin{proof}
Let $S\subseteq\Rn$ be convex, and let $\inds$ be its indicator function (see Eq.~\ref{eq:indf-defn}). Then
\[
  \indfa{\Sbar}
  =
  \indsext
  =
  \indsdub
  =
  \indsextdub
  =
  \indfadub{\Sbar}
  =
  \indfa{\ohull\Sbar}
  =
  \indfa{\ohull S},
\]
where the equalities follow, in order, from:
\Cref{pr:inds-ext};
\Cref{cor:all-red-closed-sp-cases}
(since $\inds$ is convex and $\inds\geq 0$);
\Cref{pr:fextstar-is-fstar};
\Cref{pr:inds-ext} again;
\Cref{thm:ohull-biconj};
and finally \Cref{pr:ohull:hull}.
Thus, $\Sbar=\ohull S$, so the convexity of $\Sbar$ follows by \Cref{pr:e1}(\ref{pr:e1:ohull}).
\end{proof}

In general,
\Cref{thm:e:6} does not hold for arbitrary convex sets in
$\extspace$ (rather than in $\Rn$), as will be demonstrated in
\Cref{sec:closure-convex-set}.%
\indexg{outer convex hull!sets in $\Rn$@for sets in $\Rn$|)}%
\indexg{closure, astral!convexity of|)}%
\indexg{convex sets, astral|)}

\section{Astral polytopes and polyhedral sets}
\label{sec:simplices}

\indexg{polytopes, astral|(}%
The outer convex hull, $\ohull V$, of any finite set
$V\subseteq\extspace$
is called the \emph{(astral) polytope formed by $V$}.
\indexg{polytopes, astral!sequential characterization|(}%
We next derive an alternative characterization of such sets.

To motivate our characterization, we draw an analogy with standard polytopes in~$\Rn$.
A standard polytope is the convex hull of a finite set of points in $\Rn$,
and can be characterized in two ways:
either as the intersection of all halfspaces containing the points, or
as the set of all convex combinations of the points.
Thus, there is both an ``outside'' and ``inside'' way of
characterizing convexity in this case.
The astral outer hull
has so far been described
in ``outside'' terms, as the intersection of closed astral halfspaces.
We next give an alternative ``inside'' description of this same set.

Specifically, we show that
$\ohull V$ can be characterized
in terms of sequences via a formulation
which says that a point is in the outer hull of $V$
if and only if it is the limit of points in $\Rn$ that are
themselves convex combinations of points converging to the
points in $V$.
(See the illustration in \Cref{fig:thm:e:7}.)
More precisely:

\begin{figure}
  \centering
  \includegraphics{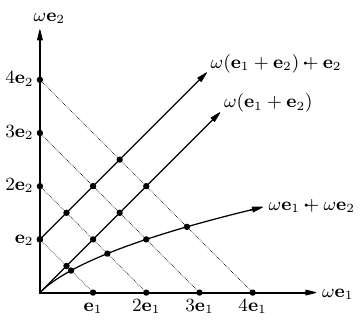}
  \mycaption{Sequential characterization of astral segments}{%
\indexf{segments, astral!sequential characterization}%
    By \Cref{thm:e:7} (and \Cref{cor:e:1}), the segment $\seg(\limray{\ee_1},\limray{\ee_2})$
    consists of points $\zbar$ that can be obtained as limits of convex combinations
    $\zz_t=(1-\lambda_t)\xx_t+\lambda_t\yy_t$ for some sequences $\seq{\xx_t}$
    and~$\seq{\yy_t}$ that converge, respectively, to $\limray{\ee_1}$ and $\limray{\ee_2}$.
    In the figure, we consider $\xx_t=t\ee_1$ and $\yy_t=t\ee_2$, and depict several
    choices of their convex combinations converging to three different points~$\zbar$,
    all of which are thus points on $\seg(\limray{\ee_1},\limray{\ee_2})$. 
    (The depicted convex combinations, from top to bottom, are obtained by setting
    $\lambda_t=1/2 + 1/(2t)$, $\lambda_t=1/2$, and
    $\lambda_t=(\sqrt{t+1}-1)/t$.)}
  \label{fig:thm:e:7}%
\end{figure}

\begin{theorem}  \label{thm:e:7}
  Let $V=\{\xbar_1,\dotsc,\xbar_m\}\subseteq\extspace$,
  and let $\zbar\in\extspace$.
  Then $\zbar\in\ohull{V}$ if
  and only if there exist sequences $\seq{\xx_{it}}$
  in $\Rn$ and
  $\seq{\lambda_{it}}$ in $\Rpos$, for $i=1,\dotsc,m$,
  such that:
  \begin{itemize}[noitemsep]
  \item
    $\xx_{it}\rightarrow\xbar_i$ for $i=1,\dotsc,m$.
  \item
    $\sum_{i=1}^m \lambda_{it} = 1$ for all $t$.
  \item
    The sequence %
    $\zz_t=\sum_{i=1}^m \lambda_{it} \xx_{it}$ converges to
    $\zbar$.
  \end{itemize}
  The same equivalence holds if we additionally require either or both
  of the following:
  \begin{letter-compact}
  \item  \label{thm:e:7:add:a}
    $\lambda_{it}\to\lambar_i$ for $i=1,\dotsc,m$, for some
    $\lambar_i\in [0,1]$ with
    $\sum_{i=1}^m \lambar_i = 1$.
  \item  \label{thm:e:7:add:b}
    Each of the sequences $\seq{\xx_{it}}$ is span-bound,
    for $i=1,\dotsc,m$.
  \end{letter-compact}
\end{theorem}

\begin{proof}
  ~

\begin{proof-parts}
\pfpart{``If'' ($\Leftarrow$):}
Suppose there exist sequences of the form given in the
theorem.
Then for each $\uu\in\Rn$,
\begin{align*}
  \zbar\cdot\uu
  = \lim %
              (\zz_t\cdot\uu)
  &= \lim %
              \BiggParens{\sum_{i=1}^m \lambda_{it} \xx_{it} \cdot\uu}
  \\
  &\leq
      \lim %
             \BigParens{\max\set{\xx_{1t}\cdot\uu,\dotsc,\xx_{mt}\cdot\uu}}
  \\
  &=
      \max\{\xbar_1\cdot\uu,\dotsc,\xbar_m\cdot\uu\}.
\end{align*}
The first equality is
by \Cref{thm:i:1}(\ref{thm:i:1c}), since
$\zz_t\rightarrow\zbar$.
The last equality (as well as the fact that the last limit exists)
is by continuity of the function which
computes the maximum of $m$ numbers in $\Rext$,
and because
$\xx_{it}\rightarrow\xbar_i$, implying
$\xx_{it}\cdot\uu\rightarrow\xbar_i\cdot\uu$, for $i=1,\dotsc,m$
(again by \Cref{thm:i:1}\ref{thm:i:1c}).
Since this holds for all $\uu\in\Rn$, we obtain $\zbar\in\simplex{V}$ by
\Cref{pr:ohull-simplify}(\ref{pr:ohull-simplify:a},\ref{pr:ohull-simplify:b}).

\pfpart{``Only if'' ($\Rightarrow$):}
Assume $\zbar\in\simplex{V}$.
We will construct the needed sequences explicitly.
For this purpose, we first prove the following lemma showing that within
any neighborhood of $\zbar$, there must exist a point that is a convex
combination (in $\Rn$) of finite points from any collection of neighborhoods of the points in $V$.

\begin{lemma}  \label{lem:e:1mod}
  Let $\zbar\in\simplex{V}$, where $V=\{\xbar_1,\dotsc,\xbar_m\}$.
  Furthermore, let $Z\subseteq\extspace$ be any neighbor\-hood of $\zbar$,
  and let $X_i\subseteq\eRn$ be any neighborhood of $\xbar_i$,
  for $i=1,\dotsc,m$.
  Then there exist $\zz\in\Rn\cap Z$ and $\xx_i\in\Rn\cap X_i$,
  for $i=1,\dotsc,m$, such that $\zz$ is a convex combination of
  $\xx_1,\dotsc,\xx_m$.
\end{lemma}

\begin{proofx}
Let $i\in\{1,\dotsc,m\}$.
Since $X_i$ is a neighborhood of $\xbar_i$,
by \Cref{pr:base-equiv-topo},
there exists a base element $X'_i$ in the astral topology
(of the form given in Eq.~\ref{eq:h:3a}) such that
$\xbar_i\in X'_i\subseteq X_i$.
Furthermore, $X'_i$ is convex by
\Cref{pr:e1}(\ref{pr:e1:d}).
Let $R_i=\Rn\cap X'_i$. Then $R_i$ is convex (by
\Cref{pr:e1}\ref{pr:e1:int-rn}),
and
$\xbar_i\in X'_i\subseteq\Xbarpi=\Rbari$,
where the last equality follows by \Cref{pr:closed-set-facts}(\ref{pr:closed-set-facts:b}).

Let
\[
   R= \conv\BiggParens{\bigcup_{i=1}^m R_i}.
\]
Then $R\subseteq\Rn$, and
for $i=1,\dotsc,m$,
$R_i\subseteq R$, so also
$\xbar_i\in\Rbari\subseteq\Rbar$.
Therefore,
\[
  \zbar\in\simplex{V}\subseteq\ohull\Rbar=\ohull R=\Rbar,
\]
where
the second inclusion follows because $\ohull$ is a hull operator and
$V\subseteq\Rbar$, and the equalities
follow respectively from \Cref{pr:ohull:hull} and \Cref{thm:e:6}. Thus, $\zbar\in\Rbar$.

Let $\zz$ be any point in $Z\cap R$;
this set cannot be empty because $Z$ is a neighborhood of $\zbar$ and $R$ intersects
every neighborhood of $\zbar$, since $\zbar\in\Rbar$
(\Cref{pr:closure:intersect}\ref{pr:closure:intersect:a}).
Since $\zz$ is in $R$, which is the convex hull of the union of the convex sets $R_i$,
\Cref{roc:thm3.3} then implies that
$\zz$ is a convex combination of some points
$\xx_1,\dotsc,\xx_m$ with
$\xx_i\in R_i\subseteq\Rn\cap X_i$
for $i=1,\dotsc,m$.
Thus, $\zz$~has the properties stated in the lemma.
\end{proofx}

To complete the proof of \Cref{thm:e:7}, let
$\countset{Z}$
be a nested countable neighborhood base for $\zbar$,
and for each $i=1,\dotsc,m$, let
$\countsetgen{X_{it}}$
be a nested countable neighborhood base for $\xbar_i$.
(These must exist by
\Cref{thm:first:local}.)

For each $t$,
by \Cref{lem:e:1mod}, applied to $Z_t$ and the $X_{it}\negKern$'s,
there exist $\zz_t\in\Rn\cap Z_t$ and $\xx_{it}\in\Rn\cap X_{it}$
such that
$\zz_t = \sum_{i=1}^m \lambda_{it} \xx_{it}$
for some
$\lambda_{it}\in [0,1]$
with $\sum_{i=1}^m \lambda_{it} = 1$.
Then $\zz_t\rightarrow \zbar$
(by \Cref{prop:nested:limit}),
and likewise, $\xx_{it}\rightarrow\xbar_i$ for each $i$.

\pfpart{Additional requirement~(\ref{thm:e:7:add:a}):}
Suppose sequences as stated in the theorem exist.
Note that the sequence of vectors
$\trans{[\lambda_{1t},\dotsc,\lambda_{mt}]}$
is in the compact set $[0,1]^m$, and therefore, there must exist a
convergent subsequence.
Discarding those $t$ outside
this subsequence yields a sequence that still satisfies all of the
properties stated in the theorem, and in addition, provides that
$\lambda_{it}\rightarrow\lambar_i$, for some $\lambar_i\in [0,1]$ with
$\sum_{i=1}^m \lambar_i = \lim (\sum_{i=1}^m \lambda_{it}) = 1$
(by continuity of addition).

\pfpart{Additional requirement~(\ref{thm:e:7:add:b}):}
Suppose sequences as stated in the theorem exist.
Then by \Cref{thm:spam-limit-seqs-exist},
for $i=1,\dotsc,m$, there exists a span-bound sequence
$\seq{\xx'_{it}}$ with
$\seqeq{\xx'_{it}}{\xx_{it}}$,
implying $\xx'_{it}\rightarrow\xbar_i$.
For all $t$, let
$\zz'_t=\sum_{i=1}^m \lambda_{it} \xx'_{it}$.
Then $\seqeq{\zz'_t}{\zz_t}$
by repeated application of
\Cref{prop:strong:eq:propties}(\ref{prop:strong:eq:propties:a},\ref{prop:strong:eq:propties:c}).
Hence, $\zz'_t\rightarrow\zbar$ by 
\Cref{pr:eq-in-lim-same-lim}
since $\zz_t\rightarrow\zbar$.%
\indexg{polytopes, astral!sequential characterization|)}%
\qedhere
\end{proof-parts}
\end{proof}

\indexg{segments, astral!sequential characterization|(}%
Since the segment joining two points is the same as the outer hull
of the two points, we immediately obtain the following corollary
which shows, in a sense, that
the segment joining points $\xbar$ and $\ybar$ in $\extspace$
is the union of all limits of sequences of line segments in $\Rn$
whose endpoints converge to $\xbar$ and $\ybar$.

\begin{corollary}  \label{cor:e:1}
  Let $\xbar,\ybar\in\extspace$,
  and let $\zbar\in\extspace$.
  Then $\zbar\in\lb{\xbar}{\ybar}$ if
  and only if there exist sequences $\seq{\xx_t}$ and
  $\seq{\yy_t}$ in $\Rn$, and
  $\seq{\lambda_t}$ in $[0,1]$
  such that
  $\xx_t\rightarrow\xbar$, $\yy_t\rightarrow\ybar$,
  and
  $ (1-\lambda_t) \xx_t + \lambda_t \yy_t \rightarrow\zbar$.

  The same equivalence holds if
  the sequence $\seq{\lambda_t}$ is additionally required
  to converge to a limit in $[0,1]$,
  and the sequences $\seq{\xx_t}$ and
  $\seq{\yy_t}$ are required to be span-bound.
\end{corollary}

\begin{example}
Returning to \Cref{ex:seg-oe1-oe2},
according to \Cref{cor:e:1}, for every point
$\zbar\in\lb{\limray{\ee_1}}{\limray{\ee_2}}$, there exist sequences
$\seq{\xx_t}$ and $\seq{\yy_t}$ in $\R^2$, and $\seq{\lambda_t}$ in
$[0,1]$ such that
$\xx_t\rightarrow\limray{\ee_1}$,
$\yy_t\rightarrow\limray{\ee_2}$,
and $\zz_t\rightarrow\zbar$ where
$ \zz_t= (1-\lambda_t) \xx_t + \lambda_t \yy_t $.
Here we give concrete examples of such sequences for each of the
points listed in \eqref{eq:ex:seg-oe1-oe2:1}.

Suppose first that $\zbar=\limray{\ee_1}\plusl\beta\ee_2$ for some
$\beta\in\Rextpos$.
In this case, we can let
$\xx_t=t \ee_1$ and $\yy_t=t \ee_2$,
and then let $\lambda_t=\min\{1,\,\beta/t\}$ if $\beta\in\Rpos$,
and $\lambda_t=t^{-1/2}$ if $\beta=\oms$.
The case $\zbar=\limray{\ee_2}\plusl\beta\ee_1$, for some
$\beta\in\Rextpos$,
is similar.

Otherwise, suppose $\zbar=\limray{\vv}\plusl\beta\ww$
for some $\vv=\trans{[v_1,v_2]}\in\R^2_{>0}$ and 
$\beta\in\eR$, where $\ww=\trans{[v_2,-v_1]}$.
In this case, we can let
$\xx_t=2 (t v_1+ \beta_t v_2) \ee_1$,
$\yy_t=2 (t v_2 -  \beta_t v_1) \ee_2$,
and $\lambda_t=1/2$
(so that $\zz_t= t \vv + \beta_t \ww$),
where $\beta_t=\beta$ if $\beta\in\R$,
and otherwise
$\beta_t=\sqrt{t}$ if $\beta=\oms$
and
$\beta_t=-\sqrt{t}$ if $\beta=-\oms$.

In all cases, it can be checked that the resulting sequences have the
stated
\indexg{segments, astral!sequential characterization|)}%
properties.
\end{example}

\indexg{caratheodorys theorem@Carath\'{e}odory's theorem!outer hull and|(}%
By combining \Cref{thm:e:7} with Carath\'{e}odory's theorem (\Cref{roc:thm17.1})
we obtain the following:

\begin{theorem}  \label{thm:carath}
  Suppose $\zbar\in\simplex{V}$ for some finite set
  $V\subseteq\extspace$.
  Then $\zbar\in\simplex{V'}$ for some $V'\subseteq V$ with
  $|V'|\leq n+1$.
\end{theorem}

\begin{proof}
Let $V=\{\xbar_1,\dotsc,\xbar_m\}$.
Since $\zbar\in\simplex{V}$, there exist sequences as given in
\Cref{thm:e:7}.
For each $t$,
let $I_t$ be the nonzero indices of the $\lambda_{it}\negKern$'s, that is,
\[
  I_t = \bigBraces{ i\in\{1,\dotsc,m\} :\: \lambda_{it}>0 }.
\]
By Carath\'{e}odory's theorem
(\Cref{roc:thm17.1}),
we can assume without loss of
generality that the $\lambda_{it}\negKern$'s have been chosen
in such a way that $|I_t|\leq n+1$, for all $t$.
Since there are only finitely many subsets of $m$ items, there must
exist some subset $I\subseteq\{1,\dotsc,m\}$ for which
$I_t=I$ for infinitely many values of $t$ (implying $|I|\leq n+1$).
On the subsequence consisting of all such values of $t$,
all conditions of the theorem are satisfied with
each $\zz_t$ a convex combination of only points $\xx_{it}$ with
$i\in I$.
Applying \Cref{thm:e:7} to this subsequence then shows that
$\zbar\in\simplex{V'}$ where
\[
  V' = \set{ \xbar_i :\: i\in I }.%
\indexg{caratheodorys theorem@Carath\'{e}odory's theorem!outer hull and|)}%
\qedhere
\]
\end{proof}

\indexg{polytopes, astral!included in convex set|(}%
The following \namecref{thm:e:2}
shows that if $V$ is a finite
subset of some convex set,
then
the polytope formed by $V$ 
must also be entirely included
in that set.
This is useful, for instance, for characterizing the convex hull, as
we will see shortly.

\begin{theorem}  \label{thm:e:2}
  Let $S\subseteq\extspace$ be convex, and let $V\subseteq S$ be a
  finite subset.
  Then $\simplex{V}\subseteq S$.
\end{theorem}

\begin{proof}
Let $V=\{\xbar_1,\dotsc,\xbar_m\}\subseteq S$.
Proof is by induction on $m=|V|$.

If $m=0$ (so that $V=\emptyset$), then
$\simplex{V}=\emptyset\subseteq S$.
And if $m=1$, then
$\simplex{V}=\{\xbar_1\}\subseteq S$,
by \Cref{pr:e1}(\ref{pr:e1:ohull}).

For the inductive step,
assume $m\geq 2$, and that the claim holds for $m-1$.
Let $\zbar\in\simplex{V}$.
Then by \Cref{thm:e:7}, there exist sequences
$\seq{\zz_t}$, $\seq{\xx_{it}}$ and $\seq{\lambda_{it}}$ as in that
\namecref{thm:e:7} with $\lambda_{it}$ converging to
some $\lambar_i\in[0,1]$ with $\sum_{i=1}^m\lambar_i=1$.
Since no more than one of these limits $\lambar_i$ can be equal to
$1$, assume without loss of generality that $\lambar_m<1$.
Then for all $t$ sufficiently large, $\lambda_{mt}<1$; by discarding
all other elements from the sequences, we can assume this holds for
all $t$.

For each $t$, let
\[
 \yy_t = \sum_{i=1}^{m-1}  \frac{\lambda_{it}}{1-\lambda_{mt}} \xx_{it},
\]
which is a convex combination of just the $\xx_{it}\negKern$'s
for $i=1,\dotsc,m-1$.
By sequential compactness, the sequence $\seq{\yy_t}$
must have a convergent subsequence;
by discarding all other elements, we can
assume the entire sequence converges to some point $\ybar\in\extspace$.
Since all the conditions of \Cref{thm:e:7} are now satisfied
for $\ybar$, it must be that
$\ybar\in\simplex\set{\xbar_1,\dotsc,\xbar_{m-1}}$,
and so that $\ybar\in S$ by inductive hypothesis.

Further,
$\zz_t = (1-\lambda_{mt}) \yy_t + \lambda_{mt} \xx_{mt}$,
which converges to $\zbar$.
Thus, the conditions of \Cref{cor:e:1} are satisfied, so
$\zbar\in\lb{\ybar}{\xbar_m}$.
Therefore, $\zbar\in S$ since $\ybar$ and $\xbar_m$ are in $S$,
and since $S$ is
\indexg{polytopes, astral!included in convex set|)}%
\indexg{polytopes, astral|)}%
convex.
\end{proof}

\indexg{polyhedral sets, astral|(}%
Astral polyhedral sets are defined analogously to the standard ones:

\begin{definition}
\indexg{polyhedral sets, astral!defined|(}%
  We say that a convex set $Q\subseteq\extspace$ is
  \emph{astral polyhedral}
  if $Q$ is the intersection of finitely many closed astral
  halfspaces.%
\indexg{polyhedral sets, astral!defined|)}%
\end{definition}

Astral polyhedral sets
are closely linked to the corresponding standard polyhedral sets in $\Rn$.
\indexg{polyhedral sets (standard)!astral closure of|(}%
\indexg{closure, astral!standard polyhedral sets@of standard polyhedral sets|(}%
To see this connection, we show first that the closure of a standard
polyhedral convex set, expressed as an intersection of closed halfspaces,
is the intersection of the corresponding closed astral halfspaces, and
is thus astral polyhedral:

\begin{theorem}   \label{thm:polyhedra-closure}
  Let $\uu_1,\dotsc,\uu_m\in\Rn$,
  $\beta_1,\dotsc,\beta_m\in\R$,
  and let
  \begin{align}
    P
    &=
    \bigBraces{
      \xx\in\Rn :\:
      \xx\cdot\uu_i \leq \beta_i
      \textup{ for }
      i=1,\dotsc,m
    },
    \label{eq:thm:polyhedra-closure:2}    
    \\
    Q
    &=
    \bigBraces{
      \xbar\in\extspace :\:
      \xbar\cdot\uu_i \leq \beta_i
      \textup{ for }
      i=1,\dotsc,m
    }.
    \label{eq:thm:polyhedra-closure:3}
  \end{align}
  Then $Q=\Pbar$.
\end{theorem}

\begin{proof}
Clearly, $P\subseteq Q$.
Also, $Q$ is closed (in $\extspace$), being the intersection of
closed sets (namely, closed astral halfspaces).
Therefore, $\Pbar\subseteq Q$.

For the reverse inclusion, suppose $\xbar\in Q$.
To show $\xbar\in\Pbar$, we construct a sequence $\seq{\xx_t}$ in $P$
converging to $\xbar$.
From \eqref{eq:thm:polyhedra-closure:3},
$\xbar\cdot\uu_i\leq\beta_i$ for all $i$.
Let $J$ be the set of indices for which this inequality holds
strictly; that is,
\[
  J
  =
  \bigBraces{
    i\in\{1,\dotsc,m\} :\:
    \xbar\cdot\uu_i < \beta_i
  }.
\]
We can write $\xbar=\limrays{\vv_1,\dotsc,\vv_k}\plusl\qq$
for some $\vv_1,\dotsc,\vv_k,\qq\in\Rn$.
For each $t$, let
$\xx_t = \sum_{j=1}^k b_{t,j} \vv_j + \qq$
where $\seq{\bb_t}$ is any sequence in $\Rk$ with entries converging
to $+\infty$ at decreasing rates (such as $\bb_t=\trans{[t^k,t^{k-1},\dotsc,t^1]}$).
Then $\xx_t\rightarrow\xbar$ (by \Cref{thm:i:seq-rep}).

Let
\[
  U
  =
  \Braces{
    \zbar\in\extspace:\:
    \zbar\cdot\uu_i<\beta_i
    \text{ for } i\in J
  },
\]  
which is open (being a base element, as in Eq.~\ref{eq:h:3a})
and includes $\xbar$.
Therefore, since $\xx_t\rightarrow\xbar$,
only finitely many $\xx_t$ can be outside $U$; by
discarding these, we can assume $\xx_t\in U$ for all $t$.

Next, let $i\in\{1,\dotsc,m\}$.
We claim that $\xx_t\cdot\uu_i\leq \beta_i$ for all $t$.
If $i\in J$, then this follows because $\xx_t\in U$.
Otherwise, suppose $i\not\in J$, implying
$\xbar\cdot\uu_i=\beta_i\in\R$,
and so, by \Cref{pr:vtransu-zero}, that
$\vv_j\cdot\uu_i=0$ for $j=1,\dotsc,k$,
and that $\xbar\cdot\uu_i=\qq\cdot\uu_i$.
Plugging these identities into the definition of $\xx_t$ yields
\[
  \xx_t\cdot\uu_i
  = \sum_{j=1}^k b_{t,j} (\vv_j\cdot\uu_i) + \qq\cdot\uu_i
  = \qq\cdot\uu_i
  = \xbar\cdot\uu_i
  = \beta_i.
\]

Thus, for all $t$,
we have 
$\xx_t\cdot\uu_i\leq\beta_i$ for $i=1,\dotsc,m$,
so
$\xx_t\in P$.
Hence,
$\xbar=\lim\xx_t$
is in
$\Pbar$, so $Q\subseteq\Pbar$.
\end{proof}

\Cref{thm:polyhedra-closure}
shows that taking the closure of a standard polyhedral convex set yields
an astral polyhedral set, and that intersecting an astral polyhedral convex set
with $\Rn$
yields a standard polyhedral set.
Moreover, as we show in the next \namecref{cor:std-ast-polyhedra},
these operations are inverses of one another, thus
defining a one-to-one correspondence between all astral and all standard
polyhedral convex sets.
This shows furthermore that an astral set is an astral polyhedral
convex set if and only if it is the closure of some standard
polyhedral convex set.

\begin{corollary}   \label{cor:std-ast-polyhedra}
  ~

  \begin{letter-compact}
  \item   \label{cor:std-ast-polyhedra:P}
    Let $P\subseteq\Rn$ be a standard polyhedral convex set.
    Then $\Pbar$ is an astral polyhedral convex set, and moreover,
    $\Pbar\cap\Rn=P$.
  \item   \label{cor:std-ast-polyhedra:Q}
    Let $Q\subseteq\extspace$ be an astral polyhedral convex set.
    Then $Q\cap\Rn$ is a standard polyhedral convex set, and moreover,
    $\clbar{Q\cap\Rn}=Q$.
  \end{letter-compact}
\end{corollary}

\begin{proof}
  ~

\begin{proof-parts}
\pfpart{Part~(\ref{cor:std-ast-polyhedra:P}):}
Since $P$ is a polyhedral convex set, it has the form given in
\eqref{eq:thm:polyhedra-closure:2}
for some $\uu_1,\dotsc,\uu_m\in\Rn\wo\{\zero\}$,
and some $\beta_1,\dotsc,\beta_m\in\R$.
By \Cref{thm:polyhedra-closure}, its closure, $\Pbar$, then has the form
given in 
\eqref{eq:thm:polyhedra-closure:3}, which evidently is an astral
polyhedral convex set, whose intersection with $\Rn$ is $P$.

\pfpart{Part~(\ref{cor:std-ast-polyhedra:Q}):}
Since $Q$ is an astral polyhedral convex set, it has the form given in
\eqref{eq:thm:polyhedra-closure:3}
for some $\uu_1,\dotsc,\uu_m\in\Rn\wo\{\zero\}$,
and some $\beta_1,\dotsc,\beta_m\in\R$.
This implies that $Q\cap\Rn$ has the form given in
\eqref{eq:thm:polyhedra-closure:2}, which evidently is a standard
polyhedral convex set, and whose closure is $Q$, by
\Cref{thm:polyhedra-closure}.%
\indexg{polyhedral sets (standard)!astral closure of|)}%
\indexg{closure, astral!standard polyhedral sets@of standard polyhedral sets|)}%
\qedhere
\end{proof-parts}
\end{proof}

Although in $\Rn$ every polytope is polyhedral
(see \Cref{roc:thm19.1}),
this is not the case in $\eRn$.
For instance, there exist astral polytopes that consist entirely
of infinite points, such as the singleton $\{\xbar\}$ for any
$\xbar\in\eRn\setminus\Rn$,
or the set $\lb{\limray{\ee_1}}{\limray{\ee_2}}$ as we saw in
\Cref{ex:seg-oe1-oe2}.
On the other hand,
by \Cref{cor:std-ast-polyhedra}(\ref{cor:std-ast-polyhedra:Q}),
every nonempty astral polyhedral set must have a nonempty intersection with $\Rn$.
Therefore, astral polytopes that are wholly contained in $\eRn\setminus\Rn$ cannot be astral polyhedral.%
\indexg{polyhedral sets, astral|)}

\section{Convex hull}
\label{sec:convex-hull}

\indexg{convex hull, astral|(}%
The astral convex hull of any set in $\extspace$ is defined in the
same way as for sets in $\Rn$:

\begin{definition}
\indexg{convex hull, astral!defined|(}%
Let $S\subseteq\extspace$.
The \emph{convex hull} of $S$, denoted $\conv{S}$,%
\indexm{conv}{$\conv S$}{convex hull}
is the intersection of all convex sets in $\extspace$ that include
$S$.%
\indexg{convex hull, astral!defined|)}
\end{definition}
By \Cref{pr:e1}(\ref{pr:e1:b}),
the convex hull $\conv{S}$ of any set $S\subseteq\extspace$ is convex.
Thus, $\conv{S}$ is the {smallest} convex set that includes $S$.

In fact, this definition and the notation $\conv{S}$ are really an
extension of the standard definition for sets in $\Rn$
(as in \Cref{sec:prelim:convex-sets}), which is possible since
astral and standard convexity have the same meaning for sets in $\Rn$
(\Cref{pr:e1}\ref{pr:e1:a}).

The astral convex hull operation is the hull operator on the set of
all astral convex sets.
As such,
basic properties of standard convex hulls in $\Rn$ carry
over easily, as stated in the next proposition:

\begin{proposition}  \label{pr:conhull-prop}
  Let $S, U\subseteq\extspace$.
  \begin{letter-compact}
  \item  \label{pr:conhull-prop:aa}
    If $S\subseteq U$ and $U$ is convex, then
    $\conv{S}\subseteq U$.
  \item  \label{pr:conhull-prop:b}
    If $S\subseteq U$, then
    $\conv{S}\subseteq\conv{U}$.
  \item  \label{pr:conhull-prop:c}
    If $S\subseteq U\subseteq\conv{S}$, then
    $\conv{U} = \conv{S}$.
  \item  \label{pr:conhull-prop:a}
    $\conv{S}\subseteq\simplex{S}$
    with equality if $|S|<+\infty$.
  \end{letter-compact}
\end{proposition}

\begin{proof}
~

\begin{proof-parts}
\pfpart{Parts~(\ref{pr:conhull-prop:aa}),~(\ref{pr:conhull-prop:b})
  and~(\ref{pr:conhull-prop:c}):}
Since the astral convex hull operation is a hull operator,
these follow from
\Cref{pr:gen-hull-ops}(\ref{pr:gen-hull-ops:b},\ref{pr:gen-hull-ops:c},\ref{pr:gen-hull-ops:d}).

\pfpart{Part~(\ref{pr:conhull-prop:a}):}
The outer hull
$\simplex{S}$ is convex and includes $S$
(\Cref{pr:e1}\ref{pr:e1:ohull}).
Therefore, $\conv{S}\subseteq\simplex{S}$
by part~(\ref{pr:conhull-prop:aa}).
On the other hand, if $S$ is finite, then
$\simplex{S}\subseteq\conv{S}$ by \Cref{thm:e:2}
since $\conv{S}$ is convex and includes $S$.
\qedhere
\end{proof-parts}
\end{proof}

By \Cref{pr:conhull-prop}(\ref{pr:conhull-prop:a}),
the convex hull, $\conv{V}$,
of any finite set $V\subseteq\extspace$ is
the same as its outer hull, $\ohull V$, and so is an astral polytope.
Nonetheless, in general, convex hull and outer hull
need not be the same. For instance, the convex hull of any open convex
set in $\eRn$ (say, an open astral halfspace) is equal to that set, whereas its outer convex hull
is necessarily a closed set.
Later
(in \Cref{sec:sep-cvx-sets}),
we will establish a more specific relationship between these
operations: whereas the convex hull $\conv{S}$
of any set $S\subseteq\extspace$ is the
smallest convex set that includes $S$, the outer hull $\ohull{S}$ is
the smallest \emph{closed} convex set that includes $S$.

\indexg{convex hull, astral!union of polytopes@as union of polytopes|(}%
\indexg{polytopes, astral!convex hull as union of|(}%
As we show next,
the convex hull of any set is the union
of all polytopes formed by all its finite subsets:

\begin{theorem}
\label{thm:convhull-of-simpices}
  Let $S\subseteq\extspace$.
  Then its convex hull is equal to the union of all polytopes formed
  by finite subsets of $S$, that is,
  \begin{align}
  \label{eq:e:6}
    \conv{S}
    &= \bigcup_{\substack{V\subseteq S:\\ \card{V}<+\infty}}  \!\!\simplex{V}
     = \bigcup_{\substack{V\subseteq S:\\ \card{V}\leq n+1}}  \!\!\simplex{V}.
  \end{align}
\end{theorem}

\begin{proof}
Consider the collection of sets appearing in the first union, that is,
\[
  \calC =\set{\simplex{V}:\: V\subseteq S,\,\card{V}<+\infty}.
\]
This collection is updirected, because for any finite sets $V_1,V_2\subseteq S$, their union is also finite and a subset of $S$, so $\ohull(V_1\cup V_2)$ is in $\calC$.
Moreover,
\[
  (\ohull V_1)\cup(\ohull V_2)\subseteq\ohull(V_1\cup V_2),
\]
since $\ohull$ is a hull operator. Thus, by \Cref{pr:e1}(\ref{pr:e1:union}), $\bigcup{\calC}$ is convex.

Also, for any $\xbar\in S$,
$\set{\xbar}=\ohull\set{\xbar}\in\calC$, so $S\subseteq\bigcup{\calC}$, and therefore
$\conv S\subseteq\bigcup{\calC}$
by \Cref{pr:conhull-prop}(\ref{pr:conhull-prop:aa}).
Furthermore, for any finite set $V\subseteq S$, $\ohull
V\subseteq\conv S$ by \Cref{thm:e:2}, so $\bigcup{\calC}\subseteq\conv
S$. Thus, $\bigcup{\calC}=\conv S$, proving the first equality
of \eqref{eq:e:6}.

The second equality then follows by \Cref{thm:carath}.%
\indexg{convex hull, astral!union of polytopes@as union of polytopes|)}%
\indexg{polytopes, astral!convex hull as union of|)}%
\indexg{convex hull, astral|)}%
\end{proof}

\chapter{Constructing and operating on convex sets}
\label{sec:convex-set-ops}  

We next study how various operations on sets in $\extspace$
interact with astral convexity.
For instance,
we will see that affine transformations preserve astral convexity in
much the same way as in standard convex analysis.
Likewise, the interior of any convex set is convex.
Nonetheless, unlike the standard setting, we will see that the closure
of an astral convex set need not be convex.

\section{Convexity under affine transformations}

\indexg{polytopes, astral!affine map@under affine map|(}%
\indexg{outer convex hull!affine map@under affine map|(}%
\indexg{affine maps, astral!applied to polytope|(}%
We first consider affine transformations and begin with the analysis
of astral polytopes. Given any affine map $F$ and a finite set of points $V$ forming
the polytope $\conv{V}=\ohull{V}$,
we show that the image of $\ohull{V}$ under $F$ is itself a polytope
formed by $F(V)$.
Among its consequences,
this will imply that the image or pre-image of a convex
set under an affine map is also convex.

\begin{theorem} \label{thm:e:9}
  Let $\A\in\R^{m\times n}$, $\bbar\in\extspac{m}$, and
  let $F:\extspace\rightarrow\extspac{m}$ be the
  affine map
  $F(\zbar)=\bbar\plusl \A\zbar$
  for $\zbar\in\extspace$.
  Let
  $V=\{\xbar_1,\dotsc,\xbar_\ell\}\subseteq\extspace$.
  Then
  \[
     \simplex{F(V)} = F(\simplex{V}).
  \]
  In particular,
  $\lb{F(\xbar)}{F(\ybar)}=F(\lb{\xbar}{\ybar})$
  for all $\xbar,\ybar\in\extspace$.
\end{theorem}

Later, in
\Cref{cor:ohull-fs-is-f-ohull-s},
we will show that
\Cref{thm:e:9} holds more generally if the finite set $V$ is
replaced by an arbitrary set $S\subseteq\extspace$.
The next lemma argues one of the needed inclusions, proved in
this broader context, as a first step in the proof
of \Cref{thm:e:9}:

\begin{lemma}   \label{lem:f-conv-S-in-conv-F-S}
  Let $\A$, $\bbar$ and $F$ be as given in
  \Cref{thm:e:9}, and let $S\subseteq\extspace$.
  Then
  \[
     F(\ohull{S}) \subseteq \ohull{F(S)}.
  \]
\end{lemma}

\begin{proof}
Let $\zbar\in\ohull{S}$; we aim to show that
$F(\zbar)\in\ohull{F(S)}$.

Let $\uu\in\R^m$.
For all $\xbar\in\extspace$, we have
\begin{equation}   \label{eq:lem:f-conv-S-in-conv-F-S:1}
  F(\xbar)\cdot\uu
  =
  (\bbar\plusl\A\xbar)\cdot\uu
  =
  \bbar\cdot \uu \plusl (\A\xbar)\cdot\uu
  =
  \bbar\cdot \uu \plusl \xbar\cdot(\transA \uu),
\end{equation}
with the second and third equalities following from
\Cref{pr:i:6} and \Cref{thm:Ax-dot-u},
respectively.
We claim that
\begin{equation}   \label{eq:lem:f-conv-S-in-conv-F-S:2}
  F(\zbar)\cdot\uu
  \leq
  \sup_{\xbar\in S} [F(\xbar)\cdot\uu].
\end{equation}
If $\bbar\cdot\uu\in\{-\infty,+\infty\}$, then,
by \eqref{eq:lem:f-conv-S-in-conv-F-S:1},
$F(\xbar)\cdot\uu=\bbar\cdot\uu$ for all $\xbar\in\extspace$,
implying that
\eqref{eq:lem:f-conv-S-in-conv-F-S:2} must hold (with equality) in
this case.

Otherwise, $\bbar\cdot\uu\in\R$, so
\begin{align*}
  F(\zbar)\cdot\uu
  &=
  \bbar\cdot\uu + \zbar\cdot(\transA \uu)
  \\
  &\leq
  \bbar\cdot\uu + \sup_{\xbar\in S} \bigBracks{\xbar\cdot(\transA \uu)}
  \\
  &=
  \sup_{\xbar\in S} \bigBracks{\bbar\cdot\uu + \xbar\cdot(\transA \uu)}
  =
  \sup_{\xbar\in S} \bigBracks{F(\xbar)\cdot\uu}.
\end{align*}
The first and last equalities are by
\eqref{eq:lem:f-conv-S-in-conv-F-S:1}.
The inequality is by
\Cref{pr:ohull-simplify}(\ref{pr:ohull-simplify:a},\ref{pr:ohull-simplify:b})
since $\zbar\in\ohull{S}$.

Thus,
\eqref{eq:lem:f-conv-S-in-conv-F-S:2} holds for all $\uu\in\R^m$.
Therefore, $F(\zbar)\in \ohull{F(S)}$, again by
\Cref{pr:ohull-simplify}(\ref{pr:ohull-simplify:a},\ref{pr:ohull-simplify:b}).
\end{proof}

To prove the reverse inclusion, we use the characterization of outer hulls via sequences of (standard) convex combinations (\Cref{thm:e:7}). Specifically, we show that given any sequence of convex combinations that characterizes a point $\zbar'\in\ohull{F(V)}$,
it is possible to exhibit a sequence of convex combinations that characterizes a point $\zbar\in\ohull{V}$ such that $F(\zbar)=\zbar'$.
In fact, we prove a somewhat stronger result which states that the point $\zbar$ can be obtained as a limit of convex combinations whose coefficients have the same limits as the coefficients of convex combinations converging to $\zbar'$. This stronger form of the result
will be required later
(in proving \Cref{thm:F:seqconv}).

\begin{lemma}
\label{lem:aff}
  Assume the setup of \Cref{thm:e:9}. For $i=1,\dotsc,\ell$, let $\xbar'_i=F(\xbar_i)$,
  let $\lambar_i\in\Rpos$ with $\sum_{i=1}^\ell\lambar_i=1$, and let $\zbar'\in\eRm$. Assume there
  exist sequences $\seq{\xx'_{it}}$ in~$\Rm$
  and $\seq{\lambda'_{it}}$ in $\Rpos$, for $i=1,\dotsc,\ell$,
  such that
  \begin{letter-compact-prime-simple}
  \item\label{i:lem:aff:ap}
    $\xx'_{it}\to\xbar'_i$ and $\lambda'_{it}\to\lambar_i$ for $i=1,\dotsc,\ell$.
  \item\label{i:lem:aff:bp}
    $\sum_{i=1}^\ell\lambda'_{it}=1$ for all $t$.
  \item\label{i:lem:aff:cp}
    $\sum_{i=1}^\ell\lambda'_{it}\xx'_{it}\to\zbar'$.
  \end{letter-compact-prime-simple}
  Then there exists $\zbar\in\eRn$ such that $F(\zbar)=\zbar'$, and for which there exist sequences $\seq{\xx_{it}}$ in $\Rn$ and $\seq{\lambda_{it}}$ in $\Rpos$, for $i=1,\dotsc,\ell$, such that
  \begin{letter-compact-prime}
  \item\label{i:lem:aff:a}
    $\xx_{it}\to\xbar_i$ and $\lambda_{it}\to\lambar_i$ for $i=1,\dotsc,\ell$.
  \item\label{i:lem:aff:b}
    $\sum_{i=1}^\ell\lambda_{it}=1$ for all $t$.
  \item\label{i:lem:aff:c}
    $\sum_{i=1}^\ell\lambda_{it}\xx_{it}\to\zbar$.
  \end{letter-compact-prime}
\end{lemma}

\begin{proof}
By \Cref{thm:icon-fin-decomp},
$\bbar=\ebar\plusl\qq$ for some icon $\ebar\in\corez{m}$ and $\qq\in\R^m$. Since $\ebar$ is an icon, we have $\limray{\ebar}=\ebar$ (\Cref{pr:i:8}\ref{pr:i:8d}), so $\xbar'_i=\limray{\ebar}\plusl\qq\plusl\A\xbar_i$ for $i=1,\dotsc,\ell$.
The next claim uses tools developed in
\Cref{sec:conv-pts-part-forms} to
show how each sequence $\seq{\xx'_{it}}$,
which converges to $\xbar'_i$, can be
written, up to strong equivalence, in terms of sequences converging to
$\ebar$ and $\xbar_i$.

\begin{claimpx}
\label{claim:aff}
   For each $i=1,\dotsc,\ell$, we have
   $\seqeq{\xx'_{it}}{M_{it}\yy_{it}+\qq+\A\xx_{it}}$
   for some span-bound sequence $\seq{\yy_{it}}$ in $\Rm$
   and some sequences
   $\seq{\xx_{it}}$ in $\Rn$ and $\seq{M_{it}}$ in $\Rstrictpos$ such that
   $\yy_{it}\rightarrow\ebar$,
   $\xx_{it}\to\xbar_i$, $M_{it}\to+\infty$, and $\norm{\A\xx_{it}}/M_{it}\to 0$.
\end{claimpx}

\begin{proofx}
  Let $i\in\set{1,\dotsc,\ell}$. Since $\xx'_{it}\to\xbar'_i$, \Cref{thm:lim-plusl-inv} implies
  $\xx'_{it}=M_{it}\yy_{it}+\rr_{it}$ for some span-bound sequence $\seq{\yy_{it}}$ in $\Rm$
and some sequences $\seq{M_{it}}$ in $\Rstrictpos$ and $\seq{\rr_{it}}$ in~$\Rm$ such that
$\yy_{it}\to\ebar$, 
  $\rr_{it}\to\qq\plusl\A\xbar_i$, $M_{it}\to+\infty$, and $\norm{\rr_{it}}/M_{it}\to 0$.

  Hence, $\rr_{it}-\qq\to\A\xbar_i$ (by \Cref{pr:i:7}\ref{pr:i:7f}), so \Cref{thm:inv-lin-seq} implies that $\seqeq{\rr_{it}-\qq}{\A\xx_{it}}$ for some sequence $\seq{\xx_{it}}$ in $\Rn$ such that $\xx_{it}\to\xbar_i$.
  By \Cref{prop:strong:eq:propties}(\ref{prop:strong:eq:propties:a})
 then
  $\seqeq{\rr_{it}}{\qq+\A\xx_{it}}$
  and so $\seqeq{\xx'_{it}}{M_{it}\yy_{it}+\qq+\A\xx_{it}}$.

  It remains to show $\norm{\A\xx_{it}}/M_{it}\to 0$.
  Writing $\rr_{it}=\qq+\A\xx_{it}+\vepsilon_{it}$ where
  $\vepsilon_{it}\to\zero$ (by strong equivalence), we have
  \[
    0\le\frac{\norm{\A\xx_{it}}}{M_{it}}
      = \frac{\norm{\rr_{it}-\qq-\vepsilon_{it}}}{M_{it}}
     \le\frac{\norm{\rr_{it}}+\norm{\qq}+\norm{\vepsilon_{it}}}{M_{it}},
  \]
  where the second inequality is by the triangle inequality.
  The expression on the right converges to $0$
  since $M_{it}\to+\infty$, $\norm{\vepsilon_{it}}\to 0$, and
  $\norm{\rr_{it}}/M_{it}\to 0$.
  Therefore, $\norm{\A\xx_{it}}/M_{it}\to 0$.  
\end{proofx}

  Let $\seq{\yy_{it}}$, $\seq{\xx_{it}}$, and $\seq{M_{it}}$, for $i=1,\dotsc,\ell$,
  be the sequences whose existence was established in \Cref{claim:aff}.
  For each $t$,
  let $\lambda_{it}=\lambda'_{it}$ for $i=1,\dotsc,\ell$,
  and
  let
\begin{equation}   \label{eq:lem:aff:r1}
   \zz_t=\sum_{i=1}^\ell\lambda_{it}\xx_{it},
   \;\text{ and }\;
   \zz'_t=
    \sum_{i=1}^\ell\lambda_{it}\BigBracks{M_{it}\yy_{it}+\qq+\A\xx_{it}}.
\end{equation}
Then $\seqeq{\zz'_t}{\sum_{i=1}^\ell\lambda'_{it}\xx'_{it}}$
by \Cref{claim:aff} and repeated application of
\Cref{prop:strong:eq:propties}(\ref{prop:strong:eq:propties:a},\ref{prop:strong:eq:propties:c}).
Therefore, $\zz'_t\rightarrow\zbar'$ by
condition~(\ref{i:lem:aff:cp}) and
\Cref{pr:eq-in-lim-same-lim}.

By sequential
  compactness of $\eRn$, the sequence $\seq{\zz_t}$ must have a
  subsequence that converges to some $\zbar\in\eRn$.
  Discarding all other sequence elements (as well as the corresponding
  elements of the sequences
  $\seq{\yy_{it}}$, $\seq{\xx_{it}}$, $\seq{M_{it}}$,
  $\seq{\lambda_{it}}$ and $\seq{\zz'_t}$),
  we obtain sequences $\seq{\xx_{it}}$ and $\seq{\lambda_{it}}$ that
  satisfy conditions~(\ref{i:lem:aff:a}),~(\ref{i:lem:aff:b}) and~(\ref{i:lem:aff:c}) of the
  \namecref{lem:aff}, with $\sum_{i=1}^\ell\lambda_{it}\xx_{it}=\zz_t\to\zbar$. 
It remains to show that $F(\zbar)=\zbar'$.
For this, since $\zz'_t\rightarrow\zbar'$, it suffices to prove that
$\zz'_t\rightarrow F(\zbar)$.
We will show this again using tools from
\Cref{sec:conv-pts-part-forms}.

For each $t$,
let $S_t=\sum_{i=1}^\ell\lambda_{it}M_{it}$,
let $\gamma_{it}=\lambda_{it} M_{it} / S_t$
for $i=1,\dotsc,\ell$
(implying $\gamma_{it}\geq 0$ and 
$\sum_{i=1}^\ell\gamma_{it}=1$),
and let $\dd_t=\sum_{i=1}^\ell\gamma_{it}\yy_{it}$.
Note that
each $S_t\in\Rstrictpos$ (since each $M_{it}\in\Rstrictpos$),
and that $S_t\to+\infty$ (by
\Cref{prop:limsup:eR:conv-comb} since
each $M_{it}\to+\infty$).
From these definitions and \eqref{eq:lem:aff:r1},
it follows, for all $t$,
that
\[\zz'_t=S_t \dd_t + \qq + \A\zz_t.\]

We claim $\dd_t\rightarrow\ebar$.
This is because, for $\uu\in\Rn$,
$\yy_{it}\cdot\uu\rightarrow\ebar\cdot\uu$,
for $i=1,\dotsc,\ell$,
by \Cref{thm:i:1}(\ref{thm:i:1c})
since $\yy_{it}\rightarrow\ebar$.
Hence,
\[
  \dd_t\cdot\uu
  =
  \Parens{\sum_{i=1}^\ell \gamma_{it}\yy_{it}}\cdot\uu
  =
  \sum_{i=1}^\ell \gamma_{it}(\yy_{it}\cdot\uu)
  \rightarrow
  \ebar\cdot\uu,
\]
with convergence by
\Cref{prop:limsup:eR:conv-comb}.
Thus, $\dd_t\rightarrow\ebar$, again by
\Cref{thm:i:1}(\ref{thm:i:1c}).
Further, since each $\yy_{it}$ is in $\rspanset{\ebar}$, $\dd_t$ is as
well, so $\seq{\dd_t}$ is span-bound.

Next, $\qq+\A\zz_t\to\qq\plusl\A\zbar$
by continuity of affine maps (\Cref{cor:aff-cont}).
Further, for each $t$,
\[
  0
  \leq
  \frac{\norm{\qq + \A\zz_t}}{S_t}
  =
  \frac
      {\bigNorm{
          \qq
          +
          \sum_{i=1}^\ell (\gamma_{it} S_t/M_{it}) \A\xx_{it}
        }
      }
      {S_t}
  \leq      
       \frac{\norm{\qq}}{S_t}
     + \sum_{i=1}^\ell \gamma_{it}\cdot\frac{\norm{\A\xx_{it}}}{M_{it}}.
\]
The equality is from the definitions of $\zz_t$ and $\gamma_{it}$.
The second inequality is by the triangle inequality.
The expression on the right converges to $0$ since
$S_t\to+\infty$ and $\norm{\A\xx_{it}}/M_{it}\to 0$ (by \Cref{claim:aff}). 
Therefore,
$\norm{\qq + \A\zz_t}/{S_t}\rightarrow 0$.

Thus, applying \Cref{thm:lim-plusl}
(with $\xx_t$, $\yy_t$, $\xbar$, $\ybar$, $\lambda_t$,
as they appear in that \namecref{thm:lim-plusl},
set to $\dd_t$, $\qq+\A\zz_t$, $\ebar$, $\qq\plusl\A\zbar$, $S_t$,
respectively),
we obtain
\[
  \zbar'
  =
  \lim\zz'_t
  =
  \limray{\ebar}\plusl\qq\plusl\A\zbar
  =
  \bbar\plusl\A\zbar
  =
  F(\zbar),
\]
proving the \namecref{lem:aff}.
\end{proof}

\begin{proof}[Proof of \Cref{thm:e:9}]
The inclusion
$F(\ohull{V})\subseteq\ohull{F(V)}$
is by \Cref{lem:f-conv-S-in-conv-F-S}.

For the reverse inclusion, suppose $\zbar'\in\ohull{F(V)}$, and let
$\xbar'_i=F(\xbar_i)$ for $i=1,\dotsc,\ell$. Then by \Cref{thm:e:7},
for $i=1,\dotsc,\ell$, there exist
sequences $\seq{\xx'_{it}}$ in $\Rm$
and $\seq{\lambda'_{it}}$ in $\Rpos$
such that conditions~(\ref{i:lem:aff:ap}),~(\ref{i:lem:aff:bp}),
and~(\ref{i:lem:aff:cp}) of
\Cref{lem:aff} are satisfied
for some $\lambar_i\in\Rpos$ with $\sum_{i=1}^\ell\lambar_i=1$.
Therefore, that \namecref{lem:aff} shows, for $i=1,\dotsc,\ell$, that there exist
sequences $\seq{\xx_{it}}$ in $\Rn$
and $\seq{\lambda_{it}}$ in $\Rpos$
satisfying conditions~(\ref{i:lem:aff:a}),~(\ref{i:lem:aff:b}),
and~(\ref{i:lem:aff:c}) of
the \namecref{lem:aff}, for some $\zbar\in\eRn$
with $F(\zbar)=\zbar'$. By \Cref{thm:e:7}, these conditions imply that
$\zbar\in\ohull{V}$.
Thus, $\zbar'=F(\zbar)\in F(\ohull{V})$, completing the proof.
\end{proof}

\Cref{thm:e:9} does not hold in general if an astral
point $\bbar$ in the definition of $F$ is leftwardly added on the right
rather than the left, which would give rise to
mappings of the form $\xbar\mapsto \A\xbar\plusl\bbar$.
For instance, in~$\eRf{1}$, suppose $F(\ex)=\ex\plusl(-\infty)$ and $V=\set{0,+\infty}$. Then $\ohull V=[0,+\infty]$, so $F(\ohull V)=\set{-\infty,+\infty}$, but $\ohull F(V)=[-\infty,+\infty]$.%
\indexg{affine maps, astral!applied to polytope|)}%
\indexg{polytopes, astral!affine map@under affine map|)}%
\indexg{outer convex hull!affine map@under affine map|)}

\indexg{affine maps, astral!convexity preserved|(}%
\indexg{convex sets, astral!affine map@under affine map|(}%
The next two corollaries
show that both the image and inverse image of a convex set under an affine map
are convex.

\begin{corollary}  \label{cor:thm:e:9}
  Let $F:\extspace\rightarrow\extspac{m}$ be an affine map
  (as in \Cref{thm:e:9}), and let $S$ be a convex subset of
  $\extspace$.
  Then $F(S)$ is also convex.
\end{corollary}

\begin{proof}
Let $F(\xbar)$ and $F(\ybar)$ be any two points of $F(S)$, where
$\xbar,\ybar\in S$.
Then
\[
   \seg\bigParens{F(\xbar),F(\ybar)}
   =
   F\bigParens{\seg(\xbar,\ybar)}
   \subseteq
   F(S),
\]
with the equality from \Cref{thm:e:9}
and the inclusion following from the convexity of~$S$.
Thus, $F(S)$ is convex.
\end{proof}

\begin{corollary}  \label{cor:inv-image-convex}
  Let $F:\extspace\rightarrow\extspac{m}$ be an affine map
  (as in \Cref{thm:e:9}), and let $S$ be a convex subset of
  $\extspac{m}$.
  Then $\finv(S)$ is also convex.
\end{corollary}

\begin{proof}
Let $\xbar,\ybar\in\finv(S)$, so that $F(\xbar),F(\ybar)\in S$.
Then
\[
   F\bigParens{\seg(\xbar,\ybar)}
   =
   \seg\bigParens{F(\xbar),F(\ybar)}
   \subseteq
   S,
\]
with the equality from \Cref{thm:e:9}
and the inclusion from the convexity of $S$.
This implies
$\lb{\xbar}{\ybar}\subseteq \finv(S)$.
Therefore, $\finv(S)$ is convex.
\end{proof}

\indexg{convex hull, astral!affine map@under affine map|(}%
The next corollary shows that \Cref{thm:e:9} holds for the
convex hull of arbitrary sets.
(\Cref{thm:e:9} is then a special case of the corollary in
which $S$ is finite.)

\begin{corollary}  \label{cor:thm:e:9b}
  Let $F:\extspace\rightarrow\extspac{m}$ be an affine map
  (as in \Cref{thm:e:9}), and let $S\subseteq\extspace$.
  Then
  \[  \conv{F(S)} = F(\conv{S}). \]
\end{corollary}

\begin{proof}
Since $F(S)$ is included in $F(\conv{S})$, and since the latter set is
convex by \Cref{cor:thm:e:9},
we must have $\conv{F(S)}\subseteq F(\conv{S})$
(by \Cref{pr:conhull-prop}\ref{pr:conhull-prop:aa}).

For the reverse inclusion, suppose
$F(\xbar)$ is any point in $F(\conv{S})$, where $\xbar\in\conv{S}$.
Then $\xbar\in\simplex{V}$ for some finite $V\subseteq S$,
by \Cref{thm:convhull-of-simpices}.
Thus,
\[
  F(\xbar)
  \in
  F(\simplex{V})
  =
  \simplex{F(V)}
  \subseteq
  \conv{F(S)}.
\]
The equality is by \Cref{thm:e:9}.
The last line again uses
\Cref{thm:convhull-of-simpices}.%
\indexg{convex hull, astral!affine map@under affine map|)}%
\indexg{affine maps, astral!convexity preserved|)}%
\indexg{convex sets, astral!affine map@under affine map|)}%
\end{proof}

\section{Segments with a finite endpoint}
\label{sec:seg-join-point-origin}

\indexg{segments, astral!finite endpoint@with finite endpoint|(}%
We next explicitly work out the segment joining any finite point in
$\Rn$ with any other point in $\extspace$, making use of some of what
was developed in the last section.
We focus first on the structure of $\lb{\zero}{\xbar}$,
the segment joining the origin and an arbitrary point
$\xbar\in\extspace$. We saw instances of this type of segment
in Examples~\ref{ex:seg-zero-e1}, \ref{ex:seg-zero-oe2},
and~\ref{ex:seg-zero-oe2-plus-e1}.

\begin{theorem}
\label{thm:lb-with-zero}
  Suppose $\xbar=\limrays{\vv_1,\dotsc,\vv_k}\plusl\qq$
  where $\vv_1,\dotsc,\vv_k,\qq\in\Rn$.
  Then
  \begin{align}
  \notag
     \lb{\zero}{\xbar}
     &=
     \bigSet{ \limrays{\vv_1,\dotsc,\vv_{j-1}}\plusl\lambda\vv_j
              :\: j\in\{1,\dotsc,k\},\, \lambda\in\R_{\ge 0} }
     \\
  \label{eq:seg:zero}
     &\qquad{}
     \cup\,
     \bigSet{ \limrays{\vv_1,\dotsc,\vv_k}\plusl\lambda\qq
              :\: \lambda\in [0,1] }.
  \end{align}
\end{theorem}

The theorem states that the segment $\seg(\zero,\xbar)$ consists of several blocks. The first block starts at the origin and passes along the ray $\set{\lambda\vv_1:\:\lambda\in\R_{\ge0}}$ to the astron $\limray{\vv_1}$, which is where the second block starts, continuing along $\set{\limray{\vv_1}\plusl\lambda\vv_2:\:\lambda\in\R_{\ge0}}$, and so on, until the final block,
$\set{\limrays{\vv_1,\dotsc,\vv_k}\plusl\lambda\qq :\: \lambda\in [0,1] }$.
This is precisely the block structure that we saw in \Cref{ex:seg-zero-oe2-plus-e1}.
Note that the representation of $\xbar$ given in the theorem need not be canonical.

\begin{proof}
Let $\VV=[\vv_1,\dotsc,\vv_k]$, so $\xbar=\VV\omm\plusl\qq$, and
let $S$ denote the set on the right-hand side of \eqref{eq:seg:zero}.

\begin{proof-parts}
  \pfpart{Part ``$\subseteq$'':}
  Let $\zbar\in\lb{\zero}{\xbar}$;
  we aim to show $\zbar\in S$.
  By \Cref{cor:e:1}, there exist a span-bound
  sequence $\seq{\xx_t}$ in $\Rn$ and a sequence $\seq{\lambda_t}$ in $[0,1]$
  such that $\xx_t\to\xbar$, $\lambda_t$ converges to some $\lamhat\in[0,1]$,
  and the sequence $\zz_t=\lambda_t\xx_t$ converges to $\zbar$.
  (This is because $\rspanset{\zero}=\{\zero\}$, so the only
  span-bound sequence converging to $\zero$ has every element equal to
  $\zero$.)
  
  If $\lamhat>0$, then
  $\zbar=\lim(\lambda_t\xx_t)=\lamhat\xbar=\VV\omm\plusl\lamhat\qq$
  (by Propositions~\ref{pr:scalar-prod-props}\ref{pr:scalar-prod-props:e}
  and~\ref{pr:i:8}\ref{pr:i:8d}),
  so $\zbar\in S$.
  We therefore assume henceforth that $\lamhat=0$, that is, $\lambda_t\to 0$.

  By \Cref{thm:seq-rep-not-lin-ind},
  we can write $\xx_t=\VV\bb_t+\qq_t$ for some sequences
  $\seq{\bb_t}$ in $\Rstrictpos^k$
  and $\seq{\qq_t}$ in $\Rn$ such that
  conditions~(\ref{thm:seq-rep:a}),~(\ref{thm:seq-rep:b}),
  and~(\ref{thm:seq-rep:c})
  of
  \Cref{thm:seq-rep} are satisfied;
  that is, the entries of $\bb_t$ converge to $+\infty$ at decreasing
  rates, and $\qq_t\to\qq$.
  
  By sequential compactness, the sequence $\seq{\lambda_t \bb_t}$
  has a convergent subsequence.
  By discarding all other elements (and corresponding elements of the
  other sequences), we can assume the entire sequence converges
  to some point $\cbar\in\extspac{k}$.
  For $i=1,\dotsc,k$, let $\alpha_i=\cbar\cdot\ee_i$,
  implying $\lambda_t b_{t,i}=\lambda_t\bb_t\cdot\ee_i\rightarrow\alpha_i$
  by \Cref{thm:i:1}(\ref{thm:i:1c}).

  If $\alpha_1=0$ then also
  $\lambda_t b_{t,i}=(\lambda_t b_{t,1})\cdot(b_{t,i}/b_{t,1})\rightarrow 0\cdot 0 = 0$
  for $i=2,\dotsc,k$ by
  \Cref{prop:dec:trans} 
  (and \Cref{prop:lim:eR}\ref{i:lim:eR:genmul}).
  Thus,
  $\lambda_t\bb_t\to\zero$.
  Since also $\lambda_t\rightarrow 0$ and $\qq_t\rightarrow\qq$,
  this implies
  $\zz_t=\lambda_t\xx_t=\VV(\lambda_t\bb_t)+\lambda_t\qq_t\to\zero$,
  so $\zbar=\lim\zz_t=\zero$, which is in $S$.

  Therefore, in the remainder, we assume $\alpha_1>0$.
  Let $j\in\{1,\dotsc,k\}$ be the last index such that
  $\alpha_j>0$.
  In particular, since $\lambda_t b_{t,1}\to\alpha_1$, this
  implies $\lambda_t>0$ for all $t$ sufficiently large;
  by discarding all other sequence elements, we can assume
  $\lambda_t\in\Rstrictpos$ for all $t$.
  
  By our choice of $j$, $\alpha_i=0$ for $i=j+1,\dotsc,k$,
  implying $\lambda_t b_{t,i}\rightarrow 0$.
  Further, for $i=1,\dotsc,j-1$, we must have
  $
    \lambda_t b_{t,i}
    =(\lambda_t b_{t,j})\cdot({b_{t,i}}/{b_{t,j}})
    \to +\infty
  $,
  with convergence following from
  \Cref{prop:lim:eR}(\ref{i:lim:eR:genmul})
  since $b_{t,i}/b_{t,j}\to+\infty$ (by \Cref{prop:dec:trans})
  and $\lambda_t b_{t,j}\to\alpha_j>0$.
  Summarizing, for $i=1,\dotsc,k$,
  \begin{equation}   \label{eq:thm:lb-with-zero:1}
     \lambda_t b_{t,i}
     \rightarrow
    \begin{cases}
      +\infty  & \text{if $i<j$,} \\
      \alpha_j & \text{if $i=j$,} \\
      0        & \text{if $i>j$.}
    \end{cases}
  \end{equation}
  Further,
  for $i=1,\dotsc,k-1$,
  \begin{equation}   \label{eq:thm:lb-with-zero:2}
    \frac{\lambda_t b_{t,i+1}}{\lambda_t b_{t,i}}
    =
    \frac{b_{t,i+1}}{b_{t,i}}
    \rightarrow 0,
  \end{equation}
  since the entries of $\bb_t$ converge to $+\infty$ at decreasing rates.

  We distinguish two cases. First, suppose $\alpha_j=+\infty$. Let
  \begin{align*}
  \qquad\quad
  &
    \VV'=[\vv_1,\dotsc,\vv_j],
  &&\VV''=[\vv_{j+1},\dotsc,\vv_k],
\intertext{%
  and define sequences
}
  &
    \bb'_t=\lambda_t\trans{[b_{t,1},\dotsc,b_{t,j}]},
  &&
    \bb''_t=\lambda_t\trans{[b_{t,{j+1}},\dotsc,b_{t,k}]}.
  \qquad\quad
  \end{align*}
  Furthermore, let $\qq'_t=\VV''\bb''_t+\lambda_t\qq_t$.
  Then, by construction, $\zz_t=\lambda_t\xx_t=\VV'\bb'_t+\qq'_t$.
  By Eqs.~(\ref{eq:thm:lb-with-zero:1})
  and~(\ref{eq:thm:lb-with-zero:2}),
  the entries of
  $\bb'_t$ converge to $+\infty$ at decreasing rates,
  and $\bb''_t\rightarrow\zero$.
  Since also $\lambda_t\rightarrow 0$ and $\qq_t\rightarrow\qq$,
  this shows that
  the sequences $(\bb'_t)$ and $(\qq'_t)$ satisfy the conditions of
  \Cref{thm:seq-rep-not-lin-ind}, with $\qq'_t\to\zero$,
  thereby showing that
  $\zbar=\lim\zz_t=\VV'\omm$, which is in $S$.

  In the remaining case, $\alpha_j\in\Rstrictpos$. Let
  \begin{align*}
  \qquad\quad
  &
    \VV'=[\vv_1,\dotsc,\vv_{j-1}],
  &&
    \VV''=[\vv_{j+1},\dotsc,\vv_k],
  \\
  &
    \bb'_t=\lambda_t\trans{[b_{t,1},\dotsc,b_{t,j-1}]},
  &&
    \bb''_t=\lambda_t\trans{[b_{t,{j+1}},\dotsc,b_{t,k}]},
  \qquad\quad
  \end{align*}
  and let
  $\qq'_t=\lambda_t b_{t,j}\vv_j+\VV''\bb''_t+\lambda_t\qq_t$.
  Then, by construction, $\zz_t=\lambda_t\xx_t=\VV'\bb'_t+\qq'_t$.
  By similar reasoning as in the last case,
  the sequences $(\bb'_t)$ and $(\qq'_t)$ satisfy the conditions of
  \Cref{thm:seq-rep-not-lin-ind}, with $\qq'_t\to\alpha_j\vv_j$,
  thereby showing that $\zbar=\lim\zz_t=\VV'\omm\plusl\alpha_j\vv_j$,
  which is again in $S$.

  Thus, in every case, $\zbar\in S$, so
  $\lb{\zero}{\xbar}\subseteq S$.

  \pfpart{Part ``$\supseteq$'':}
  Define the (polynomially graded) sequence
  $\xx_t=\VV\bb_t+\qq$, where $\bb_t=\trans{[t^k,t^{k-1},\dotsc,t]}$.
  Then $\xx_t\to\VV\omm\plusl\qq=\xbar$
  by \Cref{thm:seq-rep-not-lin-ind}.
  Let $\zbar\in S$.
  We will show that there exists a sequence $\seq{\lambda_t}$ in $[0,1]$
  such that $\lambda_t\xx_t\to\zbar$, proving that $\zbar\in\lb{\zero}{\xbar}$
  by \Cref{cor:e:1}, and so that
  $S\subseteq\lb{\zero}{\xbar}$.

  First, suppose $\zbar=\VV\omm\plusl\lambda\qq$ for some $\lambda\in[0,1]$.
  For each $t$,
  if $\lambda>0$ then set $\lambda_t=\lambda$, and if $\lambda=0$, then set $\lambda_t=t^{-1/2}$.
  In both cases, $\lambda_t\to\lambda$, and also the entries of
  $\lambda_t\bb_t$ converge to $+\infty$ at decreasing rates.
  Therefore,
  $\lambda_t\xx_t=\VV(\lambda_t\bb_t)+\lambda_t\qq\to\VV\omm\plusl\lambda\qq=\zbar$,
  with convergence by \Cref{thm:seq-rep-not-lin-ind}.

  Otherwise, we must have
  $\zbar=[\vv_1,\dotsc,\vv_{j-1}]\omm\plusl c\vv_{j}$ for some $j\in\{1,\dotsc,k\}$
  and $c\in\Rpos$.
  For each $t$, if $c>0$ then let $c_t=c$, and if
  $c=0$ then let $c_t=t^{-1/2}$.
  Let $\VV'=[\vv_1,\dotsc,\vv_{j-1}]$ and
  $\VV''=[\vv_{j+1},\dotsc,\vv_k]$,
  and define $\lambda_t=c_t\, t^{j-k-1}$.
  Then $\VV=[\VV',\vv_j,\VV'']$ and
  $\lambda_t\bb_t=\mtuple{\bb'_t,c_t,\bb''_t}$,
  where
  \begin{align*}
  \quad
  &
     \bb'_t=c_t\trans{[t^{j-1},t^{j-2},\dotsc,t]},
  &&
     \bb''_t=c_t\trans{[t^{-1},t^{-2},\dotsc,t^{-(k-j)}]}.
  \quad
  \end{align*}
  Then for either case in how $c_t$ was defined, we have
  $c_t\to c$, $\lambda_t\to 0$, $\bb''_t\to\zero$, and
  the entries of ${\bb'_t}$ converge to $+\infty$ at decreasing rates.
  Thus,
  \begin{align*}
     \lambda_t\xx_t
     =
     \VV (\lambda_t\bb_t) + \lambda_t\qq
     &=\bigBracks{\VV',\vv_{j},\VV''}\mtuple{\bb'_t, c_t,\bb''_t}+\lambda_t\qq
\\
     &=\VV'\bb'_t+\BigParens{c_t \vv_{j}+\VV''\bb''_t+\lambda_t\qq}
      \to \VV'\omm\plusl c\vv_{j}
      =\zbar,
  \end{align*}
  with convergence again by \Cref{thm:seq-rep-not-lin-ind}.
  \qedhere
\end{proof-parts}
\end{proof}

As an immediate corollary, we can also compute the segment
joining an arbitrary finite point in $\Rn$ and any other point in
$\extspace$:

\begin{corollary}  \label{cor:lb-with-finite}
  Let $\yy\in\Rn$, and
  suppose that $\xbar=\limrays{\vv_1,\dotsc,\vv_k}\plusl\qq$
  where $\vv_1,\dotsc,\vv_k,\qq\in\Rn$.
  Then
  \begin{align*}
     \lb{\yy}{\xbar}
     &=
     \bigSet{ \limrays{\vv_1,\dotsc,\vv_{j-1}}\plusl
              \parens{\lambda\vv_j+\yy}
              :\: j\in\{1,\dotsc,k\},\, \lambda\in\R_{\ge 0} }
     \\
     &\qquad{}
     \cup\,
     \bigSet{ \limrays{\vv_1,\dotsc,\vv_k}\plusl
              \parens{\lambda\qq + (1-\lambda)\yy}
              :\: \lambda\in [0,1] }.
  \end{align*}
\end{corollary}

\begin{proof}
By \Cref{thm:e:9},
\begin{equation*}  %
   \lb{\yy}{\xbar}
   =
   \yy \plusl \seg\bigParens{\zero,\xbar \plusl (-\yy)}.
\end{equation*}
The result now follows by evaluating the right-hand side using
\Cref{thm:lb-with-zero}.
\end{proof}

Here is a simple but useful consequence of
\Cref{cor:lb-with-finite}:

\begin{corollary}  \label{cor:d-in-lb-0-dplusx}
  Let $\ebar\in\corezn$,
  and let $\xbar\in\extspace$.
  Then
  \begin{equation}  \label{eq:cor:d-in-lb-0-dplusx:1}
     \ebar
     \in
     \lb{\zero}{\ebar}
     \subseteq
     \lb{\zero}{\ebar\plusl\xbar}.
  \end{equation}
\end{corollary}

\begin{proof}
The first inclusion of \eqref{eq:cor:d-in-lb-0-dplusx:1}
is immediate.

By \Crefequiv{pr:icon-equiv}{pr:icon-equiv:a}{pr:icon-equiv:c},
$\ebar=\limrays{\vv_1,\dotsc,\vv_k}$ for some
$\vv_1,\dotsc,\vv_k\in\Rn$, while
$\xbar=\limrays{\ww_1,\dotsc,\ww_\ell}\plusl\qq$ for some
$\ww_1,\dotsc,\ww_\ell,\qq\in\Rn$.
Since
\[
  \ebar\plusl\xbar=\limrays{\vv_1,\dotsc,\vv_k,
                           \ww_1,\dotsc,\ww_\ell}\plusl\qq,
\]
it follows from a direct application of
\Cref{thm:lb-with-zero}
that $\ebar \in \lb{\zero}{\ebar\plusl\xbar}$.
Since $\zero$ is also in $\lb{\zero}{\ebar\plusl\xbar}$,
this implies the second inclusion of
\eqref{eq:cor:d-in-lb-0-dplusx:1}
by definition of convexity
(and \Cref{pr:e1}\ref{pr:e1:ohull}).
\end{proof}

As another corollary, we
derive next an inductive description of $\lb{\zero}{\xbar}$. We show that
if $\xbar\in\extspace$ has dominant direction $\vv$,
and thus, if $\xbar=\limray{\vv}\plusl\xbar'$ for some $\xbar'$, then
the points in $\lb{\zero}{\xbar}$ are of two types:
the finite points in $\Rn$, which are exactly all of the nonnegative
multiples of $\vv$;
and the infinite points, which are the points of the form
$\limray{\vv}\plusl\zbar'$ where $\zbar'$ is in the segment
joining $\zero$ and $\xbar'$.

\begin{corollary}
\label{cor:seg:zero}
  Suppose $\xbar=\limray{\vv}\plusl\xbar'$, for some $\vv\in\Rn$ and
  $\xbar'\in\extspace$.
  Then
  \begin{equation}  \label{eq:cor:seg:zero:1}
     \lb{\zero}{\xbar}
     =
     \bigBraces{ \lambda \vv :\: \lambda \in \Rpos }
     \,\cup\,
     \bigBracks{\limray{\vv} \plusl \lb{\zero}{\xbar'}}.
  \end{equation}
\end{corollary}

\begin{proof}
We can write $\xbar'=\limrays{\vv_1,\dotsc,\vv_k}\plusl\qq$
for some $\vv_1,\dotsc,\vv_k,\qq\in\Rn$.
Then $\xbar=\limrays{\vv,\vv_1,\dotsc,\vv_k}\plusl\qq$.
Computing $\lb{\zero}{\xbar}$ and $\lb{\zero}{\xbar'}$
according to
\Cref{thm:lb-with-zero},
the claim then follows by direct examination of the points comprising
each side of
\eqref{eq:cor:seg:zero:1}.
\end{proof}

\indexg{segments, astral!split at a point|(}%
An ordinary segment $\lb{\xx}{\yy}$ joining finite points $\xx$ and
$\yy$ in $\Rn$ can be split into two segments at any point
$\zz\in\lb{\xx}{\yy}$.
The same is true for the segment joining a finite point and
any other astral point, as shown next.

\begin{theorem}
\label{thm:decomp-seg-at-inc-pt}
  Let $\xbar\in\extspace$, let $\yy\in\Rn$, and
  let $\zbar\in\lb{\yy}{\xbar}$.
  Then:
  \begin{letter-compact}
  \item   \label{thm:decomp-seg-at-inc-pt:a}
    $\displaystyle
    \lb{\yy}{\xbar}
    =
    \lb{\yy}{\zbar}
    \cup
    \lb{\zbar}{\xbar}
    $.
  \item   \label{thm:decomp-seg-at-inc-pt:b}
    $\displaystyle
    \lb{\yy}{\zbar}
    \cap
    \lb{\zbar}{\xbar}
    =
    \{\zbar\}
    $.
  \end{letter-compact}
\end{theorem}

\begin{proof}
We first consider the special case that $\yy=\zero$, and
then return to the general case.

Let $\xbar=[\vv_1,\dotsc,\vv_k]\omm\plusl\vv_{k+1}$ be the canonical representation of $\xbar$,
where $\vv_1,\dotsc,\vv_{k+1}\in\Rn$.
Since $\zbar\in\seg(\zero,\xbar)$, by \Cref{thm:lb-with-zero}, it can be written as
$\zbar=[\vv_1,\dotsc,\vv_{j-1}]\omm\plusl\lambda\vv_j$ for some
$j\in\{1,\dotsc,k+1\}$ and $\lambda\in\Rpos$, with $\lambda\le 1$ if $j=k+1$.
Let $\ebar=[\vv_1,\dotsc,\vv_{j-1}]\omm$ and $\vv=\vv_j$, so $\zbar=\ebar\plusl\lambda\vv$.

\begin{proof-parts}
\pfpart{%
  Part~(\ref{thm:decomp-seg-at-inc-pt:a})
  when $\yy=\zero$:
}
We distinguish two cases.
First, if $j=k+1$ then $\xbar=\ebar\plusl\vv$ and $\lambda\le 1$.
In this case,
\begin{align}
  \seg(\zbar,\xbar)
  =
  \seg(\ebar\plusl\lambda\vv,\,\ebar\plusl\vv)
  &=
  \ebar\plusl\lambda\vv\plusl \seg\bigParens{\zero,(1-\lambda)\vv}
  \notag
  \\
  &=
  \ebar\plusl\lambda\vv\plusl
  \Braces{\lambda'\vv :\: \lambda'\in [0,1-\lambda]}
  \notag
  \\
  &=
  \Braces{\ebar\plusl\lambda'\vv :\: \lambda'\in [\lambda,1]}.
  \label{eq:thm:decomp-seg-at-inc-pt:1}
\end{align}
The second equality is by
\Cref{thm:e:9},
and the third by
\Cref{pr:e1}(\ref{pr:e1:a}).
Evaluating
$\lb{\zero}{\xbar}$ and $\lb{\zero}{\zbar}$
using \Cref{thm:lb-with-zero}, and combining with
\eqref{eq:thm:decomp-seg-at-inc-pt:1}, we see that
the claim
holds in this case.

In the alternative case, $j\le k$,
so $\xbar=\ebar\plusl\limray{\vv}\plusl\xbar'$ where
$\xbar'=[\vv_{j+1},\dotsc,\vv_k]\omm\plusl\vv_{k+1}$.
Then
\begin{align}
  \seg(\zbar,\xbar)
  &=
  \seg(\ebar\plusl\lambda\vv,\,\ebar\plusl\limray{\vv}\plusl\xbar')
  \notag
  \\
  &=
  \ebar\plusl\lambda\vv\plusl \seg(\zero,\limray{\vv}\plusl\xbar')
  \notag
  \\
  &=
  \ebar\plusl\lambda\vv\plusl
  \bigBracks{
    \Braces{\lambda'\vv :\: \lambda'\in\Rpos}
    \,\cup\,
    \bigParens{\limray{\vv}\plusl \lb{\zero}{\xbar'}}
  }
  \notag
  \\
  &=
  \bigBraces{\ebar\plusl\lambda'\vv :\: \lambda'\in [\lambda,+\infty)}
  \,\cup\,
  \bigBracks{
    \ebar\plusl\limray{\vv}\plusl \lb{\zero}{\xbar'}
  }.
  \label{eq:thm:decomp-seg-at-inc-pt:4}
\end{align}
The second equality is by
\Cref{thm:e:9} (and since $\lambda\vv\plusl\limray{\vv}=\limray{\vv}$,
by \Cref{pr:i:7}\ref{pr:i:7c}).
The third is by
\Cref{cor:seg:zero}.
As before, we can then see that
the claim
holds in this case by evaluating
$\lb{\zero}{\xbar}$, $\lb{\zero}{\zbar}$ and $\lb{\zero}{\xbar'}$
using \Cref{thm:lb-with-zero}, and combining with
\eqref{eq:thm:decomp-seg-at-inc-pt:4}.

\pfpart{%
  Part~(\ref{thm:decomp-seg-at-inc-pt:b})
  when $\yy=\zero$:
}
Let $\beta=\lambda\norm{\vv}^2$, and let
$H=\set{\zbar'\in\eRn:\:\zbar'\inprod\vv\le\beta}$.
Then $\zero\in H$ since $\zero\cdot\vv=0\leq\beta$, and
$\zbar\in H$ since
$\zbar\cdot\vv=\lambda\norm{\vv}^2=\beta$
by the Case Decomposition Lemma~\ref{lemma:case},
using that $\vv_1,\dotsc,\vv_{k+1}$ are
orthogonal to one another.
Therefore, $\seg(\zero,\zbar)\subseteq H$,
since $H$ is convex
(by \Cref{pr:e1}\ref{pr:e1:c}).

We will show further that
\begin{equation}  \label{eq:thm:decomp-seg-at-inc-pt:5}
  H\cap\seg(\zbar,\xbar)\subseteq\set{\zbar},
\end{equation}
which
will imply
\[
  \set{\zbar}
  \subseteq
  \lb{\zero}{\zbar} \cap \lb{\zbar}{\xbar}
  \subseteq
  H\cap\lb{\zbar}{\xbar}
  \subseteq
  \set{\zbar},
\]
thereby proving
the claim.

Returning to the cases considered in proving
part~(\ref{thm:decomp-seg-at-inc-pt:a}),
if $j\leq k$, then
from \eqref{eq:thm:decomp-seg-at-inc-pt:4}, it is apparent that
the only point $\zbar'\in\seg(\zbar,\xbar)$ that satisfies
$\zbar'\inprod\vv\le\beta$ is $\zbar=\ebar\plusl\lambda\vv$
(noting again the orthogonality of $\vv_1,\dotsc,\vv_{k+1}$,
and also that $\norm{\vv}=1$), thus proving
\eqref{eq:thm:decomp-seg-at-inc-pt:5}.

The argument is similar if $j= k+1$.
If $\zbar'\in\seg(\zbar,\xbar)$ then by
\eqref{eq:thm:decomp-seg-at-inc-pt:1},
$\zbar'=\ebar\plusl\lambda'\vv$ for some $\lambda'\in[\lambda,1]$,
implying
$\zbar'\cdot\vv=\lambda'\norm{\vv}^2$, which is at most
$\beta=\lambda\norm{\vv}^2$ if and only if either $\lambda'=\lambda$
or $\vv=\zero$.
In either case,
$\zbar'=\ebar\plusl\lambda'\vv=\ebar\plusl\lambda\vv=\zbar$,
proving
\eqref{eq:thm:decomp-seg-at-inc-pt:5},
and completing the proof for when $\yy=\zero$.

\pfpart{General $\yy\in\Rn$:}
For the general case (when $\yy$ is not necessarily $\zero$),
let $\xbar'=-\yy\plusl\xbar$ and $\zbar'=-\yy\plusl\zbar$.
Since $\zbar\in\lb{\yy}{\xbar}$, we also have
\[
  \zbar'
  =
  -\yy\plusl\zbar
  \in
  \bigBracks{
    -\yy\plusl\lb{\yy}{\xbar}
  }
  =
  \lb{\zero}{\xbar'},
\]
where the last equality is by
\Cref{thm:e:9}.
Thus, for part~(\ref{thm:decomp-seg-at-inc-pt:a}),
\begin{align*}
  \lb{\yy}{\zbar} \,\cup\, \lb{\zbar}{\xbar}
  &=
  \bigBracks{\yy \plusl \lb{\zero}{\zbar'}}
  \,\cup\,
  \bigBracks{\yy \plusl \lb{\zbar'}{\xbar'}}
\\
  &=
  \yy \plusl \bigBracks{ \lb{\zero}{\zbar'} \,\cup\, \lb{\zbar'}{\xbar'} }
\\
  &=
  \yy \plusl \lb{\zero}{\xbar'}
  =
  \lb{\yy}{\xbar},
\end{align*}
where the first and fourth equalities are by
\Cref{thm:e:9},
and the third is by the proof above.

Similarly,
for part~(\ref{thm:decomp-seg-at-inc-pt:b}),
\begin{align*}
  \lb{\yy}{\zbar} \,\cap\, \lb{\zbar}{\xbar}
  &=
  \bigBracks{\yy \plusl \lb{\zero}{\zbar'}}
  \,\cap\,
  \bigBracks{\yy \plusl \lb{\zbar'}{\xbar'}}
  \\
  &=
  \yy \plusl \bigBracks{\lb{\zero}{\zbar'}
                     \,\cap\,
                     \lb{\zbar'}{\xbar'}
                    }
  \\
  &=
  \yy \plusl \set{\zbar'}
  =
  \set{\zbar}.
  \qedhere
\end{align*}
\end{proof-parts}
\end{proof}

\Cref{thm:decomp-seg-at-inc-pt}
does not hold in general for the segment
joining any two astral points; in other words, the
theorem is false if $\yy$ is replaced by an arbitrary point
in~$\extspace$.
For instance, in $\extspac{2}$, let
$\xbar=\limray{\ee_1}\plusl\limray{\ee_2}$.
Then $\lb{-\xbar}{\xbar}=\extspac{2}$, as we saw in
\Cref{ex:seg-negiden-iden}.
This is different from
$\lb{-\xbar}{\zero}\cup\lb{\zero}{\xbar}$,
which, for instance, does not include $\ee_2$,
as follows from
\Cref{cor:seg:zero}.%
\indexg{segments, astral!split at a point|)}%
\indexg{segments, astral!finite endpoint@with finite endpoint|)}

\section{Interior of a convex set and convex hull of an open set}

\indexg{convex sets, astral!interior of|(}%
\indexg{interior!convex set@of convex set|(}%
We prove next that the interior of a convex set in $\extspace$ is also
convex.
For the proof,
recall that a convex set in $\Rn$ is {polyhedral} if it is the
intersection of a finite collection of closed halfspaces.

\begin{theorem}  \label{thm:interior-convex}
  Let $S$ be a convex subset of $\extspace$.
  Then its interior, $\intr S$, is also convex.
\end{theorem}

\begin{proof}
We assume that $\intr S$ is neither $\emptyset$ nor $\extspace$,
since otherwise the claim holds immediately.
Let $\xbar,\ybar\in \intr S$.
We aim to prove convexity by showing that the segment joining them,
$\lb{\xbar}{\ybar}$, is entirely contained in $\intr S$.

\begin{claimpx}  \label{cl:thm:interior-convex:1}
  There exists a topological base element $X\subseteq\extspace$
  (of the form given in Eq.~\ref{eq:h:3a})
  that includes $\xbar$, and
  whose closure $\Xbar$ is included in $S$.
\end{claimpx}

\begin{proofx}
Let $C=\extspace \setminus (\intr S)$, which is closed, being the
complement of $\intr S$, and nonempty by assumption.
Furthermore, $\xbar\not\in C$.
Therefore, since astral space is regular
(\Cref{thm:i:1}\ref{thm:i:1aa}),
there exist disjoint open sets $X$ and $V$ such that
$\xbar\in X$ and $C\subseteq V$.
Without loss of generality, we can assume $X$ is a base
element (since otherwise we can replace it by a base element
containing $\xbar$ and included in $X$,
by \Cref{pr:base-equiv-topo}).
Since $X$ and $V$ are disjoint,
$X$ is included in the closed set $\extspace \setminus V$,
and therefore,
\[
  \Xbar
  \subseteq \extspace \setminus V
  \subseteq \extspace\setminus C = \intr S \subseteq S.
\qedhere
\]
\end{proofx}

Let $X$ be as in \Cref{cl:thm:interior-convex:1}, and let $Y$ be
a similar base element for $\ybar$
(so $\ybar\in Y$ and $\Ybar\subseteq S$).

Next, let
$X'=\cl(X\cap\Rn)$
and
$Y'=\cl(Y\cap\Rn)$.
The set $X\cap\Rn$ is open and nonempty since $\Rn$ is dense in
$\extspace$.
Moreover, since $X$ is a base element of the form given in
\eqref{eq:h:3a},
$X\cap\Rn$ is an intersection of (standard) open halfspaces in $\Rn$, which
means that $X'$, its closure in $\Rn$, is an intersection of closed
halfspaces in $\Rn$ by 
\Cref{roc:thm6.5}.
In other words, $X'$ is a nonempty, polyhedral convex set,
as is $Y'$ by the same argument.

Let
$R=\conv(X' \cup Y')$
be the convex hull of their union.
Since $X'$ and $Y'$ are polyhedral, $\cl R$,
the closure of $R$ in $\Rn$, is also polyhedral
by \Cref{roc:thm19.6}.
Thus,
\begin{equation}  \label{eqn:thm:interior-convex:1}
  \cl R = \Braces{ \zz\in\Rn :\:
                   \zz\cdot\uu_i \leq b_i \mbox{~for~} i=1,\dotsc,k }
\end{equation}
for some $\uu_1,\dotsc,\uu_k\in\Rn\wo\{\zero\}$ and $b_1,\dotsc,b_k\in\R$, and
some $k\geq 0$.

The closure of $R$ in $\extspace$ is exactly the convex hull of
$\Xbar\cup\Ybar$:

\begin{claimpx}  \label{cl:thm:interior-convex:2}
  $\Rbar = \conv(\Xbar\cup\Ybar)$.
\end{claimpx}

\begin{proofx}
By construction,
$X\cap\Rn \subseteq X' \subseteq R$.
Therefore, by
\Cref{pr:closed-set-facts}(\ref{pr:closed-set-facts:b}),
$\Xbar = \clbar{X\cap\Rn}\subseteq \Rbar$.
Likewise, $\Ybar\subseteq\Rbar$.
Therefore, since $\Rbar$ is convex
(by \Cref{thm:e:6}),
$\conv(\Xbar\cup\Ybar)\subseteq\Rbar$
(by \Cref{pr:conhull-prop}\ref{pr:conhull-prop:aa}).

For the reverse inclusion, suppose $\zbar\in\Rbar$, implying
there exists a sequence $\seq{\zz_t}$ in $R$ that converges to
$\zbar$.
Then by $R$'s definition, and since $X'$ and $Y'$ are convex,
for each $t$, we can write
$  \zz_t = (1-\lambda_t) \xx_t + \lambda_t \yy_t $
for some $\xx_t\in X'$, $\yy_t\in Y'$ and $\lambda_t\in [0,1]$
(by \Cref{roc:thm3.3}).
By sequential compactness,
the sequence $\seq{\xx_t}$ must have a convergent subsequence; by
discarding all other elements, we can assume that the entire sequence
converges to some point $\xbar'$, which thereby must be in $\Xbarp$.
Furthermore, $\Xbarp=\clbar{X\cap\Rn}=\Xbar$
by \Cref{pr:closed-set-facts}(\ref{pr:closed-set-facts:aa},\ref{pr:closed-set-facts:b}).
Repeating this argument, we can take the sequence $\seq{\yy_t}$ to
converge to some point $\ybar'\in\Ybar$.
Therefore, applying \Cref{cor:e:1},
$\zbar\in\lb{\xbar'}{\ybar'}\subseteq\conv(\Xbar\cup\Ybar)$.
\end{proofx}

Let $Q$ be the set
\[
  Q = \bigBraces{ \zbar\in\extspace :\:
                   \zbar\cdot\uu_i < b_i \mbox{~for~} i=1,\dotsc,k },
\]
which is the intersection of open astral halfspaces
corresponding to the closed halfspaces (in $\Rn$) whose intersection
defines $\cl R$ in \eqref{eqn:thm:interior-convex:1}.
Then $Q$ is a base element and so is open.
Moreover,
\begin{equation}
\label{eq:cl:thm:interior-convex:3}
 Q\subseteq\Qbar=\overline{Q\cap\Rn}\subseteq\clbar{\cl R}=\Rbar,
\end{equation}
where the first equality is by \Cref{pr:closed-set-facts}(\ref{pr:closed-set-facts:b}) (since $Q$ is open),
the second inclusion is because $Q\cap\Rn\subseteq\cl R$, and the final equality is by \Cref{pr:closed-set-facts}(\ref{pr:closed-set-facts:aa}).

\begin{claimpx}  \label{cl:thm:interior-convex:4}
  $\xbar\in Q$.
\end{claimpx}

\begin{proofx}
To prove the claim, we show $\xbar\cdot\uu_i < b_i$ for
each $i\in\{1,\dotsc,k\}$.

Let $H_i$ be the closed astral halfspace
\[ H_i=\{\zbar\in\extspace :\: \zbar\cdot\uu_i \leq b_i\}. \]
Then $R\subseteq \cl R \subseteq H_i$,
so $\Rbar\subseteq H_i$ since $H_i$ is closed.
Therefore, using \Cref{cl:thm:interior-convex:2},
$X \subseteq \Xbar \subseteq \Rbar \subseteq H_i$.
Since $X$ is open, we have $X\subseteq\intr H_i=\{\zbar\in\extspace :\: \zbar\cdot\uu_i < b_i\}$
with the equality following by \Cref{prop:intr:H}.
Since $\xbar\in X$, this means
$\xbar\cdot\uu_i < b_i$.
\end{proofx}

By \Cref{cl:thm:interior-convex:1}, $\Xbar\subseteq S$, and
similarly $\Ybar\subseteq S$, so
$\conv(\Xbar\cup\Ybar)\subseteq S$
by \Cref{pr:conhull-prop}(\ref{pr:conhull-prop:aa}).
Together with \eqref{eq:cl:thm:interior-convex:3} and
\Cref{cl:thm:interior-convex:2} this implies
$Q\subseteq \Rbar=\conv(\Xbar\cup\Ybar)\subseteq S$.
Since $Q$ is open, this further implies that
$Q$ is included in the interior of $S$.
By \Cref{cl:thm:interior-convex:4},
$\xbar\in Q$, and by the same argument, $\ybar\in Q$.
Since $Q$ is convex (by \Cref{pr:e1}\ref{pr:e1:d}),
it follows that
$\lb{\xbar}{\ybar}\subseteq Q \subseteq \intr S$,
completing the proof.%
\indexg{convex sets, astral!interior of|)}%
\indexg{interior!convex set@of convex set|)}%
\end{proof}

\indexg{convex hull, astral!openness preserved|(}%
As a consequence, the convex hull of an open set is also open:

\begin{corollary}  \label{cor:convhull-open}
  Let $U\subseteq\extspace$ be open.
  Then its convex hull, $\conv{U}$, is also open.
\end{corollary}

\begin{proof}
Let $S=\conv{U}$.
Then $U\subseteq S$, implying, since $U$ is open, that
$U\subseteq \intr S$.
By \Cref{thm:interior-convex}, $\intr S$ is convex.
Therefore, $S=\conv{U}\subseteq \intr S$
(by \Cref{pr:conhull-prop}\ref{pr:conhull-prop:aa}).
Thus, $S= \intr S$ (since $\intr S\subseteq S$ always),
so $S$ is open.%
\indexg{convex hull, astral!openness preserved|)}
\end{proof}

\section{Closure of a convex set and convex hull of a closed set}
\label{sec:closure-convex-set}

\indexg{closure, astral!convex set@of convex set|(}%
\indexg{convex sets, astral!closure of|(}%
In standard convex analysis, the closure of any convex set is also
convex and equal to the intersection of closed halfspaces that contain it.
In \Cref{thm:e:6},
we saw that also the {astral} closure of any
convex set in $\Rn$ is convex, and is moreover equal to the outer hull of the set,
the intersection of closed astral halfspaces that contain the set.

As we show in the next theorem, this is not true
for arbitrary convex sets in $\extspace$.
In particular, we show that for $n\geq 2$, there exist sets in $\extspace$
that are convex, but whose closures are not convex.
This also means that the closure of such a set
cannot be equal to its
outer hull, since the outer hull of any set is always convex.

\begin{theorem}  \label{thm:closure-not-always-convex}
  For $n\geq 2$, there exists a convex set $C\subseteq\extspace$ whose closure
  $\Cbar$ is not convex.
  Consequently, $C$'s closure is not equal to its outer hull;
  that is, $\Cbar\neq\ohull C$.

  Explicitly, the following set has these properties:
  \begin{equation}  \label{eq:thm:closure-not-always-convex:2}
    C
    =
    \bigBraces{ \lambda \ee_1 :\: \lambda\in\R }
    \,\cup\, \bigBraces{ \limray{\ee_1} \plusl \lambda\ee_2
                                  :\: \lambda\in\Rpos },
  \end{equation}
  where $\ee_1$ and $\ee_2$ are the first two standard
  basis vectors in $\Rn$.
\end{theorem}

\begin{proof}
Let $C$ be as given in
\eqref{eq:thm:closure-not-always-convex:2}.
We will show that $C$ is convex but that its closure $\Cbar$ is not
(implying $\Cbar\neq\ohull C$ by \Cref{pr:e1}\ref{pr:e1:ohull}).
We first derive an alternate expression for $C$.

For $\alpha\in\R$ and $\beta\in\Rpos$, let
$S_{\alpha,\beta}$ be the segment
$S_{\alpha,\beta}=\lb{\alpha \ee_1}{\,\limray{\ee_1}\plusl\beta\ee_2}$.
We can compute $S_{\alpha,\beta}$ explicitly using
\Cref{cor:lb-with-finite}:
\begin{align}
  S_{\alpha,\beta}
  &=
  \Braces{ \lambda \ee_1 + \alpha \ee_1 :\: \lambda \in \Rpos }
  \nonumber
  \\
  &\qquad{}
    \cup\, \Braces{ \limray{\ee_1} \plusl \lambda\beta\ee_2 \plusl (1-\lambda)\alpha\ee_1
                                  :\: \lambda\in [0,1] }
  \nonumber
  \\
  &=
  \{ \lambda \ee_1 :\: \lambda \in [\alpha,+\infty) \}
  \,\cup\, \bigBraces{ \limray{\ee_1} \plusl \lambda\ee_2
                                  :\: \lambda \in [0,\beta] },
  \label{eq:thm:closure-not-always-convex:1}
\end{align}
using that $\limray{\ee_1}\plusl\zbar\plusl\lambda\ee_1=\limray{\ee_1}\plusl\zbar$ for any $\zbar\in\eRn$ and any $\lambda\in\R$
(by the Projection Lemma~\ref{lemma:proj}).
The union over the expressions on the right-hand side of
\eqref{eq:thm:closure-not-always-convex:1} is evidently equal to the
expression for $C$ in
\eqref{eq:thm:closure-not-always-convex:2}.
Thus, $C$ is the union over the collection of sets
$\{S_{\alpha,\beta} :\: \alpha\in\R,\,\beta\in\Rpos\}$,
which is updirected since
$S_{\alpha,\beta}\cup S_{\alpha',\beta'}
 =S_{\min\set{\alpha,\alpha'},\max\set{\beta,\beta'}}$
for any $\alpha,\alpha'\in\R$ and $\beta,\beta'\in\Rpos$.
Therefore, $C$ is convex by \Cref{pr:e1}(\ref{pr:e1:union}).

We next show that $\Cbar$ is {not} convex.
Let $\xbar=\limray{(-\ee_1)}$ and
let $\ybar=\limray{\ee_1}\plusl\limray{\ee_2}$.
From \eqref{eq:thm:closure-not-always-convex:2}, we see that
both $\xbar$ and $\ybar$ are in $\Cbar$ since
$\xbar$ is the limit of the sequence
$\seq{-t \ee_1}$ and
$\ybar$ is the limit of the sequence
$\seq{\limray{\ee_1}\plusl t\ee_2}$, and all elements of both
sequences are in~$C$.
Let $\zz=\ee_2$.
To prove $\Cbar$ is not convex, we will show that $\zz$ is on the
segment joining $\xbar$ and $\ybar$, but that $\zz$ is not itself in
$\Cbar$.
Indeed, $\zz\not\in\Cbar$ since, for instance, the open set
$\{ \xx \in \Rn :\: \xx\cdot\ee_2 > 0 \}$
includes $\zz$ but is entirely disjoint from $C$.

It remains to show that $\zz\in\lb{\xbar}{\ybar}$.
To do so, we construct sequences satisfying the conditions of
\Cref{cor:e:1}.
Specifically, for each $t$, let
\begin{align*}
  \xx_t &= -t \ee_1, \\
  \yy_t &= t(t-1) \ee_1 + t \ee_2, \\
  \zz_t &= \biggParens{1-\frac{1}{t}} \xx_t + \frac{1}{t} \yy_t.
\end{align*}
Then $\xx_t\rightarrow\xbar$ and $\yy_t\rightarrow\ybar$
(by \Cref{thm:seq-rep}).
Also,
plugging the definitions of $\xx_t$ and $\yy_t$ into the expression
for $\zz_t$,
we obtain that
$\zz_t=\ee_2$ for all $t$, so
$\zz_t\rightarrow\zz$.
Thus, as claimed, $\zz\in\lb{\xbar}{\ybar}$
by \Cref{cor:e:1}, so $\Cbar$ is not convex.
\end{proof}

Thus, the closure of a convex set need not be 
\indexg{closure, astral!convex set@of convex set|)}%
\indexg{convex sets, astral!closure of|)}%
convex.
\indexg{convex hull, astral!closedness preserved|(}%
On the other hand, the convex hull of a closed set is always closed:

\begin{theorem}  \label{thm:cnv-hull-closed-is-closed}
  Let $S\subseteq\extspace$ be closed (in $\extspace$).
  Then $\conv{S}$ is also closed.
\end{theorem}

\begin{proof}
If $S=\emptyset$, then the claim is immediate,
so we assume $S\neq\emptyset$.

Let $\zbar\in\clbar{\conv{S}}$; we aim to show $\zbar\in\conv{S}$,
thereby proving that $\conv{S}$ is closed.
Since $\zbar\in\clbar{\conv{S}}$, there exists a sequence $\seq{\zbar_t}$ in $\conv{S}$ with
$\zbar_t\rightarrow\zbar$.
For each $t$,
by \Cref{thm:convhull-of-simpices},
$\zbar_t$ is in the outer hull of at most $n+1$ elements of $S$.
That is,
$\zbar_t\in\ohull\{\xbar_{1t},\dotsc,\xbar_{(n+1)t}\}$
for some $\xbar_{1t},\dotsc,\xbar_{(n+1)t}\in S$
(where we here allow the same point to appear repeatedly
so that each $\zbar_t$ is in the outer hull of exactly
$n+1$ not necessarily distinct points).

We can consider each sequence $\seq{\xbar_{it}}$ in turn,
for $i=1,\dotsc,n+1$.
By sequential compactness,
the sequence $\seq{\xbar_{it}}$ must have a convergent subsequence;
by discarding the other elements (and the corresponding elements
of the other sequences), we can ensure the entire sequence converges
so that $\xbar_{it}\rightarrow\xbar_i$ for some $\xbar_i\in\extspace$.
Moreover, since $S$ is closed, $\xbar_i\in S$.

Let $\uu\in\Rn$.
By
\Cref{pr:ohull-simplify}(\ref{pr:ohull-simplify:a},\ref{pr:ohull-simplify:b}),
for each $t$,
\begin{equation}  \label{eq:thm:cnv-hull-closed-is-closed:1}
  \zbar_t\cdot\uu
  \leq
  \max\{\xbar_{1t}\cdot\uu,\dotsc,\xbar_{(n+1)t}\cdot\uu\}.
\end{equation}
By \Cref{thm:i:1}(\ref{thm:i:1c}),
$\zbar_t\cdot\uu\rightarrow\zbar\cdot\uu$
and
$\xbar_{it}\cdot\uu\rightarrow\xbar_i\cdot\uu$
for $i=1,\dotsc,n+1$.
Consequently, taking limits of both sides of
\eqref{eq:thm:cnv-hull-closed-is-closed:1} yields
\[
  \zbar\cdot\uu
  \leq
  \max\{\xbar_{1}\cdot\uu,\dotsc,\xbar_{n+1}\cdot\uu\}
\]
since the maximum of $n+1$ numbers, as a function, is continuous.
Therefore,
$\zbar\in\ohull\{\xbar_1,\dotsc,\xbar_{n+1}\}\subseteq\conv{S}$
by
\Cref{pr:ohull-simplify}(\ref{pr:ohull-simplify:b},\ref{pr:ohull-simplify:a})
and
\Cref{thm:convhull-of-simpices}.%
\indexg{convex hull, astral!closedness preserved|)}%
\end{proof}

\chapter{Sequential operations}
\label{sec:sequential}

Motivated by the characterization of astral polytopes via sequences
(\Cref{thm:e:7}),
we next introduce and study a class of
operations based on sequences which will provide astral analogues of
several standard operations, namely, the standard sum of sets in
$\Rn$, scalar-vector products, and convex and conic combinations of
points in $\Rn$.
We will see that these generalized operations have numerous favorable properties,
and so will find frequent use.
For instance, the sequential astral operation for adding sets,
which generalizes ordinary addition of sets in $\Rn$,
is commutative (unlike leftward addition),
preserves both convexity and closedness,
and commutes with any linear map.

We will also see that
sequential operations relate
to other notions that are defined or characterized in terms of
sequences, such as astral polytopes, characterized in
\Cref{thm:e:7}.
In particular, we will see that
sequential convex and conic combinations
provide an algebra for expressing and working
with the detailed structure of every astral polytope, and so also
of the convex hull of any set in $\extspace$
(via \Cref{thm:convhull-of-simpices}).
Moreover, they can be used to describe the structure of the astral conic hull of a set,
which will be defined and studied later in \Cref{sec:conic:hull}.

\section{Sequential sum of sets}
\label{sec:seq-sum}

\indexg{sequential sum|(}%
Suppose $\seq{\xx_t}$ and $\seq{\yy_t}$ are two sequences in $\Rn$
converging, respectively, to astral points $\xbar$ and $\ybar$.
What can be said about the limit (if it exists) of the sum of the
sequences,
$\seq{\xx_t+\yy_t}$?
More generally, if $\seq{\xx_t}$ and $\seq{\yy_t}$ each have limits,
respectively, in some sets $X$ and $Y$ in $\extspace$, then
what can be said about the limit of
$\seq{\xx_t+\yy_t}$?

To answer this, we study the set of all such limits, called
the sequential sum:

\begin{definition}  \label{dfn:seq-sum}
\indexg{sequential sum!defined|(}%
  The \emph{sequential sum} of two sets $X,Y\subseteq\extspace$,
  denoted $X\seqsum Y$,%
\indexm{x+y900}{$X\protect\seqsum Y$}{sequential sum}
  is the set of all points $\zbar\in\extspace$ for
  which there exist sequences $\seq{\xx_t}$ and $\seq{\yy_t}$ in $\Rn$
  and points $\xbar\in X$ and $\ybar\in Y$ such that
  $\xx_t\rightarrow\xbar$,
  $\yy_t\rightarrow\ybar$,
  and
  $\xx_t+\yy_t\rightarrow \zbar$.%
\indexg{sequential sum!defined|)}%
\end{definition}

If $X$ and $Y$ are in $\Rn$, then $X\seqsum Y$ is simply the sum of
the sets, $X+Y$,
so we can view sequential addition of astral sets as a generalization
of ordinary addition of sets in~$\Rn$.

When taking sequential sums involving
singletons, as in the example below,
we usually streamline notation and omit braces, for
instance, writing
$\xbar\seqsum \ybar$ for $\{\xbar\} \seqsum \{\ybar\}$.
It is important, nonetheless, to keep in mind that sequential sum
always denotes a set in~$\extspace$,
even if that set happens to be a singleton.

\begin{example}   \label{ex:oe1-seqsum-oe2}
In $\extspac{2}$, let $\xbar=\limray{\ee_1}$ and
$\ybar=\limray{\ee_2}$.
Each line below gives example sequences $\seq{\xx_t}$ and
$\seq{\yy_t}$ converging respectively to $\xbar$ and $\ybar$,
along with the limit of their sum, $\seq{\xx_t+\yy_t}$;
every such limit must be in $\xbar\seqsum\ybar$.

{%
\squeezebin%
\begin{align*}
  \xx_t &= t \ee_1,
  &
  \yy_t &= t \ee_2,
  &
  \xx_t + \yy_t &= t (\ee_1 + \ee_2) \rightarrow \limray{(\ee_1 + \ee_2)}.
  \\
  \xx_t &= t^2 \ee_1,
  &
  \yy_t &= t \ee_2,
  &
  \xx_t + \yy_t &= t^2 \ee_1 + t\ee_2 \rightarrow \limray{\ee_1}\plusl\limray{\ee_2}.
  \\
  \xx_t &= t \ee_1,
  &
  \yy_t &= t^2 \ee_2,
  &
  \xx_t + \yy_t &= t^2 \ee_2 + t\ee_1 \rightarrow \limray{\ee_2}\plusl\limray{\ee_1}.
  \intertext{%
    Generalizing the first line, we also have, for all
    $\vv=\trans{[v_1,v_2]}\in\Rstrictpos^2$
    and all $\ww=\trans{[w_1,w_2]}\in\R^2$:
  }
  \xx_t &= (t v_1 + w_1) \ee_1,
  &
  \yy_t &= (t v_2 + w_2) \ee_2,
  &
  \xx_t + \yy_t &= t \vv + \ww \rightarrow \limray{\vv}\plusl\ww.
  \\  
  \xx_t &= (t^2 v_1 + t w_1) \ee_1,
  &
  \yy_t &= (t^2 v_2 + t w_2) \ee_2,
  &
  \xx_t + \yy_t &= t^2 \vv + t \ww \rightarrow \limray{\vv}\plusl\limray{\ww}.
\end{align*}%
}

Thus, $\xbar\seqsum\ybar$ includes
$\limray{\ee_1}\plusl\limray{\ee_2}$,
$\limray{\ee_2}\plusl\limray{\ee_1}$,
and
$\limray{\vv}\plusl\extspac{2}$, for all $\vv\in\Rstrictpos^2$.
Later, in \Cref{ex:oe1-seqsum-oe2:cont},
we will see that these points comprise the
entirety of $\xbar\seqsum\ybar$.
\end{example}

Here are some simple properties of sequential addition, whose proofs
are routine from its definition:

\begin{proposition}   \label{pr:seqsum-props}
  ~

  \begin{letter-compact}
  \item   \label{pr:seqsum-props:finite}
    For $X,Y\subseteq\Rn$,
    $X \seqsum Y = X + Y$.
  \item   \label{pr:seqsum-props:a}
    For $X,Y\subseteq\extspace$,
    $X \seqsum Y = Y \seqsum X$.
  \item   \label{pr:seqsum-props:identity}
    For $X\subseteq\extspace$,
    $X \seqsum \{\zero\}= \{\zero\} \seqsum X = X$,
    and
    $X\seqsum \emptyset = \emptyset \seqsum X = \emptyset$.
  \item   \label{pr:seqsum-props:b}
    Suppose
    $X\subseteq X'\subseteq\extspace$
    and
    $Y\subseteq Y'\subseteq\extspace$.
    Then
    $X \seqsum Y \subseteq X' \seqsum Y'$.
  \item   \label{pr:seqsum-props:c}
    Let $X_i$, for $i\in I$, and $Y_j$, for $j\in J$,
    be arbitrary families of sets in $\extspace$,
    where $I$ and $J$ are index sets.
    Then
    \[
       \BiggParens{\,\bigcup_{i\in I} X_i}
       \seqsum
       \BiggParens{\,\bigcup_{j\in J} Y_j}
       =
       \bigcup_{i\in I,\,j\in J} (X_i \seqsum Y_j).
    \]
  \end{letter-compact}
\end{proposition}

\indexg{m-ary sequential sum@$m$-ary sequential sum|(}%
\indexg{sequential sum!mary@$m$-ary|(}%
To establish several key properties of sequential addition,
we consider its generalization
to finite collections of sets:

\begin{definition}  \label{dfn:m-ary-seq-sum}
\indexg{m-ary sequential sum@$m$-ary sequential sum!defined|(}%
  For $m\ge 1$,
  the \emph{$m$-ary sequential sum} of sets $X_1,\dotsc,X_m\subseteq\eRn$,
  denoted $\Seqsum(X_1,\dotsc,X_m)$,%
\indexm{sigma+}{$\protect\Seqsum(\ldots)$}{$m$-ary sequential sum}
  is the set of all points
  $\zbar\in\extspace$ for which there exist, for $i=1,\dotsc,m$,
  a sequence $\seq{\xx_{it}}$ in $\Rn$ and a point $\xbar_i\in X_i$ 
  such that $\xx_{it}\to\xbar_i$, and such that
  $\sum_{i=1}^m\xx_{it}\to\zbar$.
  When $m=0$, the $0$-ary sequential sum $\Seqsum(\,)$ is defined to
  be
\indexg{m-ary sequential sum@$m$-ary sequential sum!defined|)}%
  $\{\zero\}$.
\end{definition}

For $m=1$ and $m=2$, this definition yields $\Seqsum(X)=X$ and
\begin{equation}   \label{eq:Seqsumxy-is-xseqsumy}
  \Seqsum(X,Y)=X\seqsum Y,
\end{equation}
for all $X,Y\subseteq\extspace$.
Later, we will show, as expected, that
\begin{equation}   \label{eq:Seqsumxm-is-seqsumxm}
  \Seqsum(X_1,\dotsc,X_m)=X_1\seqsum \cdots \seqsum X_m,
\end{equation}
for $m\geq 1$, but this will require proof.
(Indeed, we have not yet even established that sequential sum is
associative, so this latter expression is, at this point, somewhat ambiguous.)

We next characterize $m$-ary sequential sums using linear maps.
Recall that for vectors $\xx_1,\dotsc,\xx_m\in\Rn$, we write
$\mtuple{\xx_1,\dotsc,\xx_{m}}$ for the (column) vector in $\R^{mn}$
obtained by concatenating these vectors.
For $i=1,\dotsc,m$, we then define the
\indexg{component projection matrices!tuples in r mn@for tuples in $\R^{mn}$}%
\emph{$i$-th component projection matrix}
$\PPi\in\R^{n\times mn}$ as
\begin{equation}  \label{eq:comp-proj-mat-defn}
  \PPi
  =
  [
    \underbrace{\zero,\dotsc,\zero}_{i-1},
    \Iden,
    \underbrace{\zero,\dotsc,\zero}_{m-i}
  ],
\end{equation}
where $\Iden$ is the $n\times n$ identity matrix, and, in this context,
$\zero$ is the $n\times n$ all-zeros matrix. Thus,
$\PPi \mtuple{\xx_1,\dotsc,\xx_{m}} = \xx_i$, for
$i=1,\dotsc,m$. We also define the
\emph{addition matrix}%
\indexg{addition matrix}
$\A\in\R^{n\times mn}$ as $\A=[\Iden,\dotsc,\Iden]$.
Thus,
\begin{equation}  \label{eq:add-mat-defn}
  \A=\sum_{i=1}^m\PPi,
\end{equation}
so
$\A \mtuple{\xx_1,\dotsc,\xx_{m}} =\sum_{i=1}^m\xx_i$.

Using the language of linear maps, the next \lcnamecref{thm:seqsum-mat-char} shows how sequential addition can be viewed as an extension of standard vector addition:

\begin{theorem}
\label{thm:seqsum-mat-char}
Let $X_1,\dotsc,X_m\subseteq\eRn$,
let $\PPsub{1},\dotsc,\PPsub{m}$
be the component projection matrices,
and let
$\A$
be
the addition matrix
(as in Eqs.~\ref{eq:comp-proj-mat-defn}
and~\ref{eq:add-mat-defn}).
Then
\begin{align}
\notag
&\Seqsum(X_1,\dotsc,X_m)
\\
\label{eq:prop:seqsum}
&\qquad{}
  =\bigSet{\A\wbar:\:\wbar\in\eRf{mn}
          \text{ such that }\PPi\wbar\in X_i\text{ for }i=1,\dotsc,m
  }.
\end{align}
\end{theorem}
\begin{proof}
Let $N=mn$,
let $Z=\Seqsum(X_1,\dotsc,X_m)$, and let $U$ denote the right-hand side of
\eqref{eq:prop:seqsum}.
We will show that $Z=U$.

First, let $\zbar\in Z$.
Then there exist sequences $\seq{\xx_{it}}$ in $\Rn$, for
$i=1,\dotsc,m$, such that
$\xx_{it}\to\xbar_i$ for some $\xbar_i\in X_i$,
and such that
$\sum_{i=1}^m \xx_{it}\to\zbar$.
Let $\ww_t=\mtuple{\xx_{1t},\dotsc,\xx_{mt}}$.
By sequential compactness, the sequence $\seq{\ww_t}$ must have a
convergent subsequence; by discarding all other elements
as well as the corresponding elements of the sequences~$\seq{\xx_{it}}$,
we can
assume the entire sequence converges to some point
$\wbar\in\eRf{N}$. By continuity of linear maps
(\Cref{thm:linear:cont}\ref{thm:linear:cont:b}),
$\PPi \wbar = \lim(\PPi\ww_t)=\lim \xx_{it} = \xbar_i \in X_i$ for all $i$,
while $\zbar=\lim(\sum_{i=1}^m \xx_{it})=\lim(\A\ww_t)=\A\wbar$.
Therefore, $\zbar\in U$, so $Z\subseteq U$.

For the reverse inclusion, let $\zbar\in U$. Then
there exists $\wbar\in\eRf{N}$
such that $\A\wbar=\zbar$ and
$\PPi \wbar=\xbar_i$ for some $\xbar_i \in X_i$, for $i=1,\dotsc,m$.
Since $\wbar\in\eRf{N}$, there exists a sequence $\seq{\ww_t}$ in $\R^N$
such that $\ww_t\rightarrow\wbar$.
For $i=1,\dotsc,m$,
setting $\xx_{it}=\PPi \ww_t$, we obtain
$\xx_{it}=\PPi \ww_t \to \PPi \wbar = \xbar_i\in X_i$.
Also, $\sum_{i=1}^m\xx_{it}=\A\ww_t\to\A\wbar=\zbar$.
Thus, $\zbar\in Z$, so $Z=U$.
\end{proof}

Using the characterization in \Cref{thm:seqsum-mat-char}, we can show that sequential addition is associative, as an immediate corollary of the following lemma.

\begin{lemma}
\label{lemma:seqsum:assoc}
Let $X_1,\dotsc,X_m,Y\subseteq\eRn$. Then
\[
  \Seqsum(X_1,\dotsc,X_m)\seqsum Y=
  \Seqsum(X_1,\dotsc,X_m,Y).
\]
\end{lemma}
\begin{proof}

Let $Z=\Seqsum(X_1,\dotsc,X_m)$, let
$V=Z\seqsum Y$, and let $Z'=\Seqsum(X_1,\dotsc,X_m,Y)$.

We will show that $V=Z'$.
First, let $\zbar'\in Z'$.
Then there exist sequences $\seq{\xx_{it}}$ in $\Rn$, for $i=1,\dotsc,m$,
such that
$\xx_{it}\to\xbar_i$ for some $\xbar_i\in X_i$,
and a sequence $\seq{\yy_t}$ in $\Rn$,
such that $\yy_t\to\ybar$ for some $\ybar\in Y$,
and such that
$(\sum_{i=1}^m \xx_{it})+\yy_t\to\zbar'$.
Let $\zz_t=\sum_{i=1}^m \xx_{it}$.
By sequential compactness, the sequence $\seq{\zz_t}$ must have a
convergent subsequence; by discarding all other elements
as well as the corresponding elements of the sequences
$\seq{\xx_{it}}$ and~$\seq{\yy_t}$,
we can
assume the entire sequence converges to some point $\zbar\in\eRn$.
Since $\xx_{it}\to\xbar_i\in X_i$ and $(\sum_{i=1}^m \xx_{it})=\zz_t\to\zbar$,
we have $\zbar\in Z$. Moreover,
$\yy_t\to\ybar\in Y$, and $\zz_t+\yy_t=(\sum_{i=1}^m \xx_{it})+\yy_t\to\zbar'$,
so $\zbar'\in Z\seqsum Y=V$. Therefore, $Z'\subseteq V$.

For the reverse inclusion, let $\vbar\in V$. By \Cref{thm:seqsum-mat-char},
\[
  V=\bigSet{\B\ubar:\:\ubar\in\eRf{2n}\text{ such that }
    \QQz\ubar\in Z\text{ and }\QQy\ubar\in Y},
\]
where $\B=[\Iden,\Iden]$, $\QQz=[\Iden,\zero]$, and $\QQy=[\zero,\Iden]$ (writing $\Iden=\Idnn$ and $\zero=\zero_{n\times n}$ in this context).
Therefore, there exists $\ubar\in\eRf{2n}$ such that $\B\ubar=\vbar$, and also
$\QQz\ubar=\zbar$ for some $\zbar\in Z$, and $\QQy\ubar=\ybar$ for some
$\ybar\in Y$.

Let $N=mn$. Then, by \Cref{thm:seqsum-mat-char}, we also have
\[
  Z=\bigSet{\A\wbar:\:\wbar\in\eRf{N}
          \text{ such that }\PPi\wbar\in X_i\text{ for }i=1,\dotsc,m
  },
\]
with the addition matrix $\A\in\R^{n\times N}$ and the component projection matrices
$\PPi\in\R^{n\times N}$ for $i=1,\dotsc,m$.
Since $\zbar\in Z$, there exists
$\wbar\in\eRf{N}$ such that $\A\wbar=\zbar$ and $\PPi\wbar=\xbar_i$ for some $\xbar_i\in X_i$,
for $i=1,\dotsc,m$.

To finish the proof, we construct sequences $\seq{\xx_{it}}$ in $\Rn$, for $i=1,\dotsc,m$,
and $\seq{\yy_t}$ in $\Rn$ whose properties will imply that $\vbar\in Z'$.
To start, let $\seq{\uu_t}$ be a sequence in $\R^{2n}$ such that $\uu_t\to\ubar$. By continuity of linear maps
(\Cref{thm:linear:cont}\ref{thm:linear:cont:b}),
$\QQz\uu_t\to\QQz\ubar=\zbar=\A\wbar$. Therefore, by \Cref{thm:inv-lin-seq},
there exists a sequence $\seq{\ww_t}$ in $\R^N$ such that
$\ww_t\to\wbar$,
and
$\seqeq{\A\ww_t}{\QQz\uu_t}$.
Let $\veps_t=\A\ww_t-\QQz\uu_t$,
implying $\veps_t\rightarrow\zero$,
and let $\xx_{it}=\PPi\ww_t$,
for $i=1,\dotsc,m$,
and $\yy_t=\QQy\uu_t$.
By the foregoing construction and continuity
of linear maps,
\begin{align*}
&
  \xx_{it}=\PPi\ww_t\to\PPi\wbar=\xbar_i\quad\text{for $i=1,\dotsc,m$,}
\\
&
  \yy_t=\QQy\uu_t\to\QQy\ubar=\ybar,
\\
&
  \bigParens{\textstyle\sum_{i=1}^m\xx_{it}}+\yy_t
  =
  \A\ww_t+\yy_t
  =
  \QQz\uu_t+\veps_t+\QQy\uu_t
  =
  \B\uu_t+\veps_t
  \to
  \B\ubar
  =
  \vbar,
\end{align*}
where the convergence in the last line follows by \Cref{pr:i:7}(\ref{pr:i:7g}).
Thus, $\vbar\in Z'$, so $V=Z'$.
\end{proof}

\begin{corollary}
\label{cor:seqsum:assoc}
\indexg{sequential sum!associativity of|(}%
Let $X,Y,Z\subseteq\eRn$. Then
\[
  (X\seqsum Y)\seqsum Z = X\seqsum(Y\seqsum Z).
\]
\end{corollary}

\begin{proof}
We have
\[
  (X\seqsum Y)\seqsum Z
  =
  \Seqsum(X,Y,Z)
  =
  \Seqsum(Y,Z,X)
  =
  (Y\seqsum Z)\seqsum X
  =
  X\seqsum(Y\seqsum Z),
\]
The first and third equalities are by \Cref{lemma:seqsum:assoc}
(and Eq.~\ref{eq:Seqsumxy-is-xseqsumy}).
The second follows from the definition of
$m$-ary sequential sum.
The last is by commutativity of sequential addition
(\Cref{pr:seqsum-props}\ref{pr:seqsum-props:a}).
\end{proof}

Having established associativity of (binary) sequential addition, we can now omit parentheses and simply write
$X\seqsum Y\seqsum Z$, or more generally $X_1\seqsum\dotsb\seqsum X_m$
for any sets $X_1,\dotsc,X_m\subseteq\eRn$.
When $m=0$, such an expression is understood to be equal to~$\{\zero\}$.%
\indexg{sequential sum!associativity of|)}

The next theorem
shows that such a chain of binary sequential sums is equal to the
$m$-ary sequential sum,
in other words, that \eqref{eq:Seqsumxm-is-seqsumxm} holds.
It also shows that a sequential sum can be decomposed into sequential sums of singletons, and that sequential addition preserves closedness and convexity of its arguments.

\begin{theorem}
\label{prop:seqsum-multi}
Let $Z=X_1\seqsum\dotsb\seqsum X_m$ for some $X_1,\dotsc,X_m\subseteq\eRn$.
Then:
\begin{letter}
\item \label{prop:seqsum-multi:equiv}
  $Z = \Seqsum(X_1,\dotsc,X_m)$.
\item \label{prop:seqsum-multi:decomp}
$
\displaystyle
  Z = \bigcup_{\xbar_1\in X_1, \dotsc, \xbar_m\in X_m}
                   (\xbar_1 \seqsum \dotsb \seqsum \xbar_m)
$.
\item \label{prop:seqsum-multi:nonempty}
  $Z$ is nonempty if the sets $X_1,\dotsc,X_m$ are nonempty.
\item \label{prop:seqsum-multi:closed}
\indexg{sequential sum!closedness preserved|(}%
  $Z$ is closed if the sets $X_1,\dotsc,X_m$ are closed.
\item \label{prop:seqsum-multi:convex}
\indexg{sequential sum!convexity preserved|(}%
  $Z$ is convex if the sets $X_1,\dotsc,X_m$ are convex.
\end{letter}
\end{theorem}

\begin{proof}
The claims are straightforward to check when $m=0$, so we assume
$m\geq 1$.
(Note in part~(\ref{prop:seqsum-multi:decomp}) that when $m=0$,
the union on the right-hand side is over a single empty tuple, for which 
the value of $\xbar_1 \seqsum \dotsb \seqsum \xbar_m$ is equal to $\set{\zero}$
by our convention.)

\begin{proof-parts}
\pfpart{Part~(\ref{prop:seqsum-multi:equiv}):}
Proof is by induction on $m$. By definition of $m$-ary sequential sum, the claim holds for $m=1$. Let $m>1$ and assume the claim holds for $m-1$. Then
\[
 Z=(X_1\seqsum\dotsb\seqsum X_{m-1})\seqsum X_m
  =\Seqsum(X_1,\dotsc,X_{m-1})\seqsum X_m
  =\Seqsum(X_1,\dotsc,X_m),
\]
The second equality is by inductive hypothesis, and the third is by \Cref{lemma:seqsum:assoc}.%
\indexg{m-ary sequential sum@$m$-ary sequential sum|)}%
\indexg{sequential sum!mary@$m$-ary|)}

\pfpart{Part~(\ref{prop:seqsum-multi:decomp}):}
Proof is again by induction on $m$.
The base case that $m=1$ is immediate.
For the inductive step, let $m>1$ and assume that the claim holds for
$m-1$.
Then
\begin{align*}
Z&=
  \parens{ X_1 \seqsum \dotsb \seqsum X_{m-1} } \seqsum X_m
\\
&=
  \BiggParens{\,\bigcup_{\xbar_1\in X_1, \dotsc, \xbar_{m-1}\in X_{m-1}}
                (\xbar_1 \seqsum \dotsb \seqsum \xbar_{m-1}) }
  \seqsum
  \BiggParens{\,\bigcup_{\xbar_{m}\in X_{m}} \{ \xbar_m \} }
\\
&=
  \bigcup_{\xbar_1\in X_1, \dotsc, \xbar_{m-1}\in X_{m-1}, \xbar_m\in X_m}
       \BigParens{(\xbar_1 \seqsum \dotsb \seqsum \xbar_{m-1}) \seqsum \xbar_m}.
\end{align*}
The second equality is by inductive hypothesis,
and the third is by
\Cref{pr:seqsum-props}(\ref{pr:seqsum-props:c}).

\pfpart{Part~(\ref{prop:seqsum-multi:nonempty}):}
For $i=1,\dotsc,m$, suppose each set $X_i$ includes some element $\xbar_i$,
and let $\seq{\xx_{it}}$ be a sequence in $\Rn$
converging to $\xbar_i$ (which exists by
\Cref{thm:i:1}\ref{thm:i:1d}).
For each~$t$, let $\zz_t=\sum_{i=1}^m \xx_{it}$.
By sequential compactness, the sequence $\seq{\zz_t}$ has a convergent
subsequence; by discarding all other elements (and corresponding
elements of the other sequences), we can assume the entire sequence
converges to some point $\zbar\in\extspace$.
Then by \Cref{dfn:m-ary-seq-sum},
$\zbar$ is in $\Seqsum(X_1,\dotsc,X_m)$, and thus in $Z$ by
part~(\ref{prop:seqsum-multi:equiv}).
Therefore, $Z$ is nonempty.

\pfpart{Part~(\ref{prop:seqsum-multi:closed}):} By part~(\ref{prop:seqsum-multi:equiv}) and \Cref{thm:seqsum-mat-char},
\begin{equation}
\label{eq:seqsum:1}
  Z
  = \bigSet{\A\wbar:\:\wbar\in\eRf{mn}
          \text{ such that }\PPi\wbar\in X_i\text{ for }i=1,\dotsc,m
  },
\end{equation}
where $\A$
is the addition matrix and
$\PPsub{1},\dotsc,\PPsub{m}$
are the component projection matrices
(as in Eqs.~\ref{eq:comp-proj-mat-defn}
and~\ref{eq:add-mat-defn}).
Let $A:\eRf{mn}\to\eRn$ and $P_i:\eRf{mn}\to\eRn$ denote the
associated astral linear maps (so that $A(\wbar)=\A\wbar$ and
$P_i(\wbar)=\PPi\wbar$ for $\wbar\in\eRf{mn}$).
Then \eqref{eq:seqsum:1} can be rewritten as
\begin{equation}
\label{eq:seqsum:2}
  Z = A\BiggParens{\,\bigcap_{i=1}^m P_i^{-1}(X_i)}.
\end{equation}

Suppose that each $X_i$ for $i=1,\dotsc,m$ is closed. Then, by
continuity of $P_i$
(\Cref{thm:linear:cont}\ref{thm:linear:cont:b}),
each of the sets $P_i^{-1}(X_i)$ is closed as well
(\Cref{prop:cont}\ref{prop:cont:a}\ref{prop:cont:inv:closed}),
and so is their intersection, and also the image of the intersection under $A$, by
\indexg{sequential sum!closedness preserved|)}%
\Cref{cor:aff-img-closed-is-closed}(\ref{cor:aff-img-closed-is-closed:a}).

\pfpart{Part~(\ref{prop:seqsum-multi:convex}):} Write $Z$ as in \eqref{eq:seqsum:2} and suppose that each of the sets $X_i$ is convex.
Then each of the sets $P_i^{-1}(X_i)$ is also convex
(\Cref{cor:inv-image-convex}), so their intersection is convex
(\Cref{pr:e1}\ref{pr:e1:b}), implying that $Z$, the image of that
intersection under $A$, is as well (\Cref{cor:thm:e:9}).
\qedhere
\end{proof-parts}
\end{proof}

Thus, if $X$ and $Y$ are convex subsets of $\extspace$, then
$X\seqsum Y$ is also convex.
The same is not true, in general, for $X\plusl Y$, as the next example
shows.
We will say more about the relationship between $X\seqsum Y$ and
$X\plusl Y$ later in this section.

\begin{example}
In $\R^2$,
let $X=\lb{\zero}{\limray{\ee_1}}$ and
$Y=\lb{\zero}{\limray{\ee_2}}$, which are both convex.
By \Cref{thm:lb-with-zero},
$X=\{\lambda \ee_1 : \lambda\in\Rpos\} \cup \{\limray{\ee_1}\}$
and similary for $Y$,
so
$\limray{\ee_1}$ and $\limray{\ee_2}$ are in $X\plusl Y$,
but $\limray{\ee_2}\plusl\limray{\ee_1}$ is not.
On the other hand, this latter point is on the segment joining
$\limray{\ee_1}$ and $\limray{\ee_2}$,
as we saw in \Cref{ex:seg-oe1-oe2}.
Therefore, $X\plusl Y$ is not convex.%
\indexg{sequential sum!convexity preserved|)}
\end{example}

\indexg{sequential sum!characterizations of|(}%
The next theorem gives various characterizations for the sequential
sum of a tuple of astral points, $\xbar_1,\dotsc,\xbar_m$
(or, more precisely, of the singletons associated with these points).
Implicitly, it also characterizes the sequential sum of a tuple of
arbitrary sets, since these can be decomposed in terms of singletons using
\Cref{prop:seqsum-multi}(\ref{prop:seqsum-multi:decomp}).
The first two characterizations, following immediately from the
preceding development, are in terms of sequences, and in
terms of an addition matrix and component projection matrices.
The next characterization is in terms of the behavior of the
associated functions $\uu\mapsto \xbar_i\cdot\uu$.
The final characterization shows that the sequential sum
$\xbar_1 \seqsum \dotsb \seqsum \xbar_m$
is actually an astral
polytope whose vertices are formed by
permuting the
ordering of the $\xbar_i$'s, in all possible ways, and then combining
them using leftward addition.
(In the theorem, a \emph{permutation of a set $S$}
is simply a bijection
$\pi:S\rightarrow S$.)

\begin{theorem}    \label{thm:seqsum-equiv-mod}
  Let $\xbar_1,\dotsc,\xbar_m\in\extspace$,
  and let $\zbar\in\extspace$.
  Also, let $\PPsub{1},\dotsc,\PPsub{m}$
  be the component projection matrices,
  and let
  $\A$ be the addition matrix
  (as in Eqs.~\ref{eq:comp-proj-mat-defn}
  and~\ref{eq:add-mat-defn}).
  Then the following are equivalent:
  \begin{letter-compact}
  \item      \label{thm:seqsum-equiv-mod:a}
    $\zbar \in \xbar_1 \seqsum \dotsb \seqsum \xbar_m$.
  \item      \label{thm:seqsum-equiv-mod:seq}
    There exist sequences $\seq{\xx_{it}}$ in $\Rn$ such that
    $\xx_{it}\rightarrow\xbar_i$, for $i=1,\dotsc,m$,
    and such that
    $\sum_{i=1}^m \xx_{it} \rightarrow\zbar$.
  \item      \label{thm:seqsum-equiv-mod:b}
    $\zbar =\A\wbar$ for some $\wbar\in\eRf{mn}$ such that $\PPi\wbar=\xbar_i$ for $i=1,\dotsc,m$.
  \item      \label{thm:seqsum-equiv-mod:c}
    For all $\uu\in\Rn$,
    if $\xbar_1\cdot\uu,\dotsc,\xbar_m\cdot\uu$ are summable,
    then
    $\zbar\cdot\uu = \sum_{i=1}^m \xbar_i\cdot\uu$.
  \item      \label{thm:seqsum-equiv-mod:d}
    $ \zbar
      \in
      \ohull\regBraces{
        \xbar_{\pi(1)}\plusl \dotsb \plusl \xbar_{\pi(m)}
        :\:
        \pi\in \Pi_m
      }
    $
    where $\Pi_m$ is the set of all permutations of $\{1,\dotsc,m\}$.
  \end{letter-compact}
\end{theorem}

\begin{proof}
~
\begin{proof-parts}
\pfpart{%
  (\ref{thm:seqsum-equiv-mod:a})
  $\Leftrightarrow$
  (\ref{thm:seqsum-equiv-mod:seq}):
}
Immediate from \Cref{prop:seqsum-multi}(\ref{prop:seqsum-multi:equiv}),
noting that the set of points $\zbar$ satisfying
part~(\ref{thm:seqsum-equiv-mod:seq}) is precisely
$\Seqsum(\{\xbar_1\},\dotsc,\{\xbar_m\})$,
by its definition.

\pfpart{%
  (\ref{thm:seqsum-equiv-mod:seq})
  $\Leftrightarrow$
  (\ref{thm:seqsum-equiv-mod:b}):
}
Immediate from \Cref{thm:seqsum-mat-char}.

\pfpart{%
  (\ref{thm:seqsum-equiv-mod:b})
  $\Rightarrow$
  (\ref{thm:seqsum-equiv-mod:c}):
}
Suppose, for some $\wbar\in\eRf{mn}$,
that $\zbar=\A\wbar$ and $\xbar_i=\PPi\wbar$ for $i=1,\dotsc,m$.
Let $\uu\in\Rn$, and suppose that
$\xbar_1\cdot\uu,\dotsc,\xbar_m\cdot\uu$ are summable.
By
\Cref{thm:Ax-dot-u},
$\xbar_i\cdot\uu=(\PPi\wbar)\cdot\uu=\wbar\cdot(\trans{\PPi}\uu)$
for all $i$,
and $\zbar\inprod\uu=(\A\wbar)\inprod\uu=\wbar\inprod(\transA\uu)$.
Thus,
\begin{align*}
\sum_{i=1}^m \xbar_i\inprod\uu
&=
  \sum_{i=1}^m \wbar\inprod(\trans{\PPi}\uu)
 =\wbar\inprod\BiggParens{\sum_{i=1}^m \trans{\PPi}\uu}
 =\wbar\inprod(\transA\uu)
 =\zbar\inprod\uu,
\end{align*}
where
the second equality is by an iterated application of \Cref{pr:i:1}
(having assumed summability), and the third by
\eqref{eq:add-mat-defn}.

\pfpart{%
  (\ref{thm:seqsum-equiv-mod:c})
  $\Rightarrow$
  (\ref{thm:seqsum-equiv-mod:d}):
}
Suppose part~(\ref{thm:seqsum-equiv-mod:c}) holds.
Let $\uu\in\Rn$.
To prove part~(\ref{thm:seqsum-equiv-mod:d}),
by
\Cref{pr:ohull-simplify}(\ref{pr:ohull-simplify:b},\ref{pr:ohull-simplify:a}),
it suffices to show that
\begin{equation}
\label{eq:thm:seqsum-equiv-mod:1}
  \zbar\cdot\uu
  \leq
  \max\BigBraces{
    \bigParens{\xbar_{\pi(1)}\plusl \dotsb \plusl \xbar_{\pi(m)}}\cdot\uu
    :\:
    \pi\in \Pi_m
  }.
\end{equation}

Suppose first that $\xbar_i\cdot\uu=+\infty$ for some
$i\in\{1,\dotsc,m\}$.
In this case, let $\pi$ be any permutation in $\Pi_m$ for which
$\pi(1)=i$.
Then
$(\xbar_{\pi(1)}\plusl \dotsb \plusl \xbar_{\pi(m)})\cdot\uu =+\infty$,
implying \eqref{eq:thm:seqsum-equiv-mod:1} since its right-hand side must
also be $+\infty$.

Otherwise, we must have $\xbar_i\cdot\uu<+\infty$ for
$i=1,\dotsc,m$, and so
${\xbar_1\cdot\uu},\dotsc,{\xbar_m\cdot\uu}$ are summable, with their sum being equal
to $\zbar\cdot\uu$ (by assumption).
Thus, for every $\pi\in\Pi_m$, we must have
\[
  \zbar\cdot\uu
  =
  \sum_{i=1}^m \xbar_i\cdot\uu
  =
  \xbar_{\pi(1)}\cdot\uu \plusl \dotsb \plusl \xbar_{\pi(m)}\cdot\uu
  =
  \bigParens{\xbar_{\pi(1)} \plusl \dotsb \plusl \xbar_{\pi(m)}}\cdot\uu,
\]
implying
\eqref{eq:thm:seqsum-equiv-mod:1}, and completing the proof.

\pfpart{%
  (\ref{thm:seqsum-equiv-mod:d})
  $\Rightarrow$
  (\ref{thm:seqsum-equiv-mod:b}):
}
Let $Z$ be the set of all points $\zbar\in\extspace$ which satisfy
part~(\ref{thm:seqsum-equiv-mod:b}), and
let
\[
  V=\bigBraces{
        \xbar_{\pi(1)}\plusl \dotsb \plusl \xbar_{\pi(m)}
        :\:
        \pi\in\Pi_m
      }.
\]
We need to show that $\ohull V\subseteq Z$. The set $Z$ is convex by \Cref{prop:seqsum-multi}(\ref{prop:seqsum-multi:convex}), using the equivalence of parts (\ref{thm:seqsum-equiv-mod:a}) and (\ref{thm:seqsum-equiv-mod:b})
(and since all singletons are convex).
Therefore, it suffices to show that $V\subseteq Z$, which will imply
$\ohull V\subseteq Z$ by \Cref{thm:e:2}.

Let $\vbar\in V$, so $\vbar=\xbar_{\pi(1)}\plusl \dotsb \plusl \xbar_{\pi(m)}$ for some $\pi\in\Pi_m$. Let
\[
  \wbar = \trans{\PPsub{\pi(1)}}\xbar_{\pi(1)}
          \plusl\dotsb\plusl
          \trans{\PPsub{\pi(m)}}\xbar_{\pi(m)}.
\]
Note that, for $i,j\in\set{1,\dotsc,m}$, we have
\begin{equation}
\label{eq:PPiPPj}
  \PPi\trans{\PPsub{j}}=
  \begin{cases}
     \zero_{n\times n}
     &\text{if $i\ne j$,}
  \\
     \Idnn
     &\text{if $i=j$.}
  \end{cases}
\end{equation}
Thus, for each $i$,
\[
  \PPi\wbar
  =
  \PPi\trans{\PPsub{\pi(1)}}\xbar_{\pi(1)}
          \plusl\dotsb\plusl
  \PPi\trans{\PPsub{\pi(m)}}\xbar_{\pi(m)}
  =
  \xbar_i,
\]
where the first equality is by \Cref{pr:h:4}(\ref{pr:h:4c}), and the second
by \eqref{eq:PPiPPj}. Also,
\begin{align*}
  \A\wbar
  &=
  \A\trans{\PPsub{\pi(1)}}\xbar_{\pi(1)}
          \plusl\dotsb\plusl
  \A\trans{\PPsub{\pi(m)}}\xbar_{\pi(m)}
  \\
  &=
  \BiggParens{\sum_{i=1}^m \PPsub{i}\trans{\PPsub{\pi(1)}}}\xbar_{\pi(1)}
          \plusl\dotsb\plusl
  \BiggParens{\sum_{i=1}^m \PPsub{i}\trans{\PPsub{\pi(m)}}}\xbar_{\pi(m)}
  \\
  &=
  \xbar_{\pi(1)}\plusl \dotsb \plusl \xbar_{\pi(m)}
  =
  \vbar.
\end{align*}
The first equality is again by \Cref{pr:h:4}(\ref{pr:h:4c}),
the second by \eqref{eq:add-mat-defn},
and the third by \eqref{eq:PPiPPj}.
Thus,
$\vbar=\A\wbar$ and $\PPi\wbar=\xbar_i$ for all $i$, so $\vbar\in Z$.%
\indexg{sequential sum!characterizations of|)}%
\qedhere
\end{proof-parts}
\end{proof}

\indexg{sequential sum!closed graph property|(}%
As we show next,
\Cref{thm:seqsum-equiv-mod} implies that sequential addition exhibits a certain continuity property called the \emph{closed graph}%
\indexg{closed graph property}
property, which states that the set of tuples $\mtuple{\xbar,\ybar,\zbar}$ such that $\zbar\in\xbar\seqsum\ybar$ is closed in $\eRn\times\eRn\times\eRn$.

\begin{theorem}   \label{thm:seqsum-cont-prop}
  Let $\seq{\xbar_t}$, $\seq{\ybar_t}$, $\seq{\zbar_t}$ be
  sequences in $\extspace$ converging, respectively, to some points
  $\xbar,\ybar,\zbar\in\extspace$.
  Suppose
  $\zbar_t\in\xbar_t\seqsum\ybar_t$ for all $t$.
  Then $\zbar\in\xbar\seqsum\ybar$.
\end{theorem}

\begin{proof}
Let $\A\in\R^{n\times 2n}$ be the addition matrix, and
$\PPsub{1},\PPsub{2}\in\R^{n\times 2n}$ the component projection
matrices
(as in Eqs.~\ref{eq:comp-proj-mat-defn}
and~\ref{eq:add-mat-defn}).
Then, by \Cref{thm:seqsum-equiv-mod}(\ref{thm:seqsum-equiv-mod:a},\ref{thm:seqsum-equiv-mod:b}),
for each $t$, there exists $\wbar_t\in\eRf{2n}$ such that $\PPsub{1}\wbar_t=\xbar_t$, $\PPsub{2}\wbar_t=\ybar_t$, and $\A\wbar_t=\zbar_t$.
By sequential compactness, the sequence $\seq{\wbar_t}$ must have a
convergent subsequence $\seq{\wbar_{s(t)}}$,
converging to some point $\wbar\in\extspace$,
for some indices $s(t)\in\nats$ with
$s(1)<s(2)<\cdots$.
By continuity of linear maps (\Cref{thm:linear:cont}\ref{thm:linear:cont:b}),
$\PPsub{1}\wbar=\lim\PPsub{1}\wbar_{s(t)}=\lim\xbar_{s(t)}=\xbar$,
and similarly $\PPsub{2}\wbar=\ybar$ and $\A\wbar=\zbar$.
Therefore, $\zbar\in\xbar\seqsum\ybar$ by \Cref{thm:seqsum-equiv-mod}(\ref{thm:seqsum-equiv-mod:b},\ref{thm:seqsum-equiv-mod:a}).%
\indexg{sequential sum!closed graph property|)}%
\end{proof}

\indexg{sequential sum!characterizations of|(}%
Here are several additional consequences of
\Cref{thm:seqsum-equiv-mod}:

\begin{corollary}     \label{cor:seqsum-conseqs}
  Let $\xbar,\ybar,\zbar\in\extspace$, and let
  $X,Y\subseteq\extspace$.
  \begin{letter-compact}
  \item   \label{cor:seqsum-conseqs:a}
    $\xbar\seqsum \ybar = \lb{\xbar\plusl\ybar}{\,\ybar\plusl\xbar}$.
  \item   \label{cor:seqsum-conseqs:b}
    $\zbar\in\xbar\seqsum\ybar$
    if and only if,
    for all $\uu\in\Rn$,
    if $\xbar\cdot\uu$ and $\ybar\cdot\uu$ are summable
    then
    $\zbar\cdot\uu = \xbar\cdot\uu + \ybar\cdot\uu$.
  \item   \label{cor:seqsum-conseqs:c}
    If $\xbar\plusl\ybar=\ybar\plusl\xbar$
    then
    $\xbar\seqsum\ybar = \{\xbar\plusl\ybar\} = \{\ybar\plusl\xbar\}$.
  \item   \label{cor:seqsum-conseqs:d}
    If $X$ and $Y$ are convex then
    \begin{equation}    \label{eq:cor:seqsum-conseqs:1}
        X\seqsum Y = \conv\bigParens{(X\plusl Y) \cup (Y\plusl X)}.
    \end{equation}
  \item   \label{cor:seqsum-conseqs:e}
    Suppose $\xbar'\plusl\ybar'=\ybar'\plusl\xbar'$
    for all $\xbar'\in X$ and all $\ybar'\in Y$
    (as will be the case, for instance, if either $X$ or
    $Y$ is included in $\Rn$).
    Then $X\seqsum Y = X\plusl Y = Y\plusl X$.
  \end{letter-compact}
\end{corollary}

\begin{proof}
~

\begin{proof-parts}
\pfpart{Part~(\ref{cor:seqsum-conseqs:a}):}
Immediate from
\Cref{thm:seqsum-equiv-mod}(\ref{thm:seqsum-equiv-mod:a},\ref{thm:seqsum-equiv-mod:d}),
applied to $\xbar$ and $\ybar$.

\pfpart{Part~(\ref{cor:seqsum-conseqs:b}):}
Immediate from
\Cref{thm:seqsum-equiv-mod}(\ref{thm:seqsum-equiv-mod:a},\ref{thm:seqsum-equiv-mod:c}),
applied to $\xbar$ and $\ybar$.

\pfpart{Part~(\ref{cor:seqsum-conseqs:c}):}
This follows from part~(\ref{cor:seqsum-conseqs:a})
and \Cref{pr:e1}(\ref{pr:e1:ohull}).

\pfpart{Part~(\ref{cor:seqsum-conseqs:d}):}
Let $C$ denote the set on the right-hand side of
\eqref{eq:cor:seqsum-conseqs:1}.
We have
\begin{equation}    \label{eq:cor:seqsum-conseqs:2}
   X \seqsum Y
   =
   \bigcup_{\xbar\in X, \ybar\in Y}  (\xbar\seqsum\ybar)
   =
   \bigcup_{\xbar\in X, \ybar\in Y}  \lb{\xbar\plusl\ybar}{\,\ybar\plusl\xbar},
\end{equation}
where the first equality is by
\Cref{prop:seqsum-multi}(\ref{prop:seqsum-multi:decomp})
and the second by part~(\ref{cor:seqsum-conseqs:a}).
If $\xbar\in X$ and $\ybar\in Y$ then
$\xbar\plusl\ybar\in X\plusl Y\subseteq C$,
and similarly,
$\ybar\plusl\xbar\in C$.
Since $C$ is convex, it follows that
$\lb{\xbar\plusl\ybar}{\,\ybar\plusl\xbar}\subseteq C$.
Thus,
$X \seqsum Y \subseteq C$
by \eqref{eq:cor:seqsum-conseqs:2}.

For the reverse inclusion,
since $\xbar\plusl\ybar\in\lb{\xbar\plusl\ybar}{\,\ybar\plusl\xbar}$,
the rightmost expression in
\eqref{eq:cor:seqsum-conseqs:2} includes
$X\plusl Y$.
Thus
$X\plusl Y\subseteq X\seqsum Y$.
Likewise,
$Y\plusl X\subseteq X\seqsum Y$.
Since $X\seqsum Y$ is convex by
\Cref{prop:seqsum-multi}(\ref{prop:seqsum-multi:convex}),
it follows that
$C\subseteq X\seqsum Y$
(by \Cref{pr:conhull-prop}\ref{pr:conhull-prop:aa}).

\pfpart{Part~(\ref{cor:seqsum-conseqs:e}):}
This follows directly from \Cref{prop:seqsum-multi}(\ref{prop:seqsum-multi:decomp})
and part~(\ref{cor:seqsum-conseqs:c}).%
\indexg{sequential sum!characterizations of|)}%
\qedhere
\end{proof-parts}
\end{proof}

As an example,
using the tools we have been developing,
we can now complete the calculation of
$\limray{\ee_1}\seqsum\limray{\ee_2}$ that was started in
\Cref{ex:oe1-seqsum-oe2}:

\begin{example}   \label{ex:oe1-seqsum-oe2:cont}
Continuing \Cref{ex:oe1-seqsum-oe2}, we claim
\begin{equation}   \label{eq:ex:oe1-seqsum-oe2:cont:1}
  \limray{\ee_1}\seqsum\limray{\ee_2}
   =
   \Braces{%
     \limray{\ee_1}\plusl\limray{\ee_2},\,
     \limray{\ee_2}\plusl\limray{\ee_1}
   }
   \;\cup\;
   \bigcup\bigBraces{
     \limray{\vv}\plusl\extspac{2}
     :\:
     \vv\in\Rstrictpos^2
   }.
\end{equation}
Let $S$ denote the set on the right-hand side of
\eqref{eq:ex:oe1-seqsum-oe2:cont:1}.
That $S\subseteq\limray{\ee_1}\seqsum\limray{\ee_2}$ was shown in
\Cref{ex:oe1-seqsum-oe2}.
For the reverse inclusion, suppose $\zbar\in\limray{\ee_1}\seqsum\limray{\ee_2}$.
Since $\limray{\ee_1}\cdot\ee_1=+\infty$ and
$\limray{\ee_2}\cdot\ee_1=0$,
which are summable, 
\Cref{cor:seqsum-conseqs}(\ref{cor:seqsum-conseqs:b}) implies
$\zbar\cdot\ee_1=+\infty$.
Similarly, $\zbar\cdot\ee_2=+\infty$.
In particular, this means $\zbar\not\in\R^2$.
Let $\vv$ be $\zbar$'s dominant direction.
Then $v_1\geq 0$; otherwise, if $v_1<0$, we would have
$\zbar\cdot\ee_1=-\infty$
(\Cref{thm:dom-dir}\ref{thm:dom-dir:a}\ref{thm:dom-dir:b}),
a contradiction.
Similarly, $v_2\geq 0$, so $\vv\in\Rpos^2$.
If $\vv\in\Rstrictpos^2$, then $\zbar$ is evidently in $S$;
otherwise, either $\vv=\ee_1$ or $\vv=\ee_2$.

If $\vv=\ee_1$ then, in canonical form,
$\zbar=\limray{\ee_1}\plusl\alpha\ee_2$ for some $\alpha\in\Rext$.
Further, $\alpha=\zbar\cdot\ee_2=+\infty$, so
$\zbar=\limray{\ee_1}\plusl\limray{\ee_2}$.
Similarly, if $\vv=\ee_2$, then
$\zbar=\limray{\ee_2}\plusl\limray{\ee_1}$.
Thus, in all cases, $\zbar\in S$, proving
\eqref{eq:ex:oe1-seqsum-oe2:cont:1}.
\end{example}

We finish this section with several further properties of sequential addition.
\indexg{sequential sum!convex hull and|(}%
\indexg{convex hull, astral!sequential sum and|(}%
The first relates sequential sums to convex hulls:

\begin{theorem}  \label{thm:seqsum-in-union}
  Let $X_1,\dotsc,X_m\subseteq\extspace$,
  and let $\lambda_1,\dotsc,\lambda_m\in\Rpos$,
  where $m\geq 1$.
  Then
  \begin{equation}   \label{eq:thm:seqsum-in-union:2}
     (\lambda_1 X_1) \seqsum \dotsb \seqsum (\lambda_m X_m)
     \subseteq
     \BiggParens{\sum_{i=1}^m \lambda_i}
             \conv\BiggParens{\bigcup_{i=1}^m X_i}.
  \end{equation}
\end{theorem}

\begin{proof}
If any of the sets $X_i$ is empty, then the left-hand side of
\eqref{eq:thm:seqsum-in-union:2} is also empty and the claim
holds.
We therefore assume henceforth that $X_i\neq\emptyset$
for $i=1,\dotsc,m$.

Let $s=\sum_{i=1}^m \lambda_i$.
If $s=0$, so that $\lambda_i=0$ for all $i$, then
both sides of
\eqref{eq:thm:seqsum-in-union:2} are equal to $\{\zero\}$,
again implying the claim.
We therefore assume henceforth that $s>0$.
Also, by possibly rearranging the indices, we assume
without loss of generality that $\lambda_i>0$ for $i=1,\dotsc,\ell$ and $\lambda_i=0$ for $i=\ell+1,\dotsc,m$, for some $\ell\in\set{1,\dotsc,m}$.

Let $\zbar\in (\lambda_1 X_1) \seqsum \dotsb \seqsum (\lambda_m X_m)
     =(\lambda_1 X_1) \seqsum \dotsb \seqsum (\lambda_\ell X_\ell)
     =\Seqsum(\lambda_1 X_1,\dotsc,\lambda_\ell X_\ell)$.
Then,
for $i=1,\dotsc,\ell$,
there exist sequences $\seq{\xx'_{it}}$ in $\Rn$ such that $\xx'_{it}\to\lambda_i\xbar_i$
for some $\xbar_i\in X_i$, and $\sum_{i=1}^\ell \xx'_{it}\to\zbar$. Let
$\xx_{it}=\xx'_{it}/\lambda_i$ and $\lambar_i=\lambda_i/s$. Then
$\xx_{it}\to\xbar_i$ for $i=1,\dotsc,\ell$, $\sum_{i=1}^\ell\lambar_i=1$, and $\sum_{i=1}^\ell\lambar_i\xx_{it}\to\zbar/s$.
Therefore,
\[
  {\zbar}/{s}
  \in\ohull\set{\xbar_1,\dotsc,\xbar_\ell}
  \subseteq\conv\BiggParens{\bigcup_{i=1}^m X_i},
\]
where the inclusions are, respectively,
by \Cref{thm:e:7}
and \Cref{thm:convhull-of-simpices}.
Multiplying
by $s$ completes the proof.%
\indexg{sequential sum!convex hull and|)}%
\indexg{convex hull, astral!sequential sum and|)}%
\end{proof}

The next proposition relates sequential addition to an analogous
generalization of subtraction:

\begin{proposition}  \label{pr:swap-seq-sum}
  Let $\xbar,\ybar,\zbar\in\extspace$.
  Then
  $\zbar\in\xbar\seqsum\ybar$
  if and only if
  $\xbar\in\zbar\seqsum(-\ybar)$.
\end{proposition}

\begin{proof}
Suppose $\zbar\in\xbar\seqsum\ybar$.
Then there exist sequences $\seq{\xx_t}$, $\seq{\yy_t}$, and
$\seq{\zz_t}$ in $\Rn$ such that
$\xx_t\rightarrow\xbar$,
$\yy_t\rightarrow\ybar$,
$\zz_t\rightarrow\zbar$,
and $\zz_t=\xx_t+\yy_t$ for all $t$.
Since $-\yy_t\rightarrow-\ybar$ and
$\zz_t-\yy_t=\xx_t\rightarrow\xbar$,
this also shows that
$\xbar\in\zbar\seqsum(-\ybar)$.

The converse can then be derived by applying the above
with $\xbar$ and $\zbar$ swapped, and $\ybar$ and $-\ybar$ swapped.
\end{proof}

\indexg{sequential sum!linear map@under linear map|(}%
\indexg{linear maps, astral!sequential sum under|(}%
Finally, we show that sequential addition commutes with linear maps:

\begin{theorem}  \label{thm:distrib-seqsum}
  Let $\A\in\Rmn$, and let $A:\eRn\to\eRm$ be the associated astral
  linear map (so that $A(\xbar)=\A\xbar$ for $\xbar\in\extspace$).
  Let $X,Y\subseteq\extspace$.
  Then
  \[
     A (X\seqsum Y) = A(X) \seqsum A(Y).
  \]
\end{theorem}

\begin{proof}
For any $\xbar,\ybar\in\extspace$, we have
\begin{align}
\notag
  A (\xbar\seqsum\ybar)
  &=
  A\bigParens{\lb{\xbar\plusl\ybar}{\,\ybar\plusl\xbar}}
\\
\notag
  &=
  \seg\bigParens{A(\xbar\plusl\ybar),\,A(\ybar\plusl\xbar)}
\\
  &=
  \seg\bigParens{A(\xbar)\plusl A(\ybar),\, A(\ybar)\plusl A(\xbar)}
  =
  A(\xbar)\seqsum A(\ybar).
  \label{eq:thm:distrib-seqsum:1}
\end{align}
The first and fourth equalities are by
\Cref{cor:seqsum-conseqs}(\ref{cor:seqsum-conseqs:a}).
The second and third are by
\Cref{thm:e:9}
and
\Cref{pr:h:4}(\ref{pr:h:4c}),
respectively.
Thus,
\begin{align*}
  A (X \seqsum Y)
  &=
  A\BiggParens{\,\bigcup_{\xbar\in X, \ybar\in Y} (\xbar\seqsum\ybar)}
  \\
  &=
  \bigcup_{\xbar\in X, \ybar\in Y}  A(\xbar\seqsum\ybar)
  \\
  &=
  \bigcup_{\xbar\in X, \ybar\in Y} \bigBracks{A(\xbar) \seqsum A(\ybar)}
  \\
  &=
  A(X) \seqsum A(Y).
\end{align*}
The first and fourth equalities are by
\Cref{pr:seqsum-props}(\ref{pr:seqsum-props:c}).
The third is by
\eqref{eq:thm:distrib-seqsum:1}.%
\indexg{sequential sum|)}%
\indexg{sequential sum!linear map@under linear map|)}%
\indexg{linear maps, astral!sequential sum under|)}%
\end{proof}

\section{Sequential multiples}

\indexg{sequential multiples|(}%
We next develop a sequential
analogue of
multiplication of a vector by a scalar.
When combined with sequential sum, this will allow us to define
sequential generalizations of
convex and conic combinations, which will be studied in the following section.

For a point $\xbar\in\extspace$ and scalar $\alpha\in\Rextpos$,
we have already seen that the scalar product $\alpha\xbar$ provides
a natural generalization of the standard scalar-vector product.
Here,
we consider a different generalization based on
sequences, along the lines of how sequential sum was defined in
\Cref{dfn:seq-sum}:

\begin{definition}
\label{def:seq:multiple}
\indexg{sequential multiples!defined|(}%
Let $\xbar\in\extspace$,
and let
$\alpha\in\Rextpos$.
Then the \emph{sequential $\alpha$-multiple of $\xbar$} (or
simply the \emph{$\alpha$-multiple of $\xbar$}),
denoted $\mul{\alpha}{\xbar}$,%
\indexm{alpha xbar700}{$\protect\mul{\alpha}{\xbar}$}{sequential multiple}
is the set of all points $\zbar\in\extspace$
for which there exist a sequence
$\seq{\lambda_t}$ in $\Rpos$ 
and a span-bound sequence $\seq{\xx_t}$ in $\Rn$
such that
$\lambda_t\rightarrow\alpha$,
$\xx_t\rightarrow\xbar$,
and
$\lambda_t\xx_t\rightarrow\zbar$.%
\indexg{sequential multiples!defined|)}%
\end{definition}

As with sequential sum, the result of the sequential multiple
operation is always a set, even if it happens to be a singleton.

When $\alpha\in\Rpos$, the definition is unchanged if
the requirement that $\seq{\xx_t}$ be span-bound is omitted.
This is because, by \Cref{thm:spam-limit-seqs-exist}, any sequence $\seq{\xx_t}$ can
be replaced by a strongly equivalent span-bound sequence
$\seq{\xx'_t}$, always satisfying
$\lim\xx'_t=\lim\xx_t$ and
$\lim\lambda_t\xx'_t=\lim\lambda_t\xx_t$
(by \Cref{prop:strong:eq:propties}\ref{prop:strong:eq:propties:c}),
provided that the $\lambda_t\negKern$'s are bounded, as is the case
when $\lim\lambda_t\in\Rpos$.

However, when $\alpha=\oms$, the span-boundness requirement is an essential part of the definition. 
For example, in $\eR$, we have $\mul{\oms}{0}=\set{0}$ since the only
span-bound sequence converging to $0$ is the constant $0$ sequence.
On the other hand, in the absence of the span-boundness requirement,
we could pick $x_t=c/t$ and $\lambda_t=t$, for any $c\in\R$, so that
$x_t\rightarrow 0$,
$\lambda_t\rightarrow\oms$,
and
$\lambda_t x_t\rightarrow c$.
Thus, for such an alternative version of the definition, the
$\oms$-multiple of $0$ would include all of $\R$,
rather than only $\{0\}$.
We choose to require span-boundness in the definition of
$\mul{\oms}\xbar$ because doing so will yield the appropriate
generalization of a conic combination, 
compatible with astral cones, rays, and conic hulls,
as will be developed in \Cref{sec:cones}.

The next
\namecref{thm:mul-char}
characterizes $\alpha$-multiples of any point
$\xbar\in\extspace$, for all values of $\alpha\in\Rextpos$.
When $\alpha>0$, this set consists only of the single point $\alpha\xbar$, as
might be expected.
But when $\alpha=0$, the set is more subtle, and turns out to be
the segment joining $\zero$ and $\xbar$'s iconic part.
Note in particular that $\mul{1}{\xbar}=\{\xbar\}$.

\begin{theorem}   \label{thm:mul-char}
  Let $\xbar\in\extspace$, and let
  $\alpha\in\Rextpos$.
  Then:
  \begin{letter-compact}
  \item   \label{thm:mul-char:a}
    $\almul{\xbar} = \braces{\alpha\xbar}$
    if $\alpha>0$.
  \item   \label{thm:mul-char:b}
    $\zmul{\xbar} = \lb{\zero}{\ebar}$
    where $\ebar\in\corezn$ is $\xbar$'s iconic part
    (so that $\xbar\in\ebar\plusl\Rn$).
  \end{letter-compact}
  Thus, in all cases,
  $\alpha\xbar\in\almul{\xbar}$,
  $\almul{\xbar}$ is a convex set,
  and $\almul{\xx}=\set{\alpha\xx}$ for $\xx\in\Rn$.
\end{theorem}

\begin{proof}
  ~

\begin{proof-parts}
\pfpart{Part~(\ref{thm:mul-char:a}):}
Assume $\alpha>0$.
Suppose $\seq{\lambda_t}$ is any sequence in $\Rpos$ converging to~$\alpha$,
and that $\seq{\xx_t}$ is any span-bound sequence in $\Rn$ converging
to $\xbar$.
Then $\lambda_t\xx_t\rightarrow\alpha\xbar$.
This follows from
\Cref{pr:scalar-prod-props}(\ref{pr:scalar-prod-props:e})
in the case that $\alpha\in\Rstrictpos$,
and from \Cref{thm:seq-to-limray:0} in the case that
$\alpha=\oms$.
Thus, $\alpha\xbar\in\almul{\xbar}$.

Moreover, if $\zbar\in\almul{\xbar}$, then there exists a sequence
$\seq{\lambda_t}$ in $\Rpos$ and a span-bound sequence $\seq{\xx_t}$
in $\Rn$ with $\lambda_t\rightarrow\alpha$, $\xx_t\rightarrow\xbar$,
and $\lambda_t\xx_t\rightarrow\zbar$.
By the foregoing argument, these conditions further imply that
$\lambda_t\xx_t\rightarrow\alpha\xbar$.
Thus, $\zbar=\alpha\xbar$, proving the claim.

\pfpart{Part~(\ref{thm:mul-char:b}):}
Let $\VV\omm\plusl\qq$ be $\xbar$'s canonical representation
where $\VV=[\vv_1,\dotsc,\vv_k]$ and $\qq,\vv_1,\dotsc,\vv_k\in\Rn$.
This implies $\ebar=\VV\omm$ (by
\Cref{thm:icon-fin-decomp}).

We first prove $\zmul{\xbar}\subseteq\lb{\zero}{\ebar}$.
Let $\zbar\in\zmul{\xbar}$.
Then there exists a sequence $\seq{\lambda_t}$ in $\Rpos$ and
a span-bound sequence $\seq{\xx_t}$ in $\Rn$ such that $\lambda_t\rightarrow 0$,
$\xx_t\rightarrow\xbar$, and $\lambda_t\xx_t\rightarrow\zbar$.
These imply that $\lambda_t<1$ for all but finitely many values of
$t$; by discarding all other elements, we thus can assume
$\lambda_t\in[0,1]$ for all $t$.
Since
$(1-\lambda_t)\zero+\lambda_t\xx_t\rightarrow \zbar$,
it then follows that $\zbar\in\lb{\zero}{\xbar}$
by \Cref{cor:e:1}.

Suppose by way of contradiction that
$\zbar\not\in\lb{\zero}{\ebar}$.
Since $\zbar\in\lb{\zero}{\xbar}$, by \Cref{thm:lb-with-zero}
(applied to $\xbar=\VV\omm\plusl\qq$,
and again to $\ebar=\VV\omm$),
the only possibility is that
$\zbar=\VV\omm\plusl\gamma\qq$ for some $\gamma\in(0,1]$,
and that $\qq\neq\zero$.
Thus, $\lambda_t\xx_t\cdot\qq\rightarrow\zbar\cdot\qq=\gamma\norm{\qq}^2>0$
(with convergence by
\Cref{thm:i:1}\ref{thm:i:1c},
and the equality by
\Cref{pr:vtransu-zero}
since $\qq\perp\VV$).
On the other hand, by similar reasoning,
$\xx_t\cdot\qq\rightarrow\xbar\cdot\qq=\norm{\qq}^2$,
so
$\lambda_t\xx_t\cdot\qq\rightarrow 0$
by continuity of multiplication, and since $\lambda_t\rightarrow 0$.
Having reached a contradiction, we conclude that
$\zbar\in\lb{\zero}{\ebar}$.

Now suppose $\zbar\in\lb{\zero}{\ebar}$, which we aim to show is in
$\zmul{\xbar}$.
For each $t$, let
\[
  \xx_t
  =
  \sum_{i=1}^k t^{k-i+1} \vv_i + \qq.
\]
Then $\xx_t\rightarrow\xbar$ (by \Cref{thm:i:seq-rep}),
and 
furthermore, $\seq{\xx_t}$ is span-bound.
By \Cref{thm:lb-with-zero}, applied to $\ebar$,
$\zbar$ must have the form
\begin{equation}   \label{eq:thm:mul-char:1}
  \zbar
  =
  \limrays{\vv_1,\dotsc,\vv_{\ell}} \plusl \beta \vv_{\ell+1}
\end{equation}
for some $\ell\in\{0,\dotsc,k\}$ and some $\beta\in\Rpos$,
where it is understood that $\vv_{k+1}=\zero$ and that $\beta=0$ if
$\ell=k$.
For each $t$, let
\[
  \lambda_t
  =
  \begin{cases}
    t^{\ell - k - 1/2}    & \text{if $\beta=0$,} \\
    \beta t^{\ell - k}   & \text{if $\beta>0$.} \\
  \end{cases}
\]
Then $\lambda_t\rightarrow 0$, in either case
(noting in particular that $\ell<k$ if $\beta>0$).
Further,
\[
  \lambda_t \xx_t
  =
  \begin{cases}
    \sum_{i=1}^{\ell} t^{\ell - i + 1/2} \vv_i
    +
    \BigBracks{
      \sum_{i=\ell+1}^{k} t^{\ell - i + 1/2} \vv_i
      +
      \lambda_t \qq
    }
    & \text{if $\beta=0$,} \\[\medskipamount]
    \sum_{i=1}^{\ell} \beta t^{\ell - i + 1} \vv_i
    +
    \beta \vv_{\ell+1}
    +
    \BigBracks{
      \sum_{i=\ell+2}^{k} \beta t^{\ell - i + 1} \vv_i
      +
      \lambda_t \qq
    }
    & \text{if $\beta>0$.} \\
  \end{cases}
\]
Noting that, in either case, the bracketed expressions converge to
$\zero$, it then follows that $\lambda_t\xx_t\rightarrow\zbar$
(again by \Cref{thm:i:seq-rep}).
Thus, $\zbar\in\zmul{\xbar}$.

\pfpart{Consequences:}
That $\almul{\xbar}$ is convex now follows from
\Cref{pr:e1}(\ref{pr:e1:ohull}).

If $\xx\in\Rn$ then $\xx$'s iconic part is $\ebar=\zero$ so that, in all
cases, $\almul{\xx}=\set{\alpha\xx}$.
\qedhere
\end{proof-parts}
\end{proof}

We next summarize some basic properties of sequential multiples.

\begin{proposition}
\label{pr:seq:mul}
  Let $\xbar\in\eRn$. Then:
  \begin{letter-compact}
  \item \label{i:seq:mul:A}
    $\A[\mul{\alpha}{\xbar}]=\mul{\alpha}{(\A\xbar)}$
    for any $\alpha\in\Rextpos$ and any matrix $\A\in\R^{m\times n}$.
  \item \label{i:seq:mul:sum:yy}
    $\mul{\lambda}(\xbar\plusl\yy) = \mul{\lambda}{\xbar}\plusl\lambda\yy$
    for any $\lambda\in\Rpos$ and any $\yy\in\Rn$.
  \end{letter-compact}
\end{proposition}
\begin{proof}
We can write $\xbar=\ebar\plusl\qq$ for some icon $\ebar\in\corezn$
and $\qq\in\Rn$.
\begin{proof-parts}
  \pfpart{Part (\ref{i:seq:mul:A}):}
  Let $\alpha\in\Rextpos$ and $\A\in\R^{m\times n}$.
  If $\alpha>0$ then
  \[
    \A[\mul{\alpha}{\xbar}]
    =\A\set{\alpha\xbar}=\set{\A(\alpha\xbar)}=\set{\alpha(\A\xbar)}
    =\mul{\alpha}{(\A\xbar)}
  \]
  with the first and last equality by
  \Cref{thm:mul-char}(\ref{thm:mul-char:a}),
  and the third by \Cref{pr:h:4}(\ref{pr:h:4g}).

  Likewise,
  if $\alpha=0$ then
  \[
    \A[\mul{0}{\xbar}]
    =\A[\lb{\zero}{\ebar}]=\lb{\zero}{\A\ebar}=\mul{0}{(\A\xbar)}.
  \]
  The first equality is by
  \Cref{thm:mul-char}(\ref{thm:mul-char:b})
  and the second by \Cref{thm:e:9}.
  The third equality is also by
  \Cref{thm:mul-char}(\ref{thm:mul-char:b}), noting that
  $\A\ebar$ is the iconic part of $\A\xbar$ since
  $\A\xbar=\A\ebar\plusl\A\qq$ and since
  $\A\ebar$ is an icon (by \Cref{pr:i:8}\ref{pr:i:8-matprod}).

  \pfpart{Part (\ref{i:seq:mul:sum:yy}):}
  Let $\lambda\in\Rpos$ and $\yy\in\Rn$.
  If $\lambda>0$ then by \Cref{thm:mul-char}(\ref{thm:mul-char:a}) and \Cref{pr:i:7}({\ref{pr:i:7b}}),
  \[
    \mul{\lambda}(\xbar\plusl\yy)
    =\set{\lambda(\xbar\plusl\yy)}
    =\set{\lambda\xbar\plusl\lambda\yy}
    =\mul{\lambda}{\xbar}\plusl\lambda\yy.
  \]
  For $\lambda=0$, we note that $\ebar$ is the iconic part of $\xbar\plusl\yy=\ebar\plusl(\qq+\yy)$
  and apply \Cref{thm:mul-char}(\ref{thm:mul-char:b}) to obtain
  \[
    \mul{0}(\xbar\plusl\yy)
    =\lb{\zero}{\ebar}
    =\mul{0}{\xbar}
    =\mul{0}{\xbar}\plusl 0\yy.
  \qedhere
  \]
\end{proof-parts}
\end{proof}

\indexg{sequential multiples!sequences converging to|(}%
\indexg{sequence convergence!sequential multiple@to sequential multiple|(}%
Suppose $\zbar\in\almul{\xbar}$ for some $\alpha\in\Rextpos$ and
$\xbar\in\extspace$.
Then by definition, $\zbar$ is the limit of some sequence of the form
$\seq{\lambda_t\xx_t}$ where the sequences $\seq{\lambda_t}$ and
$\seq{\xx_t}$ satisfy the conditions of
\Cref{def:seq:multiple}.
The next theorem shows that every sequence that converges to
$\zbar$ must be strongly equivalent to a sequence of this form.

\begin{theorem} \label{thm:seq-mul-char}
  Let $\xbar\in\extspace$ and
  $\alpha\in\Rextpos$.
  Let $\seq{\zz_t}$ be a sequence in $\Rn$ that converges to a point
  $\zbar\in\almul{\xbar}$.
  Then $\seqeq{\zz_t}{\lambda_t\xx_t}$ for some sequence $\seq{\lambda_t}$ in $\Rpos$
  and some span-bound sequence $\seq{\xx_t}$ in $\Rn$
  such that $\lambda_t\to\alpha$ and $\xx_t\to\xbar$.
\end{theorem}

\begin{proof}
  ~

\begin{proof-parts}
\pfpart{Case $\alpha=\oms$:}
In this case, we must have $\zbar=\limray{\xbar}$ by
\Cref{thm:mul-char}(\ref{thm:mul-char:a}).
The claim then follows directly from
\Cref{thm:seq-to-limray}.

\pfpart{Case $\alpha\in\Rstrictpos$:}
By
\Cref{thm:mul-char}(\ref{thm:mul-char:a}),
$\zbar=\alpha\xbar$ in this case.
By \Cref{thm:spam-limit-seqs-exist}, there exists a span-bound
sequence $\seq{\zz'_t}$ in $\Rn$ that is strongly equivalent
to $\seq{\zz_t}$.
For each $t$, let $\xx_t=\zz'_t/\alpha$, and let $\lambda_t=\alpha$.
Then $\xx_t\rightarrow\zbar/\alpha=\xbar$
(by \Cref{pr:scalar-prod-props}\ref{pr:scalar-prod-props:e}),
$\lambda_t\rightarrow\alpha$,
and
$\seqeq{\zz_t}{\zz_t'}=\lambda_t\xx_t$.

\pfpart{Case $\alpha=0$:}
Let $\VV\omm\plusl\qq$ be $\xbar$'s canonical representation
where $\VV=[\vv_1,\dotsc,\vv_k]$ and $\vv_1,\dotsc,\vv_k,\qq\in\Rn$.
By \Cref{thm:mul-char}(\ref{thm:mul-char:b}),
$\zbar\in\lb{\zero}{\VV\omm}$.
Therefore, by \Cref{thm:lb-with-zero}, $\zbar$ must have the
form given in \eqref{eq:thm:mul-char:1}, as in the proof of
\Cref{thm:mul-char}(\ref{thm:mul-char:b}), where, as before,
$\ell\in\{0,\dotsc,k\}$ and $\beta\in\Rpos$, and
where $\vv_{k+1}=\zero$ and $\beta=0$ if $\ell=k$.
Let $\VV'=[\vv_1,\dotsc,\vv_\ell]$ and let $\rr= \beta \vv_{\ell+1}$.
For each $t$, we can decompose $\zz_t$
as
$\zz_t = \VV' \bb_t + \rr_t$
for some (unique) $\bb_t\in\R^{\ell}$ and $\rr_t\in\Rn$
with $\rr_t\perp\VV'$ (see \Cref{pr:lin-decomp-rel-vecs}).
Further, these sequences must satisfy the conditions of
\Cref{thm:seq-rep}.
In particular, $\rr_t\rightarrow\rr$.
Letting $\zz'_t=\VV'\bb_t+\rr=\zz_t+(\rr-\rr_t)$,
so that $\zz'_t-\zz_t=\rr-\rr_t\rightarrow\zero$,
it follows that $\seqeq{\zz_t}{\zz'_t}$.

Let
\[
   d=
    \begin{cases}
            \ell+1 & \text{if $\beta=0$,}\\
            \ell+2 & \text{if $\beta>0$,}
    \end{cases}
\]
and let
\begin{align}
  \xx_t
  &=
  e^t \zz'_t
  +
  \sum_{i=d}^k  t^{k-i+1} \vv_i
  +
  \qq
  \label{eq:thm:seq-mul-char:1}
  \\
  &=
  \begin{cases}
    \sum_{i=1}^\ell e^t b_{t,i} \vv_i
    +
    \sum_{i=\ell+1}^k  t^{k-i+1} \vv_i
    +
    \qq
    & \text{if $\beta=0$,}
    \\[\smallskipamount]
    \sum_{i=1}^\ell e^t b_{t,i} \vv_i
    +
    \beta e^t \vv_{\ell+1}
    +
    \sum_{i=\ell+2}^k  t^{k-i+1} \vv_i
    +
    \qq
    & \text{if $\beta>0$.}
  \end{cases}
  \nonumber
\end{align}
Considering separately when $\beta$ is zero or positive,
it can be checked that the conditions of
\Cref{thm:seq-rep} are satisfied,
implying $\xx_t\rightarrow\xbar$.
Furthermore, $\seq{\xx_t}$ is evidently span-bound.

Let $\lambda_t=e^{-t}$.
Then $\lambda_t\rightarrow 0$, and
\[
  \lambda_t \xx_t - \zz'_t
  =
  e^{-t} \BiggParens{\sum_{i=d}^k  t^{k-i+1} \vv_i + \qq}
  \rightarrow
  \zero,
\]
where the equality follows
from \eqref{eq:thm:seq-mul-char:1}.
Thus, $\seqeq{\lambda_t\xx_t}{\zz'_t}$, so
$\seqeq{\lambda_t\xx_t}{\zz_t}$ by transitivity of strong equivalence,
completing the proof for this case.

In all of the cases above, we have shown that
$\seqeq{\zz_t}{\lambda_t\xx_t}$, implying
$\lambda_t\xx_t\rightarrow\zbar$
by
\Cref{pr:eq-in-lim-same-lim}.%
\indexg{sequential multiples!sequences converging to|)}%
\indexg{sequence convergence!sequential multiple@to sequential multiple|)}%
\indexg{sequential multiples|)}%
\qedhere
\end{proof-parts}
\end{proof}

\section{Sequential convex and conic combinations}
\label{sec:seq-cvx-conic-combs}

We next develop astral analogues of standard convex and conic combinations.
As reviewed in Sections~\ref{sec:prelim:convex-sets}
and~\ref{sec:prelim:cones},
a conic combination of points $\xx_1,\ldots,\xx_m\in\Rn$ is a point
of the form
\begin{equation}  \label{eq:conv-comb}
  \lambda_1\xx_1+\cdots+\lambda_m\xx_m
\end{equation}
for some $\lambda_1,\ldots,\lambda_m\in\Rpos$,
while
a convex combination has the same form with the additional constraint
that $\sum_{i=1}^m \lambda_i=1$.

\indexg{sequential convex combinations|(}%
\indexg{sequential conic combinations|(}%
Having studied sequential addition and sequential multiples,
we can now generalize these notions to astral space from
\eqref{eq:conv-comb}
simply by replacing standard sum with sequential sum,
scalar-vector products by sequential multiples,
and each finite scalar $\lambda_i\in\Rpos$ with a possibly infinite
scalar $\alpha_i\in\Rextpos$.
This leads to the following:

\begin{definition}
\label{def:seq:conic:comb}
\indexg{sequential convex combinations!defined|(}%
\indexg{sequential conic combinations!defined|(}%
  Let $\xbar_1,\dotsc,\xbar_m\in\extspace$. A \emph{sequential conic combination}
  of $\xbar_1,\dotsc,\xbar_m$ is any set of the form
  \begin{equation}   \label{eq:def:seq:conic:comb:1}
     \mul{\alpha_1}{\xbar_1}\seqsum \cdots \seqsum\mul{\alpha_m}{\xbar_m},
  \end{equation}
  where $\alpha_1,\dotsc,\alpha_m\in\Rextpos$.
  A \emph{sequential convex combination} is a set of the same form
  that also satisfies $\sum_{i=1}^m\alpha_i=1$.%
\indexg{sequential convex combinations!defined|)}%
\indexg{sequential conic combinations!defined|)}%
\end{definition}

Note that, as with sequential sums and sequential multiples,
sequential conic or convex combinations are always sets, not
individual points.
We also define specific terminology for
the sequential convex combination of two points:

\begin{definition}  \label{def:midpoint}
\indexg{midsegments and midpoints (sequential)|(}%
  For any $\xbar,\ybar\in\extspace$,
  and any $\lambda\in[0,1]$,
  the \emph{$\lambda$-midsegment} of $\xbar$ and $\ybar$
  is the set $\mul{1-\lambda}\xbar\seqsum\mul{\lambda}\ybar$.
  Elements of this set are called \emph{$\lambda$-midpoints} of $\xbar$ and $\ybar$.%
\indexg{midsegments and midpoints (sequential)|)}%
\end{definition}

The set appearing in \eqref{eq:def:seq:conic:comb:1}
is a composition of various sequential sum and sequential multiple
operations.
\indexg{sequential conic combinations!sequential characterization|(}%
\indexg{sequential convex combinations!sequential characterization|(}%
The next theorem gives a more holistic characterization of such an
expression, showing that it
can be viewed as a generalization of a standard conic
(or convex) combination obtained by replacing
coefficients and vectors with appropriate sequences and taking the limit
of the resulting expression. The theorem is phrased with generality to allow
mixing of a convex combination (with $k$ terms) with a conic combination (with $m-k$ terms).
Setting $k=0$ corresponds to a purely conic combination, while setting $k=m$ corresponds to
a purely convex combination.

In the latter case, note that the sequences specified in the theorem
are the same as those in \Cref{thm:e:7}
which characterized the outer hull of a finite set of points, with the
additional constraint that $\lambda_{it}\rightarrow\alpha_i$
for $i=1,\ldots,m$
(so $\alpha_i\negKern$'s role is similar to that of
$\lambar_i$ in that theorem).

\begin{theorem}   \label{thm:seqsum-midrays}
  Let $\xbar_1,\dotsc,\xbar_m\in\extspace$,
  let $\alpha_1,\dotsc,\alpha_m\in\Rextpos$,
  and either let $k=0$ or let $k\in\set{1,\dotsc,m}$ be such that $\sum_{i=1}^{\smash[t]{k}}\alpha_i=1$.
  Let $\zbar\in\extspace$.
  Then
  \begin{equation}  \label{eq:thm:seqsum-midrays:1}
     \zbar
     \in
     \mul{\alpha_1}{\xbar_1}
     \seqsum \cdots \seqsum
     \mul{\alpha_m}{\xbar_m}
  \end{equation}
  if and only if, for $i=1,\dotsc,m$, there exist sequences
  $\seq{\lambda_{it}}$ in $\Rpos$
  and span-bound sequences
  $\seq{\xx_{it}}$ in $\Rn$
  such that:
  \begin{letter-compact}
  \item
    $\xx_{it}\rightarrow\xbar_i$ for $i=1,\dotsc,m$.
  \item
    $\lambda_{it}\rightarrow\alpha_i$ for $i=1,\dotsc,m$.
  \item \label{i:seqsum-midrays-sum}
    If $k\ge 1$ then $\sum_{i=1}^{\smash[t]{k}} \lambda_{it}=1$ for all $t$.
  \item
    The sequence %
    $\zz_t=\sum_{i=1}^m \lambda_{it} \xx_{it}$ converges to
    $\zbar$.
  \end{letter-compact}
  The same equivalence holds if the sequences $\seq{\xx_{it}}$ are only required to be span-bound for the indices $i$ with $\alpha_i=\omega$.

\end{theorem}

\begin{proof}
  ~

\begin{proof-parts}
\pfpart{``If'' ($\Leftarrow$):}
Suppose sequences with the stated properties exist.
By sequential compactness, for each $i=1,\dotsc,m$ considered in turn,
the sequence $\seq{\lambda_{it}\xx_{it}}$ must have a convergent
subsequence; by discarding all other elements (and the corresponding
elements of all other sequences), we can assume that this entire
sequence converges to some point $\wbar_i\in\extspace$.
Further, this point must be in $\mul{\alpha_i}{\xbar_i}$ since all
the conditions of \Cref{def:seq:multiple} are satisfied.
Since $\lambda_{it}\xx_{it}\rightarrow\wbar_i$
and $\zz_t\rightarrow\zbar$, it then follows
by
\Crefequiv{thm:seqsum-equiv-mod}{thm:seqsum-equiv-mod:a}{thm:seqsum-equiv-mod:seq}
that
$\zbar\in\wbar_1\seqsum\cdots\seqsum\wbar_m$, proving
\eqref{eq:thm:seqsum-midrays:1}
(by
\Cref{prop:seqsum-multi}\ref{prop:seqsum-multi:decomp}).

\pfpart{``Only if'' ($\Rightarrow$):}
Suppose \eqref{eq:thm:seqsum-midrays:1} holds, implying
(by
\Cref{prop:seqsum-multi}\ref{prop:seqsum-multi:decomp})
that
$\zbar\in\wbar_1\seqsum\cdots\seqsum\wbar_m$ for some
$\wbar_i\in\mul{\alpha_i}{\xbar_i}$, for $i=1,\dotsc,m$,
and so, by
\Crefequiv{thm:seqsum-equiv-mod}{thm:seqsum-equiv-mod:a}{thm:seqsum-equiv-mod:seq},
that there exist sequences $\seq{\ww_{it}}$ in $\Rn$ such that
$\ww_{it}\rightarrow\wbar_i$ for $i=1,\dotsc,m$,
and $\sum_{i=1}^m \ww_{it} \rightarrow\zbar$.
For each $i$,
\Cref{thm:seq-mul-char},
applied
to the sequence $\seq{\ww_{it}}$ (which converges to
$\wbar_i\in\mul{\alpha_i}{\xbar_i}$)
yields that there exist
a sequence $\seq{\lambda_{it}}$ in $\Rpos$ and a
span-bound sequence $\seq{\xx_{it}}$ such that
$\lambda_{it}\rightarrow\alpha_i$,
$\xx_{it}\rightarrow\xbar_i$,
and $\seqeq{\ww_{it}}{\lambda_{it}\xx_{it}}$. Since
$\sum_{i=1}^m\ww_{it}\to\zbar$, we also have
$\sum_{i=1}^m \lambda_{it} \xx_{it}\to\zbar$
(by Propositions~\ref{pr:eq-in-lim-same-lim}
and~\ref{prop:strong:eq:propties}\ref{prop:strong:eq:propties:a}).
Thus, the sequences $\seq{\lambda_{it}}$ and $\seq{\xx_{it}}$ satisfy
all the conditions of the \namecref{thm:seqsum-midrays} except possibly (\ref{i:seqsum-midrays-sum}).
If $k=0$ then condition~(\ref{i:seqsum-midrays-sum}) is vacuously satisfied. We next consider the case $k\ge 1$
and show that the sequences $\seq{\lambda_{it}}$, for $i=1,\dotsc,m$, can be replaced by suitable
sequences $\seq{\lambda'_{it}}$ to satisfy all the conditions.

Let $s_t=\sum_{i=1}^{\smash[t]{k}} \lambda_{it}$.
Then $s_t\to\sum_{i=1}^{k} \alpha_i = 1$, so
$s_t>0$ for all but finitely many $t$; by discarding all other elements
(and the corresponding elements of all other sequences),
we can assume henceforth that $s_t>0$ for all $t$.
Let $\lambda'_{it}=\lambda_{it}/s_t$ for $i=1,\dotsc,m$ and all $t$.
Then $\sum_{i=1}^{\smash[t]{k}} \lambda'_{it}=1$ for all $t$,
and, for $i=1,\dotsc,m$,
$\lambda'_{it}=\lambda_{it}/s_t\to\alpha_i$ by
\Cref{prop:lim:eR}(\ref{i:lim:eR:genmul}) since
$\lambda_{it}\rightarrow\alpha_i$ and $1/s_t\to 1$.
Moreover,
\[
  \sum_{i=1}^m \lambda'_{it}\xx_{it}
  =
  \frac{1}{s_t} \sum_{i=1}^m \lambda_{it}\xx_{it}
  \to
  \zbar
\]
by
\Cref{pr:scalar-prod-props}(\ref{pr:scalar-prod-props:e})
since
$\sum_{i=1}^m \lambda_{it}\xx_{it}\to\zbar$
and $1/s_t\to 1$.
Thus, replacing sequences $\seq{\lambda_{it}}$ by $\seq{\lambda'_{it}}$,
for $i=1,\dotsc,m$,
all conditions of the \namecref{thm:seqsum-midrays} are satisfied.

\pfpart{Relaxed span-boundness requirement.}
Suppose sequences as stated in the theorem exist, but only the
sequences $\seq{\xx_{it}}$ with $\alpha_i=\omega$ are necessarily span-bound.
Let
\[
   I=\bigBraces{i\in\{1,\ldots,m\} :\: \alpha_i\in\Rpos},
   \;\text{~and~}\;
   I^c=\bigBraces{i\in\{1,\ldots,m\}:\:\alpha_i=\omega}.
\]
By \Cref{thm:spam-limit-seqs-exist},
for $i\in I$, there exist span-bound sequences
$\seq{\xx'_{it}}$ that are strongly equivalent to
$\seq{\xx_{it}}$
so that also $\xx'_{it}\rightarrow\xbar_i$.
Moreover, for $i\in I$, the sequences $\seq{\lambda_{it}}$ are bounded
since they have limits in $\Rpos$,
so by \Cref{prop:strong:eq:propties}(\ref{prop:strong:eq:propties:a},\ref{prop:strong:eq:propties:c}),
\[
  \sum_{i\in I}\lambda_{it}\xx'_{it}
  +
  \sum_{i\in I^c}\lambda_{it}\xx_{it}
  \seqeqsymbol
  \sum_{i\in I}\lambda_{it}\xx_{it}
  +
  \sum_{i\in I^c}\lambda_{it}\xx_{it}.
\]
Therefore, since the right-hand side converges to $\zbar$,
the left-hand side does as well (by \Cref{pr:eq-in-lim-same-lim}).
Thus, we can replace $\seq{\xx_{it}}$ by $\seq{\xx'_{it}}$ for all $i\in I$, and the resulting
sequences will still satisfy the conditions of the theorem while also being span-bound.
\qedhere
\end{proof-parts}
\end{proof}

\indexg{midsegments and midpoints (sequential)!sequential characterization|(}%
As an immediate consequence of \Cref{thm:seqsum-midrays}, $\lambda$-midsegments can be characterized
in terms of sequences as in the next corollary.
These sequences are the same as those in \Cref{cor:e:1} with the
additional constraint that $\lambda_t\rightarrow\lambda$.

\begin{corollary}   \label{cor:midpoint-char}
  Let $\xbar,\ybar\in\eRn$ and let $\lambda\in[0,1]$.
  Let $\zbar\in\extspace$.
  Then $\zbar\in\flammid{\xbar}{\ybar}$ if and only if there exist sequences
  $\seq{\xx_t}$, $\seq{\yy_t}$ in $\Rn$, and
  $\seq{\lambda_t}$ in $[0,1]$
  such that $\xx_t\to\xbar$, $\yy_t\to\ybar$, $\lambda_t\to\lambda$, and
  $(1-\lambda_t)\xx_t+\lambda_t\yy_t\to\zbar$.%
\indexg{midsegments and midpoints (sequential)!sequential characterization|)}%
\indexg{sequential convex combinations!sequential characterization|)}%
\indexg{sequential conic combinations!sequential characterization|)}%
\end{corollary}

Here are some properties of sequential conic combinations (and hence also sequential convex combinations
and $\lambda$-midsegments):

\begin{proposition}
\label{pr:seq-conic:prop}
Let $\xbar_1,\dotsc,\xbar_m\in\extspace$, and 
let $\alpha_1,\dotsc,\alpha_m\in\Rextpos$.
Also, let $\xx_1,\ldots,\xx_m\in\Rn$,
and $\lambda_1,\dotsc,\lambda_m\in\Rpos$.
Then:
  \begin{letter-compact}
  \item   \label{i:seq-conic:subset}
    $\alpha_1\xbar_1\seqsum\dotsb\seqsum\alpha_m\xbar_m
       \subseteq
       \mul{\alpha_1}{\xbar_1}\seqsum \dotsb \seqsum\mul{\alpha_m}{\xbar_m}$.
  \item   \label{i:seq-conic:conv}
\indexg{sequential conic combinations!convexity of|(}%
\indexg{sequential convex combinations!convexity of|(}%
    $\mul{\alpha_1}{\xbar_1}\seqsum \dotsb \seqsum\mul{\alpha_m}{\xbar_m}$
    is a convex set.
  \item   \label{i:seq-conic:finite}
    $\mul{\lambda_1}{\xx_1}\seqsum \dotsb \seqsum\mul{\lambda_m}{\xx_m}
    = \set{\lambda_1\xx_1+\dotsb+\lambda_m\xx_m}$.
  \end{letter-compact}
\end{proposition}

\begin{proof}
By \Cref{thm:mul-char}, for $i=1,\ldots,m$,
$\mul{\alpha_i}{\xbar_i}$ is convex and includes $\alpha_i\xbar_i$,
and $\mul{\lambda_i}{\xx_i}=\set{\lambda_i\xx_i}$.
Parts~(\ref{i:seq-conic:subset}),~(\ref{i:seq-conic:conv}),
and~(\ref{i:seq-conic:finite})
then follow respectively from
\Cref{pr:seqsum-props}(\ref{pr:seqsum-props:b}),
\Cref{prop:seqsum-multi}(\ref{prop:seqsum-multi:convex}),
and
\Cref{pr:seqsum-props}(\ref{pr:seqsum-props:finite}).%
\indexg{sequential conic combinations!convexity of|)}%
\indexg{sequential convex combinations!convexity of|)}%
\end{proof}

\Cref{pr:seq-conic:prop}(\ref{i:seq-conic:finite}) shows that the only element of
a sequential conic combination of finite points
with finite nonnegative coefficients
is the corresponding standard conic combination. In this sense,
sequential conic and convex combinations generalize their standard versions.

In particular,
\Cref{pr:seq-conic:prop}(\ref{i:seq-conic:finite}) implies
that every pair of finite points has exactly one $\lambda$-midpoint
for every $\lambda\in[0,1]$. However,
in general, a pair of astral points might have
many $\lambda$-midpoints.
As an extreme example, in $\Rext$, $\flammid{(-\infty)}{(+\infty)} = \Rext$ for all $\lambda\in [0,1]$.
This is because
this set
is convex (by \Cref{pr:seq-conic:prop}\ref{i:seq-conic:conv}) and contains
both $-\infty$ and $+\infty$ (as can be shown by exhibiting suitable sequences in \Cref{cor:midpoint-char}), and hence
must contain $\lb{-\infty}{+\infty}=\Rext$.%
\indexg{sequential conic combinations|)}

\indexg{sequential convex combinations!convex hull as union of|(}%
\indexg{convex hull, astral!union of sequential convex combinations@as union of sequential convex combinations|(}%
\indexg{outer convex hull!union of sequential convex combinations@as union of sequential convex combinations|(}%
\indexg{polytopes, astral!union of sequential convex combinations@as union of sequential convex combinations|(}%
As we saw in \Cref{roc:thm2.3}, the (standard) convex hull of any set
$S\subseteq\Rn$ is equal to the union of all convex combinations of
finite sets of points in $S$.
The next
\namecref{cor:conv-as-seqsum-midrays}
uses the characterization of sequential convex
combinations in \Cref{thm:seqsum-midrays} to show that
the (outer) convex hull of a finite set of astral points is equal to
the union of all sequential convex combinations of the points.
Combined with \Cref{thm:convhull-of-simpices},
this further implies that the convex hull of any
set $S\subseteq\extspace$ is equal to the union of all
sequential convex combinations of all finite subsets of $S$,
providing an astral analogue of \Cref{roc:thm2.3}.
Also, as a special case, the
\namecref{cor:conv-as-seqsum-midrays}
shows
that the segment joining $\xbar,\ybar\in\eRn$ is the union of all
$\lambda$-midpoints of $\xbar$ and $\ybar$ across all $\lambda\in[0,1]$.

\begin{corollary}   \label{cor:conv-as-seqsum-midrays}
  Let $\xbar_1,\dotsc,\xbar_m\in\extspace$.
  Then
  \begin{align}
  \notag
     \ohull\{\xbar_1,\dotsc,\xbar_m\}
     =
     \bigcup
     \;
     \Bigl\{
     &
       \mul{\lambda_1}{\xbar_1}
       \seqsum \cdots \seqsum
       \mul{\lambda_m}{\xbar_m}:
\\
   \label{eq:cor:conv-as-seqsum-midrays:1}
     &\qquad
     \smash[t]{
       \lambda_1,\dotsc,\lambda_m\in\Rpos
       \textup{ and }
       {\textstyle \sum_{i=1}^m \lambda_i = 1}
     \Bigr\}.
     }
  \end{align}
  Consequently, for any $\xbar,\ybar\in\eRn$,
  \[
     \lb{\xbar}{\ybar}=\bigcup_{\lambda\in[0,1]}
        \bigBracks{\flammid{\xbar}{\ybar}}.
  \]
\end{corollary}

\begin{proof}
Let $\zbar\in\extspace$.
By \Cref{thm:e:7}, $\zbar$ is in
$\ohull\{\xbar_1,\dotsc,\xbar_m\}$
if and only if there exist sequences as in
\Cref{thm:seqsum-midrays} (with $k=m$)
for some
$\lambda_1,\dotsc,\lambda_m\in\Rpos$
with $\sum_{i=1}^m \lambda_i = 1$, and so,
by that same \namecref{thm:seqsum-midrays},
if and only if $\zbar$ is included in the union on the right-hand side
of \eqref{eq:cor:conv-as-seqsum-midrays:1}.%
\indexg{sequential convex combinations!convex hull as union of|)}%
\indexg{convex hull, astral!union of sequential convex combinations@as union of sequential convex combinations|)}%
\indexg{outer convex hull!union of sequential convex combinations@as union of sequential convex combinations|)}%
\indexg{polytopes, astral!union of sequential convex combinations@as union of sequential convex combinations|)}%
\end{proof}

\indexg{sequential convex combinations!affine maps@under affine maps|(}%
\indexg{affine maps, astral!sequential convex combinations under|(}%
We show next that the operation of taking a sequential convex combination commutes with the application of an affine map,
yielding
(in light of \Cref{cor:conv-as-seqsum-midrays})
a strengthening of \Cref{thm:e:9}:

\begin{theorem}
\label{thm:F:seqconv}
  Let $\A\in\R^{m\times n}$, $\bbar\in\extspac{m}$, and
  let $F:\extspace\rightarrow\extspac{m}$ be the
  affine map
  $F(\zbar)=\bbar\plusl \A\zbar$
  for $\zbar\in\extspace$.
  Let $\xbar_1,\dotsc,\xbar_\ell\in\extspace$, and let $\lambda_1,\dotsc,\lambda_\ell\in\Rpos$ with $\sum_{i=1}^\ell\lambda_i=1$.
  Then
  \begin{equation*}  %
     F\BigParens{
       \mul{\lambda_1}{\xbar_1}
       \seqsum \cdots \seqsum
       \mul{\lambda_\ell}{\xbar_\ell}
     }
     =
     \mul{\lambda_1}{F(\xbar_1)}
       \seqsum \cdots \seqsum
     \mul{\lambda_\ell}{F(\xbar_\ell)}.
  \end{equation*}
\end{theorem}

\begin{proof}
Let $\xbar'_i=F(\xbar_i)$ for $i=1,\ldots,\ell$, and
write $\bbar=\ebar\plusl\qq$ for some icon $\ebar\in\corezn$ and $\qq\in\Rn$. Since $\ebar$ is an icon, we have $\limray{\ebar}=\ebar$ (\Cref{pr:i:8}\ref{pr:i:8d}), so $\xbar'_i=\limray{\ebar}\plusl\qq\plusl\A\xbar_i$ for all $i$.

\begin{proof-parts}
\pfpart{Part ``$\subseteq$'':}
Let $\zbar\in\mul{\lambda_1}{\xbar_1}
       \seqsum \cdots \seqsum
       \mul{\lambda_\ell}{\xbar_\ell}$.
We aim to show that
\begin{equation}  \label{eq:thm:F:seqconv:3}
  F(\zbar)
  \in
  \mul{\lambda_1}\xbar'_1
  \seqsum\dotsb\seqsum
  \mul{\lambda_\ell}\xbar'_\ell.
\end{equation}
By \Cref{thm:seqsum-midrays}, for $i=1,\dotsc,\ell$, there exist
sequences $\seq{\xx_{it}}$ in $\Rn$ and $\seq{\glambdait}$ in
$\Rpos$ such that $\xx_{it}\to\xbar_i$, %
$\glambdait\to\lambda_i$, %
 $\sum_{i=1}^\ell\glambdait=1$ for all $t$, 
and $\zz_t\to\zbar$
where $\zz_t=\sum_{i=1}^\ell\glambdait\xx_{it}$ for all $t$.

We construct, for each $i$, a sequence $\seq{\xx'_{it}}$ in $\Rn$ such that $\xx'_{it}\to\xbar'_i$. 
Let $\seq{\dd_t}$ be any span-bound sequence that converges to $\ebar$
(which exists by 
Theorems~\ref{thm:i:1}\ref{thm:i:1d} 
and~\ref{thm:spam-limit-seqs-exist}),
and, for each $t$, let
\[
  M_t=t\BigParens{1+\max_{i\in\set{1,\dotsc,\ell}}\norm{\qq+\A\xx_{it}}}.
\]
This definition implies $M_t\to+\infty$, and also, for
$i=1,\ldots,\ell$, that
\begin{equation}  \label{eq:thm:F:seqconv:2}
  0
  \leq
  \frac{\norm{\qq+\A\xx_{it}}}{M_t}
  \leq
  \frac{1}{t},
\end{equation}
for all $t$, so
$\norm{\qq+\A\xx_{it}}/M_t\to 0$.
For each $t$,
let $\xx'_{it}=M_t\dd_t+\qq+\A\xx_{it}$.
Then by continuity of affine maps
(\Cref{cor:aff-cont}),
$\qq+\A\xx_{it}\to\qq\plusl\A\xbar_i$, so
\begin{equation*}   %
  \xx'_{it}=M_t\dd_t+\qq+\A\xx_{it}\to\limray{\ebar}\plusl\qq\plusl\A\xbar_i=\xbar'_i,
\end{equation*}
with convergence from \Cref{thm:lim-plusl}
(applied with $\xx_t$, $\xbar$, $\yy_t$, $\ybar$, $\lambda_t$,
as they appear in that \namecref{thm:lim-plusl},
set to
$\dd_t$, $\ebar$,
$\qq+\A\xx_{it}$, $\qq\plusl\A\xbar_i$,
$M_t$, respectively).

Similarly,
$\qq+\A\zz_t\to\qq\plusl\A\zbar$ by
continuity of affine maps,
and for all $t$,
\[
 0
 \le
 \frac{\norm{\qq+\A\zz_t}}{M_t}
 =
 \frac{\norm{\sum_{i=1}^\ell\glambdait(\qq+\A\xx_{it})}}{M_t}
 \le
 \frac{\sum_{i=1}^\ell\glambdait\norm{\qq+\A\xx_{it}}}{M_t}
 \le
 \frac{1}{t},
\]
where the equality follows from the definition of $\zz_t$,
the second inequality is by the triangle inequality,
and the third inequality is by \eqref{eq:thm:F:seqconv:2}.
Hence,
$\norm{\qq+\A\zz_t}/M_t\rightarrow 0$,
and therefore,
\begin{align*}
  \smash[b]{\sum_{i=1}^\ell\glambdait\xx'_{it}}
  =
  \smash[b]{\sum_{i=1}^\ell\glambdait(M_t\dd_t+\qq+\A\xx_{it})}
  &=
  M_t\dd_t+\qq+\A\zz_t
\\
  &
  \to
  \limray{\ebar}\plusl\qq\plusl\A\zbar = F(\zbar).
\end{align*}
The first two equalities are by definitions.
The convergence is by \Cref{thm:lim-plusl}
(applied with $\xx_t$, $\xbar$, $\yy_t$, $\ybar$, $\lambda_t$,
as in that \namecref{thm:lim-plusl},
set to
$\dd_t$, $\ebar$,
$\qq+\A\zz_t$, $\qq\plusl\A\zbar$,
$M_t$, respectively).

Hence, $\xx'_{it}\to\xbar'_i$ for all $i$, $\glambdait\to\lambda_i$ for all $i$, and $\sum_{i=1}^\ell\glambdait\xx'_{it}\to F(\zbar)$,
thus proving \eqref{eq:thm:F:seqconv:3}
by \Cref{thm:seqsum-midrays}.

\pfpart{Part ``$\supseteq$'':}
Let $\zbar'\in\mul{\lambda_1}{\xbar'_1}\seqsum\dotsb\seqsum\mul{\lambda_\ell}{\xbar'_\ell}$.
Then by \Cref{thm:seqsum-midrays}
(with $k=m=\ell$),
for $i=1,\dotsc,\ell$, there exist sequences $\seq{\xx'_{it}}$ in $\Rm$ and $\seq{\glambdapit}$ in $\Rpos$ such that $\xx'_{it}\to\xbar'_i$ and
$\glambdapit\to\lambda_i$ for all $i$,
 $\sum_{i=1}^\ell\glambdapit=1$ for all $t$, and $\sum_{i=1}^\ell\glambdapit\xx'_{it}\to\zbar'$.
By \Cref{lem:aff}, this implies that there exists $\zbar\in\extspace$ such that
$F(\zbar)=\zbar'$ and there also exist sequences
$\seq{\xx_{it}}$ in $\Rn$ and $\seq{\glambdait}$ in $\Rpos$ such that $\xx_{it}\to\xbar_i$ and
$\glambdait\to\lambda_i$ for all $i$,
 $\sum_{i=1}^\ell\glambdait=1$ for all~$t$, and $\sum_{i=1}^\ell\glambdait\xx_{it}\to\zbar$.
By \Cref{thm:seqsum-midrays}, this implies that $\zbar\in\mul{\lambda_1}{\xbar_1}
       \seqsum \cdots \seqsum
       \mul{\lambda_\ell}{\xbar_\ell}$.
Thus,
\[
  \zbar'=F(\zbar)\in
  F\BigParens{
       \mul{\lambda_1}{\xbar_1}
       \seqsum \cdots \seqsum
       \mul{\lambda_\ell}{\xbar_\ell}
     }.%
\indexg{sequential convex combinations!affine maps@under affine maps|)}%
\indexg{affine maps, astral!sequential convex combinations under|)}%
\indexg{sequential convex combinations|)}%
\qedhere
\]
\end{proof-parts}
\end{proof}

\chapter{Astral cones}
\label{sec:cones}

In standard convex analysis, convex cones play an important role;
for instance, both the recession cone and barrier cone of a function
are convex cones.
In this chapter, we extend the notion of a cone to astral
space.
Similar to standard convex cones, convex astral cones
play a critical role in understanding properties of convex functions
extended to astral space.
For instance, in \Cref{sec:minimizers},
we introduce a convex astral cone that is an analogue of the
recession cone and characterizes minimizers of such functions.
Astral cones will also play an important part in deriving some of the
astral separation theorems in \Cref{sec:sep-thms}.

\section{Definition and basic properties}
\label{sec:cones-basic}

\indexg{cones, naive|(}%
\indexg{naive cones|(}%
Recall that a set $K\subseteq\Rn$ is a cone if
$\zero\in K$ and if $K$ is closed
under positive scalar multiplication, that is, if $\lambda\xx\in K$
for all $\xx\in K$ and all $\lambda\in\Rstrictpos$.
To extend this notion to astral space, we might, as a first attempt,
consider an exactly analogous definition.
Accordingly,
we say that a set $K\subseteq\extspace$ is a
\emph{naive cone} if $\zero\in K$ and if $\lambda\xbar\in K$
for all $\xbar\in K$ and all $\lambda\in\Rstrictpos$.
Clearly, a set $K\subseteq\Rn$ is a cone if and only if it is a
naive cone.
Nonetheless, this definition is not entirely satisfactory as an astral analogue of standard cones.

To see the issue, note that a standard cone $K$ in $\Rn$ has the property that
if a point $\xx\ne\zero$ is in~$K$, then $K$ includes a ray emanating from
the origin that passes through $\xx$ and
continues on to infinity.
This is not necessarily so for naive cones in astral space.
For instance, in $\Rext=\extspac{1}$, the set
$\{0,+\infty\}$ is a naive cone, but there is no ray within this set
that would connect $0$ and $+\infty$.
More generally, in $\extspace$, every set of icons $K\subseteq\corezn$
that includes the origin is a naive cone
(as follows from
\Cref{pr:i:8}\ref{pr:i:8d}),
but every point in the set $K$ (other than $\zero$) is topologically
disconnected from the origin.
This suggests that perhaps the naive cone does not properly generalize
the standard notion of a cone.

\indexg{rays (standard)|(}%
Stepping back, we can say that a nonempty set $K\subseteq\Rn$ is a cone
if for every point $\xx\in K$,
the entire ray in the direction of
$\xx$ is also included in $K$, meaning that $K$ includes the set
\begin{equation}  \label{eq:std-ray-defn}
\indexm{ray x300}{$\ray{\xx}$}{ray (standard)}%
  \ray{\xx}
  =
  \Braces{\lambda\xx :\: \lambda\in\Rpos}
  =
  \cone\{\xx\}.
\end{equation}
This is the usual way of defining a ray in $\Rn$, which,
when generalized directly to astral space, leads to the notion of a
naive cone, as just discussed.%
\indexg{cones, naive|)}%
\indexg{naive cones|)}

Alternatively, similarly to how we generalized the notion
of a segment and convexity in~\Cref{sec:def-convexity}, we can view $\ray{\xx}$
as the intersection of
all homogeneous closed halfspaces that include $\xx$
(this is implied by
\indexg{rays (standard)|)}%
\Cref{pr:con-int-halfspaces}\ref{roc:cor11.7.2}).
\indexg{rays, astral|(}%
\indexg{cones, astral|(}%
\indexg{halfspaces, astral!homogeneous|(}%
\indexg{homogeneous (halfspace or hyperplane)!astral|(}%
The astral ray is then obtained as an intersection of homogeneous closed
astral halfspaces, and an astral cone as a nonempty set that includes
an astral ray through each of its elements. We formalize this perspective
in a sequence of definitions that mirrors our development of astral
convexity in \Cref{sec:def-convexity}.

First, we say that a closed astral halfspace $\chsua$
is \emph{homogeneous} if $\beta=0$, that is,
if it has the form
\begin{equation}   \label{eqn:homo-halfspace}
  \chsuz=\{\xbar\in\extspace :\: \xbar\cdot\uu \leq 0\}
\end{equation}
for some $\uu\in\Rn\wo\{\zero\}$.
\indexg{outer conic hull|(}%
\indexg{outer conic hull!defined|(}%
The outer conic hull of a set is then defined like the outer
convex hull, but now using homogeneous halfspaces:

\begin{definition}
Let $S\subseteq\extspace$.
The \emph{outer conic hull} of $S$,
denoted $\oconich S$, is the
intersection of all homogeneous closed astral halfspaces that include
$S$; that is,
\[
\indexm{cone s800}{$\oconich S$}{outer conic hull}
  \oconich S
  =
  \bigcap{\BigBraces{ \chsuz:\:
      \uu\in\Rn\wo\{\zero\},\, S\subseteq \chsuz}}.%
\indexg{outer conic hull!defined|)}%
\indexg{halfspaces, astral!homogeneous|)}%
\indexg{homogeneous (halfspace or hyperplane)!astral|)}
\]
\end{definition}
Similar to outer convex hull, outer conic hull is a hull operator (for the collection
consisting of all possible intersections of homogeneous closed astral halfspaces).
\indexg{rays, astral!defined|(}%
Using the notion of an outer conic hull, we define astral rays:

\begin{definition}  \label{dfn:astral-ray}
Let $\xbar\in\extspace$.
The \emph{astral ray through $\xbar$},
denoted $\aray{\xbar}$,
is the intersection of all homogeneous closed astral
halfspaces that include $\xbar$;
that is,
\[
\indexm{ray x700}{$\aray{\xbar}$}{astral ray}%
  \aray{\xbar} = \oconich\{\xbar\}.%
\indexg{rays, astral!defined|)}%
\]
\end{definition}
Finally, we can define astral cones:
\begin{definition}
\indexg{cones, astral!defined|(}%
A nonempty set $K\subseteq\extspace$ is an
\emph{astral cone} if
$K$ includes the astral ray through $\xbar$
for all $\xbar\in K$, that is, if
\[
   \aray{\xbar}\subseteq K
   \quad
   \text{for all $\xbar\in K$.}%
\indexg{cones, astral!defined|)}%
\indexg{cones, astral|)}%
  \]
\end{definition}

We study properties of astral cones and the outer conic hull.
\indexg{rays, astral!segment@as segment|(}%
\indexg{outer conic hull!outer convex hull@as outer convex hull|(}%
We begin by showing that
the outer conic hull of any set can be expressed in terms of an
outer convex hull.
Consequently, the astral ray through $\xbar\in\eRn$, $\aray{\xbar}$, is equal to the segment
joining $\zero$ and $\limray{\xbar}$.
Informally, this ray follows a path from
the origin through $\xbar$, continuing on in the ``direction'' of
$\xbar$, and reaching $\limray{\xbar}$ in the limit.
For a point $\xx\in\Rn$, $\aray{\xx}$ is thus the standard ray,
$\ray{\xx}$, together with $\limray{\xx}$, the one additional point in
its astral closure; that is, $\aray{\xx}=(\ray{\xx})\cup\{\limray{\xx}\}$.

For example, in $\Rext$, the only astral rays are
$\{ 0 \}$,
$[0,+\infty]$,
$[-\infty,0]$,
and the only astral cones are
$\{ 0 \}$,
$[0,+\infty]$,
$[-\infty,0]$,
and $\Rext$.

\begin{theorem}   \label{thm:oconichull-equals-ocvxhull}
  Let $S\subseteq\extspace$.
  Then
  \begin{equation}   \label{eq:thm:oconichull-equals-ocvxhull:2}
    \ohull{S}
    \subseteq
    \oconich{S}
    =
    \ohull(\{\zero\} \cup \lmset{S}).
  \end{equation}
  In particular, for all $\xbar\in\extspace$,
  $\aray{\xbar}=\lb{\zero}{\limray{\xbar}}$.
\end{theorem}

\begin{proof}

By definition,
the outer conic hull of $S$ is the intersection of all homogeneous
astral closed halfspaces that include $S$, while its outer convex hull is
the intersection of all astral closed halfspaces that include $S$, whether
homogeneous or not.
Therefore, the latter is included in the former;
that is, $\ohull{S}\subseteq \oconich{S}$.

It remains to prove the equality in
\eqref{eq:thm:oconichull-equals-ocvxhull:2}.
Let $U=\ohull(\{\zero\} \cup \lmset{S})$.
We aim to show that $\oconich{S}=U$.

Let $\xbar\in\oconich{S}$, and let $\uu\in\Rn$.
We claim that
\begin{equation}   \label{eq:thm:oconichull-equals-ocvxhull:1}
  \xbar\cdot\uu
  \leq
  \max\Braces{ 0, \; \sup_{\zbar\in S} \limray{\zbar}\cdot\uu }.
\end{equation}
This is immediate if $\uu=\zero$.
Further,
if $\zbar\cdot\uu>0$ for any $\zbar\in S$, then
$\limray{\zbar}\cdot\uu=+\infty$, implying that the right-hand side of
\eqref{eq:thm:oconichull-equals-ocvxhull:1} is also $+\infty$
and the inequality holds.

Otherwise, $\uu\neq\zero$ and
$\zbar\cdot\uu\leq 0$ for all $\zbar\in S$.
In this case,
$\limray{\zbar}\cdot\uu\leq 0$ for all $\zbar\in S$ so that the
right-hand side of
\eqref{eq:thm:oconichull-equals-ocvxhull:1}
is equal to~$0$.
Furthermore, in this case, $S$ is included in the homogeneous
closed astral halfspace
$\chsuz$ (Eq.~\ref{eqn:homo-halfspace}).
Therefore, by definition of outer conic hull,
$\xbar\in\oconich{S}\subseteq\chsuz$,
so $\xbar\cdot\uu\leq 0$.
Together with the preceding, this proves
\eqref{eq:thm:oconichull-equals-ocvxhull:1}
in this case as well.

Since
\eqref{eq:thm:oconichull-equals-ocvxhull:1}
holds for all $\uu\in\Rn$, it follows by
\Cref{pr:ohull-simplify}(\ref{pr:ohull-simplify:b},\ref{pr:ohull-simplify:a})
that $\xbar\in U$.
Thus, $\oconich{S}\subseteq U$.

To prove the reverse implication, suppose now that $\xbar\in U$.
Suppose $S$ is included in some homogeneous closed astral halfspace
$\chsuz$ where $\uu\in\Rn\wo\{\zero\}$.
Then for all $\zbar\in S$, $\zbar\cdot\uu\leq 0$, implying
$\limray{\zbar}\cdot\uu\leq 0$.
Since $\xbar\in U$, by
\Cref{pr:ohull-simplify}(\ref{pr:ohull-simplify:a},\ref{pr:ohull-simplify:b}),
\eqref{eq:thm:oconichull-equals-ocvxhull:1} must hold.
Furthermore, the preceding remarks imply that
its right-hand side is equal to~$0$.
Therefore, $\xbar\cdot\uu\leq 0$, so
$\xbar\in\chsuz$.
Since this holds for all homogeneous closed astral halfspaces
$\chsuz$ that include $S$, it follows that
$\xbar\in\oconich{S}$, and thus that
$U\subseteq\oconich{S}$, completing the proof of
\eqref{eq:thm:oconichull-equals-ocvxhull:2}.

Finally, for all $\xbar\in\extspace$, by the foregoing,
\[
  \aray{\xbar}
  =
  \oconich\{\xbar\}
  =
  \ohull(\{\zero\}\cup\{\limray{\xbar}\})
  =
  \ohull\{\zero,\limray{\xbar}\}
  =
  \lb{\zero}{\limray{\xbar}},
\]
proving the final claim.%
\indexg{outer conic hull|)}%
\indexg{rays, astral!segment@as segment|)}%
\indexg{outer conic hull!outer convex hull@as outer convex hull|)}%
\end{proof}

A standard ray can be expressed parametrically as $\ray\xx=\set{\lambda\xx:\:\lambda\in\Rpos}$.
\indexg{rays, astral!union of sequential multiples@as union of sequential multiples|(}%
\indexg{sequential multiples!astral rays as union of|(}%
An astral ray $\aray\xbar$ can be expressed
similarly as the union of multiples $\mul{\alpha}{\xbar}$ across all $\alpha\in\Rextpos$:

\begin{theorem}
\label{thm:aray}
Let $\xbar\in\eRn$. Then
\[
  \aray{\xbar}=\aray{\limray{\xbar}}
     =\bigcup\bigBraces{\mul{\alpha}{\xbar} :\: \alpha\in\Rextpos}.
\]
\end{theorem}

\begin{proof}
We have
\[
   \aray{\xbar}
     =\lb{\zero}{\limray{\xbar}}
     =\seg\bigParens{\zero, \limray{ (\limray{\xbar}) }}
     =\aray{\limray{\xbar}},
\]
with the first and last equalities by
\Cref{thm:oconichull-equals-ocvxhull},
and the second because
$\limray{ (\limray{\xbar}) }=\limray{\xbar}$.

To finish the proof, we show that
   $\lb{\zero}{\limray{\xbar}}
     =\bigcup_{\alpha\in\Rextpos}\mul{\alpha}{\xbar}$.
  \begin{proof-parts}
  \pfpart{Part ``$\subseteq$'':}
    Let $\zbar\in\lb{\zero}{\limray{\xbar}}$;
    we aim to show $\zbar\in\almul{\xbar}$ for some
    $\alpha\in\Rextpos$.
    By
    \Cref{cor:e:1},
    $\zbar=\lim\lambda_t\yy_t$ for some sequences $\seq{\lambda_t}$ in $[0,1]$
    and $\seq{\yy_t}$ in $\Rn$ such that $\yy_t\to\limray{\xbar}$
    (since the only span-bound sequence converging to $\zero$ is the
    all $\zero$ sequence).
    By \Cref{thm:seq-to-limray}, we can write $\seqeq{\yy_t}{\gamma_t\xx_t}$ for some sequence $\seq{\gamma_t}$ in $\Rstrictpos$ and some span-bound sequence $\seq{\xx_t}$ in $\Rn$ such that $\gamma_t\to+\infty$ and $\xx_t\to\xbar$.
    By \Cref{prop:strong:eq:propties}(\ref{prop:strong:eq:propties:c}),
     we have $\seqeq{\lambda_t\yy_t}{\lambda_t\gamma_t\xx_t}$, so
    $\lim\lambda_t\gamma_t\xx_t=\lim\lambda_t\yy_t=\zbar$
    (by \Cref{pr:eq-in-lim-same-lim}).

    Since $\lambda_t\gamma_t\ge 0$ for all $t$,
    and by compactness of $[0,+\infty]$
    (being a closed subset of $\Rext$),
    there must exist
    a subsequence of $\seq{\lambda_t\gamma_t}$ that converges to some $\alpha\in[0,+\infty]$. Discarding the remaining elements of the sequence, we thus have $\lambda_t\gamma_t\to\alpha$,
    $\xx_t\to\xbar$ and $\lambda_t\gamma_t\xx_t\to\zbar$, so
    $\zbar\in\mul{\alpha}{\xbar}$.

  \pfpart{Part ``$\supseteq$'':}
    Let $\zbar\in\mul{\alpha}{\xbar}$ for some $\alpha\in\Rextpos$;
  we aim to show $\zbar\in\lb{\zero}{\limray{\xbar}}$.
  If $\alpha=\omega$, then $\zbar=\limray{\xbar}$ by
  \Cref{thm:mul-char}(\ref{thm:mul-char:a}), implying
    $\zbar\in\lb{\zero}{\limray{\xbar}}$.
  We therefore assume henceforth that $\alpha\in\Rpos$.

  Since $\zbar\in\almul{\xbar}$, there exist a sequence $\seq{\alpha_t}$ in $\Rpos$ and a span-bound sequence $\seq{\xx_t}$ in $\Rn$ such that $\alpha_t\to\alpha$, $\xx_t\to\xbar$, and $\alpha_t\xx_t\to\zbar$.
Let $\beta=\sup\{\alpha_t :\: t\in\nats\}$, which is finite
since $\seq{\alpha_t}$ converges to finite $\alpha$.
For each $t$, let $\gamma_t=\max\{\beta,t\}$,
and let $\lambda_t=\alpha_t/\gamma_t$, implying
$\lambda_t\in[0,1]$.
By \Cref{thm:seq-to-limray:0}, $\gamma_t\xx_t\to\limray{\xbar}$
since $\gamma_t\rightarrow+\infty$.
Since also
$(1-\lambda_t)\cdot\zero+\lambda_t(\gamma_t\xx_t)=\alpha_t\xx_t\rightarrow\zbar$,
we conclude $\zbar\in\lb{\zero}{\limray{\xbar}}$
by \Cref{cor:e:1}.%
\indexg{rays, astral|)}%
\indexg{rays, astral!union of sequential multiples@as union of sequential multiples|)}%
\indexg{sequential multiples!astral rays as union of|)}%
\qedhere
  \end{proof-parts}
\end{proof}

\indexg{cones, astral|(}%
Astral cones have many of the same properties as standard cones in $\Rn$.
In particular, every astral cone is also a naive cone, and its
intersection with $\Rn$ is a standard cone:

\begin{corollary}  \label{pr:ast-cone-is-naive}
  Let $K\subseteq\extspace$ be an astral cone,
  and let $\xbar\in\extspace$.
  Suppose either $\xbar\in K$ or $\limray{\xbar}\in K$.
  Then
  $\alpha \xbar\in K$ for all $\alpha\in\Rextpos$.
  Consequently, $\zero\in K$,
  $K$ is a naive cone, and
  $K\cap\Rn$ is a cone in $\Rn$.
\end{corollary}

\begin{proof}
If $\xbar\in K$ or $\limray\xbar\in K$, then
for all $\alpha\in\Rextpos$,
\[
  \alpha\xbar
  \in
  \mul{\alpha}{\xbar}
  \subseteq
  \aray\xbar
  =
  \aray\limray\xbar
  \subseteq
  K,
\]
where the first inclusion is by 
\Cref{thm:mul-char},
the second inclusion and equality are by
\Cref{thm:aray},
and the last inclusion is because $K$ is an astral cone that includes
either $\xbar$ or $\limray{\xbar}$.

In particular, taking $\alpha=0$, this shows that $\zero\in K$
(since $K$ is nonempty),
and taking $\alpha\in\Rstrictpos$,
this shows that $K$ is a naive cone.
Applied to points in $K\cap\Rn$, it follows that this set
is a cone.
\end{proof}

\indexg{cones, astral!linear or affine maps and|(}%
\indexg{linear maps, astral!astral cone under|(}%
The image of an astral cone under a linear map is also an astral cone:

\begin{theorem}
  Let $K\subseteq\extspace$ be an astral cone,
  let $\A\in\Rmn$, and let $A:\extspace\rightarrow\extspac{m}$ be the
  associated linear map (so that $A(\xbar)=\A\xbar$ for
  $\xbar\in\extspace$).
  Then $A(K)$ is also an astral cone.
\end{theorem}

\begin{proof}
Since $K$ is nonempty, so is $A(K)$.

Let $\zbar\in A(K)$, implying $\zbar=A(\xbar)=\A\xbar$ for some
$\xbar\in K$.
Then
\begin{align*}
  \aray{\zbar}
  =
  \lb{\zero}{\limray{\zbar}}
  &=
  \seg\bigParens{\zero, \limray{(\A\xbar)}}
  \\
  &=
  \seg\bigParens{A(\zero), A(\limray{\xbar})}
  \\
  &=
  A\bigParens{\lb{\zero}{\limray{\xbar}}}
  =
  A(\aray{\xbar})
  \subseteq
  A(K).
\end{align*}
The first and last equalities are by
\Cref{thm:oconichull-equals-ocvxhull},
the third is by \Cref{pr:h:4}(\ref{pr:h:4g}),
and the fourth by \Cref{thm:e:9}.
The inclusion is because $K$ is an astral cone.
Thus, $A(K)$ is an astral cone.%
\indexg{linear maps, astral!astral cone under|)}%
\end{proof}

\indexg{affine maps, astral!astral cones and|(}%
Similarly, the inverse image of an astral cone $K$ under a linear map is an astral cone. This can be further generalized to hold for inverse images of affine maps of the form $\zbar\mapsto\ebar\plusl\A\zbar$ where $\ebar$ is an icon in $K$.

\begin{theorem}  \label{thm:affine-inv-ast-cone}
  Let $\A\in\R^{m\times n}$, $\ebar\in\corez{m}$, and
  let $F:\extspace\rightarrow\extspac{m}$ be the
  affine map
  $F(\zbar)=\ebar\plusl \A\zbar$
  for all $\zbar\in\extspace$.
  Let $K\subseteq\extspac{m}$ be an astral cone,
  and assume $\ebar\in K$.
  Then
  $\finv(K)$ is also an astral cone.
\end{theorem}

\begin{proof}
Since $F(\zero)=\ebar\in K$,
$\finv(K)$ is not empty.

Let $\zbar\in\finv(K)$, so that $F(\zbar)\in K$.
Let
\[
  S
  =
  \aray{F(\zbar)}
  =
  \seg\bigParens{\zero,\limray{F(\zbar)}}
\]
(using \Cref{thm:oconichull-equals-ocvxhull} for the second equality).
Note that
\begin{equation}   \label{eq:thm:affine-inv-ast-cone:1}
  \limray{F(\zbar)}
  =
  \limray{(\ebar\plusl \A \zbar)}
  =
  \ebar \plusl \A (\limray\zbar)
  =
  F(\limray{\zbar}),
\end{equation}
where the second equality uses
Propositions~\ref{pr:i:8}(\ref{pr:i:8d})
and~\ref{pr:h:4}(\ref{pr:h:4g}).
Then
\[
   F(\aray{\zbar})
   =
   F\bigParens{\lb{\zero}{\limray{\zbar}}}
   =
   \seg\bigParens{F(\zero),F(\limray{\zbar})}
   =
   \seg\bigParens{\ebar,\limray{F(\zbar)}}
   \subseteq
   S
   \subseteq
   K.
\]
The first equality is by \Cref{thm:oconichull-equals-ocvxhull},
the second by
\Cref{thm:e:9},
and the third by \eqref{eq:thm:affine-inv-ast-cone:1}.
The first inclusion is by
\Cref{thm:decomp-seg-at-inc-pt}(\ref{thm:decomp-seg-at-inc-pt:a})
since $\ebar\in S$
(by \Cref{cor:d-in-lb-0-dplusx},
since $\limray{F(\zbar)}=\ebar \plusl \A (\limray\zbar)$ by
Eq.~\ref{eq:thm:affine-inv-ast-cone:1}).
The last inclusion is because $K$ is an astral cone
that includes $F(\zbar)$.
Thus, $\aray{\zbar}\subseteq \finv(K)$.
Therefore, $\finv(K)$ is an astral cone.
\end{proof}

\Cref{thm:affine-inv-ast-cone} is false in general if $\ebar$
is replaced by a general astral point that is not necessarily an icon.
For instance, in $\Rext$, let $K=[0,+\infty]$, which is an astral
cone, and let $F(\barx)=\barx+1$ for $\barx\in\Rext$.
Then $\finv(K)=[-1,+\infty]$, which is not an astral cone.%
\indexg{affine maps, astral!astral cones and|)}%
\indexg{cones, astral!linear or affine maps and|)}

\section{Convex astral cones}

Of particular importance are convex astral cones, that is,
astral cones that are convex.

\begin{proposition}  \label{pr:astral-cone-props}
  Let $S\subseteq\extspace$.
  \begin{letter-compact}
  \item  \label{pr:astral-cone-props:b}
    The intersection of an arbitrary collection of astral cones is an
    astral cone.
  \item  \label{pr:astral-cone-props:c}
    Suppose $S$ is convex.
    Then $S$ is a convex astral cone if and only if $\zero\in S$ and
    $\lmset{S}\subseteq S$.
  \item  \label{pr:astral-cone-props:e}
    Any intersection of homogeneous closed astral halfspaces in $\extspace$
    is a closed convex astral cone.
    Therefore, $\oconich{S}$ is a closed convex astral cone,
    as is $\aray{\xbar}$ for all $\xbar\in\extspace$.
  \end{letter-compact}
\end{proposition}

\begin{proof}
~

\begin{proof-parts}
\pfpart{Part~(\ref{pr:astral-cone-props:b}):}
Similar to the proof of
\Cref{pr:e1}(\ref{pr:e1:b}).

\pfpart{Part~(\ref{pr:astral-cone-props:c}):}
If $S$ is an astral cone then
$\aray{\xbar}=\lb{\zero}{\limray{\xbar}}$ is included in $S$ for all
$\xbar\in S$, proving
the ``only if'' direction.

For the converse, suppose $\zero\in S$ and $\lmset{S}\subseteq S$.
Then for
all $\xbar\in S$, $\limray{\xbar}\in S$, implying
$\aray{\xbar}=\lb{\zero}{\limray{\xbar}}\subseteq S$
since $S$ is convex.
Therefore, $S$ is an astral cone.

\pfpart{Part~(\ref{pr:astral-cone-props:e}):}
Let $\chsuz$ be a homogeneous closed astral halfspace
(as in Eq.~\ref{eqn:homo-halfspace}),
for some $\uu\in\Rn\wo\{\zero\}$.
Then $\chsuz$ is convex by
\Cref{pr:e1}(\ref{pr:e1:c}).
Further, $\chsuz$ includes $\zero$, and if $\xbar\in\chsuz$ then
$\xbar\cdot\uu\leq 0$, implying $\limray{\xbar}\cdot\uu\leq 0$,
and so that $\limray{\xbar}\in\chsuz$.
Therefore, $\chsuz$ is a convex astral cone by
part~(\ref{pr:astral-cone-props:c}).

That an intersection of homogeneous closed astral halfspaces is a convex
astral cone now follows by
part~(\ref{pr:astral-cone-props:b})
and \Cref{pr:e1}(\ref{pr:e1:b}).
By definition, this includes $\oconich{S}$ and $\aray{\xbar}$ for any
$\xbar\in\extspace$.
\qedhere
\end{proof-parts}
\end{proof}

\indexg{cones, astral!convex hull of|(}%
The convex hull of an astral cone is a convex astral cone:

\begin{theorem}  \label{thm:cvxh-ast-cone-is-ast-cone}
  Let $K\subseteq\extspace$ be an astral cone.
  Then its convex hull, $\conv{K}$, is also an astral cone
  (and therefore a convex astral cone).
\end{theorem}

\begin{proof}
First, $\zero\in K$ (by \Cref{pr:ast-cone-is-naive}),
so $\zero\in\conv{K}$.
Let $\xbar\in \conv{K}$.
By
\Cref{pr:astral-cone-props}(\ref{pr:astral-cone-props:c}),
it suffices to show that $\limray{\xbar}\in\conv{K}$.

By \Cref{thm:convhull-of-simpices},
there exists a finite subset $V\subseteq K$ such that
$\xbar\in\ohull{V}$.
Thus,
\begin{equation}  \label{eq:thm:cvxh-ast-cone-is-ast-cone:1}
  \xbar
  \in
  \ohull{V}
  \subseteq
  \oconich{V}
  =
  \ohull(\{\zero\}\cup \lmset{V})
  \subseteq
  \conv{K}.
\end{equation}
The second inclusion and the equality are by
\Cref{thm:oconichull-equals-ocvxhull}.
The last inclusion is by
\Cref{thm:convhull-of-simpices}
since $\{\zero\}\cup \lmset{V}$ is of finite cardinality and,
by \Cref{pr:ast-cone-is-naive},
is included in $K$ since $K$ is an astral cone.
Since $\oconich{V}$ is an astral cone
(by
\Cref{pr:astral-cone-props}\ref{pr:astral-cone-props:e}),
it follows that $\limray{\xbar}\in\oconich{V}\subseteq\conv{K}$
by \Cref{pr:ast-cone-is-naive}
and
\indexg{cones, astral!convex hull of|)}%
\eqref{eq:thm:cvxh-ast-cone-is-ast-cone:1}.
\end{proof}

\indexg{cones, astral!sequential sum of|(}%
\indexg{sequential sum!astral cones@of astral cones|(}%
Analogous to similar results in standard convex analysis,
the sequential sum of two astral cones is also an astral cone.
Furthermore,
the sequential sum of finitely many convex astral cones is equal to the
convex hull of their union:

\begin{theorem}   \label{thm:seqsum-ast-cone}
  ~

  \begin{letter-compact}
  \item   \label{thm:seqsum-ast-cone:a}
    Let $K_1$ and $K_2$ be astral cones in $\extspace$.
    Then $K_1\seqsum K_2$ is also an astral cone.
  \item   \label{thm:seqsum-ast-cone:b}
    Let $K_1,\ldots,K_m$ be convex astral cones in $\extspace$,
    with $m\geq 1$.
    Then
    \[
       K_1\seqsum\dotsb\seqsum K_m
       =
       \conv(K_1\cup\dotsb\cup K_m),
    \]
    which is also a convex astral cone.
  \end{letter-compact}
\end{theorem}

\begin{proof}
~

\begin{proof-parts}
\pfpart{Part~(\ref{thm:seqsum-ast-cone:a}):}
Let $Z=K_1\seqsum K_2$, which is nonempty
by \Cref{prop:seqsum-multi}(\ref{prop:seqsum-multi:nonempty})
since $K_1$ and $K_2$ are nonempty.
Let $\zbar\in Z$.
We aim to show that $\aray{\zbar}\subseteq Z$, which will
prove that $Z$ is an astral cone.

Since $\zbar\in Z$, we must have $\zbar\in \xbar_1\seqsum\xbar_2$
for some $\xbar_1\in K_1$ and $\xbar_2\in K_2$.
For $i\in\{1,2\}$, let
$K'_i=\aray{\xbar_i}=\lb{\zero}{\limray{\xbar_i}}$,
and let $Z'=K'_1\seqsum K'_2$.
Note that $K'_i\subseteq K_i$ since $K_i$ is an astral cone,
so $Z'\subseteq Z$
(by \Cref{pr:seqsum-props}\ref{pr:seqsum-props:b}).
Further, each $K'_i$ is closed and convex
(\Cref{pr:e1}\ref{pr:e1:ohull}),
implying $Z'$ is as well,
by
\Cref{prop:seqsum-multi}(\ref{prop:seqsum-multi:closed},\ref{prop:seqsum-multi:convex}).
Therefore, for all $t$,
\[
  t\zbar
  \in
  t (\xbar_1\seqsum\xbar_2)
  =
  (t\xbar_1)\seqsum(t\xbar_2)
  \subseteq
  K'_1\seqsum K'_2.
\]
The equality is by
\Cref{thm:distrib-seqsum}.
The second inclusion is because,
for $i\in\{1,2\}$,
$\xbar_i\in K'_i$ and
each $K'_i$ is an astral cone
(by \Cref{pr:astral-cone-props}\ref{pr:astral-cone-props:e}),
and therefore also a naive cone
(by \Cref{pr:ast-cone-is-naive}).
Taking limits, it follows that
$\limray{\zbar}=\lim(t\zbar)\in Z'$,
since $Z'$ is closed.
In addition, $\zero\in Z'$ since $\zero$ is in both $K'_1$ and $K'_2$.
Thus, $\aray{\zbar}=\lb{\zero}{\limray{\zbar}}\subseteq Z'\subseteq Z$
since $Z'$ is convex.

\pfpart{Part~(\ref{thm:seqsum-ast-cone:b}):}
Let $J=K_1\seqsum\dotsb\seqsum K_m$
and let $U=K_1\cup \dotsb \cup K_m$.
Then $J$ is convex by \Cref{prop:seqsum-multi}(\ref{prop:seqsum-multi:convex}),
and is an astral cone by part~(\ref{thm:seqsum-ast-cone:a})
(applied repeatedly).
We aim to show that $J=\conv{U}$.

We have
\[
   J
   \subseteq
   m \conv\BiggParens{\bigcup_{i=1}^m K_i}
   =
   \conv\BiggParens{\bigcup_{i=1}^m (m K_i)}
   \subseteq
   \conv\BiggParens{\bigcup_{i=1}^m K_i}
   =
   \conv{U}.
\]
The first inclusion is by
\Cref{thm:seqsum-in-union}.
The first
equality is by
\Cref{cor:thm:e:9b}
(applied to the linear map $\xbar\mapsto m\xbar$).
The last inclusion is because each $K_i$ is an astral cone
(and by \Cref{pr:ast-cone-is-naive}),
implying
$m K_i\subseteq K_i$ for $i=1,\ldots,m$.

On the other hand, since $\zero$ is in each astral cone $K_i$,
\[
   K_1
   =
   K_1 \seqsum \underbrace{ \set{\zero} \seqsum \dotsb \seqsum \set{\zero} }_{m-1}
   \subseteq
   K_1 \seqsum (K_2\seqsum \dotsb \seqsum K_m)
   =
   J.
\]
By symmetry, this also shows that $K_i\subseteq J$
for $i=1,\ldots,m$.
Thus, $U\subseteq J$, implying
$\conv{U}\subseteq J$, since $J$ is convex.
This completes the proof.%
\indexg{cones, astral!sequential sum of|)}%
\indexg{sequential sum!astral cones@of astral cones|)}%
\qedhere
\end{proof-parts}
\end{proof}

\indexg{cones, astral!convexity of|(}%
We next show that an astral cone is convex if and only if it is closed
under the sequential sum operation.
An analogous fact from standard convex analysis is given in
\Cref{pr:scc-cone-elts}(\ref{pr:scc-cone-elts:d}).

\begin{theorem}  \label{thm:ast-cone-is-cvx-if-sum}
Let $K\subseteq\extspace$ be an astral cone.
Then $K$ is convex if and only if
$K\seqsum K\subseteq K$
(that is, if and only if
$\xbar\seqsum\ybar\subseteq K$ for all $\xbar,\ybar\in K$).
\end{theorem}

\begin{proof}
Suppose $K$ is convex.
Then
\[
   K\seqsum K
   =
   \conv(K\cup K)
   =
   K,
\]
where the first equality is by
\Cref{thm:seqsum-ast-cone}(\ref{thm:seqsum-ast-cone:b}),
and the second is because $K$ is convex.

For the converse, suppose
$K\seqsum K\subseteq K$.
Let $\xbar,\ybar\in K$.
To show that $K$ is convex, we aim to show that
$\lb{\xbar}{\ybar}\subseteq K$.
Let
$J = (\aray{\xbar})\seqsum(\aray{\ybar})$.
Then $J\subseteq K\seqsum K\subseteq K$
since $K$ is an astral cone
(so that $\aray{\xbar}\subseteq K$
and $\aray{\ybar}\subseteq K$).
Further, $\xbar\in\aray{\xbar}$ and $\zero\in\aray{\ybar}$ so
$\xbar\in\xbar\seqsum\zero\subseteq J$; similarly, $\ybar\in J$.
We also have that $\aray{\xbar}$ and $\aray{\ybar}$ are convex
(by
\Cref{pr:astral-cone-props}\ref{pr:astral-cone-props:e}),
so $J$ is convex by
\Cref{prop:seqsum-multi}(\ref{prop:seqsum-multi:convex}).
Consequently, $\lb{\xbar}{\ybar}\subseteq J \subseteq K$,
completing the proof.%
\indexg{cones, astral!convexity of|)}%
\end{proof}

\indexg{cones, astral!convex hull of icons@as convex hull of icons|(}%
As mentioned earlier, every set of icons that includes the origin
is a naive cone.
This is certainly not the case for astral cones; indeed, the only set
of icons that is an astral cone is the singleton $\{\zero\}$.
Nevertheless, as we show next,
the convex hull of any set of icons that includes the
origin is a convex astral cone.
Moreover, every convex astral cone can be expressed in this way, and
more specifically, as the convex hull of all the icons that it
includes.

\begin{theorem}  \label{thm:ast-cvx-cone-equiv}
  Let $K\subseteq\extspace$.
  Then the following are equivalent:
  \begin{letter-compact}
  \item  \label{thm:ast-cvx-cone-equiv:a}
    $K$ is a convex astral cone.
  \item  \label{thm:ast-cvx-cone-equiv:b}
    $K=\conv{E}$ for some set of icons $E\subseteq\corezn$ that
    includes the origin.
  \item  \label{thm:ast-cvx-cone-equiv:c}
    $\zero\in K$ and $K=\conv(K \cap \corezn)$.
  \end{letter-compact}
\end{theorem}

\begin{proof}
  ~

\begin{proof-parts}
\pfpart{%
  (\ref{thm:ast-cvx-cone-equiv:c})
  $\Rightarrow$
  (\ref{thm:ast-cvx-cone-equiv:b}):
}
Setting $E=K\cap \corezn$, which must include $\zero$ by assumption,
this implication follows immediately.

\pfpart{%
  (\ref{thm:ast-cvx-cone-equiv:b})
  $\Rightarrow$
  (\ref{thm:ast-cvx-cone-equiv:a}):
}
Suppose $K=\conv{E}$ for some $E\subseteq\corezn$ with $\zero\in E$.
Then $K$ is convex and includes $\zero$.
Let $\xbar\in K$.
By
\Cref{pr:astral-cone-props}(\ref{pr:astral-cone-props:c}),
to show $K$ is an astral cone, it suffices to show that
$\limray{\xbar}\in K$.
Since $\xbar\in\conv{E}$, by
\Cref{thm:convhull-of-simpices},
there exists a finite subset $V\subseteq E$ such that
$\xbar\in\ohull{V}$.
Without loss of generality, we assume $\zero\in V$.
We claim that $\limray{\xbar}$ is also in $\ohull{V}$.
To see this, let $\uu\in\Rn$.
If $\xbar\cdot\uu\leq 0$, then
\[
  \limray{\xbar}\cdot\uu
  \leq
  0
  \leq
  \max_{\ebar\in V} \ebar\cdot\uu
\]
where the second inequality is because $\zero\in V$.
Otherwise, if $\xbar\cdot\uu>0$ then
\[
  0
  <
  \xbar\cdot\uu
  \leq
  \max_{\ebar\in V} \ebar\cdot\uu,
\]
with the second inequality following from
\Cref{pr:ohull-simplify}(\ref{pr:ohull-simplify:a},\ref{pr:ohull-simplify:b}).
Since all of the points in $V$ are icons, this implies that
$\ebar\cdot\uu=+\infty$ for some $\ebar\in V$
(by
\Cref{pr:icon-equiv}\ref{pr:icon-equiv:a}\ref{pr:icon-equiv:b}).
Therefore,
\[
  \limray{\xbar}\cdot\uu
  \leq
  \max_{\ebar\in V} \ebar\cdot\uu
\]
in this case as well.
From
\Cref{pr:ohull-simplify}(\ref{pr:ohull-simplify:b},\ref{pr:ohull-simplify:a}),
it follows that $\limray{\xbar}\in\ohull{V}\subseteq\conv{E}$,
as claimed, and so that $K$ is an astral cone.

\pfpart{%
  (\ref{thm:ast-cvx-cone-equiv:a})
  $\Rightarrow$
  (\ref{thm:ast-cvx-cone-equiv:c}):
}
Suppose $K$ is a convex astral cone.
Then $\zero\in K$ by \Cref{pr:ast-cone-is-naive}.
Let $E=K\cap\corezn$.  %
Since $E\subseteq K$ and $K$ is convex,
$\conv{E}\subseteq K$
(by
\Cref{pr:conhull-prop}\ref{pr:conhull-prop:aa}).
It remains to show the reverse implication.

Let $\xbar\in K$.
Then both $\zero$ and $\limray{\xbar}$ are in $K$ by
\Cref{pr:ast-cone-is-naive}.
Since they are both icons
(\Cref{pr:i:8}\ref{pr:i:8-infprod}),
they also are both in $E$.
Therefore,
\[
  \xbar
  \in
  \aray{\xbar}
  =
  \lb{\zero}{\limray{\xbar}}
  \subseteq
  \conv{E},
\]
where the equality is by
\Cref{thm:oconichull-equals-ocvxhull},
and the last inclusion is because $\conv{E}$ is convex
and includes $\zero$ and $\limray{\xbar}$.
Thus, $K\subseteq\conv{E}$, completing the proof.%
\indexg{cones, astral|)}%
\indexg{cones, astral!convex hull of icons@as convex hull of icons|)}%
\qedhere
\end{proof-parts}
\end{proof}

\section{Conic hull operations}
\label{sec:conic:hull}

We next take a closer look at conic hull operations in astral space.
\indexg{outer conic hull|(}%
To begin, analogous to \Cref{thm:e:2}, we show that a convex astral
cone includes the outer conic hull of every finite subset:

\begin{theorem}  \label{thm:oconich-fin-subset}
  Let $K\subseteq\extspace$ be a convex astral cone, and let
  $V\subseteq K$ be a finite subset.
  Then $\oconich{V}\subseteq K$.
\end{theorem}

\begin{proof}
\mathtogether%
By \Cref{pr:ast-cone-is-naive},
$\zero\in K$, and
if $\xbar\in V\subseteq K$ then $\limray{\xbar}\in K$;
thus, $\{\zero\} \cup \lmset{V}\subseteq K$.
It follows that
\[
  \oconich{V}
  =
  \ohull(\{\zero\} \cup \lmset{V})
  \subseteq
  K,
\]
where the equality is by \Cref{thm:oconichull-equals-ocvxhull},
and the inclusion is by \Cref{thm:e:2} since
$K$ is convex.
\end{proof}

\indexg{outer conic hull!sets in rn@for sets in $\Rn$|(}%
As shown next,
the outer conic hull of a set in $\Rn$ is simply equal to the
astral closure of its (standard) conic hull.
If the set is already a convex cone, then this is simply its closure.

\begin{theorem}  \label{thm:out-conic-h-is-closure}
  Let $S\subseteq\Rn$.
  Then $\oconich{S}=\clcone{S}$.
  Consequently, if $K\subseteq\Rn$ is a convex cone in $\Rn$,
  then $\oconich{K}=\Kbar$.
\end{theorem}

\begin{proof}
Let $\chsuz$ be a homogeneous closed astral halfspace
(as in Eq.~\ref{eqn:homo-halfspace}), for some
$\uu\in\Rn\wo\{\zero\}$.
If $S\subseteq\chsuz$, then $S$ is also included in the standard
homogeneous closed halfspace $\chsuz\cap\Rn$ in $\Rn$.
By \Cref{pr:con-int-halfspaces}(\ref{roc:cor11.7.2}),
this implies that
\[
  \cone{S}
  \subseteq
  \cl(\cone{S})
  \subseteq
  \chsuz\cap\Rn
  \subseteq
  \chsuz.
\]
Since this holds for all homogeneous closed astral halfspaces that include
$S$, it follows that $\cone{S}\subseteq\oconich{S}$, by definition of
outer conic hull.
Therefore,
$\clcone{S}\subseteq\oconich{S}$
since $\oconich{S}$ is closed in $\extspace$
(being an intersection of closed astral halfspaces).

For the reverse inclusion, we show that
$\oconich{S}$ is included in $\ohull(\cone{S})$, which is the same
as $\clcone{S}$ by \Cref{thm:e:6}.
To do so,
let $\chsub$ be a closed astral halfspace
(as in Eq.~\ref{eq:chsua-defn}) that includes $\cone{S}$,
for some $\uu\in\Rn\wo\{\zero\}$ and $\beta\in\R$.
Since $\zero\in\cone{S}$, this implies that
$0=\zero\cdot\uu\leq\beta$.
Furthermore, for all $\xx\in\cone{S}$, we claim that
$\xx\cdot\uu\leq 0$.
This is because for all $\lambda\in\Rstrictpos$,
$\lambda\xx$ is also in $\cone{S}$, implying
$\lambda\xx\cdot\uu\leq \beta$.
Therefore, if $\xx\cdot\uu>0$ then the left-hand side of this
inequality could be made arbitrarily large in the limit
$\lambda\rightarrow+\infty$, contradicting that it is bounded by
$\beta$.

Thus,
\[
  S
  \subseteq
  \cone{S}
  \subseteq
  \chsuz.
\]
This implies that
\[
  \oconich{S}
  \subseteq
  \chsuz
  \subseteq
  \chsub.
\]
The first inclusion is by definition of outer conic hull since
$S\subseteq \chsuz$.
The second inclusion is because $\beta\geq 0$.
Since this holds for all closed astral halfspaces $\chsub$ that include $\cone{S}$,
it follows, by definition of outer convex hull, that
\[
  \oconich{S}
  \subseteq
  \ohull(\cone{S})
  =
  \clcone{S},
\]
with the equality from \Cref{thm:e:6}.
This completes the proof that $\oconich{S}=\clcone{S}$.

If $K$ is a convex cone in $\Rn$, then
$\cone{K}=K$, yielding the final claim.%
\indexg{outer conic hull!sets in rn@for sets in $\Rn$|)}%
\indexg{outer conic hull|)}%
\end{proof}

\indexg{conic hull, astral|(}%
\indexg{conic hull, astral!defined|(}%
Analogous to the standard conic hull operation
(\Cref{sec:prelim:cones}), we define the astral conic hull of a
set to be the smallest convex astral cone that includes the set:

\begin{definition}
  Let $S\subseteq\extspace$.
  The \emph{astral conic hull} of $S$, denoted
\indexm{cone s600}{$\acone{S}$}{astral conic hull}%
  $\acone{S}$,
  is the intersection of all convex astral cones that include $S$.%
\indexg{conic hull, astral!defined|)}%
\end{definition}

By Propositions~\ref{pr:e1}(\ref{pr:e1:b})
and~\ref{pr:astral-cone-props}(\ref{pr:astral-cone-props:b}),
the astral conic hull of any set is indeed a convex astral cone.
Further,
this operation is the hull operator for the set of all convex astral
cones, and so has properties like those given in
\Cref{pr:conhull-prop} for convex hull, as listed in the
next proposition.

\begin{proposition}  \label{pr:acone-hull-props}
  Let $S,U\subseteq\extspace$.
  \begin{letter-compact}
  \item  \label{pr:acone-hull-props:a}
    If $S\subseteq U$ and $U$ is a convex astral cone, then
    $\acone{S}\subseteq U$.
  \item  \label{pr:acone-hull-props:b}
    If $S\subseteq U$ then $\acone{S}\subseteq\acone{U}$.
  \item  \label{pr:acone-hull-props:c}
    If $S\subseteq U\subseteq\acone{S}$, then
    $\acone{U}=\acone{S}$.
  \item  \label{pr:acone-hull-props:d}
    $\acone{S}\subseteq\oconich{S}$ with equality if $|S|<+\infty$.
  \end{letter-compact}
\end{proposition}

\begin{proof}
Similar to that of \Cref{pr:conhull-prop}, where,
in the proof of part~(\ref{pr:acone-hull-props:d}), we use
\Cref{thm:oconich-fin-subset} instead of
\Cref{thm:e:2}.
\end{proof}

In standard convex analysis, the conic hull of a set $S\subseteq\Rn$
is equal to the convex hull of the origin
and all the rays $\ray{\xx}$ with $\xx\in S$.
The astral conic hull of a set $S\subseteq\extspace$ can be expressed
analogously, namely, as the convex hull of the origin
and all the astral rays $\aray{\xbar}$ with $\xbar\in S$.
\indexg{conic hull, astral!convex hull@as convex hull|(}%
More simply, as given in the next theorem,
this is the same as the convex hull of the origin
and all of the endpoints of such astral rays, that is,
all of the points $\limray{\xbar}$ with $\xbar\in S$.
This is because, by
including $\zero$ and $\limray{\xbar}$ in the convex hull, we also
implicitly include the segment between these points, which is exactly
the astral ray through $\xbar$.

Also, if the set $S$ is already an astral cone, then its astral conic
hull is the same as its convex hull.

\begin{theorem}  \label{thm:acone-char}
  Let $S\subseteq\extspace$.
  \begin{letter-compact}
  \item  \label{thm:acone-char:a}
    $\acone{S}=\conv(\{\zero\}\cup \lmset{S})$.
  \item  \label{thm:acone-char:b}
    If $S$ is an astral cone, then
    $\acone{S}=\conv{S}$.
  \end{letter-compact}
\end{theorem}

\begin{proof}
Let $E=\{\zero\}\cup \lmset{S}$, which is a set of icons
(by \Cref{pr:i:8}\ref{pr:i:8-infprod})
that includes the origin.

\begin{proof-parts}
\pfpart{Part~(\ref{thm:acone-char:a}):}
By
\Crefequiv{thm:ast-cvx-cone-equiv}{thm:ast-cvx-cone-equiv:b}{thm:ast-cvx-cone-equiv:a},
$\conv{E}$ is a convex astral cone.
Furthermore, for all $\xbar\in S$,
we have $\limray{\xbar}\in E\subseteq\conv{E}$,
implying that $\xbar$ is also in the astral cone $\conv{E}$ by
\Cref{pr:ast-cone-is-naive}.
Therefore, $S\subseteq\conv{E}$, implying
$\acone{S}\subseteq\conv{E}$ by
\Cref{pr:acone-hull-props}(\ref{pr:acone-hull-props:a}).

For the reverse inclusion, since $\acone{S}$ is an astral cone,
by \Cref{pr:ast-cone-is-naive},
it must include $\zero$ and $\limray{\xbar}$ for all
$\xbar\in\acone{S}$, and so also for all $\xbar\in S$.
Therefore, $E\subseteq\acone{S}$.
Since $\acone{S}$ is convex, it follows that
$\conv{E}\subseteq\acone{S}$
(by \Cref{pr:conhull-prop}\ref{pr:conhull-prop:aa}),
completing the proof.

\pfpart{Part~(\ref{thm:acone-char:b}):}
The set $S$ is included in $\acone{S}$, which is convex.
Therefore, $\conv{S}\subseteq\acone{S}$.

On the other hand, since $S$ is an astral cone, it includes $\zero$
and $\limray{\xbar}$ for $\xbar\in S$
(\Cref{pr:ast-cone-is-naive}).
Therefore, $E\subseteq S$, implying
$\acone{S}=\conv{E}\subseteq\conv{S}$ by
\indexg{conic hull, astral!convex hull@as convex hull|)}%
part~(\ref{thm:acone-char:a}).
\qedhere
\end{proof-parts}
\end{proof}

\indexg{conic hull, astral!union@of union|(}%
Combining
Theorems~\ref{thm:seqsum-ast-cone}(\ref{thm:seqsum-ast-cone:b})
and~\ref{thm:acone-char}(\ref{thm:acone-char:a})
yields the following identity for the astral conic hull of a union:

\begin{theorem}   \label{thm:decomp-acone}
  Let $S_1,S_2\subseteq\extspace$.
  Then
  \[
     (\acone{S_1}) \seqsum (\acone{S_2})
     =
     \acone(S_1\cup S_2).
  \]
\end{theorem}

\begin{proof}
Let $K_i=\acone{S_i}$ for $i\in\{1,2\}$.
Then
\begin{align*}
  K_1 \seqsum K_2
  &=
  \conv(K_1 \cup K_2)
  \\
  &=
  \conv\bigParens{\conv(\{\zero\}\cup\lmset{S_1})
               \,\cup\,
               \conv(\{\zero\}\cup\lmset{S_2})
              }
  \\
  &=
  \conv(\{\zero\}\cup \lmset{S_1} \cup\lmset{S_2})
  \\
  &=
  \acone(S_1\cup S_2).
\end{align*}
The first equality is by
\Cref{thm:seqsum-ast-cone}(\ref{thm:seqsum-ast-cone:b})
since $K_1$ and $K_2$ are convex astral cones.
The second and fourth equalities are both by
\Cref{thm:acone-char}(\ref{thm:acone-char:a}).
The third equality is by
\Cref{pr:conhull-prop}(\ref{pr:conhull-prop:c})
since
\begin{align*}
  \{\zero\}\cup \lmset{S_1} \cup\lmset{S_2}
  &\subseteq
  \conv(\{\zero\}\cup\lmset{S_1})
  \,\cup\,
  \conv(\{\zero\}\cup\lmset{S_2})
  \\
  &\subseteq
  \conv(\{\zero\}\cup\lmset{S_1}\cup\lmset{S_2})
\end{align*}
(with the second inclusion following from
\indexg{conic hull, astral!union@of union|)}%
\indexg{conic hull, astral|)}%
\Cref{pr:conhull-prop}\ref{pr:conhull-prop:b}).
\end{proof}

\indexg{outer conic hull!decomposed as rays|(}%
\indexg{rays, astral!sequential sum of|(}%
As corollary, we obtain that
the outer conic hull of a finite set of astral points can be decomposed
as a sequential sum of astral rays through those points, which can
also be written as a union of sequential conic combinations of those points:

\begin{corollary}
\label{cor:conic:comb}
  Let $\xbar_1,\dotsc,\xbar_m\in\extspace$.
  Then
  \begin{align*}
  &
     \oconich\{\xbar_1,\dotsc,\xbar_m\}
      =
     (\aray{\xbar_1})\seqsum\dotsb\seqsum(\aray{\xbar_m})
\\&\qquad{}
     =
     \bigcup
     \;
     \Bigl\{
       \mul{\alpha_1}{\xbar_1}
       \seqsum \cdots \seqsum
       \mul{\alpha_m}{\xbar_m}:\:
     \smash[t]{
       \alpha_1,\dotsc,\alpha_m\in
     \Rextpos
     \Bigr\}.
     }
  \end{align*}
\end{corollary}

\begin{proof}
We have
\begin{align*}
  \oconich\{\xbar_1,\dotsc,\xbar_m\}
  &=
  \acone\{\xbar_1,\dotsc,\xbar_m\}
  \\
  &=
  (\acone\{\xbar_1\})\seqsum\dotsb\seqsum(\acone\{\xbar_m\})
  \\
  &=
  (\aray{\xbar_1})\seqsum\dotsb\seqsum(\aray{\xbar_m})
  \\
  &=
  \BiggParens{\,\bigcup_{\alpha_1\in\Rextpos}\!\!\mul{\alpha_1}{\xbar_1}}
  \seqsum\cdots\seqsum
  \BiggParens{\,\bigcup_{\alpha_m\in\Rextpos}\!\!\mul{\alpha_m}{\xbar_m}}
  \\
  &=
     \bigcup
     \;
     \Bigl\{
       \mul{\alpha_1}{\xbar_1}
       \seqsum \cdots \seqsum
       \mul{\alpha_m}{\xbar_m}:\:
       \alpha_1,\dotsc,\alpha_m\in
     \Rextpos
     \Bigr\}.
\end{align*}
The first and third equalities are by
\Cref{pr:acone-hull-props}(\ref{pr:acone-hull-props:d})
(and by definition of astral rays, \Cref{dfn:astral-ray}).
The second is by
\Cref{thm:decomp-acone}.
The fourth is by \Cref{thm:aray}.
And the fifth is by
\Cref{pr:seqsum-props}(\ref{pr:seqsum-props:c}).%
\indexg{outer conic hull!decomposed as rays|)}%
\indexg{rays, astral!sequential sum of|)}%
\end{proof}

\indexg{outer conic hull!sequential characterization|(}%
Building on \Cref{cor:conic:comb}, we can also characterize
the outer conic hull of a finite set of points using limits
of sequences of standard conic combinations:

\begin{theorem}   \label{thm:oconic-hull-and-seqs}
  Let $V=\{\xbar_1,\dotsc,\xbar_m\}\subseteq\extspace$,
  and let $\zbar\in\extspace$.
  Then $\zbar\in\oconich{V}$ if
  and only if there exist sequences
  $\seq{\lambda_{it}}$ in $\Rpos$
  and span-bound sequences
  $\seq{\xx_{it}}$ in $\Rn$,
  for $i=1,\dotsc,m$,
  such that:
  \begin{item-compact}
  \item
    $\xx_{it}\rightarrow\xbar_i$
    for $i=1,\dotsc,m$.
  \item
    The sequence %
    $\zz_t=\sum_{i=1}^m \lambda_{it} \xx_{it}$ converges to
    $\zbar$.
  \end{item-compact}

  The same equivalence holds if each sequence $\seq{\lambda_{it}}$,
  for $i=1,\ldots,m$,
  is additionally required to converge to a limit in $[0,+\infty]$.
\end{theorem}

\begin{proof}
  ~
\begin{proof-parts}
\pfpart{``If'' ($\Leftarrow$):}
Suppose sequences as stated in the \namecref{thm:oconic-hull-and-seqs}
exist. By compactness of $[0,+\infty]$, for each $i=1,\dotsc,m$,
considered in turn,
the sequence $\seq{\lambda_{it}}$ must have a subsequence that converges to a point in $[0,+\infty]$; by discarding all other elements (and the corresponding elements of other sequences), we can therefore assume that for each $i=1,\dotsc,m$, we have $\lambda_{it}\to\alpha_i$ for some $\alpha_i\in[0,+\infty]$. Thus, 
\[
  \zbar
    \in \mul{\alpha_1}{\xbar_1}\seqsum\dotsb\seqsum\mul{\alpha_m}{\xbar_m}
    \subseteq \oconich{V},
\]
with the first inclusion by \Cref{thm:seqsum-midrays},
and the second by \Cref{cor:conic:comb}.

\pfpart{``Only if'' ($\Rightarrow$):}
Assume $\zbar\in\oconich{V}$. By \Cref{cor:conic:comb}, there exist $\alpha_1,\dotsc,\alpha_m\in[0,+\infty]$ such that
$\zbar\in\mul{\alpha_1}{\xbar_1}\seqsum\dotsb\seqsum\mul{\alpha_m}{\xbar_m}$,
and so by \Cref{thm:seqsum-midrays} there exist sequences that satisfy
the conditions stated in the \namecref{thm:oconic-hull-and-seqs} (including the additional convergence requirement).%
\indexg{outer conic hull!sequential characterization|)}%
\qedhere
\end{proof-parts}
\end{proof}

\indexg{rays, astral!sequential characterization|(}%
Since $\aray{\xbar}=\oconich\{\xbar\}$, we immediately obtain as
corollary:

\begin{corollary}   \label{cor:aray-and-seqs}
  Let $\xbar,\zbar\in\extspace$.
  Then $\zbar\in\aray{\xbar}$ if
  and only if there exist
  a sequence
  $\seq{\lambda_t}$ in $\Rpos$
  and a span-bound sequence $\seq{\xx_t}$ in $\Rn$
  such that
  $\xx_t\rightarrow\xbar$ and
  $\lambda_t\xx_t\rightarrow\zbar$.

  The same equivalence holds if the sequence $\seq{\lambda_t}$
  is additionally required to converge to a limit in $[0,+\infty]$.
\end{corollary}

Without the span-boundness requirement,
\Cref{thm:oconic-hull-and-seqs} and \Cref{cor:aray-and-seqs} might not hold.
For instance, in $\R$, suppose $x_t=1/t$ and $\lambda_t=t$.
Then $x_t\rightarrow 0$ and $\lambda_t x_t \rightarrow 1$, but
$1\not\in\aray{0}=\oconich\{0\}=\{0\}$.%
\indexg{rays, astral!sequential characterization|)}

\indexg{conic hull, astral!union of outer conic hulls@as union of outer conic hulls|(}%
\indexg{outer conic hull!astral conic hull as union of|(}%
Analogous to \Cref{thm:convhull-of-simpices}, the astral conic
hull of a set $S$ is equal to the union of all outer conic hulls of
all finite subsets of $S$:

\begin{theorem}   \label{thm:ast-conic-hull-union-fin}
  Let $S\subseteq\extspace$.
  Then $S$'s astral conic hull is equal to the union of all outer conic
  hulls of finite subsets of $S$.
  That is,
  \begin{equation}   \label{eq:thm:ast-conic-hull-union-fin:1}
    \acone{S}
    = \bigcup_{\scriptontop{V\subseteq S:}{|V|<+\infty}}  \oconich{V}.
  \end{equation}
\end{theorem}

\begin{proof}
Let $U$ denote the union on the right-hand side of
\eqref{eq:thm:ast-conic-hull-union-fin:1}.
First, if $V\subseteq S$ and $|V|<+\infty$ then
$\oconich{V}=\acone{V}\subseteq\acone{S}$
by
\Cref{pr:acone-hull-props}(\ref{pr:acone-hull-props:b},\ref{pr:acone-hull-props:d}).
Thus, $U\subseteq\acone{S}$.

For the reverse inclusion, we have
\begin{align*}
  \acone{S}
  =
  \conv(\{\zero\}\cup\lmset{S})
  &=
  \bigcup_{\scriptontop{V\subseteq S\cup\{\zero\}:}{|V|<+\infty}} \ohull(\lmset{V})
  \\
  &\subseteq
  \bigcup_{\scriptontop{V\subseteq S:}{|V|<+\infty}} \ohull(\{\zero\}\cup\lmset{V})
  \\
  &=
  \bigcup_{\scriptontop{V\subseteq S:}{|V|<+\infty}} \oconich{V}
  =
  U.
\end{align*}
The first equality is by
\Cref{thm:acone-char}(\ref{thm:acone-char:a}).
The second is by
\Cref{thm:convhull-of-simpices}.
The inclusion is because if $V\subseteq S\cup\{\zero\}$
then $V\subseteq \{\zero\}\cup V'$ where $V'=V\cap S\subseteq S$,
implying
$\ohull(\lmset{V})\subseteq \ohull(\{\zero\}\cup\lmset{V'})$.
The third equality is by \Cref{thm:oconichull-equals-ocvxhull}.%
\indexg{conic hull, astral!union of outer conic hulls@as union of outer conic hulls|)}%
\indexg{outer conic hull!astral conic hull as union of|)}%
\end{proof}

\section{Astral halflines}
\label{sec:hline}

\indexg{halflines (standard)|(}%
For $\xx,\vv\in\Rn$, the (standard)
\emph{halfline with endpoint $\xx$ in direction $\vv$},
denoted  $\hfline{\xx}{\vv}$,
is the set
\begin{equation}   \label{eq:std-halfline-defn}
\indexm{hline x v300}{$\hfline{\xx}{\vv}$}{halfline (standard)}%
  \hfline{\xx}{\vv}
  =
  \Braces{ \xx+\lambda \vv :\: \lambda \in \Rpos }.
\end{equation}
In this section, we study an astral analogue of this set, for an astral
``endpoint'' $\xbar\in\extspace$ and an astral ``direction''
$\vbar\in\extspace$.
This will turn out to be useful, for instance, in the study of astral
convex functions and their minimizers.

Each point $\xx+\lambda\vv$ on the halfline given in
\eqref{eq:std-halfline-defn} is a conic combination of the endpoint
$\xx$ and the direction $\vv$.
The halfline then is the union of all such conic combinations for all
\indexg{halflines (standard)|)}%
$\lambda\in\Rpos$.
\indexg{halflines, astral|(}%
In the same way, in astral space, we can form the sequential conic
combination
$\mul{1}{\xbar}\seqsum\mul{\alpha}{\vbar}=\xbar\seqsum\mul{\alpha}{\vbar}$
of the endpoint $\xbar$ with the direction $\vbar$,
and then form an astral halfline by taking the union over all
$\alpha\in\Rextpos$.
This leads to the following definition:

\begin{definition}   \label{def:ast-halfline:seg}
\indexg{halflines, astral!defined|(}%
Let $\xbar,\vbar\in\extspace$.
The
\emph{astral halfline with endpoint $\xbar$ in direction $\vbar$},
denoted $\ahfline{\xbar}{\vbar}$,
is the set
\[
\indexm{hline x v800}{$\ahfline{\xbar}{\vbar}$}{astral halfline}%
  \ahfline{\xbar}{\vbar}
  =
  \bigcup\bigBraces{
    \xbar\seqsum\mul{\alpha}{\vbar} :\:
    \alpha\in\Rextpos
  }.%
\indexg{halflines, astral!defined|)}%
\]
\end{definition}

In particular,
for $\xx,\vv\in\Rn$, this definition (combined with
\Cref{thm:mul-char}) implies that
\[
  \ahfline{\xx}{\vv}
  =
  \set{\xx+\alpha\vv:\:\alpha\in\Rpos}\,\cup\,(\xx\seqsum\limray\vv)
  =
  \hfline{\xx}{\vv}\,\cup\,\set{\limray\vv\plusl\xx}.
\]
Thus, in this case, the astral halfline is obtained by extending the
standard halfline with its endpoint at infinity,
$\limray\vv\plusl\xx$, corresponding to the limit of the points along the halfline.

The next \namecref{thm:ast-halfline:conic} gives several equivalent
formulations of an astral halfline $\ahfline{\xbar}{\vbar}$.
The first is as the sequential sum of $\xbar$ with
the astral ray $\aray{\vbar}$, analogous to how a standard halfline
can be written as $\hfline{\xx}{\vv}=\xx+(\ray{\vv})$.
The second is as the segment joining the endpoint $\xbar$ with
$\limray{\vbar}\plusl\xbar$, which could perhaps be interpreted as
the other endpoint of the halfline at infinity.
The third is as a single sequential convex combination of these two
endpoints, specifically, as their $0$-midsegment.

\begin{theorem}
\label{thm:ast-halfline:conic}
Let $\xbar,\vbar\in\extspace$. Then
\begin{equation}
  \ahfline{\xbar}{\vbar}
  =
  \xbar\seqsum (\aray{\vbar})
  =
  \lb{\xbar}{\limray{\vbar}\plusl\xbar}
  =
  \mul{1}{\xbar}\seqsum\mul{0}{(\limray{\vbar}\plusl\xbar)}.
  \label{eq:thm:ast-halfline:conic:2}
\end{equation}
\end{theorem}

We prove part of the theorem as a slightly generalized lemma for later
reference:

\begin{lemma}  \label{lem:ahfline-containments}
  Let $\xbar,\vbar,\ybar\in\extspace$.
  Then
  \begin{equation} \label{eq:lem:ahfline-containments:1}
    \xbar\seqsum (\aray{\vbar})
    \subseteq
    \mul{1}{\xbar}\seqsum\mul{0}{(\limray{\vbar}\plusl\ybar)}
    \subseteq
    \lb{\xbar}{\limray{\vbar}\plusl\ybar}.
  \end{equation}
\end{lemma}

\begin{proof}
The second inclusion in
\eqref{eq:lem:ahfline-containments:1}
is immediate from
\Cref{cor:conv-as-seqsum-midrays}.
It remains then only to prove the first inclusion.

We can write $\ybar=\ebar\plusl\qq$ for some $\ebar\in\corezn$ and
$\qq\in\Rn$.
Then
\[
  \xbar\seqsum (\aray{\vbar})
  =
  \xbar\seqsum\lb{\zero}{\limray{\vbar}}
  \subseteq
  \xbar\seqsum\lb{\zero}{\limray{\vbar}\plusl\ebar}
  =
  \mul{1}{\xbar}\seqsum\mul{0}{(\limray{\vbar}\plusl\ybar)}.
\]
The first equality is by 
\Cref{thm:oconichull-equals-ocvxhull},
The inclusion is by \Cref{cor:d-in-lb-0-dplusx}
since $\limray\vbar$ is an icon.
The last equality is by
\Cref{thm:mul-char}(\ref{thm:mul-char:b}), noting that
$\limray{\vbar}\plusl\ebar$ is the iconic part of
$\limray{\vbar}\plusl\ybar=\limray{\vbar}\plusl\ebar\plusl\qq$.
\end{proof}

\begin{proof}[Proof of \Cref{thm:ast-halfline:conic}:]
First, we have
\begin{align}
  \ahfline{\xbar}{\vbar}
  &=
  \bigcup_{\alpha\in\Rextpos}[\xbar\seqsum\mul{\alpha}{\vbar}]
  \nonumber
  \\
  &=
  \xbar\seqsum
  \BiggParens{\,
    \bigcup_{\alpha\in\Rextpos}\!\!\mul{\alpha}{\vbar}
  }
  =
  \xbar\seqsum (\aray{\vbar}),
  \nonumber
\end{align}
where the first equality is by definition of astral halfline,
the second is by \Cref{pr:seqsum-props}(\ref{pr:seqsum-props:c}),
and the third is by \Cref{thm:aray}.
This proves the first equality of
\eqref{eq:thm:ast-halfline:conic:2}.

Next, we show:

\begin{claimpx}   \label{cl:thm:ast-halfline:conic:1}
  $\lb{\xbar}{\limray{\vbar}\plusl\xbar} \subseteq \xbar\seqsum (\aray{\vbar}).$
\end{claimpx}

\begin{proofx}
Since any singleton and any astral ray are convex
(Propositions~\ref{pr:e1}\ref{pr:e1:ohull}
and~\ref{pr:astral-cone-props}\ref{pr:astral-cone-props:e}),
the set $\xbar\seqsum (\aray{\vbar})$ is also convex (by \Cref{prop:seqsum-multi}\ref{prop:seqsum-multi:convex}).
Moreover, since $\zero$ and $\limray\vbar$ are in $\aray\vbar$ (by \Cref{thm:aray}),
we have $\xbar\in\xbar\seqsum\zero\subseteq\xbar\seqsum (\aray{\vbar})$ and
    \[
      \limray\vbar\plusl\xbar
      \in\seg(\limray\vbar\plusl\xbar,\,\xbar\plusl\limray\vbar)
      =\xbar\seqsum\limray\vbar
      \subseteq\xbar\seqsum(\aray{\vbar}),
    \]
    where the equality is by \Cref{cor:seqsum-conseqs}(\ref{cor:seqsum-conseqs:a}).
Thus, $\xbar\seqsum (\aray{\vbar})$ includes both $\xbar$ and
$\limray\vbar\plusl\xbar$, and therefore, being convex, also includes
the segment joining them,
$\lb{\xbar}{\limray\vbar\plusl\xbar}$.
\end{proofx}

Combining \Cref{lem:ahfline-containments} (applied with $\ybar=\xbar$)
and \Cref{cl:thm:ast-halfline:conic:1} now yields
\begin{equation*}
  \xbar\seqsum (\aray{\vbar})
  \subseteq
  \mul{1}{\xbar}\seqsum\mul{0}{(\limray{\vbar}\plusl\xbar)}
  \subseteq
  \lb{\xbar}{\limray{\vbar}\plusl\xbar}
  \subseteq
  \xbar\seqsum (\aray{\vbar}),
\end{equation*}
proving the
remaining
equalities of
\eqref{eq:thm:ast-halfline:conic:2}, and completing the proof.
\end{proof}

\indexg{halflines, astral!sequential characterization|(}%
Much like the limit characterizations of
segments and rays (\Cref{cor:e:1,cor:aray-and-seqs}),
we next derive a limit characterization
of astral halflines showing that
a point is included in the astral halfline with
endpoint $\xbar$ in direction $\vbar$ if and only if it is the limit
of points on standard halflines with endpoints and directions
converging respectively to $\xbar$ and $\vbar$.

\begin{theorem}    \label{thm:ahfline-seq}
  Let $\xbar,\vbar,\zbar\in\extspace$.
  Then $\zbar\in\ahfline{\xbar}{\vbar}$
  if and only if there exist
  sequences
  $\seq{\xx_t}$ in $\Rn$,
  $\seq{\lambda_t}$ in $\Rpos$,
  and
  a span-bound sequence $\seq{\vv_t}$ in $\Rn$
  such that
  $\xx_t\rightarrow\xbar$,
  $\vv_t\rightarrow\vbar$,
  and
  $\xx_t+\lambda_t\vv_t\rightarrow\zbar$.
\end{theorem}

\begin{proof}
If $\zbar\in\ahfline{\xbar}{\vbar}$ then, by definition,
$\zbar\in\xbar\seqsum\mul{\alpha}{\vbar}=\mul{1}{\xbar}\seqsum\mul{\alpha}{\vbar}$
for some
$\alpha\in\Rextpos$.
Therefore, by \Cref{thm:seqsum-midrays}
(applied with $\xbar_1=\xbar$, $\alpha_1=1$, $\xbar_2=\vbar$, $\alpha_2=\alpha$, $k=1$),
there exist sequences $\seq{\xx_t},\seq{\lambda_t},\seq{\vv_t}$ with the required properties.

Conversely, suppose there exist sequences satisfying the stated conditions. By compactness of $[0,+\infty]$, we can pick a convergent subsequence of $\seq{\lambda_t}$; discarding all other terms, we assume that $\lambda_t\to\alpha$ for some $\alpha\in[0,+\infty]$. Thus, the sequences satisfy the conditions of \Cref{thm:seqsum-midrays}
(with $\xx_{1t}=\xx_t$, $\lambda_{1t}=1$, $\xx_{2t}=\vv_t$,
$\lambda_{2t}=\lambda_t$, $k=1$), implying
$\zbar\in
\mul{1}{\xbar}\seqsum\mul{\alpha}{\vbar}
=\xbar\seqsum\mul{\alpha}{\vbar}\subseteq\ahfline{\xbar}{\vbar}$.%
\indexg{halflines, astral!sequential characterization|)}%
\end{proof}

\indexg{halflines, astral!included in convex set|(}%
In standard convex analysis, it is known that, for $\vv\in\Rn$, if
the halfline $\hfline{\xx}{\vv}$ is included in some
closed (in $\Rn$) and convex set $C\subseteq\Rn$ for even one point
$\xx\in C$, then the same holds for \emph{every} point $\xx\in C$
\idxroc\citep[Theorem~8.3]{ROC}.
We next prove an analogous theorem for astral convex sets, without the closedness requirement.
Moreover, instead of requiring that an entire astral halfline be included in the set,
it suffices to require only that the set include the terminating
endpoint of one such halfline (viewed as a segment, as in
\Cref{thm:ast-halfline:conic}).

\begin{theorem}
\label{cor:gen-one-in-set-implies-all}
  Let $S\subseteq\extspace$ be convex and let $\vbar\in\extspace$.
  Suppose there exists a point $\ybar\in\extspace$
  for which $\limray{\vbar}\plusl\ybar\in S$.
  Then
  $\ahfline{\xbar}{\vbar} \subseteq S$
  for all $\xbar\in S$.
  Thus,
  \begin{equation*}  \label{eq:cor:gen-one-in-set-implies-all:1}
    S \seqsum  (\aray{\vbar}) = S.
  \end{equation*}
\end{theorem}

\begin{proof}
Let $\xbar\in S$.
Then
\[
  \ahfline{\xbar}{\vbar}
  =
  \xbar\seqsum(\aray{\vbar})
  \subseteq
  \lb{\xbar}{\limray\vbar\plusl\ybar}
  \subseteq
  S.
\]
The first equality is by \Cref{thm:ast-halfline:conic}.
The first inclusion is by \Cref{lem:ahfline-containments}.
The final inclusion is by convexity of $S$ since $\xbar$ and
$\limray\vbar\plusl\ybar$ are both in $S$.

{\mathtogether%
Since this holds for all $\xbar\in S$, by
\Cref{pr:seqsum-props}(\ref{pr:seqsum-props:c}),
it follows that
$S \seqsum  (\aray{\vbar}) \subseteq S$.
The reverse inclusion holds because $\zero\in\aray{\vbar}$.}%
\indexg{halflines, astral!included in convex set|)}%
\indexg{halflines, astral|)}%
\end{proof}

\chapter{Separation theorems}
\label{sec:sep-thms}

Separation theorems are at the foundation of convex analysis, proving,
for instance, that there must always exist a hyperplane separating
any closed convex set from any point outside the set.
Separation theorems have many implications and applications, showing,
for example, that every closed convex set is equal to the
intersection of all closed halfspaces that include it.

In this chapter, we prove analogous separation theorems for astral
space.
We show that any two closed convex sets in astral space can
be separated in a strong sense using an astral hyperplane.
We will see that this has several consequences.
For instance, it implies that the outer convex hull of any set
$S\subseteq\extspace$ is equal to the smallest closed convex set that
includes $S$.

These results use astral hyperplanes,
defined by vectors $\uu$ in
$\Rn$, to separate one set in $\extspace$ from another.
In a dual fashion, we also show how standard convex sets in
$\Rn$ can be separated using a kind of generalized hyperplane
that is
defined by an astral point in~$\extspace$.

\section{Separating convex sets}
\label{sec:sep-cvx-sets}

As defined in \Cref{sec:prelim-sep-thms}, two nonempty sets
$X,Y\subseteq\Rn$ are strongly separated by a (standard) hyperplane $J$ if
the two sets are included in the halfspaces on opposite sides of the
hyperplane with all points in $X\cup Y$
a distance at least some $\epsilon>0$ from
the separating hyperplane.
In \Cref{roc:thm11.1}, we saw that this condition is equivalent to
\begin{equation}   \label{eqn:std-str-sep-crit1}
  \sup_{\xx\in X} \xx\cdot\uu
  <
  \inf_{\yy\in Y} \yy\cdot\uu,
\end{equation}
for some $\uu\in\Rn$.
Taking a step further, it can be argued that
this is also equivalent to the condition
\begin{equation}   \label{eqn:std-str-sep-crit2}
  \sup_{\xx\in X,\yy\in Y} (\xx-\yy)\cdot\uu
  <
  0.
\end{equation}

\indexg{separation of astral sets, strong|(}%
\indexg{separation of astral sets, strong!defined|(}%
Extending these definitions,
we can now define an analogous notion of separation in which
two astral sets are strongly separated if they are included in the
astral halfspaces on opposite sides of some astral hyperplane by some
positive margin:

\begin{definition}
Let $X$ and $Y$ be nonempty subsets of $\extspace$.
Let $J=\{\xbar\in\extspace :\: \xbar\cdot\uu=\beta\}$
be an astral hyperplane
defined by some $\uu\in\Rn\wo\{\zero\}$ and $\beta\in\R$.
We say that $J$
\emph{strongly separates} $X$ and $Y$ if
there exists $\epsilon>0$ such that
$\xbar\cdot\uu < \beta-\epsilon$ for all $\xbar\in X$,
and
$\ybar\cdot\uu > \beta+\epsilon$ for all $\ybar\in Y$.
\end{definition}

More simply, analogous to \eqref{eqn:std-str-sep-crit1},
$X$ and $Y$ are strongly separated by some hyperplane
if and only if
\begin{equation}  \label{eqn:ast-strong-sep-defn}
  \sup_{\xbar\in X} \xbar\cdot \uu
  <
  \inf_{\ybar\in Y} \ybar\cdot \uu
\end{equation}
for some $\uu\in\Rn$
(where the possibility that $\uu=\zero$ is ruled out since $X$ and $Y$
are nonempty).
When \eqref{eqn:ast-strong-sep-defn} holds,
we say simply that $X$ and $Y$ are strongly separated by the vector
$\uu\in\Rn$, or that $\uu$ strongly separates them
(rather than referring to an astral hyperplane, such as the one given
above).
We say that $X$ and $Y$ are strongly separated if there exists
$\uu\in\Rn$ that strongly separates them.%
\indexg{separation of astral sets, strong!defined|)}

The condition for (standard) strong separation given in
\eqref{eqn:std-str-sep-crit2} means that $\zz\cdot\uu$ is bounded
from above by a negative constant for all $\zz\in X-Y$.
\indexg{separation of astral sets, strong!characterized|(}%
We next give an analogous criterion for strong separation in astral
space in which, for sets $X,Y\subseteq\extspace$, the set difference
$X-Y$ is replaced by $\Xbar\seqsum(-\Ybar)$,
using sequential sum after taking closures.

\begin{theorem}   \label{thm:ast-str-sep-seqsum-crit}
  Let $X,Y\subseteq\extspace$ be nonempty, and let $\uu\in\Rn$.
  Then
  \begin{equation}   \label{eq:thm:ast-str-sep-seqsum-crit:1}
    \sup_{\xbar\in X} \xbar\cdot \uu
    <
    \inf_{\ybar\in Y} \ybar\cdot \uu
  \end{equation}
  if and only if
  \begin{equation}   \label{eq:thm:ast-str-sep-seqsum-crit:2}
    \sup_{\zbar\in \Xbar\seqsum(-\Ybar)} \zbar\cdot \uu
    <
    0.
  \end{equation}
  Consequently,
  $X$ and $Y$ are strongly separated if and only if
  \eqref{eq:thm:ast-str-sep-seqsum-crit:2}
  holds for some $\uu\in\Rn$.
\end{theorem}

\begin{proof}
  ~

\begin{proof-parts}
\pfpart{``Only if'' ($\Rightarrow$):}
Suppose \eqref{eq:thm:ast-str-sep-seqsum-crit:1} holds,
implying $\uu\neq\zero$
(since $X$ and $Y$ are nonempty).
Then there exist $\beta,\gamma\in\R$ such that
\[
  \sup_{\xbar\in X} \xbar\cdot \uu
  \leq
  \beta
  <
  \gamma
  \leq
  \inf_{\ybar\in Y} \ybar\cdot \uu.
\]
Let $\zbar\in \Xbar\seqsum(-\Ybar)$, implying that
$\zbar\in\xbar\seqsum(-\ybar)$
for some $\xbar\in \Xbar$ and $\ybar\in \Ybar$.
Then $\xbar\cdot\uu\leq\beta<+\infty$ since
$X$, and therefore also $\Xbar$,
is included in the closed astral halfspace
$\chsua$ (as defined in Eq.~\ref{eq:chsua-defn}).
Similarly,
$\ybar\cdot\uu\geq\gamma$, so
$-\ybar\cdot\uu\leq-\gamma<+\infty$.
Thus, $\xbar\cdot\uu$ and $-\ybar\cdot\uu$ are summable, implying
\[
  \zbar\cdot\uu
  =
  \xbar\cdot\uu - \ybar\cdot\uu
  \leq
  \beta - \gamma,
\]
where the equality is by
\Cref{cor:seqsum-conseqs}(\ref{cor:seqsum-conseqs:b}).
Since this holds for all $\zbar\in \Xbar\seqsum(-\Ybar)$, and
since $\beta-\gamma$ is a negative constant,
this proves the claim.

\pfpart{``If'' ($\Leftarrow$):}
Suppose \eqref{eq:thm:ast-str-sep-seqsum-crit:2} holds.

Let $\xbar'\in\Xbar$ maximize $\xbar\cdot\uu$ over all $\xbar\in\Xbar$.
Such a point $\xbar'$ must exist by \Cref{pr:cont-compact-attains-max},
because $\Xbar$ is closed and therefore compact
(\Cref{prop:compact}\ref{prop:compact:closed-subset}),
and the map $\xbar\mapsto\xbar\cdot\uu$ is continuous
(\Cref{thm:i:1}\ref{thm:i:1c}).
Likewise, let $\ybar'\in\Ybar$ minimize $\ybar\cdot\uu$ over all
$\ybar\in\Ybar$.

Let $\zbar'=\xbar'\plusl(-\ybar')$, and
let $\zbar''=(-\ybar')\plusl\xbar'$.
Then $\zbar'$ and $\zbar''$ are both in
$\xbar'\seqsum(-\ybar')\subseteq  \Xbar\seqsum(-\Ybar)$
by \Cref{cor:seqsum-conseqs}(\ref{cor:seqsum-conseqs:a}).
Consequently,
\begin{equation}   \label{eq:thm:ast-str-sep-seqsum-crit:3}
  \xbar'\cdot\uu\plusl(-\ybar'\cdot\uu)=\zbar'\cdot\uu<0,
\end{equation}
implying specifically that $\xbar'\cdot\uu<+\infty$.
Likewise,
$(-\ybar'\cdot\uu)\plusl\xbar'\cdot\uu=\zbar''\cdot\uu<0$,
implying that $-\ybar'\cdot\uu<+\infty$.
Thus, $\xbar'\cdot\uu$ and $-\ybar'\cdot\uu$ are summable,
so
\begin{equation}   \label{eq:thm:ast-str-sep-seqsum-crit:4}
  \xbar'\cdot\uu\plusl(-\ybar'\cdot\uu)=\xbar'\cdot\uu - \ybar'\cdot\uu.
\end{equation}

Combining, we now have that
\[
  \sup_{\xbar\in X} \xbar\cdot\uu
  \leq
  \xbar'\cdot\uu
  <
  \ybar'\cdot\uu
  \leq
  \inf_{\ybar\in Y} \ybar\cdot\uu.
\]
The first and last inequalities are by our choice of $\xbar'$ and
$\ybar'$.
The middle (strict) inequality follows from
Eqs.~(\ref{eq:thm:ast-str-sep-seqsum-crit:3})
and~(\ref{eq:thm:ast-str-sep-seqsum-crit:4}).

\pfpart{Strong separation equivalence:}
If $X$ and $Y$ are both nonempty, then $X$ and $Y$ are strongly
separated if and only if
\eqref{eq:thm:ast-str-sep-seqsum-crit:1} holds for some $\uu\in\Rn$,
which, as shown above, holds if and only if
\eqref{eq:thm:ast-str-sep-seqsum-crit:2} holds for some $\uu\in\Rn$.
\qedhere
\end{proof-parts}
\end{proof}

The equivalence given in
\Cref{thm:ast-str-sep-seqsum-crit} is false in general if,
in \eqref{eq:thm:ast-str-sep-seqsum-crit:2}, we use the sequential sum
$X\seqsum(-Y)$, without taking closures, rather than
$\Xbar\seqsum(-\Ybar)$.
For instance, in $\Rext$, suppose $X=\{-\infty\}$, $Y=\R$, and $u=1$.
Then $X\seqsum(-Y)=\{-\infty\}$, so
$\sup_{\barz\in X\seqsum(-Y)} \barz u = -\infty<0$,
but \eqref{eq:thm:ast-str-sep-seqsum-crit:1} does not hold in this
case (since both sides of the inequality are equal to $-\infty$).%
\indexg{separation of astral sets, strong!characterized|)}

\indexg{separation of astral sets, strong!main theorem|(}%
We next give the central theorem of this section,
showing
every pair of disjoint, closed convex sets in $\extspace$ is strongly
separated.

\begin{theorem}  \label{thm:sep-cvx-sets}
  Let $X,Y\subseteq\extspace$ be nonempty, closed (in $\extspace$),
  convex, and disjoint.
  Then $X$ and $Y$ are strongly separated.
\end{theorem}

The main steps in proving this theorem are in showing that it holds in
particular special cases,
and then showing how these special cases imply the general case.
\indexg{separation of astral sets, strong!astral cone@from astral cone|(}%
To begin, we prove the next lemma showing how a closed convex astral
cone $K\subseteq\extspace$ can be separated from a single, finite
point $\zz$ in $\Rn\setminus K$,
and more specifically, that there must always exist a homogeneous
closed astral halfspace that includes $K$ but not $\zz$.
The main idea of the lemma's proof is to first separate $\zz$ from the
(standard) convex cone $K'=K\cap\Rn$, and then show that that
separation also leads to a separation of $\zz$ from the entire astral
cone $K$.
For this to work, care must be taken in how this separation is done;
in the proof, we use sharp discrimination of cones introduced in
\Cref{sec:prelim-sep-thms}.

\begin{lemma}  \label{lem:acone-sep-from-fin-pt}
  Let $K\subseteq\extspace$ be a closed convex astral cone, and let
  $\zz\in\Rn\wo K$.
  Then there exists $\uu\in\Rn\wo\{\zero\}$ such that
  \[
     \sup_{\xbar\in K} \xbar\cdot\uu
     \leq
     0
     <
     \zz\cdot\uu.
  \]
\end{lemma}

\begin{proof}
Note that $\zz\neq\zero$ since $\zero\in K$
(by \Cref{pr:ast-cone-is-naive})
and $\zz\not\in K$.

Let $K'=K\cap\Rn$.
Then $K'$ is a closed convex cone
(being the intersection of two convex sets, and by
\Cref{pr:ast-cone-is-naive}
and \Cref{prop:subspace}\ref{i:subspace:closed}).

Let $J=\cone \{\zz\}$, which is a closed convex cone.
Then $J\cap K'=\{\zero\}$, since
otherwise we would have $\lambda \zz \in K'\subseteq K$ for some
$\lambda\in\Rstrictpos$, implying $\zz\in K$
(by \Cref{pr:ast-cone-is-naive}),
a contradiction.

Thus, we can apply \Cref{pr:cones-sharp-sep} to $K'$ and $J$,
yielding their sharp discrimination.
Letting $L=K'\cap -K'$, this means
there exists $\uu\in\Rn$ such that
$\xx\cdot\uu= 0$ for $\xx\in L$;
$\xx\cdot\uu< 0$ for $\xx\in K'\wo L$;
and $\zz\cdot\uu > 0$
(noting that $J\cap-J=\{\zero\}$).
In particular, this implies $\uu\neq\zero$.

To prove the lemma, we will show that $K\subseteq\chsuz$
(where $\chsuz$ is as defined in Eq.~\ref{eqn:homo-halfspace}).
To prove this, we will see that it is sufficient to
show that all of $K$'s icons are in $\chsuz$, since $K$ is the convex
hull of its icons
(\Cref{thm:ast-cvx-cone-equiv}\ref{thm:ast-cvx-cone-equiv:a}\ref{thm:ast-cvx-cone-equiv:c}).
As such,
let $\ebar\in K\cap\corezn$ be an icon in $K$;
we aim to show that $\ebar\cdot\uu\leq 0$.
We can write $\ebar=[\vv_1,\ldots,\vv_k]\omm$ for some
$\vv_1,\ldots,\vv_k\in\Rn$.
If $\vv_i\in L$ for all $i=1,\ldots,k$, then
$\vv_i\cdot\uu=0$ for all $i$, implying in this case that
$\ebar\cdot\uu=0$
(by \Cref{pr:vtransu-zero}).

Otherwise, let $j\in\{1,\ldots,k\}$ be the smallest index for which
$\vv_j\not\in L$
(implying $\vv_i\in L$ for $i=1,\ldots,j-1$).

\begin{claimpx}
  $\vv_j\in K'$.
\end{claimpx}

\begin{proofx}
Let $\ybar=[\vv_1,\ldots,\vv_j]\omm$.
Then
\[
  \ybar
  \in
  \lb{\zero}{\ebar}
  =
  \aray{\ebar}
  \subseteq
  K.
\]
The first inclusion is by
\Cref{thm:lb-with-zero}
(or \Cref{cor:d-in-lb-0-dplusx}),
the equality is by
\Cref{thm:oconichull-equals-ocvxhull},
and the second inclusion is because $K$ is an astral cone that
includes $\ebar$.

If $j=1$, then $\limray{\vv_1}=\ybar\in K$,
implying $\vv_1\in K$ by
\Cref{pr:ast-cone-is-naive}.

So suppose henceforth that $j>1$.
Then
for $i=1,\ldots,j-1$, we have $-\vv_i\in K'$ since $\vv_i\in L$,
implying $\limray{(-\vv_i)}=-\limray{\vv_i}$ is in $K$
(by \Cref{pr:ast-cone-is-naive}).
To show $\vv_j\in K$, we will argue that it is in the outer convex
hull of $\ybar$ and $-\limray{\vv_1},\ldots,-\limray{\vv_{j-1}}$, all
of which are in $K$,
using sequences so that \Cref{thm:e:7} can be
applied.
As such, let
\[
   \yy_t
   =
   \sum_{i=1}^j t^{j-i+1} \vv_i
   \\
   \;\;\;\;\mbox{ and }\;\;\;\;
   \xx_{it}
   =
   -\frac{t^{j-i}}{\lambda_t} \vv_i,
\]
for $i=1,\ldots,j-1$, where
\[
   \gamma_t
   =
   \frac{1}{t}
   \;\;\;\;\mbox{ and }\;\;\;\;
   \lambda_t
   =
   \frac{1}{j-1}\Parens{1-\gamma_t}.
\]
Finally, let
\begin{equation}
\label{eq:acone-sep-from-fin-pt:1}
   \zz_t
   =
   \gamma_t \yy_t
   +
   \sum_{i=1}^{j-1} \lambda_t \xx_{it}.
\end{equation}
Note that $\gamma_t + (j-1)\lambda_t = 1$.
Then
$\yy_t\rightarrow\ybar$,
and
$\xx_{it}\rightarrow -\limray{\vv_i}$ for $i=1,\ldots,j-1$
(by \Cref{thm:i:seq-rep}).
Further, by evaluating \eqref{eq:acone-sep-from-fin-pt:1}, we obtain $\zz_t=\vv_j$ for all $t$, and so $\zz_t$ trivially
converges to $\vv_j$.
It follows that
\[
  \vv_j
  \in
  \ohull\braces{\ybar, -\limray{\vv_1},\ldots, -\limray{\vv_{j-1}}}
  \subseteq
  K,
\]
where the first inclusion is by \Cref{thm:e:7},
applied to the sequences above,
and the second by \Cref{thm:e:2} since
each of the points in braces are in $K$.
Thus, $\vv_j\in K'$, as claimed.
\end{proofx}

Since $\uu$ sharply discriminates $K'$ and $J$, it follows that
$\vv_j\cdot\uu<0$
since $\vv_j$ is in $K'$ but not in $L$ by our choice of $j$.
Further, $\vv_i\cdot\uu=0$ for $i=1,\ldots,j-1$ since each
$\vv_i\in L$.
Consequently, $\ebar\cdot\uu=-\infty<0$
(by the Case Decomposition Lemma~\ref{lemma:case}),
so $\ebar\in\chsuz$.

Thus, $E\subseteq\chsuz$, where $E=K\cap\corezn$.
Therefore, $K=\conv{E}\subseteq\chsuz$,
with equality by
\Cref{thm:ast-cvx-cone-equiv}(\ref{thm:ast-cvx-cone-equiv:a},\ref{thm:ast-cvx-cone-equiv:c}),
and inclusion by
\Cref{pr:conhull-prop}(\ref{pr:conhull-prop:aa})
since $\chsuz$ is convex
(\Cref{pr:e1}\ref{pr:e1:c}).
Since $\zz\not\in\chsuz$, this completes the proof.%
\indexg{separation of astral sets, strong!astral cone@from astral cone|)}%
\end{proof}

Next, we build on the preceding lemma to prove a separation result between a closed
convex set and a disjoint point.
We adapt a technique from standard convex analysis, which is sometimes
used to extend results for convex cones to general convex sets.
The idea is to take a given convex set $X\subseteq\Rn$, use it to
create a convex cone $K\subseteq\Rnp$ as the conic hull of pairs
$\rpair{\xx}{1}$ where $\xx\in X$, and then apply the result for cones to $K$
to obtain something useful about $X$.\looseness=-1

In the proof of the next lemma,
we use the same idea, adapted to astral space, and lift the separation result for a closed convex astral cone and a disjoint point (\Cref{lem:acone-sep-from-fin-pt}) to obtain a separation result for a closed convex set and a disjoint point (specifically $\zero$).

\begin{lemma}   \label{lem:sep-cvx-set-from-origin}
  Let $X\subseteq\extspace$ be nonempty, convex and closed (in $\extspace$),
  and assume $\zero\not\in X$.
  Then there exists $\uu\in\Rn\wo\{\zero\}$ such that
  \[
    \sup_{\xbar\in X} \xbar\cdot\uu
    <
    0.
  \]
\end{lemma}

\begin{proof}
As just discussed,
the main idea of the proof is to use $X$ to construct a closed convex
astral cone in $\extspacnp$, and then use
\Cref{lem:acone-sep-from-fin-pt} to find a separating hyperplane.

As in \Cref{sec:work-with-epis},
for $\xbar\in\extspace$ and $y\in\R$,
the notation $\rpair{\xbar}{y}$ denotes, depending on context,
either a pair in $\extspace\times\R$ or a point in
$\extspacnp$, namely, the image of that pair under the natural
embedding $\homf$ given in \Cref{thm:homf}.
We also make use of the matrix $\PPx=[\Idnn,\zerov{n}]$
from \eqref{eq:PPx:PPy}, which satisfies $\PPx\rpair{\xbar}{y}=\xbar$
for all $\xbar\in\extspace$ and $y\in\R$.

Let $S\subseteq\extspacnp$ be the set
\[
  S
  =
  \Braces{ \rpair{\xbar}{1} :\: \xbar\in X }.
\]
Equivalently, $S=F(X)$ where $F:\extspace\rightarrow\extspacnp$ is
the affine map
$\xbar\mapsto \rpair{\zero}{1}\plusl\trans{\PPx}\xbar$
for $\xbar\in\extspace$
(\Cref{pr:xy-pairs-props}\ref{pr:xy-pairs-props:a}).
Therefore, $S$ is convex (by \Cref{cor:thm:e:9})
and closed (by \Cref{cor:aff-img-closed-is-closed}\ref{cor:aff-img-closed-is-closed:a})
since $X$ is.

The next step is to use $S$
to construct a closed convex astral cone $K\subseteq\extspacnp$,
to which \Cref{lem:acone-sep-from-fin-pt} can be applied.
It might seem natural to choose $K$ simply to be $\acone{S}$,
the astral conic hull of $S$.
The problem, however, is that $\acone{S}$ might not be closed.

Instead, we use a variant on $\acone{S}$ that ensures that $K$ has all
the required properties.
Specifically, we define
\begin{equation}   \label{eq:lem:sep-cvx-set-from-origin:2}
  K
  =
  \conv\bigParens{ \{\zero\} \cup E }
  \;\text{ where }\;
  E
  =
  \clbar{\lmset{S}}.
\end{equation}
Thus, comparing to the expression for $\acone{S}$ given in
\Cref{thm:acone-char}(\ref{thm:acone-char:a}),
we have simply replaced $\lmset{S}$ with its closure.

The set $E$ is clearly closed, so $K$ is as well,
by \Cref{thm:cnv-hull-closed-is-closed}.
Also, since $\lmset{S}$ is included in the closed set $\coreznp$
(\Cref{pr:i:8}\ref{pr:i:8-infprod}\ref{pr:i:8e}),
$E$ is also included in $\coreznp$, and so consists only of icons.
Therefore, $K$ is a convex astral cone by
\Crefequiv{thm:ast-cvx-cone-equiv}{thm:ast-cvx-cone-equiv:b}{thm:ast-cvx-cone-equiv:a}.

Although $E$ may differ from $\lmset{S}$, every element of $E$ has
a prefix that is an icon in $\lmset{S}$, as shown in the next claim:

\begin{claimpx}   \label{cl:lem:sep-cvx-set-from-origin:1}
  Let $\ebar\in E$.
  Then $\ebar=\limray{\vbar}\plusl\ybar$
  for some $\vbar\in S$ and $\ybar\in\extspacnp$.
  That is,
  $E\subseteq\lmset{S}\plusl\extspacnp$.
\end{claimpx}

\begin{proofx}
Since $\ebar\in E$, there exists a sequence $\seq{\vbar_t}$ in $S$
such that $\limray{\vbar_t}\rightarrow\ebar$.
By sequential compactness, this sequence must have a convergent
subsequence; by discarding all other elements, we can assume the
entire sequence $\seq{\vbar_t}$ converges to some point $\vbar$,
which must also be in $S$ since $S$ is closed.
Therefore, $\ebar\in\limray{\vbar}\plusl\extspacnp$
by \Cref{thm:gen-dom-dir-converg}.
\end{proofx}

Since $\zero\not\in X$, we also have that
$\rpair{\zero}{1}\not\in S$.
The next several steps of the proof are directed toward showing
that also $\rpair{\zero}{1}\not\in K$, which
will allow us to apply
\Cref{lem:acone-sep-from-fin-pt}
to separate it from $K$.

The astral cone $K$ is largely constructed from astral rays
$\aray{\ebar}$ over $\ebar\in E$ (in a way that will be made precise
in the proof of 
\Cref{cl:lem:sep-cvx-set-from-origin:3}).
The next claim focuses on such astral rays.
As shown in
\Cref{cl:lem:sep-cvx-set-from-origin:1}, the point $\ebar$
must have the form $\limray{\zbar}\plusl\ybar$ for some $\zbar\in S$
and $\ybar\in\extspacnp$.
By \Cref{thm:aray}, every point $\wbar\in\aray{\ebar}$ must belong to
a sequential $\alpha$-multiple of $\ebar$, for some $\alpha\in\Rextpos$.
The next claim shows that if the ``last component''
of $\wbar$ is equal to some finite $\lambda\in\Rpos$, then
actually $\wbar$ must more specifically belong to
$\mul{\lambda}{\zbar}$, the sequential $\lambda$-multiple of $\zbar$
(rather than $\ebar$).
This will provide a key step in the proof of 
\Cref{cl:lem:sep-cvx-set-from-origin:3}.

In the statement of the claim,
we can think of $\xbar$ as a point in $X$ (so that $\zbar$ is in $S$),
and $\limray{\zbar}\plusl\ybar$ as a point in $E$.

\begin{claimpx}   \label{cl:lem:sep-cvx-set-from-origin:2}
  Let $\xbar\in\extspace$,
  $\zbar=\rpair{\xbar}{1}$,
  and
  let $\wbar\in\aray(\limray{\zbar}\plusl\ybar)$
  for some $\ybar\in\extspacnp$.
  Suppose $\wbar\cdot\rpair{\zero}{1}=\lambda$ where
  $\lambda\in\Rpos$.
  Then $\wbar\in\mul{\lambda}{\zbar}$.
\end{claimpx}

\begin{proofx}
We have
\begin{align}
   \aray(\limray{\zbar}\plusl\ybar)
   =
   \lb{\zero}{\limray{\zbar}\plusl\limray{\ybar}}
   &=
   \lb{\zero}{\limray{\zbar}}
   \cup
   \lb{\limray{\zbar}}{\limray{\zbar}\plusl\limray{\ybar}}
   \nonumber
   \\
   &=
   (\aray{\zbar})
   \cup
   \lb{\limray{\zbar}}{\limray{\zbar}\plusl\limray{\ybar}}.
   \label{eq:lem:sep-cvx-set-from-origin:1}
\end{align}
The first and third equalities are by
\Cref{thm:oconichull-equals-ocvxhull}
(and \Cref{pr:i:8}\ref{pr:i:8d}).
The second equality is by
\Cref{thm:decomp-seg-at-inc-pt}
since $\limray{\zbar}\in\lb{\zero}{\limray{\zbar}\plusl\limray{\ybar}}$
(by \Cref{cor:d-in-lb-0-dplusx}).

We claim that
$\wbar\not\in\lb{\limray{\zbar}}{\limray{\zbar}\plusl\limray{\ybar}}$.
This is because $\limray{\zbar}\cdot\rpair{\zero}{1}=+\infty$
(since $\zbar\cdot\rpair{\zero}{1}=1$ by
\Cref{pr:xy-pairs-props}\ref{pr:xy-pairs-props:b}), implying
$(\limray{\zbar}\plusl\limray{\ybar})\cdot\rpair{\zero}{1}=+\infty$
as well.
Thus, if $\wbar$ were on the segment joining these points,
we would have $\wbar\cdot\rpair{\zero}{1}=+\infty$
by
\Cref{pr:seg-simplify}(\ref{pr:seg-simplify:a},\ref{pr:seg-simplify:c}),
a contradiction.

Thus, in light of \eqref{eq:lem:sep-cvx-set-from-origin:1},
it follows that $\wbar\in\aray{\zbar}$.
Therefore,
by \Cref{thm:aray}, $\wbar\in\mul{\alpha}{\zbar}$ for some $\alpha\in\Rextpos$.
Let $\PPy=[\trans{\zero_n},1]$.
Then
\[
  \lambda
  =
  \PPy\wbar
  \in
  \PPy \Bracks{\mul{\alpha}{\zbar}}
  =
  \mul{\alpha}{(\PPy\zbar)}
  =
  \mul{\alpha}{1}
  =
  \set{\alpha},
\]
where the first and third equalities are by
\Cref{pr:xy-pairs-props}(\ref{pr:xy-pairs-props:b-new},\ref{pr:xy-pairs-props:c}),
the second is by \Cref{pr:seq:mul}(\ref{i:seq:mul:A}),
and the last is by \Cref{thm:mul-char}.
Thus, $\alpha=\lambda$ and $\wbar\in\mul{\lambda}{\zbar}$.
\end{proofx}

We can now show that the only points in $K$ with
``last coordinate'' equal to $1$ are the points in $S$:

\begin{claimpx}   \label{cl:lem:sep-cvx-set-from-origin:3}
  Let $\xbar\in\extspace$, and
  suppose $\zbar=\rpair{\xbar}{1}\in K$.
  Then $\xbar\in X$.
\end{claimpx}

\begin{proofx}
The set $K$, as defined in
\eqref{eq:lem:sep-cvx-set-from-origin:2},
is equal to $\acone{E}$ by
\Cref{thm:acone-char}(\ref{thm:acone-char:a})
(and since $\lmset{E}=E$).
Therefore, since $\zbar\in K$,
\Cref{thm:ast-conic-hull-union-fin}
implies that there must exist
icons $\ebar_1,\ldots,\ebar_m\in E$
for which
\[
  \zbar
  \in
  \oconich\{\ebar_1,\ldots,\ebar_m\}
  =
  (\aray{\ebar_1})\seqsum\dotsb\seqsum(\aray{\ebar_m}),
\]
where the equality is by \Cref{cor:conic:comb}.
Thus,
by \Cref{prop:seqsum-multi}(\ref{prop:seqsum-multi:decomp}),
\begin{equation}   \label{eq:lem:sep-cvx-set-from-origin:3}
  \zbar
  \in
  \wbar_1\seqsum\dotsb\seqsum\wbar_m
\end{equation}
for some $\wbar_1,\ldots,\wbar_m$ with $\wbar_i\in\aray{\ebar_i}$
for $i=1,\ldots,m$.

Let
$H = \{ \ybar\in\extspacnp :\: \ybar\cdot\rpair{\zero}{1} \geq 0 \}$.
This is a homogeneous closed astral halfspace, and therefore also a convex
astral cone
(\Cref{pr:astral-cone-props}\ref{pr:astral-cone-props:e}).
Note that $S\subseteq H$ (since $\ybar\cdot\rpair{\zero}{1}=1$ for
$\ybar\in S$), so $\lmset{S}\subseteq H$ by
\Cref{pr:ast-cone-is-naive} since $H$ is an astral cone.
Therefore, $E\subseteq H$ since $H$ is closed.
Consequently, for $i=1,\ldots,m$,
we have
$\ebar_i\in H$, implying that $\aray{\ebar_i}\subseteq H$
(by definition of astral cones),
and thus $\wbar_i\in H$.

Let $\lambda_i=\wbar_i\cdot\rpair{\zero}{1}$
for $i=1,\ldots,m$.
Then, as just argued,
$\lambda_1,\ldots,\lambda_m$
are all nonnegative and therefore summable.
Combined with
\eqref{eq:lem:sep-cvx-set-from-origin:3},
\Cref{thm:seqsum-equiv-mod}(\ref{thm:seqsum-equiv-mod:a},\ref{thm:seqsum-equiv-mod:c})
therefore implies that
\begin{equation}   \label{eq:lem:sep-cvx-set-from-origin:4}
  1
  =
  \zbar\cdot\rpair{\zero}{1}
  =
  \sum_{i=1}^m \wbar_i\cdot\rpair{\zero}{1}
  =
  \sum_{i=1}^m \lambda_i.
\end{equation}

Thus, for $i=1,\ldots,m$, each $\lambda_i\in [0,1]$.
By \Cref{cl:lem:sep-cvx-set-from-origin:1},
there exists $\xbar_i\in X$ such that
$\ebar_i\in\limray{\zbar_i}\plusl\extspacnp$
where $\zbar_i=\rpair{\xbar_i}{1}$.
Therefore, by \Cref{cl:lem:sep-cvx-set-from-origin:2},
$\wbar_i\in\mul{\lambda_i}{\zbar_i}$.

We then have that
\begin{equation}
\label{eq:lem:sep-cvx-set-from-origin:6}
  \zbar
  \in
  \wbar_1 \seqsum \dotsb \seqsum \wbar_m
  \subseteq
  \mul{\lambda_1}{\zbar_1}
  \seqsum\dotsb\seqsum
  \mul{\lambda_m}{\zbar_m}
  \subseteq
  \ohull\{\zbar_1,\ldots,\zbar_m\},
\end{equation}
where the second inclusion is by
\Cref{prop:seqsum-multi}(\ref{prop:seqsum-multi:decomp}),
and the third is by
\Cref{cor:conv-as-seqsum-midrays}
(using Eq.~\ref{eq:lem:sep-cvx-set-from-origin:4}).
It then follows that
\[
  \xbar
  =
  \PPx \zbar
  \in
  \ohull\{\PPx \zbar_1, \ldots, \PPx \zbar_m\}
  =
  \ohull\{\xbar_1, \ldots, \xbar_m\}
  \subseteq
  X.
\]
The two equalities are by
\Cref{pr:xy-pairs-props}(\ref{pr:xy-pairs-props:c}).
The first inclusion is by
\Cref{thm:e:9} applied to
\eqref{eq:lem:sep-cvx-set-from-origin:6}.
The last inclusion is by \Cref{thm:e:2}
since $X$ is convex.
\end{proofx}

In particular, since $\zero\not\in X$,
\Cref{cl:lem:sep-cvx-set-from-origin:3}
implies that $\rpair{\zero}{1}$ is not in $K$.
Therefore, we can finally apply
\Cref{lem:acone-sep-from-fin-pt} to $K$ and $\rpair{\zero}{1}$,
yielding that there exists $\uu\in\Rn$ and $v\in\R$ such that
$v=\rpair{\zero}{1}\cdot\rpair{\uu}{v} > 0$
and
$\zbar\cdot\rpair{\uu}{v}\leq 0$ for all $\zbar\in K$.
In particular, for all $\xbar\in X$, since $\rpair{\xbar}{1}\in K$,
this means that
\[ \xbar\cdot\uu + v=\rpair{\xbar}{1}\cdot\rpair{\uu}{v}\leq 0 \]
(with the equality from
\Cref{pr:xy-pairs-props}\ref{pr:xy-pairs-props:b}).
Thus,
\begin{equation}   \label{eq:lem:sep-cvx-set-from-origin:5}
  \sup_{\xbar\in X} \xbar\cdot\uu
  \leq
  -v
  <
  0,
\end{equation}
completing the proof
(after noting,
by Eq.~\ref{eq:lem:sep-cvx-set-from-origin:5},
that $\uu\neq\zero$ since $X$ is nonempty).
\end{proof}

We can now prove the general case:

\begin{proof}[Proof of \Cref{thm:sep-cvx-sets}]
Let $Z=X\seqsum(-Y)$.
Since the map $\xbar\mapsto-\xbar$ is linear,
the set $-Y$ is convex (by \Cref{cor:thm:e:9})
and closed
(by \Cref{cor:aff-img-closed-is-closed}\ref{cor:aff-img-closed-is-closed:a}).
Therefore, $Z$ is nonempty, convex and closed
(by \Cref{prop:seqsum-multi}\ref{prop:seqsum-multi:nonempty}\ref{prop:seqsum-multi:closed}\ref{prop:seqsum-multi:convex}).

We further claim $\zero\not\in Z$.
Otherwise, if $\zero$ were in $Z$,
then we would have $\zero\in\xbar\seqsum(-\ybar)$
for some $\xbar\in X$ and $\ybar\in Y$, implying
$\xbar\in\zero\seqsum\ybar=\{\ybar\}$
by \Cref{pr:swap-seq-sum}, and thus that
$\xbar=\ybar$, a contradiction since $X$ and $Y$ are disjoint.

Thus, we can apply \Cref{lem:sep-cvx-set-from-origin}
to $Z$, implying that there exists $\uu\in\Rn\wo\{\zero\}$
for which
\[
  \sup_{\zbar\in X\seqsum(-Y)} \zbar\cdot\uu
  <
  0.
\]
By \Cref{thm:ast-str-sep-seqsum-crit}
(and since $X$ and $Y$ are closed), this proves the theorem.%
\indexg{separation of astral sets, strong!main theorem|)}%
\end{proof}

\Cref{thm:sep-cvx-sets} has several direct consequences, as
we summarize in the next corollary.
\indexg{outer convex hull!smallest closed convex superset@as smallest closed convex superset|(}%
Note especially that the outer convex hull of any set
$S\subseteq\extspace$ can now be shown to be the smallest closed convex
set that includes $S$.

\begin{corollary}   \label{cor:sep-cvx-sets-conseqs}
\indexg{convex sets, astral!intersection of halfspaces@as intersection of halfspaces|(}%
\indexg{halfspaces, astral!convex sets as intersection of|(}%
  Let $S\subseteq\extspace$, and
  let $X, Y\subseteq\extspace$ be nonempty.
  \begin{letter-compact}
  \item    \label{cor:sep-cvx-sets-conseqs:d}
    The following are equivalent:
    \begin{roman-compact}
    \item    \label{cor:sep-cvx-sets-conseqs:d:1}
      $S$ is closed and convex.
    \item    \label{cor:sep-cvx-sets-conseqs:d:2}
      $S$ is equal to the intersection of some collection of closed
      astral halfspaces.
    \item    \label{cor:sep-cvx-sets-conseqs:d:3}
      $S$ is equal to the intersection of all closed
      astral halfspaces that include $S$;
      that is, $S=\ohull{S}$.
    \end{roman-compact}
  \item    \label{cor:sep-cvx-sets-conseqs:a}
    The outer convex hull of $S$, $\ohull{S}$, is the smallest closed
    convex set that includes~$S$.
    (That is, $\ohull{S}$ is equal to the intersection of all closed
    convex sets in $\extspace$ that include $S$.)
  \item    \label{cor:sep-cvx-sets-conseqs:b}
    $\ohull{S}=\conv{\Sbar}$.
  \item    \label{cor:sep-cvx-sets-conseqs:c}
    $X$ and $Y$ are strongly separated
    if and only if
    $(\ohull{X})\cap(\ohull{Y})=\emptyset$.
  \end{letter-compact}
\end{corollary}

\begin{proof}
~

\begin{proof-parts}
\pfpart{Part~(\ref{cor:sep-cvx-sets-conseqs:d}):}
That (\ref{cor:sep-cvx-sets-conseqs:d:3}) implies
(\ref{cor:sep-cvx-sets-conseqs:d:2}) is immediate.

Every closed astral halfspace is closed and convex,
so the arbitrary intersection of such sets is also closed and convex
(\Cref{pr:e1}\ref{pr:e1:b}\ref{pr:e1:c}),
proving that (\ref{cor:sep-cvx-sets-conseqs:d:2}) implies
(\ref{cor:sep-cvx-sets-conseqs:d:1}).

To prove
(\ref{cor:sep-cvx-sets-conseqs:d:1}) implies
(\ref{cor:sep-cvx-sets-conseqs:d:3}),
suppose $S$ is closed and convex.
We assume $S$ is nonempty since otherwise,
$\ohull{S}=\emptyset=S$
(since the intersection of all closed astral halfspaces is empty).
Clearly, $S\subseteq\ohull{S}$.
For the reverse inclusion, suppose $\zbar\in\extspace\setminus S$.
Then \Cref{thm:sep-cvx-sets},
applied to $S$ and $\{\zbar\}$, yields that there exist
$\uu\in\Rn\wo\{\zero\}$
and $\beta\in\R$ such that
$\xbar\cdot\uu<\beta$ for $\xbar\in S$, but $\zbar\cdot\uu>\beta$.
Thus,
$S\subseteq\chsua$ but $\zbar\not\in\chsua$,
implying that $\zbar\not\in\ohull{S}$.
Therefore, $\ohull{S}\subseteq S$.%
\indexg{halfspaces, astral!convex sets as intersection of|)}%
\indexg{convex sets, astral!intersection of halfspaces@as intersection of halfspaces|)}

\pfpart{Part~(\ref{cor:sep-cvx-sets-conseqs:a}):}
Let $U$ be equal to the intersection of all closed
convex sets in $\extspace$ that include $S$.
Then $U\subseteq\ohull{S}$
since $\ohull{S}$ is such a set.
On the other hand,
\[
  \ohull{S}
  \subseteq
  \ohull{U}
  =
  U
\]
where the inclusion is because $S\subseteq U$,
and the equality is by
part~(\ref{cor:sep-cvx-sets-conseqs:d}).%
\indexg{outer convex hull!smallest closed convex superset@as smallest closed convex superset|)}

\pfpart{Part~(\ref{cor:sep-cvx-sets-conseqs:b}):}
The set $\conv{\Sbar}$ includes $S$ and
is closed and convex by
\Cref{thm:cnv-hull-closed-is-closed}.
Therefore,
$\ohull{S}\subseteq\conv{\Sbar}$
by
part~(\ref{cor:sep-cvx-sets-conseqs:a}).
On the other hand, $\ohull{S}$ is closed and includes $S$,
so $\Sbar\subseteq\ohull{S}$.
Since $\ohull{S}$ is also convex, this implies that
$\conv{\Sbar}\subseteq\ohull{S}$
(\Cref{pr:conhull-prop}\ref{pr:conhull-prop:aa}).

\pfpart{Part~(\ref{cor:sep-cvx-sets-conseqs:c}):}
Suppose $\ohull{X}$ and $\ohull{Y}$ are disjoint.
Since they are also closed and convex,
\Cref{thm:sep-cvx-sets} then implies that
they are strongly separated.
Therefore, $X$ and~$Y$, being subsets of $\ohull{X}$ and $\ohull{Y}$,
are also strongly separated.

Conversely, suppose $X$ and $Y$ are strongly separated.
Then there exists $\uu\in\Rn\wo\{\zero\}$
and $\beta,\gamma\in\R$ such that
\[
  \sup_{\xbar\in X} \xbar\cdot\uu
  \leq
  \beta
  <
  \gamma
  \leq
  \sup_{\ybar\in Y} \ybar\cdot\uu.
\]
That is, $X\subseteq H$ and $Y\subseteq H'$, where
\[
   H=\{\zbar\in\extspace :\: \zbar\cdot\uu\leq\beta\}
   \mbox{ and }
   H'=\{\zbar\in\extspace :\: \zbar\cdot\uu\geq\gamma\}.
\]
Since $H$ and $H'$ are closed astral halfspaces, by definition of outer hull,
$\ohull{X}\subseteq H$ and $\ohull{Y}\subseteq H'$.
Since $H$ and $H'$ are disjoint, this proves that
$\ohull{X}$ and $\ohull{Y}$ are as well.
\qedhere
\end{proof-parts}
\end{proof}

\Cref{cor:sep-cvx-sets-conseqs}(\ref{cor:sep-cvx-sets-conseqs:b})
shows that the convex hull of the closure of any set
$S\subseteq\extspace$
is equal to $S$'s outer hull,
which is also the smallest closed convex set that includes $S$.
In \Cref{sec:closure-convex-set} we saw that this is not generally
true if these operations are applied in the reverse order;
that is, $\clbar{\conv{S}}$, the closure of the convex hull of $S$, need
not in general be equal to $\ohull{S}$ (nor must it even be convex).

\indexg{separation of astral sets, strong!sets in $\Rn$ and|(}%
\Cref{thm:sep-cvx-sets} shows that two convex astral sets are
strongly separated if
they are closed (in $\extspace$) and disjoint.
The same is not true when working only in~$\Rn$.
In other words, it is not true, in general, that two convex sets
$X,Y\subseteq\Rn$ must be strongly separated if
they are
disjoint and closed in $\Rn$.
Nevertheless,
\Cref{cor:sep-cvx-sets-conseqs}(\ref{cor:sep-cvx-sets-conseqs:c})
provides a simple test for determining if they are strongly separated,
specifically, if their outer hulls are disjoint.
Equivalently, in this case, $X$ and $Y$ are strongly separated if
and only if their astral closures are disjoint, that is,
if and only if $\Xbar\cap\Ybar=\emptyset$,
as follows from \Cref{thm:e:6}
since the sets are convex and in $\Rn$.

Here is an example:

\begin{example}
In $\R^2$, let $X=\{\xx\in\Rpos^2 :\: x_1 x_2 \geq 1\}$,
and let $Y=\R\times\{0\}$.
These sets are convex, closed (in $\R^2$), and disjoint.
Nonetheless, they are not strongly separated.
To prove this, note that the sequence
$\seq{\trans{[t,1/t]}}$ in $X$
and the sequence $\seq{\trans{[t,0]}}$ in $Y$
both converge to $\limray{\ee_1}$,
implying $\limray{\ee_1}\in\Xbar\cap\Ybar=(\ohull{X})\cap(\ohull{Y})$
by \Cref{thm:e:6}.
Therefore, by
\Cref{cor:sep-cvx-sets-conseqs}(\ref{cor:sep-cvx-sets-conseqs:c}),
$X$ and $Y$ are not strongly separated.
\end{example}

Thus, although the strong separation of sets in $\Rn$ involves only
notions from standard convex analysis, we see that concepts from
astral space can still provide insight.
In general, convex sets $X$ and $Y$ in $\Rn$ are strongly separated if
and only if $X-Y$ and the singleton $\{\zero\}$ are strongly separated
(as reflected in the equivalence of Eqs.~\ref{eqn:std-str-sep-crit1}
and~\ref{eqn:std-str-sep-crit2}).
Therefore, by the preceding argument, $X$ and $Y$ are strongly
separated if and only if $\{\zero\}$ and $\clbar{X-Y}$ are disjoint,
that is, if and only if
$\zero\not\in (\clbar{X-Y})\cap\Rn=\cl(X-Y)$.
This is equivalent to Theorem~11.4 of \idxroc\citet{ROC}.%
\indexg{separation of astral sets, strong!sets in $\Rn$ and|)}

\indexg{outer convex hull!affine map@under affine map|(}%
In \Cref{thm:e:9}, we saw that
$\ohull{F(S)}=F(\ohull{S})$
for any affine function $F$ and
any finite set $S\subseteq\extspace$.
Using
\Cref{cor:sep-cvx-sets-conseqs},
we can now prove that the same holds even when $S$ is not
necessarily finite.

\begin{corollary}  \label{cor:ohull-fs-is-f-ohull-s}
  Let $\A\in\R^{m\times n}$, $\bbar\in\extspac{m}$, and
  let $F:\extspace\rightarrow\extspac{m}$ be the
  affine map
  $F(\zbar)=\bbar\plusl \A\zbar$
  for $\zbar\in\extspace$.
  Let $S\subseteq\extspace$.
  Then
  \[
     \ohull{F(S)} = F(\ohull{S}).
  \]
\end{corollary}

\begin{proof}
That $F(\ohull{S}) \subseteq \ohull{F(S)}$
was proved in \Cref{lem:f-conv-S-in-conv-F-S}.
For the reverse inclusion,
$\ohull{S}$ is closed and convex, so
$F(\ohull{S})$ is also convex
(by \Cref{cor:thm:e:9})
and closed
(by
\Cref{cor:aff-img-closed-is-closed}\ref{cor:aff-img-closed-is-closed:a}).
Furthermore, $F(S)\subseteq F(\ohull{S})$ since
$S\subseteq\ohull{S}$.
Therefore, $\ohull{F(S)}\subseteq F(\ohull{S})$
since, by
\Cref{cor:sep-cvx-sets-conseqs}(\ref{cor:sep-cvx-sets-conseqs:a}),
$\ohull{F(S)}$ is the smallest closed convex set that includes $F(S)$.%
\indexg{outer convex hull!affine map@under affine map|)}%
\end{proof}

\indexg{separation of astral sets, strong!astral cone@from astral cone|(}%
We next consider separating a closed convex astral cone from a
disjoint closed convex set.
In this case, it suffices to use a homogeneous astral hyperplane
for the separation, as we show next
(thereby generalizing
\Cref{lem:acone-sep-from-fin-pt}).

\begin{theorem}   \label{thm:sep-cone-from-cvx}
  Let $K\subseteq\extspace$ be a closed convex astral cone, and let
  $Y\subseteq\extspace$ be a nonempty closed convex set
  with $K\cap Y=\emptyset$.
  Then there exists $\uu\in\Rn\wo\{\zero\}$ such that
  \[
     \sup_{\xbar\in K} \xbar\cdot\uu
     \leq
     0
     <
     \inf_{\ybar\in Y} \ybar\cdot \uu.
  \]
\end{theorem}

\begin{proof}
Since $K$ and $Y$ are nonempty, closed, convex and disjoint, we can
apply
\Cref{thm:sep-cvx-sets} yielding that there exist
$\uu\in\Rn\wo\{\zero\}$
and $\beta\in\R$ such that
\[
  \sup_{\xbar\in K} \xbar\cdot\uu
  <
  \beta
  <
  \inf_{\ybar\in Y} \ybar\cdot \uu.
\]
Since $K$ is an astral cone,
we have that $\zero\in K$
(by \Cref{pr:ast-cone-is-naive}),
implying that $0=\zero\cdot\uu<\beta$.
Also, for all $\xbar\in K$,
we claim that $\xbar\cdot\uu\leq 0$.
Otherwise, if $\xbar\cdot\uu>0$ then
$\limray{\xbar}\cdot\uu=+\infty$.
On the other hand, $\limray{\xbar}\in K$
(by \Cref{pr:ast-cone-is-naive}, since $\xbar\in K$),
implying $\limray{\xbar}\cdot\uu < \beta$, a contradiction.

Together, these prove all the parts of the theorem.%
\indexg{separation of astral sets, strong|)}%
\end{proof}

Here are some consequences of \Cref{thm:sep-cone-from-cvx}.
\indexg{outer conic hull!smallest closed convex astral cone@as smallest closed convex astral cone|(}%
In particular, we now can see that the outer conic hull of any set
$S\subseteq\extspace$ is the smallest closed convex astral cone that
includes $S$.

\begin{corollary}   \label{cor:sep-ast-cone-conseqs}
\indexg{halfspaces, astral!astral cones as intersection of|(}%
\indexg{cones, astral!intersection of halfspaces@as intersection of halfspaces|(}%
  Let $S\subseteq\extspace$.
  \begin{letter-compact}
  \item        \label{cor:sep-ast-cone-conseqs:b}
    The following are equivalent:
    \begin{roman-compact}
    \item        \label{cor:sep-ast-cone-conseqs:b:1}
      $S$ is a closed convex astral cone.
    \item        \label{cor:sep-ast-cone-conseqs:b:2}
      $S$ is equal to the intersection of some collection of
      homogeneous, closed astral halfspaces.
    \item        \label{cor:sep-ast-cone-conseqs:b:3}
      $S$ is equal to the intersection of all
      homogeneous, closed astral halfspaces
      that include $S$;
      that is, $S=\oconich{S}$.
    \end{roman-compact}
  \item        \label{cor:sep-ast-cone-conseqs:a}
    The outer conic hull of $S$, $\oconich{S}$, is the smallest closed
    convex astral cone that includes $S$.
    (That is, $\oconich{S}$ is equal to the intersection of all closed
    convex astral cones in $\extspace$ that include $S$.)
  \end{letter-compact}
\end{corollary}

\begin{proof}
~

\begin{proof-parts}
\pfpart{Part~(\ref{cor:sep-ast-cone-conseqs:b}):}
The proof is very similar to that of
\Cref{cor:sep-cvx-sets-conseqs}(\ref{cor:sep-cvx-sets-conseqs:d}),
but in proving
that (\ref{cor:sep-ast-cone-conseqs:b:1})
implies (\ref{cor:sep-ast-cone-conseqs:b:3}),
we instead use \Cref{thm:sep-cone-from-cvx} to yield
a homogeneous halfspace that includes $K$ but not any particular point
$\zbar\in\eRn\setminus K$.%
\indexg{halfspaces, astral!astral cones as intersection of|)}%
\indexg{cones, astral!intersection of halfspaces@as intersection of halfspaces|)}

\pfpart{Part~(\ref{cor:sep-ast-cone-conseqs:a}):}
Similar to the proof of
\Cref{cor:sep-cvx-sets-conseqs}(\ref{cor:sep-cvx-sets-conseqs:a}).%
\indexg{separation of astral sets, strong!astral cone@from astral cone|)}%
\indexg{outer conic hull!smallest closed convex astral cone@as smallest closed convex astral cone|)}%
\qedhere
\end{proof-parts}
\end{proof}

\section{Astral dual separation of sets in $\Rn$}
\label{sec:ast-def-sep-thms}

So far, we have focused on how convex sets in astral space can be
separated using astral halfspaces defined by a vector $\uu\in\Rn$.
Next, we consider the dual question of
how convex sets in $\Rn$ can be separated using a
generalized form of halfspaces, each defined now by an astral point
$\xbar\in\extspace$.

There are many standard results regarding separation of
convex sets in $\Rn$ (see, for example,
\Cref{sec:prelim-sep-thms}). We are specifically interested
in strong separation. In the standard setting, disjoint convex sets must
satisfy additional conditions in order to be strongly separated.
For instance, in $\R$, the sets
$\set{\zero}$ and $\Rstrictpos$ are disjoint and convex, but cannot be
strongly separated.
In this section, we use astral points to define what we call
astral dual halfspaces. Using these halfspaces,
we are able to separate any two disjoint
convex sets in $\Rn$ without further conditions. This leads to
a characterization of convex sets in $\Rn$ as intersections of
astral dual halfspaces.

Much of what follows is closely related to the previous work of
\idxmartinezsinger\citet{lexicographic_separation,hemispaces},
which used a rather different formalism.
Here, we give a development based on astral space.

We build on the standard notion of strong separation, introduced
in \Cref{sec:prelim-sep-thms}.
In \Cref{roc:thm11.1} and again in
\eqref{eqn:std-str-sep-crit1},
we saw (using slightly different notation)
that two nonempty sets $U$ and $V$ in $\Rn$ are strongly separated if,
for some $\xx\in\Rn\wo\{\zero\}$,
\begin{equation}   \label{eqn:std-str-sep-crit}
   \sup_{\uu\in U} \xx\cdot\uu
   <
   \inf_{\vv\in V} \xx\cdot\vv.
\end{equation}
In this case, the sets are separated by the
hyperplane $J=\{\uu\in\Rn :\: \xx\cdot\uu = \beta\}$,
for some $\beta\in\R$.

Geometrically, nothing changes if we translate everything by any
vector $\rr\in\Rn$.
The translations $U-\rr$ and $V-\rr$ are still
strongly separated by a translation of this same hyperplane,
namely, $J-\rr$.
Algebraically, this is saying simply that
\begin{equation}   \label{eqn:std-str-sep-crit3}
   \sup_{\uu\in U} \xx\cdot(\uu-\rr)
   <
   \inf_{\vv\in V} \xx\cdot(\vv-\rr),
\end{equation}
which is of course equivalent to
\eqref{eqn:std-str-sep-crit}
since $\xx\cdot\rr$ is a constant.
Thus, $U$ and $V$ are strongly separated if
\eqref{eqn:std-str-sep-crit3}
holds for some $\xx\in\Rn\wo\{\zero\}$
and some (or all) $\rr\in\Rn$.

\indexg{separation, strong astral dual|(}%
In this form, we can generalize to separation of sets in $\Rn$ using
separating sets defined by points in $\extspace$ simply by replacing
$\xx$ in
\eqref{eqn:std-str-sep-crit3}
with an astral point $\xbar$, leading to:

\begin{definition}
\indexg{separation, strong astral dual!defined|(}%
Let $U$ and $V$ be nonempty subsets of $\Rn$.
We say that $U$ and $V$ are
\emph{strong\-ly astrally dually separated}
(or \emph{strongly dually separated}, for short)
if there exist a point
$\xbar\in\extspace$
and a point $\rr\in\Rn$ such that
\begin{equation}   \label{eqn:ast-str-sep-defn}
   \sup_{\uu\in U} \xbar\cdot(\uu-\rr)
   <
   \inf_{\vv\in V} \xbar\cdot(\vv-\rr).
\end{equation}
(Note that $\xbar$ cannot be $\zero$ since $U$ and $V$ are nonempty.)
\end{definition}

Equivalently,
\eqref{eqn:ast-str-sep-defn}
holds if and only if, for some $\beta\in\R$ and
$\epsilon\in\Rstrictpos$,
$\xbar\cdot(\uu-\rr)<\beta-\epsilon$ for all $\uu\in U$,
and
$\xbar\cdot(\vv-\rr)>\beta+\epsilon$ for all $\vv\in V$.
Intuitively, $\rr\in\Rn$ is acting to translate the sets $U$ and $V$,
and $\xbar$ is specifying
a separation set
$\{\uu\in\Rn :\: \xbar\cdot\uu=\beta\}$,
for some $\beta\in\R$,
analogous to a standard hyperplane.%
\indexg{separation, strong astral dual!defined|)}

\begin{example}
\label{ex:dual-separable}
For instance, in $\R^2$, let $\bb\in\R^2$, and let
\begin{align*}
  U
  &=
  \Braces{\uu\in\R^2 :\: u_1 \leq b_1 \mbox{ and } u_2 \leq b_2}
  \\
  V
  &=
  \Braces{\vv\in\R^2 :\: v_1 \geq b_1} \setminus U
  \\
  &=
  \Braces{\vv\in\R^2 :\:
    v_1 > b_1
    \text{ or }
    (v_1 = b_1 \text{ and } v_2 > b_2)
  }.  
\end{align*}
Thus, $U$ includes all points in the plane that are
``below and to the left'' of $\bb$,
while $V$ includes all points ``to the right'' of $\bb$, except those
in $U$.

These sets are convex and disjoint.
The set $V$ is neither open nor closed.
These sets cannot be strongly separated
in the standard sense, using standard hyperplanes.
Indeed, suppose $\xx\in\R^2$ satisfies
\eqref{eqn:std-str-sep-crit}.
If $x_2\neq 0$ then the right-hand side of that equation
would be $-\infty$, which is impossible.
Thus, $\xx=\lambda\ee_1$ for some $\lambda\in\R$.
If $\lambda<0$ then the left-hand side would be $+\infty$, which is
impossible.
And if $\lambda\geq 0$, then both sides of the equation would be equal
to $\lambda b_1$, which also violates that strict inequality.
Thus, a contradiction is reached in all cases.

These sets are, however, strongly dually separated.
In particular, setting $\xbar=\limray{\ebar_1}\plusl\limray{\ebar_2}$
and $\rr=\bb$,
it can be checked that
\eqref{eqn:ast-str-sep-defn} holds,
specifically, that
$\xbar\cdot(\uu-\rr)\leq 0$ for all $\uu\in U$
while
$\xbar\cdot(\vv-\rr) = +\infty$ for all $\vv\in V$.

On the other hand, if $\bb=\trans{[2,2]}$,
and if we require $\rr=\zero$,
then
\eqref{eqn:ast-str-sep-defn} does not hold for any
$\xbar\in\extspac{2}$.
To see this, note first that,
as just argued, this equation cannot hold for any
$\xbar\in\R^2$.
So suppose it holds for some
$\xbar\in\extspac{2}\setminus\R^2$,
and let $\ww\in\R^2$ be $\xbar$'s dominant direction.
Then $\ww$, being a unit vector, is in $U$,
implying the left-hand side of
\eqref{eqn:ast-str-sep-defn} is equal to $+\infty$
(since $\xbar\cdot\ww=+\infty$), a contradiction.
\end{example}

This example shows why we include the translation $\rr$
in the definition of strong dual separation given in
\eqref{eqn:ast-str-sep-defn}.
\indexg{separation, strong astral dual!characterized|(}%
Nevertheless, there is an alternative formulation,
analogous to \eqref{eqn:std-str-sep-crit2},
that dispenses with $\rr$.
In this formulation, we instead require that
$\xbar\cdot(\uu-\vv)$ be bounded from above by a negative constant for all $\uu\in U$ and
$\vv\in V$, that is, that
\[
   \sup_{\uu\in U, \vv\in V} \xbar\cdot(\uu-\vv)
   <
   0.
\]
The next theorem shows that this condition is equivalent to the one
given in
\eqref{eqn:ast-str-sep-defn}, spelling out the specific relationship
between these two notions:

\begin{theorem}   \label{thm:ast-str-sep-equiv-no-trans}
  Let $U,V\subseteq\Rn$ and let $\xbar\in\extspace$.
  Then the following are equivalent:
  \begin{letter-compact}
  \item    \label{thm:ast-str-sep-equiv-no-trans:a}
    There exist $\ybar\in\lb{\zero}{\xbar}$ and $\rr\in\Rn$
    such that
    \begin{equation}    \label{eq:thm:ast-str-sep-equiv-no-trans:1}
       \sup_{\uu\in U} \ybar\cdot(\uu-\rr)
       <
       \inf_{\vv\in V} \ybar\cdot(\vv-\rr).
    \end{equation}
  \item    \label{thm:ast-str-sep-equiv-no-trans:b}
    We have
    \begin{equation}    \label{eq:thm:ast-str-sep-equiv-no-trans:3}
       \sup_{\uu\in U, \vv\in V} \xbar\cdot(\uu-\vv)
       <
       0.
    \end{equation}
  \end{letter-compact}
  Therefore, if $U$ and $V$ are nonempty, then
  they are strongly dually separated if and only if
  \eqref{eq:thm:ast-str-sep-equiv-no-trans:3} holds for some
  $\xbar\in\extspace$.
\end{theorem}

\begin{proof}
~

\begin{proof-parts}
\pfpart{%
  (\ref{thm:ast-str-sep-equiv-no-trans:a})
  $\Rightarrow$
  (\ref{thm:ast-str-sep-equiv-no-trans:b}):
}
Suppose that
\eqref{eq:thm:ast-str-sep-equiv-no-trans:1}
holds for some $\ybar\in\lb{\zero}{\xbar}$ and $\rr\in\Rn$.
Then there exist $\beta,\gamma\in\R$ such that
\[
  \sup_{\uu\in U} \ybar\cdot(\uu-\rr)
  \leq
  \beta
  <
  \gamma
  \leq
  \inf_{\vv\in V} \ybar\cdot(\vv-\rr).
\]
Suppose $\uu\in U$, $\vv\in V$.
Then $\ybar\cdot(\vv-\rr)\geq\gamma>-\infty$
and $\ybar\cdot(\uu-\rr)\leq\beta$,
implying
$\ybar\cdot(\rr-\uu)\geq -\beta>-\infty$.
Consequently, by \Cref{pr:i:1},
\begin{equation}  \label{eq:thm:ast-str-sep-equiv-no-trans:2}
  \ybar\cdot(\vv-\uu)
  =
  \ybar\cdot(\vv-\rr)
  +
  \ybar\cdot(\rr-\uu)
  \geq
  \gamma-\beta,
\end{equation}
since the expressions being added are summable.
This further implies that
\[
  0
  <
  \gamma-\beta
  \leq
  \ybar\cdot(\vv-\uu)
  \leq
  \max\braces{ 0, \xbar\cdot(\vv-\uu) }
  =
  \xbar\cdot(\vv-\uu).
\]
The second inequality is from
\eqref{eq:thm:ast-str-sep-equiv-no-trans:2}.
The third inequality is
by \Cref{pr:seg-simplify}(\ref{pr:seg-simplify:a},\ref{pr:seg-simplify:b})
since $\ybar\in\lb{\zero}{\xbar}$.
The equality is because the maximum that appears here must be equal to
one of its two arguments, but cannot be $0$ since
that would contradict the preceding inequalities.
Thus, $\xbar\cdot(\uu-\vv)\leq\beta-\gamma<0$ for all
$\uu\in U$ and $\vv\in V$, proving the claim.

\pfpart{%
  (\ref{thm:ast-str-sep-equiv-no-trans:b})
  $\Rightarrow$
  (\ref{thm:ast-str-sep-equiv-no-trans:a}):
}
Proof is by induction on the astral rank of $\xbar$.
More specifically, we prove the following by induction on
$k=0,\ldots,n$:
For all $U,V\subseteq\Rn$ and for all $\xbar\in\extspace$,
if $\xbar$ has astral rank $k$
and if \eqref{eq:thm:ast-str-sep-equiv-no-trans:3} holds,
then there exist $\ybar\in\lb{\zero}{\xbar}$ and $\rr\in\rspanxbar$
for which \eqref{eq:thm:ast-str-sep-equiv-no-trans:1} holds.

Before beginning the induction, we consider the case that either $U$
or $V$ (or both) are empty.
In this case, we can choose $\ybar=\zero$ (which is in
$\lb{\zero}{\xbar}$) and $\rr=\zero$ (which is in $\rspanxbar$).
If $V=\emptyset$, then
\eqref{eq:thm:ast-str-sep-equiv-no-trans:1}
holds since
the right-hand side
is $+\infty$ and the left-hand side is either
$-\infty$ (if $U=\emptyset$) or $0$ (otherwise).
The case $U=\emptyset$ is similar.

In the induction argument that follows, we suppose that
$\xbar$ has astral rank $k$
and that \eqref{eq:thm:ast-str-sep-equiv-no-trans:3} holds,
implying there exists $\beta\in\Rstrictpos$ such that
$\sup_{\uu\in U,\vv\in V} \xbar\cdot(\uu-\vv) \leq -\beta < 0$.

In the base case that $k=0$, we have $\xbar=\qq$ for some $\qq\in\Rn$.
The case that either $U$ or $V$ is empty was handled above, so
we assume both are nonempty.
Then
for all $\uu\in U$ and $\vv\in V$,
$\qq\cdot(\uu-\vv)\leq-\beta$,
implying
$\qq\cdot\uu \leq \qq\cdot\vv-\beta$.
Thus,
\[
  -\infty
  <
  \sup_{\uu\in U} \qq\cdot\uu
  \leq
  \inf_{\vv\in V} \qq\cdot\vv - \beta
  <
  +\infty,
\]
where the first and last inequalities are because
$U$ and $V$ are nonempty.
Letting $\ybar=\qq$ and $\rr=\zero$, the claim thereby follows in this
case.

For the inductive step, let $k>0$ and assume the claim holds for
$k-1$.
We again assume $U$ and $V$ are both nonempty, since otherwise we can
apply the argument above.
Then we can write
$\xbar=\limray{\ww}\plusl\xbarperp$
where $\ww\in\Rn$, $\norm{\ww}=1$, is $\xbar$'s dominant direction,
and $\xbarperp$ is the projection of $\xbar$ orthogonal to $\ww$
(\Cref{pr:h:6}).

For all $\uu\in U$ and $\vv\in V$, we must have
$\ww\cdot\uu\leq\ww\cdot\vv$; otherwise,
if $\ww\cdot(\uu-\vv)>0$, then we would have
$\xbar\cdot(\uu-\vv)=+\infty$, a
contradiction.
Thus,
\[
  -\infty
  <
  \sup_{\uu\in U} \ww\cdot\uu
  \leq
  \inf_{\vv\in V} \ww\cdot\vv
  <
  +\infty,
\]
where the first and last inequalities are because $U$ and $V$ are
nonempty.

Let $\lambda\in\R$ be such that
$\sup_{\uu\in U} \ww\cdot\uu \leq \lambda \leq \inf_{\vv\in V} \ww\cdot\vv$.
Also,
let
\[ J=\{\uu\in\Rn :\: \ww\cdot\uu = \lambda\}, \]
let $U'=U\cap J$, and let $V'=V\cap J$.
Then for all $\uu\in U'$ and $\vv\in V'$, we have
$\ww\cdot(\uu-\vv)=\ww\cdot\uu-\ww\cdot\vv=0$
(since $\uu,\vv\in J$).
Thus,
\[
  \xbarperp\cdot(\uu-\vv)
  =
  (\limray{\ww}\plusl\xbarperp)\cdot(\uu-\vv)
  =
  \xbar\cdot(\uu-\vv)
  \leq
  -\beta
  <
  0.
\]
Therefore, we can apply our inductive hypothesis to $U'$, $V'$ and
$\xbarperp$
(whose astral rank is $k-1$, by \Cref{pr:h:6}),
implying that there exist
$\ybar'\in\lb{\zero}{\xbarperp}$
and
$\rr'\in\rspanxbarperp$
such that
\[
  \sup_{\uu\in U'} \ybar'\cdot(\uu-\rr')
  <
  \inf_{\vv\in V'} \ybar'\cdot(\vv-\rr').
\]

Let $\ybar=\limray{\ww}\plusl\ybar'$ and let
$\rr=\lambda\ww+\rr'$.
Then
$\ybar\in\limray{\ww}\plusl\lb{\zero}{\xbarperp}\subseteq\lb{\zero}{\xbar}$,
with the second inclusion by \Cref{cor:seg:zero}.
Further, $\rr\in\rspanxbar$ since if $\xbarperp=\VV'\omm\plusl\qq'$ for
some $\VV'\in\R^{n\times k'}$, $k'\geq 0$, $\qq'\in\Rn$ then
$\rr'\in\rspanxbarperp=\colspace[\VV',\qq']$,
so
$\rr\in\colspace[\ww,\VV',\qq']=\rspanxbar$.

Also, $\xbarperp\cdot\ww=0$, implying $\rr'\cdot\ww=0$
by \Cref{pr:rspan-sing-equiv-dual}.
Moreover, $\ybar'\cdot\ww=0$ since $\zero$ and $\xbarperp$ are both in
$\{\zbar\in\extspace :\: \zbar\cdot\ww=0\}$,
which, being convex
(\Cref{pr:e1}\ref{pr:e1:c}),
must include $\lb{\zero}{\xbarperp}$, specifically, $\ybar'$.

Let $\uu\in U$.
We next argue that
\begin{equation}  \label{eq:thm:ast-str-sep-equiv-no-trans:5}
  \ybar\cdot(\uu-\rr)
  =
  \begin{cases}
    \ybar'\cdot(\uu-\rr') & \text{if $\uu\in U'$,} \\
    -\infty               & \text{otherwise.} \\
  \end{cases}
\end{equation}
We have
\begin{equation}  \label{eq:thm:ast-str-sep-equiv-no-trans:4}
  \ww\cdot(\uu-\rr)
  =
  \ww\cdot(\uu-\lambda\ww-\rr')
  =
  \ww\cdot\uu - \lambda
\end{equation}
since $\ww\cdot\rr'=0$ and $\norm{\ww}=1$.
Also, by $\lambda$'s definition, $\ww\cdot\uu\leq \lambda$.
Therefore,
if $\uu\not\in U'$ then $\uu\not\in J$ so $\ww\cdot\uu<\lambda$
and $\ww\cdot(\uu-\rr)<0$
(by Eq.~\ref{eq:thm:ast-str-sep-equiv-no-trans:4}),
implying
\[
  \ybar\cdot(\uu-\rr)
  =
  \limray{\ww}\cdot(\uu-\rr)\plusl\ybar'\cdot(\uu-\rr)
  =
  -\infty.
\]
Similarly, if $\uu\in U'$ then
$\ww\cdot\uu=\lambda$ so
$\ww\cdot(\uu-\rr)=0$
(again by Eq.~\ref{eq:thm:ast-str-sep-equiv-no-trans:4}),
implying
\[
  \ybar\cdot(\uu-\rr)
  =
  \ybar'\cdot(\uu-\rr)
  =
  \ybar'\cdot(\uu-\lambda\ww-\rr')
  =
  \ybar'\cdot(\uu-\rr')
\]
where the last equality is because $\ybar'\cdot\ww=0$
(and using \Cref{pr:i:1}).

By a similar argument, for $\vv\in V$,
\begin{equation}  \label{eq:thm:ast-str-sep-equiv-no-trans:6}
  \ybar\cdot(\vv-\rr)
  =
  \begin{cases}
    \ybar'\cdot(\vv-\rr') & \text{if $\vv\in V'$,} \\
    +\infty               & \text{otherwise.} \\
  \end{cases}
\end{equation}
Thus, combining, we have
\[
  \sup_{\uu\in U} \ybar\cdot(\uu-\rr)
  \leq
  \sup_{\uu\in U'} \ybar'\cdot(\uu-\rr')
  <
  \inf_{\vv\in V'} \ybar'\cdot(\vv-\rr')
  \leq
  \inf_{\vv\in V} \ybar\cdot(\vv-\rr).
\]
The first and third inequalities are by
Eqs.~(\ref{eq:thm:ast-str-sep-equiv-no-trans:5})
and~(\ref{eq:thm:ast-str-sep-equiv-no-trans:6}),
respectively.
The second inequality is by inductive hypothesis.

This completes the induction and the proof.%
\indexg{separation, strong astral dual!characterized|)}%
\qedhere
\end{proof-parts}
\end{proof}

The equivalence given in
\Cref{thm:ast-str-sep-equiv-no-trans}
does not hold in general if, in
part~(\ref{thm:ast-str-sep-equiv-no-trans:a}) of that 
\namecref{thm:ast-str-sep-equiv-no-trans}, we
require that $\ybar=\xbar$.
Here is an example:

\begin{example}
In $\R^2$, let $U$ be the closed halfplane
$U=\{\uu\in\R^2 : u_1 \leq 0\}$, and let
$V$ be the open halfplane that is its complement,
$V=\R^2\setminus U$.
Let $\xbar=\limray{\ee_1}\plusl\limray{\ee_2}$.
Then \eqref{eq:thm:ast-str-sep-equiv-no-trans:3} is satisfied
since
$\xbar\cdot(\uu-\vv)=-\infty$
for all $\uu\in U$ and $\vv\in V$.
However, if $\ybar=\xbar$, then
\eqref{eq:thm:ast-str-sep-equiv-no-trans:1} is not satisfied for any
$\rr\in\R^2$.
This is because if $r_1>0$ then there exists $\vv\in V$ with
$r_1>v_1>0$, implying $\xbar\cdot(\vv-\rr)=-\infty$ so that
\eqref{eq:thm:ast-str-sep-equiv-no-trans:1}
is unsatisfied since its right-hand side is $-\infty$.
Similarly, if $r_1\leq 0$ then there exists $\uu\in U$ with
$u_1=r_1$ and $u_2>r_2$ so that
$\xbar\cdot(\uu-\rr)=+\infty$, again implying that
\eqref{eq:thm:ast-str-sep-equiv-no-trans:1}
is unsatisfied.

On the other hand,
\eqref{eq:thm:ast-str-sep-equiv-no-trans:1}
is satisfied if we choose $\rr=\zero$
and $\ybar=\limray{\ee_1}$, which is in $\lb{\zero}{\xbar}$.
\end{example}

\indexg{separation, strong astral dual!main theorem|(}%
Every pair of nonempty,
disjoint convex sets in $\Rn$ is strongly dually
separated:

\begin{theorem}    \label{thm:ast-def-sep-cvx-sets}
  Let $U,V\subseteq\Rn$ be nonempty, convex and disjoint.
  Then $U$ and $V$ are strongly dually separated.
\end{theorem}

Before proving the theorem, we give a lemma handling the special case
that one set is the singleton $\{\zero\}$.
We then show how this implies the general case.

\begin{lemma}  \label{lem:ast-def-sep-from-orig}
  Let $U\subseteq\Rn\wo\{\zero\}$ be convex.
  Then there exists $\xbar\in\extspace$ for which
  \begin{equation}  \label{eq:lem:ast-def-sep-from-orig:1}
     \sup_{\uu\in U} \xbar\cdot\uu
     <
     0.
  \end{equation}
\end{lemma}

\begin{proof}
Proof is by induction on the dimension of a linear subspace that
includes $U$.
More precisely, we prove by induction on $d=0,\ldots,n$
that for all convex sets $U\subseteq\Rn$
and for all linear subspaces $L\subseteq\Rn$,
if $U\subseteq L\setminus\{\zero\}$ and $\dim{L}\leq d$
then there exists $\xbar\in\extspace$ satisfying
\eqref{eq:lem:ast-def-sep-from-orig:1}.

As a preliminary step, we note that if $U=\emptyset$, then we can
choose $\xbar=\zero$ so that
\eqref{eq:lem:ast-def-sep-from-orig:1}
holds vacuously.

In the base case that $d=0$, we must have $L=\{\zero\}$,
implying $U=\emptyset$.
Thus, as just discussed, we can choose $\xbar=\zero$ in this case.

For the inductive step, let $d>0$ and assume the claim holds for
$d-1$.
Let $U$ and $L$ be as described in the claim.
Since the case $U=\emptyset$ was handled above, we assume $U$ is
nonempty.

Since $U$ and the singleton $\{\zero\}$ are disjoint,
their relative interiors are as well.
Therefore, by \Cref{roc:thm11.3},
there exists a (standard) hyperplane
\[ J=\{ \uu\in\Rn :\: \ww\cdot\uu = \beta \} \]
that properly separates them,
for some $\ww\in\Rn$ and $\beta\in\R$.
That is,
$\ww\cdot\uu\leq \beta$ for all $\uu\in U$,
and $0=\ww\cdot\zero\geq \beta$.
Furthermore, there is some point in $U\cup \{\zero\}$
that is not in~$J$.
If $\zero\not\in J$ then $0=\ww\cdot\zero>\beta$,
implying
$\sup_{\uu\in U} \ww\cdot\uu \leq \beta < 0$,
and therefore that
\eqref{eq:lem:ast-def-sep-from-orig:1}
is satisfied if we choose $\xbar=\ww$.

Otherwise, $\zero\in J$, implying $\beta=\ww\cdot\zero=0$,
and furthermore, since the separation is proper,
that there exists a point $\uhat\in U\setminus J$.
Thus, $\ww\cdot\uu\leq 0$ for all $\uu\in U$,
and $\ww\cdot\uhat< 0$.

Let
$U'=U\cap J$, and let $L'=L\cap J$.
Then $L'$ is a linear subspace, and
$U'$ is convex and included in $L'\setminus\{\zero\}$.
Further, $\uhat\in U\subseteq L$, but
$\uhat\not\in J$, implying $\uhat\not\in L'$.
Thus, $L'$ is a proper subset of $L$, so $\dim{L'}<\dim{L}$.
Therefore, we can apply our inductive hypothesis (to $U'$ and $L'$),
yielding that there exists $\xbar'\in\extspace$ such that
$\gamma<0$ where $\gamma=\sup_{\uu\in U'}\xbar'\cdot\uu$.

Let $\xbar=\limray{\ww}\plusl\xbar'$.
Let $\uu\in U$.
If $\uu\not\in J$ then $\ww\cdot\uu<0$
so $\xbar\cdot\uu=-\infty$.
Otherwise, if $\uu\in J$ then $\uu\in U'$ and $\ww\cdot\uu=0$
so $\xbar\cdot\uu=\xbar'\cdot\uu\leq \gamma$.
Thus,
\[
  \sup_{\uu\in U} \xbar\cdot\uu \leq \gamma < 0,
\]
completing the induction and the proof.
\end{proof}

\begin{proof}[Proof of \Cref{thm:ast-def-sep-cvx-sets}]
Let $W=U-V$.
Then $-V$ is convex by \Cref{pr:aff-preserves-cvx},
so $W$ is convex by \Cref{roc:thm3.1}.
Further, $\zero\not\in W$ since $U$ and $V$ are disjoint.
Therefore, by \Cref{lem:ast-def-sep-from-orig},
applied to $W$,
there exists $\xbar\in\extspace$ such that
$\sup_{\uu\in U,\vv\in V} \xbar\cdot(\uu-\vv)<0$.
By \Cref{thm:ast-str-sep-equiv-no-trans},
it follows that
$U$ and $V$ are strongly dually separated.%
\indexg{separation, strong astral dual!main theorem|)}%
\indexg{separation, strong astral dual|)}%
\end{proof}

\indexg{dual halfspaces, astral|(}%
\indexg{dual halfspaces, astral!defined|(}%
An \emph{astral dual halfspace} (or \emph{dual halfspace}
for short) is a set of the form
\begin{equation}    \label{eqn:ast-def-hfspace-defn}
\indexm{h x r beta}{$\ahfsp{\xbar}{\rr}{\beta}$}{astral dual halfspace}%
   \ahfsp{\xbar}{\rr}{\beta}
   =
   \Braces{\uu\in\Rn :\: \xbar\cdot(\uu-\rr)\leq \beta}
\end{equation}
for some $\xbar\in\extspace\,\wo\{\zero\}$, $\rr\in\Rn$,
\indexg{dual halfspaces, astral!defined|)}%
$\beta\in\R$.
As we show next, dual halfspaces are always convex, and when $\rr=\zero$ and $\beta=0$, they are in fact convex cones.
Also,
although the definition in \eqref{eqn:ast-def-hfspace-defn} uses a nonstrict inequality, dual halfspaces
are not necessarily closed, but rather may be closed, open, or neither.
Moreover, as seen next,
it turns out that there is no need to introduce dual halfspace variants based
on a strict inequality
since such sets can be rewritten in the form of \eqref{eqn:ast-def-hfspace-defn}.
Consequently, the complement of every dual halfspace is also a
dual halfspace.

\begin{proposition}  \label{pr:ast-def-hfspace-convex}
  Let $\xbar\in\eRn\,\wo\{\zero\}$, $\rr\in\Rn$,
  $\beta\in\R$,
  and let
  \begin{align*}
    H
    &=
    \Braces{\uu\in\Rn :\: \xbar\cdot(\uu-\rr)\leq \beta}
  \intertext{be the corresponding dual halfspace. Let}
    H'
    &=
    \Braces{\uu\in\Rn :\: \xbar\cdot(\uu-\rr) < \beta},
  \end{align*}
  and let $H^c=\Rn\wo H$ be $H$'s complement.
  Then:
  \begin{letter-compact}
  \item \label{i:dual-hfspace-complement}
    $H'$ and $H^c$ are also dual halfspaces.
  \item \label{i:dual-hfspace-conv}
    $H$, $H'$, and $H^c$ are all convex.
  \item \label{i:dual-hfspace-cone}
    If $\rr=\zero$ and $\beta=0$, then $H$ is a convex cone.
  \end{letter-compact}
\end{proposition}

\begin{proof}
~

\begin{proof-parts}
\pfpart{Part (\ref{i:dual-hfspace-complement}):}
We show first that $H'$ is a dual halfspace.
If $\xbar\cdot\ww\neq\beta$ for all $\ww\in\Rn$,
then
$\xbar\cdot\ww\leq\beta$ if and only if
$\xbar\cdot\ww<\beta$, so
$H'$ is equal to $H$ and therefore a dual halfspace.

Otherwise, there exists a point $\what\in\Rn$ for which $\xbar\cdot\what=\beta$.
Let $\rr'=\rr+\what$.
Then for all $\uu\in\Rn$,
\begin{equation}  \label{eq:pr:ast-def-hfspace-convex:4}
  \xbar\cdot(\uu-\rr) - \beta
  =
  \xbar\cdot(\uu-\rr) - \xbar\cdot\what
  =
  \xbar\cdot(\uu-\rr-\what)
  =
  \xbar\cdot(\uu-\rr'),
\end{equation}
where the second equality is by \Cref{pr:i:1} (with summability ensured since $\xbar\cdot\what\in\R$).
Thus,
\begin{align}
  H'
   =
   \braces{\uu\in\Rn :\: \xbar\cdot(\uu-\rr)< \beta}
   &=
   \braces{\uu\in\Rn :\: \xbar\cdot(\uu-\rr')< 0}
   \nonumber
\\
   &=
   \braces{\uu\in\Rn :\: \limray{\xbar}\cdot(\uu-\rr')< 0}
   \nonumber
\\
   &=
   \braces{\uu\in\Rn :\: \limray{\xbar}\cdot(\uu-\rr')\le -1},
  \label{eq:pr:ast-def-hfspace-convex:3}
\end{align}
where the second equality is by \eqref{eq:pr:ast-def-hfspace-convex:4},
and the last equality is because $\limray{\xbar}\cdot(\uu-\rr')\in\set{-\infty,0,+\infty}$ for all $\uu\in\Rn$
since $\limray{\xbar}$ is an icon
(by Propositions~\ref{pr:icon-equiv}\ref{pr:icon-equiv:a}\ref{pr:icon-equiv:b}
and~\ref{pr:i:8}\ref{pr:i:8-infprod}).
Evidently, the rightmost set in
\eqref{eq:pr:ast-def-hfspace-convex:3}
is a dual halfspace; therefore, $H'$ is as well
(noting that $\limray{\xbar}\neq\zero$ since $\xbar\neq\zero$).

Since
\[
  H^c
  =
  \Braces{\uu\in\Rn :\: \xbar\cdot(\uu-\rr) > \beta}
  =
  \Braces{\uu\in\Rn :\: -\xbar\cdot(\uu-\rr) < -\beta},
\]
$H^c$ has the same form as $H'$; therefore, $H^c$ is a dual halfspace
as well.

\pfpart{Part (\ref{i:dual-hfspace-conv}):}
Setting $\ww=\uu-\rr$,
we can write $H$ as
\begin{equation}  \label{eq:ast-def-hfspace-convex:1}
  H = \rr + \Braces{\ww\in\Rn :\: \xbar\cdot\ww \leq \beta}.
\end{equation}
Thus, $H$ is the translation of a sublevel set
of the function $\ww\mapsto\xbar\cdot\ww$.
This function is convex by
\Crefequiv{thm:h:5}{thm:h:5a0}{thm:h:5b}.
Therefore, the sublevel set appearing on the far right
of \eqref{eq:ast-def-hfspace-convex:1}
is also convex
(by \Cref{roc:thm4.6}),
so $H$ is as well
(by \Cref{pr:aff-preserves-cvx}).
The convexity of $H'$ and $H^c$ then follows by
part~(\ref{i:dual-hfspace-complement}) having just proved all dual halfspaces are convex.

\pfpart{Part (\ref{i:dual-hfspace-cone}):}
If $\rr=\zero$ and $\beta=0$, then
$H$ is also a cone since it includes $\zero$,
and since if $\uu\in H$ then $\xbar\cdot\uu\leq 0$,
implying, for $\lambda\in\Rpos$, that
$\xbar\cdot(\lambda\uu)=\lambda(\xbar\cdot\uu)\leq 0$
so that $\lambda\uu$ is in $H$ as well.
\qedhere
\end{proof-parts}
\end{proof}

In this terminology,
according to
\Cref{thm:ast-def-sep-cvx-sets},
if two nonempty sets $U,V\subseteq\Rn$ are strongly dually
separated,
then they must be included respectively in a pair of disjoint
astral dual halfspaces.
More specifically, if $U,V$ satisfy \eqref{eqn:ast-str-sep-defn},
then
\begin{align*}
  U
  &\subseteq
  \Braces{\uu\in\Rn :\: \xbar\cdot(\uu-\rr)\leq \beta},
  \\
  V
  &\subseteq
  \Braces{\uu\in\Rn :\: \xbar\cdot(\uu-\rr)\geq \gamma},
\end{align*}
where $\beta,\gamma\in\R$ are such that
$\sup_{\uu\in U} \xbar\cdot(\uu-\rr)\leq\beta<\gamma\leq
 \inf_{\uu\in V} \xbar\cdot(\uu-\rr)$.
Moreover, these sets are disjoint since $\beta<\gamma$.
Thus, according to \Cref{thm:ast-def-sep-cvx-sets},
any two disjoint convex sets are included in disjoint astral dual
halfspaces.

The closure of the convex hull of any set $S\subseteq\Rn$ is equal to
the intersection of all the closed halfspaces that include it
(\Cref{pr:con-int-halfspaces}\ref{roc:cor11.5.1}).
\indexg{convex hull (standard)!intersection of dual halfspaces@as intersection of dual halfspaces|(}%
\indexg{dual halfspaces, astral!convex sets as intersection of|(}%
When working with astral dual halfspaces, the same holds simply
for the convex hull of the set, without taking its closure, as we
show next.
\indexg{convex sets (standard)!intersection of dual halfspaces@as intersection of dual halfspaces|(}%
In particular, this means that every convex set in $\Rn$ is equal to
the intersection of all the astral dual halfspaces that include
it.
This and other consequences of
\Cref{thm:ast-def-sep-cvx-sets}
are summarized next:

\begin{corollary}   \label{cor:ast-def-sep-conseqs}
  Let $S, U, V\subseteq\Rn$.
  \begin{letter-compact}
  \item   \label{cor:ast-def-sep-conseqs:a}
    The following are equivalent:
    \begin{roman-compact}
    \item   \label{cor:ast-def-sep-conseqs:a1}
      $S$ is convex.
    \item   \label{cor:ast-def-sep-conseqs:a2}
      $S$ is equal to the intersection of some collection of
      astral dual halfspaces.
    \item   \label{cor:ast-def-sep-conseqs:a3}
      $S$ is equal to the intersection of all
      astral dual halfspaces that include $S$.
    \end{roman-compact}
  \item   \label{cor:ast-def-sep-conseqs:b}
    The convex hull of $S$, $\conv{S}$, is equal to the intersection
    of all astral dual halfspaces that include $S$.
  \item   \label{cor:ast-def-sep-conseqs:c}
    $U$ and $V$ are strongly dually separated
    if and only if
    $(\conv{U})\cap(\conv{V})=\emptyset$.
  \end{letter-compact}
\end{corollary}

\begin{proof}
~

\begin{proof-parts}
\pfpart{Part~(\ref{cor:ast-def-sep-conseqs:a}):}
That (\ref{cor:ast-def-sep-conseqs:a3}) implies
(\ref{cor:ast-def-sep-conseqs:a2}) is immediate.
That (\ref{cor:ast-def-sep-conseqs:a2}) implies
(\ref{cor:ast-def-sep-conseqs:a1}) follows from
Propositions~\ref{pr:e1}(\ref{pr:e1:b})
and~\ref{pr:ast-def-hfspace-convex}(\ref{i:dual-hfspace-conv}).

To show that (\ref{cor:ast-def-sep-conseqs:a1}) implies
(\ref{cor:ast-def-sep-conseqs:a3}), suppose that $S$ is convex,
and
let $U$ be the intersection of all astral dual halfspaces that
include $S$;
we aim to show $S=U$.
We assume $S$ is not empty since otherwise
$U=\emptyset=S$.
Clearly, $S\subseteq U$.
Suppose $\ww\in\Rn\wo S$.
Then by \Cref{thm:ast-def-sep-cvx-sets},
applied to $S$ and $\{\ww\}$,
there exists $\xbar\in\extspace$, $\rr\in\Rn$ and $\beta\in\R$
such that
$\xbar\cdot(\ww-\rr)>\beta$,
while
$\xbar\cdot(\uu-\rr)\leq\beta$
for all $\uu\in S$
(implying $\xbar\neq\zero$ since $S$ is nonempty).
Thus, $S$ is in the astral dual halfspace
$\ahfsp{\xbar}{\rr}{\beta}$
(as defined in Eq.~\ref{eqn:ast-def-hfspace-defn}),
but $\ww$ is not.
Therefore, $\ww\not\in U$,
proving $S=U$.%
\indexg{convex sets (standard)!intersection of dual halfspaces@as intersection of dual halfspaces|)}

\pfpart{Part~(\ref{cor:ast-def-sep-conseqs:b}):}
For any astral dual halfspace $H\subseteq\Rn$, since $H$ is convex
(by~\Cref{pr:ast-def-hfspace-convex}\ref{i:dual-hfspace-conv}),
$S\subseteq H$ if and only if $\conv{S}\subseteq H$
(by \Cref{pr:conhull-prop}\ref{pr:conhull-prop:aa}).
Therefore, the intersection of all astral dual halfspaces that include
$S$ is the same as the intersection of all astral dual halfspaces that
include $\conv{S}$, which, by
part~(\ref{cor:ast-def-sep-conseqs:a}),
is simply $\conv{S}$.%
\indexg{dual halfspaces, astral!convex sets as intersection of|)}%
\indexg{convex hull (standard)!intersection of dual halfspaces@as intersection of dual halfspaces|)}

\pfpart{Part~(\ref{cor:ast-def-sep-conseqs:c}):}
The proof is analogous to that of
\Cref{cor:sep-cvx-sets-conseqs}(\ref{cor:sep-cvx-sets-conseqs:c}).
\qedhere
\end{proof-parts}
\end{proof}

\indexg{hemispace|(}%
A set $H\subseteq\Rn$ is said to be a \emph{hemispace}
if both $H$ and its complement, $\Rn\wo H$,
are convex.
For example, ordinary open and closed halfspaces are hemispaces.
In fact, except for $\emptyset$ and $\Rn$,
a set is a hemispace if and only if it is an astral dual
halfspace, as we show next as a
consequence
of \Cref{thm:ast-def-sep-cvx-sets}.
Closely related equivalences for hemispaces were proved by
\idxmartinezsinger\citet[Theorem~1.1]{hemispaces}.

\begin{theorem}  \label{cor:hemi-is-ast-dual-halfspace}
  Let $H\subseteq\Rn$.
  Then the following are equivalent:
  \begin{letter-compact}
  \item  \label{cor:hemi-is-ast-dual-halfspace:a}
    $H$ is a hemispace
  \item  \label{cor:hemi-is-ast-dual-halfspace:b}
    $H$ is either $\emptyset$ or $\Rn$ or an astral dual halfspace.
  \end{letter-compact}
\end{theorem}

\begin{proof}
  ~

\begin{proof-parts}
\pfpart{%
  (\ref{cor:hemi-is-ast-dual-halfspace:b})
  $\Rightarrow$
  (\ref{cor:hemi-is-ast-dual-halfspace:a}):
}
If $H$ is $\emptyset$ or $\Rn$, then it is clearly a hemispace.
If $H$ is an astral dual halfspace, then both $H$ and its complement are convex by \Cref{pr:ast-def-hfspace-convex}(\ref{i:dual-hfspace-conv}), so
$H$ is a hemispace as well.

\pfpart{%
  (\ref{cor:hemi-is-ast-dual-halfspace:a})
  $\Rightarrow$
  (\ref{cor:hemi-is-ast-dual-halfspace:b}):
}
Suppose $H$ is a hemispace and is neither $\emptyset$ nor $\Rn$.
Then both $H$ and its complement, $H^c=\Rn\wo H$,
are convex and nonempty.
Since they are also disjoint, by
\Cref{thm:ast-def-sep-cvx-sets}, they are strongly dually
separated, so that
\[
  \sup_{\uu\in H} \xbar\cdot(\uu-\rr)
  <
  \beta
  <
  \inf_{\vv\in H^c} \xbar\cdot(\vv-\rr),
\]
for some
$\xbar\in\extspace\,\wo\{\zero\}$,
$\rr\in\Rn$, and $\beta\in\R$.
It follows that $H$ is included in the astral dual halfspace
$H'=\ahfsp{\xbar}{\rr}{\beta}$,
while its complement, $H^c$, is in the complement of $H'$.
That is, $\Rn\wo H = H^c \subseteq \Rn\wo H'$,
which implies that $H'\subseteq H$.
Thus, $H=H'$, so $H$ is an astral dual halfspace.%
\indexg{hemispace|)}%
\indexg{dual halfspaces, astral|)}%
\qedhere
\end{proof-parts}
\end{proof}

Suppose $U,V\subseteq\Rn$ are strongly dually separated.
Then by
\Cref{thm:ast-str-sep-equiv-no-trans},
there exists $\xbar\in\extspace$ such that
$\xbar\cdot(\uu-\vv)<0$ for all $\uu\in U$, $\vv\in V$,
which in turn implies that
$\limray{\xbar}\cdot(\uu-\vv)=-\infty$.
Since $\limray{\xbar}$ is an icon, we can further write it in the form
$\limray{\xbar}=\WW\omm$, for some $\WW\in\Rnk$,
so that
$\WW \omm\cdot(\uu-\vv) = -\infty$
for all $\uu\in U$, $\vv\in V$.

\indexg{lexicographic separation|(}%
\indexg{separation, lexicographic|(}%
We remark that
this latter condition can be expressed in terms of
{lexicographic ordering}.
Specifically, for any two vectors $\aaa,\bb\in\Rk$,
we say that $\aaa$ is
\indexg{lexicographic order}%
\emph{lexicographically less than}
$\bb$, written $\aaa \lexless \bb$,
if $\aaa\neq\bb$ and if, in the first component where they differ, $\aaa$ is
less than $\bb$;
that is, if
for some $j\in\{1,\ldots,k\}$,
$a_j < b_j$ and
$a_i = b_i$ for $i=1,\ldots,j-1$.
Then it can be checked
(for instance, using the Case Decomposition Lemma~\ref{lemma:case})
that the condition above, that
$\WW \omm\cdot(\uu-\vv) = -\infty$,
holds if and only if
$\trans{\WW}(\uu-\vv) \lexless \zero$,
or equivalently,
$\trans{\WW}\uu \lexless \trans{\WW}\vv$.
In this sense, since this holds for all $\uu\in U$,
and all $\vv\in V$,
the sets $U$ and $V$ are
\emph{lexicographically separated}.
Thus, \Cref{thm:ast-def-sep-cvx-sets}
shows that any two disjoint, convex sets in $\Rn$
can be lexicographically separated.

This latter fact had been shown previously by
\idxmartinezsinger\citet[Theorem~2.1]{lexicographic_separation}.
Indeed, the preceding argument can be used to show that our notion of
strong dual separation is essentially equivalent to lexicographic
separation, so in this sense,
\Cref{thm:ast-def-sep-cvx-sets}
and
\Cref{cor:ast-def-sep-conseqs}(\ref{cor:ast-def-sep-conseqs:c})
were already known in their work.%
\indexg{lexicographic separation|)}%
\indexg{separation, lexicographic|)}

\indexg{dual halfspaces, astral!homogeneous|(}%
\indexg{homogeneous (halfspace or hyperplane)!dual halfspace|(}%
\indexg{separation, strong astral dual!convex cone@from convex cone|(}%
We consider next separation involving convex cones.
We say that an astral dual halfspace, as in
\eqref{eqn:ast-def-hfspace-defn},
is \emph{homogeneous} if both $\rr=\zero$ and $\beta=0$,
that is, if it has the form
\begin{equation}    \label{eqn:ast-def-homo-hfspace-defn}
   \ahfsp{\xbar}{\zero}{0}
   =
   \Braces{\uu\in\Rn : \xbar\cdot\uu\leq 0}.
\end{equation}
When separating a convex cone $K\subseteq\Rn$ from a disjoint convex
set $V\subseteq\Rn$,
we show next that there always exists a homogeneous
astral dual halfspace that includes $K$ while excluding all of $V$.%
\indexg{dual halfspaces, astral!homogeneous|)}%
\indexg{homogeneous (halfspace or hyperplane)!dual halfspace|)}

\begin{theorem}    \label{thm:ast-def-sep-cone-sing}
  Let $K\subseteq\Rn$ be a convex cone,
  and let $V\subseteq\Rn$ be nonempty and convex with
  $K\cap V=\emptyset$.
  Then there exists $\xbar\in\extspace$ such that
  \[
     \sup_{\uu\in K} \xbar\cdot\uu
     \leq
     0
     <
     \inf_{\vv\in V} \xbar\cdot\vv.
  \]
\end{theorem}

Before proving the theorem, we give a lemma that will be used in its
proof.
The lemma supposes that we are given some $\xbar\in\extspace$ and sets
$U,W\subseteq\Rn$ such
that $\xbar\cdot\uu\leq 0$ for all $\uu\in U$ and
$\xbar\cdot\ww<0$ for all $\ww\in W$.
In particular, it is possible that $\xbar\cdot\ww=-\infty$ for all
$\ww\in W$.
The lemma shows that there then must exist some point
$\ybar\in\extspace$ that preserves all these same properties with the
additional condition that $\ybar\cdot\ww\in\R$ for at least one point
$\ww\in W$.

\begin{lemma}   \label{lem:sep-one-pt-finite}
  Let $U,W\subseteq\Rn$, with $W$ nonempty,
  and let $\xbar\in\extspace$.
  Suppose $\xbar\cdot\uu\leq 0$ for all $\uu\in U$,
  and that $\xbar\cdot\ww < 0$ for all $\ww\in W$.
  Then there exists a point $\ybar\in\extspace$ such that
  $\ybar\cdot\uu\leq 0$ for all $\uu\in U$,
  $\ybar\cdot\ww < 0$ for all $\ww\in W$,
  and such that there exists some $\what\in W$
  with $\ybar\cdot\what \in \R$.
\end{lemma}

\begin{proof}
Since $\limray{\xbar}$ is an icon, we can write it as
$\limray{\xbar}=\limrays{\vv_1,\ldots,\vv_k}$
for some $\vv_1,\ldots,\vv_k\in\Rn$
(Propositions~\ref{pr:icon-equiv}\ref{pr:icon-equiv:a}\ref{pr:icon-equiv:c}
and~\ref{pr:i:8}\ref{pr:i:8-infprod}).
For $i=0,\ldots,k$, let
$ \ebar_i = \limrays{\vv_1,\ldots,\vv_i} $.

For all $\uu\in U\cup W$ and all $i\in\{0,\ldots,k\}$,
we claim $\ebar_i\cdot\uu\leq 0$.
Otherwise, if $\ebar_i\cdot\uu>0$ then
$\ebar_i\cdot\uu=+\infty$ (by
\Cref{pr:icon-equiv}\ref{pr:icon-equiv:a}\ref{pr:icon-equiv:b},
since $\ebar_i$ is an icon),
implying
$\limray{\xbar}\cdot\uu
=
(\ebar_i\plusl\limrays{\vv_{i+1},\ldots,\vv_k})\cdot\uu
=
+\infty
$,
so that $\xbar\cdot\uu>0$, a contradiction.

Let $j\in\{1,\ldots,k\}$ be the largest index for which
there exists $\what\in W$ such that $\ebar_{j-1}\cdot\what=0$.
In particular, this means that $\ebar_j\cdot\ww\neq 0$ for all
$\ww\in W$.
Note that $1\leq j\leq k$ since for all $\ww\in W$,
$\ebar_0\cdot\ww=\zero\cdot\ww=0$
and
$\ebar_k\cdot\ww=\limray{\xbar}\cdot\ww<0$
(since $\xbar\cdot\ww<0$).

Let $\ybar=\ebar\plusl\vv$ where $\ebar=\ebar_{j-1}$ and $\vv=\vv_j$.
Note that $\limray{\ybar}=\ebar\plusl\limray{\vv}=\ebar_j$
(\Cref{pr:i:8}\ref{pr:i:8d}).
Therefore, for all $\uu\in U\cup W$,
$\limray{\ybar}\cdot\uu=\ebar_j\cdot\uu\leq 0$,
as argued above, implying
$\ybar\cdot\uu\leq 0$.

Further, by $j$'s choice, for all $\ww\in W$,
$\limray{\ybar}\cdot\ww=\ebar_j\cdot\ww\neq 0$,
so $\ybar\cdot\ww\neq 0$, implying
$\ybar\cdot\ww<0$.

Finally, $\what$ and $j$ were chosen so that $\ebar\cdot\what=0$.
Therefore,
$\ybar\cdot\what=\ebar\cdot\what\plusl\vv\cdot\what=\vv\cdot\what\in\R$,
completing the proof.
\end{proof}

\begin{proof}[Proof of \Cref{thm:ast-def-sep-cone-sing}]
The sets $K$ and $V$ are nonempty, convex and disjoint,
so by Theorems~\ref{thm:ast-str-sep-equiv-no-trans}
and~\ref{thm:ast-def-sep-cvx-sets}, there exists $\zbar\in\extspace$
such that
$\zbar\cdot(\uu-\vv)<0$ for all $\uu\in K$ and $\vv\in V$.
Since $\zero\in K$, this implies $\zbar\cdot(-\vv)<0$
for all $\vv\in V$.
Therefore, by \Cref{lem:sep-one-pt-finite}
(applied with $\xbar$, $U$ and $W$, as they appear in the lemma,
set to $\zbar$, $K-V$ and $-V$),
there exists $\ybar\in\extspace$ and $\vhat\in V$ such that,
for all $\uu\in K$ and $\vv\in V$,
$\ybar\cdot(\uu-\vv)\leq 0$
and
$\ybar\cdot\vv>0$,
and furthermore,
$\ybar\cdot\vhat\in\R$.

For $\uu\in K$, we thus have that
$\ybar\cdot\uu$ and $-\ybar\cdot\vhat$ are summable,
so
$\ybar\cdot\uu-\ybar\cdot\vhat=\ybar\cdot(\uu-\vhat)\leq 0$
(by \Cref{pr:i:1} and since $\vhat\in V$).
Therefore, for all $\uu\in K$,
$\ybar\cdot\uu\leq \ybar\cdot\vhat$.
Since $K$ is a cone, this further implies, for all
$\lambda\in\Rpos$, that $\lambda\uu\in K$ so that
\begin{equation*}   \label{eq:thm:ast-def-sep-cone-sing:1}
  \lambda\ybar\cdot\uu
  =
  \ybar\cdot(\lambda\uu)
  \leq
  \ybar\cdot\vhat.
\end{equation*}
This means that $\ybar\cdot\uu\leq 0$;
otherwise,
if $\ybar\cdot\uu>0$, then
the left-hand side of
this equation
could be made arbitrarily large as $\lambda$ becomes large,
while the right-hand side remains bounded
(since $\ybar\cdot\vhat\in\R$).

Finally, let $\xbar=\limray{\ybar}$.
Then for all $\uu\in K$,
$\xbar\cdot\uu=\limray{\ybar}\cdot\uu\leq 0$
(since, as just argued, $\ybar\cdot\uu\leq 0$).
And for all $\vv\in V$,
$\ybar\cdot\vv>0$,
implying
$\xbar\cdot\vv=\limray{\ybar}\cdot\uu=+\infty$,
completing the proof.%
\indexg{separation, strong astral dual!convex cone@from convex cone|)}%
\end{proof}

Here are some consequences of
\Cref{thm:ast-def-sep-cone-sing},
analogous to
\Cref{cor:ast-def-sep-conseqs}.
\indexg{dual halfspaces, astral!cones as intersection of|(}%
\indexg{cones (standard)!intersection of dual halfspaces@as intersection of dual halfspaces|(}%
In particular, every (standard) convex cone in $\Rn$ is equal to the
intersection of all homogeneous astral dual halfspaces that include
it.

\begin{corollary}   \label{cor:ast-def-sep-cone-conseqs}
  Let $S\subseteq\Rn$.
  \begin{letter-compact}
  \item   \label{cor:ast-def-sep-cone-conseqs:a}
    The following are equivalent:
    \begin{roman-compact}
    \item   \label{cor:ast-def-sep-cone-conseqs:a1}
      $S$ is a convex cone.
    \item   \label{cor:ast-def-sep-cone-conseqs:a2}
      $S$ is equal to the intersection of some collection of
      homogeneous astral dual halfspaces.
    \item   \label{cor:ast-def-sep-cone-conseqs:a3}
      $S$ is equal to the intersection of all homogeneous
      astral dual halfspaces that include $S$.
    \end{roman-compact}
  \item   \label{cor:ast-def-sep-cone-conseqs:b}
\indexg{conic hull (standard)!intersection of dual halfspaces@as intersection of dual halfspaces|(}%
    The conic hull of $S$, $\cone{S}$, is equal to the intersection
    of all homogeneous astral dual halfspaces that include $S$.
  \end{letter-compact}
\end{corollary}

\begin{proof}
The proof is analogous to that of
\Cref{cor:ast-def-sep-conseqs},
after making the additional observation that every homogeneous astral
dual halfspace, as in
\eqref{eqn:ast-def-homo-hfspace-defn},
is itself a convex cone
(\Cref{pr:ast-def-hfspace-convex}\ref{i:dual-hfspace-cone}).%
\indexg{dual halfspaces, astral!cones as intersection of|)}%
\indexg{cones (standard)!intersection of dual halfspaces@as intersection of dual halfspaces|)}%
\indexg{conic hull (standard)!intersection of dual halfspaces@as intersection of dual halfspaces|)}%
\end{proof}

\chapter{Representational closure and astral polarity}
\chaptermark{Representational closure, astral polarity}

In this chapter, we look more closely at how to extend standard
convex cones and standard conic operations to astral space.  We
introduce an operation called representational closure, and
use it to characterize the astral conic hull of any set $S\subseteq\Rn$.
We also show how
the notion of a polar cone extends to the astral setting.
Whereas standard polarity induces a one-to-one correspondence
on the set of closed convex cones in $\Rn$, we will see that astral polarity induces a one-to-one correspondence
between the set of closed convex astral cones in $\eRn$ and the set of {all}
convex cones in $\Rn$ (not just the closed ones).

\section{Extending a standard convex cone to astral space}
\label{sec:cones:astrons}

We begin by considering natural ways in which a standard convex cone
$K\subseteq\Rn$ can be extended to form a convex astral cone.
One simple way is to take its closure, $\Kbar$,
which, by \Cref{thm:out-conic-h-is-closure},
is the same as $\oconich{K}$, its outer conic hull.
By
\Cref{cor:sep-ast-cone-conseqs}(\ref{cor:sep-ast-cone-conseqs:a}),
this is the smallest closed convex astral cone that includes $K$.

An alternative way of extending $K$
is to take its astral conic
hull, $\acone{K}$, which is the smallest convex astral cone (whether
closed or not) that includes $K$.
We next study the structure of this set and prove that every point in it can be written in a form that uses only vectors from $K$ in its representation.

To make this precise, we define the following set:

\begin{definition}
\indexg{representational closure|(}%
Let $K\subseteq\Rn$ be a convex cone.
The
\emph{representational closure of $K$}, denoted $\repcl{K}$, is
the set of all points in $\extspace$ that can be represented using
vectors in $K$;
that is,
\begin{equation}  \label{eq:rep-cl-defn}
\indexm{k delta}{$\repcl{K}$}{representational closure}%
  \repcl{K} =
    \bigBraces{
         \limrays{\vv_1,\ldots,\vv_k}\plusl \qq :\:
         \vv_1,\ldots,\vv_k,\qq\in K,\, k\geq 0
           }.
\end{equation}
\end{definition}

Note that the points appearing in
\eqref{eq:rep-cl-defn} need {not} be canonical representations.
For example, in this notation,
\Cref{cor:h:1} states exactly that
$\extspace=\repcl{(\Rn)}$.%
\indexg{representational closure|)}

\indexg{representational closure!astral conic hull as|(}%
\indexg{conic hull, astral!representational closure@as representational closure|(}%
As suggested above,
the astral conic hull of any standard convex cone
turns out to be the same as its representational closure.
More generally, as shown next,
the astral conic hull of any set in $\Rn$ is equal
to the representational closure of its (standard) conic hull.

\begin{theorem}   \label{thm:acone-is-repclos}
  Let $S\subseteq\Rn$.
  Then
  \[
     \acone{S} = \repcl{(\cone{S})}.
  \]
  Consequently, $\cone{S} = (\acone{S})\cap\Rn$.
\end{theorem}

Before proving the theorem, we give a key lemma that is central to its
proof.
The lemma shows that the topological closure of a finitely generated
cone
(meaning $\cone{V}$ for some finite set $V\subseteq\Rn$)
is always included in its representational closure.

\begin{lemma}  \label{lem:conv-lmset-char:a}
  Let $V$ be a finite subset of $\Rn$.
  Then
  $\clcone{V}\subseteq \repcl{(\cone V)}$.
\end{lemma}

\begin{proof}
If $V$ is empty, then
$\cone V = \{\zero\}$, so
$\clcone{V} = \{\zero\} = \repcl{(\cone V)}$.

For $V$ nonempty,
we prove the result by induction on astral rank.
More precisely,
we prove, by induction on $k=0,1,\ldots,n$, that for all
$\xbar\in\extspace$ and for every finite, nonempty set $V\subseteq\Rn$,
if $\xbar$ has astral rank at most $k$ and if $\xbar\in\clcone{V}$
then $\xbar\in\repcl{(\cone V)}$.

For the base case that $k=0$,
suppose $\xbar$ has astral rank $0$ and that $\xbar\in\clcone{V}$
where $V\subseteq\Rn$ is finite and nonempty.
Then $\xbar$ is equal to some $\qq\in\Rn$.
Since $\xbar\in\clcone{V}$, there exists a sequence $\seq{\xx_t}$
in $\cone V$ converging to $\xbar=\qq$.
Thus, $\qq\in \cone V$ since $\cone V$ is closed in $\Rn$
by \Cref{roc:thm19.1}(\ref{roc:thm19.1:b},\ref{roc:thm19.1:c}).
Therefore, $\xbar\in\repcl{(\cone V)}$.

For the inductive step, suppose $\xbar$ has astral rank $k>0$ and
that $\xbar\in\clcone{V}$
for some $V=\{\zz_1,\ldots,\zz_m\}\subseteq\Rn$ with $m\geq 1$.
Furthermore, assume inductively
that the claim holds for all points with astral rank strictly less
than $k$.
Then
by \Cref{pr:h:6},
we can write
$\xbar=\limray{\ww}\plusl\xbarperp$
where $\ww$ is $\xbar$'s dominant direction, and
$\xbarperp$ is the
projection of $\xbar$ orthogonal to $\ww$.

Since $\xbar\in\clcone{V}$, there exists a sequence $\seq{\xx_t}$
in $\cone V$ that converges to $\xbar$.
Since $\norm{\xx_t}\rightarrow+\infty$
(by \Cref{pr:seq-to-inf-has-inf-len}),
we can discard all elements of the sequence with $\xx_t=\zero$, of
which there can be at most finitely many.
Then by
\Crefequiv{thm:dom-dir}{thm:dom-dir:a}{thm:dom-dir:c},
\[
  \ww = \lim \frac{\xx_t}{\norm{\xx_t}}.
\]
Since $\cone V$ is a cone,
$\xx_t/\norm{\xx_t}\in\cone V$, for all $t$,
which implies that $\ww\in \cone V$
since $\cone V$ is closed in $\Rn$.

By \Cref{pr:lin-decomp-rel-vecs},
for each $t$, we have
$\xx_t = b_t \ww + \xperpt$
where $b_t=\xx_t\cdot\ww$  %
and $\xperpt$ is the projection of $\xx_t$ orthogonal to $\ww$ so that $\xperpt\cdot\ww=0$.
Since $\xx_t\rightarrow\xbar$ and since $\ww$ is $\xbar$'s dominant
direction,
$b_t=\xx_t\cdot\ww\rightarrow \xbar\cdot\ww = +\infty$
(by Theorems~\ref{thm:i:1}\ref{thm:i:1c}
and~\ref{thm:dom-dir}\ref{thm:dom-dir:a}\ref{thm:dom-dir:b}).
Thus, $b_t$ must be positive for all but finitely many values of $t$;
by discarding these, we can assume that $b_t>0$ for all $t$.

By \Cref{pr:scc-cone-elts}(\ref{pr:scc-cone-elts:b}),
each $\xx_t$, being in $\cone V$, is a conic
combination of $\zz_1,\ldots,\zz_m$, the points comprising $V$.
Therefore, each point
$\xperpt=\xx_t-b_t\ww$ is a conic combination of the
points in the expanded set
\[
  V' = V \cup \{ -\ww \} = \{ \zz_1,\ldots,\zz_m,\,-\ww \}.
\]
In other words, $\xperpt\in\cone V'$, for all $t$.
Furthermore, $\xbarperp$ must be in the astral closure of this cone,
$\clcone{V'}$,
since
$\xperpt\rightarrow\xbarperp$
(by \Cref{pr:h:5}\ref{pr:h:5b}).

Since $\xbarperp$ has strictly lower astral rank than $\xbar$
(by \Cref{pr:h:6}), we can apply our
inductive assumption which implies that
\[ \xbarperp = [\vv'_1,\ldots,\vv'_{k'}]\omm \plusl \qq' \]
for some $\qq',\vv'_1,\ldots,\vv'_{k'}\in \cone V'$.
For $i=1,\ldots,k'\negKern$,
since $\vv'_i\in \cone V'$,
we can write $\vv'_i$ as a conic combination over the
points in $V'$ so that
$\vv'_i=\vv_i - a_i \ww$
where
\[  \vv_i = \sum_{j=1}^m c_{ij} \zz_j \]
for some $a_i\geq 0$, $c_{ij}\geq 0$, $j=1,\ldots,m$.
Note that $\vv_i\in\cone V$.
Similarly we can write
$\qq'=\qq-b\ww$ for some $\qq\in \cone V$ and $b\geq 0$.

Thus,
\begin{align*}
  \xbar
  =
  \limray{\ww} \plusl \xbarperp
  =
  \limray{\ww}
  \plusl
  [\vv'_1,\ldots,\vv'_{k'}]\omm \plusl \qq'
  &=
  [\ww,\, \vv'_1,\ldots,\vv'_{k'}]\omm \plusl \qq'
  \\
  &=
  [\ww,\, \vv_1,\ldots,\vv_{k'}]\omm \plusl \qq.
\end{align*}
The last equality follows from
the Push Lemma~\ref{pr:g:1},
combined with \Cref{prop:pos-upper:descr},
since
$\vv_i=\vv'_i+a_i\ww$ for $i=1,\ldots,k'$,
and since
$\qq=\qq'+b\ww$.
Since $\ww$, $\qq$, and $\vv_1,\ldots,\vv_{k'}$ are all in
$\cone V$,
this shows that $\xbar\in\repcl{(\cone V)}$,
completing the induction and the proof.
\end{proof}

\begin{proof}[Proof of \Cref{thm:acone-is-repclos}]
Let $K=\cone{S}$.
We first prove that $\repcl{K}\subseteq\acone{S}$.

The set $(\acone{S})\cap\Rn$ includes $S$ and is a convex cone
(by \Cref{pr:ast-cone-is-naive} and since $\acone{S}$ and
$\Rn$ are both convex).
Therefore, it must include $S$'s conic hull, so
\[
  K=\cone{S}\subseteq(\acone{S})\cap\Rn\subseteq\acone{S}.
\]
This further implies that $\lmset{K}\subseteq\acone{S}$
(again by \Cref{pr:ast-cone-is-naive}).

Since $\acone{S}$ is a convex astral cone, it is closed under
sequential sum (by \Cref{thm:ast-cone-is-cvx-if-sum}).
That is, if $\xbar,\ybar\in\acone{S}$, then
$\xbar\seqsum\ybar\subseteq\acone{S}$, and so,
in particular, $\xbar\plusl\ybar\in\acone{S}$
(by \Cref{cor:seqsum-conseqs}\ref{cor:seqsum-conseqs:a}).
Therefore, $\acone{S}$ is closed under leftward addition.

Let $\xbar\in\repcl{K}$.
Then $\xbar=\limray{\vv_1}\plusl\dotsb\plusl\limray{\vv_k}\plusl\qq$
for some $\qq,\vv_1,\ldots,\vv_k\in K$.
By the foregoing, $\qq$ and
$\limray{\vv_i}$, for $i=1,\ldots,k$,
are all in $\acone{S}$,
implying $\xbar$ is as well.
Therefore, $\repcl{K}\subseteq\acone{S}$.

For the reverse inclusion, let $\xbar\in\acone{S}$.
Then by \Cref{thm:ast-conic-hull-union-fin},
$\xbar\in\oconich{V}$ for some finite subset $V\subseteq S$.
We then have
\[
  \xbar
  \in
  \oconich{V}
  =
  \clcone{V}
  \subseteq
  \repcl{(\cone V)}
  \subseteq
  \repcl{(\cone S)}
  =
  \repcl{K}.
\]
The first equality is by
\Cref{thm:out-conic-h-is-closure}.
The second inclusion is by \Cref{lem:conv-lmset-char:a}.
The third inclusion is because $V\subseteq S$, a relation that is
preserved when generating cones and taking representational closure.
Thus, $\acone{S} \subseteq \repcl{K}$,
completing the proof.
\end{proof}

Applying \Cref{thm:acone-is-repclos} to a convex cone $K\subseteq\Rn$ yields a characterization of the astral conic hull of $K$ as its representational closure: $\acone{K}=\repcl{K}$. On the other hand, as we saw earlier, the outer conic hull of $K$ is equal to its topological closure: $\oconich{K}=\Kbar$. Moreover, we always have $\acone{K}\subseteq\oconich{K}$.
We summarize these relations in the next  \namecref{cor:a:1}.
Later, in \Cref{sec:repcl-and-polyhedral},
we give necessary and sufficient conditions for
when the two types of astral conic hull (and hence the two types of closure) coincide, that is, for when $\repcl{K}=\Kbar$.

\begin{corollary}  \label{cor:a:1}
  Let $K\subseteq\Rn$ be a convex cone.
  Then $\repcl{K}$ is a convex astral cone, and
  \begin{align*}
    K
  &\subseteq
    \acone{K}
    =
    \repcl{K}
    =
    \conv(\lmset{K})
\\
  &\subseteq
    \oconich{K}
    =
    \Kbar.
  \end{align*}
\end{corollary}

\begin{proof}
Since $K$ is already a convex cone,
$\cone K = K$.
The first equality therefore follows from
\Cref{thm:acone-is-repclos}
(which also shows that $\repcl{K}$ is a convex astral cone).
The second equality is from
\Cref{thm:acone-char}(\ref{thm:acone-char:a})
(and since $\zero\in K$).
The second inclusion is by
\Cref{pr:acone-hull-props}(\ref{pr:acone-hull-props:d}),
and the final equality is by
\Cref{thm:out-conic-h-is-closure}.%
\indexg{representational closure!astral conic hull as|)}%
\indexg{conic hull, astral!representational closure@as representational closure|)}%
\end{proof}

\indexg{representational closure!convex hull of astrons and|(}%
\indexg{convex hull, astral!astrons@of astrons|(}%
In the remainder of this section, we show how the concept of  representational closure can be applied to characterize convex hulls of sets of astrons.

Any set of astrons can be expressed as $\lmset{S}$ for some $S\subseteq\Rn$.
If $\zero\in S$ then by \Cref{thm:acone-char}(\ref{thm:acone-char:a}),
$\conv(\lmset{S})$ is equal to the astral conic hull of $S$, $\acone{S}$, which is in turn
equal to the representational closure $\repcl{(\cone{S})}$ by \Cref{thm:acone-is-repclos}.
Thus, \Cref{thm:acone-is-repclos} succinctly characterizes the
representations of all points in $\conv(\lmset{S})$.

The next theorem provides a similarly simple characterization
for the convex hull of any set of astrons, whether or not that set includes the origin.
We show that if $\zero\in\conv{S}$, then
$\conv(\lmset{S})$ is equal to $\repcl{(\cone{S})}$.
Otherwise, if $\zero\not\in\conv{S}$, then
$\conv(\lmset{S})$ consists of just the infinite points in
$\repcl{(\cone{S})}$.

\begin{theorem}  \label{thm:conv-lmset-char}
  Let $S\subseteq\Rn$.
  \begin{letter-compact}
  \item      \label{thm:conv-lmset-char:a}
    If $\zero\in\conv{S}$, then
    $ \conv(\lmset{S}) = \repcl{(\cone{S})} $.
  \item      \label{thm:conv-lmset-char:b}
    If $\zero\not\in\conv{S}$, then
    $ \conv(\lmset{S}) = \repcl{(\cone{S})} \setminus \Rn$.
  \end{letter-compact}
\end{theorem}

Before proving the theorem, we first give a simple lemma to be
used in its proof
(stated in a bit more generality than is needed at this point).

\begin{lemma}    \label{lem:conv-lmset-disj-rn}
  Let $S\subseteq\extspace$.
  \begin{letter-compact}
  \item    \label{lem:conv-lmset-disj-rn:a}
    If $\zero\in\conv{S}$, then
    $\zero\in\conv(\lmset{S})$.
  \item    \label{lem:conv-lmset-disj-rn:b}
    If $\zero\not\in\conv{S}$, then
    $\conv(\lmset{S})\cap\Rn=\emptyset$.
  \end{letter-compact}
\end{lemma}

\begin{proof}
  ~

\begin{proof-parts}
\pfpart{Part~(\ref{lem:conv-lmset-disj-rn:a}):}
Suppose $\zero\in\conv{S}$.
Then by \Cref{thm:convhull-of-simpices},
$\zero\in\ohull{V}$ for some finite subset $V\subseteq S$.
Let $\uu\in\Rn$.
Then
\[
  0
  =
  \zero\cdot\uu
  \leq
  \max\,\regBraces{\xbar\cdot\uu :\: \xbar\in V}
\]
by
\Cref{pr:ohull-simplify}(\ref{pr:ohull-simplify:a},\ref{pr:ohull-simplify:b}).
Thus, for some $\xbar\in V$, $\xbar\cdot\uu\geq 0$,
implying $\limray{\xbar}\cdot\uu\geq 0$.
Therefore,
\[
  \zero\cdot\uu
  =
  0
  \leq
  \max\,\regBraces{\limray{\xbar}\cdot\uu :\: \xbar\in V},
\]
so $\zero\in\ohull(\lmset{V})\subseteq\conv(\lmset{S})$
(again by
\Cref{pr:ohull-simplify}\ref{pr:ohull-simplify:b}\ref{pr:ohull-simplify:a}
and
\Cref{thm:convhull-of-simpices}).

\pfpart{Part~(\ref{lem:conv-lmset-disj-rn:b}):}
We prove the contrapositive with a proof similar to the
preceding part.
Suppose
there exists a point $\qq\in\conv(\lmset{S})\cap\Rn$;
we aim to show this implies $\zero\in\conv{S}$.
By \Cref{thm:convhull-of-simpices},
$\qq\in \ohull(\lmset{V})$ for some finite subset $V\subseteq S$.
Let $\uu\in\Rn$.
Then
by
\Cref{pr:ohull-simplify}(\ref{pr:ohull-simplify:a},\ref{pr:ohull-simplify:b}),
\begin{equation}  \label{eq:h:8a}
   \qq\cdot\uu
   \leq
   \max\,\regBraces{\limray{\xbar}\cdot\uu :\: \xbar\in V}.
\end{equation}
This implies that
\begin{equation*}  \label{eq:h:8b}
  \zero\cdot\uu
  =
  0
  \leq
  \max\,\regBraces{\xbar\cdot\uu :\: \xbar\in V},
\end{equation*}
since otherwise, if $\xbar\cdot\uu<0$ for all $\xbar\in V$, then the
right-hand side of \eqref{eq:h:8a} would be equal to $-\infty$, which
is impossible since $\qq\cdot\uu\in\R$.
Again applying
\Cref{pr:ohull-simplify}(\ref{pr:ohull-simplify:b},\ref{pr:ohull-simplify:a})
and \Cref{thm:convhull-of-simpices},
this shows that
$\zero\in\ohull{V}\subseteq\conv{S}$.
\qedhere
\end{proof-parts}
\end{proof}

\begin{proof}[Proof of \Cref{thm:conv-lmset-char}]
  ~

\begin{proof-parts}
\pfpart{Part~(\ref{thm:conv-lmset-char:a}):}
Suppose $\zero\in\conv{S}$.
Then
\[
  \conv(\lmset{S})
  =
  \conv(\{\zero\}\cup\lmset{S})
  =
  \acone{S}
  =
  \repcl{(\cone{S})}.
\]
The first equality is by
\Cref{pr:conhull-prop}(\ref{pr:conhull-prop:c})
since, by
\Cref{lem:conv-lmset-disj-rn}(\ref{lem:conv-lmset-disj-rn:a}),
$\zero\in\conv(\lmset{S})$.
The second equality is by
\Cref{thm:acone-char}(\ref{thm:acone-char:a}),
and the third is by
\Cref{thm:acone-is-repclos}.

\pfpart{Part~(\ref{thm:conv-lmset-char:b}):}
Suppose $\zero\not\in\conv{S}$.
Then, similar to the foregoing,
\[
  \conv(\lmset{S})
  \subseteq
  \conv(\{\zero\}\cup\lmset{S})
  =
  \acone{S}
  =
  \repcl{(\cone{S})},
\]
where the equalities follow again from
Theorems~\ref{thm:acone-char}(\ref{thm:acone-char:a})
and~\ref{thm:acone-is-repclos}.
Also, by
\Cref{lem:conv-lmset-disj-rn}(\ref{lem:conv-lmset-disj-rn:b}),
$\conv(\lmset{S})\cap\Rn=\emptyset$.
Therefore,
$\conv(\lmset{S})\subseteq\repcl{(\cone{S})}\setminus\Rn$.

For the reverse inclusion,
let $\zbar\in\repcl{(\cone{S})}\setminus\Rn$,
which we aim to show is in $\conv(\lmset{S})$.
Since $\acone{S}=\repcl{(\cone{S})}$ (by \Cref{thm:acone-is-repclos}),
we must have $\zbar\in\acone{S}$, which implies, by
\Cref{thm:ast-conic-hull-union-fin},
that $\zbar\in\oconich{V}$ for some
${V=\{\xx_1,\ldots,\xx_m\}\subseteq S}$.
Note that each $\xx_i$ is in $\Rn\wo\{\zero\}$
since $S\subseteq\Rn$ and $\zero\not\in\conv{S}$.
Therefore, by
\Cref{thm:oconic-hull-and-seqs},
there exist sequences $\seq{\lambda_{it}}$ in $\Rpos$
and span-bound sequences $\seq{\xx_{it}}$ in $\Rn$
such that $\xx_{it}\rightarrow\xx_i$ for $i=1,\ldots,m$,
and $\zz_t=\sum_{i=1}^m  \lambda_{it} \xx_{it} \rightarrow \zbar$.
Since, for each $i$, the sequence $\seq{\xx_{it}}$ is span-bound, each
element $\xx_{it}$ is in $\rspanset{\xx_i}=\spn\{\xx_i\}$;
thus, $\xx_{it}=\gamma_{it}\xx_i$ for some
$\gamma_{it}\in\R$.
Furthermore, since $\xx_{it}\rightarrow\xx_i$ and $\xx_i\neq\zero$,
we must have $\gamma_{it}\rightarrow 1$.
Thus, without loss of generality, we can assume henceforth that
$\gamma_{it}>0$ for all~$i$ and all $t$, since
all other elements can be discarded.

For each $t$ and $i=1,\ldots,m$,
let $\lambda'_{it}=\lambda_{it} \gamma_{it}$, which is nonnegative
since $\lambda_{it}$ and $\gamma_{it}$ are.
We then have that
\begin{equation}   \label{eq:thm:conv-lmset-char:2}
  \zz_t = \sum_{i=1}^m \lambda'_{it} \xx_i.
\end{equation}
Also, let
$  \sigma_t = \sum_{i=1}^m  \lambda'_{it} $.
We claim that $\sigma_t\rightarrow+\infty$.
Otherwise,
if $\sigma_t$ remains bounded on any subsequence, then $\zz_t$ would
also remain in a compact subregion of $\Rn$ on that subsequence,
implying that $\zbar$, its limit, is in that same subregion,
a contradiction since $\zbar\not\in\Rn$.
As a consequence,
we can assume $\sigma_t>0$ for all $t$ since all other
elements can be discarded.

Rewriting \eqref{eq:thm:conv-lmset-char:2}, we then have that
\[
   \zz_t = \sum_{i=1}^m \biggParens{\frac{\lambda'_{it}}{\sigma_t}}
     (\sigma_t \xx_i).
\]
Since $\sigma_t \xx_i\rightarrow\limray{\xx_i}$
for all $i$
(by \Cref{thm:i:seq-rep}),
and since $\sum_{i=1}^m (\lambda'_{it}/\sigma_t)=1$ for all $t$,
this then implies that
\[
  \zbar
  \in
  \ohull\{\limray{\xx_1},\ldots,\limray{\xx_m}\}
  \subseteq
  \conv(\limray{S}),
\]
where the inclusions are by Theorems~\ref{thm:e:7}
and~\ref{thm:convhull-of-simpices},
respectively,
completing the proof.
\qedhere
\end{proof-parts}
\end{proof}

\begin{example}
\label{ex:seg-oe1-oe2:cont}
\indexg{segments, astral!examples|(}%
In \Cref{ex:seg-oe1-oe2},
we gave a somewhat involved calculation of the representations of all
points on the segment joining $\limray{\ee_1}$ and $\limray{\ee_2}$
in $\extspac{2}$.
That segment, which is the same as the convex hull of the two astrons,
can now be computed directly from
\Cref{thm:conv-lmset-char} with $S=\{\ee_1,\ee_2\}$ so that
$\cone{S}=\Rpos^2$, yielding
\[
  \lb{\limray{\ee_1}}{\limray{\ee_2}}
  =
  \conv\{\limray{\ee_1},\limray{\ee_2}\}
  =
  \repcl{(\Rpos^2)}\setminus\R^2.
\]
Thus, this segment consists of all infinite points in $\extspac{2}$
with representations that only involve vectors in $\Rpos^2$.
Although these representations are not canonical as they were in
\Cref{ex:seg-oe1-oe2}, they nonetheless represent the same set of points.%
\indexg{representational closure!convex hull of astrons and|)}%
\indexg{segments, astral!examples|)}%
\indexg{convex hull, astral!astrons@of astrons|)}%
\end{example}

\section{The astral conic hull of a polyhedral cone}
\label{sec:repcl-and-polyhedral}

\indexg{representational closure!topological closure, when equal to|(}%
\indexg{closure, astral!representational closure and|(}%
In the last section, we discussed two possible extensions of a convex
cone $K\subseteq\Rn$ to astral space:
its representational closure $\repcl{K}$, which is the
same as its astral conic hull
(by \Cref{thm:acone-is-repclos}),
and its topological closure $\Kbar$, which is the same as
its outer conic hull
(\Cref{thm:out-conic-h-is-closure}).
In this section, we consider when these extensions, both convex astral
cones, are the same.

We begin with an example showing that it is indeed possible for them
to be different,
even when $K$ is a closed (in $\Rn$) convex cone:

\begin{example}
\indexg{Sideways cone!representational and topological closures different|(}%
Let $K$ be the sideways cone in $\R^3$ introduced in
\Cref{ex:KLsets:extsum-not-sum-exts}:
\[
    K
  =\set{\xx\in\R^3:\: x_1^2\le 2 x_2 x_3\text{ and }x_2,x_3\ge 0}.
\]
It can be checked that $K$ is a closed (in $\R^3$)
convex cone.

Let $\xbar=\limray{\ee_2}\plusl\ee_1$.
In \Cref{ex:KLsets:extsum-not-sum-exts},
we argued that $\xbar\in\Kbar$.
On the other hand, $\xbar\not\in\repcl{K}$.
To see this, observe first that if $\xx\in K$ and $x_3=0$ then it must
also be the case that $x_1=0$.
Suppose $\xbar\in\repcl{K}$ so that
$\xbar=\limrays{\vv_1,\ldots,\vv_k}\plusl\qq$
for some $\vv_1,\ldots,\vv_k,\qq\in K$.
Note that $\xbar\cdot\ee_3=0$.
This implies that
$\vv_i\cdot\ee_3=0$ for $i=1,\ldots,k$, and so also that
$\qq\cdot\ee_3=0$
(by \Cref{pr:vtransu-zero}).
But then, by the preceding observation, it must also be the case that
$\qq\cdot\ee_1=0$, and that
$\vv_i\cdot\ee_1=0$ for $i=1,\ldots,k$.
These imply that $\xbar\cdot\ee_1=0$, a contradiction, since in fact,
$\xbar\cdot\ee_1=1$.
We conclude that $\repcl{K}\neq\Kbar$ in this case.%
\indexg{Sideways cone!representational and topological closures different|)}%
\end{example}

Thus, it is possible for $\repcl{K}$ and $\Kbar$ to be different, but
it is also possible for them to be the same.
\indexg{finitely generated (set)!representational closure of|(}%
Indeed, if $K$ is {finitely generated}, that is, if $K=\cone V$
for some finite set $V\subseteq\Rn$
(see
\Cref{pr:fin-gen-cvx-cone}\ref{pr:fin-gen-cvx-cone:a}\ref{pr:fin-gen-cvx-cone:b}),
then
\Cref{lem:conv-lmset-char:a}, together with
\Cref{cor:a:1}, prove that $\repcl{K}$ and $\Kbar$
must be equal, implying also that $\repcl{K}$ must be closed.
In fact, as we show in the next theorem,
these properties of $K$
and $\repcl{K}$
turn out to be equivalent:
$\repcl{K}$ is closed
and equal to $\Kbar$
if and only
if $K$ is finitely generated.
Later, in \Cref{sec:arescone},
we will see that equality between $\repcl{K}$ and $\Kbar$ (for a specific cone $K$)
is required for an important property of functions
on $\Rn$ called recessive completeness,
related to continuity of a function's extension and the structure of its minimizers.

\indexg{polyhedral sets (standard)!representational closure of|(}%
Recall that
a convex set in $\Rn$ is finitely generated
if and only if it is {polyhedral}, meaning that the set is the
intersection of finitely many closed halfspaces
(\Cref{roc:thm19.1}\ref{roc:thm19.1:a}\ref{roc:thm19.1:b}).
Thus, when discussing convex sets in $\Rn$,
particularly convex cones, we can
use these terms interchangeably.

\begin{theorem}  \label{thm:repcl-polyhedral-cone}
  Let $K\subseteq\Rn$ be a convex cone.
  Then the following are equivalent:
  \begin{letter-compact}
  \item    \label{thm:repcl-polyhedral-cone:a}
    $\repcl{K}$ is closed (in $\extspace$).
  \item    \label{thm:repcl-polyhedral-cone:b}
    $\repcl{K}=\Kbar$.
  \item    \label{thm:repcl-polyhedral-cone:c}
    $K$ is finitely generated (or equivalently, polyhedral);
    that is, $K=\cone V$ for some finite set $V\subseteq\Rn$.
  \end{letter-compact}
\end{theorem}

\begin{proof}
  ~

\begin{proof-parts}
\pfpart{%
  (\ref{thm:repcl-polyhedral-cone:b})
  $\Rightarrow$
  (\ref{thm:repcl-polyhedral-cone:a}):
}
This is immediate.

\pfpart{%
  (\ref{thm:repcl-polyhedral-cone:c})
  $\Rightarrow$
  (\ref{thm:repcl-polyhedral-cone:b}):
}
\Cref{cor:a:1} shows that $\repcl{K}\subseteq\Kbar$ in
general.
If $K$ is finitely generated, then
\Cref{lem:conv-lmset-char:a} proves $\Kbar\subseteq\repcl{K}$.

\pfpart{%
  (\ref{thm:repcl-polyhedral-cone:a})
  $\Rightarrow$
  (\ref{thm:repcl-polyhedral-cone:c}):
}
Assume $\repcl{K}$ is closed, and therefore compact
by \Cref{prop:compact}(\ref{prop:compact:closed-subset}),
being a closed subset of the compact space $\extspace$.
To prove the result, we construct an open cover of $\repcl{K}$, which,
since $\repcl{K}$ is compact, must include a finite subcover.
From this, we show that a finite set of points can be extracted that
are sufficient to generate the cone $K$.

We make use of the open sets $\Uv$ constructed in \Cref{thm:formerly-lem:h:1:new}
for all $\vv\in\Rn$. 
Each such set $\Uv$ includes the astron $\limray{\vv}$,
while being entirely contained in
$\Rn \cup [\limray{\vv}\plusl\Rn]$,
meaning all points in $\Uv$ are either in $\Rn$ or have the form
$\limray{\vv}\plusl\qq$ for some $\qq\in\Rn$.

Slightly overloading notation,
for any set $S\subseteq\Rn$, we further define $\US$ to be the convex
hull of the union of all sets $\Uv$ over $\vv\in S$:
\[
  \US = \conv\BiggParens{\,\bigcup_{\vv\in S} \Uv}.
\]
The parenthesized union is open (being the union of open sets),
implying, by \Cref{cor:convhull-open}, that its convex
hull, $\US$, is also open, for all $S\subseteq\Rn$.
Using compactness, we show next that $\repcl{K}$ is included in one of
these sets $\UV$ for some {finite} set $V\subseteq K$.

First,
for all $\vv\in\Rn$, $\limray{\vv}\in\Uv$;
therefore, for all $S\subseteq\Rn$,
$\lmset{S}\subseteq\US$, and so
\begin{equation}  \label{eq:thm:repcl-polyhedral-cone:2}
  \conv(\lmset{S})\subseteq\US
\end{equation}
by \Cref{pr:conhull-prop}(\ref{pr:conhull-prop:aa})
since $\US$ is convex.
As a result, we can cover $\repcl{K}$ using the collection of all
open sets $\UV$ for all finite subsets
$V\subseteq K$.
This is a cover because
\[
  \repcl{K}
  =
  \conv(\lmset{K})
  =
  \bigcup_{\scriptontop{V\subseteq K:}{\abs{V}<+\infty}} \!\!\conv(\lmset{V})
  \subseteq
  \bigcup_{\scriptontop{V\subseteq K:}{\abs{V}<+\infty}} \!\!\UV.
\]
Here,
the first equality is
by \Cref{cor:a:1},
and the second by
\Cref{thm:convhull-of-simpices}
(combined
with \Cref{pr:conhull-prop}\ref{pr:conhull-prop:a}).
The inclusion is by \eqref{eq:thm:repcl-polyhedral-cone:2}.

Since $\repcl{K}$ is compact, there exists a finite subcover, that is,
a finite collection of sets $V_1,\ldots,V_\ell$ where each $V_j$ is a
finite subset of $K$, for $j=1,\ldots,\ell$, and
such that
\begin{equation} \label{eqn:thm:repcl-polyhedral-cone:1}
  \repcl{K} \subseteq \bigcup_{j=1}^{\ell} \UVj.
\end{equation}
Let
\[ V=\{\zero\}\cup \bigcup_{j=1}^{\ell} V_j \]
be their union, along with the origin (added for later convenience),
which is also a finite subset of $K$.
Furthermore,
\eqref{eqn:thm:repcl-polyhedral-cone:1} implies
$\repcl{K}\subseteq \UV$ since, for $j=1,\ldots,\ell$, $V_j\subseteq V$,
implying $\UVj\subseteq\UV$
by \Cref{pr:conhull-prop}(\ref{pr:conhull-prop:b}).

Summarizing, $\repcl{K}\subseteq \UV$, where
$\zero\in V\subseteq K$ and $|V|<+\infty$.
Let $\Khat=\cone V$ be the cone generated by $V$.
To complete the proof that $K$ is finitely generated, we show that
$K=\Khat$.
We have
$\Khat=\cone V\subseteq K$, because
$V$ is included in $K$ and $K$ is a convex cone.
Thus, it remains to prove that $K\subseteq\Khat$.

Let $\ww$ be any point in $K$, which
we aim to show is in $\Khat$.
Since $\ww\in K$,
\[
   \limray{\ww}\in\repcl{K}\subseteq \UV
   =
   \conv\BiggParens{\,\bigcup_{\vv\in V} \Uv}.
\]
Therefore, by
\Cref{thm:convhull-of-simpices},
\begin{equation}  \label{eqn:thm:repcl-polyhedral-cone:5}
   \limray{\ww}
   \in
   \simplex\{ \xbar_1,\ldots,\xbar_m \}
\end{equation}
for some $\xbar_1,\ldots,\xbar_m\in\bigcup_{\vv\in V} \Uv$.
From the form of points in $\Uv$
(\Cref{thm:formerly-lem:h:1:new}\ref{thm:formerly-lem:h:1:b:new}),
this means, for $j=1,\ldots,m$, that we can write
$\xbar_j=\limray{\vv_j}\plusl\qq_j$ for some $\vv_j\in V$ and some
$\qq_j\in\Rn$.
(Note that this takes into account the possibility that $\xbar_j$
might be in $\Rn$ since in that case we can choose $\vv_j=\zero$,
which is in $V$.)

Recall from \Cref{sec:prelim:cones} that a polar of
a convex cone $J\subseteq\Rn$ is defined
as
$
  \polar{J}=\set{\uu\in\Rn:\:\xx\inprod\uu\le 0\text{ for all }\xx\in J}
$.
We claim that $\ww$ is in the polar
of the polar of $\Khat$, that is, that $\ww$ is in
$(\Khatpol)^{\circ}=\Khatdubpol$, which we then prove is equal to $\Khat$.

To see that $\ww\in\Khatdubpol$,
let $\uu$ be any point in $\Khatpol$.
Then in light of \eqref{eqn:thm:repcl-polyhedral-cone:5},
\Cref{pr:ohull-simplify}(\ref{pr:ohull-simplify:a},\ref{pr:ohull-simplify:b})
implies that
\begin{equation} \label{eqn:thm:repcl-polyhedral-cone:4}
   \limray{\ww}\cdot\uu
   \leq
   \max\{ \xbar_1\cdot\uu,\ldots,\xbar_m\cdot\uu \}.
\end{equation}
Also, for $j=1,\ldots,m$,
we have
$\vv_j\in V \subseteq \cone V = \Khat$, so
$\vv_j\cdot\uu \leq 0$.
Therefore,
\[
  \xbar_j\cdot\uu
  =
  \limray{\vv_j}\cdot\uu
  \plusl
  \qq_j\cdot\uu
  < +\infty
\]
since $\limray{\vv_j}\cdot\uu\leq 0$ and $\qq_j\cdot\uu\in\R$.
Combined with \eqref{eqn:thm:repcl-polyhedral-cone:4}, this
means $\limray{\ww}\cdot\uu<+\infty$, and therefore
$\limray{\ww}\cdot\uu\leq 0$ (since $\limray{\ww}$ is an icon)
so $\ww\cdot\uu\leq 0$.

Since this holds for all $\uu\in\Khatpol$, it follows that
$\ww\in\Khatdubpol$, and so that $K\subseteq\Khatdubpol$.
Furthermore,
because $\Khat$ is a finitely generated convex cone in $\Rn$, it
must be closed in~$\Rn$
(\Cref{roc:thm19.1}\ref{roc:thm19.1:b}\ref{roc:thm19.1:c}),
and so $\Khatdubpol=\Khat$
(by \Cref{pr:polar-props}\ref{pr:polar-props:c}).

Thus, $K=\Khat$, so $K$ is finitely generated.%
\indexg{polyhedral sets (standard)!representational closure of|)}%
\indexg{finitely generated (set)!representational closure of|)}%
\qedhere
\end{proof-parts}
\end{proof}

\indexg{closure, astral!linear subspace@of linear subspace|(}%
\indexg{representational closure!linear subspace@of linear subspace|(}%
\indexg{linear subspaces (standard)!closure same as representational closure|(}%
Applying \Cref{thm:repcl-polyhedral-cone} to linear subspaces yields the following corollary:

\begin{corollary}   \label{cor:lin-sub-bar-is-repcl}
  Let $L\subseteq\Rn$ be a linear subspace.
  Then $\Lbar=\repcl{L}$.
\end{corollary}

\begin{proof}
Let $B\subseteq\Rn$ be any basis for $L$.
Then $|B|$ is finite and
$L=\spn{B}=\cone(B\cup -B)$ (by
\Cref{pr:scc-cone-elts}\ref{pr:scc-cone-elts:span}).
Thus, $L$ is a finitely generated convex cone.
The claim then follows directly from
\Cref{thm:repcl-polyhedral-cone}(\ref{thm:repcl-polyhedral-cone:c},\ref{thm:repcl-polyhedral-cone:b}).%
\indexg{representational closure!topological closure, when equal to|)}%
\indexg{closure, astral!representational closure and|)}%
\indexg{closure, astral!linear subspace@of linear subspace|)}%
\indexg{representational closure!linear subspace@of linear subspace|)}%
\indexg{linear subspaces (standard)!closure same as representational closure|)}%
\end{proof}

\section{Astral polar cones}
\label{sec:ast-pol-cones}

We next study extensions of the standard polar cone
to astral space.
In later chapters, we will see applications of such astral polar
cones,
for instance, in the critical role they play
in the minimization and continuity properties of convex functions.

Recall from \Cref{sec:prelim:cones}
that the standard polar $\Kpol$ of a convex cone
$K\subseteq\Rn$ is the set
\begin{equation}   \label{eq:std-pol-cone-rvw}
  \Kpol
  =
  \braces{\uu\in\Rn :\: \xx\cdot\uu\leq 0 \mbox{ for all } \xx\in K}.
\end{equation}
We will consider two extensions of this notion.
\indexg{polar (primal)!redefined for naive cones|(}%
\indexg{polar (primal)|(}%
In the first of these, we simply extend this definition to apply to
any naive cone in $\extspace$
(including any astral cone or any standard cone in $\Rn$):

\begin{definition}
\label{def:polar}
Let $K\subseteq\extspace$ be a naive cone.
Then the \emph{polar} (or \emph{polar cone})
of $K$, denoted $\Kpol$, is the set
\begin{equation}   \label{eq:pol-cone-ext-defn}
\indexm{k 300}{$\Kpol$}{(primal) polar}%
  \Kpol
  =
  \braces{\uu\in\Rn :\: \xbar\cdot\uu\leq 0 \mbox{ for all } \xbar\in K}.
\end{equation}
\end{definition}
\indexg{polar (primal)!redefined for naive cones|)}%
If $K\subseteq\Rn$, then this definition matches the one
given in \eqref{eq:std-pol-cone-rvw}, so the new definition of $\Kpol$
subsumes and is consistent with the old one for standard convex cones.
We sometimes use the term \emph{primal} polar cone (or \emph{primal} polar) for $\Kpol$
to emphasize the distinction from \emph{dual} polar cones, which we define next.

\indexg{polar, dual|(}%
In
\Cref{def:polar},
the naive cone $K$ is allowed to be a subset of $\extspace$,
and its
polar~$\Kpol$ is a subset of $\Rn$.
On the other hand,
in the definition of dual polar cone,
we consider $K$ that is a subset of $\Rn$,
and define its dual polar $\apol{K}$ to be a subset of $\eRn$:

\begin{definition}
\indexg{polar, dual!defined|(}%
Let $K\subseteq\Rn$ be a cone.
Then the
\emph{dual polar}
(or \emph{dual polar cone})
of~$K$, denoted $\apol{K}$, is the set
\[
\indexm{k 800}{$\apol{K}$}{dual polar}%
  \apol{K} = \braces{\xbar\in\extspace :\:
                       \xbar\cdot\uu\leq 0 \mbox{ for all } \uu\in K }.%
\indexg{polar (primal)|)}%
\indexg{polar, dual!defined|)}%
\]
\end{definition}

In the next two \namecrefs{pr:ast-pol-props} we derive some properties of astral polarity operations.
Note that in \Cref{pr:ast-pol-props}, $J$ and $K$ are only assumed to be cones, but not necessarily
convex cones, so $K$ might be distinct from $\cone K$ in \Cref{pr:ast-pol-props}(\ref{pr:ast-pol-props:f}).
Similarly, in \Cref{pr:ext-pol-cone-props}, $J$ and $K$ are only assumed to be naive cones in $\eRn$,
but not necessarily astral cones or convex astral cones,
so $K$ might be distinct from $\acone K$ in \Cref{pr:ext-pol-cone-props}(\ref{pr:ext-pol-cone-props:b}).

\begin{proposition}  \label{pr:ast-pol-props}
  Let $J,K\subseteq\Rn$ be cones.
  Then the following hold:
  \begin{letter-compact}
  \item  \label{pr:ast-pol-props:c}
    $\apol{K}$ is a closed (in $\extspace$) convex astral cone.
  \item  \label{pr:ast-pol-props:b}
    If $J\subseteq K$, then
    $\apol{K}\subseteq\apol{J}$.
  \item  \label{pr:ast-pol-props:coneSpol}
    Let $S\subseteq\Rn$.
    Then
    $
      \apol{(\cone{S})}
      =
      \Braces{\xbar\in\extspace :\:
        \xbar\cdot\uu\leq 0 \textup{ for all } \uu\in S
      }
    $.
  \item  \label{pr:ast-pol-props:f}
    $\apol{K}=\apol{(\cone{K})}$.
  \item  \label{pr:ast-pol-props:e}
    $\apol{(J+K)}=\Japol\cap\Kapol$.
  \item  \label{pr:ast-pol-props:a}
    $\Kpol=\apol{K}\cap \Rn$.
  \end{letter-compact}
\end{proposition}

\begin{proof}
~

\begin{proof-parts}
\pfpart{Part~(\ref{pr:ast-pol-props:c}):}
The set $\apol{K}$ can be expressed as
\[
  \apol{K} = \bigcap_{\uu\in K}
                       \Braces{\xbar\in\extspace :\: \xbar\cdot\uu\leq 0}.
\]
Each set in this intersection is either a homogeneous
closed astral halfspace or all of $\extspace$ (if $\uu=\zero$).
Therefore, $\apol{K}$ is a closed convex astral cone
by
\Cref{pr:astral-cone-props}(\ref{pr:astral-cone-props:e}).

\pfpart{Part~(\ref{pr:ast-pol-props:b}):}
If $\xbar\in\apol{K}$ then $\xbar\cdot\uu\leq 0$ for all $\uu\in K$,
and therefore all $\uu\in J$; thus, $\xbar\in\apol{J}$.

\pfpart{Part~(\ref{pr:ast-pol-props:coneSpol}):}
Let
$U=\braces{\xbar\in\extspace :\: \xbar\cdot\uu\leq 0 \text{ for all } \uu\in S }$.
We aim to show a point $\xbar\in\extspace$ is in
$\apol{(\cone{S})}$ if and only if it is in $U$.
If $\xbar\in\apol{(\cone{S})}$, then $\xbar\cdot\uu\leq 0$
for all $\uu$ in $\cone{S}$, and so also for all $\uu$ in $S\subseteq\cone{S}$;
therefore, $\xbar\in U$.
For the converse, suppose $\xbar\in U$,
meaning $\xbar\cdot\uu\leq 0$ for $\uu\in S$, and so that
$S\subseteq C$ where
$C=\{\uu\in\Rn :\: \xbar\cdot\uu \leq 0\}$.
The set $C$ is a convex cone by
\Cref{pr:ast-def-hfspace-convex}(\ref{i:dual-hfspace-cone}).
Therefore, $\cone{S}\subseteq C$,
implying $\xbar\inprod\uu\le 0$ for all $\uu\in\cone{S}$, so $\xbar\in\apol{(\cone{S})}$.

\pfpart{Part~(\ref{pr:ast-pol-props:f}):}
This follows directly from part~(\ref{pr:ast-pol-props:coneSpol})
applied with $S=K$.

\pfpart{Part~(\ref{pr:ast-pol-props:e}):}
Since $\zero\in K$, $J\subseteq J+K$, so
$\apol{(J+K)}\subseteq\Japol$
by part~(\ref{pr:ast-pol-props:b}).
Similarly,
$\apol{(J+K)}\subseteq\Kapol$,
so
$\apol{(J+K)}\subseteq\Japol\cap\Kapol$.

For the reverse inclusion, let $\xbar\in\Japol\cap\Kapol$.
Let $\uu\in J$ and $\vv\in K$.
Then, by definition of polar,
$\xbar\cdot\uu\leq 0$ and $\xbar\cdot\vv\leq 0$.
Thus, $\xbar\cdot\uu$ and $\xbar\cdot\vv$ are summable,
so $\xbar\cdot(\uu+\vv) = \xbar\cdot\uu + \xbar\cdot\vv \leq 0$
by \Cref{pr:i:1}.
Since this holds for all $\uu\in J$ and $\vv\in K$, it follows that
$\xbar\in\apol{(J+K)}$, proving the claim.

\pfpart{Part~(\ref{pr:ast-pol-props:a}):}
This follows immediately from definitions.%
\indexg{polar, dual|)}%
\qedhere
\end{proof-parts}
\end{proof}

\begin{proposition}   \label{pr:ext-pol-cone-props}
\indexg{polar (primal)|(}%
  Let $J,K\subseteq\extspace$ be naive cones.
  Then:
  \begin{letter-compact}
  \item    \label{pr:ext-pol-cone-props:a}
    $\Kpol$ is a convex cone in $\Rn$.
  \item    \label{pr:ext-pol-cone-props:d}
    If $J\subseteq K$ then $\Kpol\subseteq\Jpol$.

  \item    \label{pr:ext-pol-cone-props:coneSpol}
    Let $S\subseteq\extspace$.
    Then
    $
      \polar{(\oconich{S})}
      =
      \Braces{\uu\in\Rn :\:
        \xbar\cdot\uu\leq 0 \textup{ for all } \xbar\in S
      }
    $.
  \item    \label{pr:ext-pol-cone-props:b}
    $\Kpol=\polar{(\Kbar)}=\polar{(\acone{K})}=\polar{(\oconich{K})}$.
  \item    \label{pr:ext-pol-cone-props:e}
    $\polar{(J\seqsum K)}=\Jpol\cap\Kpol$.
  \end{letter-compact}
\end{proposition}

\begin{proof}
  ~

\begin{proof-parts}
\pfpart{Part~(\ref{pr:ext-pol-cone-props:a}):}
We can write
\[
  \Kpol
  =
  \bigcap_{\xbar\in K} \Braces{\uu\in\Rn :\: \xbar\cdot\uu\leq 0}.
\]
Each of the sets appearing on the right is a convex cone by
\Cref{pr:ast-def-hfspace-convex}(\ref{i:dual-hfspace-cone}).
Therefore, $\Kpol$ is a convex cone.

\pfpart{Part~(\ref{pr:ext-pol-cone-props:d}):}
Proof is analogous to that of
\Cref{pr:ast-pol-props}(\ref{pr:ast-pol-props:b}).

\pfpart{Part~(\ref{pr:ext-pol-cone-props:coneSpol}):}
Let $U=
      \braces{\uu\in\Rn :\:
        \xbar\cdot\uu\leq 0 \text{ for all } \xbar\in S
      }
$.
For $\uu\in\Rn$, we aim to show that $\uu\in\polar{(\oconich{S})}$ if
and only if $\uu\in U$.
If $\uu\in\polar{(\oconich{S})}$ then $\xbar\cdot\uu\leq 0$ for all
$\xbar$ in $\oconich{S}$ and so also all $\xbar$ in $S$;
hence, $\uu\in U$.
For the converse, suppose $\uu\in U$.
If $\uu=\zero$ then it is also in
$\polar{(\oconich{S})}$, since the latter is a cone
by part~(\ref{pr:ext-pol-cone-props:a}).
So suppose $\uu\ne\zero$. Since
$\xbar\cdot\uu\leq 0$ for all $\xbar\in S$,
we have
$S\subseteq H$ where
$H=\{\xbar\in\extspace :\: \xbar\cdot\uu \leq 0\}$.
Since $H$
is a homogeneous closed astral halfspace that includes $S$,
this implies $\oconich{S}\subseteq H$ by definition of outer conic hull.
Therefore, $\uu\in\polar{(\oconich{S})}$.

\pfpart{Part~(\ref{pr:ext-pol-cone-props:b}):}
Let $K'\subseteq\extspace$ be any naive cone with
$K\subseteq K'\subseteq\oconich{K}$.
Then
\[
  \polar{(K')}
  =
  \polar{(\oconich{K'})}
  =
  \polar{(\oconich{K})},
\]
where the first equality is by
part~(\ref{pr:ext-pol-cone-props:coneSpol}) applied with $S=K'$,
and the second equality is because
$\oconich{K'}=\oconich{K}$
by \Cref{pr:gen-hull-ops}(\ref{pr:gen-hull-ops:d}) since
outer conic hull is a hull operator.
Applied repeatedly with $K'$ set to $K$, $\Kbar$, and $\acone{K}$,
this proves all the claimed equalities.
(We also used that $\acone{K}\subseteq\oconich{K}$ by
\Cref{pr:acone-hull-props}\ref{pr:acone-hull-props:d},
and that $\Kbar\subseteq\oconich{K}$ since $\oconich{K}$ is closed and
includes $K$.
We further note that $\Kbar$ is a naive cone since $K$ is.
This is because if $\xbar\in\Kbar$ then there exists a sequence
$\seq{\xbar_t}$ in $K$ with $\xbar_t\rightarrow\xbar$, implying for
all $\lambda\in\Rpos$ that $\lambda\xbar_t\rightarrow\lambda\xbar$ and
thus that $\lambda\xbar\in\Kbar$ since each $\lambda\xbar_t\in K$.)

\pfpart{Part~(\ref{pr:ext-pol-cone-props:e}):}
That $\polar{(J\seqsum K)}\subseteq\Jpol\cap\Kpol$
follows by an argument analogous to that given in
proving
\Cref{pr:ast-pol-props}(\ref{pr:ast-pol-props:e}).

For the reverse inclusion, suppose $\uu\in\Jpol\cap\Kpol$.
Let $\zbar\in J\seqsum K$.
Then $\zbar\in\xbar\seqsum\ybar$
for some $\xbar\in J$ and $\ybar\in K$. We have
$\xbar\cdot\uu\leq 0$ and $\ybar\cdot\uu\leq 0$
(since $\uu\in\Jpol\cap\Kpol$),
and so $\xbar\cdot\uu$ and $\ybar\cdot\uu$ are summable.
Therefore,
$\zbar\cdot\uu = \xbar\cdot\uu + \ybar\cdot\uu \leq 0$
by \Cref{cor:seqsum-conseqs}(\ref{cor:seqsum-conseqs:b}).
Thus, $\uu\in\polar{(J\seqsum K)}$, completing the proof.
\qedhere
\end{proof-parts}
\end{proof}

Applying the standard polar operation twice to a convex cone
$K\subseteq\Rn$ yields its closure; that is,
$\dubpolar{K}=\cl K$
(\Cref{pr:polar-props}\ref{pr:polar-props:c}).
More generally, it can be shown that if $K\subseteq\Rn$ is a cone
(not necessarily convex), then
$\dubpolar{K}=\cl(\cone{K})$.
We next consider the similar effect of applying the astral polarity
operations twice:

\begin{theorem}   \label{thm:dub-ast-polar}
  ~

  \begin{letter-compact}
  \item   \label{thm:dub-ast-polar:a}
\indexg{polar, dual|(}%
    Let $J\subseteq\extspace$ be a naive cone.
    Then $\Jpolapol=\oconich{J}$.
    Therefore, if $J$ is a closed convex astral cone,
    then $\Jpolapol=J$.
  \item   \label{thm:dub-ast-polar:b}
    Let $K\subseteq\Rn$ be a cone.
    Then $\Kapolpol=\cone{K}$.
    Therefore, if $K$ is a convex cone,
    then $\Kapolpol=K$.
  \end{letter-compact}
\end{theorem}

\begin{proof}
  ~

\begin{proof-parts}
\pfpart{Part~(\ref{thm:dub-ast-polar:a}):}
For $\uu\in\Rn\wo\{\zero\}$, let $\chsuz$ be a homogeneous closed
astral halfspace, as in \eqref{eqn:homo-halfspace}.
Then
\begin{align*}
  \Jpolapol
  &=
  \;\bigcap_{\mathclap{\uu\in\Jpol}}
  \;
  \bigBraces{\xbar\in\extspace :\: \xbar\cdot\uu\leq 0}
  \\
  &=
  \;\bigcap
  \;
  \bigBraces{\chsuz :\: \uu\in\Rn\wo\{\zero\},\, J\subseteq\chsuz}
  =
  \oconich{J}.
\end{align*}
The first equality is by definition of dual polar,
and the third is by definition of outer conic hull.
The second equality is because, for $\uu\in\Rn\wo\{\zero\}$,
$\uu\in\Jpol$ if and only if
$J\subseteq\chsuz$
(and because $\{\xbar\in\extspace :\: \xbar\cdot\uu\leq 0\}=\extspace$
if $\uu=\zero$).

Thus,
if $J$ is a closed convex astral cone, then
$\Jpolapol=\oconich{J}=J$ by
\Cref{cor:sep-ast-cone-conseqs}(\ref{cor:sep-ast-cone-conseqs:b}).

\pfpart{Part~(\ref{thm:dub-ast-polar:b}):}
First, if $\uu\in K$ then $\xbar\cdot\uu\leq 0$ for all
$\xbar\in\Kapol$, by its definition, so $\uu\in\Kapolpol$.
Therefore, $K$ is included in $\Kapolpol$, which is a convex cone, by
\Cref{pr:ext-pol-cone-props}(\ref{pr:ext-pol-cone-props:a}).
Thus, $\cone{K}\subseteq\Kapolpol$.

For the reverse inclusion, suppose $\vv\not\in \cone{K}$.
Then, by \Cref{thm:ast-def-sep-cone-sing},
there exists $\xbar\in\extspace$ such that
\[
  \sup_{\uu\in \cone{K}} \xbar\cdot\uu
  \leq
  0
  <
  \xbar\cdot\vv.
\]
Since $K\subseteq\cone{K}$,
the first inequality implies that $\xbar\in\Kapol$.
Combined with the second inequality, this shows that
$\vv\not\in\Kapolpol$.
Thus, $\Kapolpol\subseteq \cone{K}$.
\qedhere
\end{proof-parts}
\end{proof}

In standard convex analysis, the polarity operation $K\mapsto\Kpol$
maps the set of all closed convex cones bijectively onto itself. The operation
is its own inverse.
It is an example of a duality correspondence, much like orthogonality of linear subspaces in $\Rn$, or the conjugacy of closed proper convex functions on $\Rn$.

In an analogous way,
\Cref{thm:dub-ast-polar}
shows that astral polarity operations
define a one-to-one correspondence between the set of all closed
convex astral cones in $\extspace$,
and the set of all convex cones in $\Rn$
(whether closed or not).
The primal polarity operation $J\mapsto\Jpol$ bijectively
maps the set of all closed convex astral cones in $\eRn$ to the set of all
convex cones in $\Rn$. Its inverse is the dual polarity operation $K\mapsto\Kapol$.
These polarity operations are examples of astral duality correspondences.
Unlike dualities in standard convex analysis, astral dualities are typically between objects
defined over $\eRn$ and objects defined over $\Rn$, reflecting the asymmetry of the coupling function.

\indexg{indicator functions (astral)!polar operations and|(}%
\indexg{indicator functions (standard)!polar operations and|(}%
\indexg{polar (primal)!conjugate of indicator and|(}%
\indexg{polar, dual!conjugate of indicator and|(}%
Similar duality holds for the indicator functions
associated with convex cones in $\Rn$ and closed convex astral cones in $\eRn$,
as follows from the next \namecref{pr:ind-apol-conj-ind}.
This is analogous to a standard duality result stating that
the indicator functions of
a closed convex cone $K\subseteq\Rn$ and its standard polar
$\Kpol$ are conjugates of one another (\Cref{roc:thm14.1-conj}).

\begin{proposition}   \label{pr:ind-apol-conj-ind}
  ~

  \begin{letter-compact}
  \item    \label{pr:ind-apol-conj-ind:a}
    Let $J\subseteq\extspace$ be a naive cone.
    Then $\indaJstar=\indJpol$.
  \item    \label{pr:ind-apol-conj-ind:b}
    Let $K\subseteq\Rn$ be a cone.
    Then $\indfdstar{K} = \indfa{\apol{K}}$.
  \end{letter-compact}
\end{proposition}

\begin{proof}
We prove only part~(\ref{pr:ind-apol-conj-ind:a});
the proof for part~(\ref{pr:ind-apol-conj-ind:b}) is entirely
analogous.

Let $\uu\in\Rn$.
Then by the form of the conjugate given in
\eqref{eq:Fstar-down-def},
\[
   \indaJstar(\uu)
   =
   \sup_{\xbar\in\extspace}
       \BigBracks{ - \indaJ(\xbar) \plusd \xbar\cdot\uu }
   =
   \sup_{\xbar\in J} \xbar\cdot\uu,
\]
which we aim to show is equal to $\indJpol(\uu)$.

If $\uu\in\Jpol$, then $\xbar\cdot\uu\leq 0$ for all
$\xbar\in J$, so
$\indaJstar(\uu)\leq 0$.
On the other hand, $\zero\in J$, so
$\indaJstar(\uu)\geq 0$.
Thus,
$\indaJstar(\uu) = 0 = \indJpol(\uu)$
in this case.

Otherwise, if $\uu\not\in\Jpol$, then $\xbar\cdot\uu>0$ for some
$\xbar\in J$.
Since $J$ is a naive cone,
$\lambda\xbar$ is also in $J$
for all $\lambda\in\Rstrictpos$.
This implies that
$\indaJstar(\uu)\geq (\lambda\xbar)\cdot\uu = \lambda (\xbar\cdot\uu)$,
which tends to $+\infty$ as $\lambda\rightarrow+\infty$.
Thus, in this case,
$\indaJstar(\uu) = +\infty = \indJpol(\uu)$.%
\indexg{indicator functions (astral)!polar operations and|)}%
\indexg{indicator functions (standard)!polar operations and|)}%
\indexg{polar (primal)!conjugate of indicator and|)}%
\indexg{polar, dual!conjugate of indicator and|)}%
\indexg{polar (primal)|)}%
\end{proof}

\indexg{polar, dual!closure@of closure|(}%
A basic property of the standard polar cone is that the polar of a
convex cone $K\subseteq\Rn$
is the same as the polar of its closure, that is,
$\Kpol=\polar{(\cl K)}$.
Extending to astral space, we saw in
\Cref{pr:ext-pol-cone-props}(\ref{pr:ext-pol-cone-props:b})
that $\Jpol=\polar{(\Jbar)}$ if $J$ is a naive cone in~$\extspace$.
Nevertheless,
the analogous property does not hold, in general, for the dual polar;
in other words, it is possible that
$\apol{K}\neq\apol{(\cl K)}$.
This is because
dual polarity is a one-to-one correspondence between \emph{all} convex
cones in $\Rn$ (not just the closed ones) and closed astral convex cones in $\eRn$. So, for any standard cone $K$ such that $K\ne\cl K$, we must also have $\apol{K}\ne\apol{(\cl K)}$.

\begin{example}
In $\R^2$, let $K$ be the open left halfplane adjoined with the origin:
\[
   K
   =
   \{\zero\}
   \cup
   \{\uu \in \R^2 :\: u_1 < 0\}.
\]
The standard polar of this convex cone is the ray
$\Kpol = \{ \lambda \ee_1 :\: \lambda\in\Rpos\}$.
Its dual polar, $\apol{K}$, includes $\Kpol$, and also includes all
infinite points in $\extspac{2}$ whose dominant direction is $\ee_1$.
This is because, for every such point $\xbar$,
we have
$\xbar\cdot\uu=\limray{\ee_1}\cdot\uu=-\infty$, for all
$\uu\in K\setminus \{\zero\}$. For infinite points with dominant directions $\vv\ne\ee_1$,
we can always find $\uu\in K$ for which
$\vv\inprod\uu>0$, so such infinite points cannot be in $\apol{K}$.
Thus,
\[
  \apol{K} = \{ \lambda \ee_1 :\: \lambda \in \Rpos\}
        \cup \bracks{ \limray{\ee_1} \plusl \extspac{2} }.
\]
On the other hand,
the dual polar of
$\cl K$, which is the closed left halfplane of $\R^2$, includes only
$\Kpol=\polar{(\cl K)}$
together with the single infinite point $\limray{\ee_1}$;
that is,
\[
  \apol{(\cl K)}
  =
  \{ \lambda \ee_1 :\: \lambda \in \Rpos\}
  \cup \{ \limray{\ee_1} \}.
\]
Note that this set is exactly the closure (in $\extspac{2}$) of
$\Kpol$.
\end{example}

\indexg{polar (primal)!closure of|(}%
\indexg{closure, astral!polar@of polar|(}%
This last observation turns out to be general.
As we show next,
if $K\subseteq\Rn$ is a convex cone, then
$\apol{(\cl K)}$, the dual polar of its closure in $\Rn$,
is always equal to $\Kpolbar$, the closure in $\extspace$ of
its standard polar.
Consequently, $\apol{K}=\Kpolbar$ if and only if $K$ is closed in
$\Rn$.

\begin{theorem}    \label{thm:apol-of-closure}
  Let $K\subseteq\Rn$ be a convex cone.
  Then:
  \begin{letter-compact}
  \item     \label{thm:apol-of-closure:a}
    $\Kpolbar = \apol{(\cl K)} \subseteq \apol{K}$.
  \item     \label{thm:apol-of-closure:b}
    $\apol{K} = \Kpolbar$ if and only if
    $K$ is closed in $\Rn$.
  \end{letter-compact}
\end{theorem}

\begin{proof}
  ~

\begin{proof-parts}
\pfpart{Part~(\ref{thm:apol-of-closure:a}):}
We have
\[
   \Kpolbar
   =
   \oconich{\Kpol}
   =
   \polapol{(\Kpol)}
   =
   \apol{(\polpol{K})}
   =
   \apol{(\cl K)}
   \subseteq
   \apol{K}.
\]
The first equality is by \Cref{thm:out-conic-h-is-closure}
(since $\Kpol$ is a convex cone).
The second equality is by
\Cref{thm:dub-ast-polar}(\ref{thm:dub-ast-polar:a}).
The fourth equality is by
\Cref{pr:polar-props}(\ref{pr:polar-props:c}).
The inclusion is by
\Cref{pr:ast-pol-props}(\ref{pr:ast-pol-props:b})
since $K\subseteq \cl{K}$.

\pfpart{Part~(\ref{thm:apol-of-closure:b}):}
If $K$ is closed in $\Rn$, then
part~(\ref{thm:apol-of-closure:a})
immediately implies that $\apol{K} = \Kpolbar$.

For the converse, suppose $\apol{K} = \Kpolbar$,
implying, by part~(\ref{thm:apol-of-closure:a}),
that $\apol{K} = \apol{(\cl K)}$.
Taking the polar of both sides then yields
\[
   K
   =
   \apolpol{K}
   =
   \apolpol{(\cl K)}
   =
   \cl K,
\]
where the first and third equalities are by
\Cref{thm:dub-ast-polar}(\ref{thm:dub-ast-polar:b})
(since $K$ and $\cl K$ are both convex cones).
Thus, $K$ is closed.%
\indexg{polar, dual!closure@of closure|)}%
\indexg{polar, dual|)}%
\indexg{polar (primal)!closure of|)}%
\indexg{closure, astral!polar@of polar|)}%
\qedhere
\end{proof-parts}
\end{proof}

\chapter{Astral linear subspaces}
\label{sec:ast-lin-subspaces}

In this chapter, we study the astral analogue of linear subspaces.
There are several natural ways of extending the notion of a linear
subspace to astral space; we will see that all of these are
equivalent.
We will also see how such notions as span, basis and orthogonal
complement carry over to astral space.

\section{Definition and characterizations}

We begin with definitions.

A nonempty set $L\subseteq\Rn$ is a linear subspace if it is closed under
vector addition (so that $\xx+\yy\in L$ for all $\xx,\yy\in L$)
and scalar multiplication (so that $\lambda\xx\in L$ for all
$\xx\in L$ and $\lambda\in\R$).
However,
a direct attempt to generalize this particular definition to astral space is
likely to lead to the kind of issues that we saw in \Cref{sec:cones-basic},
when we considered applying the definition of a naive cone to astral space.

Nonetheless,
these conditions can be restated to say equivalently that $L$ is a
linear subspace if and only if $L$ is a convex cone
(which, by \Cref{pr:scc-cone-elts}\ref{pr:scc-cone-elts:d}, holds
if and only if $L$ includes the origin and is closed under
vector addition and positive scalar multiplication)
and is also closed under negation (so that $-\xx\in L$ for all $\xx\in L$).
\indexg{linear subspaces, astral|(}%
\indexg{linear subspaces, astral!defined|(}%
In this form, this definition extends straightforwardly to astral
space:

\begin{definition}
A set $M\subseteq\extspace$ is an
\emph{astral linear subspace} if $M$ is a convex astral cone
and if $M=-M$ (so that $M$ is closed under negation).%
\indexg{linear subspaces, astral!defined|)}%
\end{definition}

\indexg{linear subspaces, astral!characterizations|(}%
There are several ways in which a standard linear subspace can be
constructed, for instance, as the column space of a matrix or as
the set of solutions to a system of homogeneous linear equations.
In a similar way, it is possible to construct and characterize
astral linear subspaces, which we outline in the next theorem.

In words, the theorem states that
a set $M\subseteq\extspace$ is an astral linear subspace if
and only if:
it is the (topological or representational) closure of a standard
linear subspace;
it is the astral column space of some matrix;
it is the astral null space of some matrix;
it is the set of astral points orthogonal to some subset of $\Rn$;
it is the astral conic hull of some set that is closed under
negation;
it is the convex hull of a set of icons, closed under negation;
it is the segment joining an icon and its negation.
Furthermore, the condition that $M$ be an astral conic hull can be
restricted to require that it be the astral conic hull of only finitely
many points, all in $\Rn$,
and similarly,
the condition that $M$ be the convex hull of a set of icons can be
restricted to require that it be the convex hull of finitely many
astrons.

\begin{theorem}   \label{thm:ast-lin-sub-equivs}
  Let $M\subseteq\extspace$.
  Then the following are equivalent:
  \begin{letter-compact}
  \item   \label{thm:ast-lin-sub-equivs:a}
    $M$ is an astral linear subspace.
  \item   \label{thm:ast-lin-sub-equivs:b}
    $M=\Lbar=\repcl{L}$ for some (standard) linear subspace
    $L\subseteq\Rn$.
  \item   \label{thm:ast-lin-sub-equivs:c}
    $M=\acolspace \A = \{\A\zbar :\: \zbar\in\extspac{k}\}$
    for some matrix $\A\in\R^{n\times k}$, $k\geq 0$.
  \item   \label{thm:ast-lin-sub-equivs:d}
    $M=\{\xbar\in\extspace :\: \B\xbar = \zero\}$
    for some matrix $\B\in\R^{m\times n}$, $m\geq 0$.
  \item   \label{thm:ast-lin-sub-equivs:e}
    $M=\{\xbar\in\extspace :\: \xbar\cdot\uu=0 \textup{ for all } \uu\in U\}$,
    for some $U\subseteq\Rn$.
  \item   \label{thm:ast-lin-sub-equivs:f}
    $M=\acone{S}$ for some $S\subseteq\extspace$
    with $S=-S$.
  \item   \label{thm:ast-lin-sub-equivs:ff}
    $M=\acone{S}$ for some $S\subseteq\Rn$
    with $S=-S$ and $|S|<+\infty$.
  \item   \label{thm:ast-lin-sub-equivs:g}
    $M=\conv{E}$ for some nonempty set of icons $E\subseteq\corezn$
    with $E=-E$.
  \item   \label{thm:ast-lin-sub-equivs:gg}
    $M=\conv(\lmset{S})$ for some nonempty set $S\subseteq\Rn$
    with $S=-S$ and $|S|<+\infty$.
  \item   \label{thm:ast-lin-sub-equivs:h}
    $M=\lb{-\ebar}{\ebar}$ for some icon $\ebar\in\corezn$.
  \end{letter-compact}
\end{theorem}

Before proving the theorem, we give a lemma that
spells out in detail the structure underlying some of these equivalences.

\begin{lemma}    \label{lem:ast-subspace-equivs}
  Let $\vv_1,\ldots,\vv_k\in\Rn$, and also define
  $V=\{\vv_1,\ldots,\vv_k\}$,
  $\VV=[\vv_1,\ldots,\vv_k]$,
  $L=\spn{V}=\colspace\VV$,
  and
  $\ebar=\VV\omm$.
  Then
  \[
     \Lbar
     =
     \repcl{L}
     =
     \conv\bigParens{\{\zero\}\cup\lmset{V}\cup-\lmset{V}}
     =
     \acone(V\cup-V)
     =
     \acolspaceVV
     =
     \lb{-\ebar}{\ebar}.
  \]
\end{lemma}

\begin{proof}
We have
\begin{align*}
     \Lbar
     =
     \repcl{L}
     =
     \acone{L}
     &=
     \acone\bigParens{\cone(V\cup-V)}
\\
     &=
     \acone(V\cup-V)
     =
     \conv\bigParens{\{\zero\}\cup\lmset{V}\cup-\lmset{V}}.
\end{align*}
The first equality is by \Cref{cor:lin-sub-bar-is-repcl}. The second
is by \Cref{cor:a:1} (since $L$ is a cone).
The third is because $L=\cone(V\cup-V)$
by \Cref{pr:scc-cone-elts}(\ref{pr:scc-cone-elts:span}).
The fourth is by \Cref{pr:acone-hull-props}(\ref{pr:acone-hull-props:c}), because $\cone(V\cup-V)\subseteq\acone(V\cup-V)$
by \Cref{thm:acone-is-repclos}. And the final equality is by \Cref{thm:acone-char}(\ref{thm:acone-char:a}).

This proves the first three equalities of
the \namecref{lem:ast-subspace-equivs}.
In the remainder, we show that $\acolspaceVV=\Lbar=\lb{-\ebar}{\ebar}$
by proving three separate inclusions:
\begin{proof-parts}

\pfpart{%
  Inclusion $\acolspaceVV \subseteq \Lbar$:
}
Let $F:\extspac{k}\rightarrow\extspace$ be the astral linear map
associated with $\VV$ so that
$F(\zbar)=\VV\zbar$ for $\zbar\in\extspac{k}$.
Then
\[
   \acolspaceVV
   =
   F\bigParens{\extspac{k}}
   \subseteq
   \clbar{F(\Rk)}
   =
   \Lbar.
\]
The first equality is by definition of astral column space.
The inclusion is by
\Crefequiv{prop:cont}{prop:cont:a}{prop:cont:c}
since $F$ is continuous
(\Cref{thm:linear:cont}\ref{thm:linear:cont:b}).
The last equality is because $L$, the set of linear combinations of
columns of $\VV$, is exactly $F(\Rk)$.

\pfpart{%
  Inclusion $\Lbar \subseteq \lb{-\ebar}{\ebar}$:
}
We first show that $L \subseteq \lb{-\ebar}{\ebar}$.
Let $\xx\in L$ and let $\uu\in\Rn$.
We claim that
\begin{equation}    \label{eq:lem:ast-subspace-equivs:1}
  \xx\cdot\uu
  \leq
  |\ebar\cdot\uu|
  =
  \max\{-\ebar\cdot\uu,\ebar\cdot\uu\},
\end{equation}
which, by
\Cref{pr:seg-simplify}(\ref{pr:seg-simplify:b},\ref{pr:seg-simplify:a}),
will imply that
$\xx\in \lb{-\ebar}{\ebar}$.
If $\xx\cdot\uu=0$, then
\eqref{eq:lem:ast-subspace-equivs:1} must hold since the right-hand
side of that equation is always nonnegative.
Otherwise, suppose $\xx\cdot\uu\neq 0$.
Since $\xx$ is a linear combination of $\vv_1,\ldots,\vv_k$, this implies
that $\vv_i\cdot\uu\neq 0$, for some $i\in\{1,\ldots,k\}$.
Consequently, $\ebar\cdot\uu\in\{-\infty,+\infty\}$ by
\Cref{pr:vtransu-zero}, so
\eqref{eq:lem:ast-subspace-equivs:1} must hold since the expression on
the right-hand side is $+\infty$.
Thus, $\xx\in \lb{-\ebar}{\ebar}$,
so $L\subseteq\lb{-\ebar}{\ebar}$.
Since $\lb{-\ebar}{\ebar}$ is closed in $\extspace$
(\Cref{pr:ohull:hull}),
this further shows that
$\Lbar\subseteq\lb{-\ebar}{\ebar}$.

\pfpart{%
  Inclusion $\lb{-\ebar}{\ebar} \subseteq \acolspaceVV$:
}
We have $\ebar=\VV \omm\in\acolspaceVV$ since $\omm\in\extspac{k}$.
Likewise, $-\ebar=-\VV \omm = \VV (-\omm)\in\acolspaceVV$.
Further, with $F$ as defined above,
$\acolspaceVV =  F(\extspac{k})$ is convex
by \Cref{cor:thm:e:9},
being the image of the convex set $\extspac{k}$ under the linear map~$F$.
Therefore,
$\lb{-\ebar}{\ebar} \subseteq \acolspaceVV$.
\qedhere
\end{proof-parts}
\end{proof}

\begin{figure}
  \centering
  \includegraphics{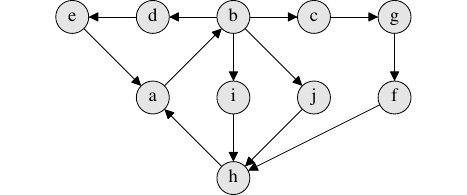}
  \mycaption{Structure of the proof of \Cref{thm:ast-lin-sub-equivs}}{%
    Statements of the theorem are depicted as nodes, proved implications as
    directed edges. Since there is a directed path from each node to every other
    node, all the statements of the theorem are equivalent.
  }%
  \label{fig:linsubspace-equiv-thm-graph}%
\end{figure}

\begin{proof}[Proof of \Cref{thm:ast-lin-sub-equivs}]
We prove several implications, each shown as a directed edge in the graph in
\Cref{fig:linsubspace-equiv-thm-graph}.
Since there is a directed path from each node to every other node,
all the statements of the theorem are equivalent.

\begin{proof-parts}
\pfpart{%
  (\ref{thm:ast-lin-sub-equivs:a})
  $\Rightarrow$
  (\ref{thm:ast-lin-sub-equivs:b}):
}
Suppose $M$ is an astral linear subspace.
Let $L=M\cap\Rn$.
We claim first that $L$ is a (standard) linear subspace.
Since $M$ is a convex astral cone, $L$ is a (standard) convex cone
by \Cref{pr:e1}(\ref{pr:e1:int-rn})
and~\Cref{pr:ast-cone-is-naive}.
And since $M$ is closed under negation, $L$ is as well.
It follows that $L$ is closed under vector addition
by \Cref{pr:scc-cone-elts}(\ref{pr:scc-cone-elts:d}).
Further, since $L$ is a cone,
if $\xx\in L$ and $\lambda\in\R$ with $\lambda\geq 0$, then
$\lambda\xx\in L$;
and if $\lambda<0$ then $\lambda\xx=(-\lambda)(-\xx)\in L$ since
$-\xx\in L$.
Thus, $L$ is closed under multiplication by any scalar.
Together, these facts show that $L$ is a linear subspace.

Next, we show that $M=\repcl{L}$.
Since $M$ is a convex astral cone that includes $L$,
we must have
$\repcl{L}=\acone{L}\subseteq M$ (with equality by
\Cref{thm:acone-is-repclos} and inclusion by \Cref{pr:acone-hull-props}\ref{pr:acone-hull-props:a}).
For the reverse inclusion, suppose %
that $\xbar\in M$.
Then $\limray{\xbar}$ is also in the astral cone~$M$
(\Cref{pr:ast-cone-is-naive}), as is
$-\limray{\xbar}=\limray{(-\xbar)}$ since $-\xbar$ is also in $M$.
Let
\[
  \xbar=\limray{\vv_1}\plusl\dotsb\plusl\limray{\vv_k}\plusl\vv_{k+1}
\]
be a representation of $\xbar$ for some $\vv_1,\dotsc,\vv_{k+1}\in\Rn$ and $k\ge 0$.
We claim that each $\vv_i$ is in $L$, implying $\xbar\in\repcl{L}$.
To see this, let $V=\{\vv_1,\ldots,\vv_{k+1}\}$
and let $\VV=[\vv_1,\ldots,\vv_{k+1}]$.
Then
\[
  V
  \subseteq
  \spn{V}
  \subseteq
  \clbar{\spn{V}}
  =
  \lb{-\VV\omm}{\VV\omm}
  =
  \lb{-\limray{\xbar}}{\limray{\xbar}}
  \subseteq
  M.
\]
The first equality is by \Cref{lem:ast-subspace-equivs}
(applied to $\vv_1,\ldots,\vv_{k+1}$).
The second equality is because
$\limray{\xbar}=\VV \omm$
(\Cref{pr:mult-rep-by-scalar}).
The last inclusion is because, as noted above, $-\limray{\xbar}$ and
$\limray{\xbar}$ are in $M$, and since $M$ is convex.
Thus, $V\subseteq M\cap\Rn=L$, so $\xbar\in\repcl{L}$.
Therefore, $M=\repcl{L}$.
By \Cref{cor:lin-sub-bar-is-repcl}, we also obtain
$\repcl{L}=\Lbar$.

\pfpart{%
  (\ref{thm:ast-lin-sub-equivs:b})
  $\Rightarrow$
  (\ref{thm:ast-lin-sub-equivs:c}),
  (\ref{thm:ast-lin-sub-equivs:b})
  $\Rightarrow$
  (\ref{thm:ast-lin-sub-equivs:gg}),
  (\ref{thm:ast-lin-sub-equivs:b})
  $\Rightarrow$
  (\ref{thm:ast-lin-sub-equivs:h}):
}
Suppose $M=\Lbar$ where $L\subseteq\Rn$ is a linear subspace.
Let $V=\{\vv_1,\ldots,\vv_k\}\subseteq\Rn$ be a basis for $L$
(so that $L=\spn{V}$),
and let $\VV=[\vv_1,\ldots,\vv_k]$.
Then by \Cref{lem:ast-subspace-equivs},
\[
  M = \acolspaceVV = \conv(\{\zero\}\cup\lmset{V}\cup-\lmset{V}) = \lb{-\VV\omm}{\VV\omm},
\]
proving (\ref{thm:ast-lin-sub-equivs:c}) with $\A=\VV$,
(\ref{thm:ast-lin-sub-equivs:gg}) with $S=\{\zero\}\cup V\cup -V$,
and (\ref{thm:ast-lin-sub-equivs:h}) with $\ebar=\VV\omm$.

\pfpart{%
  (\ref{thm:ast-lin-sub-equivs:c})
  $\Rightarrow$
  (\ref{thm:ast-lin-sub-equivs:ff}):
}
Suppose $M=\acolspace \A$ where $\A=[\vv_1,\ldots,\vv_k]$
for some $\vv_1,\ldots,\vv_k\in\Rn$.
Then by \Cref{lem:ast-subspace-equivs},
$M=\acone{S}$ where
$S=V\cup -V$ and $V=\{\vv_1,\ldots,\vv_k\}$.

\pfpart{%
  (\ref{thm:ast-lin-sub-equivs:ff})
  $\Rightarrow$
  (\ref{thm:ast-lin-sub-equivs:f}):
}
This is immediate.

\pfpart{%
  (\ref{thm:ast-lin-sub-equivs:f})
  $\Rightarrow$
  (\ref{thm:ast-lin-sub-equivs:g}):
}
Suppose $M=\acone{S}$ for some $S\subseteq\extspace$ with $S=-S$.
Then by \Cref{thm:acone-char}(\ref{thm:acone-char:a}),
$M=\conv{E}$ where $E=\{\zero\}\cup\lmset{S}$, which is a nonempty set
of icons with $E=-E$.

\pfpart{%
  (\ref{thm:ast-lin-sub-equivs:gg})
  $\Rightarrow$
  (\ref{thm:ast-lin-sub-equivs:g}):
}
This is immediate since the astrons in $\lmset{S}$, for any
$S\subseteq\Rn$, are also icons.

\pfpart{%
  (\ref{thm:ast-lin-sub-equivs:h})
  $\Rightarrow$
  (\ref{thm:ast-lin-sub-equivs:g}):
}
If $M=\lb{-\ebar}{\ebar}$ for some icon $\ebar\in\corezn$,
then $M=\conv{E}$ where $E=\{-\ebar,\ebar\}$.

\pfpart{%
  (\ref{thm:ast-lin-sub-equivs:g})
  $\Rightarrow$
  (\ref{thm:ast-lin-sub-equivs:a}):
}
Suppose $M=\conv{E}$ where $E\subseteq\corezn$ is nonempty with
$E=-E$.
We argue that
$\zero\in\conv{E}$, which will imply $M=\conv\set{\set{\zero}\cup E}$, and hence that
$M$ is a convex astral cone (by
\Cref{thm:ast-cvx-cone-equiv}\ref{thm:ast-cvx-cone-equiv:b}\ref{thm:ast-cvx-cone-equiv:a}).
Let $\ebar$ be any point in $E$ (which must exist), implying $-\ebar$
is also in $E$.
Then
$
  \zero
  \in
  \lb{-\ebar}{\ebar}
  \subseteq
  \conv{E}
$.
The first inclusion follows from
\Cref{pr:seg-simplify}(\ref{pr:seg-simplify:b},\ref{pr:seg-simplify:a})
since,
for all $\uu\in\Rn$,
$\zero\cdot\uu=0\leq \max\{-\ebar\cdot\uu,\ebar\cdot\uu\}$.
The second inclusion is because $\ebar$ and $-\ebar$ are in
$E$, and so also in $\conv{E}$, which is convex.
Thus,
by \Cref{pr:conhull-prop}(\ref{pr:conhull-prop:c}),
$M=\conv(\{\zero\}\cup E)$.
Therefore, $M$ is a convex astral cone by
\Crefequiv{thm:ast-cvx-cone-equiv}{thm:ast-cvx-cone-equiv:b}{thm:ast-cvx-cone-equiv:a}.

Next, we have
\[
  -M
  =
  -(\conv{E})
  =
  \conv(-E)
  =
  \conv{E}
  =
  M,
\]
where the second equality follows from
\Cref{cor:thm:e:9}
(applied with the linear map $\xbar\mapsto-\xbar$).
Thus, $M$ is closed under negation.
Therefore, $M$ is an astral linear subspace.

\pfpart{%
  (\ref{thm:ast-lin-sub-equivs:b})
  $\Rightarrow$
  (\ref{thm:ast-lin-sub-equivs:d}):
}
Suppose $M=\Lbar$ for some linear subspace $L\subseteq\Rn$.
Let $U=\{\uu_1,\ldots,\uu_m\}\subseteq \Rn$ be a basis for
$\Lperp$ so that $\Lperp=\spn{U}$.
Let $\B=\trans{[\uu_1,\ldots,\uu_m]}$; that is, the $i$-th row of
matrix $\B\in\Rmn$ is $\utrans_i$ for $i=1,\ldots,m$.
Then
\[
   L
   =
   \Lperperp
   =
   \Uperp
   =
   \{ \xx\in\Rn :\: \B \xx = \zero \},
\]
where the first and second equalities are by
\Cref{pr:std-perp-props}(\ref{pr:std-perp-props:c},\ref{pr:std-perp-props:d}),
and the third is by the definitions of $\B$ and $\Uperp$.
Thus, by \Cref{cor:closure-affine-set},
\[
   M
   =
   \Lbar
   =
   \{ \xbar\in\extspace :\: \B \xbar = \zero \}.
\]

\pfpart{%
  (\ref{thm:ast-lin-sub-equivs:d})
  $\Rightarrow$
  (\ref{thm:ast-lin-sub-equivs:e}):
}
Suppose $M=\{\xbar\in\extspace :\: \B\xbar = \zero\}$
for some matrix $\B\in\Rmn$.
Then $\transB = [\uu_1,\ldots,\uu_m]$
for some $\uu_1,\ldots,\uu_m\in\Rn$
(that is, $\trans{\uu}_1,\ldots,\trans{\uu}_m$ are the rows of $\B$).
Let $U=\{\uu_1,\ldots,\uu_m\}$.
Then for all $\xbar\in\extspace$,
by \Cref{pr:mat-prod-is-row-prods-if-finite},
$\B\xbar=\zero$ if and only if
$\xbar\cdot\uu_i=0$ for $i=1,\ldots,m$.
Therefore,
$M=\{\xbar\in\extspace :\: \xbar\cdot\uu=0 \mbox{ for all } \uu\in U\}$.

\pfpart{%
  (\ref{thm:ast-lin-sub-equivs:e})
  $\Rightarrow$
  (\ref{thm:ast-lin-sub-equivs:a}):
}
Suppose
$M=\{\xbar\in\extspace :\: \xbar\cdot\uu=0 \mbox{ for all } \uu\in U\}$,
for some $U\subseteq\Rn$.
Then we can write
\[
  M
  =
  \bigcap_{\uu\in U}
  \BigParens{
       \bigBraces{\xbar\in\extspace :\: \xbar\cdot\uu \leq 0}
       \cap
       \bigBraces{\xbar\in\extspace :\: \xbar\cdot\uu \geq 0}
  }.
\]
Each set appearing in this intersection is either
a homogeneous closed astral halfspace or all of $\extspace$
(if $\uu=\zero$).
Therefore, $M$ is
a convex astral cone
by \Cref{pr:astral-cone-props}(\ref{pr:astral-cone-props:e}).
Furthermore, $M=-M$ since clearly $\xbar\cdot\uu=0$ if and only if
$-\xbar\cdot\uu=0$, for all $\xbar\in\extspace$ and $\uu\in U$.
Therefore, $M$ is an astral linear subspace.%
\indexg{linear subspaces, astral!characterizations|)}%
\qedhere
\end{proof-parts}
\end{proof}

\indexg{linear subspaces, astral!correspondence with linear subspaces|(}%
\indexg{linear subspaces (standard)!correspondence with astral linear subspaces|(}%
The characterization of astral linear subspaces given in
part~(\ref{thm:ast-lin-sub-equivs:b}) of
\Cref{thm:ast-lin-sub-equivs}
implies a one-to-one correspondence between astral linear subspaces and
standard linear subspaces.
According to this correspondence,
every astral linear subspace $M\subseteq\extspace$ is
mapped to a standard linear subspace by taking its intersection with
$\Rn$, that is, via the map
$M\mapsto M\cap\Rn$.
Likewise, every linear subspace $L\subseteq\Rn$ is mapped to an astral
linear subspace via the closure operation $L\mapsto \Lbar$.
As expressed in the next corollary,
these operations are bijective and are inverses of one
another.

An analogous correspondence was proved earlier,
in \Cref{cor:std-ast-polyhedra},
for astral and standard polyhedral sets.
Indeed, these correspondences are closely related since
every linear subspace is (standard) polyhedral,
and every astral linear subspace is astral polyhedral
(as can be shown to follow from
part~(\ref{thm:ast-lin-sub-equivs:d}) of
\Cref{thm:ast-lin-sub-equivs}).

\begin{corollary}   \label{cor:std-ast-linsub-corr}
  ~

  \begin{letter-compact}
  \item   \label{cor:std-ast-linsub-corr:L}
    Let $L\subseteq\Rn$ be a linear subspace.
    Then $\Lbar$ is an astral linear subspace, and moreover,
    $\Lbar\cap\Rn=L$.
  \item   \label{cor:std-ast-linsub-corr:M}
    Let $M\subseteq\extspace$ be an astral linear subspace.
    Then $M\cap\Rn$ is a (standard) linear subspace, and moreover,
    $\clbar{M\cap\Rn}=M$.
  \end{letter-compact}
\end{corollary}

\begin{proof}
  ~

\begin{proof-parts}
\pfpart{Part~(\ref{cor:std-ast-linsub-corr:L}):}
\Crefequiv{thm:ast-lin-sub-equivs}{thm:ast-lin-sub-equivs:b}{thm:ast-lin-sub-equivs:a}
shows that $\Lbar$ is an astral linear subspace.
Further, by
\Cref{pr:closed-set-facts}(\ref{pr:closed-set-facts:a}),
$\Lbar\cap\Rn=\cl{L}=L$ since $L$ is closed in $\Rn$.

\pfpart{Part~(\ref{cor:std-ast-linsub-corr:M}):}
By \Crefequiv{thm:ast-lin-sub-equivs}{thm:ast-lin-sub-equivs:a}{thm:ast-lin-sub-equivs:b},
since $M$ is an astral linear subspace, there exists a linear subspace
$L\subseteq\Rn$ such that $M=\Lbar=\repcl{L}$.
Consequently, $M\cap\Rn=\repcl{L}\cap\Rn=L$.
Thus, $M\cap\Rn$ is a linear subspace with closure
$\clbar{M\cap\Rn}=\Lbar=M$.%
\indexg{linear subspaces, astral!correspondence with linear subspaces|)}%
\indexg{linear subspaces (standard)!correspondence with astral linear subspaces|)}%
\qedhere
\end{proof-parts}
\end{proof}

The iconic and finite parts of elements of an astral linear subspace must
be themselves in the same astral linear subspace:

\begin{proposition} \label{pr:aspan:decomp}
Let $L\subseteq\Rn$ be a linear subspace, and
let $\xbar=\ebar\plusl\qq$ for some $\ebar\in\corezn$ and $\qq\in\Rn$.
Then $\xbar\in\Lbar$
if and only if $\ebar\in\Lbar$ and $\qq\in L$.
\end{proposition}
\begin{proof}
  By \Cref{thm:ast-lin-sub-equivs}(\ref{thm:ast-lin-sub-equivs:b},\ref{thm:ast-lin-sub-equivs:d}),
  there exists a matrix $\B\in\Rmn$, $m\ge 0$, such that
  $\Lbar=\{\xbar\in\extspace:\:\B\xbar=\zero\}$.

  Suppose that $\ebar\in\Lbar$ and $\qq\in L\subseteq\Lbar$.
  Then $\B\ebar=\B\qq=\zero$, and hence also $\B\xbar=\B\ebar\plusl\B\qq=\zero$.
  Thus, $\xbar\in\Lbar$.

  For the converse, suppose $\xbar\in\Lbar$, so
  $\zero=\B\xbar=\B\ebar\plusl\B\qq$. Since $\B\ebar$ is an icon
  (\Cref{pr:i:8}\ref{pr:i:8-matprod}) and $\B\qq\in\Rn$, this implies
  that $\B\ebar=\zero$
  (by uniqueness of iconic parts, \Cref{thm:icon-fin-decomp}).
  Hence, $\B\qq=\zero$ as well.
  Thus, $\ebar\in\Lbar$ and $\qq\in\Lbar\cap\Rn=L$
  (with the last equality by \Cref{cor:std-ast-linsub-corr}\ref{cor:std-ast-linsub-corr:L}).
\end{proof}

Analogous to standard linear subspaces, intersections and sequential sums of astral linear subspaces yield astral linear subspaces:

\begin{proposition}   \label{pr:ast-lin-sub-clos-props}
  ~

  \begin{letter-compact}
  \item   \label{pr:ast-lin-sub-clos-props:a}
    The intersection of an arbitrary collection of astral linear
    subspaces is an astral linear subspace.
  \item   \label{pr:ast-lin-sub-clos-props:b}
\indexg{sequential sum!astral linear subspaces@of astral linear subspaces|(}%
    Let $M_1$ and $M_2$ be astral linear subspaces.
    Then $M_1\seqsum M_2$ is also an astral linear subspace.
  \end{letter-compact}
\end{proposition}

\begin{proof}
~

\begin{proof-parts}
\pfpart{Part~(\ref{pr:ast-lin-sub-clos-props:a}):}
Let $M=\cap_{i\in I} M_i$, where $I$ is an arbitrary index set and each $M_i\subseteq\eRn$ is an astral linear subspace. Since each $M_i$ is a convex astral cone, so must be their intersection $M$
(by Propositions~\ref{pr:e1}\ref{pr:e1:b} and~\ref{pr:astral-cone-props}\ref{pr:astral-cone-props:b}).
The set $M$ is also closed under negation, since $\xbar\in M$ if and only
if $\xbar\in M_i$ for all $i\in I$, which is the case if and only if $-\xbar\in M_i$ for all $i\in I$, which in turn holds if and only if $-\xbar\in M$.
Thus, $M$ is an astral linear subspace.

\pfpart{Part~(\ref{pr:ast-lin-sub-clos-props:b}):}
Let $M=M_1\seqsum M_2$.
Since $M_1$ and $M_2$ are convex astral cones,
$M$ is as well,
by \Cref{thm:seqsum-ast-cone}(\ref{thm:seqsum-ast-cone:b}).
Also,
\[
  -M
  =
  -(M_1\seqsum M_2)
  =
  (-M_1)\seqsum (-M_2)
  =
  M_1 \seqsum M_2
  =
  M.
\]
Here, the second equality is by
\Cref{thm:distrib-seqsum}
(applied with $\A=-\Iden$),
and the third is because $M_1$ and $M_2$ are closed under negation.
Thus, $M$ is also closed under negation, and hence $M$ is an astral linear subspace.%
\indexg{sequential sum!astral linear subspaces@of astral linear subspaces|)}%
\indexg{linear subspaces, astral|)}%
\qedhere
\end{proof-parts}
\end{proof}

\section{Astral span and bases}
\label{sec:astral-span}

\indexg{span, astral|(}%
In light of
\Cref{pr:ast-lin-sub-clos-props}(\ref{pr:ast-lin-sub-clos-props:a}),
we can define astral span as a hull operator, analogous to the standard
span in $\Rn$:

\begin{definition}
\indexg{span, astral!defined|(}%
The \emph{astral span} of a set $S\subseteq\extspace$,
denoted
\indexm{span s700}{$\aspan{S}$}{astral span}%
$\aspan{S}$,
is the smallest astral linear subspace that includes $S$, or
equivalently, the intersection of all astral linear subspaces that
include $S$.%
\indexg{span, astral!defined|)}%
\end{definition}

Astral span has the usual properties of a hull operator as stated in
\Cref{pr:gen-hull-ops}:

\begin{proposition}  \label{pr:aspan-hull-props}
  Let $S,U\subseteq\extspace$.
  \begin{letter-compact}
  \item  \label{pr:aspan-hull-props:a}
    If $S\subseteq U$ and $U$ is an astral linear subspace, then
    $\aspan{S}\subseteq U$.
  \item  \label{pr:aspan-hull-props:b}
    If $S\subseteq U$ then $\aspan{S}\subseteq\aspan{U}$.
  \item  \label{pr:aspan-hull-props:c}
    If $S\subseteq U\subseteq\aspan{S}$, then
    $\aspan{U}=\aspan{S}$.
  \end{letter-compact}
\end{proposition}

\indexg{span, astral!conic or convex hull@as conic or convex hull|(}%
The astral span of any set can always be re-expressed using
the astral conic hull or convex hull operations:

\begin{proposition}  \label{pr:aspan-ito-acone-conv}
  Let $S\subseteq\extspace$.
  Then
  \begin{equation}   \label{eq:pr:aspan-ito-acone-conv:1}
    \aspan{S}
    =
    \acone(S\cup-S)
    =
    \conv\regParens{\{\zero\}\cup\lmset{S}\cup-\lmset{S}}.
  \end{equation}
  If $S\neq\emptyset$, then also
  \begin{equation}   \label{eq:pr:aspan-ito-acone-conv:2}
    \aspan{S}
    =
    \conv\regParens{\lmset{S}\cup-\lmset{S}}.
  \end{equation}
\end{proposition}

\begin{proof}
We begin with the first equality of
\eqref{eq:pr:aspan-ito-acone-conv:1}.
The astral linear subspace $\aspan{S}$ is a convex astral cone that
includes $S$, and so also $-S$, being closed under negation.
Therefore,
$\acone(S\cup-S) \subseteq \aspan{S}$
(\Cref{pr:acone-hull-props}\ref{pr:acone-hull-props:a}).
On the other hand, $\acone(S\cup-S)$ includes $S$ and is an astral
linear subspace by
\Crefequiv{thm:ast-lin-sub-equivs}{thm:ast-lin-sub-equivs:f}{thm:ast-lin-sub-equivs:a}.
Thus,
$\aspan{S} \subseteq \acone(S\cup-S)$.

The second equality of
\eqref{eq:pr:aspan-ito-acone-conv:1}
is by
\Cref{thm:acone-char}(\ref{thm:acone-char:a}).

Suppose henceforth that $S$ is not empty.
Then $\conv(\lmset{S}\cup-\lmset{S})$ is an astral linear subspace
by
\Crefequiv{thm:ast-lin-sub-equivs}{thm:ast-lin-sub-equivs:g}{thm:ast-lin-sub-equivs:a}
and so includes the origin
(for instance, by \Cref{pr:ast-cone-is-naive}).
Therefore,
$\conv\regParens{\lmset{S}\cup-\lmset{S}} =
\conv\regParens{\{\zero\}\cup\lmset{S}\cup-\lmset{S}}$
by \Cref{pr:conhull-prop}(\ref{pr:conhull-prop:c}),
proving
\eqref{eq:pr:aspan-ito-acone-conv:2}.%
\indexg{span, astral!conic or convex hull@as conic or convex hull|)}%
\indexg{span, astral|)}%
\end{proof}

\indexg{representational span|(}%
In
\Cref{sec:strong-equiv-span-bound},
we encountered the representational span of a
point $\xbar\in\extspace$,
$\rspanxbar$, which we defined as the span (in $\Rn$)
of vectors in $\xbar$'s representation. We next extend this definition
to obtain representational span of an arbitrary set $S\subseteq\eRn$
and show how this concept relates to the astral span of $S$.

\begin{definition}   \label{def:rspan-gen}
\indexg{representational span!general definition|(}%
Let $S\subseteq\extspace$.
For each $\xbar\in S$, let
$\xbar = \Vxbar\, \omm \plusl \qqxbar$ be any
representation of $\xbar$ where
$\Vxbar\in\R^{n\times \kxbar}$, $\kxbar\geq 0$,
and $\qqxbar\in\Rn$.
Also, let
$ \Rxbar=(\columns{\Vxbar})\cup\{\qqxbar\} $
be the vectors appearing in $\xbar$'s representation,
and let $R=\bigcup_{\xbar\in S} \Rxbar$
be the union of all of these.
Then the \emph{representational span} of $S$, denoted
\indexm{rspan800}{$\rspan{S}$}{representational span (general)}%
$\rspan{S}$, is
equal to $\spn{R}$,
the (standard) span of $R$.%
\indexg{representational span!general definition|)}%
\end{definition}

Clearly, the earlier definition for singletons (\Cref{def:rep-span-snglton})
is a special case of the one given here. Moreover,
the representational span of a set can be decomposed as
\begin{equation}
\mathtogether%
\!\!
 \rspan{S}
 =\spn\BiggParens{\bigcup_{\xbar\in S}\Rxbar\!}\!
 =\spn\BiggParens{\bigcup_{\xbar\in S}\spn\Rxbar\!}\!
\label{eq:rspan:alt}
 =\spn\BiggParens{\bigcup_{\xbar\in S}\rspanxbar\!},
\end{equation}
where the second equality follows by \Cref{pr:hull:union},
using the fact that standard span is a hull operator.
As we noted earlier, the representational span of a singleton $\set{\xbar}$
is independent of
the particular representation chosen for $\xbar$
(\Cref{pr:rspan-sing-equiv-dual}). Thus, \eqref{eq:rspan:alt} shows that the
representational span of a set
is likewise independent of the particular representations chosen for individual points.

\indexg{representational span!astral span and|(}%
\indexg{span, astral!representational span and|(}%
\indexg{span, astral|(}%
We next show that the representational span of any set
$S\subseteq\extspace$ is equal to the intersection of the astral span of $S$
and $\Rn$.
In other words, $\rspan{S}$ and $\aspan{S}$ are related to each other
as the sets appearing in
\Cref{cor:std-ast-linsub-corr}, implying further that
$\aspan{S}$ is the astral closure of $\rspan{S}$.

\begin{theorem}  \label{thm:rspan-is-aspan-in-rn}
  Let $S\subseteq\extspace$, and assume the notation of
  \Cref{def:rspan-gen}.
  Then
  \begin{equation}   \label{eq:thm:rspan-is-aspan-in-rn:1}
    \rspan{S}=\spn{R}=(\aspan{S})\cap\Rn,
  \end{equation}
  and
  \begin{equation}   \label{eq:thm:rspan-is-aspan-in-rn:2}
    \clbar{\rspan{S}}
    =
    \clbar{\spn{R}}
    =
    \aspan{S}.
  \end{equation}
\end{theorem}

\begin{proof}
The first equality in \eqref{eq:thm:rspan-is-aspan-in-rn:1} is by
definition of $\rspan{S}$.
For the second equality,
let
$L=\spn{R}$
and
$L'=(\aspan{S})\cap\Rn$; we aim to prove $L=L'$.

We show first that $L'\subseteq L$.
Let $\xbar\in S$.
Then $\Rxbar\subseteq R\subseteq\spn{R}=L$, and since $\xbar$ can be written
as
$\xbar=\Vxbar\omm\plusl\qqxbar$ where the columns of $\Vxbar$ as well as the
vector $\qqxbar$ are in $\Rxbar$, we obtain $\xbar\in\repcl{L}$.
Thus, $S\subseteq\repcl{L}$.
Since $\repcl{L}$ is an astral linear subspace
(by
\Cref{thm:ast-lin-sub-equivs}\ref{thm:ast-lin-sub-equivs:b}\ref{thm:ast-lin-sub-equivs:a}),
it follows that $\aspan{S}\subseteq\repcl{L}$
(by \Cref{pr:aspan-hull-props}\ref{pr:aspan-hull-props:a}).
Thus, $L'=(\aspan{S})\cap\Rn\subseteq\repcl{L}\cap\Rn=L$.

For the reverse inclusion, let $\xbar\in S$.
Then
\begin{align*}
  \Rxbar
  \subseteq
  \spn{\Rxbar}
  \subseteq
  \clbar{\spn{\Rxbar}}
  &=
  \lb{-\limray{\xbar}}{\limray{\xbar}}
  \\
  &=
  \conv\{-\limray{\xbar},\limray{\xbar}\}
  =
  \aspan\{\xbar\}
  \subseteq
  \aspan{S}.
\end{align*}
The first equality is by \Cref{lem:ast-subspace-equivs},
applied with $V=\Rxbar$ and $\VV=[\Vxbar,\qqxbar]$ so that
$\ebar=\VV\omm=\limray{\xbar}$ by
\Cref{pr:mult-rep-by-scalar}.
The third equality is by
\Cref{pr:aspan-ito-acone-conv}.
And the final inclusion is by
\Cref{pr:aspan-hull-props}(\ref{pr:aspan-hull-props:b}).
Since this holds for all $\xbar\in S$, it follows that
$R\subseteq(\aspan{S})\cap\Rn=L'$.
Therefore,
$L=\spn{R}\subseteq L'$
since $L'$ is a linear subspace
(by
\Cref{cor:std-ast-linsub-corr}\ref{cor:std-ast-linsub-corr:M}).

This completes the proof of
\eqref{eq:thm:rspan-is-aspan-in-rn:1}.
Taking closures then yields
\eqref{eq:thm:rspan-is-aspan-in-rn:2},
again by
\Cref{cor:std-ast-linsub-corr}(\ref{cor:std-ast-linsub-corr:M}).%
\indexg{representational span!astral span and|)}%
\indexg{span, astral!representational span and|)}%
\end{proof}

\indexg{representational span!sets in rn@for sets in $\Rn$|(}%
\indexg{span, astral!sets in rn@for sets in $\Rn$|(}%
When $S$ is a subset of $\Rn$, we can express its astral span
and representational span
in terms of its standard span:

\begin{theorem}   \label{pr:aspan-for-rn}
  Let $S\subseteq\Rn$.
  Then
  \[ \aspan{S}=\clbar{\spn{S}}=\repcl{(\spn{S})}, \]
  and
  \[ \rspan{S}=\spn{S}. \]
\end{theorem}

\begin{proof}
We have
\begin{equation}
\label{eq:aspan-for-rn}
  \aspan{S}
  =
  \acone(S\cup-S)
  =
  \repcl{\bigParens{\cone(S\cup-S)}}
  =
  \repcl{(\spn{S})}
  =
  \clbar{\spn{S}}.
\end{equation}
The first equality is by
\Cref{pr:aspan-ito-acone-conv},
the second is by
\Cref{thm:acone-is-repclos},
the third is by
\Cref{pr:scc-cone-elts}(\ref{pr:scc-cone-elts:span}),
and the fourth is by \Cref{cor:lin-sub-bar-is-repcl}.
Thus,
\[
  \rspan{S}=(\aspan{S})\cap\Rn=\clbar{\spn{S}}\cap\Rn=\spn{S},
\]
where the equalities follow, respectively, by
\Cref{thm:rspan-is-aspan-in-rn},
\eqref{eq:aspan-for-rn}, and
\Cref{cor:std-ast-linsub-corr}(\ref{cor:std-ast-linsub-corr:L}).%
\indexg{representational span|)}%
\indexg{representational span!sets in rn@for sets in $\Rn$|)}%
\indexg{span, astral!sets in rn@for sets in $\Rn$|)}%
\end{proof}

There is an intuitive sense in which the (standard) span of a nonempty
set of vectors $S\subseteq\Rn$ is like the convex hull of those
vectors and their negations as they are extended to infinity.
Indeed, taking the convex hull of scalar multiples of $S$ and $-S$
will, in the limit, encompass all of the linear subspace spanned by $S$;
that is,
\[
  \spn{S}
  =
  \bigcup_{\lambda\in\Rstrictpos} \conv(\lambda S \cup -\lambda S).
\]
In astral space, this can be expressed without resorting to
limits. Specifically, as we saw in \Cref{pr:aspan-ito-acone-conv},
$\aspan{S}=\conv\regParens{\lmset{S}\cup-\lmset{S}}$,
that is, the astral span of $S$ is equal to the convex hull of all the
astrons associated with $S$ and $-S$.
Consequently,
using \Cref{pr:aspan-for-rn} combined with \Cref{thm:rspan-is-aspan-in-rn},
the standard span can be written as
\[
  \spn{S}=\conv\regParens{\lmset{S}\cup-\lmset{S}}\cap\Rn.
\]

\indexg{span, astral!union@of union|(}%
We next derive an identity for the astral span of a union using a similar
identity for the astral conic hull of a union (\Cref{thm:decomp-acone}):

\begin{theorem}   \label{thm:decomp-aspan}
  Let $S_1,S_2\subseteq\extspace$.
  Then
  \[
     (\aspan{S_1}) \seqsum (\aspan{S_2})
     =
     \aspan(S_1\cup S_2).
  \]
\end{theorem}

\begin{proof}
We have:
\begin{align*}
  (\aspan{S_1}) \seqsum (\aspan{S_2})
  &=
  \acone(S_1\cup -S_1) \seqsum \acone(S_2\cup -S_2)
  \\
  &=
  \acone(S_1\cup -S_1 \cup S_2\cup -S_2)
  \\
  &=
  \aspan(S_1\cup S_2).
\end{align*}
The first and third equalities are by
\Cref{pr:aspan-ito-acone-conv}.
The second equality is by \Cref{thm:decomp-acone}.%
\indexg{span, astral|)}%
\indexg{span, astral!union@of union|)}%
\end{proof}

\indexg{basis (for linear subspace)|(}%
A set $B\subseteq\Rn$ is a basis for a linear subspace $L\subseteq\Rn$
if $B$ spans $L$ and if $B$ is linearly
\indexg{basis (for linear subspace)|)}%
independent.
\indexg{basis, astral|(}%
\indexg{astral basis|(}%
In the same way, we say that $B\subseteq\Rn$ is an \emph{astral basis}
(or simply a \emph{basis}) for an astral linear subspace
$M\subseteq\extspace$ if $M=\aspan{B}$ and $B$ is linearly
independent.
As we show next, $B$ is an astral basis for $M$ if and only if it is a
standard basis for the linear subspace $M\cap\Rn$.

We also give another characterization for when a set is an astral basis.
In standard linear algebra, a set $B$ is a basis for the linear
subspace $L$ if the coefficients defining the linear combinations of
$B$ are in one-to-one correspondence with the elements of $L$.
More precisely,
in matrix form,
the linear combinations of $B$ are exactly the vectors
$\B\zz$ where $\B$ is a matrix with columns matching the $d$ vectors
in $B$, and $\zz\in\R^d$.
Then $B$ is a basis if and only if the map
$\zz\mapsto\B\zz$ is a bijection from $\R^d$ onto $L$.
Likewise, as we show next, $B$ is an astral basis for $M$ if and only
if the corresponding astral map $\zbar\mapsto\B\zbar$ is a bijection
from $\extspac{d}$ onto $M$.

\begin{theorem}   \label{thm:basis-equiv}
  Let $B=\{\vv_1,\ldots,\vv_d\}\subseteq\Rn$
  where $d=|B|<+\infty$, and let
  $\B=[\vv_1,\ldots,\vv_d]$.
  Let $M\subseteq\extspace$ be an astral linear subspace.
  Then the following are equivalent:
  \begin{letter-compact}
  \item   \label{thm:basis-equiv:a}
    $B$ is an astral basis for $M$.
  \item   \label{thm:basis-equiv:b}
    $B$ is a (standard) basis for the linear subspace $M\cap\Rn$.
  \item   \label{thm:basis-equiv:c}
    The map $\zbar\mapsto\B\zbar$ defines a bijection from
    $\extspac{d}$ onto $M$.
  \end{letter-compact}
  Consequently, an astral basis must exist for every astral
  linear subspace $M$.
  Furthermore, every astral basis for $M$ must have cardinality equal
  to $\dim(M\cap\Rn)$.
\end{theorem}

\begin{proof}
Let $L=M\cap\Rn$, which is a linear subspace with $M=\Lbar$
(by
\Cref{cor:std-ast-linsub-corr}\ref{cor:std-ast-linsub-corr:M}).

\begin{proof-parts}
\pfpart{%
  (\ref{thm:basis-equiv:b})
  $\Rightarrow$
  (\ref{thm:basis-equiv:a}):
}
Suppose $B$ is a standard basis for $L$.
Then $B$ is linearly independent, and
\[
  M
  =
  \Lbar
  =
  \clbar{\spn{B}}
  =
  \aspan{B},
\]
where the second equality is because $B$ is a basis for $L$, and the
third is by \Cref{pr:aspan-for-rn}.
Thus, $B$ is an astral basis for $M$.

\pfpart{%
  (\ref{thm:basis-equiv:a})
  $\Rightarrow$
  (\ref{thm:basis-equiv:c}):
}
Suppose $B$ is an astral basis for $M$.
Then
\begin{equation}   \label{eq:thm:basis-equiv:1}
  M
  =
  \aspan{B}
  =
  \clbar{\spn{B}}
  =
  \acolspace{\B},
\end{equation}
where the second equality is by
\Cref{pr:aspan-for-rn},
and the third is by
\Cref{lem:ast-subspace-equivs}.
Let $\Fbar:\extspac{d}\rightarrow M$ be the astral linear map
associated with $\B$ so that
$\Fbar(\zbar)=\B\zbar$ for $\zbar\in\extspac{d}$.
Note that the restriction of $\Fbar$'s range to $M$ is justified
by \eqref{eq:thm:basis-equiv:1}.
Indeed, $\Fbar(\extspac{d})$ is exactly the astral column space of
$\B$;
therefore, \eqref{eq:thm:basis-equiv:1} also shows that
$\Fbar$ is surjective.

To see that $\Fbar$ is injective, suppose $\Fbar(\zbar)=\Fbar(\zbar')$
for some $\zbar,\zbar'\in\extspac{d}$.
Then
\[
  \zbar
  =
  \Bpseudoinv \B\zbar
  =
  \Bpseudoinv \B\zbar'
  =
  \zbar',
\]
where $\Bpseudoinv$ is $\B$'s pseudoinverse.
The second equality is because, by assumption, $\B\zbar=\B\zbar'$.
The first and third equalities follow from the identity
$\Bpseudoinv\B=\Iden$
(where $\Iden$ is the $d\times d$ identity matrix),
which holds by
\Cref{pr:pseudoinv-props}(\ref{pr:pseudoinv-props:c}) since
$B$ is linearly independent (being a basis), implying that $\B$ has
full column rank.

Thus, $\Fbar$ is a bijection, as claimed.

\pfpart{%
  (\ref{thm:basis-equiv:c})
  $\Rightarrow$
  (\ref{thm:basis-equiv:b}):
}
Suppose $\Fbar$, as defined above, is a bijection.
Let $F:\R^d\rightarrow L$ be the restriction of $\Fbar$ to $\R^d$ so
that $F(\zz)=\Fbar(\zz)=\B\zz$ for $\zz\in\R^d$.
Note that the restriction of $F$'s range to $L$ is justified because,
for $\zz\in\R^d$,
$\Fbar(\zz)=\B\zz$ is in $M\cap\Rn=L$.

We claim $F$ is a bijection.
To see this, observe first that $F$ is injective since if
$F(\zz)=F(\zz')$ for some $\zz,\zz'\in\R^d$ then
$\Fbar(\zz)=F(\zz)=F(\zz')=\Fbar(\zz')$, implying
$\zz=\zz'$ since $\Fbar$ is injective.

To see that $F$ is surjective, let $\xx\in L$.
Then because $\Fbar$ is surjective, there exists $\zbar\in\extspac{d}$
for which $\Fbar(\zbar)=\xx$.
We can write $\zbar=\ebar\plusl\qq$ for some icon $\ebar\in\corez{d}$
and some $\qq\in\R^d$.
Thus, $\xx=\B\zbar=\B\ebar\plusl \B\qq$.
Evidently, $\B\ebar$, which is an icon
(by \Cref{pr:i:8}\ref{pr:i:8-matprod}),
is the iconic part of the
expression on the right-hand side of this equality,
while $\zero$ is the iconic part of $\xx$, the expression on the left.
Since the iconic part of any point is unique
(\Cref{thm:icon-fin-decomp}),
it follows that $\Fbar(\ebar)=\B\ebar=\zero$.
Since $\Fbar(\zero)=\zero$ and since $\Fbar$ is injective,
this implies that $\ebar=\zero$.
Thus, $\zbar=\qq$ and $F(\qq)=\Fbar(\zbar)=\xx$,
so $F$ is surjective.

Since $F$ is surjective, $\spn{B}=\colspace{\B}=F(\Rn)=L$,
so $B$ spans $L$.
And since $F$ is injective, $B$ must be linearly independent.
This is because if some linear combination of $B$ is $\zero$,
meaning $\B\zz=\zero$ for some $\zz\in\R^d$, then
$F(\zz)=\zero=F(\zero)$ implying $\zz=\zero$ since $F$ is
injective.
Thus, $B$ is a basis for $L$.

\pfpart{Existence and equal cardinality:}
The linear subspace $L=M\cap\Rn$ must have a basis, and every basis for
$L$ must have cardinality $\dim{L}$.
Since, as shown above,
a set $B$ is an astral basis for $M$ if and only if it is a standard
basis for $L$, the same properties hold for $M$.
\qedhere
\end{proof-parts}
\end{proof}

Thus, every astral linear subspace $M\subseteq\extspace$
must have an astral basis,
and the cardinality of all such bases must be
the same, namely, $\dim(M\cap\Rn)$, the usual dimension of the
standard linear subspace $M\cap\Rn$.
\indexg{dimension!astral linear subspace@of astral linear subspace|(}%
\indexg{linear subspaces, astral!dimension of|(}%
Accordingly, we define the \emph{dimension} of $M$, denoted
\indexm{dim m}{$\dim{M}$}{dimension (astral)}%
$\dim{M}$,
to be the cardinality of any basis for $M$.
Thus, in general, $\dim{M}=\dim(M\cap\Rn)$.%
\indexg{dimension!astral linear subspace@of astral linear subspace|)}%
\indexg{linear subspaces, astral!dimension of|)}%
\indexg{astral basis|)}%
\indexg{basis, astral|)}

\indexg{span, astral|(}%
We show next that every set $S\subseteq\extspace$ includes a finite
subset $V$ whose astral span is the same as that of the full set $S$;
moreover, it is always possible to choose a set $V$ whose cardinality does
not exceed the dimension of $\aspan{S}$.

\begin{theorem}   \label{thm:aspan-finite-subset}
  Let $S\subseteq\extspace$, and let $d=\dim(\aspan{S})$.
  Then there exists a subset $V\subseteq S$ with $|V|\leq d$ such that
  $\aspan{V}=\aspan{S}$.
\end{theorem}

\begin{proof}
Let $\Rxbar$ and $R$ be as in
\Cref{def:rspan-gen}
(for $\xbar\in S$), implying that
$\rspan{S}=\spn{R}$.
Let $B$ be a linearly independent subset of $R$ of maximum
cardinality among all such subsets.
Then $B$ is of a finite cardinality, because
$\card{B}\le\dim(\spn{R})\le n$,
and every vector in $R$ is equal to a linear combination of the
vectors in $B$
(otherwise, the cardinality of $B$ could be increased).
Thus, $B\subseteq R\subseteq \spn{B}$, so $\spn{R}=\spn{B}$
(by \Cref{pr:gen-hull-ops}\ref{pr:gen-hull-ops:d}
since linear span is a hull operator).

This shows that $B$, being linearly independent, is a basis for
$\spn{R}$, which, by
\Cref{thm:rspan-is-aspan-in-rn},
is the same as $(\aspan{S})\cap\Rn$.
Therefore, $B$ is an astral basis for $\aspan{S}$ as well by
\Crefequiv{thm:basis-equiv}{thm:basis-equiv:b}{thm:basis-equiv:a}.
Consequently, $|B|=d$.

Let $\bb_1,\ldots,\bb_d$ be the elements of $B$.
Since each $\bb_i\in R$, there must exist $\xbar_i\in S$ such that
$\bb_i\in\Rxbari$, for $i=1,\ldots,d$.
Let $V=\{\xbar_1,\ldots,\xbar_d\}$, whose cardinality is at most $d$,
and let
\[ R'=\bigcup_{i=1}^d \Rxbari. \]
Then $B\subseteq R'\subseteq R$ so
\begin{equation}   \label{eq:thm:aspan-finite-subset:1}
  \spn{B}
  \subseteq
  \spn{R'}
  \subseteq
  \spn{R}
  =
  \spn{B}.
\end{equation}
Therefore,
\[
  \aspan{V}
  =
  \clbar{\spn{R'}}
  =
  \clbar{\spn{R}}
  =
  \aspan{S}.
\]
The second equality is because $\spn{R'}=\spn{R}$
by \eqref{eq:thm:aspan-finite-subset:1}.
The first and third equalities are by
\Cref{thm:rspan-is-aspan-in-rn}
(applied to $V$ and $S$ respectively).
\end{proof}

We previously saw
in Theorems~\ref{thm:e:7} and~\ref{thm:oconic-hull-and-seqs}
how the outer convex hull and outer conic hull of a
finite set of points can be characterized in terms of sequences.
\indexg{span, astral!sequential characterization|(}%
The astral span of a finite set can be similarly characterized,
as shown next.
Indeed, this follows simply from the fact that the astral span of any
set can always be expressed as a conic hull.

\begin{theorem}   \label{thm:aspan-and-seqs}
  Let $V=\{\xbar_1,\dotsc,\xbar_m\}\subseteq\extspace$,
  and let $\zbar\in\extspace$.
  Then $\zbar\in\aspan{V}$ if
  and only if there exist sequences
  $\seq{\lambda_{it}}$
  and
  $\seq{\lambda'_{it}}$
  in $\Rpos$,
  and span-bound sequences
  $\seq{\xx_{it}}$
  and
  $\seq{\xx'_{it}}$
  in $\Rn$,
  for $i=1,\dotsc,m$,
  such that:
  \begin{item-compact}
  \item
    $\xx_{it}\rightarrow\xbar_i$
    and
    $\xx'_{it}\rightarrow -\xbar_i$
    for $i=1,\dotsc,m$.
  \item
    The sequence %
    $\zz_t=\sum_{i=1}^m (\lambda_{it} \xx_{it} + \lambda'_{it} \xx'_{it})$
    converges to
    $\zbar$.
  \end{item-compact}
\end{theorem}

\begin{proof}
We have
\[ \aspan{V} = \acone(V\cup -V) = \oconich(V\cup -V), \]
where the first equality is 
by \Cref{pr:aspan-ito-acone-conv}, and the second
by \Cref{pr:acone-hull-props}(\ref{pr:acone-hull-props:d}).
Therefore, applying
\Cref{thm:oconic-hull-and-seqs}
to $V\cup -V$
yields the claim.
\end{proof}

The characterization given in
\Cref{thm:aspan-and-seqs}
requires twin sequences for every point $\xbar_i$ in $V$:
one sequence, $\seq{\xx_{it}}$, converging to $\xbar_i$,
and the other, $\seq{\xx'_{it}}$, converging to its negation, $-\xbar_i$.
The sequence elements $\zz_t$ are then conic (nonnegative) combinations of
the elements of all these sequences.

It seems reasonable to wonder if such twin sequences are really necessary,
or if it would suffice to have a \emph{single} sequence
$\seq{\xx_{it}}$ for each $\xbar_i$,
and then for each $\zz_t$ to be a \emph{linear}
combination of the $\xx_{it}\negKern$'s (with possibly negative
coefficients).
Thus, the characterization would say that $\zbar\in\aspan{V}$
if and only if there exist span-bound sequences $\seq{\xx_{it}}$
converging to $\xbar_i$ and (arbitrary) sequences $\seq{\lambda_{it}}$
in $\R$ such that $\zz_t=\sum_{i=1}^m \lambda_{it} \xx_{it}$
converges to $\zbar$.
The existence of such sequences is indeed sufficient for $\zbar$ to be
in $\aspan{V}$, as can be shown to follow from
\Cref{thm:aspan-and-seqs}.
However, it is not necessary, as the next example shows:

\begin{example}
In $\extspac{2}$,
let $V=\{\xbar\}$ where $\xbar=\limray{\ee_1}\plusl\limray{\ee_2}$,
and let $\zbar=\zz=\ee_2$.
Then $\aspan{V}=\extspac{2}$
(as follows, for instance, from
\Cref{thm:rspan-is-aspan-in-rn}).
In particular, $\zz\in\aspan{V}$.
Suppose there exist sequences as described above, namely,
a sequence $\seq{\lambda_t}$ in $\R$ and
a span-bound sequence $\seq{\xx_t}$ in $\R^2$
such that $\xx_t\rightarrow\xbar$
and $\zz_t=\lambda_t \xx_t \rightarrow \zz$.
We must then have that either $\lambda_t\geq 0$ for infinitely many
values of $t$,
or that $\lambda_t\leq 0$ for infinitely many~$t$ (or both).
Suppose the first case holds.
Then by discarding all other sequence elements, we can assume
$\lambda_t\geq 0$ for all $t$.
By \Cref{cor:aray-and-seqs},
the existence of such
sequences implies that
$\zz\in\aray{\xbar}=\lb{\zero}{\limray{\xbar}}$.
However, this is a contradiction since $\zz$ is not in this set
(by \Cref{thm:lb-with-zero}).
In the alternative case that $\lambda_t\leq 0$ for infinitely many
$t$, a similar argument shows that $\zz\in\aray{(-\xbar)}$, which is
again a contradiction.

Thus, no such sequences can exist, even though $\zz\in\aspan{V}$.

To construct twin sequences as in
\Cref{thm:aspan-and-seqs},
we can choose
$\xx_t=\trans{[t^2,t+1]}$,
$\xx'_t=\trans{[-t^2,-t]}$,
and
$\lambda_t=\lambda'_t=1$,
for all $t$.
Then $\xx_t\rightarrow\xbar$,
$\xx'_t\rightarrow -\xbar$,
and
$\zz_t=\lambda_t\xx_t+\lambda'_t\xx'_t=\ee_2\rightarrow\zz$.%
\indexg{span, astral!sequential characterization|)}%
\indexg{span, astral|)}%
\end{example}

\section{Orthogonal complements}
\label{sec:ortho-compl}

\indexg{orthogonality|(}%
Orthogonality between pairs of points in $\Rn$ naturally extends to
pairs in which one point is in $\Rn$ and the other in $\extspace$.
Thus, for $\xbar\in\extspace$ and $\uu\in\Rn$, if $\xbar\cdot\uu=0$, then
we say that $\uu$ is \emph{orthogonal} to $\xbar$,
or that $\xbar$ is orthogonal to $\uu$,
written
\indexm{x s200 u}{$\xbar\perp\uu$}{orthogonal to}%
$\uu\perp\xbar$ or $\xbar\perp\uu$.%
\indexg{orthogonality|)}

Using this extended notion of orthogonality we can now also extend the notion of orthogonal complement.
Recall that,
as given in \eqref{eq:std-ortho-comp-defn},
the orthogonal complement of a set $S\subseteq\Rn$ is the set
\begin{equation}  \label{eq:std-ortho-comp-defn-rev}
   \Sperp
   =
   \regBraces{\uu\in\Rn :\: \xx\cdot\uu = 0 \mbox{ for all } \xx\in S }.
\end{equation}
As was done with the polarity operation in
\Cref{sec:ast-pol-cones},
we extend orthogonal
complement to astral space in two ways.
\indexg{orthogonal complement (primal)|(}%
\indexg{orthogonal complement (primal)!redefined for astral sets|(}%
In the first of these, we simply allow $S$ to be any subset of
$\extspace$:

\begin{definition}
The \emph{(primal) orthogonal complement} of a set
$S\subseteq\extspace$, denoted $\Sperp$,
is the set
\begin{equation}  \label{eqn:sperp-def}
\indexm{s 500}{$\Sperp$}{(primal) orthogonal complement}%
  \Sperp = \regBraces{\uu\in\Rn :\: \xbar\cdot\uu=0 \mbox{ for all } \xbar\in S}.
\end{equation}
\end{definition}
Thus, $\Sperp$ is the set of vectors in $\Rn$ that are orthogonal to
every point in $S\subseteq\extspace$.
Clearly, for $S\subseteq\Rn$, this definition is consistent with the old
one that it generalizes.%
\indexg{orthogonal complement (primal)|)}%
\indexg{orthogonal complement (primal)!redefined for astral sets|)}

\indexg{orthogonal complement, dual|(}%
In the second extension to astral space, we reverse the roles of $\Rn$
and $\extspace$:
\begin{definition}
\indexg{orthogonal complement, dual!defined|(}%
The \emph{dual orthogonal complement} of a set $U\subseteq\Rn$,
denoted $\aperp{U}$, is the set
\begin{equation}  \label{eq:aperp-defn}
  \indexm{s 550}{$\aperp{S}$}{dual orthogonal complement}%
  \aperp{U}
  =
  \regBraces{\xbar\in\extspace :\: \xbar\cdot\uu=0 \mbox{ for all } \uu\in U}.%
\indexg{orthogonal complement, dual!defined|)}%
\end{equation}
\end{definition}
Thus, $\aperp{U}$ is the
set of points in $\extspace$ that are orthogonal to all points in
$U\subseteq\Rn$.

\indexg{orthogonality|(}%
Implicitly, these two definitions extend the notion of orthogonality to sets: a vector $\uu\in\Rn$ is said to be \emph{orthogonal} to a set $S\subseteq\eRn$ if it is orthogonal to every point in~$S$, and symmetrically, $\xbar\in\eRn$ is said
to be \emph{orthogonal} to $U\subseteq\Rn$ if it is orthogonal to every point in $U$. With this terminology, $\Sperp$ is the set of vectors in $\Rn$ that are orthogonal to $S$, and $\aperp{U}$ is the set of points in $\eRn$ that are orthogonal to $U$.%
\indexg{orthogonality|)}

As with the polarity operations defined in
\Cref{sec:ast-pol-cones},
we will see that the orthogonal complement operations
form a duality relationship.
One operation maps any set
$S\subseteq\extspace$ to a set $\Sperp$ in $\Rn$, while the other maps
any set $U\subseteq\Rn$ to a set $\Uaperp$ in $\extspace$.

The next propositions summarize various properties of these operations:

\begin{proposition}  \label{pr:aperp-props}
  Let $S,U\subseteq\Rn$.
  \begin{letter-compact}
  \item  \label{pr:aperp-props:d}
    $\Uaperp$ is an astral linear subspace.
  \item  \label{pr:aperp-props:a}
    If $S\subseteq U$ then $\Uaperp\subseteq\Saperp$.
  \item  \label{pr:aperp-props:b}
    $\Uaperp=\aperp{(\cl{U})}=\aperp{(\spn{U})}$.
  \item  \label{pr:aperp-props:c}
    $\aperp{U}\subseteq\apol{U}$ with equality
    if $U$ is a linear subspace.
  \item  \label{pr:aperp-props:e}
    $\Uaperp=\clbar{\Uperp}$.
  \end{letter-compact}
\end{proposition}

\begin{proof}
  ~

\begin{proof-parts}
\pfpart{Part~(\ref{pr:aperp-props:d}):}
This is immediate from
\Crefequiv{thm:ast-lin-sub-equivs}{thm:ast-lin-sub-equivs:e}{thm:ast-lin-sub-equivs:a}
and \eqref{eq:aperp-defn}.

\pfpart{Part~(\ref{pr:aperp-props:a}):}
Let $S\subseteq U$.
If $\xbar\in\Uaperp$, then $\xbar\cdot\uu=0$ for all $\uu\in U$, and
therefore also for all $\uu\in S$.
Thus, $\uu\in\Saperp$.

\pfpart{Part~(\ref{pr:aperp-props:b}):}
We show first that $\Uaperp\subseteq\aperp{(\spn{U})}$.
Suppose $\xbar\in\Uaperp$,
and let $\uu\in\spn{U}$.
Then $\uu=\sum_{i=1}^m \lambda_i \uu_i$ for some
$\lambda_1,\ldots,\lambda_m\in\R$ and
$\uu_1,\ldots,\uu_m\in U$.
For $i=1,\ldots,m$, we then have $\xbar\cdot\uu_i=0$,
implying that
$\xbar\cdot(\lambda_i\uu_i)=\lambda_i(\xbar\cdot\uu_i)=0$,
and so that
\[
  \xbar\cdot\uu
  =
  \xbar\cdot\BiggParens{\sum_{i=1}^m \lambda_i \uu_i}
  =
  \sum_{i=1}^m \xbar\cdot(\lambda_i\uu_i)
  =
  0
\]
where the second equality is by \Cref{pr:i:1}.
Since this holds for all $\uu\in\spn{U}$, it follows that
$\xbar\in\aperp{(\spn{U})}$.

Thus,
\[
   \Uaperp
   \subseteq
   \aperp{(\spn{U})}
   \subseteq
   \aperp{(\cl{U})}
   \subseteq
   \Uaperp.
\]
The first inclusion was just shown.
The second and third inclusions follow from
part~(\ref{pr:aperp-props:a})
since $U\subseteq \cl{U} \subseteq\spn{U}$
(since the span of any set is closed).

\pfpart{Part~(\ref{pr:aperp-props:c}):}
Suppose $\xbar\in\aperp{U}$.
Then for all $\uu\in U$,
$\xbar\cdot\uu=0\leq 0$.
Thus, $\xbar\in\apol{U}$, so
$\aperp{U}\subseteq\apol{U}$.

Next, suppose $U$ is a linear subspace.
For the reverse inclusion under this additional assumption,
suppose $\xbar\in\apol{U}$.
Then for all $\uu\in U$,
$\xbar\cdot\uu\leq 0$, implying that
$-\xbar\cdot\uu=\xbar\cdot(-\uu)\leq 0$ since $-\uu$ is also in $U$.
Thus, $\xbar\cdot\uu=0$ so $\xbar\in\aperp{U}$,
proving
$\apol{U}\subseteq\aperp{U}$.

\pfpart{Part~(\ref{pr:aperp-props:e}):}
Let $L=\spn{U}$, which is a linear subspace and therefore also
a closed (in $\Rn$) convex cone.
We then have
\[
  \Uaperp
  =
  \Laperp
  =
  \Lapol
  =
  \Lpolbar
  =
  \clbar{\Lperp}
  =
  \clbar{\Uperp}.
\]
The first equality is by
part~(\ref{pr:aperp-props:b}).
The second is by
part~(\ref{pr:aperp-props:c}).
The third is by
\Cref{thm:apol-of-closure}(\ref{thm:apol-of-closure:b}).
The fourth is by
\Cref{pr:polar-props}(\ref{pr:polar-props:f}).
And the last is by
\Cref{pr:std-perp-props}(\ref{pr:std-perp-props:d}).%
\indexg{orthogonal complement, dual|)}%
\qedhere
\end{proof-parts}
\end{proof}

\begin{proposition}  \label{pr:perp-props-new}
\indexg{orthogonal complement (primal)|(}%
  Let $S, U\subseteq\extspace$.
  \begin{letter-compact}
  \item  \label{pr:perp-props-new:a}
    $\Sperp$ is a linear subspace of $\Rn$.
  \item  \label{pr:perp-props-new:b}
    If $S\subseteq U$ then $\Uperp\subseteq\Sperp$.
  \item  \label{pr:perp-props-new:c}
    $\Sperp = \Sbarperp = (\aspan{S})^\bot = (\rspan{S})^\bot$.
  \item  \label{pr:perp-props-new:d}
    $\Sperp\subseteq\polar{S}$ with equality
    if $S$ is an astral linear subspace.
  \end{letter-compact}
\end{proposition}

\begin{proof}
  ~

\begin{proof-parts}
\pfpart{Part~(\ref{pr:perp-props-new:a}):}
Let $L=\Sperp$.
Clearly, $\zero\in L$.
If $\uu,\vv\in L$, then for all $\xbar\in S$,
$\xbar\cdot\uu=\xbar\cdot\vv=0$ so
$\xbar\cdot(\uu+\vv)=0$ by \Cref{pr:i:1}; therefore,
$\uu+\vv\in L$.
Finally, if $\uu\in L$ and $\lambda\in\R$, then
for all $\xbar\in S$,
$\xbar\cdot(\lambda\uu)=\lambda(\xbar\cdot\uu)=0$,
by \Cref{pr:i:2}, so $\lambda\uu\in L$.
Thus, $L$ is a linear subspace.

\pfpart{Part~(\ref{pr:perp-props-new:b}):}
Proof is similar to that of
\Cref{pr:aperp-props}(\ref{pr:aperp-props:a}).

\pfpart{Part~(\ref{pr:perp-props-new:c}):}
We show first that $\Sperp\subseteq(\aspan{S})^\bot$.
Let $\uu\in\Sperp$.
This means $\uu$ is orthogonal to every point in $S$, and thus that
$S\subseteq C$ where $C=\aperp{\{\uu\}}$.
By \Cref{pr:aperp-props}(\ref{pr:aperp-props:d}),
$C$ is an astral linear subspace, so
$\aspan{S}\subseteq C$
(\Cref{pr:aspan-hull-props}\ref{pr:aspan-hull-props:a}),
that is,  every point in $\aspan{S}$ is orthogonal to $\uu$.
Therefore, $\uu\in(\aspan{S})^\bot$.

We thus have
\[
   \Sperp
   \subseteq
   (\aspan{S})^\bot
   \subseteq
   \Sbarperp
   \subseteq
   \Sperp
\]
with the first inclusion following from the argument above,
and the other inclusions from
part~(\ref{pr:perp-props-new:b})
since $S\subseteq \Sbar\subseteq\aspan{S}$
(since astral linear subspaces are closed,
for instance, by \Cref{thm:ast-lin-sub-equivs}\ref{thm:ast-lin-sub-equivs:a}\ref{thm:ast-lin-sub-equivs:b}).
This shows that
$\Sperp = \Sbarperp = (\aspan{S})^\bot$.

Replacing $S$ with $\rspan{S}$, this same argument shows that
$(\rspan{S})^\bot = (\clbar{\rspan{S}})^\bot=(\aspan{S})^\bot$
since $\aspan{S} = \clbar{\rspan{S}}$
(by \Cref{thm:rspan-is-aspan-in-rn}).

\pfpart{Part~(\ref{pr:perp-props-new:d}):}
The proof is analogous to that of
\Cref{pr:aperp-props}(\ref{pr:aperp-props:c}).
\qedhere
\end{proof-parts}
\end{proof}

\indexg{orthogonal complement, dual|(}%
We next consider the effect of applying orthogonal complement
operations twice, generalizing the standard result that
$\Uperperp=\spn{U}$ for $U\subseteq\Rn$
(\Cref{pr:std-perp-props}\ref{pr:std-perp-props:c}):

\begin{theorem}  \label{thm:dub-perp}
  Let $S\subseteq\extspace$ and let $U\subseteq\Rn$.
  \begin{letter-compact}
  \item    \label{thm:dub-perp:d}
    $\Uperpaperp = \clbar{\spn{U}}$.
    Thus, if $U$ is a linear subspace, then
    $\Uperpaperp = \Ubar$.
  \item    \label{thm:dub-perp:c}
    $\Uaperperp = \spn{U}$.
    Thus, if $U$ is a linear subspace, then
    $\Uaperperp = U$.
  \item    \label{thm:dub-perp:a}
    $\Sperpaperp = \aspan{S}$.
    Thus, if $S$ is an astral linear subspace, then
    $\Sperpaperp = S$.
  \item    \label{thm:dub-perp:b}
    $\Sperperp = (\aspan{S})\cap\Rn=\rspan{S}$.
  \end{letter-compact}
\end{theorem}

\begin{proof}
  ~

\begin{proof-parts}
\pfpart{Part~(\ref{thm:dub-perp:d}):}
We have
\[
   \Uperpaperp
   =
   \clbar{\Uperperp}
   =
   \clbar{\spn{U}},
\]
with equalities following by
Propositions~\ref{pr:aperp-props}(\ref{pr:aperp-props:e})
and~\ref{pr:std-perp-props}(\ref{pr:std-perp-props:c}).

\pfpart{Part~(\ref{thm:dub-perp:c}):}
We have
\[
  \Uaperperp
  =
  \aperperp{(\spn{U})}
  =
  \bigParens{\clbar{(\spn{U})^\bot}}^\bot
  =
  \perperp{(\spn{U})}
  =
  \spn{U},
\]
where the equalities follow, respectively,
by Propositions~\ref{pr:aperp-props}(\ref{pr:aperp-props:b}),
\ref{pr:aperp-props}(\ref{pr:aperp-props:e}),
\ref{pr:perp-props-new}(\ref{pr:perp-props-new:c}),
and
\ref{pr:std-perp-props}(\ref{pr:std-perp-props:c}).

\pfpart{Part~(\ref{thm:dub-perp:a}):}
Let $L=\rspan{S}$.
Then
\[
  \aspan{S}
  =
  \Lbar
  =
  \Lperpaperp
  =
  \Sperpaperp.
\]
The first equality is from
\Cref{thm:rspan-is-aspan-in-rn}.
The second is by
part~(\ref{thm:dub-perp:d})
since $L$ is a linear subspace.
And the third is by
\Cref{pr:perp-props-new}(\ref{pr:perp-props-new:c}).

\pfpart{Part~(\ref{thm:dub-perp:b}):}
We have
\[
  \Sperperp
  =
  \Sperpaperp \cap \Rn
  =
  (\aspan{S}) \cap \Rn
  =
  \rspan{S},
\]
where the first equality follows from definitions
(Eqs.~\ref{eq:std-ortho-comp-defn-rev} and~\ref{eq:aperp-defn}),
the second is by part~(\ref{thm:dub-perp:a}),
the third is by
\Cref{thm:rspan-is-aspan-in-rn}.
\qedhere
\end{proof-parts}
\end{proof}

\indexg{linear subspaces, astral!correspondence with linear subspaces|(}%
\indexg{linear subspaces (standard)!correspondence with astral linear subspaces|(}%
Parts~(\ref{thm:dub-perp:c}) and~(\ref{thm:dub-perp:a})
of \Cref{thm:dub-perp}
show that the orthogonal complement operations define a one-to-one
correspondence between astral linear subspaces and standard linear
subspaces with $M\mapsto \Mperp$
mapping astral linear subspaces $M\subseteq\extspace$
to linear subspaces in $\Rn$, and
$L\mapsto \Laperp$ mapping standard linear subspaces $L\subseteq\Rn$
to astral linear subspaces.
The theorem shows that these operations are inverses of each other,
and therefore are bijective.%
\indexg{linear subspaces, astral!correspondence with linear subspaces|)}%
\indexg{linear subspaces (standard)!correspondence with astral linear subspaces|)}%
\indexg{orthogonal complement, dual|)}

In \Cref{pr:rspan-sing-equiv-dual}, we showed that elements $\zz$
of $\rspanxbar$ are characterized by the property that for any $\uu\in\Rn$,
$\uu\perp\xbar$ implies $\uu\perp\zz$.
\Cref{thm:dub-perp}(\ref{thm:dub-perp:b}) can be rephrased to
say that $\zz\in\rspan{S}$ if and only if,
for all $\uu\in\Rn$, $\uu\perp S$ implies $\uu\perp\zz$. Thus,
\Cref{pr:rspan-sing-equiv-dual} can be viewed as a special case of this theorem.

In $\Rn$, we have $L\cap\Lperp=\set{\zero}$ for any linear subspace $L$.
The next \namecref{pr:aperp:intersect}
states an astral analogue:

\begin{theorem}
\label{pr:aperp:intersect}
  Let $S\subseteq\extspace$. Then $(\aspan S)\cap\clbar{\Sperp}=\set{\zero}$.
\end{theorem}

\begin{proof}
Let $M=\aspan{S}$, and let $L=M\cap\Rn$. Then $L$ is a linear subspace
and $M=\Lbar$ (\Cref{cor:std-ast-linsub-corr}\ref{cor:std-ast-linsub-corr:M}).
Moreover, 
$\Sperp=\Mperp=(\Lbar)^\perp=\Lperp$,
with the first and third equalities from
\Cref{pr:perp-props-new}(\ref{pr:perp-props-new:c}).
Thus, we need to show that $\Lbar\cap\clbar{\Lperp}=\set{\zero}$.

First, note that $\zero\in L$ and $\zero\in\Lperp$, so $\zero\in\Lbar\cap\clbar{\Lperp}$.
Next,
let $\xbar\in\Lbar\cap\clbar{\Lperp}$.
By \Cref{pr:aperp-props}(\ref{pr:aperp-props:e}),
$\xbar\in\clbar{\Lperp}=\Laperp$,
so
$\xbar\inprod\uu=0$ for all $\uu\in L$.
Also, by \Cref{thm:dub-perp}(\ref{thm:dub-perp:d}),
$\xbar\in\Lbar=\aperp{(\Lperp)}$,
so
$\xbar\inprod\uu=0$ for all $\uu\in\Lperp$.
Combining then yields
$\xbar\in\aperp{(L\cup\Lperp)}=\aperp{[\spn(L\cup\Lperp)]}=\aperp{(\Rn)}$
(also using \Cref{pr:aperp-props}\ref{pr:aperp-props:b}).
Thus, $\xbar\inprod\uu=0$ for all $\uu\in\Rn$, so $\xbar=\zero$.%
\indexg{orthogonal complement (primal)|)}%
\end{proof}

\part{Functions on Astral Space}
\label{part:astral-functions}

\chapter{Convex functions}
\label{sec:conv:fct}

We next define and study what it means for an astral
function $F:\extspace\rightarrow\Rext$ to be convex.
As in standard convex analysis, we define this property
in terms of the convexity of the function's epigraph, $\epi{F}$, for
whose study we are well-prepared, having closely considered astral
convex sets in previous chapters.
Convexity of standard functions is also commonly characterized in
terms of a function's values, as in \Cref{pr:stand-cvx-fcn-char}.
We provide a similar
characterization for astral convex functions.
We also look at various natural operations for constructing astral
convex functions, for instance, by adding two convex functions or by
composing a convex function with an affine function.
We
then specifically
explore functions on $\extspace$ that are derived
from a set in $\extspacnp$
(similarly to how a function on $\Rn$ can be
naturally derived from its own epigraph).

\section{Definition and properties}

As in \Cref{sec:work-with-epis}, for $\xbar\in\extspace$ and
$y\in\R$, we continue to use the notation $\rpair{\xbar}{y}$ to
denote, depending on context, either a pair in $\eRn\times\R$ or a
point in $\extspacnp$
obtained by embedding $\eRn\times\R$ in $\extspacnp$ as in \Cref{thm:homf}.
Likewise,
we regard subsets of $\extspace\times\R$, such as the
epigraph of an astral function, as subsets of $\extspacnp$.
We also make use of the component projection matrices $\PPx=[\Idnn,\zero_n]$ and
$\PPy=[\trans{\zero_n},1]$
introduced in \Cref{sec:work-with-epis}.
Their properties and the properties of the elements of $\eRn\times\R$
(viewed as elements of $\eRf{n+1}$)
were summarized in \Cref{pr:xy-pairs-props};
in particular,
$\PPx\rpair{\xbar}{y}=\xbar$
and
$\PPy\rpair{\xbar}{y}=y$.

Recall that a function $f:\Rn\rightarrow\Rext$ is convex if its
epigraph, which is a subset of $\Rn\times\R=\Rnp$, is convex.
\indexg{convex functions, astral|(}%
We define astral convex functions analogously:

\begin{definition}
\indexg{convex functions, astral!defined|(}%
  An astral function
  $F:\extspace\rightarrow\Rext$ is \emph{convex} if its epigraph,
  $\epi{F}$, is convex as a subset of $\extspacnp$.%
\indexg{convex functions, astral!defined|)}%
\end{definition}

\indexg{lower semicontinuous extension!convexity of|(}%
As a first observation,
the extension $\fext$ of any convex function
$f:\Rn\rightarrow\Rext$ is convex:

\begin{theorem}  \label{thm:fext-convex}
  Let $f:\Rn\rightarrow\Rext$ be convex.
  Then $\fext$, its lower semicontinuous extension, is also convex.
\end{theorem}

\begin{proof}
Since $f$ is convex, its epigraph, $\epi{f}$, is a convex
subset of $\R^{n+1}$,
so $\epibarbar{f}$, its closure in $\extspacnp$, is convex by
\Cref{thm:e:6}.
Furthermore, by \Cref{pr:wasthm:e:3}(\ref{pr:wasthm:e:3b}), $\epi{\fext}=\epibarbar{f}\cap\finclset$,
where, as in \eqref{eq:finclset-defn},
\begin{align*}
  \finclset &=
       \{\zbar\in\extspacnp :\: \PPy\zbar\in\R\}
       =\lmPPyinv(\R),
\end{align*}
and where $\PPy=[\trans{\zero_n},1]$, and
$\lmPPy:\extspacnp\rightarrow\Rext$ is the associated astral linear
map, $\lmPPy(\zbar)=\PPy\zbar$ for $\zbar\in\extspacnp$.
Since $\R$ is a convex subset of $\eR$, also $\finclset$ is convex
(by \Cref{cor:inv-image-convex}).
Thus, $\epi{\fext}=\epibarbar{f}\cap\finclset$ is convex as well
(by \Cref{pr:e1}\ref{pr:e1:b}). Therefore,
$\fext$ is a convex function.%
\indexg{lower semicontinuous extension!convexity of|)}%
\end{proof}

In studying standard convexity of functions, it is often convenient
to work with a characterization like the one in
\Cref{pr:stand-cvx-fcn-char}(\ref{pr:stand-cvx-fcn-char:b}),
which states that $f:\Rn\rightarrow\Rext$ is convex if and only if
\begin{equation}  \label{eqn:stand-f-cvx-ineq}
   f\bigParens{(1-\lambda)\xx_0 + \lambda \xx_1}
   \leq
   (1-\lambda) f(\xx_0) + \lambda f(\xx_1)
\end{equation}
for all $\xx_0,\xx_1\in\dom{f}$ and all $\lambda\in [0,1]$.
Indeed, when working exclusively with proper functions,
\eqref{eqn:stand-f-cvx-ineq}
is sometimes presented as a \emph{definition} of convexity.

\indexg{convex functions, astral!characterization|(}%
We can similarly characterize convex astral functions.
To do so, in
\eqref{eqn:stand-f-cvx-ineq},
we first replace $f$ with an astral
function $F:\extspace\rightarrow\Rext$,
and $\xx_0,\xx_1$ with astral points
$\xbar_0,\xbar_1\in\dom{F}\subseteq\extspace$.
In this way, the right-hand side straightforwardly generalizes to
astral functions.
On the left-hand side, a similarly direct substitution of
$(1-\lambda)\xx_0+\lambda\xx_1$ with
$(1-\lambda)\xbar_0+\lambda\xbar_1$ would be problematic, since such
addition of astral points is undefined.
Rather, we apply notions developed in
\Cref{sec:seq-cvx-conic-combs},
replacing the (standard) convex combination
$(1-\lambda)\xx_0 + \lambda \xx_1$
by any element of the sequential convex combination
$\mul{1-\lambda}\xbar_0\seqsum\mul{\lambda}\xbar_1$,
that is, any $\lambda$-midpoint of $\xbar_0$ and $\xbar_1$.
We then require that the resulting inequality hold for all such
$\lambda$-midpoints.
This generalizes \eqref{eqn:stand-f-cvx-ineq}
since
${(1-\lambda)\xx_0} + {\lambda \xx_1}$ is the sole
$\lambda$-midpoint of $\xx_0$ and $\xx_1$
when $\xx_0,\xx_1\in\Rn$ 
(\Cref{pr:seq-conic:prop}\ref{i:seq-conic:finite}).

This yields the following characterization of convexity for astral
functions:

\begin{theorem}     \label{thm:ast-F-char-fcn-vals}
  Let $F:\extspace\rightarrow\Rext$.
  Then the following are equivalent:
  \begin{letter-compact}
    \item        \label{thm:ast-F-char-fcn-vals:a}
      $F$ is convex.
    \item        \label{thm:ast-F-char-fcn-vals:b}
      For all $\xbar_0,\xbar_1\in\dom{F}$,
      for all $\lambda\in [0,1]$,
      and for all
      $\xbar\in\mul{1-\lambda}\xbar_0\seqsum\mul{\lambda}\xbar_1$,
      \begin{equation}   \label{eq:thm:ast-F-char-fcn-vals:1}
         F(\xbar)
         \leq
         (1-\lambda)  F(\xbar_0)
         +
         \lambda F(\xbar_1).
      \end{equation}
  \end{letter-compact}
\end{theorem}

\indexg{midsegments and midpoints (sequential)|(}%
Before proving this theorem, we give in the next
\namecref{pr:lam-mid-props}
some general rules for
working with $\lambda$-mipoints that will be used in its proof:

\begin{proposition}   \label{pr:lam-mid-props}
  Let $\xbar_0,\xbar_1\in\extspace$,
  let $\lambda\in[0,1]$, and
  let
  $\xbar\in\mul{1-\lambda}\xbar_0\seqsum\mul{\lambda}\xbar_1$.

  \begin{letter-compact}
  \item   \label{pr:lam-mid-props:aff}
    Let $\A\in\R^{m\times n}$, let $\bbar\in\extspac{m}$,
    and let $F:\extspace\rightarrow\extspac{m}$ be defined by
    $F(\zbar)=\bbar\plusl\A\zbar$ for $\zbar\in\extspace$.
    Then
    \[
      F(\xbar)
      \in
      \flammid{F(\xbar_0)}{F(\xbar_1)}.
    \]
  \item   \label{pr:lam-mid-props:c}
    Let $\yy_0,\yy_1\in\Rn$ and
    let $\yy=(1-\lambda) \yy_0 + \lambda \yy_1$.
    Then
    \[
       \xbar\plusl\yy
        \in
        \flammid{(\xbar_0\plusl\yy_0)}{(\xbar_1\plusl\yy_1)}.
    \]
  \item   \label{pr:lam-mid-props:d}
    Let $y_0,y_1\in\R$ and
    let $y=(1-\lambda) y_0 + \lambda y_1$.
    Then
    \[
      \rpair{\xbar}{y}
      \in
      \flammid{\rpair{\xbar_0}{y_0}}{\rpair{\xbar_1}{y_1}}.
    \]
  \end{letter-compact}
\end{proposition}

\begin{proof}
~
\begin{proof-parts}
\pfpart{Part~(\ref{pr:lam-mid-props:aff}):}
By \Cref{thm:F:seqconv},
\[
   F(\xbar)
   \in
   F\bigParens{\mul{1-\lambda}\xbar_0\seqsum\mul{\lambda}\xbar_1}
   =
   \mul{1-\lambda}F(\xbar_0)\seqsum\mul{\lambda}F(\xbar_1).
\]

\pfpart{Part~(\ref{pr:lam-mid-props:c}):}
We have
\begin{align*}
  \set{\xbar\plusl\yy}
  =
  \xbar\seqsum\yy
  &\subseteq
  \bigBracks{\mul{1-\lambda}{\xbar_0} \seqsum \mul{\lambda}{\xbar_1}}
  \seqsum
  \bigBracks{(1-\lambda)\yy_0 \seqsum \lambda\yy_1}
\\
  &=
  \bigBracks{\mul{1-\lambda}{\xbar_0} \seqsum (1-\lambda)\yy_0} \seqsum
  \bigBracks{\mul{\lambda}{\xbar_1} \seqsum \lambda\yy_1}
\\
  &=
  \bigBracks{\mul{1-\lambda}{\xbar_0} \plusl (1-\lambda)\yy_0} \seqsum
  \bigBracks{\mul{\lambda}{\xbar_1} \plusl \lambda\yy_1}
\\
  &=
  \mul{1-\lambda}{(\xbar_0\plusl\yy_0)}\seqsum\mul{\lambda}{(\xbar_1\plusl\yy_1)}.
\end{align*}
The first and third equalities are by
\Cref{cor:seqsum-conseqs}(\ref{cor:seqsum-conseqs:c}).
The inclusion is by
\Cref{pr:seqsum-props}(\ref{pr:seqsum-props:finite},\ref{pr:seqsum-props:b}).
The second equality is by commutativity and associativity of
sequential addition.
The fourth equality is by
\Cref{pr:seq:mul}(\ref{i:seq:mul:sum:yy}).

\pfpart{Part~(\ref{pr:lam-mid-props:d}):}
As earlier,
let $\PPx=[\Idnn,\zero_n]$ and $\PPy=[\trans{\zero_n},1]$.
By part~(\ref{pr:lam-mid-props:aff}),
$\trans{\PPx}\xbar \in \flammid{\regBracks{\trans{\PPx}\xbar_0}}{\regBracks{\trans{\PPx}\xbar_1}}$.
Therefore,
\begin{align*}
  \rpair{\xbar}{y}
  =
  \trans{\PPx} \xbar \plusl \trans{\PPy}y
  &\in
  \mul{1-\lambda}\Parens{\trans{\PPx} \xbar_0 \plusl \trans{\PPy}y_0}
  \seqsum
  \mul{\lambda}\Parens{\trans{\PPx} \xbar_1 \plusl \trans{\PPy}y_1}
  \\
  &=
  \flammid{\rpair{\xbar_0}{y_0}}{\rpair{\xbar_1}{y_1}},
\end{align*}
with the inclusion following from part~(\ref{pr:lam-mid-props:c})
(and since
$\trans{\PPy}y=(1-\lambda) \trans{\PPy}y_0 + \lambda \trans{\PPy}y_1
$),
and both equalities from
\Cref{pr:xy-pairs-props}(\ref{pr:xy-pairs-props:a}).%
\indexg{midsegments and midpoints (sequential)|)}%
\qedhere
\end{proof-parts}
\end{proof}

\begin{proof}[Proof of \Cref{thm:ast-F-char-fcn-vals}]
~

\begin{proof-parts}
\pfpart{%
  (\ref{thm:ast-F-char-fcn-vals:a})
  $\Rightarrow$
  (\ref{thm:ast-F-char-fcn-vals:b}):
}
Suppose $F$ is convex, and let
$\xbar_0,\xbar_1\in\dom{F}$,
$\lambda\in [0,1]$,
and $\xbar\in\flammid{\xbar_0}{\xbar_1}$.
We aim to prove \eqref{eq:thm:ast-F-char-fcn-vals:1}.
For $i\in\{0,1\}$, let $y_i\in\R$ with $y_i\geq F(\xbar_i)$
so that $\rpair{\xbar_i}{y_i}\in\epi{F}$.
Let $y=(1-\lambda) y_0 + \lambda y_1$.
Then
\[
  \rpair{\xbar}{y}
  \in
  \flammid{\rpair{\xbar_0}{y_0}}{\rpair{\xbar_1}{y_1}}
  \subseteq
  \lb{\rpair{\xbar_0}{y_0}}{\rpair{\xbar_1}{y_1}}
  \subseteq
  \epi{F}.
\]
The first inclusion is by
\Cref{pr:lam-mid-props}(\ref{pr:lam-mid-props:d}).
The second is by
\Cref{cor:conv-as-seqsum-midrays}.
The last is because $\epi{F}$ is convex
(since $F$ is).

Thus,
$\rpair{\xbar}{y}\in\epi{F}$,
so $F(\xbar)\leq y = (1-\lambda) y_0 + \lambda y_1$.
Since this holds for all $y_0\geq F(\xbar_0)$ and $y_1\geq F(\xbar_1)$,
this proves \eqref{eq:thm:ast-F-char-fcn-vals:1}.

\pfpart{%
  (\ref{thm:ast-F-char-fcn-vals:b})
  $\Rightarrow$
  (\ref{thm:ast-F-char-fcn-vals:a}):
}
Suppose statement~(\ref{thm:ast-F-char-fcn-vals:b}) holds.
Let $\rpair{\xbar_i}{y_i}\in\epi{F}$, for $i\in\{0,1\}$.
Let $\zbar\in\lb{\rpair{\xbar_0}{y_0}}{\rpair{\xbar_1}{y_1}}$.
To prove convexity,
we aim to show $\zbar\in\epi{F}$.

As usual,
let $\PPx=[\Idnn,\zero_n]$ and $\PPy=[\trans{\zero_n},1]$.
Then
\[
  \PPy\zbar
  \in\lb{\PPy\rpair{\xbar_0}{y_0}}{\PPy\rpair{\xbar_1}{y_1}}
  =\lb{y_0}{y_1}\subseteq\R,
\]
with the first inclusion by \Cref{cor:thm:e:9},
and the equality by
\Cref{pr:xy-pairs-props}(\ref{pr:xy-pairs-props:c}).
Thus, by \Cref{pr:xy-pairs-props}(\ref{pr:xy-pairs-props:d}),
$\zbar=\rpair{\xbar}{y}$ where $\xbar=\PPx\zbar\in\extspace$ and
$y=\PPy\zbar\in\R$.

Since $\zbar\in\lb{\rpair{\xbar_0}{y_0}}{\rpair{\xbar_1}{y_1}}$,
by \Cref{cor:conv-as-seqsum-midrays},
$\zbar$ must be a $\lambda$-midpoint of
$\rpair{\xbar_0}{y_0}$ and $\rpair{\xbar_1}{y_1}$,
for some $\lambda\in[0,1]$;
that is,
$\zbar\in\flammid{\rpair{\xbar_0}{y_0}}{\rpair{\xbar_1}{y_1}}$.
Thus,
\begin{align*}
  \xbar
  &=
  \PPx\zbar
  \in
  \flammid{[\PPx \rpair{\xbar_0}{y_0}]}{[\PPx \rpair{\xbar_1}{y_1}]}
  =
  \flammid{\xbar_0}{\xbar_1},
\\
  y
  &=
  \PPy\zbar
  \in
  \flammid{[\PPy \rpair{\xbar_0}{y_0}]}{[\PPy \rpair{\xbar_1}{y_1}]}
  =
  \flammid{y_0}{y_1},
\end{align*}
with each inclusion following from
\Cref{pr:lam-mid-props}(\ref{pr:lam-mid-props:aff}),
and the second equality in each line by
\Cref{pr:xy-pairs-props}(\ref{pr:xy-pairs-props:c}).
\Cref{pr:seq-conic:prop}(\ref{i:seq-conic:finite}) applied to the second line
then implies that
$y=(1-\lambda)y_0+\lambda y_1$.

Since $\xbar_0,\xbar_1\in\dom{F}$
and since $\xbar\in\flammid{\xbar_0}{\xbar_1}$,
\eqref{eq:thm:ast-F-char-fcn-vals:1} must hold so that
\[
  F(\xbar)
  \leq
  (1-\lambda) F(\xbar_0) + \lambda F(\xbar_1)
  \leq
  (1-\lambda) y_0 + \lambda y_1
  =
  y,
\]
where the second inequality is because
$\rpair{\xbar_i}{y_i}\in\epi{F}$ so that $F(\xbar_i)\le y_i$,
for $i\in\{0,1\}$.
Thus,
$\zbar=\rpair{\xbar}{y}\in\epi{F}$, completing the proof.%
\indexg{convex functions, astral!characterization|)}%
\qedhere
\end{proof-parts}
\end{proof}

Here is a simple corollary
of \Cref{thm:ast-F-char-fcn-vals}:

\begin{corollary}   \label{cor:F-conv-max-seg}
  Let $F:\extspace\rightarrow\Rext$ be convex.
  Let $\xbar,\ybar\in\extspace$, and let
  $\zbar\in\lb{\xbar}{\ybar}$.
  Then
  $F(\zbar)\leq \max\{F(\xbar),F(\ybar)\}$.
\end{corollary}

\begin{proof}
If either $F(\xbar)=+\infty$ or $F(\ybar)=+\infty$, then the
claim is immediate, so we assume $\xbar,\ybar\in\dom{F}$.
Then
by \Cref{cor:conv-as-seqsum-midrays},
$\zbar$ is a $\lambda$-midpoint of $\xbar,\ybar$, for some
$\lambda\in [0,1]$, implying, by
\Cref{thm:ast-F-char-fcn-vals}, that
\[
  F(\zbar)
  \leq
  (1-\lambda) F(\xbar)
  +
  \lambda F(\ybar)
  \leq
  \max\{ F(\xbar), F(\ybar) \}.
\qedhere
\]
\end{proof}

\indexg{convex functions, astral!sublevel sets of|(}%
\indexg{convex functions, astral!domain of|(}%
\indexg{sublevel sets!astral convex function@of astral convex function|(}%
\indexg{domain, effective!astral convex function@of astral convex function|(}%
Consequently,
as in standard convex analysis, the effective domain and all
sublevel sets of an astrally convex function are convex:

\begin{theorem}  \label{thm:f:9}
  Let $F:\extspace\rightarrow\Rext$ be convex,
  let $\alpha\in\Rext$,
  and let
  \begin{align*}
    S
    &=
    \regBraces{\xbar\in\extspace :\: F(\xbar) \leq \alpha},
    \\
    S'
    &=
    \regBraces{\xbar\in\extspace :\: F(\xbar) < \alpha}.
  \end{align*}
  Then $S$ and $S'$ are both convex.
  Consequently, $\dom F$,
  the effective domain of $F$, is convex as well.
\end{theorem}

\begin{proof}
Let $\xbar,\ybar\in S$, and let $\zbar\in\lb{\xbar}{\ybar}$.
Then, by \Cref{cor:F-conv-max-seg},
$F(\zbar)\leq\max\set{F(\xbar),F(\ybar)}\leq\alpha$,
so $\zbar\in S$.
Therefore, $S$ is convex, as claimed.
The proof that $S'$ is convex is similar.
In particular, taking $\alpha=+\infty$, this shows that
$\dom{F}$ is convex.%
\indexg{convex functions, astral!sublevel sets of|)}%
\indexg{convex functions, astral!domain of|)}%
\indexg{sublevel sets!astral convex function@of astral convex function|)}%
\indexg{domain, effective!astral convex function@of astral convex function|)}%
\end{proof}

In standard convex analysis,
if a convex, lower semicontinuous function $f:\Rn\rightarrow\Rext$
attains a limiting value less than $+\infty$ along any halfline
in the direction of some vector $\vv\in\Rn$, then the function must be nonincreasing
everywhere in that direction.
More precisely,
if $\liminf_{\lambda\rightarrow+\infty} f(\yy+\lambda\vv)<+\infty$ for any $\yy\in\Rn$,
then $f(\xx+\lambda\vv)\leq f(\xx)$ holds for all $\xx\in\Rn$ and all $\lambda\in\Rpos$
(see \Cref{pr:stan-rec-equiv}\ref{pr:stan-rec-equiv:c}\ref{pr:stan-rec-equiv:a}). The
latter condition can be restated as $f(\zz)\leq f(\xx)$ for all $\zz$ on the halfline $\hfline{\xx}{\vv}$,
or more succinctly,
\[
   \sup f\bigParens{\hfline{\xx}{\vv}} \leq f(\xx).
\]
\indexg{convex functions, astral!bounded along halfline|(}%
\indexg{halflines, astral!function bounded along|(}%
We next derive an analogue of this result for astrally convex
functions in terms of astral halflines
(as introduced in
\Cref{sec:hline}).
The standard result above can then be recovered by setting $F=\fext$, $\ybar=\yy$, $\vbar=\vv$, and $\xbar=\xx$.

\begin{theorem}   \label{thm:F-conv-res}
  Let $F:\extspace\rightarrow\Rext$ be convex,
  let $\vbar\in\extspace$,
  and assume $F(\limray{\vbar}\plusl\ybar)<+\infty$ for some
  $\ybar\in\extspace$.
  Then for all $\xbar\in\extspace$,
  \[
     \sup F\bigParens{\ahfline{\xbar}{\vbar}} \leq F(\xbar);
  \]
  that is, $F(\zbar)\leq F(\xbar)$ for all
  $\zbar\in\ahfline{\xbar}{\vbar}$.
\end{theorem}

\begin{proof}
Let $\xbar\in\extspace$,
and let $\zbar\in\ahfline{\xbar}{\vbar}$.
We aim to prove $F(\zbar)\leq F(\xbar)$.
If $F(\xbar)=+\infty$ then this holds trivially,
so we assume henceforth that $\xbar\in\dom{F}$.
We have
\[
  \zbar
  \in
  \ahfline{\xbar}{\vbar}
  =
  \xbar\seqsum (\aray{\vbar})
  \subseteq
  \fzrmid{\xbar}{(\limray{\vbar}\plusl\ybar)},
\]
where the equality is by \Cref{thm:ast-halfline:conic},
and the second inclusion is by
\Cref{lem:ahfline-containments}.
Therefore, by \Cref{thm:ast-F-char-fcn-vals}
and since both $\xbar$ and 
$\limray{\vbar}\plusl\ybar$ are in $\dom{F}$,
we obtain
\[
  F(\zbar) \leq 1\cdot F(\xbar) + 0\cdot F(\limray{\vbar}\plusl\ybar) = F(\xbar).%
\indexg{convex functions, astral|)}%
\indexg{convex functions, astral!bounded along halfline|)}%
\indexg{halflines, astral!function bounded along|)}%
\qedhere
\]
\end{proof}

\section{Constructing and operating on convex functions}
\label{sec:op-ast-cvx-fcns}

We next present various ways of constructing astral convex
functions and study operations that preserve convexity,
generalizing similar results from standard convex analysis.

\indexg{convex functions, astral!one dimension@in one dimension|(}%
As a first illustration, we characterize astral convex
functions in one dimension:

\begin{theorem}   \label{pr:1d-cvx}
  Let $F:\Rext\rightarrow\Rext$.
  Then $F$ is convex if and only if all of the following hold:
  \begin{roman-compact}
    \item      \label{pr:1d-cvx:a}
      $\resfcn{F}{\R}$ is convex.
    \item      \label{pr:1d-cvx:b}
      If $F(+\infty)<+\infty$ then $F$ is nonincreasing.
    \item      \label{pr:1d-cvx:c}
      If $F(-\infty)<+\infty$ then $F$ is nondecreasing.
  \end{roman-compact}
\end{theorem}

\begin{proof}
  ~

\begin{proof-parts}
\pfpart{``Only if'' ($\Rightarrow$):}
Suppose $F$ is convex
and let $f=\resfcn{F}{\R}$.
We aim to prove the three stated conditions.
First,
a point $\rpair{x}{y}\in\R\times\R=\R^2$
is clearly in $f$'s epigraph if and only if it is in $F$'s epigraph;
thus,
$\epi{f} = (\epi{F}) \cap \R^2$.
Since the latter two sets are convex (in $\extspac{2}$),
$\epi{f}$ is as well, proving
condition~(\ref{pr:1d-cvx:a}).

For condition~(\ref{pr:1d-cvx:b}), suppose $F(+\infty)<+\infty$.
Then $F(\limray{1}\plusl 0)=F(+\infty)<+\infty$; therefore, by
\Cref{thm:F-conv-res},
$F(\barz)\le F(\barx)$ for all $\barx\in\Rext$
and all $\barz\in\ahfline{\barx}{1}$.
Furthermore,
\[
  \ahfline{\barx}{1}
  =
  \lb{\barx}{+\infty}
  =
  [\barx,+\infty],
\]
where the first equality is by
\Cref{thm:ast-halfline:conic}, and the second 
by \Cref{ex:segs-in-Rext}.
Thus, $F(\barz)\le F(\barx)$ for all
$\barx,\barz\in\eR$ such that $\barz\ge\barx$,
so $F$ is nonincreasing.

The proof of
condition~(\ref{pr:1d-cvx:c})
is symmetric.

\pfpart{``If'' ($\Leftarrow$):}
Suppose
conditions~(\ref{pr:1d-cvx:a}),~(\ref{pr:1d-cvx:b}),
and~(\ref{pr:1d-cvx:c}) hold,
and let $f=\resfcn{F}{\R}$.
We aim to show that $F$ is convex.

If $F(-\infty)=F(+\infty)=+\infty$, then
$\epi F = \epi f$, which is convex since $f$ is convex, implying $F$
is as well.

Otherwise, let $\lambda\in [0,1]$,
let $\barx_0,\barx_1\in\dom{F}$, and let
$\barx\in\flammid{\barx_0}{\barx_1}$.
We aim to show
\begin{equation}   \label{eq:pr:1d-cvx:1}
  F(\barx) \leq (1-\lambda) F(\barx_0) + \lambda F(\barx_1),
\end{equation}
implying convexity by \Cref{thm:ast-F-char-fcn-vals}.
Without loss of generality, we assume $\barx_0\leq\barx_1$
(otherwise swapping $\barx_0$ with $\barx_1$,
and replacing $\lambda$ with $1-\lambda$).

If $F(-\infty)<+\infty$ and $F(+\infty)<+\infty$,
then conditions~(\ref{pr:1d-cvx:b})
and~(\ref{pr:1d-cvx:c}) together imply that $F$ is identically equal
to some constant in $\R\cup\{-\infty\}$,
which then implies \eqref{eq:pr:1d-cvx:1}.

Assume then that $F(-\infty)=+\infty$ and $F(+\infty)<+\infty$,
implying that $F$ is nonincreasing.
(The remaining case that
$F(-\infty)<+\infty$ and $F(+\infty)=+\infty$
can be handled symmetrically.)
In particular, this implies $\barx_0,\barx_1>-\infty$.

Suppose first that both $\barx_0$ and $\barx_1$ are in $\R$.
Then
$\barx=(1-\lambda)\barx_0+\lambda\barx_1$
by \Cref{pr:seq-conic:prop}(\ref{i:seq-conic:finite}),
so $\barx\in\R$.
Since $f$ is convex,
$f(\barx)\leq (1-\lambda) f(\barx_0) + \lambda f(\barx_1)$
(\Cref{pr:stand-cvx-fcn-char}),
implying \eqref{eq:pr:1d-cvx:1}
since $F$ equals $f$ on all points in $\R$.

In the alternative case that
$\barx_0$ and $\barx_1$ are not both finite,
we must have
$-\infty<\barx_0\leq\barx_1=+\infty$.
Then
\begin{equation}   \label{eq:pr:1d-cvx:2}
  \barx
  \in
  \flammid{\barx_0}{\barx_1}
  \subseteq
  \lb{\barx_0}{\barx_1}
  =
  [\barx_0,+\infty],
\end{equation}
with the second inclusion by \Cref{cor:conv-as-seqsum-midrays},
and the equality by \Cref{ex:segs-in-Rext}.
If $\barx_0=+\infty$, then this means also that
$\barx=\barx_0=\barx_1=+\infty$, thus implying
\eqref{eq:pr:1d-cvx:1}.
We therefore assume henceforth that $\barx_0\in\R$.

If $\lambda=0$, then because $F$ is nonincreasing and $\barx\ge\barx_0$
(by Eq.~\ref{eq:pr:1d-cvx:2}),
we obtain $F(\barx)\le F(\barx_0) = (1-\lambda) F(\barx_0) + \lambda F(\barx_1)$.

Otherwise, suppose $\lambda>0$.
Then
\begin{equation*}
  \barx
  \in
  \flammid{\barx_0}{\barx_1}
  =
  [(1-\lambda)\barx_0]\seqsum(+\infty)
  =\set{(1-\lambda)\barx_0\plusl(+\infty)}
  =\set{+\infty}.
\end{equation*}
The first equality is because, by \Cref{thm:mul-char},
$\mul{\lambda}{\barx_1}=\mul{\lambda}{(+\infty)}=\set{+\infty}$ (since $\lambda>0$)
and
$\mul{1-\lambda}{\barx_0}=\set{(1-\lambda)\barx_0}$
(since $\barx_0\in\R$).
The second equality is by
\Cref{cor:seqsum-conseqs}(\ref{cor:seqsum-conseqs:c}).
Thus, $\barx=+\infty$.
Since $F$ is nonincreasing and $\barx\ge\barx_0$,
we then have $(1-\lambda)F(\barx)\leq (1-\lambda)F(\barx_0)$.
Adding $\lambda F(\barx)=\lambda F(+\infty)=\lambda F(\barx_1)$ to both sides then yields
\eqref{eq:pr:1d-cvx:1}.%
\indexg{convex functions, astral!one dimension@in one dimension|)}%
\qedhere
\end{proof-parts}
\end{proof}

\indexg{convex functions, astral!constructing|(}%
\indexg{indicator functions (astral)!convexity of|(}%
An astral indicator function on any convex subset of $\extspace$ is
convex:

\begin{proposition}   \label{pr:ast-ind-fcn-cvx}
  Let $S\subseteq\extspace$ be convex.
  Then the astral indicator function $\indaS$ on $S$
  (as defined in Eq.~\ref{eq:indfa-defn})
  is convex.
\end{proposition}

\begin{proof}
Let $\xbar_0,\xbar_1\in S=\dom{\indaS}$, let $\lambda\in [0,1]$,
and let
$\xbar\in\flammid{\xbar_0}{\xbar_1}$.
Then
$\xbar\in\lb{\xbar_0}{\xbar_1}$
(by \Cref{cor:conv-as-seqsum-midrays}),
so
$\xbar\in S$ since $S$ is convex.
Therefore,
$\indaS(\xbar)=0=(1-\lambda)\indaS(\xbar_0)+\lambda\indaS(\xbar_1)$,
proving the claim by
\Cref{thm:ast-F-char-fcn-vals}.%
\indexg{indicator functions (astral)!convexity of|)}%
\end{proof}

\indexg{affine functions (astral)!convexity of|(}%
Astral linear or affine functions are convex:

\begin{proposition}  \label{pr:aff-fcns-cvx}
  Let $\uu\in\Rn$ and $\beta\in\R$, and
  let $F:\extspace\rightarrow\Rext$ be defined by
  $F(\xbar)=\xbar\cdot\uu+\beta$ for $\xbar\in\extspace$.
  Then $F$ is convex.
\end{proposition}

\begin{proof}
The function $\xx\mapsto\xx\cdot\uu+\beta$, for $\xx\in\Rn$,
is affine and so convex.
As we saw in \Cref{ex:ext-affine},
its extension is exactly $F$.
Therefore, $F$ is convex by \Cref{thm:fext-convex}.%
\indexg{affine functions (astral)!convexity of|)}%
\end{proof}

\indexg{upward addition!convex functions@of convex functions|(}%
\indexg{convex functions, astral!sum of|(}%
The upward sum (as defined in
Eq.~\ref{eq:up-add-def})
of two convex funtions is also convex:

\begin{theorem}   \label{thm:sum-ast-cvx-fcns}
  Let $F:\extspace\rightarrow\Rext$ and
  $G:\extspace\rightarrow\Rext$ be convex.
  Then $F\plusu G$ is also convex.
\end{theorem}

\begin{proof}
Let $H=F\plusu G$.
Let $\xbar_0,\xbar_1\in\dom{H}$, and let $\lambda\in[0,1]$.
Suppose
$\xbar\in\flammid{\xbar_0}{\xbar_1}$.
Then $\xbar_0,\xbar_1$ also are both in $\dom{F}$ and $\dom{G}$
(by definition of upward sum).
Therefore, from \Cref{thm:ast-F-char-fcn-vals},
\begin{align*}
  F(\xbar)
  &\leq
  (1-\lambda) F(\xbar_0)
  +
  \lambda F(\xbar_1),
  \\
  G(\xbar)
  &\leq
  (1-\lambda) G(\xbar_0)
  +
  \lambda G(\xbar_1).
\end{align*}
None of the terms appearing in either inequality can be $+\infty$.
Therefore, the two inequalities can be added yielding
\[
  H(\xbar)
  \leq
  (1-\lambda) H(\xbar_0)
  +
  \lambda H(\xbar_1),
\]
proving $H$ is convex
by \Cref{thm:ast-F-char-fcn-vals}.%
\indexg{upward addition!convex functions@of convex functions|)}%
\end{proof}

In general,
neither the leftward nor the downward sum of two astral convex
functions is necessarily convex, as the next example shows:

\begin{example}
Let $F:\Rext\rightarrow\Rext$ and $G:\Rext\rightarrow\Rext$ be defined
by $F(\barx)=\barx$ and $G(\barx)=\barx^2$ for
$\barx\in\Rext$
(with $(\pm\infty)^2=+\infty$).
By \Cref{pr:1d-cvx}, both $F$ and $G$ are convex.
However, if either $H=F\plusl G$ or $H=F\plusd G$, then
\[
  H(\barx)
  =
  \begin{cases}
    -\infty    & \text{if $\barx=-\infty$,} \\
    \barx^2+\barx    & \text{if $\barx\in\R$,} \\
    +\infty    & \text{if $\barx=+\infty$,}
  \end{cases}
\]
which is not convex,
because it violates condition~(\ref{pr:1d-cvx:c})
of \Cref{pr:1d-cvx}.%
\indexg{convex functions, astral!sum of|)}%
\end{example}

\indexg{affine map in composition with astral function!convexity of|(}%
\indexg{linear map in composition with astral function!convexity of|(}%
The composition of an astral convex function with a linear or affine
map is convex:

\begin{theorem}  \label{thm:cvx-compose-affine-cvx}
  Let $F:\extspac{m}\rightarrow\Rext$ be convex,
  let $\A\in\Rmn$ and let $\bbar\in\extspac{m}$.
  Let $H:\extspace\rightarrow\Rext$ be defined by
  $H(\xbar)=F(\bbar\plusl\A\xbar)$
  for $\xbar\in\extspace$.
  Then $H$ is convex.
\end{theorem}

\begin{proof}
Let $\xbar_0,\xbar_1\in\dom{H}$, let $\lambda\in[0,1]$,
and let $\xbar\in\flammid{\xbar_0}{\xbar_1}$.
Let $G(\zbar)=\bbar\plusl\A\zbar$ for $\zbar\in\extspace$.
Then by
\Cref{pr:lam-mid-props}(\ref{pr:lam-mid-props:aff}),
$G(\xbar)\in\flammid{G(\xbar_0)}{G(\xbar_1)}$.
Also, for $i\in\{0,1\}$,
$F(G(\xbar_i))=H(\xbar_i)<+\infty$ so
$G(\xbar_i)\in\dom{F}$.
Thus,
\[
   H(\xbar)
   =
   F(G(\xbar))
   \leq
   (1-\lambda) F(G(\xbar_0))
   +
   \lambda F(G(\xbar_1))
   =
   (1-\lambda) H(\xbar_0)
   +
   \lambda H(\xbar_1),
\]
with the inequality by \Cref{thm:ast-F-char-fcn-vals},
so $H$ is convex
(by that same theorem).%
\indexg{affine map in composition with astral function!convexity of|)}%
\indexg{linear map in composition with astral function!convexity of|)}%
\end{proof}

\indexg{linear image of astral function!convexity of|(}%
The linear image of a convex function (as defined in
Eq.~\ref{eqn:image-F-dfn}) is also convex:

\begin{theorem}   \label{th:lin-img-ast-fcn-cvx}
  Let $F:\extspace\rightarrow\Rext$ be convex and let $\A\in\Rmn$.
  Then $\A F$ is also convex.
\end{theorem}

\begin{proof}
Let $H=\A F$.
Let $\xbar_0,\xbar_1\in\dom{H}$, let $\lambda\in[0,1]$,
and let $\xbar\in\flammid{\xbar_0}{\xbar_1}$.
We aim to show that
\begin{equation}   \label{eq:th:lin-img-ast-fcn-cvx:1}
  H(\xbar) \leq (1-\lambda) H(\xbar_0) + \lambda H(\xbar_1).
\end{equation}

As earlier,
let $\PPx=[\Iden_n,\zero_n]$ and $\PPy=[\trans{\zero_n},1]$ be the
component projection matrices for $\eRn\times\R$,
and also let $\Qx=[\Iden_m,\zero_m]$ and $\Qy=[\trans{\zero_m},1]$ be
the analogous matrices for $\eRm\times\R$.
Let $\B\in\R^{(m+1)\times(n+1)}$ be the matrix
\begin{equation}   \label{eq:th:lin-img-ast-fcn-cvx:B-dfn}
  \B
  =
   \left[
     \begin{array}{ccc|c}
        & &       &  \\
       ~ & \A & ~ & \zerov{m} \\
        & &       &  \\
       \hline
       \rule{0pt}{2.5ex}
       & \trans{\zerov{n}} & &  1
     \end{array}
   \right]
  =
   \trans{\Qx}\A\PPx+\trans{\Qy}\PPy.
\end{equation}
Thus, the upper left $(m\times n)$-submatrix is a copy of $\A$, the
bottom right entry is $1$, and all other entries are $0$.
Then
for all $\zbar\in\extspace$ and $y\in\R$,
\begin{equation}
  \label{eq:th:lin-img-ast-fcn-cvx:2}
  \B \rpair{\zbar}{y}
  =
  \bigParens{\trans{\Qx}\A\PPx+\trans{\Qy}\PPy}
  \bigParens{\trans{\PPx}\zbar\plusl\trans{\PPy}y}
  =
  \trans{\Qx}\A\zbar\plusl\trans{\Qy}y
  =
  \rpair{\A \zbar}{y}.
\end{equation}
The first and third equalities are both
by \Cref{pr:xy-pairs-props}(\ref{pr:xy-pairs-props:a})
(and Eq.~\ref{eq:th:lin-img-ast-fcn-cvx:B-dfn}).
The second equality follows
by distributivity and associativity of astral linear maps
(\Cref{pr:h:4}\ref{pr:h:4c}\ref{pr:h:4d}),
combined with
the orthogonality identities for $\PPx$ and $\PPy$
(Eq.~\ref{eq:PPx:PPy:orth}).

Let $U=\B (\epi{F})$, the image of $\epi{F}$ under $\B$.
Since $F$ is convex, so is its epigraph, implying that $U$ is
as well (by \Cref{cor:thm:e:9}).

For $i\in\{0,1\}$, let $y_i\in\R$ with $y_i > H(\xbar_i)$.
Then by $H$'s definition, there exists $\zbar_i\in\extspace$ such that
$\A\zbar_i=\xbar_i$ and $F(\zbar_i)<y_i$.
Since $\rpair{\zbar_i}{y_i}\in\epi{F}$,
we obtain by \eqref{eq:th:lin-img-ast-fcn-cvx:2} that
$\rpair{\xbar_i}{y_i}=\B\rpair{\zbar_i}{y_i}\in U$.

Let $y=(1-\lambda)y_0 + \lambda y_1$.
Then
\[
  \rpair{\xbar}{y}
  \in
  \flammid{\rpair{\xbar_0}{y_0}}{\rpair{\xbar_1}{y_1}}
  \subseteq
  \lb{\rpair{\xbar_0}{y_0}}{\rpair{\xbar_1}{y_1}}
  \subseteq
  U.
\]
The first inclusion is by
\Cref{pr:lam-mid-props}(\ref{pr:lam-mid-props:d}).
The second is by
\Cref{cor:conv-as-seqsum-midrays}.
The last is because $U$ is convex.

Since $\rpair{\xbar}{y}\in U$, there exist $\zbar\in\extspace$ and
$y'\in\R$ such that $\rpair{\zbar}{y'}\in\epi{F}$ and
$\B\rpair{\zbar}{y'}=\rpair{\xbar}{y}$.
By \eqref{eq:th:lin-img-ast-fcn-cvx:2},
this implies
$\rpair{\A\zbar}{y'}=\rpair{\xbar}{y}$,
and so that
$y'=y$ and $\A\zbar=\xbar$
(\Cref{pr:xy-pairs-props}\ref{pr:xy-pairs-props:h}).
Therefore,
\[
  H(\xbar)
  \leq
  F(\zbar)
  \leq
  y = (1-\lambda)y_0 + \lambda y_1,
\]
with the first inequality from $H$'s definition,
and the second because $\rpair{\zbar}{y}\in\epi{F}$.
Since this holds for all $y_0>H(\xbar_0)$ and $y_1>H(\xbar_1)$,
\eqref{eq:th:lin-img-ast-fcn-cvx:1} must also hold,
proving convexity
by \Cref{thm:ast-F-char-fcn-vals}.%
\indexg{linear image of astral function!convexity of|)}%
\end{proof}

\indexg{suprema and maxima, pointwise!astral convexity of|(}%
The pointwise supremum of any collection of convex functions is convex:

\begin{theorem}  \label{thm:point-sup-is-convex}
  Let $F_i:\extspace\rightarrow\Rext$ be convex
  for all $i\in\indset$, where
  $\indset$ is any index set.
  Let $H=\sup_{i\in\indset} F_i$
  be their pointwise supremum.
  Then $H$ is convex.
\end{theorem}

\begin{proof}
If $\indset$ is empty, then $H\equiv-\infty$, which is convex (for
instance, by
\Cref{thm:ast-F-char-fcn-vals}).
So we assume henceforth that $\indset$ is not empty.

The epigraph of $H$ is exactly the intersection of the epigraphs of
the functions $F_i$; that is,
\[
   \epi{H} = \bigcap_{i\in\indset}  \epi{F_i}.
\]
This is because a pair $\rpair{\xbar}{y}$ is in $\epi{H}$, meaning
$y\geq H(\xbar)$,
if and only if
$y\geq F_i(\xbar)$
for all $i\in\indset$,
that is, if and only if
$\rpair{\xbar}{y}$ is in $\epi{F_i}$
for all $i\in\indset$.
Since $F_i$ is convex, its epigraph, $\epi{F_i}$, is convex,
for $i\in\indset$.
Thus, $\epi{H}$ is convex by
\Cref{pr:e1}(\ref{pr:e1:b}), and therefore $H$ is as well.%
\indexg{suprema and maxima, pointwise!astral convexity of|)}%
\end{proof}

We had earlier shown (\Cref{pr:conj-is-convex})
that the conjugate $\Fstar$ of any function
$F:\extspace\rightarrow\Rext$ must always be convex.
\indexg{conjugate, dual (astral)!convexity of|(}%
We can now show that the dual conjugate $\psistarb$ of any function
$\psi:\Rn\rightarrow\Rext$ must also always be convex.
In particular, this immediately implies that the biconjugate
$\Fdub$ is always convex, as well as $\fdub$, for
any function $f:\Rn\rightarrow\Rext$.

\begin{theorem}   \label{thm:dual-conj-cvx}
  Let $\psi:\Rn\rightarrow\Rext$.
  Then its dual conjugate, $\psistarb$, is convex.
\end{theorem}

\begin{proof}
As defined in \eqref{eq:psistar-def},
the dual conjugate $\psistarb$, for $\xbar\in\extspace$, is equal to
\[
   \psistarb(\xbar) = \sup_{\rpair{\uu}{v}\in\epi \psi} [\xbar\cdot\uu-v].
\]
Each function $\xbar\mapsto\xbar\cdot\uu-v$ is convex by
\Cref{pr:aff-fcns-cvx}.
Therefore, being a pointwise supremum over convex functions,
$\psistarb$ is also convex, by
\Cref{thm:point-sup-is-convex}.%
\indexg{conjugate, dual (astral)!convexity of|)}%
\end{proof}

\indexg{nondecreasing function, composition with!astral convexity of|(}%
The composition of a nondecreasing, convex astral function with a
convex astral function is convex:

\begin{theorem}     \label{thm:comp-nondec-fcn-cvx}
  Let $F:\extspace\rightarrow\Rext$ be convex,
  and let $G:\Rext\rightarrow\Rext$ be convex and nondecreasing.
  Then $G\circ F$ is convex.
\end{theorem}

\begin{proof}
Let $H=G\circ F$.
If $G(+\infty)<+\infty$, then $G$ must be nonincreasing, by
\Cref{pr:1d-cvx}.
Since $G$ is also nondecreasing, this
implies $G\equiv\barc$ for some $\barc\in\Rext$, further implying
$H\equiv\barc$, which is convex
(for instance, by \Cref{thm:ast-F-char-fcn-vals}).

Otherwise, $G(+\infty)=+\infty$.
Let $\xbar_0,\xbar_1\in\dom{H}$, let $\lambda\in[0,1]$,
and let $\xbar\in\flammid{\xbar_0}{\xbar_1}$.
For $i\in\{0,1\}$,
$G(F(\xbar_i))=H(\xbar_i)<+\infty$
so $F(\xbar_i)\in\dom{G}$,
implying also that $F(\xbar_i)<+\infty$.

By \Cref{thm:ast-F-char-fcn-vals},
\begin{equation}        \label{eq:thm:comp-nondec-fcn-cvx:1}
  F(\xbar)
  \leq
  (1-\lambda) F(\xbar_0) + \lambda F(\xbar_1).
\end{equation}
Further, the expression on the right-hand side is equal to
$(1-\lambda) F(\xbar_0) \plusl \lambda F(\xbar_1)$,
which is in
$(1-\lambda) F(\xbar_0) \seqsum \lambda F(\xbar_1)$
by \Cref{cor:seqsum-conseqs}(\ref{cor:seqsum-conseqs:a}),
and so is a $\lambda$-midpoint of $F(\xbar_0)$ and $F(\xbar_1)$ by
\Cref{pr:seq-conic:prop}(\ref{i:seq-conic:subset}).

Thus,
\begin{align*}
  H(\xbar)
  =
  G\bigParens{F(\xbar)}
  &\le
  G\bigParens{(1-\lambda) F(\xbar_0) + \lambda F(\xbar_1)}
  \\
  &\le
  (1-\lambda) G\bigParens{F(\xbar_0)} + \lambda G\bigParens{F(\xbar_1)}
  \\
  &=
  (1-\lambda) H(\xbar_0) + \lambda H(\xbar_1).
\end{align*}
The first inequality follows from
\eqref{eq:thm:comp-nondec-fcn-cvx:1} since $G$ is nondecreasing.
The second inequality is by
\Cref{thm:ast-F-char-fcn-vals}
applied to $G$.
This proves $H$ is convex
(by the same theorem).%
\indexg{convex functions, astral!constructing|)}%
\indexg{nondecreasing function, composition with!astral convexity of|)}%
\end{proof}

\section{Lower envelope and lower support functions}
\label{sec:fcns-induced-by-sets}

The properties of a function $F:\extspace\rightarrow\Rext$ are
entirely determined by its epigraph; indeed, $F$ itself can be derived
straightforwardly from $\epi{F}$ by the rule
\begin{equation}  \label{eq:F-from-epiF}
  F(\xbar)
  =
  \inf\regBraces{y\in\R :\: \rpair{\xbar}{y}\in\epi{F}}.
\end{equation}
We will see shortly that it is sometimes helpful to represent a
function by a set in $\extspacnp$ other than the epigraph, for
instance, using a set that is both convex and closed in $\extspacnp$,
which the epigraph usually is not.
\indexg{lower envelope|(}%
Generalizing \eqref{eq:F-from-epiF} along these lines leads to the
following definition:

\begin{definition}
\indexg{lower envelope!defined|(}%
Let $E\subseteq\extspacnp$.
The \emph{lower envelope} of $E$ is the
function $\settofcnE:\extspace\rightarrow\Rext$ defined,
for $\xbar\in\eRn$, as
\begin{equation}  \label{eq:settofcn-dfn}
\indexm{lenv e}{$\settofcn{E}$}{lower envelope}%
  \settofcnEf(\xbar)
  =
  \inf\regBraces{
    \PPy\zbar:\:
    \zbar\in E,\,
    \PPx\zbar = \xbar
  },
\end{equation}
where
$\PPx=[\Idnn,\zerov{n}]$ and
$\PPy=[\trans{\zerov{n}},1]$.%
\indexg{lower envelope!defined|)}%
\end{definition}
In other words, $\settofcnEf(\xbar)$ is the infimal value of the ``last
coordinate'' of any point $\zbar$ in $E$ for which
$\PPx\zbar=\xbar$.
In particular, if $F:\extspace\rightarrow\Rext$ and $E=\epi{F}$, then
$\settofcnE=F$, as evident from \eqref{eq:F-from-epiF}.

Here are other properties of lower envelopes:

\begin{proposition}   \label{pr:settofcn-props}
  Let $E\subseteq\extspacnp$, let $F:\extspace\rightarrow\Rext$,
  and
  let
  $\PPx=[\Idnn,\zerov{n}]$
  and
  $\PPy=[\trans{\zerov{n}},1]$.
  \begin{letter-compact}
  \item   \label{pr:settofcn-props:a}
    If $E$ is convex then $\settofcnE$ is convex.
  \item   \label{pr:settofcn-props:b}
    If $E$ is closed (in $\extspacnp$),
    then $\settofcnE$ is lower semicontinuous;
    furthermore,
    for each $\xbar\in\PPxE$,
    the infimum in \eqref{eq:settofcn-dfn} is attained,
    that is,
    there exists $\zbar\in E$ such that $\PPx\zbar=\xbar$
    and $\PPy\zbar=\settofcnEf(\xbar)$.
  \item   \label{pr:settofcn-props:c}
    If $E'\subseteq E$ then $\settofcn{E'}\geq \settofcnE$.
  \item   \label{pr:settofcn-props:d}
    If $E=\epi F$, then
    $\settofcnE=F$.
  \item   \label{pr:settofcn-props:e}
    If $E=\clbar{\epi F}$, then
    $\settofcnE=\lsc{F}$.
  \end{letter-compact}
\end{proposition}

\begin{proof}
The function $\settofcnE$, as given in \eqref{eq:settofcn-dfn},
can be rewritten, for $\xbar\in\extspace$, as
\begin{align*}
  \settofcnEf(\xbar)
  &=
  \inf\,\BigBraces{
    \PPy\zbar \plusu \indfa{E}(\zbar)
    :\:
    \zbar\in \extspacnp,\, \PPx\zbar = \xbar
  }
  \\
  &=
  \inf\,\BigBraces{
    G(\zbar)
    :\:
    \zbar\in \extspacnp,\, \PPx\zbar = \xbar
  },
\end{align*}
where $\indfa{E}$ is the indicator function for the set $E$,
and where we define
$G:\extspacnp\rightarrow\Rext$ by
$G(\zbar) = \PPy\zbar \plusu \indfa{E}(\zbar)$
for $\zbar\in\extspacnp$.
Thus, $\settofcnE=\PPx G$.

\begin{proof-parts}
\pfpart{Part~(\ref{pr:settofcn-props:a}):}
Suppose $E$ is convex.
Then $\indfa{E}$ and $\zbar\mapsto\PPy\zbar$
are both convex functions
(by Propositions~\ref{pr:ast-ind-fcn-cvx}
and~\ref{pr:aff-fcns-cvx}, noting that
$\PPy\zbar=\zbar\cdot\rpair{\zero}{1}$ by
\Cref{pr:xy-pairs-props}\ref{pr:xy-pairs-props:b-new}),
so $G$ is convex as well
(by \Cref{thm:sum-ast-cvx-fcns}).
Thus, $\settofcnE=\PPx G$ is convex by
\Cref{th:lin-img-ast-fcn-cvx}.

\pfpart{Part~(\ref{pr:settofcn-props:b}):}
Suppose $E$ is closed. Define functions
$G_1:\extspacnp\rightarrow\Rext$
and
$G_2:\extspacnp\rightarrow\Rext$,
for $\zbar\in\extspacnp$, as
\[
  G_1(\zbar) = \PPy\zbar,
  \qquad
  G_2(\zbar) = \begin{cases}
    -\infty &\text{if $\zbar\in E$,}
  \\
    +\infty &\text{otherwise.}
  \end{cases}
\]
Note that $G=\max\set{G_1,G_2}$.
Then $G_1$ is continuous by
\Cref{thm:linear:cont}(\ref{thm:linear:cont:b}),
and therefore lower semicontinuous as well.
Also, $G_2$ is lower semicontinuous by
\Cref{prop:lsc}(\ref{prop:lsc:b},\ref{prop:lsc:a}),
because its epigraph is $E\times\R$, which is a closed set
in $\extspacnp\times\R$ (by
\Cref{pr:prod-top-props}\ref{pr:prod-top-props:d}
since $E$ is closed in $\extspacnp$ and $\R$ is closed in $\R$).
Thus, $G=\max\set{G_1,G_2}$ is lower semicontinuous
(by \Cref{pr:lsc-sup}). Hence,
$\settofcnE=\PPx G$ is lower semicontinuous as well
(by \Cref{thm:AF-lsc}\ref{thm:AF-lsc:a}).

To show the infimum in \eqref{eq:settofcn-dfn} is attained,
suppose $\xbar$ is in $\PPxE$, and so also in the astral column
space of $\PPx$.
If $\settofcnEf(\xbar)=+\infty$, then let $\zbar$ be any point in
$E$ for which $\xbar=\PPx\zbar$, implying
$\PPy\zbar\geq \settofcnEf(\xbar)=+\infty$
(by definition of $\settofcnE$),
and thus that $\zbar$ attains the infimum in
\eqref{eq:settofcn-dfn}.

Otherwise, if $\settofcnEf(\xbar)<+\infty$, then
\Cref{thm:AF-lsc}(\ref{thm:AF-lsc:b})
implies that there exists $\zbar\in\extspacnp$ such that
$\PPx\zbar = \xbar$ and $G(\zbar)=\settofcnEf(\xbar)<+\infty$,
and so, by $G$'s definition, that $\zbar\in E$
and $\settofcnEf(\xbar)=\PPy\zbar$, proving attainment
of the infimum in this case as well.

\pfpart{Part~(\ref{pr:settofcn-props:c}):}
Let $\xbar\in\extspace$.
Then
\begin{align*}
  \settofcnEpf(\xbar)
  &=
  \inf\,\regBraces{
    \PPy\zbar:\:
    \zbar\in E',\,
    \PPx\zbar = \xbar
  }
\\
  &\ge
  \inf\,\regBraces{
    \PPy\zbar:\:
    \zbar\in E,\,
    \PPx\zbar = \xbar
  }
  =
  \settofcnEf(\xbar),
\end{align*}
where the inequality is because $E'\subseteq E$.

\pfpart{Part~(\ref{pr:settofcn-props:d}):}
This follows from the equivalence of
\eqref{eq:F-from-epiF}
and
\eqref{eq:settofcn-dfn}
when $E=\epi F$.

\pfpart{Part~(\ref{pr:settofcn-props:e}):}
Let $E=\clbar{\epi F}$.
We show first that $\settofcnE\geq\lsc{F}$.
Let $\xbar\in\extspace$.
If $\xbar\not\in\PPxE$, then
$\settofcnEf(\xbar)=+\infty$,
so $\settofcnEf(\xbar)\geq(\lsc{F})(\xbar)$.
Otherwise, if $\xbar\in\PPxE$, then
by part (\ref{pr:settofcn-props:b}), there exists $\zbar\in E$
such that $\PPx\zbar=\xbar$ and
$\PPy\zbar=\settofcnEf(\xbar)$.
Since $\zbar\in\clbar{\epi F}$, there exists a sequence
$\seq{\rpair{\xbar_t}{y_t}}$ in $\epi F$ that converges to $\zbar$.
Hence,
\[
  \settofcnEf(\xbar)
  =
  \PPy\zbar
  =
  \lim \PPy\rpair{\xbar_t}{y_t}
  =
  \lim y_t
  \geq
  \liminf F(\xbar_t)
  \geq
  (\lsc{F})(\xbar).
\]
The second equality is by
\Cref{thm:linear:cont}(\ref{thm:linear:cont:b}), since
$\rpair{\xbar_t}{y_t} \rightarrow \zbar$.
The third is by \Cref{pr:xy-pairs-props}(\ref{pr:xy-pairs-props:c}).
The first inequality is because
$\rpair{\xbar_t}{y_t}\in\epi F$.
The second inequality is because
$\xbar_t=\PPx\rpair{\xbar_t}{y_t}\rightarrow\PPx\zbar=\xbar$,
and by definition of $\lsc{F}$.

Next, since $\epi{F}\subseteq E$,
we have that $F=\settofcn{(\epi{F})}\geq\settofcnE$,
with equality by part~(\ref{pr:settofcn-props:d}),
and inequality by part~(\ref{pr:settofcn-props:c}).
Since $F\geq\settofcnE$ and since $\settofcnE$ is lower semicontinuous,
by part~(\ref{pr:settofcn-props:b}), it then follows that
$\lsc{F}\geq\settofcnE$ by
\Cref{prop:lsc:characterize}(\ref{prop:lsc:characterize:b}),
completing the proof.%
\indexg{lower envelope|)}%
\qedhere
\end{proof-parts}
\end{proof}

In \Cref{def:ast-conjugate}, the conjugate $\Fstar$
of a function $F:\extspace\rightarrow\Rext$ was defined in terms of
$F$'s epigraph.
If we
replace $\epi{F}$ with an arbitrary set $E\subseteq\extspacnp$,
we obtain the following generalization of the notion of a conjugate:

\begin{definition}
\indexg{lower support function|(}%
\indexg{lower support function!defined|(}%
Let $E\subseteq\extspacnp$.
The \emph{lower support function} of $E$ is the function
$\Esstar:\Rn\rightarrow\Rext$ defined, for $\uu\in\Rn$, as
\begin{equation}  \label{eq:esstar-dfn}
\indexm{lsup e}{$\Esstar$}{lower support function}%
   \Esstarf(\uu)
   =
   \sup_{\zbar\in E} \Bracks{\zbar\cdot\rpair{\uu}{-1}}.%
\indexg{lower support function!defined|)}%
\end{equation}
\end{definition}

Like conjugates,
the lower support function of any set is always convex:

\begin{proposition}
  Let $E\subseteq\extspacnp$.
  Then $\Esstar$ is convex.
\end{proposition}

\begin{proof}
For each $\zbar\in\extspacnp$, 
the function $\ww\mapsto\zbar\cdot\ww$, for $\ww\in\Rnp$,
is convex
by \Cref{thm:h:5}(\ref{thm:h:5a0},\ref{thm:h:5b}).
Therefore, the function
$\uu\mapsto\zbar\cdot\rpair{\uu}{-1}$, for $\uu\in\Rn$,
is convex as well by 
\Cref{roc:thm5.7:fA}, since
$\rpair{\uu}{-1}=\trans{\PPx}\uu-\trans{\PPy}$
(by 
\Cref{pr:xy-pairs-props}\ref{pr:xy-pairs-props:a},
where $\PPx=[\Idnn,\zerov{n}]$ and
$\PPy=[\trans{\zerov{n}},1]$).
Therefore, $\Esstar$ is the pointwise supremum of convex functions,
and so is convex by \Cref{roc:thm5.5}.
\end{proof}

\indexg{lower envelope|(}%
When $E=\epi{F}$, the lower support function $\Esstar$ is the same as
the conjugate $\Fstar$ defined in \eqref{eq:Fstar-def}. In this case, we also
have $\settofcnE=F$, and so $\Esstar=\settofcnEstar$. However,
this latter equality need not hold for general sets $E\subseteq\extspacnp$,
as we show in the next example.

\begin{example}
  In $\extspac{2}$, let $\zbar=\limray{\vv}$
  where $\vv=\trans{[1,1]}$, and let
  $E=\{\zbar\}$.
  Then $\PPy\zbar=+\infty$, so $\settofcnE\equiv+\infty$,
  implying $\settofcnEstar\equiv-\infty$.
  On the other hand, $\Esstarf(1)=\zbar\cdot\rpair{1}{-1}=0$.
  Thus, $\Esstar\neq\settofcnEstar$.
\end{example}

Although $\Esstar$ can differ from $\settofcnEstar$, it is always the
case that $\Esstar\geq\settofcnEstar$, as we show next:

\begin{theorem}  \label{thm:settofcnEstar-leq-Esstar}
  Let $E\subseteq\extspacnp$.
  Then $\settofcnEstar\leq\Esstar$.
\end{theorem}

\begin{proof}
Let $\uu\in\Rn$.
We aim to show that $\settofcnEstar(\uu)\leq\Esstarf(\uu)$.
Let $\xbar\in\extspace$.
We claim that
\begin{equation}  \label{eq:settofcnEstar-leq-Esstar:1}
  {-}\settofcnEf(\xbar)\plusd\xbar\cdot\uu
  \le
  \Esstarf(\uu).
\end{equation}
If either $\settofcnEf(\xbar)=+\infty$ or $\xbar\cdot\uu=-\infty$, then
the left-hand side of
\eqref{eq:settofcnEstar-leq-Esstar:1}
is $-\infty$ so that the inequality holds trivially.
We therefore assume henceforth that
$\settofcnEf(\xbar)<+\infty$ and $\xbar\cdot\uu>-\infty$.
Let $\beta\in\R$ be such that $\beta>\settofcnEf(\xbar)$.
Then by definition of $\settofcnE$,
there exists $\zbar\in E$ such that $\PPx\zbar=\xbar$
and $\PPy\zbar<\beta$,
implying that
\begin{align}
\notag
  \xbar\inprod\uu - \beta
  &\le
  (\PPx\zbar)\inprod\uu - \PPy\zbar
\\
\label{eq:settofcnEstar-leq-Esstar:3}
  &=
  \zbar\inprod\rpair{\uu}{0}
  +
  \zbar\inprod\rpair{\zero}{-1}
  =
  \zbar\inprod\rpair{\uu}{-1}
  \le\Esstarf(\uu).
\end{align}
The first inequality and the summability of the terms on its
right-hand side follow by
\Cref{pr:summable:inc},
since $\xbar\inprod\uu-\beta>-\infty$.
The first equality is by \Cref{pr:xy-pairs-props}(\ref{pr:xy-pairs-props:b-new}),
the second by \Cref{pr:i:1},
and the final inequality is by the definition of $\Esstar$
since $\zbar\in E$.
Thus,
\begin{align*}
  {-}\settofcnEf(\xbar)\plusd\xbar\cdot\uu
  &=
  \xbar\cdot\uu\plusd
  \sup\bigBraces{
    {-}\beta
    :\:
    \beta\in\R,\, \beta>\settofcnEf(\xbar)
  }
  \\
  &=
  \sup\bigBraces{
    \xbar\cdot\uu - \beta
    :\:
    \beta\in\R,\,\beta>\settofcnEf(\xbar)
  }
  \le
  \Esstarf(\uu),
\end{align*}
with the second equality by
\Cref{pr:plusd-props}(\ref{pr:plusd-props:d-gen}),
and the inequality by
\eqref{eq:settofcnEstar-leq-Esstar:3}.

Having proved \eqref{eq:settofcnEstar-leq-Esstar:1} for all
$\xbar\in\extspace$, it now follows by definition of conjugate
(Eq.~\ref{eq:Fstar-down-def})
that
\[
  \settofcnEstar(\uu)
  =
  \sup_{\xbar\in\extspace}
  \bigBracks{
    {-}\settofcnEf(\xbar)\plusd\xbar\cdot\uu
  }
  \le
  \Esstarf(\uu).
\qedhere
\]
\end{proof}

In \Cref{thm:settofcnEstar-leq-Esstar}, we considered the relationship
between $\Esstar$ and $\settofcnE$'s primal conjugate,
$\settofcnEstar$.
As a direct corollary, we obtain correspondingly an analogous
relationship between $\settofcnE$ and $\Esstar$'s dual conjugate, 
$\Esdub$, namely, that $\settofcnE\geq\Esdub$ always.
This can be viewed as a generalization of
\Cref{thm:fdub-sup-afffcns}(\ref{thm:fdub-sup-afffcns:c})'s
assertion that $F\geq\Fdub$,
for any function $F:\extspace\rightarrow\Rext$;
indeed, that fact can be shown to follow as a special case of this
corollary by setting
$E=\epi F$.

\begin{corollary}  \label{cor:ME-geq-Esdub}
  Let $E\subseteq\extspacnp$.
  Then $\settofcnE\geq\Esdub$.
\end{corollary}

\begin{proof}
By
\Cref{thm:fdub-sup-afffcns}(\ref{thm:fdub-sup-afffcns:c}),
$\settofcnE\geq\settofcnEdub$,
and by
\Cref{thm:settofcnEstar-leq-Esstar},
$\settofcnEstar\leq\Esstar$,
implying
$\settofcnEdub\geq\Esdub$
by
\Cref{pr:dual-biconj-props}(\ref{pr:dual-biconj-props:a}).
Combining yields the claim.
\end{proof}

The next proposition relates the functions $\settofcnE$ and
$\Esstar$ to the function $F$ according to how the set $E$ relates to
$F$'s epigraph.

\begin{proposition}  \label{pr:settofcn-ohull-prop}
  Let $E\subseteq\extspacnp$ and $F:\extspace\rightarrow\Rext$.
  \begin{letter-compact}
  \item   \label{pr:settofcn-ohull-prop:a}
    If $\epi{F}\subseteq E$, then
    $F\geq \settofcnE$ and $\Fstar\leq\Esstar$.
  \item   \label{pr:settofcn-ohull-prop:b}
    If $E \subseteq \ohull(\epi{F})$, then
    $\settofcnE \geq \Fdub$
    and
    $\Fstar\geq\Esstar$.
  \end{letter-compact}
\end{proposition}

\begin{proof}
  ~

\begin{proof-parts}
\pfpart{Part~(\ref{pr:settofcn-ohull-prop:a}):}
Suppose $\epi{F}\subseteq E$.
Then $F=\settofcn{(\epi{F})}\geq\settofcnE$
by
\Cref{pr:settofcn-props}(\ref{pr:settofcn-props:c},\ref{pr:settofcn-props:d}).
Consequently, $\Fstar\leq\settofcnEstar\leq\Esstar$ with the 
inequalities following respectively from
\Cref{pr:primal-biconj-props}(\ref{pr:primal-biconj-props:a})
and
\Cref{thm:settofcnEstar-leq-Esstar}.

\pfpart{Part~(\ref{pr:settofcn-ohull-prop:b}):}
Suppose $E \subseteq \ohull(\epi{F})$.
We show first that $\Fstar\geq\Esstar$.
Let $\uu\in\Rn$.
If $\zbar\in E$, then $\zbar\in\ohull(\epi F)$, implying
\[
  \zbar\cdot\rpair{\uu}{-1}
  \leq
  \sup_{\zbar'\in\epi F} [\zbar'\cdot\rpair{\uu}{-1}]
  =
  \Fstar(\uu),
\]
with the inequality from
\Crefequiv{pr:ohull-simplify}{pr:ohull-simplify:a}{pr:ohull-simplify:b}.
Since this holds for all $\zbar\in E$, we then have that
\[
  \Esstarf(\uu)
  =
  \sup_{\zbar\in E} \Bracks{\zbar\cdot\rpair{\uu}{-1}}
  \leq
  \Fstar(\uu).
\]
Thus, $\Esstar\leq\Fstar$,
implying further that
$
  \settofcnE
  \geq
  \Esdub
  \geq
  \Fdub
$
with inequalities from
\Cref{cor:ME-geq-Esdub}
and
\Cref{pr:dual-biconj-props}(\ref{pr:dual-biconj-props:a}),
respectively.
\qedhere
\end{proof-parts}
\end{proof}

We finish this section by providing general conditions for when
$\settofcnEf(\xbar)=\Esdubf(\xbar)$, where $\xbar\in\extspace$,
assuming $E\subseteq\extspacnp$ is closed and convex.
In particular, we show in the next
\namecref{thm:settofcn-eq-biconj}
that this equality holds
if $\xbar\in\PPx E$ or if $\Esdubf(\xbar)>-\infty$.
This is analogous to \Cref{thm:fext-neq-fdub},
which states that $\fext(\xbar)=\fdub(\xbar)$
if $\xbar\in\cldom{f}$ or if $\fdub(\xbar)>-\infty$
(assuming $f$ is convex);
indeed, it can be shown that that
\namecref{thm:fext-neq-fdub}
can be recovered by setting $E=\cldom{f}$ below.
Later, in \Cref{cor:ult-cvx-F-Fdub},
we will also apply the next
\namecref{thm:settofcn-eq-biconj}
to derive a generalization of
\Cref{thm:fext-neq-fdub} that holds for any convex and lower
semicontinuous astral function (not just extensions).

\begin{theorem}   \label{thm:settofcn-eq-biconj}
  Let $E\subseteq\extspacnp$ be closed (in $\extspacnp$) and convex,
  and let $\xbar\in\extspace$.
  Let $\PPx=[\Idnn,\zerov{n}]$.
  If either $\xbar\in\PPxE$ or $\Esdubf(\xbar)>-\infty$,
  then $\settofcnEf(\xbar)=\Esdubf(\xbar)$.

  Consequently, if $\settofcnEf(\xbar)\neq\Esdubf(\xbar)$
  then we must have $\settofcnEf(\xbar)=+\infty$ and $\Esdubf(\xbar)=-\infty$.
\end{theorem}

\begin{proof}
If $E=\emptyset$, then $\settofcnE=\Esdub\equiv+\infty$,
proving the claim.
We therefore assume henceforth that $E$ is not empty.

\begin{proof-parts}
\pfpart{Case $\xbar\in\PPxE$:}
Suppose first that $\xbar\in\PPxE$, that is, that there exists
$\zbar\in E$ such that $\xbar=\PPx\zbar$,
and suppose, by way of contradiction, that
$\settofcnEf(\xbar)\neq\Esdubf(\xbar)$.
In light of \Cref{cor:ME-geq-Esdub}, this last assumption
implies
$\settofcnEf(\xbar)>\Esdubf(\xbar)$, and so that there exists $y\in\R$
such that
$\settofcnEf(\xbar)>y>\Esdubf(\xbar)$.
In particular, this means that $\rpair{\xbar}{y}$ is not in $E$, since
otherwise we would have $\settofcnEf(\xbar)\leq y$.

Since the sets $E$ and $\set{\rpair{\xbar}{y}}$ are closed, convex,
nonempty and disjoint, they are strongly separated, by
\Cref{thm:sep-cvx-sets};
that is, there exist $\rpair{\uu}{v}\in\Rn\times\R$ and $\beta\in\R$
such that
$\zbar'\cdot\rpair{\uu}{v}<\beta$ for all $\zbar'\in E$,
and
\begin{equation}   \label{eq:thm:settofcn-eq-biconj:1}
  \xbar\cdot\uu+yv=\rpair{\xbar}{y}\cdot\rpair{\uu}{v}>\beta
\end{equation}
(with the equality from \Cref{pr:xy-pairs-props}\ref{pr:xy-pairs-props:b}).

Considering first the case that $v\geq 0$, we can derive a
contradiction as follows:
\begin{align*}
  \beta
  <
  \xbar\cdot\uu+yv
  &\le
  (\PPx\zbar)\cdot\uu+(\PPy\zbar)v
\\
  &=
  \zbar\cdot\rpair{\uu}{0}+\zbar\cdot\rpair{\zero}{v}
  =
  \zbar\cdot\rpair{\uu}{v}
  <
  \beta.
\end{align*}
The first inequality is by \eqref{eq:thm:settofcn-eq-biconj:1}.
The second inequality and the summability of the terms on its
right-hand side
are by \Cref{pr:summable:inc}
since $\PPx\zbar=\xbar$, $\PPy\zbar\geq\settofcnEf(\xbar)>y$
(since $\zbar\in E$ and by definition of $\settofcnE$),
and $v\ge 0$.
The first equality is by 
\Cref{pr:xy-pairs-props}(\ref{pr:xy-pairs-props:b-new}),
the second by \Cref{pr:i:1},
and the final inequality follows from strong separation,
since $\zbar\in E$.

We are left then with the alternative case that $v<0$.
In this case, we can assume without loss of generality that $v=-1$
(otherwise replacing $\rpair{\uu}{v}$ and $\beta$ with
$\rpair{\uu/|v|}{-1}$ and $\beta/|v|$).
By the foregoing assumptions and
the definition of $\Esstar$,
we then have that
\begin{equation}   \label{eq:thm:settofcn-eq-biconj:2}
  \Esstarf(\uu)
  =
  \sup_{\zbar'\in E} \Bracks{\zbar'\cdot\rpair{\uu}{-1}}
  \leq
  \beta.
\end{equation}
This in turn implies that
\[
  y
  >
  \Esdubf(\xbar)
  \geq
  \xbar\cdot\uu - \beta,
\]
with the first inequality from our choice of $y$, and the second
from the definition of dual conjugate (Eq.~\ref{eq:psistar-def})
since $\rpair{\uu}{\beta}\in\epi\Esstarf$
by \eqref{eq:thm:settofcn-eq-biconj:2}.
Thus, $\xbar\cdot\uu-y<\beta$, contradicting
\eqref{eq:thm:settofcn-eq-biconj:1}.

Having reached a contradiction in both cases, we conclude that
$\settofcnEf(\xbar)=\Esdubf(\xbar)$, as claimed.

\pfpart{Case $\Esdubf(\xbar)>-\infty$:}
Suppose next that $\Esdubf(\xbar)>-\infty$, which we aim to show implies
$\settofcnEf(\xbar)=\Esdubf(\xbar)$.
The preceding argument shows that this will be the case if
$\xbar\in\PPxE$.
We therefore assume henceforth that $\xbar\not\in\PPxE$.
Under this assumption, we will argue that $\Esdubf(\xbar)=+\infty$
which, when combined with \Cref{cor:ME-geq-Esdub}, will prove the
claim.

Since $E$ is closed and convex, $\PPxE$ is as well, by
Corollaries~\ref{cor:aff-img-closed-is-closed}(\ref{cor:aff-img-closed-is-closed:a})
and~\ref{cor:thm:e:9}.
Thus, $\PPxE$ and $\{\xbar\}$ are nonempty,
closed, convex and disjoint;
therefore, they are strongly separated, by
\Cref{thm:sep-cvx-sets}.
Thus, there exists $\uu\in\Rn$ and $\beta\in\R$ such that
$\xbar'\cdot\uu<\beta$ for all $\xbar'\in\PPxE$,
and
$\xbar\cdot\uu>\beta$.

Also, since $\Esdubf(\xbar)>-\infty$, by definition of dual conjugate
(Eq.~\ref{eq:psistar-def}),
there must exist $\rpair{\ww}{v}\in\epi\Esstarf$ such that
$\xbar\cdot\ww-v>-\infty$,
implying $\xbar\cdot\ww>-\infty$.
If $\xbar\cdot\ww=+\infty$, then
$\Esdubf(\xbar)=+\infty$ since
$\Esdubf(\xbar)\geq\xbar\cdot\ww-v$.
We therefore can assume henceforth that $\xbar\cdot\ww\in\R$.

Let $\lambda\in\Rstrictpos$, and let
$\ww'=\ww+\lambda\uu$ and $v'=v+\lambda\beta$.
We claim that $\rpair{\ww'}{v'}\in\epi\Esstarf$.
To see this, let $\zbar$ be any point in $E$.
Then $\zbar\cdot\rpair{\ww}{-1}\leq\Esstarf(\ww)\leq v<+\infty$
since $\rpair{\ww}{v}\in\epi\Esstarf$.
Also,
\[
  \zbar\cdot\rpair{\lambda\uu}{0}
  =
  \lambda\zbar\cdot\rpair{\uu}{0}
  =
  \lambda(\PPx\zbar)\cdot\uu
  <
  \lambda\beta
  <+\infty,
\]
where the second equality is by
\Cref{pr:xy-pairs-props}(\ref{pr:xy-pairs-props:b-new}),
and the first inequality follows because, by our
choice of $\uu$, we have
$(\PPx\zbar)\cdot\uu<\beta$ for all $\zbar\in E$.
Combining,
\begin{align*}
  \zbar\cdot\rpair{\ww'}{-1}
  &=
  \zbar\cdot\bigParens{\rpair{\ww}{-1}+\rpair{\lambda\uu}{0}}
  \\
  &=
  \zbar\cdot\rpair{\ww}{-1}+\zbar\cdot\rpair{\lambda\uu}{0}
  \leq
  v+\lambda\beta
  =
  v',
\end{align*}
with the inequality following from the arguments above, which also
show that the terms on the left-hand side of the second line are
summable, justifying the second equality
(by \Cref{pr:i:1}).
Since this holds for all $\zbar\in E$, it follows that
$\Esstarf(\ww')\leq v'$
(using the definition of $\Esstar$).

Thus,
\begin{align*}
\SwapAboveDisplaySkip
  \Esdubf(\xbar)
  \geq
  \xbar\cdot\ww' - v'
  &=
  (\xbar\cdot\ww + \lambda\xbar\cdot\uu) - (v + \lambda\beta)
\\
  &=
  (\xbar\cdot\ww - v) + \lambda(\xbar\cdot\uu - \beta).
\end{align*}
The inequality is by definition of dual conjugate, since
$\rpair{\ww'}{v'}\in\epi\Esstarf$.
The first equality uses
\Cref{pr:i:1}, with summability ensured since $\xbar\cdot\ww\in\R$.
Since $\xbar\cdot\uu-\beta>0$ and $\xbar\cdot\ww\in\R$,
and since this holds for all
$\lambda\in\Rstrictpos$, it follows that $\Esdubf(\xbar)=+\infty$.
Combined with \Cref{cor:ME-geq-Esdub},
this proves that $\settofcnEf(\xbar)=\Esdubf(\xbar)$.

\pfpart{When $\settofcnEf(\xbar)\neq\Esdubf(\xbar)$:}
If $\settofcnEf(\xbar)\neq\Esdubf(\xbar)$,
then the foregoing shows
that $\Esdubf(\xbar)=-\infty$
and that $\xbar\not\in\PPxE$, implying,
by the definition of $\settofcnE$,
that $\settofcnEf(\xbar)$ is vacuously
\indexg{lower support function|)}%
\indexg{lower envelope|)}%
$+\infty$.
\qedhere
\end{proof-parts}
\end{proof}

\chapter{Astral closedness of functions in general}
\label{sec:biconj-ast-close-fcns}

This chapter provides general characterizations for astral closedness,
that is, for when an
astral function $F:\extspace\rightarrow\Rext$ is equal to its own
biconjugate, $\Fdub$, and for when a (standard) function
$\psi:\Rn\rightarrow\Rext$ is equal to its own astral dual biconjugate,
$\psidub$.
Such a characterization was given
in
\Cref{sec:ent-closed-fcn} for the extension of a convex
function, crucially using properties of reductions.
The characterizations here are
based more directly on the separation theorems
from \Cref{sec:sep-thms},
and are considerably more
general.

\section{Astral primal closedness}
\label{sec:ast-close-primal-biconj}

\indexg{closedness, astral (primal)!one dimension@in one dimension|(}%
We begin by characterizing when an astral function
$F:\extspace\rightarrow\Rext$ is equal to $\Fdub$.
As we show first, when $n=1$,
this is the case when $F$ is convex and lower semicontinuous,
and when its restriction to $\R$ is closed.
We will soon see that this special case is directly relevant for
the general case.

\begin{theorem}
  \label{thm:Fdub:1d}
  Let $F:\eR\to\eR$. Then $F=\Fdub$ if and only if $F$ is convex,
  lower semicontinuous, and $\resfcn{F}{\R}$ is closed.
\end{theorem}
\begin{proof}
  Let $f=\resfcn{F}{\R}$ denote the restriction of $F$ to $\R$.
  \begin{proof-parts}
  \pfpart{%
    ``Only if'' ($\Rightarrow$):
  }
    Suppose $F=\Fdub$.
    Then $F$ is convex by \Cref{thm:dual-conj-cvx}
    and lower semicontinuous by \Cref{pr:dual-conj-lsc}.
    Moreover,
    for all $\xx\in\R$, $f(\xx)=F(\xx)=\Fdub(\xx)=(\Fstar)^*(\xx)$
    (with the last equality by \Cref{pr:psistarb:psistar} applied to $\Fstar$).
    Thus, $f=(\Fstar)^*$, so $f$ is closed
    (using \Cref{pr:conj-props}\ref{pr:conj-props:d} applied to $\Fstar$).

  \pfpart{%
    ``If'' ($\Leftarrow$):
  }
  Suppose $F$ is convex and lower semicontinuous,
  and that $f=\resfcn{F}{\R}$ is closed
  (and also convex by \Cref{pr:1d-cvx}).
  By
  \Cref{thm:fdub-sup-afffcns}(\ref{thm:fdub-sup-afffcns:d}),
  it suffices to show
    that there exists a set $E\subseteq\R^2$
    such that, for all $\barx\in\eR$,
    \begin{equation}
    \label{eq:Fdub:1d:1}
      F(\barx)=\sup_{\rpair{u}{v}\in E} [\barx u - v].
    \end{equation}
    We distinguish the following cases, based on whether $f\equiv+\infty$, 
    or $f\equiv\alpha\in\R\cup\{-\infty\}$,
    or $f$ is proper and not constant.
    (These cases are exhaustive since $f$ is closed.)

  \pfpart{%
    Case $f\equiv+\infty$:
  }
    If $F\equiv+\infty$, then $E=\set{\rpair{0}{v}:\: v\in\R}$ satisfies \eqref{eq:Fdub:1d:1}.
    Otherwise, suppose $F(+\infty)<+\infty$ (the case $F(-\infty)<+\infty$
    is symmetric). By \Cref{pr:1d-cvx}, $F$ is nonincreasing, so $F(-\infty)=+\infty$. If
    $F(+\infty)=-\infty$, then $E=\set{\rpair{-1}{v}:\:v\in\R}$ satisfies \eqref{eq:Fdub:1d:1}.
    If $F(+\infty)=c\in\R$, then $E=\set{\rpair{-1}{v}:\:v\in\R}\cup\set{\rpair{0}{-c}}$ satisfies
    \eqref{eq:Fdub:1d:1}.

  \pfpart{%
    Case $f\equiv\alpha$ for some $\alpha\in\R\cup\{-\infty\}$:
  }
    Since $F$ is lower semicontinuous,
    $F(+\infty)\le\lim f(t)=\alpha$ and
    $F(-\infty)\le\lim f(-t)=\alpha$. By \Cref{pr:1d-cvx},
    this implies that $F$ is both nonincreasing and nondecreasing.
    Thus, $F\equiv\alpha$ and
    \eqref{eq:Fdub:1d:1} is satisfied by
    $E=\set{\rpair{0}{-\alpha}}$ if $\alpha\in\R$, and by $E=\emptyset$ if $\alpha=-\infty$.

  \pfpart{%
    Case $f$ is proper, nonincreasing, and not constant.
  }
    Since $F$ is lower semicontinuous, $F(+\infty)\le\lim f(t)=\inf f<+\infty$ (since $f$ is proper). Thus, $F(+\infty)=\inf F<+\infty$.
    Moreover, $f$ cannot be nondecreasing, so
    $F(-\infty)=+\infty$ by \Cref{pr:1d-cvx}.

    Since $f$ is nonincreasing and not constant, $\resc{f}=\Rpos$,
    so by \Cref{pr:rescpol-is-con-dom-fstar}, $\cl(\cone(\dom\fstar))=\polar{(\Rpos)}=\Rneg$.
    Thus, $\dom\fstar\subseteq\Rneg$ and
    there exists $u_0\in(\dom\fstar)\cap\Rstrictneg$.
This further implies that $\ri(\dom\fstar)\subseteq\Rstrictneg$.
(Otherwise, if $0$ were in $\ri(\dom\fstar)$ then, by
\Cref{roc:thm6.4}, there would exist $\delta\in\Rstrictpos$ with
$-\delta u_0 = 0+\delta(0-u_0)\in\dom\fstar\subseteq\Rneg$,
a contradiction.)

Let $E'=\Braces{\rpair{u}{\fstar(u)} :\: u\in\ri(\dom\fstar)}$,
which is nonempty (since the convex set $\dom\fstar$ includes $u_0$
and so has a nonempty relative interior by
\Cref{pr:ri-props}\ref{pr:ri-props:roc-thm6.2b});
further, $E'\subseteq\R^2$
since $f$ is proper, and so also is $\fstar$
(\Cref{pr:conj-props-cvx}\ref{pr:conj-props-cvx:a}).
Let $G:\Rext\rightarrow\Rext$ be defined,
for $\barx\in\Rext$, by
\[ G(\barx)=\sup_{\rpair{u}{v}\in E'} [\barx u - v]. \]
Then for $x\in\R$, $G(x)=\fdubs(x)=f(x)=F(x)$,
where the first equality is by
\Cref{pr:conj-props-cvx}(\ref{pr:conj-props-cvx:c})
(applied to $\fstar$),
and the second by
\Cref{pr:conj-props-cvx}(\ref{pr:conj-props-cvx:b})
since $f$ is closed.
Since $u<0$ for all $\rpair{u}{v}\in E'$, and since $E'$ is nonempty,
we have moreover that
$G(-\infty)=+\infty$ and $G(+\infty)=-\infty$.

Letting $E=E'\cup\Braces{\rpair{0}{-v} :\: v\in\R, v\leq \inf F}$,
this then implies, for $\barx\in\Rext$, that
\[
  \sup_{\rpair{u}{v}\in E} [\barx u - v]
  =
  \max\{G(\barx),\, \inf F\}
  =
  F(\barx),
\]
where the first equality is from definitions,
and the second is because $G(x)=F(x)\geq\inf F$ for $x\in\R$,
$F(-\infty)=+\infty=G(-\infty)$, and
$F(+\infty)=\inf F\geq -\infty=G(+\infty)$.
Thus, \eqref{eq:Fdub:1d:1} holds,
proving $F=\Fdub$ by
\Cref{thm:fdub-sup-afffcns}(\ref{thm:fdub-sup-afffcns:d}).

  \pfpart{%
    Case $f$ is proper, nondecreasing, and not constant.
  }
  Symmetric to the previous case.

  \pfpart{%
    Case $f$ is proper, and neither nonincreasing nor nondecreasing:
  }%
    In this case, we have $\resc{f}=\set{0}$, so $\cl(\cone(\dom\fstar))=\polar{(\set{0})}=\R$
    (\Cref{pr:rescpol-is-con-dom-fstar}), which implies
    that there exist $u_-\in(\dom\fstar)\cap\Rstrictneg$ and
    $u_+\in(\dom\fstar)\cap\Rstrictpos$. By convexity of $\dom\fstar$,
    this means $0\in\dom\fstar$,
    so $\inf f=-\fstar(0)>-\infty$
    (by \Cref{pr:conj-props}\ref{pr:conj-props:a}).
    Since $f$ is neither
    nonincreasing nor nondecreasing, $F(-\infty)=F(+\infty)=+\infty$ by \Cref{pr:1d-cvx}.
    Thus, $F=\ftriv$, implying $\fext=\lsc\ftriv=\lsc F=F$
    with the first equality by \Cref{pr:lsc-ftriv-is-fext},
    and the third because $F$ is lower semicontinuous.
    Since $f>-\infty$, we have $\ef=\fextdub$ (by
    \Cref{cor:all-red-closed-sp-cases} and
    \Cref{pr:fextstar-is-fstar}).
    Therefore, $F=\ef=\fextdub=\Fdub$.%
\indexg{closedness, astral (primal)!one dimension@in one dimension|)}%
    \qedhere
  \end{proof-parts}

\end{proof}

\indexg{closedness, astral (primal)!general function@of general function|(}%
We next characterize when $F=\Fdub$ in general.
Unlike the characterization from \Cref{sec:ent-closed-fcn},
the one below is applicable to any astral function
${F:\eRn\to\eR}$, not just the extension
of a convex function $f$ on $\Rn$.
The characterization is stated in terms of the functions
$\utransF$, the image of $F$ under $\trans{\uu}\negKern$,
for all $\uu\in\Rn\wo\set{\zero}$.
As follows from the definition given in \eqref{eqn:image-F-dfn}
(applied with $\A=\trans{\uu}$),
such functions are defined on
the extended real line, $\Rext$; specifically, for $\alpha\in\Rext$,
\begin{equation}   \label{eq:trans-u-F-defn}
  (\utransF)(\alpha)
  =
  \inf\regBraces{F(\xbar) :\: \xbar\in\extspace,\, \xbar\cdot\uu=\alpha}
\end{equation}
(using \Cref{pr:trans-uu-xbar}).
Assuming $F$ is convex and lower semicontinuous,
the characterization given in the
next theorem states that $F$ is
astral closed if and only if every function $\utransF$
is astral closed
(for all $\uu\in\Rn\wo\set{\zero}$),
or alternatively,
if and only if
each function $\utransF$, when restricted to $\R$, is closed.

\begin{theorem}   \label{thm:F-equal-Fdub-alt}
  Let $F:\extspace\rightarrow\Rext$.
  Then the following are equivalent:
  \begin{letter-compact}
  \item   \label{thm:F-equal-Fdub-alt:a}
    $F=\Fdub$.
  \item   \label{thm:F-equal-Fdub-alt:b}
    $F$ is convex, lower semicontinuous, and
    for all $\uu\in\Rn\wo\{\zero\}$,
    $\resfcn{(\utransF)}{\R}$ is closed.
  \item   \label{thm:F-equal-Fdub-alt:c}
    $F$ is convex, lower semicontinuous, and
    for all $\uu\in\Rn\wo\{\zero\}$,
    $\utransF = (\utransF)^{*\dualstar}$.
  \end{letter-compact}
\end{theorem}

\begin{proof}
  ~
\begin{proof-parts}
\pfpart{%
  (\ref{thm:F-equal-Fdub-alt:a})
  $\Rightarrow$
  (\ref{thm:F-equal-Fdub-alt:b}):
}
Suppose $F=\Fdub$.
Then $F$ is convex by \Cref{thm:dual-conj-cvx}
and lower semicontinuous by \Cref{pr:dual-conj-lsc}.

Let $\uu\in\Rn\wo\{\zero\}$, and let
$G=\utransF$.
Let $g=\resfcn{G}{\R}$ be the restriction of $G$ to $\R$.
We need to show that $g$ is closed.
Since $F$ is convex and lower semicontinuous,
the same is true for $G$
(by Theorems~\ref{th:lin-img-ast-fcn-cvx}
and~\ref{thm:AF-lsc}\ref{thm:AF-lsc:a}),
and thus also for $g$
(by \Cref{pr:1d-cvx} and \Cref{pr:lsc-res-domain}\ref{pr:lsc-res-domain:b}).
Therefore, if $g>-\infty$, then $g$ must be closed.
Otherwise, assume
$g(z)=-\infty$ for some $z\in\R$. We need to show that
$g\equiv-\infty$; that is, letting $z'\in\R$, we must show
$g(z')=-\infty$.

By \Cref{thm:AF-lsc}(\ref{thm:AF-lsc:b}),
there must exist
$\xbar\in\extspace$ such that $\xbar\cdot\uu=\trans{\uu}\xbar=z$
and such that $F(\xbar)=G(z)=-\infty$.
Since $F=\Fdub$, this means that $\Fdub(\xbar)=-\infty$,
and so, by definition of dual conjugate
(Eq.~\ref{eq:psistar-def}), $\xbar\inprod\ww-v=-\infty$
for all $\rpair{\ww}{v}\in\epi{\Fstar}$.

Let $\lambda=(z'-z)/\norm{\uu}^2$ (noting that $\uu\neq\zero$).
Then for all $\rpair{\ww}{v}\in\epi{\Fstar}$,
\[
  (\xbar\plusl\lambda\uu)\inprod\ww-v
  =
  (\xbar\cdot\ww - v) + \lambda\uu\cdot\ww
  =
  -\infty,
\]
so, again by definition of dual conjugate,
$F(\xbar\plusl\lambda\uu)=\Fdub(\xbar\plusl\lambda\uu)=-\infty$.

Also, since $z=\xbar\cdot\uu$, we have
\[z'=z+\lambda\norm{\uu}^2=\xbar\cdot\uu+(\lambda\uu)\cdot\uu=(\xbar\plusl\lambda\uu)\cdot\uu=\trans{\uu}(\xbar\plusl\lambda\uu).\]
From the definition of~$G$, this implies
that
$G(z')\le F(\xbar\plusl\lambda\uu) = -\infty$.
Thus, $g(z')=G(z')=-\infty$, so $g\equiv-\infty$,
and therefore is closed.

\pfpart{%
  (\ref{thm:F-equal-Fdub-alt:b})
  $\Rightarrow$
  (\ref{thm:F-equal-Fdub-alt:c}):
}
Suppose (\ref{thm:F-equal-Fdub-alt:b}) holds,
let $\uu\in\Rn\wo\set{\zero}$, and let $G=\utransF$.
Then $G$ is
convex by \Cref{th:lin-img-ast-fcn-cvx}
and lower semicontinuous
by \Cref{thm:AF-lsc}(\ref{thm:AF-lsc:a}). Also, by assumption,
$\resfcn{G}{\R}$ is closed, so $G=\Gdub$ by \Cref{thm:Fdub:1d},
proving the claim.

\pfpart{%
  (\ref{thm:F-equal-Fdub-alt:c})
  $\Rightarrow$
  (\ref{thm:F-equal-Fdub-alt:a}):
}
Assume (\ref{thm:F-equal-Fdub-alt:c}) holds.
Let $\xbar\in\extspace$.
We aim to show $\Fdub(\xbar)=F(\xbar)$.
Since $\Fdub\leq F$ by
\Cref{thm:fdub-sup-afffcns}(\ref{thm:fdub-sup-afffcns:c}),
it suffices to show $\Fdub(\xbar)\geq F(\xbar)$.

We first calculate $\Fdub$ in terms of the functions $\utransF$:
\begin{align}
\notag
  \Fdub(\xbar)
  &=\sup_{\uu\in\Rn} [-\Fstar(\uu)\plusd \xbar\cdot\uu]
\\
\notag
  &=\sup_{
      \substack{
        \uu\in\Rn\wo\set{\zero}
        \\ \lambda\in\R
      }}
    \bigBracks{-\Fstar(\lambda\uu)\plusd \xbar\cdot(\lambda\uu)}
\\
\notag
  &=\sup_{\uu\in\Rn\wo\set{\zero}}
    \;
    \sup_{\lambda\in\R\vphantom{\uu\in\Rn\wo\set{\zero}}}
    \,
    \bigBracks{-(\Fstar\uu)(\lambda)\plusd (\xbar\cdot\uu)\lambda}
\\
\notag
  &=\sup_{\uu\in\Rn\wo\set{\zero}}
    (\Fstar\uu)^{\dualstar}(\xbar\cdot\uu)
\\
  &=\sup_{\uu\in\Rn\wo\set{\zero}}
    (\utransF)^{*\dualstar}(\xbar\cdot\uu)
   =\sup_{\uu\in\Rn\wo\set{\zero}}
    (\utransF)(\xbar\cdot\uu).
\label{eq:Fdub-equal-F:alt:1}
\end{align}
The first and fourth equalities are from the definition of
dual conjugate (Eq.~\ref{eq:psistar-def:2}). The fifth equality
is by \Cref{pr:prml-lin-op-conj-idens}(\ref{pr:prml-lin-op-conj-idens:a}),
and the last is by assumption.

\eqref{eq:Fdub-equal-F:alt:1} implies that $\Fdub(\xbar)\geq\inf F$
since $\utransF\geq\inf F$ for all $\uu\in\Rn\wo\{\zero\}$
(by \Cref{pr:ast-lin-img-props}\ref{pr:ast-lin-img-props:inf-AF}).
Hence,
if $F(\xbar)=\inf F$, then
$\Fdub(\xbar)\geq F(\xbar)$, completing the proof.
We therefore assume henceforth that
$F(\xbar)>\inf F$.
In that case,
the statement will follow from the next claim:
\begin{claimpx}
\label{cl:F-equal-Fdub-alt:1}
Let $\beta\in\R$ be such that $\inf F<\beta<F(\xbar)$. Then
there exists
$\uu\in\Rn\wo\set{\zero}$ such that
$(\utransF)(\xbar\cdot\uu)\geq\beta$.
\end{claimpx}

\begin{proofx}
Let $S=\set{\zbar\in\eRn:\:F(\zbar)\le\beta}$.
Then $S$ is closed (by
\Cref{prop:lsc}\ref{prop:lsc:a}\ref{prop:lsc:c}, since $F$ is lower semicontinuous),
convex (\Cref{thm:f:9}), nonempty (since $\inf F<\beta$),
and does not contain $\xbar$ (since $F(\xbar)>\beta$).
Therefore, $S$ and $\set{\xbar}$ are strongly separated
by \Cref{thm:sep-cvx-sets},
meaning there exists $\uu\in\Rn\wo\set{\zero}$
such that
\begin{equation} \label{eq:Fdub-equal-F:alt:2}
  \sup_{\zbar\in S} \zbar\inprod\uu < \xbar\inprod\uu.
\end{equation}
Let
$J=\set{\ybar\in\eRn:\:\ybar\inprod\uu=\xbar\inprod\uu}$.
Then $J$ is
disjoint from~$S$ (by Eq.~\ref{eq:Fdub-equal-F:alt:2}), so
$F(\ybar)>\beta$ for all $\ybar\in J$.
Thus,
$(\utransF)(\xbar\cdot\uu)=\inf_{\ybar\in J} F(\ybar)\ge\beta$
(with equality by Eq.~\ref{eq:trans-u-F-defn}),
so $\uu$ satisfies the claim.
\end{proofx}

\Cref{cl:F-equal-Fdub-alt:1} combined with \eqref{eq:Fdub-equal-F:alt:1} implies
that $\Fdub(\xbar)\ge\beta$ for all $\beta\in\R$ such that
$\inf F<\beta<F(\xbar)$. Taking
supremum over this (nonempty) set of $\beta$ values yields $\Fdub(\xbar)\ge F(\xbar)$,
completing the proof.%
\indexg{closedness, astral (primal)!general function@of general function|)}%
\qedhere
\end{proof-parts}
\end{proof}

\indexg{closedness, astral (primal)!extension@of extension|(}%
\indexg{lower semicontinuous extension!astral closedness of|(}%
As a special case, for a convex function $f:\Rn\rightarrow\Rext$, we
obtain the following characterization for when $\fext=\fdub$:

\begin{corollary}  \label{cor:fext-eq-fdub-1d-projs}
  Let $f:\Rn\rightarrow\Rext$ be convex.
  Then the following are equivalent:
  \begin{letter-compact}
  \item    \label{cor:fext-eq-fdub-1d-projs:a}
    $\fext=\fdub$.
  \item    \label{cor:fext-eq-fdub-1d-projs:b}
    For all $\uu\in\Rn\wo\{\zero\}$,
    $\lsc(\utransf)$ is closed.
  \item    \label{cor:fext-eq-fdub-1d-projs:c}
    For all $\uu\in\Rn\wo\{\zero\}$,
    $\trans{\uu}\fext = (\utransf)^{*\dualstar}$.
  \end{letter-compact}
\end{corollary}

\begin{proof}
Let $F=\fext$,
which is convex and lower semicontinuous
(Theorems~\ref{thm:fext-convex} and~\ref{prop:ext:F}\ref{prop:ext:F:a}).
Then $\Fstar=\fstar$ by
\Cref{pr:fextstar-is-fstar}.
Moreover,
$\utransF = \utransfext = \ufextnp$
by \Cref{thm:inf-lin-ext},
so
$\resfcn{(\utransF)}{\R}=\resfcn{\ufext}{\R}=\lsc(\utransf)$
by \Cref{pr:h:1}(\ref{pr:h:1a}), and
also
$(\utransF)^{*} = \ufextstar = \ufstar$,
again by \Cref{pr:fextstar-is-fstar}.
Combining with \Cref{thm:F-equal-Fdub-alt}, the
\namecref{cor:fext-eq-fdub-1d-projs}
now follows.%
\indexg{lower semicontinuous extension!astral closedness of|)}%
\indexg{closedness, astral (primal)!extension@of extension|)}%
\end{proof}

\indexg{closedness, astral (primal)!lower-bounded function@of lower-bounded function|(}%
As a further consequence,
a function that is bounded from below
is astral closed if and only if it is convex and lower semicontinuous:

\begin{corollary}   \label{cor:lower-bnd-ast-closed}
  Let $F:\extspace\rightarrow\Rext$ with $\inf F>-\infty$.
  Then $F = \Fdub$ if and only if $F$ is convex and lower semicontinuous.
\end{corollary}

\begin{proof}
If $F=\Fdub$ then $F$ is convex
and lower semicontinuous by
\Cref{thm:F-equal-Fdub-alt}(\ref{thm:F-equal-Fdub-alt:a},\ref{thm:F-equal-Fdub-alt:b}).

For the reverse implication, suppose that $F$ is convex and lower
semicontinuous.
Let $\uu\in\Rn\wo\{\zero\}$,
let $G=\utransF$, and let $g=\resfcn{G}{\R}$.
Then $G$ is convex and lower semicontinuous
(by Theorems~\ref{th:lin-img-ast-fcn-cvx}
and~\ref{thm:AF-lsc}\ref{thm:AF-lsc:a}),
implying the same for $g$
(by \Cref{pr:lsc-res-domain}\ref{pr:lsc-res-domain:b}
and \Cref{pr:1d-cvx}).
Also,
$\inf g\geq\inf G = \inf F>-\infty$
(with the equality by
\Cref{pr:ast-lin-img-props}\ref{pr:ast-lin-img-props:inf-AF}),
so $g>-\infty$, and thus
$g=\resfcn{(\utransF)}{\R}$ is closed.
Since this holds for all $\uu\in\Rn\wo\{\zero\}$,
we conclude that $F=\Fdub$ by
\Cref{thm:F-equal-Fdub-alt}(\ref{thm:F-equal-Fdub-alt:b},\ref{thm:F-equal-Fdub-alt:a}).%
\indexg{closedness, astral (primal)!lower-bounded function@of lower-bounded function|)}%
\end{proof}

\indexg{closedness, astral (primal)!preserved by linear image|(}%
\indexg{linear image of astral function!astral closedness preserved|(}%
In \Cref{pr:prml-lin-op-conj-idens}(\ref{pr:prml-lin-op-conj-idens:c}), we proved
that
if $F:\extspace\rightarrow\Rext$ is astral closed, then
so is~$\FA$, for any $\A\in\Rnm$.
As a final consequence of
\Cref{thm:F-equal-Fdub-alt}, we can now prove an analogous
statement for $\A F$ (with $\A\in\Rmn$), under the additional
assumption that $\A$ has full row rank:

\begin{corollary}   \label{cor:AF-ast-closed}
  Let $F:\extspace\rightarrow\Rext$ with $F=\Fdub$,
  and let
  $\A\in\Rmn$ have full row rank.
  Then $\A F = \AFdub$.
\end{corollary}

\begin{proof}
By
\Cref{thm:F-equal-Fdub-alt}(\ref{thm:F-equal-Fdub-alt:a},\ref{thm:F-equal-Fdub-alt:c}),
$F$ is convex and lower semicontinuous, and
for all $\uu\in\Rn\wo\{\zerov{n}\}$,
\begin{equation}   \label{eq:cor:AF-ast-closed:1}
  \utransF = (\utransF)^{*\dualstar}.
\end{equation}
Therefore, by Theorems~\ref{th:lin-img-ast-fcn-cvx}
and~\ref{thm:AF-lsc}(\ref{thm:AF-lsc:a}),
$\A F$ is also convex and lower semicontinuous.
Let $\vv\in\Rm\wo\{\zerov{m}\}$.
Then
\[
  \trans{\vv} (\A F)
  =
  (\trans{\vv} \kernA) F
  =
  \bigParens{(\trans{\vv} \kernA) F}^{*\dualstar}
  =
  \bigParens{\trans{\vv} (\A F)}^{*\dualstar}.
\]
The first and last equality are by
\Cref{pr:ast-lin-img-props}(\ref{pr:ast-lin-img-props:c}).
The second equality is by
\eqref{eq:cor:AF-ast-closed:1}
with $\uu=\transA \vv$, noting that $\uu\neq\zerov{n}$ since
$\vv\neq\zerov{m}$ and $\transA$ has linearly independent columns.
Thus, by
\Cref{thm:F-equal-Fdub-alt}(\ref{thm:F-equal-Fdub-alt:c},\ref{thm:F-equal-Fdub-alt:a}),
$\A F = \AFdub$.
\end{proof}

\Cref{cor:AF-ast-closed} is not true in general without the assumption
that $\A$ has full row rank.
In other words, without this assumption, it is possible that $F$ is
astral closed, but that $\A F$ is not.
Here is an example:

\begin{example}   \label{ex:AF-not-ast-closed}
For $\barx\in\Rext$, let $F(\barx)=\barx$.
It can be checked that $\Fstar=\smash[b]{\indfsing{1}}$, the indicator of the
singleton $\set{1}$, and that $\Fdub=F$.
Thus, $F$ is astral closed.

Let $\ee_1$ and $\ee_2$ be the standard basis vectors in $\R^2$,
and let $\A=[\ee_1]$.
This matrix does not have full row rank.
For $\xbar\in\extspac{2}$, it can be calculated that
\[
  (\A F)(\xbar)
  =
  \inf\regBraces{ F(\alpha) :\: \alpha\in\Rext,\, \alpha\ee_1=\xbar }
  =
  \begin{cases}
    \xbar\cdot\ee_1    & \text{if $\xbar\cdot\ee_2=0$,} \\
    +\infty            & \text{otherwise.}
  \end{cases}
\]
The conjugate of this function, $\AFstar$, can be shown to be
$\indf{L}$,
the indicator of the line $L=\{\trans{[1,v]} :\: v\in\R\}$.
Letting $\xbar=\limray{(-\ee_1)}\plusl\ee_2$,
this implies that $\AFdub(\xbar)=\indfdstar{L}(\xbar)=-\infty$,
but
$(\A F)(\xbar)=+\infty$.
Thus, $\A F$ is not astral closed.
\end{example}

The analogue of
\Cref{cor:AF-ast-closed}
does not hold for standard conjugates and standard convex functions.
In other words, if $f:\Rn\rightarrow\Rext$ is convex and closed
(so that $f=\fdubs$),
and $\A\in\Rmn$ has full row rank, $\A f$ need not be closed.
Here is an example:

\begin{example}
\label{ex:Af:not-closed}
On $\R^2$, let $f$ be the restricted linear function
from \Cref{ex:negx1-else-inf,ex:biconj:notext},
\[
   f(x_1,x_2) =
   \begin{cases}
     -x_1
     & \text{if $x_2\geq 0$,}
   \\
     +\infty
     & \text{otherwise.}
   \end{cases}
\]
This function is convex, closed and proper.
Let $\A=[0,1]$, which has full row rank.
Then,
for $z\in\R$,
\[
  (\A f)(z)
  =
  \inf\bigSet{f(x_1,z):\:x_1\in\R}
  =
  \begin{cases}
    -\infty & \text{if $z\geq 0$,} \\
    +\infty & \text{otherwise.}
  \end{cases}
\]
This function is convex, but it is not closed,
\indexg{closedness, astral (primal)!preserved by linear image|)}%
\indexg{linear image of astral function!astral closedness preserved|)}%
so $\A f\neq (\A f)^{**}$.
\end{example}

\indexg{closedness, astral dual!not preserved by linear image|(}%
\indexg{linear image of standard function!astral dual closedness not preserved|(}%
\Cref{ex:Af:not-closed} also shows that astral dual closedness of
a function $\psi:\Rn\rightarrow\Rext$
does not necessarily imply that $\A\psi$ is astral dual closed as well.
To see this, let $\psi$ be the function $f$ from the example,
which is (standard) closed, and so also astral dual closed
by \Cref{pr:std-clos-implies-ast-dual-clos}.
Nonetheless, the function $\A \psi$ (with $\A$ as in the example)
is not astral dual closed.
This is because $(\A \psi)(0)=-\infty$, implying
$(\A\psi)^{\dualstar}\equiv+\infty$
(by Eq.~\ref{eq:psistar-def:2}), and so that
$(\A\psi)^{\dualstar *}\equiv-\infty$.
Hence, $\A\psi\neq (\A\psi)^{\dualstar *}$.%
\indexg{closedness, astral dual!not preserved by linear image|)}%
\indexg{linear image of standard function!astral dual closedness not preserved|)}%

\section{Lower semicontinuity and primal biconjugacy}
\label{sec:lsc-Fdub}

In \Cref{sec:ent-closed-fcn}, we showed that 
whenever an extension $\ef$ and a biconjugate $\fdub$ of a convex function $f:\Rn\to\eR$
differ at a point $\xbar\in\eRn$, we must have
$\ef(\xbar)=+\infty$ and $\fdub(\xbar)=-\infty$. In this section, we
show that the same is true more generally for points where
any convex, lower semicontinuous function on astral space differs from its biconjugate.

\indexg{lower envelope!outer hull of epigraph@of outer hull of epigraph|(}%
\indexg{lower support function!outer hull of epigraph@of outer hull of epigraph|(}%
Our analysis builds on results from
\Cref{sec:fcns-induced-by-sets}, namely
\Cref{thm:settofcn-eq-biconj}.
To apply this
\namecref{thm:settofcn-eq-biconj} to a function, we need
to represent the function and its conjugate as the lower
envelope and lower support function (respectively)
of some closed and convex set in $\extspacnp$.
We might first try to use $\epi{F}$ itself for this purpose
since if $F:\eRn\to\eR$ is convex and lower semicontinuous then its
epigraph is convex in $\extspacnp$ (by definition) and closed in $\eRn\times\R$
(by \Cref{prop:lsc}\ref{prop:lsc:a}\ref{prop:lsc:b});
the problem,
however, is that $\epi{F}$ is not closed in $\extspacnp$ (except for the degenerate case when $F\equiv+\infty$).
Representing the function (and its conjugate) using $\epibarbar{F}$
might not work
since closures 
of convex sets in astral space are not necessarily convex (\Cref{thm:closure-not-always-convex}).
Nonetheless, as the next
\namecref{thm:cvx-lsc-implies-ult-cvx}
shows, the outer hull of $\epi F$ does always work for this purpose;
that is, if $F$ is convex and lower semicontinuous, then $F$ and its
conjugate can be represented as the lower envelope and lower support
function
of $\ohull(\epi F)$.

\begin{theorem}   \label{thm:cvx-lsc-implies-ult-cvx}
  Let $F:\extspace\rightarrow\Rext$ be convex and lower
  semicontinuous,
  and let $E=\ohull(\epi F)$.
  Then $F=\settofcnE$ and $\Fstar=\Esstar$.
\end{theorem}

\begin{proof}
By \Cref{pr:settofcn-ohull-prop}, 
$\Fstar=\Esstar$ and $F\geq \settofcnE$,
since $\epi F\subseteq E=\ohull(\epi{F})$.
Therefore, it remains only to show that $F\leq \settofcnE$.

Let 
$\PPx=[\Idnn,\zerov{n}]$ and
$\PPy=[\trans{\zerov{n}},1]$,
let $\xbar\in\extspace$, and let $\zbar\in\extspacnp$
with $\PPx\zbar=\xbar$.
We will prove that
if $\PPy\zbar < F(\xbar)$ then $\zbar\not\in E$.
In other words, we show (in the contrapositive) that
$\PPy\zbar\geq F(\xbar)$ whenever $\zbar\in E$,
proving that $\settofcnEf(\xbar)\geq F(\xbar)$
by definition of lower envelope (Eq.~\ref{eq:settofcn-dfn}).

Suppose then that $\PPy\zbar<F(\xbar)$.
Let $y\in\R$ be such that
$\PPy\zbar\leq y <F(\xbar)$, and let $c\in\R$ be such that
$c<y$.
Let $G=\max\{c,F\}$ be the pointwise maximum of $F$
and the constant function $c$.
Since both these functions are convex and lower semicontinuous,
$G$ is as well,
by
\Cref{thm:point-sup-is-convex}
and
\Cref{pr:lsc-sup}.
Note that, for $\xbar'\in\extspace$, if either
$F(\xbar')>c$ or $G(\xbar')>c$ then
$G(\xbar')=F(\xbar')$.
Further, $\inf G\geq c>-\infty$.

These facts imply that $G=\Gdub$,
by \Cref{cor:lower-bnd-ast-closed}.
Since $F(\xbar)>y>c$, $G(\xbar)=F(\xbar)$, implying
$\Gdub(\xbar)>y$.
Thus, by
\Cref{thm:fdub-sup-afffcns}(\ref{thm:fdub-sup-afffcns:b}),
there exists $\uu\in\Rn$ and $v\in\R$ such that
$G\geq A$ and $y < A(\xbar)$, where $A:\extspace\rightarrow\Rext$ is
the affine function defined by $A(\xbar')=\xbar'\cdot\uu-v$
for $\xbar'\in\extspace$.
In particular, since
$y < A(\xbar)=\xbar\cdot\uu-v$, we must have
\begin{equation}   \label{eq:thm:cvx-lsc-implies-ult-cvx:5}
  v < \xbar\cdot\uu - y.
\end{equation}

The main idea of the proof is to construct a closed halfspace $H$ in
$\extspacnp$ that includes $\epi F$ but excludes $\rpair{\xbar}{y}$.
The next claim shows that, under appropriate conditions, this is
sufficient to complete the proof:

\begin{claimpx}   \label{cl:thm:cvx-lsc-implies-ult-cvx:1}
  Let $a\in\Rstrictpos$, $b\in\Rneg$, $\beta\in\R$,
  and let
  \begin{equation}   \label{eq:thm:cvx-lsc-implies-ult-cvx:3}
    H
    =
    \Braces{
      \zbar'\in\extspacnp :\:
      \zbar'\cdot\rpair{a\uu}{b} \leq \beta
    }.
  \end{equation}
  Suppose $\epi{F}\subseteq H$ and $\rpair{\xbar}{y}\not\in H$.
  Then $\zbar\not\in E$.
\end{claimpx}

\begin{proofx}
We have
\begin{equation}   \label{eq:thm:cvx-lsc-implies-ult-cvx:1}
  \zbar\cdot\rpair{a\uu}{0}
  =
  (\PPx\zbar)\cdot(a\uu)
  =
  a \xbar\cdot\uu
  >
  -\infty,
\end{equation}
with the first equality from
\Cref{pr:xy-pairs-props}(\ref{pr:xy-pairs-props:b-new}),
and the inequality because $a > 0$ and
$\xbar\cdot\uu - v = A(\xbar)>y$,
implying $\xbar\cdot\uu>-\infty$.
Also,
\begin{equation}   \label{eq:thm:cvx-lsc-implies-ult-cvx:2}
  \zbar\cdot\rpair{\zero}{b}
  =
  b (\PPy\zbar)
  \geq
  b y
  >
  -\infty
\end{equation}
with the equality again from
\Cref{pr:xy-pairs-props}(\ref{pr:xy-pairs-props:b-new}),
and the first inequality because
$b\leq 0$ and $\PPy\zbar\leq y$.
Thus,
$\zbar\cdot\rpair{a\uu}{0}$
and
$\zbar\cdot\rpair{\zero}{b}$
are summable, so
\begin{equation*}
  \zbar\cdot\rpair{a\uu}{b}
  =
  \zbar\cdot\rpair{a\uu}{0}
  +
  \zbar\cdot\rpair{\zero}{b}
  \geq
  a \xbar\cdot\uu
  +
  b y
  =
  \rpair{\xbar}{y}\cdot\rpair{a\uu}{b}
  >
  \beta.
\end{equation*}
The first equality is by \Cref{pr:i:1}.
The first inequality is
by Eqs.~(\ref{eq:thm:cvx-lsc-implies-ult-cvx:1})
and~(\ref{eq:thm:cvx-lsc-implies-ult-cvx:2}).
The second equality is by
\Cref{pr:xy-pairs-props}(\ref{pr:xy-pairs-props:b}), and
the second inequality is because $\rpair{\xbar}{y}\not\in H$.
Thus, $\zbar$ is not in $H$.
Since $\epi{F}\subseteq H$, this shows that
$\zbar\not\in\ohull(\epi{F})$ by definition of outer hull
(\Cref{def:outer-cvx-hull}).
\end{proofx}

To construct such halfspaces, we consider two cases:

\begin{proof-parts}
\pfpart{First case:}
In the first case,
suppose for all $\xbar'\in\extspace$ that if
$-\infty<\xbar'\cdot\uu<\xbar\cdot\uu$ then $F(\xbar')>y$.
In this case, we claim that $F\geq A$.
To see this, let $\xbar'\in\extspace$; we show
$F(\xbar')\geq A(\xbar')$.
If $\xbar'\cdot\uu=-\infty$, then $A(\xbar')=-\infty$,
trivially implying this inequality.
Otherwise, if
$-\infty<\xbar'\cdot\uu<\xbar\cdot\uu$ then, by assumption,
$F(\xbar')>y>c$, implying $F(\xbar')=G(\xbar')\geq A(\xbar')$.
Finally, if $\xbar'\cdot\uu\geq\xbar\cdot\uu$, then
\[
  G(\xbar')
  \geq
  A(\xbar')
  =
  \xbar'\cdot\uu-v
  \geq
  \xbar\cdot\uu-v
  =
  A(\xbar)
  >y>c,
\]
again implying $F(\xbar')=G(\xbar')\geq A(\xbar')$.

Thus, as claimed, $F\geq A$.
Let $H=\bigBraces{\zbar'\in\extspacnp :\: \zbar'\cdot\rpair{\uu}{-1}\leq v}$.
Then ${\epi{F}\subseteq H}$ since if $\rpair{\xbar'}{y'}\in\epi{F}$ then
$y'\geq F(\xbar')\geq A(\xbar')=\xbar'\cdot\uu-v$,
implying
$v\geq\xbar'\cdot\uu-y'$.
On the other hand, $\rpair{\xbar}{y}\not\in H$
since $v < \xbar\cdot\uu - y$ by
\eqref{eq:thm:cvx-lsc-implies-ult-cvx:5}.
Therefore, $\zbar\not\in E$ by
\Cref{cl:thm:cvx-lsc-implies-ult-cvx:1}
(with $a=1$, $b=-1$, $\beta=v$),
completing this case.

\pfpart{Alternative case:}
In the alternative, suppose that the first case does not hold,
meaning that there exists $\xbar_0\in\extspace$ such that
$-\infty<\xbar_0\cdot\uu<\xbar\cdot\uu\leq+\infty$ and
$F(\xbar_0)\leq y$.
In particular, this implies $\xbar_0\cdot\uu\in\R$.

To derive a halfspace as in 
\Cref{cl:thm:cvx-lsc-implies-ult-cvx:1}
for this case, we work with
a one-variable version of $F$
similar to those appearing in the characterization
of astral closed functions in
\Cref{thm:F-equal-Fdub-alt}.
Specifically, we consider the function
$f=\resfcn{(\utransF)}{\R}$.
Since $F$ is convex, $f$ is as well,
by \Cref{th:lin-img-ast-fcn-cvx,pr:1d-cvx}.
Let $\barx=\xbar\cdot\uu$ and $x_0=\xbar_0\cdot\uu$,
implying $x_0\in\R$ and $x_0<\barx$.
We will derive a separation (in $\extspac{2}$)
of $\rpair{\barx}{y}$ from
$\epi f$, which will lead to a corresponding separation of $\epi{F}$
from $\rpair{\xbar}{y}$.

Let
\[
  U
  =
  \Braces{
    \zbar'\in\extspac{2}
    :\:
    \zbar'\cdot\rpair{1}{-1} > v
    \text{ and }
    \zbar'\cdot\rpair{0}{1} > c
  }.
\]
As an intersection of two open halfspaces, $U$ must be open.
Further, $\rpair{\barx}{y}\in U$ since
$y>c$ and
$\barx-y=\xbar\cdot\uu-y>v$
by \eqref{eq:thm:cvx-lsc-implies-ult-cvx:5}.
Thus, $U$ is a neighborhood of $\rpair{\barx}{y}$.

We claim moreover that $U$ is disjoint from $\epi f$.
To see this, suppose to the contrary that a point
$\rpair{x'}{y'}$ is both in $U$ and $\epi f$,
meaning
$x',y'\in\R$,
$x'-y'>v$, $y'>c$, and $y'\geq f(x')$.
Since $F$ is lower semicontinuous and
$(\utransF)(x')=f(x')<+\infty$,
\Cref{thm:AF-lsc}(\ref{thm:AF-lsc:b}) then implies that
there exists $\xbar'\in\extspace$ such that
$\xbar'\cdot\uu=\trans{\uu}\xbar'=x'$
and
$F(\xbar')=f(x')$.
Thus,
\[
  y'
  \geq
  \max\{c, F(\xbar')\}
  =
  G(\xbar')
  \geq
  A(\xbar')
  =
  \xbar'\cdot\uu - v
  =
  x' - v,
\]
contradicting that $x'-y'>v$.

Hence, as claimed, $U\cap (\epi f)=\emptyset$.
Therefore, $\rpair{\barx}{y}\not\in\clepif$
(by \Cref{pr:closure:intersect}\ref{pr:closure:intersect:a}).
Since $\epi f$ is a convex subset of $\R^2$,
$\clepif=\ohull(\epi f)$
by \Cref{thm:e:6}, implying
$\rpair{\barx}{y}\not\in\ohull(\epi f)$.
Therefore, by definition of outer hull
(\Cref{def:outer-cvx-hull}),
there exists $a,b,\beta\in\R$ such that the halfspace
\[
  H'
  =
  \Braces{
    \zbar'\in\extspac{2} :\:
    \zbar'\cdot\rpair{a}{b}\leq\beta
  }
\]
includes $\epi{f}$ and excludes $\rpair{\barx}{y}$.

In particular, this implies $b\leq 0$.
To see this, note first that
\begin{equation}   \label{eq:thm:cvx-lsc-implies-ult-cvx:6}
  f(x_0)
  =
  (\utransF)(\xbar_0\cdot\uu)
  =
  (\utransF)(\trans{\uu}\xbar_0)
  \leq
  F(\xbar_0)
  \leq
  y
  < +\infty,
\end{equation}
with the first inequality from
\Cref{pr:ast-lin-img-props}(\ref{pr:ast-lin-img-props:a})
and the second by assumption for this case.
Hence, for all
$y'\in [y,+\infty)$,
$\rpair{x_0}{y'}\in\epi f\subseteq H'$, implying
$a x_0 + b y' \leq \beta$.
Therefore, if $b>0$ then the left-hand side can be made
arbitrarily large (as $y'\rightarrow+\infty$),
contradicting the inequality.

We also claim $a>0$.
This is because
$a x_0 + b y\leq \beta < a \barx + b y$
since, by \eqref{eq:thm:cvx-lsc-implies-ult-cvx:6},
$\rpair{x_0}{y}\in\epi{f}\subseteq H'$
but $\rpair{\barx}{y}\not\in H'$.
Since $x_0 < \barx$, this proves $a>0$.

Let $H$ be the halfspace defined in
\eqref{eq:thm:cvx-lsc-implies-ult-cvx:3}.
We show next that $\epi{F}\subseteq H$.
Let $\rpair{\xbar'}{y'}\in\epi{F}$.
We aim to show the inequality
\begin{equation}   \label{eq:thm:cvx-lsc-implies-ult-cvx:4}
   a \xbar'\cdot\uu + b y'
   =
   \rpair{\xbar'}{y'}\cdot\rpair{a\uu}{b}
   \leq
   \beta.
\end{equation}
This is immediate if $\xbar'\cdot\uu=-\infty$
(since $a>0$).
If $\xbar'\cdot\uu\in\R$ then
\[
  y'
  \geq
  F(\xbar')
  \geq
  (\utransF)(\trans{\uu}\xbar')
  =
  f(\xbar'\cdot\uu)
\]
with the second inequality
from \Cref{pr:ast-lin-img-props}(\ref{pr:ast-lin-img-props:a}).
Thus, $\rpair{\xbar'\cdot\uu}{y'}\in\epi f\subseteq H'$,
implying
\eqref{eq:thm:cvx-lsc-implies-ult-cvx:4}.
Finally, if $\xbar'\cdot\uu=+\infty$ then
$G(\xbar')\geq A(\xbar')=\xbar'\cdot\uu-v=+\infty$,
implying $y'\geq F(\xbar')=G(\xbar')=+\infty$, a contradiction,
showing actually that this last case is impossible.

On the other hand, since $\rpair{\barx}{y}\not\in H'$,
$
  \rpair{\xbar}{y}\cdot\rpair{a\uu}{b}
  =
  a \xbar\cdot\uu + b y
  =
  a \barx + b y
  >
  \beta
$,
so $\rpair{\xbar}{y}\not\in H$.
Thus, by \Cref{cl:thm:cvx-lsc-implies-ult-cvx:1},
$\zbar\not\in E$,
completing this case and the proof.%
\indexg{lower envelope!outer hull of epigraph@of outer hull of epigraph|)}%
\indexg{lower support function!outer hull of epigraph@of outer hull of epigraph|)}%
\qedhere
\end{proof-parts}
\end{proof}

\indexg{closedness, astral (primal)!general function@of general function|(}%
Combining \Cref{thm:cvx-lsc-implies-ult-cvx}
with \Cref{thm:settofcn-eq-biconj},
we obtain the following characterization for when
a convex, lower semicontinuous function on astral space
is equal to its biconjugate.
When applied to the extension of a convex
function $f:\Rn\rightarrow\Rext$,
we recover
\Cref{thm:fext-neq-fdub} as a special case.

\begin{corollary}   \label{cor:ult-cvx-F-Fdub}
  Let $F:\extspace\rightarrow\Rext$ be convex and lower semicontinuous,
  and let $\xbar\in\extspace$.
  If either $\xbar\in\ohull(\dom{F})$ or $\Fdub(\xbar)>-\infty$,
  then $F(\xbar)=\Fdub(\xbar)$.

  Consequently, if $F(\xbar)\neq\Fdub(\xbar)$
  then $F(\xbar)=+\infty$ and $\Fdub(\xbar)=-\infty$.
\end{corollary}

\begin{proof}
By \Cref{thm:cvx-lsc-implies-ult-cvx}, there exists a closed and convex set
$E\subseteq\extspacnp$ such that $F=\settofcnE$ and $\Fstar=\Esstar$.
Let $\PPx=[\Idnn,\zerov{n}]$.
Then $\PPxE$ is also closed and convex
(by
Corollaries~\ref{cor:aff-img-closed-is-closed}\ref{cor:aff-img-closed-is-closed:a}
and~\ref{cor:thm:e:9}).
Also, if $\xbar\in\dom{F}$ then $\settofcnEf(\xbar)=F(\xbar)<+\infty$,
implying $\xbar=\PPx\zbar$ for some $\zbar\in E$, by $\settofcnE$'s
definition (Eq.~\ref{eq:settofcn-dfn}).
Thus, $\dom{F}\subseteq\PPxE$, implying
$\ohull(\dom{F})\subseteq\PPxE$ by
\Cref{cor:sep-cvx-sets-conseqs}(\ref{cor:sep-cvx-sets-conseqs:a})
since $\PPxE$ is closed and convex.
Therefore, if $\xbar$ is in $\ohull(\dom{F})$, then it is also in
$\PPxE$, so that
\begin{equation}   \label{eq:cor:ult-cvx-F-Fdub:1}
  F(\xbar)=\settofcnEf(\xbar)=\Esdubf(\xbar)=\Fdub(\xbar),
\end{equation}
with the second equality from
\Cref{thm:settofcn-eq-biconj}.

Similarly, if $\Fdub(\xbar)>-\infty$, then also
$\Esdubf(\xbar)>-\infty$,
so
\Cref{thm:settofcn-eq-biconj}
again implies
\eqref{eq:cor:ult-cvx-F-Fdub:1}.%
\indexg{closedness, astral (primal)!general function@of general function|)}%
\end{proof}

\section{Astral dual closedness}

We next consider when a function $\psi:\Rn\rightarrow\Rext$ is astral
dual closed.
\indexg{closedness, astral dual!one dimension@in one dimension|(}%
As with primal closed functions,
we begin by characterizing dual closed functions
when $n=1$.
In the \namecref{thm:psidub:1d} below,
for any function $\rho:\R\to\eR$,
we write $\rhoneg$ and $\rhopos$ to denote its restrictions to
$\Rstrictneg$ and $\Rstrictpos$, respectively.

We provide two characterizations. The first
states that dual closed functions on $\R$ are exactly those
convex functions $\psi:\R\to\eR$ that
are lower semicontinuous on $\Rstrictneg$ and $\Rstrictpos$, but not necessarily
at $0$, while also satisfying the condition that the set of points where they are
equal to $-\infty$ must be $\Rstrictneg$, or $\Rstrictpos$, or $\R$,
or $\emptyset$.
In this sense, they generalize the notion of a standard closed convex function,
which must be lower semicontinuous everywhere and equal
to $-\infty$ either on $\R$ or on $\emptyset$.
The second characterization additionally spells out that if $\psi$ is
not lower semicontinuous at $0$, then either $\psipos$ or $\psineg$ must
be identically equal to $+\infty$.

\begin{theorem}
  \label{thm:psidub:1d}
  Let $\psi:\R\to\eR$. Then the following are equivalent:
  \begin{letter}
  \item \label{thm:psidub:1d:a}
    $\psi=\psidub$.
  \item \label{thm:psidub:1d:b}
    All of the following hold:
    \begin{roman-compact}[labelsep=*]
      \item \label{thm:psidub:1d:b:i}
        $\psi$ is convex, and $\psi$ is lower semicontinuous everywhere except possibly at $0$.
      \item \label{thm:psidub:1d:b:ii}
        If $\psi(u)=-\infty$ for any $u\in\Rneg$, then $\psineg\equiv-\infty$.
      \item \label{thm:psidub:1d:b:iii}
        If $\psi(u)=-\infty$ for any $u\in\Rpos$, then $\psipos\equiv-\infty$.
    \end{roman-compact}
  \item \label{thm:psidub:1d:c}
    There exists a closed convex function $\rho:\R\to\eR$ such that
    one of the following holds:
    \begin{roman-compact}[labelsep=*]
      \newcommand{\first}{$\psineg\equiv+\infty$}
      \itemx{i$'$} \label{thm:psidub:1d:c:i}
         $\psi=\rho$.
      \itemx{ii$'$} \label{thm:psidub:1d:c:ii}
         \forcewidthof{\first,}{$\psineg=\rhoneg$,}\quad
         $\psi(0)>\rho(0)$,\quad
         $\psipos\equiv+\infty$.
      \itemx{iii$'$} \label{thm:psidub:1d:c:iii}
         $\psineg\equiv+\infty$,\quad
         $\psi(0)>\rho(0)$,\quad
         $\psipos=\rhopos$.
      \end{roman-compact}
    \end{letter}
\end{theorem}

 In proving the theorem, we use the following simple lemma:

  \begin{lemma}
  \label{cl:psidub:1d:1}
  Let $\psi:\R\to\eR$ be convex and let $u_0\in\R$ be a point
  where $\psi$ is not lower semicontinuous. Then
  either $\dom\psi\subseteq(-\infty,u_0]$ or $\dom\psi\subseteq[u_0,+\infty)$.
  \end{lemma}
  \begin{proof}
  We prove the contrapositive. Suppose
  there exist $u_{-},u_{+}\in\dom\psi$ such that $u_{-}<u_0<u_{+}$. We
  need to show that $\psi$ is lower semicontinuous at $u_0$.
  Since $\psi$ is convex, its effective domain is also convex
  (\Cref{roc:thm4.6}), so $(u_{-},u_{+})\subseteq\dom\psi$.
  Thus, $u_0\in\intr(\dom\psi)$; therefore,
  $\psi$ is continuous (and so also lower semicontinuous) at $u_0$
  (by \Cref{pr:stand-cvx-cont}).
  \end{proof}

\begin{proof}[Proof of \Cref{thm:psidub:1d}]
~
\begin{proof-parts}

\pfpart{%
  (\ref{thm:psidub:1d:a})
  $\Rightarrow$
  (\ref{thm:psidub:1d:b}):}
  Suppose $\psi=\psidub$. Then $\psi$ is convex (by \Cref{pr:conj-is-convex}). Moreover,
  from the definition of astral conjugate
  (Eq.~\ref{eq:Fstar-down-def}),
  for all $u\in\R$,
  \begin{equation}
  \label{eq:psidub:1d:1}
    \psi(u)
    =
    \psidub(u)
    =
    \sup_{\barx\in\eR} \BigBracks{-\psistarb(\barx)\plusd \barx u}.
  \end{equation}
  For each $\barx\in\Rext$, the function $u\mapsto-\psistarb(\barx)\plusd \barx u$
  is continuous on $\Rstrictpos\cup\Rstrictneg$, and hence also lower
  semicontinuous on $\Rstrictpos\cup\Rstrictneg$. By \eqref{eq:psidub:1d:1}, $\psi$ is a pointwise
  supremum of such functions, and so also lower semicontinuous on $\Rstrictpos\cup\Rstrictneg$
  (\Cref{pr:lsc-sup}).
  Thus, condition~(\ref{thm:psidub:1d:b:i}) holds.

  Next, suppose $\psi(u)=-\infty$ for some $u\in\Rneg$. If $u=0$, then \eqref{eq:psidub:1d:1} implies that
  $\psistarb\equiv+\infty$, which in turn implies (by Eq.~\ref{eq:psidub:1d:1}) that $\psi\equiv-\infty$.
  If $u<0$, then \eqref{eq:psidub:1d:1} implies that
  $\psistarb(\barx)=+\infty$ for all $\barx\in[-\infty,+\infty)$, which in turn implies (by Eq.~\ref{eq:psidub:1d:1})
  that $\psi(u)=-\infty$ for all $u\in\Rstrictneg$. Thus, in either case, condition~(\ref{thm:psidub:1d:b:ii}) holds.
  
  By a symmetric argument, if $\psi(u)=-\infty$ for some $u\in\Rpos$, then
  $\psipos\equiv-\infty$, so condition~(\ref{thm:psidub:1d:b:iii}) holds as well.

\pfpart{%
  (\ref{thm:psidub:1d:b})
  $\Rightarrow$
  (\ref{thm:psidub:1d:c}):}
  Suppose that $\psi$ satisfies the conditions of part~(\ref{thm:psidub:1d:b}).
  In particular, $\psi$ is convex, and $\psi$ is lower semicontinuous everywhere except possibly at $0$.
  We need to show that $\psi$ takes one of the forms in part~(\ref{thm:psidub:1d:c}).
  
  If $\psi(0)=-\infty$, then $\psi\equiv-\infty$ by conditions (\ref{thm:psidub:1d:b:ii})
  and (\ref{thm:psidub:1d:b:iii}) of part~(\ref{thm:psidub:1d:b}),
  so $\psi$ takes form~(\ref{thm:psidub:1d:c:i})
  with $\rho\equiv-\infty$.
  We therefore
  assume henceforth that $\psi(0)>-\infty$.

  If $\psi(u)=-\infty$ for some $u\in\Rstrictneg$, then $\psineg\equiv-\infty$ by condition~(\ref{thm:psidub:1d:b:ii}) of part~(\ref{thm:psidub:1d:b}).
  Moreover, $\psi(0)>-\infty$, so $\psi$ is not lower semicontinuous at $0$; hence,
  $\dom\psi\subseteq\Rneg$ by
  \Cref{cl:psidub:1d:1} (noting that we cannot have $\dom\psi\subseteq\Rpos$
  since $\psineg\equiv-\infty$).
  Thus, $\psi(u')=+\infty$ for all $u'\in\Rstrictpos$,
  so $\psi$ takes form~(\ref{thm:psidub:1d:c:ii}) with $\rho\equiv-\infty$.
  
  By a symmetric argument, if $\psi(u)=-\infty$ for some $u\in\Rstrictpos$, then
  $\psi$ takes form~(\ref{thm:psidub:1d:c:iii}) with $\rho\equiv-\infty$.

  It remains to consider the case $\psi>-\infty$. If $\psi$ is lower semicontinuous, it must be closed, and so takes form~(\ref{thm:psidub:1d:c:i}) with $\rho=\psi$. Otherwise, the only point where $\psi$ is not lower semicontinuous is $0$. In particular, this means that $\psi\not\equiv+\infty$, so $\psi$ is proper. Let $\rho=\lsc\psi$, which is also proper (\Cref{pr:lsc-props}\ref{pr:lsc-props:d}) and therefore closed. Moreover, $\rho(u)=\psi(u)$ for all $u\ne 0$ and $\rho(0)<\psi(0)$ (by Propositions~\ref{pr:lsc-is-idempotent}
  and~\ref{prop:lsc:characterize}\ref{prop:lsc:characterize:a}). Finally, by 
  \Cref{cl:psidub:1d:1}, we have either $\dom\psi\subseteq\Rneg$ or
  $\dom\psi\subseteq\Rpos$. Thus, $\psi$ takes either form~(\ref{thm:psidub:1d:c:ii})
  on form~(\ref{thm:psidub:1d:c:iii}).

\pfpart{%
  (\ref{thm:psidub:1d:c})
  $\Rightarrow$
  (\ref{thm:psidub:1d:a}):}
  Let $\rho:\R\to\eR$ be a closed convex function such that $\psi$ takes one of the forms in
  part~(\ref{thm:psidub:1d:c}).
  In particular, this implies $\rhodubs=\rho$ by
  \Cref{pr:conj-props-cvx}(\ref{pr:conj-props-cvx:b}).
  For each possible form, we show there exists a function
  $F:\eR\to\eR$ for which $\psi=\Fstar$.
  By \Cref{pr:dual-biconj-props}(\ref{pr:dual-biconj-props:c}), this
  will prove that $\psi=\psidub$

  First, suppose $\psi$ takes form~(\ref{thm:psidub:1d:c:i}), and let $F=\rhostarext$.
  Then
  $\Fstar=(\rhostarext)^*=\rhodubs=\rho=\psi$,
  with the second equality by \Cref{pr:fextstar-is-fstar}.
  
  Next, suppose $\psi$ takes form~(\ref{thm:psidub:1d:c:ii}),
  and define $F:\Rext\rightarrow\Rext$, for $\barx\in\Rext$, by
  \[
    F(\barx)=\begin{cases}
       +\infty
         &\text{if $\barx=-\infty$,}
       \\
       \rhostar(\barx)
         &\text{if $\barx\in\R$,}
       \\
       -\psi(0)
         &\text{if $\barx=+\infty$.}
    \end{cases}
  \]
  Then, for $u\in\Rstrictneg$,
  \[
    \Fstar(u)=\sup_{\barx\in\eR}\BigBracks{-F(\barx)\plusd\barx u}
             =\sup_{x\in\R}\BigBracks{-\rhostar(x)+xu}
             =\rhodubs(u)=\rho(u)=\psi(u).
  \]
  The second equality is from the definition of $F$ and because $u<0$. The third is
  from the definition of standard conjugate (Eq.~\ref{eq:fstar-def}).
  The last is because, by assumption, $\psineg=\rhoneg$.
  
  Furthermore, $F(+\infty)=-\psi(0)<-\rho(0)=-\rhodubs(0)=\inf\rhostar$, with the inequality
  by assumption,
  and the last equality
  by \Cref{pr:conj-props}(\ref{pr:conj-props:a}). Thus,
  $F(+\infty)=\inf F$, so
  \[
    \Fstar(0)=\sup_{\barx\in\eR} [-F(\barx)]
             =-\inf F = -F(+\infty) = \psi(0).
  \]

  Finally, if $u\in\Rstrictpos$, then
  $\Fstar(u)\ge[-F(+\infty)\plusd(+\infty)u]=+\infty$,
  with the inequality following from the definition of $\Fstar$ (Eq.~\ref{eq:Fstar-down-def}),
  and the equality because $u>0$ and
  $-F(+\infty)=\psi(0)>\rho(0)\ge-\infty$. Thus, $\Fstar(u)=+\infty=\psi(u)$
  for $u\in\Rstrictpos$.
  We conclude that $\psi=\Fstar$ when $\psi$ takes form~(\ref{thm:psidub:1d:c:ii}).

  The remaining case, when $\psi$ takes
  form~(\ref{thm:psidub:1d:c:iii}), follows by a symmetric argument.%
\indexg{closedness, astral dual!one dimension@in one dimension|)}%
\qedhere
\end{proof-parts}
\end{proof}

\indexg{closedness, astral dual!characterization|(}%
We next provide a general characterization of when
$\psi=\psidub$.
The characterization is in terms of functions
$\psi \uu$, for $\uu\in\Rn\wo\set{\zero}$, that is, functions of the form
\[
  (\psi\uu)(\lambda)
  =
  \psi(\lambda \uu)
\]
for $\lambda\in\R$.
Such functions capture $\psi$'s behavior along the line
$\{ \lambda \uu :\: \lambda\in\R \}$.
We show that $\psi$ is dual closed if and only if
$\psi$ is convex and
$\psi\uu$ is dual closed for all $\uu\in\Rn\wo\set{\zero}$.
The dual closedness of the functions $\psi\uu$ can in turn be
verified using \Cref{thm:psidub:1d}.

\begin{theorem}   \label{thm:when-psi-is-psidub}
  Let $\psi:\Rn\rightarrow\Rext$.
  Then $\psi=\psidub$ if and only if
  $\psi$ is convex and for all $\uu\in\Rn\wo\set{\zero}$, $\psi\uu=\psiuudub$.
\end{theorem}

\indexg{epigraph!separation of disjoint set from|(}%
The next lemma, which we prove first,
will be used repeatedly in the proof of the theorem.
The lemma provides a separation result for any convex function on
$\Rn$ whose epigraph is disjoint from a given convex subset of $\Rnp$.
(As usual, we here often regard points in $\Rnp$ as pairs in
$\Rn\times\R$.)

\begin{lemma}  \label{lem:sep-W-epipsi}
  Let $\psi:\Rn\rightarrow\Rext$ be convex, and
  let $W\subseteq\Rnp$ be convex and disjoint from $\epi\psi$.
  Then there exists $\xbar\in\extspace$ such that
  \begin{equation}  \label{eq:lem:sep-W-epipsi:1}
    \psi(\uu')
    \geq
    v + \xbar\cdot(\uu'-\uu)
  \end{equation}
  for all $\uu'\in\Rn$ and for all $\rpair{\uu}{v}\in W$.
\end{lemma}

\begin{proof}
If either $\psi\equiv+\infty$ or $W=\emptyset$ then
the claim holds trivially or vacuously for any $\xbar\in\extspace$.
We therefore assume henceforth that $\psi\not\equiv+\infty$ and that
$W$ is nonempty.

Let $\PPx=[\Idnn,\zerov{n}]$.
We prove the lemma separately for when $\psi$'s effective domain is or
is not disjoint from $\PPx W$.

\begin{proof-parts}
\pfpart{Case $(\dom\psi)\cap(\PPx W)\neq\emptyset$:}
Suppose first that there exists a point
$\uhat$ in $(\dom\psi)\cap(\PPx W)$, meaning that
$\psi(\uhat)<+\infty$ and that
$\rpair{\uhat}{\ving}\in W$ for some $\ving\in\R$.
Then $\ving<\psi(\uhat)$
since $W$ and $\epi\psi$ are disjoint,
implying $\psi(\uhat)\in\R$.

Since $W$ and $\epi\psi$ are convex, nonempty and disjoint,
they are strongly dually separated by
\Cref{thm:ast-def-sep-cvx-sets}.
Therefore, by 
\Cref{thm:ast-str-sep-equiv-no-trans},
there exists $\zbar\in\extspacnp$ such that
$\zbar\cdot(\ww'-\ww) < 0$
for all
$\ww'\in\epi\psi$ and $\ww\in W$.

In particular, taking $\ww'=\rpair{\uhat}{\psi(\uhat)}$
and $\ww=\rpair{\uhat}{\ving}$, this means that
\begin{align*}
  0
  >
  \zbar\cdot(\ww'-\ww)
  &=
  \zbar\cdot\bigParens{
    \rpair{\uhat}{\psi(\uhat)} - \rpair{\uhat}{\ving}
  }
  \\
  &=
  \zbar\cdot\rpair{\zero}{\psi(\uhat)-\ving}
  =
  \bigParens{\psi(\uhat)-\ving}\bigParens{\zbar\cdot\rpair{\zero}{1}}.
\end{align*} 
Thus, $\zbar\cdot\rpair{\zero}{1}<0$ since $\ving<\psi(\uhat)$.

Therefore,
since $\zbar\cdot(\ww'-\ww)\leq 0$ for all $\ww'\in\epi\psi$ and
$\ww\in W$,
we can apply \Cref{lem:sep-one-pt-finite},
with $U$, $W$, and $\xbar$, as they appear in that lemma,
set to $(\epi\psi)-W$, $\{\rpair{\zero}{1}\}$, and $\zbar$,
respectively.
This yields that there exists $\zbar'\in\extspacnp$ such that
$\zbar'\cdot(\ww'-\ww)\leq 0$ for all $\ww'\in\epi\psi$ and
$\ww\in W$,
and
$\zbar'\cdot\rpair{\zero}{1}=-y'$ for some $y'\in\Rstrictpos$.
By
\Cref{pr:xy-pairs-props}(\ref{pr:xy-pairs-props:b-new},\ref{pr:xy-pairs-props:d}),
we therefore can write $\zbar'=\rpair{\xbar'}{-y'}$ where
$\xbar'=\PPx\zbar'\in\extspace$.

Thus, for all $\rpair{\uu'}{v'}\in\epi{\psi}$ and for all
$\rpair{\uu}{v}\in W$, we have
\[
   \xbar'\cdot(\uu'-\uu) - y' (v'-v)
   =
   \zbar'\cdot\bigParens{\rpair{\uu'}{v'} - \rpair{\uu}{v}}
   \leq
   0,
\]
where the equality is by
\Cref{pr:xy-pairs-props}(\ref{pr:xy-pairs-props:b}).
Dividing the inequality above by $y'$, letting $\xbar=\xbar'/y'$,
and rearranging then yields
\[
   v'
   \geq
   v
   +
   \xbar \cdot (\uu'-\uu).
\]
Since this holds for all $v'\in\R$ with $v'\geq\psi(\uu')$,
\eqref{eq:lem:sep-W-epipsi:1} must also hold,
for all $\uu'\in\Rn$ and all $\rpair{\uu}{v}\in W$,
completing the proof in this case.

\pfpart{Case $(\dom\psi)\cap(\PPx W)=\emptyset$:}
Suppose now that $\dom\psi$ and $\PPx W$ are disjoint.
These sets are both convex
(by
Propositions~\ref{roc:thm4.6} and~\ref{pr:aff-preserves-cvx}),
and are nonempty since $\psi\not\equiv+\infty$ and $W\neq\emptyset$.
Therefore, $\dom\psi$ and $\PPx W$ are strongly dually separated, by
\Cref{thm:ast-def-sep-cvx-sets}, implying, by 
\Cref{thm:ast-str-sep-equiv-no-trans}, that
there exists $\xbar'\in\extspace$ such that
$\xbar'\cdot(\uu'-\uu)<0$
for all $\uu'\in\dom\psi$ and $\uu\in \PPx W$.

Letting $\xbar=\limray{\xbar'}$, this then implies
\eqref{eq:lem:sep-W-epipsi:1}
for all $\uu'\in\Rn$ and all $\rpair{\uu}{v}\in W$.
This is because
\eqref{eq:lem:sep-W-epipsi:1}
holds trivially if $\psi(\uu')=+\infty$;
otherwise, if $\uu'\in\dom\psi$ and $\rpair{\uu}{v}\in W$,
so that $\uu\in\PPx W$,
then $\xbar\cdot(\uu'-\uu)=\limray{\xbar'}\cdot(\uu'-\uu)=-\infty$.%
\indexg{epigraph!separation of disjoint set from|)}%
\qedhere
\end{proof-parts}
\end{proof}

\begin{proof}[Proof of \Cref{thm:when-psi-is-psidub}]
  ~
\begin{proof-parts}
\pfpart{``Only if'' ($\Rightarrow$):}
  Suppose $\psi=\psidub$.
  Then $\psi$ is convex by
  \Cref{pr:conj-is-convex},
  and for $\uu\in\Rn\wo\set{\zero}$, $\psi\uu=\psiuudub$
  by
  \Cref{pr:dual-lin-op-conj-idens}(\ref{pr:dual-lin-op-conj-idens:c}).

\pfpart{``If'' ($\Leftarrow$):}
Suppose that $\psi\uu=\psiuudub$ for all $\uu\in\Rn\wo\set{\zero}$.
Since $\psi\ge\psidub$
(\Cref{thm:psi-geq-psidub}\ref{thm:psi-geq-psidub:c}),
it suffices to show that
$\psi(\uu)\leq\psidub(\uu)$ for all $\uu$.
The main step of the proof is given by the next claim.

\begin{claimpx}
\label{cl:psidub-char:alt}
  Let $\uu\in\Rn$ and $\beta\in\R$ be such that $\beta<\psi(\uu)$.
  Then there exist $\xbar\in\extspace$ and $y\in\R$ such that:
  \begin{roman-compact}
  \item \label{cl:psidub-char:alt:i}
    $\psi(\uu')\geq y + \xbar\inprod\uu'$
    for all $\uu'\in\Rn$; and
  \item \label{cl:psidub-char:alt:ii}
    $\beta\le y+\xbar\inprod\uu$.
\end{roman-compact}
\end{claimpx}

\begin{proofx}
We will show the existence of $\xbar$ and $y$ by considering a few cases.
For each case, we
construct an appropriate set $W\subseteq\Rnp$
and invoke \Cref{lem:sep-W-epipsi}.

Suppose first that $\uu=\zero$.
In this case, we set $W=\{\rpair{\zero}{\beta}\}$,
which is convex and disjoint from $\epi\psi$.
Therefore, by \Cref{lem:sep-W-epipsi},
there exists $\xbar\in\extspace$ such that, for all $\uu'\in\Rn$,
\[
  \psi(\uu')\ge \beta+\xbar\inprod\uu',
\]
while also $\beta\le\beta+\xbar\inprod\uu$ (since $\uu=\zero$).
Thus, the given $\xbar$ and $y=\beta$ satisfy the claim
in this case.

We therefore assume henceforth that $\uu\ne\zero$.
Let $\rho=\psi\uu$.
By assumption, $\rhodub=\rho$.
We consider two cases: $\rho\equiv+\infty$ and $\rho\not\equiv+\infty$.
\begin{proof-parts}
\pfpart{
  Case $\rho\equiv+\infty$:
}
In this case, $\psi(\lambda\uu)=\rho(\lambda)=+\infty$ for all
$\lambda\in\R$.
Let
\[
   W = \bigBraces{\rpair{\lambda\uu}{\beta} :\: \lambda\in\R}.
\]
This set is convex (being a line in $\R^{n+1}$) and disjoint
from $\epi\psi$.
Therefore, by \Cref{lem:sep-W-epipsi}, there exists
$\xbar'\in\extspace$ such that, for all $\uu'\in\Rn$ and all $\lambda\in\R$,
\begin{equation}
\label{eq:psidub-char:alt:3}
  \psi(\uu')
  \geq
  \beta + \xbar'\cdot(\uu'-\lambda\uu).
\end{equation}

For any point $\ybar\in\extspace$ (including points in $\Rn$),
let $\ybarperp$ denote $\ybar$'s projection orthogonal to $\uu$.
Since each $\uu'\in\Rn$ can be decomposed as
$\uu'=\uu'^{\perp}+\lambda\uu$ for some $\lambda\in\R$
(\Cref{pr:lin-decomp-rel-vecs}),
plugging this decomposition into \eqref{eq:psidub-char:alt:3}
implies that, for all $\uu'\in\Rn$,
\[
  \psi(\uu')
  \geq
  \beta + \xbar'\cdot\uu'^{\perp}
  =
  \beta + \xbar'^{\perp} \cdot\uu',
\]
with the equality following from
\Cref{pr:h:5}(\ref{pr:h:5a}).
Thus, $y=\beta$ and $\xbar=\xbar'^{\perp}$
satisfy condition~(\ref{cl:psidub-char:alt:i}) of the claim.
Condition~(\ref{cl:psidub-char:alt:ii}) is also satisfied by these
choices since
$\xbar\cdot\uu=\xbar'\cdot\uperp=\xbar'\cdot\zero=0$
(with the first equality again by \Cref{pr:h:5}\ref{pr:h:5a}).

\pfpart{Case $\rho\not\equiv+\infty$:}
From the definition of conjugate (Eq.~\ref{eq:Fstar-down-def}), for all $\lambda\in\R$,
\begin{equation*}
  \psi(\lambda\uu)=\rho(\lambda)=\rhodub(\lambda)
  =\sup_{\rpair{\barx}{y}\in\epi\rhostarb} [ \barx\lambda-y].
\end{equation*}
Since $\beta<\psi(\uu)$, there thus exists a pair $\rpair{\barx_0}{y_0}\in\epi\rhostarb$ such that
$\beta<\barx_0-y_0$.

Let
\[
  W
  =
  \bigBraces{
    \rpair{\lambda\uu}{v}:\:
    \lambda,v\in\R,\, v<\barx_0\lambda-y_0
  }.
\]
This set is disjoint from $\epi\psi$ since if
$\zz\in W$, then $\zz=\rpair{\lambda\uu}{v}$ for some uniquely determined
$\lambda,v\in\R$ (since $\uu\neq\zero$), implying
\[ v < \barx_0\lambda-y_0 \leq \rho(\lambda)=\psi(\lambda\uu), \]
with the first inequality by $W$'s definition (and since
$\uu\neq\zero$), and the second by
\Cref{thm:psi-geq-psidub}(\ref{thm:psi-geq-psidub:a})
since $\rpair{\barx_0}{y_0}\in\epi\rhostarb$;
thus,
$\zz=\rpair{\lambda\uu}{v}\not\in\epi\psi$.

Also, $W$ is convex.
This is because we can write $W=\A W'$ where
\begin{equation*}
  W'
  =
  \Braces{
    \rpair{\lambda}{v} \in \R^2 :\:
    v < \barx_0\lambda -y_0
  }
  =
  \Braces{
    \rpair{\lambda}{v} \in \R^2 :\:
    \rpair{-\barx_0}{1} \cdot \rpair{\lambda}{v} < -y_0
  },
\end{equation*}
and where
$\A=[\rpair{\uu}{0},\rpair{\zero}{1}]\in\R^{(n+1)\times 2}$
(so that the linear map associated with $\A$ is
$\rpair{\lambda}{v}\mapsto\rpair{\lambda\uu}{v}$
for $\rpair{\lambda}{v}\in\R^2$).
The set $W'$
is an astral dual halfspace, and so is convex by
\Cref{pr:ast-def-hfspace-convex}(\ref{i:dual-hfspace-complement},\ref{i:dual-hfspace-conv}).
Therefore, $W$, its image under $\A$,
is also convex (\Cref{pr:aff-preserves-cvx}).

Hence, by \Cref{lem:sep-W-epipsi},
there exists $\xbar\in\extspace$ such that,
for all $\lambda,v\in\R$ satisfying $v<\barx_0\lambda-y_0$,
and for all $\uu'\in\Rn$,
\[
  \psi(\uu')\ge v+\xbar\cdot(\uu'-\lambda\uu).
\]
Taking supremum, we then obtain that, for each $\lambda\in\R$
and for all $\uu'\in\Rn$,
\begin{align}
  \psi(\uu')
  &\ge
  \sup\BigBraces{
    v+\xbar\cdot(\uu'-\lambda\uu) :\:
    v\in\R,\, v<\barx_0\lambda-y_0
  }
  \notag
  \\
  &=
  \sup\BigBraces{
    v\in\R :\: v<\barx_0\lambda-y_0
  }
  \plusd
  \xbar\cdot(\uu'-\lambda\uu)
  \notag
  \\
  &=
  (\barx_0\lambda-y_0)\plusd\xbar\cdot(\uu'-\lambda\uu),
\label{eq:psidub-char:alt:5}
\end{align}
with the first equality 
by \Cref{pr:plusd-props}(\ref{pr:plusd-props:d-gen}).
Setting $\lambda=0$, we thus obtain that condition~(\ref{cl:psidub-char:alt:i})
holds for the given $\xbar$ and $y=-y_0$.
It remains to show that condition~(\ref{cl:psidub-char:alt:ii}) also holds.

Since $\rho\not\equiv+\infty$, there exists $c\in\R$ such that $\psi(c\uu)=\rho(c)<+\infty$.
Setting $\uu'=c\uu$ in \eqref{eq:psidub-char:alt:5}, we thus have, for all $\lambda\in\R$,
\begin{equation}
\label{eq:psidub-char:alt:6}
  \psi(c\uu)\ge (\barx_0\lambda-y_0)\plusd\xbar\cdot(c-\lambda)\uu.
\end{equation}

If $\barx_0=x_0$ for some $x_0\in\R$, then setting
$\lambda=c+\delta$, for any $\delta\in\R$,
\eqref{eq:psidub-char:alt:6} can be rewritten as
\begin{equation*}
\notag
  \psi(c\uu)
  \ge x_0(c+\delta)-y_0+\xbar\cdot(-\delta\uu)
  = x_0 c - y_0 + \delta(x_0-\xbar\cdot\uu).
\end{equation*}
If $x_0\ne\xbar\cdot\uu$, then by an appropriate choice of $\delta$, we can
make the right-hand side arbitrarily large, contradicting that $\psi(c\uu)<+\infty$. Thus,
we must have $\barx_0=x_0=\xbar\cdot\uu$. Since $\barx_0$ and $y_0$ were chosen so that
$\beta<\barx_0-y_0$, we then have $\beta<\xbar\cdot\uu-y_0$, satisfying
condition~(\ref{cl:psidub-char:alt:ii})
(with $\xbar$ and $y$ as above).

Otherwise,
if $\barx_0\not\in\R$, then we must have $\barx_0=+\infty$ (since $\beta<\barx_0-y_0$).
Then setting $\lambda=\abs{c}+1$, \eqref{eq:psidub-char:alt:6} implies
\begin{equation}
\label{eq:psidub-char:alt:7}
  \psi(c\uu)\ge (+\infty)\plusd c'(\xbar\cdot\uu),
\end{equation}
where $c'=c-\abs{c}-1<0$. Since $\psi(c\uu)<+\infty$, \eqref{eq:psidub-char:alt:7} implies that
$c'(\xbar\cdot\uu)=-\infty$ and thus $\xbar\cdot\uu=+\infty$.
Therefore, $\beta<\xbar\cdot\uu-y_0$,
satisfying
condition~(\ref{cl:psidub-char:alt:ii})
(with $\xbar$ and $y$ as above).
\qedhere
\end{proof-parts}
\end{proofx}

To complete the proof,
let $\uu\in\Rn$. If $\psi(\uu)=-\infty$ then trivially
$\psi(\uu)\leq\psidub(\uu)$.
Otherwise, let $\beta\in\R$ with $\beta<\psi(\uu)$,
and let $\xbar\in\extspace$ and $y\in\R$ be as provided in
\Cref{cl:psidub-char:alt}.
Let $\xi:\Rn\rightarrow\Rext$
be defined by $\xi(\uu')=y+\xbar\cdot\uu'$ for $\uu'\in\Rn$.
Then
\Cref{cl:psidub-char:alt}
shows that
$\psi\geq\xi$ and that $\xi(\uu)\geq\beta$.
Since $\psi\geq\xi$,
\Cref{thm:psi-geq-psidub}(\ref{thm:psi-geq-psidub:b})
therefore implies that %
$\psidub(\uu)\geq\xi(\uu)\geq\beta$.
Since this holds for all such $\beta$, it follows that
$\psidub(\uu)\geq\psi(\uu)$.

Thus, $\psi(\uu)\le\psidub(\uu)$ for all $\uu\in\Rn$, so
$\psi=\psidub$, completing the proof.%
\indexg{closedness, astral dual!characterization|)}%
\qedhere
\end{proof-parts}
\end{proof}

\indexg{closedness, astral dual!positively homogeneous functions@of positively homogeneous functions|(}%
\indexg{positively homogeneous functions!astral dual closedness of|(}%
\indexg{positively homogeneous functions!support function@as support function|(}%
\indexg{support functions, astral!positively homogeneous functions as|(}%
As an application, the next
\namecref{thm:pos-homo-is-ast-closed}
shows that
every convex positively homogeneous function on
$\Rn$ that
vanishes at the origin is astral dual closed,
implying, moreover, that every such function is the astral support
function (as in Eq.~\ref{eqn:astral-support-fcn-def})
of some closed, convex, nonempty subset of $\extspace$:

\begin{theorem}   \label{thm:pos-homo-is-ast-closed}
  Let $\psi:\Rn\rightarrow\Rext$ be convex and positively homogeneous,
  with ${\psi(\zero)=0}$.
  Then:
  \begin{letter-compact}
  \item    \label{thm:pos-homo-is-ast-closed:a}
    $\psi=\psidub$.
  \item    \label{thm:pos-homo-is-ast-closed:b}
    $\psi=\indaSstar$ and $\psistarb=\indaS$
    where
    \begin{equation}   \label{eq:thm:pos-homo-is-ast-closed:2}
       S
       =
       \bigBraces{
         \xbar\in\extspace :\:
         \xbar\cdot\uu \leq \psi(\uu) \text{ for all } \uu\in\Rn
       },
    \end{equation}
    which is convex, closed (in $\extspace$), and nonempty.
  \end{letter-compact}
\end{theorem}

\begin{proof}
  ~

\begin{proof-parts}
\pfpart{Part~(\ref{thm:pos-homo-is-ast-closed:a}):}
Let $\uu\in\Rn\wo\set{\zero}$, and let $\rho=\psi\uu$.
We show that $\rho$ satisfies the conditions of
\Cref{thm:psidub:1d}(\ref{thm:psidub:1d:b}) and is therefore dual closed.

By \Cref{thm:cvx-compose-affine-cvx},
$\rho$ is convex since $\psi$ is.
And since $\psi$ is positively homogeneous,
for $\lambda\in\Rstrictpos$, $\rho(\lambda)=\lambda\rho(1)$,
so $\rho$ is continuous on $\Rstrictpos$,
being either a line or identically equal to $\pm\infty$.
Similarly, for $\lambda\in\Rstrictneg$,
$\rho(\lambda)=-\lambda\rho(-1)$, so $\rho$ is continuous on $\Rstrictneg$.
Thus, $\rho$ satisfies condition~(\ref{thm:psidub:1d:b:i}) of~\Cref{thm:psidub:1d}(\ref{thm:psidub:1d:b}).

Furthermore, if $\rho(u)=-\infty$ for some $u\le 0$, then actually
$u<0$ (since $\rho(0)=\psi(\zero)=0$), so by positive homogeneity,
$\rho(\lambda)=(\lambda/u)\rho(u)=-\infty$ for all $\lambda\in\Rstrictneg$. Thus,
$\rho$ satisfies condition~(\ref{thm:psidub:1d:b:ii}) of~\Cref{thm:psidub:1d}(\ref{thm:psidub:1d:b}),
and so, by a symmetric argument, condition~(\ref{thm:psidub:1d:b:iii})
as well. Hence, $\rho=\rhodub$
by
\Cref{thm:psidub:1d}(\ref{thm:psidub:1d:b},\ref{thm:psidub:1d:a}).

Since $\psi\uu$ is dual closed for all $\uu\in\Rn\wo\set{\zero}$,
the claim follows by \Cref{thm:when-psi-is-psidub}.

\pfpart{Part~(\ref{thm:pos-homo-is-ast-closed:b}):}
Let $\xbar\in\extspace$;
we first aim to prove that $\psistarb(\xbar)=\indaS(\xbar)$.
By definition of dual conjugate (Eq.~\ref{eq:psistar-def:2}),
\begin{equation}   \label{eq:thm:pos-homo-is-ast-closed:1}
  \psistarb(\xbar)
  =
  \sup_{\uu\in\Rn} \bigBracks{ - \psi(\uu) \plusd \xbar\cdot\uu }.
\end{equation}
If $\xbar\in S$, then from the definition of $S$,
for all $\uu\in\Rn$,
$0 + \xbar\cdot\uu \leq \psi(\uu)$, implying
$-\psi(\uu) \plusd \xbar\cdot\uu \leq 0$ by
\Cref{pr:plusd-props}(\ref{pr:plusd-props:e}).
Since $-\psi(\zero)\plusd\xbar\cdot\zero=0$, this implies
$\psistarb(\xbar)=0=\indaS(\xbar)$ in this case.

Otherwise, if $\xbar\not\in S$, then there exists $\uu\in\Rn$ with
$\psi(\uu)<\xbar\cdot\uu$.
This strict inequality implies that $\xbar\cdot\uu$ and
$-\psi(\uu)$ are summable, and also that
$\xbar\cdot\uu - \psi(\uu) > 0$.
For $\lambda\in\Rstrictpos$, we then have that
\[
   \xbar\cdot(\lambda\uu) - \psi(\lambda\uu)
   =
   \lambda(\xbar\cdot\uu) - \lambda \psi(\uu)
   =
   \lambda\bigParens{\xbar\cdot\uu - \psi(\uu)},
\]
where the first equality is because $\psi$ is positively homogeneous.
Since the term on the right tends to $+\infty$ as
$\lambda\rightarrow+\infty$, it follows from
\eqref{eq:thm:pos-homo-is-ast-closed:1}
that $\psistarb(\xbar)=+\infty=\indaS(\xbar)$.

Thus, $\psistarb=\indaS$.
Combined with part~(\ref{thm:pos-homo-is-ast-closed:a}),
this further implies that
$\psi=\psidub=\indaSstar$.

Since $\psistarb=\indaS$, we have $S=\set{\xbar\in\eRn:\:\psistarb(\xbar)\le 0}$,
so $S$ is convex and closed by \Cref{thm:f:9} and
\Cref{prop:lsc}(\ref{prop:lsc:a},\ref{prop:lsc:c}),
since $\psistarb$
is convex and lower semicontinuous (\Cref{thm:dual-conj-cvx} and \Cref{pr:dual-conj-lsc}).
Also, $S\ne\emptyset$, since if $S=\emptyset$ then we would have $\psistarb=\indaS\equiv+\infty$ and hence $\psi=\psidub\equiv-\infty$, contradicting $\psi(\zero)=0$.%
\indexg{closedness, astral dual!positively homogeneous functions@of positively homogeneous functions|)}%
\indexg{positively homogeneous functions!astral dual closedness of|)}%
\qedhere
\end{proof-parts}
\end{proof}

If $S\subseteq\extspace$ is convex, closed and nonempty, then
its support function $\indaSstar$, as given in
\eqref{eqn:astral-support-fcn-def}, is evidently positively
homogeneous with $\indaSstar(\zero)=0$, and is also convex by
\Cref{pr:conj-is-convex}.
\Cref{thm:pos-homo-is-ast-closed} thus shows that
a function $\psi:\Rn\rightarrow\Rext$ is convex and
positively homogeneous with $\psi(\zero)=0$ if and only if
it is the astral support function of a convex, closed, nonempty set
$S\subseteq\extspace$.
Noting further that $\indfadub{S} = \indaS$
(by \Cref{thm:ohull-biconj} and \Cref{cor:sep-cvx-sets-conseqs}\ref{cor:sep-cvx-sets-conseqs:d}),
this implies moreover that there exists a one-to-one correspondence
between all such functions $\psi$ and the indicator functions $\indaS$
of all such sets $S$.
This bijection is defined by the dual conjugation operation
$\psi\mapsto\psistarb$ (or equivalently, the mapping
$\psi\mapsto\indaS$ where $S$ is as given in
Eq.~\ref{eq:thm:pos-homo-is-ast-closed:2}), while its inverse is given
by the primal conjugation operation, $\indaS\mapsto\indaSstar$.%
\indexg{positively homogeneous functions!support function@as support function|)}%
\indexg{support functions, astral!positively homogeneous functions as|)}

\chapter{Minimizers and their structure}
\label{sec:minimizers}

From the start, our interest in convex functions has largely been
about their minimization
and the possibility that their minimizers might be at infinity.
In this chapter, we turn our attention specifically to the structural
properties of
astral minimizers of
the extension~$\fext$ of a convex function
$f:\Rn\rightarrow\Rext$.
Because of
how astral space was constructed
and the close connection between astral points and sequences
(such as \Cref{thm:seq-rep-not-lin-ind}),
understanding
the structure of minimizers of $\ef$
also tells us much about the sequences in~$\Rn$
that minimize the original function $f$.
This is especially true when $\fext$ is continuous.

As we know,
every astral point $\xbar\in\extspace$ can be decomposed as
$\xbar=\ebar\plusl\qq$ for some icon $\ebar\in\corezn$ and
finite $\qq\in\Rn$
(\Cref{thm:icon-fin-decomp}).
In the same way, the problem of minimizing the function $\fext$ decomposes
into the separate questions of how to minimize $\fext$ over the
choice of $\ebar$, and how to minimize it over $\qq$.
We will eventually study both of these in detail.
In this chapter, we will see that if $\xbar$ minimizes $\fext$, then
its iconic part $\ebar$ must belong to a particular set called the astral
recession cone, which will be our starting point.
We study the properties of this set and the structure of its elements,
leading to a procedure that, in a sense described below, enumerates
all of the minimizers of $\fext$.

\section{Astral recession cone}
\label{sec:arescone}

The standard recession cone
is the set of vectors corresponding to directions along which a function (on $\Rn$)
is never increasing.
We begin by studying an extension of this notion to astral space,
which will be centrally important
to our understanding of minimizers and continuity.

For a function $f:\Rn\rightarrow\Rext$,
the definition
from \Cref{sec:prelim:rec-cone}
states that a vector $\vv\in\Rn$ is in $f$'s standard recession cone,
$\resc{f}$, if $f(\xx+\lambda\vv)\leq f(\xx)$ for all $\xx\in\Rn$ and
all $\lambda\in\Rpos$.
This is equivalent to the condition that $f(\yy)\leq f(\xx)$ for all
points $\yy$ on the halfline $\hfline{\xx}{\vv}$ with endpoint
$\xx$ in direction $\vv$
(as defined in Eq.~\ref{eq:std-halfline-defn}).
Thus,
\[
  \resc{f}
  =
  \bigBraces{\vv\in\Rn :\:
    \forall \xx\in\Rn,\,
    \forall \yy\in \hfline{\xx}{\vv},\,
    f(\yy)\leq f(\xx)
  }.
\]
Simplifying further, this means that
\begin{equation}  \label{eq:std-rescone-w-sup}
  \resc{f}
  =
  \bigBraces{\vv\in\Rn :\:
    \forall \xx\in\Rn,\,
    \sup f(\hfline{\xx}{\vv})\leq f(\xx)
  },
\end{equation}
where, as usual, $f(\hfline{\xx}{\vv})$ is the set of all values of $f$
along the halfline $\hfline{\xx}{\vv}$.

\indexg{recession cone, astral|(}%
These definitions can be extended to astral
space simply by replacing the standard halfline with its astral
analogue, the astral halfline, as characterized in \Cref{thm:ast-halfline:conic}:
\[
  \ahfline{\xbar}{\vbar}
   =\xbar\seqsum\aray{\vbar}.
\]
This leads to the following definition of astral recession cone:

\begin{definition}
\label{def:arescone}
\indexg{recession cone, astral!defined|(}%
Let $F:\extspace\rightarrow\Rext$.
The \emph{astral recession cone} of $F$, denoted $\aresconeF$, is the set
of points $\vbar\in\extspace$ with the property that
$F(\ybar)\leq F(\xbar)$ for all $\xbar\in\extspace$ and all
$\ybar\in\ahfline{\xbar}{\vbar}$.
That is,
\begin{align}
\indexm{rec f800}{$\aresconeF$}{astral recession cone}%
  \aresconeF
  &=
  \bigBraces{ \vbar\in\extspace :\:
           \forall \xbar\in\extspace,\,
           \forall \ybar\in\ahfline{\xbar}{\vbar},\,
           F(\ybar) \leq F(\xbar) }.
  \nonumber
  \\
\intertext{Or, equivalently,}
  \aresconeF
  &=
  \bigBraces{ \vbar\in\extspace :\:
           \forall \xbar\in\extspace,\,
           \sup F(\ahfline{\xbar}{\vbar}) \leq F(\xbar) }
  \nonumber
  \\
  &=
  \bigBraces{ \vbar\in\extspace :\:
           \forall \xbar\in\extspace,\,
           \sup F(\xbar\seqsum\aray{\vbar}) \leq F(\xbar) }.%
\indexg{recession cone, astral!defined|)}%
  \label{eqn:aresconeF-def}
\end{align}
\end{definition}
The astral recession cone can be used to derive various
properties of $F$.
The next simple proposition summarizes some of the useful ones that
follow when a point $\vbar$ is in $\aresconeF$.
Shortly, we will see that, under convexity assumptions,
some of these necessary
conditions are in fact sufficient as well.

\begin{proposition}  \label{pr:arescone-def-ez-cons}
  Let $F:\extspace\rightarrow\Rext$ and let $\vbar\in\aresconeF$.
  Then:
  \begin{letter-compact}
  \item  \label{pr:arescone-def-ez-cons:a}
    $F(\zbar\plusl\xbar)\leq F(\xbar)$
    and
    $F(\xbar\plusl\zbar)\leq F(\xbar)$
    for all $\xbar\in\extspace$ and all $\zbar\in\aray{\vbar}$.
  \item  \label{pr:arescone-def-ez-cons:b}
    $F(\alpha\vbar\plusl\xbar)\leq F(\xbar)$
    and
    $F(\xbar\plusl\alpha\vbar)\leq F(\xbar)$
    for all $\xbar\in\extspace$ and all $\alpha\in\Rextpos$.
  \end{letter-compact}
\end{proposition}

\begin{proof}
  ~

\begin{proof-parts}
\pfpart{Part~(\ref{pr:arescone-def-ez-cons:a}):}
Let $\xbar\in\extspace$ and $\zbar\in\aray{\vbar}$.
Then
$\zbar\plusl\xbar \in \zbar\seqsum\xbar \subseteq \xbar\seqsum\aray{\vbar}$
by
\Cref{cor:seqsum-conseqs}(\ref{cor:seqsum-conseqs:a}).
Similarly,
$\xbar\plusl\zbar\in\xbar\seqsum\aray{\vbar}$.
The claim then follows by definition of $\aresconeF$
(Eq.~\ref{eqn:aresconeF-def}).

\pfpart{Part~(\ref{pr:arescone-def-ez-cons:b}):}
Let $\alpha\in\Rextpos$.
Since $\alpha\vbar\in\aray{\vbar}$
(by \Cref{thm:aray,thm:mul-char}),
this follows from
part~(\ref{pr:arescone-def-ez-cons:a}).
\qedhere
\end{proof-parts}
\end{proof}

\indexg{recession cone, astral!convexity of|(}%
\indexg{recession cone, astral!astral cone@as astral cone|(}%
Just as the standard recession cone of any function on $\Rn$ is a
convex cone (\Cref{pr:resc-cone-basic-props}), so
the astral recession cone of any astral function is always a convex
astral cone:

\begin{theorem}   \label{thm:arescone-is-ast-cvx-cone}
  Let $F:\extspace\rightarrow\Rext$.
  Then
  $\aresconeF$
  is a convex astral cone.
\end{theorem}

\begin{proof}
  ~

\begin{proof-parts}
\pfpart{Astral cone:}
Let $\vbar\in\aresconeF$.
We aim to show that $\aray{\vbar}\subseteq\aresconeF$.
As such, let $\ybar\in\aray{\vbar}$.
Then, for all $\xbar\in\extspace$,
\[
  \sup F(\xbar\seqsum\aray{\ybar})
  \leq
  \sup F(\xbar\seqsum\aray{\vbar})
  \leq
  F(\xbar).
\]
The first inequality is because
$\aray{\ybar}\subseteq\aray{\vbar}$ since
$\aray{\vbar}$ is an astral cone
(\Cref{pr:astral-cone-props}\ref{pr:astral-cone-props:e}).
The second inequality is because $\vbar\in\aresconeF$.
Thus, $\ybar\in\aresconeF$, so $\aresconeF$ is an astral cone.

\pfpart{Convex:}
In light of
\Cref{thm:ast-cone-is-cvx-if-sum},
to show that the astral cone $\aresconeF$ is convex, it suffices to
show that it is closed under sequential sum.
As such,
let $\vbar_0,\vbar_1\in\aresconeF$, and let
$\vbar\in\vbar_0\seqsum\vbar_1$.
We aim to show that $\vbar\in\aresconeF$.

Let $K_i=\aray{\vbar_i}$ for $i\in\{0,1\}$.
Then $K_0$ and $K_1$ are astral cones
(by
\Cref{pr:astral-cone-props}\ref{pr:astral-cone-props:e}),
so $K_0\seqsum K_1$ is also an astral cone by
\Cref{thm:seqsum-ast-cone}(\ref{thm:seqsum-ast-cone:a}).
Further, $\vbar\in\vbar_0\seqsum\vbar_1\subseteq K_0\seqsum K_1$,
so
$\aray{\vbar}\subseteq K_0\seqsum K_1$.
Thus, for all $\xbar\in\extspace$,
\begin{align*}
   \sup F(\xbar\seqsum\aray{\vbar})
   &\leq
   \sup F\bigParens{\xbar\seqsum (\aray{\vbar_0}) \seqsum (\aray{\vbar_1})}
   \\
   &\leq
   \sup F\regParens{\xbar\seqsum \aray{\vbar_0}}
   \\
   &\leq
   F(\xbar).
\end{align*}
The first inequality is because, as just argued,
$\aray{\vbar}\subseteq (\aray{\vbar_0}) \seqsum (\aray{\vbar_1})$.
The second is because $\vbar_1\in\aresconeF$, implying
$\sup F(\ybar\seqsum \aray{\vbar_1}) \leq F(\ybar)$
for all $\ybar\in\eRn$, including all $\ybar\in\xbar\seqsum \aray{\vbar_0}$.
The third is because $\vbar_0\in\aresconeF$.
Therefore, $\vbar\in\aresconeF$, so $\aresconeF$ is convex.%
\indexg{recession cone, astral!astral cone@as astral cone|)}%
\indexg{recession cone, astral!convexity of|)}%
\qedhere
\end{proof-parts}
\end{proof}

\indexg{recession cone, astral!characterizations|(}%
We next show that
if $F:\extspace\rightarrow\Rext$ is convex, then
much simpler conditions than those in \Cref{def:arescone}
can be used to characterize points
in $\aresconeF$.
First, we show that
if we merely have
$F(\limray{\vbar}\plusl\xbar)\leq F(\xbar)$
for all $\xbar\in\extspace$, then $\vbar$ must be in $\aresconeF$,
so in this sense, this one inequality is as strong
as all the others that are implicit in \eqref{eqn:aresconeF-def}.
Furthermore, if
$F(\limray{\vbar}\plusl\xbar)<+\infty$ at even a single point
$\xbar\in\extspace$, then $\vbar$ must be in $\aresconeF$.

\begin{theorem}   \label{thm:recF-equivs}
  Let $F:\extspace\rightarrow\Rext$ be convex, and let
  $\vbar\in\extspace$.
  Then the following are equivalent:
  \begin{letter-compact}
  \item   \label{thm:recF-equivs:a}
    $\vbar\in\aresconeF$.
  \item   \label{thm:recF-equivs:b}
    For all $\xbar\in\extspace$,
    $F(\limray{\vbar}\plusl\xbar)\leq F(\xbar)$.
  \item   \label{thm:recF-equivs:c}
    Either $F\equiv+\infty$ or
    there exists $\xbar\in\extspace$ such that
    $F(\limray{\vbar}\plusl\xbar) < +\infty$.
  \end{letter-compact}
\end{theorem}

\begin{proof}
  ~

\begin{proof-parts}
\pfpart{%
  (\ref{thm:recF-equivs:a})
  $\Rightarrow$
  (\ref{thm:recF-equivs:b}):
}
This is immediate from
\Cref{pr:arescone-def-ez-cons}(\ref{pr:arescone-def-ez-cons:b}).

\pfpart{%
  (\ref{thm:recF-equivs:b})
  $\Rightarrow$
  (\ref{thm:recF-equivs:c}):
}
Suppose condition~(\ref{thm:recF-equivs:b}) holds.
If $F\not\equiv+\infty$ then there exists $\xbar\in\dom{F}$
implying
$F(\limray{\vbar}\plusl\xbar)\leq F(\xbar)<+\infty$.

\pfpart{%
  (\ref{thm:recF-equivs:c})
  $\Rightarrow$
  (\ref{thm:recF-equivs:a}):
}
If $F\equiv+\infty$ then $\aresconeF=\extspace$ which trivially
includes $\vbar$.

Otherwise,
suppose there exists $\ybar\in\extspace$ such that
$F(\limray{\vbar}\plusl\ybar)<+\infty$.
Then
\Cref{thm:F-conv-res}
immediately implies that
$\vbar\in\aresconeF$.
\qedhere
\end{proof-parts}
\end{proof}

\indexg{recession cone, astral!extension of|(}%
\indexg{lower semicontinuous extension!astral recession cone of|(}%
We will be especially interested in the astral recession cone of the
extension $\fext$ of a convex function $f:\Rn\rightarrow\Rext$.
In this case, even simpler conditions
can be used to characterize points in
$\aresconef$:

\begin{theorem}   \label{thm:rec-ext-equivs}
  Let $f:\Rn\rightarrow\Rext$ be convex, and let
  $\vbar\in\extspace$.
  Then the following are equivalent:
  \begin{letter-compact}
  \item   \label{thm:rec-ext-equivs:a}
    $\vbar\in\aresconef$.
  \item   \label{thm:rec-ext-equivs:b}
    For all $\xx\in\Rn$,
    $\fext(\vbar\plusl\xx)\leq f(\xx)$.
  \item   \label{thm:rec-ext-equivs:c}
    For all $\xx\in\Rn$,
    $\fext(\limray{\vbar}\plusl\xx)\leq f(\xx)$.
  \item   \label{thm:rec-ext-equivs:d}
    Either $f\equiv+\infty$ or
    $\fshadvb\not\equiv+\infty$
    (that is, $\fext(\limray{\vbar}\plusl\yy)<+\infty$ for some
    $\yy\in\Rn$).
  \item   \label{thm:rec-ext-equivs:e}
    Either $f\equiv+\infty$ or
    $\fext(\limray{\vbar}\plusl\ybar)<+\infty$ for some
    $\ybar\in\extspace$.
  \end{letter-compact}
 \end{theorem}

\begin{proof}
  ~

\begin{proof-parts}
\pfpart{%
  (\ref{thm:rec-ext-equivs:a})
  $\Rightarrow$
  (\ref{thm:rec-ext-equivs:b}):
}
Suppose $\vbar\in\aresconef$.
Then for all $\xx\in\Rn$,
$ \fext(\vbar\plusl\xx) \leq \fext(\xx) \leq f(\xx)$
by
Propositions~\ref{pr:arescone-def-ez-cons}(\ref{pr:arescone-def-ez-cons:b})
and~\ref{pr:h:1}(\ref{pr:h:1a}).

\pfpart{%
  (\ref{thm:rec-ext-equivs:b})
  $\Rightarrow$
  (\ref{thm:rec-ext-equivs:c}):
}
Suppose condition~(\ref{thm:rec-ext-equivs:b}) holds.
Let $\xbar\in\extspace$.
Then there exists a sequence $\seq{\xx_t}$ in $\Rn$ converging to $\xbar$
and with $f(\xx_t)\rightarrow\fext(\xbar)$
(\Cref{pr:d1}).
Thus,
\begin{equation}  \label{eq:thm:rec-ext-equivs:1}
  \fext(\xbar)
  =
  \lim f(\xx_t)
  \geq
  \liminf \fext(\vbar\plusl\xx_t)
  \geq
  \fext(\vbar\plusl\xbar).
\end{equation}
The first inequality is by our assumption,
and the second is because $\fext$ is lower semicontinuous
(\Cref{prop:ext:F}\ref{prop:ext:F:a})
and $\vbar\plusl\xx_t\rightarrow\vbar\plusl\xbar$
(\Cref{pr:i:7}\ref{pr:i:7f}).

Let $\xx\in\Rn$.
Then for all $t$,
we now have
\[
  \fext(\xx\plusl t \vbar)
  =
  \fext(t\vbar \plusl \xx)
  \leq
  \fext(\xx)
  \leq
  f(\xx),
\]
where the first inequality follows by repeated application of
\eqref{eq:thm:rec-ext-equivs:1},
and the second by
\Cref{pr:h:1}(\ref{pr:h:1a}).
Thus,
\[
  \fext(\limray{\vbar}\plusl\xx)
  =
  \fext(\xx\plusl\limray{\vbar})
  \leq
  \liminf \fext(\xx\plusl t\vbar)
  \leq
  f(\xx),
\]
where the first inequality is by lower semicontinuity of $\fext$
and since $\xx\plusl t\vbar\rightarrow \xx\plusl\limray{\vbar}$
(by \Cref{pr:i:7}\ref{pr:i:7f}).

\pfpart{%
  (\ref{thm:rec-ext-equivs:c})
  $\Rightarrow$
  (\ref{thm:rec-ext-equivs:d}):
}
Suppose condition~(\ref{thm:rec-ext-equivs:c}) holds.
If $f\not\equiv+\infty$ then there exists $\yy\in\dom{f}$ so
$\fext(\limray{\vbar}\plusl\yy)\leq f(\yy)<+\infty$.

\pfpart{%
  (\ref{thm:rec-ext-equivs:d})
  $\Rightarrow$
  (\ref{thm:rec-ext-equivs:e}):
}
This is immediate.

\pfpart{%
  (\ref{thm:rec-ext-equivs:e})
  $\Rightarrow$
  (\ref{thm:rec-ext-equivs:a}):
}
Suppose condition~(\ref{thm:rec-ext-equivs:e}) holds.
Then condition~(\ref{thm:recF-equivs:c}) of
\Cref{thm:recF-equivs} must hold as well with $F=\fext$
(since if $f\equiv+\infty$ then $\fext\equiv+\infty$).
By that theorem, it then follows that $\vbar\in\aresconef$
(noting that $\fext$ is convex by
\Cref{thm:fext-convex}).%
\indexg{recession cone, astral!characterizations|)}%
\qedhere
\end{proof-parts}
\end{proof}

Here is an immediate corollary for icons:

\begin{corollary}  \label{cor:a:4}
  Let $f:\Rn\rightarrow\Rext$ be convex and not identically $+\infty$,
  and
  let $\ebar\in\corezn$.
  Then $\ebar\in\aresconef$ if and only if
  $\fext(\ebar\plusl\qq)<+\infty$
  for some $\qq\in\Rn$.
\end{corollary}

\begin{proof}
This follows immediately from
\Crefequiv{thm:rec-ext-equivs}{thm:rec-ext-equivs:a}{thm:rec-ext-equivs:d},
noting that $\limray{\ebar}=\ebar$ by
\Cref{pr:i:8}(\ref{pr:i:8d}).
\end{proof}

The following corollary provides a simple description of
$\aresconef$ when $f\not\equiv+\infty$:

\begin{corollary}
\label{cor:rec-equiv}
Let $f:\Rn\to\eR$ be convex, and suppose $\yy\in\dom f$.
Then
\[
  \aresconef = \set{\vbar\in\eRn:\: \ef(\limray{\vbar}\plusl\yy)<+\infty}.
\]
\end{corollary}
\begin{proof}
If $\vbar\in\aresconef$, then $\ef(\limray{\vbar}\plusl\yy)\le f(\yy)<+\infty$
  by
  \Cref{thm:rec-ext-equivs}(\ref{thm:rec-ext-equivs:a},\ref{thm:rec-ext-equivs:c}).
  Conversely, if $\vbar\in\eRn$ with $\ef(\limray{\vbar}\plusl\yy)<+\infty$, then
  \Cref{thm:rec-ext-equivs}(\ref{thm:rec-ext-equivs:e},\ref{thm:rec-ext-equivs:a})
  implies that $\vbar\in\aresconef$.%
\indexg{recession cone, astral!extension of|)}%
\indexg{lower semicontinuous extension!astral recession cone of|)}%
\end{proof}

\indexg{recession cone, astral!closedness of|(}%
If $F:\extspace\rightarrow\Rext$ is convex and lower semicontinuous,
then its astral recession cone must be closed:

\begin{theorem}   \label{thm:rescone-closed}
  Let $F:\extspace\rightarrow\Rext$ be convex and lower
  semicontinuous.
  Then $\aresconeF$ is a closed (in $\extspace$) convex astral cone.
\end{theorem}

\begin{proof}
That $\aresconeF$ is a convex astral cone was shown in
\Cref{thm:arescone-is-ast-cvx-cone},
so it only remains to show that it is also closed.

If $F\equiv+\infty$, then $\aresconeF=\extspace$, which is closed.
Therefore, we assume henceforth that $F\not\equiv+\infty$.

Let $\vbar$ be any point in $\clbar{\aresconeF}$, implying there
exists a sequence $\seq{\vbar_t}$ in $\aresconeF$ with
$\vbar_t\rightarrow\vbar$.
To show that $\aresconeF$ is closed,
we aim to show that $\vbar\in\aresconeF$.

Let $\ybar$ be any point in $\dom{F}$, which exists since
$F\not\equiv+\infty$.
For all $t$,
let $\xbar_t=\limray{\vbar_t}\plusl\ybar$.
Then by sequential compactness, the sequence $\seq{\xbar_t}$ must have
a subsequence converging to some point $\xbar\in\extspace$.
By discarding all other sequence elements, we can assume that
$\xbar_t\rightarrow\xbar$.
We then have
\[
  F(\xbar)
  \leq
  \liminf F(\xbar_t)
  =
  \liminf F(\limray{\vbar_t}\plusl\ybar)
  \leq
  F(\ybar)
  <
  +\infty.
\]
The first inequality is because $F$ is lower semicontinuous, and
the second is by
\Cref{pr:arescone-def-ez-cons}(\ref{pr:arescone-def-ez-cons:b})
since $\vbar_t\in\aresconeF$.
Furthermore, by \Cref{thm:gen-dom-dir-converg},
$\xbar\in\limray{\vbar}\plusl\extspace$
since $\vbar_t\rightarrow\vbar$;
that is,
$\xbar=\limray{\vbar}\plusl\zbar$ for some $\zbar\in\extspace$.
Since $F(\limray{\vbar}\plusl\zbar)<+\infty$, it now follows
by
\Crefequiv{thm:recF-equivs}{thm:recF-equivs:c}{thm:recF-equivs:a}
that $\vbar\in\aresconeF$, completing the proof.
\end{proof}

\indexg{recession cone, astral!extension of|(}%
\indexg{lower semicontinuous extension!astral recession cone of|(}%
\Cref{thm:rescone-closed} applies specifically to the extension
of any convex function on $\Rn$:

\begin{corollary}   \label{cor:res-fbar-closed}
  Let $f:\Rn\rightarrow\Rext$ be convex.
  Then $\aresconef$ is a closed (in $\extspace$) convex astral cone.
\end{corollary}

\begin{proof}
This follows immediately from
\Cref{thm:rescone-closed}
since $\fext$ is convex
by \Cref{thm:fext-convex},
and lower semicontinuous
by \Cref{prop:ext:F}(\ref{prop:ext:F:a}).%
\indexg{recession cone, astral!extension of|)}%
\indexg{lower semicontinuous extension!astral recession cone of|)}%
\indexg{recession cone, astral!closedness of|)}%
\end{proof}

\indexg{recession cone, astral!standard recession cone and|(}%
\indexg{recession cone (standard)!astral recession cone and|(}%
We next look at the relationship between standard and astral
recession cones
when working with a convex, lower semicontinuous function
$f:\Rn\rightarrow\Rext$.
We first show that
$\fext$'s astral recession cone, when restricted to $\Rn$, is
the same as $f$'s standard recession cone.

\begin{proposition}  \label{pr:f:1}
  Let $f:\Rn\rightarrow\Rext$ be convex and lower semicontinuous.
  Then $(\aresconef)\cap\Rn = \resc{f}$.
\end{proposition}

\begin{proof}
Let $\vv\in\resc{f}\subseteq\Rn$.
Then for all $\xx\in\Rn$,
$\fext(\vv+\xx)=f(\vv+\xx)\leq f(\xx)$,
where the equality is by
\Cref{pr:h:1}(\ref{pr:h:1a}) since $f$ is lower
semicontinuous, and the inequality is because $\vv\in\resc{f}$.
Therefore, $\vv\in\aresconef$ by
\Crefequiv{thm:rec-ext-equivs}{thm:rec-ext-equivs:b}{thm:rec-ext-equivs:a}.
Thus, $\resc{f}\subseteq(\aresconef)\cap\Rn$.

For the reverse inclusion, suppose now that
$\vv\in(\aresconef)\cap\Rn$.
Then for all $\xx\in\Rn$ and for all $\lambda\in\Rpos$, we have
\[
   f(\lambda\vv+\xx)
   =
   \fext(\lambda\vv+\xx)
   \leq
   \fext(\xx)
   =
   f(\xx).
\]
The two equalities are by
\Cref{pr:h:1}(\ref{pr:h:1a}).
The inequality is by
\Cref{pr:arescone-def-ez-cons}(\ref{pr:arescone-def-ez-cons:b})
since $\vv\in\aresconef$.
Therefore, $\vv\in\resc{f}$, completing the proof.
\end{proof}

Combining \Cref{pr:f:1} with general properties of astral cones yields
the following inclusions:

\begin{proposition}  \label{pr:repres-in-arescone}
  Let $f:\Rn\rightarrow\Rext$ be convex and lower semicontinuous.
  Then $\represc{f}$ and $\rescbar{f}$
  are both convex astral cones, and
  \[
     \resc{f}
     \subseteq
     \represc{f}
     \subseteq
     \rescbar{f}
     \subseteq
     \aresconef.
  \]
\end{proposition}

\begin{proof}
  By \Cref{pr:resc-cone-basic-props},
  $\resc{f}$ is a convex cone in $\Rn$.
  Therefore, \Cref{cor:a:1} immediately yields
  that
  $\resc{f} \subseteq \represc{f} \subseteq \rescbar{f}$.
  Also, \Cref{pr:f:1} implies that
  $\resc{f} \subseteq  \aresconef$.
  Since $\aresconef$ is closed in $\extspace$
  (by \Cref{cor:res-fbar-closed}),
  this further implies that this set must also include $\resc{f}$'s closure,
  $\rescbar{f}$, proving the final inclusion.

That $\represc{f}$ and $\rescbar{f}$ are convex astral cones
follows from \Cref{cor:a:1} (and
\Cref{pr:astral-cone-props}\ref{pr:astral-cone-props:e}).
\end{proof}

\indexg{recessive completeness|(}%
When the last two inclusions in \Cref{pr:repres-in-arescone}
become equalities,
$f$ is said to have an important property called recessive
completeness,
as we now define:

\begin{definition}  \label{dfn:reces-complete}
We say that 
a convex and lower semicontinuous function $f:\Rn\rightarrow\Rext$
is \emph{recessive complete}
if
$\represc{f}=\aresconef$.
\end{definition}

Thus, if $f$ is recessive complete, then
its standard recession cone,
$\resc{f}$, completely describes all the points in the astral recession cone, $\aresconef$,
in the sense that all points in that set have representations composed
from points in $\resc{f}$ 
(and all points with such representations are in
\indexg{recession cone, astral!standard recession cone and|)}%
\indexg{recession cone (standard)!astral recession cone and|)}%
$\aresconef$).
As examples of its relevance to our later development,
we will see in \Cref{subsec:rank-one-minimizers}
that recessive completeness provides a sufficient condition for an
extension $\fext$ to have minimizers of astral rank at most one.
And in \Cref{subsec:cond-for-cont},
we will see that recessive completeness of $f$ is closely related to
$\fext$'s continuity properties.

Here is
a useful sufficient condition for recessive completeness:

\begin{theorem}   \label{pr:rec-complete}
\indexg{recessive completeness!polyhedral recession cone and|(}%
\indexg{polyhedral sets, astral!recessive completeness and|(}%
  Let $f:\Rn\rightarrow\Rext$ be convex and lower semicontinuous,
  and suppose $\aresconef$ is astral polyhedral.
  Then $f$ is recessive complete.
\end{theorem}

\begin{proof}
By \Cref{pr:f:1},
$ \resc{f} = (\aresconef)\cap\Rn $, which, by
\Cref{cor:std-ast-polyhedra}(\ref{cor:std-ast-polyhedra:Q}),
is polyhedral with 
$\rescbar{f}=\aresconef$.
Since $\resc{f}$ is polyhedral and a convex cone
(\Cref{pr:resc-cone-basic-props}), it then also follows that
$\represc{f}=\rescbar{f}$ by
\Cref{thm:repcl-polyhedral-cone}(\ref{thm:repcl-polyhedral-cone:c},\ref{thm:repcl-polyhedral-cone:b}).
Thus, $f$ is recessive complete.%
\indexg{recessive completeness|)}%
\indexg{recessive completeness!polyhedral recession cone and|)}%
\indexg{polyhedral sets, astral!recessive completeness and|)}%
\end{proof}

\indexg{recession cone, astral!minimizers characterized using|(}%
\indexg{minimizers of extensions!characterizations of|(}%
Using the astral recession cone,
we can now characterize the form of all points that minimize $\fext$.
Specifically, we show that all minimizers must have the form $\ebar\plusl\qq$ where $\ebar$ is an icon in
the astral recession cone $\aresconef$,
and $\qq\in\Rn$ minimizes the reduction
$\fshadd$
(defined in \Cref{def:iconic-reduction}).
In later sections,
especially in \Cref{sec:univ-red-and-min},
we will develop a much more detailed
analysis of the minimizers of $\fext$, but this
\namecref{thm:arescone-fshadd-min}
provides a start:

\begin{theorem}  \label{thm:arescone-fshadd-min}
  Let $f:\Rn\rightarrow\Rext$ be convex.
  Let $\xbar=\ebar\plusl\qq$ where $\ebar\in\corezn$ and $\qq\in\Rn$.
  Then $\xbar$ minimizes $\fext$ if and only if $\ebar\in\aresconef$
  and $\qq$ minimizes $\fshadd$.
\end{theorem}

\begin{proof}
If $f\equiv+\infty$ then $\fext\equiv+\infty$, $\fshadd\equiv+\infty$,
and $\aresconef=\extspace$, so the claim follows trivially.
Therefore, we assume $f\not\equiv+\infty$,
so $\min \fext = \inf f < +\infty$
(by \Cref{pr:fext-min-exists}).

Suppose $\xbar$ minimizes $\fext$.
Then $\fext(\ebar\plusl\qq)<+\infty$, so $\ebar\in\aresconef$
by
\Cref{cor:a:4}.
If, contrary to the claim, $\qq$ does not minimize $\fshadd$,
then there exists $\qq'\in\Rn$ with
\[
  \fext(\ebar\plusl\qq')
  =
  \fshadd(\qq')
  <
  \fshadd(\qq)
  =
  \fext(\xbar),
\]
contradicting that $\xbar$ minimizes $\fext$.

Conversely, suppose $\ebar\in\aresconef$, and that $\qq$ minimizes
$\fshadd$.
Let $\beta\in\R$ be such that $\beta > \inf f$,
and let $\yy\in\Rn$ be such that $f(\yy) < \beta$.
Then
\[
   \inf f
   \leq
   \fext(\xbar)
   =
   \fshadd(\qq)
   \leq
   \fshadd(\yy)
   =
   \fext(\ebar\plusl\yy)
   \leq
   f(\yy)
   <
   \beta.
\]
The first inequality is by \Cref{pr:fext-min-exists},
the second is because $\qq$ minimizes $\fshadd$,
and the third is by
\Crefequiv{thm:rec-ext-equivs}{thm:rec-ext-equivs:a}{thm:rec-ext-equivs:b}.
Since this holds for all $\beta>\inf f$,
$\fext(\xbar)=\inf f$, so $\xbar$ minimizes $\fext$.%
\indexg{recession cone, astral!minimizers characterized using|)}%
\indexg{minimizers of extensions!characterizations of|)}%
\end{proof}

\indexg{recession cone, astral!examples|(}%
\indexg{recessive completeness!examples|(}%
We end this section with three examples of functions and their
astral recession cones as an illustration of some of the notions
presented above.

\begin{example}[Astral linear function]
\label{ex:astral-linear:rec}
Let $\uu\in\Rn$, and define
$f:\Rn\rightarrow\R$ by $f(\xx)=\xx\cdot\uu$ for $\xx\in\Rn$,
whose extension is the astral linear function $\fext(\xbar)=\xbar\cdot\uu$
for $\xbar\in\extspace$
(\Cref{ex:ext-affine}).
Then by \Cref{cor:rec-equiv} (setting $\yy=\zero$),
\[
  \aresconef
   =
  \{ \vbar\in\extspace :\: \limray{\vbar}\cdot\uu < +\infty \}
  =
  \{ \vbar\in\extspace :\: \vbar\cdot\uu \leq 0 \}.
\]
Intersecting with $\Rn$, \Cref{pr:f:1} then yields
\[
  \resc{f}
  =
  \{ \vv\in\Rn :\: \vv\cdot\uu \leq 0 \}.
\]
Further, $\aresconef$ is a closed astral halfspace (or all of
$\extspace$ if $\uu=\zero$), and so is astral polyhedral;
hence, $f$ is recessive complete by \Cref{pr:rec-complete}.
\end{example}

\begin{example}[Extension of restricted linear function]
\label{ex:negx1-else-inf-cont}
\indexg{Restricted linear function!recession cone of|(}%
Consider the restricted linear function
from \Cref{ex:negx1-else-inf}, defined, for $\xx\in\R^2$, as
\[
  f(\xx)=
  f(x_1,x_2) =
  \begin{cases}
    -x_1
    & \text{if $x_2\geq 0$,}
  \\
    +\infty
    & \text{otherwise.}
  \end{cases}
\]
Its extension, for $\xbar\in\extspac{2}$, can be derived by
combining
Example~\ref{ex:ext-affine} with
Corollaries~\ref{cor:ext-restricted-fcn}
and~\ref{cor:halfspace-closure}:
\begin{equation}   \label{eq:ex:negx1-else-inf-cont:2}
  \fext(\xbar)
  =
  \begin{cases}
  \xbar\cdot(-\ee_1)
     & \text{if $\xbar\cdot\ee_2\geq 0$,}
   \\
     +\infty
     & \text{otherwise.}
   \end{cases}
\end{equation}
By \Cref{cor:rec-equiv} (setting $\yy=\zero$),
\begin{align}
\notag
  \aresconef
  &=
  \bigBraces{\vbar\in\extspac{2} :\:
     \limray{\vbar}\cdot\ee_2\geq 0
     \text{ and }
     \limray{\vbar}\cdot(-\ee_1)<+\infty
  }
\\
\label{eq:ex:negx1-else-inf-cont:1}
  &=
  \bigBraces{ \vbar\in\extspac{2} :\:
    \vbar\cdot\ee_2 \geq 0
    \text{ and }
    \vbar\cdot\ee_1 \geq 0
  }.
\end{align}
\Cref{pr:f:1} then yields
\[
  \resc{f}
  =
  \bigBraces{ \vv\in\R^2 :\:
    \vv\cdot\ee_2 \geq 0
    \text{ and }
    \vv\cdot\ee_1 \geq 0
  }.
\]
Thus, $\aresconef$ is an intersection of two closed astral
halfspaces, and so is astral polyhedral.
Therefore, $f$ is recessive complete by
\Cref{pr:rec-complete}.

The extension $\fext$ is minimized, for instance, by
$\limray{\ee_2}\plusl\limray{\ee_1}$
and
$\limray{(\ee_1+\ee_2)}\plusl\ee_1$, since $\fext$ is $-\infty$ on
each of these points,
as can be seen from
\eqref{eq:ex:negx1-else-inf-cont:2}.
Consistent with \Cref{thm:arescone-fshadd-min},
the iconic part of each of these is included in $\aresconef$,
and the finite part minimizes the corresponding reduction
(since
$\fshad{\limray{\ee_2}\plusl\limray{\ee_1}}=\fshad{\limray{(\ee_1+\ee_2)}}\equiv-\infty$,
as follows from Eq.~\ref{eq:ex:negx1-else-inf-cont:2}).
Also consistent with that \namecref{thm:arescone-fshadd-min},
$\limray{\ee_2}\plusl\ee_1$ does not minimize~$\fext$, even though its
iconic part, $\limray{\ee_2}$, is included in $\aresconef$,
since $\ee_1$ does not minimize~$\fshad{\limray{\ee_2}}$,
given that $\fshad{\limray{\ee_2}}(\xx)=\xx\inprod(-\ee_1)$ for $\xx\in\R^2$
(again
from Eq.~\ref{eq:ex:negx1-else-inf-cont:2}).%
\indexg{Restricted linear function!recession cone of|)}%
\end{example}

\begin{figure}
  \centering
  \includegraphics{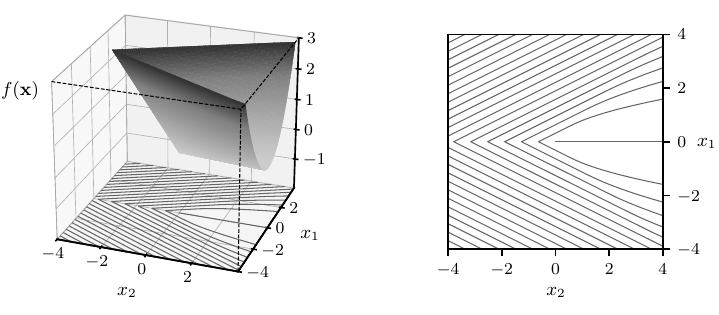}
  \mycaption{Flattening valley}{%
\indexf{Flattening valley}%
    The function $f$ from \Cref{ex:x1sq-over-x2}.%
  }%
  \label{fig:flat-valley}%
\end{figure}

\begin{example}[Flattening valley]
\label{ex:x1sq-over-x2}
\indexg{Flattening valley|(}%
Let $f:\R^2\rightarrow\R$ be defined,
for $\xx\in\R^2$, by
\begin{equation}  \label{eqn:curve-discont-finiteev-eg}
  f(\xx)
  =
  f(x_1,x_2)
  =
  \begin{cases}
    x_1^2/x_2
      & \text{if $x_2 >\regAbs{x_1}$,}
    \\
    2\regAbs{x_1}-x_2
      & \text{otherwise.}
  \end{cases}
\end{equation}
(See \Cref{fig:flat-valley}.)
This function is convex, closed, proper, finite everywhere, and
continuous everywhere.
It is also nonnegative
\indexg{Flattening valley|)}%
everywhere.
\indexg{Flattening valley!recession cone of|(}%
It can be checked that
$f$'s standard recession cone
is the closed ray
$\resc{f}=\{ \lambda \ee_2 :\: \lambda\in\Rpos \}$.
(This can be seen by noting, for fixed $x_1\in\R$, that
$f(x_1,x_2)$ is nonincreasing as a function of $x_2$.
On the other hand, if $\vv$ is not a nonnegative multiple of $\ee_2$,
then $f(\vv)>0=f(\zero)$, so $\vv\not\in\resc{f}$.)
In this case, $\resc{f}$'s representational and topological closures
are the same, namely,
\begin{equation}   \label{eq:ex:x1sq-over-x2:1}
   \represc{f}
   =
   \rescbar{f}
   =
   (\resc{f})\cup\{\limray{\ee_2}\}
   =
   \regBraces{ \alpha \ee_2 :\: \alpha\in\Rextpos },
\end{equation}
which,
by \Cref{pr:repres-in-arescone},
is included in $\aresconef$.

It can further be checked that
$\fshad{\limray{\ee_2}}(\xx)=\fext(\limray{\ee_2}\plusl\xx)=0$
for all $\xx\in\R^2$.
Let $g=\fshad{\limray{\ee_2}}$, so $\eg\equiv 0$.
Then, for all $\ybar\in\extspac{2}$ and all $\xx\in\R^2$, by \Cref{thm:d4}(\ref{thm:d4:c}),
we have
$\ef(\limray{\ee_2}\plusl\ybar\plusl\xx)=\eg(\ybar\plusl\xx)=0 \leq f(\xx)$.
This means, by
\Crefequiv{thm:rec-ext-equivs}{thm:rec-ext-equivs:b}{thm:rec-ext-equivs:a},
that any astral point of the form
$\limray{\ee_2}\plusl\ybar$ is in $\aresconef$.
Moreover, it can be argued that $\aresconef$ cannot include any
additional points without violating
\Cref{pr:f:1}.
Thus, in summary,
\[
  \aresconef = \braces{ \lambda \ee_2 :\: \lambda \in\Rpos }
        \cup
        \bracks{\limray{\ee_2} \plusl \extspac{2}}.
\]
Comparing with \eqref{eq:ex:x1sq-over-x2:1},
$\represc{f}\ne\aresconef$, so $f$ is not recessive complete.%
\indexg{Flattening valley!recession cone of|)}%
\indexg{recessive completeness!examples|)}%
\indexg{recession cone, astral!examples|)}%
\indexg{recession cone, astral|)}%
\end{example}

\section{Recursive characterization and generation of minimizers}
\label{sec:astral-cone}

From \Cref{thm:arescone-fshadd-min},
for a convex function $f:\Rn\rightarrow\Rext$,
we know that the iconic part of every
minimizer of $\fext$
must be an element of $\ef$'s astral recession cone.
Thus, to minimize $\fext$ (as well as $f$), it will be helpful to
understand the structure of points in the astral recession cone,
and how such points can be constructed.
In this section, we use methods developed earlier based on projections and reductions
to obtain a recursive characterization of the astral recession cone,
and we then use that
characterization to algorthmically describe the set of minimizers of $\fext$.

\indexg{recession cone, astral!recursive formulation of|(}%
\indexg{reductions, astronic!astral recession cone of|(}%
Let $\vv\in\Rn$, and let
$g=\fshadv$
be the reduction of $f$ at $\limray{\vv}$.
We begin by showing how points in $\aresconef$,
the astral recession cone of $\fext$,
relate to points in $\aresconeg$,
the astral recession cone of $\gext$.

\begin{theorem}   \label{thm:f:3}
  Let $f:\Rn\rightarrow\Rext$ be convex and lower semicontinuous.
  Let $\vv\in\Rn$
  and let
  $g=\fshadv$
  be the reduction of $f$ at $\limray{\vv}$.
  Then:
  \begin{letter-compact}
  \item   \label{thm:f:3a}
    $\aresconef \subseteq \aresconeg$.
    (Consequently, $\resc{f} \subseteq \resc{g}$.)
  \item   \label{thm:f:3b}
    Suppose $\ybar=\limray{\vv}\plusl\zbar$ for some $\zbar\in\extspace$.
    Then $\ybar\in\aresconef$ if and only if
    $\vv\in\resc{f}$ and
    $\zbar\in\aresconeg$.
  \end{letter-compact}
\end{theorem}

Part~(\ref{thm:f:3b}) of
this theorem provides a recursive characterization of all
the points comprising $\aresconef$:
The finite points in $\aresconef$
are exactly those in the standard recession cone $\resc{f}$,
by \Cref{pr:f:1}.
All of the other points in $\aresconef$ can be enumerated
by considering each point
$\vv\in\resc{f}$, forming the reduction $g=\fshadv$,
finding $\gext$'s astral recession cone $\aresconeg$, and then
leftwardly
adding $\limray{\vv}$ to each element in $\aresconeg$.
Thus,
\[
   \aresconef = (\resc{f}) \cup
               \bigcup_{\vv\in\resc{f}}
                 \Parens{\limray{\vv}\plusl \arescone{\fshadext{\limray{\vv}}}}.
\]

Alternatively, we can think of part~(\ref{thm:f:3b}), together with
\Cref{pr:f:1}, as providing a test for determining if a
given point $\ybar$ is in $\aresconef$:
If $\ybar$ is in $\Rn$, then it is in $\aresconef$ if and only if it
is in $\resc{f}$.
Otherwise, it is in $\aresconef$ if and only if
its dominant direction $\vv$ is in $\resc{f}$ and its projection
orthogonal to $\vv$
is in $\aresconeg$, as can be determined in a recursive
manner.

\begin{proof}[Proof of \Cref{thm:f:3}]
If $f\equiv +\infty$, then $g\equiv+\infty$, so
$\fext=\gext\equiv +\infty$,
$\aresconef=\aresconeg=\extspace$,
and $\resc{f}=\resc{g}=\Rn$,
implying all the parts of the theorem.
We therefore assume henceforth that $f\not\equiv+\infty$.

\begin{proof-parts}
\pfpart{Part~(\ref{thm:f:3a}):}
If $\vv\not\in\resc{f}$, then $\gext\equiv +\infty$ by
\Cref{thm:i:4}, again implying $\aresconeg=\extspace$ and
trivially yielding the claim.

So suppose $\vv\in\resc{f}$ and $f\not\equiv+\infty$.
Let $\ybar\in\aresconef$.
Then by
\Crefequiv{thm:rec-ext-equivs}{thm:rec-ext-equivs:a}{thm:rec-ext-equivs:d},
there exists $\rr\in\Rn$ for which
$\fext(\limray{\ybar} \plusl \rr)<+\infty$.
This implies
\[
  \gext(\limray{\ybar} \plusl \rr)
  \leq
  \fext(\limray{\ybar} \plusl \rr)
  < +\infty,
\]
where the first inequality is by
\Cref{cor:i:1}(\ref{cor:i:1c}).
Therefore, $\ybar\in\aresconeg$, again by
\Crefequiv{thm:rec-ext-equivs}{thm:rec-ext-equivs:d}{thm:rec-ext-equivs:a}.

Having proved $\aresconef \subseteq \aresconeg$,
it now follows, by \Cref{pr:f:1}, that
\[
  \resc{f}
  =
  (\aresconef)\cap\Rn
  \subseteq
  (\aresconeg)\cap\Rn
  =
  \resc{g}.
\]

\pfpart{Part~(\ref{thm:f:3b}):}
Suppose first that $\vv\in\resc{f}$ and that $\zbar\in\aresconeg$.
Then for all $\xx\in\Rn$,
\[
  \fext(\ybar\plusl \xx)
  =
  \fext(\limray{\vv}\plusl\zbar\plusl \xx)
  =
  \gext(\zbar\plusl \xx)
  \leq
  g(\xx)
  \leq
  f(\xx).
\]
The second equality is by \Cref{cor:i:1}(\ref{cor:i:1b}).
The first inequality is by
\Crefequiv{thm:rec-ext-equivs}{thm:rec-ext-equivs:a}{thm:rec-ext-equivs:b},
since $\zbar\in\aresconeg$.
The second inequality is
by \Cref{pr:d2}(\ref{pr:d2:b})
(since $\vv\in\resc{f}$).
Therefore, $\ybar\in\aresconef$,
again by
\Crefequiv{thm:rec-ext-equivs}{thm:rec-ext-equivs:b}{thm:rec-ext-equivs:a}.

For the converse, suppose for the rest of the proof that
$\ybar\in\aresconef$.
We argue separately that $\vv\in\resc{f}$ and $\zbar\in\aresconeg$.

Let $\qq\in\dom{f}$, which exists since $f\not\equiv+\infty$.
Then
\[
  \fext(\limray{\vv}\plusl\zbar\plusl\qq)
  =
  \fext(\ybar\plusl\qq)
  \leq
  f(\qq)
  <+\infty,
\]
where the first inequality is by
\Crefequiv{thm:rec-ext-equivs}{thm:rec-ext-equivs:a}{thm:rec-ext-equivs:b},
since $\ybar\in\aresconef$.
Therefore, $\vv\in(\aresconef)\cap\Rn = \resc{f}$ by
\Crefequiv{thm:rec-ext-equivs}{thm:rec-ext-equivs:e}{thm:rec-ext-equivs:a}
and \Cref{pr:f:1}.

Finally, we have
\[
  \gext(\limray{\zbar}\plusl\qq)
  =
  \fext(\limray{\vv}\plusl\limray{\zbar}\plusl\qq)
  =
  \fext(\limray{\ybar}\plusl\qq)
  \leq
  f(\qq)
  <
  +\infty.
\]
The first equality is by \Cref{cor:i:1}(\ref{cor:i:1b}).
The second equality is because
$\limray{\ybar}=\limray{\vv}\plusl\limray{\zbar}$
by \Cref{pr:i:8}(\ref{pr:i:8d},\ref{pr:i:8-infprod}).
The first inequality is by
\Crefequiv{thm:rec-ext-equivs}{thm:rec-ext-equivs:a}{thm:rec-ext-equivs:c},
since $\ybar\in\aresconef$.
Thus, $\zbar\in\aresconeg$ by
\Crefequiv{thm:rec-ext-equivs}{thm:rec-ext-equivs:d}{thm:rec-ext-equivs:a}.%
\indexg{recession cone, astral!recursive formulation of|)}%
\indexg{reductions, astronic!astral recession cone of|)}%
\qedhere
\end{proof-parts}
\end{proof}

\indexg{recession cone, astral!characterizations|(}%
Applying \Cref{thm:f:3}(\ref{thm:f:3b}) repeatedly,
we obtain the following characterization of points in the
astral recession cone, based on their representation:

\begin{theorem}  \label{thm:astral-cone-char}
  Let $f:\Rn\rightarrow\Rext$ be convex.
  Let $\xbar = \limrays{\vv_1,\ldots,\vv_k}\plusl\qq$ where
  $\vv_1,\ldots,\vv_k,\qq\in\Rn$.
  Let
  $g_0=\lsc f$ and $g_i = \gishadvi$ for $i=1,\ldots,k$.
  Then the following are equivalent:
  \begin{letter-compact}
  \item  \label{thm:astral-cone-char:a}
    $\xbar\in\aresconef$.
  \item  \label{thm:astral-cone-char:b}
    $\vv_i\in\resc{g_{i-1}}$ for $i=1,\ldots,k$,
    and $\qq\in\resc{g_k}$.
  \end{letter-compact}
\end{theorem}

\begin{proof}
As preliminary steps, we note by
\Cref{pr:icon-red-decomp-astron-red} that each
$g_i$ is convex and lower semicontinuous.
Also, for $i=0,\ldots,k$,
let $\xbar_i=\limrays{\vv_{i+1},\ldots,\vv_k}\plusl\qq$.
Then
for $i=1,\ldots,k$,
since $\xbar_{i-1}=\limray{\vv_i}\plusl\xbar_i$,
and by $g_i\negKern$'s definition,
\Cref{thm:f:3}(\ref{thm:f:3b}) implies that
$\xbar_{i-1}\in\aresconegsub{{i-1}}$
if and only if
$\vv_i\in\resc{g_{i-1}}$ and
$\xbar_{i}\in\aresconegsub{{i}}$.

Note also that $\fext=\gext_0$ by
\Cref{pr:h:1}(\ref{pr:h:1aa}).

\begin{proof-parts}
\pfpart{%
  (\ref{thm:astral-cone-char:a})
  $\Rightarrow$
  (\ref{thm:astral-cone-char:b}):
}
Suppose $\xbar\in\aresconef$, so that also
$\xbar_0=\xbar\in\aresconef=\aresconegsub{0}$.
From the preceding remarks,
it then follows by a straightforward induction that
$\vv_i\in\resc{g_{i-1}}$ and
$\xbar_{i}\in\aresconegsub{{i}}$
for $i=1,\ldots,k$. In particular, $\qq=\xbar_k\in\aresconegsub{k}$.
By assumption, $\qq\in\Rn$,
so $\qq\in(\aresconegsub{k})\cap\Rn=\resc{g_k}$ (by \Cref{pr:f:1}).

\pfpart{%
  (\ref{thm:astral-cone-char:b})
  $\Rightarrow$
  (\ref{thm:astral-cone-char:a}):
}
Suppose
$\vv_i\in\resc{g_{i-1}}$
for $i=1,\ldots,k$,
and $\qq\in\resc{g_k}$.
Then by backward induction,
we show that
$\xbar_{i}\in\aresconegsub{{i}}$
for $i=0,\ldots,k$.
The base case, when $i=k$, holds because
$\xbar_k=\qq\in\resc{g_k}\subseteq\aresconegsub{k}$
(by \Cref{pr:f:1}).
For the inductive step, when $0\leq i<k$, we have
$\xbar_{i+1}\in\aresconegsub{{i+1}}$ by inductive hypothesis, and
$\vv_{i+1}\in\resc{g_i}$ by assumption, so
$\xbar_{i}\in\aresconegsub{{i}}$
by \Cref{thm:f:3}(\ref{thm:f:3b}).
Hence, $\xbar=\xbar_0\in\aresconegsub{0}=\aresconef$.%
\indexg{recession cone, astral!characterizations|)}%
\qedhere
\end{proof-parts}
\end{proof}

\begin{algorithm}[t]
  \caption[Test whether a given point $\xbar$ minimizes $\fext$]%
          {Test whether a given point $\xbar$ minimizes $\fext$.}
    \label{fig:min-test-proc}
\begin{algorithmic}
    \Block{Input}
        \State
\indexalg{minimizers of extensions!testing for}%
            function $f:\Rn\rightarrow\Rext$ that is convex and lower
            semicontinuous,
        \State
            test point $\xbar = \limrays{\vv_1,\ldots,\vv_k} \plusl \qq$
            where $\vv_1,\ldots,\vv_k,\qq\in\Rn$
    \EndBlock
    \smallskip
    \Block{Define \Call{TestIfMinimizer}{$f,\xbar$}}
        \If{$k=0$}
            \State if $\qq$ minimizes $f$ then return \true
            \State else return \false
        \Else
            \If{$\vv_1\not\in \resc{f}$}
                \State return \false
            \Else
                \State let $g = \fshad{\limray{\vv_1}}$ and
                       $\zbar = \limrays{\vv_2,\ldots,\vv_k} \plusl \qq$
                \State return \Call{TestIfMinimizer}{$g,\zbar$}
            \EndIf
        \EndIf
    \EndBlock
    \smallskip
    \Block{Output properties}
        \State
            \true if $\xbar$ minimizes $\fext$;
            \false otherwise
    \EndBlock
\end{algorithmic}
\end{algorithm}

\indexg{minimizers of extensions!testing for|(}%
Combining the characterization of minimizers (\Cref{thm:arescone-fshadd-min}) with
the recursive characterization of the astral recession cone (\Cref{thm:f:3}\ref{thm:f:3b}),
we can now provide an algorithm
for testing if a given astral point minimizes $\fext$
(see \Cref{fig:min-test-proc}).
The input is a convex and lower semicontinuous
function $f:\Rn\rightarrow\Rext$ and an explicitly represented
test point $\xbar = \limrays{\vv_1,\ldots,\vv_k} \plusl \qq$,
where $\vv_1,\ldots,\vv_k,\qq\in\Rn$.
The algorithm determines if $\xbar$ minimizes $\fext$ using only
more basic primitives which operate on standard points and
functions over $\Rn$, specifically,
for testing if a point in $\Rn$ minimizes an ordinary
convex function on $\Rn$,
and also for testing if a
vector in $\Rn$ is in the standard recession cone of such a function.

The operation and correctness of this algorithm follow directly from
our development regarding minimizers and reductions:
If $k=0$, then $\xbar=\qq\in\Rn$, so $\xbar$ minimizes $\fext$ if and
only if $\qq$ minimizes the standard function $f$
(which we have assumed is lower semicontinuous).
Otherwise, if $k>0$, then
$\xbar=\limray{\vv_1} \plusl \zbar$
where
$\zbar = \limrays{\vv_2,\ldots,\vv_k} \plusl \qq$.
If $\vv_1\not\in\resc{f}$, then
$\limrays{\vv_1,\ldots,\vv_k}$ cannot be in $\aresconef$,
by \Cref{thm:f:3}(\ref{thm:f:3b}), and therefore
$\xbar$ cannot minimize $\fext$, by
\Cref{thm:arescone-fshadd-min}.
Otherwise, with $g = \fshad{\limray{\vv_1}}$,
if $\vv_1\in\resc{f}$, then
$\gext(\zbar)=\fext(\limray{\vv_1}\plusl\zbar)=\fext(\xbar)$
and
$\min \gext = \inf g = \inf f = \min \fext$
by \Cref{thm:d4}
and \Cref{pr:fext-min-exists}.
Therefore, $\xbar$ minimizes $\fext$
if and only if $\zbar$ minimizes~$\gext$.%
\indexg{minimizers of extensions!testing for|)}

\begin{algorithm}[t]

  \caption[Nondeterministically generate all minimizers of $\fext$]%
          {Nondeterministically generate all minimizers of $\fext$.}
\label{fig:all-min-proc}

\begin{algorithmic}
    \Block{Input}
        \State
\indexalg{minimizers of extensions!generation of all}%
            convex function $f:\Rn\rightarrow\Rext$
    \EndBlock
    \smallskip
    \Block{Define \Call{GenerateMinimizer}{$f$}}
        \State $i \leftarrow 0$
        \State $g_0 = \lsc f$
        \Loop{repeat \emph{at least} until $g_i$ has a finite minimizer}
              \Comment{nondeterministic choice}
            \State $i \leftarrow i+1$
            \State let $\vv_i$ be \emph{any} point in $\resc{g_{i-1}}$
              \Comment{nondeterministic choice}
            \State $g_i = \gishadvi$
        \EndLoop
        \State $k \leftarrow i$
        \State $\ebar=\limrays{\vv_1,\ldots,\vv_k}$
        \State let $\qq\in\Rn$ be \emph{any} finite minimizer of $g_k$
          \Comment{nondeterministic choice}
        \State $\xbar=\ebar\plusl\qq$
        \State return $(g_k,\ebar,\qq,\xbar)$
    \EndBlock
    \smallskip
    \Block{Output properties}
        \State
            $g_k = \fshadd$
        \State
            $\ebar\in (\aresconef)\cap\corezn$
        \State
            $\qq$ minimizes $\fshadd$
        \State
            $\xbar$ minimizes $\fext$
    \EndBlock
\end{algorithmic}
\end{algorithm}

\indexg{minimizers of extensions!generation of all|(}%
We next describe how to
algorithmically generate all minimizers
of $\fext$.
As a starting point, we use \Cref{thm:arescone-fshadd-min},
which states that
every minimizer $\xbar=\ebar\plusl\qq$ has an iconic part
$\ebar$ that is in the astral recession cone $\aresconef$
and a finite part $\qq$ that minimizes~$\fshadd$.
To generate an icon in $\aresconef$,
we use \Cref{thm:astral-cone-char} (applied with $\xbar=\ebar$ and $\qq=\zero$)
as follows.
We start with $g_0=\lsc f$,
and on each iteration~$i$, we generate a vector~$\vv_i$ in $\resc{g_{i-1}}$, the
standard recession cone of $g_{i-1}$.
We then define
$g_i = \gishadvi$
to be the next reduction.
In this way, we ensure (by \Cref{thm:astral-cone-char})
that the resulting icon
$\ebar=\limrays{\vv_1,\ldots,\vv_k}$
formed by the $\vv_i\negKern$'s must be in $\aresconef$.

We can continue this process until we manage to form a reduction
$g_k$ that has some finite minimizer $\qq\in\Rn$.
By such a construction, $g_k$ actually is equal to the reduction
$\fshadd$ at icon $\ebar$, so in fact, $\qq$ minimizes $\fshadd$,
which, combined with $\ebar$ being in $\aresconef$, ensures that the
point $\xbar=\ebar\plusl\qq$ minimizes $\fext$.

We summarize this process in \Cref{fig:all-min-proc}.
Although we describe the process in the form of an algorithm, we do
not literally mean to suggest that it be implemented on a computer, at
least not in this generality.
The point, rather, is to reveal the structure of the minimizers of a
function in astral space, and how that structure can be related to
standard notions from convex analysis.

\begin{example}[Two exps and a square]
   \label{ex:simple-eg-exp-exp-sq}
\indexg{Two exps and a square|(}%
Let $f:\R^3\rightarrow\Rext$ be the convex function defined, for
$\xx\in\R^3$, by
\begin{equation}  \label{eq:simple-eg-exp-exp-sq}
   f(\xx) = f(x_1,x_2,x_3) = \me^{x_3-x_1} + \me^{-x_2}
               + (2+x_2-x_3)^2.
\end{equation}
Then $f$'s extension, $\fext$, can be shown to be continuous
everywhere, and is specifically
\begin{equation}  \label{eq:simple-eg-exp-exp-sq:fext}
  \fext(\xbar) =
      \expex\bigParens{\xbar\cdot (\ee_3-\ee_1)}
      +
      \expex\bigParens{\xbar\cdot (-\ee_2)}
      +
      \bigParens{2+\xbar\cdot (\ee_2-\ee_3)}^2
\end{equation}
for $\xbar\in\extspac{3}$,
where $\expex$ is as given in
\eqref{eq:expex-defn},
$(\pm\infty)^2=+\infty$ by standard arithmetic over
$\Rext$,
and $\ee_1,\ee_2,\ee_3$ are the standard basis vectors in $\R^3$.

By \Cref{cor:rec-equiv} (setting $\yy=\zero$),
the astral recession cone of $\ef$ consists of
those points $\zbar\in\extspac{3}$ for which
      $\limray{\zbar}\cdot (\ee_3-\ee_1)<+\infty$,\;
      $\limray{\zbar}\cdot (-\ee_2)<+\infty$,
      and
      $\regParens{2+\limray{\zbar}\cdot (\ee_2-\ee_3)}^2<+\infty$.
Thus,
\begin{equation}
  \label{eq:simple-eg-exp-exp-sq:arescone}
      \resc{\ef}
    =
    \BigBraces{
      \zbar\in\R^3 :\:
      \zbar\cdot (\ee_3-\ee_1)\le 0,\;
      \zbar\cdot (-\ee_2)\le 0,\;
      \zbar\cdot (\ee_2-\ee_3) = 0
    }.
\end{equation}
Intersecting with $\Rn$, \Cref{pr:f:1} then yields
\begin{equation}
\label{eq:simple-eg-exp-exp-sq:rescone}
  \resc{f}
    =
    \BigBraces{
      \zz\in\R^3 :\:
      0\leq z_2=z_3\leq z_1
    }.%
\indexg{Two exps and a square|)}%
\end{equation}

\indexg{Two exps and a square!algorithm applied to|(}%
Suppose we apply
\Cref{fig:all-min-proc} to $f$.
On the first iteration, the algorithm chooses any vector $\vv_1$ in the
standard recession cone of $g_0=\lsc f=f$, say
$\vv_1=\trans{[1,1,1]}$.
Next, the reduction $g_1=\genshad{g_0}{\limray{\vv_1}}$
is formed, which is
\[
  g_1(\xx)
  =
  \fext(\limray{\vv_1}\plusl \xx)
  = \me^{x_3-x_1} + (2+x_2-x_3)^2
\]
for $\xx\in\R^3$.
Its recession cone is
\[
   \resc{g_1}
            = \regBraces{ \zz\in\R^3 :\:
                            z_2=z_3\leq z_1 },
\]
so on the next iteration, we can choose any $\vv_2$ in this set,
say
$\vv_2=\trans{[1,-1,-1]}$.
The next reduction $g_2=\genshad{g_1}{\limray{\vv_2}}$ is
$g_2(\xx)=(2+x_2-x_3)^2$.
This function has finite minimizers, such as
$\qq=\trans{[0,0,2]}$.
The resulting minimizer of $\fext$ is
$\xbar=\ebar\plusl\qq$
where,
as claimed in \Cref{fig:all-min-proc},
$\ebar=\limray{\vv_1}\plusl\limray{\vv_2}$
is an icon in $\aresconef$ and $\qq$ minimizes $g_2=\fshadd$.%
\indexg{Two exps and a square!algorithm applied to|)}%
\end{example}

Returning to our general discussion,
\Cref{fig:all-min-proc}
is
\indexg{nondeterminism}%
\emph{nondeterministic} in the sense
that at various points, choices are made in a way that is entirely
arbitrary.
This happens at three different points:
First, on each iteration of the main loop, an \emph{arbitrary} point
$\vv_i$ is selected from $\resc{g_{i-1}}$.
Second, this loop must iterate \emph{at least} until $g_i$ has a finite
minimizer, but can continue to iterate arbitrarily beyond that point.
Third, after terminating the loop, an \emph{arbitrary} finite
minimizer $\qq$ of $g_k$ is selected.

Clearly, the point $\xbar$ that is eventually computed by the algorithm
depends on these arbitrary choices.
Nevertheless, we show that in all cases, the resulting
point $\xbar$ must be a minimizer of~$\fext$.
We also show that if $\xbar$ minimizes $\fext$, then it must be possible for
these arbitrary choices to be made in such a way that $\xbar$ is
produced
(while still respecting the constraints imposed at each step
of the algorithm).
It is in this sense that the algorithm generates \emph{all} of the
minimizers of $\fext$.

When there exists such a sequence of choices that results in $\xbar$
as the final output of the computation, we say that $\xbar$
is a
\indexg{potential output (of algorithm)}%
\emph{potential output} of the algorithm.
Thus, we are claiming that a point $\xbar\in\extspace$ minimizes
$\fext$ if and only if $\xbar$ is a potential output of the algorithm.
This is shown formally by the next theorem, whose
condition~(\ref{thm:all-min-proc-correct:b})
captures exactly when the point $\xbar$ is a potential output.

\begin{theorem}  \label{thm:all-min-proc-correct}
\indexg{minimizers of extensions!characterizations of|(}%
  Let $f:\Rn\rightarrow\Rext$ be convex.
  Let $\xbar = \limrays{\vv_1,\ldots,\vv_k}\plusl\qq$ where
  $\vv_1,\ldots,\vv_k,\qq\in\Rn$,
  and
  let $\ebar=\limrays{\vv_1,\ldots,\vv_k}$.
  Let
  $g_0=\lsc f$, and $g_i = \gishadvi$ for $i=1,\ldots,k$.
  Then $g_k=\fshadd$, and the following are equivalent:
  \begin{letter-compact}
  \item  \label{thm:all-min-proc-correct:a}
    $\xbar$ minimizes $\fext$.
  \item  \label{thm:all-min-proc-correct:b}
    $\qq$ minimizes $g_k$,
    and
    $\vv_i\in\resc{g_{i-1}}$ for $i=1,\ldots,k$.
  \end{letter-compact}
\end{theorem}

\begin{proof}
By \Cref{pr:icon-red-decomp-astron-red}, $g_k=\fshadd$.
  
  \begin{proof-parts}
    \pfpart{%
      (\ref{thm:all-min-proc-correct:a})
      $\Rightarrow$
      (\ref{thm:all-min-proc-correct:b}):
    }
    Suppose $\xbar$ minimizes $\fext$.
    Then \Cref{thm:arescone-fshadd-min} implies that
    $\qq$ minimizes $\fshadd=g_k$, and also that
    $\ebar\in\aresconef$.
    By \Cref{thm:astral-cone-char}, applied to $\ebar$,
    we then obtain that $\vv_i\in\resc{g_{i-1}}$
    for $i=1,\ldots,k$.
    
    \pfpart{%
      (\ref{thm:all-min-proc-correct:b})
      $\Rightarrow$
      (\ref{thm:all-min-proc-correct:a}):
    }
    Suppose
    $\qq$ minimizes $g_k=\fshadd$
    and that
    $\vv_i\in\resc{g_{i-1}}$
    for $i=1,\ldots,k$.
    Then, since $\zero\in\resc{g_k}$,
    we obtain by \Cref{thm:astral-cone-char}
    that
    $\ebar\in\aresconef$. Thus, by
    \Cref{thm:arescone-fshadd-min}, $\xbar=\ebar\plusl\qq$ minimizes $\fext$.%
\indexg{minimizers of extensions!characterizations of|)}%
    \qedhere
    \end{proof-parts}
    \end{proof}

\Cref{thm:all-min-proc-correct} shows that if \Cref{fig:all-min-proc}
terminates, then the computed point $\xbar$ must minimize $\fext$.
But what if the algorithm never terminates?
Indeed, it is possible for the algorithm to never terminate,
or even to
reach the loop's (optional) termination condition.
For instance,
the same point $\vv_i=\zero$ (which is in every recession cone)
might be chosen on every iteration so
that the algorithm never makes any progress toward a solution at all.
In general, if a vector $\vv_i$ is chosen that is already in the span
of the preceding vectors $\vv_1,\ldots,\vv_{i-1}$, then no progress
is made in the sense that the icon that is being constructed has not
changed; that is,
\[
  \limrays{\vv_1,\ldots,\vv_i} = \limrays{\vv_1,\ldots,\vv_{i-1}}
\]
(by the Push Lemma~\ref{pr:g:1}).
Thus, to ensure progress, we might insist that $\vv_i$ be chosen
to be not only in the recession cone of $g_{i-1}$ but also
outside the span of $\vv_1,\ldots,\vv_{i-1}$ so that
\begin{equation}  \label{eq:sensible-vi-dfn}
   \vv_i \in (\resc{g_{i-1}}) \setminus
             \spnfin{\vv_1,\ldots,\vv_{i-1}}.
\end{equation}
We say such a choice of $\vv_i$ is \emph{nonredundant}.
In contrast, a choice $\vv_i\in(\resc{g_{i-1}})\cap\spnfin{\vv_1,\ldots,\vv_{i-1}}$
is called \emph{redundant}.

We next prove a sufficient condition for a function
to have a finite minimizer, and show that if all available
choices for $\vv_i$ are redundant then $g_i$ must have a finite minimizer
and thus the termination condition of \Cref{fig:all-min-proc} must have
already been reached.

Our sufficient condition refers to the constancy space of a function.
As discussed in \Cref{sec:prelim:rec-cone},
the constancy space of a function $f$, denoted $\conssp{f}$ and
defined in \eqref{eq:conssp-defn}, is the set of directions in
which $f$ remains constant;
it is always a linear subspace of $\Rn$,
and a subset of $\resc{f}$ (\Cref{pr:prelim:const-props}).
Our result states that if $\resc{f}=\conssp{f}$, then $f$ must have a
finite minimizer. Informally, this is because,
when this condition holds, any candidate direction in which
$f$ could possibly have a minimum at infinity is actually a direction in which $f$
is constant, 
implying that the function can be minimized without going to infinity.
This is made precise in the proof below.

Note that this proposition is entirely about notions from standard
convex analysis; nonetheless, the proof that we give is
based very much on concepts from astral space.

\begin{proposition}
\label{pr:cons-implies-fin-minimizer}
  Let $f:\Rn\rightarrow\Rext$ be convex and lower semicontinuous,
  and assume $\resc{f}=\conssp{f}$. Then $f$ must have a finite minimizer.
\end{proposition}

\begin{proof}
  If $f\equiv+\infty$, then any $\xx\in\Rn$ is a minimizer. 
  We therefore assume henceforth that $f\not\equiv+\infty$.
  Assume for the sake of contradiction that $f$ has no finite minimizer,
  implying that $\ef$ also has no finite minimizer.
(This is because any finite minimizer of $\ef$ would also minimize $f$
since $\resfcn{\ef}{\Rn}=f$ by lower semicontinuity of $f$ and
\Cref{pr:h:1}\ref{pr:h:1a}).
  
By \Cref{pr:fext-min-exists}, $\ef$ attains its minimum.
  Let $\zbar\in\eRn$ be a minimizer of $\ef$ of minimum
  astral rank $r$. Since $\zbar$ is infinite, $r\ge 1$. Let $\vv$
  be the dominant direction of $\zbar$, and let $\PP$ be the projection matrix orthogonal to $\vv$,
  so $\zbar=\limray{\vv}\plusl\PP\zbar$
  by \Cref{pr:h:6}.
  Since $\ef\not\equiv+\infty$, we have
  $\ef(\limray{\vv}\plusl\PP\zbar)=\ef(\zbar)=\min\ef<+\infty$. Thus, $\vv\in(\aresconef)\cap\Rn=\resc{f}=\conssp{f}$
  (with inclusion by
  \Cref{thm:rec-ext-equivs}\ref{thm:rec-ext-equivs:e}\ref{thm:rec-ext-equivs:a}
  and first equality by \Cref{pr:f:1}).
  Since $\vv\in\conssp{f}$, we have $f(\xx)=f(\PP\xx)$ for all
  $\xx\in\Rn$ (by \Cref{pr:cons-PP}), 
  implying $\ef(\zbar)=\ef(\PP\zbar)$ by \Cref{pr:PP:ef}.
  Hence, $\PP\zbar$ also minimizes $\ef$.
  However, this contradicts that $\zbar$ is a minimizer of minimum
  astral rank since the astral rank of $\PP\zbar$ is $r-1$
  (by \Cref{pr:h:6}).
\end{proof}

We next show that if there are no nonredundant choices for $\vv_i$,
we must have $\resc{g_i}=\conssp{g_i}$, and thus $g_i$ must have a finite minimizer.

\begin{proposition}  \label{pr:no-sense-then-done}
  Let $f:\Rn\rightarrow\Rext$ be convex.
  Let $g_0=\lsc f$, let $\vv_i\in\resc{g_{i-1}}$, and let
  $g_i = \gishadvi$,
  for $i=1,\ldots,k$.
  Suppose
  \[
      \resc{g_k} \subseteq \spnfin{\vv_1,\ldots,\vv_k}.
  \]
  Then $\resc{g_k}=\conssp{g_k}$. Consequently, $g_k$ has a finite minimizer.
\end{proposition}

\begin{proof}
  Let $\ebar=\limrays{\vv_1,\ldots,\vv_k}$. Then $g_k=\fshad{\ebar}$
(by \Cref{pr:icon-red-decomp-astron-red}), so
\[
    \conssp{g_k}
    \subseteq \resc{g_k}
    \subseteq \spnfin{\vv_1,\ldots,\vv_k}
    = \rspanset{\ebar}
    \subseteq \conssp{g_k}.
\]
  The first inclusion is
  by \Cref{pr:prelim:const-props}(\ref{pr:prelim:const-props:a}),
  and the second is by assumption.
  The equality is by definition of representational span
  (\Cref{def:rspan-gen}),
  and the last inclusion follows from
  \Cref{pr:i:9}(\ref{pr:i:9cons}).
  Thus, $\conssp{g_k}=\resc{g_k}$, and hence $g_k$ has a finite minimizer by
  \Cref{pr:cons-implies-fin-minimizer}.
\end{proof}

If, in \Cref{fig:all-min-proc},
each $\vv_i$ is chosen nonredundantly, then the dimension of the space
spanned by the $\vv_i\negKern$'s increases by one on each iteration of the
main loop.
Therefore, after at most $n$ iterations, all available choices must be
redundant, and therefore, by \Cref{pr:no-sense-then-done},
the termination condition must have been reached.
This shows that the algorithm can always be run in a way that guarantees
termination within $n$ iterations.
Furthermore, it shows that the algorithm cannot ``get stuck'' in the
sense that, no matter what preceding choices have been made,
the ensuing choices of $\vv_i$ can be made nonredundantly,
again ensuring termination within $n$ additional iterations.%
\indexg{minimizers of extensions!generation of all|)}

\section{Dual characterization}
\label{sec:astral-cone-dual}

We will often find it useful to rely on another fundamental
characterization of
the astral recession cone in terms of the function's dual
properties.
From \Cref{pr:rescpol-is-con-dom-fstar}, we know that
the standard recession cone of a
closed, proper, convex function $f$ is the
(standard) polar of $\cone(\dom{\fstar})$.
\indexg{recession cone, astral!dual characterization|(}%
We show that an analogous identity, expressed using dual polarity,
holds for astral functions and their
astral recession cones:

\begin{theorem}  \label{thm:arescone-is-conedomFstar-pol}
  Let $F:\extspace\rightarrow\Rext$ with $F\not\equiv+\infty$.
  Then the following hold:
  \begin{letter-compact}
  \item  \label{thm:arescone-is-conedomFstar-pol:a}
    $\aresconeF\subseteq\apolconedomFstar$.
  \item  \label{thm:arescone-is-conedomFstar-pol:b}
    If, in addition, $F=\Fdub$, then
    $\aresconeF=\apolconedomFstar$.
  \end{letter-compact}
\end{theorem}

\begin{proof}

~

\begin{proof-parts}
\pfpart{Part~(\ref{thm:arescone-is-conedomFstar-pol:a}):}
Let $\vbar\in\aresconeF$ and $\uu'\in\cone(\dom{\Fstar})$. We need
to show that $\vbar\cdot\uu'\leq 0$.
If $\uu'=\zero$ then this holds.
Otherwise, we must have $\uu'=\lambda\uu$ for some $\uu\in\dom{\Fstar}$
and some $\lambda\in\Rstrictpos$, by \Cref{pr:scc-cone-elts}(\ref{pr:scc-cone-elts:c:new}),
noting that $\dom{\Fstar}$ is convex, being the domain of the convex function $\Fstar$
(\Cref{pr:conj-is-convex}).

Let $\xbar$ be any point in $\dom{F}$, which exists since
$F\not\equiv+\infty$. Then:
\begin{align}
\notag
   +\infty
   >
   F(\xbar)
   &\ge
   F(\limray{\vbar}\plusl\xbar)
\\
\label{eq:dual:arescone:1}
   &\ge
   \Fdub(\limray{\vbar}\plusl\xbar)
   \ge
   -\Fstar(\uu) \plusd (\limray{\vbar}\cdot\uu\plusl\xbar\cdot\uu).
\end{align}
The first inequality is because $\xbar\in\dom{F}$, the second is
because $\vbar\in\aresconeF$
(and \Cref{pr:arescone-def-ez-cons}\ref{pr:arescone-def-ez-cons:b}),
the third is by
\Cref{thm:fdub-sup-afffcns}(\ref{thm:fdub-sup-afffcns:c}),
and the last is by definition of dual
conjugate (Eq.~\ref{eq:psistar-def:2}).

Note that $-\Fstar(\uu)>-\infty$ since $\uu\in\dom{\Fstar}$. Thus,
\eqref{eq:dual:arescone:1} implies that
$\limray{\vbar}\cdot\uu\plusl\xbar\cdot\uu<+\infty$,
which in turn implies $\limray{\vbar}\cdot\uu<+\infty$, and thus
$\vbar\cdot\uu\leq 0$. Therefore, $\vbar\cdot\uu'=\lambda(\vbar\cdot\uu)\leq 0$
(using \Cref{pr:scalar-prod-props}\ref{pr:scalar-prod-props:a}).
This shows that $\vbar\cdot\uu'\le 0$ for all
$\uu'\in\cone(\dom{\Fstar})$,
so $\vbar\in\apolconedomFstar$.

\pfpart{Part~(\ref{thm:arescone-is-conedomFstar-pol:b}):}
Suppose $F=\Fdub$, which is convex by \Cref{thm:dual-conj-cvx}.
Let $\vbar\in\apolconedomFstar$; we aim to show
$\vbar\in\aresconeF$.
Combined with part~(\ref{thm:arescone-is-conedomFstar-pol:a}), this
will prove the claim.

Let $\xbar\in\extspace$ and $\uu\in\Rn$.
We claim that
\begin{equation}  \label{eq:thm:arescone-is-conedomFstar-pol:2}
  -\Fstar(\uu) \plusd (\limray{\vbar}\cdot\uu\plusl\xbar\cdot\uu)
  \leq
  -\Fstar(\uu) \plusd \xbar\cdot\uu.
\end{equation}
If $\Fstar(\uu)=+\infty$, then this is immediate since both sides of
the inequality are $-\infty$.
Otherwise, $\uu\in\dom{\Fstar}\subseteq\cone(\dom{\Fstar})$ so
$\vbar\cdot\uu\leq 0$ (since $\vbar\in\apolconedomFstar$),
implying
$\limray{\vbar\cdot\uu}\leq 0$; hence,
\eqref{eq:thm:arescone-is-conedomFstar-pol:2}
holds in this case as well
(by Propositions~\ref{pr:i:5}\ref{pr:i:5ord} and~\ref{pr:plusd-props}\ref{pr:plusd-props:f}).
Thus, for all $\xbar\in\extspace$,
\begin{align*}
  F(\limray{\vbar}\plusl\xbar)
  =
  \Fdub(\limray{\vbar}\plusl\xbar)
  &=
  \sup_{\uu\in\Rn}
      \Bracks{ -\Fstar(\uu) \plusd (\limray{\vbar}\cdot\uu\plusl\xbar\cdot\uu) }
  \\
  &\leq
  \sup_{\uu\in\Rn}
      \Bracks{ -\Fstar(\uu) \plusd \xbar\cdot\uu }
  =
  \Fdub(\xbar)
  =
  F(\xbar).
\end{align*}
The first and last equalities are by our assumption that $F=\Fdub$.
The second and third equalities are
by definition of dual conjugate
(Eq.~\ref{eq:psistar-def:2}).
The inequality is by
\eqref{eq:thm:arescone-is-conedomFstar-pol:2}.
From
\Crefequiv{thm:recF-equivs}{thm:recF-equivs:b}{thm:recF-equivs:a},
it follows that $\vbar\in\aresconeF$, completing the proof.
\qedhere
\end{proof-parts}
\end{proof}

As corollary, we immediately obtain:

\begin{corollary}  \label{cor:res-pol-char-red-clsd}
  Let $f:\Rn\rightarrow\Rext$ be convex and reduction-closed,
  with $f\not\equiv+\infty$.
  Then $\aresconef=\apolconedomfstar$.
\end{corollary}

\begin{proof}
Since $f$ is convex and reduction-closed, we have
$\fext=\fdub$ by
\Cref{thm:dub-conj-new}.
Noting also that $\fextstar=\fstar$
(\Cref{pr:fextstar-is-fstar}),
the claim then follows directly from
\Cref{thm:arescone-is-conedomFstar-pol}(\ref{thm:arescone-is-conedomFstar-pol:b})
applied with $F=\fext$.
\end{proof}

If $f$ is not reduction-closed,
then the sets $\aresconef$ and $\apolconedomfstar$ need not be equal.
Here is an example:

\begin{example}
\indexg{Restricted linear function!astral recession cone and|(}%
Consider the restricted linear 
function $f$ from \Cref{ex:negx1-else-inf,ex:negx1-else-inf-cont}, and its conjugate,
\begin{equation}  \label{eq:restr:conj}
\fstar(u_1,u_2) =
\begin{cases}
          0   & \mbox{if $u_1=-1$ and $u_2\leq 0$,} \\
       +\infty & \mbox{otherwise.}
\end{cases}
\end{equation}
Let
$\vbar=\limray{\ee_1}\plusl\limray{(-\ee_2)}$.
Then $\vbar\in\apolconedomfstar$ since,
from \eqref{eq:restr:conj},
if $\uu\in\dom{\fstar}$ then $\ee_1\cdot\uu=-1$
implying $\vbar\cdot\uu=-\infty$.
But $\vbar$ is not in $\fext$'s astral recession cone,
$\aresconef=
  \set{ \vbar\in\extspac{2} :\:
    \vbar\cdot\ee_2 \geq 0
    \text{ and }
    \vbar\cdot\ee_1 \geq 0
  }
$,
derived in \eqref{eq:ex:negx1-else-inf-cont:1}
of \Cref{ex:negx1-else-inf-cont},
since $\vbar\cdot\ee_2=-\infty$.%
\indexg{Restricted linear function!astral recession cone and|)}%
\end{example}

\indexg{barrier cone (of function)!astral recession cone and|(}%
Nonetheless, we can generalize
\Cref{cor:res-pol-char-red-clsd} so that it
holds for all convex functions not identically $+\infty$,
even if they are not reduction-closed,
using the technique developed in \Cref{sec:exp-comp},
with $\cone(\dom{\fstar})$ replaced by $\slopes{f}$:

\begin{corollary}  \label{cor:ares-is-apolslopes}
  Let $f:\Rn\rightarrow\Rext$ be convex with $f\not\equiv+\infty$.
  Then $\aresconef=\apolslopesf$
  and $\slopes{f}=\rescpol{\ef}$.
\end{corollary}

\begin{proof}
Let $f'=\expex\circ f$.
Then
\[
   \aresconef
   =
   \aresconefp
   =
   \apolconedomfpstar
   =
   \apolslopesf.
\]
The first equality is by
\Cref{pr:j:2}(\ref{pr:j:2d}).
The second is by \Cref{cor:res-pol-char-red-clsd} applied to $f'$
(which is convex and lower-bounded
by \Cref{pr:j:2}\ref{pr:j:2a},
and therefore is reduction-closed by
\Cref{pr:j:1}\ref{pr:j:1c}).
The third is
by \Cref{thm:dom-fpstar}. This proves the first claim.

For the second claim, taking polars, we obtain
$\rescpol{\ef}=\apolpol{(\slopes{f})}=\slopes{f}$,
with the second equality by 
\Cref{thm:dub-ast-polar}(\ref{thm:dub-ast-polar:b})
since $\slopes{f}$ is a convex cone (by \Cref{thm:slopes-equiv}).%
\indexg{recession cone, astral!dual characterization|)}%
\indexg{barrier cone (of function)!astral recession cone and|)}%
\end{proof}

\indexg{recessive completeness!barrier cone and|(}%
\indexg{barrier cone (of function)!recessive completeness and|(}%
We next use \Cref{cor:ares-is-apolslopes}
to develop a dual characterization of recessive completeness.
Specifically,
we prove that a convex, lower semicontinuous function
$f:\Rn\rightarrow\Rext$ is recessive complete
if and only if its barrier cone, $\slopes{f}$, is polyhedral.
Thus, recessive completeness is fully characterized by
a simple geometric property of the barrier cone. %
We further prove that $f$ is recessive complete
if and only if $\resc{\ef}$ is astral polyhedral,
which shows that the sufficient condition given 
in \Cref{pr:rec-complete} for $f$ to be recessive complete
is actually necessary as well.

\indexg{recessive completeness!polyhedral recession cone and|(}%
\indexg{polyhedral sets, astral!recessive completeness and|(}%
\indexg{recession cone (standard)!recessive completeness when polyhedral|(}%
By definition of recessive completeness
(\Cref{dfn:reces-complete}) and \Cref{pr:repres-in-arescone},
for $f$ to be recessive complete, it must satisfy
$\represc{f}=\rescbar{f}=\aresconef$. As a first step, we characterize
when each of these two equalities holds:

\begin{theorem}  \label{thm:dual-char-specs}
  Let $f:\Rn\rightarrow\Rext$ be convex and lower semicontinuous.
  Then:
  \begin{letter-compact}
  \item  \label{thm:dual-char-specs:a}
    $\represc{f}=\rescbar{f}$ if and only if $\resc{f}$ is polyhedral.
  \item  \label{thm:dual-char-specs:b}
    $\rescbar{f}=\aresconef$ if and only if $\slopes{f}$ is closed
    in $\Rn$.
  \end{letter-compact}
\end{theorem}

\begin{proof}
  ~

\begin{proof-parts}
\pfpart{Part~(\ref{thm:dual-char-specs:a}):}
This is immediate from
\Cref{thm:repcl-polyhedral-cone}(\ref{thm:repcl-polyhedral-cone:b},\ref{thm:repcl-polyhedral-cone:c})
applied to $\resc{f}$.

\pfpart{Part~(\ref{thm:dual-char-specs:b}):}
If $f\equiv+\infty$ then $\slopes{f}=\Rn$ which is closed in $\Rn$,
and also
$\rescbar{f}=\aresconef=\extspace$, proving the claim in this case.

Otherwise, $f\not\equiv+\infty$.
Then by
\Cref{thm:apol-of-closure}(\ref{thm:apol-of-closure:b})
applied to $\slopes{f}$ (which is a convex cone
by \Cref{thm:slopes-equiv}),
$\slopes{f}$ is closed in $\Rn$ if and only if
$  \clbar{\polar{(\slopes{f})}}  =  \apol{(\slopes{f})}$.
Since also
$\polar{(\slopes{f})}=\resc{f}$ by \Cref{thm:rescpol-is-slopes},
and
$\apol{(\slopes{f})}=\aresconef$ by \Cref{cor:ares-is-apolslopes},
this proves the claim.
\qedhere
\end{proof-parts}
\end{proof}

\indexg{recession cone, astral!recessive completeness when polyhedral|(}%
Combining the two parts of \Cref{thm:dual-char-specs},
we can now prove the characterizations of recessive completeness
described above.

\begin{theorem}   \label{thm:dual-cond-char}
  Let $f:\Rn\rightarrow\Rext$ be convex and lower semicontinuous.
  Then the following are equivalent:
  \begin{letter-compact}
  \item   \label{thm:dual-cond-char:a}
    $f$ is recessive complete.
  \item   \label{thm:dual-cond-char:b}
    $\resc{f}$ is polyhedral, and also $\slopes{f}$ is closed
    in $\Rn$.
  \item   \label{thm:dual-cond-char:c}
    $\slopes{f}$ is polyhedral.
  \item   \label{thm:dual-cond-char:d}
    $\aresconef$ is astral polyhedral.
  \end{letter-compact}
\end{theorem}

\begin{proof}
If $f\equiv+\infty$, then $\resc{f}=\slopes{f}=\Rn$,
which is closed (in $\Rn$) and polyhedral,
and also
$\represc{f}=\aresconef=\extspace$,
which is astral polyhedral.
Together, these facts prove the equivalence of all parts in this case.
We therefore assume henceforth that $f\not\equiv+\infty$.

\begin{proof-parts}
\pfpart{%
  (\ref{thm:dual-cond-char:a})
  $\Rightarrow$
  (\ref{thm:dual-cond-char:b}):
}
If $\represc{f}=\aresconef$, then
$\represc{f}=\rescbar{f}=\aresconef$
by \Cref{pr:repres-in-arescone}, implying
the conditions in~(\ref{thm:dual-cond-char:b})
by \Cref{thm:dual-char-specs}.

\pfpart{%
  (\ref{thm:dual-cond-char:b})
  $\Rightarrow$
  (\ref{thm:dual-cond-char:c}):
}
Suppose $\resc{f}$ is polyhedral and that $\slopes{f}$ is
closed in $\Rn$.
Then $\rescpol{f}$
is also polyhedral
(by \Cref{roc:cor19.2.2}).
Furthermore,
$
  \rescpol{f} = \cl(\slopes{f}) = \slopes{f}
$
by \Cref{thm:rescpol-is-slopes},
and since $\slopes{f}$ is closed in $\Rn$.

\pfpart{%
  (\ref{thm:dual-cond-char:c})
  $\Rightarrow$
  (\ref{thm:dual-cond-char:d}):
}
Suppose $\slopes{f}$ is polyhedral.
Then
$\slopes{f}=\cone\set{\uu_1,\dotsc,\uu_m}$ for some
$\uu_1,\dotsc,\uu_m\in\Rn$
(by
\Cref{pr:fin-gen-cvx-cone}\ref{pr:fin-gen-cvx-cone:a}\ref{pr:fin-gen-cvx-cone:b}).
Thus,
\[
  \aresconef=\apolslopesf
  =\bigBraces{ \vbar\in\Rn :\: \vbar\cdot\uu_i \leq 0 \text{ for }i=1,\dotsc,m},
\]
with the first equality by \Cref{cor:ares-is-apolslopes},
and the second by
\Cref{pr:ast-pol-props}(\ref{pr:ast-pol-props:coneSpol}).
The expression on the right represents an intersection of closed
astral halfspaces (disregarding any $i$ for which $\uu_i=\zero$, since
for those the inequality $\vbar\cdot\uu_i\leq 0$ holds always).
Thus, $\aresconef$ is astral polyhedral.

\pfpart{%
  (\ref{thm:dual-cond-char:d})
  $\Rightarrow$
  (\ref{thm:dual-cond-char:a}):
}
Immediate by \Cref{pr:rec-complete}.%
\indexg{recession cone, astral!recessive completeness when polyhedral|)}%
\qedhere
\end{proof-parts}
\end{proof}

The condition that $\resc{f}$ is polyhedral
(which implies that its polar
$\rescpol{f}=\cl(\slopes{f})$ must also be polyhedral)
is not
in itself sufficient for $f$ to be recessive complete,
as shown by the next example:

\begin{example}
\label{ex:flat-valley-rec-poly-not-enough}
\indexg{Flattening valley!recessive completeness and|(}%
From
\Cref{ex:x1sq-over-x2},
we know that the recession cone
$\resc{f}$
of the flattening valley function
is the cone generated by
the singleton $\{\ee_2\}$, which is polyhedral.
Nevertheless,
as we argued in \Cref{ex:x1sq-over-x2},
$f$ is not recessive complete.
Indeed,
it can be calculated that the effective domain
of $\fstar$ is
\[
 \dom{\fstar} =
   \Braces{\uu \in \R^2 :\:
            -1 \leq u_2 \leq -\frac{u_1^2}{4} %
          },
\]
(and actually $\fstar$ is the indicator function for this set).
It can then be further calculated that
$\slopes{f}$ (which is the same as
$\cone(\dom{\fstar})$ in this case, by
\Cref{cor:ent-clos-is-slopes-cone}\ref{cor:ent-clos-is-slopes-cone:a}\ref{cor:ent-clos-is-slopes-cone:c}
since this function is finite everywhere and therefore reduction-closed
by \Cref{pr:j:1}\ref{pr:j:1d})
is equal to the origin adjoined to the open lower
half-plane:
\[
  \slopes{f} =
  \cone(\dom{\fstar}) = \{\zero\}
           \cup \{ \uu\in\R^2 :\: u_2 < 0 \}.
\]
Consistent with \Cref{thm:dual-cond-char}, this set is not
polyhedral or closed, even though its closure is polyhedral.%
\indexg{Flattening valley!recessive completeness and|)}%
\end{example}

In \Cref{thm:dual-cond-char},
if $f$ is additionally assumed to be reduction-closed,
then $\slopes{f}=\cone(\dom{\fstar})$ by
\Cref{thm:entclosed-implies-slopescone},
implying, in this
case, that we can replace $\slopes{f}$ by
$\cone(\dom{\fstar})$
in parts~(\ref{thm:dual-cond-char:b})
and~(\ref{thm:dual-cond-char:c})
of \Cref{thm:dual-cond-char} to obtain slightly simplified
necessary and sufficient conditions for recessive completeness.%
\indexg{recessive completeness!barrier cone and|)}%
\indexg{barrier cone (of function)!recessive completeness and|)}

Nonetheless, in general, without the additional assumption of
reduction closedness, \Cref{thm:dual-cond-char} does not hold
if $\slopes{f}$ is replaced by $\cone(\dom{\fstar})$, as shown in the
next example:

\begin{example}  \label{ex:negx1-else-inf-cont2}
\indexg{Restricted linear function!recessive completeness and|(}%
In \Cref{ex:negx1-else-inf-cont},
we saw that the restricted linear function $f$
is recessive complete, and that $\resc{f}=\Rpos^2$.
Nevertheless, $\cone(\dom{\fstar})$,
which we derived in 
\eqref{eqn:ex:1:cone-dom-fstar} of~\Cref{ex:negx1-else-inf},
is not closed in $\R^2$, and therefore not polyhedral.
(However, its closure is polyhedral, and furthermore,
$\slopes{f}=-\Rpos^2$, which is polyhedral, consistent
with \Cref{thm:dual-cond-char}.)
So $\cone(\dom{\fstar})$ being polyhedral is not necessary for
$f$ to be recessive complete.
Nor is it necessary that $\resc{f}$ be polyhedral and
$\cone(\dom{\fstar})$ be closed.%
\indexg{Restricted linear function!recessive completeness and|)}%
\end{example}

The example shows that conditions~(\ref{thm:dual-cond-char:b})
and~(\ref{thm:dual-cond-char:c}) 
of \Cref{thm:dual-cond-char}, with $\slopes{f}$ replaced by
$\cone(\dom{\fstar})$, are not necessary for recessive completeness.
On the other hand, if $f$ is convex and closed, then they are
sufficient, as shown next.
(Note that the function considered in
the last example
is convex, closed, and proper.)

\begin{theorem}   \label{thm:suff-conedomf-rec-comp}
  Let $f:\Rn\rightarrow\Rext$ be convex and closed.
  Suppose either of the following hold:
  \begin{letter-compact}
  \item   \label{thm:suff-conedomf-rec-comp:b}
    $\resc{f}$ is polyhedral and also $\cone(\dom{\fstar})$ is closed
    in $\Rn$.
  \item   \label{thm:suff-conedomf-rec-comp:c}
    $\cone(\dom{\fstar})$ is polyhedral.
  \end{letter-compact}
  Then $f$ is recessive complete.
\end{theorem}

\begin{proof}
If $f\equiv+\infty$, then $\resc{f}=\Rn$ and
$\represc{f}=\extspace=\aresconef$, proving the claim in this case.
We therefore assume henceforth that $f\not\equiv+\infty$.

\begin{proof-parts}
\pfpart{Part~(\ref{thm:suff-conedomf-rec-comp:b}):}
Suppose $\resc{f}$ is polyhedral and that $\cone(\dom{\fstar})$ is closed.
Then
\[
  \rescpol{f}
  =
  \cl\bigParens{\cone(\dom{\fstar})}
  =
  \cone(\dom{\fstar})
  \subseteq
  \slopes{f}
  \subseteq
  \cl(\slopes{f})
  =
  \rescpol{f},
\]
where the first equality is
by \Cref{pr:rescpol-is-con-dom-fstar},
the first inclusion is by
\Cref{thm:slopes-equiv},
and the last equality is by
\Cref{thm:rescpol-is-slopes}.
Thus,
$\slopes{f}=\rescpol{f}$,
which is polyhedral by \Cref{roc:cor19.2.2}, implying
$f$ is recessive complete
by
\Cref{thm:dual-cond-char}(\ref{thm:dual-cond-char:c},\ref{thm:dual-cond-char:a}).

\pfpart{Part~(\ref{thm:suff-conedomf-rec-comp:c}):}
Suppose $\cone(\dom{\fstar})$ is polyhedral, and therefore closed
(\Cref{roc:thm19.1}\ref{roc:thm19.1:a}\ref{roc:thm19.1:c}).
Then
$\resc{f}=\polar{(\cone(\dom{\fstar}))}$ by
\Cref{pr:rescpol-is-con-dom-fstar}, so $\resc{f}$ is polyhedral by
\Cref{roc:cor19.2.2}.
The claim now follows by part~(\ref{thm:suff-conedomf-rec-comp:b}).%
\indexg{recessive completeness!polyhedral recession cone and|)}%
\indexg{polyhedral sets, astral!recessive completeness and|)}%
\indexg{recession cone (standard)!recessive completeness when polyhedral|)}%
\qedhere
\end{proof-parts}
\end{proof}

\section{Reductions revisited}

We look next at how some of our earlier results regarding reductions
can be generalized using the astral recession cone.
Recall that, for a function $f:\Rn\rightarrow\Rext$ and
an icon $\ebar\in\corezn$, we defined the reduction
$\fshadd(\xx)=\fext(\ebar\plusl\xx)$ for $\xx\in\Rn$.
\indexg{reductions, iconic!astral analogue|(}%
In a similar way, for an astral function
$F:\extspace\rightarrow\Rext$, we can consider the function
$\xbar\mapsto F(\ebar\plusl\xbar)$, capturing the behavior of $F$ on
$\galcld$, the closure of $\ebar$'s galaxy.
The next theorem provides a useful representation for this function
akin to that used in \Cref{pr:gtil-is-PPfP-gen},
provided $\ebar$ is in $F$'s astral recession cone:

\begin{theorem}   \label{thm:F-reduc-PFP}
  Let $F:\extspace\rightarrow\Rext$,
  let $\VV\in\Rnk$,
  and let $\ebar=\VV\omm$.
  Let $G:\extspace\rightarrow\Rext$ be defined, for
  $\xbar\in\extspace$, by
  \begin{equation}  \label{eq:thm:F-reduc-PFP:0}
    G(\xbar) = F(\ebar\plusl\xbar),
  \end{equation}
  and let $\PP$ be the projection matrix onto
  $(\colspace{\VV})^\perp$.
  \begin{letter-compact}
  \item    \label{thm:F-reduc-PFP:a}
    If $\ebar\in\aresconeF$ then $G = (\PP F) \PP$.
  \item    \label{thm:F-reduc-PFP:b}
    If $\ebar\not\in\aresconeF$ and $F$ is convex then $G\equiv+\infty$.
  \end{letter-compact}
\end{theorem}

\begin{proof}
  ~

\begin{proof-parts}
\pfpart{Part~(\ref{thm:F-reduc-PFP:a}):}
Suppose $\ebar\in\aresconeF$.
Let $\xbar\in\extspace$, and let $H=(\PP F)\PP$.
We aim to show that $G(\xbar)=H(\xbar)$.
We have
\begin{equation}   \label{eq:thm:F-reduc-PFP:1}
  H(\xbar)
  =
  [(\PP F) \PP](\xbar)
  =
  (\PP F)(\PP\xbar)
  =
  \inf\regBraces{
    F(\zbar) :\:
    \zbar\in\extspace,\, \PP\zbar = \PP\xbar
  },
\end{equation}
with the second and third equalities following from definitions
(Eqs.~\ref{eq:FA-dfn} and~\ref{eqn:image-F-dfn}).
Note that
$\PP(\ebar\plusl\xbar)=\PP\VV\omm\plusl\PP\xbar=\PP\xbar$
since $\PP\VV=\zeromat{n}{k}$
(\Cref{pr:proj-mat-props}\ref{pr:proj-mat-props:e}).
Thus, \eqref{eq:thm:F-reduc-PFP:1}
implies that $H(\xbar)\leq F(\ebar\plusl\xbar)=G(\xbar)$.

For the reverse inequality, suppose $\zbar\in\extspace$ is such that
$\PP\zbar=\PP\xbar$.
Then
\[
   \ebar\plusl\zbar
   =
   \ebar\plusl\PP\zbar
   =
   \ebar\plusl\PP\xbar
   =
   \ebar\plusl\xbar,
\]
where the first and third equalities are by the Projection Lemma~\ref{lemma:proj}.
Thus,
\[
   G(\xbar)
   =
   F(\ebar\plusl\xbar)
   =
   F(\ebar\plusl\zbar)
   \leq
   F(\zbar),
\]
where the inequality is because $\ebar\in\aresconeF$
(\Cref{pr:arescone-def-ez-cons}\ref{pr:arescone-def-ez-cons:b}).
Since this holds for all $\zbar\in\extspace$ with $\PP\zbar=\PP\xbar$,
it follows from \eqref{eq:thm:F-reduc-PFP:1}
that $G(\xbar)\leq H(\xbar)$, completing the proof.

\pfpart{Part~(\ref{thm:F-reduc-PFP:b}):}
This follows immediately from
\Cref{thm:recF-equivs}(\ref{thm:recF-equivs:a},\ref{thm:recF-equivs:c})
(and since $\limray{\ebar}=\ebar$).
\qedhere
\end{proof-parts}
\end{proof}

Note in \Cref{thm:F-reduc-PFP}
that the projection matrix $\PP$ depends only on
$\colspace{\VV}$, in other words, not on the individual columns of
$\VV$ or their ordering within the matrix, but only on the linear
subspace that they span.
This means, assuming $\VV\omm\in\aresconeF$, that the value of the
function $F$ at a point $\VV\omm\plusl\xbar$ is determined,
according to this theorem, only by $\xbar$ and $\colspace{\VV}$,
not on any other details regarding the matrix $\VV$.%
\indexg{reductions, iconic!astral analogue|)}

We next look at how some of our earlier results for astronic
reductions can be generalized for iconic reductions.
We start by generalizing \Cref{thm:a10-nunu},
which, for a convex function $f:\Rn\rightarrow\Rext$,
earlier showed that $\fshadv=\lsc \fv$, provided
$\vv\in\resc{f}$.
\indexg{reductions, iconic!lower semicontinuous hull@as lower semicontinuous hull|(}%
\indexg{shadow (of function)!iconic reductions and|(}%
Now, for any matrix $\VV\in\Rnk$ whose associated icon $\VV\omm$ is in
$\aresconef$, we show that $\fshadVo = \lsc \fV$,
that is, that the reduction $\fshadVo$ is the lower semicontinuous
hull of the shadow $\fV$
(as was defined in \Cref{dfn:shadow-fcns}).

\begin{theorem}  \label{thm:icon-reduc-lsc-inf}
Let $f:\Rn\rightarrow\Rext$ be convex,
let $\VV\in\R^{n\times k}$, and assume
$\VV\omm\in\aresconef$.
Then $\fshadVo = \lsc \fV$.
\end{theorem}

\begin{proof}
Let $\PP$ be the projection matrix onto $(\colspace{\VV})^\perp$.
We have
\begin{equation}  \label{eq:thm:icon-reduc-lsc-inf:n1}
  (\PP\fext)\PP
  =
  (\Pfext) \PP
  =
  \PfPext.
\end{equation}
The first equality is by \Cref{thm:inf-lin-ext}.
The second is by
\Cref{thm:ext-linear-comp}
(applied with $f$ and $\A$, as they appear in that
\namecref{thm:ext-linear-comp},
set to $\PPf$ and $\PP$),
noting that
$\PPf$ is convex (\Cref{roc:thm5.7:Af}),
and that there exists $\xhat\in\Rn$ with
$\PP\xhat\in\ri(\dom \PPf)$ by
\Cref{pr:Ax-in-ri-dom-Af}.
We then have, for $\xx\in\Rn$,
\[
   \fshadVo(\xx)
   =
   \regBracks{(\PP\fext)\PP}(\xx)
   =
   \regBracks{\PfPext}(\xx)
   =
   \bigParens{\lsc [(\PPf) \PP]}(\xx)
   =
   (\lsc \fV)(\xx).
\]
The first equality is obtained by applying
\Cref{thm:F-reduc-PFP} with $F=\fext$ so that
$G(\xx)$, as given in \eqref{eq:thm:F-reduc-PFP:0},
is equal to $\fshadVo(\xx)$ as defined in
\eqref{eqn:def:iconic-reduction}.
The second equality is by \eqref{eq:thm:icon-reduc-lsc-inf:n1}.
The third is by \Cref{pr:h:1}(\ref{pr:h:1a}).
And the fourth is by \Cref{pr:gtil-is-PPfP-gen}.%
\indexg{reductions, iconic!lower semicontinuous hull@as lower semicontinuous hull|)}%
\indexg{shadow (of function)!iconic reductions and|)}%
\end{proof}

\indexg{reductions, iconic!conjugate of|(}%
Next, we
generalize \Cref{thm:e1}(\ref{thm:e1:a})
to compute the conjugate of the reduction $\fshadd$ for any
icon $\ebar$ in the astral recession cone $\aresconef$.

\begin{theorem}  \label{thm:conj-of-iconic-reduc}
Let $f:\Rn\rightarrow\Rext$ be convex, $f\not\equiv+\infty$, and
let $\ebar\in\corezn\cap (\aresconef)$.
Let
$g=\fshad{\ebar}$
be the reduction of $f$ at $\ebar$.
Then for all $\uu\in\Rn$,
\[
\gstar(\uu) =
  \begin{cases}
    \fstar(\uu)
    & \text{if $\ebar\cdot\uu=0$,}
  \\
    +\infty
    & \text{otherwise.}
  \end{cases}
\]
\end{theorem}

\begin{proof}
By \Cref{pr:icon-equiv}(\ref{pr:icon-equiv:a},\ref{pr:icon-equiv:c}),
$\ebar=\VV\omm$ for some matrix $\VV\in\Rnk$.
Let $\gtil=\fV$, and let $L=\{\uu\in\Rn :\: \uu\perp\VV\}$.
We then have:
\[
  \gstar
  =
  (\lsc \gtil)^*
  =
  \gtil^*
  =
  \fstar \plusu \indf{L}.
\]
The first equality is because $g=\lsc\gtil$
by \Cref{thm:icon-reduc-lsc-inf}.
The second is by \Cref{pr:conj-props}(\ref{pr:conj-props:e}).
The third is by \Cref{thm:fV-shad-conj}.
Further, by \Cref{pr:vtransu-zero},
$\uu\in L$ if and only if $\ebar\cdot\uu=0$,
completing the proof.%
\indexg{reductions, iconic!conjugate of|)}%
\end{proof}

\indexg{reductions, iconic!relative interior of domain|(}%
As another consequence of
\Cref{thm:icon-reduc-lsc-inf},
we can express the relative interior of the
effective domain of any reduction at an icon in $\aresconef$:

\begin{theorem}  \label{thm:ri-dom-icon-reduc}
Let $f:\Rn\rightarrow\Rext$ be convex,
let $\VV\in\R^{n\times k}$, and assume $\VV\omm\in\aresconef$.
Then
\[
  \ri(\dom{\fshadVo})
  =
  \ri(\dom{f}) + (\colspace\VV)
  \supseteq
  \ri(\dom{f}).
\]
\end{theorem}

\begin{proof}
Let $g=\fshadVo$
and let $\gtil=\fV$.
Then
$\gtil$ is convex by \Cref{pr:gtil-is-PPfP-gen},
and $g=\lsc \gtil$
by \Cref{thm:icon-reduc-lsc-inf}.
Note that, by \Cref{pr:lsc-props}(\ref{pr:lsc-props:c}),
\begin{equation}  \label{eq:thm:ri-dom-icon-reduc:1}
  \ri(\dom{g}) = \ri(\dom{\gtil}).
\end{equation}

From the form of $\gtil=\fV$
(see \Cref{dfn:shadow-fcns}), for $\xx\in\Rn$,
$\gtil(\xx)<+\infty$ if and only if there exists $\yy\in\colspace\VV$
with $f(\xx+\yy)<+\infty$.
Thus,
\[
  \dom{\gtil} = (\dom{f}) + (\colspace\VV),
\]
and therefore,
\[
  \ri(\dom{\gtil})
  =
  \ri(\dom{f}) + \ri(\colspace\VV)
  =
  \ri(\dom{f}) + (\colspace\VV)
\]
by
\Cref{pr:ri-props}(\ref{pr:ri-props:roc-cor6.6.2}),
and since $\colspace\VV$, being a linear subspace, is
relatively open.
Combined with
\eqref{eq:thm:ri-dom-icon-reduc:1}, this completes the proof.%
\indexg{reductions, iconic!relative interior of domain|)}%
\end{proof}

\chapter{Universal reduction and canonical minimizers}
\chaptermark{Universal reduction, canonical minimizers}
\label{sec:univ-red-and-min}

For a convex function $f:\Rn\rightarrow\Rext$,
we have seen so far that if $\xbar=\ebar\plusl\qq$ minimizes $\fext$,
where $\ebar\in\corezn$ and $\qq\in\Rn$, then $\xbar$'s iconic part
$\ebar$ must be in $\fext$'s
astral recession cone (\Cref{thm:arescone-fshadd-min}).
In this chapter, we delve further into the structure of $\fext$'s
minimizers, both their iconic and finite parts.
We will see that all finite parts $\qq$ of all minimizers
can be obtained as the (finite) minimizers of one
particular convex function called the universal reduction,
defined in a moment.
Furthermore, all of the minimizers of this function
are necessarily in a bounded region of $\Rn$
(up to displacements in directions where this function is constant).
This alleviates the issue of minimizers at infinity
and allows us to reduce the problem of
finding the finite parts $\qq$ of all minimizers
of $\fext$ to minimization of a standard convex
function in the most favorable setting, where finite minimizers exist
and only occur within some compact region.

Furthermore,
we will see that there exist choices for the icon $\ebar$ that
minimize
$\fext(\ebar\plusl\xx)$ in the variable $\ebar$
\emph{simultaneously} for all values of $\xx\in\Rn$.
We will discuss properties of such icons, how to find them,
and how they combine with the minimizers of the universal reduction
to give rise to minimizers of $\fext$ that are, in a sense, canonical.

\section{The universal reduction}
\label{sec:fullshad}

\indexg{universal reduction|(}%
We begin by defining the universal reduction:

\begin{definition}  \label{def:univ-reduction}
\indexg{universal reduction!defined|(}%
Let $f:\Rn\rightarrow\Rext$ be convex.
The \emph{universal reduction} of $f$ is the function
$\fullshad{f}:\Rn\rightarrow \Rext$ defined, for $\xx\in\Rn$, by
\begin{equation}  \label{eqn:fullshad-def}
\indexm{f 900}{$\fullshad{f}$}{universal reduction}%
  \fullshad{f}(\xx) = \inf_{\ebar\in\corezn} \fext(\ebar\plusl\xx)
         = \inf_{\ebar\in\corezn} \fshadd(\xx).%
\indexg{universal reduction!defined|)}%
\end{equation}
\end{definition}
Thus, $\fullshad{f}$
computes the minimum possible value
of $\fext$ when a given point $\xx\in\Rn$ is combined with any icon
$\ebar\in\corezn$.
Said differently, $\fullshad{f}$ is the pointwise infimum of all
reductions $\fshadd$ over all icons $\ebar\in\corezn$.
Intuitively, $\fullshad{f}$ shortcuts the process of minimizing
the function $f$ by jumping to the end of any astral ray, as described
by an icon $\ebar$, and focusing the minimization on the remaining finite part.
Later, in \Cref{cor:fullshad:absorb}(\ref{cor:fullshad:absorb:a}),
we will prove
that $\fullshad{f}$ can actually be obtained as
an iconic reduction, for an appropriately selected icon,
so calling it a ``reduction''
is consistent with our earlier terminology.

Here are some simple facts about universal reductions:

\begin{proposition}   \label{pr:univ-red-props}
  Let $f:\Rn\rightarrow\Rext$ be convex.
  Then:
  \begin{letter-compact}
  \item   \label{pr:univ-red-props:min}
    $\fullshad{f}$ attains its minimum, and
    $\min\fullshad{f}=\min\fext=\inf f$.
  \item   \label{pr:univ-red-props:a}
    $\fullshad{f}\leq \lsc{f}\leq f$.
  \item   \label{pr:univ-red-props:b}
    $\fullshad{(\lsc f)}=\fullshad{f}$.
  \end{letter-compact}
\end{proposition}

\begin{proof}
  ~

\begin{proof-parts}
\pfpart{Part~(\ref{pr:univ-red-props:min}):}
By \Cref{pr:fext-min-exists}, $\min\fext=\inf f$, and also $\fext$
attains its minimum
at some point $\xbar\in\extspace$.
We can write $\xbar=\ebar\plusl\qq$ for some $\ebar\in\corezn$ and
$\qq\in\Rn$.
Then
\[
  \inf\fullshad{f}
  \leq
  \fullshad{f}(\qq)
  \leq
  \fext(\ebar\plusl\qq)
  =
  \fext(\xbar)
  =
  \min\fext
  \leq
  \inf\fullshad{f}.
\]
The second inequality is by $\fullshad{f}$'s definition and since
$\ebar\in\corezn$.
The last inequality is because
$\fullshad{f}(\xx)\geq\min \fext$ for all $\xx\in\Rn$,
again by $\fullshad{f}$'s definition.
Thus, $\qq$ minimizes $\fullshad{f}$, and
$\min\fullshad{f}=\min\fext$.

\pfpart{Part~(\ref{pr:univ-red-props:a}):}
For $\xx\in\Rn$, we have
$\fullshad{f}(\xx)\leq\fext(\zero\plusl\xx) = (\lsc f)(\xx) \leq f(\xx)$.
The first inequality is from $\fullshad{f}$'s
definition, and because $\zero$ is an icon.
The equality and last inequality are by
\Cref{pr:h:1}(\ref{pr:h:1a}).

\pfpart{Part~(\ref{pr:univ-red-props:b}):}
This follows from $\fullshad{f}$'s definition
since $\lscfext=\fext$
by \Cref{pr:h:1}(\ref{pr:h:1aa}).
\qedhere
\end{proof-parts}
\end{proof}

The definition of $\fullshad{f}$ is based on minimization
of $\fext(\ebar\plusl\xx)$ over all icons $\ebar$.
\indexg{universal reduction!equivalent definitions|(}%
As shown in the next proposition,
the definition would yield the same function
if it were instead based on minimization of $\fext(\zbar\plusl\xx)$
over all $\zbar$ in the astral recession cone,
whether or not further restricting $\zbar$ to be an icon.
Furthermore, in all cases, including
\eqref{eqn:fullshad-def}, the respective infima are always attained,
which means we can state these expressions in
terms of minima rather than infima.

\begin{proposition}  \label{pr:fullshad-equivs}
  Let $f:\Rn\rightarrow\Rext$ be convex.
  Then for all $\xx\in\Rn$,
  \begin{equation}  \label{eq:pr:fullshad-equivs:8}
     \fullshad{f}(\xx)
      = \min_{\ebar\in\corezn} \fext(\ebar\plusl\xx)
      = \min_{\zbar\in\aresconef}\,\fext(\zbar\plusl\xx)
      = \min_{\ebar\in(\aresconef)\cap\corezn} \fext(\ebar\plusl\xx).
  \end{equation}
  In particular, this means that each of these minima is attained.
\end{proposition}

\begin{proof}
Let $\xx\in\Rn$.
By \Cref{cor:a:4},
if $\ebar\in\corezn\setminus(\aresconef)$,
then $\fext(\ebar\plusl\xx)=+\infty$.
Therefore,
\begin{equation}   \label{eqn:pr:fullshad-equivs:2}
  \fullshad{f}(\xx)
   = \inf_{\ebar\in(\aresconef)\cap\corezn} \fext(\ebar\plusl\xx)
   \geq
     \inf_{\zbar\in\aresconef} \fext(\zbar\plusl\xx),
\end{equation}
where the inequality is simply because
$\aresconef$ is a superset of $(\aresconef)\cap\corezn$.

Let $\zbar\in\aresconef$, implying
$\limray{\zbar}$ also is in $\aresconef$
by \Cref{pr:ast-cone-is-naive} since
$\aresconef$ is an astral cone
(by \Cref{cor:res-fbar-closed}).
Then
\[
  \ef(\zbar\plusl\xx)
  \ge
  \ef(\limray{\zbar}\plusl\zbar\plusl\xx)
  =
  \ef(\limray{\zbar}\plusl\xx)
  \ge
  \inf_{\ebar\in(\aresconef)\cap\corezn} \fext(\ebar\plusl\xx).
\]
The first inequality is
by
\Cref{thm:recF-equivs}(\ref{thm:recF-equivs:a},\ref{thm:recF-equivs:b})
since $\zbar\in\aresconef$.
The equality is because
$\limray{\zbar}\plusl\zbar=(\oms+1)\zbar=\limray{\zbar}$
by \Cref{pr:i:7}(\ref{pr:i:7c}).
The last
inequality is because $\limray{\zbar}$ is also an icon
(\Cref{pr:i:8}\ref{pr:i:8-infprod}).

Since this holds for all $\zbar\in\aresconef$, it follows that
\[
     \inf_{\zbar\in\aresconef}\,\fext(\zbar\plusl\xx)
     \ge
     \inf_{\ebar\in(\aresconef)\cap\corezn} \fext(\ebar\plusl\xx),
\]
which, combined with \eqref{eqn:pr:fullshad-equivs:2}, proves the proposition,
with infima instead of minima.

The map $\zbar\mapsto\fext(\zbar\plusl\xx)$, for $\zbar\in\extspace$,
is lower semicontinuous by
\Cref{pr:lsc-res-domain}(\ref{pr:lsc-res-domain:a})
because $\fext$ is lower semicontinuous (\Cref{prop:ext:F}\ref{prop:ext:F:a}),
and because the astral affine map $\zbar\mapsto\zbar\plusl\xx$, for
$\zbar\in\extspace$, is continuous (\Cref{cor:aff-cont}).
Also, 
the sets $\corezn$, $\aresconef$, as well as their intersection, are
all closed
(by \Cref{pr:i:8}\ref{pr:i:8e} and \Cref{cor:res-fbar-closed}), and so
are compact (\Cref{prop:compact}\ref{prop:compact:closed-subset}).
Therefore, each of the minima in \eqref{eq:pr:fullshad-equivs:8}
is attained by
\Cref{thm:weierstrass}.%
\indexg{universal reduction!equivalent definitions|)}%
\end{proof}

\indexg{universal reduction!minimizers of extension and|(}%
\indexg{minimizers of extensions!finite parts characterized|(}%
Using \Cref{pr:fullshad-equivs}, we can now show that
the set of all minimizers of $\fullshad{f}$ is exactly
equal to the set of all finite parts of all minimizers of $\fext$:

\begin{theorem}  \label{pr:min-fullshad-is-finite-min}
  Let $f:\Rn\rightarrow\Rext$ be convex,
  and let $\qq\in\Rn$.
  Then $\qq$ minimizes $\fullshad{f}$ if and only if
  there exists $\ebar\in\corezn$ such that
  $\ebar\plusl\qq$ minimizes $\fext$.
\end{theorem}

\begin{proof}
Suppose first that $\qq$ minimizes $\fullshad{f}$.
Then by \Cref{pr:fullshad-equivs}, there exists
$\ebar\in\corezn$ such that
$\fext(\ebar\plusl\qq)=\fullshad{f}(\qq)$.
Since $\fullshad{f}(\qq)=\min\fullshad{f}=\min\fext$, with the second
equality by
\Cref{pr:univ-red-props}(\ref{pr:univ-red-props:min}), it follows that
$\ebar\plusl\qq$ minimizes $\fext$.

Conversely, suppose now
that $\ebar\plusl\qq$ minimizes $\fext$, for some
$\ebar\in\corezn$.
Then
\[
  \fullshad{f}(\qq)
  \leq
  \fext(\ebar\plusl\qq)
  =
  \min \fext
  =
  \min\fullshad{f},
\]
where the inequality is from
$\fullshad{f}$'s definition
(Eq.~\ref{eqn:fullshad-def}),
and the second equality is by
\Cref{pr:univ-red-props}(\ref{pr:univ-red-props:min}).
Thus, $\qq$ minimizes $\fullshad{f}$.%
\indexg{universal reduction!minimizers of extension and|)}%
\indexg{minimizers of extensions!finite parts characterized|)}%
\end{proof}

\indexg{universal reduction!preserved under reduction|(}%
An important property of the universal reduction $\fullshad{f}$
is that it is
preserved when the function $f$ is reduced at an astron
$\limray{\vv}$ with $\vv\in\resc{f}$;
that is, if $g$ is a reduction
of $f$ at such an astron, then $\fullshad{g}=\fullshad{f}$.
This observation will be at the crux of our iterative algorithm for
constructing $\fullshad{f}$.

\begin{theorem}  \label{thm:f:5}
  Let $f:\Rn\rightarrow\Rext$ be convex,
  let $\vv\in\resc{f}$, and let
  $g=\fshadv$
  be the reduction of $f$ at $\limray{\vv}$.
  Then the universal reductions of $f$ and $g$ are identical;
  that is,
  $\fullshad{g}=\fullshad{f}$.
\end{theorem}

\begin{proof}
Let $\xx\in\Rn$.
We have
\begin{equation*}
  \fullshad{g}(\xx)
  =
  \min_{\ebar\in\corezn} \gext(\ebar\plusl\xx)
  =
  \min_{\ebar\in\corezn} \fext(\limray{\vv}\plusl \ebar \plusl \xx)
  \geq
  \min_{\ebar\in\corezn} \fext(\ebar \plusl \xx)
  =
  \fullshad{f}(\xx).
\end{equation*}
The first and last equalities are by
\Cref{pr:fullshad-equivs}.
The second equality is
by \Cref{cor:i:1}(\ref{cor:i:1b}),
and the inequality
is because
$\limray{\vv}\plusl\ebar$ is an icon since
$\ebar$ is
(by \Cref{pr:i:8}\ref{pr:i:8-leftsum}).

On the other hand,
\[
  \fullshad{f}(\xx)
  =
  \min_{\ebar\in\corezn} \fext(\ebar \plusl \xx)
  \geq
  \min_{\ebar\in\corezn} \gext(\ebar \plusl \xx)
   =
  \fullshad{g}(\xx),
\]
since $\ef\ge\eg$ by
\Cref{cor:i:1}(\ref{cor:i:1c}).%
\indexg{universal reduction!preserved under reduction|)}%
\end{proof}

\indexg{universal reduction!when same as original function|(}%
The next theorem characterizes when $f=\fullshad{f}$.
This characterization, which will soon be used as the termination
condition in our algorithm for constructing $\fullshad{f}$,
requires that
the recession cone of $f$ be the same as its constancy space.
This same condition appeared in \Cref{sec:astral-cone},
where we proved it implies the existence of a finite minimizer for $f$
(\Cref{pr:cons-implies-fin-minimizer});
that it further implies $f=\fullshad{f}$ strengthens that statement.

\begin{theorem}
\label{thm:fullshad:char}
  Let $f:\Rn\rightarrow\Rext$ be convex and lower semicontinuous.
  Then $f=\fullshad{f}$ if and only if
  $\resc{f}=\conssp{f}$.
\end{theorem}

\begin{proof}
  ~
\begin{proof-parts}
\pfpart{``Only if'' ($\Rightarrow$):}
Suppose $f=\fullshad{f}$.
By \Cref{pr:prelim:const-props}(\ref{pr:prelim:const-props:a}),
$\conssp{f}\subseteq\resc{f}$.
For the reverse inclusion, let $\vv\in\resc{f}$,
which we aim to show is in $\conssp{f}$.
Then $\fshadv\le f=\fullshad{f}\le\fshadv$, with the first inequality
by \Cref{pr:d2}(\ref{pr:d2:b}) and the second from the definition
of $\fullshad{f}$
(Eq.~\ref{eqn:fullshad-def}). Thus, $\fshadv=f$,
so
\[\vv\in\rspanset{\limray{\vv}}\subseteq\conssp{\fshadv}=\conssp{f},\]
where the first inclusion is by definition of representational span
(\Cref{def:rspan-gen}), and the second is by
\Cref{pr:i:9}(\ref{pr:i:9cons}).
Hence,
$\resc{f}\subseteq\conssp{f}$, completing the proof.

\pfpart{``If'' ($\Leftarrow$):}
For the sake of contradiction, suppose $\resc{f}=\conssp{f}$ but
$f\neq\fullshad{f}$.
Then $\fullshad{f}(\xx)\neq f(\xx)$ for some $\xx\in\Rn$.
Since also $\fullshad{f}(\xx)\leq f(\xx)$
by \Cref{pr:univ-red-props}(\ref{pr:univ-red-props:a}),
this further implies that
$\fullshad{f}(\xx) < f(\xx)$.

Let $\ebar\in\corezn$ be an icon with
$\ef(\ebar\plusl\xx)<f(\xx)$, and, among all such icons, let
$\ebar$ be one with minimum astral rank $r$.
Note that such an icon must exist 
by $\fullshad{f}$'s definition
since $\fullshad{f}(\xx) < f(\xx)$.
Also, $\ef(\zero\plusl\xx)=f(\xx)$
by \Cref{pr:h:1}(\ref{pr:h:1a})
since $f$ is lower semicontinuous;
therefore, $\ebar\neq\zero$ so $r>0$.

Let $\vv$ be
the dominant direction of $\ebar$, and
let $\PP$ be the projection matrix orthogonal to $\vv$.
Then by \Cref{pr:h:6}, $\ebar=\limray{\vv}\plusl\dbar$,
where $\dbar=\PP\ebar$ has astral rank $r-1$; further, $\dbar$ is an icon (by \Cref{pr:i:8}\ref{pr:i:8-matprod}).

Since $\ef(\ebar\plusl\xx)<f(\xx)\le+\infty$,
we have $\ebar\in\aresconef$
(by \Cref{cor:a:4}), implying
$\vv\in\resc{f}=\conssp{f}$
(by \Cref{thm:f:3}\ref{thm:f:3b}).
Therefore,
$f(\yy)=f(\PP\yy)$ for all $\yy\in\Rn$ (\Cref{pr:cons-PP}),
so also
$\ef(\ybar)=\ef(\PP\ybar)$ for all $\ybar\in\eRn$ (\Cref{pr:PP:ef}).
Hence,
      \[
        \ef(\dbar\plusl\xx)
        =
        \ef(\PP\dbar\plusl\PP\xx)
        =
        \ef(\PP\ebar\plusl\PP\xx)
        =
        \ef(\ebar\plusl\xx)
        <
        f(\xx),
      \]
where the second equality is because
$\PP\dbar=\PP\PP\ebar=\PP\ebar$
(by \Cref{pr:proj-mat-props}\ref{pr:proj-mat-props:b}).
But this contradicts that $\ebar$ has minimum astral rank among all icons
with $\fext(\ebar\plusl\xx)<f(\xx)$
(since $\dbar$ has strictly lower astral rank).
Therefore, $f=\fullshad{f}$.%
\indexg{universal reduction!when same as original function|)}%
\qedhere
\end{proof-parts}
\end{proof}

\begin{algorithm}[t]
  \caption[Generate $\fullshad{f}$ and all canonical minimizers of $\fext$]%
          {Generate the universal reduction $\fullshad{f}$ and all canonical minimizers of $\fext$.}
  \label{fig:min-proc}
\begin{algorithmic}
    \Block{Input}
        \State  convex function
\indexalg{universal reduction!construction of}%
\indexalg{canonical minimizers!generation of all}%
\indexalg{universal reducers!generation of all}%
                                 $f:\Rn\rightarrow\Rext$
    \EndBlock
    \smallskip
    \Block{Define \Call{GenerateCanonicalMinimizer}{$f$}}
        \State  $i \leftarrow 0$    
        \State  $g_0 = \lsc f$
        \Loop{repeat \emph{at least} until $\resc{g_i} = \conssp{g_i}$}
            \Comment{nondeterministic choice}
            \State  $i \leftarrow i+1$
            \State  let $\vv_i$ be \emph{any} point in $\resc{g_{i-1}}$
                \Comment{nondeterministic choice}
            \State  $g_i = \gishadvi$
        \EndLoop
        \State  $k \leftarrow i$
        \State  $\ebar=\limrays{\vv_1,\ldots,\vv_k}$
        \State  let $\qq\in\Rn$ be \emph{any} finite minimizer of $g_k$
            \Comment{nondeterministic choice}
        \State  $\xbar=\ebar\plusl\qq$
        \State  return $(g_k,\ebar,\qq,\xbar)$
    \EndBlock
    \smallskip
    \Block{Output properties}
        \State  $g_k = \fullshad{f}$
        \State  $\ebar\in \unimin{f} \subseteq (\aresconef)\cap\corezn$
           \Comment{definition and analysis in \Cref{subsec:univ-min}}
        \State  $\qq$ minimizes $\fullshad{f}$
        \State  $\xbar$ is a canonical minimizer of $\fext$
           \Comment{definition and analysis in \Cref{subsec:univ-min}}
        \EndBlock
\end{algorithmic}
\end{algorithm}

\indexg{universal reduction!construction of|(}%
\indexg{canonical minimizers!generation of all|(}%
\indexg{universal reducers!generation of all|(}%
\Cref{thm:f:5,thm:fullshad:char}
give rise to \Cref{fig:min-proc} for constructing the universal reduction
$\fullshad{f}$.
The algorithm also generates a specific type of minimizer of $\ef$
called a canonical minimizer and mentions the set $\unimin{f}$.
These concepts will be introduced in \Cref{subsec:univ-min},
and can be disregarded for now.%
\indexg{universal reducers!generation of all|)}%
\indexg{canonical minimizers!generation of all|)}

\Cref{fig:min-proc} is very similar to \Cref{fig:all-min-proc},
modified only in the termination condition for the main loop.
Both algorithms construct a sequence of functions $g_i$, each the reduction of
$g_{i-1}$ at an astron associated with a point $\vv_i\in\resc{g_{i-1}}$.
However, in contrast with \Cref{fig:all-min-proc},
the loop in \Cref{fig:min-proc} terminates
when $\resc{g_i}=\conssp{g_i}$, rather than when $g_i$ has a finite
minimizer. By \Cref{pr:cons-implies-fin-minimizer}, the new termination
condition
implies the existence of a finite minimizer, and
by \Cref{thm:f:5,thm:fullshad:char}, it also guarantees that
$\fullshad{f}=\fullshad{g_0}=\dotsb=\fullshad{g_k}=g_k$, so if the algorithm
terminates then it indeed returns the universal reduction $\fullshad{f}$.
Moreover, this universal reduction takes the form of an iconic reduction,
$g_k=\fshad{\ebar}$ where $\ebar=\limrays{\vv_1,\ldots,\vv_k}$
(by \Cref{pr:icon-red-decomp-astron-red}).

\begin{example}[Two exps and a square, continued]
\label{ex:simple-eg-exp-exp-sq-part2}
\indexg{Two exps and a square!algorithm applied to|(}%
  In \Cref{ex:simple-eg-exp-exp-sq},
  we considered a run of \Cref{fig:all-min-proc} on the function
  \begin{align*}
     f(\xx) &= f(x_1,x_2,x_3) = \me^{x_3-x_1} + \me^{-x_2}
                 + (2+x_2-x_3)^2,
  \intertext{%
    with $\vv_1=\trans{[1,1,1]}$, $\vv_2=\trans{[1,-1,-1]}$,
    and
  }
    g_1(\xx) &= \fshad{\limray{\vv_1}}(\xx) = \me^{x_3-x_1} + (2+x_2-x_3)^2,
  \\
    g_2(\xx) &= \genshad{g_1}{\limray{\vv_2}}(\xx) = (2+x_2-x_3)^2,
  \end{align*}
  terminating at $g_2$, which has a finite minimizer
  (for instance $\qq=\trans{[0,0,2]}$).

  In fact, that identical run could also have occurred using
  \Cref{fig:min-proc} since the function $g_2$ satisfies
  $\resc{g_2} = \conssp{g_2} = \set{ \zz\in\R^3 :\: z_2=z_3 }$.
  So, $\fullshad{f}=g_2$.%
\indexg{Two exps and a square!algorithm applied to|)}%
\end{example}

Similar to \Cref{fig:all-min-proc}, at any point of execution,
\Cref{fig:min-proc} can switch to making only nonredundant choices,
that is, choices
$\vv_i \in (\resc{g_{i-1}}) \setminus \spnfin{\vv_1,\ldots,\vv_{i-1}}$.
After at most $n$ nonredundant choices, the algorithm must reach a point when
no nonredundant choices are possible. By \Cref{pr:no-sense-then-done},
this means that $\resc{g_i}=\conssp{g_i}$, and so that the termination condition
has been reached. Thus, it is always possible to complete a run
of \Cref{fig:min-proc} in at most $n$ iterations.

The properties of \Cref{fig:min-proc} are summarized in the next
\namecref{cor:fig-cons-main-props}.

\begin{theorem}  \label{cor:fig-cons-main-props}
  Let $f:\Rn\rightarrow\Rext$ be convex,
  and let $g_0=\lsc f$.
  Then there exist $\vv_1,\ldots,\vv_k\in\Rn$,
  for some $k\geq 0$, so that the following hold,
  where $g_i = \gishadvi$, for $i=1,\ldots,k$:
  \begin{item-compact}
    \item
      $\vv_i\in\resc{g_{i-1}}$,
      for $i=1,\ldots,k$;
    \item
      $\resc{g_k}=\conssp{g_k}$.
  \end{item-compact}
  Whenever these hold, the following are also true:
  \begin{letter-compact}
  \item  \label{cor:fig-cons-main-props:a}
     $\fullshad{f}=\fullshad{g_k}=g_k=\fshadd$, where $\ebar=\limrays{\vv_1,\ldots,\vv_k}$.
  \item  \label{cor:fig-cons-main-props:b}
    There exists some $\qq\in\Rn$ that minimizes
    $g_k=\fullshad{f}$.
  \item  \label{cor:fig-cons-main-props:c}
    The point $\xbar=\ebar\plusl\qq$ minimizes $\fext$.
  \end{letter-compact}
\end{theorem}

\begin{proof}
Suppose that in \Cref{fig:min-proc}, we only make
nonredundant choices $\vv_i$, that is,
$\vv_i \in (\resc{g_{i-1}}) \setminus \spnfin{\vv_1,\ldots,\vv_{i-1}}$.
Then the dimension of $\spnfin{\vv_1,\ldots,\vv_i}\subseteq\Rn$ is
equal to $i$, and so, after at most $n$ itreations, the algorithm
must reach a point at which no such choice is possible, implying, by
\Cref{pr:no-sense-then-done}, that the termination
condition of the main loop has been reached.
Upon termination, all of the claimed properties hold:

\begin{proof-parts}
\pfpart{Part~(\ref{cor:fig-cons-main-props:a}):}
We have
\[
  \fullshad{f}=\fullshad{g_0}=\dotsb=\fullshad{g_k}=g_k=\fshadd.
\]
The first equality is by
\Cref{pr:univ-red-props}(\ref{pr:univ-red-props:b}),
and the last two equalities are, respectively, by
\Cref{thm:fullshad:char}
and \Cref{pr:icon-red-decomp-astron-red}.
The remaining middle equalities (that is, that
$\fullshad{g_{i-1}}=\fullshad{g_{i}}$
for $i=1,\ldots,k$)
are by \Cref{thm:f:5}.

\pfpart{Part~(\ref{cor:fig-cons-main-props:b}):}
This follows from \Cref{pr:cons-implies-fin-minimizer}
(noting that $g_k$ is convex and lower semicontinuous by
\Cref{pr:icon-red-decomp-astron-red}).

\pfpart{Part~(\ref{cor:fig-cons-main-props:c}):}
Immediate by \Cref{thm:all-min-proc-correct}.%
\indexg{universal reduction!construction of|)}%
\qedhere
\end{proof-parts}
\end{proof}

\indexg{universal reduction!iconic reduction@as iconic reduction|(}%
As corollary,
\Cref{cor:fig-cons-main-props} implies that $\fullshad{f}$ is convex
and lower semicontinuous, and that it
can always be expressed as an iconic reduction,
justifying our terminology.
We also
can show that the universal reduction $\fullshad{f}$ is unchanged
by taking an astronic reduction at $\vv\in\resc{\fullshad{f}}$
or by taking its own universal reduction. Thus,
$\fullshad{f}$ is in a sense the maximal iconic reduction of $f$,
which means
running the algorithm for further iterations would not
change the result.

\begin{corollary}
\label{cor:fullshad:absorb}
  Let $f:\Rn\rightarrow\Rext$ be convex.
  Then:
  \begin{letter-compact}
  \item   \label{cor:fullshad:absorb:b}
    $\fullshad{f}$ is convex and lower semicontinuous.
  \item   \label{cor:fullshad:absorb:a}
    $\fullshad{f}$ is an iconic reduction; that is,
    $\fullshad{f}=\fshadd$ for some $\ebar\in\corezn$.
  \item   \label{cor:fullshad:absorb:e}
    $\resc{\fullshad{f}}=\conssp{\fullshad{f}}$.
  \item   \label{cor:fullshad:absorb:c}
    $\fullshad{(\fullshad{f})}=\fullshad{f}$.
  \item   \label{cor:fullshad:absorb:d}
    $\genshad{(\fullshad{f})}{\limray{\vv}}=\fullshad{f}$
    for all $\vv\in\resc{\fullshad{f}}$.
  \end{letter-compact}
\end{corollary}

\begin{proof}
  ~

\begin{proof-parts}
\pfpart{%
  Parts~(\ref{cor:fullshad:absorb:b}),~(\ref{cor:fullshad:absorb:a}),~(\ref{cor:fullshad:absorb:e}),
  and~(\ref{cor:fullshad:absorb:c}):
}
By \Cref{cor:fig-cons-main-props}(\ref{cor:fig-cons-main-props:a}),
$ \fullshad{f}=\fullshad{g_k}=g_k=\fshadd $,
where $g_k$ and icon $\ebar$ are as were shown to exist in that
\namecref{cor:fig-cons-main-props}.
This proves
part~(\ref{cor:fullshad:absorb:a}), as well as
part~(\ref{cor:fullshad:absorb:e}) since
$\resc{g_k}=\conssp{g_k}$.
Since $\fshadd$ is convex and lower semicontinuous
by \Cref{pr:i:9}(\ref{pr:i:9a}), this also proves
part~(\ref{cor:fullshad:absorb:b}).
Finally, since $g_k=\fullshad{f}$ and $\fullshad{f}=\fullshad{g_k}$,
this proves
part~(\ref{cor:fullshad:absorb:c}).

\pfpart{Part~(\ref{cor:fullshad:absorb:d}):}
Let $\vv\in\resc{\fullshad{f}}$.
Then
\[
  \genshad{(\fullshad{f})}{\limray{\vv}}
  \leq
  \fullshad{f}
  =
  \fullshad{(\fullshad{f})}
  \leq
  \genshad{(\fullshad{f})}{\limray{\vv}}.
\]  
The first inequality is by \Cref{pr:d2}(\ref{pr:d2:b}).
The equality is by part~(\ref{cor:fullshad:absorb:c}).
The second inequality is by definition of universal reduction
(\Cref{def:univ-reduction}) since $\limray{\vv}$ is an icon.
This proves the claim.%
\indexg{universal reduction!iconic reduction@as iconic reduction|)}%
\qedhere
\end{proof-parts}
\end{proof}

\indexg{universal reduction!recessive completeness of|(}%
As a final property, we can now show that
every universal reduction is also recessive complete:

\begin{theorem}
\label{pr:fullshad-recess-complete}
  Let $f:\Rn\rightarrow\Rext$ be convex.
  Then $\fullshad{f}$ is recessive complete.
\end{theorem}

\begin{proof}
By \Cref{cor:fullshad:absorb}(\ref{cor:fullshad:absorb:b}),
$\fullshad{f}$ is convex and lower semicontinuous.
Let $\xbar\in\arescone{\fullshadfbar}$. We will show that $\xbar\in\represc{\fullshad{f}}$,
which, when combined with \Cref{pr:repres-in-arescone},
will imply
$\arescone{\fullshadfbar}=\represc{\fullshad{f}}$,
and so that $\fullshad{f}$ is recessive complete.

We can write $\xbar=\limrays{\vv_1,\ldots,\vv_k}\plusl\qq$
for some $\vv_1,\ldots,\vv_k,\qq\in\Rn$.
Let $g_0=\fullshad{f}$, and $g_i=\gishadvi$, for $i=1,\dotsc,k$. Since $\xbar\in\arescone{\fullshadfbar}$,
\Cref{thm:astral-cone-char} (applied with $f$, as it appears in that
\namecref{thm:astral-cone-char}, set to $\fullshad{f}$)
implies that $\vv_i\in\resc{g_{i-1}}$ for $i=1,\dotsc,k$,
and $\qq\in\resc{g_k}$.
By induction on $i=0,\ldots,k$, we claim that
$g_i=\fullshad{f}$:
By definition, $g_0=\fullshad{f}$ in the base case.
And for $i>0$, we have
\[
  g_i
  =
  \gishadvi
  =
  \genshad{(\fullshad{f})}{\limray{\vv_i}}
  =
  \fullshad{f},
\]
where the first equality is by definition,
the second is by inductive hypothesis,
and the third is by
\Cref{cor:fullshad:absorb}(\ref{cor:fullshad:absorb:d}).

In particular, this implies
$\vv_i\in\resc{g_{i-1}}=\resc{\fullshad{f}}$
for $i=1,\ldots,k$,
and
$\qq\in\resc{g_k}=\resc{\fullshad{f}}$.
Thus, in fact, $\vv_1,\ldots,\vv_k,\qq\in\resc{\fullshad{f}}$,
so $\xbar\in\represc{\fullshad{f}}$,
completing the proof.%
\indexg{universal reduction|)}%
\indexg{universal reduction!recessive completeness of|)}%
\end{proof}

\section{Coercivity of the universal reduction}
\label{subsec:fullshad-coercive}

Continuing our study of the structure of minimizers
of a convex function's extension,
we focus next on the set of all finite parts of all such
minimizers, which, as we proved in
\Cref{pr:min-fullshad-is-finite-min}, is identical to the set of all
minimizers of the function's universal reduction.

As an initial observation,
note that when minimizing a function $h:\Rn\to\eR$,
we can always disregard directions in which the function remains
constant, that is, in $\conssp{h}$.
Instead, we can restrict our attention exclusively to points in the
linear subspace orthogonal to that space, $L=(\conssp{h})^\perp$.
This is because the rest of the space consists of parallel affine sets
$\yy+L$, for $\yy\in\conssp{h}$, on each of which $h$'s behavior
is simply a copy of its behavior on $L$.

\indexg{universal reduction!minimizers bounded|(}%
Thus, when minimizing~$\fullshad{f}$,
where $f:\Rn\rightarrow\Rext$ is convex,
it suffices to consider only points in
$(\conssp{\fullshad{f}})^\perp$,
which so becomes a focus in the development below.
We will see that the minimization of $\fullshad{f}$ on
this set
is particularly well-behaved
because all of the sublevel sets of $\fullshad{f}$,
when restricted to $(\conssp{\fullshad{f}})^\perp$, are bounded.
As a result, 
the set of all minimizers of $\fullshad{f}$ on
$(\conssp{\fullshad{f}})^\perp$ is a nonempty, convex, compact set $C$.
The full set of minimizers of $\fullshad{f}$
is then equal to $C+\conssp{\fullshad{f}}$, and by
\Cref{pr:min-fullshad-is-finite-min}, this is
also the set of all finite parts of all minimizers of $\ef$.%
\indexg{universal reduction!minimizers bounded|)}

\begin{example}[Two exps and a square, continued]
\label{ex:simple-eg-exp-exp-sq-2}
\indexg{Two exps and a square!universal reduction's minimizers|(}%
  For instance, let $f:\R^{3}\to\eR$ be the function from
  \Cref{ex:simple-eg-exp-exp-sq,ex:simple-eg-exp-exp-sq-part2}.
  In \Cref{ex:simple-eg-exp-exp-sq-part2}, we saw that $f$'s
  universal reduction, for $\xx\in\R^3$, is
  \[
    \fullshad{f}(\xx)=(2+x_2-x_3)^2,
  \]
  with constancy space $\conssp{\fullshad{f}}
  =\set{\xx\in\R^{3}:\:x_2=x_3}$. Setting $\vv=\trans{[0,1,-1]}$,
  the constancy space
  can be rewritten as
  $\conssp{\fullshad{f}}=\set{\xx\in\R^{3}:\:\xx\inprod\vv=0}=\set{\vv}^\perp$.
  Thus,
  \begin{equation}
  \label{eq:simple-eg-exp-exp-sq:cons-perp}
    (\conssp{\fullshad{f}})^\perp
    =\spn\set{\vv}
    =\bigBraces{\trans{[0,\lambda,-\lambda]}:\:\lambda\in\R}
  \end{equation}
  (by \Cref{pr:std-perp-props}\ref{pr:std-perp-props:c}).
  Furthermore, $\fullshad{f}$ has bounded sublevel sets
  when restricted to $(\conssp{\fullshad{f}})^\perp$, because
  \[
    \fullshad{f}\bigParens{\trans{[0,\lambda,-\lambda]}}
    =
    (2+2\lambda)^2
  \]
  for $\lambda\in\R$.
  The only minimizer of $\fullshad{f}$ on $(\conssp{\fullshad{f}})^\perp$ is
  $\trans{[0,-1,1]}$. 
  Therefore, the entire set of minimizers of $\fullshad{f}$
  (and hence the set of all finite parts of
  all minimizers of $\fext$) is
  exactly
  \[
     \trans{[0,-1,1]} + \conssp{\fullshad{f}}
     =
     \bigBraces{ \trans{[\alpha,\beta-1,\beta+1]} : \alpha,\beta \in \R },
  \]
  where the expression on the right is
  because $\conssp{\fullshad{f}}$ is the span of
  $\trans{[1,0,0]}$
  and
  $\trans{[0,1,1]}$ (which are linearly independent and orthogonal to $\vv$,
  so they form a basis for $\conssp{\fullshad{f}}=\set{\vv}^\perp$).%
\indexg{Two exps and a square!universal reduction's minimizers|)}%
  \end{example}

\indexg{coercivity|(}%
\indexg{sublevel sets!all bounded|(}%
Returning to the general case, we next consider
what it means for all of a function's sublevel sets to be compact.
For lower semicontinuous functions,
we show that this is equivalent to a property called coercivity, which
says that the function's values must grow to infinity on every
sequence whose elements' norms grow without bound:

\begin{definition}
\label{def:coercive}
\indexg{coercivity!defined|(}%
Let $f:\Rn\rightarrow\Rext$, and let $L$ be a linear subspace of $\Rn$.
We say that $f$ is \emph{coercive} on $L$ if for every sequence $\seq{\xx_t}$ in $L$,
if $\norm{\xx_t}\to+\infty$, then
\indexg{coercivity!defined|)}%
$f(\xx_t)\to+\infty$.
\end{definition}

\begin{theorem}
\label{thm:coercive:equiv}
Let $f:\Rn\rightarrow\Rext$ be lower semicontinuous,
and let $L$ be a linear subspace of $\Rn$.
Then $f$ is coercive on $L$ if and only if, for all $\beta\in\R$,
the set $\set{\xx\in L : f(\xx)\le \beta}$ is a compact subset of $\Rn$.
\end{theorem}

\begin{proof}
  ~
  \begin{proof-parts}
  \pfpart{``Only if'' ($\Rightarrow$):}
      We prove the contrapositive. Suppose
      there exists $\beta\in\R$ such that the set $S=\set{\xx\in L :\: f(\xx)\le \beta}$
      is not compact.
      The set $S$ is closed in $\Rn$, being
      the intersection of~$L$, which is closed in $\Rn$, and
      $\set{\xx\in\Rn:\: f(\xx)\le\beta}$, which is a sublevel set of a lower semicontinuous
      function and therefore also closed in $\Rn$
      (\Cref{prop:lsc}\ref{prop:lsc:a}\ref{prop:lsc:c}).
      Since $S$ is not compact, it therefore must
      be unbounded
      (by \Cref{pr:compact-in-rn}).
      Hence, there exists a sequence $\seq{\xx_t}$ in $S\subseteq L$ such that
      $\norm{\xx_t}\to+\infty$. However, $f(\xx_t)\le\beta$ for all~$t$, so we cannot have
      $f(\xx_t)\to+\infty$.
      Thus, $f$ is not coercive on $L$.\looseness=-1
  \pfpart{``If'' ($\Leftarrow$):}
      Again we prove the contrapositive. Suppose $f$ is not coercive on $L$.
      Then there exists a sequence $\seq{\xx_t}$ in $L$ such that $\norm{\xx_t}\to+\infty$,
      but $f(\xx_t)\not\rightarrow+\infty$.
This latter condition means that
there exists $\beta\in\R$ such that $f(\xx_t)\leq\beta$ for
infinitely many sequence elements.
By discarding all other elements, we can assume $f(\xx_t)\leq\beta$
for all~$t$.
Consequently, the set $\set{\xx\in L:\: f(\xx)\le\beta}$ includes the
entire sequence $\seq{\xx_t}$ and is therefore unbounded;
thus, it also is not compact
(again by \Cref{pr:compact-in-rn}),
finishing the proof.%
\indexg{sublevel sets!all bounded|)}%
      \qedhere
  \end{proof-parts}
\end{proof}

In a moment, we will show that $\fullshad{f}$ is coercive on
$(\conssp{\fullshad{f}})^\perp$,
implying, by \Cref{thm:coercive:equiv},
that all its minimizers in this linear subspace are included in
some compact region.
For this, we will rely on the following general result:

\begin{theorem}
\label{thm:coercive}
   Let $f:\Rn\to\eR$ be convex and lower semicontinuous. Then $f$
   is coercive on $(\resc{f})^\perp$.
\end{theorem}

Like \Cref{pr:cons-implies-fin-minimizer},
this theorem only concerns concepts from standard convex
analysis, and could be proved using standard techniques
(indeed, closely related results are given, for instance, in
\idxroc\citealp[Section~8]{ROC}).
Here, we instead give a more direct proof as an illustration of astral
techniques applied to standard convex analysis.

\begin{proof}
  For the sake of contradiction, suppose $f$ is not coercive on
  $(\resc{f})^\perp$. Then there exists a sequence $\seq{\xx_t}$ in
  $(\resc{f})^\perp$ such that
  $\norm{\xx_t}\to+\infty$, but $f(\xx_t)\not\rightarrow+\infty$.
  This latter condition implies, for some $\beta\in\R$, that
  $f(\xx_t)\leq\beta$ for infinitely many sequence elements;
  by discarding all other elements, we can assume $f(\xx_t)\leq\beta$
  for all $t$.
  Moreover,
  by sequential compactness of $\eRn$, the sequence $\seq{\xx_t}$ has
  a convergent subsequence;
  by again discarding all other elements, we can assume
  the entire sequence converges to some $\xbar\in\extspace$.
  
  Since $\norm{\xx_t}\to+\infty$,
  $\xbar$ is infinite
  (by \Cref{pr:seq-to-inf-has-inf-len}), so $\xbar=\limray{\vv}\plusl\zbar$ 
  for some $\zbar\in\eRn$, 
  where $\vv$ is $\xbar$'s dominant direction.
  Then
  \[
    \ef(\limray{\vv}\plusl\zbar)
    =
    \ef(\xbar)
    \le
    \liminf f(\xx_t)\le\beta<+\infty,
  \]
  where the first inequality is by $\ef$'s definition since
  $\xx_t\rightarrow\xbar$.
  Therefore, $\vv\in(\aresconef)\cap\Rn=\resc{f}$ (by
  \Cref{thm:rec-ext-equivs}\ref{thm:rec-ext-equivs:e}\ref{thm:rec-ext-equivs:a}
  and
  \Cref{pr:f:1}).
  
  Since $\vv$ is $\xbar$'s dominant direction, we have
  $\xx_t/\norm{\xx_t}\to\vv$ (\Cref{thm:dom-dir}\ref{thm:dom-dir:a}\ref{thm:dom-dir:c}).
  Moreover, for all $t$, $\xx_t$ is in $(\resc{f})^\perp$, which
  is a linear subspace, so
  $\xx_t/\norm{\xx_t}\in(\resc{f})^\perp$, and hence also
  $\lim(\xx_t/\norm{\xx_t})=\vv\in(\resc{f})^\perp$ since $(\resc{f})^\perp$ is closed.
  Thus, $\vv\in(\resc{f})^\perp\cap\resc{f}$, implying, by definition
  of orthogonal complement (Eq.~\ref{eq:std-ortho-comp-defn}),
  that $\vv\cdot\vv=0$,
  so $\vv=\zero$,
  contradicting that
  $\vv$ is $\xbar$'s dominant direction
  (which would mean $\vv\ne\zero$).
This completes the proof.%
\indexg{coercivity|)}%
\end{proof}

\indexg{universal reduction!coercivity of|(}%
\indexg{coercivity!universal reduction@of universal reduction|(}%
Since $\conssp{\fullshad{f}}=\resc{\fullshad{f}}$, we immediately obtain that
$\fullshad{f}$ is coercive on $(\conssp{\fullshad{f}})^\perp$.

\begin{corollary}
\label{cor:fullshad-coercive}
  Let $f:\Rn\rightarrow\Rext$ be convex.
  Then $\fullshad{f}$ is
  coercive on $(\conssp{\fullshad{f}})^\perp=(\resc{\fullshad{f}})^\perp$.
\end{corollary}

\begin{proof}
  By \Cref{cor:fullshad:absorb}(\ref{cor:fullshad:absorb:b},\ref{cor:fullshad:absorb:e}),
  $\fullshad{f}$ is convex and lower semicontinuous,
  while also satisfying $\conssp{\fullshad{f}}=\resc{\fullshad{f}}$.
  By \Cref{thm:coercive}, $\fullshad{f}$ is therefore coercive on
  $(\resc{\fullshad{f}})^\perp=(\conssp{\fullshad{f}})^\perp$.%
\indexg{universal reduction!coercivity of|)}%
\indexg{coercivity!universal reduction@of universal reduction|)}%
\end{proof}

We next relate the set $(\resc{\fullshad{f}})^\perp$,
where $\fullshad{f}$ is coercive,
to properties of $f$ and its extension $\fext$.
\indexg{recession cone, astral!orthogonal complement of|(}%
As a first step,
we show in the next \namecref{thm:f:6}
that the set $(\resc{\ef})^\perp$
is preserved under the operation of taking reduction at any astron associated
with a vector in $\resc{f}$. Since $\fullshad{f}$ is obtained from $f$ by iteratively
taking such reductions, this will imply that $(\resc{\ef})^\perp=(\resc{\fullshadfbar})^\perp$.
We then show that if $f$ is recessive complete, then
$(\resc{\fext})^\perp=(\resc{f})^\perp$; in particular, this applies
to $\fullshad{f}$, by \Cref{pr:fullshad-recess-complete}.
Taken together, these results will imply that
$(\conssp{\fullshad{f}})^\perp=(\resc{\fullshad{f}})^\perp=(\resc{\ef})^\perp$,
in other words, that $(\conssp{\fullshad{f}})^\perp$, where
$\fullshad{f}$'s minimizers are to be found, is the same as
$(\resc{\ef})^\perp$, the orthogonal complement of $f$'s astral
recession cone.

\begin{theorem}  \label{thm:f:6}
  Let $f:\Rn\rightarrow\Rext$ be convex and lower semicontinuous,
  let $\vv\in\resc{f}$, and let
  $g=\fshadv$
  be the reduction of $f$ at $\limray{\vv}$.
  Then $\perpresg=\perpresf$.
\end{theorem}

\begin{proof}
By \Cref{thm:f:3}(\ref{thm:f:3a}),
$\aresconef\subseteq\aresconeg$,
implying
$\perpresg\subseteq\perpresf$
by
\Cref{pr:perp-props-new}(\ref{pr:perp-props-new:b}).
To prove the reverse inclusion,
suppose $\uu\in\perpresf$.
Let $\ybar\in\aresconeg$.
Then $\limray{\vv}\plusl\ybar\in\aresconef$
by \Cref{thm:f:3}(\ref{thm:f:3b}), so
$\limray{\vv}\cdot\uu\plusl\ybar\cdot\uu=(\limray{\vv}\plusl\ybar)\cdot\uu=0$.
This is only possible if $\vv\cdot\uu=0$
and $\ybar\cdot\uu=0$.
Since this holds for all $\ybar\in\aresconeg$, this implies
$\uu\in\perpresg$, proving
$\perpresf\subseteq\perpresg$.
\end{proof}

\begin{proposition}  \label{pr:perpres-is-rescperp}
\indexg{recession cone, astral!standard recession cone and|(}%
\indexg{recession cone (standard)!astral recession cone and|(}%
  Let $f:\Rn\rightarrow\Rext$ be convex and lower semicontinuous.
  Then $\perpresf\subseteq\rescperp{f}$.
  If, in addition,
  $f$ is recessive complete,
  then
  $\perpresf = \rescperp{f}$.
\end{proposition}

\begin{proof}
The
inclusion $\perpresf\subseteq\rescperp{f}$ follows from
\Cref{pr:perp-props-new}(\ref{pr:perp-props-new:b})
since
$\resc{f} \subseteq  \aresconef$
(by \Cref{pr:f:1}).

If, in addition,
$f$ is recessive complete, then $\aresconef=\clbar{\resc{f}}$
(by \Cref{pr:repres-in-arescone});
hence, $(\resc{f})^\perp=(\clbar{\resc{f}})^\perp=(\aresconef)^\perp$,
with the first equality by
\indexg{recession cone, astral!standard recession cone and|)}%
\indexg{recession cone (standard)!astral recession cone and|)}%
\Cref{pr:perp-props-new}(\ref{pr:perp-props-new:c}).
\end{proof}

\indexg{universal reduction!coercivity of|(}%
\indexg{coercivity!universal reduction@of universal reduction|(}%
\indexg{universal reduction!minimizers bounded|(}%
\indexg{universal reduction!sublevel sets bounded|(}%
Combining the above yields that
$\fullshad{f}$ is coercive on
$(\conssp{\fullshad{f}})^\perp =(\resc{\fullshad{f}})^\perp =(\resc{\ef})^\perp$.
In minimizing $\fullshad{f}$,
as discussed earlier, this means we can
effectively restrict attention only to $\perpresf$,
the linear subspace orthogonal to the astral recession cone.
Moreover, within $\perpresf$, all of the sublevel sets of
$\fullshad{f}$ are convex and compact,
implying all such minimizers are restricted
to such a convex, compact region,
and furthermore must themselves comprise such a set.
We summarize these facts next:

\begin{theorem}  \label{thm:f:4x}
    Let $f:\Rn\rightarrow\Rext$ be convex.
    Then:
    \begin{letter-compact}
    \item  \label{thm:f:4xa}
      $\conssp{\fullshad{f}}=\resc{\fullshad{f}}=\rspan(\aresconef)$,
      and also
      $(\conssp{\fullshad{f}})^\perp=(\resc{\fullshad{f}})^\perp=(\aresconef)^\perp$.
    \item  \label{thm:f:4xb}
      $\fullshad{f}$ is coercive on $(\conssp{\fullshad{f}})^\perp=(\aresconef)^\perp$.
      Consequently,
      for all $\beta\in\R$, the set
      \begin{equation}  \label{eqn:thm:f:4x:1}
        \set{ \xx \in \perpresf :\: \fullshad{f}(\xx) \leq \beta }
      \end{equation}
      is convex and compact.
    \item  \label{thm:f:4xd}
      The set of all minimizers of $\fullshad{f}$ on $\perpresf$,
      \begin{equation}  \label{eqn:thm:f:4x:2}
         C
         =
         \set{
           \xx\in\perpresf :\:
           \fullshad{f}(\xx) = \min\fullshad{f}
         },         
      \end{equation}
      is nonempty, convex, and compact.
      The set of all minimizers of $\fullshad{f}$ on all of $\Rn$ is
      then $C + \conssp{\fullshad{f}}$.
    \end{letter-compact}
\end{theorem}

\begin{proof}
  ~

\begin{proof-parts}
\pfpart{Part~(\ref{thm:f:4xa}):}
By \Cref{cor:fig-cons-main-props}, there exist $\vv_1,\dotsc,\vv_k\in\Rn$
such that $\vv_i\in\resc{g_{i-1}}$ and $g_k=\fullshad{f}$,
where $g_0=\lsc f$, and
$g_i=\gishadvi$, for $i=1,\dotsc,k$.
Then
\begin{equation}  \label{eq:fullshad:perpresf:1}
  \perpresf
  =
  \perpresgsub{0}=\dotsb=\perpresgsub{k}
  =
  (\resc{\fullshadfbar})^\perp
  =
  (\resc{\fullshad{f}})^\perp.
\end{equation}
The first equality is because $\fext=\gext_0$ 
(by \Cref{pr:h:1}\ref{pr:h:1aa}).
The last equality is by
\Cref{pr:perpres-is-rescperp} since
$\fullshad{f}$ is recessive complete
(by \Cref{pr:fullshad-recess-complete}).
The middle equalities (that $\perpresgsub{i-1}=\perpresgsub{i}$
for $i=1,\ldots,k$) are by
\Cref{thm:f:6} (noting that each $g_i$ is convex and lower
semicontinuous by \Cref{pr:icon-red-decomp-astron-red}).
    
Moreover, by
\Cref{cor:fullshad:absorb}(\ref{cor:fullshad:absorb:e}),
$\conssp{\fullshad{f}}=\resc{\fullshad{f}}$, which
is a linear subspace of $\Rn$
(\Cref{pr:prelim:const-props}\ref{pr:prelim:const-props:b}),
so
    \begin{equation*}
        \conssp{\fullshad{f}}
      =
      \resc{\fullshad{f}}
      =
      (\resc{\fullshad{f}})^{\perp\perp}
      =
      \perperpresf
      =
      \rspan(\aresconef),
    \end{equation*}
    where the second equality is by \Cref{pr:std-perp-props}(\ref{pr:std-perp-props:c}),
    the third is by \eqref{eq:fullshad:perpresf:1}, and the fourth is by
    \Cref{pr:perp-props-new}(\ref{pr:perp-props-new:c}).
Together with \eqref{eq:fullshad:perpresf:1},
this proves all the parts of the claim.

\pfpart{Part~(\ref{thm:f:4xb}):}
That $\fullshad{f}$ is coercive on
$(\conssp{\fullshad{f}})^\perp=(\aresconef)^\perp$
is immediate from
\Cref{cor:fullshad-coercive} and part~(\ref{thm:f:4xa}).%
\indexg{universal reduction!coercivity of|)}%
\indexg{coercivity!universal reduction@of universal reduction|)}

For the second claim, let $S$ denote the set in
\eqref{eqn:thm:f:4x:1},
for some $\beta\in\R$.
By \Cref{cor:fullshad:absorb}(\ref{cor:fullshad:absorb:b}),
$\fullshad{f}$ is convex and lower semicontinuous.
Therefore, $S$ is compact by \Cref{thm:coercive:equiv}.
Further, $S$ is convex, being the intersection of
$\perpresf$, which is a linear subspace
(by \Cref{pr:perp-props-new}\ref{pr:perp-props-new:a})
and so convex,
with the sublevel set
$\set{\xx\in\Rn : \fullshad{f}(\xx)\leq\beta}$,
which is also convex by \Cref{roc:thm4.6}.

\pfpart{Part~(\ref{thm:f:4xd}):}
We claim first that $C$ is convex and compact.
If $f\equiv+\infty$, then
$\fullshad{f}\equiv+\infty$ (by
\Cref{pr:univ-red-props}\ref{pr:univ-red-props:min}),
implying
$C=\perpresf=(\conssp{\fullshad{f}})^\perp=\{\zero\}$,
which is convex and compact.
So assume $f\not\equiv+\infty$.
Then we can write
\[
   C
   =
   \bigcap_{\scriptontop{\beta\in\R:}{\min\fullshad{f}\leq\beta}}
      \set{ \xx \in \perpresf :\: \fullshad{f}(\xx) \leq \beta }.
\]
Each set on the right is convex and compact, by
part~(\ref{thm:f:4xb}), and so also closed (in $\Rn$)
and bounded (by \Cref{pr:compact-in-rn}).
Therefore, $C$, their intersection, is also convex, closed,
and bounded (noting that the intersection cannot be vacuous since
$\min\fullshad{f}<+\infty$
by \Cref{pr:univ-red-props}\ref{pr:univ-red-props:min}).
Thus, $C$ is convex and compact.

Let
$C'= \set{\xx\in\Rn :\: \fullshad{f}(\xx) = \min\fullshad{f}}$
be the set of all minimizers of $\fullshad{f}$ over all of~$\Rn$.
We claim next that $C'=C+\conssp{\fullshad{f}}$.
Let $\xx\in C+\conssp{\fullshad{f}}$,
meaning $\xx=\qq+\rr$ for some $\qq\in C$
and $\rr\in\conssp{\fullshad{f}}$.
Then
$\fullshad{f}(\xx)=\fullshad{f}(\qq+\rr)=\fullshad{f}(\qq)=\min\fullshad{f}$,
where the second equality is because 
$\rr\in\conssp{\fullshad{f}}$, and the third is because
$\qq\in C$.
Thus, $\xx\in C'$, so
$C+\conssp{\fullshad{f}}\subseteq C'$.

For the reverse inclusion, let $\xx\in C'$.
Then we can write $\xx=\qq+\rr$ for some
$\qq\in(\conssp{\fullshad{f}})^\perp=\perpresf$
and $\rr\in\conssp{\fullshad{f}}$
(by \Cref{pr:lin-decomp}), implying
$\fullshad{f}(\qq)=\fullshad{f}(\qq+\rr)=\fullshad{f}(\xx)=\min\fullshad{f}$,
where the first equality is because 
$\rr\in\conssp{\fullshad{f}}$, and the third is because
$\xx\in C'$.
Thus, $\qq\in C$, so $\xx\in C+\conssp{\fullshad{f}}$.
Hence,
$C'= C+\conssp{\fullshad{f}}$.

By \Cref{pr:univ-red-props}(\ref{pr:univ-red-props:min}),
$\fullshad{f}$ attains its minimum, so
$C'$ is nonempty.
Therefore, since $C'=C+\conssp{\fullshad{f}}$,
and since $\zero\in\conssp{\fullshad{f}}$,
$C$ must also be nonempty.%
\indexg{universal reduction!minimizers bounded|)}%
\indexg{universal reduction!sublevel sets bounded|)}%
\indexg{recession cone, astral!orthogonal complement of|)}%
\qedhere
\end{proof-parts}
\end{proof}

Thus, if $\xbar=\ebar\plusl\qq$ minimizes $\fext$,
where $\ebar\in\corezn$ and $\qq\in\Rn$, then
$\ebar$ must be in $\aresconef$, by
\Cref{thm:arescone-fshadd-min}, and $\qq$ must minimize
$\fullshad{f}$, by
\Cref{pr:min-fullshad-is-finite-min}.
Furthermore, as a consequence of \Cref{thm:f:4x}(\ref{thm:f:4xd}),
$\qq=\qq'+\rr$,
where
$\rr\in\conssp{\fullshad{f}}$,
and
$\qq'$ is an element of the convex and compact subset $C$ of
$\perpresf$, as given in \eqref{eqn:thm:f:4x:2}.

\indexg{universal reduction!minimizing sequence@on minimizing sequence|(}%
Regarding sequences, \Cref{thm:f:4x} implies that if a convex
function $f$ is minimized by some sequence, then
that sequence must
also minimize $\fullshad{f}$, as does also the projection of that
sequence onto the linear subspace $\perpresf$.
Further, that projected sequence cannot be unbounded.

\begin{proposition}  \label{pr:proj-mins-fullshad}
  Let $f:\Rn\rightarrow\Rext$ be convex.
  Let $\seq{\xx_t}$ be a sequence in $\Rn$, and for each $t$, let
  $\qq_t$ be the projection of $\xx_t$ onto the linear subspace
  $\perpresf$.
  Assume $f(\xx_t)\rightarrow\inf f$.
  Then
  $\fullshad{f}(\xx_t)=\fullshad{f}(\qq_t)\rightarrow \min\fullshad{f}$.
  Furthermore, the entire sequence $\seq{\qq_t}$ is included in a
  compact subset of $\Rn$.
\end{proposition}

\begin{proof}
If $f\equiv+\infty$, then $\fullshad{f}\equiv+\infty$ and
$\perpresf=\{\zero\}$, so $\qq_t=\zero$ for all $t$, implying the
claim.
Therefore, we assume henceforth that $f\not\equiv+\infty$.

For all $t$,
\begin{equation*}  %
  \inf f
  =
  \min \fullshad{f}
  \leq
  \fullshad{f}(\qq_t)
  =
  \fullshad{f}(\xx_t)
  \leq
  f(\xx_t).
\end{equation*}
The first equality is by
\Cref{pr:univ-red-props}(\ref{pr:univ-red-props:min}),
the second equality is by
\Cref{pr:cons-PP}
since $(\conssp{\fullshad{f}})^\perp=(\resc{\ef})^\perp$
(by \Cref{thm:f:4x}\ref{thm:f:4xa}),
and the second inequality is by
\Cref{pr:univ-red-props}(\ref{pr:univ-red-props:a}).

Since $f(\xx_t)\rightarrow\inf f$,
this
implies
$\fullshad{f}(\xx_t)=\fullshad{f}(\qq_t) \rightarrow \min \fullshad{f}$,
as claimed.
This also shows,
for any $\beta\in\R$ with $\beta>\min\fullshad{f}$,
that all but finitely many of the $\qq_t\negKern$'s are
in some sublevel set of $\fullshad{f}$, as in
\eqref{eqn:thm:f:4x:1}.
By \Cref{thm:f:4x}(\ref{thm:f:4xb}), every such sublevel set is
compact.
Therefore, there exists a (possibly larger) compact subset of $\Rn$
that includes
the entire sequence $\seq{\qq_t}$.%
\indexg{universal reduction!minimizing sequence@on minimizing sequence|)}%
\end{proof}

\indexg{universal reduction!when same as original function|(}%
For a convex, lower semicontinuous function $f:\Rn\rightarrow\Rext$,
we saw in \Cref{thm:fullshad:char} that
the condition $\conssp{f}=\resc{f}$
characterizes when $f=\fullshad{f}$.
We next provide an additional dual characterization of this property in terms of
the barrier cone of $f$.

In general, assuming $f\not\equiv+\infty$, we have
by \Cref{pr:perp-props-new}(\ref{pr:perp-props-new:d})
and \Cref{cor:ares-is-apolslopes} that
$(\resc{\ef})^\perp\subseteq\rescpol{\ef}=\slopes{f}$.
We show that if the opposite inclusion also holds, we must have $\resc{f}=\conssp{f}$
(and vice versa). Thus, the inclusion $\slopes{f}\subseteq\perpresf$
characterizes when $f=\fullshad{f}$.
Furthermore, if $f$ is
closed,
then this condition simplifies to $\dom{f^*}\subseteq\rescperp{\ef}$.

\begin{theorem}  \label{cor:thm:f:4:1}
    Let $f:\Rn\rightarrow\Rext$ be convex and lower semicontinuous.
    Then the following are equivalent:
    \begin{letter-compact}
    \item  \label{cor:thm:f:4:1c}
      $\fullshad{f} = f$.
    \item  \label{cor:thm:f:4:1a}
      $\resc{f}=\conssp{f}$.
    \item  \label{cor:thm:f:4:1b}
      Either $f\equiv+\infty$ or
      $\slopes{f}\subseteq\perpresf$.
    \end{letter-compact}
    If, in addition, $f$ is
    closed,
    then the above are also equivalent to:
    \begin{letter-compact}[resume]
    \item \label{cor:thm:f:4:1d}
      Either $f\equiv+\infty$ or
      $\dom{f^*}\subseteq\perpresf$.
    \end{letter-compact}
\end{theorem}

\begin{proof}
    If $f\equiv +\infty$,
    then the \namecref{cor:thm:f:4:1}
    holds since then
    $\resc{f}=\conssp{f}=\Rn$,
    and
    $\fullshad{f}=f$.
    We therefore assume henceforth that
    $f\not\equiv +\infty$.
  \begin{proof-parts}
  \pfpart{(\ref{cor:thm:f:4:1a}) $\Rightarrow$ (\ref{cor:thm:f:4:1c}):}
    This was proved in \Cref{thm:fullshad:char}.

  \pfpart{(\ref{cor:thm:f:4:1c}) $\Rightarrow$ (\ref{cor:thm:f:4:1b}):}
    Suppose $\fullshad{f}=f$. Then $\resc{f}=\conssp{f}$
    (by \Cref{thm:fullshad:char}), so $\resc{f}$ is a linear subspace
    of $\Rn$
    (\Cref{pr:prelim:const-props}\ref{pr:prelim:const-props:b}).
    Since
    $f=\fullshad{f}$ is recessive complete
    (\Cref{pr:fullshad-recess-complete}),
    we obtain
    $\resc{\ef}=\clbar{\resc{f}}$
    (by \Cref{pr:repres-in-arescone}), which
    is an astral linear subspace
    (by \Cref{cor:std-ast-linsub-corr}\ref{cor:std-ast-linsub-corr:L}).
    Thus, $\slopes{f}=\rescpol{\ef}=\perpresf$,
    where the first equality is by \Cref{cor:ares-is-apolslopes},
    and the second by
    \Cref{pr:perp-props-new}(\ref{pr:perp-props-new:d}).

\pfpart{(\ref{cor:thm:f:4:1b}) $\Rightarrow$ (\ref{cor:thm:f:4:1a}):}
Suppose $\slopes{f}\subseteq\perpresf$.
Then
\[
  \aresconef
  \subseteq
  \aspan(\aresconef)
  =
  \perpaperp{(\aresconef)}
  \subseteq
  \apol{\Parens{(\aresconef)^\perp}}
  \subseteq
  \apolslopesf  
  =
  \aresconef.
\]
The first equality is by \Cref{thm:dub-perp}(\ref{thm:dub-perp:a}).
The second inclusion is by
\Cref{pr:aperp-props}(\ref{pr:aperp-props:c}).
The third inclusion is by
\Cref{pr:ast-pol-props}(\ref{pr:ast-pol-props:b})
since $\slopes{f}\subseteq\perpresf$.
The last equality is by \Cref{cor:ares-is-apolslopes}.
Thus, $\aresconef$ is its own astral span, and so is
an astral linear subspace.
Hence, $\resc{f}=(\aresconef)\cap\Rn$ is
a linear subspace of $\Rn$
(using \Cref{pr:f:1} and
\Cref{cor:std-ast-linsub-corr}\ref{cor:std-ast-linsub-corr:M}), and so
is closed under negation.
Thus, $\conssp{f}=(\resc{f})\cap(-\resc{f})=\resc{f}$
(with first equality by
\Cref{pr:prelim:const-props}\ref{pr:prelim:const-props:a}).

\pfpart{%
  (\ref{cor:thm:f:4:1b})
  $\Leftrightarrow$
  (\ref{cor:thm:f:4:1d})
  if $f$ is closed:
}
Suppose $f$ is closed.
We show that
(\ref{cor:thm:f:4:1b})~$\Leftrightarrow$~(\ref{cor:thm:f:4:1d}).

If $\slopes{f}\subseteq\perpresf$, then
    \[
      \dom{f^*}\subseteq\cone(\dom{f^*})\subseteq\slopes{f}\subseteq\perpresf,
    \]
    where the second inclusion is by
    \Cref{thm:slopes-equiv}.

    Conversely, if $\dom{f^*}\subseteq\perpresf$,
    then
\[
  \slopes{f}
  =
  \polar{(\aresconef)}
  \subseteq
  \polar{(\resc{f})}
  =
  \cl\bigParens{\cone(\dom{\fstar})}
  \subseteq
  \perpresf.
\]  
The first equality is by 
\Cref{cor:ares-is-apolslopes}.
The first inclusion is by 
\Cref{pr:ext-pol-cone-props}(\ref{pr:ext-pol-cone-props:d})
since $\resc{f}\subseteq\aresconef$
by \Cref{pr:f:1}.
The second equality is by \Cref{pr:rescpol-is-con-dom-fstar}.
The last inclusion is because $\perpresf$ is a linear subspace
(\Cref{pr:perp-props-new}\ref{pr:perp-props-new:a}), and so is a closed
(in $\Rn$) convex cone; since this space includes $\dom{\fstar}$, it
therefore also includes $\cone(\dom{f^*})$,
and so too $\cl\regParens{\cone(\dom{\fstar})}$.%
\indexg{universal reduction!when same as original function|)}%
\qedhere
\end{proof-parts}    
\end{proof}

\section{Universal reducers and canonical minimizers}
\label{subsec:univ-min}

\indexg{universal reducers|(}%
From \Cref{cor:fullshad:absorb}(\ref{cor:fullshad:absorb:a}),
we know that the universal reduction
$\fullshad{f}$
of a convex function $f:\Rn\rightarrow\Rext$
can always be expressed as an iconic reduction
$\fshadd$
for some icon $\ebar$.
By \Cref{pr:fullshad-equivs},
this means, for all $\xx\in\Rn$, that
\begin{equation*}  %
\fext(\ebar\plusl\xx)
=
\fshadd(\xx)
=
\fullshad{f}(\xx)
=
\min_{\ebar'\in\corezn} \fext(\ebar'\plusl\xx).
\end{equation*}
Thus, $\ebar$ realizes
the minimum that appears on the right
for \emph{all} $\xx$ \emph{simultaneously}.
Furthermore, since $\ebar$ arises from \Cref{fig:min-proc},
there is a whole set of
points with this same property, since the construction and proof hold
for a whole range of nondeterministic choices in \Cref{fig:min-proc}.
Here, we study some of the properties of that set.

\begin{definition}  \label{def:univ-reducer}
\indexg{universal reducers!defined|(}%
Let $f:\Rn\rightarrow\Rext$ be convex.
We say that an icon $\ebar\in\corezn$ is a
\emph{universal reducer} for $f$
if $\fshadd=\fullshad{f}$, that is, if
$ \fext(\ebar\plusl\xx)=\fullshad{f}(\xx) $
for all $\xx\in\Rn$.
We write $\unimin{f}$ for the set of all such universal reducers:
\[
\indexm{univ f}{$\unimin{f}$}{set of universal reducers}%
  \unimin{f}
    = \BigBraces{ \ebar\in\corezn :\: \fshadd=\fullshad{f} }.%
\indexg{universal reducers!defined|)}%
\]
\end{definition}
We call such icons ``universal'' because they give rise to
universal reductions.
The next proposition gives some simple properties of $\unimin{f}$:

\begin{proposition}  \label{pr:new:thm:f:8a}
  Let $f:\Rn\rightarrow\Rext$ be convex.
  Then:
  \begin{letter-compact}
  \item     \label{pr:new:thm:f:8a:nonemp-closed}
    $\unimin{f}$ is nonempty and closed (in $\extspace$).
  \item     \label{pr:new:thm:f:8a:a}
    $\unimin{f}\subseteq\conv(\unimin{f})\subseteq\aresconef$.
  \end{letter-compact}
\end{proposition}

\begin{proof}
~

\begin{proof-parts}
\pfpart{Part~(\ref{pr:new:thm:f:8a:nonemp-closed}):}
That $\unimin{f}$ is nonempty follows directly from
\Cref{cor:fullshad:absorb}(\ref{cor:fullshad:absorb:a}).

Let $\ebar\in\clbar{\unimin{f}}$, which we aim to show is in
$\unimin{f}$.
Then there exists a sequence
$\seq{\ebar_t}$ in $\unimin{f}$ which converges to $\ebar$.
This implies $\ebar$ is an icon since $\corezn$ is closed in $\extspace$
(\Cref{pr:i:8}\ref{pr:i:8e}).
Also,
for all $\xx\in\Rn$,
\[
  \fullshad{f}(\xx)
  \leq
  \fext(\ebar\plusl\xx)
  \leq
  \liminf \fext(\ebar_t\plusl\xx)
  =
  \fullshad{f}(\xx).
\]
The first inequality is by definition of $\fullshad{f}$
(\Cref{def:univ-reduction}).
The second inequality is by
lower semicontinuity of $\fext$
(\Cref{prop:ext:F}\ref{prop:ext:F:a}),
and because
$\ebar_t\plusl\xx\rightarrow\ebar\plusl\xx$
(by \Cref{pr:i:7}\ref{pr:i:7g}).
The equality is because each $\ebar_t\in\unimin{f}$.
Therefore, $\ebar\in\unimin{f}$,
so $\unimin{f}$ is closed.

\pfpart{Part~(\ref{pr:new:thm:f:8a:a}):}
We argue that $\unimin{f}\subseteq\aresconef$.
Since $\aresconef$ is convex
(by \Cref{cor:res-fbar-closed}), it will then follow that
$\conv(\unimin{f})\subseteq\aresconef$
(by \Cref{pr:conhull-prop}\ref{pr:conhull-prop:aa}).

Let $\ebar\in\unimin{f}$.
Then for all $\xx\in\Rn$,
$\fext(\ebar\plusl\xx)=\fullshad{f}(\xx)\leq f(\xx)$,
where the equality is because $\ebar$ is a universal reducer,
and the inequality is by
\Cref{pr:univ-red-props}(\ref{pr:univ-red-props:a}).
Therefore, $\ebar\in\aresconef$ by
\Crefequiv{thm:rec-ext-equivs}{thm:rec-ext-equivs:b}{thm:rec-ext-equivs:a}.
\qedhere
\end{proof-parts}
\end{proof}

Suppose $\xbar=\ebar\plusl\qq$ where $\ebar\in\corezn$ and
$\qq\in\Rn$.
Previously, we saw that if $\xbar$ minimizes $\fext$, then
$\ebar$ must be in $\aresconef$
(by \Cref{thm:arescone-fshadd-min}),
and $\qq$ must be a minimizer of $\fullshad{f}$
(by \Cref{pr:min-fullshad-is-finite-min}).
The converse is false.
In other words, in general,
it is not the case that $\xbar$ minimizes $\fext$
for every choice of $\ebar\in\aresconef$ and every $\qq$ that minimizes
$\fullshad{f}$:

\begin{example}[Two exps and a square, continued]
\label{ex:simple-eg-exp-exp-sq-not-min}
\indexg{Two exps and a square!nonminimizer counterexample|(}%
As in Examples~\ref{ex:simple-eg-exp-exp-sq}
and~\ref{ex:simple-eg-exp-exp-sq-2}, consider the function
$
  f(\xx) = \me^{x_3-x_1} + \me^{-x_2} + (2+x_2-x_3)^2
$
over $\xx\in\R^3$. In those examples, we derived
\begin{gather*}
\SwapAboveDisplaySkip
  \fullshad{f}(\xx)
    =(2+x_2-x_3)^2,
\\
  \resc{\ef}
    =
    \regBraces{
      \zbar\in\R^3 :\:
      \zbar\cdot (\ee_3-\ee_1)\le 0,\;
      \zbar\cdot (-\ee_2)\le 0,\;
      \zbar\cdot (\ee_2-\ee_3) = 0
    }.
\end{gather*}
Thus,
$\qq=\trans{[0,0,2]}$ minimizes $\fullshad{f}$,
and
$\ebar=\limray{\ee_1}\in\aresconef$,
but
$\xbar=\ebar\plusl\qq$ does not minimize $\fext$ since
$\fext(\xbar)=1>0=\inf f$.%
\indexg{Two exps and a square!nonminimizer counterexample|)}%
\end{example}

\indexg{canonical minimizers|(}%
Nevertheless,
\Cref{cor:fig-cons-main-props} shows that
\Cref{fig:min-proc} yields a minimizer
$\xbar=\ebar\plusl\qq$ of a particular form, namely, with
$\qq$ a finite minimizer of $\fullshad{f}$,
and icon $\ebar$ not only in $\aresconef$, but also a universal
reducer, as shown in part~(\ref{cor:fig-cons-main-props:a}) of that
\namecref{cor:fig-cons-main-props}.
Such points are called canonical minimizers:

\begin{definition}   \label{dfn:canon-minimizer}
\indexg{canonical minimizers!defined|(}%
  Let $f:\Rn\rightarrow\Rext$ be convex.
  We say that a point $\xbar\in\extspace$ is a
  \emph{canonical minimizer} of $\fext$ if $\xbar=\ebar\plusl\qq$ for
  some universal reducer $\ebar\in\unimin{f}$ and some $\qq\in\Rn$
  that minimizes $\fullshad{f}$.%
\indexg{canonical minimizers!defined|)}%
\end{definition}

Every such point is indeed a minimizer of $\fext$,
as follows from the next proposition.
Later, we will see that
\Cref{fig:min-proc} finds all of the canonical minimizers
(and thereby all of the universal reducers as well).

\begin{proposition}  \label{pr:unimin-to-global-min}
  Let $f:\Rn\rightarrow\Rext$ be convex.
  Suppose $\ebar\in\unimin{f}$ and that $\qq\in\Rn$ minimizes
  $\fullshad{f}$.
  Then $\ebar\plusl\qq$ minimizes $\fext$.
\end{proposition}

\begin{proof}
We have
\[
   \fext(\ebar\plusl\qq)
   =
   \fullshad{f}(\qq)
   =
   \min\fullshad{f}
   =
   \min\fext,
\]
where the first equality is because $\ebar\in\unimin{f}$,
the second is because $\qq$ minimizes $\fullshad{f}$,
and the third is by
\Cref{pr:univ-red-props}(\ref{pr:univ-red-props:min}).
Therefore, $\ebar\plusl\qq$ also minimizes $\fext$.
\end{proof}

Not all minimizers of $\fext$ are canonical minimizers.
\indexg{Product of hyperbolas!canonical minimizers and|(}%
For example, for the product of hyperbolas function $f$
from \Cref{ex:recip-fcn-eg,ex:recip-fcn-eg-cont},
$\fext$ is minimized by $\limray{\ee_1}\plusl\ee_2$,
but $\limray{\ee_1}$ is not a universal reducer,
since, for instance,
$\fext(\limray{\ee_1}\plusl(-\ee_2))=+\infty$,
but
$\fullshad{f}\equiv 0$, as can be seen by noting that
$\ef\ge 0$ and $\fshad{\limray{\ee_1}\plusl\limray{\ee_2}}\equiv 0$,
so $\fullshad{f}=\fshad{\limray{\ee_1}\plusl\limray{\ee_2}}\equiv 0$.%
\indexg{Product of hyperbolas!canonical minimizers and|)}

However, by
\Cref{cor:fig-cons-main-props},
every minimizer returned by \Cref{fig:min-proc} is a canonical minimizer.
\indexg{canonical minimizers!generation of all|(}%
\indexg{universal reducers!generation of all|(}%
We next show
that in fact all canonical minimizers, and so also all universal reducers,
can be generated by \Cref{fig:min-proc}.

\begin{theorem}  \label{thm:min-proc-all-can-min}
  Let $f:\Rn\rightarrow\Rext$ be convex.
  Let $\vv_1,\dotsc,\vv_k,\linebreak[0]\qq\in\Rn$,
  and let $\ebar = \limrays{\vv_1,\ldots,\vv_k}$.
  Let $g_0=\lsc f$, and $g_i = \gishadvi$ for $i=1,\ldots,k$.
  Then the following are equivalent:
  \begin{letter-compact}
  \item  \label{thm:min-proc-all-can-min:a}
    $\ebar\in\unimin{f}$ and $\qq$ minimizes $\fullshad{f}$
    (so that $\ebar\plusl\qq$ is a canonical minimizer of $\fext$).
  \item  \label{thm:min-proc-all-can-min:b}
    $\resc{g_k}=\conssp{g_k}$;
    $\qq$ minimizes $g_k$;
    and
    $\vv_i\in\resc{g_{i-1}}$ for $i=1,\ldots,k$.
  \end{letter-compact}
\end{theorem}

\begin{proof}
~

\begin{proof-parts}
\pfpart{%
  (\ref{thm:min-proc-all-can-min:a})
  $\Rightarrow$
  (\ref{thm:min-proc-all-can-min:b}):
}
Suppose $\ebar\in\unimin{f}$ and
$\qq$ minimizes $\fullshad{f}$.
Then $\ebar\plusl\qq$ minimizes $\fext$, by
\Cref{pr:unimin-to-global-min}.
Therefore,
by \Cref{thm:all-min-proc-correct},
$\qq$ minimizes $g_k$,
and
$\vv_i\in\resc{g_{i-1}}$ for $i=1,\ldots,k$.

Also, $\fshadd=\fullshad{f}$ (since $\ebar\in\unimin{f}$),
and $g_k=\fshadd$
by \Cref{pr:icon-red-decomp-astron-red}.
Thus, $g_k=\fullshad{f}$, so
$\resc{g_k}=\conssp{g_k}$
by \Cref{cor:fullshad:absorb}(\ref{cor:fullshad:absorb:e}).

\pfpart{%
  (\ref{thm:min-proc-all-can-min:b})
  $\Rightarrow$
  (\ref{thm:min-proc-all-can-min:a}):
}
Suppose
$\resc{g_k}=\conssp{g_k}$,
$\qq$ minimizes $g_k$,
and
$\vv_i\in\resc{g_{i-1}}$ for $i=1,\ldots,k$.
These conditions imply $g_k=\fshadd=\fullshad{f}$,
by \Cref{cor:fig-cons-main-props}(\ref{cor:fig-cons-main-props:a}).
Therefore, $\ebar\in\unimin{f}$.
This also shows $\qq$ minimizes $\fullshad{f}$ since it minimizes
$g_k$.%
\indexg{canonical minimizers!generation of all|)}%
\indexg{universal reducers!generation of all|)}%
\indexg{canonical minimizers|)}%
\qedhere
\end{proof-parts}
\end{proof}

\indexg{universal reducers!recursive characterization|(}%
We show next how
the set $\unimin{f}$ can be characterized recursively,
much like the characterization of the astral recession cone
in \Cref{thm:f:3}(\ref{thm:f:3b}). For the astral recession cone,
we saw that
the recursion ``bottoms out'' with vectors
$\zz\in\resc{f}=(\resc{\ef})\cap\Rn$ (\Cref{pr:f:1}).
For universal reducers,
the recursion instead bottoms out with the icon $\ebar=\zero$,
which is the only icon in $\Rn$, and which is in $\unimin{f}$
if and only if $f=\fullshad{f}$, as shown below.

\begin{theorem}  \label{thm:new:f:8.2}
    Let $f:\Rn\rightarrow\Rext$ be convex and lower semicontinuous.
    Let $\vv\in\Rn$
    and let
    $g=\fshadv$
    be the reduction of $f$ at $\limray{\vv}$.
    Suppose $\dbar=\limray{\vv}\plusl\ebar$ for some $\ebar\in\corezn$.
    Then $\dbar\in\unimin{f}$ if and only if
    $\vv\in\resc{f}$ and
    $\ebar\in\unimin{g}$.
  \end{theorem}
  
  \begin{proof}
  Suppose $\dbar\in\unimin{f}$.
  Then $\dbar\in\aresconef$ by
  \Cref{pr:new:thm:f:8a}(\ref{pr:new:thm:f:8a:a}),
  so $\vv\in\resc{f}$ by \Cref{thm:f:3}(\ref{thm:f:3b}).
  For all $\xx\in\Rn$,
  \[
    \gext(\ebar\plusl\xx)
    =
    \fext(\limray{\vv}\plusl\ebar\plusl\xx)
    =
    \fext(\dbar\plusl\xx)
    =
    \fullshad{f}(\xx)
    =
    \fullshad{g}(\xx).
  \]
  The first equality is by \Cref{thm:d4}(\ref{thm:d4:c}).
  The third is because $\dbar\in\unimin{f}$.
  The last equality is by \Cref{thm:f:5}.
  Thus, $\ebar\in\unimin{g}$.
  
  The converse is similar.
  Suppose that $\vv\in\resc{f}$ and $\ebar\in\unimin{g}$.
  Then for all $\xx\in\Rn$,
  \[
    \fext(\dbar\plusl\xx)
    =
    \fext(\limray{\vv}\plusl\ebar\plusl\xx)
    =
    \gext(\ebar\plusl\xx)
    =
    \fullshad{g}(\xx)
    =
    \fullshad{f}(\xx).
  \]
  The second equality is by \Cref{thm:d4}(\ref{thm:d4:c}).
  The third is because $\ebar\in\unimin{g}$.
  The last equality is by \Cref{thm:f:5}.
  Thus, $\dbar\in\unimin{f}$.%
\indexg{universal reducers!recursive characterization|)}%
\end{proof}

\begin{proposition}    \label{pr:zero-in-univf}
\indexg{universal reduction!when same as original function|(}%
  Let $f:\Rn\rightarrow\Rext$ be convex and lower semicontinuous.
  Then the following are equivalent:
  \begin{letter-compact}
  \item     \label{pr:zero-in-univf:a}
    $\zero\in\unimin{f}$.
  \item     \label{pr:zero-in-univf:b}
    $f=\fullshad{f}$.
  \item     \label{pr:zero-in-univf:c}
    $\resc{f}=\conssp{f}$.
  \end{letter-compact}
\end{proposition}

\begin{proof}
~

\begin{proof-parts}
\pfpart{%
  (\ref{pr:zero-in-univf:a})
  $\Leftrightarrow$
  (\ref{pr:zero-in-univf:b}):
}
For $\xx\in\Rn$, $\fshad{\zero}(\xx)=\fext(\xx)=f(\xx)$,
by \Cref{pr:h:1}(\ref{pr:h:1a}).
Thus, $\fullshad{f}=f$ if and only if $\fullshad{f}=\fshad{\zero}$,
and so if and only if $\zero\in\unimin{f}$.

\pfpart{%
  (\ref{pr:zero-in-univf:b})
  $\Leftrightarrow$
  (\ref{pr:zero-in-univf:c}):
}
This is immediate from
\Cref{thm:fullshad:char}.%
\indexg{universal reduction!when same as original function|)}%
\qedhere
\end{proof-parts}
\end{proof}

We saw in
\Cref{pr:new:thm:f:8a}(\ref{pr:new:thm:f:8a:a})
that
$\unimin{f}\subseteq\aresconef$.
In fact, there is a much more precise relationship that exists between
these two sets.
\indexg{universal reducers!astral recession cone generated by|(}%
\indexg{recession cone, astral!conic hull of universal reducers@as conic hull of universal reducers|(}%
Specifically, the astral recession cone, $\aresconef$, is
the astral conic hull of the set $\unimin{f}$ of universal reducers,
or equivalently,
the convex hull of $\unimin{f}$ adjoined with the origin, as stated in
the next theorem.
Thus, in this way,
the astral recession cone is
generated by the universal reducers.

\begin{theorem}  \label{thm:res-convhull-unimin}
  Let $f:\Rn\rightarrow\Rext$ be convex.
  Then
  \begin{equation}     \label{eq:thm:res-convhull-unimin:1}
    \aresconef
    =
    \acone(\unimin{f})
    =
    \conv\bigParens{(\unimin{f}) \cup \{\zero\}}.
  \end{equation}
\end{theorem}

\begin{proof}
Let $K=\acone(\unimin{f})$.
We first argue that
\begin{equation}  \label{eqn:thm:res-convhull-unimin:a}
  (\aresconef)\cap\corezn   \subseteq  K.
\end{equation}
Let $\ebar\in(\aresconef)\cap\corezn$, which we aim to show is in
$K$.
Let $\dbar$ be any point in $\unimin{f}$
(which exists by
\Cref{pr:new:thm:f:8a}\ref{pr:new:thm:f:8a:nonemp-closed}).
Then for all $\xx\in\Rn$,
\[
  \fullshad{f}(\xx)
  \leq
  \fext(\ebar\plusl\dbar\plusl\xx)
  \leq
  \fext(\dbar\plusl\xx)
  =
  \fullshad{f}(\xx).
\]
The first inequality is by defintion of $\fullshad{f}$
(\Cref{def:univ-reduction}).
The second is by
\Cref{pr:arescone-def-ez-cons}(\ref{pr:arescone-def-ez-cons:b})
since $\ebar\in\aresconef$.
The equality is because $\dbar\in\unimin{f}$.
Thus, $\ebar\plusl\dbar$
(which is an icon
by \Cref{pr:i:8}\ref{pr:i:8-leftsum})
is a universal reducer,
so $\ebar\plusl\dbar\in\unimin{f}\subseteq K$.
We then have that
\[
  \ebar
  \in
  \lb{\zero}{\ebar\plusl\dbar}
  =
  \aray(\ebar\plusl\dbar)
  \subseteq
  K.
\]
The first inclusion is by
\Cref{cor:d-in-lb-0-dplusx}.
The equality is by
\Cref{thm:oconichull-equals-ocvxhull}
(and since
$\ebar\plusl\dbar$ is an icon).
The last inclusion is because
$K$
is an astral cone that
includes $\ebar\plusl\dbar$.

Having proved \eqref{eqn:thm:res-convhull-unimin:a},
we now have
\[
  \aresconef
  =
  \conv\bigParens{(\aresconef) \cap \corezn}
  \subseteq
  K
  =
  \acone(\unimin{f})
  \subseteq
  \aresconef.
\]
The first equality is by
\Crefequiv{thm:ast-cvx-cone-equiv}{thm:ast-cvx-cone-equiv:a}{thm:ast-cvx-cone-equiv:c}
since $\aresconef$ is a convex astral cone
(\Cref{cor:res-fbar-closed}).
The first inclusion is by
\eqref{eqn:thm:res-convhull-unimin:a},
and by \Cref{pr:conhull-prop}\ref{pr:conhull-prop:aa}
since $K$ is convex.
The last inclusion is because, by
\Cref{pr:new:thm:f:8a}(\ref{pr:new:thm:f:8a:a}),
$\unimin{f}$ is included in the convex astral cone
$\aresconef$
(and using
\Cref{pr:acone-hull-props}\ref{pr:acone-hull-props:a}).
This proves the first equality of
\eqref{eq:thm:res-convhull-unimin:1}.
The second equality then follows
from \Cref{thm:acone-char}(\ref{thm:acone-char:a})
since $\unimin{f}$ consists only of icons.%
\indexg{universal reducers!astral recession cone generated by|)}%
\indexg{recession cone, astral!conic hull of universal reducers@as conic hull of universal reducers|)}%
\end{proof}

\Cref{thm:res-convhull-unimin} shows that the convex hull of
$\unimin{f}$, adjoined with the origin, yields the astral recession
cone.
\indexg{universal reducers!convex hull of|(}%
We next characterize the convex hull of just
$\unimin{f}$, without adjoining the origin.
This set turns out to consist of all points $\zbar$ in $\extspace$,
not just the icons, for which
$\fext(\zbar\plusl\xx)=\fullshad{f}(\xx)$ for all $\xx\in\Rn$.
Equivalently, it consists of all points whose iconic part is a
universal reducer and whose finite part is in the constancy space of
$\fullshad{f}$.

\begin{theorem}   \label{thm:conv-univ-equiv}
  Let $f:\Rn\rightarrow\Rext$ be convex, and let
  $\zbar=\ebar\plusl\qq$ where $\ebar\in\corezn$ and $\qq\in\Rn$.
  Then the following are equivalent:
  \begin{letter-compact}
  \item    \label{thm:conv-univ-equiv:a}
    $\zbar\in\conv(\unimin{f})$.
  \item    \label{thm:conv-univ-equiv:b}
    For all $\xx\in\Rn$, $\fext(\zbar\plusl\xx)\leq \fullshad{f}(\xx)$.
  \item    \label{thm:conv-univ-equiv:c}
    For all $\xx\in\Rn$, $\fext(\zbar\plusl\xx)= \fullshad{f}(\xx)$.
  \item    \label{thm:conv-univ-equiv:d}
    $\ebar\in\unimin{f}$ and $\qq\in\conssp{\fullshad{f}}$.
  \end{letter-compact}
\end{theorem}

\begin{proof}
~

\begin{proof-parts}
\pfpart{%
  (\ref{thm:conv-univ-equiv:a})
  $\Rightarrow$
  (\ref{thm:conv-univ-equiv:b}):
}
Let $\zbar\in\conv(\unimin{f})$,
let $\xx\in\Rn$, and let $G:\extspace\rightarrow\Rext$ be defined by
$G(\ybar)=\fext(\ybar\plusl\xx)$ for $\ybar\in\extspace$.
Then $G$ is convex by \Cref{thm:cvx-compose-affine-cvx}
since $\fext$ is convex
(\Cref{thm:fext-convex}).
Therefore, the sublevel set
$S=\set{\ybar\in\extspace :\: G(\ybar)\leq \fullshad{f}(\xx)}$
is convex by \Cref{thm:f:9}.
Further, $\unimin{f}\subseteq S$ (since
$G(\ebar)=\fext(\ebar\plusl\xx)=\fullshad{f}(\xx)$
for $\ebar\in\unimin{f}$), implying
$\conv(\unimin{f})\subseteq S$
(\Cref{pr:conhull-prop}\ref{pr:conhull-prop:aa}).
In particular, $\zbar\in S$, so
$\fext(\zbar\plusl\xx)=G(\zbar)\leq \fullshad{f}(\xx)$.

\pfpart{%
  (\ref{thm:conv-univ-equiv:b})
  $\Rightarrow$
  (\ref{thm:conv-univ-equiv:d}):
}
Suppose (\ref{thm:conv-univ-equiv:b}) holds.
Then $\ef(\zbar\plusl\xx)\leq \fullshad{f}(\xx)\le f(\xx)$
for all $\xx\in\Rn$
(with second inequality by
\Cref{pr:univ-red-props}\ref{pr:univ-red-props:a}).
Thus, $\zbar\in\aresconef$
by \Cref{thm:rec-ext-equivs}(\ref{thm:rec-ext-equivs:b},\ref{thm:rec-ext-equivs:a}).
Therefore,
\[
  \qq\in\rspanset{\zbar}\subseteq\rspan(\aresconef)=\conssp{\fullshad{f}},
\]
with equality
by \Cref{thm:f:4x}(\ref{thm:f:4xa}).
Further, for all $\xx\in\Rn$,
\[
  \fullshad{f}(\xx)
  \le
  \ef\bigParens{\ebar\plusl\xx}
  =
  \ef\bigParens{\ebar\plusl\qq\plusl(\xx-\qq)}
  =
  \ef\bigParens{\zbar\plusl(\xx-\qq)}
  \le
  \fullshad{f}\bigParens{\xx-\qq}
  =
  \fullshad{f}(\xx).
\]
The first inequality is from the definition of $\fullshad{f}$,
the second inequality is by assumption, and the final equality is
because $\qq\in\conssp{\fullshad{f}}$. 
Thus, $\fshadd=\fullshad{f}$,
so $\ebar\in\unimin{f}$.

\pfpart{%
  (\ref{thm:conv-univ-equiv:d})
  $\Rightarrow$
  (\ref{thm:conv-univ-equiv:a}):
}
Suppose $\ebar\in\unimin{f}$ and $\qq\in\conssp{\fullshad{f}}$.
Then
\[
  \fshad{\ebar\plusl\limray{\qq}}
  =
  \genshad{(\fshad{\ebar})}{\limray{\qq}}
  =
  \genshad{(\fullshad{f})}{\limray{\qq}}
  =
  \fullshad{f}.
\]
The first equality is by
\Cref{pr:j:1}(\ref{pr:j:1a}),
and the second because $\ebar\in\unimin{f}$.
The third equality is by
\Cref{cor:fullshad:absorb}(\ref{cor:fullshad:absorb:e},\ref{cor:fullshad:absorb:d}),
since
$\qq\in\conssp{\fullshad{f}}$.
Thus, $\ebar\plusl\limray{\qq}$ is also a universal reducer.
Therefore,
\[
  \zbar
  =
  \ebar\plusl\qq
  \in
  \ebar\plusl \lb{\zero}{\limray{\qq}}
  =
  \lb{\ebar}{\ebar\plusl\limray{\qq}}
  \subseteq
  \conv(\unimin{f}).
\]
The first inclusion is by
\Cref{thm:lb-with-zero}.
The equality is by \Cref{thm:e:9}.
The last inclusion is because $\ebar$ and $\ebar\plusl\limray{\qq}$
are both in $\unimin{f}$, and so also in its convex hull.

\pfpart{%
  (\ref{thm:conv-univ-equiv:c})
  $\Rightarrow$
  (\ref{thm:conv-univ-equiv:b}):
}
This is immediate.

\pfpart{%
  (\ref{thm:conv-univ-equiv:d})
  $\Rightarrow$
  (\ref{thm:conv-univ-equiv:c}):
}
Suppose $\ebar\in\unimin{f}$ and $\qq\in\conssp{\fullshad{f}}$.
Then for all $\xx\in\Rn$,
\[
  \fext(\zbar\plusl\xx)
  =
  \fext(\ebar\plusl\qq\plusl\xx)
  =
  \fullshad{f}(\qq+\xx)
  =
  \fullshad{f}(\xx),
\]
where the second equality is because $\ebar\in\unimin{f}$,
and the last because $\qq\in\conssp{\fullshad{f}}$.
\qedhere
\end{proof-parts}
\end{proof}

As corollary, the only icons in $\conv(\unimin{f})$ are the universal
reducers:

\begin{corollary}  \label{cor:conv-univ-icons}
  Let $f:\Rn\rightarrow\Rext$ be convex.
  Then
  $\unimin{f}=\conv(\unimin{f})\cap\corezn$.
\end{corollary}

\begin{proof}
That
$\unimin{f}\subseteq\conv(\unimin{f})\cap\corezn$
is immediate.
The reverse inclusion follows because if
$\ebar\in\conv(\unimin{f})\cap\corezn$,
then
$\fshadd=\fullshad{f}$ by
\Crefequiv{thm:conv-univ-equiv}{thm:conv-univ-equiv:a}{thm:conv-univ-equiv:c}.%
\indexg{universal reducers|)}%
\indexg{universal reducers!convex hull of|)}%
\end{proof}

\chapter{The structure of minimizers in some particular cases}
\chaptermark{The structure of minimizers in some cases}
\label{sec:minimizers-examples}

We next study properties of minimizers of an extension $\fext$ in
some particular cases, focusing on the astral rank of minimizers, and
also on a natural class of minimization problems commonly encountered in
statistics and machine learning.

\section{Astral rank of minimizers}
\label{sec:max-rank-minimizers}
\label{subsec:rank-one-minimizers}

When a convex function has no finite minimizer, but can only be
minimized by an unbounded sequence of points, it seems natural to
wonder if the function can in fact 
be minimized by
a sequence of points that go to infinity along some
halfline.
We have seen this in many examples,
\indexg{Product of hyperbolas!minimized along halfline|(}%
like \Cref{ex:recip-fcn-eg}, where we observed that the product of hyperbolas
can be minimized along a halfline from $\ee_2$ in the direction of $\ee_1$,
by the sequence of points $\xx_t=t\ee_1+\ee_2$,
implying
that its extension is minimized by
\indexg{Product of hyperbolas!minimized along halfline|)}%
$\limray{\ee_1}\plusl\ee_2$.
Can every convex function be minimized in this way?
Or is it sometimes
necessary to use a more complicated sequence of points?
In astral terms, we are asking
if the extension $\fext$
of a convex function can always be minimized by a point in $\extspace$
whose astral rank is at most one.\looseness=-1

\indexg{minimizers of extensions!maximum rank@of maximum rank|(}%
\indexg{rank, astral!minimizer of maximum possible|(}%
We begin this section by showing
that this is not always possible.
Specifically,
we study an example of a convex function
$f:\Rn\rightarrow\R$ whose extension $\fext$
can only be minimized by a point with astral rank $n$, the maximum
possible.
Thus, it is not just that the function cannot be minimized by following
a standard halfline; rather,
the only way to minimize
the function is by pursuing a trajectory which must span and grow
to infinity
in all $n$ dimensions.
We saw the same behavior in
\Cref{ex:two-speed-exp}, and indeed, the function presented
below is a variant of that example, generalized to $\Rn$.

We define the function $f:\Rn\rightarrow\R$ described above as
follows, for $\xx\in\Rn$:
\begin{align}
  f(\xx)
  &=
  \exp\regParens{x_2^2 - 2 x_1}
  +
  \exp\regParens{x_3^2 - 2 x_2}
  + \dotsb +
  \exp\regParens{x_n^2 - 2 x_{n-1}}
  +
  \exp({-x_n})
  \nonumber
  \\
  &=
  \sum_{i=1}^n h_i(\xx),
  \label{eq:h:9}
\end{align}
where, as usual,
$x_i$ is the $i$-th component of vector $\xx$,
and $h_i:\Rn\rightarrow\R$, for $i=1,\ldots,n$,
is defined by
\[
  h_i(\xx) =
  \begin{cases}
      \exp\regParens{x_{i+1}^2 - 2 x_i}
         & \text{if $i<n$,}
      \\
      \exp({- x_n})
         & \text{if $i=n$.}
  \end{cases}
\]
Each function $h_i$ is convex
(since $x_{i+1}^2 - 2 x_i$ and $-x_n$, as functions of $\xx\in\Rn$,
are both convex, implying $h_i$
is too, by \Cref{pr:j:2}\ref{pr:j:2a}).
Therefore, $f$ is convex as well.
Clearly, $f$ and the $h_i\negKern$'s are all also continuous, closed, proper,
finite everywhere, and strictly positive everywhere.

Intuitively, to minimize $h_i$, for $i<n$,
we need $x_i$ to be growing to
$+\infty$ faster than~$x_{i+1}^2$.
Thus, to minimize $f$, we need every variable $x_i$ to tend to
$+\infty$, with $x_1$ growing faster than $x_2$,
which is growing faster than $x_3$,
and so on.
In other words, the entries of a minimizing sequence need to converge
to $+\infty$ at decreasing rates, implying
(by \Cref{cor:omm:seq}) that the sequence
itself converges to the point $\omm=\ommsub{n}$,
as was defined in \eqref{eq:i:7},
and whose astral rank is $n$.
In accord with this intuition, we show in the next
\namecref{pr:hi-rank-f-eg} that $\omm$ is indeed the only minimizer of
$\fext$:

\begin{proposition}  \label{pr:hi-rank-f-eg}
\indexg{canonical minimizers!maximum rank@of maximum rank|(}%
\indexg{universal reducers!maximum rank@of maximum rank|(}%
  Let $f$ be as defined in \eqref{eq:h:9}.
  Then $\fext$ is uniquely minimized at~$\omm$,
  with $\fext(\omm)=0$.
  Consequently, $\omm$ is also $\fext$'s only canonical minimizer, and
  $f$'s only universal reducer.
\end{proposition}

\begin{proof}
  ~
  
\begin{proof-parts}
\pfpart{Only minimizer:}
First, to see that $\omm$ minimizes $\fext$, 
let
\begin{equation}  \label{eq:pr:hi-rank-f-eg:2}
  \zz_t
  =
  \sum_{i=1}^n t^{2^{n-i}} \ee_i
  =
  \trans{\bigBracks{t^{2^{n-1}}\!,\ldots,t^4,t^2,t^1}},
\end{equation}
for all $t$.
The resulting sequence $\seq{\zz_t}$ has entries that converge to
$+\infty$ at decreasing rates, so $\zz_t\rightarrow\omm$ by
\Cref{cor:omm:seq}.
Further, by straightforward algebra,
\[
  f(\zz_t)
  =
  \sum_{i=1}^{n} \exp\bigParens{-t^{2^{n-i}}}
  \rightarrow
  0.
\]
Thus,
\[ 0=\lim f(\zz_t)\geq \fext(\omm)\geq\min\fext=\inf f\geq 0, \]
with first inequality by $\fext$'s definition,
and second equality by \Cref{pr:fext-min-exists}.
Therefore, $\fext(\omm)=\min\fext=0$.

To show $\fext$ has no other minimizers, let $\xbar\in\extspace$ be
any minimizer of $\fext$, so that $\fext(\xbar)=0$.
We will show this implies $\xbar=\omm$.
Let $\seq{\xx_t}$ be a sequence in $\Rn$
such that $\xx_t\rightarrow\xbar$
and $f(\xx_t)\rightarrow\fext(\xbar)=0$
(which exists by \Cref{pr:d1}).
Then for $i=1,\ldots,n$, we claim $x_{t,i}\rightarrow+\infty$.
When $i<n$,
this is because, for all $t$, 
\begin{equation}  \label{eq:pr:hi-rank-f-eg:3}
  f(\xx_t)
  \geq
  h_i(\xx_t)
  =
  \exp\bigParens{x_{t,i+1}^2 - 2 x_{t,i}}
  \geq
  e^{-2x_{t,i}}
  >
  0.
\end{equation}
Since also $f(\xx_t)\rightarrow 0$, this implies
$e^{-2x_{t,i}}\rightarrow 0$, and therefore that
$x_{t,i}\rightarrow+\infty$.
Similarly, for $i=n$, we have
$f(\xx_t)\geq h_n(\xx_t) = e^{- x_{t,n}}>0$,
implying
$x_{t,n}\rightarrow+\infty$.

Consequently, there can only be finitely many elements with
$x_{t,i}\leq 0$ for any $i$; by discarding these, we assume henceforth that
$x_{t,i}>0$ for all $t$ and all $i$.

Next, for $i=1,\ldots,n-1$, we show
$x_{t,i+1}/x_{t,i}\rightarrow 0$.
Since $f(\xx_t)\rightarrow 0$, we must have
$f(\xx_t)<1$ for all $t$ sufficiently large,
implying, when combined with
\eqref{eq:pr:hi-rank-f-eg:3}, that
$x_{t,i+1}^2 - 2 x_{t,i} < 0$,
and so that
$0 < x_{t,i+1}/x_{t,i} < 2/x_{t,i+1}$.
Since $x_{t,i+1}\rightarrow +\infty$, this proves
$x_{t,i+1}/x_{t,i}\rightarrow 0$.

Thus, the entries of $\seq{\xx_t}$ converge to $+\infty$ at decreasing
rates, so $\xx_t\rightarrow\omm$ by \Cref{cor:omm:seq}.
Therefore, $\xbar=\omm$ (since also $\xx_t\rightarrow\xbar$).

\pfpart{Only canonical minimizer and only universal reducer:}
Suppose $\ebar\in\unimin{f}$ and $\qq\in\Rn$ minimizes
$\fullshad{f}$.
Then $\ebar\plusl\qq$ is a canonical minimizer and so also a minimizer
of $\fext$ (by \Cref{pr:unimin-to-global-min}).
Therefore, $\ebar\plusl\qq=\omm$ since $\omm$ is $\fext$'s only
minimizer.
Further, since every point's iconic part is unique
(\Cref{thm:icon-fin-decomp}), this
also shows that $\ebar=\omm$.
Thus, every canonical minimizer and universal reducer is equal to
$\omm$.

Moreover, by Propositions~\ref{pr:univ-red-props}(\ref{pr:univ-red-props:min})
and~\ref{pr:new:thm:f:8a}(\ref{pr:new:thm:f:8a:nonemp-closed}),
there must exist such points $\ebar$ and
$\qq$, implying by the foregoing that $\omm$ actually is a canonical
minimizer and universal reducer.%
\indexg{canonical minimizers!maximum rank@of maximum rank|)}%
\indexg{universal reducers!maximum rank@of maximum rank|)}%
\qedhere
\end{proof-parts}
\end{proof}

Thus,
for $n\ge 2$,
the function $f$ in \eqref{eq:h:9}
cannot be minimized by a sequence following a
halfline, nor converging asymptotically to a halfline,
nor converging to any linear or affine subspace of dimension less than $n$.
On the contrary, the function can only be minimized by a sequence that
grows to infinity across all $n$ dimensions.
For example, $f$ is minimized by the sequence $\seq{\zz_t}$ from
\eqref{eq:pr:hi-rank-f-eg:2} which was used in the proof of
\Cref{pr:hi-rank-f-eg}, and 
which converges to $\omm$.
However,
$f$ need not be minimized by {every} sequence
converging to $\omm$;
for instance,
the sequence $\zz'_t=3\zz_t$
converges to~$\omm$, but
$f(\zz'_t)\rightarrow +\infty$.
In other words, $\fext$ is not continuous at its unique minimizer
$\omm$.
Still, \Cref{pr:hi-rank-f-eg}
does imply that convergence
to $\omm$ is a {necessary} condition for a sequence to minimize
$f$, meaning $f$ cannot be minimized by any sequence that does
\emph{not} converge in $\extspace$ to $\omm$.%
\indexg{minimizers of extensions!maximum rank@of maximum rank|)}%
\indexg{rank, astral!minimizer of maximum possible|)}%

\indexg{minimizers of extensions!functions with any given|(}%
Generalizing \Cref{pr:hi-rank-f-eg},
for any point $\zbar\in\extspace$,
we next show how to construct a
function whose extension is minimized only at $\zbar$:

\begin{theorem}  \label{thm:gen-uniq-min}
  Let $\zbar\in\extspace$.
  Then there exists a convex function $f:\Rn\rightarrow\R$ whose
  extension $\fext$ is uniquely minimized at $\zbar$.
\end{theorem}

\begin{proof}
Let $\zbar=\VV\omm\plusl\qq$ be the canonical representation of
$\zbar$ where $\VV\in\R^{n\times k}$ is column-orthogonal and $\qq\in\Rn$ with
${\VV}\perp\qq$.

Let $\fomin:\Rk\rightarrow\R$ denote the function given in
\eqref{eq:h:9} but with $n$ replaced by~$k$ so that its extension,
$\fominext$, is uniquely minimized at $\omm=\omm_k$,
with $\fominext(\omm)=0$,
by
\Cref{pr:hi-rank-f-eg}.
Let
$\normfcn:\Rn\rightarrow\R$ be the norm function,
$\normfcn(\xx)=\norm{\xx}$ for $\xx\in\Rn$,
and
let $\PP=\Iden - \VV \trans{\VV}$,
where $\Iden$ is the $n\times n$ identity matrix.
Then $\PP$ is the projection matrix onto $(\colspace{\VV})^\perp$
(by \Cref{pr:basis-to-proj-mat}),
since $\VV$ is column-orthogonal.

To construct $f$, we define
\[
  f(\xx)
  =
  \fomin\regParens{\trans{\VV} \xx}
  +
  \normfcn\regParens{\PP \xx - \qq}
\]
for $\xx\in\Rn$.
Both $\fomin$ and $\normfcn$ are
convex, closed, proper, finite everywhere, and
nonnegative everywhere, so $f$ is as well
(\Cref{pr:std-sum-fcns-cvx,roc:thm5.7:fA}).
In particular, this implies that their extensions are nonnegative as
well (by \Cref{pr:fext-min-exists}).

By \Cref{pr:seq-to-inf-has-inf-len},
$\normfcn$'s extension is
  \[
    \normfcnext(\xbar)
    =
    \begin{cases}
      \norm{\xbar} &
        \text{if $\xbar\in\Rn$,}
    \\
      +\infty      &
        \text{otherwise,}
    \end{cases}
  \]
for $\xbar\in\extspace$,
and moreover, that proposition shows that $\normfcn$ is extensibly
continuous everywhere.
Clearly, $\normfcnext$ is uniquely minimized by $\zero$.

By \Cref{pr:ext-affine-comp}(\ref{pr:ext-affine-comp:c}),
the extension of
$\xx\mapsto\normfcn(\PP \xx - \qq)$
is
$\xbar\mapsto\normfcnext(-\qq \plusl \PP\xbar)$,
and
by \Cref{thm:ext-linear-comp},
the extension of
$\xx\mapsto \fomin\regParens{\trans{\VV} \xx}$
is
$\xbar\mapsto \fominext\regParens{\trans{\VV} \xbar}$.
Combining then yields
\begin{equation}  \label{eqn:thm:gen-uniq-min:1}
  \fext(\xbar)
  =
  \fominext\regParens{\trans{\VV} \xbar}
  +
  \normfcnext(-\qq \plusl \PP\xbar)
\end{equation}
by \Cref{thm:ext-sum-fcns-w-duality}
since $\fomin$ and $\normfcn$ are both nonnegative and finite
everywhere.

To see that $\zbar$ minimizes $\ef$, note that
\begin{align*}
  \trans{\VV}\zbar
  &=
  \trans{\VV} \VV\omm \plusl \trans{\VV}\qq
  =
  \omm,
\end{align*}
since $\VV$ is column-orthogonal and $\qq\perp\VV$,
and
\[
  -\qq\plusl\PP\zbar
  =
  -\qq\plusl\PP \VV \omm \plusl \PP \qq
  =
  -\qq\plusl\qq
  =
  \zero,
\]
since $\PP\VV=\zeromat{n}{k}$ and $\PP\qq=\qq$
(by
\Cref{pr:proj-mat-props}\ref{pr:proj-mat-props:e}\ref{pr:proj-mat-props:d},
and since $\qq\perp\VV$).
Thus, by \eqref{eqn:thm:gen-uniq-min:1},
$\fext(\zbar)=\fominext(\omm)+\normfcnext(\zero)=0$,
so $\zbar$ minimizes $\fext$
(since $\fext\geq 0$).

It remains to show that $\zbar$ is the only minimizer of $\fext$.
Let $\xbar\in\extspace$ be any minimizer of~$\fext$.
Then $\fext(\xbar)=\min\fext=0$, implying,
by \eqref{eqn:thm:gen-uniq-min:1} and since
$\fominext$ and $\normfcnext$ are nonnegative
with unique minimizers,
that $\trans{\VV}\xbar=\omm$ and $-\qq\plusl\PP\xbar=\zero$.
The latter implies that $\PP\xbar=\qq$, so
\[
  \zbar
  =
  \VV\omm \plusl \qq
  =
  \VV\trans{\VV}\xbar \plusl \PP\xbar
  =
  (\VV\trans{\VV} + \PP)\xbar
  =
  \Iden \xbar = \xbar.
\]
The third equality follows from
\Cref{prop:commute:AB} since
$\PP\xbar=\qq\in\Rn$, which therefore commutes
with any astral point (by \Cref{pr:i:7}\ref{pr:i:7d}).
The fourth equality is because
$\PP=\Iden-\VV\trans{\VV}$.
Thus, as claimed, $\zbar$ is $\fext$'s only minimizer.%
\indexg{minimizers of extensions!functions with any given|)}%
\end{proof}

\Cref{pr:hi-rank-f-eg} shows
that the extension of a convex function
need not have a minimizer
of astral rank at most one, nor even strictly less than $n$.
\indexg{rank, astral!minimizers with rank at most one|(}%
\indexg{minimizers of extensions!rank at most one@of rank at most one|(}%
\indexg{canonical minimizers!rank at most one@of rank at most one|(}%
We next
provide a contrasting result, a general sufficient condition for such
minimizers to exist, showing
that the extension $\fext$ of every convex, lower semicontinuous
function $f:\Rn\rightarrow\Rext$ that is also recessive complete
must in fact have a (canonical) minimizer
of astral rank at most one.
In other words, such functions can be minimized
either at a finite point or by a sequence that (asyptotically) follows
a halfline towards infinity.
This is dramatically different from what we saw in
\Cref{pr:hi-rank-f-eg}.

\indexg{universal reducers!rank at most one@of rank at most one|(}%
\indexg{recessive completeness!ensuring rank-one minimizer|(}%
Specifically, we show that for such a recessive complete function $f$,
if $\vv$ is any point in the relative interior of
$f$'s standard recession cone, $\resc{f}$,
then the associated astron $\limray{\vv}$ must be a universal reducer.
As a
consequence, $\limray{\vv}\plusl\qq$
(whose astral rank is~$0$ or~$1$)
must be a canonical minimizer of $\fext$
for every $\qq\in\Rn$ that minimizes $\fullshad{f}$.
Since such points $\vv$ and $\qq$ must exist, this shows that such
a universal reducer and canonical minimizer must exist as well.

\begin{theorem}\label{thm:unimin-can-be-rankone}
  Let $f:\Rn\rightarrow\Rext$ be convex, lower semicontinuous,
  and recessive complete.
  Let $\vv\in\ric{f}$, and let $\qq\in\Rn$ minimize
  $\fullshad{f}$ (which both must exist).
  Then $\limray{\vv}$ is a universal reducer of $f$,
  and $\limray{\vv}\plusl\qq$ is a canonical
  minimizer of $\fext$.
\end{theorem}

\begin{proof}
First, $\resc{f}$ is convex and nonempty
(\Cref{pr:resc-cone-basic-props}), so its relative interior is
nonempty as well (\Cref{pr:ri-props}\ref{pr:ri-props:roc-thm6.2b}).
Also, $\fullshad{f}$ attains its minimum
(by \Cref{pr:univ-red-props}\ref{pr:univ-red-props:min}).
Therefore, $\vv$ and $\qq$, as in the 
\namecref{thm:unimin-can-be-rankone}'s statement, must exist.

It suffices to show that
$\limray{\vv}$ is a universal reducer, that is,
that $\fshadv=\fullshad{f}$,
since this will imply that
$\limray{\vv}\plusl\qq$ is a canonical minimizer as well,
by definition.

If $f\equiv+\infty$ then the claim holds trivially since then
$\fext\equiv+\infty$, implying
$\fshadv=\fullshad{f}\equiv+\infty$.
We therefore assume henceforth that $f\not\equiv+\infty$.

Let $g=\fshadv$, which is convex and lower semicontinuous
(\Cref{thm:a10-nunu}).
We claim first that $\slopes{g}\subseteq\rescperp{f}$.
To prove this, suppose $\uu\in\slopes{g}$ and that
$\ww\in\resc{f}$.
We aim to show $\uu\cdot\ww=0$, and thus that
$\uu\in\rescperp{f}$.

Since $\uu\in\slopes{g}$ and
$\vv\in\resc{f}$,
\Cref{thm:e1}(\ref{thm:e1:d}) implies that
$\uu\in\slopes{f}$ and $\uu\inprod\vv=0$.
Since $\resc{f}=\polar{(\slopes{f})}$ by 
\Cref{thm:rescpol-is-slopes}, this further implies
that $\uu\cdot\ww\leq 0$.

Next, because $\vv\in\ric{f}$ and 
$\ww\in\resc{f}$,
there must exist $\delta>0$ such that the point
$\vv+\delta(\vv-\ww)$
is in $\resc{f}$ as well
(by \Cref{roc:thm6.4}).
This then implies that
\[  -\delta \uu\cdot\ww
   = \uu \cdot [\vv+\delta(\vv-\ww)]
  \leq 0,
\]
where the equality is because $\uu\cdot\vv=0$,
and the inequality is again because
$\resc{f}=\polar{(\slopes{f})}$.
Since $\delta>0$,
it follows that $\uu\cdot\ww\geq 0$, and so that
$\uu\cdot\ww=0$.
Because this holds for all $\ww\in\resc{f}$,
we conclude that $\uu\in\rescperp{f}$.

Hence,
\begin{equation}   \label{eq:thm:unimin-can-be-rankone:1}
  \slopes{g}
  \subseteq
  \rescperp{f}
  =
  \rescperp{\ef}
  =\rescperp{\eg},
\end{equation}
where the inclusion is by the preceding argument,
the first equality is by
\Cref{pr:perpres-is-rescperp} since $f$ is recessive complete,
and the second equality is by \Cref{thm:f:6}.
Therefore,
\[
  \fshadv
  =
  g
  =
  \fullshad{g}
  =\fullshad{f},
\]
where the second equality is by 
\Cref{cor:thm:f:4:1}(\ref{cor:thm:f:4:1b},\ref{cor:thm:f:4:1c})
as a result of \eqref{eq:thm:unimin-can-be-rankone:1},
and the third is by \Cref{thm:f:5}.
Thus, $\limray{\vv}$ is a universal reducer for $f$, as claimed.
\end{proof}

\begin{example}[Two exps and a square, continued]
   \label{ex:two-exps-sq-cont}
\indexg{Two exps and a square!rank-one minimizer|(}%
In \Cref{ex:simple-eg-exp-exp-sq,ex:simple-eg-exp-exp-sq-part2},
we studied the minimization of
\[
  f(\xx) = \me^{x_3-x_1} + \me^{-x_2}
  + (2+x_2-x_3)^2.
\]
We derived the astral
recession cone $\resc{\ef}$, writing it as an intersection
of four closed astral halfspaces. Thus,
by \Cref{pr:rec-complete}, $f$ must be recessive complete.

The minimizer of $\ef$ that we exhibited in \Cref{ex:simple-eg-exp-exp-sq,ex:simple-eg-exp-exp-sq-part2}
had astral rank two.
However, by \Cref{thm:unimin-can-be-rankone},
since $f$ is recessive complete, $\ef$ must have a
(canonical) minimizer of astral rank at most one.
For instance,
letting $\vv=\trans{[2,1,1]}$, which is in $\resc{f}$'s relative
interior,
it can be checked that
$\limray{\vv}$ is a universal reducer.
Therefore, combining with $\qq=\trans{[0,0,2]}$, which minimizes
$\fullshad{f}$, yields the canonical minimizer
$\limray{\vv}\plusl\qq$
of astral rank one.%
\indexg{Two exps and a square!rank-one minimizer|)}%
\end{example}

Later, in \Cref{subsec:cond-for-cont},
we study the close connection between recessive completeness and
continuity.
In particular,
\Cref{thm:cont-at-min-implies-ares} shows that if
$\fext$ is continuous at all its minimizers
(and so also if it is continuous everywhere),
then $f$ is recessive
complete, and therefore, by \Cref{thm:unimin-can-be-rankone},
$\fext$ must have a universal reducer and canonical minimizer of
astral rank at most one.
Thus, every continuous extension can be minimized by some point of astral
rank at most one.
In the contrapositive, this means (for $n\geq 2$)
that because the function from
\Cref{pr:hi-rank-f-eg} has no minimizer of astral rank at most one,
it cannot be recessive complete, and it therefore
(by \Cref{thm:cont-at-min-implies-ares})
must be discontinuous at its unique minimizer, $\omm$;
indeed, we earlier saw that this was so (see discussion following the
proof of \Cref{pr:hi-rank-f-eg}).%
\indexg{universal reducers!rank at most one@of rank at most one|)}%
\indexg{rank, astral!minimizers with rank at most one|)}%
\indexg{minimizers of extensions!rank at most one@of rank at most one|)}%
\indexg{canonical minimizers!rank at most one@of rank at most one|)}%
\indexg{recessive completeness!ensuring rank-one minimizer|)}%

\section{Empirical risk minimization}
\label{sec:emp-loss-min}

\indexg{empirical risk functions|(}%
\indexg{empirical risk functions!general form|(}%
We next look closely at
functions $f:\Rn\rightarrow\R$ of the form
\begin{equation}   \label{eqn:loss-sum-form}
  f(\xx) = \sum_{i\in\indset} \ell_i(\xx\cdot\uu_i),
\end{equation}
for $\xx\in\Rn$, where $\indset$ is a finite index set,
each $\uu_i\in\Rn$, and each function
$\ell_i:\R\rightarrow\R$ is convex, lower-bounded, and
\indexg{empirical risk functions!general form|)}%
nondecreasing.
We focus on the minimizers of such functions, and especially how these
relate to concepts developed earlier.

\indexg{empirical risk functions!machine learning and statistics@in machine learning and statistics|(}%
Minimizing functions of the form given in \eqref{eqn:loss-sum-form}
is a fundamental problem
in machine learning and statistics.
Very briefly, in a typical setting, a learning algorithm might be
given random ``training examples'' $(\zz_i,y_i)$, for $i=1,\ldots,m$,
where $\zz_i\in\Rn$ is an
``instance'' or ``pattern'' (such as an image or photograph, treated
as a vector in $\Rn$ of pixel intensities),
and $y_i\in\{-1,+1\}$ is a ``label''
(that might indicate, for example, if the photograph is or is not of a
person's face).
The goal then is to find a rule for predicting if a new instance
$\zz\in\Rn$ should be labeled $-1$ or $+1$.
\indexg{logistic regression|(}%
As an example, in logistic regression,
the learner finds a vector $\ww\in\Rn$, based on the
training examples, and then predicts that a new instance $\zz$ should
be labeled according to the sign of $\ww\cdot\zz$.
Specifically, $\ww$ is chosen to minimize the
``logistic loss'' on the
training examples, that is,
\begin{equation}  \label{eqn:logistic-reg-obj}
  f(\ww) = \sum_{i=1}^m \ln\bigParens{1+\exp(-y_i \ww\cdot\zz_i)}.
\end{equation}
This kind of function, which is more generally called the
\emph{empirical risk}, has the same form as in
\eqref{eqn:loss-sum-form}
(with $\xx=\ww$, $\uu_i = -y_i \zz_i$, and $\ell_i(z)=\ln(1+\me^z)$
for $z\in\R$).%
\indexg{logistic regression|)}%
\indexg{empirical risk functions!machine learning and statistics@in machine learning and statistics|)}%

Returning to the general case in \eqref{eqn:loss-sum-form},
for $i\in\indset$,
we will assume
that $\inf \ell_i = 0$, and that $\ell_i$ is not constant
(in addition to the other assumptions mentioned above).
For the purposes of minimization, this is without loss of generality
since each function $\ell_i$ can be replaced by $\ell_i-(\inf\ell_i)$,
and any constant function $\ell_i$ can be discarded;
these modifications only change $f$ by finite additive constants.

Since $\ell_i$ is nondecreasing,
these conditions imply that
$\lim_{x\rightarrow -\infty} \ell_i(x) = 0$ and
$\lim_{x\rightarrow +\infty} \ell_i(x) = +\infty$
(by \Cref{pr:conv-inc:prop}\ref{pr:conv-inc:infsup}\ref{pr:conv-inc:nonconst}).
Each $\ell_i$ is convex and finite everywhere, and therefore
continuous everywhere
(\Cref{pr:stand-cvx-cont});
the same is also true of~$f$.
So $\ell_i$'s extension is
\begin{equation}  \label{eqn:hard-core:1}
  \ellbar_i(\barx) =
  \begin{cases}
            0              & \text{if $\barx=-\infty$,} \\
            \ell_i(\barx)  & \mbox{if $\barx\in\R$,} \\
            +\infty        & \mbox{if $\barx=+\infty$,}
  \end{cases}
\end{equation}
for $\barx\in\Rext$
(by \Cref{pr:conv-inc:prop}\ref{pr:conv-inc:infsup}).

Throughout this section, for a set $J\subseteq\indset$, we write
\indexm{u j}{$\uset{J}$}{points with given indices}%
$\uset{J}$ for the set of points $\uu_i$ with indices in $J$;
that is,
$  \uset{J} = \Braces{\uu_i :\: i\in J} $.

\indexg{empirical risk functions!extension of|(}%
\indexg{empirical risk functions!continuity of|(}%
\indexg{empirical risk functions!recessive completeness of|(}%
\indexg{empirical risk functions!recession cone of|(}%
The next \namecref{pr:hard-core:1}
gives the form of $f$'s extension, $\fext$, and shows that
$\fext$ is continuous everywhere, and that $f$ is recessive complete.
We also give explicit expressions for some of the sets that were
studied in preceding chapters, including $f$'s barrier cone and
recession cone, and $\fext$'s astral recession cone.
Part~(\ref{pr:hard-core:1:conedomfstar}) provides another example of a
fact whose statement is only in terms of notions from standard convex
analysis, but for which we give an astral-based proof (in this
case, a particularly simple one).

\begin{theorem}  \label{pr:hard-core:1}
  Let $f:\Rn\rightarrow\R$ have the form given in
  \eqref{eqn:loss-sum-form}, where,
  for $i\in\indset$,
  $\uu_i\in\Rn$  %
  and
  $\ell_i:\R\rightarrow\R$ is convex, nondecreasing, not constant,
  with $\inf \ell_i = 0$.
  Then:
  \begin{letter-compact}
  \item  \label{pr:hard-core:1:a-ext}
    For $\xbar\in\extspace$,
    \[
       \fext(\xbar) = \sum_{i\in\indset} \ellbar_i(\xbar\cdot\uu_i).
    \]
  \item  \label{pr:hard-core:1:a-cnt}
    $\fext$ is continuous everywhere.
  \item  \label{pr:hard-core:1:b-cones}
    The astral recession cone of $\fext$ and standard recession cone
    of $f$ are:
    \begin{align*}
      \aresconef &= \{ \ybar\in\extspace :\: \ybar\cdot\uu_i \leq 0
                            \text{ for } i\in\indset \}, \\
      \resc{f} &= \{ \yy\in\Rn :\: \yy\cdot\uu_i \leq 0
                            \text{ for } i\in\indset \}.
    \end{align*}
  \item  \label{pr:hard-core:1:b-resc1}
    $f$ is recessive complete.
  \item  \label{pr:hard-core:1:b-resc2}
    $\rescperp{f}=\perpresf$.
  \item  \label{pr:hard-core:1:conedomfstar}
    $\barr{f}=\cone(\dom\fstar)=\cone\uset{\indset}$.
  \end{letter-compact}
\end{theorem}

\begin{proof}
~

\begin{proof-parts}
\pfpart{%
  Parts~(\ref{pr:hard-core:1:a-ext})
  and~(\ref{pr:hard-core:1:a-cnt}):
}
For $i\in\indset$, let $h_i(\xx)=\ell_i(\xx\cdot\uu_i)$ for
$\xx\in\Rn$.
As noted above, $\ell_i$ is continuous everywhere, implying
$\ellbar_i$ is continuous everywhere
by \Cref{pr:conv-inc:prop}(\ref{pr:conv-inc:infsup}).
In turn,
\Cref{pr:Gf-cont}(\ref{pr:Gf-cont:b})
(applied with $G$ and $f$, as they appear in that
\namecref{pr:Gf-cont}, set respectively to
$\ellbar_i$ and $\xx\mapsto\xx\cdot\uu_i$)
implies that $\hext_i(\xbar)=\ellbar_i(\xbar\cdot\uu_i)$ for $\xbar\in\extspace$,
and that $\hext_i$ is continuous everywhere.
The form and continuity of $\fext$ now follow from
\Cref{pr:ext-sum-fcns}(\ref{pr:ext-sum-fcns:b},\ref{pr:ext-sum-fcns:c})
(with summability following from $\hext_i\geq 0$ since $h_i\geq 0$).
(Note that in this reasoning, we are also using the equivalence of continuity and extensible
continuity for convex functions, \Cref{thm:ext-cont-f}\ref{thm:ext-cont-f:b}.)

\pfpart{Part~(\ref{pr:hard-core:1:b-cones}):}
By \Cref{cor:rec-equiv}
(with $\yy=\zero$),
$\resc{\ef}=\set{\ybar\in\eRn:\:\ef(\limray{\ybar})<+\infty}$.
Since $\ellbar_i\ge 0$, part~(\ref{pr:hard-core:1:a-ext})
implies that 
$\ef(\limray{\ybar})<+\infty$
if and only if
$\ellbar_i(\limray{\ybar})<+\infty$ for all $i$,
which is the case if and only if
$\limray{\ybar}\inprod\uu_i<+\infty$ for all $i$
(by Eq.~\ref{eqn:hard-core:1}),
which in turn holds if and only if
$\ybar\inprod\uu_i\le 0$ for all $i$
(\Cref{pr:astrons-exist}).
This proves the identity for $\resc{\ef}$.

The expression for $\resc{f}$ now follows immediately from
\Cref{pr:f:1}.

\pfpart{Part~(\ref{pr:hard-core:1:b-resc1}):}
From part~(\ref{pr:hard-core:1:b-cones}), $\aresconef$ is an
intersection of closed astral halfspaces
(noting that inequalities $\ybar\cdot\uu_i\leq 0$ with
$\uu_i=\zero$ always hold and so can be disregarded).
Therefore, $\aresconef$ is astral polyhedral, so $f$ is recessive
complete by \Cref{pr:rec-complete}.

\pfpart{Part~(\ref{pr:hard-core:1:b-resc2}):}
Since $f$ is recessive complete by part~(\ref{pr:hard-core:1:b-resc1}),
this follows from
\Cref{pr:perpres-is-rescperp}.

\pfpart{Part~(\ref{pr:hard-core:1:conedomfstar}):}
We have
\[
  \cone(\dom\fstar)
  =
  \barr{f}
  =
  \polar{(\aresconef)}
  =
  \apolpol{(\cone\uset{\indset})}
  =
  \cone\uset{\indset}.
\]
The first equality is by
\Cref{cor:ent-clos-is-slopes-cone}(\ref{cor:ent-clos-is-slopes-cone:a},\ref{cor:ent-clos-is-slopes-cone:c})
since $f$ is finite everywhere and therefore reduction-closed
(by \Cref{pr:j:1}\ref{pr:j:1d}).
The second equality is by \Cref{cor:ares-is-apolslopes}.
The third equality is because the expression for $\aresconef$ given in 
part~(\ref{pr:hard-core:1:b-cones})
matches that obtained by applying 
\Cref{pr:ast-pol-props}(\ref{pr:ast-pol-props:coneSpol})
to $\uset{\indset}$,
thus implying that $\aresconef = \apol{(\cone\uset{\indset})}$.
The final equality is by
\Cref{thm:dub-ast-polar}(\ref{thm:dub-ast-polar:b}).%
\indexg{empirical risk functions!extension of|)}%
\indexg{empirical risk functions!continuity of|)}%
\indexg{empirical risk functions!recessive completeness of|)}%
\indexg{empirical risk functions!recession cone of|)}%
\qedhere
\end{proof-parts}
\end{proof}

\idxmatus\citet{primal_dual_boosting}
studied functions of the form we are considering,
and showed that the index set $\indset$
can be usefully partitioned
into a so-called
\indexg{easy set}%
easy set and a
\indexg{hard core}%
hard core.
\indexg{empirical risk functions!persistent and reducible indices of|(}%
\indexg{persistent indices|(}%
\indexg{reducible indices|(}%
In this book, we will instead refer to these two kinds of indices as reducible and persistent.
We will define these formally in a moment, but first pause to provide
motivation.
Consider for a moment the special
case that each $\ell_i$ is not just nondecreasing, but is in fact
strictly increasing (like in Eq.~\ref{eqn:logistic-reg-obj}).
In that case, roughly speaking,
reducible indices are those $i\in\indset$
for which the corresponding term $\ell_i(\xx\cdot\uu_i)$
can be reduced to its minimum value, achieved by driving
$\xx\cdot\uu_i$ to~$-\infty$, without making any other
term go to $+\infty$.
The remaining indices are persistent.
Their corresponding terms cannot be driven to their minimum values
while minimizing $f$; rather, for these indices~$i$,
$\xx\cdot\uu_i$ must converge to some finite value.

In astral terms,
an index $i$ is considered reducible
if there exists a point
$\xbar\in\extspace$ with $\xbar\cdot\uu_i=-\infty$ (so that
$\ellbar_i(\xbar\cdot\uu_i)$ equals $\ellbar_i\negKern$'s minimum value)
and for which $\fext(\xbar)<+\infty$.
This is equivalent to there existing a sequence $\seq{\xx_t}$
in~$\Rn$ for
which $\xx_t\cdot\uu_i\rightarrow-\infty$,
without $f(\xx_t)$ becoming unboundedly large.
Otherwise, if there exists no such $\xbar$, then $i$ is persistent.

We can write any $\xbar\in\extspace$ as $\xbar=\ebar\plusl\qq$ for some
$\ebar\in\corezn$ and $\qq\in\Rn$.
The condition that $\fext(\xbar)<+\infty$ implies that
$\ebar\in\aresconef$,
by \Cref{cor:a:4}.
In that case, for any $\uu_i\in(\resc{\ef})^\perp$,
we have $\ebar\cdot\uu_i=0$, implying that $\xbar\cdot\uu_i=\qq\cdot\uu_i$,
which is in $\R$, and therefore $\xbar$ cannot minimize $\ellbar_i(\xbar\inprod\uu_i)$.
This shows that if $\uu_i\in\perpresf$ then $i$ must be persistent,
because
for all $\xbar\in\extspace$, either $\fext(\xbar)=+\infty$ or
$\xbar\cdot\uu_i>-\infty$.

In accord with this motivation, we can now state our formal definitions
(again assuming each $\ell_i$ is nondecreasing, but not necessarily strictly increasing):

\begin{definition}   \label{def:hard-core}
  Let $f:\Rn\rightarrow\R$ have the form given in
  \eqref{eqn:loss-sum-form}, where,
  for $i\in\indset$,
  $\uu_i\in\Rn$
  and
  $\ell_i:\R\rightarrow\R$ is convex, nondecreasing, not constant,
  with
  $\inf \ell_i = 0$.
  We say that $i\in\indset$ is a \emph{persistent index} for $f$ if
  $\uu_i\in\rescperp{f}$. Otherwise, $i$ is a \emph{reducible index} for $f$.
  The set of \irredIndicesFor~$f$ is denoted
  \begin{equation}  \label{eq:hard-core:3}
\indexm{pers f}{$\hardcore{f}$}{persistent indices}%
    \hardcore{f} = \bigBraces{ i \in \indset :\: \uu_i\in\rescperp{f} }.
  \end{equation}
\end{definition}
  Note that the definition of persistent and reducible indices
  depends on both $f$ and the
  specific tuple
  $\tupset{\uu_i}{i\in\indset}$;
  however, for simplicity,
  we suppress the dependence on the $\uu_i\negKern$'s
  from the notation and assume
  that these are provided implicitly
  as part of the definition of $f$.

Also,
  under the conditions of the definition,
  $\rescperp{f}=\perpresf$ (by \Cref{pr:hard-core:1}\ref{pr:hard-core:1:b-resc2}), so
  persistent and reducible indices could be
  equivalently defined instead using
\indexg{empirical risk functions!persistent and reducible indices of|)}%
\indexg{persistent indices|)}%
\indexg{reducible indices|)}%
  $\perpresf$.

\begin{example}[Two exps and a square, continued]
\label{ex:erm-running-eg1}
\indexg{Two exps and a square!persistent and reducible indices of|(}%
The function $f$
from \Cref{ex:simple-eg-exp-exp-sq}
can be put in the form of
\eqref{eqn:loss-sum-form}.
To see this, let $\sqp:\R\rightarrow\R$ be defined,
for $z\in\R$,
by $\sqp(z)=\paren{\max\{0,z\}}^2$,
the square of the positive part of $z$.
This function is convex, nondecreasing,
not constant, with $\inf \sqp=0$.
Also, $z^2=\sqp(z)+\sqp(-z)$ for all $z\in\R$.
Thus, we can write $f$, for $\xx\in\R^3$, as
\begin{equation}
\label{eq:simple-eg-exp-exp-sq:modform}
  f(\xx) = \me^{x_3-x_1} + \me^{-x_2}
               + \sqp(2+x_2-x_3)
               + \sqp(-2-x_2+x_3),
\end{equation}
which satisfies the conditions of
\Cref{pr:hard-core:1}, with
\begin{align*}
  &
    \uu_1=\trans{[-1,0,1]},
    \ell_1(z)=e^z,
  &&
    \uu_3=\trans{[0,1,-1]},
    \ell_3(z)=\sqp(2+z),
  \\
  &
    \uu_2=\trans{[0,-1,0]},
    \ell_2(z)=e^z,
  &&
    \uu_4=\trans{[0,-1,1]},
    \ell_4(z)=\sqp(-2+z),
\end{align*}
for $z\in\R$.
As a result,
that proposition confirms various previously determined facts
about $f$.
The expressions for $\ef$, $\resc{f}$, and $\resc{\ef}$
implied by the proposition match the expressions we derived
in \Cref{ex:simple-eg-exp-exp-sq}, and as the proposition states,
$\ef$ is continuous everywhere and $f$ is recessive complete.

Since $f$ is recessive complete and by
\Cref{thm:f:4x}(\ref{thm:f:4xa}),
$\rescperp{f}=\perpresf=(\conssp{\fullshad{f}})^\perp$.
Using the expression for $(\conssp{\fullshad{f}})^\perp$,
from
\Cref{ex:simple-eg-exp-exp-sq-2}, we thus obtain
\[
  \rescperp{f}
  =\perpresf=(\conssp{\fullshad{f}})^\perp
  =\set{\trans{[0,\lambda,-\lambda]}:\:\lambda\in\R}.
\]
This set includes
$\uu_3$ and $\uu_4$, but
not $\uu_1$ or $\uu_2$.
Thus indices $1$ and $2$ are reducible, and $3$ and $4$ are
persistent, forming the \irredIndexSet
$\hardcore{f}=\{3,4\}$.%
\indexg{Two exps and a square!persistent and reducible indices of|)}%
\end{example}

\indexg{persistent indices!universal reducers and reduction and|(}%
\indexg{empirical risk functions!universal reducers of|(}%
\indexg{empirical risk functions!universal reduction of|(}%
We have seen that
the set $\unimin{f}$ of universal reducers
together with the universal reduction $\fullshad{f}$
play important roles in the theory of
astral minimizers.
We next show how to express them in terms of $f$'s \irredIndices.

\begin{theorem}  \label{thm:hard-core:3}
  Let $f:\Rn\rightarrow\R$ have the form given in
  \eqref{eqn:loss-sum-form}, where,
  for $i\in\indset$,
  $\uu_i\in\Rn$
  and
  $\ell_i:\R\rightarrow\R$ is convex, nondecreasing, not constant,
  with $\inf \ell_i = 0$.
  \begin{letter-compact}
  \item  \label{thm:hard-core:3:a}
    If $\ybar\in\aresconef$, then $\ybar\cdot\uu_i=0$ for
    all $i\in\hardcore{f}$.
  \item  \label{thm:hard-core:3:b:conv}
    Let $\ybar\in\extspace$.
    Then $\ybar\in\conv(\unimin{f})$ if and only if,
    for all $i\in\indset$,
    \begin{equation}  \label{eq:hard-core:conv:2}
       \ybar\cdot\uu_i =
       \begin{cases}
              0       & \text{if $i\in\hardcore{f}$,} \\
              -\infty & \text{otherwise.}
       \end{cases} 
    \end{equation}
  \item  \label{thm:hard-core:3:b:univ}
    Let $\ebar\in\corezn$.
    Then $\ebar\in\unimin{f}$ if and only if,
    for all $i\in\indset$,
    \[
       \ebar\cdot\uu_i =
       \begin{cases}
              0       & \text{if $i\in\hardcore{f}$,} \\
              -\infty & \text{otherwise.}
       \end{cases}
    \]
  \item  \label{thm:hard-core:3:c}
    Let $\ybar\in\extspace$, and
    suppose, for some $i\in\hardcore{f}$, that
    $\ybar\cdot\uu_i<0$.
    Then there exists $j\in\hardcore{f}$ for which
    $\ybar\cdot\uu_j>0$.
  \item  \label{thm:hard-core:3:d}
    For $\xx\in\Rn$,
    \[
       \fullshad{f}(\xx) = \sum_{i\in\hardcore{f}} \ell_i(\xx\cdot\uu_i).
    \]
  \item  \label{thm:hard-core:3:e}
    $\rescperp{f}=\spn\uset{\hardcore{f}}$.
  \end{letter-compact}
\end{theorem}

\begin{proof}
~

\begin{proof-parts}
\pfpart{Part~(\ref{thm:hard-core:3:a}):}
If $i\in\hardcore{f}$, then $\uu_i\in\rescperp{f}=\perpresf$
(with equality by \Cref{pr:hard-core:1}\ref{pr:hard-core:1:b-resc2}),
which implies that $\ybar\cdot\uu_i=0$ for all $\ybar\in\aresconef$.

\pfpart{Part~(\ref{thm:hard-core:3:b:conv}):}
Suppose $\ybar\in\conv(\unimin{f})$,
implying $\ybar\in\aresconef$
by \Cref{pr:new:thm:f:8a}(\ref{pr:new:thm:f:8a:a}).
Then
$\ybar\cdot\uu_i\leq 0$
for $i\in\indset$,
by \Cref{pr:hard-core:1}(\ref{pr:hard-core:1:b-cones}).
Specifically,
$\ybar\cdot\uu_i=0$ for $i\in \hardcore{f}$,
by part~(\ref{thm:hard-core:3:a}).
It remains then only to show that if $i\not\in\hardcore{f}$ then
$\ybar\cdot\uu_i=-\infty$, which would be implied by showing that
$\ybar\cdot\uu_i\not\in\R$.

Suppose then, by way of contradiction, that there exists
$j\in\indset\setminus(\hardcore{f})$
with $\ybar\cdot\uu_j\in\R$.
By definition,
since $j\not\in\hardcore{f}$,
there must exist $\vv\in\resc{f}$ with $\vv\cdot\uu_j\neq 0$,
implying, by \Cref{pr:hard-core:1}(\ref{pr:hard-core:1:b-cones}),
that actually $\vv\cdot\uu_j< 0$.

Let $\lambda\in\R$.
To derive a contradiction,
we compare function values at
$\ybar\plusl\lambda\uu_j$
and
$\limray{\vv}\plusl\ybar\plusl\lambda\uu_j$.
For $i\in\indset$,
by application of
\Cref{pr:hard-core:1}(\ref{pr:hard-core:1:b-cones}),
$\vv\cdot\uu_i\leq 0$, and so
$\limray{\vv}\cdot\uu_i\leq 0$,
implying $\limray{\vv}\in\aresconef$.
Thus,
\begin{align}
  -\infty
  <
  \ellbar_i\bigParens{(\limray{\vv}\plusl\ybar\plusl\lambda\uu_j)\cdot\uu_i}
  &=
  \ellbar_i\bigParens{\limray{\vv}\cdot\uu_i
              \plusl(\ybar\plusl\lambda\uu_j)\cdot\uu_i}
  \nonumber
  \\
  &\leq
  \ellbar_i\bigParens{(\ybar\plusl\lambda\uu_j)\cdot\uu_i}
  <+\infty.
  \label{eqn:thm:hard-core:3:2}
\end{align}
The first inequality is because $\ellbar_i\geq 0$
since $\ell_i\geq 0$.
The second inequality is because
$\limray{\vv}\cdot\uu_i\leq 0$
and $\ellbar_i$ is nondecreasing.
The last inequality is because $\ybar\cdot\uu_i\in\R$,
implying $(\ybar\plusl\lambda\uu_j)\cdot\uu_i\in\R$
(and using Eq.~\ref{eqn:hard-core:1}).

In particular, when $i=j$,
\[
  \ellbar_j\bigParens{(\limray{\vv}\plusl\ybar\plusl\lambda\uu_j)\cdot\uu_j}
  =
  \ellbar_j\bigParens{\limray{\vv}\cdot\uu_j
              \plusl(\ybar\plusl\lambda\uu_j)\cdot\uu_j}
  =
  \ellbar_j(-\infty)
  = 0,
\]
since $\limray{\vv}\cdot\uu_j=-\infty$
(and by Eq.~\ref{eqn:hard-core:1}).
On the other hand,
\[
  \ellbar_j\bigParens{(\ybar\plusl\lambda\uu_j)\cdot\uu_j}
  =
  \ellbar_j(\ybar\cdot\uu_j \plusl \lambda\uu_j\cdot\uu_j)
  \rightarrow +\infty
\]
as $\lambda\rightarrow+\infty$,
since $\ybar\cdot\uu_j\in\R$ and $\uu_j\neq\zero$ (since $\vv\cdot\uu_j<0$),
and by
\Cref{pr:conv-inc:prop}(\ref{pr:conv-inc:infsup},\ref{pr:conv-inc:nonconst}).
Thus, \eqref{eqn:thm:hard-core:3:2} holds for all
$i\in\indset$,
and furthermore
the inequality is strict when $i=j$ and
when $\lambda$ is sufficiently large.

Therefore,
for $\lambda$ sufficiently large, we have shown that
\begin{align*}
  \fullshad{f}(\lambda\uu_j)
  \leq
   \fext(\limray{\vv}\plusl \ybar \plusl \lambda\uu_j)
   &=
   \sum_{i\in\indset}
   \ellbar_i\bigParens{(\limray{\vv}\plusl\ybar\plusl\lambda\uu_j)\cdot\uu_i}
   \nonumber
   \\
   &<
   \sum_{i\in\indset}
   \ellbar_i\bigParens{(\ybar\plusl\lambda\uu_j)\cdot\uu_i}
   \nonumber
   \\
   &=
   \fext(\ybar \plusl \lambda\uu_j)
    \leq
   \fullshad{f}(\lambda\uu_j).
\end{align*}
The equalities are by
\Cref{pr:hard-core:1}(\ref{pr:hard-core:1:a-ext}).
The second inequality is by the argument above.
The final inequality is by
\Crefequiv{thm:conv-univ-equiv}{thm:conv-univ-equiv:a}{thm:conv-univ-equiv:b}
since $\ybar\in\conv(\unimin{f})$.
For the first inequality, note that
because $\limray{\vv}$ and $\ybar$ are both in
$\aresconef$, their leftward sum $\limray{\vv}\plusl\ybar$ is as well;
this is because $\aresconef$ is a convex astral cone
(by \Cref{cor:res-fbar-closed}),
which is therefore closed under sequential sum
(by \Cref{thm:ast-cone-is-cvx-if-sum}), and so also under leftward
addition
(by \Cref{cor:seqsum-conseqs}\ref{cor:seqsum-conseqs:a}).
The first inequality therefore follows from
\Cref{pr:fullshad-equivs}.

Thus, having reached a contradiction, we conclude that
if $\ybar\in\conv(\unimin{f})$, then
\eqref{eq:hard-core:conv:2} is satisfied for
$i\in\indset$.

For the converse, suppose
\eqref{eq:hard-core:conv:2} is satisfied for
all $i\in\indset$.
Let $\ebar$ be any point in $\unimin{f}$ (which is nonempty by
\Cref{pr:new:thm:f:8a}\ref{pr:new:thm:f:8a:nonemp-closed}).
Then, as just argued, $\ebar$ satisfies \eqref{eq:hard-core:conv:2} as
well, so $\ebar\cdot\uu_i=\ybar\cdot\uu_i$
for $i\in\indset$.
So for all $\xx\in\Rn$,
$\fext(\ybar\plusl\xx)=\fext(\ebar\plusl\xx)=\fullshad{f}(\xx)$
with the first equality from
\Cref{pr:hard-core:1}(\ref{pr:hard-core:1:a-ext}),
and the second because $\ebar\in\unimin{f}$.
Therefore, $\ybar\in\conv(\unimin{f})$
by
\Crefequiv{thm:conv-univ-equiv}{thm:conv-univ-equiv:c}{thm:conv-univ-equiv:a}.

\pfpart{Part~(\ref{thm:hard-core:3:b:univ}):}
By \Cref{cor:conv-univ-icons}, the icon $\ebar$ is in
$\unimin{f}$ if
and only if it is in $\conv(\unimin{f})$.
Combining with part~(\ref{thm:hard-core:3:b:conv}), this proves the
claim.

\pfpart{Part~(\ref{thm:hard-core:3:c}):}
Let $\ebar$ be any point in $\unimin{f}$
(which exists by
\Cref{pr:new:thm:f:8a}\ref{pr:new:thm:f:8a:nonemp-closed}),
and let
$\zbar=\ebar\plusl\ybar$.
Then $\zbar\cdot\uu_i<0$ since
$\ebar\cdot\uu_i=0$
by part~(\ref{thm:hard-core:3:b:univ}), so
$\zbar\not\in\aresconef$,
by part~(\ref{thm:hard-core:3:a}).
Therefore, for some
$j\in\indset$,
$\zbar\cdot\uu_j>0$, by
\Cref{pr:hard-core:1}(\ref{pr:hard-core:1:b-cones}).
Further, it must be that $j\in\hardcore{f}$ since
otherwise
part~(\ref{thm:hard-core:3:b:univ}) would imply
$\ebar\cdot\uu_j=-\infty$, so that also
$\zbar\cdot\uu_j=-\infty$.
Thus,
$\ybar\cdot\uu_j>0$
since $\ebar\cdot\uu_j=0$
by part~(\ref{thm:hard-core:3:b:univ}).

\pfpart{Part~(\ref{thm:hard-core:3:d}):}
Let $\ebar$ be any point in $\unimin{f}$
(which exists by
\Cref{pr:new:thm:f:8a}\ref{pr:new:thm:f:8a:nonemp-closed}).
Then for $\xx\in\Rn$,
\[
  \fullshad{f}(\xx)
  =
  \fshadd(\xx)
  =
  \fext(\ebar\plusl\xx)
  =
  \sum_{i\in\indset} \ellbar_i(\ebar\cdot\uu_i\plusl\xx\cdot\uu_i)
  =
  \sum_{i\in\hardcore{f}} \ell_i(\xx\cdot\uu_i).
\]
The first equality is because $\ebar\in\unimin{f}$.
The third is by
\Cref{pr:hard-core:1}(\ref{pr:hard-core:1:a-ext}).
The fourth is by part~(\ref{thm:hard-core:3:b:univ})
and \eqref{eqn:hard-core:1}.

\pfpart{Part~(\ref{thm:hard-core:3:e}):}
Let $U=\uset{\hardcore{f}}$.
We aim to show $\rescperp{f}=\spn U$.

If $i\in\hardcore{f}$, then $\uu_i\in\rescperp{f}$
by definition, %
so $U\subseteq\rescperp{f}$, implying
$\spn U\subseteq\rescperp{f}$ since $\rescperp{f}$ is a linear
subspace
(\Cref{pr:std-perp-props}\ref{pr:std-perp-props:a}).

For the reverse inclusion, suppose $\yy\in \Uperp$,
meaning $\yy\cdot\uu_i=0$ for all $i\in\hardcore{f}$.
By applying
\Cref{pr:hard-core:1}(\ref{pr:hard-core:1:b-cones})
to $\fullshad{f}$, whose form is given in
part~(\ref{thm:hard-core:3:d}),
this shows that $\yy\in\resc{\fullshad{f}}$.
Thus, $\Uperp\subseteq\resc{\fullshad{f}}$, so
\[
  \rescperp{f}
  =
  \perpresf
  =
  \rescperp{\fullshad{f}}
  \subseteq
  \Uperperp
  =
  \spn U.
\]
The first equality is by
\Cref{pr:hard-core:1}(\ref{pr:hard-core:1:b-resc2}),
and the second by
\Cref{thm:f:4x}(\ref{thm:f:4xa}).
The inclusion and last equality are by
\Cref{pr:std-perp-props}(\ref{pr:std-perp-props:b},\ref{pr:std-perp-props:c}).
\qedhere
\end{proof-parts}
\end{proof}

\begin{example}[Two exps and a square, continued]
\label{ex:erm-running-eg2}
\indexg{Two exps and a square!universal reduction of|(}%
Continuing \Cref{ex:erm-running-eg1},
\Cref{thm:hard-core:3}(\ref{thm:hard-core:3:d}) implies
that $\fullshad{f}$, the universal reduction of $f$, is
\[
  \fullshad{f}(\xx)
  =
  \sum_{i\in\hardcore{f}} \ell_i(\xx\cdot\uu_i)
  =
  \sqp(2+x_2-x_3) + \sqp(-2-x_2+x_3)
  =
  (2+x_2-x_3)^2,
\]
as was previously noted in
\Cref{ex:simple-eg-exp-exp-sq-part2,ex:simple-eg-exp-exp-sq-2}.
Moreover, \Cref{thm:hard-core:3}(\ref{thm:hard-core:3:e}) implies
that $(\conssp{\fullshad{f}})^\perp=\perpresf=\rescperp{f}$
is the linear subspace (in this case a line)
spanned by the set of points $\uu_i$ with \irredIndices, that is, by $\{\uu_3,\uu_4\}$,
as we saw in
\eqref{eq:simple-eg-exp-exp-sq:cons-perp} of
\Cref{ex:simple-eg-exp-exp-sq-2}.%
\indexg{Two exps and a square!universal reduction of|)}%
\end{example}

\indexg{empirical risk functions!canonical minimizers of|(}%
For functions of the form given in \eqref{eqn:loss-sum-form},
\Cref{thm:hard-core:3}
provides a characterization of all canonical minimizers
in terms of \irredIndices,
namely, as points $\xbar=\ebar\plusl\qq$
whose finite part $\qq\in\Rn$ minimizes $\fullshad{f}$
in part~(\ref{thm:hard-core:3:d}),
and
whose iconic part $\ebar\in\corezn$
satisfies the conditions for a universal reducer
in part~(\ref{thm:hard-core:3:b:univ}).
The theorem
shows that all universal reducers $\ebar$ are identical in terms of
the values of $\ebar\cdot\uu_i$,
for $i\in\indset$,
as determined by whether each $i$ is reducible or persistent.
Every universal reducer has the effect of causing the terms
with reducible indices
in \eqref{eqn:loss-sum-form} to vanish.
The remaining terms, those with persistent indices, constitute
exactly the universal reduction $\fullshad{f}$.
As is generally the case (\Cref{thm:f:4x}\ref{thm:f:4xa}\ref{thm:f:4xb}),
all sublevel sets of this function are compact when restricted to
$(\conssp{\fullshad{f}})^\perp=\perpresf=\rescperp{f}$.%
\indexg{empirical risk functions!canonical minimizers of|)}%
\indexg{empirical risk functions!universal reducers of|)}%
\indexg{empirical risk functions!universal reduction of|)}%
\indexg{persistent indices!universal reducers and reduction and|)}%

\indexg{empirical risk functions!rank of minimizers of|(}%
We also remark that since $f$ is recessive complete
(\Cref{pr:hard-core:1}\ref{pr:hard-core:1:b-resc1}),
\Cref{thm:unimin-can-be-rankone}
implies that $\ef$
can always be minimized at some point of astral rank at most one,
specifically, at any point
$\xbar=\limray{\vv}\plusl\qq$ where $\vv\in\ric{f}$ and $\qq$ minimizes
$\fullshad{f}$.%
\indexg{empirical risk functions!rank of minimizers of|)}%

\indexg{empirical risk functions!faces of associated convex hull|(}%
\indexg{persistent indices!faces of associated convex hull|(}%
We look next at a geometric characterization of \irredIndices,
focusing on the faces of
$S=\conv(\uset{\indset})$,
the convex hull of all the points $\uu_i$.
(See \Cref{sec:prelim:faces} for a brief introduction to the
faces of a convex set.)
The next theorems show that the \irredIndexSet
is fully determined by the location of the origin relative to $S$.
In particular, if the origin is not included in $S$, then the
\irredIndexSet must be empty (implying that $\fext$ is minimized by
a point $\xbar$ for which $\xbar\cdot\uu_i=-\infty$ for all $i\in\indset$).
Otherwise, the origin must be in $\ri{C}$ for exactly one face $C$,
and the \irredIndexSet consists of indices of all
points $\uu_i$ included in $C$.
Alternatively, we can say that $\conv(\uset{\hardcore{f}})$ is a face
of $S$, and is specifically the smallest face that includes the origin
(meaning that it is included in all other faces that include the origin).

\begin{figure}
  \centering
  \includegraphics{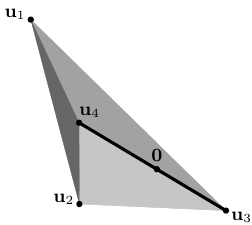}
  \mycaption{Convex hull of points $\uu_1,\uu_2,\uu_3,\uu_4$ from \Cref{ex:erm-running-eg1}}{%
\indexf{Two exps and a square!associated convex hull}%
    These points alongside loss functions $\ell_1,\ell_2,\ell_3,\ell_4$ describe the convex function
    $f$ in \eqref{eq:simple-eg-exp-exp-sq:modform}.
    Since the origin is in the relative interior of the edge connecting $\uu_3$ and $\uu_4$,
    the persistent index set for $f$ is $\set{3,4}$ (see \Cref{thm:erm-faces-hardcore2}).
  }%
  \label{fig:tetra}%
\end{figure}

\begin{example}[Two exps and a square, continued]
\label{ex:erm-running-eg1:cont}
\indexg{Two exps and a square!associated convex hull|(}%
In the setup of \Cref{ex:erm-running-eg1},
the convex hull of the points
$\uu_1,\ldots,\uu_4$ is an (irregular) tetrahedron in $\R^3$
depicted in \Cref{fig:tetra}.
Its faces consist of the tetrahedron itself, its four triangular
faces, six edges, four vertices, and the empty set.
As shown in the figure,
the origin is in the relative interior of the edge connecting $\uu_3$
and $\uu_4$ (since $\zero=\sfrac{1}{2}\uu_3+\sfrac{1}{2}\uu_4$),
corresponding to the \irredIndexSet being $\{3,4\}$ in this case.
That edge is indeed the smallest face that includes the origin.%
\indexg{Two exps and a square!associated convex hull|)}%
\end{example}

\begin{theorem}  \label{thm:erm-faces-hardcore}
  Let $f:\Rn\rightarrow\R$ have the form given in
  \eqref{eqn:loss-sum-form}, where,
  for $i\in\indset$,
  $\uu_i\in\Rn$
  and
  $\ell_i:\R\rightarrow\R$ is convex, nondecreasing, not constant,
  with $\inf \ell_i = 0$.
  Let $S=\conv{\uset{I}}$.
  Then:
  \begin{letter-compact}
  \item  \label{thm:erm-faces-hardcore:b}
    $\conv(\uset{\hardcore{f}})$ is a face of $S$.
  \item  \label{thm:erm-faces-hardcore:a}
    Let $J\subseteq I$, and suppose
    $\zero\in\ri(\conv{\uset{J}})$.
    Then $J\subseteq\hardcore{f}$.
  \item  \label{thm:erm-faces-hardcore:aa}
    Let $C$ be a face of $S$, and suppose $\zero\in C$.
    Then
    $\zero\in\conv(\uset{\hardcore{f}})\subseteq C$.
  \item  \label{thm:erm-faces-hardcore:c}
    $\zero\in S$ if and only if $\hardcore{f}\neq\emptyset$.
  \end{letter-compact}
\end{theorem}

\begin{proof}
~

\begin{proof-parts}
\pfpart{Part~(\ref{thm:erm-faces-hardcore:b}):}
Let $U=\uset{\hardcore{f}}$, and
let $C=S \cap (\spn U)$.
We show first that $C$ is a face of $S$,
and later show $C=\conv U$.

Let $\xx,\zz\in S$ and $\lambda\in (0,1)$.
Let $\ww=(1-\lambda)\xx+\lambda\zz$.
Assume $\ww\in C$, which we aim to show implies that $\xx$ and $\zz$
are also in $C$,
since this will prove that $C$ is a face of $S$.
Let $\yy\in\resc{f}$, implying
$\yy\cdot\uu_i\leq 0$ for $i\in\indset$,
by \Cref{pr:hard-core:1}(\ref{pr:hard-core:1:b-cones}).
Since $\xx$ and $\zz$ are in $S$, the convex hull of the
$\uu_i$'s, this implies
$\yy\cdot\xx\leq 0$ and
$\yy\cdot\zz\leq 0$.
Also, for all $i\in\hardcore{f}$,
$\uu_i\in\rescperp{f}$ so $\yy\cdot\uu_i=0$.
Since $\ww\in\spn U$, this also means $\yy\cdot\ww=0$.
Thus,
\[
  0 = \yy\cdot\ww = (1-\lambda)(\yy\cdot\xx) + \lambda(\yy\cdot\zz).
\]
Since $\lambda\in (0,1)$ and the two terms on the right are
nonpositive, we must have $\yy\cdot\xx=\yy\cdot\zz=0$.
Therefore, $\xx,\zz\in\rescperp{f}$, since
this holds for all $\yy\in\resc{f}$.
Thus,
\[
   \xx,\zz
   \in S\cap\rescperp{f} = S\cap(\spn{U})=C,
\]
where the first equality is because
$\rescperp{f}=\spn U$,
by \Cref{thm:hard-core:3}(\ref{thm:hard-core:3:e}).

We have shown $C$ is a face of $S$.
As such, $C$ is equal to the convex hull of the points in
$\uset{\indset}$ that are included in $C$,
by \Cref{roc:thm18.3}.
Moreover,
a point $\uu_i$, for $i\in\indset$, is included in $C$
if and only if it is in $\spn U=\rescperp{f}$,
that is, if and only if $i\in\hardcore{f}$.
We conclude that $C=\conv{\uset{\hardcore{f}}}$, completing the proof.

\pfpart{Part~(\ref{thm:erm-faces-hardcore:a}):}
Let $j\in J$, and let $C=\conv(\uset{J})$.
Since $\zero\in\ri{C}$,
and since $\uu_j\in C$,
there exists $\delta>0$ for which the point
$\ww=(1+\delta)\zero-\delta\uu_j=-\delta\uu_j$
is also in $C$
(\Cref{roc:thm6.4}).

Let $\yy\in\resc{f}$, implying
$\yy\cdot\uu_i\leq 0$ for all $i\in\indset$,
by \Cref{pr:hard-core:1}(\ref{pr:hard-core:1:b-cones});
in particular, $\yy\cdot\uu_j\leq 0$.
Also, $\ww$ is in $C$ and therefore a convex combination of points in
$\uset{J}$.
Thus, $\yy\cdot(-\delta\uu_j)=\yy\cdot\ww\leq 0$ as well.
Together, these imply $\yy\cdot\uu_j=0$ since $\delta>0$.
Since this holds for all $\yy\in\resc{f}$, we have shown that
$\uu_j\in\rescperp{f}$, that is, $j\in\hardcore{f}$.

\pfpart{Part~(\ref{thm:erm-faces-hardcore:aa}):}
Let $F=\conv(\uset{\hardcore{f}})$, which is a face of $S$ by
part~(\ref{thm:erm-faces-hardcore:b}).
We show first that $\zero\in F$.
Since $\zero\in C\subseteq S$,
and since the relative interiors of faces of $S$ form a partition
(\Cref{roc:thm18.2}),
there must exist a face $D$ of $S$ for which
$\zero\in\ri{D}$.
Let $J=\{i\in\indset : \uu_i\in D\}$.
Then $D=\conv(\uset{J})$,
by \Cref{roc:thm18.3}.
From
part~(\ref{thm:erm-faces-hardcore:a}),
$J\subseteq\hardcore{f}$, so
$\zero\in D\subseteq F$, as claimed.

We next show $F\subseteq C$.
Suppose not.
Let $C'=F\cap C$, which is a face of $S$
since both $F$ and $C$ are faces
(\Cref{pr:face-props}\ref{pr:face-props:intersect}).
Also, $\zero\in C'$, but $F\not\subseteq C'$
since we have assumed $F\not\subseteq C$.

Because $F$ and $C'$ are distinct faces of $S$, their relative
interiors are disjoint
so that $(\ri{F})\cap(\ri{C'})=\emptyset$
(by \Cref{pr:face-props}\ref{pr:face-props:cor18.1.2}).
As a result, there exists a hyperplane properly separating
$F$ and $C'$
(by \Cref{roc:thm11.3}).
That is, there exist $\vv\in\Rn$ and $\beta\in\R$ for which
$\vv\cdot\ww\leq \beta$ for all $\ww\in F$
and
$\vv\cdot\ww\geq \beta$ for all $\ww\in C'$.
Since $C'\subseteq F$, this actually implies
$\vv\cdot\ww = \beta$ for all $\ww\in C'$.
Moreover, because $\zero\in C'$, we must have $\beta=0$.

Furthermore, this hyperplane \emph{properly} separates these sets,
meaning there must exist a point in $F\cup C'$ not in the separating hyperplane itself.
Since $C'$ is entirely included in the hyperplane, this implies there
must be a point $\zz\in F$ for which $\vv\cdot\zz<0$.
Since $\zz\in\conv{\uset{\hardcore{f}}}$,
it must be a convex combination of points $\uu_i$, for
$i\in\hardcore{f}$.
Therefore, there must exist some $i\in\hardcore{f}$ with
$\vv\cdot\uu_i<0$.
By
\Cref{thm:hard-core:3}(\ref{thm:hard-core:3:c}),
this implies there exists a point $j\in\hardcore{f}$ with
$\vv\cdot\uu_j>0$, contradicting that
$\vv\cdot\ww\leq 0$ for all $\ww\in F$.

\pfpart{Part~(\ref{thm:erm-faces-hardcore:c}):}
Suppose $\zero\in S$, and that, contrary to the claim,
$\hardcore{f}=\emptyset$.
Let $\ebar\in\unimin{f}$
(which must exist by
\Cref{pr:new:thm:f:8a}\ref{pr:new:thm:f:8a:nonemp-closed}).
Then $\ebar\cdot\uu_i=-\infty$ for all $i\in\indset$,
by \Cref{thm:hard-core:3}(\ref{thm:hard-core:3:b:univ}).
Since $\zero\in S$,
$\zero$ is a convex combination of the $\uu_i\negKern$'s,
implying
$\ebar\cdot\zero=-\infty$ using
\Cref{pr:i:1}, a contradiction.

For the converse, suppose $\zero\not\in S$, and, contrary to the
claim, that $\hardcore{f}\neq\emptyset$.
Then because both $S$ and $\{\zero\}$ are convex,
closed (in $\Rn$), and bounded,
there exists a hyperplane strongly separating them
by \Cref{roc:cor11.4.2}.
That is, there exists $\vv\in\Rn$ for which
\[
  \sup_{\ww\in S} \vv\cdot\ww < \vv\cdot\zero = 0
\]
by \Cref{roc:thm11.1}.
In particular, this means $\vv\cdot\uu_i<0$ for all $i\in\indset$.
Let $i\in\hardcore{f}$, which we have assumed is not empty.
Then by
\Cref{thm:hard-core:3}(\ref{thm:hard-core:3:c}),
because $\vv\cdot\uu_i<0$,
there also must exist $j\in\hardcore{f}$ with $\vv\cdot\uu_j>0$, a
contradiction.
\qedhere
\end{proof-parts}
\end{proof}

\indexg{faces!smallest face including set|(}%
For a convex set $S\subseteq\Rn$ and a set $A\subseteq S$, we say that
a face $C$ of $S$ is the
\emph{smallest face of $S$ that includes $A$}
if $A\subseteq C$ and if for all faces $C'$ of $S$,
if $A\subseteq C'$ then $C\subseteq C'$.
Equivalently, the smallest face of $S$ that includes $A$ is the
intersection of all of the faces of $S$ that include $A$, which is itself a
face
(\Cref{pr:face-props}\ref{pr:face-props:intersect}).%
\indexg{faces!smallest face including set|)}%

\begin{theorem}  \label{thm:erm-faces-hardcore2}
  Let $f:\Rn\rightarrow\R$ have the form given in
  \eqref{eqn:loss-sum-form}, where,
  for $i\in\indset$,
  $\uu_i\in\Rn$
  and
  $\ell_i:\R\rightarrow\R$ is convex, nondecreasing, not constant,
  with $\inf \ell_i = 0$.
  Let $C$ be any nonempty face of $S=\conv(\uset{I})$.
  Then the following are equivalent:
  \begin{letter-compact}
  \item  \label{thm:erm-faces-hardcore2:a}
    $\zero\in\ri{C}$.
  \item  \label{thm:erm-faces-hardcore2:b}
    $C$ is the smallest face of $S$ that includes $\{\zero\}$.
  \item  \label{thm:erm-faces-hardcore2:c}
    $\hardcore{f}=\{i\in\indset :\: \uu_i\in C\}$.
  \item  \label{thm:erm-faces-hardcore2:d}
    $C = \conv(\uset{\hardcore{f}})$.
  \end{letter-compact}
\end{theorem}

\begin{proof}
~

\begin{proof-parts}
\pfpart{%
  (\ref{thm:erm-faces-hardcore2:a})
  $\Rightarrow$
  (\ref{thm:erm-faces-hardcore2:b}):}
Suppose $\zero\in\ri{C}$.
Then clearly $\zero\in C$.
If $\zero$ is included in some face $C'$ of $S$,
then $C'$ and $\ri{C}$ are not disjoint,
implying $C\subseteq C'$
by \Cref{pr:face-props}(\ref{pr:face-props:thm18.1}).

\pfpart{%
  (\ref{thm:erm-faces-hardcore2:b})
  $\Rightarrow$
  (\ref{thm:erm-faces-hardcore2:c}):}
Suppose (\ref{thm:erm-faces-hardcore2:b}) holds.
Let $J=\{i\in\indset : \uu_i\in C\}$.
By
\Cref{thm:erm-faces-hardcore}(\ref{thm:erm-faces-hardcore:aa}),
$\zero\in\conv(\uset{\hardcore{f}})\subseteq C$,
since $\zero\in C$.
Therefore, $\hardcore{f}\subseteq J$.

For the reverse inclusion, since $\zero\in S$,
by \Cref{roc:thm18.2},
there exists a face $C'$ of $S$ with
$\zero\in\ri{C'}$.
By assumption, this implies
$C\subseteq C'$, and so
$J\subseteq J'$
where $J'=\{i\in\indset : \uu_i\in C'\}$.
Furthermore, $J'\subseteq \hardcore{f}$ by
\Cref{thm:erm-faces-hardcore}(\ref{thm:erm-faces-hardcore:a}).
Combining yields $J=\hardcore{f}$ as claimed.

\pfpart{%
  (\ref{thm:erm-faces-hardcore2:c})
  $\Rightarrow$
  (\ref{thm:erm-faces-hardcore2:d}):}
This is immediate
by \Cref{roc:thm18.3}.

\pfpart{%
  (\ref{thm:erm-faces-hardcore2:d})
  $\Rightarrow$
  (\ref{thm:erm-faces-hardcore2:a}):}
Suppose $C = \conv(\uset{\hardcore{f}})$.
Since $C$ is not empty, $\hardcore{f}\neq\emptyset$, so
$\zero\in S$ by
\Cref{thm:erm-faces-hardcore}(\ref{thm:erm-faces-hardcore:c}).
Therefore,
by \Cref{roc:thm18.2},
there exists a face $C'$ of $S$ for which
$\zero\in\ri{C'}$.
That is, $C'$ satisfies
(\ref{thm:erm-faces-hardcore2:a}),
and so also
(\ref{thm:erm-faces-hardcore2:d}),
by the implications proved above.
Thus,
$C'=\conv(\uset{\hardcore{f}})=C$,
so $\zero\in\ri{C}$.%
\indexg{empirical risk functions!faces of associated convex hull|)}%
\indexg{persistent indices!faces of associated convex hull|)}%
\qedhere
\end{proof-parts}
\end{proof}

Notice that the sets we have been considering, namely, the standard
and astral recession cones, the set of universal reducers, as well
as the set of \irredIndices, all depend exclusively
on the $\uu_i\negKern$'s,
and are entirely {independent} of the specific functions~$\ell_i$.
In other words, suppose we form a new function $f'$ as in
\eqref{eqn:loss-sum-form}
with the $\uu_i\negKern$'s
unchanged, but
with each $\ell_i$ replaced by some other function $\ell'_i$ (though
still satisfying the same assumed properties).
Then the sets listed above are unchanged.
That is,
$\resc{f}=\resc{f'}$ and
$\aresconef=\aresconefp$
by \Cref{pr:hard-core:1}(\ref{pr:hard-core:1:b-cones});
$\hardcore{f}=\hardcore{f'}$ by the definition in
\eqref{eq:hard-core:3}; and
$\unimin{f}=\unimin{f'}$ by
\Cref{thm:hard-core:3}(\ref{thm:hard-core:3:b:univ}).

\indexg{empirical risk functions!minimizers of|(}%
From \Cref{pr:unimin-to-global-min},
$\fext$ is minimized at a point
$\xbar=\ebar\plusl\qq$
if $\ebar\in\unimin{f}$ and $\qq$ minimizes $\fullshad{f}$.
If, in addition to our preceding assumptions,
each $\ell_i$ is \emph{strictly} increasing, then these
conditions are not only sufficient but also necessary for $\xbar$ to
minimize $\fext$.
Furthermore, if each function $\ell_i$ is \emph{strictly} convex,
then $\fullshad{f}$
is uniquely minimized over the linear subspace
$\perpresf = \rescperp{f}$, implying also that,
for all $i\in\indset$,
the value of $\xbar\cdot\uu_i$ will be the same for all of $\fext$'s
minimizers $\xbar$.
We show these in the next two theorems.

\begin{theorem}  \label{thm:waspr:hard-core:2}
  Let $f:\Rn\rightarrow\R$ have the form given in
  \eqref{eqn:loss-sum-form}, where,
  for $i\in\indset$,
  $\uu_i\in\Rn$  %
  and
  $\ell_i:\R\rightarrow\R$ is convex,
  with $\inf \ell_i = 0$.
  Suppose further that each $\ell_i$
  is strictly increasing
  (as will be the case if each $\ell_i$ is nondecreasing and strictly convex).
  Let $\xbar=\ebar\plusl\qq$ where $\ebar\in\corezn$ and $\qq\in\Rn$.
  Then $\xbar$ minimizes $\fext$ if and only if $\ebar\in\unimin{f}$
  and $\qq$ minimizes $\fullshad{f}$.
\end{theorem}

\begin{proof}
The ``if'' direction follows from
\Cref{pr:unimin-to-global-min}.
For the converse, suppose $\xbar=\ebar\plusl\qq$ minimizes $\fext$.
By \Cref{pr:min-fullshad-is-finite-min}, this implies that
$\qq$ minimizes $\fullshad{f}$.

To show $\ebar\in\unimin{f}$, note first
that $\ebar\in\aresconef$
by \Cref{thm:arescone-fshadd-min},
since $\xbar$ minimizes $\fext$.
Therefore,
by
\Crefequiv{pr:icon-equiv}{pr:icon-equiv:a}{pr:icon-equiv:b}
and \Cref{pr:hard-core:1}(\ref{pr:hard-core:1:b-cones}),
$\ebar\cdot\uu_i\in\{-\infty,0\}$
for $i\in\indset$
since
$\ebar\in\corezn$.
Let $J=\{ i \in\indset :\: \ebar\cdot\uu_i = 0\}$.

We claim that $J=\hardcore{f}$.
The inclusion $\hardcore{f}\subseteq J$ follows directly from
\Cref{thm:hard-core:3}(\ref{thm:hard-core:3:a})
since $\ebar\in\aresconef$.
To show these sets are actually equal,
suppose by way of contradiction
that %
$\hardcore{f}$ is a proper subset of $J$.
Also, let $\ebar'$ be any point in $\unimin{f}$
(which exists by
\Cref{pr:new:thm:f:8a}\ref{pr:new:thm:f:8a:nonemp-closed}).
Then
\begin{align*}
  \fext(\ebar\plusl\qq)
  =
  \sum_{i\in\indset}
     \ellbar_i(\ebar\cdot\uu_i\plusl\qq\cdot\uu_i)
  &=
  \sum_{i\in J} \ell_i(\qq\cdot\uu_i)
  \\
  &>
  \sum_{\mathclap{i\in \hardcore{f}}} \ell_i(\qq\cdot\uu_i)
  =
  \fullshad{f}(\qq)
  =
  \fext(\ebar'\plusl\qq).
\end{align*}
The first two equalities are from
\Cref{pr:hard-core:1}(\ref{pr:hard-core:1:a-ext})
and \eqref{eqn:hard-core:1}.
The inequality is because
$\hardcore{f}$ is a proper subset of $J$, and because
$\ell_j(\qq\cdot\uu_j)>\ellbar_j(-\infty)=0$
for all $j\in J\setminus(\hardcore{f})$
since $\ell_j$ is strictly increasing
(and also by Eq.~\ref{eqn:hard-core:1}).
The third equality is by
\Cref{thm:hard-core:3}(\ref{thm:hard-core:3:d}),
and the last equality is because $\ebar'\in\unimin{f}$.
This contradicts that $\ebar\plusl\qq$ minimizes
$\fext$.

Thus, $J=\hardcore{f}$, and therefore $\ebar\in\unimin{f}$ by
\Cref{thm:hard-core:3}(\ref{thm:hard-core:3:b:univ}).

Finally, we note that if $\ell_i$ is nondecreasing and strictly
convex, then it is also strictly increasing.
Otherwise, there would exist real numbers $x<y$ for which
$\ell_i(x)=\ell_i(y)$.
Letting $z=(x+y)/2$, this implies, by strict convexity, that
$\ell_i(z)<(\ell_i(x)+\ell_i(y))/2=\ell_i(x)$.
Thus, $x<z$, but $\ell_i(x)>\ell_i(z)$, contradicting that $\ell_i$ is
nondecreasing.
\end{proof}

\begin{theorem}  \label{thm:waspr:hard-core:4}
  Let $f:\Rn\rightarrow\R$ have the form given in
  \eqref{eqn:loss-sum-form}, where,
  for $i\in\indset$,
  $\uu_i\in\Rn$  %
  and
  $\ell_i:\R\rightarrow\R$ is nondecreasing and strictly convex,
  with $\inf \ell_i = 0$.
  Then~$\fullshad{f}$, if restricted to
  $(\conssp{\fullshad{f}})^\perp=\perpresf = \rescperp{f}$,
  has a unique minimizer $\qq$.
  Furthermore, the following are equivalent, for $\xbar\in\extspace$:
  \begin{letter-compact}
  \item  \label{pr:hard-core:4:a}
    $\xbar$ minimizes $\fext$.
  \item  \label{pr:hard-core:4:b}
    $\xbar=\zbar\plusl\qq$ for some $\zbar\in\conv(\unimin{f})$.
  \item  \label{pr:hard-core:4:c}
    For $i\in\indset$,
    \[
       \xbar\cdot\uu_i
       =
       \begin{cases}
         \qq\cdot\uu_i   & \text{if $i\in\hardcore{f}$,} \\
         -\infty & \text{otherwise.} \\
       \end{cases}
    \]
  \end{letter-compact}
\end{theorem}

\begin{proof}
That
$(\conssp{\fullshad{f}})^\perp=\perpresf = \rescperp{f}$
follows from
\Cref{thm:f:4x}(\ref{thm:f:4xa})
and
\Cref{pr:hard-core:1}(\ref{pr:hard-core:1:b-resc2}).

Let $\qq$ be a minimizer of $\fullshad{f}$ in $\rescperp{f}$,
which must exist by
\Cref{thm:f:4x}(\ref{thm:f:4xd}).
Suppose, by way of contradiction, that some other point
$\qq'\in\rescperp{f}$ also minimizes $\fullshad{f}$,
with $\qq\neq\qq'$.

We claim first that
$\qq\cdot\uu_i\neq\qq'\cdot\uu_i$
for some $i\in\hardcore{f}$.
Suppose to the contrary that $\dd\cdot\uu_i=0$ for all
$i\in\hardcore{f}$, where $\dd=\qq'-\qq\neq\zero$.
Then because $\rescperp{f}$ is a linear subspace
(\Cref{pr:std-perp-props}\ref{pr:std-perp-props:a}),
for all $\lambda\in\R$,
$\lambda\dd\in\rescperp{f}$.
Furthermore, by
\Cref{thm:hard-core:3}(\ref{thm:hard-core:3:d}),
$\fullshad{f}(\lambda\dd)=\fullshad{f}(\zero)\in\R$.
In other words, the entire line
$\{ \lambda\dd : \lambda\in\R\}$
is included in the set
$\{\xx\in\rescperp{f} :\: \fullshad{f}(\xx)\leq\fullshad{f}(\zero)\}$.
This, however, is a contradiction since the latter set is compact
by \Cref{thm:f:4x}(\ref{thm:f:4xb}), and therefore bounded.

So let $i\in\hardcore{f}$ be such that
$\qq\cdot\uu_i\neq\qq'\cdot\uu_i$.
Let $\zz=(\qq+\qq')/2$.
Since each $\ell_j$ is convex,
$\ell_j(\zz\cdot\uu_j)\leq(\ell_j(\qq\cdot\uu_j)+\ell_j(\qq'\cdot\uu_j))/2$.
Furthermore,
when $j=i$, by strict convexity of $\ell_i$, this inequality is strict.
Therefore, applying
\Cref{thm:hard-core:3}(\ref{thm:hard-core:3:d})
yields
$\fullshad{f}(\zz) < (\fullshad{f}(\qq) + \fullshad{f}(\qq'))/2$,
contradicting the assumption that $\qq$ and $\qq'$ both
minimize~$\fullshad{f}$.

Thus, $\qq$ is the only minimizer of $\fullshad{f}$ in $\rescperp{f}$.
We next prove the stated equivalences:

\begin{proof-parts}

\pfpart{%
  (\ref{pr:hard-core:4:a})
  $\Rightarrow$
  (\ref{pr:hard-core:4:b})
}:
Suppose $\xbar$ minimizes $\fext$.
Then
by \Cref{thm:waspr:hard-core:2},
$\xbar=\ebar\plusl\yy$ for some $\ebar\in\unimin{f}$
and some $\yy\in\Rn$ that minimizes $\fullshad{f}$.
Since $\conssp{\fullshad{f}}$ is a linear subspace
(\Cref{pr:prelim:const-props}\ref{pr:prelim:const-props:b}),
we can write
$\yy=\yy'+\yy''$
for some
$\yy'\in(\conssp{\fullshad{f}})^\perp$
and $\yy''\in\conssp{\fullshad{f}}$
(\Cref{pr:lin-decomp}).
This last fact implies
$\fullshad{f}(\yy)=\fullshad{f}(\yy')$,
so $\yy'$ also minimizes $\fullshad{f}$.
Therefore, $\yy'=\qq$ since, as already shown, $\qq$ is the only
minimizer of $\fullshad{f}$ in $\rescperp{f}=(\conssp{\fullshad{f}})^\perp$.
Thus, $\xbar=\zbar\plusl\qq$ where $\zbar=\ebar\plusl\yy''$,
which is in $\conv(\unimin{f})$ by
\Crefequiv{thm:conv-univ-equiv}{thm:conv-univ-equiv:d}{thm:conv-univ-equiv:a}.

\pfpart{%
  (\ref{pr:hard-core:4:b})
  $\Rightarrow$
  (\ref{pr:hard-core:4:c})
}:
Suppose $\xbar=\zbar\plusl\qq$ for some $\zbar\in\conv(\unimin{f})$.
Then for each $i\in\indset$,
$\xbar\cdot\uu_i=\zbar\cdot\uu_i\plusl\qq\cdot\uu_i$.
That these values take the form given
in~(\ref{pr:hard-core:4:c})
therefore follows directly from
\Cref{thm:hard-core:3}(\ref{thm:hard-core:3:b:conv}).

\pfpart{%
  (\ref{pr:hard-core:4:c})
  $\Rightarrow$
  (\ref{pr:hard-core:4:a})
}:
Suppose $\xbar\in\extspace$ has the property stated
in~(\ref{pr:hard-core:4:c}).
Let $\xbar'\in\extspace$ be any minimizer of $\fext$
(which exists by \Cref{pr:fext-min-exists}).
Since $\xbar'$ satisfies~(\ref{pr:hard-core:4:a}),
by the foregoing implications, it must satisfy
(\ref{pr:hard-core:4:c}) as well.
Thus, $\xbar\cdot\uu_i=\xbar'\cdot\uu_i$ for all $i\in\indset$,
implying $\fext(\xbar)=\fext(\xbar')$ by
\Cref{pr:hard-core:1}(\ref{pr:hard-core:1:a-ext}).
Therefore, $\xbar$ also minimizes $\fext$.%
\indexg{empirical risk functions!minimizers of|)}%
\qedhere
\end{proof-parts}
\end{proof}

\indexg{empirical risk functions!minimizing sequences|(}%
Finally, we mention some implications for minimizing sequences.
Let $\seq{\xx_t}$ be a sequence in $\Rn$ that minimizes some function
$f$ satisfying the assumptions of
\Cref{thm:waspr:hard-core:2}.
Then it can be argued using
\Cref{thm:waspr:hard-core:2}
that
$\xx_t\cdot\uu_i\rightarrow-\infty$ for all $i\not\in\hardcore{f}$.
Also, let $\qq_t$ be the projection of $\xx_t$ onto the linear
subspace $\rescperp{f}=\perpresf$.
Then by
\Cref{pr:proj-mins-fullshad},
$\fullshad{f}(\xx_t)=\fullshad{f}(\qq_t) \rightarrow \min \fullshad{f}$,
and furthermore, the $\qq_t\negKern$'s are all
in a compact region of $\rescperp{f}$.
Thus,
if $i\in\hardcore{f}$, then $\uu_i\in\rescperp{f}$,
by $\hardcore{f}$'s definition,
implying $\xx_t\cdot\uu_i=\qq_t\cdot\uu_i$ for all $t$
(since $\xx_t-\qq_t$ is orthogonal to $\rescperp{f}$).
Therefore, $\xx_t\cdot\uu_i$ remains always in some bounded interval of
$\R$.
If, in addition, each $\ell_i$ is strictly convex, then
$\fullshad{f}$ has a unique minimum $\qq$ in $\rescperp{f}$,
by
\Cref{thm:waspr:hard-core:4};
in this case, it can be further argued that $\qq_t\rightarrow\qq$,
and $\xx_t\cdot\uu_i\rightarrow\qq\cdot\uu_i$ for $i\in\hardcore{f}$.%
\indexg{empirical risk functions!minimizing sequences|)}%
\indexg{empirical risk functions|)}%

\chapter{Continuity}
\label{sec:continuity}

As we saw in earlier examples, the extension $\fext$ of a convex function
$f:\Rn\rightarrow\Rext$
may not be continuous at a particular
point, even if the function $f$ being extended is continuous everywhere and
well-behaved in other ways.
In this chapter, we characterize precisely the set of points where
$\fext$ is continuous in terms of properties of the original function
$f$.
We also give necessary and sufficient conditions for
$\fext$ to be continuous everywhere.

\indexg{continuity of extensions!sources of discontinuity|(}%
We begin with some examples.
In these, as well as in the development to follow, we focus mainly on
the extensible continuity of $f$,
as was defined in \Cref{dfn:extens-cont},
which, because $f$ is convex, is equivalent to the ordinary continuity
of $\fext$ (by \Cref{thm:ext-cont-f}\ref{thm:ext-cont-f:b}).

\begin{example}[Continuity of the product of hyperbolas]
\label{ex:recip-fcn-eg:cont}
\indexg{Product of hyperbolas!continuity of|(}%
In \Cref{ex:recip-fcn-eg}, we considered the closed, proper, convex and
everywhere continuous function
$f:\R^2\to\eR$ defined, for $\xx\in\R^2$, by
\[
f(\xx) =
\begin{cases}
  \dfrac{1}{x_1 x_2}
  & \text{if $x_1>0$ and $x_2>0$,}
\\[2ex]
  +\infty
  & \text{otherwise.}
\end{cases}
\]
We argued
that if $\seq{\xx_t}$ is a sequence in $\R^2$
converging to some $\xbar=\limray{\ee_1}\plusl\beta\ee_2$ with $\beta\ne 0$,
then $f(\xx_t)\to\ef(\xbar)$, meaning $f$ is extensibly continuous
at such points $\xbar$, and implying $\ef$ is continuous at such
points as well
(\Cref{thm:ext-cont-f}\ref{thm:ext-cont-f:a}).
However, $\fext$ is not continuous at $\xbar=\limray{\ee_1}$.
For instance, if $\xx_t=\trans{[t^2,1/t]}$,
then
$\ef(\xx_t)=f(\xx_t)=1/t\to 0$,
whereas if $\xx_t=\trans{[t,0]}$,
then
$\ef(\xx_t)=f(\xx_t)=+\infty$ for all $t$.
Since in both cases $\xx_t\to\xbar$, this means
$\ef$ is not continuous at $\xbar$
(and also $f$ is not extensibly continuous at $\xbar$).%
\indexg{Product of hyperbolas!continuity of|)}%
\end{example}

\begin{example}[Continuity of the flattening valley]
  \label{ex:x1sq-over-x2:cont}
\indexg{Flattening valley!continuity of|(}%
  In \Cref{ex:x1sq-over-x2}, we studied the function
  \[
    f(\xx)
    =
    \begin{cases}
      x_1^2/x_2
        & \text{if $x_2 >\regAbs{x_1}$,}
      \\
      2\regAbs{x_1}-x_2
        & \text{otherwise,}
    \end{cases}
  \]
  for $\xx\in\R^2$.
We noted that
this function is convex, closed, proper, finite everywhere and
continuous everywhere.
Nevertheless,
$f$'s extension, $\fext$, is not continuous, for instance, at
$\xbar=\limray{\ee_2}\plusl\limray{\ee_1}$.
For example, if $\xx_t=\trans{[t,t^{3}]}$, then
$\fext(\xx_t)=f(\xx_t)=1/t\rightarrow 0$,
whereas
if $\xx_t=\trans{[t^2,t^3]}$, then
$\fext(\xx_t)=f(\xx_t)={t}\rightarrow +\infty$.
Since $\xx_t\rightarrow\xbar$ in both cases,
$\fext$ is not continuous at $\xbar$
(nor is $f$ extensibly continuous at $\xbar$).%
\indexg{Flattening valley!continuity of|)}%
\end{example}

These examples suggest two different ways in
which $\fext$ can be discontinuous.
In \Cref{ex:recip-fcn-eg:cont},
the discontinuity
appeared as a result of
reaching the boundary between where the function $f$ is finite
(namely, all points with $x_1>0$ and $x_2>0$), and where it is
infinite, in other words, the boundary of $\dom{f}$.
On the other hand, 
in \Cref{ex:x1sq-over-x2:cont},
the function $f$ is finite
everywhere so there is no such boundary to its effective domain.
Instead, the discontinuity seemed to arise as a result of the
variety of ways in which it is possible to construct a sequence
reaching the same
astral point at infinity, but on which the function takes very
different values.
We will soon see how our characterization of continuity
captures these two different kinds of discontinuity.%
\indexg{continuity of extensions!sources of discontinuity|)}%

\section{Characterizing where \texorpdfstring{$\fext$}{the extension of f} is continuous}
\label{subsec:charact:cont}

We turn now to
analyzing the continuity of $\ef$ for
a convex function $f:\Rn\rightarrow\Rext$.
Note first that
if $\fext(\xbar)=+\infty$, then,
by $\fext$'s definition (\Cref{def:lsc-ext}),
$f$ must converge to $+\infty$ on every sequence in $\Rn$ converging to
$\xbar$, so $f$ is extensibly continuous at $\xbar$,
and $\fext$ is continuous at $\xbar$ (by
\Cref{thm:ext-cont-f}\ref{thm:ext-cont-f:a}).
Therefore, we focus on understanding continuity at points $\xbar$
where $\fext(\xbar)< +\infty$, that is, astral points in $\dom{\fext}$.
Let
\indexm{cont f}{$\contsetf$}{points where continuous and not $+\infty$}%
$\contsetf$ denote the set of all points $\xbar\in\extspace$
where $\fext$ is continuous and $\fext(\xbar)<+\infty$.

\indexg{continuity of extensions!characterizations|(}%
In this section,
we provide characterization of the points in $\contsetf$.
We do this in two different ways.
First, we will see that $\contsetf$ is equal to the
interior of the effective domain of $\fext$, that is,
$\contsetf=\intdom{\fext}$.
This means that $\fext$ is continuous everywhere except for points
that are in $\dom{\fext}$, but not its interior.
This provides a close analogue to the continuity properties of standard
convex functions on $\Rn$ (\Cref{pr:stand-cvx-cont}).

In addition, we further characterize $\contsetf$ in terms of the
original function $f$.
In particular, we will see that $\contsetf$ consists of all
points $\xbar\in\extspace$ of a specific form
$\xbar=\VV\omm\plusl\qq$
where $\qq\in\Rn$ is in the interior of the effective domain
of $f$, and
all of the columns of $\VV$ are in $\resc{f}$ so that
$\VV\omm\in\represc{f}$.
Thus,
\[
  \contsetf = \paren{\represc{f}\cap\corezn} \plusl \intdom{f}.
\]

Moreover, we will see that the value of $\fext$ at any point
$\xbar$ of this form is equal to $\fV(\qq)$ where
$\fV$ is the shadow of $f$ along $\VV$, as was defined in
\Cref{dfn:shadow-fcns}.
Since $f$ is extensibly continuous at $\xbar$
(by \Cref{thm:ext-cont-f}\ref{thm:ext-cont-f:b}),
this means that
for every sequence $\seq{\xx_t}$ in $\Rn$ with
$\xx_t\rightarrow\xbar$, we must also have
$f(\xx_t)\rightarrow\fV(\qq)$.

We prove all this in a series of theorems.
The first regards convergence of
sequences of a particular form, and
will imply as corollary that $\fext$ must be continuous at every point
of the form just described.
Note in the next
\namecref{thm:recf-seq-cont} that the sequence $\seq{\xx_t}$ need not
have a limit, and that we also make no assumption about the relative
speed at which each sequence $\seq{b_{t,i}}$ converges to $+\infty$.

In what follows,
for vectors $\bb,\,\cc\in\Rk$, we write $\bb\geq\cc$ to mean
$b_i\geq c_i$ for all $i\in\set{1,\dotsc,k}$.

\begin{theorem}  \label{thm:recf-seq-cont}
  Let $f:\Rn\rightarrow\Rext$ be convex and lower semicontinuous.
  Let $\VV=[\vv_1,\ldots,\vv_k]$
  where $\vv_1,\ldots,\vv_k\in\resc{f}$,
  and let
  $\qq=\qhat+\VV\chat$ for some $\qhat\in\intdom{f}$
  and $\chat\in\R^k$.
  Suppose
  $\xx_t = \VV \bb_t + \qq_t$
  for all $t$,
  where $\seq{\bb_t}$ is a sequence in $\R^k$
  such that  $b_{t,i}\rightarrow+\infty$ for $i=1,\ldots,k$,
  and $\seq{\qq_t}$ is a sequence in $\Rn$
  such that $\qq_t\rightarrow\qq$.
  Then $f(\xx_t)\rightarrow \fV(\qq)=\fV(\qhat)$.
\end{theorem}

\begin{proof}
Since $\qhat\in\intdom{f}$, there exists an open set $U\subseteq\Rn$
that includes $\zero$ and
such that $\qhat+U\subseteq\dom{f}$
(\Cref{pr:closure:intersect}\ref{pr:closure:intersect:b}).

By \Cref{pr:gtil-is-PPfP-gen}, $\fV$ is convex.
Also, by its definition (Eq.~\ref{eqn:fV-init-defn}),
\begin{equation} \label{eqn:h-prop:1}
\fV(\VV\bb+\xx)=\fV(\xx)\leq f(\xx)
\end{equation}
for all $\xx\in\Rn$ and $\bb\in\Rk$.
In particular,
\[\fV(\qq)=\fV(\VV\chat+\qhat)=\fV(\qhat)\leq f(\qhat)<+\infty,\]
since $\qhat\in\dom{f}$.

\begin{claimpx}   \label{cl:a4}
  $\fV$ is continuous at $\qq$.
\end{claimpx}

\begin{proofx}
For all $\sS\in U$,
\[
  \fV(\qq+\sS)
  =\fV(\VV\chat + \qhat+\sS)
  =\fV(\qhat +\sS)
 \leq f(\qhat+\sS)<+\infty.
\]
The second equality and first inequality are both
from \eqref{eqn:h-prop:1}.
The last inequality is because $\qhat+U\subseteq\dom{f}$.
Thus, $\qq+U\subseteq\dom{h}$, so
$\qq\in\intdom{h}$.
Since $\fV$ is convex, this implies that $\fV$ is continuous at $\qq$
(by \Cref{pr:stand-cvx-cont}).
\end{proofx}

Next,
we have that
\begin{equation}   \label{eq:thm:recf-seq-cont:2}
  \liminf f(\xx_t)
  =
  \liminf f(\VV \bb_t + \qq_t)
  \geq
  \liminf \fV(\qq_t)
  =
  \fV(\qq).
\end{equation}
The inequality is by \eqref{eqn:h-prop:1},
and the second equality is
by \Cref{cl:a4}, since {$\qq_t\rightarrow\qq$}.
It remains then only to show that
$\limsup f(\xx_t)\leq \fV(\qq)$.

Let $\beta\in\R$ be such that $\beta>\fV(\qq)=\fV(\qhat)$.
Then by $\fV$'s definition (Eq.~\ref{eqn:fV-init-defn}),
there exists $\dd\in\Rk$ for which
$f(\VV\dd+\qhat) < \beta$.
Without loss of generality, we can assume that $\dd\geq\zero$;
otherwise, if this is not the case, we can replace $\dd$ by
$\dd'\in\Rk$ where $d'_i=\max\{0,d_i\}$ for $i=1,\ldots,k$,
so that
\[
   f(\VV \dd' + \qhat)
   =
   f\bigParens{\VV \dd + \qhat + \VV (\dd' - \dd)}
   \leq
   f(\VV \dd + \qhat)
   <
   \beta,
\]
where the first inequality is because $\dd' - \dd \geq \zero$,
implying
$\VV (\dd' - \dd) \in \resc{f}$.
(This is because $\resc{f}$ is a convex cone,
by \Cref{pr:resc-cone-basic-props},
so any conic combination of the columns of $\VV$,
written as $\VV\bb$ for some $\bb\ge\zero$, must be in
$\resc{f}$.)

\begin{claimpx}   \label{cl:thm:recf-seq-cont:1}
  $f$ is continuous at $\VV \dd + \qhat$.
\end{claimpx}

\begin{proofx}
For all $\sS\in U$,
\[ f(\VV \dd + \qhat + \sS) \leq f(\qhat + \sS) < +\infty. \]
The first inequality is because $\dd\geq\zero$, implying
$\VV\dd\in\resc{f}$, and the second is because
$\qhat+U\subseteq \dom{f}$.
Thus,
$\VV\dd+\qhat+U\subseteq\dom{f}$, implying
$\VV\dd+\qhat\in\intdom{f}$, and proving the claim
(by \Cref{pr:stand-cvx-cont}).
\end{proofx}

From \Cref{cl:thm:recf-seq-cont:1}
and since $f(\VV\dd + \qhat) < \beta$,
it follows that there exists $\delta\in\Rstrictpos$ such that
for all $\yy\in\Rn$,
if $\norm{\yy}<\delta$ then
  $f(\VV\dd + \qhat + \yy) < \beta$.

Since $b_{t,i}\rightarrow+\infty$ and $\qq_t\rightarrow\qq$,
we must have that
for all $t$ sufficiently large,
$b_{t,i}\geq d_i - \hatc_i$ for $i=1,\ldots,k$
(that is, $\bb_t \geq \dd - \chat$),
and also
$\norm{\qq_t - \qq}<\delta$.
When these conditions hold, we have
\begin{align}
  f(\xx_t)
  =
  f(\VV \bb_t + \qq_t)
  &=
  f\bigParens{\VV \dd + \qhat + (\qq_t-\qq) + \VV(\bb_t - \dd + \chat)}
  \nonumber
  \\
  &\leq
  f\bigParens{\VV \dd + \qhat + (\qq_t-\qq)}
  <
  \beta.
  \label{eq:thm:recf-seq-cont:3}
\end{align}
The second equality is by algebra and since
$\qq=\qhat+\VV\chat$.
The first inequality is because
$\bb_t\geq\dd-\chat$ implying
$\VV(\bb_t - \dd + \chat)\in\resc{f}$.
The last inequality is
because $\norm{\qq_t-\qq}<\delta$.

Thus, \eqref{eq:thm:recf-seq-cont:3} holds for all $t$ sufficiently
large, so 
$\limsup f(\xx_t)\leq \beta$.
Since this holds for all $\beta>\fV(\qq)$,
it follows, with \eqref{eq:thm:recf-seq-cont:2}, that
\[  \fV(\qq) \leq \liminf f(\xx_t) \leq \limsup f(\xx_t)\leq \fV(\qq), \]
completing the proof.
\end{proof}

As corollary, we obtain the continuity of $\fext$ at points
of the form discussed above:

\begin{corollary}  \label{thm:g:2}
  Let $f:\Rn\rightarrow\Rext$ be convex and lower semicontinuous.
  Let $\VV=[\vv_1,\ldots,\vv_k]$
  where $\vv_1,\ldots,\vv_k\in\resc{f}$,
  and let
  $\qq\in\intdom{f}$.
  Let $\xbar=\VV\omm\plusl\qq$.
  Then $\fext(\xbar)=\fV(\qq)<+\infty$, and
  $\fext$ is continuous at $\xbar$.
  Thus,
  \[
      \paren{\represc{f}\cap\corezn} \plusl \intdom{f}
                   \subseteq \contsetf.
  \]
\end{corollary}

\begin{proof}
Let $\seq{\xx_t}$ be any sequence in $\Rn$ that converges to
$\xbar$.
Then by \Cref{thm:seq-rep-not-lin-ind}, there exist sequences
$\seq{\bb_t}$ in $\Rk$ and $\seq{\qq_t}$ in $\Rn$ such that
$\xx_t = \VV \bb_t + \qq_t$ for all $t$,
with
$b_{t,i}\rightarrow+\infty$ for $i=1,\ldots,k$,
and $\qq_t\rightarrow\qq$.
Therefore, by
\Cref{thm:recf-seq-cont},
$f(\xx_t)\rightarrow\fV(\qq)$.
Since this holds for every such sequence, this implies that
$\fext(\xbar)=\fV(\qq)$ by definition of $\fext$
(Eq.~\ref{eq:e:7}), and furthermore that $f$ is extensibly continuous
at $\xbar$, so also $\fext$ is continuous at
$\xbar$ (by \Cref{thm:ext-cont-f}\ref{thm:ext-cont-f:a}).
In addition,
$\fV(\qq)\leq f(\qq)<+\infty$
by $\fV$'s definition and
since $\qq\in\dom{f}$, so
$\xbar\in\contsetf$.
\end{proof}

Next, we show that if $\xbar$ is in $\dom{\fext}$,
and if $\fext$ is continuous at $\xbar$, then actually
$\xbar$ must be in the interior of $\dom{\fext}$.

\begin{theorem}  \label{thm:cont-in-intdom}
  Let $f:\Rn\rightarrow\Rext$ be convex.
  Suppose $\fext$ is continuous at some point $\xbar\in\extspace$, and
  that $\fext(\xbar)<+\infty$.
  Then $\xbar\in\intdom{\fext}$.
  In other words,
  \[   \contsetf \subseteq \intdom{\fext}. \]
\end{theorem}

\begin{proof}
Suppose, by way of contradiction, that $\xbar\not\in\intdom{\fext}$.
Then
\[
  \xbar\in\eRn\setminus[\intdom{\ef}]=\clbar{\eRn\setminus(\dom{\ef})},
\]
with equality from
\Cref{pr:closure:intersect}(\ref{pr:closure:intersect:comp}).
Therefore, there exists a sequence $\seq{\xbar_t}$ in $\eRn\setminus(\dom{\ef})$
such that $\xbar_t\rightarrow\xbar$. Since $\ef(\xbar_t)=+\infty$ for all $t$,
the continuity of $\ef$ at $\xbar$ implies that $\ef(\xbar)=+\infty$,
which is a contradiction.
\end{proof}

As a final step in our characterization,
we show that every point in $\intdom{\fext}$ must have the
form given in \Cref{thm:g:2}:

\begin{theorem}   \label{thm:intdom-then-resc-form}
  Let $f:\Rn\rightarrow\Rext$ be convex and lower semicontinuous.
  Then
  \[   \intdom{\fext} \subseteq
      \paren{\represc{f}\cap\corezn} \plusl \intdom{f}.
  \]
  That is, if $\xbar\in\intdom{\fext}$, then
  $\xbar=\ebar\plusl\qq$ for some $\ebar\in\represc{f}\cap\corezn$
  and some $\qq\in\intdom{f}$.
\end{theorem}

\begin{proof}
Let $\xbar\in\intdom{\fext}$,
implying there exists a neighborhood $U$ of
$\xbar$ that is included in $\dom{\fext}$.
Also,
we can write $\xbar=\ebar\plusl\qq'$,
where $\ebar=[\vv_1,\ldots,\vv_k]\omm$,
for some
$\vv_1,\ldots,\vv_k,\qq'\in\Rn$.

We prove the \namecref{thm:intdom-then-resc-form}
in two parts:

\begin{claimpx}
  There exists some $\qq\in\intdom{f}$ for which $\xbar=\ebar\plusl\qq$.
\end{claimpx}

\begin{proofx}
For each $t$, let
$  \dd_t = \sum_{i=1}^{k}  t^{k+1-i} \vv_i $,
and
$ \xx_t=\dd_t+\qq'$.
Then $\dd_t\rightarrow\ebar$ and
$\xx_t\rightarrow\xbar$
(by \Cref{thm:i:seq-rep}).

We claim there must exist some $t_0\in\nats$ and some $\epsilon\in\Rstrictpos$ for
which $\ball(\xx_{t_0},\epsilon)\subseteq \dom{f}$
(where $\ball(\cdot,\cdot)$ denotes an open ball in $\Rn$,
as in Eq.~\ref{eqn:open-ball-defn}).
Suppose this claim is false.
Then for each $t$, there must exist a point $\xx'_t$ that is in
$\ball(\xx_t,1/t)$ but not $\dom{f}$.
Then $\xx'_t\sim\xx_t$, so $\xx'_t\to\xbar$
(by \Cref{pr:eq-in-lim-same-lim}).
Since $U$ is a neighborhood of $\xbar$, it follows that,
for all $t$ sufficiently large,
$\xx'_t\in U\subseteq\dom{\fext}$, implying
$f(\xx'_t)=\fext(\xx'_t)<+\infty$
(with equality by \Cref{pr:h:1}\ref{pr:h:1a}
since $f$ is lower semicontinuous).
But this is a contradiction since $\xx'_t\not\in\dom{f}$ for all $t$.

So let $t_0\in\nats$ and $\epsilon\in\Rstrictpos$ be such that
$\ball(\xx_{t_0},\epsilon)\subseteq \dom{f}$, and let
$\qq=\xx_{t_0} = \dd_{t_0}+\qq'$.
Then
$\qq\in\ball(\xx_{t_0},\epsilon)\subseteq \dom{f}$, so
$\qq\in\intdom{f}$.
Furthermore, $\xbar=\ebar\plusl\qq$, by the Push Lemma~\ref{pr:g:1},
since $\qq-\qq'=\dd_{t_0}$ is in the span of $\vv_1,\ldots,\vv_k$.
\end{proofx}

\begin{claimpx}
  $\ebar\in\represc{f}\cap\corezn$.
\end{claimpx}

\begin{proofx}
Let
$\sbar_j = \limrays{\vv_1,\ldots,\vv_j}$
for $j=0,1,\ldots,k$; in particular, $\ebar=\sbar_k$.
We prove by induction on $j=0,1,\ldots,k$ that
$\sbar_j=\limrays{\ww_1,\ldots,\ww_j}$
for some vectors $\ww_1,\ldots,\ww_j\in\resc{f}$.
When $j=k$, the claim is proved.

The base case, when $j=0$, holds vacuously.
For the inductive step, let $j\geq 1$, and assume
$\sbar_{j-1}=\limrays{\ww_1,\ldots,\ww_{j-1}}$
for some $\ww_1,\ldots,\ww_{j-1}\in\resc{f}$.

Let
\[
   \ybar = \sbar_j
   = \sbar_{j-1}\plusl \limray{\vv_j}
   = \limrays{\ww_1,\ldots,\ww_{j-1},\vv_j},
\]
and let
$\zbar=\limrays{\vv_{j+1},\ldots,\vv_k}\plusl\qq'$
so that $\xbar=\ybar\plusl\zbar$.
For each $t$, let
\begin{gather*}
   \yy_t = \BiggParens{\sum_{i=1}^{j-1} t^{k+1-i} \ww_i} + t^{k+1-j} \vv_j,
   \qquad
   \zz_t = \sum_{i=j+1}^{k} t^{k+1-i} \vv_i + \qq',
\\
  \xx_t = \yy_t + \zz_t,
  \qquad
  \ybar_t=\limray{\yy_t},
  \qquad
  \xbar_t=\ybar_t\plusl\zz_t.
\end{gather*}
Then $\xx_t\rightarrow\xbar$ and $\yy_t\rightarrow\ybar$
(by \Cref{thm:i:seq-rep}).

We claim furthermore that $\xbar_t\rightarrow\xbar$.
To see this, let $\uu\in\Rn$;
we aim to show that $\xbar_t\cdot\uu\rightarrow\xbar\cdot\uu$.
First, suppose $\ybar\cdot\uu=+\infty$.
Then $\xbar\cdot\uu=\ybar\cdot\uu\plusl\zbar\cdot\uu=+\infty$.
Since $\yy_t\cdot\uu\rightarrow\ybar\cdot\uu$,
we must have $\yy_t\cdot\uu>0$ for all $t$ sufficiently large, implying,
for all such $t$, that
$\ybar_t\cdot\uu=+\infty$, and so also
$\xbar_t\cdot\uu=\ybar_t\cdot\uu\plusl\zz_t\cdot\uu=+\infty$.
Thus, $\xbar_t\cdot\uu\rightarrow\xbar\cdot\uu$ in this case.

The case $\ybar\cdot\uu=-\infty$ follows symmetrically.

The only remaining case is that $\ybar\cdot\uu=0$
(by
\Cref{pr:icon-equiv}\ref{pr:icon-equiv:c}\ref{pr:icon-equiv:b}).
In this case, we must have $\vv_j\cdot\uu=0$ and
$\ww_i\cdot\uu=0$ for $i=1,\ldots,j-1$, by
\Cref{pr:vtransu-zero}.
These imply, for all $t$, that $\yy_t\cdot\uu=0$, so also
$\ybar_t\cdot\uu=0$.
Thus,
\[
  \xbar_t\cdot\uu
  =
  \ybar_t\cdot\uu \plusl \zz_t\cdot\uu
  =
  \zz_t\cdot\uu
  =
  \yy_t\cdot\uu + \zz_t\cdot\uu
  =
  \xx_t\cdot\uu \rightarrow \xbar\cdot\uu.
\]
We conclude $\xbar_t\rightarrow\xbar$
(by \Cref{thm:i:1}\ref{thm:i:1c}).

Therefore,
for all sufficiently large $t$, $\xbar_t$ must be in the neighborhood
$U\subseteq\dom{\fext}$.
For the rest of the proof,
let $t$ be any such index so that $\fext(\xbar_{t})<+\infty$.
Then
\[
   \fext(\limray{\yy_t}\plusl\zz_t)
   =
   \fext(\ybar_t\plusl\zz_t)
   =
   \fext(\xbar_t)
   <
   +\infty,
\]
implying $\yy_t\in\resc{f}$ (by \Cref{thm:rec-ext-equivs}\ref{thm:rec-ext-equivs:d}\ref{thm:rec-ext-equivs:a} and \Cref{pr:f:1}).

Furthermore, 
\[
  \sbar_j = \ybar = \limrays{\ww_1,\ldots,\ww_{j-1},\yy_t}
\]
by the Push Lemma~\ref{pr:g:1} and
\Cref{prop:pos-upper:descr},
since
$\yy_t$ is a linear combination of
$\ww_1,\dotsc,\ww_{j-1},\vv_j$, with a positive coefficient on $\vv_j$.
Setting $\ww_j=\yy_t\in\resc{f}$ now completes the induction.
\end{proofx}

Together, the two claims prove the
\namecref{thm:intdom-then-resc-form}.
\end{proof}

Combining the foregoing, we obtain a full characterization of all
points where $\fext$ is continuous (and not $+\infty$):

\begin{corollary}  \label{cor:cont-gen-char}
  Let $f:\Rn\rightarrow\Rext$ be convex and lower semicontinuous,
  and let $\xbar\in\extspace$. Then the following are equivalent:
  \begin{letter-compact}
    \item  \label{cor:cont-gen-char:a}
      $\fext(\xbar)<+\infty$ and $\fext$ is continuous at $\xbar$.
    \item  \label{cor:cont-gen-char:b}
      $\xbar\in\intdom{\fext}$.
    \item  \label{cor:cont-gen-char:c}
      $\xbar=\limrays{\vv_1,\ldots,\vv_k}\plusl\qq$ for some
      $\qq\in\intdom{f}$ and some $\vv_1,\ldots,\vv_k\in\resc{f}$.
  \end{letter-compact}
  That is,
  \begin{align}
    \contsetf
    =
    \intdom{\fext}
    &=
    \bigParens{\represc{f}\cap\corezn} \plusl \intdom{f}
    \nonumber
    \\
    &=
    \represc{f} \plusl \intdom{f}.
    \label{eq:cor:cont-gen-char:1}
  \end{align}
\end{corollary}

\begin{proof}
  ~

\begin{proof-parts}

\pfpart{%
  (\ref{cor:cont-gen-char:a})
  $\Rightarrow$
  (\ref{cor:cont-gen-char:b}):
}
This was shown in \Cref{thm:cont-in-intdom}.

\pfpart{%
  (\ref{cor:cont-gen-char:b})
  $\Rightarrow$
  (\ref{cor:cont-gen-char:c}):
}
This was shown in \Cref{thm:intdom-then-resc-form}.

\pfpart{%
  (\ref{cor:cont-gen-char:c})
  $\Rightarrow$
  (\ref{cor:cont-gen-char:a}):
}
This was shown in \Cref{thm:g:2}.

\pfpart{\eqref{eq:cor:cont-gen-char:1}:}
Let $K=\represc{f}$ and $D=\intdom{f}$.
Corresponding to the three implications above, we have
\[
  \contsetf
  \subseteq
  \intdom{\fext}
  \subseteq
  \paren{K\cap\corezn} \plusl D
  \subseteq
  \contsetf,
\]
proving the first two equalities of \eqref{eq:cor:cont-gen-char:1}.

For the last equality,
it suffices to show
$K\plusl D\subseteq\regParens{K\cap\corezn}\plusl D$ (since the reverse
inclusion is immediate).
Let $\xbar=\ybar\plusl\zz$ where
$\ybar\in K$ and $\zz\in D$.
Then $\ybar=\ebar\plusl\qq$ for some
$\ebar\in K\cap\corezn$
and
$\qq\in\resc{f}$
(since $\ybar$ can be represented using only vectors in $\resc{f}$).
Since $\zz\in\intdom{f}$, there exists an open set $U\subseteq\Rn$
including $\zero$ such that $\zz+U\subseteq\dom{f}$.
Since $\qq\in\resc{f}$, for all $\sS\in U$,
$f(\qq+\zz+\sS)\leq f(\zz+\sS)<+\infty$,
so $\qq+\zz+U\subseteq\dom{f}$, implying $\qq+\zz\in\intdom{f}=D$.
Thus,
$\xbar=\ebar\plusl(\qq+\zz)\in \paren{K\cap\corezn} \plusl D$,
completing the proof.
\qedhere
\end{proof-parts}
\end{proof}

Let $\xbar\in\extspace$, and
suppose $\fext(\xbar)<+\infty$.
We can write $\xbar=\ebar\plusl\qq$ for some $\ebar\in\corezn$ and
$\qq\in\Rn$; furthermore, $\ebar\in\aresconef$ by
\Cref{cor:a:4}.
\Cref{cor:cont-gen-char}(\ref{cor:cont-gen-char:a},\ref{cor:cont-gen-char:c})
makes explicit the
conditions under which $\fext$ is or is not continuous at
$\xbar$, namely,
$\fext$ is continuous at $\xbar$
if $\ebar\in\represc{f}$ and if also $\qq$ can be chosen
to be in $\intdom{f}$.
Otherwise, if
$\ebar\not\in\represc{f}$ or if there is no way of choosing $\qq$ so
that $\xbar=\ebar\plusl\qq$ still holds and also $\qq\in\intdom{f}$,
then $\fext$ is discontinuous at $\xbar$.

These two cases of how a discontinuity can arise are also
demonstrated in the two examples from the beginning of this chapter:

\begin{example}
\indexg{Product of hyperbolas!continuity of|(}%
In \Cref{ex:recip-fcn-eg:cont}, we saw that the extension $\ef$
of the product of
hyperbolas function is not continuous at $\xbar=\limray{\ee_1}$.
We mentioned earlier that $\resc{f}=\Rpos^2$ for this function, so
$\limray{\ee_1}\in\represc{f}$, which means the first condition for
continuity is satisfied.
However, we can only write $\xbar=\limray{\ee_1}\plusl\qq$ if
$\qq=\beta\ee_1$ for some $\beta\in\R$.
Since no such point is in the effective domain of $f$, let alone its
interior, the extension $\ef$ is discontinuous at $\xbar$.%
\indexg{Product of hyperbolas!continuity of|)}%
\end{example}

\begin{example}
\indexg{Flattening valley!continuity of|(}%
For the flattening valley function, we showed in
\Cref{ex:x1sq-over-x2:cont} that $\ef$ is not continuous
at the point
$\xbar=\ebar\plusl\zero$ where
$\ebar=\limray{\ee_2}\plusl\limray{\ee_1}$.
In this case,
the function $f$ is finite everywhere, so all points in $\R^2$,
including the origin, are in the
interior of $\dom{f}=\R^2$, thereby satisfying the second condition for
continuity.
However, from \Cref{ex:x1sq-over-x2}, we know that
$\resc{f}=\{ \beta \ee_2 :\: \beta \in \Rpos \}$,
which implies that the only icons in $\represc{f}$ are
$\zero$ and $\limray{\ee_2}$.
In particular, this means
$\ebar\not\in\represc{f}$, yielding the
discontinuity at $\xbar$.%
\indexg{Flattening valley!continuity of|)}%
\end{example}

The next theorem uses sequential sums to generalize \Cref{thm:g:2} by
broadening the set of conditions on a point $\xbar$ under which
continuity of $\fext$ is ensured, while also providing an explicit expression for
$\fext(\xbar)$.
\Cref{thm:g:2} then follows as a special case since,
in the notation that follows,
$\VV\omm\plusl\qq$ is always included in the sequential sum appearing
in \eqref{eq:thm:seqsum-cont:1}.

\begin{theorem}   \label{thm:seqsum-cont}
  Let $f:\Rn\rightarrow\Rext$ be convex and lower semicontinuous.
  Let $\vv_1,\ldots,\vv_k\in\resc{f}$,
  let $\VV=[\vv_1,\ldots,\vv_k]$,
  and let $\qq\in\intdom{f}$.
  Let
  \begin{equation}  \label{eq:thm:seqsum-cont:1}
    \xbar
    \in
    \limray{\vv_1}\seqsum\dotsb\seqsum\limray{\vv_k}\seqsum\qq.
  \end{equation}
  Then $\fext(\xbar)=\fV(\qq)<+\infty$, and
  $\fext$ is continuous at $\xbar$.
\end{theorem}

\begin{proof}
  Since $\limray{\vv_i}\in\represc{f}$ for $i=1,\dotsc,k$, and $\represc{f}$ is a convex astral
  cone (\Cref{pr:repres-in-arescone}),
  \[
    (\limray{\vv_1}\seqsum\dotsb\seqsum\limray{\vv_k})\seqsum\qq
    \subseteq
    \represc{f}\seqsum\qq
    =
    \represc{f}\plusl\qq
    \subseteq
    \represc{f}\plusl\intdom{f},
  \]
  where the first inclusion is by \Cref{thm:ast-cone-is-cvx-if-sum},
  and the equality is by \Cref{cor:seqsum-conseqs}(\ref{cor:seqsum-conseqs:e}).
  Therefore, $\xbar\in\represc{f}\plusl\intdom{f}$, so
  $\ef(\xbar)<+\infty$ and $\ef$ is continuous at $\xbar$
  (by
  \Cref{cor:cont-gen-char}\ref{cor:cont-gen-char:c}\ref{cor:cont-gen-char:a}).
  It remains to show that $\ef(\xbar)=\fV(\qq)$.

By \Cref{thm:mul-char}(\ref{thm:mul-char:a}),
we have
$\xbar\in\mul{\oms}{\vv_1}\seqsum\dotsb\seqsum\mul{\oms}{\vv_k}\seqsum\mul{1}{\qq}$.
Therefore,
by \Cref{thm:seqsum-midrays},
for $i=1,\ldots,k$,
there exist sequences
$\seq{\lambda_{it}}$ in $\Rpos$
and span-bound sequences
$\seq{\ww_{it}}$ in $\Rn$
and
$\seq{\qq_t}$ in $\Rn$
such that
$\ww_{it}\rightarrow\vv_i$,
$\qq_t\rightarrow\qq$,
$\lambda_{it}\rightarrow+\infty$,
and
$\xx_t\rightarrow\xbar$ where
$\xx_t=\sum_{i=1}^k \lambda_{it}\ww_{it} + \qq_t$
for all $t$.
(Note that \Cref{thm:seqsum-midrays} allows us in this case to choose
sequences in which $\qq_t\negKern$'s coefficient is always~$1$.)

Further, for $i=1,\ldots,k$,
since $\seq{\ww_{it}}$ is span-bound and converges to
$\vv_i$, we must have $\ww_{it}=\gamma_{it}\vv_i$ for some
$\gamma_{it}\in\R$ for all $t$; moreover, $\gamma_{it}\rightarrow 1$.
(This assumes $\vv_i\neq\zero$; if $\vv_i=\zero$, then also
$\ww_{it}=\zero$, so we can choose $\gamma_{it}=1$ for all $t$.)
Thus, for all $t$,
$\xx_t=\sum_{i=1}^k \lambda_{it}\gamma_{it}\vv_i + \qq_t$.
Setting $\bb_t=\trans{[b_{t,1},\ldots,b_{t,k}]}$
with $b_{t,i}=\lambda_{it}\gamma_{it}$ for $i=1,\dotsc,k$,
we obtain
sequences that satisfy the conditions of \Cref{thm:recf-seq-cont}
(that is,
$\xx_t=\VV\bb_t + \qq_t$
with $b_{t,i}\rightarrow+\infty$ for $i=1,\ldots,k$,
and $\qq_t\rightarrow\qq$).
Since also $\ef$ is continuous at $\xbar$, that
\namecref{thm:recf-seq-cont} then implies
$\ef(\xbar)=\lim f(\xx_t)=\fV(\qq)$.%
\indexg{continuity of extensions!characterizations|)}%
\end{proof}

\section{Conditions for continuity}
\label{subsec:cond-for-cont}

We next explore general conditions for continuity, especially for
$\fext$ to be continuous everywhere, and especially as
implications of the characterization
given in
\Cref{cor:cont-gen-char}(\ref{cor:cont-gen-char:a},\ref{cor:cont-gen-char:c}).

\indexg{recessive completeness!characterizations|(}%
From \Cref{cor:a:4}, we know that
if $\fext(\xbar)<+\infty$, where $\xbar=\ebar\plusl\qq$, with
$\ebar\in\corezn$ and $\qq\in\Rn$, then $\ebar\in\aresconef$.
Moreover, if $\ebar\not\in\represc{f}$, then $\fext$ cannot be continuous
at~$\xbar$
(by \Cref{cor:cont-gen-char}\ref{cor:cont-gen-char:c}\ref{cor:cont-gen-char:a}).
Thus, for $\fext$ to be continuous everywhere, it is necessary that
$(\aresconef)\cap\corezn\subseteq\represc{f}$.
As we show in the next
\namecref{pr:r-equals-c-equiv-forms},
this latter condition is equivalent to
$\aresconef$ being equal to $\represc{f}$,
that is, to $f$ being recessive complete.

\begin{proposition}  \label{pr:r-equals-c-equiv-forms}
  Let $f:\Rn\rightarrow\Rext$ be convex and lower semicontinuous.
  Then the following are equivalent:
  \begin{letter-compact}
  \item  \label{pr:r-equals-c-equiv-forms:a}
    $(\aresconef) \cap \corezn = \represc{f} \cap \corezn$.
  \item  \label{pr:r-equals-c-equiv-forms:b}
    $(\aresconef) \cap \corezn \subseteq \represc{f}$.
  \item  \label{pr:r-equals-c-equiv-forms:c}
    $\aresconef = \represc{f}$.
    (That is,    $f$ is recessive complete.)
  \end{letter-compact}
\end{proposition}

\begin{proof}
That
(\ref{pr:r-equals-c-equiv-forms:a}) $\Rightarrow$
(\ref{pr:r-equals-c-equiv-forms:b}),
and
(\ref{pr:r-equals-c-equiv-forms:c}) $\Rightarrow$
(\ref{pr:r-equals-c-equiv-forms:a})
are both immediate.

To see
(\ref{pr:r-equals-c-equiv-forms:b}) $\Rightarrow$
(\ref{pr:r-equals-c-equiv-forms:c}),
suppose
$(\aresconef) \cap \corezn \subseteq \represc{f}$.
Then
\[
  \aresconef
  =
  \conv\bigParens{(\aresconef) \cap \corezn}
  \subseteq
  \represc{f}
  \subseteq
  \aresconef.
\]
The equality follows from
\Crefequiv{thm:ast-cvx-cone-equiv}{thm:ast-cvx-cone-equiv:a}{thm:ast-cvx-cone-equiv:c}
since $\aresconef$ is a convex astral cone
(\Cref{cor:res-fbar-closed}).
The first inclusion is by our assumption, and since
$\represc{f}$ is convex
by
\Cref{pr:repres-in-arescone}
(and also
\Cref{pr:conhull-prop}\ref{pr:conhull-prop:aa}).
The second inclusion is 
by \Cref{pr:repres-in-arescone} as well.%
\indexg{recessive completeness!characterizations|)}%
\end{proof}

Expanding on the discussion above, we prove several direct
consequences of the characterization given in
\Cref{cor:cont-gen-char}(\ref{cor:cont-gen-char:a},\ref{cor:cont-gen-char:c}).

\indexg{recessive completeness!continuity and|(}%
\indexg{continuity of extensions!recessive completeness and|(}%
\indexg{continuity of extensions!characterizations|(}%
First, if $f$ is recessive complete, then the next
\namecref{thm:res-comp-cont}
shows that $\fext$ is continuous
exactly at those points where either $\fext$ is $+\infty$, or
which can be expressed with finite part in $\intdom{f}$.
Thus, for establishing continuity at a given point in this case,
we can largely disregard the point's iconic part.

\begin{corollary}  \label{thm:res-comp-cont}
  Let $f:\Rn\rightarrow\Rext$ be convex, lower semicontinuous, and
  recessive complete.
  Let $\xbar\in\extspace$.
  Then $\fext$ is continuous at $\xbar$ if and only if
  either $\fext(\xbar)=+\infty$
  or $\xbar\in\corezn\plusl\intdom{f}$.
\end{corollary}

\begin{proof}
Suppose first that $\fext$ is continuous at $\xbar$.
If $\fext(\xbar)<+\infty$, then by
\Cref{cor:cont-gen-char}(\ref{cor:cont-gen-char:a},\ref{cor:cont-gen-char:c}),
$\xbar\in(\represc{f}\cap\corezn) \plusl \intdom{f}\subseteq\corezn\plusl\intdom{f}$.

For the converse,
if $\fext(\xbar)=+\infty$, then because $\fext$ is lower semicontinuous
(\Cref{prop:ext:F}\ref{prop:ext:F:a}),
it must converge to $+\infty$ on every sequence in $\extspace$
converging to $\xbar$;
thus, $\fext$ is continuous at $\xbar$.

Otherwise, suppose $\fext(\xbar)<+\infty$ and
$\xbar=\ebar\plusl\qq$ for some $\ebar\in\corezn$ and
$\qq\in\intdom{f}$.
Then $\fext(\ebar\plusl\qq)<+\infty$, so
$\ebar\in\aresconef=\represc{f}$
by \Cref{cor:a:4}.
Thus,
$\xbar\in(\represc{f}\cap\corezn) \plusl \intdom{f}$,
implying $\fext$ is continuous at $\xbar$ by
\indexg{continuity of extensions!characterizations|)}%
\Cref{cor:cont-gen-char}(\ref{cor:cont-gen-char:c},\ref{cor:cont-gen-char:a}).
\end{proof}

Along similar lines,
if $\ebar$ is an icon in
$\aresconef$
but not in $\represc{f}$, then no point $\zbar$
in $\ebar\plusl\extspace$, the closure of $\ebar$'s galaxy,
can be in $\contsetf$;
in other words, it is not possible that both $\fext(\zbar)<+\infty$ and
that $\fext$ is continuous at $\zbar$:

\begin{theorem}  \label{thm:d5}
  Let $f:\Rn\rightarrow\Rext$ be convex and lower semicontinuous.
  Suppose $\ebar\in \bigParens{(\aresconef)\cap\corezn} \setminus\represc{f}$.
  Then for all $\xbar\in\extspace$, either
  $\fext(\ebar\plusl\xbar)=\fext(\xbar)=+\infty$,
  or
  $\fext$ is not continuous at $\ebar\plusl\xbar$.
\end{theorem}

\begin{proof}
Let $\xbar\in\extspace$.
If $\fext(\ebar\plusl\xbar)=+\infty$ then $\fext(\xbar)=+\infty$
by \Cref{pr:arescone-def-ez-cons}(\ref{pr:arescone-def-ez-cons:b})
since $\ebar\in\aresconef$.

So assume $\fext(\ebar\plusl\xbar)<+\infty$, and suppose, by way of
contradiction, that $\fext$ is continuous at $\ebar\plusl\xbar$.
Then by
\Cref{cor:cont-gen-char}(\ref{cor:cont-gen-char:a},\ref{cor:cont-gen-char:c}),
there exists
$\dbar\in\represc{f}\cap\corezn$ and $\qq\in\intdom{f}$ such that
$\ebar\plusl\xbar=\dbar\plusl\qq$.
This implies
$\dbar=\ebar\plusl\xbar\plusl (-\qq)$,
and so that
\[\ebar\in\lb{\zero}{\dbar}\subseteq\represc{f},\]
where
the first inclusion is by
\Cref{cor:d-in-lb-0-dplusx},
and the second is because both $\zero$ and $\dbar$ are in $\represc{f}$, which is convex
by \Cref{pr:repres-in-arescone}.
However, this contradicts our assumption that
$\ebar\not\in\represc{f}$.
\end{proof}

We previously remarked that if $\fext$ is continuous everywhere then
$f$ is recessive complete.
\indexg{continuity of extensions!all minimizers@at all minimizers|(}%
\indexg{minimizers of extensions!continuity at all|(}%
Actually, as shown next,
we can make a somewhat stronger statement, namely, that if
$\fext$ is continuous \emph{at all of its minimizers}, then
$f$ is recessive complete.
Clearly, this implies the former assertion.

Combined with
\Cref{thm:unimin-can-be-rankone},
this shows furthermore that if $\fext$ is continuous at all of its
minimizers, then it also must have
canonical minimizers with astral rank at most one.

\begin{theorem}  \label{thm:cont-at-min-implies-ares}
  Let $f:\Rn\rightarrow\Rext$ be convex and lower semicontinuous.
  If $\fext$ is continuous at all of its minimizers,
  then
  $f$ is recessive complete.
\end{theorem}

\begin{proof}
If $f\equiv +\infty$ then
$\aresconef=\represc{f}=\extspace$, proving the claim in this case.
We therefore assume henceforth that $f\not\equiv +\infty$.

We prove the theorem in the contrapositive.
Suppose
$f$ is not recessive complete.
Then by
\Cref{pr:r-equals-c-equiv-forms}(\ref{pr:r-equals-c-equiv-forms:c},\ref{pr:r-equals-c-equiv-forms:b}),
there exists a point
$\ebar \in ((\aresconef) \cap \corezn) \setminus \represc{f}$.
Let $\ybar\in\extspace$ be any point that minimizes
$\fext$, implying $\fext(\ybar)<+\infty$.
Furthermore,
$\fext(\ebar\plusl\ybar)\leq\fext(\ybar)$ by
\Cref{pr:arescone-def-ez-cons}(\ref{pr:arescone-def-ez-cons:b}) since $\ebar\in\aresconef$, so
$\ebar\plusl\ybar$ must also minimize $\fext$.
It now follows immediately from
\Cref{thm:d5}
that $\fext$ is not continuous at $\ebar\plusl\ybar$, which is one of
its minimizers.
\end{proof}

So
recessive completeness
is a necessary condition for $\fext$ to
be continuous everywhere, or even for it to be continuous at all its minimizers.
\indexg{continuity of extensions!when finite everywhere|(}%
When $f$ is convex and finite everywhere,
these conditions all turn out to be
equivalent.
In other words, in this case, $\fext$ is continuous everywhere if and
only if
$f$ is recessive complete.
Furthermore, and quite remarkably, if $\fext$ is continuous at all its
minimizers, then it must actually be continuous everywhere.
Equivalently, if $\fext$ is discontinuous anywhere, then it must be
discontinuous at one or more of its minimizers
(as was the case for the flattening valley function in
\Cref{ex:x1sq-over-x2:cont}, which is finite everywhere).

We prove these equivalences in the next
\namecref{thm:cont-conds-finiteev}.
We also
use the dual characterization of recessive completeness
from \Cref{thm:dual-cond-char} to expand the list of
conditions that are necessary and sufficient for $\ef$ to be
continuous everywhere.
Since we here assume that $f$ is finite everywhere, $f$
must be reduction-closed (\Cref{pr:j:1}\ref{pr:j:1d}),
allowing us, by \Cref{cor:ent-clos-is-slopes-cone},
to replace $\slopes{f}$ by $\cone(\dom{\fstar})$
in the conditions
from \Cref{thm:dual-cond-char}.

Parts~(\ref{thm:cont-conds-finiteev:d})
and~(\ref{thm:cont-conds-finiteev:e})
provide particularly simple
geometric criteria, based only on notions from standard convex
analysis,
for determining if the astral extension $\fext$ is continuous everywhere.

\begin{theorem}  \label{thm:cont-conds-finiteev}
\indexg{polyhedral sets, astral!continuity and|(}%
\indexg{recession cone (standard)!continuity when polyhedral|(}%
\indexg{recession cone, astral!continuity when polyhedral|(}%
\indexg{continuity of extensions!polyhedral recession cones and|(}%
  Let $f:\Rn\rightarrow \R$ be convex.
  Then the following are equivalent:
  \begin{letter-compact}
  \item  \label{thm:cont-conds-finiteev:a}
    $\fext$ is continuous everywhere.
  \item  \label{thm:cont-conds-finiteev:b}
    $\fext$ is continuous at all its minimizers.
  \item  \label{thm:cont-conds-finiteev:c}
    $f$ is recessive complete.
  \item  \label{thm:cont-conds-finiteev:d}
    $\resc{f}$ is polyhedral, and also $\cone(\dom{\fstar})$ is closed
    in $\Rn$.
  \item  \label{thm:cont-conds-finiteev:e}
    $\cone(\dom{\fstar})$ is polyhedral.
  \item  \label{thm:cont-conds-finiteev:f}
    $\aresconef$ is astral polyhedral.
  \end{letter-compact}
\end{theorem}

\begin{proof}
Since $f$ is convex and finite everywhere, it is continuous
everywhere as well
(\Cref{pr:stand-cvx-cont}), and so also lower semicontinuous.

\begin{proof-parts}
\pfpart{%
  (\ref{thm:cont-conds-finiteev:a})
  $\Rightarrow$
  (\ref{thm:cont-conds-finiteev:b}):
}
Immediate.

\pfpart{%
  (\ref{thm:cont-conds-finiteev:b})
  $\Rightarrow$
  (\ref{thm:cont-conds-finiteev:c}):
}
Immediate from \Cref{thm:cont-at-min-implies-ares}.

\pfpart{%
  (\ref{thm:cont-conds-finiteev:c})
  $\Rightarrow$
  (\ref{thm:cont-conds-finiteev:a}):
}
Suppose $f$ is recessive complete.
Since $\dom{f}=\Rn$, its interior is also all of $\Rn$.
By \Cref{thm:res-comp-cont}, $\fext$ is continuous at every point in
$\corezn\plusl\intdom{f}=\extspace$,
and thus is continuous everywhere.

\pfpart{%
  (\ref{thm:cont-conds-finiteev:c}) $\Leftrightarrow$
  (\ref{thm:cont-conds-finiteev:d}) $\Leftrightarrow$
  (\ref{thm:cont-conds-finiteev:e}) $\Leftrightarrow$
  (\ref{thm:cont-conds-finiteev:f}):
}
Since $f$ is finite everywhere, it is also reduction-closed
(\Cref{pr:j:1}\ref{pr:j:1d}),
implying
$\slopes{f}=\cone(\dom{\fstar})$ (by
\Cref{thm:entclosed-implies-slopescone}).
These equivalences then follow by \Cref{thm:dual-cond-char}.%
\indexg{recessive completeness!continuity and|)}%
\indexg{continuity of extensions!recessive completeness and|)}%
\indexg{polyhedral sets, astral!continuity and|)}%
\indexg{recession cone (standard)!continuity when polyhedral|)}%
\indexg{recession cone, astral!continuity when polyhedral|)}%
\indexg{continuity of extensions!polyhedral recession cones and|)}%
\qedhere
\end{proof-parts}
\end{proof}

If $f$ is convex but not finite everywhere,
then it is possible that $\fext$ is
continuous at all of its minimizers, but discontinuous somewhere else.
This is possible even if the function $f$ itself is continuous
everywhere, as the next example shows:

\begin{example}
\indexg{Product of hyperbolas!variant continuous at all minimizers|(}%
Let $f:\R^2\rightarrow\Rext$ be
the following variation on the product of hyperbolas
function from \Cref{ex:recip-fcn-eg}:
\begin{equation}  \label{eqn:recip-fcn-eg-mod}
  f(\xx) = f(x_1, x_2)
  = \begin{cases}
      \dfrac{1}{x_1 x_2}
      + \me^{- x_1}
      + \me^{- x_2}
      & \text{if $x_1>0$ and $x_2>0$,}
    \\[2ex]
      +\infty
      & \text{otherwise,}
  \end{cases}
\end{equation}
for $\xx\in\R^2$.
It can be checked that
$f$ converges to~$0$, and is thereby minimized, on just
those sequences $\seq{\xx_t}$ in $\R^2$ for which
$\xx_t\cdot\ee_1\rightarrow+\infty$
and
$\xx_t\cdot\ee_2\rightarrow+\infty$.
Thus, $\fext$ is minimized, with $\fext(\xbar)=0$,
exactly at those points
$\xbar\in\extspac{2}$ for which
$\xbar\cdot\ee_1=\xbar\cdot\ee_2=+\infty$.
Moreover, $\fext$ is continuous at all such points.
On the other hand, $\fext$ is not continuous at $\limray{\ee_1}$ by
the same reasoning as in \Cref{ex:recip-fcn-eg:cont},
but this point is not a
minimizer since $\fext(\limray{\ee_1})=1$.%
\indexg{Product of hyperbolas!variant continuous at all minimizers|)}%
\indexg{continuity of extensions!when finite everywhere|)}%
\indexg{continuity of extensions!all minimizers@at all minimizers|)}%
\indexg{minimizers of extensions!continuity at all|)}%
\end{example}

Finally, we relate the equivalences given in
\Cref{thm:cont-conds-finiteev}
to functions of the form studied in
\Cref{sec:emp-loss-min}:

\begin{example}  \label{ex:erm-cont}
\indexg{empirical risk functions|(}%
Suppose $f$ is a function of the form given in
\eqref{eqn:loss-sum-form}
(under the same assumptions as in \Cref{pr:hard-core:1} and
throughout \Cref{sec:emp-loss-min}).
Every such function satisfies all of the equivalent
conditions given in \Cref{thm:cont-conds-finiteev}, as we now confirm.
First, $\fext$ is continuous everywhere, by
\Cref{pr:hard-core:1}(\ref{pr:hard-core:1:a-cnt}),
and so also at all its minimizers,
confirming parts~(\ref{thm:cont-conds-finiteev:a})
and~(\ref{thm:cont-conds-finiteev:b})
of \Cref{thm:cont-conds-finiteev}.
That $f$ is recessive complete was shown in
\Cref{pr:hard-core:1}(\ref{pr:hard-core:1:b-resc1}),
confirming part~(\ref{thm:cont-conds-finiteev:c}).
It follows from 
\Cref{pr:hard-core:1}(\ref{pr:hard-core:1:b-cones})
that $\resc{f}$ is polyhedral, and that $\aresconef$ is astral
polyhedral,
confirming part~(\ref{thm:cont-conds-finiteev:f}).
Finally, 
\Cref{pr:hard-core:1}(\ref{pr:hard-core:1:conedomfstar})
shows that $\cone(\dom\fstar)$ is finitely generated,
and therefore polyhedral and closed in $\Rn$
(by \Cref{roc:thm19.1}),
confirming parts~(\ref{thm:cont-conds-finiteev:e})
and~(\ref{thm:cont-conds-finiteev:d}) (since also $\resc{f}$ is polyhedral).%
\indexg{empirical risk functions|)}%
\end{example}

\part{Differential Theory}
\label{part:differential-theory}

\chapter{Defining astral subgradients}
\label{sec:gradients}

Gradients and subgradients are centrally important to convex analysis
because of the power they afford in tackling and understanding
optimization problems.
In this chapter, we introduce a theory of subdifferentials for functions
defined over astral space.
We study two different generalizations of the standard
subgradient.
The first allows us to assign meaningful subgradients for points at
infinity.
The other allows us to assign subgradients  for finite points where the
standard subgradient would otherwise be undefined, for instance, due
to a tangent slope that is infinite.
In such cases, the subgradients we assign are themselves astral
points, and thus potentially infinite.
Later,
we will see that these two forms of astral subgradients are dually
related, and that they are, to a degree,
inverses of one another via conjugacy.

The notions of subgradient that we are about to present differ
from the
\indexg{horizon subgradients}%
``horizon subgradients''
of \idxroc\citet{rockafellar1985extensions}
and \idxroc\idxwets\citet[Chapter 8]{rock_wets}.
They also differ
from those introduced by
\idxsinger\citet[Chapter 10]{singer_book},
although our development of conjugacy largely matches
his abstract framework.

\section{Astral primal subgradients}
\label{sec:gradients-def}

As briefly discussed in \Cref{sec:prelim:subgrads},
a vector $\uu\in\Rn$ is said to be a {subgradient} of a
function $f:\Rn\rightarrow\Rext$ at $\xx\in\Rn$
if
\begin{equation}  \label{eqn:standard-subgrad-ineq}
   f(\xx') \geq f(\xx) + (\xx'-\xx)\cdot\uu
\end{equation}
for all $\xx'\in\Rn$, so that the affine function (in $\xx'$) on the
right-hand side of this inequality is supporting $f$ at $\xx$,
by which we mean that it lower bounds $f$ everywhere
and matches its value at $\xx$.
The {subdifferential} of $f$ at $\xx$, denoted
$\partial f(\xx)$, is the set of all subgradients of $f$ at $\xx$.
(Although these definitions are intended for convex functions, we also
apply them in what follows to general functions
$f$ that are not necessarily convex.)
If $f$ is convex, then its gradient, $\gradf(\xx)$,
if it exists, is always also a subgradient.
Subgradients are central to optimization of convex functions since
$\zero\in\partial f(\xx)$ for a convex function $f$ if and only if $f$ is
minimized at $\xx$.

\begin{figure}
  \centering
  \includegraphics{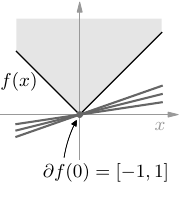}\hfill%
  \includegraphics{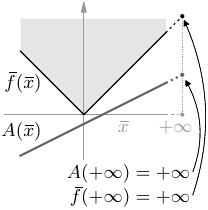}\hfill%
  \includegraphics{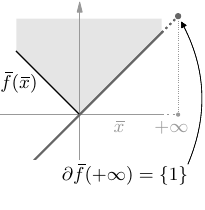}
  \mycaption{Generalizing subgradients of the absolute value function}{%
\indexf{absolute value function, subgradients of}%
\indexf{subdifferentials of extensions!examples}%
\indexf{subdifferentials, astral (primal)!examples}%
  \emph{Left:}
  Several standard subgradients of $f(x)=\abs{x}$ at $x=0$,
  each corresponding to an affine lower bound meeting the graph of $f$ at $x=0$.
  \emph{Center:}
  The shown astral affine function $A$ lower bounds the extension $\ef$ everywhere
  and matches its value at $\barx=+\infty$, but
  the gap
  between the graphs of $\ef$ and $A$ is infinitely large at $\barx=+\infty$.
  \emph{Right:}
  The astral affine function $A$ corresponding to the sole astral subgradient (equal to $1$)
  at $\barx=+\infty$.
  The affine function $A$ lower bounds $\ef$ and also satisfies
  $\ef(x_t)-A(x_t)\to 0$ for the sequence $x_t=t$, which converges to $\barx=+\infty$.
  }%
  \label{fig:1d-abs}%
\end{figure}

\begin{example}[Subgradients of absolute value]
\label{ex:standard-abs-val-subgrad}
\indexg{absolute value function, subgradients of|(}%
Consider the absolute value function,
$f(x)=|x|$ for $x\in\R$.
Although $f$ is not differentiable at $0$, it does have a meaningful
subdifferential at that point, namely, $\partial f(0) = [-1,1]$.
This is because, for $u\in [-1,1]$, $f(x)\geq x u$
for all $x\in\R$, with equality at $x=0$
(see \Cref{fig:1d-abs}, left).
By similar reasoning, for $x\in\R$,
\[
  \partial f(x)
  =
  \begin{cases}
    \{-1\}  & \text{if $x<0$,} \\
    [-1,1]  & \text{if $x=0$,} \\
    \{1\}   & \text{if $x>0$.}
  \end{cases}
\]
The standard subdifferential is not defined at $+\infty$, but if it
were, we might reasonably expect that $1$ should be included as a
subgradient at this point since $f(x)\geq x$ for all $x\in\R$, with
equality holding asymptotically in the limit as $x\rightarrow+\infty$.
Our definition of astral subgradient introduced
below captures this intutition.%
\indexg{absolute value function, subgradients of|)}%
\end{example}

We introduce two generalizations of the standard subdifferential
to astral space, analogous to
two kinds of astral conjugates introduced in \Cref{sec:conjugacy}.
\indexg{subdifferentials, astral (primal)|(}%
The first kind of astral subdifferential maps astral
points in $\extspace$ to subgradients which are real vectors in~$\Rn$;
the second (dual) kind of astral subdifferential maps in the reverse
direction from~$\Rn$ to subsets of $\extspace$.
As with conjugates, this asymmetry
is a consequence of the asymmetry of the coupling function
$\xbar\cdot\uu$, defined over $\xbar\in\extspace$ but $\uu\in\Rn$.

In this section, we define the first generalization, that is,
we define finite subgradients at astral points in
$\extspace$.
As a starting point, we use the definition of a standard
subgradient in \eqref{eqn:standard-subgrad-ineq}.
Note that for standard convex functions $f$,
only $f(\xx')$ and $f(\xx)$ can be infinite
in \eqref{eqn:standard-subgrad-ineq}, so there is no possibility
of adding $-\infty$ and $+\infty$ in that expression.
However, when extending to astral space,
other quantities,
particularly those involving inner products, may become infinite.
Furthermore, there is no operation for directly adding or subtracting
astral points analogous to the difference of vectors, $\xx'-\xx$, that
appears in
\eqref{eqn:standard-subgrad-ineq}.
As a result,
\eqref{eqn:standard-subgrad-ineq} cannot be simply generalized
by replacing each variable and function by its
astral counterpart.
Instead, our approach builds on the correspondence of
subgradients to affine functions that support the graph of $f$.

For the purposes of this discussion, suppose that
$f:\Rn\rightarrow\Rext$ is proper. Then $f$ can have
a subgradient at a point $\xx\in\Rn$ only if $f(\xx)\in\R$,
in which case the right-hand
side of \eqref{eqn:standard-subgrad-ineq}
is equal to an affine function $a(\xx')=\xx'\inprod\uu-\beta$ for some $\beta\in\R$.
The condition given in that inequality requires that
$f(\xx')\ge a(\xx')$ for all~$\xx'$, with equality at~$\xx$.
Thus, for $f$ proper,
$\uu\in\partial f(\xx)$ if and only if there exists $\beta\in\R$ such
that
\begin{letter-compact}
\item     \label{it:std-sub-mod-dfn:a}
  $f(\xx')\ge\xx'\cdot\uu - \beta$
  for all $\xx'\in\Rn$; and
\item     \label{it:std-sub-mod-dfn:b}
  $f(\xx)=\xx\cdot\uu-\beta$.
\end{letter-compact}
Restated in these terms, we can more readily extend subgradients
to astral space.

Let $F:\extspace\rightarrow\Rext$.
We aim to define what it means for some vector $\uu\in\Rn$ to be an astral
subgradient of $F$ at some point $\xbar\in\extspace$.
We begin by
requiring that there exist some $\beta\in\R$ such that
$F(\xbar')\ge\xbar'\cdot\uu - \beta$
for all $\xbar'\in\eRn$
exactly as in condition~(\ref{it:std-sub-mod-dfn:a})
above for standard subgradients.
That is, we require $F\geq A$, where $A:\extspace\rightarrow\Rext$ is
the astral affine function $A(\xbar')=\xbar'\cdot\uu-\beta$
for $\xbar'\in\extspace$.

To generalize condition~(\ref{it:std-sub-mod-dfn:b}),
we could simply require that $F(\xbar)=\xbar\inprod\uu-\beta$,
that is, that $F(\xbar)=A(\xbar)$.
However, when $F(\xbar)$ is infinite, this condition 
does not fully capture the geometric intuition of the graph of $F$
meeting the graph of the astral affine function $A$
at~$\xbar$.
For instance, suppose $F:\Rext\rightarrow\Rext$ is the (extended)
absolute value function $F(\barx')=\abs{\barx'}$ for $\barx'\in\Rext$
(where $\abs{-\infty}=\abs{+\infty}=+\infty$),
and that $A:\Rext\rightarrow\Rext$ is the affine function
$A(\barx')=\barx'/2-1$.
Also, let $\barx=+\infty$.
Then $F\geq A$, and 
$A$ satisfies
$F(\barx)=+\infty=A(\barx)$, naively suggesting, in this view,
that $1/2$ should be a subgradient of $F$ at $\barx$.
However,
$F(x_t)-A(x_t)\to+\infty$ for any sequence $\seq{x_t}$ in $\R$
converging to $\barx$. In this sense,
the gap between the graphs of $F$ and $A$ is infinitely large at
$\barx$, even though the values of $F(\barx)$ and $A(\barx)$ are equal
(see \Cref{fig:1d-abs}, center).

We therefore impose a stronger condition in formulating a general
definition.
Specifically, we require that there exist
a sequence $\seq{\xbar_t}$ in $\eRn$ converging to $\xbar$
on which both $A$ and $F$ are finite, and such that
$F(\xbar_t)\to F(\xbar)$
and
$F(\xbar_t)-A(\xbar_t)\to 0$.
We thus arrive at the following definition:

\begin{definition}   \label{def:ast-subgrad}
\indexg{subdifferentials, astral (primal)!defined|(}%
Let $F:\extspace\rightarrow\Rext$, and
let $\uu\in\Rn$ and $\xbar\in\extspace$.
We say that $\uu$ is an \emph{astral subgradient}
of
$F$ at $\xbar$ if there exists $\beta\in\R$ such that:
\begin{letter-compact}
\item  \label{en:ast-sub-defn-cond-1}
  $F(\xbar')\ge\xbar'\cdot\uu-\beta$
  for all $\xbar'\in\eRn$; and
\item  \label{en:ast-sub-defn-cond-2}
  there exists a sequence $\seq{\xbar_t}$ in $\eRn$ with
  $F(\xbar_t)\in\R$ and $\xbar_t\inprod\uu\in\R$ for all $t$,
  and such that 
  $\xbar_t\to\xbar$, $F(\xbar_t)\to F(\xbar)$,
  and
  $F(\xbar_t)-[\xbar_t\inprod\uu-\beta]\to 0$.
\end{letter-compact}
The
\emph{astral subdifferential}
of $F$ at
$\xbar\in\extspace$, denoted
\indexm{d fx600}{$\asubdifF{\xbar}$}{astral (primal) subdifferential}%
$\asubdifF{\xbar}$,
is the set of all such astral subgradients of $F$ at $\xbar$.
We say that $F$ is
\emph{astral subdifferentiable}
at $\xbar$ if
$\asubdifF{\xbar}$ is not empty.\looseness=-1
\end{definition}

We also sometimes use the term astral \emph{primal} subgradient
(or subdifferential) to emphasize distinction from astral \emph{dual}
subgradients which will be introduced in
\Cref{sec:astral-dual-subgrad}.

In some of the proofs,
it will be convenient to work with the following
alternate pair of conditions,
which are equivalent to those in the definition above:
\begin{letter-compact-prime-simple}
\item  \label{eq:astral:subgrad:alt:1}
  $\beta\ge-F(\xbar')\plusd\xbar'\cdot\uu$
  for all $\xbar'\in\eRn$; and
\item  \label{eq:astral:subgrad:alt:2}
  there exists a sequence $\seq{\xbar_t}$ in $\eRn$ with
  $F(\xbar_t)\in\R$ and $\xbar_t\inprod\uu\in\R$ for all $t$,
  and such that 
  $\xbar_t\to\xbar$, $F(\xbar_t)\to F(\xbar)$,
  and
  $-F(\xbar_t)+\xbar_t\inprod\uu\to\beta$.
\end{letter-compact-prime-simple}
The equivalence of (\ref{en:ast-sub-defn-cond-1}) and
(\ref{eq:astral:subgrad:alt:1}) follows by \Cref{pr:plusd-props}(\ref{pr:plusd-props:e}), and
the equivalence of (\ref{en:ast-sub-defn-cond-2}) and
(\ref{eq:astral:subgrad:alt:2}) follows from properties of limits in $\R$.%
\indexg{subdifferentials, astral (primal)!defined|)}%

Here are some examples of astral subgradients:

\begin{example}  \label{ex:standard-abs-val-subgrad-cont}
\indexg{subdifferentials of extensions!examples|(}%
\indexg{subdifferentials, astral (primal)!examples|(}%
\indexg{absolute value function, subgradients of|(}%
Let $f$ be the absolute value function, as in
Example~\ref{ex:standard-abs-val-subgrad},
whose
extension $\fext$ is the same as $f$ with
$\fext(-\infty)=\fext(+\infty)=+\infty$.
At $\barx=+\infty$,
we can see that $u=1$ is an astral subgradient of $F=\fext$ according to
\Cref{def:ast-subgrad}, as witnessed by
$\beta=0$ and the sequence $x_t=t$
(see \Cref{fig:1d-abs}, right).
There are no other astral subgradients at $+\infty$ since,
for any finite $\beta\in\R$, if $u>1$ then
condition~(\ref{en:ast-sub-defn-cond-1}) of
\Cref{def:ast-subgrad} is violated for sufficiently large $x'$,
and if $u<1$ then $F(x_t)-[x_tu-\beta]\to+\infty$ for
any sequence $\seq{x_t}$ in~$\R$ with $x_t\to+\infty$,
so condition~(\ref{en:ast-sub-defn-cond-2}) cannot be satisfied.
Thus, $\partial \fext(+\infty)=\{1\}$. Similarly,
$\partial \fext(-\infty)=\{-1\}$.
It can also be checked that $\partial\fext(x)=\partial f(x)$
for $x\in\R$.%
\indexg{absolute value function, subgradients of|)}%
\end{example}

\begin{figure}
  \centering
  \includegraphics{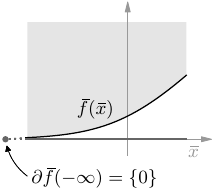}\hfill%
  \includegraphics{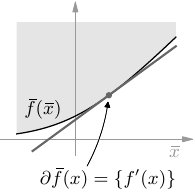}\hfill%
  \includegraphics{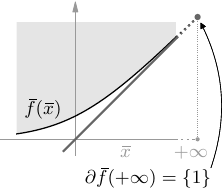}
  \caption[Astral subgradients of the function $f$ from \Cref{ex:subgrad-log1+ex}]{%
\indexf{subdifferentials of extensions!examples}%
\indexf{subdifferentials, astral (primal)!examples}%
\emph{Astral subgradients of the function $f$ from \Cref{ex:subgrad-log1+ex}.}}
  \label{fig:1d-log1exp}%
\end{figure}

\begin{example}  \label{ex:subgrad-log1+ex}
Let $f:\R\rightarrow\Rext$ be defined, for $x\in\R$, by
\begin{equation} \label{eq:ex:ln1plusexp}
  f(x) = \ln(1+e^x),
\end{equation}
and let $\fext$ be the extension of $f$.
The standard subgradients of this function at points $x\in\R$ are
simply given by its derivative $f'$.
Consistent with this derivative tending to~$1$ as
$x\rightarrow+\infty$, it can be checked that $u=1$ is an astral subgradient
of $\fext$ at $\barx=+\infty$ (witnessed by $\beta=0$ and the sequence $x_t=t$).
Indeed,
\begin{equation} \label{eq:ex:ln1plusexp-subgrad}
  \partial \fext(\barx)
  =
  \begin{cases}
           \{0\}   & \mbox{if $\barx = -\infty$,} \\
           \{f'(\barx)\}   & \mbox{if $\barx\in\R$,} \\
           \{1\}   & \mbox{if $\barx = +\infty$,}
  \end{cases}
\end{equation}
as depicted in \Cref{fig:1d-log1exp}.
Note in particular that $0$ is an astral subgradient of $\fext$
at $-\infty$,
consistent with this point being $\fext$'s unique minimizer.
\end{example}

\begin{figure}
  \centering
  \includegraphics{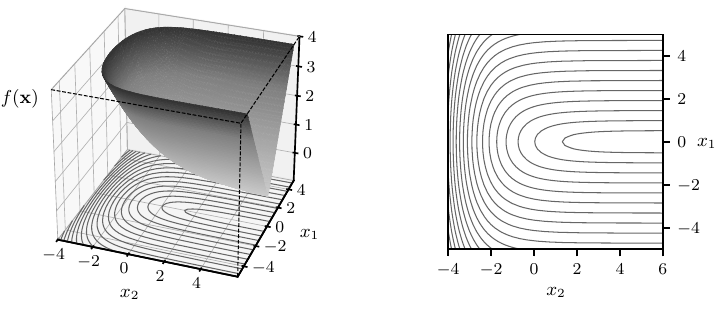}
  \mycaption{Absolute value at infinity}{%
\indexf{Absolute value at infinity}%
      The function $f$ from \Cref{ex:subgrad-sqrt-approaches-abs-val}.
    The surface plot on the left illustrates that as $x_2\to+\infty$,
    cross-sections of $f$'s graph approach the absolute value function in the variable $x_1$.
  }%
  \label{fig:abs-at-inf}%
\end{figure}

\begin{example}[Absolute value at infinity]
   \label{ex:subgrad-sqrt-approaches-abs-val}
\indexg{Absolute value at infinity|(}%
Let $f:\R^2\rightarrow\Rext$ be defined, for $\xx\in\R^2$, by
\begin{equation}   \label{eq:ex:subgrad-sqrt-approaches-abs-val:2}
   f(\xx)
   =
   f(x_1,x_2)
   =
   \sqrt{x_1^2 + e^{-x_2}}
\end{equation}
(see \Cref{fig:abs-at-inf}).
Then $f$ is convex and differentiable
\indexg{Absolute value at infinity|)}%
everywhere.
Its extension is continuous
everywhere and, for $\xbar\in\extspac{2}$, is given by
\[
  \ef(\xbar)=\sqrt{\strut(\xbar\inprod\ee_1)^2 + \expex(-\xbar\inprod\ee_2)}\,.
\]
Therefore, its value at the point $\xbar_\alpha= \limray{\ee_2} \plusl \alpha \ee_1$,
for $\alpha\in\R$, is
$\fext(\xbar_\alpha) = |\alpha|$.
\indexg{Absolute value at infinity!astral subgradients of extension|(}%
Intuitively, we therefore might expect $\fext$'s astral subgradients at
$\xbar_\alpha$ to be akin to the standard subgradients of the absolute value
function, $\alpha\mapsto|\alpha|$, which were considered in
Example~\ref{ex:standard-abs-val-subgrad}.
Indeed, these astral subgradients can be shown to be
\begin{equation}   \label{eq:ex:subgrad-sqrt-approaches-abs-val:1}
  \asubdiffext{\xbar_\alpha}
  =
  \begin{cases}
    \{-\ee_1\}    & \text{if $\alpha < 0$,} \\
    \Braces{\lambda\ee_1 :\: \lambda\in[-1,1]}
                  & \text{if $\alpha = 0$,} \\
    \{\ee_1\}     & \text{if $\alpha > 0$.}
  \end{cases}
\end{equation}
(These can be derived directly from \Cref{def:ast-subgrad},
but we will instead postpone the proof until
\Cref{ex:subgrad-sqrt-approaches-abs-val:cont},
by which point we will have developed tools
to simplify the arguments.)

Intuitively, $f(x_1,x_2)$ flattens out in the direction of $\ee_2$ as
$x_2\rightarrow+\infty$, so all the astral subgradients for points
$\xbar_\alpha$ are orthogonal to $\ee_2$, and therefore scalar
multiples of~$\ee_1$.
Moreover, their coordinate corresponding to $\ee_1$
matches the standard subgradients of the absolute value function.
Note in particular that $\zero\in\asubdiffext{\limray{\ee_2}}$,
corresponding to $\xbar_0=\limray{\ee_2}$ being $\fext$'s unique minimizer.%
\indexg{subdifferentials of extensions!examples|)}%
\indexg{subdifferentials, astral (primal)!examples|)}%
\indexg{Absolute value at infinity!astral subgradients of extension|)}%
\end{example}

For a convex and proper function $f:\Rn\rightarrow\Rext$, if
$\uu\in\partial f(\xx)$ then
Propositions~\ref{roc:thm23.4}(\ref{roc:thm23.4:b})
and~\refequiv{pr:stan-subgrad-equiv-props}{pr:stan-subgrad-equiv-props:a}{pr:stan-subgrad-equiv-props:b}
imply that $\xx\in\dom{f}$ and that $f^*(\uu)\in\R$.
The next proposition gives somewhat analogous properties for astral
subgradients, both for general astral functions on $\extspace$,
and specifically for extensions of functions on $\Rn$:

\begin{proposition}  \label{pr:subgrad-imp-in-cldom}
  Let $\xbar\in\extspace$ and $\uu\in\Rn$.
  \begin{letter-compact}
  \item  \label{pr:subgrad-imp-in-cldom:a}
    For $F:\extspace\rightarrow\Rext$,
    if $\uu\in\asubdifF{\xbar}$
    then $\xbar\in\cldom{F}$ and $\Fstar(\uu)\in\R$.
  \item  \label{pr:subgrad-imp-in-cldom:c}
    For $f:\Rn\rightarrow\Rext$,
    if $\uu\in\asubdiffext{\xbar}$
    or $\uu\in\asubdifftriv{\xbar}$
    then $f$ is proper,
    $\xbar\in\cldom{f}$, and $\fstar(\uu)\in\R$.
  \end{letter-compact}
\end{proposition}

\begin{proof}
~

\begin{proof-parts}
\pfpart{Part~(\ref{pr:subgrad-imp-in-cldom:a}):}
Let $F:\extspace\rightarrow\Rext$.
Suppose $\uu\in\asubdifF{\xbar}$, implying that
there exist $\beta\in\R$
and a sequence $\seq{\xbar_t}$ in $\eRn$
satisfying \Cref{def:ast-subgrad}.
Since $F(\xbar_t)\in\R$ for all $t$,
we have $\xbar_t\in\dom{F}$ for all $t$,
and so $\xbar=\lim\xbar_t\in\cldom{F}$.

Furthermore, 
\[
  \Fstar(\uu)
  =
  \sup_{\xbar'\in\eRn}
  [-F(\xbar')\plusd\xbar'\inprod\uu]
  \le
  \beta
  =
  \lim
  [-F(\xbar_t) + \xbar_t\inprod\uu]
  \le
  \Fstar(\uu).
\]
The first equality is by definition of astral conjugate
(Eq.~\ref{eq:Fstar-down-def}).
The first inequality and the second equality follow from the
alternate conditions~(\ref{eq:astral:subgrad:alt:1})
and~(\ref{eq:astral:subgrad:alt:2}) of \Cref{def:ast-subgrad}.
The last inequality is again by
\eqref{eq:Fstar-down-def}.
Thus, $\Fstar(\uu)=\beta\in\R$.

\pfpart{Part~(\ref{pr:subgrad-imp-in-cldom:c}):}
Let $f:\Rn\rightarrow\Rext$,
and suppose either $\uu\in\asubdiffext{\xbar}$ or
$\uu\in\asubdifftriv{\xbar}$.
Then applying
part~(\ref{pr:subgrad-imp-in-cldom:a})
with $F=\fext$ or $F=\ftriv$ yields that
$\fstar(\uu)\in\R$ since
$\fextstar=\ftrivstar=\fstar$
(\Cref{pr:fextstar-is-fstar}),
and that
$\xbar\in\cldom{f}$ since
$\dom{\ftriv}=\dom{f}$ and
$\cldomfext=\cldom{f}$
(by \Cref{pr:h:1}\ref{pr:h:1c}).

Furthermore, if $f$ is not proper, then
either $\fstar\equiv-\infty$ or
$\fstar\equiv+\infty$
(by
\Cref{pr:conj-props}\ref{pr:conj-props:c1}\ref{pr:conj-props:c2}),
contradicting that $\fstar(\uu)\in\R$, as just argued.
Therefore, $f$ is proper.
\qedhere
\end{proof-parts}
\end{proof}

Thus, if
$\asubdifF{\xbar}\neq\emptyset$ then $\xbar$ must be in the closure of
$\dom{F}$.
Note, however, that $\xbar$ need not be in $\dom{F}$ itself, as
we saw, for instance, in \Cref{ex:standard-abs-val-subgrad-cont}.

\indexg{subdifferentials, astral (primal)!characterizations|(}%
In the next \namecref{pr:equiv-ast-subdif-defn},
we give several equivalent formulations of what
it means for $\uu$ to be an astral subgradient of $F$ at $\xbar$.
Here, we assume $\Fstar(\uu)\in\R$, which, as we just proved,
is a necessary condition for $\uu$ to be an astral subgradient.

Our first equivalent formulation
(in \Cref{pr:equiv-ast-subdif-defn}\ref{pr:equiv-ast-subdif-defn:b})
captures the observation
from the proof of \Cref{pr:subgrad-imp-in-cldom}(\ref{pr:subgrad-imp-in-cldom:a})
that we can always take $\beta$
in the definition of astral subgradient (\Cref{def:ast-subgrad})
to be equal to $\Fstar(\uu)$. In that case,
condition~(\ref{en:ast-sub-defn-cond-1}) of
that definition
is satisfied automatically (and so can be omitted).

In our next equivalent formulation
(in \Cref{pr:equiv-ast-subdif-defn}\ref{pr:equiv-ast-subdif-defn:plusd}),
we further relax the requirement that
$F(\xbar_t)\in\R$ and $\xbar_t\inprod\uu\in\R$ for all $t$, and replace
ordinary addition in the requirement that $-F(\xbar_t)+\xbar_t\inprod\uu\to\Fstar(\uu)$
with downward addition, obtaining the requirement
$-F(\xbar_t)\plusd \xbar_t\inprod\uu\to\Fstar(\uu)$, much like the expression
in the definition of astral conjugate (Eq.~\ref{eq:Fstar-down-def}).

Building on this characterization, we next obtain a formulation
in terms of epigraphs in a way that mirrors the relationship between the two definitions
of astral conjugate given in Eqs.~(\ref{eq:Fstar-def}) and~(\ref{eq:Fstar-down-def}). Similar
to the characterization with downward sum, this epigraph-based characterization
(in \Cref{pr:equiv-ast-subdif-defn}\ref{pr:equiv-ast-subdif-defn:epi})
does
not require the finiteness of $F(\xbar_t)$ or $\xbar_t\inprod\uu$.

Finally, using epigraphs, we can restate the definition of astral subgradient
without using sequences at all, instead using properties of points in the closure
of $F$'s epigraph. Two such formulations are presented in
\Cref{pr:equiv-ast-subdif-defn}(\ref{pr:equiv-ast-subdif-defn:c},\ref{pr:equiv-ast-subdif-defn:d}).

\begin{proposition}  \label{pr:equiv-ast-subdif-defn}
  Let $F:\extspace\rightarrow\Rext$,
  let $\xbar\in\extspace$ and $\uu\in\Rn$.
  Also,
  let $\PPx=[\Idnn,\zerov{n}]$ and
  $\PPy=[\trans{\zerov{n}},1]$.
  Assume $\Fstar(\uu)\in\R$.
  Then the following are equivalent:
  \begin{letter-compact}
  \item  \label{pr:equiv-ast-subdif-defn:a}
    $\uu\in\partial F(\xbar)$.

  \item  \label{pr:equiv-ast-subdif-defn:b}
    There exists a sequence $\seq{\xbar_t}$ in $\eRn$ with
    $F(\xbar_t)\in\R$ and $\xbar_t\inprod\uu\in\R$ for all $t$,
    and such that
    $\xbar_t\to\xbar$,
    $F(\xbar_t)\to F(\xbar)$,
    and
    $\xbar_t\inprod\uu-F(\xbar_t)\to\Fstar(\uu)$.
    
  \item  \label{pr:equiv-ast-subdif-defn:plusd}
    There exists a sequence $\seq{\xbar_t}$ in $\eRn$
    such that
    $\xbar_t\to\xbar$,
    $F(\xbar_t)\to F(\xbar)$,
    and
    $-F(\xbar_t) \plusd \xbar_t\inprod\uu\to\Fstar(\uu)$.
    
  \item  \label{pr:equiv-ast-subdif-defn:epi}
    There exists a sequence $\seq{\rpair{\xbar_t}{y_t}}$ in $\epi F$
    such that
    $\xbar_t\rightarrow\xbar$,
    $y_t \rightarrow F(\xbar)$,
    and
    $\xbar_t\cdot\uu - y_t \rightarrow \Fstar(\uu)$.

  \item  \label{pr:equiv-ast-subdif-defn:c}
    There exists $\zbar\in \clepi{F}$ such that
    $\PPx\zbar = \xbar$,
    $\PPy\zbar = F(\xbar)$,
    and
    $\zbar \cdot \rpair{\uu}{-1} = \Fstar(\uu)$.
  \item  \label{pr:equiv-ast-subdif-defn:d}
    $\zbar'\cdot \rpair{\uu}{-1}$, as a function of $\zbar'$, is
    maximized over $\clepi{F}$ by some $\zbar\in\clepi{F}$ with
    $\PPx\zbar = \xbar$
    and
    $\PPy\zbar = F(\xbar)$.
  \end{letter-compact}
\end{proposition}

\begin{proof}
~

\begin{proof-parts}
\pfpart{%
  (\ref{pr:equiv-ast-subdif-defn:b})
  $\Rightarrow$
  (\ref{pr:equiv-ast-subdif-defn:a}):
}
Suppose $\seq{\xbar_t}$ is a sequence as in
part~(\ref{pr:equiv-ast-subdif-defn:b}).
Let $\beta=\Fstar(\uu)$.
Then from the definition of the conjugate (Eq.~\ref{eq:Fstar-down-def}),
$\beta\ge-F(\xbar')\plusd\xbar'\inprod\uu$ for all $\xbar'\in\eRn$,
so $\beta$ and $\seq{\xbar_t}$ satisfy the alternate
conditions~(\ref{eq:astral:subgrad:alt:1})
and~(\ref{eq:astral:subgrad:alt:2})
of \Cref{def:ast-subgrad}, proving that $\uu\in\partial F(\xbar)$.

\pfpart{%
  (\ref{pr:equiv-ast-subdif-defn:plusd})
  $\Rightarrow$
  (\ref{pr:equiv-ast-subdif-defn:b}):
}
Suppose $\seq{\xbar_t}$ is a sequence as in
part~(\ref{pr:equiv-ast-subdif-defn:plusd}).
Since $\Fstar(\uu)\in\R$, this implies, for all $t$ sufficiently
large, that $-F(\xbar_t) \plusd \xbar_t\cdot\uu$ is finite,
implying further that $F(\xbar_t)$ and $\xbar_t\cdot\uu$ are finite as
well.
By discarding all other elements, we thus obtain a sequence as in
part~(\ref{pr:equiv-ast-subdif-defn:b}).

\pfpart{%
  (\ref{pr:equiv-ast-subdif-defn:epi})
  $\Rightarrow$
  (\ref{pr:equiv-ast-subdif-defn:plusd}):
}
Suppose $\seq{\rpair{\xbar_t}{y_t}}$ is a sequence as in
part~(\ref{pr:equiv-ast-subdif-defn:epi}).
Then
\begin{align}
  \Fstar(\uu)
  =
  \lim[-y_t\plusd\xbar_t\inprod\uu]
  &\le
  \liminf[-F(\xbar_t)\plusd\xbar_t\inprod\uu]
  \nonumber
  \\
  &\leq
  \limsup[-F(\xbar_t)\plusd\xbar_t\inprod\uu]
  \le
  \Fstar(\uu),
  \label{eq:pr:equiv-ast-subdif-defn:1}
\end{align}
where
the first inequality is because $y_t\ge F(\xbar_t)$,
and the last is by definition of conjugate (Eq.~\ref{eq:Fstar-down-def}).
Thus, $-F(\xbar_t)\plusd\xbar_t\inprod\uu\to\Fstar(\uu)$.
Furthermore, by the same reasoning as above,
since $\Fstar(\uu)\in\R$, we must have $F(\xbar_t)\in\R$ and $\xbar_t\inprod\uu\in\R$
for all $t$ sufficiently large; by discarding all other sequence elements,
we assume these hold for all $t$.

It remains to argue
that $F(\xbar_t)\to F(\xbar)$.
We have
\begin{align*}
    \lim F(\xbar_t)
    &=\lim\bigParens{[F(\xbar_t)-\xbar_t\inprod\uu]+\xbar_t\inprod\uu}
\nonumber
\\
    &=-\Fstar(\uu)+\xbar\inprod\uu
\nonumber
\\
    &=\lim\bigParens{[y_t-\xbar_t\inprod\uu]+\xbar_t\inprod\uu}
    =\lim y_t = F(\xbar).
\end{align*}
The second and third equalities both follow from
continuity of addition
(\Cref{prop:lim:eR}\ref{i:lim:eR:sum}, noting that summability of the
limits is ensured since $\Fstar(\uu)\in\R$),
and because
$\xbar_t\inprod\uu-F(\xx_t)\to\Fstar(\uu)$
(by Eq.~\ref{eq:pr:equiv-ast-subdif-defn:1})
and
$\xbar_t\inprod\uu-y_t\to\Fstar(\uu)$
(by assumption),
and finally because $\xbar_t\inprod\uu\to\xbar\inprod\uu$
(\Cref{thm:i:1}\ref{thm:i:1c}).

\pfpart{%
  (\ref{pr:equiv-ast-subdif-defn:c})
  $\Rightarrow$
  (\ref{pr:equiv-ast-subdif-defn:epi}):
}
Suppose $\zbar$ is as specified in
part~(\ref{pr:equiv-ast-subdif-defn:c}).
Because $\zbar\in\clepi{F}$, there exists a sequence
$\seq{\rpair{\xbar_t}{y_t}}$ in $\epi F$ that converges to $\zbar$.
Then
\begin{align*}
  \xbar_t
  &=
  \PPx \rpair{\xbar_t}{y_t}
  \rightarrow
  \PPx\zbar = \xbar,
  \\
  y_t
  &=
  \PPy \rpair{\xbar_t}{y_t}
  \rightarrow
  \PPy \zbar = F(\xbar),
\end{align*}
with convergence in each line following from continuity
of linear maps (\Cref{thm:linear:cont}\ref{thm:linear:cont:b}),
and the first equality in each line from
\Cref{pr:xy-pairs-props}(\ref{pr:xy-pairs-props:c}).
Moreover,
\begin{equation*}
  \xbar_t \cdot \uu - y_t
  =
  \rpair{\xbar_t}{y_t} \cdot \rpair{\uu}{-1}
  \rightarrow
  \zbar \cdot \rpair{\uu}{-1} = \Fstar(\uu)
\end{equation*}
(using
\Cref{pr:xy-pairs-props}\ref{pr:xy-pairs-props:b}
and
\Cref{thm:i:1}\ref{thm:i:1c}).
Thus, the sequence
$\seq{\rpair{\xbar_t}{y_t}}$ satisfies
part~(\ref{pr:equiv-ast-subdif-defn:epi}).

\pfpart{%
  (\ref{pr:equiv-ast-subdif-defn:d})
  $\Rightarrow$
  (\ref{pr:equiv-ast-subdif-defn:c}):
}
Suppose $\zbar$ is as specified in
part~(\ref{pr:equiv-ast-subdif-defn:d}).
Then
\begin{align*}
  \zbar\cdot \rpair{\uu}{-1}
  =
  \max_{\zbar'\in\clepi{F}} \zbar'\cdot \rpair{\uu}{-1}
  &=
  \sup_{\rpair{\xbar'}{y'}\in\epi F} \rpair{\xbar'}{y'}\cdot \rpair{\uu}{-1}
  \\
  &=
  \sup_{\rpair{\xbar'}{y'}\in\epi F} [\xbar'\cdot\uu - y']
  = \Fstar(\uu).
\end{align*}
The second equality is
by \Cref{pr:sup-f-Abar} since the function
$\zbar'\mapsto\zbar'\cdot\rpair{\uu}{-1}$, for $\zbar'\in\extspacnp$,
is continuous
(\Cref{thm:i:1}\ref{thm:i:1c}).
The third equality is
by
\Cref{pr:xy-pairs-props}(\ref{pr:xy-pairs-props:b}).
And the fourth is by
definition of astral conjugate
(Eq.~\ref{eq:Fstar-def}).
Thus, $\zbar$ satisfies
part~(\ref{pr:equiv-ast-subdif-defn:c}).

\pfpart{%
  (\ref{pr:equiv-ast-subdif-defn:a})
  $\Rightarrow$
  (\ref{pr:equiv-ast-subdif-defn:d}):
}
Suppose $\beta$
and $\seq{\xbar_t}$ are as in
\Cref{def:ast-subgrad},
and let $\zbar_t=\rpair{\xbar_t}{F(\xbar_t)}$ for all~$t$.
Then $\seq{\zbar_t}$ is a sequence in $\epi F$, which,
by sequential compactness,
must have a convergent subsequence; by discarding
all other elements, we can assume the entire
sequence converges to some point $\zbar\in\extspacnp$, which must then
be in
$\clepi F$.
Similar to previous arguments,
$\xbar_t = \PPx\zbar_t \rightarrow \PPx\zbar$,
so $\PPx\zbar = \xbar$ since $\xbar_t\rightarrow\xbar$.
Likewise,
$F(\xbar_t) = \PPy\zbar_t \rightarrow \PPy\zbar$,
implying
$\PPy\zbar = F(\xbar)$,
and
$\xbar_t\cdot\uu - F(\xx_t) = \zbar_t \cdot \rpair{\uu}{-1} \rightarrow \zbar\cdot\rpair{\uu}{-1}$,
implying
$\zbar\cdot\rpair{\uu}{-1} = \beta$,
using condition~(\ref{eq:astral:subgrad:alt:2}) of
\Cref{def:ast-subgrad}.

For all $\rpair{\xbar'}{y'}\in\epi{F}$,
we have
\begin{equation}   \label{eq:pr:equiv-ast-subdif-defn:2}
  \rpair{\xbar'}{y'} \cdot \rpair{\uu}{-1}
  =
  \xbar'\cdot\uu-y'
  =
  -y'\plusd\xbar'\cdot\uu
  \leq
  -F(\xbar')\plusd\xbar'\cdot\uu
  \leq
  \beta,
\end{equation}
where the first equality is by
\Cref{pr:xy-pairs-props}(\ref{pr:xy-pairs-props:b}),
the first inequality is  because ${F(\xbar')\leq y'}$,
and the second inequality is by
condition~(\ref{eq:astral:subgrad:alt:1}) in the definition of astral
subgradient.
Thus,
\[
  \sup_{\zbar'\in\clepi{F}} \zbar'\cdot\rpair{\uu}{-1}
  =
  \sup_{\zbar'\in\epi{F}} \zbar'\cdot\rpair{\uu}{-1}
  \leq
  \beta
  =
  \zbar\cdot\rpair{\uu}{-1},
\]
where the first equality is
by \Cref{pr:sup-f-Abar}
(since $\zbar'\mapsto\zbar'\cdot\rpair{\uu}{-1}$ is continuous),
and the inequality is by
\eqref{eq:pr:equiv-ast-subdif-defn:2}.
Therefore,
$ \zbar' \cdot \rpair{\uu}{-1} $ is maximized over $\clepi{F}$
when $\zbar'=\zbar$,
so all conditions of part~(\ref{pr:equiv-ast-subdif-defn:d}) are
satisfied.%
\indexg{subdifferentials, astral (primal)!characterizations|)}%
\qedhere

\end{proof-parts}
\end{proof}

\indexg{subdifferentials, astral (primal)!astral closedness and|(}%
\indexg{subdifferentials, astral (primal)!lower semicontinuity and|(}%
We next show that at any point $\xbar$ where an astral function is
astral subdifferentiable, the function must also be astral closed,
and therefore lower semicontinuous.
Moreover, any subgradient $\uu$ at $\xbar$
achieves the supremum in the definition of the biconjugate $\Fdub(\xbar)$
(as in Eq.~\ref{eq:psistar-def:2} with $\psi=\Fstar$).

\begin{theorem}    \label{thm:subgrad-then-lsc}
  Let $F:\extspace\rightarrow\Rext$, let $\xbar\in\extspace$ and $\uu\in\Rn$.
  If $\uu\in\asubdifF{\xbar}$,
  then
  \[
    F(\xbar)=(\lsc F)(\xbar)=\Fdub(\xbar)=-\Fstar(\uu)+\xbar\inprod\uu.
  \]
  (Consequently,
  the first two equalities hold
  whenever $\asubdifF{\xbar}\neq\emptyset$.)
\end{theorem}

\begin{proof}
Suppose $\uu\in\asubdifF{\xbar}$,
implying $\Fstar(\uu)\in\R$ by
\Cref{pr:subgrad-imp-in-cldom}(\ref{pr:subgrad-imp-in-cldom:a}),
and, by
\Cref{pr:equiv-ast-subdif-defn},
that there exists a sequence $\seq{\xbar_t}$ satisfying
part~(\ref{pr:equiv-ast-subdif-defn:b}) of that
\namecref{pr:equiv-ast-subdif-defn}.
Then
\begin{align*}
  F(\xbar)
  =
  \lim F(\xbar_t)
  &=
  \lim \bigParens{[F(\xbar_t) - \xbar_t\cdot\uu] + \xbar_t\cdot\uu}
\\
  &=
  -\Fstar(\uu)+\xbar\inprod\uu
\\
  &\le
  \Fdub(\xbar)
  \le
  (\lsc F)(\xbar)
  \le
  F(\xbar).
\end{align*}
The third equality is because
$\xbar_t\cdot\uu - F(\xbar_t) \rightarrow \Fstar(\uu)$
(a defining property of the sequence $\seq{\xbar_t}$)
and ${\xbar_t\cdot\uu\rightarrow\xbar\cdot\uu}$
(\Cref{thm:i:1}\ref{thm:i:1c}),
and by continuity of addition
(\Cref{prop:lim:eR}\ref{i:lim:eR:sum}, noting $\Fstar(\uu)\in\R$ for
summability of the limits).
The first inequality is by definition of dual conjugate
(Eq.~\ref{eq:psistar-def:2}).
The last two inequalities are by
\Cref{pr:conj-lsc-props}.
This proves the
\namecref{thm:subgrad-then-lsc}.%
\indexg{subdifferentials, astral (primal)!astral closedness and|)}%
\end{proof}

\Cref{thm:subgrad-then-lsc} shows that at any point
$\xbar\in\extspace$ where a function $F:\extspace\rightarrow\Rext$ is
astral subdifferentiable, it is also lower semicontinuous.
The next theorem shows moreover that at all such points,
the astral subdifferentials of $F$ and $\lsc F$ are identical:

\begin{theorem}    \label{thm:subgrad-F-lsc-F}
  Let $F:\extspace\rightarrow\Rext$, and let $\xbar\in\extspace$.
  Then
  \[
     \asubdifF{\xbar}
     =
     \begin{cases}
       \asubdiflscF{\xbar}  & \text{if $(\lsc F)(\xbar)=F(\xbar)$,} \\
       \emptyset            & \text{otherwise.}
     \end{cases}
  \]
\end{theorem}
\begin{proof}
If $(\lsc F)(\xbar)\neq F(\xbar)$, then
$\asubdifF{\xbar}=\emptyset$ by \Cref{thm:subgrad-then-lsc}.
We therefore assume henceforth that
$(\lsc F)(\xbar) = F(\xbar)$.

Let $G=\lsc F$, so $\Gstar=\Fstar$
by \Cref{pr:conj-lsc-props}.
Let $\uu\in\Rn$.
We aim to show that $\uu\in\asubdifF{\xbar}$ if and only if
$\uu\in\asubdifG{\xbar}$.
If $\Fstar(\uu)=\Gstar(\uu)\not\in\R$ then
$\uu\not\in\asubdifF{\xbar}$ and
$\uu\not\in\asubdifG{\xbar}$
by
\Cref{pr:subgrad-imp-in-cldom}(\ref{pr:subgrad-imp-in-cldom:a}).
So assume henceforth that $\Fstar(\uu)=\Gstar(\uu)\in\R$.
Then, by \Cref{pr:equiv-ast-subdif-defn}(\ref{pr:equiv-ast-subdif-defn:a},\ref{pr:equiv-ast-subdif-defn:d}),
$\uu\in\partial F(\xbar)$
if and only if $\uu$ satisfies the conditions of 
part~(\ref{pr:equiv-ast-subdif-defn:d}) of that
\namecref{pr:equiv-ast-subdif-defn},
which is the case if and only if $\uu$ satisfies those
same conditions for the function $G$, since $G(\xbar)=F(\xbar)$ by assumption,
and
$\clepi{G}=\clepi{F}$ by
\Cref{pr:epiF-epi-lscF}(\ref{pr:epiF-epi-lscF:c}).
By \Cref{pr:equiv-ast-subdif-defn}(\ref{pr:equiv-ast-subdif-defn:d},\ref{pr:equiv-ast-subdif-defn:a}), this is in turn equivalent
to $\uu\in\asubdifG{\xbar}$. Thus,
$\asubdifF{\xbar}=\asubdifG{\xbar}$.%
\indexg{subdifferentials, astral (primal)!lower semicontinuity and|)}%
\indexg{subdifferentials, astral (primal)|)}%
\end{proof}

\indexg{subdifferentials, standard!generalized by astral subgradients|(}%
\indexg{subdifferentials of extensions!finite points@at finite points|(}%
Using
\Cref{thm:subgrad-F-lsc-F},
we can characterize the relationship
between the standard
subgradients of a proper function $f:\Rn\rightarrow\Rext$ and the
astral subgradients of the function's lower semicontinuous and trivial
extensions, $\fext$ and $\ftriv$.
In particular, the standard subdifferential of $f$ is always the same
as the astral subdifferential of $\ftriv$,
and is also the same as the astral subdifferential of
$\fext$ at points where $f$ is lower semicontinuous.
In this sense, for proper functions,
the astral subdifferentials of extensions generalize standard
subdifferentials.

\begin{theorem}  \label{pr:asubdiffext-at-x-in-rn}
  Let $f:\Rn\rightarrow\Rext$ be proper,
  and let $\xx\in\Rn$.
  Then:
  \begin{letter}
  \item  \label{pr:asubdiffext-at-x-in-rn:a}
    $\asubdifftriv{\xx} = \partial f(\xx)$.
  \item  \label{pr:asubdiffext-at-x-in-rn:b}
    $\displaystyle
       \partial f(\xx)
       =
       \begin{cases}
         \asubdiffext{\xx}
                           & \text{if $\fext(\xx) = f(\xx)$,} \\
         \emptyset         & \mbox{otherwise.}
       \end{cases}
    $
  \item  \label{pr:asubdiffext-at-x-in-rn:c}
    $\asubdiffext{\xx} = \partial (\lsc f)(\xx)$.
  \end{letter}
\end{theorem}

\begin{proof}
  ~

\begin{proof-parts}
\pfpart{Part~(\ref{pr:asubdiffext-at-x-in-rn:a}):}
Suppose $\uu\in\partial\ftriv(\xx)$. Then
$\ftriv(\xx) = \xx\cdot\uu - \ftrivstar(\uu)$ (by \Cref{thm:subgrad-then-lsc}). Since $\ftriv(\xx)=f(\xx)$ and $\ftrivstar(\uu)=\fstar(\uu)$ (by \Cref{pr:fextstar-is-fstar}),
we also have $f(\xx) = \xx\cdot\uu - \fstar(\uu)$, so
$\uu\in\partial f(\xx)$ by
\Cref{pr:stan-subgrad-equiv-props}(\ref{pr:stan-subgrad-equiv-props:b},\ref{pr:stan-subgrad-equiv-props:a}).

Conversely, suppose $\uu\in\partial f(\xx)$.
Then $\fstar(\uu) = \xx\cdot\uu - f(\xx)$ by
\Cref{pr:stan-subgrad-equiv-props}(\ref{pr:stan-subgrad-equiv-props:a},\ref{pr:stan-subgrad-equiv-props:b}).
Since $f$ is proper, $f(\xx)>-\infty$,
and hence $\fstar(\uu)<+\infty$. Also, $f\not\equiv+\infty$, so
$\fstar(\uu)>-\infty$ (by \Cref{pr:conj-props}\ref{pr:conj-props:c1}).
Thus, $\ftrivstar(\uu)=\fstar(\uu)\in\R$,
implying also $f(\xx)\in\R$.
Let $\zz=\rpair{\xx}{f(\xx)}$. Then $\zz\in\epi f=\epi\ftriv$. Also, $\PPx\zz=\xx$,
$\PPy\zz=f(\xx)=\ftriv(\xx)$, and
$\zz\cdot\rpair{\uu}{-1}=\xx\inprod\uu-f(\xx)=\fstar(\uu)=\ftrivstar(\uu)$.
Thus, $\uu\in\partial\ftriv(\xx)$
by
\Cref{pr:equiv-ast-subdif-defn}(\ref{pr:equiv-ast-subdif-defn:c},\ref{pr:equiv-ast-subdif-defn:a}),
since $\zz$ satisfies part~(\ref{pr:equiv-ast-subdif-defn:c}) of that
\namecref{pr:equiv-ast-subdif-defn}
(applied to $\ftriv$).

\pfpart{Part~(\ref{pr:asubdiffext-at-x-in-rn:b}):}
We have
\begin{align*}
  \partial f(\xx)
  =
  \asubdifftriv{\xx}
  &=
  \begin{cases}
    \asubdiflscftriv{\xx} & \text{if $(\lsc\ftriv)(\xx) = \ftriv(\xx)$,} \\
    \emptyset             & \text{otherwise}
  \end{cases}
  \\
  &=
  \begin{cases}
    \asubdiffext{\xx} &
                    \text{if $\fext(\xx) = f(\xx)$,}
    \\
    \emptyset         & \text{otherwise.}
  \end{cases}
\end{align*}
The first equality is by part~(\ref{pr:asubdiffext-at-x-in-rn:a}).
The second is by \Cref{thm:subgrad-F-lsc-F}
(applied to $\ftriv$).
The last equality is because $\fext=\lsc\ftriv$
(\Cref{pr:lsc-ftriv-is-fext})
and $\ftriv(\xx)=f(\xx)$.

\pfpart{Part~(\ref{pr:asubdiffext-at-x-in-rn:c}):}
Let $g=\lsc f$.
Then
by \Cref{pr:h:1}(\ref{pr:h:1a},\ref{pr:h:1aa}),
$\fext=\gext$ so $\gext(\xx)=\fext(\xx)=(\lsc f)(\xx)=g(\xx)$.
Consequently,
applying part~(\ref{pr:asubdiffext-at-x-in-rn:b}) to $g$
yields
$\asubdiffext{\xx} = \asubdifgext{\xx} = \partial g(\xx) = \partial(\lsc f)(\xx)$.
\qedhere
\end{proof-parts}
\end{proof}

\Cref{pr:asubdiffext-at-x-in-rn}
is false for all improper functions.
Indeed, if $f$ is improper then
$\asubdifftriv{\xx} = \asubdiffext{\xx} = \emptyset$
for all $\xx\in\Rn$ by
\Cref{pr:subgrad-imp-in-cldom}(\ref{pr:subgrad-imp-in-cldom:c}).
On the other hand, 
$\partial f(\xx)=\Rn$ at every point where $f(\xx)=-\infty$,
and at every point in $\Rn$ if $f\equiv+\infty$.
(Note further that $f$ is lower semicontinuous at all such points.)%
\indexg{subdifferentials of extensions!finite points@at finite points|)}%
\indexg{subdifferentials, standard!generalized by astral subgradients|)}%

\section{Astral dual subgradients}
\label{sec:astral-dual-subgrad}

\begin{figure}
  \centering
  \includegraphics{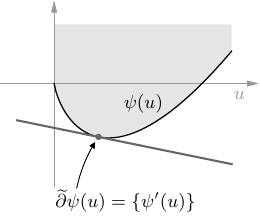}\hspace{0.15\textwidth}%
  \includegraphics{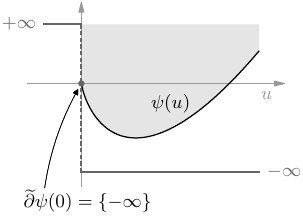}
  \caption[Astral dual subgradients of the function $\psi$ from \Cref{ex:entropy-1d}]{%
  \indexf{subdifferentials, astral dual!examples}%
\emph{Astral dual subgradients of the function $\psi$ from \Cref{ex:entropy-1d}.}}
  \label{fig:xlogx}%
\end{figure}

The
astral subdifferential introduced in the previous section
allows us to assign finite subgradients to points at infinity
for functions that are defined on $\extspace$.
\indexg{subdifferentials, astral dual|(}%
We turn next to a different kind of subdifferential,
called the astral
dual subdifferential, which, for a function $\psi$ defined on $\Rn$, assigns
subgradients to finite points in $\Rn$, but the subgradients themselves are
allowed to be infinite points in $\eRn$.
As with astral primal subgradients, astral dual subgradients correspond
to generalized affine functions, which support the graph of $\psi$, but by allowing infinite subgradients we can now
support $\psi$'s graph at points where only vertical supporting hyperplanes
exist. Thus, astral dual subgradients allow us to conceptualize
``infinite slopes,'' which are not captured by standard subgradients.
\indexg{subdifferentials, astral dual!examples|(}%
Here is an example:

\begin{example}   \label{ex:entropy-1d}
Let $\psi:\R\rightarrow\Rext$ be defined, for $u\in\R$, by
\[
   \psi(u)
       =
       \begin{cases}
               u \ln u    & \text{if $u\geq 0$,} \\
               +\infty    & \mbox{otherwise,}
       \end{cases}
\]
(where, as usual, $0\ln 0 = 0$).
As depicted in \Cref{fig:xlogx},
for $u > 0$,
the subgradients of this function are the same as its
derivative and they also coincide with astral dual subgradients
(as we will see momentarily).
But for $u\leq 0$, the function has no standard subgradients, even at
$u=0$, which is in the function's effective domain.
This is because, for any $x\in\partial \psi(0)$,
we must have
$\psistar(x) = 0\inprod x - \psi(0) = 0$
(by \Cref{pr:stan-subgrad-equiv-props}\ref{pr:stan-subgrad-equiv-props:a}\ref{pr:stan-subgrad-equiv-props:b}).
However, this is impossible since it can be calculated that
$\psistar(x)=e^{x-1}>0$ for all $x\in\R$.
Nevertheless, at $u=0$, we will see that
$-\infty$ is an astral dual subgradient, corresponding to the infinite
rate at which the function is decreasing at this point.%
\indexg{subdifferentials, astral dual!examples|)}%
\end{example}

To define astral dual subgradients,
let $\psi:\Rn\rightarrow\Rext$ be any function, changing
notation to emphasize the switch to dual space;
later, we will often take $\psi$ to be the conjugate of some
function on astral space.
Replacing $f$ by $\psi$ and swapping variable names in
\eqref{eqn:standard-subgrad-ineq},
we can say that $\xx\in\Rn$ is a standard
subgradient of $\psi$ at $\uu\in\Rn$ if
\begin{equation}  \label{eqn:psi-subgrad:1}
   \psi(\uu') \geq \psi(\uu) + \xx\cdot(\uu'-\uu)
\end{equation}
for all $\uu'\in\Rn$.
We extend this notion to astral space by allowing infinite subgradients
$\xbar\in\eRn$, and by replacing ordinary addition with downward addition
on the right-hand side of \eqref{eqn:psi-subgrad:1}.
This is analogous to the use of downward addition in the definition of astral conjugate (Eq.~\ref{eq:Fstar-down-def}),
and is in a sense the most permissive way to ``resolve a tie,'' as it results in the lower bound evaluating to $-\infty$ whenever $+\infty$ and $-\infty$ are being added.

\begin{definition}
\label{def:dual-subgrad}
\indexg{subdifferentials, astral dual!defined|(}%
  Let $\psi:\Rn\rightarrow\Rext$, and
  let $\xbar\in\extspace$ and $\uu\in\Rn$.
  We say that $\xbar$ is an \emph{astral dual subgradient} of
  $\psi$ at $\uu$ if
  \begin{equation}
  \label{eqn:psi-subgrad:3-alt}
    \psi(\uu') \geq \psi(\uu) \plusd \xbar\cdot(\uu'-\uu)
  \end{equation}
  for all $\uu'\in\Rn$.
  The
  \emph{astral dual subdifferential} of $\psi$ at $\uu$,
  denoted
\indexm{d psiu}{$\adsubdifpsi{\uu}$}{astral dual subdifferential}%
  $\adsubdifpsi{\uu}$,
  is the set of all such astral dual subgradients
  of $\psi$ at $\uu$.
  We say that $\psi$ is \emph{astral dual subdifferentiable} at $\uu$ if
  $\adsubdifpsi{\uu}$ is not empty.
  \end{definition}

  Equivalently, 
  $\xbar$ is an astral dual subgradient of $\psi$ at $\uu$ if and only if
  \begin{equation}
  \label{eqn:psi-subgrad:3}
    {-\psi(\uu)} \geq -\psi(\uu') \plusd \xbar\cdot(\uu'-\uu)
  \end{equation}
  for all $\uu'\in\Rn$ (as follows from Eq.~\ref{eqn:psi-subgrad:3-alt} by
  \Cref{pr:plusd-props}\ref{pr:plusd-props:e}),
or also
if and only if
\begin{equation}  \label{eqn:psi-subgrad:3-plusu}
  {-\psi}(\uu) \plusu \psi(\uu') \geq \xbar\cdot(\uu'-\uu)
\end{equation}
for all $\uu'\in\Rn$ (by
\Cref{pr:plusd-plusu-props}\ref{pr:plusd-plusu-props:b}).%
\indexg{subdifferentials, astral dual!defined|)}%

When $\psi(\uu)\in\R$, astral dual subgradients
coincide with the ``extended subgradients'' of
\idxwaggoner\citet{waggoner25}.
Since our definition also captures when
$\psi(\uu)\in\set{\pm\infty}$, it is slightly more general.

  Compared with the definition of astral primal subgradient (\Cref{def:ast-subgrad}), we do not impose the condition 
  that the gap between the function $\psi$ and its
  lower bound should go to zero along some sequence converging to $\uu$.
  Instead, the lower bound is required to exactly match
  the value of $\psi$ at $\uu$ (which is ensured by the form of the lower bound).

We use the notation $\adsubdifplain{\psi}$ to distinguish the astral
dual subdifferential from the standard subdifferential
$\partial \psi$, since either operation can be applied to an ordinary
function $\psi$ over $\Rn$.
(In contrast, the astral subdifferential $\asubdifplain{F}$ is
only applied to functions $F$ over $\extspace$.)

As discussed above,
the astral primal subdifferential
$\asubdifF{\xbar}$ captures finite subgradients $\uu\in\Rn$
at astral (and therefore potentially infinite) points
$\xbar\in\extspace$,
while the astral dual subdifferential
$\adsubdifpsi{\uu}$ captures astral (and so potentially infinite)
subgradients $\xbar\in\extspace$
at finite points $\uu\in\Rn$
(as we saw in \Cref{ex:entropy-1d}).
\indexg{subdifferentials, astral dual!examples|(}%
The next example shows further how it is possible for an astral dual
subgradient to encode an infinite slope (or rate of change) in one
direction, while simultaneously encoding an ordinary finite slope in
another direction:

\begin{figure}[t]
  \centering
  \includegraphics{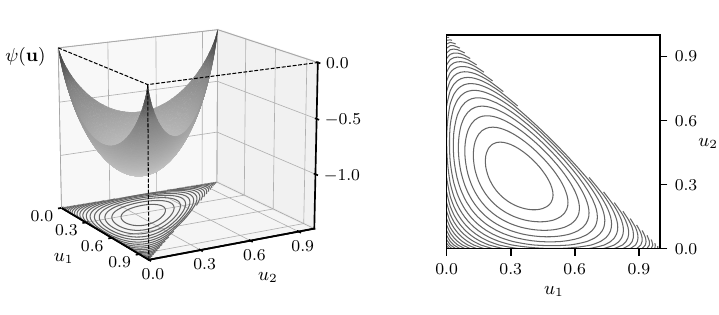}
  \mycaption{Entropy}{%
\indexf{Entropy function}%
    The function $\psi$ from \Cref{ex:entropy-ast-dual-subgrad}.
  }%
  \label{fig:entropy}%
\end{figure}

\begin{example}
[Astral dual subgradients of the entropy function]
\label{ex:entropy-ast-dual-subgrad}
\indexg{Entropy function|(}%
Let
\[
  \Delta
  =
  \regBraces{\uu\in\Rpos^2 :\:  u_1+u_2\leq 1},
\]
so that, for any $\uu\in\Delta$, the triple $\mtuple{u_1,u_2,1-u_1-u_2}$ describes a probability distribution over three possible outcomes. Let
$\psi:\R^2\rightarrow\Rext$ be the associated negative entropy
function, defined for $\uu\in\R^2$ as
\begin{equation*} %
   \psi(\uu)
   =
   \begin{cases}
     u_1 \ln u_1
     +
     u_2 \ln u_2
     +
     (1-u_1-u_2) \ln(1-u_1-u_2)
     & \text{if $\uu\in\Delta$,}
     \\
     +\infty
     & \text{otherwise.}
   \end{cases}
\end{equation*}
This function, depicted in \Cref{fig:entropy},
is closed, proper, convex, with effective domain
$\dom{\psi}=\Delta$.
It is also differentiable at all points in $\Delta$'s interior,
so its standard and astral dual subgradients are the same
as its gradient at all such points.
In particular, for $\uu\in\intr{\Delta}$,
\begin{equation}  \label{eq:ex:entropy-ast-dual-subgrad:3}
  \nabla\psi(\uu)
  =
  \trans{\biggBracks{
    \ln\biggParens{\frac{u_1}{1-u_1-u_2}},\,
    \ln\biggParens{\frac{u_2}{1-u_1-u_2}}
  }}.%
\indexg{Entropy function|)}%
\end{equation}
\indexg{Entropy function!astral dual subgradients of|(}%
Nonetheless, let us consider a point $\tuu=\trans{[\tu,1-\tu]}$, for some
$\tu\in(0,1)$, which is on
$\Delta$'s boundary.
At such a point, even though $\psi(\tuu)\in\R$,
$\psi$ has no standard subgradient;
however, it does have an astral dual subgradient, as we will now see.

The boundary of $\Delta$ consists of three line segments
connecting the points $\trans{[0,0]}$, $\trans{[0,1]}$, and $\trans{[1,0]}$.
Let $L=\set{\trans{[u,1-u]}\!:\:u\in\R}$ be the line in $\R^2$
containing the segment $\lb{\trans{[0,1]}\!}{\trans{[1,0]}}$, as well
as the point $\tuu$.
Let $\rho:\R\rightarrow\Rext$ be the function evaluating $\psi$ along
the line $L$; that is, for $u\in\R$,
\begin{equation}  \label{eq:ex:entropy-ast-dual-subgrad:6}
  \rho(u)=\psi(u,1-u)=
  \begin{cases}
      u\ln u + (1-u)\ln(1-u)
      &\text{if $u\in[0,1]$,}
  \\
      +\infty
      &\text{otherwise,}
  \end{cases}
\end{equation}
as shown in \Cref{fig:entropy-subgrad} (right).
This function is differentiable at all points in $(0,1)$, the interior
of its domain, with derivative $\rho'(u)=\ln[u/(1-u)]$ for $u\in(0,1)$.
Let $\ww_1=\trans{[1,-1]}$, which is parallel to $L$, and let
$\ww_2=\trans{[1,1]}$, which is orthogonal to both $\ww_1$ and $L$.
Then $L=\set{\uu\in\Rn:\:\ww_2\inprod\uu=1}$.

We claim that $\psi$ has an astral dual subgradient at $\tuu$
equal to
\begin{equation}  \label{eq:ex:entropy-ast-dual-subgrad:1}
  \xbar
  =
  \limray{\ww_2}\plusl\rho'(\tu)\frac{\ww_1}{2}
  =
  \limray{\ww_2}\plusl\ln\biggParens{\frac{\tu}{1-\tu}}\frac{\ww_1}{2}.
\end{equation}
Intuitively, this astral point
$\xbar$ is capturing both $\psi$'s infinite rate of change in
the direction of $\ww_2$, as well as its
finite derivative
information in the direction of $\ww_1$, that is, along the line $L$.
Below, this connection between astral dual subgradients and the function's
rate of change in various directions will be made rigorous.

\begin{figure}[t]
  \centering
  \figentsubgrad{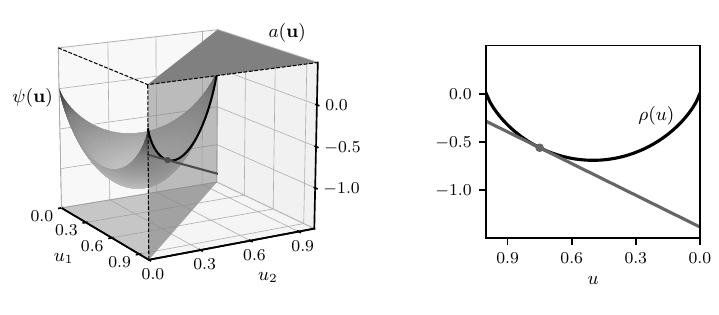}
  \mycaption{Astral dual subgradient of the entropy function}{%
\indexf{Entropy function!astral dual subgradients of}%
\indexf{subdifferentials, astral dual!examples}%
    On the left, we show the affine function $a$ from
    \Cref{ex:entropy-ast-dual-subgrad}
    (Eq.~\ref{eq:entropy-ast-dual-subgrad:1}),
    corresponding to the astral dual subgradient of $\psi$ at $\tuu=\trans{[\tu,1-\tu]}$ (gray dot).
    The function $a$ matches the value of $\psi$ at $\tuu$ and bounds $\psi$ from below
    everywhere else:
    It is equal to $-\infty$
    over the region of points $\uu$ with $u_1+u_2<1$, and to $+\infty$
    for points with $u_1+u_2>1$, where also $\psi(\uu)=+\infty$.
    For points with $u_1+u_2=1$ (comprising a line),
    the function $a$ is equal to the affine function corresponding to
    the subgradient of $\psi$'s restriction to that line.
    That restriction of $\psi$, denoted $\rho$, and the corresponding
    restriction of $a$ are depicted on the right.}
  \label{fig:entropy-subgrad}%
\end{figure}

To see that $\xbar\in\adsubdifpsi{\tuu}$, we need to show that,
for all $\uu'\in\R^2$,
$\psi(\uu')$ is bounded below by the function
\begin{align}
\notag
   a(\uu')
   &=\xbar\inprod(\uu'-\tuu)+\psi(\tuu)
\\
\notag
   &=\limray{\ww_2}\inprod(\uu'-\tuu) + \rho'(\tu)
   \frac{\ww_1\inprod(\uu'-\tuu)}{2} + \rho(\tu)
\\
\label{eq:entropy-ast-dual-subgrad:1}
  &=\omega(\ww_2\inprod\uu'-1)
     + \rho'(\tu)\frac{\ww_1\inprod(\uu'-\tuu)}{2}
     + \rho(\tu),
\end{align}
where the second equality is by $\xbar$'s definition, and the third is because
$\ww_2\inprod\tuu=1$. 
(We also are using that both $\psi(\tuu)=\rho(\tu)$ and $\rho'(\tu)$ are
finite.)
The function $a$ is depicted in \Cref{fig:entropy-subgrad} (left),
which also illustrates how $a$ bounds $\psi$ from below.
To prove this formally,
we distinguish three cases based on the value of $\ww_2\inprod\uu'$:
If $\ww_2\inprod\uu'<1$, then $a(\uu')=-\infty$, so $\psi(\uu')\ge a(\uu')$. If $\ww_2\inprod\uu'>1$, that is, $u'_1+u'_2>1$, then $\uu'\not\in\Delta$, so $\psi(\uu')=+\infty$ and the inequality $\psi(\uu')\ge a(\uu')$ holds as well.

The only remaining case is $\ww_2\inprod\uu'=1$, that is, $\uu'=[u',1-u']$ for some $u'\in\R$. Then
\begin{equation}
\label{eq:entropy-ast-dual-subgrad:2}
  \ww_1\inprod(\uu'-\tuu)
  =
  \trans{[1,-1]}\inprod\trans{[u'-\tu,\,-(u'-\tu)]}
  =
  2(u'-\tu),
\end{equation}
so
\begin{equation*}
   \psi(\uu')
   =
   \rho(u')
   \ge
   \rho'(\tu)(u'-\tu)+\rho(\tu)
   =
   \rho'(\tu)\frac{\ww_1\inprod(\uu'-\tuu)}{2}+\rho(\tu)
   =
   a(\uu').
\end{equation*}
The inequality is because $\rho'(\tu)$ is a subgradient of $\rho$ at $\tu$
(\Cref{roc:thm25.1}\ref{roc:thm25.1:a}). The second equality is by
\eqref{eq:entropy-ast-dual-subgrad:2}. And the final equality is by
\eqref{eq:entropy-ast-dual-subgrad:1}, since $\ww_2\inprod\uu'=1$. Thus, in all cases,
$\psi(\uu')\ge a(\uu')$, so $\xbar\in\adsubdifpsi{\tuu}$.
Below, in
\Cref{ex:entropy-ast-dual-subgrad-cont},
we show that $\xbar$ is 
the only astral dual subgradient of $\psi$ at $\tuu$.%
\indexg{subdifferentials, astral dual!examples|)}%
\indexg{Entropy function!astral dual subgradients of|)}%
\end{example}

\indexg{optimality conditions!astral|(}%
\indexg{subdifferentials, astral dual!0 included in@$\zero$ included in|(}%
For any function $\psi:\Rn\rightarrow\Rext$,
the definition of astral dual subgradient
immediately implies that
$\uu\in\Rn$ minimizes $\psi$ if and only if
$\zero$ is an astral dual subgradient of $\psi$ at
\indexg{subdifferentials, astral dual!0 included in@$\zero$ included in|)}%
\indexg{optimality conditions!astral|)}%
$\uu$.
\indexg{subdifferentials, standard!generalized by astral subgradients|(}%
\indexg{subdifferentials, astral dual!generalizing standard subgradients|(}%
Also,
as we show next,
astral dual subgradients generalize standard subgradients in the
sense that
$\psi$'s standard subdifferential at $\uu$
is exactly equal to the
finite points included in $\psi$'s astral dual subdifferential at $\uu$:

\begin{proposition}  \label{pr:adsubdif-int-rn}
  Let $\psi:\Rn\rightarrow\Rext$, and let $\uu\in\Rn$.
  Then $\partial\psi(\uu) = \adsubdifpsi{\uu}\cap\Rn$.
\end{proposition}

\begin{proof}
When $\xbar=\xx\in\Rn$,
the definitions of standard and astral dual
subgradients (Eqs.~\ref{eqn:psi-subgrad:1} and \ref{eqn:psi-subgrad:3-alt})
are equivalent, proving the claim.%
\indexg{subdifferentials, astral dual!generalizing standard subgradients|)}%
\indexg{subdifferentials, standard!generalized by astral subgradients|)}%
\end{proof}

\indexg{one-sided directional derivatives!astral dual subgradients and|(}%
\indexg{subdifferentials, astral dual!directional derivatives and|(}%
In
\Cref{sec:prelim:subgrads},
we defined
the one-sided directional derivative of a function as a way
of capturing the instantaneous rate at which the function is changing
in a given direction. We then presented several results that relate
directional derivatives to standard subgradients. We next extend
those results to astral dual subgradients.

As given
in \eqref{eq:direc-deriv-dfn}, the directional derivative
of a function $\psi:\Rn\rightarrow\Rext$ at a point $\uu\in\Rn$
with respect to a vector $\vv\in\Rn$,
assuming $\psi(\uu)\in\R$,
was defined to be
\begin{equation}   \label{eq:direc-deriv-dfn-rpt}
  \dderpsi{\uu}{\vv}=\lim_{\rightlim{\lambda}{0}}
  \frac{\psi(\uu+\lambda\vv)-\psi(\uu)}{\lambda}.
\end{equation}
We begin by generalizing
\Cref{roc:thm23.2}, which used directional derivatives
to characterize the standard subdifferentials of a convex
function at points where the function is finite.
We prove next that the same characterization
can be applied more generally to astral dual subdifferentials:

\begin{theorem}  \label{thm:dderiv-gives-adsubdif}
  Let $\psi:\Rn\rightarrow\Rext$ be convex, let $\uu\in\Rn$,
  and let $\xbar\in\extspace$.
  Assume $\psi(\uu)\in\R$.
  Then $\xbar\in\adsubdifpsi{\uu}$ if and only if
  $\xbar\cdot\vv\leq \dderpsi{\uu}{\vv}$ for all $\vv\in\Rn$.
\end{theorem}

\begin{proof}
  ~

\begin{proof-parts}
\pfpart{``Only if'' ($\Rightarrow$):}
Suppose $\xbar\in\adsubdifpsi{\uu}$, and let $\vv\in\Rn$.
Then the definition of astral dual subgradient (Eq.~\ref{eqn:psi-subgrad:3-alt})
yields, for all $\lambda\in\Rstrictpos$,
that
$\psi(\uu+\lambda\vv)\geq \psi(\uu) + \lambda \xbar\cdot\vv$,
and thus that
\[
  \frac{\psi(\uu+\lambda\vv) - \psi(\uu)}{\lambda}
  \geq
  \xbar\cdot\vv.
\]
By \Cref{roc:thm23.1}, this implies
$\dderpsi{\uu}{\vv}\geq \xbar\cdot\vv$.

\pfpart{``If'' ($\Leftarrow$):}
Suppose $\xbar\cdot\vv\leq \dderpsi{\uu}{\vv}$ for all $\vv\in\Rn$.
Let $\uu'\in\Rn$, and let $\vv=\uu'-\uu$.
Then
\[
  \psi(\uu') - \psi(\uu)
  =
  \psi(\uu+\vv) - \psi(\uu)
  \geq
  \dderpsi{\uu}{\vv}
  \geq
  \xbar\cdot\vv
  =
  \xbar\cdot(\uu'-\uu),
\]
where the first inequality is by
\Cref{roc:thm23.1}, and the second by assumption.
After rearranging, this is the same as
\eqref{eqn:psi-subgrad:3-alt},
proving $\xbar\in\adsubdifpsi{\uu}$.
\qedhere
\end{proof-parts}
\end{proof}

\Cref{thm:dderiv-gives-adsubdif} shows
how directional derivatives characterize astral dual subgradients.
In the opposite direction, we saw in \Cref{roc:thm23.4-dd}
that standard subgradients characterize
the directional derivatives at a point, as long as that point is in the
relative interior of the function's effective domain
(and the function is convex and proper).
As we show next,
for astral dual subgradients, this holds more generally for all points
where the function is finite (and without requiring properness).

For instance, in \Cref{ex:entropy-ast-dual-subgrad}, we
noted informally that the function $\psi$'s
astral dual subgradient $\xbar$
at the point $\tuu$ (Eq.~\ref{eq:ex:entropy-ast-dual-subgrad:1})
captures both $\psi$'s infinite slope in the direction of $\ww_2$,
as well as the function's finite
slope in the perpendicular direction of $\ww_1$.
The next \namecref{thm:adsub-gives-dderiv} makes this precise,
showing generally how astral dual subdifferentials capture
the instantaneous rate of increase in
all directions at all points where a convex function is finite.

\begin{theorem}   \label{thm:adsub-gives-dderiv}
  Let $\psi:\Rn\rightarrow\Rext$ be convex, and let $\uu\in\Rn$.
  Assume $\psi(\uu)\in\R$.
  Then for all $\vv\in\Rn$,
  \begin{equation}   \label{eq:thm:adsub-gives-dderiv:1}
     \dderpsi{\uu}{\vv}     
     =
     \indfastar{\adsubdifpsi{\uu}}(\vv)
     =
     \sup\bigBraces{
       \xbar\cdot\vv
       :\:
       \xbar\in\adsubdifpsi{\uu}
     }.
  \end{equation}
\end{theorem}

\begin{proof}
Let $\rho:\Rn\rightarrow\Rext$ be defined by
$\rho(\vv)=\dderpsi{\uu}{\vv}$ for $\vv\in\Rn$.
Then, by \Cref{roc:thm23.1},
$\rho$ is convex and positively homogeneous
with $\rho(\zero)=0$.
Applying
\Cref{thm:pos-homo-is-ast-closed}(\ref{thm:pos-homo-is-ast-closed:b})
to the function $\rho$ then yields $\rho=\indaSstar$
with
\begin{equation}   \label{eq:thm:adsub-gives-dderiv:2}
  S
  =
  \Braces{
    \xbar\in\extspace
    :\:
    \xbar\cdot\vv \leq \rho(\vv)
    \text{ for all } \vv\in\Rn
  }
  =
  \adsubdifpsi{\uu},
\end{equation}
where the second equality is by
\Cref{thm:dderiv-gives-adsubdif}.
Thus, for $\vv\in\Rn$, we have
\[
  \dderpsi{\uu}{\vv}
  =
  \rho(\vv)
  =
  \indaSstar(\vv)
  =
  \indfastar{\adsubdifpsi{\uu}}(\vv),
\]
proving the first equality of
\eqref{eq:thm:adsub-gives-dderiv:1}.
(The second equality is by
Eq.~\ref{eqn:astral-support-fcn-def}.)
\end{proof}

In fact, the use of
\Cref{thm:pos-homo-is-ast-closed}(\ref{thm:pos-homo-is-ast-closed:b})
in the proof above
further implies
that $\adsubdifpsi{\uu}$ is nonempty, closed, and convex,
assuming
that $\psi(\uu)$ is finite.
Later, 
we will prove that $\adsubdifpsi{\uu}$ has these same properties even if
$\psi(\uu)$ is infinite
(see \Cref{thm:adsubdiff-nonempty,thm:ast-dual-subdif-is-convex}).

\indexg{subdifferentials, astral dual!singleton|(}%
In standard convex analysis, the subdifferential of a convex function
at a point is a singleton if and only if the function is
differentiable at that point, in which case, its only subgradient
is the function's gradient (\Cref{roc:thm25.1}).
The case that the astral dual subdifferential
$\adsubdifpsi{\uu}$ is a singleton is therefore also of particular
interest as a generalization of standard differentiability.
The next \namecref{thm:adsub-singleton}
uses directional derivatives to characterize when
this is so.

Note that condition~(\ref{thm:adsub-singleton:b}) of the
\namecref{thm:adsub-singleton} is equivalent to requiring
that the {two}-sided limit
\[
  \lim_{\lambda\rightarrow 0}
  \frac{\psi(\uu+\lambda\vv)-\psi(\uu)}{\lambda}
\]
must exist for all $\vv\in\Rn$, and so be equal to
$\dderpsi{\uu}{\vv}$,
as can be argued from the definition of directional derivative
(Eq.~\ref{eq:direc-deriv-dfn-rpt}).

\begin{theorem}  \label{thm:adsub-singleton}
  Let $\psi:\Rn\rightarrow\Rext$ be convex, and let $\uu\in\Rn$.
  Assume $\psi(\uu)\in\R$.
  Then the following are equivalent:
  \begin{letter-compact}
  \item    \label{thm:adsub-singleton:a}
    $\adsubdifpsi{\uu}$ is a singleton.
  \item    \label{thm:adsub-singleton:b}
    $\dderpsi{\uu}{\vv} = -\dderpsi{\uu}{-\vv}$
    for all $\vv\in\Rn$.
  \item    \label{thm:adsub-singleton:c}
    There exists $\xbar\in\extspace$ such that
    $\dderpsi{\uu}{\vv} = \xbar\cdot\vv$
    for all $\vv\in\Rn$.
  \end{letter-compact}
  Moreover, if any (and therefore all) of these conditions hold, then
  the single element of $\adsubdifpsi{\uu}$ in
  condition~(\ref{thm:adsub-singleton:a}) is also the unique point
  $\xbar$ satisfying
  condition~(\ref{thm:adsub-singleton:c}).
\end{theorem}

\begin{proof}
  ~
  
\begin{proof-parts}
\pfpart{%
  (\ref{thm:adsub-singleton:a})
  $\Rightarrow$
  (\ref{thm:adsub-singleton:c}):
}
If $\adsubdifpsi{\uu}=\{\xbar\}$ for some $\xbar\in\extspace$,
then \Cref{thm:adsub-gives-dderiv} immediately implies that
$\dderpsi{\uu}{\vv}=\xbar\cdot\vv$ for all $\vv\in\Rn$.

\pfpart{%
  (\ref{thm:adsub-singleton:c})
  $\Rightarrow$
  (\ref{thm:adsub-singleton:b}):
}
Suppose, for some $\xbar\in\extspace$, that
$\dderpsi{\uu}{\vv}=\xbar\cdot\vv$ for all $\vv\in\Rn$.
Then for $\vv\in\Rn$, we have
\[
   -\dderpsi{\uu}{-\vv}
   =
   -\xbar\cdot(-\vv)
   =
   \xbar\cdot\vv
   =
   \dderpsi{\uu}{\vv}.
\]

\pfpart{%
  (\ref{thm:adsub-singleton:b})
  $\Rightarrow$
  (\ref{thm:adsub-singleton:a}):
}
Suppose condition~(\ref{thm:adsub-singleton:b}) holds.
Note that $\adsubdifpsi{\uu}$ cannot be empty since otherwise
\Cref{thm:adsub-gives-dderiv} would imply
$\dderpsi{\uu}{\zero}=-\infty$, contradicting
that
$\dderpsi{\uu}{\zero}=0$
(\Cref{roc:thm23.1}).

Suppose both $\xbar$ and $\xbar'$ are elements of $\adsubdifpsi{\uu}$.
Then for all $\vv\in\Rn$,
\[
  \xbar\cdot\vv
  \leq
  \dderpsi{\uu}{\vv}
  =
  - \dderpsi{\uu}{-\vv}
  \leq
  - \xbar'\cdot(-\vv)
  =
  \xbar'\cdot\vv,
\]
where both inequalities are by
\Cref{thm:dderiv-gives-adsubdif}.
Applied
to $-\vv$,
this also shows that $\xbar\cdot\vv\geq\xbar'\cdot\vv$.
Thus, $\xbar\cdot\vv=\xbar'\cdot\vv$ for all $\vv\in\Rn$, so
$\xbar=\xbar'$ (\Cref{pr:i:4}).
Therefore, $\adsubdifpsi{\uu}$ is a singleton.

\pfpart{Final claim:}
Suppose any, and therefore all, of the conditions hold,
implying $\adsubdifpsi{\uu}=\{\xbar\}$ for some $\xbar\in\extspace$.
Then, as shown above, $\xbar$ satisfies
condition~(\ref{thm:adsub-singleton:c}).
If some point $\xbar'\in\extspace$ also satisfies this condition,
then this would imply that
$\xbar'\cdot\vv=\dderpsi{\uu}{\vv}=\xbar\cdot\vv$
for all $\vv\in\Rn$, and so that
$\xbar'=\xbar$.
Thus, $\xbar$ is the unique point satisfying
condition~(\ref{thm:adsub-singleton:c}).%
\indexg{subdifferentials, astral dual!singleton|)}%
\qedhere
\end{proof-parts}
\end{proof}

\indexg{Entropy function!astral dual subgradients of|(}%
\indexg{Entropy function!directional derivatives of|(}%
As an illustration, we apply these tools to the entropy function from
\Cref{ex:entropy-ast-dual-subgrad}:

\begin{example}
  [Entropy function continued, with directional derivatives]
  \label{ex:entropy-ast-dual-subgrad-cont}
We continue \Cref{ex:entropy-ast-dual-subgrad},
with
all notation
as defined earlier in that example.
We saw already that $\xbar$ is an astral dual subgradient of $\psi$ at
$\tuu=\trans{[\tu,1-\tu]}$, where $\tu\in(0,1)$.
Here, we show that $\xbar$ is the only such subgradient, while also
deriving all of the directional derivatives at $\tuu$.
As before, $\ww_1=\trans{[1,-1]}$ and $\ww_2=\trans{[1,1]}$. Hence, $\ww_2\inprod\tuu=1$, and
all points $\uu\in\Delta$ satisfy $\ww_2\inprod\uu\le 1$.

Let
$h:\Rpos\to\R$ be the function defined by $h(u)=u\ln u$,
which is differentiable, with derivative $h'$, for $u\in\Rstrictpos$.
Then
$\psi(\uu)=h(u_1)+h(u_2)+h(1-u_1-u_2)$ for all $\uu\in\Delta$.
Let $\vv\in\Rn$.
We seek the directional derivative
\begin{equation}   \label{eq:entropy-direct-der:1}
  \dderpsi{\tuu}{\vv}
  =
  \lim_{\rightlim{\lambda}{0}}\frac{\psi(\tuu+\lambda\vv)-\psi(\tuu)}{\lambda}.
\end{equation}

If $\ww_2\inprod\vv>0$, then for all $\lambda\in\Rstrictpos$,
$\ww_2\inprod(\tuu+\lambda\vv)>1$, implying
$\tuu+\lambda\vv\not\in\Delta$,
so $\psi(\tuu+\lambda\vv)=+\infty$;
hence, the limit in \eqref{eq:entropy-direct-der:1} is $+\infty$.

Otherwise, $\ww_2\inprod\vv\le 0$.
The components of $\tuu$ are $\tu_1=\tu$ and $\tu_2=1-\tu$, which are
both in $(0,1)$.
Therefore, for sufficiently small $\lambda\in\Rstrictpos$,
$\tuu+\lambda\vv\in\Delta$,
in which case,
\begin{align*}
  \psi(\tuu+\lambda\vv)
  &=
    h(\tu_1+\lambda v_1)
    +h(\tu_2+\lambda v_2)
    +h(1-\tu_1-\tu_2-\lambda v_1-\lambda v_2)
\\ 
  &=
    h(\tu_1+\lambda v_1)
    +h(\tu_2+\lambda v_2)
    +h(-\lambda v_1-\lambda v_2),
\\
  \psi(\tuu)
  &=
    h(\tu_1)
    +h(\tu_2)
    +h(1-\tu_1-\tu_2)
  =
    h(\tu_1)
    +h(\tu_2),
\end{align*}
using $\tu_1+\tu_2=1$.
Plugging these into \eqref{eq:entropy-direct-der:1} then yields
\begin{align*}
  &
  \dderpsi{\tuu}{\vv}
\\
  &\qquad{}
=
\lim_{\rightlim{\lambda}{0}}
\smash[t]{\biggBracks{
  \frac{h(\tu_1+\lambda v_1)-h(\tu_1)}{\lambda}
  + \frac{h(\tu_2+\lambda v_2)-h(\tu_2)}{\lambda}
  + \frac{h(-\lambda v_1-\lambda v_2)}{\lambda}
}}
\\
&\qquad{}
=
v_1 h'(\tu_1)+v_2 h'(\tu_2)+
\lim_{\rightlim{\lambda}{0}}
  (-v_1-v_2)\ln(-\lambda v_1-\lambda v_2)
,
\end{align*}
where the second equality follows
because $h$ is differentiable at $\tu_1,\tu_2\in(0,1)$
(and by continuity of addition).
If $v_1+v_2<0$, then the final limit is $-\infty$, and
if $v_1+v_2=0$, then this limit is~$0$.

Summarizing, we have shown that
\begin{equation}
\label{eq:entropy-direct-der:2}
  \dderpsi{\tuu}{\vv}=
  \begin{cases}
    +\infty & \text{if $\ww_2\inprod\vv>0$,}
  \\
    v_1 h'(\tu_1)+v_2 h'(\tu_2) & \text{if $\ww_2\inprod\vv=0$,}
  \\
    -\infty & \text{if $\ww_2\inprod\vv<0$.}
  \end{cases}
\end{equation}
In particular, this implies that
$\dderpsi{\tuu}{\vv}=-\dderpsi{\tuu}{-\vv}$ for all $\vv\in\Rn$.
Therefore, by
\Crefequiv{thm:adsub-singleton}{thm:adsub-singleton:b}{thm:adsub-singleton:a},
$\adsubdifpsi{\tuu}$ is a singleton, whose only element must be
\[
  \xbar=\limray{\ww_2}\plusl\frac{\rho'(\tu)}{2}\ww_1,
\]
as given in \eqref{eq:ex:entropy-ast-dual-subgrad:1},
which was shown in \Cref{ex:entropy-ast-dual-subgrad} to be
in $\adsubdifpsi{\tuu}$.
Thus, $\adsubdifpsi{\tuu}=\{\xbar\}$, and
$\dderpsi{\tuu}{\vv}=\xbar\cdot\vv$ for all $\vv\in\Rn$,
as follows from \Cref{thm:adsub-gives-dderiv}
(or \Cref{thm:adsub-singleton}).%
\indexg{subdifferentials, astral dual|)}%
\indexg{one-sided directional derivatives!astral dual subgradients and|)}%
\indexg{subdifferentials, astral dual!directional derivatives and|)}%
\indexg{Entropy function!astral dual subgradients of|)}%
\indexg{Entropy function!directional derivatives of|)}%
\end{example}

\chapter{%
  Characterizations
  of astral subgradients
}
\label{chp:char-astral-subgrads}

\indexg{Fenchel-Young inequality|(}%
In standard convex analysis, the Fenchel-Young inequality states that,
for any
function $f:\Rn\rightarrow\Rext$,
\begin{equation}  \label{eqn:fenchel-stand}
  \fstar(\uu) \geq -f(\xx) +\xx\cdot\uu
\end{equation}
for all $\xx\in\Rn$ and all
\indexg{Fenchel-Young inequality|)}%
$\uu\in\Rn$.
Furthermore,
this holds with equality if and only if
${\uu\in\partial f(\xx)}$
(\Cref{pr:stan-subgrad-equiv-props}\ref{pr:stan-subgrad-equiv-props:a}\ref{pr:stan-subgrad-equiv-props:b}),
providing a very useful characterization of subgradients.
It is also known that if $f$ is closed, proper and convex, then
$\partial f$ and $\partial \fstar$ act as inverses of one another in
the sense that,
for $\xx\in\Rn$ and $\uu\in\Rn$,
$\uu\in\partial f(\xx)$ if and only if
$\xx\in\partial \fstar(\uu)$
(\Cref{pr:stan-subgrad-equiv-props}\ref{pr:stan-subgrad-equiv-props:a}\ref{pr:stan-subgrad-equiv-props:c}).
Thus, if $f$ is closed, proper and convex, then the following are equivalent:
\begin{letter-compact}
\item
  $\fstar(\uu) = -f(\xx)+\xx\cdot\uu$.
\item
  $\uu\in\partial f(\xx)$.
\item
  $\xx\in\partial \fstar(\uu)$.
\end{letter-compact}
This equivalence holds under
slightly weaker conditions as well:
Instead of assuming that $f$ is closed, meaning $\fdubs=f$, it suffices to assume that
$f$ is just closed at~$\xx$, meaning $\fdubs(\xx)=f(\xx)$,
while still assuming that $f$ is convex and proper
(see Propositions~\ref{pr:stan-subgrad-equiv-props}
and~\ref{pr:conj-props-cvx}\ref{pr:conj-props-cvx:b}).

We turn now to exploring in detail analogous connections in the astral
setting, showing how astral subdifferentials can be characterized
using astral versions of the Fenchel-Young inequality, and how astral
primal and dual subdifferentials can act as inverses of each other.

\section{Conditions for astral subgradients}

\indexg{Fenchel-Young inequality|(}%
To begin, we note that
the Fenchel-Young inequality generalizes directly to functions
$F:\extspace\rightarrow\Rext$ over astral space
since the definition of astral conjugate in \eqref{eq:Fstar-down-def}
immediately implies that
\begin{equation}  \label{eqn:ast-fenchel}
  \Fstar(\uu)
  \geq
  - F(\xbar) \plusd \xbar\cdot\uu
\end{equation}
for all $\xbar\in\extspace$ and all $\uu\in\Rn$.
Moreover,
by \Cref{pr:plusd-props}(\ref{pr:plusd-props:e}), this is equivalent
to the alternative form
\begin{equation}  \label{eqn:ast-fenchel-alt}
  F(\xbar)
  \geq
  - \Fstar(\uu) \plusd \xbar\cdot\uu,
\end{equation}
and, by \Cref{pr:plusd-plusu-props}(\ref{pr:plusd-plusu-props:b}),
is also equivalent to
\begin{equation}  \label{eqn:ast-fenchel-alt-plusu}
  F(\xbar)
  \plusu
  \Fstar(\uu)
  \geq
  \xbar\cdot\uu.%
\indexg{Fenchel-Young inequality|)}%
\end{equation}

\indexg{Fenchel-Young inequality!astral subgradients and|(}%
\indexg{subdifferentials, astral (primal)!Fenchel-Young and|(}%
\indexg{subdifferentials, astral (primal)!dual subgradients of conjugate and|(}%
\indexg{subdifferentials, astral dual!conjugate@of conjugate|(}%
Analogous to the conditions above for the standard setting,
the next theorem gives a chain of implications connecting
astral versions of the Fenchel-Young inequality
with a function's astral (primal) subdifferentials
and its conjugate's astral dual subdifferentials.
More specifically,
assuming $\Fstar(\uu)\in\R$,
the theorem shows that if the version of the Fenchel-Young inequality
in \eqref{eqn:ast-fenchel} holds with equality, then $\uu$ must be an
astral subgradient of $F$ at $\xbar$.
This in turn implies that $F$ is astral closed at $\xbar$ (as we already
showed), and that $\xbar$ is an astral dual subgradient of
$\Fstar$ at $\uu$, recovering part of the inverse relationship
discussed above for standard subgradients.
These conditions further imply
that the alternative form of the Fenchel-Young inequality
in \eqref{eqn:ast-fenchel-alt}
(and so also the one in Eq.~\ref{eqn:ast-fenchel-alt-plusu})
must hold with equality.

Finally, if
$-F(\xbar)$ and $\xbar\cdot\uu$ are summable,
then all four statements are equivalent.

\begin{theorem}  \label{thm:fenchel-subgrad}
  Let $F:\extspace\rightarrow\Rext$,
  $\xbar\in\extspace$, and $\uu\in\Rn$.
  Assume $\Fstar(\uu)\in\R$.
  Consider the following statements:
  \begin{letter-compact}
  \item  \label{thm:fenchel-subgrad:a}
    $\Fstar(\uu) = -F(\xbar) \plusd \xbar\cdot\uu$.
  \item  \label{thm:fenchel-subgrad:b}
    $\uu\in\asubdifF{\xbar}$.
  \item  \label{thm:fenchel-subgrad:b-dual}
    $\xbar\in\adsubdifFstar{\uu}$
    and
    $\Fdub(\xbar)=F(\xbar)$.
  \item  \label{thm:fenchel-subgrad:c}
    $F(\xbar) = -\Fstar(\uu) + \xbar\cdot\uu$.
  \end{letter-compact}
  Then
  (\ref{thm:fenchel-subgrad:a})
  $\Rightarrow$
  (\ref{thm:fenchel-subgrad:b}),
  and
  (\ref{thm:fenchel-subgrad:b})
  $\Rightarrow$
  (\ref{thm:fenchel-subgrad:b-dual}),
  and
  (\ref{thm:fenchel-subgrad:b-dual})
  $\Rightarrow$
  (\ref{thm:fenchel-subgrad:c}).
  
  Furthermore, if
  $-F(\xbar)$ and $\xbar\cdot\uu$ are summable, then
  the four statements are equivalent.
\end{theorem}

Before proving the theorem, we give a lemma that effectively proves
the implication
(\ref{thm:fenchel-subgrad:b-dual})~$\Rightarrow$~(\ref{thm:fenchel-subgrad:c})
in a somewhat more general form.
This lemma will be used again shortly in proving
\Cref{thm:psi-subgrad-conds}.

\begin{lemma}   \label{lem:adsub-implies-fenchel}
\indexg{subdifferentials, astral dual!Fenchel-Young and|(}%
  Let $\psi:\Rn\rightarrow\Rext$, $\uu\in\Rn$, and
  $\xbar\in\extspace$.
  Assume $\psi(\uu)\in\R$ and that $\xbar\in\adsubdifpsi{\uu}$.
  Then $\psistarb(\xbar)=\xbar\cdot\uu - \psi(\uu)$.
\end{lemma}

\begin{proof}
That $\psistarb(\xbar)\geq\xbar\cdot\uu - \psi(\uu)$
follows from the definition of dual conjugate
(Eq.~\ref{eq:psistar-def:2}).
For the reverse inequality, let $\ww\in\Rn$.  
We claim that
\begin{equation}
\label{eq:lem:adsub-implies-fenchel:1}
  {-\psi(\ww)}\plusd\xbar\cdot\ww
  \leq
  \xbar\cdot\uu - \psi(\uu).
\end{equation}
This is immediate if $\psi(\ww)=+\infty$ or $\xbar\cdot\ww=-\infty$ or
$\xbar\cdot\uu=+\infty$.
Therefore, we assume henceforth that none of these conditions hold.

We then have that
\begin{align*}
  -\psi(\uu)
  &\geq
  -\psi(\ww)\plusd\xbar\cdot(\ww-\uu)
  \\
  &=
  -\psi(\ww)\plusd(\xbar\cdot\ww - \xbar\cdot\uu)
  \\
  &=
  \bigParens{-\psi(\ww) + \xbar\cdot\ww} - \xbar\cdot\uu.
\end{align*}
The inequality is by definition of astral dual
subgradient
(Eq.~\ref{eqn:psi-subgrad:3})
since $\xbar\in\adsubdifpsi{\uu}$.
The first equality is by summability of $\xbar\inprod\ww$ and $-\xbar\inprod\uu$ (using \Cref{pr:i:1}) since
$\xbar\cdot\ww>-\infty$ and $-\xbar\cdot\uu>-\infty$.
Since also $-\psi(\ww)>-\infty$, this also yields the second equality.
By \Cref{pr:plusd-props}(\ref{pr:plusd-props:e}), this implies
$\psi(\ww)-\xbar\cdot\ww\geq\psi(\uu)-\xbar\cdot\uu$, which is
equivalent to 
\eqref{eq:lem:adsub-implies-fenchel:1}.

Thus, \eqref{eq:lem:adsub-implies-fenchel:1} holds for all
$\ww\in\Rn$.
Taking supremum of the left-hand side, which is equal
to $\psistarb(\xbar)$ (by definition of astral dual conjugate in Eq.~\ref{eq:psistar-def:2}),
it follows that
$\psistarb(\xbar)\leq\xbar\cdot\uu - \psi(\uu)$,
completing the proof.%
\indexg{subdifferentials, astral dual!Fenchel-Young and|)}%
\end{proof}

\begin{proof}[Proof of \Cref{thm:fenchel-subgrad}]
~

\begin{proof-parts}
\pfpart{%
  (\ref{thm:fenchel-subgrad:a})
  $\Rightarrow$
  (\ref{thm:fenchel-subgrad:b}):
}
Suppose statement~(\ref{thm:fenchel-subgrad:a}) holds.
Then the sequence $\xbar_t=\xbar$ satisfies the conditions of
\Cref{pr:equiv-ast-subdif-defn}(\ref{pr:equiv-ast-subdif-defn:plusd}),
so $\uu\in\asubdifF{\xbar}$ by the equivalence of
parts~(\ref{pr:equiv-ast-subdif-defn:plusd})
and~(\ref{pr:equiv-ast-subdif-defn:a}) of
that \namecref{pr:equiv-ast-subdif-defn}.

\pfpart{%
  (\ref{thm:fenchel-subgrad:b})
  $\Rightarrow$
  (\ref{thm:fenchel-subgrad:b-dual}):
}
Suppose $\uu\in\asubdifF{\xbar}$.
Then $F(\xbar)=\Fdub(\xbar)$ by
\Cref{thm:subgrad-then-lsc}.
By
\Crefequiv{pr:equiv-ast-subdif-defn}{pr:equiv-ast-subdif-defn:a}{pr:equiv-ast-subdif-defn:b},
there exists a sequence $\seq{\xbar_t}$ in $\eRn$
with $F(\xbar_t)\in\R$ and $\xbar_t\inprod\uu\in\R$
for all $t$, and such that
$\xbar_t \to \xbar$,
$F(\xbar_t) \to F(\xbar)$,
and
$\xbar_t\cdot\uu - F(\xbar_t) \to \Fstar(\uu)$.

Let $\ww\in\Rn$.
Then for all $t$,
\[
  \Fstar(\ww)
  \geq
  \xbar_t\cdot\ww - F(\xbar_t)
  =
  \xbar_t\cdot\uu + \xbar_t\cdot(\ww-\uu) - F(\xbar_t).
\]
The inequality is by definition of conjugate
(Eq.~\ref{eq:Fstar-down-def}),
and the equality is by \Cref{pr:i:1}
since $\xbar_t\cdot\uu\in\R$.
Therefore, taking limits yields
\begin{align*}
  \Fstar(\ww)
  &\geq
  \lim\BigParens{\bigBracks{\xbar_t\cdot\uu - F(\xbar_t)}
  +              \xbar_t\cdot(\ww-\uu)}
  =
  \Fstar(\uu) + \xbar\cdot(\ww-\uu).
\end{align*}
The equality is because
$\xbar_t\cdot\uu-F(\xbar_t)\rightarrow\Fstar(\uu)$ and
$\xbar_t\cdot(\ww-\uu)\rightarrow\xbar\cdot(\ww-\uu)$
(by \Cref{thm:i:1}\ref{thm:i:1c}),
and also using that $\Fstar(\uu)\in\R$, so that the
limit of the sum is equal to the sum of the limits
(\Cref{prop:lim:eR}\ref{i:lim:eR:sum}).
Hence, $\xbar\in\adsubdifFstar{\uu}$ by
definition of astral dual subgradient
(Eq.~\ref{eqn:psi-subgrad:3-alt}).

\pfpart{%
  (\ref{thm:fenchel-subgrad:b-dual})
  $\Rightarrow$
  (\ref{thm:fenchel-subgrad:c}):
}
If
$\xbar\in\adsubdifFstar{\uu}$
and
$\Fdub(\xbar)=F(\xbar)$
then
\[
  F(\xbar)
  =
  \Fdub(\xbar)
  =
  \xbar\cdot\uu-\Fstar(\uu),
\]
where the second equality is by
\Cref{lem:adsub-implies-fenchel}
(with $\psi=\Fstar$).

\pfpart{Under summability,
(\ref{thm:fenchel-subgrad:c})
   $\Rightarrow$ (\ref{thm:fenchel-subgrad:a}):}   
Suppose $-F(\xbar)$ and $\xbar\cdot\uu$ are summable
and
that statement~(\ref{thm:fenchel-subgrad:c}) holds.
Then $\xbar\cdot\uu\in\R$, since otherwise, if
$\xbar\cdot\uu\in\{-\infty,+\infty\}$, then
statement~(\ref{thm:fenchel-subgrad:c}) would imply
$F(\xbar)=\xbar\cdot\uu$, violating the assumed summability of
$-F(\xbar)$ and $\xbar\cdot\uu$.
Thus, $F(\xbar)$, $\xbar\cdot\uu$ and $\Fstar(\uu)$ are all in $\R$,
so (\ref{thm:fenchel-subgrad:c})
implies (\ref{thm:fenchel-subgrad:a}).
\qedhere
\end{proof-parts}
\end{proof}

Without the summability assumption in the last part of
\Cref{thm:fenchel-subgrad},
the four statements appearing in the \namecref{thm:fenchel-subgrad}
are increasingly more permissive,
allowing more pairs of $\xbar$ and $\uu$ to satisfy them.
In other words, as we show in the next example,
none of the \namecref{thm:fenchel-subgrad}'s
statements generally imply the preceding
statement, even if we also assume that $F$ is convex:

\begin{example}  \label{ex:fenchel-subgrad-thm-strict}
\indexg{absolute value function, subgradients of|(}%
Let $F:\Rext\rightarrow\Rext$
be defined by $F(\barx)=\abs{\barx}$ for $\barx\in\eR$.
Then $\Fstar(u)=\indf{[-1,1]}(u)$ for $u\in\R$, and $\Fdub=F$
(by \Cref{thm:Fdub:1d}).
For $\barx=+\infty$ and $u=1$, we have $u\in \asubdifF{\barx}$
(see \Cref{ex:standard-abs-val-subgrad-cont}),
so statement~(\ref{thm:fenchel-subgrad:b})
of \Cref{thm:fenchel-subgrad}
holds, but~(\ref{thm:fenchel-subgrad:a})
is false.
And for $\barx=+\infty$ and $u=1/2$, we have $F(\barx) = -\Fstar(u) + \barx\cdot u = +\infty$,
so~(\ref{thm:fenchel-subgrad:c}) holds, but~(\ref{thm:fenchel-subgrad:b-dual})
is false since $\barx\not\in\adsubdifFstar{u}=\set{0}$.
We defer to \Cref{ex:cvx-not-subgrad-inv} a more involved construction
in which~(\ref{thm:fenchel-subgrad:b-dual}) holds
but~(\ref{thm:fenchel-subgrad:b}) does not.%
\indexg{Fenchel-Young inequality!astral subgradients and|)}%
\indexg{subdifferentials, astral (primal)!Fenchel-Young and|)}%
\indexg{subdifferentials, astral (primal)!dual subgradients of conjugate and|)}%
\indexg{subdifferentials, astral dual!conjugate@of conjugate|)}%
\indexg{absolute value function, subgradients of|)}%
\end{example}

\indexg{Absolute value at infinity!astral subgradients of extension|(}%
With \Cref{thm:fenchel-subgrad} at hand, we can now complete our discussion of
\Cref{ex:subgrad-sqrt-approaches-abs-val}:

\begin{example}[Absolute value at infinity, continued]
   \label{ex:subgrad-sqrt-approaches-abs-val:cont}
Continuing \Cref{ex:subgrad-sqrt-approaches-abs-val},
we consider $f(\xx)=\sqrt{x_1^2 + e^{-x_2}}$ for ${\xx\in\R^2}$,
and let $\xbar_\alpha=\limray{\ee_2} \plusl \alpha \ee_1$
for $\alpha\in\R$, so $\fext(\xbar_\alpha)=|\alpha|$ for all $\alpha\in\R$.
We aim to prove that
\begin{equation}
\label{eq:ex:subgrad-sqrt-approaches-abs-val:1:again}
\asubdiffext{\xbar_\alpha}
=
\begin{cases}
  \{-\ee_1\}    & \text{if $\alpha < 0$,} \\
  \Braces{\lambda\ee_1 :\: \lambda\in[-1,1]}
                & \text{if $\alpha = 0$,} \\
  \{\ee_1\}     & \text{if $\alpha > 0$.}
\end{cases}
\end{equation}

First, suppose $\uu\in\asubdiffext{\xbar_\alpha}$, for some
$\alpha\in\R$.
Then
$\fstar(\uu)\in\R$
by \Cref{pr:subgrad-imp-in-cldom}(\ref{pr:subgrad-imp-in-cldom:c}),
and
$\fext(\xbar_\alpha)=\xbar_\alpha\cdot\uu - \fstar(\uu)$
by
\Crefequiv{thm:fenchel-subgrad}{thm:fenchel-subgrad:b}{thm:fenchel-subgrad:c}.
Since $\fext(\xbar_\alpha)\in\R$, this implies that
$\xbar_\alpha\cdot\uu$ is in $\R$, and therefore that
$\uu\cdot\ee_2=0$ by \Cref{pr:vtransu-zero}.
Consequently, $\asubdiffext{\xbar_\alpha}$ can only consist of points
of the form $\lambda\ee_1$ for some $\lambda\in\R$.

Using the definition of the standard conjugate
(Eq.~\ref{eq:fstar-def}), we can compute, for $\lambda\in\R$, that
$\fstar(\lambda\ee_1)=\indf{[-1,1]}(\lambda)$
(that is, $0$ if $\lambda\in[-1,1]$ and $+\infty$ otherwise).
Thus, if $\lambda\not\in[-1,1]$ then $\lambda\ee_1$ cannot be an
astral subgradient of $\fext$ at any point
(again by
\Cref{pr:subgrad-imp-in-cldom}\ref{pr:subgrad-imp-in-cldom:c}),
so $\uu=\lambda\ee_1$ for some $\lambda\in[-1,1]$.

Since $\fext(\xbar_\alpha)\in\R$,
\Cref{thm:fenchel-subgrad}
(applied to $\fext$, $\xbar_\alpha$, and $\uu$)
holds with equivalence among its various statements.
Thus, by \Cref{thm:fenchel-subgrad}(\ref{thm:fenchel-subgrad:b},\ref{thm:fenchel-subgrad:a}),
we obtain,
for $\lambda\in[-1,1]$, that
$\lambda\ee_1\in\asubdiffext{\xbar_\alpha}$ if and only if
\[
  0
  =
  \fstar(\lambda\ee_1)
  =
  \xbar_\alpha\cdot(\lambda\ee_1) - \fext(\xbar_\alpha)
  =
  \alpha\lambda - |\alpha|.
\]
This equality holds if and only if
one of the following holds:
$\alpha=0$;
$\alpha>0$ and $\lambda=1$;
$\alpha<0$ and $\lambda=-1$.
This completes the proof of
\eqref{eq:ex:subgrad-sqrt-approaches-abs-val:1:again}.%
\indexg{Absolute value at infinity!astral subgradients of extension|)}%
\end{example}

The summability assumption in the last part of
\Cref{thm:fenchel-subgrad}
always holds when
$F(\xbar)\in\R$ or
$\xbar\cdot\uu\in\R$, including when $\xbar=\xx\in\R$.
\indexg{subdifferentials, astral (primal)!0 included in@$\zero$ included in|(}%
\indexg{optimality conditions!astral|(}%
In the important case that $\uu=\zero$, the \namecref{thm:fenchel-subgrad}
thus yields that
$\zero$ is a subgradient of $F$ at $\xbar$ if and only if
$\xbar$ minimizes $F$ (provided $F$ has finite infimum):

\begin{proposition}  \label{pr:asub-zero-is-min}
  Let $F:\extspace\rightarrow\Rext$, and
  let $\xbar\in\extspace$.
  Then:
  \begin{letter-compact}
  \item  \label{pr:asub-zero-is-min:a}
    $\zero\in\asubdifF{\xbar}$
    if and only if
    $\inf F\in\R$ and $\xbar$ minimizes $F$.
  \item  \label{pr:asub-zero-is-min:b}
    Consequently,
    $\xbar$ minimizes $F$ if and only if one of the following holds:
    \begin{roman-compact}
    \item  \label{pr:asub-zero-is-min:b:1}
      $\zero\in\asubdifF{\xbar}$;
    \item  \label{pr:asub-zero-is-min:b:2}
      $F(\xbar)=-\infty$; or
    \item  \label{pr:asub-zero-is-min:b:3}
      $F\equiv+\infty$.
    \end{roman-compact}
  \end{letter-compact}
\end{proposition}

\begin{proof}
~

\begin{proof-parts}
\pfpart{Part~(\ref{pr:asub-zero-is-min:a}):}
Note that, by
\eqref{eq:Fstar-down-def},
\begin{equation}  \label{eq:pr:asub-zero-is-min:1}
  \Fstar(\zero)
  =
  \sup_{\xbar'\in\extspace} [-F(\xbar')]
  =
  -\inf F.
\end{equation}
Thus, if $\zero\in\asubdifF{\xbar}$ then
$\inf F = -\Fstar(\zero) \in\R$ by
\Cref{pr:subgrad-imp-in-cldom}(\ref{pr:subgrad-imp-in-cldom:a}),
and $F(\xbar)=-\Fstar(\zero)=\inf F$ by
\Crefequiv{thm:fenchel-subgrad}{thm:fenchel-subgrad:b}{thm:fenchel-subgrad:c}.

Conversely, if $\inf F\in\R$ and $F(\xbar)=\inf F$
then, by
\eqref{eq:pr:asub-zero-is-min:1},
$\Fstar(\zero)\in\R$ and $\Fstar(\zero)=-F(\xbar)$,
implying that $\zero\in\asubdifF{\xbar}$ by
\Crefequiv{thm:fenchel-subgrad}{thm:fenchel-subgrad:a}{thm:fenchel-subgrad:b}.

\pfpart{Part~(\ref{pr:asub-zero-is-min:b}):}
If $F\equiv+\infty$ or $F(\xbar)=-\infty$ then $\xbar$ clearly
minimizes $F$.
And if $\zero\in\asubdifF{\xbar}$ then $\xbar$ minimizes $F$ by
part~(\ref{pr:asub-zero-is-min:a}).

Conversely, suppose $\xbar$ minimizes $F$, and consider cases:
If $\inf F\in\R$ then $\zero\in\asubdifF{\xbar}$
by part~(\ref{pr:asub-zero-is-min:a}).
Otherwise,
if $\inf F = -\infty$ then $F(\xbar)=\inf F=-\infty$.
And
if $\inf F = +\infty$ then $F\equiv+\infty$.%
\indexg{optimality conditions!astral|)}%
\indexg{subdifferentials, astral (primal)!0 included in@$\zero$ included in|)}%
\qedhere
\end{proof-parts}
\end{proof}

\indexg{Fenchel-Young inequality!astral subgradients and|(}%
\indexg{subdifferentials, astral dual!Fenchel-Young and|(}%
\indexg{subdifferentials, astral dual!primal subgradients of conjugate and|(}%
\indexg{subdifferentials, astral (primal)!dual conjugate@of dual conjugate|(}%
We next take a dual perspective. Starting with a function
$\psi:\Rn\rightarrow\Rext$, we derive relationships between
$\psi$'s astral dual subdifferentials,
the astral primal subdifferentials of its dual conjugate $\psistarb$, and
analogues of the Fenchel-Young inequalities involving $\psi$ and $\psistarb$.
Similar to \Cref{thm:fenchel-subgrad}, these relationships are presented
as a chain of implications, which become equivalences under a summability
condition.

\begin{theorem}   \label{thm:psi-subgrad-conds}
  Let $\psi:\Rn\rightarrow\Rext$, $\uu\in\Rn$, and
  $\xbar\in\extspace$.
  Assume $\psi(\uu)\in\R$.
  Consider the following statements:
  \begin{letter-compact}
  \item    \label{thm:psi-subgrad-conds:a}
    $\psi(\uu) = -\psistarb(\xbar) \plusd \xbar\cdot\uu$.
  \item    \label{thm:psi-subgrad-conds:b}
    $\uu\in\asubdifpsistarb{\xbar}$
    and
    $\psidub(\uu) = \psi(\uu)$.
  \item    \label{thm:psi-subgrad-conds:c}
    $\xbar\in\adsubdifpsi{\uu}$.
  \item    \label{thm:psi-subgrad-conds:d}
    $\psistarb(\xbar) =- \psi(\uu)+ \xbar\cdot\uu$.
  \end{letter-compact}
  Then
  (\ref{thm:psi-subgrad-conds:a})
  $\Rightarrow$
  (\ref{thm:psi-subgrad-conds:b}),
  and
  (\ref{thm:psi-subgrad-conds:b})
  $\Rightarrow$
  (\ref{thm:psi-subgrad-conds:c}),
  and
  (\ref{thm:psi-subgrad-conds:c})
  $\Rightarrow$
  (\ref{thm:psi-subgrad-conds:d}).

  Furthermore, if
  $-\psistarb(\xbar)$ and $\xbar\cdot\uu$ are summable, then
  the four
  statements
  are equivalent.
\end{theorem}

\begin{proof}
~

\begin{proof-parts}
\pfpart{%
  (\ref{thm:psi-subgrad-conds:a})
  $\Rightarrow$
  (\ref{thm:psi-subgrad-conds:b}):
}
Suppose statement~(\ref{thm:psi-subgrad-conds:a}) holds.
Then
\[
  \psi(\uu)
  \geq
  \psidub(\uu)
  \geq
  -\psistarb(\xbar) \plusd \xbar\cdot\uu
  =
  \psi(\uu).
\]
The equality is by assumption.
The first inequality is by
\Cref{thm:psi-geq-psidub}(\ref{thm:psi-geq-psidub:c}),
and the second by
definition of astral conjugate
(Eq.~\ref{eq:Fstar-down-def}).
Thus,
\[
  \psi(\uu)
  =
  \psidub(\uu)
  =
  -\psistarb(\xbar) \plusd \xbar\cdot\uu,
\]
implying $\uu\in\asubdifpsistarb{\xbar}$
by
\Crefequiv{thm:fenchel-subgrad}{thm:fenchel-subgrad:a}{thm:fenchel-subgrad:b}
(applied to $F=\psistarb$).

\pfpart{%
  (\ref{thm:psi-subgrad-conds:b})
  $\Rightarrow$
  (\ref{thm:psi-subgrad-conds:c}):
}
Suppose statement~(\ref{thm:psi-subgrad-conds:b}) holds.
Then
$\xbar\in\adsubdifpsidub{\uu}$
by
\Crefequiv{thm:fenchel-subgrad}{thm:fenchel-subgrad:b}{thm:fenchel-subgrad:b-dual}
(applied to $F=\psistarb$).
For $\uu'\in\Rn$, we then have
\[
   \psi(\uu')
   \geq
   \psidub(\uu')
   \geq
   \psidub(\uu) + \xbar\cdot(\uu'-\uu)
   =
   \psi(\uu) + \xbar\cdot(\uu'-\uu).
\]   
The first inequality is by
\Cref{thm:psi-geq-psidub}(\ref{thm:psi-geq-psidub:c}).
The second is by definition of astral dual subgradient
(Eq.~\ref{eqn:psi-subgrad:3-alt})
since
$\xbar\in\adsubdifpsidub{\uu}$.
Thus, $\xbar\in\adsubdifpsi{\uu}$ according to that same definition.

\pfpart{%
  (\ref{thm:psi-subgrad-conds:c})
  $\Rightarrow$
  (\ref{thm:psi-subgrad-conds:d}):
}
This is \Cref{lem:adsub-implies-fenchel}.

\pfpart{%
  Under summability,
  (\ref{thm:psi-subgrad-conds:d})
  $\Rightarrow$
  (\ref{thm:psi-subgrad-conds:a}):
}
Suppose $-\psistarb(\xbar)$ and $\xbar\cdot\uu$ are summable
and
that statement~(\ref{thm:psi-subgrad-conds:d}) holds.
Then, as in the proof of \Cref{thm:fenchel-subgrad},
$\xbar\cdot\uu\in\R$; otherwise, if
$\xbar\cdot\uu\in\{-\infty,+\infty\}$, 
then also $\psistarb(\xbar)=\xbar\cdot\uu$, 
contradicting that
$-\psistarb(\xbar)$ and $\xbar\cdot\uu$ are summable.
Thus, $\psistarb(\xbar)$, $\xbar\cdot\uu$ and $\psi(\uu)$ are all finite,
so (\ref{thm:psi-subgrad-conds:d})
implies (\ref{thm:psi-subgrad-conds:a}).
\qedhere
\end{proof-parts}
\end{proof}

Similar to \Cref{thm:fenchel-subgrad},
without the summability assumption, the four statements
of \Cref{thm:psi-subgrad-conds} are increasingly more permissive,
with no statement generally implying the preceding one,
as we show next:

\begin{example}
First, let $F$ be as in \Cref{ex:fenchel-subgrad-thm-strict},
and let $\psi=\Fstar=\indf{[-1,1]}$, implying
$\psistarb=\Fdub=F$ and $\psidub=\psi$.
Then, as in that example, if $\barx=+\infty$ and $u=1$, then
statement~(\ref{thm:psi-subgrad-conds:b}) of
\Cref{thm:psi-subgrad-conds} holds, but
not~(\ref{thm:psi-subgrad-conds:a}).
And if $\barx=+\infty$ and $u=1/2$,
then~(\ref{thm:psi-subgrad-conds:d}) holds, but
not~(\ref{thm:psi-subgrad-conds:c}).

To see that~(\ref{thm:psi-subgrad-conds:c}) can hold
without~(\ref{thm:psi-subgrad-conds:b}),
let $\psi:\R\rightarrow\Rext$ be defined by
$\psi(u)=\omega(u-1)$ for $u\in\R$.
Then $\psistarb\equiv+\infty$,
because
$\psistarb(\barx)\ge-\psi(0)=+\infty$ for all $\barx\in\eR$
from the definition of dual conjugate (Eq.~\ref{eq:psistar-def:2}).
For $\barx=+\infty$ and $u=1$, we have
$\barx\in\adsubdifpsi{u}$
since, for all $u'\in\R$,
$\psi(u')=\omega(u'-1)=\psi(u)+\barx\inprod(u'-u)$,
so~(\ref{thm:psi-subgrad-conds:c})
holds. However, $\asubdif{\psistarb}{\barx}=\emptyset$
by \Cref{pr:subgrad-imp-in-cldom}(\ref{pr:subgrad-imp-in-cldom:a})
since $\psistarb\equiv+\infty$, so~(\ref{thm:psi-subgrad-conds:b})
does not hold.%
\indexg{Fenchel-Young inequality!astral subgradients and|)}%
\indexg{subdifferentials, astral dual!Fenchel-Young and|)}%
\indexg{subdifferentials, astral dual!primal subgradients of conjugate and|)}%
\indexg{subdifferentials, astral (primal)!dual conjugate@of dual conjugate|)}%
\end{example}

\section{Characterizing primal subgradients with linear tilts}

As we saw in \Cref{ex:fenchel-subgrad-thm-strict},
the equivalence between
statements~(\ref{thm:fenchel-subgrad:a}) and~(\ref{thm:fenchel-subgrad:b})
of \Cref{thm:fenchel-subgrad}
does not hold in general without the additional assumption that
$-F(\xbar)$ and $\xbar\cdot\uu$ are summable.
In other words, without the summability assumption,
it is possible
that $\uu$ is a subgradient of $F$ at $\xbar$ but that
$\Fstar(\uu) \neq - F(\xbar) \plusd \xbar\cdot\uu$.
Nonetheless, there is an equality
related to the Fenchel-Young inequality that does always
characterize when a vector is an astral subgradient of a function.
To explain this, we start with a definition:

\begin{definition}  \label{def:ast-lin-tilt}
\indexg{linear tilt of astral function|(}%
\indexg{linear tilt of astral function!defined|(}%
Let $F:\extspace\rightarrow\Rext$ and let $\uu\in\Rn$.
The \emph{linear tilt} of $F$ by $\uu$
is the function $\Fminusu:\extspace\rightarrow\Rext$ given,
for $\xbar\in\extspace$, by
\begin{equation} \label{eqn:Fminusu-defn}
\indexm{f u 700}{$\Fminusu$}{linear tilt of astral function}%
  \Fminusu(\xbar) = - \xbar\cdot\uu \plusu F(\xbar).
\indexg{linear tilt of astral function!defined|)}%
\end{equation}
\end{definition}

\indexg{Fenchel-Young inequality!linear tilts and|(}%
For $\uu\in\Rn$,
the generalized form of the Fenchel-Young inequality
given in \eqref{eqn:ast-fenchel} can then be rewritten as
\[
   -\Fstar(\uu) \leq \Fminusu(\xbar),
\]
which holds for all $\xbar\in\extspace$.
As a result,
by \Cref{pr:lsc-min-exists-equals-sup},
it also holds for $\Fminusu\negKern$'s lower semicontinuous hull
so that
\begin{equation}   \label{eq:Fstar-leq-lscFminusu}
   {-\Fstar(\uu)} \leq (\lsc\Fminusu)(\xbar)
\end{equation}
for all $\xbar\in\extspace$. This can be viewed as a lower semicontinuous
variant of the Fenchel-Young inequality.%
\indexg{Fenchel-Young inequality!linear tilts and|)}%

In a moment, we will show that
the inequality in \eqref{eq:Fstar-leq-lscFminusu}
holds with equality if and only if $\uu$ is an astral subgradient of
$F$ at $\xbar$ (provided $\Fstar(\uu)\in\R$ and
that $F$ is lower semicontinuous at $\xbar$).
First, we prove
some general properties of linear tilts.

\begin{proposition}   \label{pr:Fminusu-basic-props}
  Let $F:\extspace\rightarrow\Rext$,
  let $\xbar\in\extspace$,
  and let $\uu\in\Rn$.
  Then:
  \begin{letter-compact}
  \item   \label{pr:Fminusu-basic-props:a}
    $\lsc\Fminusu$ attains its minimum at some point in $\extspace$.
    Moreover,
    \begin{equation}   \label{eq:pr:Fminusu-basic-props:1}
      {-\Fstar(\uu)}
      =
      \inf\Fminusu
      =
      \min(\lsc\Fminusu)
      \leq
      (\lsc\Fminusu)(\xbar)
      \leq
      \Fminusu(\xbar).
    \end{equation}
  \item   \label{pr:Fminusu-basic-props:d}
    If $(\lsc\Fminusu)(\xbar)\in\R$, then there exists a sequence
    $\seq{\xbar_t}$ in $\extspace$ with $F(\xbar_t)\in \R$ and
    $\xbar_t\cdot\uu\in\R$ for all $t$, and such that
    $\xbar_t\rightarrow\xbar$ and
    $F(\xbar_t) - \xbar_t\cdot\uu \rightarrow (\lsc\Fminusu)(\xbar)$.
  \item   \label{pr:Fminusu-basic-props:b}
    If $(\lsc F)(\xbar)$ and $-\xbar\cdot\uu$ are summable, then
    \[
      (\lsc\Fminusu)(\xbar)
      =
      (\lsc F)(\xbar)
      -
      \xbar\cdot\uu.
    \]
  \item   \label{pr:Fminusu-basic-props:c}
    Let $G=\lsc F$.
    Then $\lsc\Gminusu=\lsc\Fminusu$.
  \item   \label{pr:Fminusu-basic-props:convex}
    If $F$ is convex then so is $\Fminusu$.
  \end{letter-compact}
\end{proposition}

\begin{proof}
  ~

\begin{proof-parts}
\pfpart{Part~(\ref{pr:Fminusu-basic-props:a}):}
The first equality of
\eqref{eq:pr:Fminusu-basic-props:1}
is because
\[
  -\Fstar(\uu)
  =
  - \sup_{\xbar'\in\extspace} \Bracks{-F(\xbar')\plusd\xbar'\cdot\uu}
  =
  \inf_{\xbar'\in\extspace} \Bracks{F(\xbar')\plusu(-\xbar'\cdot\uu)}
  =
  \inf\Fminusu,
\]
with the first equality by definition of astral conjugate
(Eq.~\ref{eq:Fstar-down-def}).

The second equality of
\eqref{eq:pr:Fminusu-basic-props:1}, as well as
$\lsc\Fminusu\negKern$'s attainment of its minimum,
are by \Cref{pr:lsc-min-exists-equals-sup}.
The second inequality of
\eqref{eq:pr:Fminusu-basic-props:1}
is by
\Cref{prop:lsc:characterize}(\ref{prop:lsc:characterize:a}).

\pfpart{Part~(\ref{pr:Fminusu-basic-props:d}):}
Suppose $(\lsc\Fminusu)(\xbar)\in\R$.
Then by \Cref{pr:lsc-seq-lim-exists}, there exists a sequence
$\seq{\xbar_t}$ in $\extspace$ that converges to $\xbar$ and with
$\Fminusu(\xbar_t)\rightarrow(\lsc\Fminusu)(\xbar)$.
Since $(\lsc\Fminusu)(\xbar)\in\R$, $\Fminusu(\xbar_t)$
must also be finite for all $t$ sufficiently large.
By discarding all other sequence elements, we can assume
$\Fminusu(\xbar_t)=-\xbar_t\cdot\uu\plusu F(\xbar_t)\in\R$ for all $t$,
implying that both $F(\xbar_t)$ and $\xbar_t\cdot\uu$ are
in $\R$ as well, thereby proving the claim.

\pfpart{Part~(\ref{pr:Fminusu-basic-props:b}):}
Suppose $(\lsc F)(\xbar)$ and $-\xbar\cdot\uu$ are summable.
By \Cref{pr:lsc-seq-lim-exists}, there exists a sequence
$\seq{\xbar_t}$ in $\extspace$ converging to $\xbar$ and with
$F(\xbar_t)\rightarrow(\lsc F)(\xbar)$,
further implying $\xbar_t\cdot\uu\rightarrow\xbar\cdot\uu$
(by \Cref{thm:i:1}\ref{thm:i:1c}).
Since $(\lsc F)(\xbar)$ and $-\xbar\cdot\uu$ are summable,
it must also hold that $F(\xbar_t)$ and $-\xbar_t\cdot\uu$ are
summable as well for all sufficiently large $t$
(by \Cref{prop:lim:eR}\ref{i:lim:eR:sum});
by discarding all other
sequence elements, we can assume this holds for all $t$.
Also, the sequence $\seq{\Fminusu(\xbar_t)}$ in $\Rext$ must have a
convergent subsequence; by discarding all other elements, we can
assume that this sequence has a limit.

Thus,
\begin{align*}
  (\lsc\Fminusu)(\xbar)
  \leq
  \lim \Fminusu(\xbar_t)
  &=
  \lim \bigParens{F(\xbar_t) - \xbar_t\cdot\uu}
  \\
  &=
  \lim F(\xbar_t)
  -
  \lim (\xbar_t\cdot\uu)
  =
  (\lsc F)(\xbar) - \xbar\cdot\uu.
\end{align*}
The inequality is by definition of lower semicontinuous hull
(Eq.~\ref{eq:lsc:liminf:X:prelims}).
The second equality is by
\Cref{prop:lim:eR}(\ref{i:lim:eR:sum}) and our summability assumption.

It remains to prove the reverse inequality,
$(\lsc\Fminusu)(\xbar)\geq  (\lsc F)(\xbar) - \xbar\cdot\uu$.
We assume $(\lsc\Fminusu)(\xbar)<+\infty$,
since otherwise the inequality holds trivially.
By \Cref{pr:lsc-seq-lim-exists}, there exists a sequence
$\seq{\xbar_t}$ in $\extspace$ converging to $\xbar$ with
$\Fminusu(\xbar_t)\rightarrow(\lsc\Fminusu)(\xbar)$.
Since $(\lsc\Fminusu)(\xbar)<+\infty$, we can discard all elements for
which $\Fminusu(\xbar_t)=+\infty$, of which there can be only finitely
many, and so assume henceforth that $\Fminusu(\xbar_t)<+\infty$ for
all $t$.
This implies that $F(\xbar_t)<+\infty$ and
$-\xbar_t\cdot\uu<+\infty$ for all $t$, so these terms are summable.
Thus,
\begin{align*}
  (\lsc\Fminusu)(\xbar)
  =  
  \lim \Fminusu(\xbar_t)
  &=
  \liminf\bigParens{F(\xbar_t) - \xbar_t\cdot\uu}
  \\
  &\geq
  \liminf F(\xbar_t)
  \plusd
  \liminf(-\xbar_t\cdot\uu)
  \\
  &\geq
  (\lsc F)(\xbar) - \xbar\cdot\uu.
\end{align*}
The first inequality is by \Cref{prop:lim:eR}(\ref{i:liminf:eR:sum}).
The second inequality is by definition of lower semicontinuous hull
(Eq.~\ref{eq:lsc:liminf:X:prelims}), and since
$\xbar_t\cdot\uu\rightarrow\xbar\cdot\uu$
(by \Cref{thm:i:1}\ref{thm:i:1c}),
also using the fact that
$(\lsc F)(\xbar)$ and $-\xbar\cdot\uu$ are summable.
This completes the proof.

\pfpart{Part~(\ref{pr:Fminusu-basic-props:c}):}
First, $G\leq F$
(by \Cref{prop:lsc:characterize}\ref{prop:lsc:characterize:a}),
implying, from definitions
(Eqs.~\ref{eqn:Fminusu-defn} and~\ref{eq:lsc:liminf:X:prelims}),
that $\Gminusu\leq\Fminusu$
and so that $\lsc\Gminusu\leq\lsc\Fminusu$.

For the reverse inequality,
we first claim that
$(\lsc\Fminusu)(\xbar)\leq\Gminusu(\xbar)$ for all $\xbar\in\extspace$.
This is immediate if $\Gminusu(\xbar)=+\infty$, so
suppose $\Gminusu(\xbar)<+\infty$.
Then also $G(\xbar)<+\infty$ and $\xbar\cdot\uu>-\infty$.
Thus, $G(\xbar)$ and $-\xbar\cdot\uu$ are summable,
implying, by
part~(\ref{pr:Fminusu-basic-props:b}),
that
  $(\lsc\Fminusu)(\xbar)
   =(\lsc F)(\xbar)-\xbar\inprod\uu
   =G(\xbar)-\xbar\cdot\uu
   =\Gminusu(\xbar)$.

Hence,
$\lsc\Fminusu\leq\Gminusu$.
Since $\lsc\Fminusu$ is lower semicontinuous
(by \Cref{prop:lsc:characterize}\ref{prop:lsc:characterize:a}),
this implies that
$\lsc\Fminusu\leq\lsc\Gminusu$
by \Cref{prop:lsc:characterize}(\ref{prop:lsc:characterize:b}),
completing the proof.

\pfpart{Part~(\ref{pr:Fminusu-basic-props:convex}):}
This follows from
\Cref{thm:sum-ast-cvx-fcns}
since the function $\xbar\mapsto\xbar\cdot(-\uu)$ is convex by
\Cref{pr:aff-fcns-cvx}.%
\indexg{linear tilt of astral function|)}%
\qedhere
\end{proof-parts}
\end{proof}

\indexg{linear tilt of astral function!subgradients and|(}%
\indexg{subdifferentials, astral (primal)!characterizations|(}%
\indexg{subdifferentials, astral (primal)!linear tilts and|(}%
We now prove that $\uu\in\asubdifF{\xbar}$
if and only if
the lower semicontinuous variant of the Fenchel-Young inequality
in \eqref{eq:Fstar-leq-lscFminusu}
holds with equality.
We only assume that
$\Fstar(\uu)\in\R$ and that $F$ is lower semicontinuous
at $\xbar$.
Under these assumptions,
we derive additional characterizations for when
$\uu\in\asubdifF{\xbar}$ that are based on sequences
similar to
\Cref{pr:equiv-ast-subdif-defn}(\ref{pr:equiv-ast-subdif-defn:b}).
Note that since $F$ is assumed to be lower semicontinuous at $\xbar$,
the requirement that $F(\xbar_t)\to F(\xbar)$ can be dropped,
leading to simpler conditions.

\begin{theorem}   \label{thm:subgrad-eq-fenchel}
  Let $F:\extspace\rightarrow\Rext$, let $\xbar\in\extspace$, and let
  $\uu\in\Rn$.
  Assume $\Fstar(\uu)\in\R$ and that $(\lsc F)(\xbar)=F(\xbar)$.
  Then the following are equivalent:
  \begin{letter-compact}
  \item   \label{thm:subgrad-eq-fenchel:a}
    $\uu\in\asubdifF{\xbar}$.
  \item   \label{thm:subgrad-eq-fenchel:b}
    $(\lsc\Fminusu)(\xbar) = -\Fstar(\uu)$
  \item   \label{thm:subgrad-eq-fenchel:c}
    There exists a sequence $\seq{\xbar_t}$ in $\extspace$ with
    $F(\xbar_t)\in\R$ and $\xbar_t\cdot\uu\in\R$ for all $t$,
    and such that
    $\xbar_t\rightarrow\xbar$
    and
    $\xbar_t\cdot\uu - F(\xbar_t) \rightarrow \Fstar(\uu)$.
  \item   \label{thm:subgrad-eq-fenchel:d}
    There exists a sequence $\seq{\xbar_t}$ in $\extspace$ with
    $\xbar_t\rightarrow\xbar$ and
    \begin{equation}   \label{eq:thm:subgrad-eq-fenchel:3}
       \limsup\bigBracks{-F(\xbar_t) \plusd \xbar_t\cdot\uu}
           \geq \Fstar(\uu).
    \end{equation}
  \item   \label{thm:subgrad-eq-fenchel:e}
    $\xbar$ minimizes $\lsc\Fminusu$.
  \end{letter-compact}
\end{theorem}

\begin{proof}
~

\begin{proof-parts}
\pfpart{%
  (\ref{thm:subgrad-eq-fenchel:a})
  $\Rightarrow$
  (\ref{thm:subgrad-eq-fenchel:d}):
}
Immediate by
\Cref{pr:equiv-ast-subdif-defn}(\ref{pr:equiv-ast-subdif-defn:a},\ref{pr:equiv-ast-subdif-defn:b}).
\pfpart{%
  (\ref{thm:subgrad-eq-fenchel:d})
  $\Rightarrow$
  (\ref{thm:subgrad-eq-fenchel:b}):
}
Suppose there exists a sequence as in 
statement~(\ref{thm:subgrad-eq-fenchel:d}).
Then
\[
  (\lsc\Fminusu)(\xbar)
  \leq
  \liminf \Fminusu(\xbar_t)
  \leq
  -\Fstar(\uu)
  \leq
  (\lsc\Fminusu)(\xbar).
\]
The first inequality is by definition of lower semicontinuous hull
(Eq.~\ref{eq:lsc:liminf:X:prelims}).
The second inequality is equivalent to the assumption in
\eqref{eq:thm:subgrad-eq-fenchel:3}.
The last inequality is by
\Cref{pr:Fminusu-basic-props}(\ref{pr:Fminusu-basic-props:a}).
This proves the claimed equality.

\pfpart{%
  (\ref{thm:subgrad-eq-fenchel:b})
  $\Rightarrow$
  (\ref{thm:subgrad-eq-fenchel:c}):
}
This follows directly from
\Cref{pr:Fminusu-basic-props}(\ref{pr:Fminusu-basic-props:d})
since $\Fstar(\uu)\in\R$.

\pfpart{%
  (\ref{thm:subgrad-eq-fenchel:c})
  $\Rightarrow$
  (\ref{thm:subgrad-eq-fenchel:a}):
}
Suppose there exists a sequence $\seq{\xbar_t}$ as in 
statement~(\ref{thm:subgrad-eq-fenchel:c}).
The sequence $(F(\xbar_t))$ in $\eR$ must have a convergent subsequence;
by discarding all other elements, we can assume that this
sequence has a limit.
Then
\begin{align*}
  F(\xbar)
  =
  (\lsc F)(\xbar)
  \le
  \lim F(\xbar_t)
  &=
  \lim\bigParens{
    [F(\xbar_t) - \xbar_t\cdot\uu] + \xbar_t\cdot\uu
  }
\\
  &=
  -\Fstar(\uu)+\xbar\inprod\uu
  \le
  F(\xbar).
\end{align*}
The first equality is by assumption.
The first inequality is from the definition of lower
semicontinuous hull (Eq.~\ref{eq:lsc:liminf:X:prelims}).
The third equality is
because 
$F(\xbar_t) - \xbar_t\cdot\uu\rightarrow -\Fstar(\uu)$
(by assumption)
and $\xbar_t\cdot\uu\rightarrow\xbar\cdot\uu$
(by \Cref{thm:i:1}\ref{thm:i:1c}), and also using
that ${\Fstar(\uu)\in\R}$, so that the limit of the sum is the sum of the limits
(\Cref{prop:lim:eR}\ref{i:lim:eR:sum}). The final inequality
is from the alternative form of the Fenchel-Young inequality in
\eqref{eqn:ast-fenchel-alt}.

Thus, $\lim F(\xbar_t)=F(\xbar)$.
Therefore, $\uu\in\asubdifF{\xbar}$ by
\Cref{pr:equiv-ast-subdif-defn}(\ref{pr:equiv-ast-subdif-defn:b},\ref{pr:equiv-ast-subdif-defn:a})
since the sequence $\seq{\xbar_t}$ satisfies all the conditions of
part~(\ref{pr:equiv-ast-subdif-defn:b}) of that
\namecref{pr:equiv-ast-subdif-defn}.

\pfpart{%
  (\ref{thm:subgrad-eq-fenchel:b})
  $\Leftrightarrow$
  (\ref{thm:subgrad-eq-fenchel:e}):
}
Immediate, since
$\min(\lsc\Fminusu)=-\Fstar(\uu)$
by \Cref{pr:Fminusu-basic-props}(\ref{pr:Fminusu-basic-props:a}).
\qedhere
\end{proof-parts}
\end{proof}

The assumptions in
\Cref{thm:subgrad-eq-fenchel}
that $\Fstar(\uu)\in\R$ and that
$(\lsc F)(\xbar)=F(\xbar)$ are both necessary for $\uu$ to be an
astral subgradient of $F$ at $\xbar$
(by
\Cref{pr:subgrad-imp-in-cldom}\ref{pr:subgrad-imp-in-cldom:a}
and~\Cref{thm:subgrad-then-lsc}).
Thus, even without these assumptions, if 
$\uu\in\asubdifF{\xbar}$, then the theorem implies that each of the
other conditions in
parts~(\ref{thm:subgrad-eq-fenchel:b}),~(\ref{thm:subgrad-eq-fenchel:c}),~(\ref{thm:subgrad-eq-fenchel:d})
and~(\ref{thm:subgrad-eq-fenchel:e})
hold as well.
However, without these assumptions, the converse is false.
In other words, it is possible, say, that
$(\lsc\Fminusu)(\xbar) = -\Fstar(\uu)$, but that
$\uu$ is not a subgradient at $\xbar$.
Here is an example:

\begin{example}
On $\Rext$,
let $F$ be the indicator function $\indfa{[-\infty,0)}$.
Let $\barx=u=0$.
Then $(\lsc\FminusU)(\barx)=0=-\Fstar(u)$.
However,
by \Cref{thm:subgrad-then-lsc},
$u$ cannot be an astral subgradient of $F$ at $\barx$
since $F$ is not lower semicontinuous
at $0$.%
\indexg{linear tilt of astral function!subgradients and|)}%
\indexg{subdifferentials, astral (primal)!characterizations|)}%
\indexg{subdifferentials, astral (primal)!linear tilts and|)}%
\end{example}

\section{Linear-tilt characterization for an extension}
\label{sec:lin-tilt-ast-subgrad-ext}

Next, we focus on how characterizations of astral subgradients
based on linear tilts
apply in the important case that $F$ is the extension $\fext$ of a
function $f:\Rn\rightarrow\Rext$.

In \Cref{def:ast-lin-tilt}, we defined the linear tilt of
an astral function on $\extspace$.
\indexg{linear tilt of standard function|(}%
\indexg{linear tilt of standard function!defined|(}%
In studying extensions, it will be useful to consider a form of linear
tilt for functions on $\Rn$, defined as follows:
\begin{definition}  \label{dfn:lin-tilt-ext}
Let $f:\Rn\rightarrow\Rext$ and let $\uu\in\Rn$.
The \emph{linear tilt} of $f$ by $\uu$
is the function $\fminusu:\Rn\rightarrow\Rext$ given,
for $\xx\in\Rn$, by
\begin{equation} \label{eqn:fminusu-defn}
\indexm{f u 300}{$\fminusu$}{linear tilt of standard function}%
  \fminusu(\xx) = f(\xx) - \xx\cdot\uu.%
\indexg{linear tilt of standard function!defined|)}%
\end{equation}
\end{definition}
Thus, if $F$ is defined on $\extspace$, then
$\Fminusu$ denotes the function on $\extspace$ defined in
\eqref{eqn:Fminusu-defn}, while
if $f$ is defined on $\Rn$, then
$\fminusu$ denotes the function on $\Rn$ defined in
\eqref{eqn:fminusu-defn}.

The standard Fenchel-Young inequality (Eq.~\ref{eqn:fenchel-stand})
can then be rewritten as
\[
   -\fstar(\uu) \leq \fminusu(\xx),
\]
which holds for all $\xx\in\Rn$
(by Eq.~\ref{eq:fstar-def}).
As a result, it also holds for $\fminusu\negKern$'s extension
so that
\begin{equation}   \label{eq:fstar-leq-fminusuext}
   {-\fstar(\uu)} \leq \fminusuext(\xbar)
\end{equation}
for all $\xbar\in\extspace$
(by \Cref{pr:fext-min-exists}).
As a special case of
\Cref{thm:subgrad-eq-fenchel}, we will see below that
$\uu$ is an astral subgradient of $\fext$ at $\xbar$ if and only if
\eqref{eq:fstar-leq-fminusuext} holds with equality.
We will also give characterizations using sequences, as in
\Cref{thm:subgrad-eq-fenchel}(\ref{thm:subgrad-eq-fenchel:c},\ref{thm:subgrad-eq-fenchel:d}).

\indexg{linear tilt of standard function!notational convention|(}%
Before continuing, we draw attention to the potential ambiguity of
notation such as~$\fminusuext$, which, without clarification, might mean
the extension of the linear tilt $\fminusu$, or the linear tilt of
the extension $\fext$.
To resolve this ambiguity, we adopt the convention that the linear
tilt operation always takes highest precedence.
Thus, according to this convention, the notation $\fminusuext$
unambiguously denotes
the lower semicontinuous extension of the linear tilt $\fminusu$.
The linear tilt of $\fext$, the lower semicontinuous
extension of $f$, is then instead denoted $\fextminusu$.
(Indeed, these can be different.
For instance, if $f(x)=x$ for $x\in\R$ and $u=1$ then
$\fminusU\equiv 0$ so
$\fminusUext(+\infty)=0$,
but $\fext(+\infty)=+\infty$ so
$\fextminusU(+\infty)=+\infty$.)

Likewise, we write $\fminusutriv$ for the trivial extension of the
linear tilt $\fminusu$, and $\ftrivminusu$ for the linear tilt of
the trivial extension $\ftriv$.
(Here, however, the ambiguity is unimportant, other than for
presentation, since we will see in a moment that
$\fminusutriv=\ftrivminusu$ always.)
The same convention also applies to conjugates, so that
$\fminusustar$ denotes the conjugate
\indexg{linear tilt of standard function!notational convention|)}%
of~$\fminusu$.

The next proposition delineates
how
$\fminusuext$ and $\fminusutriv$ relate to
$\fextminusu$ and $\ftrivminusu$:

\begin{proposition}   \label{pr:minusuext-vs-extminusu}
  Let $f:\Rn\rightarrow\Rext$, and let $\uu\in\Rn$.
  Then:
  \begin{letter-compact}
  \item   \label{pr:minusuext-vs-extminusu:a}
    $\fminusutriv = \ftrivminusu$.
  \item   \label{pr:minusuext-vs-extminusu:b}
    $\fminusuext = \lsc{\fextminusu} = \lsc{\ftrivminusu}$.
  \end{letter-compact}
\end{proposition}

\begin{proof}
~

\begin{proof-parts}
\pfpart{Part~(\ref{pr:minusuext-vs-extminusu:a}):}
Let $\xbar\in\extspace$.
By definition of trivial extension
(Eq.~\ref{eq:def:triv-ext:1}),
if $\xbar\not\in\Rn$ then
$\ftriv(\xbar)=\fminusutriv(\xbar)=+\infty$,
implying
$\ftrivminusu(\xbar)=-\xbar\cdot\uu\plusu\ftriv(\xbar)=+\infty=\fminusutriv(\xbar)$.
Otherwise, if $\xbar=\xx\in\Rn$, then
\[
  \ftrivminusu(\xx)
  =
  \ftriv(\xx) - \xx\cdot\uu
  =
  f(\xx) - \xx\cdot\uu
  =
  \fminusu(\xx)
  =
  \fminusutriv(\xx).
\]
Thus, in all cases,
$\ftrivminusu(\xbar)=\fminusutriv(\xbar)$, proving the claim.

\pfpart{Part~(\ref{pr:minusuext-vs-extminusu:b}):}
We have
\[
  \fminusuext
  =
  \lsc{\fminusutriv}
  =
  \lsc{\ftrivminusu}
  =
  \lsc{\lscftrivminusu}
  =
  \lsc{\fextminusu},
\]
where $\lscftrivminusu$ denotes the linear tilt (by $\uu$) of
$\lsc\ftriv$.
The first and last equalities are by
\Cref{pr:lsc-ftriv-is-fext}
(applied respectively to $\fminusu$ and to $f$).
The second equality is by
part~(\ref{pr:minusuext-vs-extminusu:a}).
The third equality is by
\Cref{pr:Fminusu-basic-props}(\ref{pr:Fminusu-basic-props:c})
(applied to $\ftriv$).
\qedhere
\end{proof-parts}
\end{proof}
  
We next give some simple facts about linear tilts:

\begin{proposition}  \label{pr:fminusu-props}
  Let $f:\Rn\rightarrow\Rext$,
  $\xbar\in\extspace$,
  and
  $\uu\in\Rn$.
  Then: %
  \begin{letter-compact}
  \item     \label{pr:fminusu-props:a}
    $\dom{\fminusu}=\dom{f}$.
  \item     \label{pr:fminusu-props:b}
    $\fminusustar(\ww)=\fstar(\ww+\uu)$ for $\ww\in\Rn$.
  \item     \label{pr:fminusu-props:d}
    $-\fstar(\uu)=\inf \fminusu=\min\fminusuext \leq \fminusuext(\xbar)$.
  \item     \label{pr:fminusu-props:e}
    If $\fext(\xbar)$ and $-\xbar\cdot\uu$ are summable
    then $\fminusuext(\xbar)=\fext(\xbar)-\xbar\cdot\uu$.
  \item     \label{pr:fminusu-props:convex}
    If $f$ is convex then so is $\fminusu$.
  \end{letter-compact}
\end{proposition}

\begin{proof}
~

\begin{proof-parts}
\pfpart{Part~(\ref{pr:fminusu-props:a}):}
For $\xx\in\Rn$, $\fminusu(\xx)=+\infty$ if and only if
$f(\xx)=+\infty$, proving the claim.

\pfpart{Part~(\ref{pr:fminusu-props:b}):}
Immediate from \Cref{pr:conj-shift}(\ref{pr:conj-shift:b}).

\pfpart{Part~(\ref{pr:fminusu-props:d}):}
From %
\eqref{eq:fstar-def},
$\fstar(\uu)=\sup (-\fminusu)=-\inf \fminusu=-\min\fminusuext$,
with the last equality from \Cref{pr:fext-min-exists}.

\pfpart{Part~(\ref{pr:fminusu-props:e}):}
This follows from
\Cref{pr:ext-sum-fcns}(\ref{pr:ext-sum-fcns:b})
applied to $f$ and $g(\xx)=-\xx\cdot\uu$,
since $g$'s extension is $\gext(\xbar)=-\xbar\cdot\uu$
(\Cref{ex:ext-affine}),
and since $g$ is extensibly continuous everywhere
(\Cref{thm:i:1}\ref{thm:i:1c}).

\pfpart{Part~(\ref{pr:fminusu-props:convex}):}
This follows from \Cref{pr:std-sum-fcns-cvx}
since every linear function is convex.
\qedhere
\end{proof-parts}
\end{proof}

The astral subgradients of the extension of a linear tilt,
as well as the
\indexg{subdifferentials, astral dual!conjugate@of conjugate|(}%
astral dual subgradients of its conjugate,
are both straightforward to compute, as we show next:

\begin{proposition}  \label{pr:fminusu-subgrad}
  Let $f:\Rn\rightarrow\Rext$,
  and
  let $\uu\in\Rn$.
  Then:
  \begin{letter-compact}
  \item  \label{pr:fminusu-subgrad:ext}
    $\asubdiffminusuext{\xbar} = \asubdiffext{\xbar} - \uu$
    for $\xbar\in\extspace$.
  \item  \label{pr:fminusu-subgrad:conj}
    $\adsubdiffminusustar{\ww} = \adsubdiffstar{\ww+\uu}$
    for $\ww\in\Rn$.
  \end{letter-compact}
\end{proposition}

\begin{proof}
~

\begin{proof-parts}
\pfpart{Part~(\ref{pr:fminusu-subgrad:ext}):}
Let $\xbar\in\extspace$.
Suppose $\ww\in\asubdiffext{\xbar}$.
Then $\fstar(\ww)\in\R$ by
\Cref{pr:subgrad-imp-in-cldom}(\ref{pr:subgrad-imp-in-cldom:c}),
so
$\fminusustar(\ww-\uu)=\fstar(\ww)\in\R$ by
\Cref{pr:fminusu-props}(\ref{pr:fminusu-props:b}).
Further, by
\Crefequiv{thm:fminus-subgrad-char}{thm:fminus-subgrad-char:a}{thm:fminus-subgrad-char:b},
there exists a sequence $\seq{\xx_t}$ in $\Rn$
for which $\xx_t\rightarrow\xbar$ and
\begin{equation}  \label{eqn:pr:fminusu-subgrad:1}
  \xx_t\cdot\ww - f(\xx_t) \rightarrow \fstar(\ww).
\end{equation}
By adding and subtracting $\xx_t\inprod\uu$ on the left-hand side, and applying
\Cref{pr:fminusu-props}(\ref{pr:fminusu-props:b})
on the right-hand side, this is equivalent to
\[
  \xx_t\cdot(\ww-\uu) - \fminusu(\xx_t)
  \rightarrow
  \fminusustar(\ww-\uu),
\]
implying
$\ww-\uu \in \asubdiffminusuext{\xbar}$
by
\Crefequiv{thm:fminus-subgrad-char}{thm:fminus-subgrad-char:c}{thm:fminus-subgrad-char:a}.
Thus, 
$\asubdiffext{\xbar} - \uu \subseteq\asubdiffminusuext{\xbar}$.

For the reverse inclusion,
suppose
$\ww \in \asubdiffminusuext{\xbar}$.
Let $f'=\fminusu$ so that
$\ww\in\asubdif{\efp}{\xbar}$.
Writing $\fpminusminusu$ for the linear tilt of $f'$ by $-\uu$,
it then follows by the first part of the proof that
$\ww+\uu\in\asubdif{\fpminusminusuext}{\xbar}$.
Also, $\fpminusminusu=f$ since, for all $\xx\in\Rn$,
\[
  \fpminusminusu(\xx)
  =
  f'(\xx) + \xx\cdot\uu
  =
  \bigParens{f(\xx) - \xx\cdot\uu} + \xx\cdot\uu
  =
  f(\xx).
\]
Thus, $\ww+\uu\in\asubdiffext{\xbar}$ so
$\asubdiffminusuext{\xbar} \subseteq \asubdiffext{\xbar} - \uu$.

\pfpart{Part~(\ref{pr:fminusu-subgrad:conj}):}
Let $\xbar\in\extspace$ and $\ww\in\Rn$.
Then
\begin{alignat*}{3}
  \xbar\in\adsubdiffminusustar{\ww}
  &\;\Leftrightarrow\;
  \forall \ww'\in\Rn,\;\;
  &
  \forcewidthof[c]{$\fstar(\ww+\uu+\ww')$}{%
     $\fminusustar(\ww+\ww')$}
  &\geq
  \forcewidthof[c]{$\fstar(\ww+\uu)$}{%
     $\fminusustar(\ww)$}
  &&\plusd
  \xbar\cdot\ww'
  \\
  &\;\Leftrightarrow\;
  \forall \ww'\in\Rn,\;\;
  &\fstar(\ww+\uu+\ww')
  &\geq
  \fstar(\ww+\uu)
  &&\plusd
  \xbar\cdot\ww'
  \\
  &\;\Leftrightarrow\;
  \makebox[0pt][l]{$\xbar\in\adsubdiffstar{\ww+\uu}$.}
\end{alignat*}
The first and last equivalence are by definition of astral dual
subdifferential (Eq.~\ref{eqn:psi-subgrad:1}).
The second equivalence is by
\Cref{pr:fminusu-props}(\ref{pr:fminusu-props:b}).%
\indexg{linear tilt of standard function|)}%
\indexg{subdifferentials, astral dual!conjugate@of conjugate|)}%
\qedhere
\end{proof-parts}
\end{proof}

\indexg{linear tilt of standard function!subgradients and|(}%
\indexg{subdifferentials of extensions!characterizations|(}%
\indexg{subdifferentials of extensions!linear tilts and|(}%
As discussed above,
the next theorem characterizes the subgradients of an extension along
the lines of
\Cref{thm:subgrad-eq-fenchel}.

\begin{theorem}  \label{thm:fminus-subgrad-char}
  Let $f:\Rn\rightarrow\Rext$,
  $\xbar\in\extspace$, and $\uu\in\Rn$.
  Assume $\fstar(\uu)\in\R$.
  Then the following are equivalent:
  \begin{letter-compact}
  \item  \label{thm:fminus-subgrad-char:a}
    $\uu\in\asubdiffext{\xbar}$.
  \item  \label{thm:fminus-subgrad-char:d}
    $\fminusuext(\xbar) = -\fstar(\uu)$.
  \item  \label{thm:fminus-subgrad-char:b}
    There exists a sequence $\seq{\xx_t}$ in $\Rn$
    with
    $\xx_t\rightarrow\xbar$,
    $f(\xx_t)\in\R$ for all $t$,
    and
    $\xx_t\cdot\uu - f(\xx_t) \rightarrow \fstar(\uu)$.
  \item  \label{thm:fminus-subgrad-char:c}
    There exists a sequence $\seq{\xx_t}$ in $\Rn$
    such that
    $\xx_t\rightarrow\xbar$
    and
    \begin{equation}  \label{eqn:thm:fminus-subgrad-char:2}
      \limsup\bigBracks{\xx_t\cdot\uu - f(\xx_t)} \geq \fstar(\uu).
    \end{equation}
  \item  \label{thm:fminus-subgrad-char:e}
    $\xbar$ minimizes $\fminusuext$;
    that is, $\fminusuext(\xbar) = \inf \fminusu$.
  \end{letter-compact}
\end{theorem}

\begin{proof}
~

\begin{proof-parts}
\pfpart{%
  (\ref{thm:fminus-subgrad-char:a})
  $\Rightarrow$
  (\ref{thm:fminus-subgrad-char:d}):
}
If $\uu\in\asubdiffext{\xbar}$ then
\[
  \fminusuext(\xbar)
  =
  \bigBracks{\lsc{\fextminusu}}(\xbar)
  =
  -\fextstar(\uu)
  =
  -\fstar(\uu).
\]
The first equality is by
\Cref{pr:minusuext-vs-extminusu}(\ref{pr:minusuext-vs-extminusu:b}),
the second by
\Crefequiv{thm:subgrad-eq-fenchel}{thm:subgrad-eq-fenchel:a}{thm:subgrad-eq-fenchel:b}
(applied with $F=\fext$),
and 
the third by
\Cref{pr:fextstar-is-fstar}.

\pfpart{%
  (\ref{thm:fminus-subgrad-char:d})
  $\Rightarrow$
  (\ref{thm:fminus-subgrad-char:b}):
}
Suppose statement~(\ref{thm:fminus-subgrad-char:d}) holds.
Then
\[
   \bigBracks{\lsc{\ftrivminusu}}(\xbar)
   =
   \fminusuext(\xbar)
   =
   -\fstar(\uu),
\]
where the first equality is by
\Cref{pr:minusuext-vs-extminusu}(\ref{pr:minusuext-vs-extminusu:b}),
and the second by assumption.
Therefore, by
\Cref{pr:Fminusu-basic-props}(\ref{pr:Fminusu-basic-props:d})
(applied with $F=\ftriv$)
and since $\fstar(\uu)\in\R$,
there exists a sequence $\seq{\xbar_t}$ in $\extspace$ with
$\xbar_t\rightarrow\xbar$,
$\ftriv(\xbar_t)\in\R$ for all $t$,
and
$\xbar_t\cdot\uu - \ftriv(\xbar_t)\rightarrow\fstar(\uu)$.
For all $t$, 
since $\ftriv(\xbar_t)\in\R$, we must have $\xbar_t=\xx_t$ for some
$\xx_t\in\Rn$ so 
that $\ftriv(\xbar_t)=f(\xx_t)$.
This proves all the parts of the claim.

\pfpart{%
  (\ref{thm:fminus-subgrad-char:b})
  $\Rightarrow$
  (\ref{thm:fminus-subgrad-char:c}):
}
This is immediate.

\pfpart{%
  (\ref{thm:fminus-subgrad-char:c})
  $\Rightarrow$
  (\ref{thm:fminus-subgrad-char:a}):
}
Suppose there exists a sequence $\seq{\xx_t}$ as in
statement~(\ref{thm:fminus-subgrad-char:c}).
Then
\[
  \limsup\bigBracks{\xx_t\cdot\uu - \fext(\xx_t)}
  \geq
  \limsup\bigBracks{\xx_t\cdot\uu - f(\xx_t)}
  \geq
  \fstar(\uu)
  =
  \fextstar(\uu),
\]
where the first inequality is because $\fext(\xx_t)\leq f(\xx_t)$ for
all $t$ (by \Cref{pr:h:1}\ref{pr:h:1a}),
and the second is by assumption.
Note also that $\fext$ is lower semicontinuous
(by \Cref{prop:ext:F}\ref{prop:ext:F:a}).
By
\Crefequiv{thm:subgrad-eq-fenchel}{thm:subgrad-eq-fenchel:d}{thm:subgrad-eq-fenchel:a},
it therefore follows that
$\uu\in\asubdiffext{\xbar}$.

\pfpart{%
  (\ref{thm:fminus-subgrad-char:d})
  $\Leftrightarrow$
  (\ref{thm:fminus-subgrad-char:e}):
}
This is immediate from
\Cref{pr:fminusu-props}(\ref{pr:fminusu-props:d}).
\qedhere
\end{proof-parts}
\end{proof}

Here is an example:

\begin{example}[Astral subgradients of an affine function]
\label{ex:affine-subgrad-new}
\indexg{subdifferentials of extensions!affine function@of affine function|(}%
\indexg{subdifferentials of extensions!examples|(}%
\indexg{affine functions (standard)!subgradients of extension|(}%
We compute the astral subgradients of the extension of an affine function
$f(\xx)=\xx\cdot\aaa+\beta$ for $\xx\in\Rn$,
where $\aaa\in\Rn$ and $\beta\in\R$.
This function's conjugate is
\[
  \fstar(\uu)
  =
  \begin{cases}
    -\beta          & \text{if $\uu=\aaa$,} \\
    +\infty         & \text{otherwise.}
  \end{cases}
\]
Therefore, by
\Cref{pr:subgrad-imp-in-cldom}(\ref{pr:subgrad-imp-in-cldom:c}),
$\fext$ can have no subgradients other than $\aaa$ at any point.
Moreover, the function's linear tilt by $\aaa$ is the constant function
$\fminusa\equiv\beta$,
whose extension is also constant,
$\fminusaext\equiv\beta$.
Thus, for all $\xbar\in\extspace$,
$\fminusaext(\xbar)=\beta=-\fstar(\aaa)$, so
$\aaa\in\asubdiffext{\xbar}$ by
\Cref{thm:fminus-subgrad-char}(\ref{thm:fminus-subgrad-char:d},\ref{thm:fminus-subgrad-char:a}).
Therefore, $\asubdiffext{\xbar}=\{\aaa\}$ for all $\xbar\in\extspace$.%
\indexg{subdifferentials of extensions!affine function@of affine function|)}%
\indexg{subdifferentials of extensions!examples|)}%
\indexg{affine functions (standard)!subgradients of extension|)}%
\end{example}

\Cref{thm:fminus-subgrad-char}(\ref{thm:fminus-subgrad-char:a},\ref{thm:fminus-subgrad-char:d})     
characterizes astral primal
subgradients in terms of the analogue
of the Fenchel-Young inequality given in
\eqref{eq:fstar-leq-fminusuext} holding with
\indexg{subdifferentials of extensions!characterizations|)}%
\indexg{subdifferentials of extensions!linear tilts and|)}%
equality.
\indexg{subdifferentials, astral dual!conjugate@of conjugate|(}%
In a similar but slightly different way,
for a convex function $f:\Rn\rightarrow\Rext$,
we can further use this latter condition to characterize the
astral dual subgradients of $f$'s conjugate, as shown next:\looseness=-1

\begin{theorem}   \label{thm:dual-subgrad-fenchel-tilt}
  Let $f:\Rn\rightarrow\Rext$ be convex with $f\not\equiv+\infty$,
  let $\xbar\in\extspace$ and $\uu\in\Rn$.
  Then the following are equivalent:
  \begin{letter-compact}
  \item   \label{thm:dual-subgrad-fenchel-tilt:a}
    $\fminusuext(\xbar) = -\fstar(\uu)$.
  \item   \label{thm:dual-subgrad-fenchel-tilt:b}
    $\xbar\in\adsubdiffstar{\uu}$
    and
    $\xbar\in\cldom{f}$.
  \end{letter-compact}
\end{theorem}

\begin{proof}
~

\begin{proof-parts}
\pfpart{%
  (\ref{thm:dual-subgrad-fenchel-tilt:a})
  $\Rightarrow$
  (\ref{thm:dual-subgrad-fenchel-tilt:b}):
}
Suppose $\fminusuext(\xbar) = -\fstar(\uu)$.
Since $f\not\equiv+\infty$, $\fstar(\uu)>-\infty$
(by \Cref{pr:conj-props}\ref{pr:conj-props:c1})
so $\fminusuext(\xbar) < +\infty$.
Thus,
\[
  \xbar
  \in
  \dom\fminusuext
  \subseteq
  \cldom\fminusuext
  =
  \cldom\fminusu
  =
  \cldom f,
\]
where the first equality is by
\Cref{pr:h:1}(\ref{pr:h:1c})
and the second by
\Cref{pr:fminusu-props}(\ref{pr:fminusu-props:a}).

Further,
\[
  -\fminusustar(\zero)
  =
  -\fstar(\uu)
  =
  \fminusuext(\xbar)
  \geq
  \fminusudub(\xbar)
  =
  \sup_{\ww\in\Rn}
    \bigBracks{-\fminusustar(\ww)\plusd\xbar\cdot\ww}.
\]
The first equality is by
\Cref{pr:fminusu-props}(\ref{pr:fminusu-props:b}),
and the second by assumption.
The inequality and third equality are by
\Cref{thm:fext-dub-sum}(\ref{thm:fext-dub-sum:a},\ref{thm:fext-dub-sum:b}).
Thus,
$-\fminusustar(\zero)\geq-\fminusustar(\ww)\plusd\xbar\cdot\ww$
for all $\ww\in\Rn$, so by definition of dual subgradient
(Eq.~\ref{eqn:psi-subgrad:3}),
$\xbar\in\adsubdiffminusustar{\zero}=\adsubdiffstar{\uu}$,
with equality by
\Cref{pr:fminusu-subgrad}(\ref{pr:fminusu-subgrad:conj}).

\pfpart{%
  (\ref{thm:dual-subgrad-fenchel-tilt:b})
  $\Rightarrow$
  (\ref{thm:dual-subgrad-fenchel-tilt:a}):
}
Suppose
    $\xbar\in\adsubdiffstar{\uu}$
    and
    $\xbar\in\cldom{f}$.
Then by similar reasoning as above,
\[
  -\fstar(\uu)
  =
  -\fminusustar(\zero)
  =
  \sup_{\ww\in\Rn}
    \bigBracks{-\fminusustar(\ww)\plusd\xbar\cdot\ww}
  =
  \fminusudub(\xbar)
  =
  \fminusuext(\xbar).
\]
The second equality is because
$\xbar\in\adsubdiffminusustar{\zero}$
and by definition of dual subgradient
(and since the supremum is attained when $\ww=\zero$).
The third equality is by
\Cref{thm:fext-dub-sum}(\ref{thm:fext-dub-sum:b}).
The fourth equality is by
\Cref{thm:fext-neq-fdub} since
$\xbar\in\cldom{f}=\cldom\fminusu$
(using \Cref{pr:fminusu-props}\ref{pr:fminusu-props:a}).
\qedhere
\end{proof-parts}
\end{proof}

The equivalence given in \Cref{thm:dual-subgrad-fenchel-tilt} does not
hold if the condition that $\xbar\in\cldom{f}$ is omitted.
For instance, if $f:\R^2\rightarrow\Rext$
is the restricted linear function from \Cref{ex:biconj:notext}
with $\xbar=\limray{\ee_1}\plusl(-\ee_2)$ and $\uu=\zero$,
then it can be checked that
$\xbar\in\adsubdiffstar{\uu}$
but
$\fminusuext(\xbar) = +\infty\neq-\infty = -\fstar(\uu)$.%
\indexg{linear tilt of standard function!subgradients and|)}%
\indexg{subdifferentials, astral dual!conjugate@of conjugate|)}%

\section{Inverse relation between primal and dual subgradients}

\indexg{subdifferentials, astral (primal)!dual subgradients of conjugate and|(}%
\indexg{subdifferentials, astral dual!conjugate@of conjugate|(}%
We turn now to conditions under which
the astral primal subdifferential
$\asubdifplain{F}$ of a function $F:\extspace\rightarrow\Rext$ and the
astral dual subdifferential $\adsubdifplain{\Fstar}$ of its conjugate
are inverses of one another.
\Cref{thm:fenchel-subgrad}(\ref{thm:fenchel-subgrad:b},\ref{thm:fenchel-subgrad:b-dual})
shows that $\uu\in\asubdifF{\xbar}$ always implies
$\xbar\in\adsubdifFstar{\uu}$ (assuming throughout this discussion
that $\Fstar(\uu)\in\R$), and that the converse also holds
if $\Fdub(\xbar)=F(\xbar)$ and if
$-F(\xbar)$ and $\xbar\cdot\uu$ are summable.
This summability condition will always hold, for instance, if
$F(\xbar)\in\R$ or if $\xbar=\xx\in\Rn$. The latter condition, when
applied to an extension $\ef$ of a convex and proper function $f:\Rn\to\eRn$,
recovers the standard result from
\Crefequiv{pr:stan-subgrad-equiv-props}{pr:stan-subgrad-equiv-props:a}{pr:stan-subgrad-equiv-props:c}
(using that
$\smash{\cl f=\fdubs=\fextdub}$ by \Cref{thm:fext-dub-sum}\ref{thm:fext-dub-sum:a}\ref{thm:fext-dub-sum:fdubs}).

Our next example shows that without the summability condition,
the inverse relationship between the subdifferentials of $F$ and $\Fstar$ need not hold, even if $F$ is convex.
In other words,
it is possible that
$\xbar\in\adsubdifFstar{\uu}$ but that
$\uu\not\in\asubdifF{\xbar}$,
even if $\Fstar(\uu)\in\R$,
$F(\xbar)=(\lsc F)(\xbar)=\Fdub(\xbar)$,
and $F$ is convex.

\begin{example}   \label{ex:cvx-not-subgrad-inv}
As in
\Cref{thm:closure-not-always-convex} (with $n=2$),
let
\[
    C
    =
    \bigBraces{ \lambda_1 \ee_1 :\: \lambda_1\in\R }
    \,\cup\, \bigBraces{ \limray{\ee_1} \plusl \lambda_2\ee_2
                                  :\: \lambda_2\in\Rpos }.
\]
As shown in that theorem, $C$ is convex, but its closure is not.
For $\xbar\in\extspac{2}$,
let
\[
  F(\xbar)
  =
  \indfa{C}(\xbar) \plusu \xbar\cdot\ee_2
  =
  \begin{cases}
    \xbar\cdot\ee_2 & \text{if $\xbar\in C$,} \\
    +\infty         & \text{otherwise.} \\
  \end{cases}
\]
This function is convex
(by
Propositions~\ref{pr:ast-ind-fcn-cvx}
and~\ref{pr:aff-fcns-cvx},
and \Cref{thm:sum-ast-cvx-fcns}).
It is equal to $0$ at all points 
$\lambda_1\ee_1$ for $\lambda_1\in\R$,
and to $\lambda_2$ at every point
$\limray{\ee_1}\plusl\lambda_2\ee_2$, for $\lambda_2\in\Rpos$;
the function is $+\infty$ everywhere else.

It can be calculated that $F$'s conjugate is the indicator function
$\Fstar=\inds$ where
$S=\{\lambda \ee_2 :\: \lambda\in(-\infty,1]\}$.
Let $\uu=\ee_2$ and $\xbar=\limray{\ee_2}$.
Let $\uu'=\lambda\ee_2$ with $\lambda\in(-\infty,1]$.
Then
\[
   \Fstar(\uu')
   =
   0
   \geq
   \omsf{(\lambda-1)}
   =
   0 + \limray{\ee_2}\cdot(\lambda\ee_2 - \ee_2)
   =
   \Fstar(\uu) + \xbar\cdot(\uu' - \uu).
\]
Therefore, $\xbar\in\adsubdifFstar{\uu}$ by definition
of astral dual subdifferential
(Eq.~\ref{eqn:psi-subgrad:3-alt}).

By
\Cref{thm:psi-subgrad-conds}(\ref{thm:psi-subgrad-conds:c},\ref{thm:psi-subgrad-conds:d}),
this further implies that
$\Fdub(\xbar) = \xbar\cdot\uu-\Fstar(\uu) = +\infty$,
and therefore that $F(\xbar)=(\lsc F)(\xbar)=\Fdub(\xbar)$
by
\Cref{pr:conj-lsc-props}.

Nonetheless, $\xbar\not\in\cldom{F}$
(since, for instance, $\xbar$ is in the open set
$\{\xbar'\in\extspac{2} :\: \xbar'\cdot\ee_2>0,\, \xbar'\cdot\ee_1<1\}$,
which is disjoint from $C=\dom{F}$).
Therefore, $\uu\not\in\asubdifF{\xbar}$ by
\Cref{pr:subgrad-imp-in-cldom}(\ref{pr:subgrad-imp-in-cldom:a}).
\end{example}

Thus, convexity of $F$ and its astral closedness at $\xbar$
do not ensure an inverse relationship between $\asubdifplain{F}$ and
$\adsubdifplain{\Fstar}$.
Nonetheless, as we show next,
a strengthening of the convexity condition does turn out to be sufficient.
In particular, if instead of merely requiring that $F$'s epigraph be convex, we require that the closure of $F$'s epigraph be convex, then
$\uu\in\asubdifF{\xbar}$ if and only if $\xbar\in\adsubdifFstar{\uu}$,
provided $\Fstar(\uu)\in\R$ and that $F$ is lower semicontinuous at
$\xbar$.

\begin{theorem}   \label{thm:conv-implies-subgrad-equiv}
  Let $F:\extspace\rightarrow\Rext$, let $\xbar\in\extspace$, and let
  $\uu\in\Rn$.
  Assume $\clepi{F}$ is convex.
  Then the following are equivalent:
  \begin{letter-compact}
  \item   \label{thm:conv-implies-subgrad-equiv:a}
    $\uu\in\asubdifF{\xbar}$.
  \item   \label{thm:conv-implies-subgrad-equiv:b}
    $\xbar\in\adsubdifFstar{\uu}$,
    $\Fstar(\uu)\in\R$,
    and
    $(\lsc F)(\xbar)=F(\xbar)$.
  \end{letter-compact}
\end{theorem}

\begin{proof}
That (\ref{thm:conv-implies-subgrad-equiv:a}) implies
(\ref{thm:conv-implies-subgrad-equiv:b})
is immediate from
\Cref{pr:subgrad-imp-in-cldom}(\ref{pr:subgrad-imp-in-cldom:a}),
\Cref{thm:fenchel-subgrad}(\ref{thm:fenchel-subgrad:b},\ref{thm:fenchel-subgrad:b-dual}),
and
\Cref{thm:subgrad-then-lsc}.

For the converse, suppose
statement~(\ref{thm:conv-implies-subgrad-equiv:b}) holds.
We aim to prove that $\uu\in\asubdifF{\xbar}$.
The main idea of the proof is to transform the problem
in a way that makes it easier to reason about.
This is not done by direcly transforming the function $F$, but rather
by applying a linear transformation to $\clepi{F}$, the closure of its
epigraph, which we can then work with using the tools developed in
\Cref{sec:fcns-induced-by-sets} for functions induced by sets.
The chosen transformation specifically has the effect of
transforming an astral dual subgradient at $\uu$ to one at $\zero$, a
much easier point to work with.
Translating back to the original setting then completes the proof.

In what follows,
for a set $E\subseteq\extspacnp$,
we make use of the lower envelope function $\settofcnE$
and lower support function $\Esstar$, defined in
Eqs.~(\ref{eq:settofcn-dfn}) and~(\ref{eq:esstar-dfn}).

In particular, let $E=\clepi{F}$, implying that
$\settofcnE=\lsc F$ by 
\Cref{pr:settofcn-props}(\ref{pr:settofcn-props:e}),
and that $\Esstar=\Fstar$ by
\Cref{pr:settofcn-ohull-prop}
since $\epi{F}\subseteq E\subseteq\ohull(\epi{F})$.

For the linear transformation discussed above, as well as its inverse,
we define matrices $\amatu,\amatuinv\in\Rnpnp$ 
which are each identical to the 
$(n+1)\times(n+1)$ identity matrix, except that the first $n$ entries
of the bottom row are equal to $\trans{\uu}$ for $\amatu$,
and to $-\trans{\uu}$ for $\amatuinv$.
Thus, in block form,
\[
   \amatu
   =
   \left[
     \begin{array}{ccc|c}
        & &       &  \\
       ~ & \Idnn & ~ & \zerov{n} \\
        & &       &  \\
       \hline
       \rule{0pt}{2.5ex}
       & \trans{\uu} & &  1
     \end{array}
   \right]
   \quad
   \text{and}
   \quad
   \amatuinv
   =
   \left[
     \begin{array}{ccc|c}
        & &       &  \\
       ~ & \Idnn & ~ & \zerov{n} \\
        & &       &  \\
       \hline
       \rule{0pt}{2.5ex}
       & -\trans{\uu} & &  1
     \end{array}
   \right]
   ,
\]
or equivalently
\begin{equation}
\label{eq:conv-implies-subgrad-equiv:1}
  \amatu
  =
  \Idn{n+1}+\trans{\PPy}\trans{\uu}\PPx
  \quad
  \text{and}
  \quad
  \amatuinv
  =
  \Idn{n+1}-\trans{\PPy}\trans{\uu}\PPx,
\end{equation}
where $\PPx=[\Idnn,\zero_n]$ and $\PPy=[\trans{\zero_n},1]$.
Using \eqref{eq:conv-implies-subgrad-equiv:1} and orthogonality
identities for the matrices $\PPx$ and $\PPy$ (see Eq.~\ref{eq:PPx:PPy:orth}),
we obtain the following identities:
\begin{equation}
  \label{eq:conv-implies-subgrad-equiv:2}
    \PPx\amatu
    =
    \PPx,
    \quad
    \PPy\amatu
    =
    \PPy+\trans{\uu}\PPx,
    \quad
    \PPy\amatuinv
    =
    \PPy-\trans{\uu}\PPx,
    \quad
    \amatu\amatuinv
    =
    \Idn{n+1}.
\end{equation}
In particular, $\amatuinv=\amatu^{-1}$.

Let $E'=\amatuinv E$.
Then $E'$ is closed (in $\extspacnp$) and convex by
Corollaries~\ref{cor:aff-img-closed-is-closed}(\ref{cor:aff-img-closed-is-closed:a})
and~\ref{cor:thm:e:9}.

\begin{claimpx}   \label{cl:thm:conv-implies-subgrad-equiv:1}
  Let $\ww\in\Rn$.
  Then
  $\Epsstarf(\ww)=\Esstarf(\ww+\uu)=\Fstar(\ww+\uu)$.
\end{claimpx}

\begin{proofx}
That $\Esstar=\Fstar$ was noted earlier.
For the claim's first equality, we have
\begin{align*}
  \Epsstarf(\ww)
  =
  \sup_{\zbar'\in\amatuinv E} \bigBracks{\zbar'\cdot\rpair{\ww}{-1}}
  &=
  \sup_{\zbar\in E} \bigBracks{(\amatuinv\zbar)\cdot\rpair{\ww}{-1}}
  \\
  &=
  \sup_{\zbar\in E} \bigBracks{\zbar\cdot\regParens{\transamatuinv\rpair{\ww}{-1}}}
  \\
  &=
  \sup_{\zbar\in E} \bigBracks{\zbar\cdot\rpair{\ww+\uu}{-1}}
  =
  \Esstarf(\ww+\uu).
\end{align*}
The first and last equalities are by definition of
lower support function
(Eq.~\ref{eq:esstar-dfn}).
The third equality is by \Cref{thm:Ax-dot-u}.
The fourth equality is by matrix algebra.
\end{proofx}

From this claim, it follows that, for all $\ww\in\Rn$,
\[
  \Epsstarf(\ww)
  =
  \Fstar(\ww+\uu)
  \geq
  \Fstar(\uu) \plusd \xbar\cdot\ww
  =
  \Epsstarf(\zero) \plusd \xbar\cdot\ww,
\]
where the inequality is from the definition in
\eqref{eqn:psi-subgrad:3-alt} since $\xbar\in\adsubdifFstar{\uu}$.
By that same definition, this proves that
$\xbar\in\adsubdifEpstar{\zero}$.

Next, we have that
\begin{equation}   \label{eq:thm:conv-implies-subgrad-equiv:1}
  \Fstar(\uu)
  =
  \Epsstarf(\zero)
  =
  - \Epsdubf(\xbar)
  =
  - \settofcnEpf(\xbar).
\end{equation}
The first equality is by
\Cref{cl:thm:conv-implies-subgrad-equiv:1}.
The second equality is by
\Cref{thm:psi-subgrad-conds}(\ref{thm:psi-subgrad-conds:c},\ref{thm:psi-subgrad-conds:d})
since $\xbar\in\adsubdifEpstar{\zero}$.
The last equality is by
\Cref{thm:settofcn-eq-biconj} since
the first two equalities imply that
$\Epsdubf(\xbar)=-\Fstar(\uu)\in\R$.

From \eqref{eq:thm:conv-implies-subgrad-equiv:1} combined with
\Cref{pr:settofcn-props}(\ref{pr:settofcn-props:b})
(applicable since $E'$ is closed), there exists $\zbar'\in E'$
such that
\begin{equation}   \label{eq:thm:conv-implies-subgrad-equiv:2}
  \PPx\zbar'=\xbar
  \quad\text{and}\quad
  \PPy\zbar'=-\Fstar(\uu).
\end{equation}
Since $\zbar'\in E'$, there exists $\zbar\in E$ such that
$\zbar'=\amatuinv\zbar$, implying
$\zbar=(\amatu\amatuinv)\zbar=\amatu\zbar'$.
To prove $\uu\in\asubdifF{\xbar}$,
we will argue that $\zbar$ satisfies the conditions of
\Cref{pr:equiv-ast-subdif-defn}(\ref{pr:equiv-ast-subdif-defn:c}).

First,
\begin{equation}   \label{eq:thm:conv-implies-subgrad-equiv:4}
  \PPx\zbar
  = \PPx\amatu\zbar'
  = \PPx\zbar'
  = \xbar,
\end{equation}
where the second equality is by \Cref{eq:conv-implies-subgrad-equiv:2}.
This in turn implies that
\begin{equation}   \label{eq:thm:conv-implies-subgrad-equiv:3}
   F(\xbar)
   =
   (\lsc F)(\xbar)
   =
   \settofcnEf(\xbar)
   =
   \Esdubf(\xbar)
   =
   \Fdub(\xbar).
\end{equation}
The first equality is by assumption.
The third equality is by
\Cref{thm:settofcn-eq-biconj}
since $E$ is closed and convex, and since
$\xbar\in\PPxE$ by
\eqref{eq:thm:conv-implies-subgrad-equiv:4}.
The identities $\settofcnE=\lsc F$ and $\Esstar=\Fstar$ were derived earlier
(from the definition of $E$).

Next, we have
\begin{align}
\notag
  \PPy\zbar
  =
  \PPy\amatu\zbar'
  &=
  (\PPy+\trans{\uu}\PPx)\zbar'
\\
\label{eq:thm:conv-implies-subgrad-equiv:5}
  &=
  \PPy\zbar'\plusl\trans{\uu}\PPx\zbar'
  =
  -\Fstar(\uu) + \trans{\uu}\xbar
  =
  F(\xbar).
\end{align}
The second equality is by \eqref{eq:conv-implies-subgrad-equiv:2}.
The third is by \Cref{prop:commute:AB},
noting that $\PPy\zbar'=-\Fstar(\uu)\in\R$.
The fourth is by \eqref{eq:thm:conv-implies-subgrad-equiv:2}
(and since $\Fstar(\uu)\in\R$).
The last is by
\Cref{thm:fenchel-subgrad}(\ref{thm:fenchel-subgrad:b-dual},\ref{thm:fenchel-subgrad:c})
since $\xbar\in\adsubdifFstar{\uu}$,
using
\eqref{eq:thm:conv-implies-subgrad-equiv:3}
(as well as \Cref{pr:trans-uu-xbar}).

Also,
\begin{align}
  \notag
    \Fstar(\uu)
  =
  -\PPy\zbar'
  =
  -\PPy\amatuinv\zbar
  &=
  -(\PPy-\trans{\uu}\PPx)\zbar
\\
\label{eq:thm:conv-implies-subgrad-equiv:6}
  &=
  \zbar\inprod(\trans{\PPx}\uu-\trans{\PPy})
  =
  \zbar\cdot\rpair{\uu}{-1},
\end{align}
with the first equality from
\eqref{eq:thm:conv-implies-subgrad-equiv:2},
the third from
\eqref{eq:conv-implies-subgrad-equiv:2},
the fourth from \Cref{pr:trans-uu-xbar},
and the fifth by matrix algebra.

Combining
Eqs.~(\ref{eq:thm:conv-implies-subgrad-equiv:4}),~(\ref{eq:thm:conv-implies-subgrad-equiv:5})
and~(\ref{eq:thm:conv-implies-subgrad-equiv:6}),
and since $\zbar\in\clepi{F}$, it now follows that
$\uu\in\asubdifF{\xbar}$ by
\Cref{pr:equiv-ast-subdif-defn}(\ref{pr:equiv-ast-subdif-defn:c},\ref{pr:equiv-ast-subdif-defn:a}),
completing the proof.
\end{proof}

\begin{example}  \label{ex:cvx-not-subgrad-inv-cont}
  \Cref{thm:conv-implies-subgrad-equiv}
  is consistent with
  \Cref{ex:cvx-not-subgrad-inv}.
  Recall that the function $F$ used in that example satisfies
  $\dom F=C$, where $C$ is a convex set in $\extspac{2}$ whose closure, $\Cbar$,
  is not convex. As we argued in \Cref{ex:cvx-not-subgrad-inv},
  $F$~is convex, so $\epi{F}$ is convex.
  Nevertheless, the closure of its epigraph, $\clepi{F}$, is not convex.
  This is because
  \begin{equation}  \label{eq:ex:cvx-not-subgrad-inv-cont:1}
    \PPx(\clepi{F})=\clbar{\PPx(\epi{F})}=\cldom{F}=\Cbar,
  \end{equation}
  where $\PPx=[\Idnn,\zerov{n}]$,
  and the first equality follows by
  \Cref{cor:aff-img-closed-is-closed}(\ref{cor:aff-img-closed-is-closed:b}).
  Since $\Cbar$ is not convex by
  \Cref{thm:closure-not-always-convex},
  \eqref{eq:ex:cvx-not-subgrad-inv-cont:1} implies that
  $\clepi{F}$ also cannot be convex
  (by \Cref{cor:thm:e:9}).
  Thus,
  \Cref{thm:conv-implies-subgrad-equiv}
  is not applicable to this function.%
\indexg{subdifferentials, astral (primal)!dual subgradients of conjugate and|)}%
\end{example}

\indexg{subdifferentials of extensions!dual subgradients of conjugate and|(}%
The condition that $\clepi{F}$ is convex is always satisfied if $F$ is
the extension of a convex function on $\Rn$.
Consequently, we have:

\begin{corollary}   \label{cor:strict-adif-fext-inverses}
  Let $f:\Rn\rightarrow\Rext$ be convex,
  and let $\xbar\in\extspace$ and
  $\uu\in\Rn$.
  Then
  $\uu\in\asubdiffext{\xbar}$
  if and only if
  $\xbar\in\adsubdiffstar{\uu}$
  and
  $\fstar(\uu)\in\R$.
\end{corollary}

\begin{proof}
By \Cref{pr:wasthm:e:3}(\ref{pr:wasthm:e:3c}),
$\clepifext=\clepi{f}$.
Since $f$ is convex, $\epi{f}$ is a convex set, implying
$\clepi{f}$ is as well, by \Cref{thm:e:6}.
Thus, $\clepi{\fext}$ is convex.

Noting that $\fextstar=\fstar$ (\Cref{pr:fextstar-is-fstar})
and that $\fext$ is lower semicontinuous
(\Cref{prop:ext:F}\ref{prop:ext:F:a}),
the corollary now follows directly from
\Cref{thm:conv-implies-subgrad-equiv}.%
\indexg{subdifferentials of extensions!dual subgradients of conjugate and|)}%
\indexg{subdifferentials, astral dual!conjugate@of conjugate|)}%
\end{proof}

\indexg{subdifferentials, astral dual!one dimension@in one dimension|(}%
As a simple example, we can now straightforwardly
characterize astral primal and dual subdifferentials
of proper
convex functions on $\R$.
We first characterize astral dual subgradients of any proper function
$\psi:\R\rightarrow\Rext$.
We already saw in
\Cref{pr:adsubdif-int-rn}
that any finite $x\in\R$ is an astral dual subgradients of $\psi$
at $u\in\R$ if and only if it is a standard subgradient of $\psi$ at $u$.
The next proposition completes the picture,
characterizing when $\barx=\pm\infty$ is an astral dual
subgradient of $\psi$.
Building on the characterization of astral dual subgradients
we then use
\Cref{cor:strict-adif-fext-inverses}
to fully characterize astral primal subgradients
of the extension of any convex proper function on $\R$.

\begin{proposition}   \label{pr:adsubgrad-in-1d}
  Let $\psi:\R\rightarrow\Rext$ be proper, and let $u\in\R$.
  Then:
  \begin{letter-compact}
  \item   \label{pr:adsubgrad-in-1d:a}
    $-\infty\in\adsubdifpsi{u}$
    if and only if $u\leq\inf(\dom\psi)$.
  \item   \label{pr:adsubgrad-in-1d:b}
    $+\infty\in\adsubdifpsi{u}$
    if and only if $u\geq\sup(\dom\psi)$.
  \end{letter-compact}
\end{proposition}

\begin{proof}
We only prove part~(\ref{pr:adsubgrad-in-1d:a});
the proof of part~(\ref{pr:adsubgrad-in-1d:b}) is symmetric.

\begin{proof-parts}
\pfpart{``Only if'' ($\Rightarrow$):}
Suppose $-\infty\in\adsubdifpsi{u}$, and
let $u'\in\R$. Then, from the definition of astral dual subgradient
(Eq.~\ref{eqn:psi-subgrad:3-alt}),
if $u'<u$, then
\[
  \psi(u')\ge\psi(u)\plusd (-\infty)(u'-u)=\psi(u)\plusd(+\infty)=+\infty,
\]
where the last equality is because $\psi(u)>-\infty$ since
$\psi$ is proper.
Thus, $u'<u$ implies that $u'\not\in\dom\psi$,
so all $u''\in\dom\psi$ satisfy $u''\ge u$. Hence, $\inf(\dom\psi)\ge u$.

\pfpart{``If'' ($\Leftarrow$):}
Suppose $u\leq\inf(\dom\psi)$.
Letting $u'\in\R$, we need to show that
\[
  \psi(u')
  \geq
  \psi(u)
  \plusd
  (-\infty)\cdot(u'-u).
\]
If $u'<u$ then $u'\not\in\dom\psi$, so $\psi(u')=+\infty$ and the inequality holds.
If $u'>u$ then the right-hand side is equal to $-\infty$, so the
inequality holds as well. In the remaining case, when
$u'=u$, the inequality also holds (with equality), proving the claim.%
\indexg{subdifferentials, astral dual!one dimension@in one dimension|)}%
\qedhere
\end{proof-parts}
\end{proof}

\indexg{subdifferentials of extensions!one dimension@in one dimension|(}%
We next
compute astral primal subgradients of an extension $\ef$
of a proper convex function $f:\R\to\eR$.
By
\Cref{pr:asubdiffext-at-x-in-rn}(\ref{pr:asubdiffext-at-x-in-rn:c}),
for $x\in\R$, we have
$\asubdiffext{x} = \partial (\lsc f)(x)$.
The remaining subgradients at $\pm\infty$ are given by the next
proposition:

\begin{proposition}  \label{pr:subdif-in-1d}
 Let $f:\R\rightarrow\Rext$ be convex and proper.
 Let $\alpha=\inf(\dom{\fstar})$ and $\beta=\sup(\dom{\fstar})$.
 \begin{letter-compact}
 \item  \label{pr:subdif-in-1d:a}
   $\displaystyle
   \asubdiffext{-\infty}
   =
   \{\alpha\} \cap (\dom\fstar)
   =
   \begin{cases}
     \{\alpha\}   & \text{if $\alpha\in\dom\fstar$,} \\
     \emptyset    & \text{otherwise.}
   \end{cases}
   $
 \item  \label{pr:subdif-in-1d:b}
   $\displaystyle
   \asubdiffext{+\infty}
   =
   \{\beta\} \cap (\dom\fstar)
   =
   \begin{cases}
     \{\beta\}    & \text{if $\beta\in\dom\fstar$,} \\
     \emptyset    & \text{otherwise.}
   \end{cases}
   $
 \end{letter-compact}
\end{proposition}

\begin{proof}
As in \Cref{pr:adsubgrad-in-1d},
we prove only part~(\ref{pr:subdif-in-1d:a}).

Let $u\in\R$.
We have
\begin{align*}
  u\in\asubdiffext{-\infty}
  &\;\;\Leftrightarrow\;\;
  -\infty\in\adsubdiffstar{u}
  \;\;\text{and}\;\;
  \fstar(u)\in\R
  \\
  &\;\;\Leftrightarrow\;\;
  u\leq\alpha
  \;\;\text{and}\;\;
  u\in\dom\fstar.
\end{align*}
The first equivalence is  
by
\Cref{cor:strict-adif-fext-inverses}.
The second equivalence is by
\Cref{pr:adsubgrad-in-1d} (with $\psi=\fstar$),
and because $\fstar$ is proper
(by
\Cref{pr:conj-props-cvx}\ref{pr:conj-props-cvx:a}).
Since $\dom\fstar\subseteq[\alpha,+\infty)$,
this shows that
$u\in\asubdiffext{-\infty}$ if and only if $u=\alpha$ and
$u\in\dom\fstar$, proving the claim.
\end{proof}

It can be checked that
Examples~\ref{ex:standard-abs-val-subgrad-cont}
and~\ref{ex:subgrad-log1+ex}
are consistent with
\indexg{subdifferentials of extensions!one dimension@in one dimension|)}%
\Cref{pr:subdif-in-1d},
and
Example~\ref{ex:entropy-1d}
with
\Cref{pr:adsubgrad-in-1d}.

\chapter{More properties of astral subdifferentials}

In this chapter,
we continue to explore fundamental properties of astral primal and dual
subdifferentials.
We show that every convex function on $\Rn$ is astral
dual subdifferentiable at every point, even at points outside its
effective domain.
Indeed, this property exactly characterizes the convexity of a
function.
We further show how strict convexity is characterized by the disjointness
of astral dual subdifferentials.
We also study monotonicity properties of astral subdifferentials,
and show that astral primal and dual subdifferentials at every point
must be convex sets, even for nonconvex functions.
Finally, we characterize astral subdifferentials of indicator
functions.

\section{Nonemptiness and disjointness of dual subdifferentials}
\label{sec:dual-subdif-not-empty}

In standard convex analysis, a proper convex function has a nonempty
subdifferential at every point in the relative interior of its
effective domain
(\Cref{roc:thm23.4}\ref{roc:thm23.4:a}).
Nevertheless,
it is possible for the function to have no
subgradients at other points, as indeed will be the case at all points
outside its effective domain, and possibly at some or all of its
relative boundary points.

\indexg{subdifferentials, astral dual!nonemptiness of|(}%
In contrast, as we prove next, every convex function on $\Rn$ 
has an astral dual subgradient at \emph{all}
points in $\Rn$, even those outside its effective domain.
Thus, for convex functions, astral dual subdifferentials can never be
empty at any point.
This result was also proved by
\idxwaggoner\citet[Proposition~3.2]{waggoner25}
under the additional assumption that
$\psi(\uu)\in\R$.

\begin{theorem} \label{thm:adsubdiff-nonempty}
  Let $\psi:\Rn\rightarrow\Rext$ be convex, and let $\uu\in\Rn$.
  Then $\psi$ has an astral dual subgradient at $\uu$;
  that is, $\adsubdifpsi{\uu}\neq\emptyset$.
\end{theorem}

\begin{proof}
If either $\psi\equiv+\infty$ or
$\psi(\uu)=-\infty$, then actually $\adsubdifpsi{\uu}=\extspace$,
which is obviously nonempty,
so we assume henceforth that $\psi\not\equiv+\infty$
and that $\psi(\uu)>-\infty$.

The proof will follow from the separation result for the epigraph
of a convex function given in
\Cref{lem:sep-W-epipsi}.
As such, let
\[
   W
   =
   \bigBraces{\rpair{\uu}{v} :\: v\in\R,\, v  < \psi(\uu)}.
\]
This set is convex, being either a line or an open halfline.
It is also disjoint from $\epi\psi$.
Therefore, by \Cref{lem:sep-W-epipsi},
applied to $\psi$ and $W$,
there exists $\xbar\in\extspace$ such that
\[
  \psi(\uu')
  \geq
  v + \xbar\cdot(\uu'-\uu)
\]
for all $\uu'\in\Rn$ and all $v \in\regParens{-\infty,\psi(\uu)}$.
Hence,
for all $\uu'\in\Rn$,
\begin{align*}
  \psi(\uu')
  &\geq
  \sup \bigBraces{
    \xbar\cdot(\uu'-\uu) + v
    :\:
    v\in \parens{-\infty, \psi(\uu)}
  }
  \\
  &=
  \xbar\cdot(\uu'-\uu)
  \plusd
  \sup \bigBraces{
    v
    :\:
    v\in \parens{-\infty, \psi(\uu)}
  }
  \\
  &=
  \xbar\cdot(\uu'-\uu)
  \plusd
  \psi(\uu),
\end{align*}
where the first equality is by
\Cref{pr:plusd-props}(\ref{pr:plusd-props:d-gen}).
Thus, $\xbar\in\adsubdifpsi{\uu}$.
\end{proof}

As noted earlier,
if $\uu\not\in\dom\psi$, then $\psi$ has no standard subgradient at
$\uu$, but nevertheless has an astral dual subgradient, by
\Cref{thm:adsubdiff-nonempty}.
The next proposition characterizes these dual subgradients:

\begin{proposition}  \label{pr:dual-subgrad-outside-dom}
  Let $\psi:\Rn\rightarrow\Rext$,
  let $\uu\in\Rn$, and suppose $\psi(\uu)=+\infty$.
  Let $\xbar\in\extspace$.
  Then
  $\xbar\in\adsubdifpsi{\uu}$
  if and only if
  for all $\ww\in\dom\psi$,
  $\xbar\cdot (\ww-\uu)=-\infty$.
\end{proposition}

\begin{proof}
By the definition of astral dual subgradient given in
\eqref{eqn:psi-subgrad:3},
$\xbar\in\adsubdifpsi{\uu}$
if and only if
for all $\ww\in\Rn$,
\[
  -\infty
   =
   -\psi(\uu)
   \geq
   -\psi(\ww)\plusd \xbar\cdot(\ww-\uu).
\]
This equation holds if and only if either of the terms on the
right-hand side is equal to $-\infty$, that is,
if and only if either $\psi(\ww)=+\infty$
or
$\xbar\cdot(\ww-\uu)=-\infty$.
This proves the proposition.
\end{proof}

In particular,
\Cref{pr:dual-subgrad-outside-dom}
shows that astral dual subgradients of $\psi$
at any point $\uu$ that is outside $\dom{\psi}$
correspond to a form of
astral dual separation, specifically, a strong dual
separation, between $\set{\uu}$ and $\dom{\psi}$
(using the characterization in \Cref{thm:ast-str-sep-equiv-no-trans}\ref{thm:ast-str-sep-equiv-no-trans:b}).

\indexg{convex functions (standard)!nonemptiness of dual subdifferentials and|(}%
\Cref{thm:adsubdiff-nonempty}
shows that if $\psi:\Rn\rightarrow\Rext$ is convex,
then it has an astral dual subgradient at every point.
The next theorem shows that the converse holds as well, that is, that
this latter property implies convexity.
The theorem also proves another equivalence, which we
explain next, before stating it formally.

In standard convex analysis,
as discussed in \Cref{sec:prelim:conjugate},
the biconjugate $\psidubs$ of a function $\psi:\Rn\rightarrow\Rext$
is exactly the pointwise supremum
over all affine functions majorized by $\psi$,
that is, over functions of the form
$\uu\mapsto \xx\cdot\uu+\beta$,
for $\uu\in\Rn$, where $\xx\in\Rn$ and $\beta\in\R$.
Furthermore, it is known that
$\psi$ is equal to the supremum over
the majorized affine functions (that is, $\psi=\psidubs$), 
if and only if $\psi$ is convex and closed.
Astral analogues of such results were given in
Theorems~\ref{thm:fdub-sup-afffcns}
and~\ref{thm:psi-geq-psidub}, and in
Chapters~\ref{sec:ast-close-extensions}
and~\ref{sec:biconj-ast-close-fcns}.

The next theorem considers expressing $\psi$ in a similar fashion
using instead functions of the form
$\uu\mapsto \xbar\cdot(\uu-\uu_0)+\beta$,
for $\uu\in\Rn$, where $\xbar\in\extspace$,
$\uu_0\in\Rn$, and $\beta\in\R$.
These can be viewed
as astral dual variants of standard affine functions, which have already appeared
in our dual separation results (most notably \Cref{lem:sep-W-epipsi}).
As we show next,
\Cref{thm:adsubdiff-nonempty}
implies that every convex function is equal to the pointwise supremum
over all majorized functions of this form.
Importantly, no other conditions are required beyond
convexity.
In particular, the epigraph is not required to be closed,
and the function is not required to be proper.

\Cref{thm:conv-equiv-dualsubgrad} is based on a very similar result
proved originally by
\idxwaggoner\citet[Propositions~3.7 and~3.10]{waggoner25},
though differing somewhat in the details.
(Specifically, instead of statement~(\ref{thm:conv-equiv-dualsubgrad:b})
as in \Cref{thm:conv-equiv-dualsubgrad},
his version requires that $\dom\psi$ be convex
and that $\adsubdifpsi{\uu}\neq\emptyset$ for all $\uu\in\Rn$ with
$\psi(\uu)\in\R$.)

\begin{theorem}  \label{thm:conv-equiv-dualsubgrad}
  Let $\psi:\Rn\rightarrow\Rext$.
  For $\xbar\in\extspace$, $\uu_0\in\Rn$, and $\beta\in\R$,
  define the function $\dlinxbuz:\Rn\rightarrow\Rext$ as
  \[\dlinxbuz(\uu)=\xbar\cdot(\uu-\uu_0)+\beta\]
  for $\uu\in\Rn$.
  Then the following are equivalent:
  \begin{letter-compact}
  \item  \label{thm:conv-equiv-dualsubgrad:a}
    $\psi$ is convex.
  \item  \label{thm:conv-equiv-dualsubgrad:b}
    $\adsubdifpsi{\uu}\neq\emptyset$
    for all $\uu\in\Rn$.
  \item  \label{thm:conv-equiv-dualsubgrad:c}
    For all $\uu\in\Rn$,
    \begin{equation}  \label{eqn:thm:conv-equiv-dualsubgrad:1}
      \psi(\uu)
      =
      \sup\,\bigBraces{ \dlinxbuz(\uu) :\:
               \xbar\in\extspace,\, \uu_0\in\Rn,\, \beta\in\R,\,
               \dlinxbuz \leq \psi
                 }.
    \end{equation}
  \end{letter-compact}
\end{theorem}

\begin{proof}
Throughout this proof,
let $\sigma(\uu)$ denote the supremum appearing on the right-hand side
of \eqref{eqn:thm:conv-equiv-dualsubgrad:1}.

\begin{proof-parts}
\pfpart{%
  (\ref{thm:conv-equiv-dualsubgrad:a})
  $\Rightarrow$
  (\ref{thm:conv-equiv-dualsubgrad:b}):
}
This is exactly \Cref{thm:adsubdiff-nonempty}.

\pfpart{%
  (\ref{thm:conv-equiv-dualsubgrad:b})
  $\Rightarrow$
  (\ref{thm:conv-equiv-dualsubgrad:c}):
}
Assume $\adsubdifpsi{\uu}\neq\emptyset$
for all $\uu\in\Rn$.
Since each $\dlinxbuz$ appearing in the
supremum defining $\sigma$ is majorized by $\psi$,
it follows immediately that $\sigma\leq\psi$.

To show the reverse inequality, let $\uu\in\Rn$.
We aim to show $\psi(\uu)\leq\sigma(\uu)$.
This is immediate if $\psi(\uu)=-\infty$, so we assume henceforth that
$\psi(\uu)>-\infty$.

By assumption, $\psi$ has an astral dual subgradient
$\xbar\in\extspace$ at $\uu$.
Let $\beta\in\R$ with $\beta\leq\psi(\uu)$.
Then for all $\uu'\in\Rn$,
\[
  \psi(\uu')
  \geq
  \psi(\uu) \plusd \xbar\cdot(\uu'-\uu)
  \geq
  \xbar\cdot(\uu'-\uu) + \beta
  =
  \dlinxbu(\uu').
\]
The first inequality is
\eqref{eqn:psi-subgrad:3-alt}, which holds since
$\xbar\in\adsubdifpsi{\uu}$.
The second inequality is
because $\psi(\uu)\geq\beta$ and $\beta\in\R$.
Thus,
$\psi\geq\dlinxbu$
so $\dlinxbu$ is included in the supremum
defining $\sigma$.
Therefore,
$\sigma(\uu)\geq \dlinxbu(\uu) = \beta$.
Since this holds for all $\beta\leq \psi(\uu)$, it follows that
$\sigma(\uu)\geq \psi(\uu)$, completing the proof.

\pfpart{%
  (\ref{thm:conv-equiv-dualsubgrad:c})
  $\Rightarrow$
  (\ref{thm:conv-equiv-dualsubgrad:a}):
}
Suppose $\psi=\sigma$.
For all $\xbar\in\extspace$,
the function $\ph{\xbar}$, as defined in \eqref{eq:ph-xbar-defn},
is convex,
by \Crefequiv{thm:h:5}{thm:h:5a0}{thm:h:5b}.
Therefore,
for all $\xbar\in\extspace$, $\uu_0\in\Rn$, and $\beta\in\R$,
the function $\dlinxbuz(\uu)=\ph{\xbar}(\uu-\uu_0)+\beta$
is also convex.
Hence, $\psi$,
being a pointwise supremum over such functions, is
convex as well
(\Cref{roc:thm5.5}).%
\indexg{subdifferentials, astral dual!nonemptiness of|)}%
\indexg{convex functions (standard)!nonemptiness of dual subdifferentials and|)}%
\qedhere
\end{proof-parts}
\end{proof}

\indexg{subdifferentials, astral dual!overlapping|(}%
We next consider what it means for two astral dual subdifferentials
of a convex function $\psi:\Rn\rightarrow\Rext$
to overlap.
The next theorem shows that
if $\adsubdifpsi{\uu}\cap\adsubdifpsi{\vv}\neq\emptyset$
at any two points $\uu,\vv\in\Rn$ on which $\psi$ is finite,
then $\psi$ must behave like an affine function along the segment
joining $\uu$ and $\vv$.
Moreover, the astral dual subdifferential for every point along
this segment (except possibly for the endpoints) is equal to
the intersection
$\adsubdifpsi{\uu}\cap\adsubdifpsi{\vv}$.
The converses of these implications hold as well.

\begin{theorem}   \label{thm:disj-dsubgrad-eq-affine}
  Let $\psi:\Rn\rightarrow\Rext$ be convex, let
  $\uu,\vv\in\Rn$, $\lambda\in(0,1)$,
  and let
  $\ww=(1-\lambda)\uu+\lambda\vv$.  
  Assume $\psi(\uu),\psi(\vv)\in\R$.
  Then the following are equivalent:
  \begin{letter-compact}
  \item   \label{thm:disj-dsubgrad-eq-affine:a}
    $\adsubdifpsi{\uu}\cap\adsubdifpsi{\vv}\neq\emptyset$.
  \item   \label{thm:disj-dsubgrad-eq-affine:b}
    $\adsubdifpsi{\uu}\cap\adsubdifpsi{\vv}=\adsubdifpsi{\ww}$.
  \item   \label{thm:disj-dsubgrad-eq-affine:c}
    $\psi(\ww)\geq (1-\lambda)\psi(\uu) + \lambda\psi(\vv)$.
  \item   \label{thm:disj-dsubgrad-eq-affine:d}
    $\psi(\ww) = (1-\lambda)\psi(\uu) + \lambda\psi(\vv)$.
  \end{letter-compact}
\end{theorem}

\begin{proof}
  ~

\begin{proof-parts}
\pfpart{%
  (\ref{thm:disj-dsubgrad-eq-affine:a})
  $\Rightarrow$
  (\ref{thm:disj-dsubgrad-eq-affine:c}):
}
Suppose there exists a point
$\xbar\in\adsubdifpsi{\uu}\cap\adsubdifpsi{\vv}$.
Since $\xbar\in\adsubdifpsi{\vv}$,
the definition of astral dual subgradient (specifically,
the negation of Eq.~\ref{eqn:psi-subgrad:3}) implies that
\begin{equation}   \label{eq:thm:disj-dsubgrad-eq-affine:4}
  \psi(\vv) \leq \psi(\uu) + \xbar\cdot(\vv-\uu).
\end{equation}
And since $\xbar\in\adsubdifpsi{\uu}$,
\begin{align*}
  \psi(\ww)
  &\geq
  \psi(\uu) + \xbar\cdot(\ww-\uu)
  \\
  &=
  \psi(\uu) + \lambda \xbar\cdot(\vv-\uu)
  \\
  &=
  (1-\lambda)\psi(\uu) +
  \lambda\bigParens{\psi(\uu)+ \xbar\cdot(\vv-\uu)}
  \\
  &\geq
  (1-\lambda)\psi(\uu) + \lambda\psi(\vv).
\end{align*}
The first equality is because
$\ww-\uu=\lambda(\vv-\uu)$.
The final inequality is by \eqref{eq:thm:disj-dsubgrad-eq-affine:4}.

\pfpart{%
  (\ref{thm:disj-dsubgrad-eq-affine:c})
  $\Rightarrow$
  (\ref{thm:disj-dsubgrad-eq-affine:d}):
}
Since $\psi$ is convex and $\uu,\vv\in\dom\psi$,
this follows from
\Cref{pr:stand-cvx-fcn-char}.

\pfpart{%
  (\ref{thm:disj-dsubgrad-eq-affine:d})
  $\Rightarrow$
  (\ref{thm:disj-dsubgrad-eq-affine:b}):
}
Suppose
\begin{equation}   \label{eq:thm:disj-dsubgrad-eq-affine:1}
  \psi(\ww) = (1-\lambda)\psi(\uu) + \lambda\psi(\vv).
\end{equation}
Since $\psi(\uu)$ and $\psi(\vv)$ are both in $\R$, this implies
$\psi(\ww)$ is as well.

We show first that
$\adsubdifpsi{\ww}\subseteq\adsubdifpsi{\uu}$.
Suppose $\xbar\in\adsubdifpsi{\ww}$.
Then
\[
  \psi(\vv)
  \geq
  \psi(\ww) + \xbar\cdot(\vv-\ww)
  =
  (1-\lambda)\psi(\uu)
  + \lambda\psi(\vv)
  + (1-\lambda) \xbar\cdot(\vv-\uu),
\]  
where the equality is by
\eqref{eq:thm:disj-dsubgrad-eq-affine:1}
and because
$\vv-\ww = (1-\lambda)(\vv-\uu)$.
Rearranging and simplifying, this shows that
\begin{equation}   \label{eq:thm:disj-dsubgrad-eq-affine:2}
  \psi(\vv)
  \geq
  \psi(\uu) + \xbar\cdot(\vv-\uu).
\end{equation}

Letting $\uu'\in\Rn$, we show
$\xbar\in\adsubdifpsi{\uu}$ by proving
that $\psi(\uu')\geq \psi(\uu) + \xbar\cdot(\uu'-\uu)$.
We assume $\psi(\uu')<+\infty$ since otherwise
this inequality is immediate.
We have
\begin{align}
  \psi(\uu')
  &\geq
  \psi(\ww) + \xbar\cdot(\uu'-\ww)
  \label{eq:thm:disj-dsubgrad-eq-affine:3}
  \\
  &=
  (1-\lambda)\psi(\uu)+\lambda\psi(\vv)
  + \xbar\cdot(\uu'-\ww)
  \nonumber
  \\
  &\geq
  (1-\lambda)\psi(\uu)
  + \lambda\bigParens{
    \psi(\uu) + \xbar\cdot(\vv-\uu)
  }
  + \xbar\cdot(\uu'-\ww)
  \nonumber
  \\
  \label{eq:thm:disj-dsubgrad-eq-affine:3:final}
  &=
  \psi(\uu)
  + \xbar\cdot\bigParens{\lambda(\vv-\uu) + (\uu'-\ww)}
  =
  \psi(\uu) + \xbar\cdot(\uu'-\uu).
\end{align}
The first inequality is
because $\xbar\in\adsubdifpsi{\ww}$.
The first equality is by
\eqref{eq:thm:disj-dsubgrad-eq-affine:1}.
The second inequality is by
\eqref{eq:thm:disj-dsubgrad-eq-affine:2}.
The second equality is by
\Cref{pr:i:1}, noting that
$\xbar\cdot(\vv-\uu)<+\infty$
by \eqref{eq:thm:disj-dsubgrad-eq-affine:2},
and that
$\xbar\cdot(\uu'-\ww)<+\infty$
by
\eqref{eq:thm:disj-dsubgrad-eq-affine:3},
implying these terms are summable.
Thus,
\eqref{eq:thm:disj-dsubgrad-eq-affine:3:final}
holds for all $\uu'\in\Rn$, so
$\xbar\in\adsubdifpsi{\uu}$.

That also
$\xbar\in\adsubdifpsi{\vv}$ follows by a symmetric argument
(or by
applying the argument above with $\uu$ and $\vv$ swapped,
and $\lambda$ replaced with $1-\lambda$).
Thus, 
$\adsubdifpsi{\ww}\subseteq\adsubdifpsi{\uu}\cap\adsubdifpsi{\vv}$.

For the reverse inclusion, suppose
$\xbar\in\adsubdifpsi{\uu}\cap\adsubdifpsi{\vv}$.
Let $\ww'\in\Rn$.
We aim to prove
$ \psi(\ww') \geq \psi(\ww) + \xbar\cdot(\ww'-\ww) $,
which will imply $\xbar\in\adsubdifpsi{\ww}$.
We assume $\psi(\ww')<+\infty$, since otherwise this inequality is
immediate.

Since $\xbar\in\adsubdifpsi{\uu}\cap\adsubdifpsi{\vv}$, we have
\begin{align*}
  \psi(\ww')
  &\geq
  \psi(\uu) + \xbar\cdot(\ww'-\uu),
  \\
  \psi(\ww')
  &\geq
  \psi(\vv) + \xbar\cdot(\ww'-\vv).
\end{align*}
These imply that
$\xbar\cdot(\ww'-\uu)$ and
$\xbar\cdot(\ww'-\vv)$ are summable since neither can be equal to
$+\infty$.
Multiplying the first line above by $1-\lambda$ and the second
by $\lambda$, and then adding yields
\begin{align*}
  \psi(\ww')
  &\geq
  (1-\lambda)\psi(\uu) + \lambda\psi(\vv)
  + (1-\lambda)\xbar\cdot(\ww'-\uu)
  + \lambda\xbar\cdot(\ww'-\vv)
  \\
  &=
  \psi(\ww) + \xbar\cdot(\ww'-\ww).
\end{align*}
The equality is by
\eqref{eq:thm:disj-dsubgrad-eq-affine:1},
and by
\Cref{pr:i:1}
since $\xbar\cdot(\ww'-\uu)$ and
$\xbar\cdot(\ww'-\vv)$ are summable
(and noting that
$\ww'-\ww=(1-\lambda)(\ww'-\uu)+\lambda(\ww'-\vv)$).

This completes the proof.

\pfpart{%
  (\ref{thm:disj-dsubgrad-eq-affine:b})
  $\Rightarrow$
  (\ref{thm:disj-dsubgrad-eq-affine:a}):
}
By \Cref{thm:adsubdiff-nonempty}, 
$\adsubdifpsi{\ww}$ is not empty, proving the claim.%
\indexg{subdifferentials, astral dual!overlapping|)}%
\qedhere
\end{proof-parts}
\end{proof}

\indexg{subdifferentials, astral dual!disjoint|(}%
\indexg{strictly convex functions!astral dual subdifferentials and|(}%
As a consequence of \Cref{thm:disj-dsubgrad-eq-affine},
a convex function $\psi$ is strictly convex 
over a convex region in its domain
if and only if
all of the astral dual subdifferentials at points in that region are
disjoint from one another.
This was originally proved by
\idxwaggoner\citet[Lemma~3.11]{waggoner25}.

\begin{theorem}   \label{thm:strict-cvx-subdiff-disjoint}
  Let $\psi:\Rn\rightarrow\Rext$ be convex, and let
  $C$ be a convex subset of $\dom\psi$.
  Then the following are equivalent:
  \begin{letter-compact}
  \item   \label{thm:strict-cvx-subdiff-disjoint:a}
    $\psi$ is strictly convex on $C$.
  \item   \label{thm:strict-cvx-subdiff-disjoint:b}
    For all $\uu,\vv\in C$, if $\uu\neq\vv$ then
    $\adsubdifpsi{\uu}\cap\adsubdifpsi{\vv}=\emptyset$.
  \end{letter-compact}
\end{theorem}

\begin{proof}
  ~

\begin{proof-parts}
\pfpart{%
  (\ref{thm:strict-cvx-subdiff-disjoint:a})
  $\Rightarrow$
  (\ref{thm:strict-cvx-subdiff-disjoint:b}):
}
We prove the contrapositive.
Suppose $\adsubdifpsi{\uu}\cap\adsubdifpsi{\vv}\neq\emptyset$
for some distinct points $\uu,\vv\in C$.
Let $\ww=(\uu+\vv)/2$, which is in $C$, being a convex set,
implying $\uu,\vv,\ww$ are all also in $\dom\psi$.
We claim that
$\psi(\ww)\geq(\psi(\uu)+\psi(\vv))/2$.
This inequality is immediate if either $\psi(\uu)$ or $\psi(\vv)$ is
$-\infty$.
Otherwise, if both
$\psi(\uu)$ and $\psi(\vv)$ are in $\R$,
then the inequality follows from
\Cref{thm:disj-dsubgrad-eq-affine}(\ref{thm:disj-dsubgrad-eq-affine:a},\ref{thm:disj-dsubgrad-eq-affine:c}).
Thus, $\psi$ is not strictly convex on $C$.

\pfpart{%
  (\ref{thm:strict-cvx-subdiff-disjoint:b})
  $\Rightarrow$
  (\ref{thm:strict-cvx-subdiff-disjoint:a}):
}
We again prove the contrapositive.
Suppose $\psi$ is not strictly convex on $C$,
though convex nonetheless.
In light of
\Cref{pr:stand-cvx-fcn-char},
this means that there exist distinct points $\uu,\vv\in C$ and
$\lambda\in(0,1)$ such that
\begin{equation}   \label{eq:thm:strict-cvx-subdiff-disjoint:3}
  \psi(\ww)=(1-\lambda)\psi(\uu)+\lambda\psi(\vv),
\end{equation}
where $\ww=(1-\lambda)\uu+\lambda\vv$.
Note that $\ww$ is distinct from both $\uu$ and $\vv$,
and that $\ww$ is in $C$, being convex.
We aim to show that
$\adsubdifpsi{\uu}$,
$\adsubdifpsi{\vv}$ and
$\adsubdifpsi{\ww}$
are not all disjoint from one another.

If $\psi(\uu)=-\infty$, then
\eqref{eq:thm:strict-cvx-subdiff-disjoint:3}
implies that also $\psi(\ww)=-\infty$
so that
$\adsubdifpsi{\uu}=\adsubdifpsi{\ww}=\extspace$,
proving the claim in this case.
Likewise if $\psi(\vv)=-\infty$.
Otherwise, if $\psi(\uu)$ and $\psi(\vv)$ are both in $\R$,
then the claim follows
from
\Cref{thm:disj-dsubgrad-eq-affine}(\ref{thm:disj-dsubgrad-eq-affine:d},\ref{thm:disj-dsubgrad-eq-affine:a}).%
\indexg{subdifferentials, astral dual!disjoint|)}%
\indexg{strictly convex functions!astral dual subdifferentials and|)}%
\qedhere
\end{proof-parts}
\end{proof}

\section{Monotonicity, convexity and closedness of subdifferentials}

\indexg{monotonicity of subgradients|(}%
Standard subdifferentials are said to be \emph{monotone}
in the sense that,
for a closed, proper, convex function $f:\Rn\rightarrow\Rext$,
if $\uu\in\partial f(\xx)$ and $\vv\in\partial f(\yy)$,
then $(\xx-\yy)\cdot(\uu-\vv)\geq 0$
\idxroc\citep[Corollary~31.5.2]{ROC}.
\indexg{subdifferentials, astral dual!monotonicity of|(}%
Astral dual subdifferentials have an analogous property:

\begin{theorem}   \label{thm:ast-dual-subgrad-monotone}
Let $\psi:\Rn\rightarrow\Rext$,
let $\uu,\vv\in\Rn$, and
assume $\psi(\uu)$ and $-\psi(\vv)$ are summable.
Let $\xbar\in\adsubdifpsi{\uu}$
and $\ybar\in\adsubdifpsi{\vv}$.
Then
\begin{equation*}  %
  \xbar\cdot (\uu-\vv)
  \,\geq\,
  \psi(\uu) - \psi(\vv)
  \,\geq\,
  \ybar\cdot(\uu-\vv).
\end{equation*}
In particular, if $\xbar=\ybar$ (meaning
$\xbar\in\adsubdifpsi{\uu}\cap\adsubdifpsi{\vv}$),
then
$\psi(\uu)-\psi(\vv)=\xbar\cdot(\uu-\vv)$.
\end{theorem}

\begin{proof}
We have  
\begin{align*}
  \psi(\vv) - \psi(\uu)
  =
  -\psi(\uu) \plusu \psi(\vv)
  &\geq
  \xbar\cdot(\vv - \uu),
  \\
  \psi(\uu) - \psi(\vv)
  =
  -\psi(\vv) \plusu \psi(\uu)
  &\geq
  \ybar\cdot(\uu - \vv).
\end{align*}
The equalities are because $\psi(\uu)$ and $-\psi(\vv)$ are
summable.
The first inequality is by \eqref{eqn:psi-subgrad:3-plusu}
since $\xbar\in\adsubdifpsi{\uu}$,
and the second likewise since $\ybar\in\adsubdifpsi{\vv}$.
Negating the first line and combining then proves the claim.
\end{proof}

\Cref{thm:ast-dual-subgrad-monotone}'s conclusion that
$\xbar\cdot(\uu-\vv) \geq \ybar\cdot(\uu-\vv)$
is not true in general
without the summability assumption, as the next example shows:

\begin{example}
Let $\psi:\R^2\rightarrow\Rext$ be the indicator
of the singleton $\set{-\ee_1}$,
$\psi=\indf{\{-\ee_1\}}$.
Let $\uu=\ee_2$, $\xbar=\limray{\ee_1}$,
$\vv=\zero$, and $\ybar=\limray{\ee_1}\plusl\ee_2$.
Then $\psi(\uu)=\psi(\vv)=+\infty$.
Using the definition of astral dual subgradient,
it can be checked that
$\xbar\in\adsubdifpsi{\uu}$
and
$\ybar\in\adsubdifpsi{\vv}$.
However,
$\xbar\cdot (\uu-\vv) =0 < 1 = \ybar\cdot(\uu-\vv)$.%
\indexg{subdifferentials, astral dual!monotonicity of|)}%
\end{example}

\indexg{subdifferentials, astral (primal)!monotonicity of|(}%
\Cref{thm:ast-dual-subgrad-monotone}
yields as corollary a corresponding property for astral primal
subgradients:

\begin{corollary}   \label{cor:ast-subgrad-monotone}
Let $F:\extspace\rightarrow\Rext$ and
$\xbar,\ybar\in\extspace$, and
let $\uu\in\asubdifF{\xbar}$
and $\vv\in\asubdifF{\ybar}$.
Then
\[
  \xbar\cdot (\uu-\vv)
  \,\geq\,
  \Fstar(\uu) - \Fstar(\vv)
  \,\geq\,
  \ybar\cdot(\uu-\vv).
\]
In particular, if $\xbar=\ybar$ (meaning
$\uu,\vv\in\asubdifF{\xbar}$),
then
$\Fstar(\uu)-\Fstar(\vv)=\xbar\cdot(\uu-\vv)$.
\end{corollary}

\begin{proof}
By
\Cref{thm:fenchel-subgrad}(\ref{thm:fenchel-subgrad:b},\ref{thm:fenchel-subgrad:b-dual}),
$\xbar\in\adsubdifFstar{\uu}$
and $\ybar\in\adsubdifFstar{\vv}$,
and by
\Cref{pr:subgrad-imp-in-cldom}(\ref{pr:subgrad-imp-in-cldom:a}),
$\Fstar(\uu)$ and $\Fstar(\vv)$ are both in $\R$.
Applying \Cref{thm:ast-dual-subgrad-monotone}
(with $\psi=\Fstar$)
therefore proves the claim.%
\indexg{subdifferentials, astral (primal)!monotonicity of|)}%
\indexg{monotonicity of subgradients|)}%
\end{proof}

\indexg{subdifferentials, astral dual!convexity of|(}%
\indexg{subdifferentials, astral dual!closedness of|(}%
The
astral dual subdifferential at any point $\uu\in\Rn$
is a closed and convex set, even for
functions that are not convex:

\begin{theorem}  \label{thm:ast-dual-subdif-is-convex}
  Let $\psi:\Rn\rightarrow\Rext$, and let $\uu\in\Rn$.
  Then $\adsubdifpsi{\uu}$ is convex and closed (in~$\extspace$).
\end{theorem}

\begin{proof}
We have
\begin{align*}
  \adsubdifpsi{\uu}
  &=
  \quad\;
  \bigcap_{\mathclap{\uu'\in\Rn}}
  \quad\;
  \bigBraces{
    \xbar\in\extspace
    :\:
    \xbar\cdot(\uu'-\uu) \leq -\psi(\uu)\plusu\psi(\uu')
  }
  \\[\smallskipamount]
  &=
  \quad\;
  \bigcap_{
    \mathclap{\substack{
        \uu'\in\Rn,\, \beta\in\R:
    \\
        \beta\geq -\psi(\uu)\plusu\psi(\uu')
    }}}
  \quad\;
  \bigBraces{
    \xbar\in\extspace
    :\:
    \xbar\cdot(\uu'-\uu) \leq \beta
  },
\end{align*}
where
the first equality is by definition of astral dual subdifferential as
given in \eqref{eqn:psi-subgrad:3-plusu}.
Each set in braces on the last line is a closed astral halfspace
(if $\uu'\ne\uu$) or $\eRn$ (if $\uu'=\uu$),
and thus in both cases
closed and convex (\Cref{pr:e1}\ref{pr:e1:c}).
Therefore, their intersection, $\adsubdifpsi{\uu}$,
is as well (\Cref{pr:e1}\ref{pr:e1:b}).%
\indexg{subdifferentials, astral dual!convexity of|)}%
\indexg{subdifferentials, astral dual!closedness of|)}%
\end{proof}

\indexg{subdifferentials, astral (primal)!convexity of|(}%
The
astral primal subdifferential at any point $\xbar\in\eRn$
is also always convex, even for astral functions that are not convex:

\begin{theorem}   \label{thm:asubdifF-convex}
  Let $F:\extspace\rightarrow\Rext$ and let $\xbar\in\extspace$.
  Then $\asubdifF{\xbar}$ is convex.
\end{theorem}

\begin{proof}
We assume $\asubdifF{\xbar}$ is not empty since otherwise the claim is
immediate.
Let $\uu,\vv\in\asubdifF{\xbar}$,
let $\lambda\in(0,1)$, and let
$\ww=(1-\lambda)\uu+\lambda\vv$.
To prove convexity,
we aim to show that $\ww$ is also in $\asubdifF{\xbar}$.

Since $\uu,\vv\in\asubdifF{\xbar}$,
\Cref{thm:subgrad-then-lsc} implies $(\lsc F)(\xbar)=F(\xbar)$,
\Cref{pr:subgrad-imp-in-cldom}(\ref{pr:subgrad-imp-in-cldom:a})
implies $\Fstar(\uu)$ and $\Fstar(\vv)$ are in $\R$,
and \Cref{cor:ast-subgrad-monotone}
implies that
\begin{equation}   \label{eq:thm:asubdifF-convex:2}
  \Fstar(\uu) - \Fstar(\vv) = \xbar\cdot(\uu-\vv).
\end{equation}
In addition, by
\Cref{thm:fenchel-subgrad}(\ref{thm:fenchel-subgrad:b},\ref{thm:fenchel-subgrad:b-dual}),
$\xbar\in\adsubdifFstar{\uu}\cap\adsubdifFstar{\vv}$,
so
\begin{equation}   \label{eq:thm:asubdifF-convex:3}
  \Fstar(\ww)=(1-\lambda)\Fstar(\uu)+\lambda\Fstar(\vv)
\end{equation}
by
\Cref{thm:disj-dsubgrad-eq-affine}(\ref{thm:disj-dsubgrad-eq-affine:a},\ref{thm:disj-dsubgrad-eq-affine:d})
(applied with $\psi=\Fstar$).
Thus, $\Fstar(\ww)\in\R$.

By
\Cref{thm:subgrad-eq-fenchel}(\ref{thm:subgrad-eq-fenchel:a},\ref{thm:subgrad-eq-fenchel:c}),
since $\uu\in\asubdifF{\xbar}$,
there exists a sequence $\seq{\xbar_t}$ in $\extspace$ with
$F(\xbar_t)\in\R$ and $\xbar_t\cdot\uu\in\R$ for all $t$,
and such that
$\xbar_t\rightarrow\xbar$
and
\begin{equation}   \label{eq:thm:asubdifF-convex:1}
  \xbar_t\cdot\uu - F(\xbar_t) \rightarrow \Fstar(\uu).
\end{equation}
This implies that
\begin{align*}
  -F(\xbar_t) + \xbar_t\cdot\ww
  &=
  -F(\xbar_t)
  +
  \xbar_t\cdot\bigParens{\uu - \lambda(\uu - \vv)}
  \\
  &=
  \bigParens{
    -F(\xbar_t)
    +
    \xbar_t\cdot\uu
  }
  - \lambda\xbar_t\cdot(\uu - \vv)
  \\
  &\rightarrow
  \Fstar(\uu) - \lambda\xbar\cdot(\uu-\vv)
  \\
  &=
  \Fstar(\uu) - \lambda\bigParens{\Fstar(\uu) - \Fstar(\vv)}
  =
  \Fstar(\ww).
\end{align*}
The second equality is by
\Cref{pr:i:1}
since $\xbar_t\cdot\uu\in\R$, implying summability.
The convergence is by
\eqref{eq:thm:asubdifF-convex:1}
and because $\xbar_t\cdot(\uu-\vv)\rightarrow\xbar\cdot(\uu-\vv)$
(\Cref{thm:i:1}\ref{thm:i:1c}),
also using that
$\Fstar(\uu)\in\R$, so the limit of the sum is the sum of the limits
(\Cref{prop:lim:eR}\ref{i:lim:eR:sum}).
The third and fourth equalities are by
Eqs.~(\ref{eq:thm:asubdifF-convex:2})
and~(\ref{eq:thm:asubdifF-convex:3}),
respectively.

Applying
\Cref{thm:subgrad-eq-fenchel}(\ref{thm:subgrad-eq-fenchel:d},\ref{thm:subgrad-eq-fenchel:a})
to the same sequence $\seq{\xbar_t}$,
this implies that ${\ww\in\asubdifF{\xbar}}$, completing the proof.%
\indexg{subdifferentials, astral (primal)!convexity of|)}%
\end{proof}

\indexg{subdifferentials, astral (primal)!closedness of|(}%
Unlike astral dual subdifferentials, the astral primal subdifferential
at a given point~$\xbar$ is not necessarily closed (in $\Rn$).
Here is an example:

\begin{example}
Define $F:\Rext\rightarrow\Rext$,
for $\barx\in\Rext$, by
\[
   F(\barx)
   =
   \begin{cases}
     -\infty  &  \text{if $\barx=-\infty$,}  \\
     0        &  \text{if $\barx\in\Rneg$,}  \\
     +\infty  &  \text{if $\barx>0$.}  \\
   \end{cases}
\]
This function is convex by
\Cref{pr:1d-cvx}.
It can be calculated that $\Fstar=\indf{>0}$, the indicator
function of $\Rstrictpos$.
In particular, this implies that
$\asubdifF{0}\subseteq\Rstrictpos$ by
\Cref{pr:subgrad-imp-in-cldom}(\ref{pr:subgrad-imp-in-cldom:a}).
On the other hand, if $u\in\Rstrictpos$, then
$\Fstar(u)=0=0\cdot u - F(0)$
so $u\in\asubdifF{0}$ by
\Cref{thm:fenchel-subgrad}(\ref{thm:fenchel-subgrad:a},\ref{thm:fenchel-subgrad:b}).
Thus, $\asubdifF{0}=\Rstrictpos$, which is not closed in $\R$.%
\indexg{subdifferentials, astral (primal)!closedness of|)}%
\end{example}

\section{Subdifferentials of indicator functions and normal cones}
\label{sec:subdif-ind-fcn-norm-cone}

We look next at astral subdifferentials of indicator functions.
\indexg{normal cone (standard)|(}%
\indexg{subdifferentials, standard!indicator function@of indicator function|(}%
\indexg{indicator functions (standard)!standard subgradients of|(}%
For the indicator function~$\inds$ of a nonempty set $S\subseteq\Rn$,
the standard subdifferential $\partial \inds(\xx)$ at a point
$\xx\in S$ is simply equal to a set called the \emph{normal cone}
to $S$ at $\xx$, denoted $\ncone{\xx}{S}$, defined as
\begin{align}
  \ncone{\xx}{S}
  &=
  \bigBraces{
    \uu\in\Rn :\:
    (\yy-\xx)\cdot\uu \leq 0
    \text{ for all }
    \yy\in S
  }
  \label{eqn:ncone-dfn-1}
  \\
  &=
  \bigBraces{
    \uu\in\Rn :\:
    \yy\cdot\uu \leq \xx\cdot\uu
    \text{ for all }
    \yy\in S
  }.
  \label{eqn:ncone-dfn-2}
\end{align}
Geometrically, $\ncone{\xx}{S}$ is the set of vectors
$\uu\in\Rn$ that make an angle of at least $90^{\circ}$ with $\yy-\xx$
for every $\yy\in S$, or equivalently, for which the scalar projection
$\zz\mapsto\zz\cdot\uu$ is maximized over $S$ at $\xx$.
That $\partial \inds(\xx)=\ncone{\xx}{S}$ (for $\xx\in S$)
follows straightforwardly from definitions.
\indexg{normal cone (standard)|)}%
\indexg{subdifferentials, standard!indicator function@of indicator function|)}%
\indexg{indicator functions (standard)!standard subgradients of|)}%

As with other operations we have earlier considered
(such as polarity in \Cref{sec:ast-pol-cones}
and orthogonal complement in \Cref{sec:ortho-compl}),
we can generalize normal cones to astral space in two different ways.
\indexg{normal cone (primal)|(}%
First, we can simply extend the definition in
\eqref{eqn:ncone-dfn-2} to sets in $\extspace$:

\begin{definition}   \label{dfn:ncone-ext}
\indexg{normal cone (primal)!defined|(}%
  Let $S\subseteq\extspace$, and let $\xbar\in\extspace$.
  Then the \emph{(primal) normal cone} to $S$ at~$\xbar$, denoted
  $\ncone{\xbar}{S}$, is the set
  \begin{equation}  \label{eq:dfn:ncone-ext}
\indexm{n s 300}{$\ncone{\xbar}{S}$}{(primal) normal cone}%
     \ncone{\xbar}{S}
     =
     \bigBraces{
       \uu\in\Rn :\:
       \ybar\cdot\uu \leq \xbar\cdot\uu
       \text{ for all }
       \ybar\in S
     }.%
\indexg{normal cone (primal)!defined|)}%
  \end{equation}
\end{definition}
Note that if $S\subseteq\Rn$, then this definition is the same as the
standard one.
Note also that other authors would more typically define the normal cone
only at points that are actually in the set $S$.
In our definition, we nonetheless allow $\xbar$ to be outside $S$,
although we are generally still thinking of $\xbar$ as a point in $S$.

For a set $S\subseteq\extspace$ and point $\xbar\in\extspace$,
$\ncone{\xbar}{S}$ is a set in
\indexg{normal cone (primal)|)}%
$\Rn$.
\indexg{normal cone, dual|(}%
Alternatively, 
for a set and point in $\Rn$, we can consider a dual notion of normal
cone that is instead a subset of~$\extspace$, and which
straightforwardly generalizes the definition given in
\eqref{eqn:ncone-dfn-1}:

\begin{definition}   \label{dfn:ast-ncone}
\indexg{normal cone, dual!defined|(}%
  Let $S\subseteq\Rn$, and let $\uu\in\Rn$.
  Then the \emph{dual normal cone} to $S$ at $\uu$, denoted
  $\ancone{\uu}{S}$, is the set
  \begin{equation}  \label{eq:dfn:ast-ncone}
\indexm{n s 700}{$\ancone{\uu}{S}$}{dual normal cone}%
     \ancone{\uu}{S}
     =
     \bigBraces{
       \xbar\in\extspace :\:
       \xbar\cdot(\vv-\uu) \leq 0
       \text{ for all }
       \vv\in S
     }.
\indexg{normal cone, dual!defined|)}%
\indexg{normal cone, dual|)}%
  \end{equation}
\end{definition}

\indexg{normal cone (primal)|(}%
Here are some facts about primal normal cones:

\begin{proposition}   \label{pr:ncone-props}
  Let $S\subseteq\extspace$, and let $\xbar\in\extspace$.
  Then:
  \begin{letter-compact}
  \item    \label{pr:ncone-props:a}
    $\ncone{\xbar}{S}$ is a cone.
  \item    \label{pr:ncone-props:b}
    If $S\cap\Rn\neq\emptyset$ then
    $\ncone{\xbar}{S}$ is a convex cone.
  \end{letter-compact}
\end{proposition}

\begin{proof}
~

\begin{proof-parts}
\pfpart{Part~(\ref{pr:ncone-props:a}):}
This is straightforward from the definition in
\eqref{eq:dfn:ncone-ext}.

\pfpart{Part~(\ref{pr:ncone-props:b}):}
Suppose there exists a point $\zz\in S\cap\Rn$.
Let $N=\ncone{\xbar}{S}$, and let $\uu,\vv\in N$.
Since $N$ is a cone,
to show it is also convex,
by \Cref{pr:scc-cone-elts}(\ref{pr:scc-cone-elts:d}),
it suffices to prove that
$\uu+\vv$ is in $N$.
In particular, letting $\ybar\in S$, we have
\begin{equation}   \label{eq:pr:ncone-props:2}
  \ybar\cdot\uu\leq\xbar\cdot\uu
  \;\;\text{ and }\;\;
  \ybar\cdot\vv\leq\xbar\cdot\vv,
\end{equation}
and we aim to prove that
\begin{equation}   \label{eq:pr:ncone-props:1}
  \ybar\cdot(\uu+\vv)\leq\xbar\cdot(\uu+\vv).
\end{equation}
Further, we can assume
$\xbar\cdot(\uu+\vv)<+\infty$ since otherwise
this last inequality holds trivially.
Since $\zz\in S$ and $\uu,\vv\in N$, we have
$-\infty<\zz\cdot\uu\leq\xbar\cdot\uu$
and
$-\infty<\zz\cdot\vv\leq\xbar\cdot\vv$.
Thus, $\xbar\cdot\uu$ and $\xbar\cdot\vv$ are summable,
so $\xbar\cdot(\uu+\vv)=\xbar\cdot\uu+\xbar\cdot\vv$
by \Cref{pr:i:1}.
Since $\xbar\cdot(\uu+\vv)<+\infty$, it follows further that
$\xbar\cdot\uu$ and $\xbar\cdot\vv$ are both in $\R$.
Combined with
\eqref{eq:pr:ncone-props:2}, this implies
$\ybar\cdot\uu<+\infty$
and
$\ybar\cdot\vv<+\infty$,
so $\ybar\cdot\uu$ and $\ybar\cdot\vv$ are summable, and therefore,
$\ybar\cdot(\uu+\vv)=\ybar\cdot\uu+\ybar\cdot\vv$
(again by \Cref{pr:i:1}).
Adding the inequalities in \eqref{eq:pr:ncone-props:2}
now yields \eqref{eq:pr:ncone-props:1},
proving the claim.
\qedhere
\end{proof-parts}
\end{proof}

Although $\ncone{\xbar}{S}$ is always a cone, it need not be convex if
$S$ is disjoint from $\Rn$, even if $S$ is convex.
Here is an example:

\begin{example}
In $\extspac{2}$,
let $\xbar=\limray{(-\ee_1)}\plusl\limray{\ee_2}$,
let $\xbar'=\limray{\ee_2}\plusl\limray{(-\ee_1)}$,
and let $S=\lb{\xbar}{\xbar'}$.
Let $N=\ncone{\xbar}{S}$.
Then $\ee_1\in N$ since $\xbar\cdot\ee_1=\xbar'\cdot\ee_1=-\infty$,
implying $\ybar\cdot\ee_1=-\infty$ for $\ybar\in S$ by
\Cref{pr:seg-simplify}(\ref{pr:seg-simplify:a},\ref{pr:seg-simplify:c}).
Likewise, $\ee_2\in N$.
Nonetheless, the point $\uu=\frac{1}{2}(\ee_1+\ee_2)$ is not in $N$
since
$\xbar'\cdot\uu=+\infty$ but
$\xbar\cdot\uu=-\infty$.
Therefore, $N$ is not convex.
\end{example}

\indexg{support functions, astral!dual subgradients of|(}%
\indexg{subdifferentials, astral dual!support function@of support function|(}%
The next proposition shows how inclusion in a normal cone relates
to the
astral support function and its astral dual subdifferential:

\begin{proposition}  \label{pr:ncone-ast-sup-fcn}
  Let $S\subseteq\extspace$, 
  let $\uu\in\Rn$, and
  suppose $\xbar\in S$.
  Then the following are equivalent:
  \begin{letter-compact}
  \item    \label{pr:ncone-ast-sup-fcn:a}
    $\uu\in\ncone{\xbar}{S}$.
  \item    \label{pr:ncone-ast-sup-fcn:b}
    $\xbar\cdot\uu = \indaSstar(\uu)$.
  \item    \label{pr:ncone-ast-sup-fcn:c}
    $\xbar\in\adsubdif{\indaSstar}{\uu}$.
  \end{letter-compact}
\end{proposition}

\begin{proof}
~

\begin{proof-parts}
\pfpart{%
  (\ref{pr:ncone-ast-sup-fcn:a})
  $\Rightarrow$
  (\ref{pr:ncone-ast-sup-fcn:b}):
}
Suppose $\uu\in\ncone{\xbar}{S}$.
Then
\begin{equation*}  \label{eq:pr:ncone-ast-sup-fcn:1}
  \xbar\cdot\uu
  \leq
  \sup_{\ybar\in S} \ybar\cdot\uu
  \leq
  \xbar\cdot\uu.
\end{equation*}
The first inequality is because $\xbar\in S$,
and the second is from the definition in \eqref{eq:dfn:ncone-ext}
since $\uu\in\ncone{\xbar}{S}$.
Since the middle expression is equal to
$\indaSstar(\uu)$, this proves
statement~(\ref{pr:ncone-ast-sup-fcn:b}).

\pfpart{%
  (\ref{pr:ncone-ast-sup-fcn:b})
  $\Rightarrow$
  (\ref{pr:ncone-ast-sup-fcn:c}):
}
Suppose
$\xbar\cdot\uu = \indaSstar(\uu)$.
Let $\ww\in\Rn$.
We claim that
\begin{equation}  \label{eq:pr:ncone-ast-sup-fcn:2}
  \xbar\cdot\ww
  \geq
  \xbar\cdot\uu
  \plusd
  \xbar\cdot(\ww-\uu).
\end{equation}
This is immediate if the right-hand side is equal to $-\infty$, and so
also if either $\xbar\cdot\uu$ or $\xbar\cdot(\ww-\uu)$
is equal to $-\infty$.
Otherwise, if neither of these terms is $-\infty$, then they must be
summable, implying that their (downward) sum is equal to
$\xbar\cdot\uu+\xbar\cdot(\ww-\uu)=\xbar\cdot\ww$
by \Cref{pr:i:1}, proving the claim in this case as well.
Thus,
\[
  \indaSstar(\ww)
  \geq
  \xbar\cdot\ww
  \geq
  \xbar\cdot\uu
  \plusd
  \xbar\cdot(\ww-\uu)
  =
  \indaSstar(\uu)
  \plusd
  \xbar\cdot(\ww-\uu).
\]
The first inequality is because $\xbar\in S$
and by definition of astral support function
(Eq.~\ref{eqn:astral-support-fcn-def}).
The second inequality is by
\eqref{eq:pr:ncone-ast-sup-fcn:2}.
The equality is by assumption.
Since this holds for all $\ww\in\Rn$,
this proves $\xbar\in\adsubdif{\indaSstar}{\uu}$.

\pfpart{%
  (\ref{pr:ncone-ast-sup-fcn:c})
  $\Rightarrow$
  (\ref{pr:ncone-ast-sup-fcn:a}):
}
Suppose $\xbar\in\adsubdif{\indaSstar}{\uu}$.
Then
\[
   0
   =
   \indaSstar(\zero)
   \geq
   \indaSstar(\uu)
   \plusd
   \xbar\cdot(-\uu).
\]
The equality is by definition of astral support function
(and since $S$, which includes $\xbar$, cannot be empty).
The inequality is by
the definition of astral dual subgradient.
By \Cref{pr:plusd-props}(\ref{pr:plusd-props:e}), this implies
\[
  \xbar\cdot\uu
  \geq
  \indaSstar(\uu)
  \plusd
  0
  =
  \sup_{\ybar\in S} \ybar\cdot\uu,
\]
proving $\uu\in\ncone{\xbar}{S}$
(by Eq.~\ref{eq:dfn:ncone-ext}).%
\indexg{normal cone (primal)|)}%
\indexg{support functions, astral!dual subgradients of|)}%
\indexg{subdifferentials, astral dual!support function@of support function|)}%
\qedhere
\end{proof-parts}
\end{proof}

\indexg{normal cone, dual|(}%
Here are some simple facts about dual normal cones:

\begin{proposition}   \label{pr:ancone-props}
  Let $S\subseteq\Rn$, and let $\uu\in\Rn$.
  Then:
  \begin{letter-compact}
  \item    \label{pr:ancone-props:a}
    $\ancone{\uu}{S}$ is a closed convex astral cone.
  \item    \label{pr:ancone-props:b}
    $\ancone{\uu}{S}=\apol{\bigParens{\cone(S-\uu)}}$.
  \end{letter-compact}
\end{proposition}

\begin{proof}
~

\begin{proof-parts}
\pfpart{Part~(\ref{pr:ancone-props:b}):}
This follows directly from comparison of the definition for
$\ancone{\uu}{S}$ given in
\eqref{eq:dfn:ast-ncone}
with the expression for 
$\apol{\regParens{\cone(S-\uu)}}$ obtained by applying
\Cref{pr:ast-pol-props}(\ref{pr:ast-pol-props:coneSpol}).

\pfpart{Part~(\ref{pr:ancone-props:a}):}
This follows from part~(\ref{pr:ancone-props:b})
combined with
\Cref{pr:ast-pol-props}(\ref{pr:ast-pol-props:c}).%
\indexg{normal cone, dual|)}%
\qedhere
\end{proof-parts}
\end{proof}

\indexg{subdifferentials, astral (primal)!indicator function@of indicator function|(}%
\indexg{indicator functions (astral)!subgradients of|(}%
We can now compute the astral subdifferential of the indicator
function of any set $S\subseteq\extspace$ at a point
$\xbar\in\extspace$, stated both in terms of the
normal cone $\ncone{\xbar}{S}$ and the support function $\indaSstar$:

\begin{theorem}  \label{thm:subdif-ind-fcn-gen}
  Let $S\subseteq\extspace$,
  $\xbar\in\extspace$, and $\uu\in\Rn$.
  Then the following are equivalent:
  \begin{letter-compact}
  \item  \label{thm:subdif-ind-fcn-gen:a}
    $\uu\in\asubdifindaS{\xbar}$.
  \item  \label{thm:subdif-ind-fcn-gen:b}
    $\xbar\in S$,
    $\xbar\cdot\uu=\indfastar{S}(\uu)$,
    and
    $\indfastar{S}(\uu)\in\R$.
  \item  \label{thm:subdif-ind-fcn-gen:c}
    $\xbar\in S$,
    $\uu\in\ncone{\xbar}{S}$,
    and
    $\xbar\cdot\uu\in\R$.
  \end{letter-compact}
\end{theorem}

\begin{proof}
~

\begin{proof-parts}
\pfpart{%
  (\ref{thm:subdif-ind-fcn-gen:a})
  $\Rightarrow$
  (\ref{thm:subdif-ind-fcn-gen:b}):
}
Suppose $\uu\in\asubdifindaS{\xbar}$.
Then by
\Cref{pr:subgrad-imp-in-cldom}(\ref{pr:subgrad-imp-in-cldom:a}),
$\indfastar{S}(\uu)\in\R$
and
$\xbar\in\cldom{\indaS}=\Sbar$.
We further claim that $\xbar\in S$:
Since $\xbar\in\Sbar$, there exists a sequence
$\seq{\xbar_t}$ in $S$ with $\xbar_t\rightarrow\xbar$, implying
\[
  0
  \leq
  \indaS(\xbar)
  =
  (\lsc\indaS)(\xbar)
  \leq
  \lim \indaS(\xbar_t)
  =
  0.
\]
The first equality is 
by \Cref{thm:subgrad-then-lsc} since
$\uu\in\asubdifindaS{\xbar}$, and the second inequality is by definition of
lower semicontinuous hull
(Eq.~\ref{eq:lsc:liminf:X:prelims}).
Thus, $\indaS(\xbar)=0$ so $\xbar\in S$.

By
\Cref{thm:fenchel-subgrad}(\ref{thm:fenchel-subgrad:b},\ref{thm:fenchel-subgrad:c}),
it now follows that
$0=\indaS(\xbar)=\xbar\cdot\uu - \indfastar{S}(\uu)$,
and so that
$\xbar\cdot\uu = \indfastar{S}(\uu)$,
completing the proof.

\pfpart{%
  (\ref{thm:subdif-ind-fcn-gen:b})
  $\Rightarrow$
  (\ref{thm:subdif-ind-fcn-gen:a}):
}
Suppose 
$\xbar\in S$,
$\xbar\cdot\uu=\indfastar{S}(\uu)$,
and
$\indfastar{S}(\uu)\in\R$.
Then
$\indaS(\xbar)=0$ so
$\indfastar{S}(\uu)=\xbar\cdot\uu-\indaS(\xbar)$, implying
$\uu\in\asubdifindaS{\xbar}$
by
\Cref{thm:fenchel-subgrad}(\ref{thm:fenchel-subgrad:a},\ref{thm:fenchel-subgrad:b}).

\pfpart{%
  (\ref{thm:subdif-ind-fcn-gen:b})
  $\Leftrightarrow$
  (\ref{thm:subdif-ind-fcn-gen:c}):
}
This follows directly from
\Cref{pr:ncone-ast-sup-fcn}(\ref{pr:ncone-ast-sup-fcn:a},\ref{pr:ncone-ast-sup-fcn:b}).%
\indexg{subdifferentials, astral (primal)!indicator function@of indicator function|)}%
\indexg{indicator functions (astral)!subgradients of|)}%
\qedhere
\end{proof-parts}
\end{proof}

\indexg{subdifferentials of extensions!indicator function@of indicator function|(}%
\indexg{indicator functions (standard)!subgradients of extension|(}%
Specializing to extensions of indicators of sets in $\Rn$ then yields:

\begin{corollary}  \label{cor:subdif-ind-fcn-rn}
  Let $S\subseteq\Rn$,
  $\xbar\in\extspace$, and $\uu\in\Rn$.
  Then the following are equivalent:
  \begin{letter-compact}
  \item  \label{cor:subdif-ind-fcn-rn:a}
    $\uu\in\asubdifindsext{\xbar}$.
  \item  \label{cor:subdif-ind-fcn-rn:b}
    $\xbar\in\Sbar$,
    $\xbar\cdot\uu=\indstars(\uu)$,
    and 
    $\indstars(\uu)\in\R$.
  \item  \label{cor:subdif-ind-fcn-rn:c}
    $\xbar\in \Sbar$,
    $\uu\in\ncone{\xbar}{\mkern-0.5mu\Sbar}$,
    and
    $\xbar\cdot\uu\in\R$.
  \end{letter-compact}
\end{corollary}

\begin{proof}
We have
$\indfa{\mkern-0.5mu\Sbar}=\indsext$ by
\Cref{pr:inds-ext},
implying $\indfastar{\mkern-0.5mu\Sbar}=\indsextstar=\indstars$
by \Cref{pr:fextstar-is-fstar}.
The claim therefore follows directly from
\Cref{thm:subdif-ind-fcn-gen}
(applied to~$\Sbar$).%
\indexg{subdifferentials of extensions!indicator function@of indicator function|)}%
\indexg{indicator functions (standard)!subgradients of extension|)}%
\end{proof}

\indexg{subdifferentials, astral dual!indicator function@of indicator function|(}%
\indexg{indicator functions (standard)!dual subgradients of|(}%
For astral dual subdifferentials of indicator functions, we obtain:

\begin{proposition}  \label{pr:adsubdif-ind-fcn}
  Let $S\subseteq\Rn$, and suppose $\uu\in S$.
  Then
  \[
     \adsubdif{\inds}{\uu}
     =
     \ancone{\uu}{S}.
  \]
\end{proposition}

\begin{proof}
This is straightforward from the definition of astral dual
subdifferential
(Eq.~\ref{eqn:psi-subgrad:3-alt})
and the definition of dual normal cone
(Eq.~\ref{eq:dfn:ast-ncone}).
\end{proof}

This proposition only characterizes
$\adsubdif{\inds}{\uu}$ at points $\uu\in S$.
For points outside of $S$, this dual subdifferential can instead be
characterized as a special case of
\Cref{pr:dual-subgrad-outside-dom}.%
\indexg{subdifferentials, astral dual!indicator function@of indicator function|)}%
\indexg{indicator functions (standard)!dual subgradients of|)}%

\chapter{Calculus rules for astral subgradients}
\label{sec:calc-subgrads}

We next develop rules for computing astral subgradients analogous to
the standard differential calculus and to rules for standard
subgradients, which we reviewed in \Cref{sec:prelim:subgrads}.
We focus particularly on the astral subdifferential
$\asubdifplain{\ef}$ of the extension of a function
$f:\Rn\rightarrow\Rext$.
We derive the astral subdifferential for a scalar multiple of a function,
the sum
of two functions, the composition of a function with a linear
map, and for various other functional operations.
The resulting rules can be used,
for instance, to derive optimality conditions, as
we will do in
Chapters~\ref{sec:fenchel-duality} and~\ref{sec:KKT}.

\section{Scalar multiple}

\indexg{scalar multiples (astral)!subgradients of extension|(}%
\indexg{subdifferentials of extensions!scalar multiple@of scalar multiple|(}%
Let $f:\Rn\rightarrow\Rext$, and
suppose $h=\lambda f$ for some $\lambda\in\Rstrictpos$.
Then the astral subgradients of $\hext$ can be calculated from
the astral subgradients of $\fext$ simply by multiplying by $\lambda$,
as in standard calculus:

\begin{proposition}  \label{pr:subgrad-scal-mult}
  Let $f:\Rn\rightarrow\Rext$, let $\xbar\in\extspace$,
  and let $\lambda\in\Rstrictpos$.
  Then
  \begin{equation}
  \label{eq:pr:subgrad-scal-mult:2}
    \basubdif{\lamfext}{\xbar} = \lambda[\basubdiffext{\xbar}].
  \end{equation}
  Moreover, if $\basubdiffext{\xbar}\ne\emptyset$,
  then \eqref{eq:pr:subgrad-scal-mult:2} holds for $\lambda\in\Rpos$.
\end{proposition}

\begin{proof}
First consider the case $\lambda\in\Rstrictpos$.
Let $h=\lambda f$,
let $\uu\in\Rn$, and let $\ww=\lambda\uu$.
We aim to prove that
$\uu\in\basubdiffext{\xbar}$
if and only if
$\ww\in\basubdifhext{\xbar}$.

Let $\fminusu$ and $\hminusw$ be linear tilts of $f$ and $h$.
Then for $\xx\in\Rn$,
\[
  \hminusw(\xx)
  =
  h(\xx) - \xx\cdot \ww
  =
  \lambda f(\xx) - \xx\cdot (\lambda \uu)
  =
  \lambda \fminusu(\xx).
\]
Consequently,
\begin{equation}    \label{eq:pr:subgrad-scal-mult:3}
  \inf \hminusw = \lambda \inf \fminusu,
\end{equation}
and
\begin{equation}    \label{eq:pr:subgrad-scal-mult:4}
  \hminuswext(\xbar)=\lambda \fminusuext(\xbar)
\end{equation}
by \Cref{pr:scal-mult-ext}.
Also,
by \Cref{pr:fminusu-props}(\ref{pr:fminusu-props:d})
and
\eqref{eq:pr:subgrad-scal-mult:3},
\begin{equation}    \label{eq:pr:subgrad-scal-mult:5}
  \hstar(\ww)
  =
  -\inf\hminusw
  =
  -\lambda \inf\fminusu
  =
  \lambda \fstar(\uu).
\end{equation}

Combining, we then have that
\begin{alignat*}{3}
  \uu\in\basubdiffext{\xbar}\quad
    {\Leftrightarrow}\quad
    && \fstar(\uu)\in\R
    &  \quad\text{and}\quad
    &  \fminusuext(\xbar)&=-\fstar(\uu)
  \\
    {\Leftrightarrow}\quad
    && \lambda\fstar(\uu)\in\R
    &  \quad\text{and}\quad
    &  \lambda\fminusuext(\xbar)&=-\lambda\fstar(\uu)
  \\
    {\Leftrightarrow}\quad
    && \hstar(\ww)\in\R
    &  \quad\text{and}\quad
    &  \hminuswext(\xbar)&=-\hstar(\uu)
  \\
    {\Leftrightarrow}\quad
    && \ww\in\basubdifhext{\xbar}.
\end{alignat*}
The first and last equivalences are both by
\Cref{pr:subgrad-imp-in-cldom}(\ref{pr:subgrad-imp-in-cldom:c})
and
\Cref{thm:fminus-subgrad-char}(\ref{thm:fminus-subgrad-char:a},\ref{thm:fminus-subgrad-char:d}).
The third equivalence is
by
Eqs.~(\ref{eq:pr:subgrad-scal-mult:4})
and~(\ref{eq:pr:subgrad-scal-mult:5}).
This completes the proof for $\lambda\in\Rstrictpos$.

Now suppose $\lambda=0$ and $\basubdiffext{\xbar}\ne\emptyset$.
Then both the left-hand side and right-hand side of
\eqref{eq:pr:subgrad-scal-mult:2}
are equal to $\set{\zero}$
(for instance, by \Cref{ex:affine-subgrad-new}),
so the equation holds in this case as well.%
\indexg{scalar multiples (astral)!subgradients of extension|)}%
\indexg{subdifferentials of extensions!scalar multiple@of scalar multiple|)}%
\end{proof}

\section{Sums of functions}

\indexg{sum of functions (standard)!subgradients of extension|(}%
\indexg{subdifferentials of extensions!sum of functions@of sum of functions|(}%
We next consider the astral subgradients of the extension of
a sum of functions.
For two convex functions
$f:\Rn\rightarrow\Rext$ and $g:\Rn\rightarrow\Rext$,
we expect
$\asubdiffplusgext{\xbar}$ to be equal to
$\asubdiffext{\xbar}+\asubdifgext{\xbar}$,
as holds analogously for ordinary gradients in standard calculus,
as well as for standard subgradients,
provided that $f$ and $g$ are proper, and that
the relative interiors of their effective
domains overlap
(\Cref{roc:thm23.8}\ref{roc:thm23.8:b}).
The next theorem shows that this identity for astral subgradients does
indeed hold under the same topological condition as for standard
subgradients:

\begin{theorem}   \label{thm:subgrad-sum-fcns}
  Let $f_i:\Rn\rightarrow\Rext$ be convex,
  for $i=1,\ldots,m$.
  Assume $f_1,\ldots,f_m$ are summable,
  and that
  $\bigcap_{i=1}^m \ri(\dom f_i) \neq \emptyset$.
  Let $h=f_1+\dotsb+f_m$, and let $\xbar\in\extspace$.
  Then
  \begin{equation*}
    \basubdifhext{\xbar}
    =
    \basubdiffextsub{1}{\xbar}
    + \dotsb +
    \basubdiffextsub{m}{\xbar}.
  \end{equation*}
\end{theorem}

\begin{proof}
We begin by proving that
\begin{equation}  \label{eq:thm:subgrad-sum-fcns:2}
      \basubdiffextsub{1}{\xbar}
      + \dotsb +
      \basubdiffextsub{m}{\xbar}
      \subseteq
      \basubdifhext{\xbar}.
\end{equation}
Suppose
$\uu_i\in \asubdiffextsub{i}{\xbar}$,
for $i=1,\ldots,m$, and
let $\uu=\sum_{i=1}^m \uu_i$.
We aim to prove that
$\uu\in\basubdifhext{\xbar}$.

Let
$\hminusu$ and $\fminususubd{i}$ be linear tilts of $h$ and $f_i$
(that is, $\fminususubd{i}=\fminusgen{(f_i)}{\uu_i}$).
Then for all $\xx\in\Rn$,
\begin{equation}   \label{eqn:thm:subgrad-sum-fcns:4}
   \hminusu(\xx)
   =
   h(\xx) - \xx\cdot\uu
   =
   \sum_{i=1}^m \bigParens{f_i(\xx) - \xx\cdot\uu_i}
   =
   \sum_{i=1}^m \fminususubd{i}(\xx),
\end{equation}
noting that, since $f_1(\xx),\ldots,f_m(\xx)$ are summable,
$\fminususubd{1}(\xx),\ldots,\fminususubd{m}(\xx)$ are as well.
Since the rightmost expression of this equation is at least
$\sum_{i=1}^m \inf \fminususubd{i}$
for all $\xx\in\Rn$, this implies
\begin{equation}   \label{eqn:thm:subgrad-sum-fcns:6}
  \inf \hminusu
  \geq
   \sum_{i=1}^m \inf \fminususubd{i}.
\end{equation}

From \Cref{pr:fminusu-props}(\ref{pr:fminusu-props:a}),
we also have that
\[
  \bigcap_{i=1}^m \ri(\dom{\fminususubd{i}})
  =
  \bigcap_{i=1}^m \ri(\dom{f_i})
  \neq
  \emptyset.
\]

For each $i$,
since
$\uu_i\in \asubdiffextsub{i}{\xbar}$,
we have
\begin{equation}   \label{eqn:thm:subgrad-sum-fcns:3}
  \fminususubdext{i}(\xbar)
  =
  \inf \fminususubd{i}
  =
  -\fistar(\uu_i)
  \in\R,
\end{equation}
with both equalities from
\Cref{thm:fminus-subgrad-char}(
\ref{thm:fminus-subgrad-char:a},%
\ref{thm:fminus-subgrad-char:d},%
\ref{thm:fminus-subgrad-char:e}),
and the finiteness of $\fistar(\uu_i)$ from
\Cref{pr:subgrad-imp-in-cldom}(\ref{pr:subgrad-imp-in-cldom:c}).
In particular, this implies that
$\fminususubdext{1}(\xbar),\dotsc,\fminususubdext{m}(\xbar)$
are summable.

Thus, combining, we have that
\begin{equation}  \label{eqn:thm:subgrad-sum-fcns:5}
  \inf \hminusu
  \leq
  \hminusuext(\xbar)
  =
  \sum_{i=1}^m \fminususubdext{i}(\xbar)
  =
  \sum_{i=1}^m \inf \fminususubd{i}
  \leq
  \inf \hminusu.
\end{equation}
The first inequality is by
\Cref{pr:fext-min-exists}.
The first equality follows from
\Cref{thm:ext-sum-fcns-w-duality}
(with $h$ and $f_i$, as they appear in that theorem, set to
$\hminusu$ and $\fminususubd{i}$), all of whose conditions are
satisfied, as argued above.
The second equality and final inequality
are from
Eqs.~(\ref{eqn:thm:subgrad-sum-fcns:3})
and~(\ref{eqn:thm:subgrad-sum-fcns:6}),
respectively.
This shows that
\[
  \hstar(\uu)
  =
  -\inf \hminusu
  =
  -  \sum_{i=1}^m \inf \fminususubd{i}
  =
  \sum_{i=1}^m \fistar(\uu_i),
\]
where the first and last equalities are by
\Cref{pr:fminusu-props}(\ref{pr:fminusu-props:d}),
and the second equality is by
\eqref{eqn:thm:subgrad-sum-fcns:5}.
Combined with
\eqref{eqn:thm:subgrad-sum-fcns:3},
this implies that $\hstar(\uu)\in\R$.

Further, from
\eqref{eqn:thm:subgrad-sum-fcns:5}, it follows that
$\hminusuext(\xbar)=  \inf \hminusu$.
Therefore, $\uu\in\asubdifhext{\xbar}$ by
\Cref{thm:fminus-subgrad-char}(\ref{thm:fminus-subgrad-char:e},\ref{thm:fminus-subgrad-char:a}),
proving
\eqref{eq:thm:subgrad-sum-fcns:2}.

We next prove the reverse inclusion.
Suppose $\uu\in\basubdifhext{\xbar}$.
Then $\hstar(\uu)\in\R$ and $h$ is proper
(by \Cref{pr:subgrad-imp-in-cldom}\ref{pr:subgrad-imp-in-cldom:c}).
This further implies that each $f_i$ is proper (since
if $f_i\equiv+\infty$ then $h\equiv+\infty$, and
if $f_i(\xx)=-\infty$ for any $\xx\in\Rn$, then also $h(\xx)=-\infty$).

We can therefore apply a standard result for the conjugate of a sum
(\Cref{roc:thm16.4})
from which it follows that there exist $\uu_1,\ldots,\uu_m\in\Rn$
with $\uu=\sum_{i=1}^m \uu_i$
and
\begin{equation}   \label{eqn:thm:subgrad-sum-fcns:8}
 \hstar(\uu)=\sum_{i=1}^m \fistar(\uu_i).
\end{equation}
In particular, this implies $\fistar(\uu_i)\in\R$ for all $i$.

As above, let $\hminusu$ and $\fminususubd{i}$ be linear tilts of $h$
and $f_i$, implying
\eqref{eqn:thm:subgrad-sum-fcns:4},
as before.
Then for each $i$,
\begin{equation}   \label{eqn:thm:subgrad-sum-fcns:9}
  \fminususubdext{i}(\xbar)
  \geq
  \inf \fminususubd{i}
  =
  -\fistar(\uu_i)
  \in \R,
\end{equation}
where the inequality and equality are from
\Cref{pr:fext-min-exists}
and
\Cref{pr:fminusu-props}(\ref{pr:fminusu-props:d}).
Thus, $\fminususubdext{i}(\xbar)>-\infty$ for all $i$, so
$\fminususubdext{1}(\xbar),\ldots,\fminususubdext{m}(\xbar)$
are summable.

Consequently, we have that
\begin{align} 
\notag
  \smash[b]{\sum_{i=1}^m \fminususubdext{i}(\xbar)}  %
  =
  \hminusuext(\xbar)
  &=
  \inf \hminusu
\\[\medskipamount]
\label{eqn:thm:subgrad-sum-fcns:7}
  &=
  -\hstar(\uu)
  =
  \smash[t]{
  -\sum_{i=1}^m \fistar(\uu_i)
  =
  \sum_{i=1}^m \inf \fminususubd{i}
  }.
\end{align}
The first equality is from
\Cref{thm:ext-sum-fcns-w-duality},
applied as in the first part of this proof, noting again
that all of that theorem's conditions are
satisfied.
The second equality is from
\Cref{thm:fminus-subgrad-char}(\ref{thm:fminus-subgrad-char:a},\ref{thm:fminus-subgrad-char:e})
since $\uu\in\asubdifhext{\xbar}$.
The third and fifth equalities are both from
\Cref{pr:fminusu-props}(\ref{pr:fminusu-props:d}).
And the fourth equality is from
\eqref{eqn:thm:subgrad-sum-fcns:8}.

\eqref{eqn:thm:subgrad-sum-fcns:7},
combined with
\eqref{eqn:thm:subgrad-sum-fcns:9},
implies
$\fminususubdext{i}(\xbar)=\inf \fminususubd{i}$
for all $i$,
since otherwise, if
$\fminususubdext{i}(\xbar)>\inf \fminususubd{i}$
for some $i$,
then the far left-hand side of
\eqref{eqn:thm:subgrad-sum-fcns:7}
would be strictly greater than the far right-hand side.
Therefore, $\uu_i\in \basubdiffextsub{i}{\xbar}$
by
\Cref{thm:fminus-subgrad-char}(\ref{thm:fminus-subgrad-char:e},\ref{thm:fminus-subgrad-char:a}),
for $i=1,\ldots,m$.
That is,
\[
  \uu
  =
  \uu_1+\dotsb+\uu_m
  \in
  \basubdiffextsub{1}{\xbar}
  + \dotsb +
  \basubdiffextsub{m}{\xbar},
\]
completing the proof.
\end{proof}

As we saw in \Cref{roc:thm23.8}(\ref{roc:thm23.8:a}),
the analogue of the inclusion given in
\eqref{eq:thm:subgrad-sum-fcns:2}
holds for standard subgradients \emph{without} assuming the
relative interiors of the effective domains of the functions have a
point in common.
For astral subgradients, however, that inclusion does not hold in
general if this assumption is omitted:

\begin{example}
\label{ex:subgrad-no-inclusion}
Let $f$, $g$ and $\xbar$ be as in
Example~\ref{ex:KLsets:extsum-not-sum-exts},
and let $h=f+g$.
Then, as shown in that example,
$\fext(\xbar)=\gext(\xbar)=0$, but
$\hext(\xbar)=+\infty$.
Moreover, $\inf f = \inf g = \inf h = 0$.
Therefore, by
\Cref{pr:asub-zero-is-min}(\ref{pr:asub-zero-is-min:a}),
$\zero\in\asubdiffext{\xbar}$ and
$\zero\in\asubdifgext{\xbar}$,
but
$\zero+\zero=\zero\not\in\asubdifhext{\xbar}$.
Thus,
$\asubdiffext{\xbar}+\asubdifgext{\xbar}\not\subseteq\asubdifhext{\xbar}$.%
\indexg{sum of functions (standard)!subgradients of extension|)}%
\indexg{subdifferentials of extensions!sum of functions@of sum of functions|)}%
\end{example}

\section{Composition with an affine map}
\label{sec:subgrad-comp-aff-map}

\indexg{linear map in composition with standard function!subgradients of extension|(}%
\indexg{affine map in composition with standard function!subgradients of extension|(}%
\indexg{subdifferentials of extensions!composition with affine map@of composition with affine map|(}%
We next study the astral subgradients of the extension of
a composition of a convex function with a linear or affine function.

{\mathtogether%
Let $f:\Rm\rightarrow\Rext$, let $\A\in\Rmn$, and let $h=\fA$, that
is, $h(\xx)=f(\A\xx)$ for $\xx\in\Rn$.
By the chain rule from calculus,
if $f$ is differentiable, then
$\gradh(\xx)=\transA \gradf(\A\xx)$.
More generally, if $f$ is convex and proper, and if the relative
interior of its effective domain overlaps with the column space of
$\A$, then the analogous relation holds for standard subgradients, as
we saw in
\Cref{roc:thm23.9}(\ref{roc:thm23.9:b}).
The next theorem shows that the analogous identity holds also for
astral subgradients, assuming this same topological condition.}

\begin{theorem}   \label{thm:subgrad-fA}
  Let $f:\Rm\rightarrow\Rext$ be convex, and
  let $\A\in\Rmn$.
  Assume there exists $\xing\in\Rn$ such that
  $\A\xing\in\ri(\dom f)$.
  Let $\xbar\in\extspace$.
  Then
  \begin{equation*}  %
      \basubdiffAext{\xbar}
      =
      \transA\basubdiffext{\A \xbar}.
  \end{equation*}
\end{theorem}

\begin{proof}
Let $h=\fA$.
Suppose $\uu\in\asubdiffext{\A\xbar}$, and
let $\ww=\transAk\uu$.
We aim to show $\ww\in\asubdifhext{\xbar}$,
which will prove that
$\transA\basubdiffext{\A \xbar} \subseteq \asubdifhext{\xbar}$.

Let $\fminusu$ and $\hminusw$ be linear tilts of $f$ and $h$.
Then for all $\xx\in\Rn$,
\begin{equation}   \label{eq:thm:subgrad-fA:1}
  \hminusw(\xx)
  =
  h(\xx) - \xx\cdot (\transAk\uu)
  =
  f(\A \xx) - (\A\xx)\cdot\uu
  =
  \fminusu(\A \xx).
\end{equation}
By \Cref{thm:ext-linear-comp}, it follows that
\begin{equation}   \label{eq:thm:subgrad-fA:3}
  \hminuswext(\xbar)=\fminusuext(\A\xbar)
\end{equation}
since
$\A\xing\in\ri(\dom{f})=\ri(\dom{\fminusu})$
(by \Cref{pr:fminusu-props}\ref{pr:fminusu-props:a}).

Thus,
\begin{equation}     \label{eq:thm:subgrad-fA:5}
\mathtogether%
  \inf \hminusw
  =
  -\hstar(\ww)
  \leq
  \hminuswext(\xbar)
  =
  \fminusuext(\A\xbar)
  =
  -\fstar(\uu)
  =
  \inf \fminusu
  \leq
  \inf \hminusw.
\end{equation}
The first equality and first inequality are by
\Cref{pr:fminusu-props}(\ref{pr:fminusu-props:d}).
The second equality is \eqref{eq:thm:subgrad-fA:3}.
The third and fourth equalities are by
\Cref{thm:fminus-subgrad-char}(%
\ref{thm:fminus-subgrad-char:a},%
\ref{thm:fminus-subgrad-char:d},%
\ref{thm:fminus-subgrad-char:e})
since $\uu\in\asubdiffext{\A\xbar}$.
And the last inequality follows from
\eqref{eq:thm:subgrad-fA:1}
(since the rightmost expression of that equation is at least
$\inf\fminusu$ for all $\xx\in\Rn$).

\eqref{eq:thm:subgrad-fA:5}
implies that
$\hstar(\ww)=\fstar(\uu)\in\R$
by
\Cref{pr:subgrad-imp-in-cldom}(\ref{pr:subgrad-imp-in-cldom:c})
since $\uu\in\asubdiffext{\A\xbar}$.
It also shows that
$\hminuswext(\xbar)=-\hstar(\ww)$.
Therefore,
$\transAk\uu=\ww\in\asubdifhext{\xbar}$
by
\Cref{thm:fminus-subgrad-char}(\ref{thm:fminus-subgrad-char:d},\ref{thm:fminus-subgrad-char:a}),
proving that
$\transA\basubdiffext{\A \xbar} \subseteq \basubdiffAext{\xbar}$.

For the reverse inclusion,
suppose $\ww\in \basubdifhext{\xbar}$.
We aim to prove that
$\ww\in\transA\asubdiffext{\A\xbar}$.
By \Cref{roc:thm16.3:fA}, which gives
a rule for the conjugate of a function
$\fA$ under the assumptions of this
\namecref{thm:subgrad-fA},
\begin{equation}   \label{eq:thm:subgrad-fA:2}
  \hstar(\ww)
  =
  \inf\Braces{ \fstar(\uu) :\: \uu\in\Rn,\,
                       \transAk\uu = \ww },
\end{equation}
and furthermore, this infimum must be realized
since it cannot be vacuous
since ${\hstar(\ww)\in\R}$
(by \Cref{pr:subgrad-imp-in-cldom}\ref{pr:subgrad-imp-in-cldom:c}).
Therefore, there exists $\uu\in\Rn$ for which
$\transAk\uu = \ww$ and $\hstar(\ww)=\fstar(\uu)$.
In particular, this implies
that $\fstar(\uu)\in\R$, and also that
Eqs.~(\ref{eq:thm:subgrad-fA:1}) and~(\ref{eq:thm:subgrad-fA:3})
continue to hold.
Thus,
\begin{equation}   \label{eq:thm:subgrad-fA:6}
  \fminusuext(\A\xbar)
  =
  \hminuswext(\xbar)
  =
  -\hstar(\ww)
  =
  -\fstar(\uu),
\end{equation}
with the first equality from
\eqref{eq:thm:subgrad-fA:3}, and
the second from
\Cref{thm:fminus-subgrad-char}(\ref{thm:fminus-subgrad-char:a},\ref{thm:fminus-subgrad-char:d})
since $\ww\in \asubdifhext{\xbar}$.
By that same 
\namecref{thm:fminus-subgrad-char},
\eqref{eq:thm:subgrad-fA:6} implies
$\uu\in\asubdiffext{\A\xbar}$, and thus that
$\ww\in\transA\asubdiffext{\A\xbar}$.
This completes the proof.
\end{proof}

In general,
without the assumption that $\colspace \A$ overlaps
$\ri(\dom{f})$,
\Cref{thm:subgrad-fA} need not hold:

\begin{example}
\indexg{Sideways cone!subgradients of composition with linear map and|(}%
Let $f$, $\A$ and $\xbar$ be as in
Example~\ref{ex:fA-ext-countereg},
and let $h=\fA$.
Then, as shown in that example,
$\fext(\A\xbar)=0$ and $\hext(\xbar)=+\infty$.
Since $\inf f=\inf h=0$, this implies
that
$\zero\in\basubdiffext{\A\xbar}$,
but that $\zero=\transA\zero\not\in\asubdifhext{\xbar}$,
both by \Cref{pr:asub-zero-is-min}(\ref{pr:asub-zero-is-min:a}).
Thus,
\indexg{linear map in composition with standard function!subgradients of extension|)}%
\indexg{Sideways cone!subgradients of composition with linear map and|)}%
$\transA\basubdiffext{\A \xbar} \not\subseteq \basubdiffAext{\xbar}$.
\end{example}

When a fixed vector $\bb\in\Rn$ is added to the argument of a function
$f:\Rn\rightarrow\Rext$, yielding a function
$\xx\mapsto f(\bb+\xx)$, the astral subgradients are
unaffected, as expected.
This is in accord with analogous results for standard gradients and
subgradients from calculus (based on the chain rule)
and standard convex analysis:

\begin{proposition}   \label{pr:subgrad-shift-arg}
  Let $f:\Rn\rightarrow\Rext$, let $\bb\in\Rn$,
  and let $h:\Rn\rightarrow\Rext$ be defined by
  $h(\xx)=f(\bb+\xx)$ for $\xx\in\Rn$.
  Let $\xbar\in\extspace$.
  Then
  \[
    \asubdifhext{\xbar} = \asubdiffext{\bb\plusl\xbar}.
  \]
\end{proposition}

\begin{proof}
~
\begin{proof-parts}
\pfpart{``$\subseteq$'':}
Let $\uu\in\asubdifhext{\xbar}$. By \Cref{thm:fminus-subgrad-char}(\ref{thm:fminus-subgrad-char:a},\ref{thm:fminus-subgrad-char:b}),
there exists a sequence $\seq{\xx_t}$ in $\Rn$ with $\xx_t\to\xbar$,
$h(\xx_t)\in\R$ for all $t$,
and $\xx_t\inprod\uu-h(\xx_t)\to\hstar(\uu)$.
Let $\ybar=\bb\plusl\xbar$ and $\yy_t=\bb+\xx_t$ for all $t$.
Then $\yy_t\to\ybar$
(by \Cref{pr:i:7}\ref{pr:i:7f})
and
\[
  \yy_t\inprod\uu-f(\yy_t)
  =
  \bb\inprod\uu+[\xx_t\inprod\uu-h(\xx_t)]
  \to
  \bb\inprod\uu+\hstar(\uu)
  =
  \fstar(\uu).
\]
The first equality is because $\yy_t=\bb+\xx_t$ and $f(\yy_t)=h(\xx_t)$.
The last equality is by \Cref{pr:conj-shift}(\ref{pr:conj-shift:a}).
The last equality also shows that $\fstar(\uu)\in\R$,
because $\hstar(\uu)\in\R$ by \Cref{pr:subgrad-imp-in-cldom}(\ref{pr:subgrad-imp-in-cldom:c}),
since $\uu\in\asubdifhext{\xbar}$.
Thus, $\uu\in\asubdiffext{\ybar}$
by \Cref{thm:fminus-subgrad-char}(\ref{thm:fminus-subgrad-char:c},\ref{thm:fminus-subgrad-char:a}).

\pfpart{``$\supseteq$'':}
Let $f'=h$, $\bb'=-\bb$, and 
$h'(\xx)=f'(\bb'+\xx)=h(-\bb+\xx)=f(\xx)$ for $\xx\in\Rn$. The first part of the proof
then implies
$\asubdiffext{\xbar'}\subseteq \asubdifplain{\hext}\bigParens{(-\bb)\plusl\xbar'}$
for $\xbar'\in\eRn$. Letting $\xbar'=\bb\plusl\xbar$ then
yields $\asubdiffext{\bb\plusl\xbar}\subseteq \asubdifplain{\hext}\bigParens{\xbar}$.
\qedhere
\end{proof-parts}
\end{proof}

Combining
\Cref{thm:subgrad-fA} and
\Cref{pr:subgrad-shift-arg}
yields a corollary for astral subgradients of the extension of
a composition of a convex function with an affine function:

\begin{corollary}   \label{cor:subgrad-comp-w-affine}
  Let $f:\Rm\rightarrow\Rext$ be convex,
  let $\A\in\Rmn$, and let $\bb\in\Rm$.
  Assume there exists $\xing\in\Rn$ such that
  $\bb+\A\xing\in\ri(\dom f)$.
  Define $h:\Rn\rightarrow\Rext$ by
  $h(\xx)=f(\bb+\A\xx)$ for $\xx\in\Rn$.
  Let $\xbar\in\extspace$.
  Then
  \begin{equation*}
    \basubdifhext{\xbar}
    =
    \transA \basubdiffext{\bb \plusl \A \xbar}.
  \end{equation*}
\end{corollary}

\begin{proof}
Let $g:\Rm\rightarrow\Rext$ be defined by
$g(\yy)=f(\bb+\yy)$ for $\yy\in\Rm$.
Then $h=g\A$.
Also,
$\bb+\A\xing\in\ri(\dom f)$
implies $\A\xing\in\ri(\dom g)$.
Thus,
\[
  \asubdifhext{\xbar}
  =
  \transA \asubdifgext{\A \xbar}
  =
  \transA \asubdiffext{\bb\plusl\A\xbar},
\]
where the first equality is by
\Cref{thm:subgrad-fA}, and the second is by
\Cref{pr:subgrad-shift-arg}.%
\indexg{affine map in composition with standard function!subgradients of extension|)}%
\indexg{subdifferentials of extensions!composition with affine map@of composition with affine map|)}%
\end{proof}

\section{Reductions}

\indexg{reductions, iconic!subgradients of extension|(}%
\indexg{subdifferentials of extensions!reduction@of reduction|(}%
We next study astral subgradients of
reductions of convex functions.
Given a convex function $f:\Rn\rightarrow\Rext$ and an icon
$\ebar\in\corezn$, we saw
in \Cref{pr:i:9}(\ref{pr:i:9b}) that
the extension $\eg$ of the reduction
$g=\fshadd$
satisfies
$\eg(\xbar) = \ef(\ebar\plusl\xbar)$
for $\xbar\in\extspace$.
This is similar to the extension of the function
$h(\xx)=f(\bb+\xx)$
with $\bb\in\Rn$,
considered in \Cref{pr:subgrad-shift-arg}, whose extension
is
$\eh(\xbar)=\ef(\bb\plusl\xbar)$ (by \Cref{pr:ext-affine-comp}\ref{pr:ext-affine-comp:b}).
Thus,
the results of this section can be seen as a generalization of
\Cref{pr:subgrad-shift-arg},
which showed that the astral subdifferential of $\eh$ at $\xbar$
coincides with the astral subdifferential of $\fext$ at $\bb\plusl\xbar$,
similar to analogous rules from calculus and standard convex analysis.

We expect similar behavior for reductions.
That is, supposing $g(\xx)=\fshadd(\xx)=\fext(\ebar\plusl\xx)$,
we might expect $\asubdifgext{\xbar}=\asubdiffext{\ebar\plusl\xbar}$.
As shown in the next theorem,
this turns out to be the case, provided that $g$ is proper.
Otherwise, if $g$ is improper, then $\gext$ is not astral
subdifferentiable anywhere:

\begin{theorem}  \label{thm:subgrad-icon-reduc-proper}
  Let $f:\Rn\rightarrow\Rext$ be convex,
  let $\ebar\in\corezn$,
  and
  let $g=\fshadd$.
  Let $\xbar\in\extspace$.
  If $g$ is proper, then
  \begin{equation}  \label{eq:thm:subgrad-icon-reduc-proper:statement}
     \basubdifgext{\xbar}
     =
     \basubdiffext{\ebar\plusl\xbar}.
  \end{equation}
  Otherwise, if $g$ is improper, then
  $\basubdifgext{\xbar} = \emptyset$.
\end{theorem}

Before proving the theorem, we show that
a vector
$\uu\in\Rn$ that is orthogonal to icon
$\ebar\in\aresconef$
will either be
in both or neither of the astral subdifferentials
$\basubdifgext{\xbar}$ and $\basubdiffext{\ebar\plusl\xbar}$:

\begin{lemma}   \label{lem:subgrad-reduc-equiv-eeu-z}
  Let $f:\Rn\rightarrow\Rext$ be convex,
  let $\ebar\in\corezn\cap(\aresconef)$,
  and
  let $g=\fshadd$.
  Let $\xbar\in\extspace$ and $\uu\in\Rn$.
  Assume $\ebar\cdot\uu=0$.
  Then $\uu\in\basubdifgext{\xbar}$
  if and only if
  $\uu\in\basubdiffext{\ebar\plusl\xbar}$.
\end{lemma}

\begin{proof}
If $f\equiv+\infty$ then also $g\equiv+\infty$
(\Cref{pr:d2}\ref{pr:d2:c}), implying
$\basubdiffext{\ebar\plusl\xbar}=\basubdifgext{\xbar}=\emptyset$
(by \Cref{pr:subgrad-imp-in-cldom}\ref{pr:subgrad-imp-in-cldom:c}),
thereby proving the claim in this case.
We therefore assume henceforth that $f\not\equiv+\infty$.

Let $\fminusu$ and $\gminusu$ be linear tilts of $f$ and $g$ by
$\uu$.
Then for $\xx\in\Rn$,
\[
  \gminusu(\xx)
  =
  g(\xx) - \xx\cdot\uu
  =
  \fext(\ebar\plusl\xx) - (\ebar\plusl\xx)\cdot\uu
  =
  \fminusuext(\ebar\plusl\xx).
\]
The second equality is because $\ebar\cdot\uu=0$ so that
$(\ebar\plusl\xx)\cdot\uu=\xx\cdot\uu\in\R$.
The last equality then follows from
\Cref{pr:fminusu-props}(\ref{pr:fminusu-props:e})
since $\fext(\ebar\plusl\xx)$
and $-(\ebar\plusl\xx)\cdot\uu$ are summable.
Thus, $\gminusu=\fminusushadd$, so
\begin{equation}  \label{eq:thm:subgrad-icon-reduc:1}
  \gminusuext(\xbar)
  =
  \fminusuext(\ebar\plusl\xbar)
\end{equation}
by
\Cref{pr:i:9}(\ref{pr:i:9b}).

We then have that
\begin{alignat*}{3}
  \uu\in\asubdifgext{\xbar}\quad
    {\Leftrightarrow}\quad
    && \gstar(\uu)\in\R
    &  \quad\text{and}\quad
    &  \gminusuext(\xbar) &= -\gstar(\uu)
  \\
    {\Leftrightarrow}\quad
    && \fstar(\uu)\in\R
    &  \quad\text{and}\quad
    &  \fminusuext(\ebar\plusl\xbar) &= -\fstar(\uu)
  \\
    {\Leftrightarrow}\quad
    && \uu\in\asubdiffext{\ebar\plusl\xbar}\makebox[0pt][l]{.}
\end{alignat*}
The first and last equivalences are both by
\Cref{pr:subgrad-imp-in-cldom}(\ref{pr:subgrad-imp-in-cldom:c})
and
\Cref{thm:fminus-subgrad-char}(\ref{thm:fminus-subgrad-char:a},\ref{thm:fminus-subgrad-char:d}).
The second equivalence is by
\eqref{eq:thm:subgrad-icon-reduc:1}
and because $\gstar(\uu)=\fstar(\uu)$
by \Cref{thm:conj-of-iconic-reduc}
since $\ebar\cdot\uu=0$.
This proves the claim.
\end{proof}

We also note that all astral subgradients of
the extension of a reduction $\fshadd$ must be orthogonal to $\ebar$:

\begin{proposition}   \label{pr:subg-reduc-perp}
  Let $f:\Rn\rightarrow\Rext$ be convex,
  let $\ebar\in\corezn$,
  and
  let $g=\fshadd$.
  Let $\xbar\in\extspace$.
  If $\uu\in\asubdifgext{\xbar}$ then $\ebar\cdot\uu=0$.
\end{proposition}

\begin{proof}
If either $f\equiv+\infty$ or $\ebar\not\in\aresconef$ then
$g\equiv+\infty$ by
\Cref{pr:d2}(\ref{pr:d2:c}) and \Cref{cor:a:4},
implying that $\gstar\equiv-\infty$ and so that
$\asubdifgext{\xbar}=\emptyset$ by
\Cref{pr:subgrad-imp-in-cldom}(\ref{pr:subgrad-imp-in-cldom:c}).
Thus, in this case, the claim holds vacuously.

Otherwise, if $f\not\equiv+\infty$ and $\ebar\in\aresconef$,
and if $\uu\in\asubdifgext{\xbar}$,
then $\gstar(\uu)\in\R$ by
\Cref{pr:subgrad-imp-in-cldom}(\ref{pr:subgrad-imp-in-cldom:c}),
implying $\ebar\cdot\uu=0$
by \Cref{thm:conj-of-iconic-reduc}.
\end{proof}

\begin{proof}[Proof of \Cref{thm:subgrad-icon-reduc-proper}]
If $g$ is not proper then
$\asubdifgext{\xbar}=\emptyset$ by
\Cref{pr:subgrad-imp-in-cldom}(\ref{pr:subgrad-imp-in-cldom:c}),
proving the theorem for this case.
We therefore assume henceforth that $g$ is proper.
Consequently,
$f\not\equiv+\infty$ and $\ebar\in\aresconef$
since otherwise we would have
$g\equiv+\infty$ by
\Cref{pr:d2}(\ref{pr:d2:c}) and \Cref{cor:a:4}.

Let $\uu\in\Rn$.
We aim to show that 
$\uu\in\basubdifgext{\xbar}$
if and only if
$\uu\in\basubdiffext{\ebar\plusl\xbar}$.

Suppose first that $\uu\in\basubdifgext{\xbar}$.
Then $\ebar\cdot\uu=0$ by
\Cref{pr:subg-reduc-perp}.
Therefore, by
\Cref{lem:subgrad-reduc-equiv-eeu-z},
$\uu\in\asubdiffext{\ebar\plusl\xbar}$.

For the reverse implication, suppose
$\uu\in \basubdiffext{\ebar\plusl\xbar}$.
By \Cref{pr:i:9}(\ref{pr:i:9a}),
$g$~is convex and lower semicontinuous,
and it is also proper by assumption.
Therefore,
its effective domain,
$\dom{g}$, is convex and nonempty, so its relative interior
is also nonempty
(\Cref{pr:ri-props}\ref{pr:ri-props:roc-thm6.2b}).
Thus, there exists a point $\yy\in\ri(\dom{g})$.
Every such point must have a standard subgradient
$\ww\in\partial g(\yy)$
(\Cref{roc:thm23.4}\ref{roc:thm23.4:a}),
implying
$\ww\in\basubdifgext{\yy}$
by
\Cref{pr:asubdiffext-at-x-in-rn}(\ref{pr:asubdiffext-at-x-in-rn:c}).
It then follows that
$\ebar\cdot\ww=0$
by \Cref{pr:subg-reduc-perp},
and so that
$\ww\in\basubdiffext{\ebar\plusl\yy}$
by \Cref{lem:subgrad-reduc-equiv-eeu-z}
(applied with $\uu$ and $\xbar$, as they appear there,
set to $\ww$ and $\yy$, respectively).

We claim that $\ebar\cdot\uu=0$.
To see this, suppose to the contrary that
$\ebar\cdot\uu\neq 0$, implying, since $\ebar$ is an icon, that
$\ebar\cdot\uu\in\{-\infty,+\infty\}$
(by
\Cref{pr:icon-equiv}\ref{pr:icon-equiv:a}\ref{pr:icon-equiv:b}).
Note that
\begin{equation}   \label{eq:thm:subgrad-icon-reduc-proper:5}
  \ebar\cdot(\uu-\ww)=\ebar\cdot\uu-\ebar\cdot\ww=\ebar\cdot\uu
\end{equation}
since $\ebar\cdot\ww=0$
(and by \Cref{pr:i:1} in the first equality).
Since
$\uu\in \basubdiffext{\ebar\plusl\xbar}$
and
$\ww\in\basubdiffext{\ebar\plusl\yy}$,
we can apply monotonicity, yielding
\begin{align}
  \ebar\cdot\uu\plusl\xbar\cdot(\uu-\ww)
  &=
  (\ebar\plusl\xbar)\cdot(\uu-\ww)
  \notag
  \\
  &\geq
  \fstar(\uu) - \fstar(\ww)
  \notag
  \\
  &\geq
  (\ebar\plusl\yy)\cdot(\uu-\ww)
  =
  \ebar\cdot\uu\plusl\yy\cdot(\uu-\ww).
  \label{eq:thm:subgrad-icon-reduc-proper:6}
\end{align}
The two inequalities are by
\Cref{cor:ast-subgrad-monotone}
and because $\fextstar=\fstar$
(by \Cref{pr:fextstar-is-fstar}).
The two equalities both follow from
\eqref{eq:thm:subgrad-icon-reduc-proper:5}.
Thus, if $\ebar\cdot\uu\in\set{-\infty,+\infty}$,
then the far left-hand side and far right-hand side
of \eqref{eq:thm:subgrad-icon-reduc-proper:6} both must be
equal to $\ebar\cdot\uu$,
implying
$\fstar(\uu)-\fstar(\ww)=\ebar\cdot\uu\in\set{-\infty,+\infty}$.
But this
is a contradiction since
$\fstar(\uu)$ and $\fstar(\ww)$ are both in $\R$
(by
\Cref{pr:subgrad-imp-in-cldom}\ref{pr:subgrad-imp-in-cldom:c}).

We conclude that $\ebar\cdot\uu=0$.
Since also
$\uu\in \basubdiffext{\ebar\plusl\xbar}$,
this implies that
$\uu\in\basubdifgext{\xbar}$ as well,
by \Cref{lem:subgrad-reduc-equiv-eeu-z},
completing the proof.
\end{proof}

Here is a concrete example showing how
\eqref{eq:thm:subgrad-icon-reduc-proper:statement}
from \Cref{thm:subgrad-icon-reduc-proper}
can fail to hold when $g$ is improper:

\begin{example}  \label{ex:simple-neg-lin}
Let $f(x)=-x$ for $x\in\R$,
let $\bare=+\infty$,
and let $g=\fshad{\bare}$.
Then $g\equiv-\infty$, which is improper.
For all $\barx\in\Rext$, this implies that
$\asubdifgext{\barx}=\emptyset$.
On the other hand,
$\asubdiffext{\bare\plusl\barx}=\{-1\}$,
as we saw in \Cref{ex:affine-subgrad-new}.
\end{example}

The next example shows how the tools we have developed can be
applied to calculate astral subdifferentials:

\begin{example}[Absolute value at infinity, alternate derivation]
     \label{ex:subgrad-sqrt-approaches-abs-val:cont2}
\indexg{Absolute value at infinity!astral subgradients of extension|(}%
We return to the function $f$
from
Examples~\ref{ex:subgrad-sqrt-approaches-abs-val}
and~\ref{ex:subgrad-sqrt-approaches-abs-val:cont},
which was defined as
$f(\xx)=\sqrt{x_1^2 + e^{-x_2}}$ for ${\xx\in\R^2}$,
with
$\ef(\xbar)=\sqrt{(\xbar\inprod\ee_1)^2 + \expex(-\xbar\inprod\ee_2)}$
for $\xbar\in\extspac{2}$.
In
Example~\ref{ex:subgrad-sqrt-approaches-abs-val:cont},
we showed how to compute
the astral subdifferentials of this function at points
$\xbar_\alpha= \limray{\ee_2} \plusl \alpha \ee_1$,
for $\alpha\in\R$.
Here, we give an alternative and somewhat more general derivation.
Let $g=\fshad{\limray{\ee_2}}$, and let $h:\R\rightarrow\R$ be the
absolute value function, $h(x)=|x|$ for $x\in\R$.
Then for $\zz\in\R^2$,
\[
  g(\zz)
  =
  \fext(\limray{\ee_2}\plusl\zz)
  =
  |\zz\cdot\ee_1|
  =
  h(\zz\cdot\ee_1)
  =
  h(\trans{\ee_1}\zz).
\]
Thus, $g=h \trans{\ee_1}$.
We therefore have, for $\zbar\in\extspac{2}$, that
\begin{equation}
\label{eq:ex:subgrad-sqrt-approaches-abs-val:cont2}
  \asubdiffext{\limray{\ee_2}\plusl\zbar}
  =
  \asubdifgext{\zbar}
  =
  \ee_1 \asubdifhext{\trans{\ee_1} \zbar}
  =
  \ee_1 \asubdifhext{\zbar\cdot\ee_1}.
\end{equation}
The first equality is by
\Cref{thm:subgrad-icon-reduc-proper}.
The second is by
\Cref{thm:subgrad-fA}
(applicable since $\dom h = \R$).
In particular, instantiating \eqref{eq:ex:subgrad-sqrt-approaches-abs-val:cont2} with $\zbar=\alpha\ee_1$, we obtain
\[
  \asubdiffext{\xbar_\alpha}
  =
  \asubdiffext{\limray{\ee_2}\plusl\alpha\ee_1}
  =
  \ee_1 \asubdifhext{\alpha}
  =
  \ee_1 \partial h(\alpha),
\]
where the last equality is by \Cref{pr:asubdiffext-at-x-in-rn}(\ref{pr:asubdiffext-at-x-in-rn:c}).
Plugging in the subdifferential of $h$ from \Cref{ex:standard-abs-val-subgrad}
then recovers the expression for $\asubdiffext{\xbar_\alpha}$ from
\indexg{Absolute value at infinity!astral subgradients of extension|)}%
\Cref{ex:subgrad-sqrt-approaches-abs-val:cont}.
\end{example}

As shown in the next \namecref{cor:subgrad-arb-reduc},
\Cref{thm:subgrad-icon-reduc-proper}
can be straightforwardly combined with
\Cref{pr:subgrad-shift-arg} to obtain rules for the astral
subgradients of the extension of a function of the form
$\xx\mapsto\fext(\zbar\plusl\xx)$, for any $\zbar\in\extspace$.
Here, we assume there is at least one point in
$\zbar$'s galaxy where $\fext$ is finite.
This condition is equivalent to the function $h(\xx)=\fext(\zbar\plusl\xx)$
being proper (by
\Cref{pr:improper-vals}\ref{pr:improper-vals:cor7.2.1}).
Similar assumptions are made in several of the results that follow.

\begin{corollary}  \label{cor:subgrad-arb-reduc}
  Let $f:\Rn\rightarrow\Rext$ be convex,
  let $\zbar\in\extspace$, and let
  $h:\Rn\rightarrow\Rext$ be defined by
  $h(\xx)=\fext(\zbar\plusl\xx)$
  for $\xx\in\Rn$.
  Assume there exists $\xhat\in\Rn$
  for which $h(\xhat)\in\R$.
  Then:
  \begin{letter-compact}
  \item  \label{cor:subgrad-arb-reduc:a}
    $\basubdiffext{\zbar\plusl\xbar}=\basubdifhext{\xbar}$
    for all $\xbar\in\extspace$.
  \item  \label{cor:subgrad-arb-reduc:b}
    $\basubdiffext{\zbar\plusl\xx}=\basubdifhext{\xx}=\partial h(\xx)$
    for all $\xx\in\Rn$.
  \end{letter-compact}
\end{corollary}

\begin{proof}
We can write $\zbar=\ebar\plusl\qq$ for some $\ebar\in\corezn$ and
$\qq\in\Rn$
(\Cref{thm:icon-fin-decomp}).
Let $g=\fshadd$.
Then,
for $\xx\in\Rn$,
\begin{equation}   \label{eq:cor:subgrad-arb-reduc:1}
  h(\xx) = \fext(\ebar\plusl\qq\plusl\xx) = g(\qq+\xx).
\end{equation}
Since $g$ is convex and lower semicontinuous
(by \Cref{pr:i:9}\ref{pr:i:9a}), $h$ is as well, being the
composition of the convex, lower semicontinuous function $g$ with the
affine function $\yy\mapsto\qq+\yy$
(Propositions~\ref{pr:lsc-res-domain}\ref{pr:lsc-res-domain:a}
and~\ref{roc:thm5.7:fA}).
Therefore, since $h(\xhat)\in\R$, $h$ is proper
(by
\Cref{pr:improper-vals}\ref{pr:improper-vals:cor7.2.1}).
Thus, $g$ is proper as well.

\begin{proof-parts}
\pfpart{Part~(\ref{cor:subgrad-arb-reduc:a}):}
Let $\xbar\in\extspace$.
Then
\[
  \basubdifhext{\xbar}
  =
  \basubdifgext{\qq\plusl\xbar}
  =
  \basubdiffext{\ebar\plusl\qq\plusl\xbar}
  =
  \basubdiffext{\zbar\plusl\xbar},
\]
where the first equality is
by \Cref{pr:subgrad-shift-arg}
(combined with Eq.~\ref{eq:cor:subgrad-arb-reduc:1}),
and the second by
\Cref{thm:subgrad-icon-reduc-proper}.

\pfpart{Part~(\ref{cor:subgrad-arb-reduc:b}):}
Let $\xx\in\Rn$.
Then
$\basubdifhext{\xx}=\basubdiffext{\zbar\plusl\xx}$
by part~(\ref{cor:subgrad-arb-reduc:a}),
and
$\partial h(\xx)=\basubdifhext{\xx}$
by
\Cref{pr:asubdiffext-at-x-in-rn}(\ref{pr:asubdiffext-at-x-in-rn:c})
(since $h$ is lower semicontinuous).%
\indexg{reductions, iconic!subgradients of extension|)}%
\indexg{subdifferentials of extensions!reduction@of reduction|)}%
\qedhere
\end{proof-parts}
\end{proof}

For the remainder of this section, we consider
various consequences of the previous results.
\indexg{subdifferentials of extensions!characterizations|(}%
We begin with
a characterization of astral subgradients
that is analogous to the standard definition of subgradient as given in
\eqref{eqn:standard-subgrad-ineq}:

\begin{proposition}  \label{pr:ast-subgrad-char-std-dfn:alt}
  Let $f:\Rn\rightarrow\Rext$ be convex, 
  let $\xbar\in\extspace$, and let $\uu\in\Rn$.
  Assume there exists $\yhat\in\Rn$ for which
  $\fext(\xbar\plusl\yhat)\in\R$.
  Then $\uu\in \basubdiffext{\xbar}$
  if and only if
  \begin{equation}
  \label{eq:ast-subgrad-char-std-dfn:alt}
    \fext(\xbar\plusl\yy)\geq\fext(\xbar) + \yy\cdot\uu
  \end{equation}
  for all $\yy\in\Rn$.
\end{proposition}

\begin{proof}
Similar to the proof of \Cref{cor:subgrad-arb-reduc},
we can write
$\xbar=\ebar\plusl\qq$ for some $\ebar\in\corezn$ and
$\qq\in\Rn$.
Let $g=\fshadd$, and define
$h:\Rn\rightarrow\Rext$ by
$h(\yy)=g(\qq+\yy)=\fext(\ebar\plusl\qq\plusl\yy)=\fext(\xbar\plusl\yy)$
for $\yy\in\Rn$.

By assumption,
$h(\yhat)\in\R$, so
$\asubdiffext{\xbar}=\partial h(\zero)$ by
\Cref{cor:subgrad-arb-reduc}(\ref{cor:subgrad-arb-reduc:b}).
Further,
by definition of standard subgradient
(Eq.~\ref{eqn:standard-subgrad-ineq}),
$\uu\in\basubdiffext{\xbar}=\partial h(\zero)$ if and only if
$h(\yy)\geq h(\zero) + \yy\cdot\uu$
for all $\yy\in\Rn$,
proving the \namecref{pr:ast-subgrad-char-std-dfn:alt}.
\end{proof}

Without the assumption that $\fext$ is finite at some point
$\xbar\plusl\yhat$,
\Cref{pr:ast-subgrad-char-std-dfn:alt} does not hold in general.
For instance, if $f(x)=-x$ for $x\in\R$
and $\barx=+\infty$, then $\asubdiffext{\barx}=\{-1\}$
(by \Cref{ex:affine-subgrad-new}), but
the condition given in \eqref{eq:ast-subgrad-char-std-dfn:alt}
is satisfied by all points
$u\in\Rext$ since $\fext(\barx)=-\infty$.%
\indexg{subdifferentials of extensions!characterizations|)}%

As we saw in \Cref{sec:prelim:subgrads},
the one-sided directional derivative $\dderf{\xx}{\yy}$ of a function
$f:\Rn\rightarrow\Rext$ provides the instantaneous rate
at which $f$ is increasing at a point $\xx\in\Rn$
in the direction of any vector $\yy\in\Rn$.
Directional derivatives are closely related to standard
subdifferentials, and, as we saw in
\Cref{sec:astral-dual-subgrad}, to astral dual subdifferentials
as well.\looseness=-1

\indexg{one-sided directional derivatives!defined for astral functions|(}%
We can straightforwardly
extend the definition of directional derivative
given in \eqref{eq:direc-deriv-dfn} to astral functions.
Thus, we define the one-sided directional derivative of
a function $F:\extspace\rightarrow\Rext$ at a point
$\xbar\in\extspace$ with respect to a vector $\yy\in\Rn$
to be
\[
  \dderF{\xbar}{\yy}
  =
  \lim_{\rightlim{\lambda}{0}} \frac{F(\xbar\plusl\lambda\yy) - F(\xbar)}
                                      {\lambda},
\]
provided $F(\xbar)\in\R$.
As is the case for standard functions, such an astral directional
derivative captures the instantaneous rate at which $F$ is increasing
at $\xbar$ in the direction of the (finite) vector $\yy$.%
\indexg{one-sided directional derivatives!defined for astral functions|)}%

\indexg{one-sided directional derivatives!an extension@of an extension|(}%
\indexg{subdifferentials of extensions!directional derivatives and|(}%
For a convex function $f:\Rn\rightarrow\Rext$, the next
\namecref{pr:ast-direc-deriv}
relates the directional derivatives of $\fext$ to those of a standard
function on $\Rn$, as well as to the astral subgradients of $\fext$,
analogous to \Cref{roc:thm23.2}:

\begin{proposition}  \label{pr:ast-direc-deriv}
  Let $f:\Rn\rightarrow\Rext$ be convex,
  let $\xbar\in\extspace$ and $\uu\in\Rn$,
  and let
  $h:\Rn\rightarrow\Rext$ be defined by
  $h(\yy)=\fext(\xbar\plusl\yy)$
  for $\yy\in\Rn$.
  Assume $\fext(\xbar)\in\R$.
  Then:
  \begin{letter-compact}
  \item  \label{pr:ast-direc-deriv:a}
    For $\yy\in\Rn$, $\dderfext{\xbar}{\yy}=\dderh{\zero}{\yy}$.
  \item  \label{pr:ast-direc-deriv:b}
    $\uu\in\asubdiffext{\xbar}$
    if and only if
    $\yy\cdot\uu\leq\dderfext{\xbar}{\yy}$
    for all $\yy\in\Rn$.
  \end{letter-compact}
\end{proposition}

\begin{proof}
  ~

\begin{proof-parts}
\pfpart{Part~(\ref{pr:ast-direc-deriv:a}):}
From definitions,
for $\yy\in\Rn$,
\[
  \dderfext{\xbar}{\yy}
  =
  \lim_{\rightlim{\lambda}{0}}
                  \frac{\fext(\xbar\plusl\lambda\yy) - \fext(\xbar)}
                                      {\lambda}
  =
  \lim_{\rightlim{\lambda}{0}} \frac{h(\lambda\yy) - h(\zero)}
                                      {\lambda}
  =                                      
  \dderh{\zero}{\yy}.
\]

\pfpart{Part~(\ref{pr:ast-direc-deriv:b}):}
By \Cref{cor:subgrad-arb-reduc}(\ref{cor:subgrad-arb-reduc:b}),
$\asubdiffext{\xbar}=\partial h(\zero)$
(noting, by assumption, that $h(\zero)\in\R$).
By \Cref{roc:thm23.2},
$\uu\in\partial h(\zero)$ if and only if
$\yy\cdot\uu\leq\dderh{\zero}{\yy}$ for all $\yy\in\Rn$.
Since, by part~(\ref{pr:ast-direc-deriv:a}),
$\dderh{\zero}{\yy}=\dderfext{\xbar}{\yy}$ for all $\yy$,
this proves the claim.%
\indexg{subdifferentials of extensions!directional derivatives and|)}%
\indexg{one-sided directional derivatives!an extension@of an extension|)}%
\qedhere
\end{proof-parts}
\end{proof}

\indexg{subdifferentials of extensions!nonemptiness of|(}%
In \Cref{sec:dual-subdif-not-empty}, we saw that every convex
function $\psi:\Rn\rightarrow\Rext$ is astral dual subdifferentiable
at every point, even points outside its domain.
The same is not true for the astral primal subdifferentials
of extensions, which certainly can be empty; for instance,
\Cref{pr:subgrad-imp-in-cldom}(\ref{pr:subgrad-imp-in-cldom:c}) shows
that they are necessarily empty at all points outside $\cldom{f}$.

Applying
\Cref{thm:subgrad-icon-reduc-proper}
and
\Cref{cor:subgrad-arb-reduc},
we provide next some
specific conditions for when an extension is astral
subdifferentiable,
beginning with the necessary conditions given in the following
theorem.
For a convex function $f:\extspace\rightarrow\Rext$ and $\xbar\in\extspace$,
these show that if $\fext$ is finite at any point in $\xbar$'s galaxy
(which is equivalent to the function $\zz\mapsto\fext(\xbar\plusl\zz)$
being proper)
and if $\fext$ is astral subdiffentiable at $\xbar$,
then $\fext$ must necessarily be finite at $\xbar$, and all astral
subgradients $\uu\in\asubdiffext{\xbar}$ must be such that
$\xbar\cdot\uu$ is also finite.
Further, if this properness condition on $\fext$ does not hold, then
both $\xbar\cdot\uu$ and $\fext(\xbar)$ must be infinite.

\begin{theorem}   \label{thm:nec-fext-subdif}
  Let $f:\Rn\rightarrow\Rext$ be convex and
  let $\xbar\in\extspace$.
  Suppose $\uu\in\asubdiffext{\xbar}$.
  Then the following are equivalent:
  \begin{letter-compact}
  \item   \label{thm:nec-fext-subdif:a}
    There exists $\zhat\in\Rn$ such that $\fext(\xbar\plusl\zhat)\in\R$.
  \item   \label{thm:nec-fext-subdif:b}
    $\xbar\cdot\uu\in\R$.
  \item   \label{thm:nec-fext-subdif:c}
    $\fext(\xbar)\in\R$.
  \end{letter-compact}
\end{theorem}

\begin{proof}
  ~

\begin{proof-parts}
\pfpart{%
  (\ref{thm:nec-fext-subdif:a})
  $\Rightarrow$
  (\ref{thm:nec-fext-subdif:b}):
}
Suppose $\fext(\xbar\plusl\zhat)\in\R$ for some $\zhat\in\Rn$.
Similar to the proof of \Cref{cor:subgrad-arb-reduc}, we can write
$\xbar=\ebar\plusl\qq$ for some $\ebar\in\corezn$ and $\qq\in\Rn$.
Let $g=\fshadd$.
Then
$g(\qq+\zhat)=\fext(\ebar\plusl\qq\plusl\zhat)=\fext(\xbar\plusl\zhat)\in\R$.
Therefore, $g$, which is convex and lower semicontinuous
(by \Cref{pr:i:9}\ref{pr:i:9a}),
must also be proper by
\Cref{pr:improper-vals}(\ref{pr:improper-vals:cor7.2.1}).
Hence,
$\uu\in\asubdiffext{\xbar}=\asubdiffext{\ebar\plusl\qq}=\asubdifgext{\qq}$
by \Cref{thm:subgrad-icon-reduc-proper},
implying $\ebar\cdot\uu=0$
by \Cref{pr:subg-reduc-perp}.
Thus, $\xbar\cdot\uu=\ebar\cdot\uu\plusl\qq\cdot\uu=\qq\cdot\uu\in\R$.

\pfpart{%
  (\ref{thm:nec-fext-subdif:b})
  $\Rightarrow$
  (\ref{thm:nec-fext-subdif:c}):
}
Since $\uu\in\asubdiffext{\xbar}$,
$\fstar(\uu)\in\R$
by
\Cref{pr:subgrad-imp-in-cldom}(\ref{pr:subgrad-imp-in-cldom:c}),
so
$\fext(\xbar)=\xbar\cdot\uu - \fstar(\uu)$
by
\Cref{thm:fenchel-subgrad}(\ref{thm:fenchel-subgrad:b},\ref{thm:fenchel-subgrad:c})
(and \Cref{pr:fextstar-is-fstar}).
Therefore, if $\xbar\cdot\uu$ is in $\R$, then so is $\fext(\xbar)$.

\pfpart{%
  (\ref{thm:nec-fext-subdif:c})
  $\Rightarrow$
  (\ref{thm:nec-fext-subdif:a}):
}
This is immediate.
\qedhere
\end{proof-parts}
\end{proof}

We next give sufficient conditions for when an extension $\fext$ is
astral subdifferentiable.
We show that if $\xbar\in\extspace$
can be expressed with finite part in $\ri(\dom{f})$,
then $\basubdiffext{\xbar}$ is nonempty if $\fext(\xbar)\in\R$.
This implies that the same holds
at all points where $\fext$ is continuous,
and holds everywhere if $f$ is finite everywhere.
\indexg{subdifferentials of extensions!directional derivatives and|(}%
\indexg{one-sided directional derivatives!an extension@of an extension|(}%
The theorem also shows that, under these same conditions,
$\fext$'s directional derivatives at $\xbar$ can be expressed in terms
of its astral subgradient at that point, analogous to
\Cref{roc:thm23.4-dd}.

\begin{theorem}   \label{thm:subgrad-not-empty-fin-ri}
  Let $f:\Rn\rightarrow\Rext$ be convex, and let $\xbar\in\extspace$.
  Suppose $\fext(\xbar)\in\R$ and that at least one of the following
  holds:
  \begin{roman-compact}
  \item   \label{thm:subgrad-not-empty-fin-ri:a}
    $\xbar\in\corezn\plusl\ri(\dom{f})$;
    that is, $\xbar=\ebar\plusl\qq$ for some
    $\ebar\in\corezn$ and $\qq\in\ri(\dom{f})$.
  \item   \label{thm:subgrad-not-empty-fin-ri:b}
    $\dom{f}=\Rn$.
  \item   \label{thm:subgrad-not-empty-fin-ri:c}
    $\fext$ is continuous at $\xbar$.
  \end{roman-compact}
  Then
    $\basubdiffext{\xbar}\neq\emptyset$.
  Furthermore,
    for $\yy\in\Rn$,
    \begin{equation}   \label{eq:thm:subgrad-not-empty-fin-ri:1}
       \dderfext{\xbar}{\yy}
       =
       \sup\bigBraces{\yy\cdot\uu :\: \uu\in\asubdiffext{\xbar}}.
    \end{equation}
\end{theorem}

\begin{proof}
  ~

\begin{proof-parts}
\pfpart{Condition~(\ref{thm:subgrad-not-empty-fin-ri:a}):}
Suppose $\xbar=\ebar\plusl\qq$ where
$\ebar\in\corezn$ and $\qq\in\ri(\dom{f})$.
Let $g=\fshadd$, which is convex and lower semicontinuous
(\Cref{pr:i:9}\ref{pr:i:9a}).
Let $h:\Rn\rightarrow\Rext$ be defined by
$h(\yy)=g(\qq+\yy)=\fext(\ebar\plusl\qq\plusl\yy)=\fext(\xbar\plusl\yy)$
for $\yy\in\Rn$,
which is also convex (\Cref{roc:thm5.7:fA}).
Then $h(\zero)=g(\qq)=\fext(\xbar)\in\R$, so
\begin{equation}   \label{eq:thm:subgrad-not-empty-fin-ri:2}
  \partial h(\zero) = \basubdiffext{\xbar}
\end{equation}
by
\Cref{cor:subgrad-arb-reduc}(\ref{cor:subgrad-arb-reduc:b}).
It further follows that $g$ is proper (by
\Cref{pr:improper-vals}\ref{pr:improper-vals:cor7.2.1}),
so $h$ is as well.

Also, we have
\[
  \qq
  \in
  \ri(\dom{f})
  \subseteq
  \ri(\dom{g})
  =
  \ri(\dom{h}) + \qq,
\]
where the second inclusion is
by \Cref{thm:ri-dom-icon-reduc},
and the equality is by
\Cref{pr:ri-props}(\ref{pr:ri-props:roc-cor6.6.2})
since $\dom{g}= (\dom{h}) + \qq$
(and since $\ri\{\qq\}=\{\qq\}$).
Thus, $\zero\in\ri(\dom{h})$.
Therefore, 
$\partial h(\zero)\neq\emptyset$
by \Cref{roc:thm23.4}(\ref{roc:thm23.4:a}),
proving
$\basubdiffext{\xbar}\neq\emptyset$
by \eqref{eq:thm:subgrad-not-empty-fin-ri:2}.

To prove \eqref{eq:thm:subgrad-not-empty-fin-ri:1},
let $\yy\in\Rn$.
Then
\[
  \dderfext{\xbar}{\yy}
  =
  \dderh{\zero}{\yy}
  =
  \sup\bigBraces{\yy\cdot\uu :\: \uu\in\partial h(0) }
  =
  \sup\bigBraces{\yy\cdot\uu :\: \uu\in\asubdiffext{\xbar}}.
\]
The first equality is by
\Cref{pr:ast-direc-deriv}(\ref{pr:ast-direc-deriv:a}),
the second by \Cref{roc:thm23.4-dd}
(since $\zero\in\ri(\dom h)$),
the third by
\eqref{eq:thm:subgrad-not-empty-fin-ri:2}.

\pfpart{Condition~(\ref{thm:subgrad-not-empty-fin-ri:b}):}
If $\dom{f}=\Rn$,
then
$\xbar\in\extspace=\corezn\plusl\ri(\dom{f})$
by \Cref{thm:icon-fin-decomp},
so the claim follows from the preceding argument for
condition~(\ref{thm:subgrad-not-empty-fin-ri:a}).

\pfpart{Condition~(\ref{thm:subgrad-not-empty-fin-ri:c}):}
Suppose $\fext$ is continuous at $\xbar$.
Let $\ell=\lsc f$, implying $\ellbar=\fext$ by
\Cref{pr:h:1}(\ref{pr:h:1aa}).
Then by
\Cref{cor:cont-gen-char}(\ref{cor:cont-gen-char:a},\ref{cor:cont-gen-char:c}),
$\xbar=\ebar\plusl\qq$ for some
$\ebar\in\represc{\ell}\cap\corezn \subseteq \corezn$
and some
$\qq\in\intdom{\ell}\subseteq\ri(\dom{\ell})=\ri(\dom{f})$,
where the equality is by
\Cref{pr:lsc-props}(\ref{pr:lsc-props:c}).
The claim 
therefore follows again by the argument above for
condition~(\ref{thm:subgrad-not-empty-fin-ri:a}).%
\indexg{subdifferentials of extensions!nonemptiness of|)}%
\indexg{subdifferentials of extensions!directional derivatives and|)}%
\indexg{one-sided directional derivatives!an extension@of an extension|)}%
\qedhere
\end{proof-parts}
\end{proof}

\section{Linear image of a function}

\indexg{linear image of standard function!subgradients of extension|(}%
\indexg{subdifferentials of extensions!linear image of function@of linear image of function|(}%
We next identify the astral subgradients of the extension of $\A f$,
the image of a function $f:\Rn\rightarrow\Rext$ under $\A\in\Rmn$,
as defined in
\eqref{eq:lin-image-fcn-defn}.
Let $h=\A f$.
As we saw in \Cref{thm:inf-lin-ext},
if a point $\xbar$ is in the astral column space of $\A$,
then there exists a point~$\zbar$ with $\A\zbar=\xbar$
and $\hext(\xbar)=\fext(\zbar)$, that is,
$\fext(\zbar)=\hext(\A\zbar)$.
Informally then,
in light of \Cref{thm:subgrad-fA},
we might expect $\uu$ to be an astral subgradient of
$\hext$ at $\xbar=\A\zbar$ when $\transAk\uu$
is an astral subgradient of $\fext$ at $\zbar$.
This heuristic reasoning about the astral subgradients of $\A f$
is made precise in the next theorem:

\begin{theorem}  \label{thm:subgrad-Af}
  Let $f:\Rn\rightarrow\Rext$ and $\A\in\Rmn$.
  Let $\xbar\in\extspac{m}$ and $\uu\in\Rm$.
  Then
    $\uu\in\basubdifAfext{\xbar}$
    if and only if
    there exists $\zbar\in\extspace$ such that
    $\A\zbar=\xbar$
    and
    $\transAk\uu \in \basubdiffext{\zbar}$.
\end{theorem}

\begin{proof}
Let $h=\A f$ and let
$\ww=\transAk\uu$.

We first make some preliminary observations.
Let $\fminusw$ and $\hminusu$ be linear tilts of $f$ and~$h$.
Then for $\xx\in\Rm$,
\begin{align}
  \hminusu(\xx)
  &=
  \inf {\bigBraces{ f(\zz) :\: \zz\in\Rn,\, \A\zz=\xx }}
        - \xx\cdot\uu
  \nonumber
  \\
  &=
  \inf {\bigBraces{ f(\zz) - \xx\cdot\uu :\: \zz\in\Rn,\, \A\zz=\xx }}
  \nonumber
  \\
  &=
  \inf {\bigBraces{ f(\zz) - \zz\cdot\ww :\: \zz\in\Rn,\, \A\zz=\xx }}
  \nonumber
  \\
  &=
  \inf {\bigBraces{ \fminusw(\zz) :\: \zz\in\Rn,\, \A\zz=\xx }}.
  \label{eq:thm:subgrad-Af:3}
\end{align}
The first and last equalities are by definition of linear tilt and
linear image of a function
(Eq.~\ref{eq:lin-image-fcn-defn}).
The third equality is because if $\A\zz=\xx$ then
\[
  \xx\cdot\uu = (\A\zz)\cdot\uu = \zz\cdot(\transAk\uu) = \zz\cdot\ww.
\]
Thus,
$\hminusu = \A \fminusw$.
This implies
\begin{equation}  \label{eq:thm:subgrad-Af:1}
   \hminusuext(\xbar)
   =
   \inf{\bigBraces{\fminuswext(\zbar) :\: \zbar\in\extspace,\, \A\zbar=\xbar}},
\end{equation}
from \Cref{thm:inf-lin-ext},
and furthermore that
\begin{equation}  \label{eq:thm:subgrad-Af:2}
  {-\hstar(\uu)}
  =
  \inf \hminusu
  =
  \inf \fminusw
  =
  -\fstar(\ww),
\end{equation}
where the first and last equalities are by
\Cref{pr:fminusu-props}(\ref{pr:fminusu-props:d}),
and the second is by
\Cref{pr:std-lin-img-props}(\ref{pr:std-lin-img-props:inf-Af}).

We now prove the theorem.
Suppose there exists $\zbar\in\extspace$ with
$\A\zbar=\xbar$ and
$\ww=\transAk\uu \in \asubdiffext{\zbar}$.
Then
$\hstar(\uu)=\fstar(\ww)\in\R$
by \eqref{eq:thm:subgrad-Af:2}
and
\Cref{pr:subgrad-imp-in-cldom}(\ref{pr:subgrad-imp-in-cldom:c}).
Also,
\[
  \inf \hminusu
  \leq
  \hminusuext(\xbar)
  \leq
  \fminuswext(\zbar)
  =
  \inf \fminusw
  =
  \inf \hminusu.
\]
The first inequality is from \Cref{pr:fext-min-exists}.
The second inequality
is from \eqref{eq:thm:subgrad-Af:1}, since $\A\zbar=\xbar$.
The first equality is from
\Cref{thm:fminus-subgrad-char}(\ref{thm:fminus-subgrad-char:a},\ref{thm:fminus-subgrad-char:e}),
since
$\ww\in\asubdiffext{\zbar}$.
And the last equality is from
\eqref{eq:thm:subgrad-Af:2}.
Thus,
$ \hminusuext(\xbar) =  \inf \hminusu$,
so $\uu\in\asubdifhext{\xbar}$ by
\Cref{thm:fminus-subgrad-char}(\ref{thm:fminus-subgrad-char:e},\ref{thm:fminus-subgrad-char:a}).

For the converse, suppose now that
$\uu\in\asubdifhext{\xbar}$.
Then
$\fstar(\ww)=\hstar(\uu)\in\R$
by \eqref{eq:thm:subgrad-Af:2}
and
\Cref{pr:subgrad-imp-in-cldom}(\ref{pr:subgrad-imp-in-cldom:c}).
In addition,
\Cref{thm:fminus-subgrad-char}(\ref{thm:fminus-subgrad-char:a},\ref{thm:fminus-subgrad-char:d})
implies $\hminusuext(\xbar)=-\hstar(\uu)$.
Since $\hstar(\uu)\in\R$,
$\xbar$ must be in the astral column space of $\A$
(since otherwise the infimum in Eq.~\ref{eq:thm:subgrad-Af:1} would be
vacuous).
Therefore,
by \Cref{thm:inf-lin-ext},
there exists $\zbar$ attaining the infimum in
\eqref{eq:thm:subgrad-Af:1} so that $\A\zbar=\xbar$
and $\hminusuext(\xbar)=\fminuswext(\zbar)$.
Combining, we then have
\[
  \fminuswext(\zbar)
  =
  \hminusuext(\xbar)
  =
  -\hstar(\uu)
  =
  -\fstar(\ww).
\]
Since $\fstar(\ww)\in\R$, this implies
$\transAk\uu = \ww\in\asubdiffext{\zbar}$
by
\Cref{thm:fminus-subgrad-char}(\ref{thm:fminus-subgrad-char:d},\ref{thm:fminus-subgrad-char:a}),
completing the proof.
\end{proof}

\indexg{subdifferentials, standard!linear image of function@of linear image of function|(}%
\indexg{linear image of standard function!standard subgradients of|(}%
An analogue of \Cref{thm:subgrad-Af} for standard subgradients
was given in
\Cref{pr:stan-subgrad-lin-img}.
However, that result requires an additional condition that the infimum
defining $(\A f)(\xx)$ in \eqref{eq:lin-image-fcn-defn}
be attained.
That condition is unnecessary in the astral case since the analogous
infimum is always attained (unless vacuous), as shown in
\Cref{thm:inf-lin-ext}.
Without this additional condition, the analogue of
\Cref{thm:subgrad-Af} for standard subgradients
does not hold in general.
In other words, with $f$ and $\A$ as in the theorem, and
for $\xx,\uu\in\Rm$, it is not the case that
$\uu\in\partial \Af(\xx)$
if and only if there exists $\zz\in\Rn$ such that
$\A\zz=\xx$ and $\transAk\uu \in \partial f(\zz)$,
as shown in the next example.

\begin{example}
Let $f:\R^2\rightarrow\R$ be defined by
$f(\zz)=f(z_1,z_2)=e^{z_1}$ for $\zz\in\R^2$, and let
$\A = [0,1] = \trans{\ee_2}$
(so $m=1$ and $n=2$).
Then
\[\Af(x)=\inf\{e^{z_1} :\: \zz\in\R^2, z_2 = x\} = 0\]
for all $x\in\R$,
so $\A f\equiv 0$.
Now let $x=u=0$.
Then $\nabla\Af(0)=0$, so
$u\in\partial \Af(x)$.
On the other hand, for all $\zz\in\R^2$,
$\gradf(\zz)=\trans{[e^{z_1},0]}\neq \zero$.
Therefore,
$\transAk u = \zero$ cannot be in $\partial f(\zz)$
for any $\zz\in\R^2$.

Although its analogue for standard subgradients is false for this
example,
\Cref{thm:subgrad-Af} is of course true in the astral setting.
In this case, $u\in\asubdifAfext{x}$, and, setting
$\zbar=\limray{(-\ee_1)}$, we have
$\A\zbar = 0 = x$ and
$\transAk u = \zero \in \asubdiffext{\zbar}$
(by \Cref{pr:asub-zero-is-min}\ref{pr:asub-zero-is-min:a}
since $\fext$ is minimized by $\zbar$).%
\indexg{linear image of standard function!subgradients of extension|)}%
\indexg{subdifferentials of extensions!linear image of function@of linear image of function|)}%
\indexg{subdifferentials, standard!linear image of function@of linear image of function|)}%
\indexg{linear image of standard function!standard subgradients of|)}%
\end{example}

\section{Pointwise supremum of a collection of functions}

Next, we study the astral subgradients of the extension of
a function $h$ that is the
pointwise supremum of a collection of functions, that is,
for $\xx\in\Rn$,
\begin{equation}  \label{eqn:point-sup-def}
  h(\xx) = \sup_{i\in\indset} f_i(\xx),
\end{equation}
where $f_i:\Rn\rightarrow\Rext$ for $i$ in some index
set $\indset$.

\indexg{subdifferentials, standard!pointwise supremum@of pointwise supremum|(}%
\indexg{suprema and maxima, pointwise!standard subgradients of|(}%
In standard convex analysis and under suitable conditions,
many or all of
the subgradients of $h$ at a point $\xx\in\Rn$ can be determined from
the subgradients of all those functions $f_i$ at which
the supremum in \eqref{eqn:point-sup-def} is attained.
More specifically, if $h(\xx)=f_i(\xx)$, then
every subgradient of $f_i$ at $\xx$ is also
a subgradient of $h$ at $\xx$ so that
$\partial f_i(\xx) \subseteq \partial h(\xx)$.
Furthermore, since $\partial h(\xx)$ is closed and convex, this
implies
\begin{equation*}  %
  \cl \conv\!\Parens{\bigcup_{i\in I:\: h(\xx)=f_i(\xx)}
                     \mkern-20mu
                     \partial f_i(\xx)
                    }
  \subseteq
  \partial h(\xx),
\end{equation*}
as we saw in
\Cref{pr:stnd-sup-subgrad}.

Under some conditions, this can account for all of $h$'s
subgradients;
for instance, this is the case if each $f_i$ is convex and finite
everywhere, and if $\indset$ is finite,
as we saw in
\Cref{pr:std-max-fin-subgrad}.

\begin{example}  \label{ex:subgrad-max-linear}
  Suppose $h=\max\{f,g\}$
  where,
  for $x\in\R$,
  $f(x)=x$ and $g(x)=2x$.
  Then $\partial f(x)=\{1\}$, $\partial g(x)=\{2\}$
  and, for $x\in\R$,
  \[
    \partial h(x)
    =
    \begin{cases}
        \{1\}  &  \mbox{if $x<0$,} \\
        {[1,2]}  &  \mbox{if $x=0$,} \\
        \{2\}  &  \mbox{if $x>0$.}
    \end{cases}
  \]
  Thus, for points $x<0$, $h(x)=f(x)$ and $h$'s derivative
  (unique subgradient) is equal to
  $f$'s; similarly at points $x>0$.
  At $x=0$, $h(x)=f(x)=g(x)$, and
  $\partial h(0)$ is equal to the convex hull of $f$'s derivative and $g$'s
  derivative.%
\indexg{subdifferentials, standard!pointwise supremum@of pointwise supremum|)}%
\indexg{suprema and maxima, pointwise!standard subgradients of|)}%
\end{example}

\indexg{subdifferentials of extensions!pointwise supremum@of pointwise supremum|(}%
\indexg{suprema and maxima, pointwise!subgradients of extension|(}%
For astral subgradients, we expect a similar relationship between
$\asubdifhext{\xbar}$ and $\asubdiffextsub{i}{\xbar}$ for those
$i\in\indset$ with $\hext(\xbar)=\fextsub{i}(\xbar)$.
Indeed, if $\hext(\xbar)=\fextsub{i}(\xbar)$ and if this value is
finite, then
$\asubdiffextsub{i}{\xbar} \subseteq \asubdifhext{\xbar}$, as
expected.
But if $\hext(\xbar)=\fextsub{i}(\xbar)\in\{-\infty,+\infty\}$, then
this might not be true.
For instance, for the functions in
Example~\ref{ex:subgrad-max-linear},
$\hext(+\infty)=\fext(+\infty)=+\infty$
and $\asubdiffext{+\infty}=\{1\}$,
but $\asubdifhext{+\infty}=\{2\}$,
so $\asubdiffext{+\infty}\not\subseteq\asubdifhext{+\infty}$.
This demonstrates that
when their common value is infinite,
requiring merely that
$\hext(\xbar)=\fextsub{i}(\xbar)$ is too coarse
to ensure that
$\asubdiffextsub{i}{\xbar} \subseteq \asubdifhext{\xbar}$.

Nonetheless, although
$\hext(+\infty)=\fext(+\infty)=\gext(+\infty)=+\infty$
in this example,
it is also true that the difference between $h(x)=g(x)$ and $f(x)$
grows to $+\infty$ as $x\to+\infty$. Therefore, it feels natural
that we should only consider $\gext$, but not $\fext$, as matching
the value of $\hext$ at $\barx=+\infty$. This is analogous to
the definition of astral primal subgradients, where
an affine function is considered to meet an astral function at 
a point $\xbar$ only when
the gap between the two vanishes along some sequence converging to $\xbar$.

To formalize the intuition above for the analysis of the
pointwise supremum of a collection of functions, we
leverage the framework of linear tilts.
Let $\hminusu$ and $\fminususub{i}$ be linear tilts
of $h$ and $f_i$ by $\uu\in\Rn$, for $i\in\indset$.
As just discussed, for
$\asubdiffextsub{i}{\xbar} \subseteq \asubdifhext{\xbar}$
to hold,
we will prove it is sufficient that
$\hext(\xbar)=\fextsub{i}(\xbar)\in\R$.
Much more generally, we will see in fact
that it is also sufficient if this holds for any
linear tilt, that is, if
$\hminusuext(\xbar)=\fminususubext{i}(\xbar)\in\R$,
for any $\uu\in\Rn$.

For instance, for the functions from
Example~\ref{ex:subgrad-max-linear},
with $u=2$,
$\fminusU(x)=-x$,
$\gminusU(x)=0$,
and
$\hminusU(x)=\max\{-x,0\}$,
for $x\in\R$.
So
$\hminusUext(+\infty)=\gminusUext(+\infty)=0\in\R$,
and indeed,
$\asubdifgext{+\infty}=\{2\}\subseteq\asubdifhext{+\infty}$.
On the other hand, we cannot have
$\hminusUext(+\infty)=\fminusUext(+\infty)\in\R$
for any $u\in\R$.
Together, these facts capture the
intuition discussed above regarding the attainment of the
pointwise maximum at those points where $\eh$'s value is infinite.

The next theorem brings these ideas together:

\begin{theorem}   \label{thm:sup-subgrad-subset}
  Let $f_i:\Rn\rightarrow\Rext$ for all $i\in\indset$, and let
  $h = \sup_{i\in\indset} f_i$.
  Let $\hminusu$ and $\fminususub{i}$ be linear tilts of $h$ and
  $f_i$ by $\uu\in\Rn$.
  Let $\xbar\in\extspace$, and let
  \begin{equation}   \label{eq:thm:sup-subgrad-subset:J-dfn}
    J =
    \bigBraces{i \in \indset :\:
               \exists \uu\in\Rn,\,
                   \hminusuext(\xbar)=\fminususubext{i}(\xbar)\in\R }.
  \end{equation}
  Then
  \begin{equation}   \label{eq:thm:sup-subgrad-subset:0}
         \conv\BiggParens{ \bigcup_{i\in J} \basubdiffextsub{i}{\xbar} }
           \subseteq
             \basubdifhext{\xbar}.
  \end{equation}
\end{theorem}

\begin{proof}
{\mathtogether%
Let $\ww\in\asubdiffextsub{i}{\xbar}$ for some $i\in J$.
We first aim to show that $\ww\in\asubdifhext{\xbar}$.
From $J$'s definition, there exists $\uu\in\Rn$ for which
\begin{equation}   \label{eq:thm:sup-subgrad-subset:1}
  \hminusuext(\xbar)=\fminususubext{i}(\xbar)\in\R.
\end{equation}
Note that
$\ww-\uu\in\asubdiffminususubext{i}{\xbar}$
by \Cref{pr:fminusu-subgrad}(\ref{pr:fminusu-subgrad:ext}),
further implying that
$\fiminusuextstar(\ww-\uu)\in\R$ by
\Cref{pr:subgrad-imp-in-cldom}(\ref{pr:subgrad-imp-in-cldom:a}).}

We then have
\begin{align}
  \hstar(\ww)
  \leq
  \fistar(\ww)
  =
  \fiminusustar(\ww-\uu)
  &=
  \fiminusuextstar(\ww-\uu)
  \nonumber
  \\
  &=
  \xbar\cdot(\ww-\uu) - \fminususubext{i}(\xbar)
  \nonumber
  \\
  &=
  \xbar\cdot(\ww-\uu) - \hminusuext(\xbar)
  \nonumber
  \\
  &\leq
  \hminusuextstar(\ww-\uu)
  =
  \hminusustar(\ww-\uu)
  =
  \hstar(\ww).
  \label{eq:thm:sup-subgrad-subset:2}
\end{align}
The first inequality is because,
by $h$'s definition, $f_i\leq h$, implying
$\hstar\leq\fistar$
(by \Cref{pr:conj-props}\ref{pr:conj-props:b}).
The first and last equalities are by
\Cref{pr:fminusu-props}(\ref{pr:fminusu-props:b}).
The second and fifth equalities are by
\Cref{pr:fextstar-is-fstar}.
The summability of the terms appearing on the second and third lines,
as well as the equality of these expressions, follow from
\eqref{eq:thm:sup-subgrad-subset:1}.
Since $\ww-\uu\in\asubdiffminususubext{i}{\xbar}$,
the third equality follows
from
\Cref{thm:fenchel-subgrad}(\ref{thm:fenchel-subgrad:b},\ref{thm:fenchel-subgrad:a}),
noting that
$\fiminusuextstar(\ww-\uu)$ and 
$\fminususubext{i}(\xbar)$ are both in $\R$.
The last inequality is by definition of astral conjugate
(Eq.~\ref{eq:Fstar-down-def}).

Thus,
\eqref{eq:thm:sup-subgrad-subset:2}
implies that $\hminusuextstar(\ww-\uu)=\fiminusuextstar(\ww-\uu)\in\R$,
and also that
\[
  \hminusuextstar(\ww-\uu)
  =
  \xbar\cdot(\ww-\uu) - \hminusuext(\xbar).
\]
Therefore,
$\ww-\uu\in\asubdifhminusuext{\xbar}$
by
\Cref{thm:fenchel-subgrad}(\ref{thm:fenchel-subgrad:a},\ref{thm:fenchel-subgrad:b}).
In turn, this implies
$\ww\in\asubdifhext{\xbar}$
by Proposi\-tion~\ref{pr:fminusu-subgrad}(\ref{pr:fminusu-subgrad:ext}).

We conclude that
$\asubdiffextsub{i}{\xbar}\subseteq\asubdifhext{\xbar}$
for all $i\in J$, that is,
$\bigcup_{i\in J} \asubdiffextsub{i}{\xbar}\subseteq\asubdifhext{\xbar}$.
Since
$\asubdifhext{\xbar}$ is convex by
\Cref{thm:asubdifF-convex},
\eqref{eq:thm:sup-subgrad-subset:0}
now follows by
\Cref{pr:conhull-prop}(\ref{pr:conhull-prop:aa}).
\end{proof}

Note that, in the notation of \Cref{thm:sup-subgrad-subset},
if $\xbar=\xx\in\Rn$, then 
$\fminususubext{i}(\xx)=\efsub{i}(\xx)-\xx\cdot\uu$
and
$\hminusuext(\xbar)=\hext(\xx)-\xx\cdot\uu$
by \Cref{pr:fminusu-props}(\ref{pr:fminusu-props:e}).
Therefore, in this case, the set $J$
given in \eqref{eq:thm:sup-subgrad-subset:J-dfn}
is simply equal to
$J = \Braces{i \in \indset :\: \hext(\xx)=\fext_i(\xx)\in\R}$.

In general, under the conditions of \Cref{thm:sup-subgrad-subset},
\eqref{eq:thm:sup-subgrad-subset:0} need not hold with equality.
\indexg{subdifferentials, standard!pointwise supremum@of pointwise supremum|(}%
\indexg{suprema and maxima, pointwise!standard subgradients of|(}%
For standard subgradients as well, the analogous inclusion given in
\Cref{pr:stnd-sup-subgrad}
need not generally hold with equality.
Here is an example:

\begin{example}
Let $I=\{0,1,2,\ldots\}$, and
for $i\in I$, define $f_i:\R\rightarrow\Rext$
to be the
constant function $f_i\equiv -1/i$ if $i>0$, and otherwise,
$f_0(x)=x$ for $x\in\R$.
Let
$h=\sup_{i\in I} f_i$.
Then $h(x)=\max\{0,x\}$ for $x\in\R$.
Let $x=0$.
Then the set $J$ in
\eqref{eq:thm:sup-subgrad-subset:J-dfn}
is the singleton $\{0\}$ since
$\hext(0)=\efsub{0}(0)=0$, but for all $u\in\R$ and all $i\in I\posKern\wo\{0\}$,
$\fminusUsubext{i}(0)=-1/i$
but $\hminusUext(0)=0$.
Thus, the left-hand side of \eqref{eq:thm:sup-subgrad-subset:0} is
equal to $\asubdiffextsub{0}{0}=\partial f_0(0)=\{1\}$,
whereas $\asubdifhext{0}=\partial h(0)=[0,1]$
(by \Cref{pr:asubdiffext-at-x-in-rn}\ref{pr:asubdiffext-at-x-in-rn:c}
and direct calculation of standard subgradients).
Thus, \eqref{eq:thm:sup-subgrad-subset:0} does not hold with equality
in this case.
This reasoning shows moreover that
the inclusion appearing in
\Cref{pr:stnd-sup-subgrad}
for standard subgradients also
does not hold with equality in this
\indexg{subdifferentials, standard!pointwise supremum@of pointwise supremum|)}%
\indexg{suprema and maxima, pointwise!standard subgradients of|)}%
case.
\end{example}

Nonetheless,
analogous to \Cref{pr:std-max-fin-subgrad},
when $I$ is finite and each $f_i$ is convex and finite everywhere, the
set given in
\eqref{eq:thm:sup-subgrad-subset:0}
does exactly capture all of $\hext$'s astral subgradients:

\begin{theorem}   \label{thm:ben-subgrad-max}
  Let $f_i:\Rn\rightarrow\R$ be convex, for $i=1,\dotsc,m$,
  and let
  \[  h=\max\bigBraces{f_1,\dotsc,f_m}. \]
  Let $\xbar\in\extspace$, and let
  \[
    J =
    \bigBraces{i \in \{1,\ldots,m\} :\:
               \exists \uu\in\Rn,\,
                   \hminusuext(\xbar)=\fminususubext{i}(\xbar)\in\R }.
  \]
  Then
  \[
     \basubdifhext{\xbar}
     =
     \conv\BiggParens{ \bigcup_{i\in J} \basubdiffextsub{i}{\xbar} }.
  \]
\end{theorem}

\begin{proof}
It suffices to prove
\begin{equation}   \label{eq:thm:ben-subgrad-max:1}
   \basubdifhext{\xbar}
   \subseteq
   \conv\BiggParens{ \bigcup_{i\in J} \basubdiffextsub{i}{\xbar} }
\end{equation}
since the reverse inclusion follows from
\Cref{thm:sup-subgrad-subset}.
Let $\uu\in\basubdifhext{\xbar}$, which we aim to show is in the set
on the right-hand side of \eqref{eq:thm:ben-subgrad-max:1}.

By a standard result regarding the conjugate of a maximum
(\Cref{pr:prelim-conj-max}),
there exist a subset $M\subseteq\{1,\ldots,m\}$, and numbers
$\lambda_i\in[0,1]$ and vectors $\uu_i\in\Rn$, for $i\in M$, such that
\begin{align}
  &\sum_{i\in M} \lambda_i
  =
  1,
  \label{eq:thm:ben-subgrad-max:2a}
  \\
  &\sum_{i\in M} \lambda_i \uu_i
  =
  \uu,
  \label{eq:thm:ben-subgrad-max:2}
  \\
  &\sum_{i\in M} \lambda_i \fistar(\uu_i)
  =
  \hstar(\uu).
  \label{eq:thm:ben-subgrad-max:8}
\end{align}
Furthermore, without loss of generality, we can assume $\lambda_i>0$
for all $i\in M$ (since discarding elements with $\lambda_i=0$ does not
affect these properties).
Since $\hstar(\uu)\in\R$
(by \Cref{pr:subgrad-imp-in-cldom}\ref{pr:subgrad-imp-in-cldom:c}),
this implies $\fistar(\uu_i)\in\R$ for
$i\in M$.

The idea of the proof is to show that
$\uu_i\in\basubdiffextsub{i}{\xbar}$
for $i\in M$, and also that $M\subseteq J$.
Together, these will prove the theorem.

To show all this, we define a function $g:\Rn\rightarrow\R$ whose
extension will be computed in two different ways.
For $\xx\in\Rn$, let
\[
  g(\xx)
  =
  \sum_{i\in M} \lambda_i \fminususubd{i}(\xx).
\]
Note that, for $i\in M$,
\begin{equation}   \label{eq:thm:ben-subgrad-max:4}
 \fminususubdext{i}(\xbar) \geq -\fistar(\uu_i) > -\infty
\end{equation}
by \Cref{pr:fminusu-props}(\ref{pr:fminusu-props:d}),
and since $\fistar(\uu_i)\in\R$.
Therefore, the values $\fminususubdext{i}(\xbar)$, for $i\in M$, are
summable.
By \Cref{thm:ext-sum-fcns-w-duality}
and \Cref{pr:scal-mult-ext}, this implies
\begin{equation}   \label{eq:thm:ben-subgrad-max:5}
  \gext(\xbar)
  =
  \sum_{i\in M} \lambda_i \fminususubdext{i}(\xbar)
\end{equation}
(using $\dom{\fminususubd{i}}=\Rn$, for $i\in M$).

We can re-express $g$ as
\begin{equation}  \label{eq:thm:ben-subgrad-max:3}
  g(\xx)
  =
  \sum_{i\in M} \lambda_i \bigParens{f_i(\xx) - \xx\cdot\uu_i}
  =
  \sum_{i\in M} \lambda_i \bigParens{f_i(\xx) - \xx\cdot\uu}
  =
  \sum_{i\in M} \lambda_i \fminususub{i}(\xx),
\end{equation}
for $\xx\in\Rn$, where the first and last equality are by definition
of linear tilt, and the second follows from
Eqs.~(\ref{eq:thm:ben-subgrad-max:2a})
and~(\ref{eq:thm:ben-subgrad-max:2}).
For $i\in M$,
\begin{equation}   \label{eq:thm:ben-subgrad-max:7}
  \fminususubext{i}(\xbar)
  \leq
  \hminusuext(\xbar)
  =
  -\hstar(\uu)
  <
  +\infty.
\end{equation}
The first inequality is because $f_i\leq h$, so
$\fminususub{i}(\xx)
 = f_i(\xx) - \xx\cdot\uu
 \leq h(\xx) - \xx\cdot\uu
 = \hminusu(\xx)$
for $\xx\in\Rn$, and hence also
$\fminususubext{i}\le \hminusuext$ by
\Cref{pr:h:1}(\ref{pr:h:1:geq}).
The equality is by
\Cref{thm:fminus-subgrad-char}(\ref{thm:fminus-subgrad-char:a},\ref{thm:fminus-subgrad-char:d}),
since $\uu\in\asubdifhext{\xbar}$.
The last inequality is because $\hstar(\uu)\in\R$.

Thus, the values $\fminususubext{i}(\xbar)$, for $i\in M$, are
summable, so we can again apply
\Cref{thm:ext-sum-fcns-w-duality}
and \Cref{pr:scal-mult-ext} to
\eqref{eq:thm:ben-subgrad-max:3}
yielding
\begin{equation}   \label{eq:thm:ben-subgrad-max:6}
  \gext(\xbar)
  =
  \sum_{i\in M} \lambda_i \fminususubext{i}(\xbar).
\end{equation}

Combining, we now have:
\begin{align}
  -\hstar(\uu)
  =
  -\sum_{i\in M} \lambda_i \fistar(\uu_i)
  &\leq
  \sum_{i\in M} \lambda_i \fminususubdext{i}(\xbar)
  \nonumber
  \\
  &=
  \gext(\xbar)
  \nonumber
  \\
  &=
  \sum_{i\in M} \lambda_i \fminususubext{i}(\xbar)
  \leq
  \hminusuext(\xbar)
  =
  -\hstar(\uu).
  \label{eq:thm:ben-subgrad-max:9}
\end{align}
The first equality is
\eqref{eq:thm:ben-subgrad-max:8}.
The first inequality is by \eqref{eq:thm:ben-subgrad-max:4}.
The second and third equalities are by
Eqs.~(\ref{eq:thm:ben-subgrad-max:5})
and~(\ref{eq:thm:ben-subgrad-max:6}).
And the last inequality and last equality are both by
\eqref{eq:thm:ben-subgrad-max:7}.

Thus,
\[
  \sum_{i\in M} \lambda_i \fminususubdext{i}(\xbar)
  =
  -\hstar(\uu)
  =
  -\sum_{i\in M} \lambda_i \fistar(\uu_i),
\]
with first equality from
\eqref{eq:thm:ben-subgrad-max:9}, and second from
\eqref{eq:thm:ben-subgrad-max:8}.
In view of \eqref{eq:thm:ben-subgrad-max:4}, and since
$\hstar(\uu)\in\R$ and every $\lambda_i>0$, this implies that
$\fminususubdext{i}(\xbar) = -\fistar(\uu_i)$ for $i\in M$.
Therefore,
$\uu_i\in\basubdiffextsub{i}{\xbar}$,
by
\Cref{thm:fminus-subgrad-char}(\ref{thm:fminus-subgrad-char:d},\ref{thm:fminus-subgrad-char:a})
and since $\fistar(\uu_i)\in\R$.

\eqref{eq:thm:ben-subgrad-max:9} also implies that
\[
  \sum_{i\in M} \lambda_i \fminususubext{i}(\xbar)
  =
  \hminusuext(\xbar)
  =
  -\hstar(\uu).
\]
In view of \eqref{eq:thm:ben-subgrad-max:7},
and again since $\hstar(\uu)\in\R$ and every $\lambda_i>0$, 
this implies
$ \fminususubext{i}(\xbar) = \hminusuext(\xbar)$
for $i\in M$, and thus that $M\subseteq J$.

By the form of $\uu$ given in \eqref{eq:thm:ben-subgrad-max:2},
we can therefore conclude that
\[
  \uu
  \in
  \conv\BiggParens{ \bigcup_{i\in M} \basubdiffextsub{i}{\xbar} }
  \subseteq
  \conv\BiggParens{ \bigcup_{i\in J} \basubdiffextsub{i}{\xbar} },
\]
completing the proof.%
\indexg{subdifferentials of extensions!pointwise supremum@of pointwise supremum|)}%
\indexg{suprema and maxima, pointwise!subgradients of extension|)}%
\end{proof}

\section{Composition with a nondecreasing function}
\label{sec:subgrad-comp-inc-fcn}

\indexg{subdifferentials of extensions!composition with nondecreasing function@of composition with nondecreasing function|(}%
\indexg{nondecreasing function, composition with!astral subgradients of|(}%
We next consider the astral subgradients of the extension of the
composition of a convex, nondecreasing function
$g:\R\rightarrow\Rext$
with another convex function $f:\Rn\rightarrow\Rext$.
Since $f$ may take values of $\pm\infty$, for this to make sense, we
compose $f$ with the extension $\eg$.
Letting $h=\eg\circk f$,
our aim then is to determine the astral subgradients of $\hext$.

In standard calculus, assuming differentiability,
the gradient of $h$ can be
computed using the chain rule,
  $\gradh(\xx)=g'(f(\xx))\posKern\gradf(\xx)$,
where $g'$ is the
derivative of~$g$.
Analogous results hold for standard subgradients, as we saw in
\Cref{pr:std-subgrad-comp-inc}.
Correspondingly, for astral subgradients, if
$v\in\asubdifgext{\fext(\xbar)}$ and
$\uu\in\asubdiffext{\xbar}$, 
we expect that the product $v \uu$ should be an astral
subgradient of $\hext$ at $\xbar$.
Indeed, this is so, as shown in the next theorem under the stated
conditions:

\begin{theorem}  \label{thm:uv-cond-suff-comp-subgrad}
  Let $f:\Rn\rightarrow\Rext$ be convex, and let
  $g:\R\rightarrow\Rext$ be
  convex, proper, closed, nondecreasing, and not constant.
  Assume there exists $\xing\in\Rn$ such that
  $f(\xing)<\sup(\dom g)$.
  Let $h=\eg\circk f$,
  and let $\xbar\in\extspace$.

  Suppose
  $\uu\in\basubdiffext{\xbar}$
  and $v\in\basubdifgext{\fext(\xbar)}$.
  Then $v \uu \in\basubdifhext{\xbar}$.
\end{theorem}

As a preliminary step to proving the theorem,
we first observe
that the astral
subgradients of any nondecreasing function on $\Rext$ are all
nonnegative:

\begin{proposition}  \label{pr:subgrad-inc-not-neg}
  Let
  $G:\Rext\rightarrow\Rext$ be nondecreasing.
  Suppose $v\in\asubdifG{\barx}$ for some $\barx\in\Rext$.
  Then $v\geq 0$.
\end{proposition}

\begin{proof}
\mathtogether%
By
\Cref{pr:subgrad-imp-in-cldom}(\ref{pr:subgrad-imp-in-cldom:a}),
$\Gstar(v)\in\R$ 
since $v\in\asubdifG{\barx}$.
This implies $G\not\equiv+\infty$,
since otherwise we would have $\Gstar\equiv-\infty$.
Thus, $G(-\infty)<+\infty$, since $G$ is nondecreasing.
Hence, by definition of astral conjugate
(Eq.~\ref{eq:Fstar-down-def}),
if $v<0$, then we would have
$\Gstar(v)\ge-G(-\infty)\plusd (-\infty)v=+\infty$
(since $-G(-\infty)>-\infty$),
contradicting that $\Gstar(v)\in\R$.
Therefore, $v\geq 0$.
\end{proof}

\begin{proof}[Proof of \Cref{thm:uv-cond-suff-comp-subgrad}]
Let $\ww=v\uu$.
We aim to prove that $\ww\in\basubdifhext{\xbar}$.

Note first, for $\bary\in\Rext$, that
\begin{align}
\SwapAboveDisplaySkip
\label{eq:thm:uv-cond-suff-comp-subgrad:3}
  \gext(\bary)
  &=
  \begin{cases}
    \inf g     & \text{if $\bary=-\infty$,} \\
    g(\bary)   & \text{if $\bary\in\R$,} \\
    +\infty    & \text{if $\bary=+\infty$,}
  \end{cases}
\end{align}
by
Propositions~\ref{pr:h:1}(\ref{pr:h:1a})
and~\ref{pr:conv-inc:prop}(\ref{pr:conv-inc:infsup},\ref{pr:conv-inc:nonconst}).
Thus, $\gext$ is nondecreasing.

Since
$v\in\basubdifgext{\fext(\xbar)}$ and since $\gext$ is nondecreasing,
$v$ must be nonnegative by 
\Cref{pr:subgrad-inc-not-neg}.
We consider separately the cases that $v=0$ or $v>0$:

\begin{proof-parts}
\pfpart{Case $v=0$:}
Since $0\in\basubdifgext{\fext(\xbar)}$,
by \Cref{pr:asub-zero-is-min}(\ref{pr:asub-zero-is-min:a}),
\begin{equation}  \label{eq:thm:uv-cond-suff-comp-subgrad:2}
  \gext(\fext(\xbar))=\inf \gext\in\R.
\end{equation}
We then have
\[
  \inf h
  \leq
  \hext(\xbar)
  =
  \gext(\fext(\xbar))
  =
  \inf \gext
  \leq
  \inf h.
\]
The first inequality is by \Cref{pr:fext-min-exists}.
The first equality is
by \Cref{thm:Gf-conv}, due to $\xing$'s assumed existence.
The second equality is by
\eqref{eq:thm:uv-cond-suff-comp-subgrad:2}.
The last inequality is because, for all $\xx\in\Rn$,
$h(\xx)=\gext(f(\xx))\geq\inf\gext$.
Thus, $\hext(\xbar)=\inf h$, and moreover,
$\inf h=\inf\gext\in\R$ by
\eqref{eq:thm:uv-cond-suff-comp-subgrad:2}.
Therefore, $\ww=\zero\in \asubdifhext{\xbar}$
by \Cref{pr:asub-zero-is-min}(\ref{pr:asub-zero-is-min:a}).

\pfpart{Case $v>0$:}
We adopt a similar setup to the one used in proving
\Cref{thm:conj-compose-our-version}.
As such, let $\PPx=[\Idnn,\zerov{n}]$ and
$\PPy=[\trans{\zerov{n}},1]$,
and define functions
$r=g\PPy$ and $s=r+\indepif$.
Then $r$ is convex (by \Cref{roc:thm5.7:fA}),
and $s$ is convex (by \Cref{pr:std-sum-fcns-cvx},
and since $f$ is
convex, implying $\epi{f}$ is convex).

These definitions imply that if $\zz=\rpair{\xx}{y}\in\Rn\times\R$,
then
$\PPx \zz = \xx$,
$\PPy \zz = y$,
$r(\zz)=g(y)$,
and
$s(\zz)$ equals $g(y)$ if $\zz\in\epi{f}$ and $+\infty$ otherwise.
We can then rewrite~$h$, for $\xx\in\Rn$, as
\begin{align}
  \notag
  h(\xx)
    &=
    \inf{\bigBraces{\eg(\ey) :\: \ey\in\eR,\, \ey\geq f(\xx) }}
  \\
  \notag
    &=
    \inf{\bigBraces{ g(y) :\: y\in\R,\, y\geq f(\xx) }}
  \\
  \notag
    &=
    \inf{\bigBraces{ g(y) + \indepif({\xx},{y}) :\: y\in\R }}
  \\
  \notag
    &=
    \inf{\bigBraces{ r(\zz) + \indepif(\zz) :\: \zz\in\R^{n+1},\, \PPx\zz=\xx }}
  \\
  \label{eq:thm:subgrad-comp-inc-fcn:1}
    &=
    \inf{\bigBraces{ s(\zz) :\: \zz\in\R^{n+1},\, \PPx\zz=\xx }}.
\end{align}
The first equality is because $\gext$ is nondecreasing.
The second equality follows from
\eqref{eq:thm:uv-cond-suff-comp-subgrad:3}; note that this holds even
if $f(\xx)$ is $\pm\infty$.
Once written in this form, $h$'s astral subgradients can largely be
computed using rules developed in the preceding sections, as we show
now.

Since $\uu\in\asubdiffext{\xbar}$,
$\fstar(\uu)\in\R$ by
\Cref{pr:subgrad-imp-in-cldom}(\ref{pr:subgrad-imp-in-cldom:c}),
and furthermore,
there exists $\zbar\in \clepi{f}$ such that
$\PPx\zbar = \xbar$,
$\PPy\zbar = \fext(\xbar)$,
and
$\zbar \cdot \rpair{\uu}{-1} = \fstar(\uu)$,
as follows from
\Cref{pr:equiv-ast-subdif-defn}(\ref{pr:equiv-ast-subdif-defn:a},\ref{pr:equiv-ast-subdif-defn:c})
(noting that $\clepifext = \clepi{f}$
by
\Cref{pr:wasthm:e:3}\ref{pr:wasthm:e:3c}
and that $\fextstar=\fstar$ by
\Cref{pr:fextstar-is-fstar}).
Multiplying by $v$, this last fact implies
\begin{equation}   \label{eq:thm:uv-cond-suff-comp-subgrad:4}
  \zbar \cdot \rpair{v\uu}{-v}
  =
  v \fstar(\uu)
  =
  \indepifstar(v\uu,-v),
\end{equation}
with the second equality from
\Cref{pr:support-epi-f-conjugate}.
Noting that $\zbar\in \clepi{f}$ and that
$\indepifstar(v\uu,-v)\in\R$
(from Eq.~\ref{eq:thm:uv-cond-suff-comp-subgrad:4}
since $\fstar(\uu)\in\R$),
it now follows from
\Cref{cor:subdif-ind-fcn-rn}(\ref{cor:subdif-ind-fcn-rn:b},\ref{cor:subdif-ind-fcn-rn:a})
that
\begin{equation}  \label{eq:thm:subgrad-comp-inc-fcn:6}
  \rpair{\ww}{-v}
  =
  \rpair{v\uu}{-v}
  \in
  \asubdifindepifext{\zbar}
\end{equation}
(here applied with $S$, $\xbar$ and $\uu$, as they appear in that
\namecref{cor:subdif-ind-fcn-rn}, set to $\epi{f}$, $\zbar$ and
$\rpair{\ww}{-v}$).

By definition,
$r=g\PPy$.
Therefore,
\begin{equation}  \label{eq:thm:subgrad-comp-inc-fcn:7}
  \rpair{\zero}{v}
  =
  \trans{\PPy} v
  \in
  \trans{\PPy} \asubdifgext{\PPy\zbar}
  =
  \asubdifrext{\zbar}.
\end{equation}
The inclusion is because
$v\in\asubdifgext{\fext(\xbar)}=\asubdifgext{\PPy\zbar}$.
The last equality is by \Cref{thm:subgrad-fA}
(applied with $f$ and $\A$, as they appear in that theorem, set to $g$
and $\PPy$),
noting that, by
\Cref{lem:thm:conj-compose-our-version:1}(\ref{lem:thm:conj-compose-our-version:1:a}),
$\PPy\zing\in\ri(\dom{g})$
for some $\zing\in\R^{n+1}$.

Thus,
\begin{equation}  \label{eq:thm:subgrad-comp-inc-fcn:3}
  \trans{\PPx} \ww
  =
  \rpair{\ww}{0}
  =
  \rpair{\zero}{v}
  +
  \rpair{\ww}{-v}
  \in
  \asubdifrext{\zbar}
  +
  \asubdifindepifext{\zbar}
  =
  \asubdifsext{\zbar}.
\end{equation}
The inclusion is by Eqs.~(\ref{eq:thm:subgrad-comp-inc-fcn:6})
and~(\ref{eq:thm:subgrad-comp-inc-fcn:7}).
The last equality is by
\Cref{thm:subgrad-sum-fcns}
(with $f_1=r$ and $f_2=\indepif$),
noting that
$\ri(\dom{r})\cap\ri(\epi{f})\neq\emptyset$
by
\Cref{lem:thm:conj-compose-our-version:1}(\ref{lem:thm:conj-compose-our-version:1:b}).

From \eqref{eq:thm:subgrad-comp-inc-fcn:1},
$h=\PPx s$.
Since $\PPx\zbar=\xbar$ and by
\eqref{eq:thm:subgrad-comp-inc-fcn:3},
\Cref{thm:subgrad-Af}
(applied with $f$, $\A$, $\xbar$, $\uu$, $\zbar$,
as they appear in that \namecref{thm:subgrad-Af},
set to $s$, $\PPx$, $\xbar$, $\ww$, $\zbar$)
now yields
$\ww\in\asubdifhext{\xbar}$, as claimed.
\qedhere
\end{proof-parts}
\end{proof}

\Cref{thm:uv-cond-suff-comp-subgrad} gives a sufficient condition for
a vector to be an astral subgradient analogous to the chain rule from
calculus.
Nonetheless, as the next example shows,
there can be other subgradients besides these.
\indexg{subdifferentials, standard!composition with nondecreasing function@of composition with nondecreasing function|(}%
\indexg{nondecreasing function, composition with!standard subgradients of|(}%
This is true for both astral and standard subgradients.
(In the latter case, this does not contradict
\Cref{pr:std-subgrad-comp-inc}
which assumes $f$ is finite everywhere.)

\begin{example}   \label{ex:uv-not-suf-for-subgrad-g-of-f}
Define $f:\R\rightarrow\Rext$ and $g:\R\rightarrow\Rext$,
for $x\in\R$, by
\[
  f(x)
  =
  \begin{cases}
         x \ln x   & \mbox{if $x\in [0,1]$,} \\
         +\infty   & \mbox{otherwise,}
  \end{cases}
\]
and by $g(x)=\max\{0,x\}$.
Then $f$ and $g$ are closed, proper and convex, and $g$ is also
nondecreasing and not constant.
By \Cref{pr:conv-inc:prop}(\ref{pr:conv-inc:infsup}),
$g$'s extension is
$\gext(\barx)=\max\{0,\barx\}$ for $\barx\in\Rext$.
Thus, $h=\eg\circk f$ is the indicator function $\indf{[0,1]}$.
It can be checked that
$\partial h(0) = (-\infty,0]$.
Nonetheless,
$\partial f(0) = \emptyset$
(by a similar argument as in
Example~\ref{ex:entropy-1d}).

Therefore, letting $x=0$ and $w=-1$, this shows that
$w\in\partial h(x)$, and yet there does not exist
$v\in\partial g(f(x))$ and $u\in\partial f(x)$
with $v u = w$.
By
\Cref{pr:asubdiffext-at-x-in-rn}(\ref{pr:asubdiffext-at-x-in-rn:c}),
the same statement holds for the corresponding astral subgradients of
these functions' extensions; that is,
$w\in\asubdifhext{x}$,
but there does not exist
$v\in\asubdifgext{\fext(x)}$ and $u\in\asubdiffext{x}$
with $v u = w$.
\end{example}

Thus, in general, not all (astral or standard) subgradients of the
composition $h$ can be computed using the most direct generalization
of the chain rule.
In this example, the exception occurred for a subgradient $w=-1$ at a
point $x=0$ at the boundary of $h$'s effective domain.
At such points, subgradients can ``wrap around''
the boundary of the effective domain;
for instance, in this case,
$h$'s subdifferential at $0$ is $\partial h(0)=(-\infty,0]$.

In terms of $f$ and $g$, the following facts seem
plausibly relevant to such behavior at the edges:
  First, $g$ was minimized by $f(x)$, that is,
  $0\in\partial g(f(x))$.
  And second, $x=0$ was at the ``far end''
  of $\dom{f}=[0,1]$, the effective domain of $f$,
  in the direction of $w$.
  In other words,
  \[ x w = \sup_{z\in\dom{f}} z w = \inddomfstar(w). \]
  For points in $\dom f$, this condition is equivalent
  to $w$ belonging to $\partial \inddomf(x)$,
  or alternatively, to $w$ being in the normal cone to $\dom f$ at
  $x$
  (see \Cref{sec:subdif-ind-fcn-norm-cone}).%
\indexg{subdifferentials, standard!composition with nondecreasing function@of composition with nondecreasing function|)}%
\indexg{nondecreasing function, composition with!standard subgradients of|)}%

In fact, in astral space, the analogous forms of
these conditions are generally sufficient
for $\hext$ to have a particular astral subgradient.
Specifically, as shown next,
if $\ww\in\asubdifinddomfext{\xbar}$ and
if
$\zero\in\asubdifgext{\fext(\xbar)}$, then
$\ww\in\asubdifhext{\xbar}$.

Geometrically, as shown in
\Cref{cor:subdif-ind-fcn-rn}(\ref{cor:subdif-ind-fcn-rn:a},\ref{cor:subdif-ind-fcn-rn:c}),
for points
$\xbar\in\cldom{f}$, the condition 
$\ww\in\asubdifinddomfext{\xbar}$ is equivalent
to $\ww$ being in the normal cone to $\cldom{f}$ at $\xbar$,
with $\xbar\cdot\ww\in\R$.

\begin{theorem}  \label{thm:idomf-cond-suff-comp-subgrad}
  Let $f:\Rn\rightarrow\Rext$ be convex, and let
  $g:\R\rightarrow\Rext$ be
  convex, proper, closed, nondecreasing, and not constant.
  Assume there exists $\xing\in\Rn$ such that
  $f(\xing)<\sup(\dom g)$.
  Let $h=\eg\circk f$,
  and let $\xbar\in\extspace$.

  Suppose $\ww\in\basubdifinddomfext{\xbar}$ and
  $0\in\basubdifgext{\fext(\xbar)}$.
  Then $\ww\in\basubdifhext{\xbar}$.
\end{theorem}

\begin{proof}
We have
\begin{equation}
\label{eq:thm:subgrad-comp-inc-fcn:4}
  \hext(\xbar)
  =
  \gext(\fext(\xbar))
  =
  -\gstar(0)
  \in
  \R.
\end{equation}
The first equality is by
\Cref{thm:Gf-conv}, due to $\xing$'s existence, and
the second is by
\Cref{thm:fminus-subgrad-char}(\ref{thm:fminus-subgrad-char:a},\ref{thm:fminus-subgrad-char:d}) since
$0\in\asubdifgext{\fext(\xbar)}$,
noting that $\gstar(0)\in\R$ by
\Cref{pr:subgrad-imp-in-cldom}(\ref{pr:subgrad-imp-in-cldom:c}).
Furthermore, since
$\ww\in\asubdifinddomfext{\xbar}$,
\begin{equation}  \label{eq:thm:subgrad-comp-inc-fcn:5}
  \xbar\cdot\ww=\inddomfstar(\ww)\in\R,
\end{equation}
by
\Cref{cor:subdif-ind-fcn-rn}(\ref{cor:subdif-ind-fcn-rn:a},\ref{cor:subdif-ind-fcn-rn:b}).

Thus,
\begin{equation}  \label{eq:lem:idomf-cond-suff-comp-subgrad:1}
  {-\hstar(\ww)}
  \leq
  \hminuswext(\xbar)
  =
  \hext(\xbar) - \xbar\cdot\ww
  =
  -\gstar(0) - \inddomfstar(\ww)
  \leq
  -\hstar(\ww).
\end{equation}
The first inequality is by
\Cref{pr:fminusu-props}(\ref{pr:fminusu-props:d}).
The first equality is by
\Cref{pr:fminusu-props}(\ref{pr:fminusu-props:e})
since the arguments above imply that
$\hext(\xbar)$ and $-\xbar\cdot\ww$ are summable.
The second equality is from
Eqs.~(\ref{eq:thm:subgrad-comp-inc-fcn:4})
and~(\ref{eq:thm:subgrad-comp-inc-fcn:5}).
And the last inequality is from the formula for the conjugate of a
composition given in
\Cref{thm:conj-compose-our-version},
combined with
\Cref{pr:support-epi-f-conjugate}
(with $\uu=\ww$ and $v=0$).

\eqref{eq:lem:idomf-cond-suff-comp-subgrad:1} implies
that $\hminuswext(\xbar) = -\hstar(\ww)$,
and also (when combined with
Eqs.~\ref{eq:thm:subgrad-comp-inc-fcn:4}
and~\ref{eq:thm:subgrad-comp-inc-fcn:5})
that $\hstar(\ww)\in\R$.
Therefore, $\ww\in\asubdifhext{\xbar}$ by
\Cref{thm:fminus-subgrad-char}(\ref{thm:fminus-subgrad-char:d},\ref{thm:fminus-subgrad-char:a}).
\end{proof}

It turns out that
the two sufficient conditions given in
\Cref{thm:uv-cond-suff-comp-subgrad,thm:idomf-cond-suff-comp-subgrad}
account for all astral subgradients of
the extension of $h=\eg\circk f$
(under the same conditions as before).
This is shown below in
\Cref{thm:subgrad-comp-inc-fcn-equiv}, with these two conditions
appearing in part~(\ref{thm:subgrad-comp-inc-fcn-equiv:c}).
This \namecref{thm:subgrad-comp-inc-fcn-equiv} also gives another
characterization of $\hext$'s subgradients,
as we now explain.

Compared to the familiar chain rule from calculus,
the sufficient condition given in
\Cref{thm:idomf-cond-suff-comp-subgrad}
may seem odd, as is the apparent need to
be using two seemingly unconnected conditions
rather than a single unified rule.
In fact, this condition
can be expressed in a form that better reveals its connection to the
chain rule, while also unifying the two conditions into a single rule.

As discussed earlier, the chain rule
states that
$\gradh(\xx)=g'(f(\xx))\posKern\gradf(\xx)$,
where $h=g\circk f$, and
$g'$ is the
derivative of $g$.
Letting $v=g'(f(\xx))$,
we can restate this as
\begin{equation}
\label{eq:gradh:vf}
  \gradh(\xx)
  =
  v \gradf(\xx)
  =
  \nabla (v f)(\xx).
\end{equation}

In an analogous way,
in \Cref{thm:subgrad-comp-inc-fcn-equiv},
we can re-express the
conditions of
\Cref{thm:uv-cond-suff-comp-subgrad,thm:idomf-cond-suff-comp-subgrad}
using the operation $\sfprod{v}{f}$
introduced in \Cref{sec:prelim:cvx-fcns} (Eqs.~\ref{eq:sfprod-defn} and~\ref{eq:sfprod-identity}).
In particular, by \eqref{eq:sfprod-identity},
for any $f:\Rn\to\eR$ and any $v\in\Rpos$,
\begin{equation}
  \label{eq:sfprod-vf}
      \sfprod{v}{f}
      =
      \begin{cases}
        v f
          & \text{if $v>0$,} \\
        \inddomf
          & \text{if $v=0$.}
      \end{cases}
\end{equation}
Analyzing $h=\eg\circk f$,
we show that $\ww\in\asubdifhext{\xbar}$
if and only if $\ww\in\asubdifsfprodext{v}{f}{\xbar}$ for some
$v\in\asubdifgext{\fext(\xbar)}$,
thus generalizing \eqref{eq:gradh:vf}.
The case $v>0$, with
$\asubdifsfprodext{v}{f}{\xbar}=\asubdifvfext{\xbar}=v\asubdiffext{\xbar}$,
corresponds to \Cref{thm:uv-cond-suff-comp-subgrad},
expressed as
\Cref{thm:subgrad-comp-inc-fcn-equiv}(\ref{thm:subgrad-comp-inc-fcn-equiv:c})(\ref{thm:subgrad-comp-inc-fcn-equiv:c:1}).
The case $v=0$, with
$\asubdifsfprodext{v}{f}{\xbar}=\basubdifinddomfext{\xbar}$,
corresponds to \Cref{thm:idomf-cond-suff-comp-subgrad},
expressed as
\Cref{thm:subgrad-comp-inc-fcn-equiv}(\ref{thm:subgrad-comp-inc-fcn-equiv:c})(\ref{thm:subgrad-comp-inc-fcn-equiv:c:2}).
Note that we only need to consider $v\ge 0$, because
any $v\in\basubdifgext{\fext(\xbar)}$ must satisfy $v\ge 0$,
by \Cref{pr:subgrad-inc-not-neg}, since
$\eg$ is nondecreasing.
Thus, the following equivalences hold:

\begin{theorem}  \label{thm:subgrad-comp-inc-fcn-equiv}
  Let $f:\Rn\rightarrow\Rext$ be convex, and let
  $g:\R\rightarrow\Rext$ be
  convex, proper, closed, nondecreasing, and not constant.
  Assume there exists $\xing\in\Rn$ such that
  $f(\xing)<\sup(\dom g)$.
  Let $h=\eg\circk f$,
  and let $\xbar\in\extspace$ and $\ww\in\Rn$.
  Then the following are equivalent:
  \begin{letter-compact}
  \item  \label{thm:subgrad-comp-inc-fcn-equiv:a}
    $\ww\in\asubdifhext{\xbar}$.    
  \item  \label{thm:subgrad-comp-inc-fcn-equiv:c}
    Either of the following hold:
    \begin{roman-compact}
    \item  \label{thm:subgrad-comp-inc-fcn-equiv:c:1}
      There exists $\uu\in\basubdiffext{\xbar}$
      and $v\in\basubdifgext{\fext(\xbar)}$
      such that $\ww=v \uu$; or
    \item  \label{thm:subgrad-comp-inc-fcn-equiv:c:2}
      $\ww\in\basubdifinddomfext{\xbar}$ and
      $0\in\basubdifgext{\fext(\xbar)}$.
    \end{roman-compact}
  \item  \label{thm:subgrad-comp-inc-fcn-equiv:b}
    There exists $v\in\asubdifgext{\fext(\xbar)}$
    such that
    $\ww\in\asubdifsfprodext{v}{f}{\xbar}$.
  \end{letter-compact}
  Moreover,
  the same equivalence holds if, in
  statement~(\ref{thm:subgrad-comp-inc-fcn-equiv:c})(\ref{thm:subgrad-comp-inc-fcn-equiv:c:1}),
  we further require that $v>0$,
  and also if, in
  statement~(\ref{thm:subgrad-comp-inc-fcn-equiv:b}),
  we further require that $v\geq 0$.
\end{theorem}

In the proof, we will use the following simple proposition:

\begin{proposition}
\label{pr:a+b}
Let $\alpha,\beta\in\R$,
and let $\seq{\alpha_t}$ and $\seq{\beta_t}$
be sequences in $\eR$ such that $\alpha_t\le\alpha$ and $\beta_t\le\beta$ for
all $t$, and $\alpha_t+\beta_t\to\alpha+\beta$. Then $\alpha_t\to\alpha$ and $\beta_t\to\beta$.
\end{proposition}

\begin{proof}
Since $\alpha-\alpha_t\ge 0$ and $\beta-\beta_t\ge 0$,
the sequences $\seq{\alpha-\alpha_t}$
and $\seq{\beta-\beta_t}$ are both nonnegative and majorized by
$\seq{(\alpha-\alpha_t)+(\beta-\beta_t)}$,
which, by assumption, converges to $0$.
This implies the claim.
\end{proof}

\begin{proof}[Proof of \Cref{thm:subgrad-comp-inc-fcn-equiv}]

~

\begin{proof-parts}
\pfpart{%
  (\ref{thm:subgrad-comp-inc-fcn-equiv:c})
  $\Rightarrow$
  (\ref{thm:subgrad-comp-inc-fcn-equiv:a}):
}
This is immediate from
\Cref{thm:uv-cond-suff-comp-subgrad,thm:idomf-cond-suff-comp-subgrad}.

\pfpart{%
  (\ref{thm:subgrad-comp-inc-fcn-equiv:b})
  $\Rightarrow$
  (\ref{thm:subgrad-comp-inc-fcn-equiv:c}):
}
Suppose $\ww\in\asubdifsfprodext{v}{f}{\xbar}$
for some $v\in\asubdifgext{\fext(\xbar)}$.
Then $v\geq 0$
by \Cref{pr:subgrad-inc-not-neg},
since $\gext$ is nondecreasing
(\Cref{pr:conv-inc:prop}\ref{pr:conv-inc:nondec}).

If $v>0$, then $\sfprod{v}{f}=v f$ 
(by Eq.~\ref{eq:sfprod-vf}),
so
$\ww\in\asubdifvfext{\xbar}=v\asubdiffext{\xbar}$
(with equality by \Cref{pr:subgrad-scal-mult}),
implying
$\ww/v \in \basubdiffext{\xbar}$.
Thus,
condition~\partsubref{thm:subgrad-comp-inc-fcn-equiv:c}{thm:subgrad-comp-inc-fcn-equiv:c:1}
holds (with $\uu=\ww/v$ and $v>0$).

Otherwise, if $v=0$, then
$\sfprod{0}{f}=\inddomf$
(by Eq.~\ref{eq:sfprod-vf}),
implying
condition~\partsubref{thm:subgrad-comp-inc-fcn-equiv:c}{thm:subgrad-comp-inc-fcn-equiv:c:2}.

\pfpart{%
  (\ref{thm:subgrad-comp-inc-fcn-equiv:a})
  $\Rightarrow$
  (\ref{thm:subgrad-comp-inc-fcn-equiv:b}):
}
Suppose $\ww\in\basubdifhext{\xbar}$.
Then by \Cref{thm:fminus-subgrad-char}(\ref{thm:fminus-subgrad-char:a},\ref{thm:fminus-subgrad-char:c}),
  there exists a sequence $\seq{\xx_t}$ in $\Rn$ with $\xx_t\to\xbar$, $h(\xx_t)\in\R$ for all $t$,
  and $\xx_t\inprod\ww-h(\xx_t)\to\hstar(\ww)$.
  Also, since $h(\xx_t)\in\R$, we must have $f(\xx_t)<+\infty$
  for all $t$.
  
  By \Cref{thm:conj-compose-our-version}
  (with $G=\eg$,
   as follows from
  Propositions~\ref{pr:h:1}\ref{pr:h:1a}
  and~\ref{pr:conv-inc:prop}\ref{pr:conv-inc:infsup}\ref{pr:conv-inc:nonconst}),
  there exists
  $v\in\R$ such that
  \begin{equation*}
    \hstar(\ww)
    =
    g^*(v) + \indepifstar(\rpairf{\ww}{-v}).
  \end{equation*}
  Further,
  by \Cref{pr:subgrad-imp-in-cldom}(\ref{pr:subgrad-imp-in-cldom:c}), $\hstar(\ww)\in\R$,
  so also $\indepifstar(\rpairf{\ww}{-v})\in\R$.
\Cref{pr:support-epi-f-conjugate},
combined with this last fact,
then implies
$v\geq 0$ and that
$\indepifstar(\rpairf{\ww}{-v})=\sfprodstar{v}{f}(\ww)$.
Thus,
\begin{equation} \label{eq:thm:subgrad-comp-inc-fcn-equiv:5}
      \hstar(\ww)
      =
      g^*(v) + \sfprodstar{v}{f}(\ww).
\end{equation}
Since $\hstar(\ww)\in\R$,
also $\gstar(v)\in\R$ and $\sfprodstar{v}{f}(\ww)\in\R$.
In particular, this means that $\sfprod{v}{f}>-\infty$ (\Cref{pr:conj-props}\ref{pr:conj-props:c2}).
Since, for all $t$, $f(\xx_t)<+\infty$, 
this implies 
\begin{equation} \label{eq:thm:subgrad-comp-inc-fcn-equiv:a1}
  v f(\xx_t) = (\sfprod{v}{f})(\xx_t)\in\R ,
\end{equation}
with the equality from $\sfprod{v}{f}$'s definition
(Eq.~\ref{eq:sfprod-defn}).

For each $t$, let
\begin{align*}
  \alpha_t  &=  \xx_t\cdot\ww - (\sfprod{v}{f})(\xx_t),
  &
  \alpha    &=  \sfprodstar{v}{f}(\ww),
  \\
  \beta_t   &=  v f(\xx_t) - \gext\bigParens{f(\xx_t)},
  &
  \beta     &=  \gstar(v).
\end{align*}
As noted earlier, $\alpha$ and $\beta$ are both finite,
as is $\alpha_t$.
Moreover,
$\alpha_t\leq\alpha$
by definition of standard conjugate
(Eq.~\ref{eq:fstar-def:intro}),
and
$\beta_t\leq\beta$
by definition of astral conjugate
(Eq.~\ref{eq:Fstar-down-def}),
and also using \Cref{pr:fextstar-is-fstar}
and that $v f(\xx_t)\in\R$.
Furthermore,
\begin{equation} \label{eq:thm:subgrad-comp-inc-fcn-equiv:n3}
       \alpha_t+\beta_t
       =
       \xx_t\inprod\ww
       -\gext\bigParens{f(\xx_t)}
       =\xx_t\inprod\ww -h(\xx_t)
       \to\hstar(\ww) = \alpha+\beta,
\end{equation}
where the first equality uses
\eqref{eq:thm:subgrad-comp-inc-fcn-equiv:a1},
the convergence is a defining property of the sequence $\seq{\xx_t}$,
and the last equality is by 
\eqref{eq:thm:subgrad-comp-inc-fcn-equiv:5}.
Therefore,
by \Cref{pr:a+b}, $\alpha_t\to\alpha$ and $\beta_t\to\beta$;
that is,
\begin{align}
  \xx_t\cdot\ww - (\sfprod{v}{f})(\xx_t)
  &\rightarrow
  \sfprodstar{v}{f}(\ww),
  \label{eq:thm:subgrad-comp-inc-fcn-equiv:8}
  \\
  v f(\xx_t)-\gext\bigParens{f(\xx_t)}
  &\rightarrow
  \gstar(v).
  \label{eq:thm:subgrad-comp-inc-fcn-equiv:7}
\end{align}
    By \Cref{thm:fminus-subgrad-char}(\ref{thm:fminus-subgrad-char:c},\ref{thm:fminus-subgrad-char:a}),
and since $\sfprodstar{v}{f}(\ww)\in\R$,
\eqref{eq:thm:subgrad-comp-inc-fcn-equiv:8} implies that
$\ww\in\asubdifsfprodext{v}{f}{\xbar}$.
It remains to show that also $v\in\basubdifgext{\fext(\xbar)}$,
which we prove in cases.

Suppose first that $v>0$.
Then
\begin{align*}
  v f(\xx_t)
  =
  (\sfprod{v}{f})(\xx_t)
  &=
  \bigBracks{(\sfprod{v}{f})(\xx_t)-\xx_t\inprod\ww}
  +\xx_t\inprod\ww
  \\
  &\rightarrow
  -\sfprodstar{v}{f}(\ww) + \xbar\inprod\ww
  =
  \sfprodext{v}{f}(\xbar)
  =
  \vfext(\xbar)
  =
  v \ef(\xbar).
\end{align*}
{\mathtogether%
The first and fourth equalities are by
\eqref{eq:sfprod-vf}.
The convergence follows from
\eqref{eq:thm:subgrad-comp-inc-fcn-equiv:8}
and because
$\xx_t\cdot\ww\rightarrow\xbar\cdot\ww$
(\Cref{thm:i:1}\ref{thm:i:1c}),
and also using that
$\sfprodstar{v}{f}(\ww)\in\R$
so that the sum of the
limits is the limit of the sums
(\Cref{prop:lim:eR}\ref{i:lim:eR:sum}).
The third equality is by
\Cref{thm:fenchel-subgrad}(\ref{thm:fenchel-subgrad:b},\ref{thm:fenchel-subgrad:c}),
since
$\ww\in\asubdifsfprodext{v}{f}{\xbar}$
(and by \Cref{pr:fextstar-is-fstar}).
The last equality is by
\Cref{pr:scal-mult-ext}.}

Thus, $f(\xx_t)\to\ef(\xbar)$.
Combined with
\eqref{eq:thm:subgrad-comp-inc-fcn-equiv:7}
and since
$\gext$ is lower semicontinuous (\Cref{prop:ext:F}\ref{prop:ext:F:a})
and $\gextstar(v)=\gstar(v)\in\R$,
it follows that
$v\in\basubdifgext{\fext(\xbar)}$
by
\Cref{thm:subgrad-eq-fenchel}(\ref{thm:subgrad-eq-fenchel:d},\ref{thm:subgrad-eq-fenchel:a})
(applied with $F$, $\xbar_t$, $\xbar$, $\uu$, as they appear in that
\namecref{thm:subgrad-eq-fenchel}, set to
$\gext$, $f(\xx_t)$, $\fext(\xbar)$, $v$).
This completes the proof in this case.

In the alternative case that $v=0$,
we have
    \[
      -\gstar(0)
      \leq
      \eg\bigParens{\ef(\xbar)}
      \leq
      \eh(\xbar)
      \le
      \lim h(\xx_t)
      =
      \lim \eg\bigParens{f(\xx_t)}
      =
      -\gstar(0).
    \]
The first inequality is by
\Cref{pr:fminusu-props}(\ref{pr:fminusu-props:d}).
The second inequality is by
\Cref{pr:Gf-cont}(\ref{pr:Gf-cont:a})
(and since $\gext$ is lower semicontinuous).
The third inequality is from $\hext$'s definition
(Eq.~\ref{eq:e:7}).
The final equality is by \eqref{eq:thm:subgrad-comp-inc-fcn-equiv:7}.
    Thus,
    $\eg\regParens{\ef(\xbar)}= -\gstar(0)$,
    so $0\in\asubdifgext{\ef(\xbar)}$
    by \Cref{thm:fminus-subgrad-char}(\ref{thm:fminus-subgrad-char:d},\ref{thm:fminus-subgrad-char:a}),
    completing the proof in this case as well.
\qedhere
\end{proof-parts}
\end{proof}

Thus, under the conditions of
\Cref{thm:subgrad-comp-inc-fcn-equiv}, the equivalence between
statements~(\ref{thm:subgrad-comp-inc-fcn-equiv:a})
and~(\ref{thm:subgrad-comp-inc-fcn-equiv:b})
implies that $\basubdifhext{\xbar}$ can be written succinctly as
\begin{equation*}
    \partial\eh(\xbar)
    \;=\;
    \bigcup\;\BigBraces{%
      \basubdifsfprodext{v}{f}{\xbar}
      :\:
      v\in \partial\eg\bigParens{\ef(\xbar)}
    }.
\end{equation*}

As discussed earlier and demonstrated in
\Cref{ex:uv-not-suf-for-subgrad-g-of-f},
\Cref{thm:subgrad-comp-inc-fcn-equiv}
cannot in general be simplified
by omitting the less intuitive
condition~(\ref{thm:subgrad-comp-inc-fcn-equiv:c:2})
in statement~(\ref{thm:subgrad-comp-inc-fcn-equiv:c}).
However, there are some special cases when this is possible,
as enumerated in the next
\namecref{cor:suff-cond-subgrad-comp-chain-rule}:

\begin{theorem}   \label{cor:suff-cond-subgrad-comp-chain-rule}
  Let $f:\Rn\rightarrow\Rext$ be convex, and let
  $g:\R\rightarrow\Rext$ be
  convex, proper, closed, nondecreasing, and not constant.
  Assume there exists $\xing\in\Rn$ such that
  $f(\xing)<\sup(\dom g)$.
  Let $h=\eg\circk f$,
  and let $\xbar\in\extspace$.
  Assume further that at least one of the following holds:
  \begin{roman-compact}
  \item   \label{cor:suff-cond-subgrad-comp-chain-rule:b}
    $\inf g = -\infty$.
  \item   \label{cor:suff-cond-subgrad-comp-chain-rule:c}
    $\fext(\xbar)$ does not minimize $\gext$; that is,
    $\gext(\fext(\xbar))>\inf g$.
  \item   \label{cor:suff-cond-subgrad-comp-chain-rule:a}
    $\dom{f}=\Rn$ and $\fext(\xbar)\in\R$.
  \end{roman-compact}
  Then
  \begin{equation}   \label{eq:cor:suff-cond-subgrad-comp-chain-rule:1}
     \basubdifhext{\xbar}
     =
     \Braces{ v \uu :\: \uu\in\basubdiffext{\xbar},\,
                      v\in\basubdifgext{\fext(\xbar)}
     }.
  \end{equation}
\end{theorem}

\begin{proof}
Let $U$ denote the set on the right-hand side of
\eqref{eq:cor:suff-cond-subgrad-comp-chain-rule:1}.

If $\uu\in\basubdiffext{\xbar}$
and $v\in\basubdifgext{\fext(\xbar)}$,
then $v\uu\in\basubdifhext{\xbar}$
by
\Cref{thm:uv-cond-suff-comp-subgrad}.
Thus,
$U\subseteq\basubdifhext{\xbar}$.

For the reverse inclusion, let $\ww\in\basubdifhext{\xbar}$.
By
\Cref{thm:subgrad-comp-inc-fcn-equiv}(\ref{thm:subgrad-comp-inc-fcn-equiv:a},\ref{thm:subgrad-comp-inc-fcn-equiv:c}),
either
condition~\partsubref{thm:subgrad-comp-inc-fcn-equiv:c}{thm:subgrad-comp-inc-fcn-equiv:c:1}
or
condition~\partsubref{thm:subgrad-comp-inc-fcn-equiv:c}{thm:subgrad-comp-inc-fcn-equiv:c:2}
of that
\namecref{thm:subgrad-comp-inc-fcn-equiv} must hold.

\begin{proof-parts}
\pfpart{%
  Conditions~(\ref{cor:suff-cond-subgrad-comp-chain-rule:b})
  or~(\ref{cor:suff-cond-subgrad-comp-chain-rule:c}):
}
Suppose either
$\inf g = -\infty$
or
$\gext(\fext(\xbar))>\inf g$.
Then
$0\not\in\asubdifgext{\fext(\xbar)}$
by \Cref{pr:asub-zero-is-min}(\ref{pr:asub-zero-is-min:a})
(and \Cref{pr:fext-min-exists}).
Therefore, 
condition~\partsubref{thm:subgrad-comp-inc-fcn-equiv:c}{thm:subgrad-comp-inc-fcn-equiv:c:2}
of
\Cref{thm:subgrad-comp-inc-fcn-equiv}
cannot hold, which means instead that that
\namecref{thm:subgrad-comp-inc-fcn-equiv}'s
condition~\partsubref{thm:subgrad-comp-inc-fcn-equiv:c}{thm:subgrad-comp-inc-fcn-equiv:c:1}
must instead hold, implying $\ww\in U$.

\pfpart{Condition~(\ref{cor:suff-cond-subgrad-comp-chain-rule:a}):}
Suppose
$\dom{f}=\Rn$ and $\fext(\xbar)\in\R$.
As a first case,
suppose $\ww\neq\zero$.
Then
\[
  \inddomfstar(\ww) = \sup_{\xx\in\dom{f}} \xx\cdot\ww = +\infty
\]
since $\dom{f}=\Rn$.
Thus,
by
\Cref{cor:subdif-ind-fcn-rn}(\ref{cor:subdif-ind-fcn-rn:a},\ref{cor:subdif-ind-fcn-rn:b}),
$\ww\not\in\basubdifinddomfext{\xbar}$ so,
as in the preceding case,
condition~\partsubref{thm:subgrad-comp-inc-fcn-equiv:c}{thm:subgrad-comp-inc-fcn-equiv:c:2}
of
\Cref{thm:subgrad-comp-inc-fcn-equiv}
cannot hold, implying instead that
condition~\partsubref{thm:subgrad-comp-inc-fcn-equiv:c}{thm:subgrad-comp-inc-fcn-equiv:c:1}
of that
\namecref{thm:subgrad-comp-inc-fcn-equiv}
holds so that $\ww\in U$.

In the alternative case, $\ww=\zero$.
Suppose
condition~\partsubref{thm:subgrad-comp-inc-fcn-equiv:c}{thm:subgrad-comp-inc-fcn-equiv:c:2}
of
\Cref{thm:subgrad-comp-inc-fcn-equiv}
holds
(since otherwise, if
condition~\partsubref{thm:subgrad-comp-inc-fcn-equiv:c}{thm:subgrad-comp-inc-fcn-equiv:c:1}
holds, then $\ww\in U$, as in the other cases).
Then
$0\in\basubdifgext{\fext(\xbar)}$.
Also, because $\fext(\xbar)\in\R$ and $\dom{f}=\Rn$,
there exists a point $\uu\in\basubdiffext{\xbar}$
by
\Cref{thm:subgrad-not-empty-fin-ri}.
Therefore, $\ww\in U$
since $\ww=\zero=0\cdot\uu$.%
\indexg{subdifferentials of extensions!composition with nondecreasing function@of composition with nondecreasing function|)}%
\indexg{nondecreasing function, composition with!astral subgradients of|)}%
\qedhere
\end{proof-parts}
\end{proof}

\part{Optimality and Optimization}
\label{part:optimality}

\chapter{Fenchel duality}
\label{sec:fenchel-duality}

In this and the next chapter,
we consider convex optimization problems
with special structure.
In the cases we study,
there are already classical conditions characterizing the
solution.
However, these characterizations can fail, for instance, if the
solution is infinite or if it is at a point where the function of
interest has no standard subgradients.
Here, we will see how these classical conditions can be
naturally generalized to astral space, yielding astral conditions that
apply more broadly than the standard ones, for instance,
when no finite solution exists.
These generalizations build heavily on the calculus developed in
\Cref{sec:calc-subgrads}.

\indexg{Fenchel duality|(}%
We begin in this chapter
with an astral generalization of Fenchel's duality theorem.
Here, the goal is to minimize the sum of two
convex functions,
or, in slightly more generality, the sum of a convex function with
another convex function that has been composed with a linear map.
Thus, the goal is to minimize
$f(\xx)+g(\A\xx)$ over $\xx\in\Rn$,
where $f:\Rn\rightarrow\Rext$ and
$g:\Rm\rightarrow\Rext$ are convex and proper,
and $\A\in\Rmn$.
For instance, $g$ might be an indicator function that encodes various
constraints that the solution must satisfy, so that the goal is to
minimize $f$ subject to those constraints.

We will see how the machinery of Fenchel duality generalizes to astral
space, yielding astral characterizations of the solutions to such
problems.
We then apply these to specific constrained optimization problems, as
well as to the derivation of an astral version of von Neumann's minimax
theorem.%

\section{Generalizing Fenchel duality to astral space}

\indexg{Fenchel duality!standard theorem|(}%
Here is a version of Fenchel's duality theorem
(similar, for instance, to Corollary~31.2.1 of \idxroc\citealp{ROC}):

\begin{theorem}   \label{thm:std-fenchel-duality}
  Let $f:\Rn\rightarrow\Rext$ and
  $g:\Rm\rightarrow\Rext$ be convex and proper,
  and let $\A\in\Rmn$.
  Assume there exists $\xing\in\ri(\dom f)$ such that
  $\A\xing\in\ri(\dom g)$.
  Then
  \begin{equation}   \label{eq:thm:std-fenchel-duality:1}
     \inf_{\xx\in\Rn} [f(\xx) + g(\A\xx)]
     =
     - \min_{\uu\in\Rm} [\fstar(-\transAk\uu) + \gstar(\uu)],
  \end{equation}
  meaning, in particular, that the minimum on the right is attained.
\end{theorem}

\begin{proof}
  First, note that $g\A$ is convex (by \Cref{roc:thm5.7:fA})
  and proper, since $g>-\infty$ and 
  $g(\A\xing)<+\infty$.
  Thus, $f+g\A$ is convex and proper (by \Cref{pr:std-sum-fcns-cvx}).
  Also, $\xing\in\ri(\dom(g\A))$ by
  \Cref{pr:Ax-ridomf-ridomfA},  so $\ri(\dom f)\cap\ri(\dom(g\A))\ne\emptyset$.
  
  Given these conditions, we have the following equalities:
  \begin{align*}
    &\inf_{\xx\in\Rn} (f + g\A)(\xx)
    =
    -(f+g\A)^*(\zero)
  \\&\qquad{}
    =
    -\min
    \BigBraces{\fstar(-\ww) + (g\A)^*(\ww):\:\ww\in\Rn}
  \\&\qquad{}
    =
    -\min
    \BigBraces{\fstar(-\ww) + (g\A)^*(\ww):\:\ww\in\colspace\transA}
  \\&\qquad{}
    =
    -\min
    \BigBraces{\fstar(-\ww) + g^*(\uu):\:
               \ww\in\colspace\transA,\,
               \uu\in\Rm,\,\transAk\uu=\ww}
  \\&\qquad{}
    =
    -\min
    \BigBraces{\fstar(-\transAk\uu) + g^*(\uu):\:
               \uu\in\Rm}.
\end{align*}
The first equality is by \Cref{pr:conj-props}(\ref{pr:conj-props:a}).
The second is by \Cref{roc:thm16.4}, which also establishes that the
minimum in $\ww$ is attained.
The third equality is because by \Cref{roc:thm16.3:fA},
$(g\A)^*(\ww)=+\infty$ if $\ww\not\in\colspace \transA$,
so we can assume without loss of generality that the minimum is attained at some $\ww\in\colspace\transA$.
The fourth equality is by \Cref{roc:thm16.3:fA}, which also establishes
attainment of the minimum in $\uu$. Finally, the fifth equality follows by noting that $\transAk\uu\in\colspace\transA$ for all $\uu\in\Rm$.
\end{proof}

Working in standard Euclidean space,
\Cref{thm:std-fenchel-duality}
only guarantees attainment for the minimum on the right-hand side
of \eqref{eq:thm:std-fenchel-duality:1}, but not necessarily
for the infimum on the
\indexg{Fenchel duality!standard theorem|)}%
left.
\indexg{Fenchel duality!astral version|(}%
Nevertheless, when extended in a natural way to astral space, this
infimum is also always attained, as shown next:

\begin{theorem}   \label{thm:ast-fenchel-duality}
  Let $f:\Rn\rightarrow\Rext$ and
  $g:\Rm\rightarrow\Rext$ be convex and proper,
  and let $\A\in\Rmn$.
  Assume there exists $\xing\in\ri(\dom f)$ such that
  $\A\xing\in\ri(\dom g)$.
  Then for all $\xbar\in\extspace$ and $\uu\in\Rm$, the following are
  equivalent:
  \begin{letter}
  \item   \label{thm:ast-fenchel-duality:a}
    $\xbar$ minimizes $\fgAext$,
    and $\uu$ minimizes $\fAgstar$.
  \item   \label{thm:ast-fenchel-duality:b}
    $\fgAext(\xbar) = -\bigBracks{\fstar(-\transAk\uu) + \gstar(\uu)}$.
  \end{letter}
  Furthermore, there exists such a pair $\xbar,\uu$ satisfying both of
  these conditions.
  Thus,
  \[
     \min_{\xbar\in\extspace} \fgAext(\xbar)
     =
     - \min_{\uu\in\Rm} \bigBracks{\fstar(-\transAk\uu) + \gstar(\uu)},
  \]
  with both minima attained.
\end{theorem}

\begin{proof}
Let $p=f+g\A$.

\begin{proof-parts}
\pfpart{%
  (\ref{thm:ast-fenchel-duality:a})
  $\Rightarrow$
  (\ref{thm:ast-fenchel-duality:b}):
}
Suppose
statement~(\ref{thm:ast-fenchel-duality:a}) holds for some
$\xbar\in\extspace$ and $\uu\in\Rn$.
Then
\[
  \pext(\xbar)
  =
  \inf p
  =
  -\inf\bigParens{ \fAgstar }
  =
  -\bigBracks{\fstar(-\transAk\uu) + \gstar(\uu)}.
\]
The first equality is by assumption and
\Cref{pr:fext-min-exists}.
The second equality is by
\Cref{thm:std-fenchel-duality}.
And the last equality is also by assumption.

\pfpart{%
  (\ref{thm:ast-fenchel-duality:b})
  $\Rightarrow$
  (\ref{thm:ast-fenchel-duality:a}):
}
Suppose now that
statement~(\ref{thm:ast-fenchel-duality:b}) holds for some
$\xbar\in\extspace$ and $\uu\in\Rn$.
Then
\[
   \inf p
   \leq
   \pext(\xbar)
   =
   -\bigBracks{\fstar(-\transAk\uu) + \gstar(\uu)}
   \leq
   -\inf\bigParens{ \fAgstar }
   =
   \inf p.
\]
The first inequality is by
\Cref{pr:fext-min-exists}.
The first equality is by assumption.
And the last equality is by
\Cref{thm:std-fenchel-duality}.
Thus, $\pext(\xbar)=\inf p$ so $\xbar$ minimizes $\pext$,
and similarly, $\uu$ minimizes
$\fAgstar$.

\pfpart{Existence:}
By
\Cref{thm:std-fenchel-duality},
there exists $\uu\in\Rn$ minimizing
$\fAgstar$.
And by
\Cref{pr:fext-min-exists},
there exists $\xbar\in\extspace$ minimizing $\pext$.
Together, these thus satisfy
statement~(\ref{thm:ast-fenchel-duality:a})
(and so
statement~(\ref{thm:ast-fenchel-duality:b})
as well).%
\indexg{Fenchel duality!astral version|)}%
\qedhere
\end{proof-parts}
\end{proof}

\indexg{Fenchel duality!minima characterized by subdifferentials|(}%
The next theorem gives necessary and sufficient optimality conditions
based on
astral primal and dual subgradients for when
a pair $\xbar,\uu$ is a solution in the sense of
satisfying the conditions of
\Cref{thm:ast-fenchel-duality},
assuming the objective function,
$f+g\A$, is lower-bounded:

\begin{theorem}   \label{thm:subgrad-equiv-fenchel}
  Let $f:\Rn\rightarrow\Rext$ and
  $g:\Rm\rightarrow\Rext$ be convex and proper,
  and let $\A\in\Rmn$.
  Assume there exists $\xing\in\ri(\dom f)$ such that
  $\A\xing\in\ri(\dom g)$,
  and assume also that
  $\inf (f+g\A) > -\infty$.
  Then for all $\xbar\in\extspace$ and $\uu\in\Rm$, the following are
  equivalent:
  \begin{letter-compact}
  \item   \label{thm:subgrad-equiv-fenchel:a}
    $\xbar$ minimizes $\fgAext$,
    and $\uu$ minimizes $\fAgstar$.
  \item   \label{thm:subgrad-equiv-fenchel:b}
    $\fgAext(\xbar) = -\bigBracks{\fstar(-\transAk\uu) + \gstar(\uu)}$.
  \item   \label{thm:subgrad-equiv-fenchel:c}
    \fencheldiffpair
        {-\transAk\uu\in\basubdiffext{\xbar}}
        {\uu\in\basubdifgext{\A\xbar}}
  \item     \label{thm:subgrad-equiv-fenchel:c-mixed}
    \fencheldiffpair
        {-\transAk\uu\in\basubdiffext{\xbar}}
        {\A\xbar\in\adsubdifgstar{\uu}}
  \item      \label{thm:subgrad-equiv-fenchel:c-dual}
    \fencheldiffpair
        {\xbar\in\adsubdiffstar{-\transAk\uu}}
        {\A\xbar\in\adsubdifgstar{\uu}}
  \end{letter-compact}
  Furthermore, there exists such a pair $\xbar,\uu$ satisfying all
  these conditions.
\end{theorem}

\begin{proof}
Let $p=f+g\A$.
Note that $\fstar$ and $\gstar$ are proper
(by
\Cref{pr:conj-props-cvx}\ref{pr:conj-props-cvx:a},
since $f$ and $g$ are proper).

\begin{proof-parts}
\pfpart{%
  (\ref{thm:subgrad-equiv-fenchel:a})
  $\Rightarrow$
  (\ref{thm:subgrad-equiv-fenchel:b}):
}
This is immediate from \Cref{thm:ast-fenchel-duality}.

\pfpart{%
  (\ref{thm:subgrad-equiv-fenchel:b})
  $\Rightarrow$
  (\ref{thm:subgrad-equiv-fenchel:c}):
}
Suppose
statement~(\ref{thm:subgrad-equiv-fenchel:b})
holds for some $\xbar\in\extspace$ and $\uu\in\Rm$.
Since $\fstar$ and $\gstar$ are proper,
this implies
$\pext(\xbar)<+\infty$.
Also, $\pext(\xbar)\geq \inf p>-\infty$ by assumption
(and \Cref{pr:fext-min-exists}),
so $\pext(\xbar)\in\R$.
Consequently, $\fstar(-\transAk\uu)$ and $\gstar(\uu)$ are both
finite.
Let $\alpha=\fstar(-\transAk\uu)$
and $\beta=\gstar(\uu)$.

Next,
let $\seq{\xx_t}$ be a sequence in $\Rn$ converging to $\xbar$ and
with $p(\xx_t)\rightarrow\pext(\xbar)$
(which exists by \Cref{pr:d1}).
For each $t$, let
\begin{align*}
  \alpha_t
  &=
  \xx_t\cdot(-\transAk\uu) - f(\xx_t)
  =
  -(\A\xx_t)\cdot\uu - f(\xx_t),
  \\
  \beta_t &= (\A\xx_t)\cdot\uu - g(\A\xx_t).
\end{align*}
Note that $-p(\xx_t)=\alpha_t+\beta_t$,
and that
statement~(\ref{thm:subgrad-equiv-fenchel:b})
means that $-\pext(\xbar)=\alpha+\beta$.
Thus,
$\alpha_t+\beta_t\rightarrow\alpha+\beta$.
Also, for all $t$,
$\alpha_t\leq\alpha$
and
$\beta_t\leq\beta$
by definition of conjugate
(Eq.~\ref{eq:fstar-def}).
Therefore, by \Cref{pr:a+b}, $\alpha_t\to\alpha$ and $\beta_t\to\beta$.
By
\Cref{thm:fminus-subgrad-char}(\ref{thm:fminus-subgrad-char:c},\ref{thm:fminus-subgrad-char:a}),
these in turn imply, respectively,
that $-\transAk\uu\in\asubdiffext{\xbar}$
and $\uu\in\basubdifgext{\A\xbar}$, where we also used that
$\fstar(-\transAk\uu)\in\R$, $\gstar(\uu)\in\R$, and that $\A\xx_t\to\A\xbar$
(by \Cref{thm:linear:cont}\ref{thm:linear:cont:b}).

\pfpart{%
  (\ref{thm:subgrad-equiv-fenchel:c})
  $\Rightarrow$
  (\ref{thm:subgrad-equiv-fenchel:c-mixed}),
  and
  (\ref{thm:subgrad-equiv-fenchel:c-mixed})
  $\Rightarrow$
  (\ref{thm:subgrad-equiv-fenchel:c-dual}):
}
These both follow directly from
\Cref{cor:strict-adif-fext-inverses}.

\pfpart{%
  (\ref{thm:subgrad-equiv-fenchel:c-dual})
  $\Rightarrow$
  (\ref{thm:subgrad-equiv-fenchel:a}):
}
Suppose statement~(\ref{thm:subgrad-equiv-fenchel:c-dual}) holds.
We first argue that $\uu$ minimizes 
$\fAgstar$.
By \Cref{thm:std-fenchel-duality},
there exists some point $\uu'\in\Rm$ minimizing
$\fAgstar$ and for which we have
\[
  -\bigBracks{\fstar(-\transAk\uu') + \gstar(\uu')}
  =
  \inf p
  >
  -\infty.
\]
Since $\fstar$ and $\gstar$ are proper, this implies that
$\fstar(-\transAk\uu')$ and $\gstar(\uu')$ are both
finite.

By the definition of astral dual subgradient given in
\eqref{eqn:psi-subgrad:3-alt},
statement~(\ref{thm:subgrad-equiv-fenchel:c-dual}) implies that
\begin{align}
  \label{eq:thm:subgrad-suff-fenchel:2}
  \gstar(\uu')
  &\geq
  \gstar(\uu) \plusd (\A\xbar)\cdot(\uu'-\uu),
  \;\;\text{and}
  \\
  \notag
  \fstar(-\transAk\uu')
  &\geq
  \fstar(-\transAk\uu)
  \plusd
  \xbar\cdot\regParens{-\transAk\uu'-\transAk\uu}
  \\
  \label{eq:thm:subgrad-suff-fenchel:3}
  &=
  \fstar(-\transAk\uu)
  \plusd
  (-\A\xbar)\cdot(\uu'-\uu),
\end{align}
where the last equality follows by
\Cref{thm:Ax-dot-u}.

We claim $(\A\xbar)\cdot(\uu'-\uu)\in\R$.
Otherwise, if
$(\A\xbar)\cdot(\uu'-\uu)=+\infty$, then the right-hand side of
\eqref{eq:thm:subgrad-suff-fenchel:2} must also be $+\infty$
(since $\gstar>-\infty$), contradicting that
$\gstar(\uu')\in\R$.
By similar reasoning, if
$(\A\xbar)\cdot(\uu'-\uu)=-\infty$, then the right-hand side of
\eqref{eq:thm:subgrad-suff-fenchel:3} must be $+\infty$,
contradicting that
$\fstar(-\transAk\uu')\in\R$.

Thus, $(\A\xbar)\cdot(\uu'-\uu)\in\R$.
Since $\gstar$ is proper and $\gstar(\uu')\in\R$,
\eqref{eq:thm:subgrad-suff-fenchel:2} therefore implies
that $\gstar(\uu)\in\R$.
By similar reasoning,
\eqref{eq:thm:subgrad-suff-fenchel:3}
now implies that
$\fstar(-\transAk\uu)\in\R$.

Thus, all of the terms appearing in
Eqs.~(\ref{eq:thm:subgrad-suff-fenchel:2})
and~(\ref{eq:thm:subgrad-suff-fenchel:3})
must in fact be finite, implying that these two inequalities can be
added (with the downward addition replaced by ordinary addition).
This yields that
\[
  \fstar(-\transAk\uu') + \gstar(\uu')
  \geq
  \fstar(-\transAk\uu) + \gstar(\uu).
\]
Since $\uu'$ minimizes $\fAgstar$, this shows that $\uu$ does as well.

It remains to prove that $\xbar$ minimizes $\pext$.
Applying
\Cref{cor:strict-adif-fext-inverses} to the astral dual subgradients
in statement~(\ref{thm:subgrad-equiv-fenchel:c-dual}),
and noting that
$\fstar(-\transAk\uu)$
and
$\gstar(\uu)$ are both in $\R$, as shown above,
we obtain that
$-\transAk\uu\in\asubdiffext{\xbar}$
and that
$\uu\in\asubdifgext{\A\xbar}$.
Consequently, we have
\begin{equation*}
  \zero
  =
  -\transAk\uu + \transAk\uu
  \in
  \asubdiffext{\xbar}
  +
  \transA \asubdifgext{\A\xbar}
  =
  \asubdiffext{\xbar}
  +
  \asubdifgAext{\xbar}
  =
  \asubdifpext{\xbar}.
\end{equation*}  
The inclusion is by the astral subgradient relations just noted.
The second equality
is by \Cref{thm:subgrad-fA},
noting that $\A\xing\in\ri(\dom g)$.
The third equality is by
\Cref{thm:subgrad-sum-fcns}
(with $f_1=f$ and $f_2=g\A$),
noting that 
$\xing\in\ri(\dom(g\A))$
by \Cref{pr:Ax-ridomf-ridomfA} since
$\A\xing\in\ri(\dom g)$, so that
$\xing\in\ri(\dom{f})\cap\ri(\dom(g\A))$.
Therefore, $\xbar$ minimizes $\pext$ by
\Cref{pr:asub-zero-is-min}(\ref{pr:asub-zero-is-min:a})
since $\zero\in\asubdifpext{\xbar}$,
completing the proof.

\pfpart{Existence:}
\Cref{thm:ast-fenchel-duality} proves that there exists a pair
$\xbar,\uu$ satisfying statement~(\ref{thm:subgrad-equiv-fenchel:a}),
and so all of the other statements as well.
\qedhere
\end{proof-parts}
\end{proof}

Without the assumption that $f+g\A$ is lower-bounded,
the equivalence given in
\Cref{thm:subgrad-equiv-fenchel} is no longer true in general.
More specifically,
as shown in the next example, if $\inf(f+g\A)=-\infty$, then it is
possible that
statements~(\ref{thm:subgrad-equiv-fenchel:a})
and~(\ref{thm:subgrad-equiv-fenchel:b})
of that theorem are true, but that the other
statements,~(\ref{thm:subgrad-equiv-fenchel:c}),~(\ref{thm:subgrad-equiv-fenchel:c-mixed})
and~(\ref{thm:subgrad-equiv-fenchel:c-dual}),
are false.

\begin{example}
Let $f(x)=-2x$ and $g(x)=x$ for $x\in\R$, and
let $\Ascalar=[1]$ (the $1\times 1$ identity matrix).
Both $f$ and $g$ are convex, proper, closed,
and finite everywhere, implying $0\in\ri(\dom{f})$
and $\Ascalar 0 = 0 \in\ri(\dom{g})$.
Their conjugates are the indicator functions
$\fstar=\indf{\{-2\}}$ and $\gstar=\indf{\{1\}}$.
Also,
let $\barx=+\infty$ and $u=0$, and
let $p=f+g\Ascalar$, implying $p(x)=-x$ for $x\in\R$.
Then $\pext(\barx)=\inf p=-\infty$, and
$\fstar(-\transAscalar\negKern u) = \fstar(0) = +\infty$
and
$\gstar(u)=\gstar(0)=+\infty$.
Thus, statement~(\ref{thm:subgrad-equiv-fenchel:b})
of \Cref{thm:subgrad-equiv-fenchel} holds,
and so statement~(\ref{thm:subgrad-equiv-fenchel:a}) does as well,
by \Cref{thm:ast-fenchel-duality}.

On the other hand,
$\asubdiffext{+\infty}=\{-2\}$
(by \Cref{ex:affine-subgrad-new}),
so $-\transAscalar\negKern u = 0 \not\in\asubdiffext{\barx}$.
Also,
$\Ascalar\barx=\barx\not\in\adsubdifgstar{u}$
since $\barx$ fails to satisfy the definition of astral dual
subgradient at $u$
(Eq.~\ref{eqn:psi-subgrad:3-alt});
specifically,
$\gstar(1)=0\not\geq+\infty=\gstar(u)\plusd\barx(1-u)$.
Thus, none of the
statements~(\ref{thm:subgrad-equiv-fenchel:c}),~(\ref{thm:subgrad-equiv-fenchel:c-mixed})
and~(\ref{thm:subgrad-equiv-fenchel:c-dual})
of \Cref{thm:subgrad-equiv-fenchel} hold for this example.
\end{example}

Nevertheless, even if the assumption that
$\inf(f+g\A)>-\infty$ is omitted from
\Cref{thm:subgrad-equiv-fenchel},
statement~(\ref{thm:subgrad-equiv-fenchel:c})
still implies each of the other statements.
This is because if 
statement~(\ref{thm:subgrad-equiv-fenchel:c})
holds, then $\fstar(-\transAk\uu)$ and $\gstar(\uu)$
are both in $\R$ by
\Cref{pr:subgrad-imp-in-cldom}(\ref{pr:subgrad-imp-in-cldom:c}),
implying by \Cref{thm:std-fenchel-duality} that
\[
  \inf (f+g\A)
  =
  -\inf (\fAgstar)
  \geq
  -\bigBracks{\fstar(-\transAk\uu) + \gstar(\uu)}
  >
  -\infty,
\]
and so further implying each of the other statements of
\Cref{thm:subgrad-equiv-fenchel}
as well.

With the additional assumption that $f$ and $g$ are closed,
an analogue of 
\Cref{thm:subgrad-equiv-fenchel}
for the standard (non-astral) setting can be obtained as a special
case. This is because, under these assumptions,
$f$ and $g$ are lower semicontinuous,
and hence the function $f+g\A$ is also lower semicontinuous
(by \Cref{pr:lsc-res-domain}\ref{pr:lsc-res-domain:b}\ref{pr:lsc-res-domain:c}). Thus,
if we only consider $\xx\in\Rn$ instead of $\xbar\in\eRn$,
then
$\fgAext(\xx)=(f+g\A)(\xx)$ (\Cref{pr:h:1}\ref{pr:h:1a}),
and also all the astral primal and dual
subdifferentials appearing in various conditions in \Cref{thm:subgrad-equiv-fenchel} can be replaced by standard subdifferentials
(by
\Cref{pr:asubdiffext-at-x-in-rn}\ref{pr:asubdiffext-at-x-in-rn:c}
and
\Cref{pr:adsubdif-int-rn}).
Thus, the equivalence in \Cref{thm:subgrad-equiv-fenchel}
holds in a standard setting for pairs $\xx\in\Rn$ and $\uu\in\Rn$,
as was previously known
\idxroc\citep[Theorem~31.3]{ROC}.
Although the equivalence holds in this standard setting, in general,
there need not exist a solution pair $\xx,\uu$ satisfying all (or any)
of the statements of that equivalence.
This is unlike the astral setting where such a pair $\xbar,\uu$ (with
$\xbar\in\extspace$) must always exist, as was proved in
\Cref{thm:subgrad-equiv-fenchel}.%
\indexg{Fenchel duality!minima characterized by subdifferentials|)}%

Here is an example:

\begin{figure}
  \centering
  \includegraphics{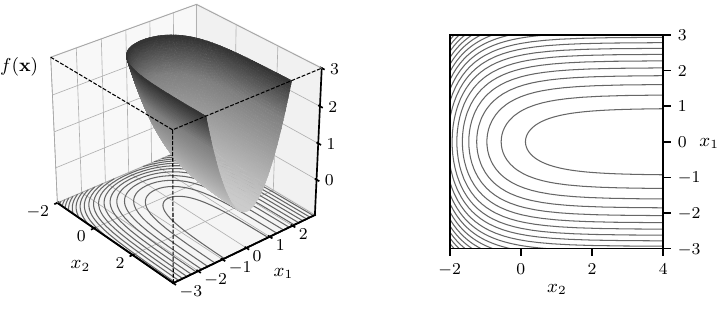}
  \mycaption{Gliding parabola}{%
\indexf{Gliding parabola}%
    The function $f$ from \Cref{ex:std-fenchel-equiv-not-exist}.
  }%
  \label{fig:para-exp}%
\end{figure}

\begin{example}[Gliding parabola]
\label{ex:std-fenchel-equiv-not-exist}
\indexg{Gliding parabola|(}%
\indexg{Gliding parabola!Fenchel duality and|(}%
Let $f:\R^2\rightarrow\Rext$ be defined, for $\xx\in\R^2$, by
\begin{equation}    \label{eq:ex:std-fenchel-equiv-not-exist:1}
  f(\xx)=f(x_1,x_2)=x_1^2+e^{-x_2}%
\indexg{Gliding parabola|)}%
\end{equation}
(see \Cref{fig:para-exp}),
and
let $g:\R\rightarrow\Rext$ be the indicator function
$g=\indf{[1,+\infty)}$.
Let $\ee_1,\ee_2$ be the standard basis vectors in $\R^2$, and let
$\A=[1,0]=\trans{\ee_1}$.
It can be checked that all the assumptions of 
\Cref{thm:subgrad-equiv-fenchel}
are satisfied.
Let $p=f+g\A$.
Then the problem of minimizing $p$ is equivalent to minimizing
$f(\xx)$ subject to $x_1\geq 1$.
Evidently, this minimum is attained by setting $x_1=1$ and taking the
limit as $x_2\rightarrow+\infty$; thus, $\inf p=1$.\looseness=-1

It can be calculated that, for $\uu\in\R^2$,
\begin{equation*}
  \fstar(\uu)
  =
  \fstar(u_1,u_2)
  =
  \begin{cases}
    \frac{1}{4} u_1^2 + (-u_2)\ln(-u_2) + u_2
              & \text{if $u_2 \leq 0$,} \\
    +\infty   & \text{otherwise,}
  \end{cases}
\end{equation*}
and for $u\in\R$,
\begin{equation*}
  \gstar(u)
  =
  \begin{cases}
    u         & \text{if $u \leq 0$,} \\
    +\infty   & \text{otherwise.}
  \end{cases}
\end{equation*}
Let $\psi=\fAgstar$.
Then, for $u\in\R$,
\[
  \psi(u)
  =
  \fstar(-u\ee_1) + \gstar(u)
  =
  \begin{cases}
    \frac{1}{4} u^2 + u
              & \text{if $u \leq 0$,} \\
    +\infty   & \text{otherwise.} \\
  \end{cases}
\]
Thus, by calculus, $\psi$ is minimized when $u=-2$,
so $\inf\psi=\psi(-2)=-1$, consistent with
\Cref{thm:std-fenchel-duality}.

We now argue that, for this example, no pair
$\xx\in\R^2$, $u\in\R$ can satisfy any of the statements
of a standard analogue of
\Cref{thm:subgrad-equiv-fenchel} (of the form discussed above).
Specifically, let $\xx$ be any point in $\R^2$ and $u$ any point in $\R$.
Then statement~(\ref{thm:subgrad-equiv-fenchel:a})
cannot hold because we always have $p(\xx+\ee_2)<p(\xx)$,
so $\xx$ cannot minimize $p$.
Statement~(\ref{thm:subgrad-equiv-fenchel:b}) cannot hold since
$p(\xx)>1=-\inf\psi\geq -\psi(u)$.
Statements~(\ref{thm:subgrad-equiv-fenchel:c})
and~(\ref{thm:subgrad-equiv-fenchel:c-mixed})
cannot hold since this would imply
$-u\ee_1=-\transA u\in\partial f(\xx)$, and so that
$-u\ee_1=\gradf(\xx)$ (by
\Cref{roc:thm25.1}\ref{roc:thm25.1:a}), which is impossible
since $\gradf(\xx)=\trans{[2x_1, -e^{-x_2}]}$.
Finally, 
statement~(\ref{thm:subgrad-equiv-fenchel:c-dual}) cannot hold since
$\partial \fstar(-\transA u)=\partial \fstar(-u\ee_1)=\emptyset$
by a similar argument as in
Example~\ref{ex:entropy-1d}.

On the other hand, in the astral setting, there does exist a solution
pair, namely $\xbar=\limray{\ee_2}\plusl\ee_1$ and $u=-2$,
which, for illustration, we verify satisfies each statement
of \Cref{thm:subgrad-equiv-fenchel}:
First, $\pext(\xbar)=\inf p=1=-\inf\psi=-\psi(u)$,
satisfying
statements~(\ref{thm:subgrad-equiv-fenchel:a})
and~(\ref{thm:subgrad-equiv-fenchel:b}).
Next, $\asubdiffext{\xbar}=\{2\ee_1\}$ by an argument similar to
\Cref{ex:subgrad-sqrt-approaches-abs-val:cont2},
and $\asubdifgext{\A\xbar}=\partial g(1)=(-\infty,0]$;
thus, 
statement~(\ref{thm:subgrad-equiv-fenchel:c})
is satisfied.
Finally, that
$\xbar\in\adsubdiffstar{-u\ee_1}=\adsubdiffstar{-\transA u}$
can be shown as in 
Example~\ref{ex:entropy-ast-dual-subgrad},
while
$\A\xbar = 1 \in \adsubdifgstar{u}$
follows straightforwardly from the definition of astral dual
subgradient (Eq.~\ref{eqn:psi-subgrad:3-alt}).
Thus,
statements~(\ref{thm:subgrad-equiv-fenchel:c-mixed})
and~(\ref{thm:subgrad-equiv-fenchel:c-dual})
are satisfied as well.%
\indexg{Gliding parabola!Fenchel duality and|)}%
\indexg{Fenchel duality|)}%
\end{example}

\section{Constrained optimization}
\label{sec:constrained-optimization}

We next look at how
\Cref{thm:subgrad-equiv-fenchel}
can be applied to characterize the solutions of constrained
optimization problems in which the goal is
to minimize a convex function $f:\Rn\rightarrow\Rext$ over some
convex set $C\subseteq\Rn$.
Said differently, the goal is to minimize the function
$p=f\plusu\indC$.

\indexg{constrained optimality conditions, standard|(}%
If $f$ is proper and differentiable everywhere, then a standard first-order
optimality condition states that a point $\xx\in C$ minimizes $f$ over
$C$ if and only
$-\gradf(\xx) \cdot (\yy-\xx) \leq 0$
for all $\yy\in C$, in other words, if and only if
$-\gradf(\xx)$ is in the normal cone $\ncone{\xx}{C}$.
More generally, even if $f$ is not everywhere differentiable,
but assuming $\ri(\dom f)$ and $\ri{C}$ overlap,
then
$\xx\in C$ minimizes $f$ over $C$ if and only if
$[-\partial f(\xx)] \cap \ncone{\xx}{C} \neq \emptyset$
(see, for example, \idxroc\citealp[Theorem 27.4]{ROC}).%
\indexg{constrained optimality conditions, standard|)}%

\indexg{constrained optimality conditions, astral|(}%
In this section, we show how more general optimality conditions can be
obtained using astral techniques which are applicable in cases where
standard subgradients are inadequate, for instance,
because a solution of the
optimization problem is at infinity, or because standard subgradients
do not exist at a solution.
Examples are given below.
We consider first the former setting in which the solution we seek may
be at infinity.

\indexg{constrained optimality conditions, astral!equivalent definitions|(}%
To begin, we need to define what it means for an astral point
$\xbar\in\extspace$ to solve the problem at hand of minimizing $f$
over $C$.
We consider three ways of defining what this means.
First, we might say $\xbar$ solves the optimization problem if
it minimizes $\pext$, the extension of $p=f\plusu\indC$, so that
$\pext(\xbar)=\min\pext=\inf p$.
Alternatively, considering $\fext$ and the set $C$ more directly,
we might say $\xbar$ solves the optimization problem if
$\xbar\in\Cbar$ and if $\xbar$
minimizes $\fext$ over $\Cbar$.
Finally, using sequences,
we might say that $\xbar$ solves the optimization problem if
there exists a sequence $\seq{\xx_t}$ in $C$ such that
$\xx_t\rightarrow\xbar$ and
$f(\xx_t)\rightarrow\inf p$.

In fact, under the conditions of the next proposition, these three
definitions are equivalent:

\begin{proposition}   \label{pr:cons-opt-equiv-soln}
  Let $f:\Rn\rightarrow\Rext$ be convex,
  and let $C\subseteq\Rn$ be convex.
  Assume $\ri(\dom f)\cap\ri{C}\neq\emptyset$.
  Let $p=f\plusu\indC$, and let $\xbar\in\extspace$.
  Then the following are equivalent:
  \begin{letter-compact}
  \item   \label{pr:cons-opt-equiv-soln:a}
    $\xbar$ minimizes $\pext$; that is,
    $\pext(\xbar)=\inf p$.
  \item   \label{pr:cons-opt-equiv-soln:b}
    $\xbar\in\Cbar$ and $\xbar$ minimizes $\fext$ over all points in
    $\Cbar$;
    that is, $\xbar$ minimizes $\fext\plusu\indfa{\Cbar}$.
  \item   \label{pr:cons-opt-equiv-soln:c}
    There exists a sequence $\seq{\xx_t}$ in $C$ such that
    $\xx_t\rightarrow\xbar$ and
    $f(\xx_t)\rightarrow\inf p$.
  \end{letter-compact}
\end{proposition}

\begin{proof}
We note that $\inf p<+\infty$ since, by assumption, there exists a
point in $\ri(\dom f)\cap\ri C\subseteq (\dom f)\cap C=\dom p$.

\begin{proof-parts}
\pfpart{%
  (\ref{pr:cons-opt-equiv-soln:a})
  $\Leftrightarrow$
  (\ref{pr:cons-opt-equiv-soln:b}):
}
By \Cref{cor:ext-restricted-fcn},
$\pext=\fext\plusu\indfa{\Cbar}$.
From this, and since $\inf p<+\infty$,
the stated equivalence is immediate.

\pfpart{%
  (\ref{pr:cons-opt-equiv-soln:a})
  $\Rightarrow$
  (\ref{pr:cons-opt-equiv-soln:c}):
}
Suppose $\pext(\xbar)=\inf p$.
By \Cref{pr:d1}, there exists a sequence
$\seq{\xx_t}$ in $\Rn$ that converges to $\xbar$ and with
$p(\xx_t)\rightarrow\pext(\xbar)$.
Since $\pext(\xbar)<+\infty$,
we can have $p(\xx_t)=+\infty$ for at most finitely many $t$.
By discarding these, we can assume for all $t$ that
$p(\xx_t)<+\infty$, and therefore that $\xx_t\in C$ so that
$p(\xx_t)=f(\xx_t)$.
Thus, $f(\xx_t)=p(\xx_t)\rightarrow\pext(\xbar)=\inf p$,
as claimed.

\pfpart{%
  (\ref{pr:cons-opt-equiv-soln:c})
  $\Rightarrow$
  (\ref{pr:cons-opt-equiv-soln:a}):
}
Suppose $\seq{\xx_t}$ is a sequence as in
statement~(\ref{pr:cons-opt-equiv-soln:c}).
Since $\xx_t\in C$, $p(\xx_t)=f(\xx_t)$
for all $t$.
Thus,
$p(\xx_t)=f(\xx_t)\rightarrow \inf p$,
and therefore,
\[
  \inf p
  =
  \lim p(\xx_t)
  \geq
  \pext(\xbar)
  \geq
  \inf p,
\]
where the inequalities are from the definition of
an extension (Eq.~\ref{eq:e:7})
and by \Cref{pr:fext-min-exists}.
This proves the claim.%
\indexg{constrained optimality conditions, astral!equivalent definitions|)}%
\qedhere
\end{proof-parts}
\end{proof}

\indexg{constrained optimality conditions, astral!general|(}%
Using \Cref{thm:subgrad-equiv-fenchel}, we can now generalize
the standard first-order optimality conditions mentioned earlier to
problems whose solution may be at infinity.
In particular,
as shown next under the stated conditions,
for the problem of minimizing a convex function
$f:\Rn\rightarrow\Rext$ over a convex set $C\subseteq\Rn$,
an astral point $\xbar\in\extspace$ is a solution in any and all
of the senses of \Cref{pr:cons-opt-equiv-soln} if and only if
$\xbar\in\Cbar$ and
$\fext$ has a negative astral subgradient at $\xbar$ that is also in the
normal cone to $\Cbar$ at $\xbar$.

\begin{theorem}  \label{thm:ast-const-opt-primal}
  Let $f:\Rn\rightarrow\Rext$ be convex and proper,
  and let $C\subseteq\Rn$ be convex.
  Let $p=f + \indC$, and let $\xbar\in\extspace$.
  Assume $\ri(\dom f)\cap\ri{C}\neq\emptyset$
  and that $\inf p > -\infty$.
  Then the following are equivalent:
  \begin{letter-compact}
  \item   \label{thm:ast-const-opt-primal:a}
    $\xbar$ minimizes $\pext$.
  \item   \label{thm:ast-const-opt-primal:b}
    $\xbar\in\Cbar$
    and
    $[-\asubdiffext{\xbar}] \cap \ncone{\xbar}{\Cbar} \neq \emptyset$.
  \end{letter-compact}
\end{theorem}

\begin{proof}
By assumption, there exists a point $\xing\in\ri(\dom f)\cap\ri{C}$.
  
We will shortly apply
\Cref{thm:subgrad-equiv-fenchel}.
To do so, we let $\A=\Iden$ (the $n\times n$ identity matrix)
and $g=\indf{C}$, the indicator function for $C$.
Then $g$ is convex since $C$ is, and is also proper
since $\xing\in\ri C\subseteq C$.
Further, $p=f+g\A$, and the assumptions of
\Cref{thm:subgrad-equiv-fenchel} are satisfied since
$\inf p>-\infty$,
$\xing\in\ri(\dom f)$, and $\A\xing=\xing\in\ri(\dom g)$.

\begin{proof-parts}
\pfpart{%
  (\ref{thm:ast-const-opt-primal:a})
  $\Rightarrow$
  (\ref{thm:ast-const-opt-primal:b}):
}
Suppose $\xbar$ minimizes $\pext$.
Then by
\Cref{thm:std-fenchel-duality},
there exists a point $\uu\in\Rn$ that minimizes
$\fstar[-\A]+\gstar$.
Therefore, by
\Cref{thm:subgrad-equiv-fenchel}(\ref{thm:subgrad-equiv-fenchel:a},\ref{thm:subgrad-equiv-fenchel:c}),
$-\uu=-\transAk\uu\in\asubdiffext{\xbar}$
and $\uu\in\asubdifgext{\A\xbar}=\asubdifgext{\xbar}$.
Since $g=\indf{C}$,
\Cref{cor:subdif-ind-fcn-rn}(\ref{cor:subdif-ind-fcn-rn:a},\ref{cor:subdif-ind-fcn-rn:c})
then implies that
$\xbar\in\Cbar$ and $\uu\in\ncone{\xbar}{\Cbar}$, completing the
proof.

\pfpart{%
  (\ref{thm:ast-const-opt-primal:b})
  $\Rightarrow$
  (\ref{thm:ast-const-opt-primal:a}):
}
Suppose $\xbar\in\Cbar$ and that there exists a point
$\uu\in [-\asubdiffext{\xbar}] \cap \ncone{\xbar}{\Cbar}$.
Then
$ \xbar \in \adsubdif{\indfastar{\Cbar}}{\uu} $
by
\Cref{pr:ncone-ast-sup-fcn}(\ref{pr:ncone-ast-sup-fcn:a},\ref{pr:ncone-ast-sup-fcn:c}).
Also, $\indfastar{\Cbar}=\indgenextstar{C}=\indfstar{C}=\gstar$,
with the first two equalities following respectively from
Propositions~\ref{pr:inds-ext} and~\ref{pr:fextstar-is-fstar}.
Thus, $\A\xbar=\xbar\in\adsubdifgstar{\uu}$ and $-\transAk\uu=-\uu\in\asubdiffext{\xbar}$,
so $\xbar$ minimizes $\pext$ by
\Cref{thm:subgrad-equiv-fenchel}(\ref{thm:subgrad-equiv-fenchel:c-mixed},\ref{thm:subgrad-equiv-fenchel:a}).
\qedhere
\end{proof-parts}
\end{proof}

\begin{example}    \label{ex:std-fenchel-equiv-not-exist:cont1}
\indexg{Gliding parabola!constrained optimality and|(}%
Continuing \Cref{ex:std-fenchel-equiv-not-exist},
we consider $f(\xx)=x_1^2+e^{-x_2}$ for $\xx\in\R^2$,
and we aim to minimize $f$ over the (standard) halfspace
$C=\{\xx\in\R^2 :\: \xx\cdot\ee_1\geq 1\}$.
Let $p=f + \indC$.
As noted earlier, this problem has no finite solution, but does have an
astral solution, namely, $\xbar=\limray{\ee_2}\plusl\ee_1$.
We can check this using
\Cref{thm:ast-const-opt-primal}:
The conditions of that theorem are satisfied since $\dom f=\R^2$ and
since $p\geq f\geq 0$.
Further, $\Cbar=\{\ybar\in\extspac{2} :\: \ybar\cdot\ee_1\geq 1\}$
by \Cref{cor:halfspace-closure}.
Let $\uu=-2 \ee_1$.
Then $\uu\in-\asubdiffext{\xbar}$, as shown in
\Cref{ex:std-fenchel-equiv-not-exist}.
Also, $\uu\in\ncone{\xbar}{\Cbar}$ by definition of normal cone
(Eq.~\ref{eq:dfn:ncone-ext})
since, for all $\ybar\in\Cbar$,
$\ybar\cdot\uu=-2\ybar\cdot\ee_1\leq -2=\xbar\cdot\uu$.
Thus, 
$\uu\in[-\asubdiffext{\xbar}]\cap\ncone{\xbar}{\Cbar}$,
so $\xbar$ minimizes $\pext$ by
\Cref{thm:ast-const-opt-primal}.%
\indexg{Gliding parabola!constrained optimality and|)}%
\indexg{constrained optimality conditions, astral!general|)}%
\end{example}

\indexg{constrained optimality conditions, astral!affine constraints@with affine constraints|(}%
The next theorem specializes
\Cref{pr:cons-opt-equiv-soln} and
\Cref{thm:ast-const-opt-primal}
to the problem of minimization
of a convex function subject to affine constraints,
in other words,
to the case
when the constraint set $C$ is an affine set.
The theorem could be proved as a corollary of 
\Cref{thm:ast-const-opt-primal}
by computing the normal cone $\ncone{\xbar}{\Cbar}$ for this special
case.
We instead give a direct proof using
\Cref{thm:subgrad-equiv-fenchel}.

\begin{theorem}    \label{thm:ast-opt-lin-cons-primal}
  Let $f:\Rn\rightarrow\Rext$ be convex and proper,
  let $\A\in\Rmn$,
  $\bb\in\Rm$,
  and $C=\{\xx\in\Rn :\: \A\xx=\bb\}$.
  Let $p=f + \indC$, and let $\xbar\in\extspace$.
  Assume $\inf p > -\infty$, and that
  there exists $\xing\in\ri(\dom{f})$ such that $\A\xing=\bb$.
  Then the following are equivalent:
  \begin{letter-compact}
  \item   \label{thm:ast-opt-lin-cons-primal:a}
    $\xbar$ minimizes $\pext$.
  \item   \label{thm:ast-opt-lin-cons-primal:b}
    $\A\xbar=\bb$ and
    $\xbar$ minimizes $\fext$ over all
    $\xbar'\in\extspace$ with $\A\xbar'=\bb$.
  \item   \label{thm:ast-opt-lin-cons-primal:c}
    There exists a sequence $\seq{\xx_t}$ in $\Rn$ converging to
    $\xbar$ with $\A\xx_t=\bb$ for all $t$ and
    $f(\xx_t)\rightarrow \inf p$.
  \item   \label{thm:ast-opt-lin-cons-primal:d}
    $\A\xbar=\bb$ and
    $\asubdiffext{\xbar} \cap (\colspace \transA)\neq\emptyset$.
  \end{letter-compact}
\end{theorem}

\begin{proof}
In a moment,
we will apply
\Cref{thm:subgrad-equiv-fenchel} with $f$ and $\A$ as given, and
letting
$g=\indf{\{\bb\}}$, the indicator function on $\Rm$ for the
singleton $\{\bb\}$.
Then for $\xx\in\Rn$,
\[
  p(\xx)
  =
  f(\xx)+\indC(\xx)
  =
  f(\xx)+\indf{\{\bb\}}(\A\xx)
  =
  f(\xx)+g(\A\xx),
\]
where the second equality is because
$\xx\in C$ if and only if $\A\xx\in\{\bb\}$.
Thus, $p=f+g\A$.
Furthermore,
$\A\xing\in\{\bb\}=\ri\,\{\bb\}=\ri(\dom g)$.

\begin{proof-parts}
\pfpart{%
  (\ref{thm:ast-opt-lin-cons-primal:a})
  $\Rightarrow$
  (\ref{thm:ast-opt-lin-cons-primal:d}):
}
Suppose $\xbar$ minimizes $\pext$.
Let $\uu\in\Rm$ be a minimizer of $\fstar[-\transA]+\gstar$,
which exists by \Cref{thm:std-fenchel-duality}.
Then by
\Cref{thm:subgrad-equiv-fenchel}(\ref{thm:subgrad-equiv-fenchel:a},\ref{thm:subgrad-equiv-fenchel:c}),
it follows that $-\transAk\uu\in\asubdiffext{\xbar}$
and $\uu\in\asubdifgext{\A\xbar}$.
The latter inclusion implies, by
\Cref{pr:subgrad-imp-in-cldom}(\ref{pr:subgrad-imp-in-cldom:c}),
that $\A\xbar\in\cldom{g}=\{\bb\}$, so
$\A\xbar=\bb$.
And
since $-\transAk\uu$ is evidently in $\colspace\transA$,
the former inclusion implies that
$\asubdiffext{\xbar} \cap (\colspace \transA)\neq\emptyset$.

\pfpart{%
  (\ref{thm:ast-opt-lin-cons-primal:d})
  $\Rightarrow$
  (\ref{thm:ast-opt-lin-cons-primal:a}):
}
Suppose $\A\xbar=\bb$ and that there exists a point
$\ww\in\asubdiffext{\xbar} \cap (\colspace \transA)$.
Since $\ww$ is in the column space of $\transA$, we can write it as
$\ww=-\transAk\uu$ for some $\uu\in\Rm$, so
$-\transAk\uu\in\asubdiffext{\xbar}$.
Also, $\indfstar{\{\bb\}}(\uu)=\bb\cdot\uu\in\R$
(by Eq.~\ref{eq:e:2}), so we can apply
\Cref{cor:subdif-ind-fcn-rn}(\ref{cor:subdif-ind-fcn-rn:b},\ref{cor:subdif-ind-fcn-rn:a})
(with $\xbar$, $S$, $\uu$, as they appear in that
\namecref{cor:subdif-ind-fcn-rn}, set to $\bb$, $\{\bb\}$, $\uu$), yielding
$\uu\in\asubdifindext{\{\bb\}}{\bb}=\asubdifgext{\bb}=\asubdifgext{\A\xbar}$.
By
\Cref{thm:subgrad-equiv-fenchel}(\ref{thm:subgrad-equiv-fenchel:c},\ref{thm:subgrad-equiv-fenchel:a}),
it now follows that $\xbar$ minimizes $\pext$.

\pfpart{%
  (\ref{thm:ast-opt-lin-cons-primal:a})
  $\Leftrightarrow$
  (\ref{thm:ast-opt-lin-cons-primal:b})
  $\Leftrightarrow$
  (\ref{thm:ast-opt-lin-cons-primal:c}):
}
Since $C$ is affine, it is also relatively open, so
$\xing\in C=\ri C$, and hence,
$\ri(\dom f)\cap\ri{C}\neq\emptyset$.
Also, by \Cref{cor:closure-affine-set},
$C$'s closure in $\extspace$ is
$\Cbar=\{\zbar\in\extspace :\: \A\zbar=\bb\}$.
With these observations,
these equivalences now follow directly from
\Cref{pr:cons-opt-equiv-soln}.
\qedhere
\end{proof-parts}
\end{proof}

\begin{example}    \label{ex:std-fenchel-equiv-not-exist:cont2}
\indexg{Gliding parabola!constrained optimality and|(}%
Suppose, as in Example~\ref{ex:std-fenchel-equiv-not-exist:cont1},
that $f(\xx)=x_1^2+e^{-x_2}$ for $\xx\in\R^2$,
but that our aim now is to minimize $f$ over the affine set
$C=\{\xx\in\R^2 :\: \xx\cdot\ee_1 = 1\}$.
Evidently, this modified problem is also solved by the same point
$\xbar=\limray{\ee_2}\plusl\ee_1$.
We can check this using
\Cref{thm:ast-opt-lin-cons-primal}.
To do so, we let
$\A=\trans{\ee_1}$ and $b=1$, and let
$p=f + \indC$.
The conditions of
\Cref{thm:ast-opt-lin-cons-primal}(\ref{thm:ast-opt-lin-cons-primal:d})
are satisfied:
We have $\A\xbar=\xbar\cdot\ee_1=1$.
Further, $2\ee_1$ is in $\asubdiffext{\xbar}$, as shown in
\Cref{ex:std-fenchel-equiv-not-exist:cont1},
and is also in $\colspace \ee_1=\colspace\transA$.
Therefore, $\xbar$ minimizes $\pext$ by
\Cref{thm:ast-opt-lin-cons-primal}(\ref{thm:ast-opt-lin-cons-primal:d},\ref{thm:ast-opt-lin-cons-primal:a}).%
\indexg{constrained optimality conditions, astral!affine constraints@with affine constraints|)}%
\indexg{Gliding parabola!constrained optimality and|)}%
\end{example}

So far, we have seen how an astral version of Fenchel-duality theory
can be used to characterize constrained minima attained at points
at infinity, where standard first-order criteria are inapplicable.
But there are also cases in which a constrained optimization problem
might have a finite minimizer and yet such standard criteria are
inapplicable because the objective function is not
(standard) subdifferentiable at its minimizer.
We will see such an example shortly
(\Cref{ex:entropy-ast-dual-subgrad:min-C}).
This particular difficulty can never be a problem when working with
astral dual subgradients since, as shown in \Cref{thm:adsubdiff-nonempty},
every convex function is astral dual subdifferentiable at every point.

\indexg{constrained optimality conditions, astral!general|(}%
The next theorem provides a tool for taking advantage of this fact in
such situations.
Under the stated assumptions, the theorem characterizes the
minimizers of a convex function $\psi:\Rn\rightarrow\Rext$ over a
convex set $C\subseteq\Rn$, providing an astral dual analogue of
standard first-order optimality criteria.

\begin{theorem}    \label{thm:fenchel-cons-opt-dual}
  Let $\psi:\Rn\rightarrow\Rext$ be closed, proper and convex,
  and let $C\subseteq\Rn$ be closed and convex.
  Assume that:
  \begin{roman-compact}
  \item    \label{thm:fenchel-cons-opt-dual:asm:1}
    $[-\ri(\dom \psistar)] \cap \ri(\barr C) \neq\emptyset$; and
  \item    \label{thm:fenchel-cons-opt-dual:asm:2}
    $C \cap (\dom \psi) \neq\emptyset$.
  \end{roman-compact}
  Let $\uu\in C$.
  Then the following are equivalent:
  \begin{letter-compact}
  \item    \label{thm:fenchel-cons-opt-dual:a}
    $\uu$ minimizes $\psi$ over $C$.
  \item    \label{thm:fenchel-cons-opt-dual:b}
    $[-\adsubdifpsi{\uu}]\cap\ancone{\uu}{C}\neq\emptyset$.
  \end{letter-compact}
\end{theorem}

\begin{proof}
The proof is based on an application of
\Cref{thm:subgrad-equiv-fenchel}.
As such,
let $f=\indfstar{C}$, $g=\psistar$, and $\A=-\Iden$
(where $\Iden$ is the $n\times n$ identity matrix).
Then $f$ and $g$ are convex and proper
(by Propositions~\ref{pr:conj-props}\ref{pr:conj-props:d}
and~\ref{pr:conj-props-cvx}\ref{pr:conj-props-cvx:a},
and since $C$ is nonempty by
assumption~(\ref{thm:fenchel-cons-opt-dual:asm:2})).
Further, $\fstar=\indC$ and $\gstar=\psi$
(by
\Cref{pr:conj-props-cvx}\ref{pr:conj-props-cvx:b}
since
$\psi$ and $\indC$ are both convex and closed).
Therefore,
$\fstar[-\transA]+\gstar=\psi+\indC$.

We claim the assumptions of
\Cref{thm:subgrad-equiv-fenchel} are satisfied.
This is because, by
assumption~(\ref{thm:fenchel-cons-opt-dual:asm:1}),
there exists a point
$\xing\in [-\ri(\dom \psistar)] \cap \ri(\barr C)$.
Since $\psistar=g$ and
$\barr C = \dom \indfstar{C} = \dom f$
(by Eq.~\ref{eqn:bar-cone-defn}),
this implies $\xing\in\ri(\dom f)$ and
$\A\xing=-\xing\in\ri(\dom g)$.
Also, by
assumption~(\ref{thm:fenchel-cons-opt-dual:asm:2}),
there exists $\uhat\in C$ with $\psi(\uhat)<+\infty$,
so
\begin{align*}
  \inf (f+g\A)
  &=
  -\min \bigParens{\fstar[-\transA]+\gstar}
  \\
  &\geq
  - \bigParens{\fstar(\uhat) + \gstar(\uhat)}
  =
  - \bigParens{\psi(\uhat) + \indC(\uhat)}
  >
  -\infty,
\end{align*}
where the first equality is by
\Cref{thm:std-fenchel-duality}.

Note also that, for $\uu\in\Rn$,
\begin{equation}   \label{eq:thm:fenchel-cons-opt-dual:1}
  \adsubdiffstar{\uu}=\adsubdif{\indC}{\uu}=\ancone{\uu}{C}
\end{equation}
with the second equality from 
\Cref{pr:adsubdif-ind-fcn}.

\begin{proof-parts}
\pfpart{%
  (\ref{thm:fenchel-cons-opt-dual:a})
  $\Rightarrow$
  (\ref{thm:fenchel-cons-opt-dual:b}):
}
Suppose $\uu$ minimizes $\psi$ over $C$, implying
$\uu$ minimizes $\psi+\indC=\fstar[-\transA]+\gstar$.
Let $\xbar\in\extspace$ be a minimizer of
$\fgAext$ (which exists by \Cref{pr:fext-min-exists}).
Then by
\Cref{thm:subgrad-equiv-fenchel}(\ref{thm:subgrad-equiv-fenchel:a},\ref{thm:subgrad-equiv-fenchel:c-dual}),
$\xbar\in\adsubdiffstar{-\transAk\uu}=\adsubdiffstar{\uu}$
and
$-\xbar=\A\xbar\in\adsubdifgstar{\uu}=\adsubdifpsi{\uu}$.
Combined with
\eqref{eq:thm:fenchel-cons-opt-dual:1},
this proves the claim.

\pfpart{%
  (\ref{thm:fenchel-cons-opt-dual:b})
  $\Rightarrow$
  (\ref{thm:fenchel-cons-opt-dual:a}):
}
Suppose there exists a point
$\xbar\in [-\adsubdifpsi{\uu}]\cap\ancone{\uu}{C}$.
Then $\xbar\in\adsubdiffstar{\uu}=\adsubdiffstar{-\transAk\uu}$
by \eqref{eq:thm:fenchel-cons-opt-dual:1},
and
$\A\xbar=-\xbar\in\adsubdifgstar{\uu}$.
Therefore, by
\Cref{thm:subgrad-equiv-fenchel}(\ref{thm:subgrad-equiv-fenchel:c-dual},\ref{thm:subgrad-equiv-fenchel:a}),
$\uu$ minimizes $\fstar[-\transA]+\gstar=\psi+\indC$,
so $\uu$ minimizes $\psi$ over~$C$.
\qedhere
\end{proof-parts}
\end{proof}

\begin{example}   \label{ex:entropy-ast-dual-subgrad:min-C}
\indexg{Entropy function!constrained optimality and|(}%
Let $\psi:\R^2\rightarrow\Rext$ be the entropy function from
\Cref{ex:entropy-ast-dual-subgrad}, and consider the problem of
minimizing $\psi$ over the set
$C=\{\uu\in\R^2 :\: u_1+u_2\geq 1\}$.
Note that $\psi$ is not standard subdifferentiable at any point in
$C$, so standard first-order optimality criteria are inapplicable
here.
These criteria are also inapplicable because $C$ intersects $\dom\psi$
only at the relative boundaries of the two sets.

Nonetheless,
\Cref{thm:fenchel-cons-opt-dual}
can be used here to confirm that
$\tuu=\trans{\regBracks{\frac{1}{2},\frac{1}{2}}}$ minimizes $\psi$ over $C$.
It can be checked that the two assumptions of the theorem are satisfied:
the first because
$\psistar(\xx)=\ln(1+e^{x_1}+e^{x_2})$ for $\xx\in\R^2$,
so that $\dom\psistar=\R^2$,
and the second because, for instance,
$\ee_1\in C\cap(\dom\psi)$.
Let $\ww_2=\trans{[1,1]}$.
Then, as shown
in \Cref{ex:entropy-ast-dual-subgrad} (Eq.~\ref{eq:ex:entropy-ast-dual-subgrad:1}),
$\limray{\ww_2}\in\adsubdifpsi{\tuu}$.
Further, $-\limray{\ww_2}\in\ancone{\tuu}{C}$ since
for all $\uu\in C$,
$\ww_2\cdot(\uu-\tuu)=u_1+u_2-1\geq 0$, implying
$-\limray{\ww_2}\cdot(\uu-\tuu)\leq 0$.
Thus,
$-\limray{\ww_2}\in[-\adsubdifpsi{\tuu}]\cap\ancone{\tuu}{C}$,
so
\Cref{thm:fenchel-cons-opt-dual} proves that
$\tuu$ minimizes $\psi$ over $C$.%
\indexg{Entropy function!constrained optimality and|)}%
\indexg{constrained optimality conditions, astral!general|)}%
\end{example}

\indexg{constrained optimality conditions, astral!affine constraints@with affine constraints|(}%
The next theorem specializes
\Cref{thm:fenchel-cons-opt-dual} to characterize the minimizers of a
convex function subject to affine equality constraints,
that is, $C=\set{\uu\in\Rn :\: \transAk\uu=\bb}$ for some $\A\in\Rnm$ and
$\bb\in\Rm$.
This could be proved as a corollary of
\Cref{thm:fenchel-cons-opt-dual} by calculating the barrier cone and dual
normal cone for an affine set.
We instead give a direct proof
based on \Cref{thm:subgrad-equiv-fenchel}.

\begin{theorem}    \label{cor:fenchel-aff-constraints}
  Let $\psi:\Rn\rightarrow\Rext$ be closed, proper and convex,
  let $\A\in\R^{n\times m}$ and $\bb\in\Rm$.
  Assume that:
  \begin{roman-compact}
  \item    \label{cor:fenchel-aff-constraints:asm:1}
    $\ri(\dom{\psistar})\cap(\colspace{\A})\neq\emptyset$;
    and
  \item    \label{cor:fenchel-aff-constraints:asm:2}
    there exists $\uhat\in\Rn$ such that
    $\transAk\uhat=\bb$ and $\psi(\uhat)<+\infty$.
  \end{roman-compact}
  Let $\uu\in\Rn$ be such that $\transAk\uu=\bb$.
  Then the following are equivalent:
  \begin{letter-compact}
  \item    \label{cor:fenchel-aff-constraints:a}
    $\uu$ minimizes $\psi$ over all
    $\uu'\in\Rn$ with $\transAk\uu'=\bb$.
  \item    \label{cor:fenchel-aff-constraints:b}
    $\adsubdifpsi{\uu}\cap(\acolspace{\A})\neq\emptyset$.
  \end{letter-compact}
\end{theorem}

\begin{proof}
The proof is similar to that of \Cref{thm:fenchel-cons-opt-dual}.
Let
$g=\psistar$, and let $f:\Rm\rightarrow\R$ be defined by
$f(\xx)=-\xx\cdot\bb$ for $\xx\in\Rm$.
Then $f$ and $g$ are convex and proper
(by Propositions~\ref{pr:conj-props}\ref{pr:conj-props:d}
and~\ref{pr:conj-props-cvx}\ref{pr:conj-props-cvx:a}).
Note also that
$f=\indfsingstar{-\bb}$ (by Eq.~\ref{eq:e:2}),
the support function for the singleton $\{-\bb\}$.
Further, $\fstar=\indfsing{-\bb}$
and $\gstar=\psi$
(by
\Cref{pr:conj-props-cvx}\ref{pr:conj-props-cvx:b}
since $\psi$ and $\indfsing{-\bb}$ are both closed).
Therefore, for $\ww\in\Rn$,
\begin{equation}    \label{eq:cor:fenchel-aff-constraints:1}
  \fstar\regParens{-\transA \ww} + \gstar(\ww)
  =
  \indfsing{-\bb}(-\transA\ww) + \psi(\ww)
  =
  \begin{cases}
    \psi(\ww) & \text{if $\transA \ww = \bb$,} \\
    +\infty   & \text{otherwise.}
  \end{cases}
\end{equation}

Below, we apply \Cref{thm:subgrad-equiv-fenchel} with
$f$ and $g$ as
above, and with $\A$ as given.
For these choices, we claim that the assumptions
of \Cref{thm:subgrad-equiv-fenchel} hold.
First, by assumption~(\ref{cor:fenchel-aff-constraints:asm:1}), there
exists a point $\xhat\in\Rm$ such that
$\A\xhat\in\ri(\dom{\psistar})=\ri(\dom{g})$.
Furthermore, $\xhat\in\ri(\dom{f})$ since $\dom{f}=\Rm$.

Next, for $\uhat$ as in
assumption~(\ref{cor:fenchel-aff-constraints:asm:2}),
we have that
$\fstar(-\transAk\uhat) + \gstar(\uhat) < +\infty$
by
\eqref{eq:cor:fenchel-aff-constraints:1}.
Combined with \Cref{thm:std-fenchel-duality},
this implies that
$\inf (f+g\A) > -\infty$.

\begin{proof-parts}
\pfpart{%
  (\ref{cor:fenchel-aff-constraints:a})
  $\Rightarrow$
  (\ref{cor:fenchel-aff-constraints:b}):
}
Suppose statement~(\ref{cor:fenchel-aff-constraints:a}) holds,
implying, by
\eqref{eq:cor:fenchel-aff-constraints:1},
that $\uu$ minimizes
$\fAgstar$.
Furthermore, there exists $\xbar\in\extspac{m}$ that minimizes
$\fgAext$ (by \Cref{pr:fext-min-exists}).
By
\Cref{thm:subgrad-equiv-fenchel}(\ref{thm:subgrad-equiv-fenchel:a},\ref{thm:subgrad-equiv-fenchel:c-dual}),
applied as described above,
it therefore follows that
$\A\xbar\in\adsubdifgstar{\uu}=\adsubdifpsi{\uu}$.
Since $\A\xbar$ is also in $\acolspace{\A}$,
this proves the claim.

\pfpart{%
  (\ref{cor:fenchel-aff-constraints:b})
  $\Rightarrow$
  (\ref{cor:fenchel-aff-constraints:a}):
}
Suppose statement~(\ref{cor:fenchel-aff-constraints:b}) holds.
Then there exists $\xbar\in\extspac{m}$ such that
$\A\xbar\in\adsubdifpsi{\uu}=\adsubdifgstar{\uu}$.
In addition, $\transAk\uu=\bb$, so
$-\transAk\uu\in\{-\bb\}=\basubdiffext{\xbar}$
by Example~\ref{ex:affine-subgrad-new}.
Together, by
\Cref{thm:subgrad-equiv-fenchel}(\ref{thm:subgrad-equiv-fenchel:c-mixed},\ref{thm:subgrad-equiv-fenchel:a}),
these imply that
$\uu$ minimizes $\fAgstar$.
By
\eqref{eq:cor:fenchel-aff-constraints:1}, this means
statement~(\ref{cor:fenchel-aff-constraints:a}) holds.
\qedhere
\end{proof-parts}
\end{proof}

\begin{example}   \label{ex:entropy-ast-dual-subgrad:min-aff}
\indexg{Entropy function!constrained optimality and|(}%
Continuing Example~\ref{ex:entropy-ast-dual-subgrad:min-C},
let $\psi:\R^2\rightarrow\Rext$ be the entropy function
as defined in \Cref{ex:entropy-ast-dual-subgrad},
and consider minimizing $\psi(\uu)$ over $\uu\in\R^2$,
now subject to the affine constraint that $u_1+u_2=1$.
We can apply
\Cref{cor:fenchel-aff-constraints} to show that
$\tuu=\trans{[\frac{1}{2},\frac{1}{2}]}$ minimizes $\psi$ subject to
this constraint.
To do so,
we let $\A=[\ww_2]$ and $b=1$,
where $\ww_2=\trans{[1,1]}$.
The assumptions of \Cref{cor:fenchel-aff-constraints} are then satisfied since
$\dom\psistar=\R^2$, and, for instance,
$\psi(\ee_1)<+\infty$.
As shown
in Example~\ref{ex:entropy-ast-dual-subgrad}
(Eq.~\ref{eq:ex:entropy-ast-dual-subgrad:1}),
$\limray{\ww_2}\in\adsubdifpsi{\tuu}$.
Since $\limray{\ww_2}=\A\omm\in\acolspace{\A}$,
it therefore follows from
\Cref{cor:fenchel-aff-constraints}
that $\tuu$ minimizes $\psi$ subject to this affine constraint.%
\indexg{Entropy function!constrained optimality and|)}%
\indexg{constrained optimality conditions, astral|)}%
\indexg{constrained optimality conditions, astral!affine constraints@with affine constraints|)}%
\end{example}

\section{Application to astral matrix games}
\label{sec:astral-matrix-games}

\indexg{games, matrix!standard|(}%
As another application, we can apply Fenchel duality to obtain an astral
generalization of
von Neumann's minimax theorem for two-person,
zero-sum matrix games
\idxneumann%
\citep{vonNeumann28alt}.
In one version of this setting, there are two players,
Minnie and Max, playing a game described by a matrix $\A\in\Rmn$.
To play the game, Minnie chooses a probability vector $\xx\in\probsim$
(where $\probsim$ is as defined in Eq.~\ref{eqn:probsim-defn}),
and simultaneously, Max chooses a probability vector
$\uu\in\probsimm$.
In this setting, such vectors are often called
\indexg{strategies}%
\emph{strategies}.
The outcome of the game is then
\[
   \utransA \xx
   =
   (\A \xx)\cdot\uu
   =
   \xx\cdot(\transAk\uu),
\]
which represents the ``loss'' to Minnie that she aims to minimize,
and the ``gain'' to Max that he aims to
\indexg{games, matrix!standard|)}%
maximize.
\indexg{minimax theorems!standard|(}%
Von Neumann's minimax theorem states that
\begin{equation}       \label{eq:std-von-neumann}
   \min_{\xx\in\probsim} \max_{\uu\in\probsimm} \Bracks{(\A\xx)\cdot\uu}
   =
   \max_{\uu\in\probsimm} \min_{\xx\in\probsim} \Bracks{(\A\xx)\cdot\uu},
\end{equation}
with all minima and maxima attained.
The common value $v$ is called the
\indexg{value (of game)}%
\emph{value} of the game.
This theorem means that there exists a
\indexg{minimax strategy}%
minimax strategy for Minnie such
that, no matter how Max plays, her loss will never exceed $v$;
similarly, there exists a
\indexg{maximin strategy}%
maximin strategy for Max that ensures his
gain can never be less than $v$, regardless of how Minnie plays.%

We can also interpret \eqref{eq:std-von-neumann}
in terms of sequential, rather than simultaneous, play of the game:
Suppose Minnie chooses (and reveals) her strategy $\xx$ first.
Then knowing $\xx$, Max will choose $\uu$ to maximize his
gain $(\A\xx)\cdot\uu$.
Thus, if Minnie plays $\xx$, her gain will be
$\max_{\uu\in\probsimm} \Bracks{(\A\xx)\cdot\uu}$.
Knowing this, she will therefore want to choose $\xx$ to minimize this
expression, yielding the outcome on the left-hand side of
\eqref{eq:std-von-neumann}.
Likewise, the expression on the right-hand side represents the outcome
if instead Max plays first.
Intuitively, if Minnie plays first, then Max should have the
advantage, knowing Minnie's chosen strategy.
The fact that the expressions in
\eqref{eq:std-von-neumann}
are equal means that actually there is no such advantage for such
games, and that if Minnie blindly plays her minimax strategy (which
must exist),
then she will do just as well as if she had known Max's chosen
strategy
before choosing her own.

We can straightforwardly generalize this setup, replacing
$\probsim$, the set of strategies from which Minnie can choose,
with an arbitrary closed, convex set $X\subseteq\Rn$, and similarly
replace $\probsimm$ with some closed, convex set $U\subseteq\Rm$.
Now, however, if either set is unbounded, there might not exist a
minimax (or maximin) strategy, as the next example shows.

\begin{example}   \label{ex:game-no-minmax-w-strat}
Let
\[
  X = \Braces{
    \xx\in\R^3 :\:
    e^{x_1} + 1 \leq x_2
  },
\]
and 
let $U = \probsimgen{2} = \Braces{ \ulam :\: \gamep\in[0,1] }$,
where $\ulam=\trans{[1-\gamep, \gamep]}$ for $\gamep\in[0,1]$.
Both $X$ and $U$ are closed and convex.
Finally, let
\begin{equation}   \label{eq:ex:game-no-minmax-w-strat:2}
  \A =
  \begin{bmatrix*}[r]
      0 & 1 &  1 \\
      0 & 1 & -1
  \end{bmatrix*}.
\end{equation}
Thus,
if Minnie plays $\xx\in X$ and Max plays $\ulam$, then the outcome
will be
\begin{equation}   \label{eq:ex:game-no-minmax-w-strat:1}
  (\A\xx)\cdot\ulam = x_2 + (1-2\gamep) x_3.
\end{equation}

In this game, there exists a maximin strategy:
If Max plays first and chooses $\gamep\neq 1/2$, then
the infimum of the outcome in
\eqref{eq:ex:game-no-minmax-w-strat:1} over Minnie's strategy $\xx\in X$
will be $-\infty$
(taking the limit $x_3\rightarrow\pm\infty$).
Otherwise, if Max chooses $\gamep=1/2$, then
the infimum over $X$ of the outcome in
\eqref{eq:ex:game-no-minmax-w-strat:1} is $1$.
Thus, $\ulamhalf$ is a maximin strategy
(or really, ``maxi-inf'' strategy).

On the other hand, there is no minimax strategy for this game:
If Minnie plays first with strategy $\xx\in X$,
then to maximize the outcome in
\eqref{eq:ex:game-no-minmax-w-strat:1},
Max will choose $\gamep=0$ if $x_3>0$ and $\gamep=1$ if $x_3<0$
(and arbitrarily if $x_3=0$).
Thus, in all cases, the outcome will be
$x_2+|x_3|$.
Hence, Minnie can achieve an outcome of $1+\epsilon$, for all
$\epsilon\in\Rstrictpos$,
by choosing
$\xx=\trans{[\ln\epsilon,\,1+\epsilon,\,0]}$;
therefore, the infimum
of the outcome over her choices is $1$.
Nonetheless, there is no strategy $\xx\in X$ that allows her to
achieve this infimum, so no minimax strategy exists for this game.

Summarizing, the arguments above show that
\[
   \inf_{\xx\in X} \max_{\uu\in U} \Bracks{(\A\xx)\cdot\uu}
   =
   1
   =
   \max_{\uu\in U} \inf_{\xx\in X} \Bracks{(\A\xx)\cdot\uu},
\]
with the maxima attained, but not the infima.
\end{example}

Thus, in this example, there exists a maximin strategy for Max, but no
minimax strategy for
\indexg{minimax theorems!standard|)}%
Minnie.
\indexg{games, matrix!astral|(}%
It seems natural to extend such a case to astral space where
Minnie would have a minimax strategy.
To do so, we allow Minnie to choose her astral strategy $\xbar$ from
$\Xbar$, the astral closure of $X$.
The resulting outcome then becomes
\indexg{games, matrix!astral|)}%
$(\A\xbar)\cdot\uu$.
\indexg{minimax theorems!astral|(}%
In this astral version of the problem, there must exist minimax and
maximin strategies for both sides under a topological
condition, as shown in the next
\namecref{thm:ast-von-neumann}:

\begin{theorem}   \label{thm:ast-von-neumann}
  Let $X\subseteq\Rn$ and $U\subseteq\Rm$ be convex and nonempty, and
  assume also that $U$ is closed (in $\Rm$).
  Let $\A\in\Rmn$.
  Assume there exists $\xing\in\ri X$ such that
  $\A\xing\in\ri(\barr U)$.
  Then
  \begin{equation}   \label{eqn:thm:ast-von-neumann:8}
    \adjustlimits\min_{\xbar\in\Xbar} \sup_{\uu\in U} [(\A\xbar)\cdot\uu]
    =
    \adjustlimits\max_{\uu\in U} \min_{\xbar\in\Xbar} [(\A\xbar)\cdot\uu]
  \end{equation}
  with the minima and the maximum attained.
\end{theorem}

\begin{proof}
{\mathtogether%
The key step of the proof is in applying
Fenchel duality (\Cref{thm:std-fenchel-duality})
with $f=\indf{X}$ and $g=\indstar{U}$.
This is justified
because both $\indf{X}$ and
$\indf{U}$ are convex and proper
(since $X$ and $U$ are convex and nonempty),
so $\indstar{U}$ is as well
(Propositions~\ref{pr:conj-props}\ref{pr:conj-props:d}
and~\ref{pr:conj-props-cvx}\ref{pr:conj-props-cvx:a}).
Also, by assumption, 
$\xing\in\ri X=\ri(\dom\indf{X})$, and 
$\A\xing\in\ri(\barr U)=\ri(\dom\indstar{U})$ (using
Eq.~\ref{eqn:bar-cone-defn}).
Note further that
$\indf{U}$ is closed since $U$ is a closed set, so
$\inddubstar{U}=\indf{U}$
(\Cref{pr:conj-props-cvx}\ref{pr:conj-props-cvx:b}).}

We first prove \eqref{eqn:thm:ast-von-neumann:8} with infima rather
than minima, and then show below that each is attained.
We have
\begin{align}
\notag
  \adjustlimits\inf_{\xbar\in\Xbar} \sup_{\uu\in U} {[(\A\xbar)\inprod\uu]}
  &\le
  \adjustlimits\inf_{\xx\in X} \sup_{\uu\in U} {[(\A\xx)\cdot\uu]}
  =
  \inf_{\xx\in X} \indstar{U}(\A\xx)
\\
\notag
  &=
  \inf_{\xx\in\Rn} \bigBracks{\indf{X}(\xx) + \indstar{U}(\A\xx)}
\\
\notag
  &=
  -\min_{\uu\in\Rm} \bigBracks{\indstar{X}(-\transAk\uu) + \inddubstar{U}(\uu)}
\\
\notag
  &=
  -\min_{\uu\in\Rm} \bigBracks{\indfastar{\Xbar}(-\transAk\uu) + \indf{U}(\uu)}
\\
\notag
  &
  =
  -\min_{\uu\in U}  \indfastar{\Xbar}(-\transAk\uu)
\\
\notag  
  &
  =
  -\adjustlimits\min_{\uu\in U} \sup_{\xbar\in \Xbar} {[\xbar\inprod(-\transAk\uu)]}
\\
\label{eq:vonNeumann:4}
  &=
  \adjustlimits\max_{\uu\in U} \inf_{\xbar\in \Xbar} {[(\A\xbar)\inprod\uu]}.
\end{align}
The inequality is because $X\subseteq\Xbar$.
The first equality is by $\indstar{U}$'s definition
(Eq.~\ref{eq:e:2}).
The third equality is by \Cref{thm:std-fenchel-duality},
whose application with $f=\indf{X}$ and $g=\indstar{U}$
was justified above.
The fourth equality is because
$\inddubstar{U}=\indf{U}$,
and
$\indstar{X}=\indgenextstar{X}=\indfastar{\Xbar}$ (by Propositions~\ref{pr:fextstar-is-fstar}
and~\ref{pr:inds-ext}).
The fifth equality is because
the minimum
on the fourth line must be attained at a point in $U$
since if it were attained at a point in $\Rm\wo U$,
then the argument to the minimum would be $+\infty$
for all $\uu\in\Rm$, implying the minimum is attained also at any point
in $U$.
The sixth equality is by $\indfastar{\Xbar}$'s definition
(Eq.~\ref{eqn:astral-support-fcn-def}).
The last equality uses \Cref{thm:Ax-dot-u}.

For the opposite inequality, let $\uu_0\in U$
attain the maximum in the
final expression of \eqref{eq:vonNeumann:4}. Then
\[
    \adjustlimits\max_{\uu\in U} \inf_{\xbar\in\Xbar} {[(\A\xbar)\inprod\uu]}
    =
    \inf_{\xbar\in\Xbar} {[(\A\xbar)\inprod\uu_0]}
    \le
    \adjustlimits\inf_{\xbar\in\Xbar} \sup_{\uu\in U} {[(\A\xbar)\inprod\uu]}.
\]
Combined with \eqref{eq:vonNeumann:4}, we thus have
\begin{equation}   \label{eq:thm:ast-von-neumann:9}
  \adjustlimits\inf_{\xbar\in\Xbar} \sup_{\uu\in U} {[(\A\xbar)\inprod\uu]}
  =
  \adjustlimits\max_{\uu\in U} \inf_{\xbar\in\Xbar} {[(\A\xbar)\inprod\uu]}.
\end{equation}

It remains to show that both infima are attained.
Note first that
$\Xbar$, being closed, is also compact
(\Cref{prop:compact}\ref{prop:compact:closed-subset}).
Further,
for each $\uu\in U$, the astral linear function
$\xbar\mapsto(\A\xbar)\inprod\uu=\xbar\cdot(\transAk\uu)$,
for $\xbar\in\extspace$, is continuous
(\Cref{thm:i:1}\ref{thm:i:1c}), and so also lower semicontinuous.
Thus, the infimum on the right-hand side of
\eqref{eq:thm:ast-von-neumann:9}
is attained, by \Cref{thm:weierstrass}.
The continuity of such astral linear functions also implies,
by \Cref{pr:lsc-sup}, that the function
$\xbar\mapsto\sup_{\uu\in U} {[(\A\xbar)\inprod\uu]}$,
for $\xbar\in\extspace$, is lower semicontinuous.
Hence,
the infimum on the left-hand side of
\eqref{eq:thm:ast-von-neumann:9}
is attained as well (again by \Cref{thm:weierstrass}).
\end{proof}

When $U$ is also bounded (and therefore compact),
a simplified form of
\Cref{thm:ast-von-neumann}
holds:

\begin{theorem}   \label{cor:ast-von-neumann-bounded-u}
  Let $X\subseteq\Rn$ and $U\subseteq\Rm$ be convex and nonempty, and
  assume also that $U$ is compact
  (that is, bounded and closed in $\Rm$).
  Let $\A\in\Rmn$.
  Then
  \begin{equation}    \label{eq:cor:ast-von-neumann-bounded-u:0}
    \adjustlimits
    \min_{\xbar\in\Xbar} \max_{\uu\in U\vphantom{\Xbar}} [(\A\xbar)\cdot\uu]
    =
    \adjustlimits
    \max_{\uu\in U} \min_{\xbar\in\Xbar}  [(\A\xbar)\cdot\uu]
  \end{equation}
  with all minima and maxima attained.
\end{theorem}

\begin{proof}

Since $U$ is bounded, $\barr U=\Rm$
(\Cref{pr:bar-cone-props}\ref{pr:bar-cone-props:c}).
Moreover, $X$ is convex and nonempty, so
there exists $\xing\in\ri X$
(\Cref{pr:ri-props}\ref{pr:ri-props:roc-thm6.2b}),
and necessarily,
$\A\xing\in\Rm=\ri(\barr U)$.
Thus, we can apply
\Cref{thm:ast-von-neumann},
yielding
\[
    \adjustlimits\min_{\xbar\in\Xbar} \sup_{\uu\in U} [(\A\xbar)\cdot\uu]
    =
    \adjustlimits\max_{\uu\in U} \min_{\xbar\in\Xbar} [(\A\xbar)\cdot\uu],
\]
with the minima and the maximum attained.
It remains to show that the supremum on the
left-hand side is also attained.

To prove this, let $\zbar\in\extspac{m}$.
We claim the supremum of $\zbar\cdot\uu$ over $\uu\in U$
must be attained.
If $\zbar\cdot\uu=+\infty$ for some $\uu\in U$, then this $\uu$ must
attain that supremum.
Likewise, if $\zbar\cdot\uu=-\infty$ for all $\uu\in U$, then any
$\uu\in U$ attains the supremum.
Therefore, we assume henceforth that $\zbar\cdot\uu<+\infty$ for all
$\uu\in U$, and that $\zbar\cdot\uu\in\R$ for some $\uu\in U$.\looseness=-1

We can write $\zbar=\VV\omm\plusl\qq$ for some $\VV\in\R^{m\times k}$
and some $\qq\in\Rm$.
Let $L$ be the linear subspace orthogonal to the columns of $\VV$,
that is,
$L=\{\uu\in\Rm :\: \uu\perp\VV\}$.
Let $U'=U\cap L$, which, by
\Cref{pr:vtransu-zero},
is exactly the set of points $\uu\in U$ for which $\zbar\cdot\uu\in\R$.
By the above assumptions, $U'$ is nonempty, and
if $\uu\in U\setminus L$, then $\zbar\cdot\uu=-\infty$.
Thus,
\begin{equation}   \label{eq:cor:ast-von-neumann-bounded-u:1}
  \sup_{\uu\in U} \zbar\cdot\uu
  =
  \sup_{\uu\in U'} \zbar\cdot\uu
  =
  \sup_{\uu\in U'} \qq\cdot\uu,
\end{equation}
with the last equality following from
\Cref{pr:vtransu-zero}.
Since $U$ and $L$ are closed (in~$\Rm$), and since $U$ is bounded,
$U'$ is compact in $\Rm$.
Since also the function $\uu\mapsto\qq\cdot\uu$ is continuous,
the supremum in
\eqref{eq:cor:ast-von-neumann-bounded-u:1}
must be attained
(by \Cref{thm:weierstrass}), proving the claim.
\end{proof}

We return to the matrix game of \Cref{ex:game-no-minmax-w-strat}, where
we can now show that in the astral
setting there exists a minimax strategy for Minnie:

\begin{example}
\label{ex:game-no-minmax-w-strat-cont}
As in \Cref{ex:game-no-minmax-w-strat},
let
\[
  X = \Braces{
    \xx\in\R^3 :\:
    e^{x_1} + 1 \leq x_2
  },
\]
let $\ulam=\trans{[1-\gamep, \gamep]}$ for $\gamep\in[0,1]$,
and let $\A$ be the matrix from
\eqref{eq:ex:game-no-minmax-w-strat:2}.
In the astral version of the game from that example,
Minnie plays $\xbar\in\Xbar$, Max plays
$\ulam$ for some $\gamep\in[0,1]$, and the outcome is
\begin{equation}   \label{eq:ex:game-no-minmax-w-strat-cont:1}
  (\A\xbar)\cdot\ulam
  =
  \xbar\cdot (\transA\ulam)
  =
  \xbar\cdot \trans{[0,\, 1,\, 1-2\gamep]},
\end{equation}
with the first equality from \Cref{thm:Ax-dot-u},
and the second 
from the definition of $\A$.
Consequently, setting $\xbar=\limray{(-\ee_1)}\plusl\ee_2$
(where $\ee_1,\ee_2,\ee_3$ are standard basis vectors in
$\R^3$)
attains an outcome of $1$, the
value of the game, regardless of Max's play.
Furthermore, this point is in $\Xbar$, being the
limit of the sequence
$\seq{\trans{[-t,\,1+e^{-t},\,0]}}$, whose elements are all
in $X$.
Thus, $\xbar$ is an astral minimax strategy for the game, attaining
the minimum on the left-hand side of
\eqref{eq:cor:ast-von-neumann-bounded-u:0}.

The minimum on the right-hand side of that equation is attained as
well:
If $\gamep<1/2$ or $\gamep>1/2$, then
\eqref{eq:ex:game-no-minmax-w-strat-cont:1}
is minimized by setting
$\xbar=\limray{(-\ee_3)}\plusl 2\ee_2$
or
$\xbar=\limray{\ee_3}\plusl 2\ee_2$,
respectively, in both cases yielding an outcome of $-\infty$.
And if $\gamep=1/2$, then
\eqref{eq:ex:game-no-minmax-w-strat-cont:1}
is minimized by setting
$\xbar=\limray{(-\ee_1)}\plusl\ee_2$ yielding an outcome of~1.
It can be checked that each of these strategies $\xbar$ is in $\Xbar$.

The maxima in
\eqref{eq:cor:ast-von-neumann-bounded-u:0}
are attained as before:
The foregoing shows that
the maximin strategy on the right-hand side is attained by
$\ulamhalf$.
And for any $\xbar\in\Xbar$,
the outcome in
\eqref{eq:ex:game-no-minmax-w-strat-cont:1}
is maximized by
choosing $\gamep=0$ if $\xbar\cdot\ee_3>0$
and $\gamep=1$ if $\xbar\cdot\ee_3<0$
(and arbitrarily if $\xbar\cdot\ee_3=0$).%
\indexg{minimax theorems!astral|)}%
\end{example}

\chapter{KKT conditions and Lagrange multipliers}
\label{sec:KKT}

We next consider constrained optimization problems in which the goal
is to minimize a convex function subject to explicitly given convex
inequality constraints and affine equality constraints.
In standard convex analysis, such optimization problems are often
handled using the Karush-Kuhn-Tucker (KKT) conditions
(originally due to
\indexa{Karush, W.}\citealp{Karush39},
and
\indexa{Kuhn, H. W.}\indexa{Tucker, A. W.}%
\citealp{KuhnTu51}%
),
which generalize and justify the method of Lagrange multipliers.
In this chapter, we will see how these notions extend to astral space,
thereby allowing them to be used for
convex optimization problems that only have solutions at infinity.\looseness=-1

\section{The standard setting}
\label{sec:kkt-stand-set}

We begin by defining the kind of optimization problem we aim to solve
throughout this chapter:

\begin{definition}  \label{dfn:ord-cvx-opt-prog}
\indexg{ordinary convex programs|(}%
\indexg{ordinary convex programs!defined|(}%
  An \emph{ordinary convex program} $P$ is defined by:
  \begin{item-compact}
  \item
    an
\indexg{objective function (of convex program)}%
    \emph{objective function}
    $f:\Rn\rightarrow\Rext$ that is convex and proper;
  \item
\indexg{inequality constraints (of convex program)}%
    \emph{inequality constraints} given by functions
    $g_i:\Rn\rightarrow\Rext$, for $i=1,\ldots,r$,
    that are each convex and proper
    with $\dom f\subseteq\dom g_i$
    and $\ri(\dom f)\subseteq\ri(\dom g_i)$;
    and
  \item
\indexg{equality constraints (of convex program)}%
    \emph{equality constraints} given by affine functions
    $h_j:\Rn\rightarrow\R$,
    for $j=1,\ldots,s$.
  \end{item-compact}
  In such a program, the goal is to minimize $f(\xx)$ over $\xx\in\Rn$
  subject to the constraints that
  $g_i(\xx)\leq 0$ for $i=1,\ldots,r$,
  and that
  $h_j(\xx) = 0$ for $j=1,\ldots,s$.%
\indexg{ordinary convex programs!defined|)}%
\end{definition}
We assume the notation of this definition throughout this chapter.
Also, if there is only one inequality or one equality constraint, 
we simply write $g$ or $h$ for its defining function.\looseness=-1

We can combine the various components of such a program into
a single convex function to be minimized.
To do so, we use the shorthand notation
$\negf=\indfAlt{(-\infty,0]}$ and
$\indfa{\le 0\vphantom{[}}=\indfa{[-\infty,0]}$
for the indicator functions of the
sets of nonpositive numbers in $\R$ and $\Rext$, respectively.
As usual,
we also write $\indz$ and $\indaz$ for the indicator functions of the
singleton $\{0\}$ in $\R$ and $\Rext$.
Note that
$\Negf=\negfext$
and
$\indazAlt=\indzext$
(by \Cref{pr:inds-ext}).
Solving the ordinary convex program $P$ is then equivalent to
minimizing the convex function $p:\Rn\rightarrow\Rext$ given by
\begin{equation}  \label{eqn:kkt-p-defn}
  p(\xx)
  =
  f(\xx)
  +
  \sumitor (\Negf\circ g_i)(\xx)
  +
  \sumjtos (\indzohj)(\xx),
\end{equation}
since this function is $+\infty$ if any constraint of $P$ is violated,
and otherwise is equal to $f$.

\indexg{feasibility (of convex program)!standard (primal)|(}%
A point $\xx\in\Rn$ is said to be
\emph{(standard) feasible}
for $P$ if $\xx$
satisfies all of the constraints so that
$g_i(\xx)\leq 0$ for $i=1,\ldots,r$,
and
$h_j(\xx)=0$ for $j=1,\ldots,s$. Note that we do not require
that $\xx\in\dom{f}$ (unlike some other authors, including
\idxroc\citealp[Section~28]{ROC}).
We let $C$ denote the set of all standard feasible points:
\begin{equation}  \label{eqn:feas-set-defn}
  C
  =
  \Braces{\xx\in\Rn :\: g_i(\xx)\leq 0 \mbox{ for $i=1,\ldots,r$};\;
                      h_j(\xx) = 0 \mbox{ for $j=1,\ldots,s$}
         }.%
\indexg{feasibility (of convex program)!standard (primal)|)}%
\end{equation}
A feasible point $\xx\in C$ is a
\indexg{solution (of convex program)!standard}%
\emph{(standard) solution} of $P$ if
$f$ is minimized by $\xx$ among all feasible points, or equivalently,
if $\xx$ minimizes $p$.
Thus, the program $P$ encodes a specific form of constrained
optimization of the kind studied in
\Cref{sec:constrained-optimization}.
\indexg{value (of convex program)|(}%
We call $\inf p$ the \emph{value} of program $P$;
it is the value of $f(\xx)$ at any solution $\xx$, and more generally,
the infimum of $f$ over all feasible points $\xx\in C$.
\indexg{value (of convex program)|)}%

\indexg{multipliers (of convex program)|(}%
A classical approach to solving such a problem is to introduce new
variables which are then used to form
a simpler function which can then be optimized or reasoned about.
The new variables or \emph{multipliers}, one for each constraint,
are denoted
$\lambda_i$ for $i=1,\ldots,r$
and
$\mu_j$ for $j=1,\ldots,s$.
We also write these in vector form as $\lamvec\in\R^r$ and
\indexg{multipliers (of convex program)|)}%
$\muvec\in\R^s$.
\indexg{Lagrangian|(}%
\indexg{Lagrangian!defined|(}%
For such vectors $\lamvec,\muvec$, we define the function
$\lagrangflm:\Rn\rightarrow\Rext$ by
\begin{equation}  \label{eqn:lagrang-defn}
  \lagranglm{\xx}
  =
  f(\xx)
  +
  \sumitor \lambda_i g_i(\xx)
  +
  \sumjtos \mu_j h_j(\xx),
\end{equation}
for $\xx\in\Rn$.
We often view $\lagranglm{\xx}$ as a function of $\xx$ only, with
$\lammupair$ fixed, but sometimes we regard it as
a function of both $\xx$ and $\lammupair$,
called the
\indexg{Lagrangian!defined|)}%
\emph{Lagrangian}.
Note importantly that we use the notation
$\lagranglmext{\xbar}$ and $\lagranglmstar{\uu}$
to denote the extension and conjugate of
$\lagranglm{\xx}$ regarded as a function of $\xx$ only, with
$\lammupair$ fixed.
That is, if
$h(\xx)=\lagranglm{\xx}$ for $\xx\in\Rn$, then
$\lagranglmext{\xbar}=\hext(\xbar)$ and
$\lagranglmstar{\uu}=\hstar(\uu)$;
in particular,
$\asubdiflagranglmext{\xbar}=\asubdifhext{\xbar}$.

\indexg{feasibility (of convex program)!dual|(}%
Generally, we require $\lambda_i\geq 0$ for $i=1,\ldots,r$
so that
$\lammupair\in\Lampairs$
where
\indexm{gammars}{$\Lampairs$}{all pairs in $\Rpos^r\times\R^s$}%
$\Lampairs=\Rpos^r\times\R^s$.
This condition is called
\emph{dual feasibility}.%
\indexg{feasibility (of convex program)!dual|)}%

\indexg{Slater's condition|(}%
In addition to the assumptions defining an ordinary convex program,
we assume throughout this chapter
that the program $P$ satisfies
\emph{Slater's condition}
requiring
that there exist a point $\xing\in\ri(\dom{f})$
that is ``strictly feasible'' in the sense that
$g_i(\xing)<0$ for $i=1,\ldots,r$,
and
$h_j(\xing)=0$ for $j=1,\ldots,s$.
We refer to such a point $\xing$
as a
\emph{witness for Slater's condition}.%
\indexg{Slater's condition|)}%

Here are some simple facts regarding these various notions:

\begin{proposition}  \label{pr:lagran-facts}
  Let $P$ be an ordinary convex program, as in
  \Cref{dfn:ord-cvx-opt-prog},
  and let $p$ and $\lagrangfcn$ be as in
  Eqs.~(\ref{eqn:kkt-p-defn}) and (\ref{eqn:lagrang-defn}).
  Let $\lammupair\in\Lampairs$.
  Then:
  \begin{letter-compact}
  \item  \label{pr:lagran-facts:b}
    $\dom{\lagrangflm} = \dom{f}$.
  \item  \label{pr:lagran-facts:cvx}
    $\lagrangflm$ is convex and proper.
  \item  \label{pr:lagran-facts:sup}
    For $\xx\in\Rn$,
    $p(\xx) = \sup\regBraces{\lagranglmp{\xx} :\: \lammupairp\in\Lampairs}$.
  \item  \label{pr:lagran-facts:a}
    $\lagrangflm \leq p$.
  \item  \label{pr:lagran-facts:c}
    If $P$ satisfies Slater's condition, then
    $\inf p < +\infty$.
  \end{letter-compact}
\end{proposition}

\begin{proof}
~

\begin{proof-parts}
\pfpart{Part~(\ref{pr:lagran-facts:b}):}
By assumption, $\dom{f}\subseteq\dom{g_i}$
for $i=1,\ldots,r$,
and also $\dom{h_j}=\Rn$ for $j=1,\ldots,s$.
Therefore, $\dom{\lagrangflm}=\dom{f}$.

\pfpart{Part~(\ref{pr:lagran-facts:cvx}):}
That $\lagrangflm$ is convex is straightforward using
\Cref{pr:std-sum-fcns-cvx}.
Also, $\lagrangflm>-\infty$ from its definition in terms of only
proper functions.
Moreover,
since $f$ is proper, 
$\lagrangflm\not\equiv+\infty$
by part~(\ref{pr:lagran-facts:b}).
Therefore, $\lagrangflm$ is proper.

\pfpart{Part~(\ref{pr:lagran-facts:sup}):}
Let $\xx\in\Rn$, and 
let $s(\xx)=\sup\{ \lagranglmp{\xx} :\: \lammupairp\in\Lampairs \}$.
We aim to prove that $p(\xx)=s(\xx)$.

Suppose $g_i(\xx)>0$ for some $i\in\{1,\ldots,r\}$,
implying $p(\xx)=+\infty$.
Then ${\lagranglmp{\xx}}\rightarrow+\infty$ in the limit
$\lambda'_i\rightarrow+\infty$ (with the other components of
$\lammupairp$ arbitrary but fixed).
Thus, $s(\xx)=+\infty=p(\xx)$.

Similarly,
suppose next that $h_j(\xx)\neq 0$ for some $j\in\{1,\ldots,s\}$,
again implying that $p(\xx)=+\infty$.
Then $\lagranglmp{\xx}\rightarrow+\infty$ in the limit
as $\mu'_j\rightarrow+\infty$ if $h_j(\xx)>0$,
or
as $\mu'_j\rightarrow-\infty$ if $h_j(\xx)<0$.
Thus, in this case as well, $s(\xx)=+\infty=p(\xx)$.

If neither of the above cases hold, then $\xx$ is feasible and
$p(\xx)=f(\xx)$.
For all $\lammupairp\in\Lampairs$,
we then have
$\lambda'_i\geq 0$ and $g_i(\xx)\leq 0$
for $i=1,\ldots,r$,
and $h_j(\xx)=0$ for $j=1,\ldots,s$, implying,
by \eqref{eqn:lagrang-defn}, that
$\lagranglmp{\xx}\leq f(\xx)$.
Also, $\lagrang{\xx}{\zerov{r},\zerov{s}}=f(\xx)$.
Thus, $s(\xx)=f(\xx)=p(\xx)$.

\pfpart{Part~(\ref{pr:lagran-facts:a}):}
This is immediate from
part~(\ref{pr:lagran-facts:sup}).

\pfpart{Part~(\ref{pr:lagran-facts:c}):}
Let $\xing$ be a witness for Slater's condition.
By definition, $\xing$ is feasible and in $\dom{f}$, so
$\inf p \leq p(\xing) = f(\xing) < +\infty$.
\qedhere
\end{proof-parts}
\end{proof}

\indexg{KKT vectors|(}%
\Cref{pr:lagran-facts}(\ref{pr:lagran-facts:a}) implies that, for all
$\lammupair\in\Lampairs$,
the infimum of $\lagrangflm$ cannot exceed the value of the
program $P$; that is,
\begin{equation}      \label{eq:inf-lagrang-leq-infp}
  \inf \lagrangflm
  \leq
  \inf p.
\end{equation}
We say that a vector of multipliers
$\lammupair\in\Lampairs$ is a \emph{KKT vector} for $P$
if this equation holds with equality so that
$\lagrangflm$ has the same infimum as $p$, that is,
if $\inf\lagrangflm=\inf p$.
Under Slater's condition, and if additionally
$\inf p>-\infty$, it is known that such a vector must exist
\idxroc\citep[Theorem~28.2]{ROC}.
Combined with \eqref{eq:inf-lagrang-leq-infp}, this means that
\begin{equation}   \label{eq:kkt-infsup-supinf}
  \adjustlimits\sup_{\lammupair\in\Lampairs} \inf_{\;\xx\in\Rn\;} \lagranglm{\xx}
  =
  \inf p
  =
  \adjustlimits\inf_{\;\xx\in\Rn\;} \sup_{\lammupair\in\Lampairs} \lagranglm{\xx}.
\end{equation}
The second equality here follows directly from
\Cref{pr:lagran-facts}(\ref{pr:lagran-facts:sup}),
which shows that solving the program $P$ is equivalent to solving
(that is, finding $\xx$ attaining the infimum for)
the expression on the right-hand side of \eqref{eq:kkt-infsup-supinf}.
Thus, swapping the $\inf$ and $\sup$ in the expression
on the left-hand side of \eqref{eq:kkt-infsup-supinf}
yields an expression with the same value.
Moreover,
KKT vectors, which must exist
(under the conditions above),
are exactly those vectors attaining the
supremum on the left-hand side,
whereas
solutions for $P$ are exactly the points attaining the infimum on the
right, which might or might not exist.%
\indexg{KKT vectors|)}%

\indexg{saddle points (of Lagrangian)!standard|(}%
\indexg{Lagrangian!standard saddle points of|(}%
For $\xx\in\Rn$ and $\lammupair\in\Lampairs$, we say that the pair
$\xx,\lammupair$ is a
\emph{saddle point} of the Lagrangian if
$\xx$ minimizes $\lagrangflm$ and simultaneously
$\lammupair$ maximizes $\lagrang{\xx}{\cdot\posSkip}$, that is, if
\begin{equation}
\label{eq:saddle-point}
   \lagrang{\xx}{\lamvec',\muvec'}
   \leq
   \lagranglm{\xx}
   \leq
   \lagranglm{\xx'}
\end{equation}
for all $\lammupairp\in\Lampairs$,
and
for all $\xx'\in\Rn$.
When \eqref{eq:kkt-infsup-supinf} holds (and therefore
under Slater's condition, with $\inf p>-\infty$),
this saddle-point condition is known to correspond
to
$\xx$ attaining the infimum on the right of that equation
(thereby being a solution of $P$)
and $\lammupair$ attaining the supremum on the left
(and so being a KKT vector).%
\indexg{saddle points (of Lagrangian)!standard|)}%
\indexg{Lagrangian!standard saddle points of|)}%
\indexg{Lagrangian|)}%
\looseness=-1

\indexg{KKT conditions!standard|(}%
These same conditions are further known to be equivalent to
the pair $\xx,\lammupair$ satisfying the \emph{KKT conditions},
namely:
\begin{item-compact}
\item
\indexg{feasibility (of convex program)!standard (primal)}%
  \emph{primal feasibility:}
  $\xx$ is standard feasible;
\item
\indexg{feasibility (of convex program)!dual}%
  \emph{dual feasibility:}
  $\lammupair\in\Lampairs$;
\item
\indexg{complementary slackness!standard}%
  \emph{complementary slackness:}
  $\lambda_i g_i(\xx)=0$ for $i=1,\ldots,r$;
\item
\indexg{first-order optimality (of convex program)!standard}%
  \emph{first-order optimality:}
\begin{equation}  \label{eqn:kkt-lagrang-opt}
  \zero
  \in
  \partial f(\xx)
  +
  \sum_{\substack{1\le i\le r:\\[1pt] \lambda_i>0}}
  \lambda_i \partial g_i(\xx)
  +
  \sum_{1\le j\le s}
  \mu_j \partial h_j(\xx).
\end{equation}
\end{item-compact}

Summarizing:

\begin{theorem}  \label{roc:thm28.3}
\indexg{constrained optimality conditions, standard!convex program@for convex program|(}%
\indexg{ordinary convex programs!standard optimality conditions|(}%
  Let $P$ be an ordinary convex program, as in
  \Cref{dfn:ord-cvx-opt-prog},
  that satisfies Slater's condition,
  and let $p$ and $\lagrangfcn$ be as in
  Eqs.~(\ref{eqn:kkt-p-defn}) and (\ref{eqn:lagrang-defn}).
  Assume $\inf p>-\infty$.
  Let $\xx\in\Rn$ and $\lammupair\in\Lampairs$.
  Then the following are equivalent:
  \begin{letter-compact}
  \item
    $\xx$ is a (standard) solution,
    and $\lammupair$ is a KKT vector.
  \item
    The pair $\xx,\lammupair$ satisfies the KKT conditions.
  \item
    The pair $\xx,\lammupair$ is a
    {saddle point} of the Lagrangian $\lagrangfcn$.
  \end{letter-compact}
\end{theorem}

\begin{proof}
  See \idxroc\citet[Theorem~28.3]{ROC}.
\end{proof}

As an immediate consequence,
since KKT vectors always exist under the conditions
of \Cref{roc:thm28.3},
the program $P$'s solutions can be characterized
as follows:

\begin{corollary}  \label{roc:cor28.3.1}
  Let $P$ be an ordinary convex program, as in
  \Cref{dfn:ord-cvx-opt-prog},
  that satisfies Slater's condition,
  and let $p$ and $\lagrangfcn$ be as in
  Eqs.~(\ref{eqn:kkt-p-defn}) and (\ref{eqn:lagrang-defn}).
  Assume $\inf p>-\infty$.
  Let $\xx\in\Rn$.
  Then the following are equivalent:
  \begin{letter-compact}
  \item  \label{roc:cor28.3.1:sol}
    $\xx$ is a (standard) solution of $P$.
  \item  \label{roc:cor28.3.1:kkt}
    There exists $\lammupair\in\Lampairs$
    such that the pair $\xx,\lammupair$ satisfies the KKT conditions.
  \item  \label{roc:cor28.3.1:sdl}
    There exists $\lammupair\in\Lampairs$
    such that the pair $\xx,\lammupair$ is a
    {saddle point} of the Lagrangian $\lagrangfcn$.
  \end{letter-compact}
\end{corollary}

\begin{proof}
  See \idxroc\citet[Corollary~28.3.1]{ROC}.%
\indexg{constrained optimality conditions, standard!convex program@for convex program|)}%
\indexg{ordinary convex programs!standard optimality conditions|)}%
\end{proof}

\indexg{first-order optimality (of convex program)!standard|(}%
In the first-order optimality condition given in
\eqref{eqn:kkt-lagrang-opt},
notice that the first sum is only over those
indices $i$ for which $\lambda_i>0$; in other words, the sum
explicitly omits terms for which $\lambda_i=0$.
Naively, it might seem that this condition
should be unchanged if
we instead took the sum over all terms
(not just those for which $\lambda_i>0$),
since if $\lambda_i=0$, this should ``zero out'' the product
$\lambda_i \partial g_i(\xx)$ that appears in that sum.
The problem is that if $\lambda_i=0$ and
also $\partial g_i(\xx)$ is empty,
then the corresponding term of the sum would become an empty set,
and so the entire sum would be an empty set;
consequently, the inclusion given in such a modified version
of \eqref{eqn:kkt-lagrang-opt} would never hold.

We will instead find it convenient to rewrite
\eqref{eqn:kkt-lagrang-opt}
in a slightly different form, namely, as
\begin{equation}  \label{eqn:kkt-lagrang-opt-mod}
  \zero
  \in
  \partial f(\xx)
  +
  \sumitor
  \partial (\lambda_i g_i)(\xx)
  +
  \sumjtos
  \partial (\mu_j h_j)(\xx).
\end{equation}
In this revised form,
both summations are over all indices
(rather than omitting terms with $\lambda_i=0$).
Also, the scalars $\lambda_i$ and $\mu_j$ have been moved inside the
subdifferentials so that
$\lambda_i \partial g_i(\xx)$ is replaced with
$\partial (\lambda_i g_i)(\xx)$,
and $\mu_j \partial h_j(\xx)$
with $\partial (\mu_j h_j)(\xx)$. 
This does not change the value of the expression:
if $\lambda_i>0$ then
$\lambda_i \partial g_i(\xx)=\partial (\lambda_i g_i)(\xx)$,
and if $\lambda_i=0$ then $\partial (\lambda_i g_i)(\xx)=\set{\zero}$
(since $0 g_i\equiv 0$).
Also,
$\mu_j \partial h_j(\xx)=\partial (\mu_j h_j)(\xx)$ always by ordinary
calculus since $h_j$ is affine and so differentiable everywhere.%
\indexg{first-order optimality (of convex program)!standard|)}%
\indexg{KKT conditions!standard|)}%

\indexg{ordinary convex programs!no finite solution|(}%
As we saw in \Cref{roc:cor28.3.1},
the KKT conditions characterize when a finite point
$\xx\in\Rn$ is a solution of the convex program $P$.
Nevertheless,
even under the favorable conditions of that
\namecref{roc:cor28.3.1},
it is possible that no finite solution exists.
Here is an example, which we will later analyze using astral tools.

\begin{example}   \label{ex:kkt-running-ex-standard}
  Working in $\R^2$,
  consider minimizing
  $x_2$
  subject to the constraints that
  $e^{x_1}+1\le x_2$,
  and that
  $x_1\le 0$.
  In other words, we aim to minimize
  $f(\xx)$ subject to $g_i(\xx)\leq 0$ for $i=1,2$,
  where
  \begin{equation}
  \label{eq:ex:kkt-running-ex-standard:1}
      f(\xx) = x_2,
    \qquad
      g_1(\xx) = e^{x_1}+1-x_2,
    \qquad
      g_2(\xx) = x_1.
  \end{equation}
  (There are no equality constraints.)
  This convex program satisfies Slater's condition, as witnessed,
  for instance, by $\xing=\trans{[-1,2]}$.
  Also,
  $\inf p>-\infty$,
  because at any feasible point $\xx$
  we have $f(\xx)=x_2\ge e^{x_1}+1\ge 1$.
  However, this program has no finite solution.
  This is because, for any feasible point $\xx\in\R^2$,
  we can decrease $x_1$ to any value $x'_1<x_1\le 0$,
  and then decrease $x_2$ to $x'_2=e^{x'_1}+1<e^{x_1}+1\le x_2$,
  thus obtaining
  another feasible point $\xx'\in\R^2$ that is strictly
  better than $\xx$, meaning $f(\xx')<f(\xx)$.%
\indexg{ordinary convex programs!no finite solution|)}%
\end{example}

\section{Defining astral optimality}

The standard KKT theory only characterizes a convex program's finite
solutions.
However, as we saw in \Cref{ex:kkt-running-ex-standard},
the program's solutions might only exist at infinity.
We show now how that theory can be extended to astral
space, yielding an astral version of the KKT conditions that
characterizes the solutions of the original convex program, even if
those solutions are infinite.

We first need to define what it might mean for an astral point
$\xbar\in\extspace$
to solve the convex program $P$.
We saw that solving $P$ is equivalent to minimizing the
function $p$ given in \eqref{eqn:kkt-p-defn}.
\indexg{ordinary convex programs!astral solution definitions|(}%
\indexg{solution (of convex program)!astral|(}%
Therefore, we say that $\xbar$ is an \emph{(astral) solution}
of $P$ if it minimizes $p$'s extension, $\pext$, that is,
if $\pext(\xbar)=\inf p$.

Nonetheless, as was the case for constrained optimization problems
considered in \Cref{sec:constrained-optimization}, there are
other reasonable definitions of what it might mean for an astral point
to solve $P$.
As one alternative, we might instead consider more directly the astral
extensions of the objective function and individual constraints.
Specifically, we define an astral point $\xbar\in\extspace$ to be
\indexg{feasibility (of convex program)!astral (primal)}%
\emph{astral feasible} for $P$ if
$\gext_i(\xbar)\leq 0$ for $i=1,\ldots,r$,
and
$\hext_j(\xbar)=0$ for $j=1,\ldots,s$.
We then might define $\xbar$ to solve $P$ if
it minimizes $\fext(\xbar)$ over all astral feasible points
$\xbar$.

In yet another alternative, based on sequences,
we might say that $\xbar$ is a solution of~$P$
if there exists a sequence $\seq{\xx_t}$ of standard feasible
points (that is, in $C$) that converges to $\xbar$ and such that
the objective function values $f(\xx_t)$ converge to the program's
value, $\inf p$.

In fact, all three of these alternative definitions are equivalent, as
we state formally below in
\indexg{ordinary convex programs!astral solution definitions|)}%
\indexg{solution (of convex program)!astral|)}%
\Cref{pr:equiv-astral-kkt}.
First, we make some preliminary observations that will be needed for
that and other proofs.
We start with some consequences of Slater's
condition, and then give expressions for both the extension $\pext$
and the set of astral feasible points.

\begin{proposition}   \label{pr:slater-conseqs}
  Let $P$ be an ordinary convex program, as in
  \Cref{dfn:ord-cvx-opt-prog},
  and let $p$, $C$ and $\lagrangfcn$ be as in
  Eqs.~(\ref{eqn:kkt-p-defn}),~(\ref{eqn:feas-set-defn})
  and (\ref{eqn:lagrang-defn}).
  Assume Slater's condition holds, and let
  $\xing\in\Rn$ be a witness for it.
  Then:
  \begin{letter-compact}
  \item   \label{pr:slater-conseqs:b}
    $\xing\in\ri(\dom(\Negf\circ g_i))$
    for $i=1,\ldots,r$.
  \item   \label{pr:slater-conseqs:c}
    $\xing\in\ri(\dom(\indzohj))$
    for $j=1,\ldots,s$.
  \item   \label{pr:slater-conseqs:riC}
    $\xing\in\ri C$.
  \end{letter-compact}
\end{proposition}

\begin{proof}
~

\begin{proof-parts}
\pfpart{Part~(\ref{pr:slater-conseqs:b}):}
Let $i\in\{1,\ldots,r\}$.
By assumption, $\xing\in\ri(\dom{f})\subseteq\ri(\dom{g_i})$.
Also, $g_i(\xing)<0$ so
\[
  \xing
  \in
  \bigBraces{\xx\in\ri(\dom{g_i}) :\: g_i(\xx) < 0}
  =
  \ri {\bigBraces{\xx\in\Rn :\: g_i(\xx) \leq 0}}
  =
  \ri\bigParens{\dom(\Negf\circ g_i)},
\]
where
the first equality is by
\Cref{roc:thm7.6-mod}.

\pfpart{Part~(\ref{pr:slater-conseqs:c}):}
Let $j\in\{1,\ldots,s\}$, and
let
\[ M=\dom(\indzohj)=\{\xx\in\Rn :\: h_j(\xx)=0\}. \]
Since $h_j$ is affine, $M$ is an affine set, which is therefore
relatively open.
Since $h_j(\xing)=0$, it follows that
$\xing\in M = \ri M$.

\pfpart{Part~(\ref{pr:slater-conseqs:riC}):}
We have
\begin{align*}
  \xing
  &\in
  \BiggParens{\,
     \bigcap_{i=1}^r \ri\bigParens{\dom(\Negf\circ g_i)}
  }
  \cap
  \BiggParens{\,
     \bigcap_{j=1}^s \ri\bigParens{\dom(\indzohj)}
  }
  \\
  &=
  \ri\BiggBracks{
    \BiggParens{\,
    \bigcap_{i=1}^r \dom(\Negf\circ g_i)
    }
    \cap
    \BiggParens{\,
    \bigcap_{j=1}^s \dom(\indzohj)
    }
  }
  \\
  &=
  \ri C.
\end{align*}
The inclusion is by
parts~(\ref{pr:slater-conseqs:b}) and~(\ref{pr:slater-conseqs:c}).
The first equality is by
\Cref{roc:thm6.5}.
The last equality is by $C$'s definition.
\qedhere
\end{proof-parts}
\end{proof}

\begin{proposition}   \label{pr:kkt-p-C-exps}
  Let $P$ be an ordinary convex program, as in
  \Cref{dfn:ord-cvx-opt-prog},
  that satisfies Slater's condition.
  Let $p$ and $C$ be as in
  Eqs.~(\ref{eqn:kkt-p-defn})
  and~(\ref{eqn:feas-set-defn}),
  and let $\xbar\in\extspace$.
  Then:
  \begin{letter-compact}
  \item   \label{pr:kkt-p-C-exps:a}
\indexg{feasibility (of convex program)!astral (primal)|(}%
    $\xbar$ is astral feasible if and only if $\xbar\in\Cbar$.
  \item   \label{pr:kkt-p-C-exps:b}
\indexg{ordinary convex programs!extension of|(}%
    $p$'s extension is
    \begin{align*}
      \pext(\xbar)
      &=
      \BiggBracks{\,\sumitor \Negf(\gext_i(\xbar))
              +
              \sumjtos \indaz(\hext_j(\xbar))
             }
      \plusu
      \fext(\xbar)
      \\
      &=
      \indfa{\Cbar}(\xbar)\plusu \fext(\xbar)
      \\
      &=
      \begin{cases}
          \fext(\xbar)    & \text{if $\xbar$ is astral feasible,} \\
          +\infty         & \text{otherwise.}
      \end{cases}
    \end{align*}
  \item   \label{pr:kkt-p-C-exps:c}
\indexg{ordinary convex programs!subgradients of extension|(}%
    $\displaystyle
      \asubdifpext{\xbar}
      =
      \asubdiffext{\xbar}
      +
      \sumitor \asubdifNegogiext{\xbar}
      +
      \sumjtos \asubdifindzohjext{\xbar}
    $.
  \end{letter-compact}
\end{proposition}

\begin{proof}
The proofs are mainly based on applications of previously developed
tools.
Let $\xing$ be a witness for Slater's condition.

\begin{proof-parts}
\pfpart{Part~(\ref{pr:kkt-p-C-exps:a}):}
Let $i\in\{1,\ldots,r\}$.
Since $g_i(\xing)<0=\sup(\dom{\Negf})$,
we can apply \Cref{thm:Gf-conv},
yielding
\begin{equation}   \label{eq:pr:kkt-p-C-exps:1}
  \Negogiext(\xbar)
  =
  \negfext(\gext_i(\xbar))
  =
  \Negf(\gext_i(\xbar)).
\end{equation}

Next, we compute the extension of $\indz$ composed with any affine
function:

\begin{claimpx}   \label{cl:pr:kkt-p-C-exps:1}
  Let $h:\Rn\rightarrow\Rext$ be affine,
  and let $\xbar\in\extspace$.
  Assume $h(\xing)=0$.
  Then
  $\indzohext(\xbar)=\indaz(\hext(\xbar))$.
\end{claimpx}

\begin{proofx}
Since $h$ is affine, we can write it as
$h(\xx)=\xx\cdot\ww-\beta$ for some $\ww\in\Rn$ and $\beta\in\R$.
Then
$\indz(h(\xx)) = \indbeta(\xx\cdot\ww) = \indbeta(\trans{\ww}\xx)$
for $\xx\in\Rn$;
that is,
$\indzoh=\indbeta\trans{\ww}$.
Therefore,
\[
  \indzohext(\xbar)
  =
  \indbetaext(\trans{\ww} \xbar)
  =
  \indabetaAlt(\xbar\cdot\ww)
  =
  \indazAlt(\xbar\cdot\ww - \beta)
  =
  \indazAlt(\hext(\xbar)).
\]
The first equality is by \Cref{thm:ext-linear-comp}
(with $\A=\trans{\ww}$ and $f=\indbeta$),
which we can apply since $h(\xing)=0$, so
$\trans{\ww}{\xing}\in\{\beta\}=\ri(\dom{\indbeta})$.
The second equality is by
Propositions~\ref{pr:trans-uu-xbar}
and~\ref{pr:inds-ext}.
The last equality uses
Example~\ref{ex:ext-affine}.
\end{proofx}

By $C$'s definition, we can express its indicator function as
\[
  \indC(\xx)
  =
  \sumitor \Negf(g_i(\xx))
  +
  \sumjtos \indz(h_j(\xx)),
\]
for $\xx\in\Rn$.
We can then compute $\indC\negKern$'s extension to be
\begin{align}
  \indfa{\Cbar}(\xbar)
  =
  \indCext(\xbar)
  &=
  \sumitor \Negogiext(\xbar)
  +
  \sumjtos \indzohjext(\xbar)
  \nonumber
  \\
  &=
  \sumitor \Negf(\gext_i(\xbar))
  +
  \sumjtos \indaz(\hext_j(\xbar))
  \label{eq:pr:kkt-p-C-exps:3}
  \\
  &=
   \begin{cases}
        0         & \mbox{if $\xbar$ is astral feasible,} \\
      +\infty     & \mbox{otherwise.}
   \end{cases}
  \label{eq:pr:kkt-p-C-exps:2}
\end{align}
The first equality is by
\Cref{pr:inds-ext}.
The second is by
\Cref{thm:ext-sum-fcns-w-duality},
which we can apply here because the extensions
$\Negogiext(\xbar)$ and $\indzohjext(\xbar)$
are all nonnegative and therefore summable, and also because
the relative interiors of the domains of the functions
being added have a point in common (namely, $\xing$),
by
\Cref{pr:slater-conseqs}(\ref{pr:slater-conseqs:b},\ref{pr:slater-conseqs:c}).
The third equality above is by
\eqref{eq:pr:kkt-p-C-exps:1}
and
\Cref{cl:pr:kkt-p-C-exps:1}.
The last equality is by definition of astral feasibility.
This proves the claim.%
\indexg{feasibility (of convex program)!astral (primal)|)}%

\pfpart{Part~(\ref{pr:kkt-p-C-exps:b}):}
Since $p=\indC\plusu f$,
it follows from 
\Cref{cor:ext-restricted-fcn}
that $\pext(\xbar) = \indfa{\Cbar}(\xbar)\plusu \fext(\xbar)$,
noting that $\xing\in\ri(\dom f)\cap\ri C$
by \Cref{pr:slater-conseqs}(\ref{pr:slater-conseqs:riC}).
Combined with
Eqs.~(\ref{eq:pr:kkt-p-C-exps:3})
and~(\ref{eq:pr:kkt-p-C-exps:2}),
this gives all three stated expressions for $\pext(\xbar)$.%
\indexg{ordinary convex programs!extension of|)}%

\pfpart{Part~(\ref{pr:kkt-p-C-exps:c}):}
This follows directly from
\Cref{thm:subgrad-sum-fcns}
applied to $f$, $\Negf\circ g_i$ for $i=1,\ldots,r$,
and $\indzohj$ for $j=1,\ldots,s$,
noting that these functions have domains with overlapping relative
interiors by
\Cref{pr:slater-conseqs}(\ref{pr:slater-conseqs:b},\ref{pr:slater-conseqs:c}).%
\indexg{ordinary convex programs!subgradients of extension|)}%
\qedhere
\end{proof-parts}
\end{proof}

\indexg{ordinary convex programs!astral solution definitions|(}%
\indexg{solution (of convex program)!astral|(}%
We can now prove the equivalence of the astral solution concepts
discussed above, which is really just a restatement of 
\Cref{pr:cons-opt-equiv-soln} for the present setting:

\begin{proposition}   \label{pr:equiv-astral-kkt}
  Let $P$ be an ordinary convex program, as in
  \Cref{dfn:ord-cvx-opt-prog},
  that satisfies Slater's condition.
  Let $p$ be as in
  \eqref{eqn:kkt-p-defn},
  and let $\xbar\in\extspace$.
  Then the following are equivalent:
  \begin{letter-compact}
  \item   \label{pr:equiv-astral-kkt:a}
    $\xbar$ is an astral solution of $P$;
    that is, $\xbar$ minimizes $\pext$.
  \item   \label{pr:equiv-astral-kkt:b}
    $\xbar$ is astral feasible and
    minimizes $\fext$ over all astral feasible points.
  \item   \label{pr:equiv-astral-kkt:c}
    There exists a sequence
    $\seq{\xx_t}$ of standard feasible points in $\Rn$
    such that $\xx_t\rightarrow\xbar$
    and $f(\xx_t)\rightarrow \inf p$.
  \end{letter-compact}
\end{proposition}

\begin{proof}
Let $C$ be the set of all standard feasible points,
as in \eqref{eqn:feas-set-defn}.
Then by \Cref{pr:kkt-p-C-exps}(\ref{pr:kkt-p-C-exps:a}),
an astral point is astral feasible if and only if it is in $\Cbar$.
Also, any witness for Slater's condition
must be both in $\ri(\dom f)$ and $\ri C$ by
\Cref{pr:slater-conseqs}(\ref{pr:slater-conseqs:riC}),
so $\ri(\dom f)\cap\ri C\neq\emptyset$.
In light of these facts, and since $p=f\plusu \indC$,
the equivalences stated in this proposition
now follow directly from
\Cref{pr:cons-opt-equiv-soln}.%
\indexg{ordinary convex programs|)}%
\indexg{ordinary convex programs!astral solution definitions|)}%
\indexg{solution (of convex program)!astral|)}%
\end{proof}

\section{Astral optimality conditions}

\indexg{KKT conditions!astral|(}%
We define next a natural astral analogue of the standard
KKT conditions.
We will soon see that, assuming Slater's condition,
the astral KKT conditions are always sufficient to ensure that an
astral point 
is an astral solution of the convex program $P$,
and are necessary as well if the program's value is finite.
As with the standard KKT conditions, the astral analogue captures
astral versions of
primal feasibility, dual feasibility, complementary slackness and
first-order optimality.

\begin{definition}  \label{dfn:ast-kkt-defn}
  Let $P$ be an ordinary convex program, as in
  \Cref{dfn:ord-cvx-opt-prog},
  let $\xbar\in\extspace$ and let $\lammupair\in\Lampairs$.
  We say that the pair $\xbar,\lammupair$ satisfies the
  \emph{astral KKT conditions} if
  all of the following hold:
  \begin{item-compact}
  \item
\indexg{feasibility (of convex program)!astral (primal)}%
    \emph{primal feasibility:}
  $\xbar$ is astral feasible;
  \item
\indexg{complementary slackness!astral}%
    \emph{complementary slackness:}
  $\lambda_i \gext_i(\xbar)=0$ for $i=1,\ldots,r$;
  \item
\indexg{first-order optimality (of convex program)!astral|(}%
    \emph{first-order optimality:}
  \begin{equation}  \label{eq:dfn:ast-kkt-defn:1}
       \zero
       \in
       \basubdiffext{\xbar}
       +
       \sumitor \asubsym\lamgiext(\xbar)
       +
       \sumjtos \basubdifmuhjext{\xbar}.
  \end{equation}
  \end{item-compact}
\end{definition}

Dual feasibility is included as part
of the assumptions.
The other three conditions
are derived from the corresponding standard
conditions by replacing standard functions with their astral
extensions.
The first-order optimality condition
is based on the modified version in \eqref{eqn:kkt-lagrang-opt-mod}.%
\indexg{first-order optimality (of convex program)!astral|)}%
\indexg{KKT conditions!astral|)}%

\indexg{saddle points (of Lagrangian)!astral|(}%
\indexg{Lagrangian!astral saddle points of|(}%
\indexg{Lagrangian!extension of|(}%
We next define an astral saddle-point condition based on the extended
Lagrangian~$\lagrangextsym$.
Whereas, under Slater's condition, the KKT conditions are
sufficient to ensure that $\xbar$ is an astral solution of $P$,
the saddle-point condition turns out to be necessary.
Similar to the KKT conditions,
it becomes both necessary and sufficient
if the program's value is finite.\looseness=-1

\begin{definition}  \label{dfn:ast-saddle-defn}
  Let $P$ be an ordinary convex program, as in
  \Cref{dfn:ord-cvx-opt-prog},
  let $\xbar\in\extspace$ and $\lammupair\in\Lampairs$,
  and let $\lagrangfcn$ be as in \eqref{eqn:lagrang-defn}.
  We say that the pair $\xbar,\lammupair$ is a
  \emph{saddle point} of the extended Lagrangian
  $\lagrangextsym$ if
  \begin{equation}  \label{eq:dfn:ast-saddle-defn:1}
      \lagrangext{\xbar}{\lamvec',\muvec'}
      \leq
       \lagranglmext{\xbar}
       \leq
       \lagranglmext{\xbar'}
  \end{equation}
  for all $\lammupairp\in\Lampairs$
  and all
\indexg{saddle points (of Lagrangian)!astral|)}%
\indexg{Lagrangian!astral saddle points of|)}%
  $\xbar'\in\extspace$.
\end{definition}

\indexg{Lagrangian!subgradients of extension|(}%
As a first step in
deriving relationships between the conditions introduced
above and astral solutions of the convex program $P$,
we next give expressions for the extended Lagrangian and for its
astral subdifferential.
In particular, part~(\ref{pr:lagrang-ext-subdif:a2:new}) implies that the
first-order optimality condition
in
the astral KKT conditions (Eq.~\ref{eq:dfn:ast-kkt-defn:1})
is equivalent to
$\zero\in\basubdiflagranglmext{\xbar}$, and so to
$\xbar$ minimizing $\lagrangflmext$ with
$\inf\lagrangflmext\in\R$.
(The function $\sfprodlamgi$ used in
part~(\ref{pr:lagrang-ext-subdif:a2}) is as defined in
Eq.~\ref{eq:sfprod-defn}.)

\begin{proposition}   \label{pr:lagrang-ext-subdif}
  Let $P$ be an ordinary convex program, as in
  \Cref{dfn:ord-cvx-opt-prog},
  that satisfies Slater's condition.
  Let $\lagrangfcn$ be as in
  \eqref{eqn:lagrang-defn}.
  Let $\xbar\in\extspace$ and $\lammupair\in\Lampairs$.
  Then:
  \begin{letter}
    \item   \label{pr:lagrang-ext-subdif:a2:new}
         $\displaystyle
           \basubdiflagranglmext{\xbar}
           =
           \basubdiffext{\xbar}
           +
           \sumitor \asubsym\lamgiext(\xbar)
           +
           \sumjtos \basubdifmuhjext{\xbar}
         $.
    \item   \label{pr:lagrang-ext-subdif:a2}
     $\displaystyle
       \basubdiflagranglmext{\xbar}
       =
       \basubdiffext{\xbar}
       +
       \sumitor \basubdiflamgiext{\xbar}
       +
       \sumjtos \basubdifmuhjext{\xbar}
     $.
  \item   \label{pr:lagrang-ext-subdif:b}
    Suppose
    $\fext(\xbar),
     \lambda_1 \gext_1(\xbar),\dotsc, \lambda_r \gext_r(\xbar),
     \mu_1 \hext_1(\xbar),\dotsc, \mu_s \hext_s(\xbar)$
    are summable.
    Then
    \[
    \lagranglmext{\xbar}
    =
    \fext(\xbar)
    +
    \sumitor \lambda_i \gext_i(\xbar)
    +
    \sumjtos \mu_j \hext_j(\xbar).
    \]
  \end{letter}
\end{proposition}

\begin{proof}
Let $\xing$ be
a witness for Slater's condition.
Then
\begin{equation}   \label{eq:pr:lagrang-ext-subdif:3}
  \xing
  \in
  \ri(\dom{f})
  \subseteq
  \ri(\dom{g_i})
  =
  \ri(\dom(\sfprodlamgi))
  \subseteq
  \ri(\dom(\lambda_i g_i))
\end{equation}
for $i=1,\ldots,r$,
where
the equality is by \eqref{eq:sfprod-defn},
and the last inclusion is because if
$\lambda_i>0$ then $\sfprodlamgi=\lambda_i g_i$,
and if $\lambda_i=0$ then $\dom(\lambda_i g_i)=\Rn$.

\begin{proof-parts}
\pfpart{%
  Part~(\ref{pr:lagrang-ext-subdif:a2:new}):
}
This follows from the definition of the Lagrangian
(Eq.~\ref{eqn:lagrang-defn})
by \Cref{thm:subgrad-sum-fcns}, whose assumptions
are satisfied by \eqref{eq:pr:lagrang-ext-subdif:3}
and because
$\dom(\mu_j h_j)=\Rn$ for all $j=1,\ldots,s$.

\pfpart{%
  Part~(\ref{pr:lagrang-ext-subdif:a2}):
}
We can slightly rewrite the Lagrangian from
\eqref{eqn:lagrang-defn}
as
\[
  \lagranglm{\xx}
  =
  f(\xx)
  +
  \sumitor (\sfprodlamgi)(\xx)
  +
  \sumjtos (\mu_j h_j)(\xx)
\]
since
$(\sfprodlamgi)(\xx)= \lambda_i g_i(\xx)$
for all $\xx\in\dom{g_i}$, and therefore for all
$\xx\in\dom{f}$.
From this, the claimed expression follows
from \Cref{thm:subgrad-sum-fcns} by the same reasoning
as in part~(\ref{pr:lagrang-ext-subdif:a2:new}).%
\indexg{Lagrangian!subgradients of extension|)}%

\pfpart{Part~(\ref{pr:lagrang-ext-subdif:b}):}
Suppose the values
    $\fext(\xbar),
     \lambda_1 \gext_1(\xbar),\dotsc, \lambda_r \gext_r(\xbar),
     \mu_1 \hext_1(\xbar),\dotsc, \mu_s \hext_s(\xbar) $
are summable.
By \Cref{pr:scal-mult-ext},
\begin{equation}   \label{eq:pr:lagrang-ext-subdif:1}
  \lambda_i \gext_i(\xbar)
  =
  \lamgiext(\xbar).
\end{equation}
Also, for $j=1,\ldots,s$, since $h_j$ is affine, we can write
$h_j(\xx)=\xx\cdot\ww_j-\beta_j$ for some $\ww_j\in\Rn$ and $\beta_j\in\R$.
Then
\begin{equation}   \label{eq:pr:lagrang-ext-subdif:2}
  \mu_j \hext_j(\xbar)
  =
  \mu_j (\xbar\cdot\ww_j - \beta_j)
  =
  \xbar\cdot (\mu_j \ww_j) - \mu_j \beta_j
  =
  \muhjext(\xbar),
\end{equation}
where the first and third equalities are by
Example~\ref{ex:ext-affine}.

Thus,
$\fext(\xbar),
 \lamgextgen{1}(\xbar),\dotsc, \lamgextgen{r}(\xbar),
 \muhextgen{1}(\xbar),\dotsc, \muhextgen{s}(\xbar) $
are also summable.
By \eqref{eq:pr:lagrang-ext-subdif:3} and since, as already noted,
$\dom(\mu_j h_j)=\Rn$ for all $j$,
we can apply
\Cref{thm:ext-sum-fcns-w-duality}
yielding
\begin{align*}
  \lagranglmext{\xbar}
  &=
  \fext(\xbar)
  +
  \sumitor \lamgiext(\xbar)
  +
  \sumjtos \muhjext(\xbar)
  \\
  &=
  \fext(\xbar)
  +
  \sumitor \lambda_i \gext_i(\xbar)
  +
  \sumjtos \mu_j \hext_j(\xbar),
\end{align*}
where the second equality follows from
Eqs.~(\ref{eq:pr:lagrang-ext-subdif:1})
and~(\ref{eq:pr:lagrang-ext-subdif:2}).%
\indexg{Lagrangian!extension of|)}%
\qedhere
\end{proof-parts}
\end{proof}

\indexg{KKT vectors!existence of|(}%
\indexg{ordinary convex programs!existence of KKT vector|(}%
As discussed in \Cref{sec:kkt-stand-set},
under Slater's condition, a KKT vector must always exist.
We prove this in the next theorem.
Notice that
this theorem concerns only notions from standard convex analysis.
Indeed, for the case $\inf p>-\infty$,
it is taken directly from \idxroc\citet[Theorem~28.2]{ROC} where it is
proved using standard (non-astral) methods.
Here, we give a different proof based centrally on astral notions and
techniques as an illustration of how these can be used even for
proving standard results.

\begin{theorem}  \label{roc:thm28.2}
  Let $P$ be an ordinary convex program, as in
  \Cref{dfn:ord-cvx-opt-prog},
  that satisfies Slater's condition,
  and let $p$ be as in
  \eqref{eqn:kkt-p-defn}.
  Then there exists a KKT vector in $\Lampairs$ for $P$.
\end{theorem}

The proof is based largely on relating the astral subdifferentials
appearing in the expression for $\asubdiflagranglmext{\xbar}$
in \Cref{pr:lagrang-ext-subdif}(\ref{pr:lagrang-ext-subdif:a2})
to the astral subdifferentials
appearing in the expression for $\asubdifpext{\xbar}$ in
\Cref{pr:kkt-p-C-exps}(\ref{pr:kkt-p-C-exps:c}).
As such, we first prove the following two lemmas, which will be used
both in proving \Cref{roc:thm28.2} and
again later in proving \Cref{thm:lammu-equivs}.

\begin{lemma}  \label{lem:subdif-neg-o-g}
  Let $g:\Rn\rightarrow\Rext$ be convex and proper,
  and let $\xbar\in\extspace$.
  Assume $\gext(\xbar)\leq 0$, and that there exists
  $\xing\in\Rn$ with $g(\xing)<0$.
  Then
  \[
    \basubdifNegogext{\xbar}
    \;=\;
    \bigcup\;\BigBraces{%
      \basubdiflamgext{\xbar}
      :\:
      \lambda\in\Rpos,\,
      \lambda\gext(\xbar)=0
    }.
  \]
\end{lemma}

\begin{proof}
By a straightforward calculation
(for instance, using
\Cref{cor:subdif-ind-fcn-rn}\ref{cor:subdif-ind-fcn-rn:a}\ref{cor:subdif-ind-fcn-rn:b}),
\[
  \basubdifnegext{\bary}
  =
  \begin{cases}
        \{0\}       & \text{if $\bary < 0$,} \\
        [0,+\infty) & \text{if $\bary = 0$,} \\
        \emptyset   & \text{if $\bary > 0$,}
  \end{cases}
\]
for $\bary\in\Rext$.
Thus, since $\gext(\xbar)\leq 0$, for $\lambda\in\R$,
it follows that
$\lambda\in\basubdifnegext{\gext(\xbar)}$
if and only if
$\lambda\geq 0$
and
$\lambda\gext(\xbar)=0$.
With this observation, the claim now follows directly
from
\Crefequiv{thm:subgrad-comp-inc-fcn-equiv}{thm:subgrad-comp-inc-fcn-equiv:a}{thm:subgrad-comp-inc-fcn-equiv:b}
(applied with $f$ and $g$, as they appear in that
\namecref{thm:subgrad-comp-inc-fcn-equiv},
set to $g$ and $\negf$, and noting that
$g(\xing)<0=\sup(\dom{\negf})$).
\end{proof}

\begin{lemma}  \label{lem:subdif-indz-o-h}
  Let $h:\Rn\rightarrow\R$ be affine,
  and let $\xbar\in\extspace$.
  Assume $\hext(\xbar) = 0$.
  Then
  \[
    \basubdifindzohext{\xbar}
    \;=\;
    \bigcup\;\BigBraces{%
      \basubdifmuhext{\xbar}
      :\:
      \mu\in\R
    }.
  \]
\end{lemma}

\begin{proof}
Since $h$ is affine, it can be written in the form
$h(\xx)=\xx\cdot\ww-\beta$ for $\xx\in\Rn$,
for some $\ww\in\Rn$ and $\beta\in\R$.
This implies that $\hext(\xbar)=\xbar\cdot\ww-\beta$
(Example~\ref{ex:ext-affine}), and so that
$\xbar\cdot\ww=\beta$.
Let
$M=\{\xx\in\Rn :\: \xx\cdot\ww=\beta\}$.
Then $\indzoh$ is exactly the indicator function for $M$;
that is, $\indzoh=\indf{M}$.
Also, $\xbar\in\Mbar$
by \Cref{cor:closure-affine-set}
(applied with $\A=\trans{\ww}$ and $\bb=\beta$), implying
$M\neq\emptyset$ and so that there exists $\xing\in M$.

\Cref{cor:subdif-ind-fcn-rn}(\ref{cor:subdif-ind-fcn-rn:b},\ref{cor:subdif-ind-fcn-rn:a})
(applied with $\xbar$ and $S$, as they appear in that corollary, set
to $0$ and $\{0\}$)
implies $\asubdifindzext{0}=\R$ since $\indzstar\equiv 0$.
Also,
$\xing\cdot\ww-\beta\in\{0\}=\ri(\dom \indz)$.
Therefore, 
\begin{equation}  \label{eq:lem:subdif-indz-o-h:1}
  \basubdifindzohextTight{\xbar}
  =
  \ww\basubdifindzext{\trans{\ww\negKern}\xbar-\beta}
  =
  \ww\basubdifindzext{\hext(\xbar)}
  =
  \ww\basubdifindzext{0}
  =
  \R\ww,
\end{equation}
where the first equality is by
\Cref{cor:subgrad-comp-w-affine}
(applied with $\A=\trans{\ww}$, $\bb=-\beta$, and $f=\indz$).

On the other hand, by
Example~\ref{ex:affine-subgrad-new},
for $\mu\in\R$,
$\basubdifmuhext{\xbar} = \{\mu \ww\}$,
since
$(\mu h)(\xx)=\xx\cdot(\mu\ww) - \mu \beta$.
Taking the union over $\mu\in\R$ then yields exactly
$\basubdifindzohextTight{\xbar}$ by \eqref{eq:lem:subdif-indz-o-h:1},
completing the proof.
\end{proof}

\begin{proof}[Proof of \Cref{roc:thm28.2}:]
By \Cref{pr:lagran-facts}(\ref{pr:lagran-facts:a}),
$\inf \lagrangflm\leq \inf p$
for all $\lammupair\in\Lampairs$.
Consequently, if $\inf p=-\infty$, then every vector
$\lammupair$ is a KKT vector, proving the claim in this case.
We therefore assume henceforth that $\inf p\in\R$
(with the case $\inf p=+\infty$ ruled out by
\Cref{pr:lagran-facts}\ref{pr:lagran-facts:c}).

By \Cref{pr:fext-min-exists}, there exists a point $\xbar$ that
minimizes $\pext$ so that $\pext(\xbar)=\min\pext=\inf p$,
which is finite.
Hence, 
$\zero\in\asubdifpext{\xbar}$
by \Cref{pr:asub-zero-is-min}(\ref{pr:asub-zero-is-min:a}).
Moreover, by \Cref{pr:kkt-p-C-exps}(\ref{pr:kkt-p-C-exps:b}),
$\xbar$ is astral feasible and
$\pext(\xbar)=\fext(\xbar)$.

We thus have that
\[
  \zero
  \in
  \basubdifpext{\xbar}
  =
  \basubdiffext{\xbar}
  +
  \sumitor \basubdifNegogiext{\xbar}
  +
  \sumjtos \basubdifindzohjext{\xbar},
\]
with the equality from
\Cref{pr:kkt-p-C-exps}(\ref{pr:kkt-p-C-exps:c}).
Consequently, there exists
$\uu\in\basubdiffext{\xbar}$,
$\vv_i\in\basubdifNegogiext{\xbar}$ for $i=1,\ldots,r$,
and
$\ww_j\in\basubdifindzohjext{\xbar}$ for $j=1,\ldots,s$,
such that
$\uu+\sumitor \vv_i + \sumjtos \ww_j = \zero$.

For each $i$,
since
$\vv_i\in\basubdifNegogiext{\xbar}$,
it follows from
\Cref{lem:subdif-neg-o-g}
that there exists $\lambda_i\in\Rpos$ with
$\lambda_i \gext_i(\xbar)=0$
such that
$\vv_i\in\basubdiflamgiext{\xbar}$.
Likewise, for each~$j$, since
$\ww_j\in\basubdifindzohjext{\xbar}$,
from
\Cref{lem:subdif-indz-o-h},
there exists $\mu_j\in\R$ such that
$\ww_j\in \basubdifmuhjext{\xbar}$.
Thus,
\begin{align*}
  \zero
  &=
  \uu + \sumitor \vv_i + \sumjtos \ww_j
  \\
  &\in
  \basubdiffext{\xbar}
  +
  \sumitor \basubdiflamgiext{\xbar}
  +
  \sumjtos \basubdifmuhjext{\xbar}
  =
  \asubdiflagranglmext{\xbar},
\end{align*}
where the last equality is by
\Cref{pr:lagrang-ext-subdif}(\ref{pr:lagrang-ext-subdif:a2}).
By \Cref{pr:asub-zero-is-min}(\ref{pr:asub-zero-is-min:a}),
this implies that $\xbar$ minimizes $\lagrangflmext$,
so $\lagranglmext{\xbar}=\inf\lagrangflm$
(by  \Cref{pr:fext-min-exists}).
Combining facts then yields that
\[
  \inf\lagrangflm
  =
  \lagranglmext{\xbar}
  =
  \fext(\xbar)
  =
  \pext(\xbar)
  =
  \inf p,
\]
where the second equality is by
\Cref{pr:lagrang-ext-subdif}(\ref{pr:lagrang-ext-subdif:b})
(since $\lambda_i\gext_i(\xbar)=0$ for all $i$,
and since $\xbar$ is astral feasible, implying
$\hext_j(\xbar)=0$ for all $j$).
This completes the proof.%
\indexg{KKT vectors!existence of|)}%
\indexg{ordinary convex programs!existence of KKT vector|)}%
\end{proof}

\indexg{constrained optimality conditions, astral!convex program@for convex program|(}%
\indexg{ordinary convex programs!astral optimality conditions|(}%
We next examine, as discussed above,
how the astral conditions we have been considering
relate to one another and to a point being a solution of program $P$.
Under Slater's condition, the next theorem shows that, for a pair
$\xbar,\lammupair$, the astral KKT conditions from
\Cref{dfn:ast-kkt-defn} are sufficient for
$\xbar$ to be an astral solution and for $\lammupair$ to be a KKT
vector, while the saddle-point condition from
\Cref{dfn:ast-saddle-defn} is necessary for these to hold.
We will then see in
\Cref{thm:lammu-equivs-fin-val} that these three conditions are in
fact equivalent if the program's value is finite.

\begin{theorem}   \label{thm:lammu-equivs}
  Let $P$ be an ordinary convex program, as in
  \Cref{dfn:ord-cvx-opt-prog},
  that satisfies Slater's condition.
  Let $p$ and $\lagrangfcn$ be as in
  Eqs.~(\ref{eqn:kkt-p-defn}) and (\ref{eqn:lagrang-defn}).
  Let $\xbar\in\extspace$ and let $\lammupair\in\Lampairs$.
  Consider the following statements:
  \begin{letter-compact}
  \item   \label{thm:lammu-equivs:kkt}
    The pair $\xbar,\lammupair$ satisfies the astral KKT
    conditions.
  \item   \label{thm:lammu-equivs:sln}
    $\xbar$ is an astral solution,
    and $\lammupair$ is a KKT vector.
  \item   \label{thm:lammu-equivs:sdl}
    The pair $\xbar,\lammupair$ is a
    saddle point of the extended Lagrangian $\lagrangextsym$.
  \end{letter-compact}
  Then
  (\ref{thm:lammu-equivs:kkt})~$\Rightarrow$~(\ref{thm:lammu-equivs:sln})
  and
  (\ref{thm:lammu-equivs:sln})~$\Rightarrow$~(\ref{thm:lammu-equivs:sdl}).
  Furthermore, if the pair $\xbar,\lammupair$ satisfies any one of these
  conditions, then
  $\lagranglmext{\xbar}=\inf p$.
\end{theorem}

\begin{proof}
  ~

\begin{proof-parts}
\pfpart{%
  (\ref{thm:lammu-equivs:kkt})
  $\Rightarrow$
  (\ref{thm:lammu-equivs:sln}):
}
{\mathtogether%
Suppose the pair $\xbar,\lammupair$ satisfies the astral KKT
conditions from \Cref{dfn:ast-kkt-defn};
we aim to prove $\pext(\xbar)=\inf p = \inf\lagrangflm$.
The first-order optimality condition
combined with \Cref{pr:lagrang-ext-subdif}(\ref{pr:lagrang-ext-subdif:a2:new})
implies that $\zero\in\asubdiflagranglmext{\xbar}$.
Consequently, by
\Cref{pr:asub-zero-is-min}(\ref{pr:asub-zero-is-min:a}),
$\xbar$~minimizes $\lagrangflmext$, 
so
$\lagranglmext{\xbar}=\inf {\lagrangflm}$
(by \Cref{pr:fext-min-exists}).}

Primal feasibility and complementary slackness
imply that
$\gext_i(\xbar)\leq 0$ and
$\lambda_i \gext_i(\xbar)=0$ for $i=1,\ldots,r$,
and that
$\hext_j(\xbar)=0$ for $j=1,\ldots,s$.
Thus, by
\Cref{pr:lagrang-ext-subdif}(\ref{pr:lagrang-ext-subdif:b}),
$ \lagranglmext{\xbar} = \fext(\xbar) $.
These conditions further imply
by \Cref{lem:subdif-neg-o-g} that
\begin{equation}   \label{eq:thm:lammu-equivs:6}
  \basubdiflamgiext{\xbar}\subseteq \basubdifNegogiext{\xbar}
\end{equation}
for $i=1,\ldots,r$,
noting that $g_i(\xing)<0$
and $\lambda_i\in\Rpos$,
where $\xing$ is a witness for Slater's condition.
Likewise, by \Cref{lem:subdif-indz-o-h},
\begin{equation}   \label{eq:thm:lammu-equivs:7}
  \basubdifmuhjext{\xbar}\subseteq \basubdifindzohjext{\xbar}
\end{equation}
for $j=1,\ldots,s$,
since $h_j$ is affine and $\mu_j\in\R$.

Thus,
\begin{align*}
  \zero
  \in
  \asubdiflagranglmext{\xbar}
  &=
  \basubdiffext{\xbar}
  +
  \sumitor \basubdiflamgiext{\xbar}
  +
  \sumjtos \basubdifmuhjext{\xbar}
  \\
  &\subseteq
  \basubdiffext{\xbar}
  +
  \sumitor \basubdifNegogiext{\xbar}
  +
  \sumjtos \basubdifindzohjext{\xbar}
  \nonumber
  =
  \basubdifpext{\xbar}.
\end{align*}
The first inclusion was observed above,
and the first equality is by
\Cref{pr:lagrang-ext-subdif}(\ref{pr:lagrang-ext-subdif:a2}).
The second inclusion follows from
Eqs.~(\ref{eq:thm:lammu-equivs:6})
and~(\ref{eq:thm:lammu-equivs:7}).
The second equality is by
\Cref{pr:kkt-p-C-exps}(\ref{pr:kkt-p-C-exps:c}).

Since $\zero\in\asubdifpext{\xbar}$, it follows from
\Cref{pr:asub-zero-is-min}(\ref{pr:asub-zero-is-min:a})
(and \Cref{pr:fext-min-exists})
that $\pext(\xbar)=\inf p$.
Combining, we have
\[
  \inf\lagrangflm
  =
  \lagranglmext{\xbar}
  =
  \fext(\xbar)
  =
  \pext(\xbar)
  =
  \inf p,
\]
where the third equality is by
\Cref{pr:kkt-p-C-exps}(\ref{pr:kkt-p-C-exps:b})
since $\xbar$ is astral feasible.
This shows both that $\xbar$ is an astral solution and that
$\lammupair$ is a KKT vector, completing the proof.

\pfpart{%
  (\ref{thm:lammu-equivs:sln})
  $\Rightarrow$
  (\ref{thm:lammu-equivs:sdl}):
}
Suppose statement~(\ref{thm:lammu-equivs:sln}) holds;
we aim to show $\xbar,\lammupair$ is a saddle point.
We have
\[
   \inf p
   =
   \inf\lagrangflm
   =
   \min\lagrangflmext
   \leq
   \lagranglmext{\xbar}
   \leq
   \pext(\xbar)
   =
   \inf p.
\]
The first equality is because $\lammupair$ is a KKT vector.
The second equality is by \Cref{pr:fext-min-exists}.
The second inequality is by 
Propositions~\ref{pr:lagran-facts}(\ref{pr:lagran-facts:a})
and~\ref{pr:h:1}(\ref{pr:h:1:geq}).
The last equality is because $\xbar$ is an astral solution.
Thus,
\begin{equation}   \label{eq:thm:lammu-equivs:1}
   \min\lagrangflmext
   =
   \lagranglmext{\xbar}
   =
   \pext(\xbar)
   =
   \inf p.
\end{equation}
This proves the second inequality of
\eqref{eq:dfn:ast-saddle-defn:1}.

Let $\lammupairp\in\Lampairs$.
Then
\[
  \lagranglmpext{\xbar}
  \leq
  \pext(\xbar)
  =
  \lagranglmext{\xbar},
\]
where the inequality is again by
Propositions~\ref{pr:lagran-facts}(\ref{pr:lagran-facts:a})
and~\ref{pr:h:1}(\ref{pr:h:1:geq}),
and the equality is by
\eqref{eq:thm:lammu-equivs:1}.
Thus, the first inequality of
\eqref{eq:dfn:ast-saddle-defn:1}
holds as well, completing the proof.

\pfpart{%
  Additional claim if any of
  (\ref{thm:lammu-equivs:kkt}),~(\ref{thm:lammu-equivs:sln}),
  or~(\ref{thm:lammu-equivs:sdl}) holds:
}
Suppose that any one of the statements
(\ref{thm:lammu-equivs:kkt}),
(\ref{thm:lammu-equivs:sln}),
or~(\ref{thm:lammu-equivs:sdl})
holds, implying, by the foregoing, that in every case,
statement~(\ref{thm:lammu-equivs:sdl}) must hold,
so $\xbar,\lammupair$ is a saddle point of $\lagrangextsym$.
We aim to show
$\lagranglmext{\xbar}=\inf p$.

By \Cref{roc:thm28.2},
there exists a KKT vector
$\lammupairz\in\Lampairs$.
We then have
\begin{align*}
  \inf p
  =
  \inf \lagranglmz{\posSkip\cdot}
  &\leq
  \lagranglmzext{\xbar}
  \\
  &\leq
  \lagranglmext{\xbar}
  =
  \min \lagrangflmext
  =
  \inf \lagrangflm
  \leq
  \inf p.
\end{align*}
The first equality is because $\lammupairz$ is a KKT vector.
The first inequality and third equality are both
by \Cref{pr:fext-min-exists}.
The second inequality and second equality are both because
$\xbar,\lammupair$ is a saddle point.
The last inequality is by
\Cref{pr:lagran-facts}(\ref{pr:lagran-facts:a}).
This proves the claim.
\qedhere
\end{proof-parts}
\end{proof}

With the additional assumption that the program's
value $\inf p$ is finite, 
the next theorem shows that the three conditions of
\Cref{thm:lammu-equivs} are equivalent, thus providing an astral
analogue of \Cref{roc:thm28.3}.
\indexg{ordinary convex programs!existence of astral solution|(}%
Importantly, whereas the program $P$ might have no standard solution,
it always has an astral solution.
That
solution, together with
any KKT vector (which exists by \Cref{roc:thm28.2}), then satisfies
all three conditions of the theorem.

\begin{theorem}   \label{thm:lammu-equivs-fin-val}
  Let $P$, $p$, $\lagrangfcn$, $\xbar$, and $\lammupair$ be as in
  \Cref{thm:lammu-equivs}, and additionally assume that
  $\inf p > -\infty$.
  Then that theorem's
  statements~(\ref{thm:lammu-equivs:kkt}),~(\ref{thm:lammu-equivs:sln}),
  and~(\ref{thm:lammu-equivs:sdl})
  are equivalent.

  Also, there exists a pair
  $\xbar,\lammupair \in\extspace\times\Lampairs$
  that satisfies all three
  conditions.
\end{theorem}

\begin{proof}
  ~

\begin{proof-parts}
\pfpart{%
  (\ref{thm:lammu-equivs:sdl})
  $\Rightarrow$
  (\ref{thm:lammu-equivs:kkt}):
}
Suppose that
$\xbar,\lammupair$ is a saddle point of $\lagrangextsym$;
we aim to show that all three conditions of
\Cref{dfn:ast-kkt-defn} hold.
We begin by showing that $\lagranglmext{\xbar}$ and $\fext(\xbar)$ are
finite and equal.

\begin{claimpx}   \label{cl:thm:lammu-equivs:o1}
  $\lagranglmext{\xbar}=\fext(\xbar)\in\R$.
\end{claimpx}

\begin{proofx}
Since $\xbar,\lammupair$ is a saddle point,
$\lagranglmext{\xbar}=\inf p$ by \Cref{thm:lammu-equivs}.
Further, $\inf p\in\R$
(by assumption and
\Cref{pr:lagran-facts}\ref{pr:lagran-facts:c}).

Observe next that
\begin{equation}  \label{eq:thm:lammu-equivs:o3}
  \fext(\xbar)
  =
  \lagrangext{\xbar}{\zerov{r},\zerov{s}}
  \leq
  \lagranglmext{\xbar}
  <
  +\infty.
\end{equation}
The equality is because $\lagrang{\posSkip\cdot}{\zerov{r},\zerov{s}}=f$,
and the first inequality is because
$\xbar,\lammupair$ is a saddle point.
It remains to show that $\fext(\xbar)\geq\lagranglmext{\xbar}$.

Let $d:\Rn\rightarrow\Rext$ be defined by
$d(\xx)=2\lagranglm{\xx}$
for $\xx\in\Rn$,
which can be rewritten as
\[ d(\xx)=f(\xx)+\lagranglmd{\xx}. \]
By Slater's condition and
\Cref{pr:lagran-facts}(\ref{pr:lagran-facts:b}), there exists a point
in $\ri(\dom f)=\ri(\dom{\lagranglmd{\posSkip\cdot}})$.
Also,
\[ \lagranglmdext{\xbar}\leq\lagranglmext{\xbar}<+\infty, \]
since $\xbar,\lammupair$ is a saddle point,
so $\fext(\xbar)$ and $\lagranglmdext{\xbar}$ are summable
(using Eq.~\ref{eq:thm:lammu-equivs:o3}).
Therefore,
\[
  2 \lagranglmext{\xbar}
  =
  \dext(\xbar)
  =
  \fext(\xbar)
  +
  \lagranglmdext{\xbar}
  \leq
  \fext(\xbar)
  +
  \lagranglmext{\xbar}.
\]
The first equality is by
\Cref{pr:scal-mult-ext},
and the second by
\Cref{thm:ext-sum-fcns-w-duality}.
It follows that $\lagranglmext{\xbar}\leq\fext(\xbar)$
(since $\lagranglmext{\xbar}\in\R$).
Combined with \eqref{eq:thm:lammu-equivs:o3},
this completes the proof.
\end{proofx}

Since $\xbar$ minimizes $\lagrangflmext$
with a minimum value $\lagranglmext{\xbar}$ that is finite,
\Cref{pr:asub-zero-is-min}(\ref{pr:asub-zero-is-min:a})
yields $\zero\in\asubdiflagranglmext{\xbar}$.
By \Cref{pr:lagrang-ext-subdif}(\ref{pr:lagrang-ext-subdif:a2:new}),
this further implies that the first-order optimality condition
(Eq.~\ref{eq:dfn:ast-kkt-defn:1}) holds.

For each $i\in\{1,\ldots,r\}$,
we show next that $\gext_i(\xbar)\leq 0$.
Let $\ee_i$ be a standard basis vector in $\R^r$, implying
$\lagpair{\ee_i}{\zerov{s}}\in\Lampairs$.
Then
\begin{equation*}  %
  \fext(\xbar)
  =
  \lagranglmext{\xbar}
  \geq
  \lagrangext{\xbar}{\ee_i,\zerov{s}}
  =
  \fext(\xbar)
  +
  \gext_i(\xbar).
\end{equation*}
The first equality is by
\Cref{cl:thm:lammu-equivs:o1}.
The inequality is because
$\xbar,\lammupair$ is a saddle point.
The second equality is by
\Cref{pr:lagrang-ext-subdif}(\ref{pr:lagrang-ext-subdif:b}),
whose summability assumption is satisfied since $\fext(\xbar)\in\R$.
As claimed, this implies that $\gext_i(\xbar)\leq 0$.

Likewise, for each $j\in\{1,\ldots,s\}$,
we argue that $\hext_j(\xbar)=0$.
Let $\ee_j$ now be a standard basis vector in $\R^s$,
implying $\lagpair{\zerov{r}}{\pm\ee_j}\in\Lampairs$.
Then by similar reasoning as above,
\begin{align*}
  \fext(\xbar)
  &=
  \lagranglmext{\xbar}
  \geq
  \lagrangext{\xbar}{\zerov{r},+\ee_j}
  =
  \fext(\xbar)
  +
  \hext_j(\xbar),
\\  
\fext(\xbar)
&=
\lagranglmext{\xbar}
\geq
\lagrangext{\xbar}{\zerov{r},-\ee_j}
=
\fext(\xbar)
-
\hext_j(\xbar).
\end{align*}
Thus, $\ef(\xbar)\ge\ef(\xbar)+\regAbs{\hext_j(\xbar)}$,
which, since $\fext(\xbar)\in\R$,
is only possible if $\hext_j(\xbar)=0$.
Hence, primal feasibility holds.

Furthermore,
\begin{align}
\SwapAboveDisplaySkip
\label{eq:thm:lammu-equivs:8}
  \fext(\xbar)
  &=
  \lagranglmext{\xbar}
  =
  \fext(\xbar)
  +
  \sumitor \lambda_i \gext_i(\xbar).
\end{align}
The first equality is by
\Cref{cl:thm:lammu-equivs:o1}.
The second is because $\hext_j(\xbar)=0$ for all $j$,
and by
\Cref{pr:lagrang-ext-subdif}(\ref{pr:lagrang-ext-subdif:b}),
whose summability assumption is satisfied because
$\fext(\xbar)\in\R$
and because $\lambda_i \gext_i(\xbar)\leq 0$ for all $i$
since $\gext_i(\xbar)\leq 0$ and $\lambda_i\in\Rpos$.
Combined with these last facts,
\eqref{eq:thm:lammu-equivs:8} then implies that actually
$\lambda_i \gext_i(\xbar) = 0$ for $i=1,\ldots,r$.
Hence, complementary slackness also holds,
completing the proof.

\pfpart{Existence:}
By \Cref{pr:fext-min-exists}, there exists $\xbar\in\extspace$ that
minimizes $\pext$, and so is an astral solution.
By \Cref{roc:thm28.2}, there exists a KKT vector $\lammupair\in\Lampairs$.
Together, the pair $\xbar,\lammupair$ satisfies
statement~(\ref{thm:lammu-equivs:sln})
of \Cref{thm:lammu-equivs}, and therefore the other
statements as well.
\qedhere
\end{proof-parts}
\end{proof}

As immediate corollary, we obtain the following optimality conditions
for an astral point to be an astral solution of a convex program,
provided the program's value is finite.
These are exactly analogous to the standard conditions given in
\Cref{roc:cor28.3.1} with the important difference, noted earlier,
that whereas a convex program might not have any standard solutions,
it always must have an astral solution:

\begin{corollary}  \label{cor:ast-soln-kkt-saddle}
  Let $P$ be an ordinary convex program, as in
  \Cref{dfn:ord-cvx-opt-prog},
  that satisfies Slater's condition,
  and let $p$ and $\lagrangfcn$ be as in
  Eqs.~(\ref{eqn:kkt-p-defn}) and (\ref{eqn:lagrang-defn}).
  Assume $\inf p>-\infty$.
  Let $\xbar\in\extspace$.
  Then the following are equivalent:
  \begin{letter-compact}
  \item  \label{cor:ast-soln-kkt-saddle:kkt}
    There exists $\lammupair\in\Lampairs$
    such that the pair $\xbar,\lammupair$ satisfies the astral KKT
    conditions.
  \item  \label{cor:ast-soln-kkt-saddle:soln}
    $\xbar$ is an astral solution of $P$.
  \item  \label{cor:ast-soln-kkt-saddle:saddle}
    There exists $\lammupair\in\Lampairs$
    such that the pair $\xbar,\lammupair$ is a
    saddle point of the extended Lagrangian $\lagrangextsym$.
  \end{letter-compact}
  Furthermore, there exists such a point $\xbar\in\extspace$
  satisfying all of these conditions.
\end{corollary}

\begin{proof}
Implications
  (\ref{cor:ast-soln-kkt-saddle:saddle})
  $\Rightarrow$
  (\ref{cor:ast-soln-kkt-saddle:kkt}),
  (\ref{cor:ast-soln-kkt-saddle:kkt})
  $\Rightarrow$
  (\ref{cor:ast-soln-kkt-saddle:soln}),
  and the
  existence of a point satisfying the three conditions
all follow immediately
from \Cref{thm:lammu-equivs-fin-val}.

To show that
  (\ref{cor:ast-soln-kkt-saddle:soln})
  $\Rightarrow$
  (\ref{cor:ast-soln-kkt-saddle:saddle}),
suppose $\xbar$ is an astral solution.
By \Cref{roc:thm28.2}, there also exists a KKT vector
$\lammupair\in\Lampairs$.
Together, these satisfy statement~(\ref{thm:lammu-equivs:sln}) of 
\Cref{thm:lammu-equivs}, thereby implying
statement~(\ref{thm:lammu-equivs:sdl}) of 
that \namecref{thm:lammu-equivs}, which in turn
implies
this \namecref{cor:ast-soln-kkt-saddle}'s
statement~(\ref{cor:ast-soln-kkt-saddle:saddle}).%
\indexg{ordinary convex programs!existence of astral solution|)}%
\end{proof}

Without the assumption that $\inf p > -\infty$
in \Cref{cor:ast-soln-kkt-saddle},
statement~(\ref{cor:ast-soln-kkt-saddle:kkt}) is still sufficient
and statement~(\ref{cor:ast-soln-kkt-saddle:saddle}) is still necessary
for $\xbar$ to be an astral solution of $P$
(as implied by \Cref{thm:lammu-equivs}). However, the equivalence among the
three statements no longer holds.
As we show in the next example, if $\inf p=-\infty$, then
a point $\xbar$ can be an astral solution without satisfying the astral
KKT conditions for any $\lammupair$,
and also a point $\xbar$ can fail to be an astral solution,
while satisfying the saddle point condition for some $\lammupair$.

\begin{example}
In $\R^2$, consider the problem of minimizing $x_1$ subject to the
constraint that $x_2\leq 0$.
In other words, we aim to minimize $f(\xx)$ subject to
$g(\xx)\leq 0$
where
$f(\xx)=\xx\cdot\ee_1$
and
$g(\xx)=\xx\cdot\ee_2$.
(There are no equality constraints.)
The Lagrangian is
\[
   \lagrang{\xx}{\lambda}
   =
   f(\xx)+\lambda g(\xx)
   =
   \xx\cdot(\ee_1 + \lambda\ee_2).
\]
It can be checked that Slater's condition is satisfied, but
$\inf p=-\infty$ (where $p=f + \Negf\circ g$).

Let $\xbar_0=\limray{(-\ee_1)}$.
Then $\xbar_0$ is an astral solution since $\gext(\xbar_0)=0$
and $\fext(\xbar_0)=-\infty$.
On the other hand, for all $\lambda\in\Rpos$,
$\asubdiflagrangext{\xbar_0}{\lambda}=\{\ee_1 + \lambda \ee_2\}$
by Example~\ref{ex:affine-subgrad-new}.
Since this astral subdifferential does not include $\zero$,
the first-order optimality condition
(Eq.~\ref{eq:dfn:ast-kkt-defn:1}) cannot be satisfied
(by \Cref{pr:lagrang-ext-subdif}\ref{pr:lagrang-ext-subdif:a2:new}),
so the astral KKT conditions are not satisfied for any
$\lambda\in\Rpos$.

Next,
let $\xbar_1=\limray{(-\ee_1)}\plusl\ee_2$ and let $\lambda=1$.
This pair is a saddle point
because, for all $\lambda'\in\Rpos$,
$\lagrangext{\xbar_1}{\lambda'}=\xbar_1\cdot(\ee_1 + \lambda'\ee_2)=-\infty$,
so $\xbar_1$ minimizes $\lagrangext{\posSkip\cdot}{\lambda}$
over $\extspac{2}$,
and $\lambda$ maximizes $\lagrangext{\xbar_1}{\cdot\posSkip}$
over $\Rpos$.
On the other hand,
$\gext(\xbar_1)=1>0$, so
$\xbar_1$ is not astral feasible and therefore also is not an astral
solution.%
\indexg{constrained optimality conditions, astral!convex program@for convex program|)}%
\indexg{ordinary convex programs!astral optimality conditions|)}%
\end{example}

\indexg{Slater's condition!necessity for theory|(}%
Throughout this chapter, we have repeatedly assumed that Slater's
condition holds, as is standard in this area.
Without this assumption, much of the theory we have developed is no
longer valid, as we show in the next example:

\begin{example}    \label{ex:no-slater}
\indexg{Sideways cone!Slater's condition and|(}%
In $\R^3$, let $K$ be the cone and $L$ the plane given in
Eqs.~(\ref{eqn:bad-set-eg-K}) and~(\ref{eqn:bad-set-eg-L})
of
Example~\ref{ex:KLsets:extsum-not-sum-exts},
that is,
\begin{equation}
\label{eq:no-slater:KL}
  \begin{aligned}
    K
    &=\regBraces{\xx\in\R^3:\: x_1^2\le 2 x_2 x_3\text{ and }x_2,x_3\ge 0},
  \\
    L
    &=
    \regBraces{
      \xx\in\R^3 :\: x_3 = 0
    }.
  \end{aligned}
\end{equation}
Let $P$ be the ordinary convex program with objective function
$f:\R^3\to\eR$ defined by
\begin{equation*}
  f(\xx)
  =
  \begin{cases}
    \abs{x_1 - 1}   & \text{if $\xx\in K$,} \\
    +\infty     & \text{otherwise,}
  \end{cases}
\end{equation*}
and with a single equality constraint given by
$h(\xx)=\xx\cdot\ee_3=x_3$, for $\xx\in\R^3$.
(There are no inequality constraints.)

The function $f$ is convex, closed, proper and nonnegative
with effective domain
$\dom{f}=K$.
The set $C$ of feasible points is equal to
$L=\{\xx \in\R^3 :\: h(\xx)=0\}$.
Thus, Slater's condition does not hold since $\ri K$ and $L$ are
disjoint.

The aim of the program $P$ is to minimize
the function $p:\R^3\to\eR$, given by
\begin{align}
  p(\xx)
  =
  f(\xx) + \indL(\xx)
  &=
  \abs{x_1 - 1} + \indK(\xx) + \indL(\xx)
  \nonumber
  \\
  &=
  \abs{x_1 - 1} + \indKL(\xx)
  =
  1 + \indKL(\xx),
  \label{eq:ex:no-slater:1}
\end{align}
where the last equality is because
$x_1=0$
for all
$\xx\in K\cap L$ (see Eq.~\ref{eq:no-slater:KL}).
Thus, $\inf p = 1$.

For $\mu\in\R$, this program's Lagrangian is
\[
  \lagrang{\xx}{\mu}
  =
  f(\xx) + \mu h(\xx)
  =
  \abs{x_1 - 1} + \indK(\xx) + \mu x_3.
\]
In particular, this implies that $\inf \lagrang{\posSkip\cdot}{\mu}\leq 0$
for all $\mu\in\R$ since the sequence
of points
$\xx_t=\transKern{[1,\,t,\,1/(2t)]}$
is entirely in $K$ with $\lagrang{\xx_t}{\mu}\rightarrow 0$.
Thus,
\[
  \adjustlimits\sup_{\mu\in\R} \inf_{\xx\in\R^3} \lagrang{\xx}{\mu}
  \leq
  0
  <
  1
  =
  \inf p
  =
  \adjustlimits\inf_{\xx\in\R^3} \sup_{\mu\in\R} \lagrang{\xx}{\mu},
\]
so \eqref{eq:kkt-infsup-supinf} does not hold in this case.
This also shows that there does not exist a KKT vector for $P$.

The astral feasible set is $\Cbar=\Lbar$.
The function $f$'s extension is
\[
  \fext(\xbar)
  =
  \abs{\xbar\cdot\ee_1 - 1} + \indfa{\Kbar}(\xbar),
\]
for $\xbar\in\extspac{3}$,
using
\Cref{pr:ext-affine-comp}(\ref{pr:ext-affine-comp:c}),
\Cref{thm:ext-sum-fcns-w-duality},
and
\Cref{pr:inds-ext}.
Also, from \eqref{eq:ex:no-slater:1},
$\pext=1+\indfa{\KLbar}$,
using \Cref{pr:inds-ext}.

Let $\xbar_0 = \limray{\ee_2}$ and $\xbar_1=\limray{\ee_2}\plusl\ee_1$.
Then $\xbar_0$ is in $\KLbar$, and so also in $\Kbar$ and $\Lbar$,
since $\xbar_0$ is the limit of
$\seq{t\ee_2}$ whose every element is in $K\cap L$.
Also, $\xbar_1$ is in $\Kbar$ and~$\Lbar$, but not $\KLbar$,
as argued in \Cref{ex:KLsets:extsum-not-sum-exts}.
Thus, $\fext(\xbar_0)=\pext(\xbar_0)=1$,
$\fext(\xbar_1)=0$,
and
$\pext(\xbar_1)=+\infty$.
Furthermore, both $\xbar_0$ and $\xbar_1$ are astral feasible,
being in $\Cbar=\Lbar$.
Consequently, $\xbar_0$ is an astral solution of $P$, but does not
minimize $\fext$ over $\Cbar$.
On the other hand, $\xbar_1$~is not an astral solution, but does
minimize $\fext$ over $\Cbar$.
Thus, in neither case does the equivalence over astral solution
concepts given in
\Cref{pr:equiv-astral-kkt}(\ref{pr:equiv-astral-kkt:a},\ref{pr:equiv-astral-kkt:b})
hold.\looseness=-1

Turning to the characterization of astral solutions given in
\Cref{cor:ast-soln-kkt-saddle},
as already noted, $\xbar_1$ is not an astral solution.
Nonetheless, the pair $\xbar_1,0$ satisfies
the astral KKT conditions since $\xbar_1$ is astral feasible;
complementary
slackness trivially holds;
and $\xbar_1$ minimizes $\fext$ with $\fext(\xbar_1)$ finite,
so $\zero\in\asubdiffext{\xbar_1}$
(by \Cref{pr:asub-zero-is-min}\ref{pr:asub-zero-is-min:a}),
thereby satisfying first-order optimality (Eq.~\ref{eq:dfn:ast-kkt-defn:1}).
The pair $\xbar_1,0$ is also a saddle point.
This is because
\[
  \lagrangext{\xbar_1}{0}
  =
  \fext(\xbar_1)
  =
  0
  =
  \inf f
  =
  \min \fext
  =
  \min\lagrangext{\posSkip\cdot}{0},
\]  
and because, for all $\mu\in\R$,
and for the sequence $\seq{\xx_t}$ defined above,
$\xx_t\rightarrow\xbar_1$ and $\lagrang{\xx_t}{\mu}\rightarrow 0$;
hence,
$\lagrangext{\xbar_1}{\mu}\leq 0=\lagrangext{\xbar_1}{0}$.
Thus, $\xbar_1$ is not an astral solution despite satisfying
statements~(\ref{cor:ast-soln-kkt-saddle:kkt})
and~(\ref{cor:ast-soln-kkt-saddle:saddle})
of
\Cref{cor:ast-soln-kkt-saddle}.

On the other hand, $\xbar_0$ is an astral solution, but
there is no $\mu\in\R$ for which $\xbar_0,\mu$ would satisfy the astral KKT conditions
since $\asubdiffext{\xbar_0}=\{-\ee_1\}$ by
\Cref{cor:subgrad-comp-w-affine}
and $\asubdifmuhext{\xbar_0}=\{\mu \ee_3\}$ by
Example~\ref{ex:affine-subgrad-new},
so
$\zero\not\in\{-\ee_1+\mu\ee_3\}=\asubdiffext{\xbar_0}+\asubdifmuhext{\xbar_0}$,
implying that first-order optimality
(Eq.~\ref{eq:dfn:ast-kkt-defn:1}) cannot hold.
Also, $\xbar_0,\mu$ is not a saddle point since
$\lagrangext{\xbar_0}{\mu}=\fext(\xbar_0)+\mu \xbar\cdot\ee_3=1$
by \Cref{thm:ext-sum-fcns-w-duality}, but on the other hand,
$\min \lagrangext{\posSkip\cdot}{\mu}=\inf \lagrang{\posSkip\cdot}{\mu}\leq 0$,
as argued above.
Consequently, $\xbar_0$ is an astral solution, but does not satisfy
either of the
statements~(\ref{cor:ast-soln-kkt-saddle:kkt})
or~(\ref{cor:ast-soln-kkt-saddle:saddle})
of
\Cref{cor:ast-soln-kkt-saddle}.%
\indexg{Sideways cone!Slater's condition and|)}%
\indexg{Slater's condition!necessity for theory|)}%
\end{example}

\indexg{equality constraints (of convex program)!matrix form@in matrix form|(}%
Throughout this section,
we have chosen to work with affine equality constraints,
$h_j(\xx)=0$, separately for each $j$.
Nonetheless,
when convenient, these can be combined into a single matrix
constraint, both for standard and astral feasibility.
More precisely,
suppose, for $j=1,\ldots,s$, that $h_j(\xx)=\xx\cdot\aaa_j-b_j$
for $\xx\in\Rn$ where $\aaa_j\in\Rn$ and $b_j\in\R$.
Then the equality constraints appearing in $P$
mean that $\xx\cdot\aaa_j=b_j$ for all~$j$, or, in matrix form,
that $\A\xx=\bb$ where $\A\in\R^{s\times n}$ is a matrix whose $j$-th
row is equal to $\trans{\aaa}_j$
(and, as usual, $\bb=\trans{[b_1,\ldots,b_s]}$).
For astral feasibility, these equality constraints require,
for all $j$,
that $\hext_j(\xbar)=0$, that is, that
$\xbar\cdot\aaa_j=b_j$ (by Example~\ref{ex:ext-affine}).
By \Cref{pr:mat-prod-is-row-prods-if-finite},
this condition is equivalent to the astral matrix constraint
$\A\xbar=\bb$.%
\indexg{equality constraints (of convex program)!matrix form@in matrix form|)}%

\indexg{ordinary convex programs!only equality contraints|(}%
\indexg{constrained optimality conditions, astral!affine constraints@with affine constraints|(}%
For the special case that there are no inequality constraints but only
affine equality constraints, the convex program $P$ is equivalent to
the problem of minimizing a convex function $f(\xx)$ subject to
the affine constraint $\A\xx=\bb$, for some
$\A\in\R^{s\times n}$ and $\bb\in\R^s$.
We considered exactly this problem in
\Cref{thm:ast-opt-lin-cons-primal}(\ref{thm:ast-opt-lin-cons-primal:a},\ref{thm:ast-opt-lin-cons-primal:d}),
where we showed that $\xbar$ is an astral solution of $P$ if and only if
$\A\xbar=\bb$ and $\basubdiffext{\xbar}\cap(\colspace \transA)\neq\emptyset$
(assuming Slater's condition and $\inf p>-\infty$).
We next give an alternate proof using astral KKT conditions and
\Cref{cor:ast-soln-kkt-saddle}(\ref{cor:ast-soln-kkt-saddle:soln},\ref{cor:ast-soln-kkt-saddle:kkt}):

\begin{proof}[Alternate proof of
    \Cref{thm:ast-opt-lin-cons-primal}
    (\ref{thm:ast-opt-lin-cons-primal:a})
    $\Leftrightarrow$
    (\ref{thm:ast-opt-lin-cons-primal:d}):%
  ]
~  

Assume the setup of \Cref{thm:ast-opt-lin-cons-primal}.
For $j=1,\ldots,m$, let $\trans{\aaa_j}$ be the $j$-th row of~$\A$
(hence $\aaa_j\in\Rn$),
and, for $\xx\in\Rn$, let $h_j(\xx)=\xx\cdot\aaa_j - b_j$
(where $\bb=\trans{[b_1,\ldots,b_m]}$).
Let $P$ be the resulting convex program defined by $f$ and
$h_1,\ldots,h_m$ (with no inequality constraints,
so $r=0$ and $s=m$).
Then the resulting function $p$ as in
\eqref{eqn:kkt-p-defn} is the same as in
\Cref{thm:ast-opt-lin-cons-primal}, and
the assumption of $\xing$'s existence ensures that Slater's condition
is satisfied.
The set $C$ of standard feasible points in
\eqref{eqn:feas-set-defn} 
is also the same as in \Cref{thm:ast-opt-lin-cons-primal},
with closure
$\Cbar=\{\xbar'\in\extspace :\: \A\xbar'=\bb\}$
by \Cref{cor:closure-affine-set}.

For all $\muvec\in\R^m$,
\begin{equation}  \label{eq:cor:kkt-no-ineq-consts:2}
  \sumjtom \basubdifmuhjext{\xbar} = \BiggBraces{\,\sumjtom \mu_j \aaa_j }
\end{equation}
since,
by Example~\ref{ex:affine-subgrad-new},
$\basubdifmuhjext{\xbar}=\{\mu_j \aaa_j\}$
for all $j$.
By \Cref{cor:ast-soln-kkt-saddle}(\ref{cor:ast-soln-kkt-saddle:soln},\ref{cor:ast-soln-kkt-saddle:kkt}),
$\xbar$ is an astral solution of $P$
if and only if
it satisfies the astral KKT conditions,
which in this case require that
$\A\xbar=\bb$ (so that $\xbar$ is in $\Cbar$ and so astral feasible)
and that
there exist $\uu\in\basubdiffext{\xbar}$
such that
$-\uu \in \sumjtom \basubdifmuhjext{\xbar}$
for some $\muvec\in\R^m$
(so that $\zero$ is in the sum of these astral subdifferentials).
By \eqref{eq:cor:kkt-no-ineq-consts:2},
this last condition holds if and only if
$-\uu\in\colspace \transA$, and so also if and only if
$\uu\in\colspace \transA$.
Thus, as claimed in \Cref{thm:ast-opt-lin-cons-primal},
$\xbar$ is an astral solution of $P$
if and only if
$\A\xbar=\bb$
and
$\basubdiffext{\xbar} \cap (\colspace \transA)\neq\emptyset$.%
\indexg{ordinary convex programs!only equality contraints|)}%
\indexg{constrained optimality conditions, astral!affine constraints@with affine constraints|)}%
\end{proof}

\indexg{ordinary convex programs!no finite solution|(}%
We earlier saw that the convex program in
Example~\ref{ex:kkt-running-ex-standard}
had no finite solution.
We now return to that example to see in detail
how the same program can
be analyzed using the astral tools developed above.

\begin{example}   \label{ex:kkt-running-ex-astral}
  We consider the convex program of
  \Cref{ex:kkt-running-ex-standard}
  with $f$, $g_1$, and $g_2$ defined as in
  \eqref{eq:ex:kkt-running-ex-standard:1}.
  We can rewrite these functions, for $\xx\in\R^2$, as
  \[
    f(\xx)=\xx\inprod\ee_2,
    \qquad
    g_1(\xx) = \exp(\xx\inprod\ee_1) + 1 - \xx\inprod\ee_2,
    \qquad
    g_2(\xx) = \xx\inprod\ee_1.
  \]
  We aim to find a pair
  $\xbar\in\extspac{2}$ and $\lamvec\in\Rpos^2$
  which satisfies the
  astral KKT conditions, thereby implying,
  by
  \Cref{cor:ast-soln-kkt-saddle}(\ref{cor:ast-soln-kkt-saddle:kkt},\ref{cor:ast-soln-kkt-saddle:soln}),
 that $\xbar$ is an astral solution.
  To do so, we suppose that the pair $\xbar,\lamvec$
  already satisfies the astral KKT conditions, and analyze the necessary
  properties of such a pair, thus deriving a solution.

  These functions' extensions are
  \begin{equation}
  \label{eq:kkt-running:3}
    \fext(\xbar)=\xbar\inprod\ee_2,
    \qquad
    \gext_1(\xbar) = \expex(\xbar\inprod\ee_1) + 1 - \xbar\inprod\ee_2,
    \qquad
    \gext_2(\xbar) = \xbar\inprod\ee_1.
  \end{equation}
  The first and third expressions are from \Cref{ex:ext-affine}.
  Note, by primal feasibility, that this implies
  $\xbar\inprod\ee_1=\gext_2(\xbar)\le 0$.
  The expression for $\gext_1$ then follows from
  \Cref{thm:ext-sum-fcns-w-duality} applied to
  $s_1(\xx)=\exp(\xx\inprod\ee_1)$
  and $s_2(\xx)=1 - \xx\inprod\ee_2$,
  for $\xx\in\R^2$,
  with extensions
  $\sext_1(\xbar)=\expex(\xbar\inprod\ee_1)$
  and $\sext_2(\xbar)=1 - \xbar\inprod\ee_2$,
  by
  \Cref{pr:j:2}(\ref{pr:j:2c})
  and again \Cref{ex:ext-affine}.
  For summability, we note that $\sext_1(\xbar)\in\R$
  since $\xbar\cdot\ee_1\leq 0$.

  We next calculate astral subdifferentials:
  \begin{equation}
  \label{eq:kkt-running:1}
      \asubsym\ef(\xbar)=\set{\ee_2},
      \quad
      \asubsym\eg_1(\xbar)=\regBraces{\ee_1\expex(\xbar\inprod\ee_1)-\ee_2},
      \quad
      \asubsym\eg_2(\xbar)=\set{\ee_1}.
  \end{equation}
  The first and third expressions are
  from \Cref{ex:affine-subgrad-new}.
  The second follows
  by \Cref{thm:subgrad-sum-fcns} applied to
  $s_1$ and $s_2$, as above, whose astral subdifferentials are
  $\asubdifsextsub{2}{\xbar}=\set{-\ee_2}$ by \Cref{ex:affine-subgrad-new},
  and
  $\asubdifsextsub{1}{\xbar}=\set{\ee_1\expex(\xbar\inprod\ee_1)}$ by
  \Cref{thm:subgrad-fA}
  (applied with $\A$ and $f$, as they appear in that
  \namecref{thm:subgrad-fA}, set to $\trans{\ee_1}$ and $\exp$).
  We also use that
  $\asubsym\expex(\bary)=\set{\expex(\bary)}$ for
  $\bary\in[-\infty,+\infty)$
  (by
  \Cref{pr:asubdiffext-at-x-in-rn}\ref{pr:asubdiffext-at-x-in-rn:c}
  and~\Cref{pr:subdif-in-1d}),
  and that $\xbar\inprod\ee_1\le 0$.
  
  Thus, by first-order optimality,
  \begin{align*}
    \zero &\in
           \asubdiffext{\xbar}
           +\asubsym(\overline{\lambda_1 g_1})(\xbar)
           +\asubsym(\overline{\lambda_2 g_2})(\xbar)
  \\
    &=\bigBraces{
              \ee_2
              + \lambda_1\bigParens{\ee_1\expex(\xbar\inprod\ee_1)-\ee_2}
              + \lambda_2\ee_1
    }
  \\
    &=\bigBraces{\ee_2(1-\lambda_1)
              + \ee_1\bigParens{\lambda_1\expex(\xbar\inprod\ee_1) + \lambda_2}
    },
  \end{align*}
  where the first equality follows from \eqref{eq:kkt-running:1}
  and \Cref{pr:subgrad-scal-mult} (noting that $\lambda_1,\lambda_2\in\Rpos$
  by dual feasibility, and that $\asubsym\eg_1(\xbar)\ne\emptyset$
  and $\asubsym\eg_2(\xbar)\ne\emptyset$).
  Hence,
  \[
      \zero = \ee_2(1-\lambda_1)
      + \ee_1\bigParens{\lambda_1\expex(\xbar\inprod\ee_1) + \lambda_2}.
  \]
  This is only possible if $\lambda_1=1$, and so,
  since $\lambda_2\ge 0$, if also $\lambda_2=0$
  and $\xbar\inprod\ee_1=-\infty$.

  Complementary slackness then implies that
  \[
     0=\lambda_1 \eg_1(\xbar)
      =\expex(\xbar\inprod\ee_1)+1-\xbar\inprod\ee_2
      =1-\xbar\inprod\ee_2,
  \]
  with the second equality by \eqref{eq:kkt-running:3}
  and since $\lambda_1=1$, 
  and the third because
  $\xbar\inprod\ee_1=-\infty$.
  Thus, $\xbar\inprod\ee_2=1$.
  
  The two conditions $\xbar\inprod\ee_1=-\infty$
  and $\xbar\inprod\ee_2=1$ are satisfied by $\xbar=\limray{(-\ee_1)}\plusl\ee_2$
  (which, it can be checked, is also the only point satisfying these two conditions).
  By construction, the pair of $\xbar=\limray{(-\ee_1)}\plusl\ee_2$ and $\lamvec=\trans{[1,0]}$
  satisfies first-order optimality. It can be checked that it
  also satisfies primal and dual feasibility and complementary slackness. Thus,
  the pair $\xbar,\lamvec$ satisfies the astral KKT conditions, and therefore $\xbar$
  is an astral solution of the studied convex program
  with value $\fext(\xbar)=1$.%
\indexg{ordinary convex programs!no finite solution|)}%
\end{example}

\chapter{Differential continuity and descent methods}
\label{sec:diff-cont:descent}

Many methods for minimizing a convex function are iterative, meaning
they compute a sequence of iterates $\seq{\xx_t}$ in $\Rn$
with the goal of minimizing
the function in the limit.
Often, such methods are related to subgradients, for instance,
operating in a way to drive the subgradients at the computed iterates
to zero, or by taking steps in the direction of a subgradient.
This chapter applies astral differential theory to the analysis of
such methods.
We first consider continuity properties of subdifferentials.
These are then applied to derive general conditions for when a class of
iterative methods based on subgradients will converge,
leading to a convergence analysis of many standard algorithms
on a wide range of convex optimization problems.

\section{Convergence of subgradient sequences}
\label{sec:subdif:cont}

\indexg{subgradient sequences (standard)!converging to 0@converging to $\zero$|(}%
A convex function $f:\Rn\rightarrow\Rext$ is minimized at a
point $\xx\in\Rn$ if and only if $\zero$ is a subgradient of $f$ at
$\xx$, and in particular, if the gradient $\gradf(\xx)$ exists and is
equal to~$\zero$.
Therefore, to minimize~$f$ numerically, it is natural to construct a
sequence $\seq{\xx_t}$ in~$\Rn$ whose gradients converge to $\zero$,
that is, for which $\gradf(\xx_t)\rightarrow\zero$.
If this is possible, and if the sequence converges to some point
$\xx\in\Rn$, then indeed, $\xx$ must minimize~$f$,
as follows from
\Cref{roc:thm24.4},
assuming $f$ is closed and proper;
in that case, $f(\xx)=\inf f$.
Moreover, even if the sequence $\seq{\xx_t}$ does not converge but
nevertheless remains within a bounded region of $\Rn$, then an
argument can again be made that $f(\xx_t)\rightarrow \inf f$.

Since driving the gradient of a function to $\zero$ seems so closely
connected to minimization, especially for convex functions,
one might expect that it should also be effective as a means of
minimizing the function for any sequence, not just
bounded sequences.
In other words, for a convex function
$f:\Rn\rightarrow\Rext$ and sequence $\seq{\xx_t}$ in $\Rn$, we might
expect that if $\gradf(\xx_t)\rightarrow\zero$ then
$f(\xx_t)\rightarrow \inf f$.
However, this is false in general, even for a convex function with
many favorable properties, as shown in the next example:

\begin{example}   \label{ex:x1sq-over-x2:grad}
\indexg{subgradient sequences (standard)!continuity and|(}%
\indexg{continuity of extensions!subgradient sequences and|(}%
\indexg{continuity of extensions!minimizing sequences and|(}%
\indexg{Flattening valley!not minimized with vanishing gradients|(}%
Let $f:\R^2\rightarrow\Rext$ be the flattening valley function from
Examples~\ref{ex:x1sq-over-x2} and~\ref{ex:x1sq-over-x2:cont},
that is, for $\xx\in\R^2$,
  \[
    f(\xx)
    =
    \begin{cases}
      x_1^2/x_2
        & \text{if $x_2 >\regAbs{x_1}$,}
      \\
      2\regAbs{x_1}-x_2
        & \text{otherwise.}
    \end{cases}
  \]
As already noted,
this function is convex,
closed, proper, continuous everywhere, finite everywhere, and
nonnegative everywhere.
It is also continuously differentiable everywhere except along the ray
$\{\trans{[0,x_2]} :\: x_2 \in \Rneg \}$,
a part of the space that is
far from the sequences we will be considering.
For each $t$, let
$\xx_t=\trans{[t^2,\,t^3]}$.
Then it can be calculated
that $f$'s gradient at each $\xx_t$ is
$ \gradf(\xx_t) = \trans{[ 2/{t},\, -1/t^2]} $,
so $\gradf(\xx_t)\rightarrow\zero$.
Nevertheless,
$f(\xx_t)={t}\rightarrow +\infty$.
Thus, the gradients are converging to $\zero$, but the function values
are growing to $+\infty$, whereas $\inf f = 0$.%
\indexg{Flattening valley!not minimized with vanishing gradients|)}%
\end{example}

It is no coincidence that this same function $f$ was used earlier
in Example~\ref{ex:x1sq-over-x2:cont}
as an example of a function whose extension $\fext$ is discontinuous.
The same sequence was also used in that discussion,
where we showed that
$\xx_t\rightarrow \xbar=\limray{\ee_2}\plusl\limray{\ee_1}$,
and that $\fext$ is discontinuous at $\xbar$.
Indeed, there is a close connection between continuity in astral
space and convergence of gradients.
We prove below that for any convex function $f:\Rn\rightarrow\Rext$,
if $\seq{\xx_t}$ is a sequence in $\Rn$ converging to a point
$\xbar\in\extspace$ where $\fext$ is continuous
(and also with $\fext(\xbar)<+\infty$),
and if $\gradf(\xx_t)\rightarrow\zero$ then
$f(\xx_t)\rightarrow \inf f$.
If $\fext$ is not continuous at $\xbar$, then this statement need not
hold, as demonstrated in the preceding example.%
\indexg{subgradient sequences (standard)!converging to 0@converging to $\zero$|)}%
\indexg{subgradient sequences (standard)!continuity and|)}%
\indexg{continuity of extensions!subgradient sequences and|)}%
\indexg{continuity of extensions!minimizing sequences and|)}%

\indexg{subgradient sequences, astral (primal)|(}%
\indexg{continuity!astral subgradients@of astral subgradients|(}%
In more generality,
for a function $f:\Rn\rightarrow\Rext$,
we study next the convergence properties of a sequence of astral subgradients of
the extension $\fext$.
More specifically,
suppose $\seq{\xbar_t}$ is a sequence in $\extspace$ and $\seq{\uu_t}$
a sequence in $\Rn$,
where each $\uu_t$ is an astral subgradient of $\fext$ at $\xbar_t$
so that $\uu_t\in\asubdiffext{\xbar_t}$.
Suppose further that $\xbar_t\rightarrow\xbar$ and
$\uu_t\rightarrow\uu$ for some $\xbar\in\extspace$ and $\uu\in\Rn$.
We seek conditions that guarantee $\uu\in\asubdiffext{\xbar}$
so that $\uu$, the limit of the subgradients $\uu_t$,
will itself be an astral subgradient of $\fext$ at $\xbar$,
the limit of the points $\xbar_t$.%
\indexg{continuity!astral subgradients@of astral subgradients|)}%

\indexg{subgradient sequences (standard)!gradients and|(}%
This setting includes, when $f$ is proper, the special case that
$\xbar_t=\xx_t\in\Rn$ for all $t$ and that
$\uu_t\in\partial f(\xx_t)$ so that
each $\uu_t$ is a standard subgradient of $f$ at $\xx_t$
and therefore also an astral subgradient of $\fext$ by
\Cref{pr:asubdiffext-at-x-in-rn}(\ref{pr:asubdiffext-at-x-in-rn:b}).
Indeed, from standard convex analysis, it is known in this case that if
$\xx_t\rightarrow\xbar=\xx\in\Rn$ and if $f$ is convex, closed and proper, then
$\uu\in\partial f(\xx)=\asubdiffext{\xx}$
(by
\Cref{roc:thm24.4}).%
\indexg{subgradient sequences (standard)!gradients and|)}%

\indexg{subgradient sequences, astral (primal)!coupling function and|(}%
\indexg{coupling function!subgradient sequences@on subgradient sequences|(}%
The sequence of values $\xbar_t\cdot\uu_t$ will play a
central role in characterizing convergence of the kind just described.
In general, this sequence need {not} converge to $\xbar\cdot\uu$.
Here is an example:

\begin{example}
Let $f$ and $\seq{\xx_t}$ be as in
Example~\ref{ex:x1sq-over-x2:grad},
and let
$\xbar=\limray{\ee_2}\plusl\limray{\ee_1}$, $\uu=\zero$, and
$\uu_t=\gradf(\xx_t)\in\partial f(\xx_t)=\asubdiffext{\xx_t}$
(by \Cref{roc:thm25.1}\ref{roc:thm25.1:a} and
\Cref{pr:asubdiffext-at-x-in-rn}\ref{pr:asubdiffext-at-x-in-rn:b}).
Then $\xx_t\rightarrow\xbar$ and $\uu_t\rightarrow\uu$.
It can be checked that $\xx_t\cdot\uu_t={t}\rightarrow+\infty$,
but $\xbar\cdot\uu=0$;
thus, $\xx_t\cdot\uu_t\not\rightarrow\xbar\cdot\uu$.
\end{example}

The convergence properties of $\seq{\xbar_t\cdot\uu_t}$
turn out to be closely connected both
to the convergence of subgradients $\uu_t$ to a subgradient at $\xbar$
and to the continuity of $\fext$ at
$\xbar$, as will be proved in the next few theorems.
We show first in the next lemma that the closer $\xbar_t\cdot\uu_t$ comes to
$\xbar\cdot\uu$, the closer will the function values $\fext(\xbar_t)$ get to
$\fext(\xbar)$.
Here and for most of the results that follow, we assume
$\xbar\cdot\uu\in\R$, thus including the case 
$\uu=\zero$, which is especially relevant for minimization.

\begin{lemma}   \label{thm:limsup-beta-bnd}
  Let $f:\Rn\rightarrow\Rext$.
  Let $\seq{\xbar_t}$ be a sequence in $\extspace$
  and $\seq{\uu_t}$ a sequence in $\Rn$
  with each $\uu_t\in\asubdiffext{\xbar_t}$, and with
  $\xbar_t\rightarrow\xbar$ and $\uu_t\rightarrow\uu$ for some
  $\xbar\in\extspace$ and $\uu\in\Rn$.
  Assume $\xbar\cdot\uu\in\R$, and
  suppose $\limsup(\xbar_t\cdot\uu_t) \leq \beta$
  for some $\beta\in\R$.
  Then
  \begin{equation}    \label{eq:thm:limsup-beta-bnd:2}
    \fext(\xbar)
    \leq
    \liminf \fext(\xbar_t)
    \leq
    \limsup \fext(\xbar_t)
    \leq
    \beta - \fstar(\uu)
    \leq
    \fext(\xbar) + \beta - \xbar\cdot\uu.
  \end{equation}
\end{lemma}

\begin{proof}
The first inequality is because $\fext$ is lower semicontinuous by
\Cref{prop:ext:F}(\ref{prop:ext:F:a}),
and since $\xbar_t\rightarrow\xbar$.
The last inequality
follows from the definition of astral conjugate
(Eq.~\ref{eq:Fstar-down-def},
applied to $F=\fext$,
and using \Cref{pr:fextstar-is-fstar}).
It remains to prove the third inequality.

Let $\epsilon\in\Rstrictpos$.
We then have
\begin{align*}
  \limsup \fext(\xbar_t)
  &=
  \limsup \bigBracks{\xbar_t\cdot\uu_t - \fstar(\uu_t)}
  \\
  &\leq
  \limsup \bigBracks{\beta+\epsilon - \fstar(\uu_t)}
  \leq
  \beta + \epsilon - \fstar(\uu).
\end{align*}
The equality is
by
\Cref{thm:fenchel-subgrad}(\ref{thm:fenchel-subgrad:b},\ref{thm:fenchel-subgrad:c})
(and \Cref{pr:fextstar-is-fstar}) since
$\uu_t\in\asubdiffext{\xbar_t}$ for all $t$,
noting that $\fstar(\uu_t)\in\R$ by 
\Cref{pr:subgrad-imp-in-cldom}(\ref{pr:subgrad-imp-in-cldom:c}).
The first inequality is because
$\xbar_t\cdot\uu_t\leq\beta+\epsilon$
for all $t$ sufficiently large
since $\limsup\xbar_t\cdot\uu_t\leq\beta$.
The last inequality is because $\fstar$
is closed and therefore
lower semicontinuous
(\Cref{pr:conj-props}\ref{pr:conj-props:d}),
implying
$\liminf \fstar(\uu_t)\geq \fstar(\uu)$
since $\uu_t\rightarrow\uu$.
Since this holds for all $\epsilon\in\Rstrictpos$,
this completes the proof.
\end{proof}

\indexg{subgradient sequences, astral (primal)!convergence of function values|(}%
\indexg{subgradient sequences, astral (primal)!convergence to subgradient|(}%
\indexg{continuity!astral subgradients@of astral subgradients|(}%
As we show next,
under the same assumptions and 
taking $\beta=\xbar\cdot\uu$, \Cref{thm:limsup-beta-bnd}
immediately implies that
if $\limsup \xbar_t\cdot\uu_t \leq \xbar\cdot\uu$,
then all the following must also hold:
the sequence of values $\fext(\xbar_t)$ must converge to $\fext(\xbar)$;
the astral form of the Fenchel-Young inequality
(Eq.~\ref{eqn:ast-fenchel-alt})
must hold with equality for $\fext$;
and, if $\fstar(\uu)$ is finite, then
$\uu$ must be an astral subgradient of $\fext$ at $\xbar$.

\begin{theorem}  \label{cor:xuleqxbaru-then-subgrad}
  Let $f:\Rn\rightarrow\Rext$.
  Let $\seq{\xbar_t}$ be a sequence in $\extspace$
  and $\seq{\uu_t}$ a sequence in $\Rn$
  with each $\uu_t\in\asubdiffext{\xbar_t}$, and with
  $\xbar_t\rightarrow\xbar$ and $\uu_t\rightarrow\uu$ for some
  $\xbar\in\extspace$ and $\uu\in\Rn$.
  Assume $\xbar\cdot\uu\in\R$ and
  that $\limsup(\xbar_t\cdot\uu_t) \leq \xbar\cdot\uu$.
  Then:
  \begin{letter-compact}
  \item  \label{cor:xuleqxbaru-then-subgrad:a}
    $\fext(\xbar_t)\rightarrow \fext(\xbar)$.
  \item  \label{cor:xuleqxbaru-then-subgrad:b}
    $\fext(\xbar)=\xbar\cdot\uu - \fstar(\uu)$.
  \item  \label{cor:xuleqxbaru-then-subgrad:c}
    If $\fstar(\uu)\in\R$, then $\uu\in\asubdiffext{\xbar}$.
  \end{letter-compact}
\end{theorem}

\begin{proof}
Given the stated assumptions,
we can immediately apply \Cref{thm:limsup-beta-bnd}
with $\beta=\xbar\cdot\uu$, yielding
$ \fext(\xbar) = \lim \fext(\xbar_t) = \xbar\cdot\uu - \fstar(\uu)$,
thus proving
parts~(\ref{cor:xuleqxbaru-then-subgrad:a})
and~(\ref{cor:xuleqxbaru-then-subgrad:b}).

Assuming $\fstar(\uu)\in\R$,
part~(\ref{cor:xuleqxbaru-then-subgrad:b})
then implies $\uu\in\asubdiffext{\xbar}$
by
\Cref{thm:fenchel-subgrad}(\ref{thm:fenchel-subgrad:c},\ref{thm:fenchel-subgrad:b})
(and \Cref{pr:fextstar-is-fstar})
since $\xbar\cdot\uu$ is finite and therefore summable with $-\fext(\xbar)$.%
\indexg{subgradient sequences, astral (primal)!convergence of function values|)}%
\indexg{subgradient sequences, astral (primal)!convergence to subgradient|)}%
\indexg{continuity!astral subgradients@of astral subgradients|)}%
\end{proof}

\Cref{thm:limsup-beta-bnd} also provides
asymptotic lower bounds on the sequence $\seq{\xbar_t\cdot\uu_t}$,
showing that
$\liminf(\xbar_t\cdot\uu_t)\geq\xbar\cdot\uu$ if $\fext(\xbar)>-\infty$,
and that $\xbar_t\cdot\uu_t\rightarrow+\infty$ if
$\fext(\xbar)=+\infty$:

\begin{theorem}  \label{cor:xugeqxbaru}
  Let $f:\Rn\rightarrow\Rext$, and assume $f\not\equiv+\infty$.
  Let $\seq{\xbar_t}$ be a sequence in $\extspace$
  and $\seq{\uu_t}$ a sequence in $\Rn$
  with each $\uu_t\in\asubdiffext{\xbar_t}$, and with
  $\xbar_t\rightarrow\xbar$ and $\uu_t\rightarrow\uu$ for some
  $\xbar\in\extspace$ and $\uu\in\Rn$.
  Assume $\xbar\cdot\uu\in\R$ and $\fext(\xbar)>-\infty$.
  Then $\liminf \xbar_t\cdot\uu_t \geq \xbar\cdot\uu$.
  If, in addition, $\fext(\xbar)=+\infty$, then
  $\xbar_t\cdot\uu_t\rightarrow+\infty$.
\end{theorem}

\begin{proof}
We consider first the case that $\fext(\xbar)\in\R$.
Suppose, by way of contradiction,
that $\liminf \xbar_t\cdot\uu_t < \xbar\cdot\uu$.
Then there exists $\epsilon\in\Rstrictpos$
and infinitely many values of~$t$ for
which
$\xbar_t\cdot\uu_t \leq \xbar\cdot\uu - \epsilon$.
By discarding all other sequence elements, we can assume that this
holds for all values of $t$.
We can therefore apply
\Cref{thm:limsup-beta-bnd}
with $\beta=\xbar\cdot\uu - \epsilon$.
However, this yields
$\fext(\xbar) \leq \fext(\xbar) - \epsilon$,
an obvious contradiction.
Therefore,
$\liminf \xbar_t\cdot\uu_t \geq \xbar\cdot\uu$
if $\fext(\xbar)\in\R$.

Next,
consider the case that $\fext(\xbar)=+\infty$.
Suppose $\xbar_t\cdot\uu_t \not\rightarrow +\infty$.
Then there exists $\beta\in\R$ such that
$\xbar_t\cdot\uu_t\leq\beta$ for
infinitely many values of $t$.
As before, we can discard all other sequence elements, so that this
holds for all values of $t$.
We can again apply
\Cref{thm:limsup-beta-bnd}
with this choice of $\beta$, yielding
$\fext(\xbar) \leq \beta - \fstar(\uu)$.
Since $\fext(\xbar)=+\infty$, this implies $\fstar(\uu)=-\infty$.
But this is a contradiction since we assumed $f\not\equiv+\infty$,
implying $\fstar>-\infty$
(\Cref{pr:conj-props}\ref{pr:conj-props:c1}).
Therefore, $\lim(\xbar_t\cdot\uu_t) = +\infty \geq \xbar\cdot\uu$
if $\fext(\xbar)=+\infty$.
\end{proof}

\indexg{subgradient sequences, astral (primal)!convergence to subgradient|(}%
\indexg{continuity!astral subgradients@of astral subgradients|(}%
The next example shows that
\Cref{cor:xuleqxbaru-then-subgrad}(\ref{cor:xuleqxbaru-then-subgrad:c})
is false in general if the condition $\fstar(\uu)\in\R$ is omitted.
In other words, if all of the conditions of that theorem are
satisfied,
it is possible that $\fstar(\uu)\not\in\R$ and that
$\uu\not\in\asubdiffext{\xbar}$.
This example also shows that
if $\fext(\xbar)=-\infty$, then
\Cref{cor:xugeqxbaru} need not hold;
that is, it need not be the case that
$\liminf(\xbar_t\cdot\uu_t) \geq \xbar\cdot\uu$.

\begin{example}
Let $f:\R\rightarrow\Rext$ be the convex function
\[
   f(x)
   =
   \begin{cases}
               -\ln x         & \mbox{if $x > 0$,} \\
               +\infty         & \mbox{otherwise,}
   \end{cases}
\]
for $x\in\R$.
For all $t$,
let $x_t=t$, and $u_t=f'(x_t)=-1/t$, where $f'$ is the
derivative of~$f$ (where defined),
so that $u_t\in\partial f(x_t)=\asubdiffext{x_t}$.
Then $x_t \rightarrow \barx$ and $u_t \rightarrow u$
where $\barx=+\infty$ and $u=0$.
Also, $x_t u_t = -1$ for all $t$, so
$x_t u_t \rightarrow -1 < 0 = \barx u$;
thus, $\liminf x_t u_t < \barx u$.

Further,
all the assumptions of
\Cref{cor:xuleqxbaru-then-subgrad}
are satisfied, but $\fstar(u)=-\inf f =+\infty$, implying
$u\not\in\asubdiffext{\barx}$ by 
\Cref{pr:subgrad-imp-in-cldom}(\ref{pr:subgrad-imp-in-cldom:c}).%
\indexg{subgradient sequences, astral (primal)|)}%
\indexg{subgradient sequences, astral (primal)!convergence to subgradient|)}%
\indexg{subgradient sequences, astral (primal)!coupling function and|)}%
\indexg{coupling function!subgradient sequences@on subgradient sequences|)}%
\end{example}

\indexg{subgradient sequences (standard)!convergence of function values|(}%
\indexg{subgradient sequences (standard)!continuity and|(}%
\indexg{continuity of extensions!subgradient sequences and|(}%
\indexg{minimizing sequences!vanishing subgradients@with vanishing subgradients|(}%
We now focus on deriving more specific conditions for
convergence of subgradients
in the special case that $\xbar_t=\xx_t\in\Rn$ for all $t$,
and that $\uu_t\in\partial f(\xx_t)$.
\Cref{cor:xuleqxbaru-then-subgrad} shows
that if $\xbar\inprod\uu\in\R$ and $\fstar(\uu)\in\R$, then
to prove
$\uu$ is an astral subgradient of $\fext$ at~$\xbar$,
it suffices to show
$\limsup(\xx_t\cdot\uu_t)\leq \xbar\cdot\uu$.
Indeed, assuming $f$ is convex and proper,
this will be the case if $\fext$ is continuous at $\xbar$
and if $\fext(\xbar)<+\infty$, as we show next, along with some other
more general sufficient conditions.
Thus, continuity in astral space, together with our other assumptions,
provides a sufficient condition for a
sequence of standard subgradients to converge to an astral subgradient.
Our earlier counterexample
(Example~\ref{ex:x1sq-over-x2:grad})
shows that this need not be true in general
without continuity.

\begin{theorem}   \label{thm:cont-subgrad-converg}
  Let $f:\Rn\rightarrow\Rext$ be convex and proper.
  Let $\seq{\xx_t}$ and $\seq{\uu_t}$ be sequences in $\Rn$
  with each $\uu_t\in\partial f(\xx_t)$, and with
  $\xx_t\rightarrow\xbar$ and $\uu_t\rightarrow\uu$ for some
  $\xbar\in\extspace$ and $\uu\in\Rn$.
  Assume $\xbar\cdot\uu\in\R$, and that at least one of the following
  conditions also holds:
  \begin{letter}
    \item      \label{thm:cont-subgrad-converg:a}
      $\xbar=\ebar\plusl\qq$ for some $\ebar\in\corezn$ and
      $\qq\in\Rn$,
      and for each $t$, $\xx_t=\ww_t+\qq_t$ for some
      $\ww_t\in\resc{f}$ and $\qq_t\in\Rn$ with $\qq_t\rightarrow\qq$.
    \item      \label{thm:cont-subgrad-converg:b}
      $\xbar\in\represc{f}\plusl\Rn$.
    \item      \label{thm:cont-subgrad-converg:c}
      $\xbar\in\contsetf$; that is, $\fext(\xbar)<+\infty$ and $\fext$
      is continuous at $\xbar$.
  \end{letter}
  Then $\limsup(\xx_t\cdot\uu_t) \leq \xbar\cdot\uu$.
  Consequently,
  $f(\xx_t)\rightarrow\fext(\xbar)$,
  $\fext(\xbar)=\xbar\cdot\uu - \fstar(\uu)$,
  and
  if $\fstar(\uu)\in\R$, then $\uu\in\asubdiffext{\xbar}$.
\end{theorem}

\begin{proof}
~

\begin{proof-parts}
\pfpart{Condition~(\ref{thm:cont-subgrad-converg:a}):}
Suppose this condition holds.
For each $t$,
by
\Cref{pr:asubdiffext-at-x-in-rn}(\ref{pr:asubdiffext-at-x-in-rn:b}),
$\uu_t\in\partial f(\xx_t)=\asubdiffext{\xx_t}$,
and $\fext(\xx_t)=f(\xx_t)$.
Also,
$f(\xx_t)>-\infty$ since $f$ is proper, and
$f(\xx_t)<+\infty$
by \Cref{roc:thm23.4}(\ref{roc:thm23.4:b})
since $\uu_t\in\partial f(\xx_t)$;
thus, $f(\xx_t)\in\R$.

Note that
$\xbar\cdot\uu=\ebar\cdot\uu\plusl\qq\cdot\uu=\qq\cdot\uu$,
since if $\ebar\cdot\uu\neq 0$, then $\xbar\cdot\uu$ would be infinite
(by \Cref{pr:icon-equiv}\ref{pr:icon-equiv:a}\ref{pr:icon-equiv:b}),
contradicting our assumption that it is finite.

Next we have, for all $t$,
\[
   f(\xx_t)
   \geq
   f(\xx_t+\ww_t)
   \geq
   f(\xx_t) + \ww_t\cdot\uu_t.
\]
The first inequality is because $\ww_t\in\resc{f}$, and the second is
because $\uu_t\in\partial f(\xx_t)$.
Since $f(\xx_t)\in\R$, this implies $\ww_t\cdot\uu_t\leq 0$.

Thus,
$
  \xx_t\cdot\uu_t
  =
  \ww_t\cdot\uu_t + \qq_t\cdot\uu_t
  \leq
  \qq_t\cdot\uu_t
$.
Since $\qq_t\cdot\uu_t \rightarrow \qq\cdot\uu = \xbar\cdot\uu$,
it follows that $\limsup(\xx_t\cdot\uu_t) \leq \xbar\cdot\uu$,
as claimed.

Hence, by \Cref{cor:xuleqxbaru-then-subgrad},
  $f(\xx_t)=\fext(\xx_t)\rightarrow\fext(\xbar)$,
  $\fext(\xbar)=\xbar\cdot\uu - \fstar(\uu)$,
  and
  if $\fstar(\uu)\in\R$ then $\uu\in\asubdiffext{\xbar}$.

\pfpart{Condition~(\ref{thm:cont-subgrad-converg:b}):}
Under this condition, $\xbar=(\VV\omm\plusl\rr)\plusl\qq'$
for some $\rr\in\resc{f}$ and $\qq'\in\Rn$,
and
for some $\VV\in\Rnk$ with each column in $\resc{f}$.
More simply, $\xbar=\VV\omm\plusl\qq$ where $\qq=\rr+\qq'$.

Since $\xx_t\rightarrow\xbar$,
by \Cref{thm:seq-rep-not-lin-ind},
there exist sequences $\seq{\bb_t}$ in $\Rstrictpos^k$ and $\seq{\qq_t}$ in
$\Rn$ such that $\xx_t=\VV\bb_t+\qq_t$ for all $t$,
and with
$\qq_t\rightarrow\qq$,
and
$b_{t,i}\rightarrow+\infty$ for $i=1,\ldots,k$.
Then
$\VV\bb_t\in\resc{f}$ since each column of $\VV$ is in
$\resc{f}$ and since $\resc{f}$ is a convex cone
(\Cref{pr:resc-cone-basic-props}).
Therefore, condition~(\ref{thm:cont-subgrad-converg:a}) is satisfied
(with $\ww_t=\VV\bb_t$ and $\ebar=\VV\omm$),
implying the claim, as argued above.

\pfpart{Condition~(\ref{thm:cont-subgrad-converg:c}):}
Suppose this condition holds.
Let $h=\lsc f$, which implies $h=\cl f$, since $f$ is convex and
proper.
Then for each $t$, because $\uu_t\in\partial f(\xx_t)$, it also holds that
$h(\xx_t)=f(\xx_t)$, and that $\uu_t\in\partial h(\xx_t)$
(by \Cref{pr:stan-subgrad-equiv-props}\ref{pr:stan-subgrad-equiv-props:a}\ref{pr:stan-subgrad-equiv-props:d}).

Also,
$\hstar=\fstar$ by \Cref{pr:conj-props}(\ref{pr:conj-props:e}),
and
$\hext=\fext$ by \Cref{pr:h:1}(\ref{pr:h:1aa}).
So in particular, $\contsetf=\contset{\hext}$.
Thus,
\[
  \xbar
  \in
  \contsetf
  =
  \contset{\hext}
  =
  \represc{h}\plusl\intdom{h}
  \subseteq
  \represc{h}\plusl\Rn,
\]
with the second equality following from
\Cref{cor:cont-gen-char}
(applied to $h$).
Therefore, condition~(\ref{thm:cont-subgrad-converg:b}) is satisfied
(with $f$ replaced by $h$).
As argued above, this implies
$\limsup(\xx_t\cdot\uu_t)\leq\xbar\cdot\uu$, and that
  $f(\xx_t)=h(\xx_t)\rightarrow\hext(\xbar)=\fext(\xbar)$,
  $\fext(\xbar)=\hext(\xbar)=\xbar\cdot\uu - \hstar(\uu)=\xbar\cdot\uu-\fstar(\uu)$,
  and
  if $\fstar(\uu)=\hstar(\uu)\in\R$ then
  $\uu\in\asubdifhext{\xbar}=\asubdiffext{\xbar}$.%
\qedhere
\end{proof-parts}
\end{proof}

\indexg{subgradient sequences (standard)!converging to 0@converging to $\zero$|(}%
\indexg{continuity of extensions!minimizing sequences and|(}%
In particular, when $\uu=\zero$,
\Cref{thm:cont-subgrad-converg} implies that
for a proper convex function~$f$,
if $\xx_t\rightarrow\xbar$,
with one of the three conditions given in the theorem satisfied
(for instance, $\xbar\in\contsetf$),
then it is sufficient for a sequence of subgradients
$\uu_t\in\partial f(\xx_t)$ to converge to $\zero$
to guarantee that $\xbar$ minimizes $\fext$,
and that $f(\xx_t)\rightarrow\inf f$.

Further, as shown next, even if the sequence $\seq{\xx_t}$
does not converge in $\eRn$,
if the subgradients $\uu_t\in\partial f(\xx_t)$ are converging to $\zero$,
then the function values must be converging to the infimum of $f$,
that is, $f(\xx_t)\rightarrow\inf f$,
under conditions similar to those in
\Cref{thm:cont-subgrad-converg},
but now assumed to hold for any
subsequential limit of the sequence
(that is, for the limit of any convergent subsequence).
As discussed earlier
(Example~\ref{ex:x1sq-over-x2:grad}),
in general,
such convergence to $f$'s infimum
cannot be guaranteed
without additional conditions beyond convergence
of the subgradients to $\zero$.

\begin{theorem}  \label{thm:subdiff-min-cont}
  Let $f:\Rn\rightarrow\Rext$ be convex and proper.
  Let $\seq{\xx_t}$ and $\seq{\uu_t}$ be sequences in $\Rn$
  with each $\uu_t\in\partial f(\xx_t)$.
  Assume $\uu_t\rightarrow\zero$,
  and also that every subsequential limit of
  $\seq{\xx_t}$ is either in $\represc{f}\plusl\Rn$
  or in $\contsetf$
  (as will be the case, for instance, if
  $\limsup f(\xx_t)<+\infty$ and $\fext$ is continuous everywhere).
  Then:
  \begin{letter-compact}
  \item    \label{thm:subdiff-min-cont:a}
    $f(\xx_t)\rightarrow \inf f$.
  \item    \label{thm:subdiff-min-cont:b}
    Every subsequential limit of $\seq{\xx_t}$ minimizes $\fext$.
    In particular, if $\xx_t\rightarrow\xbar$ for some
    $\xbar\in\extspace$, then $\fext(\xbar)=\inf f$.
  \end{letter-compact}
\end{theorem}

\begin{proof}
Suppose first that some subsequence $\seq{\xx_{s(t)}}$, with
indices $s(1)<s(2)<\dotsb$, converges to some point
$\xbar\in\extspace$.
By assumption, $\xbar$ is either in $\represc{f}\plusl\Rn$
or $\contsetf$.
In either case,
applying
\Cref{thm:cont-subgrad-converg}(\ref{thm:cont-subgrad-converg:b},\ref{thm:cont-subgrad-converg:c})
with $\uu=\zero$ therefore yields that
$f(\xx_t)\rightarrow\fext(\xbar)$, and that
$\fext(\xbar)=-\fstar(\zero)=\inf f$, so $\xbar$ minimizes $\fext$ (using \Cref{pr:fext-min-exists}).
This proves 
part~(\ref{thm:subdiff-min-cont:b}).

For
part~(\ref{thm:subdiff-min-cont:a}),  
suppose, contrary to the theorem's claim, that
$f(\xx_t)\not\rightarrow\inf f$.
Then there exists $\beta\in\R$ with $\beta>\inf f$ such that
$f(\xx_t)\geq \beta$ for infinitely many values of $t$.
By discarding all other sequence elements, we can assume that this
holds for all $t$.

By sequential compactness, the sequence $\seq{\xx_t}$
must have a subsequence converging to some point
$\xbar\in\extspace$.
By again discarding all elements not in this subsequence, we can
assume $\xx_t\rightarrow\xbar$.
Part~(\ref{thm:subdiff-min-cont:b}),
proved above, then implies that
$f(\xx_t)\rightarrow\inf f$,
a contradiction since $f(\xx_t)\geq\beta>\inf f$
for all $t$.
We conclude that $f(\xx_t)\rightarrow \inf f$.

Finally, for an arbitrary sequence $\seq{\xx_t}$ in $\Rn$,
suppose $\limsup f(\xx_t)<+\infty$ and that $\fext$ is
continuous everywhere.
We argue, under these assumptions,
that every subsequential limit of
$\seq{\xx_t}$ is in $\contsetf$.
Since $\limsup f(\xx_t)<+\infty$, there exists $\beta\in\R$ such
that $f(\xx_t)\leq \beta$ for all $t$ sufficiently large.
Therefore, if $\xbar\in\extspace$ is the limit of some convergent
subsequence, then we must have $\fext(\xbar)\leq\beta<+\infty$.
Further, since $\fext$ is continuous everywhere, it must be continuous
at $\xbar$, implying $\xbar\in\contsetf$.%
\indexg{subgradient sequences (standard)!convergence of function values|)}%
\indexg{subgradient sequences (standard)!continuity and|)}%
\indexg{continuity of extensions!subgradient sequences and|)}%
\indexg{continuity of extensions!minimizing sequences and|)}%
\indexg{subgradient sequences (standard)!converging to 0@converging to $\zero$|)}%
\indexg{minimizing sequences!vanishing subgradients@with vanishing subgradients|)}%
\end{proof}

\indexg{subgradient sequences, astral dual!convergence to subgradient|(}%
Having explored astral primal subdifferentials,
we next turn our attention to astral dual subdifferentials
and their continuity properties.
The following theorem
is somewhat analogous to \Cref{cor:xuleqxbaru-then-subgrad}
(as well as \Cref{roc:thm24.4}).
In this setting, we suppose that each $\xbar_t$ is an astral dual
subgradient of some function $\psi:\Rn\rightarrow\Rext$ at $\uu_t$,
and we aim to prove that $\xbar$, the limit of the $\xbar_t\negKern$'s, is
likewise an astral dual subgradient of~$\uu$, the limit of the
$\uu_t\negKern$'s.
Note that this setting includes the case that each $\xbar_t=\xx_t\in\Rn$, and
that $\xx_t$ is a standard subgradient of $\psi$ at $\uu_t$
(by \Cref{pr:adsubdif-int-rn}).

\begin{theorem}   \label{thm:dual-subgrad-cont}
  Let $\psi:\Rn\rightarrow\Rext$.
  Let
  $\seq{\uu_t}$ be a sequence in $\Rn$
  and $\seq{\xbar_t}$ a sequence in $\extspace$
  with each $\xbar_t\in\adsubdifpsi{\uu_t}$, and with
  $\uu_t\rightarrow\uu$
  and
  $\xbar_t\rightarrow\xbar$
  for some
  $\uu\in\Rn$
  and
  $\xbar\in\extspace$.
  Assume $\psi$ is lower semicontinuous at $\uu$,
  $\xbar\cdot\uu\in\R$, and
  that $\liminf(\xbar_t\cdot\uu_t)\leq \xbar\cdot\uu$.
  Then
  $\xbar\in\adsubdifpsi{\uu}$.
\end{theorem}

\begin{proof}
Let $\alpha=\liminf(\xbar_t\cdot\uu_t)$; by assumption,
$\alpha\leq\xbar\cdot\uu<+\infty$.
Then there exists a subsequence of $\seq{\xbar_t\cdot\uu_t}$
converging to $\alpha$.
By discarding all elements of $\seq{\xbar_t}$ and $\seq{\uu_t}$ not
corresponding to this subsequence, we can assume henceforth that
$\xbar_t\cdot\uu_t\rightarrow\alpha$.
Furthermore, since $\alpha<+\infty$, $\xbar_t\cdot\uu_t$ can only be
equal to $+\infty$ for finitely many $t$; by again discarding these,
we can assume henceforth that $\xbar_t\cdot\uu_t<+\infty$ for all $t$.

Let $\uu'\in\Rn$.
We need to
show that $\psi(\uu')\ge\psi(\uu)\plusd\xbar\cdot(\uu'-\uu)$.
A key step is the following claim:

\begin{claimpx}   \label{cl:thm:dual-subgrad-cont:1}
  $\liminf [\xbar_t\cdot(\uu'-\uu_t)]\geq\xbar\cdot(\uu'-\uu)$.
\end{claimpx}

\begin{proofx}
If $\xbar\cdot\uu'=-\infty$, then
$\xbar\cdot(\uu'-\uu)=\xbar\cdot\uu'-\xbar\cdot\uu=-\infty$
by
\Cref{pr:i:1} since $\xbar\cdot\uu\in\R$, so the claim holds
trivially in this case.
Therefore,
in the remainder, assume $\xbar\cdot\uu'>-\infty$.

We have $\xbar_t\cdot\uu'\rightarrow\xbar\cdot\uu'$ by
\Cref{thm:i:1}(\ref{thm:i:1c}),
and $\xbar_t\cdot\uu_t\rightarrow\alpha$.
Further, $\xbar\cdot\uu'$ and $-\alpha$ are summable since neither is
$-\infty$.
Therefore,
\begin{equation}   \label{eq:thm:dual-subgrad-cont:2}
  \xbar_t\cdot\uu' - \xbar_t\cdot\uu_t
  \rightarrow
  \xbar\cdot\uu' - \alpha
  \geq
  \xbar\cdot\uu' - \xbar\cdot\uu
  =
  \xbar\cdot(\uu'-\uu).
\end{equation}
The equality is by 
\Cref{pr:i:1} (since $\xbar\cdot\uu\in\R$).
The convergence is by 
\Cref{prop:lim:eR}(\ref{i:lim:eR:sum}), which also shows that
$\xbar_t\cdot\uu'$ and $-\xbar_t\cdot\uu_t$ are summable for all $t$
sufficiently large.
This last fact implies further that
$\xbar_t\cdot(\uu'-\uu_t) = \xbar_t\cdot\uu' - \xbar_t\cdot\uu_t$,
again by \Cref{pr:i:1},
for all $t$ sufficiently large.
Combined with
\eqref{eq:thm:dual-subgrad-cont:2}, this completes the proof.
\end{proofx}

With this claim, we now have
\begin{align*}
  \psi(\uu')
  &\geq
  \liminf \bigBracks{
    \psi(\uu_t)
    \plusd
    \xbar_t\cdot(\uu'-\uu_t)
  }
  \\
  &\geq
  \liminf \psi(\uu_t)
  \plusd
  \liminf\bigBracks{\xbar_t\cdot(\uu'-\uu_t)}
  \\
  &\geq
  \psi(\uu)
  \plusd
  \xbar\cdot(\uu'-\uu).
\end{align*}  
The first inequality is because
$\psi(\uu')\geq\psi(\uu_t)\plusd\xbar_t\cdot(\uu'-\uu_t)$
for all $t$
by definition of astral dual subgradient (Eq.~\ref{eqn:psi-subgrad:3-alt})
since $\xbar_t\in\adsubdifpsi{\uu_t}$.
The second inequality is by
\Cref{prop:lim:eR}(\ref{i:liminf:eR:sum}).
The last inequality is because $\psi$ is
lower semicontinuous at~$\uu$,
meaning $\liminf\psi(\uu_t)\geq\psi(\uu)$,
and by
\Cref{cl:thm:dual-subgrad-cont:1}.
Thus, $\xbar\in\adsubdifpsi{\uu}$.%
\indexg{continuity!astral subgradients@of astral subgradients|)}%
\indexg{subgradient sequences, astral dual!convergence to subgradient|)}%
\end{proof}

\section{Existence of convergent subgradient sequences}

In the previous section,
we studied conditions under which,
for a given $f:\Rn\to\eR$,
a sequence of pairs $\rpair{\xx_t}{\uu_t}$, with
$\uu_t\in\partial f(\xx_t)$,
converges to a pair $\rpair{\xbar}{\uu}$ in $\eRn\times\Rn$
with $\uu\in\asubdiffext{\xbar}$.
A key application of these results is for the case when
the subgradients $\uu_t$ converge to $\zero$, so that the
pairs
$\rpair{\xx_t}{\uu_t}$ converge to $\rpair{\xbar}{\zero}$
with $\zero\in\asubdiffext{\xbar}$,
implying that $\xbar$ is a minimizer of $\fext$
(\Cref{pr:asub-zero-is-min}\ref{pr:asub-zero-is-min:a}).
Thus, convergence conditions from the previous section
can be used to prove convergence guarantees for minimization algorithms
(as we will do in \Cref{sec:iterative}).

In this section, we study a complementary question. Given a
convex function $f$, is it possible to construct
a sequence of pairs $\rpair{\xx_t}{\uu_t}$, with
$\uu_t\in\partial f(\xx_t)$, that 
converges to a pair $\rpair{\xbar}{\zero}$
for some $\xbar\in\eRn$
with $\zero\in\asubdiffext{\xbar}$?
In other words, is
it possible to find a minimizer of $\ef$ by some algorithm that drives
the subgradients of $f$ to zero?\looseness=-1

Similar to the previous section,
we view this question of convergence to a point with subgradient $\zero$
as a special case of a more general question of convergence to a point
with subgradient equal to any value $\uu\in\Rn$.
\indexg{subgradient sequences, astral dual!existence of convergent|(}%
To answer this, we begin with a dual result
which, for a convex
function $\psi:\Rn\to\eR$, shows
under suitable conditions that
if $\uu\in\cl(\dom\psi)$, then
there must exist a sequence of pairs $\rpair{\uu_t}{\xx_t}$, with
$\xx_t\in\partial\psi(\uu_t)$, converging to a pair
$\rpair{\uu}{\xbar}$, for some $\xbar\in\eRn$ with
$\xbar\in\adsubdifpsi{\uu}$.

\begin{theorem}   \label{thm:ast-dual-subgrad-limit}
  Let $\psi:\Rn\rightarrow\Rext$ be convex with $\psi\not\equiv+\infty$.
  Let $\uu\in\Rn$, and assume $\psi$ is lower semicontinuous at $\uu$.
  Then the following are equivalent:
  \begin{letter-compact}
  \item    \label{thm:ast-dual-subgrad-limit:a}
    $\uu\in\cl(\dom\psi)$.
  \item    \label{thm:ast-dual-subgrad-limit:b}
    There exist sequences
    $\seq{\uu_t}$ and $\seq{\xx_t}$ in $\Rn$
    with $\xx_t\in\partial \psi(\uu_t)$ such that
    $\rpair{\xx_t}{\uu_t}\rightarrow\rpair{\xbar}{\uu}$
    for some $\xbar\in\adsubdifpsi{\uu}$.
  \end{letter-compact}
\end{theorem}

\begin{proof}
~

\begin{proof-parts}
\pfpart{%
  (\ref{thm:ast-dual-subgrad-limit:b})
  $\Rightarrow$
  (\ref{thm:ast-dual-subgrad-limit:a}):
}
Suppose sequences as given in
part~(\ref{thm:ast-dual-subgrad-limit:b}) exist.
Since $\psi\not\equiv+\infty$, there exists a point $\uu'\in\dom\psi$.
For each $t$, since $\xx_t\in\partial\psi(\uu_t)$,
\[
  +\infty
  >
  \psi(\uu')
  \geq
  \psi(\uu_t) + \xx_t\cdot(\uu'-\uu_t),
\]
implying $\uu_t\in\dom\psi$.
Since $\uu_t\rightarrow\uu$, it follows that
$\uu\in\cl(\dom\psi)$.

\pfpart{%
  (\ref{thm:ast-dual-subgrad-limit:a})
  $\Rightarrow$
  (\ref{thm:ast-dual-subgrad-limit:b}):
}
Suppose $\uu\in\cl(\dom\psi)$.

We consider first the case that $\psi$ is improper.
Since $\uu\in\cl(\dom\psi)=\cl(\ri(\dom\psi))$
(with equality by \Cref{pr:ri-props}\ref{pr:ri-props:roc-thm6.3}),
there exists a
sequence $\seq{\uu_t}$ in $\ri(\dom\psi)$ with $\uu_t\rightarrow\uu$.
Since $\psi$ is improper, this implies $\psi(\uu_t)=-\infty$ for all
$t$
(by \Cref{pr:improper-vals}\ref{pr:improper-vals:thm7.2}),
and so that $\zero\in\partial\psi(\uu_t)$.
Further, since $\psi$ is lower semicontinuous at $\uu$,
$\psi(\uu)\leq\liminf\psi(\uu_t)=-\infty$, so $\psi(\uu)=-\infty$ and
$\zero\in\adsubdifpsi{\uu}$.
Thus, setting $\xbar=\zero$ and $\xx_t=\zero$ for all $t$,
this proves the claim in this case.

Henceforth, we assume that $\psi$ is proper.
Let $\ww$ be any point in $\ri(\dom\psi)$ (which is nonempty since
$\psi$ is convex and proper).
Let $\seq{\lambda_t}$ be any decreasing sequence in $(0,1)$ converging
to $0$ (so that $\lambda_{t+1}<\lambda_t$ for all $t$, and
$\lambda_t\rightarrow 0$).
For each $t$, let
$\uu_t=(1-\lambda_t)\uu+\lambda_t\ww$, which
converges to $\uu$.
Since $\uu\in\cl(\dom\psi)$ and
$\ww\in\ri(\dom\psi)$,
each $\uu_t$ is in $\ri(\dom\psi)$
(\Cref{roc:thm6.1}).
Therefore,
by \Cref{roc:thm23.4}(\ref{roc:thm23.4:a}),
$\psi$ has a subgradient $\xx_t$ at each point $\uu_t$ so
that $\xx_t\in\partial\psi(\uu_t)$.

By sequential compactness, the sequence $\seq{\xx_t}$ has a
subsequence converging to some point $\xbar\in\extspace$.
By discarding all other elements (and the corresponding
elements of the $\seq{\lambda_t}$ and $\seq{\uu_t}$ sequences),
we can assume henceforth that the entire sequence converges
so that $\xx_t\rightarrow\xbar$.

By monotonicity of subgradients, for all $t$,
\[
   0
   \leq
   (\xx_{t+1}-\xx_t)\cdot(\uu_{t+1}-\uu_t)
   =
   (\lambda_{t+1} - \lambda_t)(\xx_{t+1}-\xx_t)\cdot(\ww-\uu).
\]
The inequality follows from
\Cref{thm:ast-dual-subgrad-monotone}
(since $\xx_t\in\partial\psi(\uu_t)\subseteq\adsubdifpsi{\uu_t}$ and
$\psi(\uu_t)\in\R$ for all $t$), and also using
\Cref{pr:adsubdif-int-rn}.
The equality is from the definition of each~$\uu_t$.
Since $\lambda_{t+1}<\lambda_t$, it follows that
$\xx_{t+1}\cdot(\ww-\uu)\leq\xx_t\cdot(\ww-\uu)$.
Consequently,
for all $t$,
\begin{equation}   \label{eq:thm:ast-dual-subgrad-limit:1}
  \xx_t\cdot(\ww-\uu)\leq\xx_1\cdot(\ww-\uu).
\end{equation}

Let $\uu'\in\Rn$.
Then for all $t$,
\begin{align}
\notag
  \psi(\uu')
  &\geq
  \psi(\uu_t) + \xx_t \cdot (\uu'-\uu_t)
\\
\notag
  &=
  \psi(\uu_t) + \xx_t\cdot(\uu'-\uu) - \lambda_t \xx_t\cdot(\ww-\uu)
\\
\label{eq:thm:ast-dual-subgrad-limit:3}
  &\geq
  \psi(\uu_t) + \xx_t\cdot(\uu'-\uu) - \lambda_t \xx_1\cdot(\ww-\uu).
\end{align}
The first inequality is because $\xx_t\in\partial\psi(\uu_t)$.
The equality is from the definition of~$\uu_t$.
The last inequality is from
\eqref{eq:thm:ast-dual-subgrad-limit:1} and since $\lambda_t>0$.

We claim that
\begin{equation}   \label{eq:thm:ast-dual-subgrad-limit:2}
  \psi(\uu') \geq \psi(\uu) \plusd \xbar\cdot(\uu'-\uu).
\end{equation}
This is immediate if $\xbar\cdot(\uu'-\uu)=-\infty$, so we assume
henceforth that $\xbar\cdot(\uu'-\uu)>-\infty$.

Since \eqref{eq:thm:ast-dual-subgrad-limit:3} holds for all $t$, we
can take limits, yielding
\begin{align*}
  \psi(\uu')
  &\geq
  \liminf \bigBracks{
    \psi(\uu_t) + \xx_t\cdot(\uu'-\uu) - \lambda_t \xx_1\cdot(\ww-\uu)
  }
\\
  &\geq
  \liminf \psi(\uu_t)
  \plusd
  \liminf \bigBracks{
    \xx_t\cdot(\uu'-\uu) - \lambda_t \xx_1\cdot(\ww-\uu)
  }
\\
  &\geq
  \psi(\uu) + \xbar\cdot(\uu'-\uu).
\end{align*}
The second inequality is by superadditivity of $\liminf$
(\Cref{prop:lim:eR}\ref{i:liminf:eR:sum}).
The third inequality is because $\psi$ is assumed to be lower
semicontinuous at $\uu$, and because
$\xx_t\cdot(\uu'-\uu)\rightarrow\xbar\cdot(\uu'-\uu)$
(\Cref{thm:i:1}\ref{thm:i:1c}),
and $\lambda_t\rightarrow 0$. The summability in the last line
follows because $\xbar\cdot(\uu'-\uu)>-\infty$ and also
$\psi(\uu)>-\infty$ since $\psi$ is proper.

Thus, \eqref{eq:thm:ast-dual-subgrad-limit:2} holds for all
$\uu'\in\Rn$,  so $\xbar\in\adsubdifpsi{\uu}$
(by Eq.~\ref{eqn:psi-subgrad:3-alt}).%
\indexg{subgradient sequences, astral dual!existence of convergent|)}%
\qedhere
\end{proof-parts}
\end{proof}

\indexg{subgradient sequences (standard)!existence of convergent|(}%
Returning to primal subgradients, we next apply
\Cref{thm:ast-dual-subgrad-limit} to the conjugate of any
convex, closed, proper function $f:\Rn\rightarrow\Rext$
to prove the existence of subgradient sequences
with particular convergence properties.
Specifically,
the next theorem shows that for every $\uu\in\cl(\dom\fstar)$,
there exists a sequence of pairs $\rpair{\xx_t}{\uu_t}$, with
${\uu_t\in\partial f(\xx_t)}$, converging to a pair
$\rpair{\xbar}{\uu}$
with the property that
$\fminusuext(\xbar) = -\fstar(\uu)$, so that
the variant of the Fenchel-Young inequality
from \eqref{eq:fstar-leq-fminusuext} holds with equality.
In particular, when $\uu=\zero$, this equality means that $\xbar$
minimizes $\fext$.
We also show in
\Cref{cor:ast-subgrad-limit} that
under a stronger
assumption, namely, that $\uu\in\dom\fstar$, the limit pair
$\rpair{\xbar}{\uu}$
additionally satisfies $\uu\in\asubdiffext{\xbar}$.
Thus, every $\uu\in\dom\fstar$
is an astral subgradient of $\ef$ at some point $\xbar$ such that the
pair $\rpair{\xbar}{\uu}$ is the limit of some sequence of pairs
$\rpair{\xx_t}{\uu_t}$ where $\uu_t$ is a standard subgradient of $f$
at $\xx_t$.

\begin{theorem}   \label{thm:ast-subgrad-fench-lim}
  Let $f:\Rn\rightarrow\Rext$ be convex, closed and proper.
  Let $\uu\in\Rn$.
  Then the following are equivalent:
  \begin{letter-compact}
  \item   \label{thm:ast-subgrad-fench-lim:a}
    $\uu\in\cl(\dom\fstar)$.
  \item   \label{thm:ast-subgrad-fench-lim:b}
    There exist sequences $\seq{\xx_t}$ and $\seq{\uu_t}$ in $\Rn$
    with $\xx_t\in\partial f(\xx_t)$
    such that
    $\rpair{\xx_t}{\uu_t}\rightarrow\rpair{\xbar}{\uu}$
    for some $\xbar\in\extspace$
    which satisfies
    $\fminusuext(\xbar) = -\fstar(\uu)$.
\end{letter-compact}
\end{theorem}

\begin{proof}
~

\begin{proof-parts}
\pfpart{%
  (\ref{thm:ast-subgrad-fench-lim:b})
  $\Rightarrow$
  (\ref{thm:ast-subgrad-fench-lim:a}):
}
Suppose sequences as in
statement~(\ref{thm:ast-subgrad-fench-lim:b}) exist.
Then for each $t$, $\xx_t\in\partial\fstar(\uu_t)$,
by \Cref{pr:stan-subgrad-equiv-props}(\ref{pr:stan-subgrad-equiv-props:a},\ref{pr:stan-subgrad-equiv-props:c}),
implying
$\uu_t\in\dom{\fstar}$
by \Cref{roc:thm23.4}(\ref{roc:thm23.4:b}).
Therefore, $\uu\in\cl(\dom{\fstar})$.

\pfpart{%
  (\ref{thm:ast-subgrad-fench-lim:a})
  $\Rightarrow$
  (\ref{thm:ast-subgrad-fench-lim:b}):
}
Suppose $\uu\in\cl(\dom{\fstar})$.
Note that $\fstar$ is convex, closed and proper since $f$ is
(Propositions~\ref{pr:conj-props}\ref{pr:conj-props:d}
and~\ref{pr:conj-props-cvx}\ref{pr:conj-props-cvx:a}).
We therefore can apply \Cref{thm:ast-dual-subgrad-limit} (with $\psi=\fstar$),
yielding sequences $\seq{\xx_t}$ and $\seq{\uu_t}$ in $\Rn$
and $\xbar\in\extspace$ such that
$\xx_t\rightarrow\xbar$, $\uu_t\rightarrow\uu$,
$\xx_t\in\partial \fstar(\uu_t)$ for all $t$,
and
$\xbar\in\adsubdiffstar{\uu}$.
This implies, for all $t$, that $\uu_t\in\partial f(\xx_t)$
by \Cref{pr:stan-subgrad-equiv-props}(\ref{pr:stan-subgrad-equiv-props:a},\ref{pr:stan-subgrad-equiv-props:c}),
and also
that $\xx_t\in\dom{f}$
by \Cref{roc:thm23.4}(\ref{roc:thm23.4:b}).
Thus, $\xbar\in\cldom{f}$.
Since also $\xbar\in\adsubdiffstar{\uu}$, it hence follows that
$\fminusuext(\xbar) = -\fstar(\uu)$
by \Cref{thm:dual-subgrad-fenchel-tilt}.
\qedhere
\end{proof-parts}
\end{proof}

\begin{theorem}   \label{cor:ast-subgrad-limit}
  Let $f:\Rn\rightarrow\Rext$ be convex, closed and proper.
  Let $\uu\in\Rn$.
  Then the following are equivalent:
  \begin{letter-compact}
  \item   \label{cor:ast-subgrad-limit:a}
    $\uu\in\dom\fstar$.
  \item   \label{cor:ast-subgrad-limit:b}
    There exist sequences $\seq{\xx_t}$ and $\seq{\uu_t}$ in $\Rn$
    with $\xx_t\in\partial f(\xx_t)$
    such that
    $\rpair{\xx_t}{\uu_t}\rightarrow\rpair{\xbar}{\uu}$
    for some $\xbar\in\extspace$
    with
    $\uu\in\asubdiffext{\xbar}$.
  \end{letter-compact}
\end{theorem}

\begin{proof}
~

\begin{proof-parts}
\pfpart{%
  (\ref{cor:ast-subgrad-limit:b})
  $\Rightarrow$
  (\ref{cor:ast-subgrad-limit:a}):
}
If statement~(\ref{cor:ast-subgrad-limit:b}) holds, then
$\uu\in\asubdiffext{\xbar}$, implying $\uu\in\dom\fstar$ by
\Cref{pr:subgrad-imp-in-cldom}(\ref{pr:subgrad-imp-in-cldom:c}).

\pfpart{%
  (\ref{cor:ast-subgrad-limit:a})
  $\Rightarrow$
  (\ref{cor:ast-subgrad-limit:b}):
}
Suppose $\uu\in\dom{\fstar}$, implying $\fstar(\uu)\in\R$ since
$\fstar$ is proper
(\Cref{pr:conj-props-cvx}\ref{pr:conj-props-cvx:a}).
Then by \Cref{thm:ast-subgrad-fench-lim},
there exist sequences $\seq{\xx_t}$ and $\seq{\uu_t}$ in $\Rn$
and $\xbar\in\extspace$ such that
$\xx_t\rightarrow\xbar$, $\uu_t\rightarrow\uu$,
$\uu_t\in\partial f(\xx_t)$ for all $t$,
and
$\fminusuext(\xbar) = -\fstar(\uu)$.
Since also $\fstar(\uu)\in\R$, this implies
$\uu\in\asubdiffext{\xbar}$ by
\Cref{thm:fminus-subgrad-char}(\ref{thm:fminus-subgrad-char:d},\ref{thm:fminus-subgrad-char:a}).
\qedhere
\end{proof-parts}
\end{proof}

\indexg{subgradient sequences (standard)!converging to 0@converging to $\zero$|(}%
\indexg{minimizing sequences!existence of|(}%
In the special case $\uu=\zero$, 
\Cref{thm:ast-subgrad-fench-lim} shows that
if $\zero\in\cl(\dom\fstar)$, then there must exist
a sequence $\seq{\xx_t}$ in $\Rn$ converging to some minimizer of
$\fext$ and with subgradients converging to $\zero$.
In particular, this will be the case whenever $f$ is lower-bounded
(so that $\fstar(\zero)=-\inf f<+\infty$).
On the other hand, if $\zero\not\in\cl(\dom\fstar)$, then no such
sequence can exist, which means $\fext$'s minimizers cannot be
attained by attempting to drive $f$'s subgradients to zero.
Here are examples:

\begin{example}
Suppose, for $x\in\R$, that
\[
  f(x) = \begin{cases}
      -2\sqrt{x}    & \mbox{if $x\geq 0$,} \\
        +\infty     & \mbox{otherwise.}
               \end{cases}
\]
Then it can be calculated that, for $u\in\R$,
\[
  \fstar(u) = \begin{cases}
      -1/u    & \mbox{if $u < 0$,} \\
        +\infty     & \mbox{otherwise.}
               \end{cases}
\]
Thus, $0$ is in $\cl(\dom\fstar)$, though not in $\dom\fstar$.
In this case, we can construct sequences as in
\Cref{thm:ast-subgrad-fench-lim}
by letting
$x_t=t$ and $u_t=f'(x_t)=-1/\sqrt{t}$ (where $f'$ is the
derivative of $f$).
Then $u_t\rightarrow 0$ while $x_t\rightarrow+\infty$, which minimizes
$\fext$.

Suppose now instead that $f(x)=-x$ for $x\in\R$.
Then $\fstar=\indf{\{-1\}}$, the indicator function for the single
point $-1$, so $\dom\fstar=\{-1\}$.
Thus, $0\not\in\cl(\dom\fstar)$.
Indeed, for every sequence $\seq{x_t}$ in $\Rn$, we must have
$u_t=f'(x_t)=-1$, so no such sequence of subgradients can converge to
$0$.%
\indexg{subgradient sequences (standard)!existence of convergent|)}%
\indexg{subgradient sequences (standard)!converging to 0@converging to $\zero$|)}%
\indexg{minimizing sequences!existence of|)}%
\end{example}

\section{Convergence of iterative methods}
\label{sec:iterative}

The preceding results,
especially \Cref{thm:subdiff-min-cont},
can be applied to prove the convergence of
iterative methods for minimizing a function, as we now
illustrate.
Let $f:\Rn\rightarrow\Rext$ be convex.
\indexg{minimization methods|(}%
\indexg{algorithms, minimization|(}%
We consider methods that compute a sequence of iterates $\seq{\xx_t}$
in $\Rn$
with the goal of minimizing $f$ in the limit.
\indexg{gradient descent|(}%
A classic example is \emph{gradient descent} in which $\xx_1\in\Rn$ is
arbitrary, and each successive iterate is defined by
\begin{equation}  \label{eqn:grad-desc-defn}
  \xx_{t+1} = \xx_t - \eta_t \gradf(\xx_t),
\end{equation}
for some step size
\indexg{gradient descent|)}%
$\eta_t\in\Rstrictpos$.
Although certainly an important example,
there are many well-studied optimization techniques beyond
gradient descent. For instance, line-search methods descend
along any direction
that has a positive inner product with the negative gradient
\indexa{Nocedal, J.}\indexa{Wright, S. J.}%
\citep[Chapter 3]{nocedal_wright}, and coordinate-descent methods update only one of
the coordinates in each iteration.
Our aim is to develop
techniques that are broadly applicable to this wide range of methods.

In analyzing the convergence of such iterative methods, it is quite
common to assume that $f$ has a finite minimizer in $\Rn$,
and often also that we are effectively searching for a minimizer over
only a compact subset of $\Rn$.
(See, for example,
\idxboydvand\citealp[Chapters~9,~10,~11]{boyd_vandenberghe}.)
Depending on the problem setting, such assumptions may or may not be
reasonable.
A primary purpose of the current work, of course,
has been to develop a foundation that overcomes such difficulties and
that can be
applied without relying on such assumptions.
Indeed, as we have seen,
astral space is itself compact, and
the extension $\fext$ of any convex function $f:\Rn\rightarrow\Rext$
always has a minimizer that is attained
at some astral point in $\extspace$.

\indexg{subgradient descent|(}%
\indexg{minimizing sequences!using subgradient descent|(}%
Before exploring how astral methods can be used to prove general
convergence results, we first give an example
of how standard, non-astral techniques can sometimes be
applied to prove specialized convergence for particular algorithms, in
this case, subgradient descent, a generalization of gradient descent
in which
$\xx_{t+1} = \xx_t - \eta_t \uu_t$,
where $\eta_t\in\Rstrictpos$,
and $\uu_t$ is any subgradient of $f$ at $\xx_t$.
The next
\namecref{thm:grad-desc-converges:alt}
proves convergence to the function's infimum
assuming a particular lower bound on how much $f(\xx_t)$
decreases on each iteration. We will discuss this assumption
further below. The proof of the
\namecref{thm:grad-desc-converges:alt}
is based on a proof
of
\idxziwei\idxmiro\idxmatus\idxschapire%
\citet[Lemma~2]{ziwei_poly_tail_losses}.

\begin{theorem}  \label{thm:grad-desc-converges:alt}
  Let $f:\Rn\rightarrow\Rext$ be convex and proper.
  Let $\seq{\xx_t}$ and $\seq{\uu_t}$ be sequences in $\Rn$,
  where, for all $t$,
  $f(\xx_t)\in\R$,
  $\uu_t\in\partial f(\xx_t)$,
  $\xx_{t+1} = \xx_t - \eta_t \uu_t$,
  and $\eta_t\in\Rstrictpos$.
  Assume that, for all $t$,
  \begin{equation}  \label{eqn:thm:grad-desc-converges:alt:1}
    f(\xx_{t+1}) \leq f(\xx_t) - \frac{\eta_t}{2} \norm{\uu_t}^2,
  \end{equation}
  and that
  $\sum_{t=1}^\infty \eta_t = +\infty$.
  Then $f(\xx_t) \rightarrow \inf f$.
\end{theorem}

\begin{proof}
Our proof relies on the following observation.
For any $\zz\in\Rn$, and for all $t$,
\begin{align*}
  \|\xx_{t+1} - \zz\|^2
  =
  \norm{(\xx_{t} - \zz) - \eta_t \uu_t}^2
  &=
  \|\xx_{t} - \zz\|^2
  - 2 \eta_t \uu_t\cdot (\xx_t - \zz) + \eta_t^2 \norm{\uu_t}^2
  \\
  &\leq
  \|\xx_{t} - \zz\|^2
  - 2\eta_t \bigParens{f(\xx_t) - f(\zz)}
  + \eta_t^2 \norm{\uu_t}^2.
\end{align*}
The first equality is by $\xx_{t+1}\negKern$'s definition.
The inequality is because
$\uu_t\in\partial f(\xx_t)$
(by definition of subgradient, Eq.~\ref{eqn:prelim-standard-subgrad-ineq}).
Rearranging and applying this inequality repeatedly then yields
\begin{align}
\notag
  2\sum_{s=1}^{t} \eta_{s} \bigParens{ f(\xx_{s}) - f(\zz) }
  &\le
  \sum_{s=1}^{t}
  \BigParens{
    \norm{\xx_{s} - \zz}^2
    -
    \norm{\xx_{s+1} - \zz}^2
    +
    \eta_{s}^2 \norm{\uu_{s}}^2
  }
\\[\medskipamount]
\notag
  &=
  \norm{\xx_{1} - \zz}^2
  -
  \norm{\xx_{t+1} - \zz}^2
  +
  \smash{\sum_{s=1}^{t} \eta_{s}^2 \norm{\uu_{s}}^2}
\\[\medskipamount]
\label{eq:sgd:magic:alt}
&\leq
  \norm{\xx_1 - \zz}^2
  +
  \sum_{s=1}^{t} \eta_{s}^2 \norm{\uu_{s}}^2.
\end{align}

To prove the theorem, suppose, by way of contradiction, that
$f(\xx_t)\not\rightarrow \inf f$.
\eqref{eqn:thm:grad-desc-converges:alt:1} implies
that the sequence of values $f(\xx_t)$ is nonincreasing, which means
that they must have a limit, which is also equal to their infimum.
Let
  $ \gamma=\lim f(\xx_t) =\inf \set{f(\xx_t) :\: t\in\nats} $.
By our assumption, $\gamma>\inf f$, which also implies $\gamma\in\R$.
Thus, there exists a point $\zz\in\Rn$ with
$\gamma > f(\zz) \geq \inf f$
(and with $f(\zz)\in\R$ since $f$ is proper).

Thus, for all $t$,
\begin{align*}
  2\BiggParens{\,\sum_{s=1}^t \eta_{s}} \bigParens{\gamma - f(\zz)}
  &\leq
  2 \sum_{s=1}^t \eta_{s} \bigParens{f(\xx_{s+1}) - f(\zz)}
  \\
  &\leq
  2 \sum_{s=1}^t \eta_{s} \bigParens{f(\xx_{s}) - f(\zz)}
  -
  \sum_{s=1}^t \eta_{s}^2 \norm{\uu_{s}}^2
  \leq
  \norm{\xx_1 - \zz}^2.
\end{align*}
The first inequality is
by $\gamma$'s definition.
The second and third inequalities are by
Eqs.~(\ref{eqn:thm:grad-desc-converges:alt:1})
and~(\ref{eq:sgd:magic:alt}), respectively.
The left-hand side of this inequality converges to $+\infty$
as $t\rightarrow+\infty$, since $\gamma>f(\zz)$
and since $\sum_{t=1}^\infty \eta_t = +\infty$.
But this is a contradiction since
the right-hand side is constant and finite.
\end{proof}

\Cref{thm:grad-desc-converges:alt}
proves convergence assuming
a lower bound on how much $f(\xx_t)$ decreases on each iteration
as a function of its gradient.
\indexg{gradient descent|(}%
As an example of when such a bound is possible,
\indexg{beta-smoothness@$\beta$-smoothness|(}%
suppose $f$ is \emph{$\beta$-smooth}, meaning
\begin{equation}  \label{eqn:grad-bnd-smooth}
  f(\xx')
  \leq
  f(\xx) + \gradf(\xx)\cdot(\xx'-\xx)
          + \frac{\beta}{2} \norm{\xx'-\xx}^2
\end{equation}
for all $\xx,\,\xx'\in\Rn$, for some constant $\beta\in\Rstrictpos$
(and assuming $f$ is
\indexg{beta-smoothness@$\beta$-smoothness|)}%
differentiable).
Then if $\xx_{t+1}$ is computed as in \eqref{eqn:grad-desc-defn} with
$\eta_t=1/\beta$, then the smoothness condition in \eqref{eqn:grad-bnd-smooth} implies
\begin{equation}  \label{eqn:grad-iter-bnd-smooth}
  f(\xx_{t+1})
  \leq
  f(\xx_t)
  - \frac{1}{2\beta} \norm{\gradf(\xx_t)}^2.
\end{equation}
Thus, this update is guaranteed to decrease
the function values of the iterates from $f(\xx_t)$ to $f(\xx_{t+1})$
by at least a constant times $\norm{\gradf(\xx_t)}^2$.%
\indexg{algorithms, minimization|)}%
\indexg{minimization methods|)}%
\indexg{minimizing sequences!using subgradient descent|)}%
\indexg{gradient descent|)}%
\indexg{subgradient descent|)}%

Once established, such a guarantee of progress can sometimes
be sufficient to
ensure $f(\xx_t)\rightarrow \inf f$, for instance,
as we just showed in \Cref{thm:grad-desc-converges:alt}.
Intuitively, if $\gradf(\xx_t)$  %
is getting close to $\zero$, then we
should be approaching $f$'s minimum; on the other hand, as long as
$\norm{\gradf(\xx_t)}$ remains large, we are assured of
significant progress in reducing $f(\xx_t)$ on each iteration.
As such, we might expect that a progress guarantee of this kind
should suffice to ensure the convergence of a broad family of methods,
not just (sub)gradient descent.
Indeed, we will see soon that this can be the case.
But on the other hand,
we will also see that such a progress guarantee does not always
ensure convergence to the function's infimum
(just as driving subgradients to zero does not guarantee such
convergence, as we saw in \Cref{ex:x1sq-over-x2:grad}).

\indexg{minimizing sequences!progress guarantee@with progress guarantee|(}%
For the remainder of this section, we apply astral methods
to study in substantially greater generality
when such convergence is ensured for an arbitrary iterative method,
given a generalized form of a lower bound on per-iteration progress in terms
of gradients or subgradients.
In particular,
we will see that continuity in astral space is sufficient to
ensure convergence to a minimum.
We will also see that such a result is not possible, in general,
when astral continuity is not assumed.

Let $\seq{\xx_t}$ be any sequence in $\Rn$ (not necessarily computed
using gradient descent), and let
$\seq{\uu_t}$ in $\Rn$ be a corresponding sequence of subgradients
so that $\uu_t\in\partial f(\xx_t)$ for all $t$.
Generalizing the kind of progress bounds considered above
(such as Eq.~\ref{eqn:thm:grad-desc-converges:alt:1}),
we suppose
that
\begin{equation}  \label{eqn:aux-fcn-prog-bnd}
   f(\xx_{t+1})
   \leq
   f(\xx_t) - \alpha_t h(\uu_t)
\end{equation}
for some $\alpha_t\in\Rpos$, and some function
$h:\Rn\rightarrow\Rpos$.
We assume $h$ satisfies the following two conditions:
\begin{roman-compact}
\item $h(\zero)=0$;
\item $
  \inf \regBraces{
    h(\uu)
    :\:
    \uu\in\Rn,\, \norm{\uu}\geq\epsilon
  }
  > 0
$
for all
$\epsilon\in\Rstrictpos$.
\end{roman-compact}
We refer to such a function as a
\indexg{magnitude function}%
\emph{\auxiliary function}.
Intuitively, if $h(\uu)$ is small, these properties force $\uu$ to be
close to $\zero$.
For example,
\eqref{eqn:grad-iter-bnd-smooth}
satisfies
\eqref{eqn:aux-fcn-prog-bnd}
with $\alpha_t=1/(2\beta)$ and $h(\uu)=\norm{\uu}^2$, which clearly
is a \auxiliary function.
In general, if $h$ is continuous, strictly positive except at $\zero$,
and radially nondecreasing
(meaning, for $\uu\in\Rn$, that $h(\lambda\uu)$ is nondecreasing as a
function of $\lambda\in\Rpos$),
then it must be
a \auxiliary function
(by \Cref{thm:weierstrass}
and compactness of
$\{ \uu\in\Rn :\: \norm{\uu}=\epsilon \}$,
for $\epsilon\in\Rstrictpos$).

\indexg{continuity of extensions!minimizing sequences and|(}%
Given the progress bound in \eqref{eqn:aux-fcn-prog-bnd},
if $\fext$ is continuous
(either everywhere, or just at the subsequential limits of
the sequence of iterates), then
the next theorem shows how we can use our
previous results to prove convergence to $f$'s minimum,
without requiring $f$ to have a finite minimizer, nor the sequence of
iterates $\seq{\xx_t}$ to remain bounded.
Unlike
\Cref{thm:grad-desc-converges:alt},
this theorem can be applied to any sequence $\seq{\xx_t}$, regardless
of how it is computed or constructed (provided, of course, that it
satisfies the stated conditions).
Furthermore, the theorem relies on
an assumed progress bound of a much weaker and more general form.

\begin{theorem}\label{thm:fact:gd}
  Let $f:\Rn\rightarrow\Rext$ be convex.
  Let $\seq{\xx_t}$ and $\seq{\uu_t}$ be sequences in $\Rn$
  with each $\uu_t\in\partial f(\xx_t)$.
  Assume that $\fext$ is continuous at every subsequential limit
  of the sequence $\seq{\xx_t}$.
  Also assume that
  \begin{equation}  \label{thm:fact:gd:1}
    f(\xx_{t+1}) \leq  f(\xx_t) - \alpha_t h(\uu_t)
  \end{equation}
  for each $t$, where $\alpha_t\in\Rpos$ and
  $\sum_{t=1}^\infty \alpha_t = +\infty$,
  and where
  $h:\Rn\rightarrow\Rpos$ is
  a \auxiliary function.
  Then $f(\xx_t)\rightarrow \inf f$.
\end{theorem}

\begin{proof}
We first consider when $f$ is improper.
If $f\equiv+\infty$, then the claim holds trivially.
If $f(\zz)=-\infty$ for some $\zz\in\Rn$, then
$-\infty=f(\zz)\geq f(\xx_t)+\uu_t\cdot(\zz-\xx_t)$
for all $t$,
since $\uu_t\in\partial f(\xx_t)$,
so $f(\xx_t)=-\infty$.
Thus, $f(\xx_t)\rightarrow-\infty=\inf f$,
proving the claim in this case too.

Henceforth, we assume $f$ is proper.
In that case,
by \Cref{roc:thm23.4}(\ref{roc:thm23.4:b}),
we must have $f(\xx_t)\in\R$
for all $t$,
since $\uu_t\in\partial f(\xx_t)$.

Suppose that $\liminf h(\uu_t) = 0$.
Then in this case,
there exists a subsequence $\seq{\xx_{s(t)}}$, with indices
$s(1)<s(2)<\dotsb$, such that $h(\uu_{s(t)})\rightarrow 0$.
We claim further that $\uu_{s(t)}\rightarrow\zero$.
Suppose otherwise.
Then for some $\epsilon\in\Rstrictpos$,
$\norm{\uu_{s(t)}}\geq\epsilon$ for
infinitely many values of $t$.
Let
$ \delta = \inf \{ h(\uu) :\: \uu\in\Rn,\, \norm{\uu}\geq\epsilon \} $.
Then $\delta>0$, since $h$ is
a \auxiliary function.
Thus, $h(\uu_{s(t)})\geq\delta>0$ for infinitely many values of $t$.
But this contradicts that $h(\uu_{s(t)})\rightarrow 0$.

Since $\uu_{s(t)}\rightarrow\zero$,
we can apply \Cref{thm:subdiff-min-cont} to the
extracted subsequences
$\seq{\xx_{s(t)}}$ and $\seq{\uu_{s(t)}}$,
noting that, by assumption, $\fext$ is continuous at
all the subsequential limits of $\seq{\xx_{s(t)}}$,
and that $f(\xx_t)\leq f(\xx_1)<+\infty$ for all $t$ since
the function values $f(\xx_t)$ are nonincreasing
(from Eq.~\ref{thm:fact:gd:1},
since $\alpha_t\geq 0$ and $h\geq 0$).
Together, these imply that
all the subsequential limits of $\seq{\xx_{s(t)}}$
are in $\contsetf$.
\Cref{thm:subdiff-min-cont}(\ref{thm:subdiff-min-cont:a}) thus
yields $f(\xx_{s(t)})\rightarrow \inf f$.
Furthermore, this shows the entire sequence, which is nonincreasing,
converges to $f$'s infimum as well,
so that $f(\xx_t)\rightarrow\inf f$.

In the alternative case, suppose
$\liminf h(\uu_t) > 0$.
Then there exists $t_0\in\nats$ and $\epsilon\in\Rstrictpos$ such that
$h(\uu_t)\geq\epsilon$ for all $t\geq t_0$.
Summing \eqref{thm:fact:gd:1}
yields, for $t > t_0$,
\[
   f(\xx_t)
   \leq
   f(\xx_1)
   -
   \sum_{s=1}^{t-1} \alpha_{s} h(\uu_{s})
   \leq
   f(\xx_1)
   -
   \epsilon \sum_{s=t_0}^{t-1} \alpha_{s}.
\]
As $t\rightarrow+\infty$, the sum on the right converges to $+\infty$
(by assumption, even disregarding finitely many terms), implying
$f(\xx_t)\rightarrow -\infty$.
Thus, $\inf f=-\infty$ and
$f(\xx_t)\rightarrow\inf f$ in this case as well.
\end{proof}

The next example shows that
if we drop the assumption regarding $\fext$'s continuity, then
the convergence proved in \Cref{thm:fact:gd} is no longer
ensured, in general, even given the progress bound in
\eqref{eqn:aux-fcn-prog-bnd}
with the $\alpha_t\negKern$'s all equal to a positive constant,
and even when $h(\uu)=\norm{\uu}^2$, the
most standard case:

\begin{example}   \label{ex:x1sq-over-x2:progress}
Consider again the flattening valley function $f:\R^2\rightarrow\Rext$
from Examples~\ref{ex:x1sq-over-x2}, \ref{ex:x1sq-over-x2:cont}
and~\ref{ex:x1sq-over-x2:grad},
which satisfies
$f(x_1,x_2)=x_1^2/x_2$ if $x_2\ge x_1>0$,
and is also continuously differentiable over this region.
Let $\xx_t=\trans{[t+1,\, t(t+1)]}$.
Then $f(\xx_t)=(t+1)/t$
and
$\gradf(\xx_t)=\trans{[2/t,\, -1/t^2]}$.
It can be checked that, for all $t$,
\[
  f(\xx_{t+1}) - f(\xx_t)
  =
  -\frac{1}{t(t+1)}
  \leq
  -\frac{1}{10}
  \cdot
  \frac{4t^2+1}{t^4}
  =
  -\frac{1}{10}
  \norm{\gradf(\xx_t)}^2.
\]
In other words,
\eqref{eqn:aux-fcn-prog-bnd} is satisfied with $\alpha_t=1/10$, for
all $t$, and $h(\uu)=\norm{\uu}^2$ (and with $\uu_t=\gradf(\xx_t)$).
Thus, all of the conditions of \Cref{thm:fact:gd} are
satisfied, except that $\fext$ is not continuous everywhere, nor at
the limit of this particular sequence
(namely, $\limray{\ee_2}\plusl\limray{\ee_1}$).
And indeed, the theorem's conclusion is false in this case since
$f(\xx_t)\rightarrow 1 > 0 = \inf f$.%
\indexg{minimizing sequences!progress guarantee@with progress guarantee|)}%
\end{example}

This shows that \Cref{thm:fact:gd} is false, in general,
if continuity is not
\indexg{continuity of extensions!minimizing sequences and|)}%
assumed.
Nevertheless, this does not rule out the possibility that particular
algorithms might be effective at minimizing a convex function, even
without a continuity assumption;
indeed, this was shown to be true of subgradient descent in
\Cref{thm:grad-desc-converges:alt}.

Although \Cref{thm:fact:gd} establishes convergence of the function
values $f(\xx_t)$ to $\inf f$,
it does not guarantee convergence of the iterates $\xx_t$ to some minimizer~$\xbar$.
\indexg{implicit regularization|(}%
This
is similar to other results that study minimization
at infinity (under specialized assumptions), such as the results from the
literature on implicit regularization
(%
\idxziwei\idxmatus%
\citealp{riskparam};
\idxschapire\idxfreund\indexa{Bartlett, P.}\indexa{Lee, W. S.}%
\citealp{boosting_margin};
\indexa{Soudry, D.}\indexa{Hoffer, E}\indexa{Nacson, M. S.}\indexa{Gunasekar, S.}\indexa{Srebro, N.}%
\citealp{nati_logistic};
\idxmatus%
\citealp{mjt_margins};
\idxtongbin%
\citealp{zhang_yu_boosting}%
),
which only guarantee convergence to a specific cone or
to a particular dominant direction. Compared with those works, our approach is
much more general since it makes only generic assumptions on the minimized function
and the optimization algorithm.%
\indexg{implicit regularization|)}%

\indexg{algorithms, minimization|(}%
\indexg{minimization methods|(}%
To showcase this generality, we now apply \Cref{thm:fact:gd}
to prove the effectiveness of a number of standard algorithms,
even when no finite minimizer exists.
\indexg{coordinate descent|(}%
\emph{Coordinate descent} is one such method in which, at each
iteration, just one of the coordinates %
of~$\xx_t$ is chosen and
updated, not modifying any of the others; %
thus,
$\xx_{t+1}=\xx_t + \eta_t \ee_{i_t}$ for some standard basis vector $\ee_{i_t}$ and
some $\eta_t\in\R$.
In a \emph{gradient-based} version of coordinate descent, $i_t$ is
chosen to be the largest coordinate of the gradient $\gradf(\xx_t)$.
In a \emph{fully greedy} version, both $i_t$ and $\eta_t$ are chosen to
effect the maximum possible decrease in the function value, that is,
to minimize $f(\xx_t+\eta\ee_i)$ among all choices of
$i\in\{1,\ldots,n\}$ and $\eta\in\R$.
Many other variations are possible.%
\indexg{coordinate descent|)}%

\indexg{steepest descent|(}%
\emph{Steepest descent}
\idxboydvand\citep[Section~9.4]{boyd_vandenberghe}
is a technique generalizing both gradient
descent and one form of coordinate descent.
In this method, for some $p\geq 1$, on each iteration $t$,
$\vv_t\in\Rn$ is chosen to
maximize $\gradf(\xx_t)\cdot\vv_t$ subject to $\normp{\vv_t}\leq 1$,
and then subtracted from $\xx_t$ after scaling by some step size
$\eta_t \in \Rstrictpos$;
thus,
$\xx_{t+1} = \xx_t - \eta_t \vv_t$.
Here,
$\normp{\xx}=\regParens{\sum_{i=1}^n |x_i|^p}^{1/p}$
is the $\ell_p$-norm of $\xx\in\Rn$; we also allow $p=+\infty$,
in which case $\norm{\xx}_{\infty}=\max_{i=1,\ldots,n} |x_i|$.
When $p=2$, we obtain gradient descent;
when $p=1$, we obtain (gradient-based) coordinate descent.
The
next example shows that for suitable step sizes,
steepest descent is guaranteed to minimize
any smooth convex function $f$ whose extension $\ef$ is continuous.

\begin{example}[Convergence of steepest descent]
\label{ex:steepest-descent}
Let $f:\Rn\to\R$ be convex and $\beta$-smooth
(as in Eq.~\ref{eqn:grad-bnd-smooth}),
for some $\beta\in\Rstrictpos$. Suppose $\ef$ is continuous,
and let $p\ge 1$.
We will argue that the
steepest-descent algorithm (for a suitable choice of $\eta_t$) satisfies the progress bound in \eqref{eqn:aux-fcn-prog-bnd}, and so minimizes $f$ in the limit
(by \Cref{thm:fact:gd}).

To start, note that there exists a constant $C_p\in\Rstrictpos$ for which
$\norms{\xx}\leq C_p \normp{\xx}$ for all $\xx\in\Rn$; for instance,
it can be argued that
$\norms{\xx}\le\sqrt{n}\norm{\xx}_{\infty}\le\sqrt{n}\norm{\xx}_p$.
Let $\beta'=C_p^2 \beta$, and let $q\geq 1$ be such that
$1/p + 1/q = 1$ (allowing either $p$ or $q$ to be $+\infty$).
Then for all $t$,
\begin{align*}
  f(\xx_{t+1})
  &\leq
  f(\xx_t)
  + \gradf(\xx_t)\cdot\parens{\xx_{t+1} - \xx_t}
  + \frac{\beta}{2} \norms{\xx_{t+1} - \xx_{t}}^2
  \\
  &\leq
  f(\xx_t)
  + \gradf(\xx_t)\cdot\parens{\xx_{t+1} - \xx_t}
  + \frac{\beta'}{2} \normp{\xx_{t+1} - \xx_{t}}^2
  \\
  &=
  f(\xx_t)
  - \eta_t \gradf(\xx_t)\cdot\vv_t
  + \frac{\beta'\eta_t^2}{2}
  \\
  &=
  f(\xx_t)
  - \eta_t \normq{\gradf(\xx_t)}
  + \frac{\beta'\eta_t^2}{2}
  \\
  &=
  f(\xx_t) - \frac{1}{2\beta'}  \normq{\gradf(\xx_t)}^2.
\end{align*}
The first inequality is by $f$'s smoothness, and
the second is by our choice of $\beta'$.
The first equality is because $\xx_{t+1}-\xx_t=-\eta_t \vv_t$
and $\normp{\vv_t}=1$.
The second equality is by the algorithm's choice of $\vv_t$ since,
for all $\zz\in\Rn$, it is known that
\[
  \max\braces{\zz\cdot\vv :\: \vv\in\Rn, \normp{\vv}\leq 1} = \normq{\zz}.
\]
The last equality holds if we set
$\eta_t = \normq{\gradf(\xx_t)} / \beta'$.
Thus, the progress bound in \eqref{eqn:aux-fcn-prog-bnd} is satisfied with
$\alpha_t = 1/(2\beta')$
and
$\smash[b]{h(\uu) = \normq{\uu}^2}$, which is
a \auxiliary function. Since $\ef$ is continuous, it must be continuous
at every subsequential limit of the sequence $\seq{\xx_t}$, so
the conditions of \Cref{thm:fact:gd} are satisfied,
and hence the steepest-descent algorithm converges to the infimum of $f$.%
\indexg{steepest descent|)}%
\end{example}

\indexg{coordinate descent|(}%
In particular,
for smooth, convex functions $f$ as above, gradient-based coordinate
descent effectively minimizes $f$, for appropriate step sizes.
This also shows that the same holds for fully-greedy coordinate
descent since this version makes at least as much progress
on each iteration at decreasing the function values as the
gradient-based version.%
\indexg{coordinate descent|)}%

As a consequence, any of the algorithms just discussed can be applied
to a range of commonly encountered convex optimization problems which,
in general, might have no finite
\indexg{algorithms, minimization|)}%
\indexg{minimization methods|)}%
minimizers.
Here are some examples:

\begin{example}[Logistic regression]
\label{ex:logistic}
\indexg{logistic regression|(}%
In \Cref{sec:emp-loss-min} (Eq.~\ref{eqn:logistic-reg-obj}),
we saw that logistic regression
is based on minimization of a function $f:\Rn\rightarrow\R$
of the form
\[
   f(\xx) = \sum_{i=1}^m \ln\bigParens{1+\exp(\xx\cdot\uu_i)}
\]
over $\xx\in\Rn$, for some vectors $\uu_1,\ldots,\uu_m\in\Rn$.
This function is convex and proper, but its minimizers might be at infinity.
It turns out that $f$ is also smooth.
We sketch the main ideas of why this is so.

In general,
it is known that a convex, twice-differentiable function $f$ is
$\beta$-smooth if its Hessian
$\nabla^2 f$ satisfies
$\trans{\zz} [\nabla^2 f(\xx)] \zz \leq \beta \norm{\zz}^2$
for all $\xx,\zz\in\Rn$
\idxnesterov\citep[Theorem~2.1.6]{nesterov-lectures-cvx-opt}.
To compute the Hessian,
let $g(x)=\ln(1+e^x)$ for $x\in\R$,
whose second derivative is
$g''(x) = e^x/(1+e^x)^2$. 
Since $f(\xx)=\sum_{i=1}^m g(\xx\inprod\uu_i)$,
we then obtain that $f$'s Hessian is
$
   \nabla^2 f(\xx) = \sum_{i=1}^m g''(\xx\inprod\uu_i) \uu_i\trans{\uu_i}
$.
Hence, for $\xx,\zz\in\Rn$,
\[   
  \trans{\zz} [\nabla^2 f(\xx)] \zz
  =
  \sum_{i=1}^m g''(\xx\inprod\uu_i) \, (\zz\cdot\uu_i)^2
  \leq
  \norm{\zz}^2 \sum_{i=1}^m \norm{\uu_i}^2,
\]   
where the inequality is by
the Cauchy-Schwarz inequality,
and because $g''(x)\in[0,1]$ for all $x\in\R$.
Thus, $f$ is $\beta$-smooth, with
$\beta=\sum_{i=1}^m \norm{\uu_i}^2$.

Moreover, $f$'s extension $\fext$ is continuous everywhere
(\Cref{pr:hard-core:1}\ref{pr:hard-core:1:a-cnt}),
so $f$ satisfies the conditions of \Cref{ex:steepest-descent}.
Hence, any of the methods discussed above can be applied to minimize
this function.%
\indexg{logistic regression|)}%
\end{example}

\begin{example}[Exponential loss, AdaBoost]
\indexg{exponential loss|(}%
\indexg{AdaBoost|(}%
For $\xx\in\Rn$, let
\[
   f(\xx) = {\sum_{i=1}^m \exp(\xx\cdot\uu_i) },
\]
for some vectors $\uu_1,\ldots,\uu_m\in\Rn$, and let
$g(\xx)=\ln f(\xx)$.
\indexg{logarithm, extension of|(}%
Also, let $\logex:\Rext\to\eR$ denote
the continuous extension of the natural logarithm function
defined, for $\barx\in\Rext$, by
\begin{equation}  \label{eq:logex-dfn}
\indexm{ln}{$\logex$}{natural logarithm's extension}%
  \logex(\barx)=\begin{cases}
    -\infty
      &\text{if $\barx\leq 0$,}
    \\
    \ln \barx
      &\text{if $0<\barx<+\infty$,}
    \\
    +\infty
      &\text{if $\barx=+\infty$.}%
\indexg{logarithm, extension of|)}%
  \end{cases}
\end{equation}
Both functions $f$ and $g$ are convex and proper.
The extension $\fext$ is continuous everywhere
(by \Cref{pr:hard-core:1}\ref{pr:hard-core:1:a-cnt}),
which also implies that $\gext=\logex\circ\fext$ is continuous everywhere
(using \Cref{pr:Gf-cont}\ref{pr:Gf-cont:b}).
The function $g$ can be shown to be smooth, using
similar techniques as in \Cref{ex:logistic}.
Since $f(\xx)\in\Rstrictpos$ for all $\xx\in\Rn$ and since natural logarithm
is strictly increasing on $\Rstrictpos$, minimizing $g$ is equivalent
to minimizing $f$.
Thus, either $f$ or $g$ can be minimized by applying any of the
methods above to $g$ (even if $\inf g = -\infty$).
Moreover, whether applied to $f$ or~$g$,
fully-greedy coordinate descent is identical in its
updates (as are the other methods,
for appropriate choices of step size).
In particular, the AdaBoost
algorithm~\idxschapire\idxfreund\citep{schapire_freund_book_final}
can be
viewed as minimizing an \emph{exponential loss} function of exactly
the same form as $f$ using fully-greedy coordinate descent;
therefore, the arguments above prove that AdaBoost effectively
minimizes exponential loss
(as had previously been proved by
\idxcollinsetal\citealp{collins_schapire_singer_adaboost_bregman},
using more specialized techniques based on Bregman distances).%
\indexg{exponential loss|)}%
\indexg{AdaBoost|)}%
\end{example}

\begin{example}[Maximum likelihood]
\label{ex:MLE}
\indexg{maximum-likelihood distribution!methods for computing|(}%
\indexg{exponential-family distributions!estimation from data|(}%
Let $\distset$ be a finite and nonempty set with elements referred to as
\emph{outcomes}, described by a map
$\featmap:\distset\rightarrow\Rn$,
referred to as a \emph{feature map}.
For any $\vparam\in\Rn$, we can form a probability distribution
$\qx$ over outcomes $\distelt\in\distset$ of the form
\[
  \qx(\distelt) = \frac{e^{\vparam\cdot\featmap(\distelt)}}
                  {\sum_{\distalt\in\distset} e^{\vparam\cdot\featmap(\distalt)}}.
\]
The set of distributions $\qx$ across all $\vparam\in\Rn$ is an
example of an
\emph{exponential family}.

Our goal is to estimate an unknown distribution $\pi$ over $\distset$
based on $m$ independent random draws from $\pi$. If each
outcome $\distelt$ is observed $m_{\distelt}$ times,
we can estimate the underlying distribution by $\qx$, where $\vparam$
maximizes the likelihood of the data $\prod_{\distelt\in\distset}[\qx(\distelt)]^{m_{\distelt}}$. 
This maximization problem is equivalent to minimizing
the
\indexg{negative log-likelihood}%
\emph{negative log likelihood}
\[
  f(\vparam)
  =-\sum_{i\in I}\frac{m_i}{m} \ln \qx(\distelt),
\]
over $\vparam\in\Rn$.
Similar to the preceding examples,
this function
is convex, proper, smooth, and can be shown to have an extension
that is continuous everywhere.
Therefore, once again, any of the methods above can be applied to
minimize it, even if none of its minimizers are finite.
In the next chapter,
we study
this family of distributions and
the maximum likelihood estimation approach in detail.%
\indexg{exponential-family distributions!estimation from data|)}%
\indexg{maximum-likelihood distribution!methods for computing|)}%
\end{example}

\chapter{Exponential-family distributions}
\label{sec:exp-fam}

\indexg{exponential-family distributions|(}%
We next study in detail a broad and well-established family of probability
distributions,
briefly introduced in \Cref{ex:MLE},
called the \emph{exponential family}.
Typically, the practical problem of interest is to estimate the
parameters of such a distribution from data, which can be formulated
as minimization of a particular convex function.
In this chapter, we apply astral theory to explore the
common situation in which the minimizers, corresponding to parameter
values, might be at infinity.

The development we give here for the
standard setting is quite well-studied.
Many of the astral results presented below are generalizations or
extensions of known results for the standard case.
See, for instance, \idxwainjord\citet[Chapter~3]{WainwrightJo08}
for further background.

Throughout this chapter,
in accord with the exponential-family literature,
we use $\vparam$ to denote the vector of parameters of an
exponential-family distribution.
This entails a slight deviation from our usual notational conventions, with
the vector $\vparam$ and its astral
version $\parambar$ corresponding
to the primal variables $\xx$ and $\xbar$ from the rest of the book.
We continue to use $\uu$ as the dual variable.%

\section{The standard setting}
\label{sec:exp-fam:standard-setting}

For simplicity,
we focus on probability distributions defined over some nonempty,
finite set
$\distset$.
Elements of $\distset$ are denoted as $\distelt$ and referred to as
\indexg{outcomes}%
\emph{outcomes}.
In this chapter, a distribution $p$ is a function
$p:\distset\rightarrow [0,1]$
with $\sum_{\distelt\in\distset} p(\distelt)=1$.
We let $\Delta$ denote the set of all such distributions.

\indexg{feature map|(}%
We suppose we are given a \emph{feature map}
\indexm{phi}{$\featmap$}{feature map}%
$\featmap:\distset\rightarrow \Rn$,
for some $n\geq 1$.
For each $\distelt\in\distset$,
$\featmap(\distelt)$ can be regarded as a kind of description of
outcome $\distelt$, with each of the $n$ components $\featmapj$
providing one feature or descriptor.
For example, if each $\distelt\in\distset$ corresponds to a person, then the features
might provide that person's height, weight, age, and so
\indexg{feature map|)}%
on.\looseness=-1

We can use the feature map to construct probability distributions
over $\distset$.
\indexg{exponential-family distributions!defined|(}%
An \emph{exponential-family distribution} is defined by parameters
\indexm{theta}{$\vparam$}{exponential-family parameters}%
$\vparam\in\Rn$, and denoted $\qx$, placing probability mass on
$\distelt\in\distset$ proportional to
$e^{\vparam\cdot\featmap(\distelt)}$.
Thus,
\begin{equation}  \label{eqn:exp-fam-defn}
\indexm{p theta}{$\qx$}{exponential-family distribution}%
  \qx(\distelt)
  =
  \frac{e^{\vparam\cdot\featmap(\distelt)}}
       {\sumex(\vparam)}
  =
  \exp\bigParens{\vparam\cdot\featmap(\distelt) - \lpart(\vparam)},
\end{equation}
where $\sumex:\Rn\rightarrow \R$ provides normalization,
\indexg{log-partition function|(}%
and $\lpart:\Rn\rightarrow \R$ is its logarithm, called
the
\indexg{exponential-family distributions!defined|)}%
\emph{log-partition function}.
That is, for $\vparam\in\Rn$,
\begin{equation}  \label{eqn:sumex-dfn}
\indexm{z theta}{$\sumex(\vparam)$}{normalization factor}%
  \sumex(\vparam) =
    {\sum_{\distelt\in\distset} e^{\vparam\cdot\featmap(\distelt)}},
\end{equation}
and
\[
\indexm{a theta}{$\lpart(\vparam)$}{log-partition function}%
  \lpart(\vparam)
    = \ln \sumex(\vparam)
    = \ln\BiggParens{\sum_{\distelt\in\distset} e^{\vparam\cdot\featmap(\distelt)}}.%
\indexg{log-partition function|)}%
\]
Both of these functions are convex, proper, closed, and finite
everywhere;
they will play a central role in the development to follow.

\indexg{exponential-family distributions!estimation from data|(}%
Exponential-family distributions are commonly used to estimate an
unknown distribution from data.
In such a setting, we suppose access to $\numsamp$ independent samples from
some unknown distribution $\popdist$ over $\distset$.
Let $\Nelt$ be the number of times that outcome $\distelt$ is
observed.
From such information, we aim to estimate $\popdist$.
A standard approach is to posit that some exponential-family
distribution $\qx$ is a reasonable approximation to $\popdist$, and to then
select the parameter vector $\vparam\in\Rn$ yielding the best fit to the
data.
\indexg{likelihood|(}%
In particular, as a measure of fit,
we can use the \emph{likelihood} of the data, that is,
the probability of
observing the sequence of outcomes that actually were observed if
we suppose that the unknown distribution $\popdist$ is in fact $\qx$.
This likelihood
can be computed to be
\[
  \prod_{\distelt\in\distset} [\qx(\distelt)]^{\Nelt}.
\]
\indexg{maximum-likelihood distribution|(}%
According to the \emph{maximum-likelihood principle},
the parameters $\vparam\in\Rn$ should be chosen to maximize this likelihood.
\indexg{negative log-likelihood|(}%
Taking the negative logarithm and dividing by $\numsamp$, this is equivalent to
choosing $\vparam$ to minimize the negative log-likelihood,
\begin{equation}  \label{eqn:neg-log-like:1}
  \loglik(\vparam)
  =
  -\sum_{\distelt\in\distset} \frac{\Nelt}{\numsamp} \ln \qx(\distelt).%
\indexg{maximum-likelihood distribution|)}%
\indexg{exponential-family distributions!estimation from data|)}%
\end{equation}

The negative log-likelihood function is convex in $\vparam$;
however, it might not have a finite
minimizer.
Here are some examples:

\begin{example}
Suppose $\distset=\{1,\ldots,n\}$, and that
$\featmap(\distelt)=\ee_{\distelt}$, for $\distelt\in\distset$,
where $\ee_1,\ldots,\ee_n$ are the standard basis vectors in $\Rn$.
In this case, $\qx(\distelt)=e^{\param_\distelt}/\sumex(\vparam)$
for all $\distelt$.
Therefore,
every distribution $p\in\Delta$ that assigns positive probability
to each outcome
can be represented in this form by
setting $\param_\distelt=\ln p(\distelt)$
for $\distelt\in\distset$.
Further,
$\loglik(\vparam)$
has exactly the form of the function given in
\eqref{eqn:log-sum-exp-fcn} as part of Example~\ref{ex:log-sum-exp}
(with $\lsePh_i=\Nelt/\numsamp$).
Nonetheless,
as we saw in that example, 
if one or more of the counts $\Nelt$ is zero, then
it may happen that $\loglik$ will have no
finite minimizer (for instance, if $n=3$, $\numsamp=6$,
$\Nj{1}=0$, $\Nj{2}=2$, and $\Nj{3}=4$).
\end{example}

\begin{figure}
  \centering
  \includegraphics{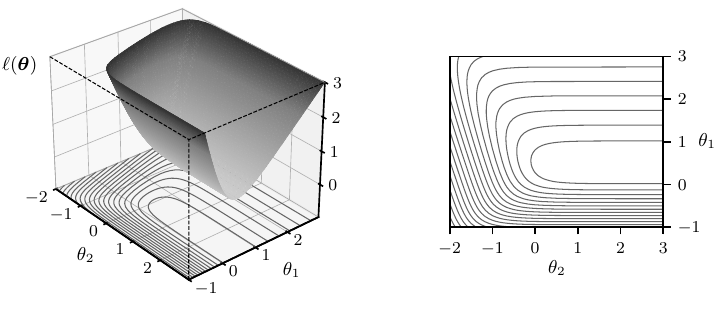}
  \mycaption{Negative log-likelihood}{%
\indexf{Four-outcome exponential family!negative log-likelihood of}%
    The negative log-likelihood function $\ell$ for the
    four-outcome exponential family from
    \Cref{ex:exp-fam-eg1}.
    The contour plot on the right indicates that the function
    is decreasing in the direction of $\ee_2$.
  }%
  \label{fig:neg-log-lik}%
\end{figure}

\begin{example}[Four-outcome exponential family]   \label{ex:exp-fam-eg1}
\indexg{Four-outcome exponential family|(}%
Suppose $\distset=\{\outA,\outB,\outC,\outD\}$, $n=2$, and
$\featmap$ is defined as follows:
\[
  \featmap(\outA)
  = \begin{bmatrix*}[r] 3 \\ 1 \end{bmatrix*},
  \qquad
  \featmap(\outB)
  = \begin{bmatrix*}[r] 0 \\ 1 \end{bmatrix*},
  \qquad
  \featmap(\outC)
  = \begin{bmatrix*}[r] -3 \\ 1 \end{bmatrix*},
  \qquad
  \featmap(\outD)
  = \begin{bmatrix*}[r] 1 \\ -2 \end{bmatrix*}.
\]
Suppose further that
$\numsamp=3$,
$\Nj{\outA}=2$,
$\Nj{\outB}=1$,
and
\indexg{Four-outcome exponential family|)}%
$\Nj{\outC}=\Nj{\outD}=0$.
\indexg{Four-outcome exponential family!no finite likelihood maximizer|(}%
Then the negative log-likelihood can be computed,
for $\vparam\in\R^2$,
to be
\begin{align}
  \loglik(\vparam)
  =
  \loglik(\param_1,\param_2)
  &=
  \regBracks{\ln\sumex(\vparam)}
  - \sfrac{2}{3} (3\param_1+\param_2)
  - \sfrac{1}{3} \param_2
  \nonumber
  \\
  &=
  \ln\bigParens{
      e^{3\param_1 + \param_2}
      +
      e^{\param_2}
      +
      e^{-3\param_1 + \param_2}
      +
      e^{\param_1 - 2\param_2}
  }
  -
  (2\param_1 + \param_2)
  \nonumber
  \\
  &=
  \ln\bigParens{
      e^{\param_1}
      +
      e^{-2\param_1}
      +
      e^{-5\param_1}
      +
      e^{-\param_1 - 3\param_2}
  }.
\label{eq:log-like-no-max-eg:2}
\end{align}
The second equality is by direct evaluation of 
\eqref{eqn:neg-log-like:1},
and the fourth by writing
$-(2\param_1 + \param_2)=\ln(e^{-2\param_1 - \param_2})$
and simplifying.
The negative-likelihood function $\ell$ is shown in \Cref{fig:neg-log-lik}.
This function cannot have a finite minimizer since
adding $\ee_2$ to any point $\vparam\in\R^2$
(that is, incrementing $\param_2$)
leaves the first three terms inside the logarithm in the last line unchanged
while strictly diminishing the last term; thus,
$\loglik(\vparam+\ee_2) < \loglik(\vparam)$ for all $\vparam\in\R^2$.%
\indexg{Four-outcome exponential family!no finite likelihood maximizer|)}%
\indexg{negative log-likelihood|)}%
\indexg{likelihood|)}%
\end{example}

These examples show that even in simple cases, when working with an
exponential family of distributions, there may be
no value $\vparam\in\Rn$ that maximizes
the likelihood of the observed data.
Nevertheless, in this chapter, we will see how the parameter space can
be extended from $\Rn$ to astral space, $\extspace$, in a way that
preserves and even enhances some of the key properties
of exponential-family distributions, and in particular, ensures that
there always exists a parameter vector that maximizes the likelihood of the data.

The fact that there might not exist a finite parameter vector maximizing the likelihood
has been noted and studied in the previous literature,
for instance, by
\indexa{Haberman, S. J.}\citet[Chapter~2]{haberman74},
\indexa{Wedderburn, R. W. M.}\citet{wedderburn76},
\indexa{Silvapulle, M. J.}\citet{silvapulle81},
\indexa{Albert, A.}\indexa{Anderson, J. A.}\citet{albert-anderson},
and
\indexa{Geyer, C. J.}\citet{geyer09}.
This further led to the study of the topological closure of
the space of exponential-family distributions, thereby extending
the space via sequential limits; see
\indexa{Barndorff-Nielsen, O. E.}\citet[\ppcite{154--156}]{barndorff-nielsen},
\indexa{Brown, L. D.}\citet[\ppcite{191--202}]{brown86},
and
\indexa{csiszar, i.@Csisz\'{a}r, I.}\indexa{matus, f.@Mat\'{u}\v{s}, F.}%
\citet{csiszar-matus03,csiszar-matus05}.
In this chapter, we address this difficulty by directly
extending both the parameter space and the space of distributions
using notions from astral
\indexg{exponential-family distributions|)}%
space.

\section{Extending to astral space}

In what follows, it will be helpful to consider simple translations of
the feature map in which some fixed vector $\uu\in\Rn$ is
subtracted from all its values.
\indexg{feature map!translated|(}%
Thus, we define a modified feature map
$\featmapu:\distset\rightarrow\Rn$, for $\distelt\in\distset$,
as
\begin{equation}   \label{eq:featmapu-dfn}
\indexm{phi u}{$\featmapu$}{translated feature map}%
  \featmapu(\distelt)=\featmap(\distelt)-\uu.
\end{equation}
\indexg{log-partition function!translated feature map@for translated feature map|(}%
Likewise, we define variants of $\lpart$ and $\sumex$, denoted
$\lpartu$ and $\sumexu$, that are associated with $\featmapu$ rather
than $\featmap$.
That is, for $\vparam\in\Rn$, 
\begin{equation}   \label{eq:sumexu-dfn}
\indexm{zu theta}{$\sumexu(\vparam)$}{normalization factor (for $\featmapu$)}%
  \sumexu(\vparam) =
    \sum_{\distelt\in\distset} e^{\vparam\cdot\featmapu(\distelt)}
    =
    \sum_{\distelt\in\distset} e^{\vparam\cdot(\featmap(\distelt)-\uu)}
    =
    \sumex(\vparam)e^{-\vparam\cdot\uu},
\end{equation}
and
\begin{equation}  \label{eqn:lpart-defn}
\indexm{a u theta}{$\lpartu(\vparam)$}{log-partition function (for $\featmapu$)}%
  \lpartu(\vparam)
    = \ln \sumexu(\vparam)
    = \lpart(\vparam) - \vparam\cdot\uu,
\end{equation}
where the last line follows from \eqref{eq:sumexu-dfn}.
Thus, $\lpartu$ is the same as $\aminusu$,
the linear tilt of~$a$ by~$\uu$,
as defined in
\indexg{log-partition function!translated feature map@for translated feature map|)}%
\Cref{dfn:lin-tilt-ext}.
(However,
the same is not true for $\featmapu$ and~$\sumexu$, nor for
the map $\meanmapu$ defined below in
Eq.~\ref{eqn:meanmapu-defn}.)

Observe that the exponential-family distributions associated with
$\featmapu$ are the same as the original distributions $\qx$
associated with $\featmap$
since, for $\distelt\in\distset$,
\[
  \frac{e^{\vparam\cdot\featmapu(\distelt)}}
       {\sumexu(\vparam)}
  =
  \frac{e^{\vparam\cdot\featmap(\distelt)}}
       {\sumex(\vparam)}
  =
  \qx(\distelt).%
\indexg{feature map!translated|)}%
\]

Finally, we note that $\sumexu$, as given in the first sum in
\eqref{eq:sumexu-dfn},
has exactly the form of the functions
studied in \Cref{sec:emp-loss-min}, specifically
\eqref{eqn:loss-sum-form}
(with the variables appearing in that equation set to
$\ell_\distelt=\exp$ and $\uu_\distelt=\featmapu(\distelt)$
for $\distelt\in\distset$).
This connection will make it possible later to apply much of what was
developed in that section to the present setting.

\indexg{log-partition function!extension of|(}%
As we show next,
the functions $\sumexu$ and $\lpartu$ extend continuously
to astral space, for any $\uu\in\Rn$.
In what follows,
$\expex$ and $\logex$
denote the continuous, nondecreasing extensions of the
exponential and logarithm functions
that were defined in
Eqs.~(\ref{eq:expex-defn})
and~(\ref{eq:logex-dfn}).

\begin{proposition}  \label{pr:sumex-lpart-cont}
  Let $\featmap:\distset\rightarrow\Rn$, where $\distset$ is finite
  and nonempty.
  Let $\uu\in\Rn$, $\parambar\in\extspace$,
  and let
  $\featmapu$, $\sumex$, $\lpart$, $\sumexu$, $\lpartu$ be as defined
  above.
  Then:
  \begin{letter}
  \item    \label{pr:sumex-lpart-cont:a}
    The extensions $\sumexextu$ and $\lpartextu$
    are continuous everywhere.
  \item    \label{pr:sumex-lpart-cont:b}
    $\displaystyle
      \sumexextu(\parambar)
      =
      \sum_{\distelt\in\distset}
      \expex\bigParens{\parambar\cdot\featmapu(\distelt)}
    $.
  \item    \label{pr:sumex-lpart-cont:c}
    $\displaystyle
      \lpartextu(\parambar)
      =
      \logex\bigParens{\sumexextu(\parambar)}
    $.
  \end{letter}
\end{proposition}

\begin{proof}
As noted earlier, $\sumexu$ has exactly the form
of the functions considered in \Cref{sec:emp-loss-min}
(as in Eq.~\ref{eqn:loss-sum-form}).
Therefore, the form and continuity of $\sumexextu$ follow directly from
\Cref{pr:hard-core:1}(\ref{pr:hard-core:1:a-ext},\ref{pr:hard-core:1:a-cnt}).
The form and continuity of $\lpartextu$ then follow
from \Cref{pr:Gf-cont}(\ref{pr:Gf-cont:b})
(and \Cref{thm:ext-cont-f}),
since $\sumexextu$ and $\logex$
are continuous.%
\indexg{log-partition function!extension of|)}%
\end{proof}

\indexg{extended exponential-family distributions|(}%
We next extend the exponential family $\set{\qx:\:\vparam\in\Rn}$ to include
distributions with parameters $\parambar\in\extspace$.
We do this by continuously extending
the map $\vparam\mapsto\qx$ to astral space. 
From Eqs.~(\ref{eqn:exp-fam-defn}) and~(\ref{eqn:lpart-defn}),
for $\distelt\in\distset$ and $\vparam\in\Rn$,
\begin{equation}  \label{eqn:qx-alt-form}
  \qx(\distelt)
  =
  \exp\bigParens{\vparam\cdot\featmap(\distelt) - \lpart(\vparam)}
  =
  \exp\bigParens{-\lpartfi(\vparam)}.
\end{equation}
In this form, extending to astral space is straightforward
since all of the functions involved extend continuously.
\indexg{extended exponential-family distributions!defined|(}%
Thus, for $\parambar\in\extspace$, we define a distribution
$\qxbar\in\Delta$ with astral parameters $\parambar$ as
\begin{equation}  \label{eqn:qxbar-defn}
  \qxbar(\distelt)
  =
  \expex\bigParens{-\lpartextfi(\parambar)}.
\end{equation}
The set of distributions $\set{\qxbar:\:\parambar\in\extspace}$
is referred to as the
\indexg{extended exponential-family distributions!defined|)}%
\emph{extended exponential family}.

Note that the notation $\qxbar$
defined in \eqref{eqn:qxbar-defn}
is compatible with our earlier
notation $\qx$ from
\eqref{eqn:exp-fam-defn}.
In other words, if $\parambar=\vparam\in\Rn$, then
$\qxbar$, as defined in 
\eqref{eqn:qxbar-defn}, is the same as 
$\qx$
in Eqs.~(\ref{eqn:exp-fam-defn}) and~(\ref{eqn:qx-alt-form}).
This is because, for $\distelt\in\distset$,
$\lpartextfi(\vparam)$ is equal to $\lpartfi(\vparam)$
(by \Cref{pr:h:1}\ref{pr:h:1a},
since $\lpartfi$ is continuous),
and so is also finite.

Applying \Cref{pr:sumex-lpart-cont}(\ref{pr:sumex-lpart-cont:b}),
each distribution $\qxbar$ can be expressed more explicitly as
\begin{equation}  \label{eqn:qxbar-expr}
  \qxbar(\distelt)
  =
  \invex\bigParens{\sumexextfi(\parambar)}
  =
  \invex\BiggParens{\,\sum_{\distalt\in\distset}
    \expex\bigParens{\parambar\cdot[\featmap(\distalt)-\featmap(\distelt)]}}.
\end{equation}
\indexg{reciprocal function, extension of|(}%
Here, $\invex:[0,+\infty]\rightarrow\Rext$ extends the
reciprocal
function
$x\mapsto 1/x$;
that is,
for $\barx\in [0,+\infty]$,
\[
\indexm{inv}{$\invex$}{reciprocal function's extension}%
   \invex(\barx) =
   \begin{cases}
                         +\infty       & \text{if $\barx=0$,} \\
                     {1}/{\barx}   & \text{if $0<\barx<+\infty$,} \\
                         0             & \text{if $\barx=+\infty$.}%
\indexg{reciprocal function, extension of|)}%
   \end{cases}
\]

\begin{proposition}   \label{pr:qxbar-cont}
  Let $\qxbar$ be as defined in \eqref{eqn:qxbar-defn}
  in terms of some feature map
  $\featmap:\distset\rightarrow\Rn$, where $\distset$ is finite
  and nonempty.
  Let $\distelt\in\distset$.
  Then:
  \begin{letter-compact}
  \item   \label{pr:qxbar-cont:a}
    As a function of $\parambar\in\extspace$,
    $\qxbar(\distelt)$ is continuous everywhere.

  \item   \label{pr:qxbar-cont:b}
    For $\parambar\in\extspace$,
    $\qxbar$ is a probability distribution in $\Delta$.
  \end{letter-compact}
\end{proposition}

\begin{proof}
  ~

\begin{proof-parts}
\pfpart{Part~(\ref{pr:qxbar-cont:a}):}
Continuity of $\qxbar(\distelt)$ follows directly from
continuity of $\lpartextfi(\parambar)$
(\Cref{pr:sumex-lpart-cont}\ref{pr:sumex-lpart-cont:a})
and of $\expex$.

\pfpart{Part~(\ref{pr:qxbar-cont:b}):}
Let $\parambar\in\extspace$, and
let
$\seq{\vparam_t}$ be a sequence in $\Rn$ converging to $\parambar$
(which exists by \Cref{thm:i:1}\ref{thm:i:1d}).
By part~(\ref{pr:qxbar-cont:a}),
this implies
$\qxt(\distelt)\rightarrow\qxbar(\distelt)$
for each $\distelt\in\distset$,
and so that
$\qxbar(\distelt)\geq 0$
since
$\qxt(\distelt)\geq 0$
for all $t$.
Likewise, $\sum_{\distelt\in\distset}\qxbar(\distelt)=1$
since
\[
  1 = \sum_{\distelt\in\distset}\qxt(\distelt)
    \rightarrow
 \sum_{\distelt\in\distset}\qxbar(\distelt),
\]
by
continuity of addition (\Cref{prop:lim:eR}\ref{i:lim:eR:sum}).
Thus, $\qxbar\in\Delta$.
\qedhere
\end{proof-parts}
\end{proof}

For a distribution $p\in\Delta$, and any (scalar- or vector-valued)
function $f$ on $\distset$,
we write $\regExp{p}{f}$ for the expected value of $f$ with respect to
$p$:
\begin{equation}  \label{eqn:expect-defn}
\indexm{ep f}{$\regExp{p}{f}$}{expected value}%
  \regExp{p}{f}
  =
  \regExp{\distelt\sim p}{f(\distelt)}
  =
  \sum_{\distelt\in\distset} p(\distelt) f(\distelt).
\end{equation}

\indexg{mean map|(}%
\indexg{mean map!defined|(}%
We then define the \emph{mean map}
$\meanmap:\extspace\rightarrow\Rn$ which maps the
parameters $\parambar\in\extspace$ to the mean of the feature map $\featmap$
under the distribution $\qxbar$ defined by $\parambar$:
\begin{equation}  \label{eq:mean-map-dfn}
\indexm{m theta 300}{$\meanmap(\parambar)$}{mean map}%
  \meanmap(\parambar)
  =
  \regExp{\qxbar}{\featmap}.%
\indexg{mean map!defined|)}%
\end{equation}
We will see that this map plays an important role in the development
to follow, as it does generally in the study of standard
exponential-family distributions.

We also define $\meanmapu$ to be the same as $\meanmap$ with
$\featmap$ replaced by $\featmapu$; thus,
\begin{equation}  \label{eqn:meanmapu-defn}
\indexm{m theta 700}{$\meanmapu(\parambar)$}{mean map (for $\featmapu$)}%
  \meanmapu(\parambar)
  =
  \regExp{\qxbar}{\featmapu}
  =
  \meanmap(\parambar) - \uu,
\end{equation}
since, as noted earlier,
$\qxbar$ is unaffected when $\featmap$ is shifted
by a fixed vector $\uu$.

\indexg{log-partition function!gradient of|(}%
Straightforward calculus shows that, for $\vparam\in\Rn$,
$\meanmap(\vparam)$ is in fact the gradient of the log-partition function
$\lpart$ at $\vparam$; that is,
\begin{equation}  \label{eqn:grad-is-meanmap}
  \nabla\lpart(\vparam) = \meanmap(\vparam).
\end{equation}
We will see below how this fact extends to astral parameters.%
\indexg{log-partition function!gradient of|)}%

We first show that the function $\meanmap$ is continuous and also
commutes with the operation of taking a topological closure.
Here and for the rest of this chapter,
we refer to the foregoing definitions of
$\featmapu$, $\qx$, $\qxbar$, $\sumex$, $\lpart$, $\sumexu$, $\lpartu$,
$\meanmap$ and $\meanmapu$, all in terms of
the feature map $\featmap:\distset\rightarrow\Rn$,
where $\distset$ is finite and nonempty,
as the
\emph{general setup of this chapter}.

\begin{proposition}  \label{pr:meanmap-cont}
  Assume the general setup of this chapter.
  Then:
  \begin{letter-compact}
  \item   \label{pr:meanmap-cont:a}
    $\meanmap$ is continuous.
  \item   \label{pr:meanmap-cont:b}
    For every set $S\subseteq\extspace$,
    $\meanmap(\Sbar)=\cl(\meanmap(S))$.
  \end{letter-compact}
\end{proposition}

\begin{proof}
From \Cref{pr:qxbar-cont}(\ref{pr:qxbar-cont:a}),
for each $\distelt\in\distset$,
the map $\parambar\mapsto\qxbar(\distelt)$ is continuous.
The continuity of $\meanmap$ then follows by continuity of
the linear map $p\mapsto\regExp{p}{\featmap}$,
proving part~(\ref{pr:meanmap-cont:a}).
Since $M$ is continuous,
part~(\ref{pr:meanmap-cont:b}) follows by
\Cref{pr:cont-from-compact}(\ref{pr:cont-from-compact:b}).%
\indexg{mean map|)}%
\indexg{extended exponential-family distributions|)}%
\end{proof}

\section{Conjugate and astral subgradients}

\indexg{extended exponential-family distributions|(}%
In the development to follow,
\indexg{log-partition function!conjugate of}%
$\lpartstar$, the conjugate of $\lpart$,
will play an important role.
\indexg{entropy (of distribution)|(}%
\indexg{extended exponential-family distributions!entropy of|(}%
The next lemma shows that if $\meanmap(\parambar)=\uu$,
then $\lpartstar(\uu)$ is equal to $-\entropy(\qxbar)$,
where $\entropy(p)$ denotes the \emph{entropy} of any distribution $p\in\Delta$:
\[
\indexm{h p}{$\entropy(p)$}{entropy}%
  \entropy(p)
  =
  -\sum_{\distelt\in\distset} p(\distelt) \ln p(\distelt)
  =
  -\regExp{p}{\ln p},
\]
where,
in expressions like the one on the right,
we use $\ln p$ as shorthand for the
function $\distelt\mapsto \ln p(\distelt)$. Also, we
define $\ln 0=-\infty$, so using our arithmetic, $0\ln 0=0$.
Note that $\entropy(p)$ is finite for all $p\in\Delta$.

First, we state some simple facts that will be used here and
elsewhere:

\begin{proposition}   \label{pr:exp-log-props}
  Let $p\in\Delta$.
  Then:
  \begin{letter-compact}
  \item   \label{pr:exp-log-props:a}
    For all $q\in\Delta$,
    \[
      \regExp{p}{\ln q} \leq \regExp{p}{\ln p},
    \]
    with equality if and only if $q=p$.
    In other words, $\regExp{p}{\ln q}$ is uniquely maximized over
    $q\in\Delta$ when $q=p$.
  \item   \label{pr:exp-log-props:b}
    For all $\vparam\in\Rn$,
    \[
      \regExp{p}{\ln \qx}
      =
      \vparam\cdot\regExp{p}{\featmap} - \lpart(\vparam).
    \]
  \item   \label{pr:exp-log-props:c}
    Let $\seq{\vparam_t}$ be a sequence in $\Rn$ that converges to some
    point $\parambar\in\extspace$.
    Then
    \[
      \bigExp{p}{\ln \qxt} \rightarrow \bigExp{p}{\ln \qxbar}.
    \]
  \end{letter-compact}
\end{proposition}

\begin{proof}
~

\begin{proof-parts}
\pfpart{Part~(\ref{pr:exp-log-props:a}):}
See
\indexa{Cover, T. M.}\indexa{Thomas, J. A.}\citet[Theorem~2.6.3]{cover-thomas-2nd}.

\pfpart{Part~(\ref{pr:exp-log-props:b}):}
By \eqref{eqn:exp-fam-defn} and
linearity of expectations,
for all $\vparam\in\Rn$,
\[
  \regExp{p}{\ln \qx}
  =
  \regExp{\distelt\sim p}{\vparam\cdot\featmap(\distelt) - \lpart(\vparam)}
  =
  \vparam\cdot\regExp{p}{\featmap} - \lpart(\vparam).
\]

\pfpart{Part~(\ref{pr:exp-log-props:c}):}
By
\Cref{pr:qxbar-cont}(\ref{pr:qxbar-cont:a}),
for $\distelt\in\distset$,
$\qxt(\distelt)\rightarrow\qxbar(\distelt)$,
implying
$\ln \qxt(\distelt)\rightarrow \ln \qxbar(\distelt)$,
by continuity of the logarthm function.
Therefore,
$\regExp{p}{\ln \qxt}\rightarrow \regExp{p}{\ln \qxbar}$
by continuity of scalar multiplication
and addition
(\Cref{prop:lim:eR}\ref{i:lim:eR:sum}\ref{i:lim:eR:mul}),
noting for summability that
$p(\distelt)\ln \qxbar(\distelt) \leq 0<+\infty$
for $\distelt\in\distset$.%
\indexg{entropy (of distribution)|)}%
\qedhere
\end{proof-parts}
\end{proof}

\begin{lemma} \label{lem:lpart-conj}
\indexg{log-partition function!conjugate of}%
\indexg{mean map|(}%
  Assume the general setup of this chapter.
  Let $\parambar\in\extspace$ and $\uu\in\Rn$, and
  suppose $\meanmap(\parambar)=\uu$.
  Then
  $\lpartstar(\uu)=-\entropy(\qxbar)$.
\end{lemma}

\begin{proof}
By
\Cref{pr:exp-log-props}(\ref{pr:exp-log-props:b}),
for $\vparam\in\Rn$,
$
  \regExp{\qxbar}{\ln \qx}
  =
  \vparam\cdot\uu - \lpart(\vparam)
$
since $\regExp{\qxbar}{\featmap}=\meanmap(\parambar)=\uu$.
As a result, by definition of conjugate,
\begin{equation}  \label{eqn:lem:lpart-conj:1}
  \lpartstar(\uu)
  =
  \sup_{\vparam\in\Rn} [\vparam\cdot\uu - \lpart(\vparam)]
  =
  \sup_{\vparam\in\Rn}  \regExp{\qxbar}{\ln \qx} .
\end{equation}
By
\Cref{pr:exp-log-props}(\ref{pr:exp-log-props:a}),
for all $\vparam\in\Rn$,
$\regExp{\qxbar}{\ln \qx} \leq \regExp{\qxbar}{\ln \qxbar} = -\entropy(\qxbar)$.
Therefore, $\lpartstar(\uu)\leq -\entropy(\qxbar)$.

For the reverse inequality,
there must exist a sequence $\seq{\vparam_t}$ in $\Rn$ that
converges to~$\parambar$
(by \Cref{thm:i:1}\ref{thm:i:1d}).
By \Cref{pr:exp-log-props}(\ref{pr:exp-log-props:c}),
$\regExp{\qxbar}{\ln \qxt}\rightarrow \regExp{\qxbar}{\ln \qxbar}$.
Combined with \eqref{eqn:lem:lpart-conj:1}, this implies
$\lpartstar(\uu)\geq \regExp{\qxbar}{\ln \qxbar}$, completing the proof.%
\indexg{extended exponential-family distributions!entropy of|)}%
\end{proof}

\begin{figure}
  \centering
  \includegraphics{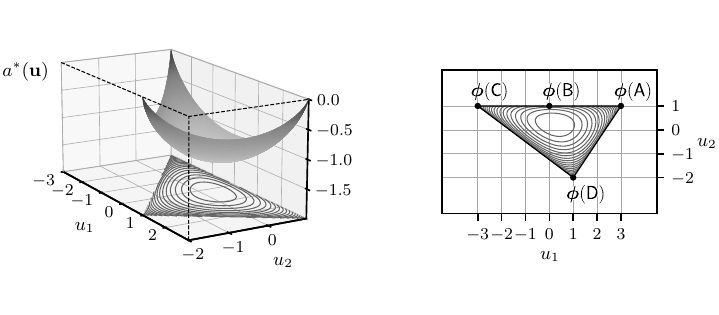}
  \mycaption{Conjugate of log-partition function and marginal polytope}{%
\indexf{Four-outcome exponential family!conjugate of log-partition function}%
    \emph{Left:}
    The conjugate $\lpartstar$ of the log-partition function
    for the four-outcome exponential family from \Cref{ex:exp-fam-eg1}.
    \emph{Right:}
    The effective domain of $\lpartstar$ coincides with the marginal polytope,
    which is the convex hull of the four feature vectors
    $\featmap(\outA)$, $\featmap(\outB)$, $\featmap(\outC)$,
    $\featmap(\outD)$.}
  \label{fig:astar}%
\end{figure}

\indexg{marginal polytope|(}%
The convex hull of the set $\featmap(\distset)$
is called the \emph{marginal polytope},
consisting of all convex combinations of the points
$\featmap(\distelt)$ for $\distelt\in\distset$
(\Cref{roc:thm2.3}).   %
This is exactly the set of means
$\regExp{p}{\featmap}$ that can be realized
by any distribution $p\in\Delta$
(not necessarily in the exponential family).
The next theorem shows that for every point
$\uu\in\conv{\featmap(\distset)}$, which is to say
every point for which there exists \emph{some} distribution
$p\in\Delta$ with $\regExp{p}{\featmap}=\uu$,
there must also exist an exponential-family distribution
with parameters $\parambar\in\extspace$ realizing the same mean,
so that $\meanmap(\parambar)=\regExp{\qxbar}{\featmap}=\uu$.
Thus, $\meanmap(\extspace)=\convfeat$, which,
in light of \Cref{lem:lpart-conj},
is also equal to the effective domain of $\lpartstar$.
(The conjugate $\lpartstar$ and the marginal polytope
for the four-outcome exponential family from
\Cref{ex:exp-fam-eg1}
are shown in \Cref{fig:astar}.)
The theorem further shows that the image of $\Rn$ under
$\meanmap$ is equal to the relative interior of this same set.

In the proof of the \namecref{thm:meanmap-onto} and throughout the rest of
this chapter,
we use the notion of
\irredIndexSet
from \Cref{sec:emp-loss-min}.
Specifically, we consider the
\irredIndexSetFor~$\sumexu$,
which is, following \Cref{def:hard-core}, defined as
\begin{equation}
\label{eq:hc:sumexu}
  \hardcore{\sumexu} =
  \set{\distelt\in\distset:\: \featmapu(\distelt)\in(\resc\sumexu)^\perp}.
\end{equation}

\begin{theorem}  \label{thm:meanmap-onto}
  Assume the general setup of this chapter.
  Then:
  \begin{letter-compact}
  \item    \label{thm:meanmap-onto:a}
    $\meanmap(\Rn) = \ri\regParens{\convfeat}$.
  \item    \label{thm:meanmap-onto:b}
    $\meanmap(\extspace) = \convfeat = \dom \lpartstar$.
  \item    \label{thm:meanmap-onto:c}
    For each $\uu\in\convfeat$,
    there exists $\parambar\in\extspace$ for which
    $\meanmap(\parambar)=\uu$, implying
    $\lpartstar(\uu)=-\entropy(\qxbar)$.
  \end{letter-compact}
\end{theorem}

\begin{proof}
  ~
  
\begin{proof-parts}
\pfpart{Part~(\ref{thm:meanmap-onto:a}):}
For $\vparam\in\Rn$,
note that
$\qx$ is a distribution in $\Delta$
with $\qx>0$.
Therefore, $\meanmap(\vparam)\in\ri(\convfeat)$ by
\Cref{pr:ri-conv-finite}.
Thus,
$\meanmap(\Rn)\subseteq\ri(\convfeat)$.

For the reverse inclusion,
let $\uu\in\ri(\convfeat)$;
we aim to prove $\uu\in\meanmap(\Rn)$.

First, $\zero\in\ri(\convfeatu)$ since subtracting $\uu$ from
$\featmap$ simply translates all points and sets, including
$\convfeat$, by $-\uu$.
Since $\sumexu$ has the form of functions in
\Cref{sec:emp-loss-min},
we can apply
\Cref{thm:erm-faces-hardcore}(\ref{thm:erm-faces-hardcore:a})
yielding that $\distset\subseteq \hardcore{\sumexu}$,
so $\hardcore{\sumexu}=\distset$.
This in turn implies, by
\Cref{thm:hard-core:3}(\ref{thm:hard-core:3:b:univ}),
that  %
$\zero\in\unimin{\sumexu}$.
Letting $\vparam\in\Rn$ minimize $\fullshadsumexu$
(which exists by
\Cref{pr:univ-red-props}\ref{pr:univ-red-props:min}),
it therefore follows by
\Cref{pr:unimin-to-global-min} that
$\zero\plusl\vparam=\vparam$ minimizes $\sumexextu$.

Since $\vparam$ minimizes $\sumexextu$, it also minimizes $\sumexu$
(by \Cref{pr:h:1}\ref{pr:h:1a}, since $\sumexu$ is continuous
everywhere),
and therefore $\lpartu$ as well
since logarithm is strictly increasing.
Being a differentiable function, this implies
$\nabla\lpartu(\vparam)=\zero$,
so
$\meanmapu(\vparam)=\zero$ and ${\meanmap(\vparam)=\uu}$
(by
Eqs.~\ref{eqn:meanmapu-defn}
and~\ref{eqn:grad-is-meanmap}).
Therefore,
$\uu\in\meanmap(\Rn)$, completing the proof.

\pfpart{Part~(\ref{thm:meanmap-onto:b}):}
We first prove two inclusions:

\begin{claimpx}  \label{cl:thm:meanmap-onto:3}
  $\meanmap(\extspace)\subseteq \dom \lpartstar$.
\end{claimpx}

\begin{proofx}
Suppose $\uu\in\meanmap(\extspace)$, so that $\meanmap(\parambar)=\uu$ for
some $\parambar\in\extspace$.
Then by \Cref{lem:lpart-conj},
$\lpartstar(\uu)=-\entropy(\qxbar)\leq 0$.
Thus, $\uu\in \dom \lpartstar$.
\end{proofx}

\begin{claimpx}  \label{cl:thm:meanmap-onto:1}
  $\ri(\dom \lpartstar)\subseteq \meanmap(\Rn)$.
\end{claimpx}

\begin{proofx}
Suppose $\uu\in\ri(\dom\lpartstar)$.
Since $\lpartstar$ is convex and proper
(by Propositions~\ref{pr:conj-props}\ref{pr:conj-props:d}
and~\ref{pr:conj-props-cvx}\ref{pr:conj-props-cvx:a},
since $\lpart$ is convex and proper),
it has a subgradient at every
point in $\ri(\dom \lpartstar)$
(by \Cref{roc:thm23.4}\ref{roc:thm23.4:a}),
implying there exists $\vparam\in\subdiflpartstar(\uu)$.
Since $\lpart$ is closed and proper,
this further implies that $\uu\in\subdiflpart(\vparam)$
by \Cref{pr:stan-subgrad-equiv-props}(\ref{pr:stan-subgrad-equiv-props:c}\ref{pr:stan-subgrad-equiv-props:a}).
Since $\lpart$ is differentiable everywhere,
the only element of
$\subdiflpart(\vparam)$ is $\nabla\lpart(\vparam)$
(by \Cref{roc:thm25.1}\ref{roc:thm25.1:a}).
Therefore, $\uu=\nabla\lpart(\vparam)=\meanmap(\vparam)$
by \eqref{eqn:grad-is-meanmap},
so $\uu\in\meanmap(\Rn)$.
\end{proofx}

Combining now yields
\begin{align}
\notag
  \meanmap(\extspace)
  \subseteq
  \dom \lpartstar
  &\subseteq
  \cl\regParens{\dom \lpartstar}
\\
\label{eq:thm:meanmap-onto:1}
  &=
  \cl\bigParens{\ri(\dom \lpartstar)}
  \subseteq
  \cl\bigParens{\meanmap(\Rn)}
  =
  \meanmap(\extspace).
\end{align}
The first inclusion is by
\Cref{cl:thm:meanmap-onto:3}.
The first equality is by 
\Cref{pr:ri-props}(\ref{pr:ri-props:roc-thm6.3})
(and since $\lpartstar$ is convex).
The last inclusion is by
\Cref{cl:thm:meanmap-onto:1}.
The last equality is from
\Cref{pr:meanmap-cont}(\ref{pr:meanmap-cont:b}).

Also,
\begin{equation}    \label{eq:thm:meanmap-onto:2}
  \cl\bigParens{\meanmap(\Rn)}
  =
  \cl\bigParens{\ri(\convfeat)}
  =
  \cl\bigParens{\convfeat}
  =
  \convfeat,
\end{equation}
where the first equality is by
part~(\ref{thm:meanmap-onto:a}), the second
by
\Cref{pr:ri-props}(\ref{pr:ri-props:roc-thm6.3}),
and the third because $\convfeat$ is closed,
being the convex hull of finitely many points
(\Cref{roc:thm19.1}\ref{roc:thm19.1:b}\ref{roc:thm19.1:c}).
Combining Eqs.~(\ref{eq:thm:meanmap-onto:1})
and~(\ref{eq:thm:meanmap-onto:2})
completes the proof.

\pfpart{Part~(\ref{thm:meanmap-onto:c}):}
Suppose $\uu\in\convfeat$.
Then by part~(\ref{thm:meanmap-onto:b}), there exists
$\parambar\in\extspace$ with $\meanmap(\parambar)=\uu$,
and by \Cref{lem:lpart-conj},
$\lpartstar(\uu)=-\entropy(\qxbar)$.%
\indexg{marginal polytope|)}%
\indexg{log-partition function!conjugate of}%
\qedhere
\end{proof-parts}
\end{proof}

\indexg{log-partition function!subgradients of extension|(}%
As we saw in \eqref{eqn:grad-is-meanmap}, for $\vparam\in\Rn$,
$\meanmap(\vparam)$ is exactly the gradient of the log-partition function
at $\vparam$, which means it is the only standard
subgradient of $\lpart$ at $\vparam$
so that $\subdiflpart(\vparam)=\{\meanmap(\vparam)\}$.
This fact generalizes to astral parameters,
as we show next.
In particular, for all $\parambar\in\extspace$,
$\lpartext$ has one and only one astral subgradient at $\parambar$,
namely, $\meanmap(\parambar)$:

\begin{theorem}  \label{thm:subgrad-lpart}
  Assume the general setup of this chapter.
  Let $\parambar\in\extspace$.
  Then
  \[ \basubdiflpart(\parambar) = \regBraces{\meanmap(\parambar)}. \]
\end{theorem}

\begin{proof}
Let $\uu=\meanmap(\parambar)$;
we aim to first show
$\uu\in\asubdiflpart(\parambar)$.
Let $\seq{\vparam_t}$ be any sequence in~$\Rn$ converging to
$\parambar$ (which exists by \Cref{thm:i:1}\ref{thm:i:1d}).
Then
\[
  \vparam_t\cdot\uu - \lpart(\vparam_t)
  =
  \bigExp{\qxbar}{\ln \qxt}
  \rightarrow
  \bigExp{\qxbar}{\ln \qxbar}
  =
  \lpartstar(\uu).
\]
The first equality is by \Cref{pr:exp-log-props}(\ref{pr:exp-log-props:b}),
because $\uu=\regExp{\qxbar}{\featmap}$.
The convergence is by
\Cref{pr:exp-log-props}(\ref{pr:exp-log-props:c}).
The final equality is from \Cref{lem:lpart-conj}, which also shows
that $\lpartstar(\uu)\in\R$.
Thus, $\uu\in\asubdiflpart(\parambar)$ by
\Cref{thm:fminus-subgrad-char}(\ref{thm:fminus-subgrad-char:c},\ref{thm:fminus-subgrad-char:a}).

Now let $\uu'\in\basubdiflpart(\parambar)$; we will show that $\uu'=\meanmap(\parambar)$.
Since $\uu'\in\basubdiflpart(\parambar)$,
\Cref{pr:subgrad-imp-in-cldom}(\ref{pr:subgrad-imp-in-cldom:c})
implies that
$\lpartstar(\uu')\in\R$.
Thus, by 
\Cref{thm:meanmap-onto}(\ref{thm:meanmap-onto:b}),
there exists $\parambar'\in\extspace$
for which $\meanmap(\parambar')=\uu'$.
Also, since $\uu'\in\basubdiflpart(\parambar)$,
\Cref{thm:fminus-subgrad-char}(\ref{thm:fminus-subgrad-char:a},\ref{thm:fminus-subgrad-char:b})
implies
there must exist a sequence
$\seq{\vparam_t}$ in $\Rn$ with
$\vparam_t\rightarrow\parambar$
and
$\vparam_t\cdot\uu'-\lpart(\vparam_t)\rightarrow\lpartstar(\uu')$.
Thus,
\begin{align}
\notag
  \lpartstar(\uu')
  =
  \lim \bigParens{ \vparam_t\cdot\uu' - \lpart(\vparam_t) }
  &=
  \lim \bigParens{ \bigExp{\qxbarp}{\ln \qxt} }
\\
\label{eq:subgrad-lpart}
  &=
  \bigExp{\qxbarp}{\ln \qxbar}
  \le
  \bigExp{\qxbarp}{\ln \qxbarp}
  =\lpartstar(\uu').
\end{align}
As before,
the second equality is by \Cref{pr:exp-log-props}(\ref{pr:exp-log-props:b}),
since $\uu'=\regExp{\qxbarp}{\featmap}$,
and
the third equality is by
\Cref{pr:exp-log-props}(\ref{pr:exp-log-props:c}).
The inequality is by
\Cref{pr:exp-log-props}(\ref{pr:exp-log-props:a}).
The final equality is by
\Cref{lem:lpart-conj}.
Thus, the inequality in \eqref{eq:subgrad-lpart}
must hold with equality, so
$\qxbar=\qxbarp$
by \Cref{pr:exp-log-props}(\ref{pr:exp-log-props:a}).
Therefore,
\indexg{log-partition function!subgradients of extension|)}%
\indexg{mean map|)}%
$\uu'=\meanmap(\parambar')=\regExp{\qxbarp}{\featmap}=\regExp{\qxbar}{\featmap}=\meanmap(\parambar)$.
\end{proof}

\begin{example}[Bernoulli distributions]
\indexg{Bernoulli distributions|(}%
\indexg{extended exponential-family distributions!Bernoulli|(}%
As perhaps the simplest possible example,
suppose $n=1$, $\distset=\{0,1\}$, and
$\featmapsc(\distelt)=\distelt$
for $\distelt\in\distset$.
Then $\sumex(\param)=1+e^\param$,
and $\lpart(\param)=\ln(1+e^\param)$,
for $\param\in\R$,
so $\qsx(1)=e^\param/(1+e^\param)=1/(1+e^{-\param})$,
and $\qsx(0)=1-\qsx(1)$.
Thus, this exponential family of distributions consists of all
Bernoulli distributions
whose mean (which is equal to the probability of outcome $1$) is
in the interval $(0,1)$.
When extended to astral space,
$\lpartext(\barparam)=\logex\bigParens{1+\expex(-\barparam)}$,
for $\barparam\in\Rext$,
and
\[
  \qbarx(1)
  =
  \invex\bigParens{1+\expex(-\barparam)}
  =
  \begin{cases}
        0 & \text{if $\barparam=-\infty$,}
    \\
        1/\bigParens{1+e^{-\barparam}}
          & \text{if $\barparam\in\R$,}
    \\
        1 & \text{if $\barparam=+\infty$,}
  \end{cases}
\]
and $\qbarx(0)=1-\qbarx(1)$.
In this way, Bernoulli distributions with mean equal to $0$ and $1$ are now
also included in the extended exponential family.

The map $\meanmap$ is simply
$\meanmap(\barparam)=\qbarx(1)\cdot 1 + \qbarx(0)\cdot 0 = \qbarx(1)$ 
for $\barparam\in\Rext$.
Thus,
consistent with
\Cref{thm:meanmap-onto}(\ref{thm:meanmap-onto:a},\ref{thm:meanmap-onto:b}),
$\meanmap(\R)= \ri(\conv \featmapsc(\distset))=(0,1)$,
and
$\meanmap(\Rext)= \conv \featmapsc(\distset)=[0,1]$.

As we previously saw in \Cref{ex:subgrad-log1+ex}
(or as we could derive using
\Cref{pr:asubdiffext-at-x-in-rn}\ref{pr:asubdiffext-at-x-in-rn:c}
and
\Cref{pr:subdif-in-1d}),
the astral subdifferentials of $\lpartext$ are,
for $\barparam\in\Rext$,
\[
  \asubdiflpart(\barparam)
  =
  \begin{cases}
           \{0\}   & \text{if $\barparam = -\infty$,} \\
           \{\lpart'(\barparam)\}
                   & \text{if $\barparam\in\R$,} \\
           \{1\}   & \text{if $\barparam = +\infty$,}
  \end{cases}
\]
where $\lpart'(\param)=1/(1+e^{-\param})$ is the derivative of $\lpart$.
Thus,
$\asubdiflpart(\barparam)=\{\meanmap(\barparam)\}$,
consistent with
\Cref{thm:subgrad-lpart}.

Let $u\in [0,1]$.
Then by straightforward calculation,
$\meanmap(\barparam)=u$ holds
if and only if
$\barparam = \ln(u) - \ln(1-u)$,
in which case, $\qbarx(1)=u$.
Thus, by \Cref{thm:meanmap-onto}(\ref{thm:meanmap-onto:c}),
\[
  \lpartstar(u)
  =
  -\entropy(\qbarx)
  =
  u \ln(u) + (1-u)\ln(1-u).%
\indexg{Bernoulli distributions|)}%
\indexg{extended exponential-family distributions!Bernoulli|)}%
\indexg{extended exponential-family distributions|)}%
  \qedhere
\]
\end{example}

\section{Maximum likelihood and maximum entropy}

\indexg{extended exponential-family distributions!estimating from data|(}%
\indexg{maximum-likelihood distribution|(}%
As discussed in
\Cref{sec:exp-fam:standard-setting},
given a random sample drawn from an unknown
distribution~$\popdist$, it is common to estimate $\popdist$ by finding the
exponential-family distribution~$\qx$, parameterized by $\vparam\in\Rn$,
with maximum likelihood, or
equivalently, minimum negative log-likelihood
as given in
\eqref{eqn:neg-log-like:1}.
Having extended the parameter space to all of astral space, we can
now more generally seek parameters $\parambar\in\extspace$ for which the
associated distribution $\qxbar$ has maximum likelihood.

Note that the fractions $\Nelt/\numsamp$ appearing in
\eqref{eqn:neg-log-like:1}
themselves form
a distribution $\phat\in\Delta$, called the
\indexg{empirical distribution}%
\emph{empirical distribution}.
Thus, in slightly more generic terms, given a distribution
$\phat\in\Delta$ (which may or may not have this fractional form),
we take the log-likelihood to be
\begin{equation}  \label{eqn:gen-log-like}
   \regExp{\phat}{\ln \qxbar},
\end{equation}
which, in this approach, we aim to maximize over
\indexg{maximum-likelihood distribution|)}%
$\parambar\in\extspace$.
\indexg{maximum-likelihood distribution!equivalent formulations|(}%
As we show next,
the negative of this log-likelihood is equal to the
extended
log-partition function, re-centered at $\regExp{\phat}{\featmap}$.
As such, maximizing log-likelihood has several equivalent
formulations, which we now summarize:

\begin{proposition}  \label{pr:ml-is-lpartu}
  Assume the general setup of this chapter.
  Let $\uu\in\convfeat$,
  and let $\parambar\in\extspace$.
  Let $\phat\in\Delta$ be any distribution
  for which $\uu=\regExp{\phat}{\featmap}$
  (which must exist).
  Then
  \begin{equation}  \label{eqn:pr:ml-is-lpartu:1}
    {-\regExp{\phat}{\ln \qxbar}} = \lpartextu(\parambar).
  \end{equation}
  Consequently, the following are equivalent:

  \begin{letter-compact-prime}
    \item  \label{pr:ml-is-lpartu:a}
      $\parambar$ maximizes $\regExp{\phat}{\ln \qxbarp}$
      over $\parambar'\in\extspace$.
    \medskip
    \item  \label{pr:ml-is-lpartu:b}
      $\parambar$ minimizes $\lpartextu$.
    \itemprime \label{pr:ml-is-lpartu:bz}
      $\parambar$ minimizes $\sumexextu$.
    \medskip
    \item  \label{pr:ml-is-lpartu:cp}
      $\uu\in\asubdiflpart(\parambar)$.
    \itemprime  \label{pr:ml-is-lpartu:c}
      $\zero\in\asubdiflpartu(\parambar)$.
    \medskip
    \item  \label{pr:ml-is-lpartu:dp}
      $\parambar\in\adsubdiflpartstar(\uu)$.
    \itemprime  \label{pr:ml-is-lpartu:d}
      $\parambar\in\adsubdiflpartustar(\zero)$.
    \medskip
    \item  \label{pr:ml-is-lpartu:ep}
      $\meanmap(\parambar)=\uu$.
    \itemprime  \label{pr:ml-is-lpartu:e}
      $\meanmapu(\parambar)=\zero$.
    \end{letter-compact-prime}

\end{proposition}

\begin{proof}
  ~
  
\begin{proof-parts}
\pfpart{Eq.~(\ref{eqn:pr:ml-is-lpartu:1}):}
Let $\seq{\vparam_t}$ be any sequence in $\Rn$ converging to $\parambar$
(which exists by \Cref{thm:i:1}\ref{thm:i:1d}).
Then for all $t$,
\[
  -\bigExp{\phat}{\ln \qxt}
  =
  \lpart(\vparam_t) - \vparam_t\cdot\uu
  =
  \lpartu(\vparam_t),
\]
from
\Cref{pr:exp-log-props}(\ref{pr:exp-log-props:b})
and
\eqref{eqn:lpart-defn}.
By \Cref{pr:sumex-lpart-cont}(\ref{pr:sumex-lpart-cont:a})
(and \Cref{thm:ext-cont-f}\ref{thm:ext-cont-f:b}),
$\lpartu(\vparam_t)\rightarrow\lpartextu(\parambar)$.
On the other hand,
$-\regExp{\phat}{\ln \qxt}\rightarrow-\regExp{\phat}{\ln \qxbar}$
by \Cref{pr:exp-log-props}(\ref{pr:exp-log-props:c}).
These two limits, which are both for the same sequence,
must be equal, proving the claim.

\medskip

\newcommand{\wideitem}{(d) $\Leftrightarrow$ (d$'$):}

\pfpart{\forcewidthof{\wideitem}{%
  (\ref{pr:ml-is-lpartu:a})
  $\Leftrightarrow$
  (\ref{pr:ml-is-lpartu:b}):}
}
This is immediate from
\eqref{eqn:pr:ml-is-lpartu:1}.

\pfpart{\forcewidthof{\wideitem}{%
  (\ref{pr:ml-is-lpartu:b})
  $\Leftrightarrow$
  (\ref{pr:ml-is-lpartu:cp}):}
}
This is
by
\Cref{thm:fminus-subgrad-char}(\ref{thm:fminus-subgrad-char:a},\ref{thm:fminus-subgrad-char:e}),
noting that $\fminusgen{\lpart}{\uu}=\lpartu$
and $\lpartstar(\uu)\in\R$
(by \Cref{thm:meanmap-onto}\ref{thm:meanmap-onto:c}).

\pfpart{\forcewidthof{\wideitem}{%
  (\ref{pr:ml-is-lpartu:cp})
  $\Leftrightarrow$
  (\ref{pr:ml-is-lpartu:dp}):}
}
This follows from
\Cref{cor:strict-adif-fext-inverses}
(and since $\lpartstar(\uu)\in\R$).

\pfpart{\forcewidthof{\wideitem}{%
  (\ref{pr:ml-is-lpartu:cp})
  $\Leftrightarrow$
  (\ref{pr:ml-is-lpartu:ep}):}
}
This is because
$\asubdiflpart(\parambar)=\{\meanmap(\parambar)\}$
by 
\Cref{thm:subgrad-lpart}.

\pfpart{\forcewidthof{\wideitem}{%
  (\ref{pr:ml-is-lpartu:b})
  $\Leftrightarrow$
  (\ref{pr:ml-is-lpartu:bz}):}
}
This is by \Cref{pr:sumex-lpart-cont}(\ref{pr:sumex-lpart-cont:c}),
and since $\logex$ is a strictly increasing bijection
when restricted to $[0,+\infty]$.

\pfpart{\forcewidthof{\wideitem}{%
  (\ref{pr:ml-is-lpartu:cp})
  $\Leftrightarrow$
  (\ref{pr:ml-is-lpartu:c}):}
}
This is because
$\asubdiflpartu(\parambar)=\asubdif{\fminusgen{\lpartext}{\uu}}{\parambar}=\asubdiflpart(\parambar)-\uu$,
with the second equality
by
\Cref{pr:fminusu-subgrad}(\ref{pr:fminusu-subgrad:ext}).

\pfpart{\forcewidthof{\wideitem}{%
  (\ref{pr:ml-is-lpartu:dp})
  $\Leftrightarrow$
  (\ref{pr:ml-is-lpartu:d}):}
}
This is because
$\adsubdiflpartustar(\zero)=\adsubdif{\fminusgenstar{\lpart}{\uu}}{\zero}=\adsubdiflpartstar(\uu)$,
with the second equality
by
\Cref{pr:fminusu-subgrad}(\ref{pr:fminusu-subgrad:conj}).

\pfpart{\forcewidthof{\wideitem}{%
  (\ref{pr:ml-is-lpartu:ep})
  $\Leftrightarrow$
  (\ref{pr:ml-is-lpartu:e}):}
}
This is by \eqref{eqn:meanmapu-defn}.%
\indexg{maximum-likelihood distribution!equivalent formulations|)}%
\qedhere
\end{proof-parts}
\end{proof}

\begin{example}[Four-outcome exponential family, continued]  \label{ex:exp-fam-eg2}
\indexg{Four-outcome exponential family!maximum-likelihood distribution|(}%
\indexg{maximum-likelihood distribution|(}%
Continuing \Cref{ex:exp-fam-eg1},
we show how the maximum-likelihood distribution
$\qxbar$ can be calculated in this case.

Let $\phat(\distelt)=\Nj{\distelt}/\numsamp$ for
$\distelt\in\distset$.
Then the negative log-likelihood for a distribution $\qx$ with
parameter vector $\vparam\in\R^2$
is $-\regExp{\phat}{\ln \qx}$, which is the same as what was
earlier denoted $\loglik(\vparam)$ in
\eqref{eq:log-like-no-max-eg:2}.
As argued earlier, this function has no finite minimizer in~$\R^2$.

Also, let
\begin{equation}  \label{eq:log-like-no-max-eg:3}
  \uu
  =
  \regExp{\phat}{\featmap}
  =
  \sfrac{2}{3} \featmap(\outA)
  +
  \sfrac{1}{3} \featmap(\outB)
  =
  \trans{[2,1]}.
\end{equation}
Then it can be checked that
the expressions given in
\eqref{eq:log-like-no-max-eg:2} for
$\loglik(\vparam)$ are equal to
$\lpart(\vparam)-\vparam\cdot\uu=\lpartu(\vparam)$,
consistent with
\Cref{pr:exp-log-props}(\ref{pr:exp-log-props:b}).

As we saw in
\Cref{pr:ml-is-lpartu}(\ref{pr:ml-is-lpartu:a},\ref{pr:ml-is-lpartu:bz}),
maximizing the log-likelihood
$\regExp{\phat}{\ln \qxbar}$ over astral parameters $\parambar\in\extspace$
is equivalent to minimizing $\sumexextu$,
for which we have developed extensive tools.
Especially applicable are those from \Cref{sec:emp-loss-min}
since,
for $\vparam\in\R^2$,
\begin{equation}   \label{eq:ex:exp-fam-eg2:1}
  \sumexu(\vparam)
  =
      e^{\param_1}
      +
      e^{-2\param_1}
      +
      e^{-5\param_1}
      +
      e^{-\param_1 - 3\param_2},
\end{equation}
which has the form given in \eqref{eqn:loss-sum-form},
with the variables appearing in that equation set to
$\uu_{\outA}=\trans{[1,0]}$,
$\uu_{\outB}=\trans{[-2,0]}$,
$\uu_{\outC}=\trans{[-5,0]}$,
$\uu_{\outD}=\trans{[-1,-3]}$.
Using \Cref{pr:hard-core:1}(\ref{pr:hard-core:1:b-cones}), we can calculate
$\sumexu\negKern$'s recession cone, yielding
\begin{equation}   \label{eq:ex:exp-fam-eg2:2}
  \resc{\sumexu}
  =
  \{ \lambda \ee_2 :\: \lambda\in\Rpos \}.
\end{equation}
Thus, 
$(\resc\sumexu)^\perp$ includes
$\uu_{\outA}$, $\uu_{\outB}$, and $\uu_{\outC}$,
but not $\uu_{\outD}$;
hence,
by definition of the \irredIndexSet (Eq.~\ref{eq:hc:sumexu}),
$\hardcore{\sumexu} = \{\outA,\outB,\outC\}$.
Consequently, 
\Cref{thm:hard-core:3}(\ref{thm:hard-core:3:b:univ}) implies that the
icon $\ebar=\limray{\ee_2}$ is a universal reducer for
$\sumexu$.
(In fact, it is the only one.)

By \Cref{thm:hard-core:3}(\ref{thm:hard-core:3:d}),
$\sumexu\negKern$'s universal reduction is, for $\vparam\in\R^2$,
\[
  \fullshadsumexu(\vparam)
  =
  \sumexextu(\ebar\plusl\vparam)
  =
      e^{\param_1}
      +
      e^{-2\param_1}
      +
      e^{-5\param_1}.
\]
Setting this function's gradient to zero and solving for $\vparam$,
we find that
$\fullshadsumexu$
is minimized, for instance, at $\qq=\alpha\ee_1$ where
\[
   \alpha
   =
   \frac{\ln\bigParens{1+\sqrt{6}}}{3}
   \approx
   0.41274
\]
(as well as at all points $\trans{[\alpha,\lambda]}$
for $\lambda\in\R$).

Since $\ebar\in\unimin{\sumexu}$ and $\qq$ minimizes
$\fullshadsumexu$,
\Cref{thm:waspr:hard-core:2}
implies that the resulting astral point $\parambar=\ebar\plusl\qq$
minimizes $\sumexextu$, and so also maximizes the log-likelihood
$\regExp{\phat}{\ln \qxbar}$.
The resulting probability distribution $\qxbar$ can then be calculated
from \eqref{eqn:qxbar-expr} to be:
\begin{equation}  \label{eq:log-like-no-max-eg:4}
  \qxbar(\outA)
  \approx
  0.72783,
  \quad
  \qxbar(\outB)
  \approx
  0.21100,
  \quad
  \qxbar(\outC)
  \approx
  0.06117,
  \quad
  \qxbar(\outD)
  =
  0.
\end{equation}
It can be checked that $\meanmap(\parambar)=\uu$, consistent with
\Cref{pr:ml-is-lpartu}(\ref{pr:ml-is-lpartu:a},\ref{pr:ml-is-lpartu:ep}).%
\indexg{Four-outcome exponential family!maximum-likelihood distribution|)}%
\indexg{maximum-likelihood distribution|)}%
\end{example}

\indexg{maximum entropy|(}%
\indexg{entropy (of distribution)!maximum|(}%
We next consider an alternative approach for estimating the unknown
distribution $\popdist$.
The idea, first, is to find a distribution $q\in\Delta$
(not necessarily in the exponential family)
under which
the expectation of the features matches what was observed on the data,
that is, for which $\regExp{q}{\featmap}=\uu$, where
$\uu=\regExp{\phat}{\featmap}$
(and where $\phat\in\Delta$ is given or observed, as above).
Typically, there will be many distributions satisfying this
property.
Among these, we choose the one of \emph{maximum entropy}, that is,
for which $\entropy(q)$ is largest.
Such a distribution is, in an information-theoretic sense,
closest
to the uniform distribution, which is viewed as the default, ``uninformed''
estimate we would choose in the absence of any data
(%
\indexa{Jaynes, E. T.}%
\citealp{Jaynes1957};
\indexa{Kullback, S.}%
\citealp{Kullback1959}%
).
\looseness=-1

\indexg{maximum-likelihood distribution!maximum entropy and|(}%
Thus,
the maximum-entropy approach estimates $\popdist$
by that distribution $q\in\medists$ for which $\entropy(q)$
is maximized, where
$\medists$ is the set of all distributions $q\in\Delta$ for which
$\regExp{q}{\featmap} = \uu$.
By comparison, the maximum-likelihood approach estimates $\popdist$
by that exponential-family distribution $q\in\mldists$ which maximizes
the log-likelihood
$\regExp{\phat}{\ln q}$,
where $\mldists$ is the set of all (extended) exponential-family
distributions $\qxbar$ for $\parambar\in\extspace$.

Remarkably, these two approaches always yield the same distribution,
as we show next.
Moreover, that distribution is always the unique point at the
intersection of the two sets $\mldists\cap\medists$, which is to say,
that exponential-family distribution $\qxbar$ for which
$\regExp{\qxbar}{\featmap} = \uu$.
This extends a similar formulation for the standard setting given by
\indexa{Della Pietra, S.}\indexa{Della Pietra, V.}\indexa{Lafferty, J.}%
\citet{DellaDeLa97}
to the current astral setting.

\begin{theorem}  \label{thm:ml-equal-maxent}
  Assume the general setup of this chapter.
  Let $\phat\in\Delta$, let $\uu=\regExp{\phat}{\featmap}$, and let
  \begin{align*}
  \SwapAboveDisplaySkip
    \mldists &= \bigBraces{ \qxbar :\: \parambar\in\extspace }, \\
    \medists &= \bigBraces{ q\in\Delta :\: \regExp{q}{\featmap} = \uu }.
  \end{align*}
  Let $q\in\Delta$.
  Then the following are equivalent:
  \begin{letter}
  \item  \label{thm:ml-equal-maxent:a}
    $\displaystyle q = \argmax*_{q'\in\mldists} \regExp{\phat}{\ln q'}$.
  \item  \label{thm:ml-equal-maxent:b}
    $q \in \mldists\cap\medists$.
  \item  \label{thm:ml-equal-maxent:c}
    $\displaystyle q = \argmax*_{q'\in\medists} \entropy(q')$.
  \end{letter}
  Furthermore, there exists a unique distribution $q$ satisfying all
  of these.
\end{theorem}

\begin{proof}
~

\begin{proof-parts}
\pfpart{%
  (\ref{thm:ml-equal-maxent:a})
  $\Rightarrow$
  (\ref{thm:ml-equal-maxent:b}):
}
Suppose $q$
satisfies statement~(\ref{thm:ml-equal-maxent:a}).
Then $q\in\mldists$, so $q=\qxbar$ for some ${\parambar\in\extspace}$,
and furthermore,
$\parambar$ maximizes $\regExp{\phat}{\ln \qxbarp}$
over $\parambar'\in\extspace$.
By \Cref{pr:ml-is-lpartu}(\ref{pr:ml-is-lpartu:a},\ref{pr:ml-is-lpartu:ep}),
this implies that
$\regExp{\qxbar}{\featmap}=\meanmap(\parambar)=\uu$.
Therefore, $q=\qxbar$
satisfies statement~(\ref{thm:ml-equal-maxent:b}).

\pfpart{%
  (\ref{thm:ml-equal-maxent:b})
  $\Rightarrow$
  (\ref{thm:ml-equal-maxent:c}):
}
Suppose statement~(\ref{thm:ml-equal-maxent:b}) holds, so that
$q=\qxbar\in\medists$, for some $\parambar\in\extspace$.
Let $q'$ be any distribution that is also in
$\medists$, so that $\regExp{q'}{\featmap}=\uu=\regExp{q}{\featmap}$.

Let $\seq{\vparam_t}$ in $\Rn$ be a sequence converging to $\parambar$
(which exists by \Cref{thm:i:1}\ref{thm:i:1d}).
Then for each $t$,
\[
  \bigExp{q}{\ln \qxt}
  =
  \vparam_t\cdot\bigExp{q}{\featmap} - \lpart(\vparam_t)
  =
  \vparam_t\cdot\bigExp{q'\negKern}{\featmap} - \lpart(\vparam_t)
  =
  \bigExp{q'\negKern}{\ln \qxt}
  \leq
  \bigExp{q'\negKern}{\ln q'}.
\]
The first and third equalities are both by
\Cref{pr:exp-log-props}(\ref{pr:exp-log-props:b}),
and the inequality is
by
\Cref{pr:exp-log-props}(\ref{pr:exp-log-props:a}).
In the limit,
\[
  \bigExp{q}{\ln \qxt}
  \rightarrow
  \bigExp{q}{\ln \qxbar}
  =
  \bigExp{q}{\ln q},
\]
with convergence by
\Cref{pr:exp-log-props}(\ref{pr:exp-log-props:c}).
Hence,
$\regExp{q}{\ln q} \leq \regExp{q'}{\ln q'}$,
so
$\entropy(q')\leq\entropy(q)$.
Since this holds for all $q'\in\medists$,
$q$ must have maximum entropy among all such distributions.

\pfpart{Existence:}
By \Cref{pr:fext-min-exists},
$\lpartextu$ must attain its minimum at some point
$\parambar\in\extspace$.
By
\Cref{pr:ml-is-lpartu}(\ref{pr:ml-is-lpartu:b},\ref{pr:ml-is-lpartu:a}),
the resulting distribution $q=\qxbar$ must
satisfy statement~(\ref{thm:ml-equal-maxent:a}),
and so statements~(\ref{thm:ml-equal-maxent:b})
and~(\ref{thm:ml-equal-maxent:c}) as well, as just shown.

\pfpart{Uniqueness:}
We show that at most one distribution can satisfy
statement~(\ref{thm:ml-equal-maxent:c}), implying, by the foregoing, that
only one distribution can satisfy statements~(\ref{thm:ml-equal-maxent:a})
or~(\ref{thm:ml-equal-maxent:b}) as well.\looseness=-1

Suppose, by way of contradiction, that two distributions $p$ and $p'$
in $\medists$ both have maximum entropy among all such distributions,
and that $p\neq p'$.
Let $r=(p+p')/2$, which is also in $\medists$
(using linearity of Eq.~\ref{eqn:expect-defn} in $p$),
and which is also distinct from both $p$ and $p'$.
Then
\[
  \regExp{p}{\ln p} + \regExp{p'\negKern}{\ln p'}
  >
  \regExp{p}{\ln r} + \regExp{p'\negKern}{\ln r}
  =
  2 \regExp{r}{\ln r},
\]
with the strict inequality following from
\Cref{pr:exp-log-props}(\ref{pr:exp-log-props:a}),
and the equality by direct calculation from \eqref{eqn:expect-defn}.
Thus,
\[
  \frac{\entropy(p)+\entropy(p')}{2}
  <
  \entropy(r),
\]
contradicting that $\entropy(p)=\entropy(p')$ is maximum among all
distributions in $\medists$.

\pfpart{%
  (\ref{thm:ml-equal-maxent:c})
  $\Rightarrow$
  (\ref{thm:ml-equal-maxent:a}):
}
Suppose $q$ satisfies statement~(\ref{thm:ml-equal-maxent:c}).
As shown above, there exists a distribution $p$ satisfying
statement~(\ref{thm:ml-equal-maxent:a}), and so also 
statement~(\ref{thm:ml-equal-maxent:c}), by the implications proved already.
Having proved uniqueness, this implies $p=q$,
so $q$ satisfies statement~(\ref{thm:ml-equal-maxent:a})
\indexg{entropy (of distribution)!maximum|)}%
\indexg{maximum entropy|)}%
\indexg{maximum-likelihood distribution!maximum entropy and|)}%
as well.
\qedhere
\end{proof-parts}
\end{proof}

The maximum-likelihood problem,
or the equivalent maximum-entropy problem
from \Cref{thm:ml-equal-maxent},
are typically solved by iterative
optimization algorithms, like those considered in
\Cref{sec:iterative}.
By \Cref{pr:ml-is-lpartu}, negative log likelihood
takes the form~$\lpartextu$, for some $\uu\in\conv\featmap(\distset)$, which is
continuous on astral space (by
\Cref{pr:sumex-lpart-cont}\ref{pr:sumex-lpart-cont:a}).
Therefore,
as we noted in \Cref{ex:MLE},
any standard
algorithm that satisfies the conditions of \Cref{thm:fact:gd} will maximize
the likelihood in the limit, even when the maximizer is at infinity.%
\indexg{extended exponential-family distributions!estimating from data|)}%

\section{Astral galaxies and faces of the marginal polytope}

Earlier, we studied the mean map $\meanmap$ and its relationship with
the marginal polytope, $\convfeat$,
showing
in
\Cref{thm:meanmap-onto}(\ref{thm:meanmap-onto:a},\ref{thm:meanmap-onto:b})
that $\meanmap$ maps $\extspace$ onto that polytope,
and maps $\Rn$ onto its relative interior.
Below, we fill in more detail regarding this relationship.
Previously, in \Cref{sec:galaxies}, we saw how astral space
can be naturally partitioned into galaxies.
In standard convex analysis,
convex sets in $\Rn$, including polytopes, can be partitioned in a
different way, into the relative interiors of their faces,
as we saw in \Cref{roc:thm18.2}.
As we will soon show, $\meanmap$ directly links these two partitions by
continuously mapping every
galaxy $\galaxd$, for any icon $\ebar$,
onto the relative interior of one face $C$ of the marginal polytope,
and also mapping
the closure of that galaxy, $\galcld$, onto the entire face $C$.
We will also give various characterizations of which galaxies map
onto relative interiors of which faces.

\indexg{marginal polytope!faces of|(}%
To begin, the next proposition shows that there is a direct
correspondence between each face of the marginal polytope
and the convex hull of the feature vectors associated with
\irredIndicesFor~$\sumexu$,
where $\uu$ is any point in the face's
relative interior:

\begin{proposition}   \label{pr:hardcore-face-margpoly}
  Assume the general setup of this chapter.
  Let $C$ be a nonempty face of $\convfeat$, and
  let $\uu\in\ri C$.
  Then
  $\hardcore{\sumexu}=\Braces{\distelt\in\distset :\: \featmap(\distelt)\in C}$,
  and
  $C=\conv{\featmap(\hardcore{\sumexu})}$.
\end{proposition}

\begin{proof}
Let $S=\convfeat$, and
let $C'=C-\uu$, and
$S'=S-\uu=\conv{\featmapu(\distset)}$
be translations of $C$ and $S$ by $\uu$.
Then $C'$ is a face of $S'$,
and $\zero=\uu-\uu\in(\ri{C}) - \uu = \ri{C'}$.
Thus, 
\begin{equation*}  %
  \hardcore{\sumexu}
  =
  \Braces{\distelt\in\distset :\: \featmapu(\distelt)\in C'}
  =
  \Braces{\distelt\in\distset :\: \featmap(\distelt)\in C},
\end{equation*}
where the first equality follows
directly from
\Cref{thm:erm-faces-hardcore2}(\ref{thm:erm-faces-hardcore2:a},\ref{thm:erm-faces-hardcore2:c})
(applied to $\sumexu$).
Likewise, by
\Cref{thm:erm-faces-hardcore2}(\ref{thm:erm-faces-hardcore2:a},\ref{thm:erm-faces-hardcore2:d}),
$C'=\conv{\featmapu(\hardcore{\sumexu})}$,
implying
\indexg{marginal polytope!faces of|)}%
$C=\conv{\featmap(\hardcore{\sumexu})}$.
\end{proof}

\indexg{support (of distribution)|(}%
The \emph{support} of a distribution $p\in\Delta$,
denoted $\support p$, is the set of points assigned nonzero
probability:
\[
\indexm{supp p}{$\support p$}{support of distribution}%
  \support p = \bigBraces{ \distelt\in\distset :\: p(\distelt)>0 }.%
\indexg{support (of distribution)|)}%
\]
\indexg{extended exponential-family distributions!support of|(}%
The next proposition characterizes the support of $\qxbar$ for any
$\parambar\in\extspace$:

\begin{proposition}   \label{pr:support-face-margpoly}
  Assume the general setup of this chapter.
  Suppose $\parambar=\ebar\plusl\qq$ for some $\ebar\in\corezn$
  and $\qq\in\Rn$.
  Let $\distelt\in\distset$.
  Then $\qxbar(\distelt)>0$ if and only if
  $\ebar\cdot (\featmap(\distalt) - \featmap(\distelt)) \leq 0$
  for all $\distalt\in\distset$.
\end{proposition}

\begin{proof}
The form of $\qxbar(\distelt)$ is given in \eqref{eqn:qxbar-expr}.
As such, $\qxbar(\distelt)=0$ if and only if
$\parambar\cdot(\featmap(\distalt)-\featmap(\distelt))=+\infty$
for some $\distalt\in\distset$,
which holds if and only if
$\ebar\cdot(\featmap(\distalt)-\featmap(\distelt))=+\infty$
for some $\distalt\in\distset$,
since $\parambar=\ebar\plusl\qq$.
Since $\ebar$ is an icon, this proves the claim
(using
\Cref{pr:icon-equiv}\ref{pr:icon-equiv:a}\ref{pr:icon-equiv:b}).%
\indexg{extended exponential-family distributions!support of|)}%
\end{proof}

\indexg{mean map!image of galaxies under|(}%
\indexg{galaxies!marginal polytope faces and|(}%
\indexg{marginal polytope!faces of|(}%
\indexg{persistent indices!marginal polytope faces and|(}%
\indexg{extended exponential-family distributions!support of|(}%
As discussed above, the next theorem shows that every
galaxy $\galaxd$ is mapped by the mean map to the relative interior of
a face of the marginal polytope, as specified in the theorem.
In addition,
the theorem characterizes when a particular galaxy $\galaxd$  is
mapped to the relative interior of a particular face $C$, that is,
when $\meanmap(\galaxd)=\ri C$.
For any (and every) point $\uu\in\ri C$,
and for any (and every) point $\parambar\in\galaxd$,
this relationship holds if and only if the support of $\qxbar$ matches
the \irredIndexSetFor~$\sumexu$, or equivalently,
if and only if $\ebar$ is a univeral reducer of the function
$\lpartu$.

\begin{theorem}  \label{thm:galaxy-mapto-face}
  Assume the general setup of this chapter.
  Let $\ebar\in\corezn$ and let
  \begin{equation}   \label{eq:thm:galaxy-mapto-face:0}
     Z
     =
     \bigBraces{
       \distelt \in\distset :\:
       \forall \distalt\in\distset,\;
          \ebar\cdot \bigParens{\featmap(\distalt) - \featmap(\distelt)}
          \leq 0
     }.
  \end{equation}
  Also, let $D=\conv{\featmap(Z)}$, and $S=\convfeat$.
  Then $D$ is a nonempty face of $S$,
  and
  $\meanmap(\galaxd)=\ri{D}$.

  Moreover, for every nonempty face $C$ of $S$,
  for every $\uu\in\ri{C}$, and for every
  $\parambar\in\galaxd$,
  the following are equivalent:
  \begin{letter-compact}
  \item  \label{thm:galaxy-mapto-face:d:1}
    $\meanmap(\galaxd) = \ri{C}$.
  \item  \label{thm:galaxy-mapto-face:d:2}
    $\meanmap(\galcld) = C$.
  \item  \label{thm:galaxy-mapto-face:d:0}
    $\meanmap(\galaxd) \cap (\ri{C}) \neq \emptyset$.
  \medskip
  \item  \label{thm:galaxy-mapto-face:d:5}
    $C=\conv{\featmap(\support\qxbar)}$.
  \item  \label{thm:galaxy-mapto-face:d:4}
    $\support\qxbar=\{\distelt\in\distset :\: \featmap(\distelt)\in C\}$.
  \medskip
  \item  \label{thm:galaxy-mapto-face:d:n1}
    $\hardcore{\sumexu} = \support\qxbar$.
  \item  \label{thm:galaxy-mapto-face:d:n2}
    $\hardcore{\sumexu} = Z$.
  \medskip
  \item  \label{thm:galaxy-mapto-face:d:3z}
      $\ebar\in\unimin{\sumexu}$.
  \item  \label{thm:galaxy-mapto-face:d:3a}
      $\ebar\in\unimin{\lpartu}$.
  \end{letter-compact}
\end{theorem}

Note that the set $Z$
in \eqref{eq:thm:galaxy-mapto-face:0} consists of the points in
$\distset$ that satisfy the characterization given in
\Cref{pr:support-face-margpoly}, and therefore is equal, by that
proposition, to the support of $\qxbar$ for any $\parambar\in\galaxd$.
Before proving this theorem, we look more closely in
the next lemma at the set $Z$,
as well as the convex hull $D$ of
its image under $\featmap$.
The lemma gives a more explicit expression for
$\ebar\cdot (\featmap(\distalt)-\featmap(\distelt))$
where $\distelt\in Z$ and $\distalt\in\distset$,
and more generally of $\ebar\cdot(\ww-\uu)$
for $\uu\in D$ and $\ww\in S$.

\begin{lemma}   \label{lem:Z-set-char}
  Assume the general setup of this chapter,
  and let $\ebar$, $Z$, $D$, and $S$ be as in
  \Cref{thm:galaxy-mapto-face}.
  \begin{letter-compact}
  \item   \label{lem:Z-set-char:a}
    Suppose $\distelt\in Z$ and $\distalt\in\distset$.
    Then
    \[
      \ebar\cdot \bigParens{\featmap(\distalt)-\featmap(\distelt)}
      =
      \begin{cases}
        0           & \text{if $\distalt\in Z$,} \\
        -\infty     & \text{otherwise.}
      \end{cases}
    \]
  \item   \label{lem:Z-set-char:b}
    Suppose $\uu\in D$ and $\ww\in S$.
    Then
    \[
       \ebar\cdot (\ww-\uu)
       =
       \begin{cases}
         0           & \text{if $\ww\in D$,} \\
         -\infty     & \text{otherwise.}
       \end{cases}
    \]
  \item   \label{lem:Z-set-char:c}
    Suppose $\distalt\in\distset$.
    Then $\distalt\in Z$ if and only if
    $\featmap(\distalt)\in D$.
  \end{letter-compact}
\end{lemma}

\begin{proof}
  ~

\begin{proof-parts}
\pfpart{Part~(\ref{lem:Z-set-char:a}):}
Since $\distelt\in Z$,
$\ebar\cdot (\featmap(\distalt)-\featmap(\distelt))\leq 0$.
Since $\ebar$ is an icon,
$\ebar\cdot (\featmap(\distalt)-\featmap(\distelt))$
must therefore be equal to $-\infty$ or $0$
(by
\Cref{pr:icon-equiv}\ref{pr:icon-equiv:a}\ref{pr:icon-equiv:b}).
Thus, to prove the claim, it suffices to show that
$\distalt\in Z$ if and only if
\begin{equation}      \label{eq:lem:Z-set-char:1}
  \ebar\cdot \bigParens{\featmap(\distalt)-\featmap(\distelt)} = 0.
\end{equation}

If $\distalt\in Z$, then we can apply $Z$'s definition to
$\distalt$, yielding
$\ebar\cdot (\featmap(\distelt)-\featmap(\distalt))\leq 0$,
thereby proving \eqref{eq:lem:Z-set-char:1}.

Conversely, suppose
\eqref{eq:lem:Z-set-char:1} holds.
Let $\distaltii\in\distset$.
Then
\[
  \ebar\cdot\bigParens{\featmap(\distaltii)-\featmap(\distalt)}
  =
  \ebar\cdot\bigParens{\featmap(\distaltii)-\featmap(\distelt)}
  +
  \ebar\cdot\bigParens{\featmap(\distelt)-\featmap(\distalt)}
  \le
  0.
\]
The equality is by
\Cref{pr:i:1}
combined with
\eqref{eq:lem:Z-set-char:1}.
The inequality follows by \eqref{eq:lem:Z-set-char:1}
and because $\distelt\in Z$.
Since this holds for all $\distaltii\in\distset$, this proves
$\distalt\in Z$.

\pfpart{Part~(\ref{lem:Z-set-char:b}):}
Since $\uu$ and $\ww$ are in $S$, there exist distributions $p$ and
$q$ in $\Delta$ for which
$\regExp{p}{\featmap}=\uu$
and
$\regExp{q}{\featmap}=\ww$.
Furthermore, since $\uu\in D$, we can assume
without loss of generality that
$\support p\subseteq Z$.

We then have
\begin{align}
  \ebar\cdot(\ww-\uu)
  &=
  \ebar\cdot\BiggParens{\,
    \sum_{\distelt\in Z} \sum_{\distalt\in\distset}
         p(\distelt) q(\distalt)
         \bigParens{\featmap(\distalt)-\featmap(\distelt)}
    }
  \nonumber
  \\
  &=
  \sum_{\distelt\in Z} \sum_{\distalt\in\distset}
     p(\distelt) q(\distalt)
     \BigBracks{\ebar\cdot\bigParens{\featmap(\distalt)-\featmap(\distelt)}}.
  \label{eq:lem:Z-set-char:2} %
\end{align}
The first equality is because
$\regExp{q}{\featmap}=\ww$
and
$\regExp{p}{\featmap}=\uu$.
The second equality follows because
$\ebar\cdot(\featmap(\distalt)-\featmap(\distelt))\leq 0$
for $\distelt\in Z$ and $\distalt\in\distset$,
by $Z$'s definition,
so the terms in the brackets are summable,
allowing us to repeatedly apply
\Cref{pr:i:1}
(as well as \Cref{pr:i:2}).

In the case that $\ww\in D$,
we can choose $q$ with
$\support q\subseteq Z$.
Combined with
part~(\ref{lem:Z-set-char:a}),
it then follows that every term of
\eqref{eq:lem:Z-set-char:2} is equal to $0$, so
$\ebar\cdot(\ww-\uu)=0$, as claimed.

In the remaining case that $\ww\not\in D$,
we must have $q(\distalt)>0$ for some $\distalt\not\in Z$.
Further, $p$ cannot have empty support, so $p(\distelt)>0$ for some
$\distelt\in Z$.
By part~(\ref{lem:Z-set-char:a}),
$\ebar\cdot(\featmap(\distalt)-\featmap(\distelt))=-\infty$
for this pair, so at least this one term of
\eqref{eq:lem:Z-set-char:2} is $-\infty$.
Since every other term
is nonpositive
(by $Z$'s definition),
we conclude $\ebar\cdot(\ww-\uu)=-\infty$.

\pfpart{Part~(\ref{lem:Z-set-char:c}):}
If $\distalt\in Z$, then clearly
$\featmap(\distalt)\in\conv{\featmap(Z)}=D$.
Conversely, if $\featmap(\distalt)\in D$, then
for any $\distelt\in Z$,
$\ebar\cdot(\featmap(\distalt)-\featmap(\distelt))=0$
by part~(\ref{lem:Z-set-char:b})
(with $\ww=\featmap(\distalt)\in D$
and $\uu=\featmap(\distelt)\in D$).
Therefore, $\distalt\in Z$
by part~(\ref{lem:Z-set-char:a}).
\qedhere
\end{proof-parts}
\end{proof}

\begin{proof}[Proof of \Cref{thm:galaxy-mapto-face}]
  ~
  
\begin{proof-parts}
\pfpart{That $D$ is a nonempty face of $S$:}
By \Cref{pr:support-face-margpoly}, the set $Z$
is equal to the support of $\qebar$, which cannot be empty.
Therefore, $Z$ is nonempty, so $D$ is as well.

Let $\ww_0,\ww_1\in S$, let $\lambda\in (0,1)$,
and suppose that
$\uu=(1-\lambda)\ww_0+\lambda\ww_1$
is in $D$.
To prove $D$ is a face of $S$,
we aim to show that $\ww_0$ and $\ww_1$ are also in $D$.

\Cref{lem:Z-set-char}(\ref{lem:Z-set-char:b})
implies
$\ebar\cdot(\ww_b-\uu)\leq 0$,
for $b\in\{0,1\}$.
Also, by rearranging the definition of $\uu$,
we have
$\zero=(1-\lambda)(\ww_0-\uu)+\lambda(\ww_1-\uu)$,
so
\[
  0
  =
  \ebar\cdot\zero
  =
  (1-\lambda)\regBracks{\ebar\cdot(\ww_0-\uu)}
  +
  \lambda\regBracks{\ebar\cdot(\ww_1-\uu)},
\]
where the second equality
uses
\Cref{pr:i:1}.
Since neither term is positive on the right, and
since $\lambda\in (0,1)$,
this implies
$\ebar\cdot(\ww_b-\uu)=0$,
for $b\in\{0,1\}$.
Therefore, $\ww_b\in D$
by
\Cref{lem:Z-set-char}(\ref{lem:Z-set-char:b}),
proving that $D$ is a face
(as defined in \Cref{sec:prelim:faces}).

\pfpart{That $\meanmap(\galaxd)=\ri{D}$:}
We next give two claims
related to $\ebar$ being
a universal
reducer for some $\sumexu$, which will be used together in proving
$\meanmap(\galaxd)=\ri{D}$:

\begin{claimpx}    \label{cl:lem:Z-set-char:1}
  Suppose $\uu\in\ri{D}$.
  Then
  $\ebar\in\unimin{\sumexu}$.
\end{claimpx}

\begin{proofx}
We have
\begin{equation}  \label{eq:thm:galaxy-mapto-face:1}
  \hardcore{\sumexu}
  =
  \bigBraces{\distelt\in\distset :\: \featmap(\distelt)\in D}
  =
  Z,
\end{equation}
where
the first equality is by
\Cref{pr:hardcore-face-margpoly},
and the second is by
\Cref{lem:Z-set-char}(\ref{lem:Z-set-char:c}).
Also, for $\distelt\in\distset$,
\begin{equation*}
  \ebar\cdot\featmapu(\distelt)
  =
  \ebar\cdot\bigParens{\featmap(\distelt)-\uu}
  =
  \begin{cases}
    0           & \text{if $\distelt\in Z$,} \\
    -\infty     & \text{otherwise,}
  \end{cases}
\end{equation*}
where the second equality follows
from \Cref{lem:Z-set-char}(\ref{lem:Z-set-char:b}),
noting that $\distelt\in Z$ if and only if $\featmap(\distelt)\in D$,
by \Cref{lem:Z-set-char}(\ref{lem:Z-set-char:c}).
Since $Z=\hardcore{\sumexu}$
(by Eq.~\ref{eq:thm:galaxy-mapto-face:1}),
this implies by
\Cref{thm:hard-core:3}(\ref{thm:hard-core:3:b:univ})
that $\ebar\in\unimin{\sumexu}$, completing the proof.
\end{proofx}

\begin{claimpx}    \label{cl:lem:Z-set-char:2}
  Let $\uu\in S$, and suppose
  $\ebar\in\unimin{\sumexu}$.
  Then
  $\uu\in\meanmap(\galaxd)$.
\end{claimpx}

\begin{proofx}
By
\Cref{pr:univ-red-props}(\ref{pr:univ-red-props:min}),
$\fullshadsumexu$ attains its minimum at some point $\qq\in\Rn$,
implying
by \Cref{pr:unimin-to-global-min} that
$\sumexextu$ is minimized by
$\zbar=\ebar\plusl\qq$, which is in $\galaxd$.
Then $\meanmap(\zbar)=\uu$ by
\Cref{pr:ml-is-lpartu}(\ref{pr:ml-is-lpartu:bz},\ref{pr:ml-is-lpartu:ep}),
proving the claim.
\end{proofx}

To prove $\meanmap(\galaxd)=\ri{D}$,
suppose first that $\parambar\in\galaxd$.
By \Cref{pr:support-face-margpoly},
$\qxbar(\distelt)>0$ if and only if $\distelt\in Z$.
Therefore,
$\meanmap(\parambar)=\regExp{\qxbar}{\featmap}\in\ri(\conv{\featmap(Z)})=\ri{D}$,
with the inclusion by \Cref{pr:ri-conv-finite}.
Thus, $\meanmap(\galaxd)\subseteq\ri{D}$.

For the reverse inclusion, if $\uu\in\ri{D}$,
then
$\ebar\in\unimin{\sumexu}$ by
\Cref{cl:lem:Z-set-char:1},
implying
$\uu\in\meanmap(\galaxd)$
by \Cref{cl:lem:Z-set-char:2}.
Thus,
$\ri{D}\subseteq\meanmap(\galaxd)$,
so
$\meanmap(\galaxd)=\ri{D}$.

\medskip

Finally, we prove the stated equivalences.
Let $C$ be a nonempty face of $S$,
let $\uu\in\ri{C}$, and let $\parambar\in\galaxd$.

\pfpart{%
  (\ref{thm:galaxy-mapto-face:d:0})
  $\Rightarrow$
  (\ref{thm:galaxy-mapto-face:d:1}):
}
As just shown, $\meanmap(\galaxd)=\ri D$,
where $D$ is a nonempty face of $S$.
Therefore, if $(\ri D)\cap(\ri C)\neq\emptyset$,
then $D=C$
by \Cref{pr:face-props}(\ref{pr:face-props:cor18.1.2}),
implying
$\meanmap(\galaxd)=\ri D = \ri C$.

\pfpart{%
  (\ref{thm:galaxy-mapto-face:d:1})
  $\Rightarrow$
  (\ref{thm:galaxy-mapto-face:d:2}):
}
Suppose $\meanmap(\galaxd) = \ri{C}$.
Then taking closures yields
\[
  \meanmap(\galcld)
  =
  \cl\bigParens{\meanmap(\galaxd)}
  =
  \cl\bigParens{\ri{C}}
  =
  \cl{C}
  =
  C.
\]
The first equality is by
\Cref{pr:meanmap-cont}(\ref{pr:meanmap-cont:b}).
The third equality is by
\Cref{pr:ri-props}(\ref{pr:ri-props:roc-thm6.3}).
The last equality is because $S$ is a polytope and therefore closed
(\Cref{roc:thm19.1}\ref{roc:thm19.1:b}\ref{roc:thm19.1:c}),
implying that its face $C$ is also closed
(\Cref{pr:face-props}\ref{pr:face-props:cor18.1.1}).

\pfpart{%
  (\ref{thm:galaxy-mapto-face:d:2})
  $\Rightarrow$
  (\ref{thm:galaxy-mapto-face:d:5}):
}
Suppose $\meanmap(\galcld) = C$.
We showed above that $\meanmap(\galaxd)=\ri D$.
Applied with $C$ replaced by $D$, the last implication therefore
shows that
$\meanmap(\galcld) = D = \conv{\featmap(Z)}$.
Since $Z=\support\qxbar$
(by \Cref{pr:support-face-margpoly}),
this proves the claim.

\pfpart{%
  (\ref{thm:galaxy-mapto-face:d:5})
  $\Rightarrow$
  (\ref{thm:galaxy-mapto-face:d:4}):
}
Suppose $C=\conv{\featmap(\support\qxbar)}$.
Noting that $Z=\support\qxbar$, this
implies $C=\conv{\featmap(Z)}=D$,
and so,
by \Cref{lem:Z-set-char}(\ref{lem:Z-set-char:c}),
that $Z=\{\distelt\in\distset :\: \featmap(\distelt)\in C\}$,
proving the claim.

\pfpart{%
  (\ref{thm:galaxy-mapto-face:d:4})
  $\Rightarrow$
  (\ref{thm:galaxy-mapto-face:d:n1}):
}
This follows directly from \Cref{pr:hardcore-face-margpoly}.

\pfpart{%
  (\ref{thm:galaxy-mapto-face:d:n1})
  $\Rightarrow$
  (\ref{thm:galaxy-mapto-face:d:n2}):
}
This is immediate since $Z = \support\qxbar$ (by \Cref{pr:support-face-margpoly}).

\pfpart{%
  (\ref{thm:galaxy-mapto-face:d:n2})
  $\Rightarrow$
  (\ref{thm:galaxy-mapto-face:d:3z}):
}
Suppose $\hardcore{\sumexu} = Z$.
Then
\[
  C
  =
  \conv\featmap(\hardcore{\sumexu})
  =
  \conv\featmap(Z)
  =
  D,
\]
where the first equality is by
\Cref{pr:hardcore-face-margpoly}.
Thus, $\uu\in\ri{C}=\ri{D}$,
so $\ebar\in\unimin{\sumexu}$ by \Cref{cl:lem:Z-set-char:1}.

\pfpart{%
  (\ref{thm:galaxy-mapto-face:d:3z})
  $\Leftrightarrow$
  (\ref{thm:galaxy-mapto-face:d:3a}):
}
We claim
$\unimin{\sumexu}=\unimin{\lpartu}$.
To see this, let $\ebar\in\corezn$.
Then
$\ebar\in\unimin{\sumexu}$
if and only if
$\sumexextu(\ebar\plusl\vparam)\leq\sumexextu(\ebar'\plusl\vparam)$
for all $\ebar'\in\corezn$,
by definition of universal reducer and universal reduction
(Definitions~\ref{def:univ-reducer}
and~\ref{def:univ-reduction}).
Since
$\logex$
is a strictly increasing bijection when
restricted to $[0,+\infty]$,
and by
\Cref{pr:sumex-lpart-cont}(\ref{pr:sumex-lpart-cont:c}),
this in turn holds if and only if
$\lpartextu(\ebar\plusl\vparam)\leq\lpartextu(\ebar'\plusl\vparam)$
for all $\ebar'\in\corezn$,
that is,
if and only if
$\ebar\in\unimin{\lpartu}$.

\pfpart{%
  (\ref{thm:galaxy-mapto-face:d:3z})
  $\Rightarrow$
  (\ref{thm:galaxy-mapto-face:d:0}):
}
Suppose $\ebar\in\unimin{\sumexu}$.
Then $\uu\in\meanmap(\galaxd)$
by \Cref{cl:lem:Z-set-char:2}, so
$\meanmap(\galaxd)$ and $\ri{C}$ are not disjoint
(since $\uu$ is in both).
\qedhere
\end{proof-parts}
\end{proof}

\begin{figure}
  \centering
  \includegraphics{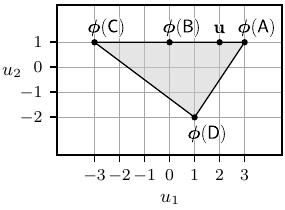}
  \mycaption{Maximum likelihood and faces of the marginal polytope}{
\indexf{Four-outcome exponential family!marginal polytope}%
    The mar\-gi\-nal polytope for the four-outcome exponential family
    and the empirical feature mean $\uu=\regExp{\phat}{\featmap}$, where $\phat$
    is the empirical distribution from \Cref{ex:exp-fam-eg1,ex:exp-fam-eg2}.
    The empirical mean $\uu$ lies in the relative interior of the edge connecting
    $\featmap(\outA)$ and $\featmap(\outC)$, which includes an additional feature vector
    $\featmap(\outB)$. As argued in \Cref{ex:exp-fam-eg3}, this means that the support
    of the maximum-likelihood distribution is the set $\set{\outA,\outB,\outC}$.}
  \label{fig:marginal-u}%
\end{figure}

\begin{example}[Four-outcome exponential family, continued]
\label{ex:exp-fam-eg3}%
\indexg{Four-outcome exponential family!marginal polytope|(}%
Continuing our running example, let $\distset$, $\featmap$, $\uu$,
$\ebar$, and $\parambar$ be as in \Cref{ex:exp-fam-eg2}.
Then the marginal polytope,
$\convfeat$, is a triangle in $\R^2$
with vertices
$\featmap(\outA),\featmap(\outC),\featmap(\outD)$
(see \Cref{fig:marginal-u}).
Its faces are the triangle itself,
its three edges, three vertices,
and the empty set.
Let $C$ be the edge
connecting
$\featmap(\outA)$ and $\featmap(\outC)$, which also includes $\featmap(\outB)$.
As shown in \Cref{fig:marginal-u}, $\uu$ is in the relative interior of this edge.

Earlier, we noted that $\hardcore{\sumexu}=\{\outA,\outB,\outC\}$, which is
consistent, according to \Cref{pr:hardcore-face-margpoly},
with the corresponding points
$\featmap(\outA),\featmap(\outB),\featmap(\outC)$
being the only ones in $C$.
The support of $\qxbar$
(given in Eq.~\ref{eq:log-like-no-max-eg:4})
is also equal to this same set
$\{\outA,\outB,\outC\}$; it can be checked that this is the same as the set $Z$
defined in \eqref{eq:thm:galaxy-mapto-face:0}, and so consistent
with \Cref{pr:support-face-margpoly}.
Thus, the set $D=\conv{\featmap(Z)}$ is the same as edge $C$.

We also earlier established that
$\ebar\in\unimin{\sumexu}=\unimin{\lpartu}$.
From these facts, we see that
statements~(\ref{thm:galaxy-mapto-face:d:5}),
(\ref{thm:galaxy-mapto-face:d:4}),
(\ref{thm:galaxy-mapto-face:d:n1}),
(\ref{thm:galaxy-mapto-face:d:n2}),
(\ref{thm:galaxy-mapto-face:d:3z}),
(\ref{thm:galaxy-mapto-face:d:3a})
of \Cref{thm:galaxy-mapto-face} all hold.
That theorem therefore also
implies that
$\meanmap(\galaxd) = \ri{C}$
and
$\meanmap(\galcld) = C$.%
\indexg{mean map!image of galaxies under|)}%
\indexg{galaxies!marginal polytope faces and|)}%
\end{example}

\indexg{maximum-likelihood distribution!support of|(}%
In this example,
although $\phat(\outC)=\phat(\outD)=0$, the resulting
maximum-likelihood distribution $\qxbar$
in \eqref{eq:log-like-no-max-eg:4}
includes outcome~$\outC$ in its support, but not outcome~$\outD$.
This shows that, in general, the support of the maximum-likelihood
distribution $\qxbar$ need not match the support of a given
empirical distribution $\phat$.
There is, nonetheless, a strong relationship between these.
The next theorem shows that the point $\uu=\regExp{\phat}{\featmap}$ must
be in the relative interior of the smallest face $C$ of the marginal
polytope that includes all points
$\featmap(\distelt)$ for which $\phat(\distelt)>0$.
The \irredIndexSet $\hardcore{\sumexu}$ then consists of those
$\distelt$ for which $\featmap(\distelt)\in C$.
Further, the support of $\qxbar$ is exactly
$\hardcore{\sumexu}$.

For instance, in the example just discussed,
the smallest face that includes
$\featmap(\outA)$ and $\featmap(\outB)$
(corresponding to the points in the support of $\phat$)
is the edge $C$ identified above,
which also includes $\featmap(\outC)$,
yielding the \irredIndexSet $\{\outA,\outB,\outC\}$.
The support of the resulting maximum-likelihood
distribution $\qxbar$ is exactly the \irredIndexSet, which thus
includes outcome~$\outC$ (since it is in the smallest face $C$ containing
$\support\phat$), but not outcome~$\outD$.%
\indexg{Four-outcome exponential family!marginal polytope|)}%

\begin{theorem}  \label{thm:comp-supports-phat-ml}
  Assume the general setup of this chapter.
  Let $\phat\in\Delta$,
  let $\uu=\regExp{\phat}{\featmap}$,
  and let $C$ be the smallest face of $\convfeat$ that includes
  $\featmap(\support\phat)$.
  Let $\parambar\in\extspace$ maximize
  $\regExp{\phat}{\ln \qxbar}$.
  Then $\uu\in\ri C$, and
  \begin{equation}    \label{eqn:thm:comp-supports-phat-ml:1}
     \support\phat
     \subseteq
     \hardcore{\sumexu}
     =
     \support\qxbar
     =
     \{\distelt\in\distset :\: \featmap(\distelt)\in C\}.
  \end{equation}
\end{theorem}

\begin{proof}
Let $E=\support\phat$.
Then $\uu\in\ri(\conv{\featmap(E)})$
by \Cref{pr:ri-conv-finite},
implying, after translation, that
$\zero\in\ri(\conv{\featmapu(E)})$.
Therefore, $E\subseteq\hardcore{\sumexu}$
by
\Cref{thm:erm-faces-hardcore}(\ref{thm:erm-faces-hardcore:a})
(applied to $\sumexu$),
proving the inclusion in
\eqref{eqn:thm:comp-supports-phat-ml:1}.

Next, let $F$ be the unique face of $\conv\featmap(\distset)$
for which $\uu\in\ri F$
(which exists by \Cref{roc:thm18.2}). Since $\uu\in\ri(\conv{\featmap(E)})$
and $\uu\in F$, we have $\conv{\featmap(E)}\subseteq F$
by \Cref{pr:face-props}(\ref{pr:face-props:thm18.1}), so $C\subseteq F$
by minimality of $C$.
Since $\uu\in(\ri F)\cap C$, also $F\subseteq C$ (again by \Cref{pr:face-props}\ref{pr:face-props:thm18.1}). Thus, $F=C$, so $\uu\in\ri C$.

We can write $\parambar=\ebar\plusl\qq$ for some icon $\ebar\in\corezn$
and some $\qq\in\Rn$.
By \Cref{pr:ml-is-lpartu}(\ref{pr:ml-is-lpartu:a},\ref{pr:ml-is-lpartu:bz}),
$\parambar$ minimizes $\sumexextu$, so
$\ebar\in\unimin{\sumexu}$
by \Cref{thm:waspr:hard-core:2}
(since $\exp$ is strictly increasing).
By \Cref{thm:galaxy-mapto-face}(\ref{thm:galaxy-mapto-face:d:3z},\ref{thm:galaxy-mapto-face:d:n1},\ref{thm:galaxy-mapto-face:d:4}),
this in turn implies the two equalities of
\eqref{eqn:thm:comp-supports-phat-ml:1},
\indexg{extended exponential-family distributions!support of|)}%
\indexg{persistent indices!marginal polytope faces and|)}%
\indexg{maximum-likelihood distribution!support of|)}%
completing the proof.
\end{proof}

\indexg{mean map!yielding homeomorphism|(}%
\Cref{thm:galaxy-mapto-face}
shows that for any icon $\ebar\in\corezn$,
$\meanmap$ maps the galaxy $\galaxd$ surjectively onto $\ri{D}$, for
some face $D$ of $\convfeat$.
In general, this mapping is not injective.
However, as we show next, if we further restrict the mapping to just
part of the galaxy, then it becomes a bijection, and in fact induces
a homeomorphism.

In more detail, suppose $\uu\in\ri{D}$.
By
\Cref{pr:ml-is-lpartu}(\ref{pr:ml-is-lpartu:ep},\ref{pr:ml-is-lpartu:bz}),
$\meanmap(\parambar)=\uu$ if and only if $\parambar$ minimizes
$\sumexextu$.
From \Crefequiv{thm:waspr:hard-core:4}{pr:hard-core:4:a}{pr:hard-core:4:b},
all such minimizers have the form
$\zbar\plusl\qq$ where $\zbar\in\conv(\unimin{\sumexu})$,
and $\qq\in\Rn$ is a uniquely determined point in
$\rescperp{\sumexu}$.
We will see next that the set
$\rescperp{\sumexu}$
is the same linear subspace $L$ for all $\uu\in\ri{D}$.
We then restrict $\meanmap$ to the part of the galaxy
$\galaxd$ corresponding to $L$, resulting in the mapping
$\qq\mapsto \meanmap(\ebar\plusl\qq)$, for $\qq\in L$.
This mapping is a homeomorphism.

To state the formal result,
recall from \Cref{sec:prelim:affine-sets}
that the affine hull of a set
$S\subseteq\Rn$, denoted $\affh{S}$,
is the smallest affine set (in $\Rn$) that includes $S$,
and that
for every nonempty affine set $A\subseteq\Rn$,
there exists a unique linear subspace
that is parallel to $A$
(\Cref{roc:thm1.2}).

\begin{theorem}  \label{thm:aff-homeo}
  Assume the general setup of this chapter,
  and let $\ebar$, $Z$, $D$, and $S$ be as in
  \Cref{thm:galaxy-mapto-face}
  (implying, by that theorem, that $D$ is a face of $S$
  with $\meanmap(\galaxd)=\ri{D}$).
  Also, let
  $L$ be the linear subspace of $\Rn$ that is parallel to
  $\affh{\featmap(Z)}$,
  the affine hull of~$\featmap(Z)$.
  Then:
  \begin{letter-compact}
  \item    \label{thm:aff-homeo:a}
    $L=\rescperp{\sumexu}=\rescperp{\lpartu}$,
    for all $\uu\in\ri{D}$.
  \item    \label{thm:aff-homeo:b}
    Let $\rho:L\rightarrow\ri{D}$
    be defined by $\rho(\vparam)=\meanmap(\ebar\plusl\vparam)$,
    for $\vparam\in L$.
    Then $\rho$ is a homeomorphism.
  \end{letter-compact}
\end{theorem}

\begin{proof}
  ~

\begin{proof-parts}
\pfpart{Part~(\ref{thm:aff-homeo:a}):}
Let $\uu\in\ri{D}$.
That $\resc{\sumexu}=\resc{\lpartu}$
follows from the definition of standard recession cone
(Eq.~\ref{eqn:resc-cone-def})
and
$\ln$
being strictly increasing on $\Rstrictpos$.
Furthermore,
\[
  \rescperp{\sumexu}
  =
  \spn\featmapu(\hardcore{\sumexu})
  =
  \spn\featmapu(Z),
\]
where the first equality is by
\Cref{thm:hard-core:3}(\ref{thm:hard-core:3:e}),
and the second is because
$\hardcore{\sumexu}=Z$
by \Cref{thm:galaxy-mapto-face}(\ref{thm:galaxy-mapto-face:d:1},\ref{thm:galaxy-mapto-face:d:n2}).
Also, $\uu$ is in $D=\conv{\featmap(Z)}$, and so in
$\affh{\featmap(Z)}$ as well.
Therefore, by
\Cref{pr:lin-aff-par},
$\spn\featmapu(Z)=\spn(\featmap(Z) - \uu)$
is the linear subspace $L$ parallel to
$\affh{\featmap(Z)}$, as claimed.

\pfpart{Part~(\ref{thm:aff-homeo:b}):}
The map $\zbar\mapsto\ebar\plusl\zbar$,
for $\zbar\in\extspace$,
is continuous by
\Cref{cor:aff-cont}, and $\meanmap$ is continuous by
\Cref{pr:meanmap-cont}(\ref{pr:meanmap-cont:a}).
Therefore, $\rho$, which is constructed by composing these and then
restricting to a subspace, is as well.

Next, we argue that $\rho$ is a bijection.
Let $\uu\in\ri{D}$; we aim to show there exists a unique point
$\vparam\in L$ with $\rho(\vparam)=\uu$.

By \Cref{thm:galaxy-mapto-face}(\ref{thm:galaxy-mapto-face:d:1},\ref{thm:galaxy-mapto-face:d:3z}),
$\ebar\in\unimin{\sumexu}$.
We apply results from \Cref{sec:emp-loss-min}
to the function $\sumexu$, noting that
$\exp$
is strictly
increasing and strictly convex.
By \Cref{thm:waspr:hard-core:2},
a point $\vparam\in\Rn$
minimizes $\fullshadsumexu$
if and only if
$\ebar\plusl\vparam$ minimizes $\sumexextu$,
which, by
\Cref{pr:ml-is-lpartu}(\ref{pr:ml-is-lpartu:bz},\ref{pr:ml-is-lpartu:ep}),
holds if and only if
$\meanmap(\ebar\plusl\vparam)=\uu$.
By \Cref{thm:waspr:hard-core:4},
there exists
a unique point $\vparam\in \rescperp{\sumexu}=L$
minimizing $\fullshadsumexu$,
and therefore
a unique point $\vparam\in L$
with
$\rho(\vparam)=\meanmap(\ebar\plusl\vparam)=\uu$,
proving $\rho$ is a bijection.

It remains then only to argue that
$\rhoinv$, the functional inverse of $\rho$,
is continuous.
Suppose by way of contradiction that $\rhoinv$ is not continuous.
Then there exists a sequence $\seq{\uu_t}$ in $\ri{D}$ converging to
some point $\uu\in\ri{D}$ such that
$\rhoinv(\uu_t)\not\rightarrow\rhoinv(\uu)$.
Let $\vparam_t=\rhoinv(\uu_t)$ for each $t$,
and let $\vparam=\rhoinv(\uu)$.
Since $\vparam_t\not\rightarrow\vparam$,
there exists a neighborhood $U\subseteq\extspace$ of $\vparam$ for which
$\vparam_t\not\in U$ for infinitely many values of $t$.
By discarding all other sequence elements
and corresponding elements of $\seq{\uu_t}$,
we assume henceforth that
$\vparam_t\not\in U$ for all $t$.

Further, by sequential compactness, the sequence $\seq{\vparam_t}$
must have a
convergent subsequence in $\extspace$;
by again discarding all other sequence elements, we can assume the
entire sequence converges to some point $\parambar'\in\extspace$.
Since each $\vparam_t$ is in the closed set
$\extspace\setminus U$, their limit, $\parambar'$,
must be as well.
In particular, since $\vparam\in U$,
this implies $\parambar'\neq \vparam$.

By continuity of $\meanmap$
(and of the affine map $\zbar\mapsto\ebar\plusl\zbar$),
we then have
\[
  \uu_t=\rho(\vparam_t)=\meanmap(\ebar\plusl\vparam_t)
  \to
  \meanmap(\ebar\plusl\parambar').
\]
Hence, $\meanmap(\ebar\plusl\parambar')=\uu\in\ri{D}$,
since also $\uu_t\rightarrow\uu$.

We can write
$\parambar'=\ebar'\plusl\qq'$ for some $\ebar'\in\corezn$ and $\qq'\in\Rn$,
so
$\meanmap(\ebar\plusl\ebar'\plusl\qq')=\uu\in\ri{D}$.
Therefore,
$\ebar\plusl\ebar'\in\unimin{\sumexu}\subseteq\resc\sumexextu$,
with first inclusion
by \Cref{thm:galaxy-mapto-face}(\ref{thm:galaxy-mapto-face:d:0},\ref{thm:galaxy-mapto-face:d:3z}),
and the second by 
\Cref{pr:new:thm:f:8a}(\ref{pr:new:thm:f:8a:a}).
Hence,
\begin{align}
\notag
  \ebar'\in\aspan\set{\ebar'}
  =
  \clbar{\rspanset{\ebar'}}
  &\subseteq
  \clbar{\rspanset{\ebar\plusl\ebar'}}
\\
\label{eq:thm:aff-homeo:2}
  &\subseteq
  \clbar{\rspan(\resc\sumexextu)}
  =
  \aspan(\resc\sumexextu),
\end{align}
where
the second and third inclusions are by definition of restricted span,
and both equalities are from \Cref{thm:rspan-is-aspan-in-rn}.

On the other hand, $\parambar'\in\Lbar$ since each $\vparam_t\in L$.
Thus,
\[
  \ebar'\in\Lbar=\clbar{\rescperp{\sumexu}}=\clbar{\rescperp{\sumexextu}},
\]
where the inclusion is by \Cref{pr:aspan:decomp},
the first equality is by part~(\ref{thm:aff-homeo:a}),
and the second equality is by
\Cref{pr:hard-core:1}(\ref{pr:hard-core:1:b-resc2}).
Combined with \eqref{eq:thm:aff-homeo:2}, this implies
$\ebar'=\zero$ by \Cref{pr:aperp:intersect}.
\Cref{pr:aspan:decomp} also implies that $\qq'\in L$,
so $\parambar'=\qq'\in L$. Since
$\rho(\parambar')=\meanmap(\ebar\plusl\parambar')=\uu=\rho(\vparam)$,
and $\rho$ is a bijection, this contradicts that
$\parambar'\neq\vparam$.

Having reached a contradiction, we conclude that $\rhoinv$ is
continuous, completing the proof.
\indexg{marginal polytope!faces of|)}%
\qedhere
\end{proof-parts}
\end{proof}

\indexg{marginal polytope!when feature map minimal|(}%
The feature map $\featmap$ is said to be
\indexg{minimal (feature map)}%
\indexg{feature map!minimal}%
\emph{minimal} if
there does
not exist any $\vparam\in\Rn\wo\set{\zero}$ for which
$\vparam\inprod\featmap(\distelt)$ is
constant as a function of $\distelt\in\distset$.
We saw in
\Cref{thm:meanmap-onto}(\ref{thm:meanmap-onto:a})
that $\meanmap$ maps $\Rn$ onto $\ri(\convfeat)$.
We show next that
when also $\featmap$ is minimal,
the affine hull of $\featmap(\distset)$ is all of $\Rn$, implying,
as a consequence of \Cref{thm:aff-homeo},
that $\meanmap$ induces a homeomorphism from
$\Rn$ to the interior of the marginal polytope $\convfeat$,
which is the same as its relative interior in this case.

\indexg{Four-outcome exponential family!minimality|(}%
For instance, the feature map $\featmap$ given in
\Cref{ex:exp-fam-eg1}
is minimal, and as such, $\meanmap$, in this case,
maps $\R^2$ homeomorphically onto the interior of the
marginal polytope, the triangle with vertices
$\featmap(\outA),\featmap(\outB),\featmap(\outD)$.%
\indexg{Four-outcome exponential family!minimality|)}%

\begin{corollary}
  Assume the general setup of this chapter,
  and also that $\featmap$ is minimal.
  Let $S=\conv{\featmap(\distset)}$.
  Then $\meanmap(\Rn)=\intr S$.
  Moreover, the function $\rho:\Rn\rightarrow\intr S$
  defined by $\rho(\vparam)=\meanmap(\vparam)$,
  for $\vparam\in\Rn$,
  is a homeomorphism.
\end{corollary}

\begin{proof}
Let $A=\affh{\featmap(\distset)}$, and let $L$ be the linear subspace parallel to $A$.
We claim first that $A=L=\Rn$.
By possibly renaming the elements of $\distset$, we
assume that $\distset=\set{1,\dotsc,m}$ for some $m\ge 1$. Let
$\ww$ be any element of $A$, and let
$\B\in\R^{n\times m}$ be the matrix whose $i$-th column is
$\featmap(i)-\ww$, for $i=1,\ldots,m$.
Then, by \Cref{pr:lin-aff-par}, $L=\colspace\B$.

Let $\vparam\in\Rn\wo\set{\zero}$. Then $\trans{\B}\vparam\ne\zero$,
because $\trans{\B}\vparam=\zero$ would imply,
for $i=1,\ldots,m$, that
$\trans{(\featmap(i)-\ww)}\vparam=0$,
and so that $\vparam\cdot\featmap(i)=\vparam\cdot\ww$,
contradicting that $\featmap$ is minimal. Thus, the $n$ columns of
$\trans{\B}$
are linearly independent, so the rank of $\trans{\B}$ is $n$, implying
also that the rank of $\B$ is $n$
(\Cref{pr:row-rank-is-col-rank}).
Since $\colspace\B\subseteq\Rn$,
this implies that $L=\colspace\B=\Rn$.
Since $A$ is parallel to $L$, we therefore obtain $A=\Rn$.

Since
$\affh{\featmap(\distset)}=A=\Rn$,
we also have $\affh{S}=\Rn$ since
$S\supseteq \featmap(\distset)$.
Thus, $\ri{S}=\intr{S}$
(by definition of relative interior),
implying $\meanmap(\Rn)=\intr{S}$ by
\Cref{thm:meanmap-onto}(\ref{thm:meanmap-onto:a}).

Finally, we can apply
\Cref{thm:aff-homeo}(\ref{thm:aff-homeo:b})
with $\ebar=\zero$,
$Z=\distset$, $D=S$, and $L=\Rn$,
since,
by the foregoing arguments, these settings satisfy
the assumptions of that \namecref{thm:aff-homeo}.
This yields that $\rho$ is a homeomorphism, as claimed.%
\indexg{mean map!yielding homeomorphism|)}%
\indexg{marginal polytope!when feature map minimal|)}%
\end{proof}

\backmatter
\bibliographystyle{plainnat}
\clearpage
\phantomsection
\bookmarksetup{startatroot}
\addtocontents{toc}{\vspace{1.25\baselineskip}}
\addcontentsline{toc}{chapter}{\bibname}

\DeclareRobustCommand{\VAN}[3]{#3}
\bibliography{ab}

\index{astral space!topology|see{astral topology}}
\index{topology!astral|see{astral topology}}
\index{upper addition|see{upward addition}}
\index{extension (of function)|see{lower semicontinuous extension}}
\index{astral rank|see{rank, astral}}
\index{effective domain|see{domain, effective}}
\index{subspaces, linear|see{linear subspaces}}
\index{compactness!sequential|see{sequential compactness}}
\index{countability axioms|see{first-countability; second-countability}}
\index{convex cone generated by set|see{conic hull}}
\index{linear map, image of function under|see{linear image of astral function; linear image of standard function}}
\index{image of function under linear map|see{linear image of astral function; linear image of standard function}}
\index{subgradients|see{subdifferentials}}
\index{derivative, directional|see{one-sided directional derivative}}
\index{directional derivative|see{one-sided directional derivative}}
\index{astral points!functions@as functions|see{functional representation of astral points}}
\index{decompositions, astral|seealso{representations of astral points}}
\index{astral points!representations of|see{representations of astral points}}
\index{matrix representation|see{representations of astral points}}
\index{neighborhood base|seealso{countable neighborhood base}}
\index{astral linear maps|see{linear maps, astral}}
\index{astral affine maps|see{affine maps, astral}}

\index{astral points!sequences involving|seeunder{sequence convergence}}
\index{convergence|see{sequential convergence}}
\index{representations of astral points!canonical|see{canonical representation}}
\index{continuity!extensible|see{extensible continuity}}
\index{primal conjugate|see{conjugate, primal (astral)}}
\index{dual conjugate|see{conjugate, dual (astral)}}
\index{primal biconjugate|see{biconjugate, primal (astral)}}
\index{dual biconjugate|see{biconjugate, dual (astral)}}
\index{increasing function|see{nondecreasing function; strictly increasing function}}
\index{astral topology!closure in|see{closure, astral}}
\index{maxima, pointwise|see{suprema and maxima, pointwise}}
\index{outer hull|see{outer convex hull}}
\index{convex hull, outer|see{outer convex hull}}
\index{multiples|see{$\oms$-multiples; scalar multiples (astral); sequential multiples}}
\index{convex combinations!sequential|see{sequential convex combinations}}
\index{conic combinations!sequential|see{sequential conic combinations}}
\index{astral cones|see{cones, astral}}
\index{astral rays|see{rays, astral}}
\index{conic hull, outer|see{outer conic hull}}
\index{linear maps, astral|seealso{affine maps}}
\index{astral conic hull|see{conic hull, astral}}
\index{hull|see{astral conic hull; convex hull; outer conic hull; outer convex hull}}
\index{outer convex hull!finite set@of finite set|see{polytopes, astral}}
\index{dual separation|see{separation, strong astral dual}}
\index{closure, representational|see{representational closure}}
\index{dual polar|see{polar, dual}}
\index{astral linear subspaces|see{linear subspaces, astral}}
\index{astral span|see{span, astral}}
\index{linear functions|see{affine functions}}
\index{astral closedness|see{closedness, astral (primal)}}
\index{astral dual closedness|see{closedness, astral dual}}
\index{astral recession cone|see{recession cone, astral}}
\index{lower semicontinuous extension!minimizers of|see{minimizers of extensions}}
\index{reduction, universal|see{universal reduction}}
\index{reducers, universal|see{universal reducers}}
\index{minimizers of extensions!canonical|see{canonical minimizers}}
\index{easy set|seealso{reducible indices}}
\index{hard core|seealso{persistent indices}}
\index{continuity!extensions@of extensions|see{continuity of extensions}}
\index{lower semicontinuous extension!continuity of|see{continuity of extensions}}
\index{lower semicontinuous extension!subdifferentials of|see{subdifferentials of extensions}}
\index{composition (with function)|see{affine map in composition with; linear map in composition with; nondecreasing function, composition with}}
\index{optimality conditions!constrained|see{constrained optimality conditions}}
\index{games, matrix|seealso{minimax theorems}}
\index{convex programs|see{ordinary convex programs}}
\index{Karush-Kuhn-Tucker|see{KKT}}
\index{gradients (standard)|seealso{subdifferentials}}
\index{gradients (standard)!sequences of|see{subgradient sequences}}
\index{subdifferentials of extensions!sequences in|seeunder{subgradient sequences}}
\index{subdifferentials, astral dual!sequences in|seeunder{subgradient sequences}}
\index{sequence convergence!subgradients@of subgradients|see{subgradient sequences}}
\index{sequence convergence!minimizers@to minimizers|see{minimizing sequences}}
\index{minimizers of extensions!sequences converging to|see{minimizing sequences}}
\index{exponential-family distributions!astral extension|see{extended exponential-family distributions}}
\index{likelihood|seealso{maximum likelihood}}
\index{dual orthogonal complement|see{orthogonal complement, dual}}
\index{upper triangular matrices|seealso{positive upper triangular matrices}}
\index{differential continuity|seeunder{continuity}}
\index{scalar multiples (astral)!omega@by $\oms$|see{$\oms$-multiples}}
\index{projection matrices!component|see{component projection matrices}}
\index{matrices|seealso{component projection matrices; positive upper triangular matrices}}
\index{basis (topological)|see{base}}

\index{addition|see{downward addition; leftward addition; sequential sum; set sum; upward addition}}

\addtocontents{toc}{\vspace{1.25\baselineskip}}
\printindex[notind]
\printindex[default][\indexintrotext]
\printindex[autind]
\end{document}